\documentclass[11pt]{report}

\usepackage{t-angles}
\usepackage[all]{xy}
\usepackage{amsmath,amsthm, amscd, amssymb, amsfonts}

\usepackage{hhline}

\title{
BRAIDED HOPF ALGEBRAS }
\author{
Shouchuan Zhang \\ Department  of Mathematics, Hunan University,
410082
\\ P.R.China\\
 }
\date{}
\begin{document}
% \tableofcontents
\newtheorem{Theorem}{\quad Theorem}[section]
\newtheorem{Proposition}[Theorem]{\quad Proposition}
\newtheorem{Definition}[Theorem]{\quad Definition}
\newtheorem{Corollary}[Theorem]{\quad Corollary}
\newtheorem{Lemma}[Theorem]{\quad Lemma}
\newtheorem{Example}[Theorem]{\quad Example}
\newtheorem{Remark}[Theorem]{\quad Remark}

\maketitle \addtocounter{chapter}{0}

 \newpage
 {\ }
\vskip16cm {\ }\\
Hunan Normal University Press \\
First edition  \  1999 \\
Second edition  \ 2005 \\
Third edition  \ 2007 \\
ISBN7-81031-812-8/0.036

 \newpage
%$$\hbox { \bf \large PREFACE}$$
\chapter *{Preface}
The term ``quantum group'' was popularized by Drinfeld in his
address to the International Congress of Mathematicians in Berkeley
. However, the concepts of quantum groups and quasitriangular Hopf
algebras
 are the same.  Therefore, Hopf algebras have close
 connections
with various areas of mathematics and physics.

      The development of  Hopf algebras can be divided into five stages  .

      The first stage is that integral and Maschke's theorem are found.
  Maschke's theorem gives a nice  criterion  of semisimplicity, which
  is due to M.E. Sweedler.  The second
  stage is that the Lagrange's theorem is  proved, which is due to
  W.D. Nichols and M.B. Zoeller. The third stage is
  the research of the actions of Hopf algebras, which unifies  actions
  of groups and Lie algebras. S. Montgomery's ``Hopf algebras and their actions on rings'' contains the main results
    in this area. S. Montgomery  and R.J. Blattner show the duality theorem.
     S. Montegomery,
     M. Cohen, H.J. Schneider, W. Chin, J.R. Fisher and the author
     study the relation between algebra $R$ and the crossed
     product $R \#_\sigma H$ of the semisimplicity, semiprimeness,
     Morita equivalence and radical.
     S. Montgomery and M. Cohen ask whether $R \#_\sigma H$  is semiprime
     when $R$ is semiprime for a finite-dimensional semiprime Hopf algebra
     $H$. This is a famous semiprime problem. The answer is sure when
     the action of $H$ is inner or  $H $  is commutative or cocommutative,
 which is due to S. Montgomery, H.J. Schneider and the author.
 This is the best
 result in the study of the semiprime problem up to now. However, this problem is still
 open. Y. Doi  and M. Takeachi show that the crossed product
 is a cleft extension.

    The research of quantum groups is the fourth stage. The concepts
    of quantum groups and (co)quasitriangular Hopf algebras
    are the same.
The Yang-Baxter equation first came up in a paper by Yang as
factorization condition of the scattering S-matrix in the many-body
problem in one dimension and in work of Baxter on exactly solvable
models in statistical mechanics. It has been playing an important
role in mathematics and physics ( see \cite {BD82}, \cite {Ma93a} ).
Attempts to find solutions of the Yang-Baxter equation (YBE) in a
systematic way have led to the theory of quantum groups.
    In other words, we can obtain a solution of  the Yang-Baxter equation
    by a (co)quasitriangular Hopf algebra. It is well-known  that
    the universal enveloping algebras of Lie algebras  are Hopf algebras.
    But they are not quasitriangular  in general. Drinfeld  obtained
    a quasitriangular Hopf algebra by means of constructing
    a quantum enveloping
    algebras of simple Lie algebra. Drinfeld  also obtained
    a quasitriangular Hopf algebra $D(H)$, quantum double,
    for any Hopf algebra $H$.
    D.E. Radford and R.G. Larson study the algebraic
    structure of (co)quasitriangular Hopf algebras. Kassel's ``Quantum Group''
    contains the main results in this area.

    The research of braided Hopf algebras and classification of finite-dimensional Hopf
    algebras  are  the fifth
    stage.
Supersymmetry has attracted a great deal of interest from physicists
and mathematicians (see \cite {Ka77}, \cite {Sc79}, \cite {MR94} ).
It finds numerous applications in particle physics and statistical
mechanics ( see \cite {CNS75} ) . Of central importance in
supersymmetric theories is the  $Z_2$-graded structure which
permeates them. Superalgebras and super-Hopf algebras are most
naturally very important.

Braided tensor categories were introduced by Joyal and  Street in
1986 \cite {JS86}. It is a generalization of super case. Majid,
Joyal, Street and Lyubasheko have reached many interesting
conclusions in braided tensor categories, for example, the braided
reconstruction theorem, transmutation and bosonisation, integral,
q-Fourier transform, q-Mikowski space, random walk and so on  (see
\cite {Ma93a}, \cite {Ma95a}, \cite {Ma95b}, \cite {MR94}, \cite
{Ly95} ).

   The main results of my research   consist  of the following:

   In braided tensor categories we show the Maschke's theorem and give the
   necessary and sufficient conditions for double cross biproducts and
   crossbiproducts and biproducts to be
   bialgebras. We  obtain the factorization theorem for braided Hopf algebras;

   In symmetric tensor categories, we show the duality theorem and construct
   quantum double;

   In ordinary vector space category with ordinary twist,
       we obtain  the relation between the global dimension,
       the  weak dimension,
       the Jacbson radical, Baer radical of algebra $R$ and its
       crossed product
$R \# _\sigma H$. We also give the relation between the
decompositions of comodules and coalgebras. We obtain all  solutions
of constant classical Yang-Baxter
 equation (CYBE)  in Lie algebra $L$ with dim $L \le 3$. We also give
   the sufficient and  necessary conditions for  $(L, \hbox {[ \ ]},
 \Delta _r, r)$  to be a coboundary
 (or triangular ) Lie bialgebra.

 These are also the points of  originality of this book.

This book is organized into three parts. In the three parts we study
Hopf algebras in braided tensor categories, symmetric braided tensor
categories and the category of vector spaces over field $k$ with
ordinary twist braiding respectively. Furthermore, we divide the
three parts  into twelve  chapters and they are briefly described as
follows.

In Chapter \ref {c1}, we recall the basic concepts and conclusions
of Hopf algebras living in braided tensor categories.

In  Chapter \ref {c2},  we obtain the fundamental theorem of Hopf
modules living in braided tensor   categories. In Yetter-Drinfeld
category, we obtain that the antipode $S$  of $H$ is invertible and
the integral $\int _H^l$  of $H$  is a one-dimensional space for any
finite-dimensional braided  Hopf algebra $H$. We give Maschke's
Theorem for YD Hopf algebras.

In Chapter \ref {c3}, the double biproduct
 $D = A \stackrel {b}{\bowtie } H $ of two bialgebras $A$  and $H$ is
 constructed in a braided tensor category and the necessary  and
 sufficient conditions for  $D$ to be a  bialgebra are given.
 The universal property of
 double biproduct is  obtained. In particular,
 the necessary and sufficient  conditions for
 double biproducts, double products (or bicrossed products ),
  double coproducts, Drinfeld double and smash products
  of two anyonic (or super ) Hopf algebras $A$  and $H$  to be
   anyonic Hopf algebras
  (or super-Hopf algebras)
  are obtained.

In  Chapter \ref {c4}, we construct the bicrossproducts and
biproducts in braided tensor category  $({\cal C}, C)$ and give the
necessary and sufficient conditions for them to be bialgebras. We
obtain that  for any left module $(A, \alpha )$  of bialgebra $H$
the smash product $A \#H$  is a bialgebra iff $A$  is an $H$-module
bialgebra and $H$ is cocommutative with respect to $(A, \alpha ).$

In  Chapter \ref {c5}, We show that the braided group analogue of
double cross (co)product is the  double cross (co)product of braided
group analogues. we give the factorization theorem of Hopf algebras
living in braided tensor categories and
 the relation between
 Hopf algebra  $H$ and  its factors.

In Chapter \ref {c4'}, braided m-Lie algebras induced by
multiplication are introduced, which generalize Lie algebras, Lie
color algebras and quantum Lie algebras. The necessary and
sufficient conditions for the braided m-Lie  algebras to be strict
Jacobi braided Lie algebras are given. Two classes of braided m-Lie
algebras are given, which are generalized matrix braided m-Lie
algebras and braided m-Lie subalgebras of $End _F M$, where $M$ is a
Yetter-Drinfeld module over $B$ with dim $B< \infty $ . In
particular, generalized classical braided m-Lie algebras $sl_{q, f}(
GM_G(A),  F)$  and $osp_{q, t} (GM_G(A), M, F)$  of generalized
matrix algebra $GM_G(A)$ are constructed  and their connection with
special generalized matrix Lie superalgebra $sl_{s, f}( GM_{{\bf
Z}_2}(A^s), F)$  and orthosymplectic generalized matrix Lie super
algebra $osp_{s, t} (GM_{{\bf Z}_2}(A^s), M^s, F)$  are established.
The relationship between representations of braided m-Lie algebras
and their  associated algebras are established.

In Chapter \ref {c6}, we construct the Drinfeld double in symmetric
tensor category.

In Chapter \ref {c7}, we give the duality theorem
 for  Hopf
 algebras living in a symmetric braided tensor category. Using the result,
 we show that
    $$
(R \# H)\# H^{\hat *}   \cong M_n(R)  \hbox { \ \ \  as algebras }
$$
when $H$ is a finite-dimensional Hopf algebra living in a symmetric
braided vector category. In particular, the duality theorem holds
for any
  finite-dimensional
super-Hopf algebra $H.$

In chapter \ref {c8}, we give the relation between the decomposition
of coalgebras and comodules.

In Chapter \ref {c9}, we obtain that the global dimensions of $R$
and
 the crossed product
$R \# _\sigma H$ are the same; in the meantime,
 their weak dimensions are also the same  when
 $H$ is finite-dimensional semisimple and cosemisimple Hopf algebra.

In Chapter \ref {c10},
 we obtain that     the relation between $H$-radical of $H$-module algebra $R$ and
radical of $R \# H$.

In Chapter \ref {c11},
 we give that     the relation between $H$-radical of twisted
 $H$-module algebra $R$ and
radical of crossed product $R \# _\sigma  H$. We also obtain that if
$H$ is a finite-dimensional semisimple, cosemisimple, and either
commutative
   or cocommutative Hopf algebra, then $R$ is $H$-semiprime iff
      $R$ is semiprime iff
$R\#_{\sigma}H$ is semiprime.

In Chapter \ref {c12}, all  solutions of constant classical
Yang-Baxter
 equation (CYBE)  in Lie algebra $L$ with dim $L \le 3$  are obtained
  and
 the sufficient and  necessary conditions for  $(L, \hbox {[ \ ]},
 \Delta _r, r)$ to be a coboundary
 (or triangular ) Lie bialgebra  are given.

The strongly symmetric elements in $L\otimes L $ are found and they
are  all
 solutions of CYBE.

In Chapter \ref {c14}, we classify quiver Hopf algebras. This leads
to the classification of multiple Taft algebras as well as pointed
Yetter-Drinfeld modules and their corresponding Nichols algebras. In
particular, when the ground-field $k$ is the complex field and $G$
is a finite abelian group, we classify quiver Hopf algebras over
$G$, multiple Taft algebras over $G$ and Nichols algebras in
$^{kG}_{kG} {\cal YD}$. We show that the quantum enveloping algebra
of a complex semisimple Lie algebra is a quotient of a semi-path
Hopf algebra.

I would like to express my gratitude to my advisors Professor
Wenting Tong and Professor Yonghua Xu.

  I would like to express my gratitude to Professor Boxun Zhou
for his support.

I am very grateful to Professor Nanqing Ding,
 Professor Fanggui Wang, Professor Jianlong Chen, Professor
 Huixiang Chen, Professor Guoqiang Hu,
 Professor Jinqing Li and  Professor Minyi Wang  for
their help.

Finally, I wish to thank my parents, my wife and all of my friends
for their support.

 \tableofcontents

\part {Hopf Algebras Living in Braided Tensor Categories}

\chapter { Preliminaries} \label {c1}

We begin with the tensor category. We define the product ${\cal C}
\times {\cal D}$ of two category ${\cal C}$  and ${\cal D}$ whose
objects are pairs of objects $(U,V) \in (ob {\cal C}, ob {\cal D})$
and whose morphisms are given by
$$Hom _{{\cal C}\times {\cal D}}((V,W) (V',W'))  = Hom _{\cal C} (V,V')
\times Hom _{\cal D} (W, W'). $$ Let ${\cal C}$  be a category and
$\otimes $   be a  functor from ${\cal C} \times {\cal C} $  to
${\cal C}$. This means

(i) we have object $V \otimes W$  for any $V, W \in ob {\cal C};$

(ii) we have morphism $f \otimes g$  from $U\otimes V$  to $X\otimes
Y$ for any morphisms $f$  and $g$  from $U$ to $X$  and from $V$  to
$Y$;

(iii)  we have
$$(f \otimes g) (f' \otimes g') = (ff') \otimes (gg') $$
for any morphisms $f: U\rightarrow X, g : V\rightarrow Y, f' :
U'\rightarrow V$ and  $g' : V' \rightarrow V;$

(iv) $id_{U\otimes V} = id_U \otimes id _V.$

Let $\otimes \tau $ denote the functor from ${\cal C} \times {\cal
C} $  to ${\cal C}$ such that $(\otimes \tau)(U, V) =( V \otimes U)$
and  $(\otimes \tau )(f, g) = g \otimes f$,
 for any objects $U, V, X, Y$ in ${\cal C},$
and for any morphisms $f : U \rightarrow X$ and $g : V \rightarrow
Y.$

 An associativity constraint
$a$ for tensor $\otimes $ is a natural isomorphism
$$a: \otimes (\otimes \times id ) \rightarrow \otimes (id \times \otimes ).$$
This means that, for any triple $(U,V,W)$  of  objects of ${\cal
C}$, there is a morphism $a_{U,V,W} : (U\otimes V) \otimes W
\rightarrow U\otimes (V\otimes W)$  such that
$$a_{U',V',W'} ((f\otimes g) \otimes h) = (f\otimes (g\otimes h) ) a_{U,V,W} $$
for any morphisms $f, g$ and $h$  from $U$ to $U'$, from $V$  to
$V'$ and from $W$ to $W'$ respectively.

Let $I$ be an object of ${\cal C}$. If there exist  natural
isomorphisms
  $$ l : \otimes (I \times id )  \rightarrow  id  \hbox { \ \ \ and \ \ \ }
   r : \otimes (id \times I )  \rightarrow  id \ , $$
   then $I$  is called the unit object of ${\cal C}$
   with left unit constraint $l$
   and right unit constraint $r$.

\begin {Definition}  \label {0.1}
$({\cal C}, \otimes, I,  a, l, r )$ is called a tensor category if
${\cal C}$  is equipped with a tensor product $\otimes $, with a
unit object $I$, an associativity constraint   $a$ , a left unit
constraint $l$  and a right unit constraint $r$ such that the
Pentagon Axiom and the Triangle Axiom are satisfied, i.e.
$$(id _U \otimes a_{V,W,X}) a_{U, V\otimes W, X} (a_{U,V,W} \otimes
id _X) = a_{U,V,W\otimes X}a_{U\otimes V,W,X}$$ and
$$(id _V \otimes l_W) a_{V,I,W} = r_{V} \otimes id _W$$
for any $U, V, W, X \in ob {\cal C}.$

Furthermore, if  there exists a  natural isomorphism
$$C: \otimes \rightarrow \otimes \tau$$
such that the Hexagon Axiom holds, i.e.

$$a_{V,W,U} C_{U,V\otimes W} a_{U,V,W} =
(id _V \otimes C_{U,W}) a_{V,U,W}(C_{U,V}\otimes id _W)$$ and
$$a^{-1}_{W,U,V} C_{U\otimes V, W} a^{-1}_{U,V,W} =
( C_{U,W} \otimes id_ V ) a^{-1}_{U,W,V}(id_U \otimes C_{V,W}),$$
for any $U,V,W \in ob {\cal C},$ then  $({\cal C}, \otimes, I,  a,
l, r, C )$ is called a braided
 tensor category.
In this case, $C$ is called a braiding of ${\cal C}.$  If $C_{U,V} =
C_{V,U} ^{-1}$ for any $U, V \in ob {\cal C},$ then $({\cal C}, C)$
is called a symmetric braided tensor category or a symmetric tensor
category. Here functor $\tau: {\cal C} \times {\cal C} \rightarrow
{\cal C} \times {\cal C}$ is the flip functor defined by $\tau
(U\times V) = V\times U$ and $\tau (f \times g) = g \times f,$ for
any $U,V \in {\cal C},$  and morphisms $f$ and $g.$ Note that we
denote the braiding $C$ in braided tensor category $({\cal C},
\otimes, I,  a, l, r, C )$
 by $^{\cal C}C$ sometimes.

A tensor category $({\cal C}, \otimes , I,  a, l, r )$ is said to be
strict if $a, r, l$ are all identities in the category.
\end {Definition}

By Definition, the following Lemmas are clear
\begin {Lemma}  \label {0.-1} Let $\otimes $ be a functor from
${\cal C} \times {\cal C}$ to  ${\cal C}$ and $I$ an object in
${\cal C}$. Then  $({\cal C}, \otimes, I, a=id, l = id, r= id)$ is a
strick tensor category iff

(ST1):  $(U \otimes V)\otimes W =U \otimes (V\otimes W ) $ for any
$U,V,W \in ob {\cal C}$;

(ST2):  $I \otimes U = U \otimes I = U$ for any $U \in ob {\cal C}$.

 \end {Lemma}

 \begin {Lemma}  \label {0.0} Let $\otimes $ be a functor from
${\cal C} \times {\cal C}$ to  ${\cal C}$ and $I$ an object in
${\cal C}$. Assume that there exists a  natural isomorphism $C:
\otimes \rightarrow \otimes \tau$. Then  $({\cal C}, \otimes, I,
a=id, l = id, r= id, C)$ is a strick braided tensor category iff
$({\cal C}, \otimes, I, a=id, l = id, r= id)$ is a strick  tensor
category and

(SBT3):  $C_{U,V\otimes W}  = (id _V \otimes C_{U,W})
(C_{U,V}\otimes id _W)$ for any $U,V,W \in ob {\cal C}$;

(SBT4):
 $C_{U\otimes V, W}  =
( C_{U,W} \otimes id_ V ) (id_U \otimes C_{V,W})$ for any $U,V,W \in
ob {\cal C}.$

 \end {Lemma}

Let $({\cal C}, \otimes , I,  a, l, r )$ and  $({\cal D}, \otimes ,
I,  a, l, r )$ be two tensor categories. A tensor functor from
${\cal C}$  to ${\cal D}$  is a triple $(F, \mu _0, \mu ),$ where
$F$  is a functor from ${\cal C}$  to ${\cal D}$, $\mu _0$  is an
isomorphism  from $I$  to $F(I)$, and $$\mu (U,V): F(U) \otimes F(V)
\rightarrow F(U\otimes V)$$ is a family  of natural isomorphisms
indexed by all couples $(U, V)$  of objects of ${\cal C}$ such that
$$\mu (U, V \otimes W) (id _{F(U)} \otimes
\mu (V, W))a_{F(U), F(V), F(W)}$$
$$=F(a_{U,V,W})\mu (U\otimes V,W)
(\mu (U,V) \otimes id _{F(W)}),$$
$$l_{F(U)}= F(l_U) \mu (I,U)(\varphi _0 \otimes id _{F(U)})$$
and
$$r_{F(U)} = F(r_U)\mu (U, I) (id _{F(U)}\otimes \mu _0)$$
for any $U, V, W \in ob {\cal C}.$ Furthermore, if $({\cal C}, C)$
and  $({\cal D}, C')$  be two braided tensor categories, and
$$F(C_{V,V'})\mu (F(V), F(V')) =
\mu (F(V'), F(V))C_{F(V), F(V')}$$ for any $V , V' \in ob {\cal C},$
then $(F, \mu _0 , \mu )$  is called a braided tensor functor.

  We now construct a strict tensor category $ ( { \cal C}^{str},
  \underline {\otimes } )$ by a general braided tensor category
   $( {\cal  C }, \otimes, I, a,l, r, C )$
  as follows. The objects  of  $ { \cal C}^{str}$
  consist of all finite sequences $S = (V_1, V_2 , \cdots V_n)$ of $ { \cal C}$,
  including the empty sequence $\emptyset$. The integer $n$ is called the length
  of the sequence    $S = (V_1, V_2 , \cdots V_n)$. The length of empty sequence is zero.
 As in \cite [P 288 -- 291]{Ka95}, let  $$F ((V_1, V_2, \cdots ,V_n)) =  (( \cdots ( V_1 \otimes V_2) \otimes
  \cdots ) \otimes V_{n-1}) \otimes V_{n} , \hbox { \ \ \ \ \ \ \ }
   F(\emptyset ) = I, $$
       $$Hom _{{\cal C } ^{str}} (S'', S''') = Hom _{{\cal C } } (F(S''), F(S''')),$$
       $$\mu (\emptyset, \emptyset ) = l_{I} , \hbox {\ \ \ \ \ } \mu (S,\emptyset) = r _{F(S)}
   \hbox { \ \ \ \ \  } \mu (\emptyset , S) = l_{F(S)},$$
  $$\mu (S, (V)) = id _{F(S) \otimes V},$$
  \begin {eqnarray}\label {e15} \mu (S, S' \underline {\otimes } (V)) = (\mu (S,S') \otimes id_ V)
      a ^{-1}_{F(S), F(S'), V} \end {eqnarray}
for  $S, S', S'', S''' \in ob ( {\cal C } ^{str} )$ with $S \not=
\emptyset$ and $ S' \not= \emptyset$. It is shown that the two
tensor categories
 $ {\cal C }^{str}$ and ${\cal C }$ are equivalent  in \cite [Proposition XI.5.1]{Ka95} .
 If $S= (V)$, then we
                 write $S=V$ in short.

If, for any   pair of objects   $(S  , S')$  in ${\cal C} ^{str}$,
we  define
  $$F(\overline {C}_{S , S' }) = \mu (S',S)C_{F(S), F(S')} \mu(S,S')^{-1},$$
we shall show that $({\cal C}^{str}, \bar C)$ is a strict braided
tensor category.

\begin {Definition}  \label {0.0.1}
let $F:{ \mathcal C} \rightarrow { \mathcal  D}$be a functor.Then
$F$ is an equivalence of categories if there exist a functor
$G:\mathcal{D}\rightarrow \mathcal{C}$ and natural isomorphisms\\
  $$ \eta:id_{\mathcal{D}}\rightarrow FG\ \  \hbox{and} \ \ \theta:GF\rightarrow
id_{\mathcal{C}}.$$
\end {Definition}
we now give a useful criterion for a functor $F:
\mathcal{C}\rightarrow\mathcal{D}$ to be equivalence of categories.
let us first say that a functor
$F:\mathcal{C}\rightarrow\mathcal{D}$ is $essentially$ $subjective$
if, for any object $W$ of $\mathcal{D}$, there exists an object $V$
of $\mathcal{C}$ such that $F(V)\cong W$ in $\mathcal{D}$. It is
said to be $faithful$  ( resp.fully faithful ) if, for any couple
($V,V'$) of object of $\mathcal{C}$, the map $$ F: Hom
 (V, V') \rightarrow Hom  (F(V),
F(V'))$$ on morphisms is injective (resp.bijective).
\begin {Proposition} \label {0.0.2}
A functor $F:\mathcal{C}\rightarrow\mathcal{D}$ is an equivalence of
categories if  and only if $F$ is essentially surjective and fully
faithful.
 \end {Proposition}
{\bf Proof. }(a)Suppose that $F$ is an equivalence. Then there exist
a functor $G:\mathcal{D}\rightarrow\mathcal{C}$ and natural
isomorphisms $\eta:id_{\mathcal{D}}\rightarrow FG$ and
$\theta:GF\rightarrow id_{\mathcal{C}}$. The first isomorphism shows
that $W\cong F(G(W))$ for any object $W$ of $\mathcal{C}$. In other
words. $F$ is essentially surjective.Now consider a
morphism $f:V\rightarrow V'$ in $\mathcal{C}$. The square\\
$$\begin {array}{lcr}
&    \theta(V) & \\
GF(V)&\longrightarrow &V \\
\downarrow GF(f) & & \downarrow f \\
 GF(V')&\longrightarrow &V'\\
&    \theta(V') &
\end {array}  $$
commutes. It results that if $F(f)=F(f')$, hence $GF(f)=GF(f')$,
then we have $f=f'$. Therefore, the functor $F$ is faithful.Using
the natural isomorphism $\eta$ in a similar way,we prove that $G $
is faithful too. Now consider a morphism $g: F(V \rightarrow F(v')$.
Let us show
that $g=F(f)$ where $f=\theta(v')\circ G(g)\circ\theta(V)^{-1}$. Indeed,\\
$$\theta(v')\circ GF(g)\circ \theta(V)^{-1}=f=\theta(v')\circ G(g)\circ \theta(V)^{-1}.$$
Therefore $GF(f)=G(g)$. As $G$ is faithful, we get $g=F(f)$. This
proves that $F$ is fully faithful.\\

(b) Let $F$ is essentially surjective and fully faithful functor.
For any object $W$ in $\mathcal{D}$, we choose an object $G(W)$ of
$\mathcal{C}$ and an isomorphism $\eta(W):W\rightarrow FG(W)$ in
$\mathcal{D}$. If $g:W\rightarrow W'$ is a morphism of
$\mathcal{D}$. We may
consider\\
$$\eta(W')\circ g\circ \eta(W)^{-1}:FG(W)\rightarrow FG(W').$$\\
Since $F$ is fully faithful, there exists a unique morphism $G(g)$ from $G(W)$ to $G(W')$ such that \\
$$FG(g)=\eta(W')\circ g\circ \eta(W)^{-1}:FG(W)\rightarrow FG(W').$$\\
One checks easily that this define a functor $G$ from $\mathcal{D}$
into $ \mathcal{C}$ and that $\eta:id_{\mathcal{D}}\rightarrow FG$
is a natural isomorphism. In order to show that $F$ and $G$ are
equivalences of categories, we have only to find a natural
isomorphism $\theta:GF\rightarrow id_{\mathcal{C}}$. We define
$\theta(V):GF(V)\rightarrow V$ for any object $V\in Ob(\mathcal{C}$)
as the unique morphism such that $F(\theta(V))=\eta(F(V))^{-1}$. It
is easily checked that this
formula defines a natural isomorphism. $\Box $\\
\begin {Proposition}\label {0.0.3}
The categories $\mathcal{C}^{str}$ and $\mathcal{C}$ are equivalent.
\end {Proposition}
{\bf Proof. }The map $F$ defined above on objects of
$\mathcal{C}^{str}$ extends to a functor
$F:\mathcal{C}^{str}\rightarrow\mathcal{C}$ which is the identity on
morphisms, hence fully faithful. As any object in $\mathcal{C}$ is
clearly isomorphic to the image under $F$ of a sequence of length
one, we see that $F$ is essentially surjective. This proves the
proposition in view of Proposition 1.0.5. Observe that $G(V)=V$
defines a functor $G:\mathcal{C}\rightarrow\mathcal{C}^{str}$ which
is the inverse equivalence to $F$. Indeed,we have$FG=
id_{\mathcal{C}}$. $\theta:GF\rightarrow id_{\mathcal{C}^{str}}$ via
the natural isomorphism $$\theta(S)=id_{F(S)}:GF(S)\rightarrow S. \
\Box$$

\begin {Lemma}\label {0.0.4}
If $S,S',S''$ are objects on $\mathcal{C}^{str}$, we have \\
$$\mu(S,S'\underline{\otimes} S'')\circ(id_{S}\otimes\mu(S',S''))\circ
a_{F(S),F(S'),F(S'')}=\mu(S\underline{\otimes} S',S'')\circ
(\mu(S,S')\otimes id_{S''}).$$
\end {Lemma}
{\bf Proof. }We proceed by induction on the length of $S''$. If
$S''=\emptyset$, we have \\
\begin {eqnarray*}
\mu(S,S')(id_{S}\otimes\mu(S',\emptyset)) a_{F(S),F(S'),I}
&=&\mu(S,S')(id_{F(S)}\otimes r_{F(S'')})a_{F(S),F(S'),I}\\
&=&\mu(S,S')r_{F(S)\otimes F(S')} \\
&=&r_{F(S\underline{\otimes} S')}(\mu(S,S')\otimes id_{I}) \\
&=&\mu(S\underline{\otimes}S',\emptyset)(\mu(S,S')\otimes id_{I}).
\end {eqnarray*}
 The first and last equalities are by definition,the second
one by \cite [Lemma XI.2.2]{Ka95}
 and the third one by naturality  of r.\\
Let $V$ be an object of the category. Let us prove that the equality
of \cite [Lemma XI. 5.2]{Ka95} for the triple($S,S',S''$) implies
the
equality for ($S,S',S''\underline{\otimes}(V)$). We have \\
\begin {eqnarray*}
&&\mu(S,S'\underline{\otimes}S''\underline{\otimes}(V))(id_{S}\otimes\mu(S',S''\underline{\otimes}(V)))
a_{F(S),F(S'),F(S''\underline{\otimes}(V))}\\
&=&(\mu(S,S'\underline{\otimes} S'')\otimes
id_{V})a^{-1}_{F(S),F(S'\underline{\otimes} S''),V}
(id_{S}\otimes\mu(S',S'')\otimes id_{V})\\
&&(id_{S}\otimes a^{-1}_{F(S'),F(S''),V})
a_{F(S),F(S'),F(S''\underline{\otimes}(V))}\\
&=&(\mu(S,S'\underline{\otimes} S'')\otimes
id_{V})(id_{S}\otimes\mu(S',S''\otimes id_{V}))
a^{-1}_{F(S),F(S')\underline{\otimes} F(S''),V}\\
&&(id_{S}\otimes a^{-1}_{F(S'),F(S''),V})
a_{F(S),F(S'),F(S''\underline{\otimes}(V))}\\
&=&(\mu(S,S'\underline{\otimes} S'')\otimes
id_{V})(id_{S}\otimes\mu(S',S''\otimes id_{V}))\\
&& (a_{F(S),F(S'), F(S'')}\otimes id_{V})
a^{-1}_{F(S)\otimes F(S'),F(S''),V}\\
&=&(\mu(S\underline{\otimes} S',S'')\otimes id_{V})(\mu(S,S')\otimes
id_{S''}\otimes
id_{V}))a^{-1}_{F(S)\otimes F(S'),F(S''),V}\\
&=&(\mu(S\underline{\otimes} S',S'')\otimes
id_{V})a^{-1}_{F(S)\otimes F(S\underline{\otimes} S'),F(S''),V}
(\mu(S,S')\otimes id_{S''}\otimes id_{V})\\
&=&(\mu(S\underline{\otimes} S',S''\underline{\otimes} (V))
(\mu(S,S')\otimes id_{S''\underline{\otimes}(V)}).
\end {eqnarray*}
The first and last equalities follow from (\ref {e15}), the second
and fifth ones from the naturality of the associativity constraint,
 the third from the Pentagon Axiom,
and the fourth one from the induction hypothesis.$\Box$
\begin {Theorem} \label {0.0.5}Equipped with this tensor
product $\mathcal{C}^{str}$ is a strict tensor category. The
categories $\mathcal{C}$ and $\mathcal{C}^{str}$ are tensor
equivalent.
\end {Theorem}
{\bf Proof. }It is easy to check that $\underline{\otimes}$ is a
functor. This functor is strictly associative by construction.
Therefore $\mathcal{C}^{str}$ is a strict tensor
category.\\

In order to prove that it is tensor equivalent to $\mathcal{C}$, we
have to exhibit tensor functors and natural isomorphisms. We first
claim that the triple($F,id_{I},\mu $) is a tensor functor from
$\mathcal{C}^{str}$ to $\mathcal{C}$ where $\mu $ is the natural
isomorphism defined above. By  \cite [Lemma XI.5.2]{Ka95} and \cite
[Proposition XI.5.2] {Ka95}, the functor $G$ of the proof of  is a
tensor functor. Finally, the natural isomorphism $\theta$ is a
natural tensor isomorphism.$\Box$

\begin  {Theorem} \label {0.2}  If   $({\cal C}, C)$  is a  braided tensor
 category, then  so is $({\cal C}^{str}, \bar C)$, and the two braided
 tensor categories are equivalent.
\end {Theorem}

{\bf Proof. }  We  first show that
 $$  \overline {C}_{S,S' \underline {\otimes} S'' } =
     (id_ {S'} \underline {\otimes}
\overline  C _{S , S''}) (\overline C _{S,S'} \underline {\otimes}
id _{S''})$$ and
$$ \overline {C}_{S \underline {\otimes} S' , S'' } =
      (\overline C _{S,S''} \underline {\otimes} id _{S'})
                      (id_ S \underline {\otimes}
\overline  C _{S' , S''}). $$

We only check the first relation. The second one can similarly be
shown.

We see that

\begin {eqnarray*}
 F &(& \hbox { the right side of the first relation } ) \\
  &=& \mu (S', S'' \underline \otimes S)(id _{F(S')} \otimes F( \bar C_{S,S''}))
\mu (S', S\underline \otimes S'')^{-1}  \\
&{}&   \mu (S' \underline \otimes S, S'')(F(\bar C_{S,S'} )\otimes
id _{F(S'')})
 \mu ( S\underline \otimes S', S'')^{-1}   \\
&=& \mu (S', S'' \underline \otimes S)(id _{F(S')} \otimes \mu
(S'',S) C_{F(S),F(S'')}) \mu (S,  S'')^{-1})
\mu (S', S \underline \otimes S'')^{-1} \mu (S'\underline { \otimes } S,S'') \\
&{}&  ( \mu (S', S) C _ {F(S), F(S')} \mu (S,S')^{-1}  \otimes id _{
F(S'')})
\mu ( S \underline \otimes S', S'')^{-1}                  \\
&=& \mu  (S', S'' \underline \otimes S)( id _{F(S')} \otimes \mu
(S'', S) ) ( id _{F(S')}  \otimes C_{F(S),F(S'')})
 ( id _{F(S')} \otimes \mu ( S,  S'')^{-1} ) \\
 &{}&   \mu (S', S \underline \otimes S'')^{-1}
\mu ( S'\underline \otimes S, S'') \\
&{}& (\mu (S', S)  \otimes id _{F(S'')})
 ( C_{F(S), F(S')} \otimes id _{F(S'')}) (\mu (S,S')^{-1} \otimes id _{F(S'')})
\mu ( S \underline \otimes S' , S'')^{-1}                   \\
&=&  \mu ( S'\underline \otimes S'',S)(\mu (S',S'') \otimes id
_{F(S)}) a_{F(S'),F(S''),F(S)}^{-1}(id _{F(S')} \otimes C_{F(S),
F(S'')})
a_{F(S'),F(S),F(S'')} \\
&{}& (C_{F(S),F(S')} \otimes id_{F(S'')}) a
_{F(S),F(S'),F(S'')}^{-1}(id_{F(S)} \otimes \mu (S',S'')^{-1})
\mu (S, S' \underline \otimes S'')^{-1} \\
&{}&
{ \hbox { \ \ \ by  \cite [ lemma XI.5.2] {Ka95} } }                        \\
&=&   \mu ( S' \underline \otimes S'',S)(\mu (S',S'')\otimes id
_{F(S)})
 C_{F(S),F(S') \otimes F(S'')}(id _{F(S)} \otimes \mu (S',S'')^{-1})
 \mu (S, S'\underline \otimes S'')^{-1}                         \\
&=& \mu ( S' \underline \otimes S'',S)C_{F(S), F(S' \underline
{\otimes } S'')}
\mu (S, S'\underline \otimes S'')^{-1}                            \\
&=&  F (  \hbox { the left side of the first relation. })
\end {eqnarray*}

Therefore the first relation holds.

Next we can obtain that $\bar C$ is  natural  by the diagram below:
$$\begin {array}{lcccccr}
{}&{}&{}&            C_{F(S),F(S')} &{}&{}&{} \\
F(S) \otimes  F(S')&{ \ }& { \ }&\longrightarrow& {} & { \ } &
F(S') \otimes  F(S) \\
{}  & \searrow \mu (S,S') & {} &  {} &{} & \swarrow  \mu (S',S) &{} \\
{}&{}& F(S \underline {\otimes} S')& \stackrel {F(\bar C_{S,S'})}{
 \longrightarrow }&   F(S' \underline {\otimes} S)&  { \ } \\
 f\otimes g \downarrow & {}& \downarrow  f \underline {\otimes} g & {}& \downarrow
 g \underline{\otimes} f &{}& \downarrow  g\otimes f \\
{} &{}& F(T \underline {\otimes} T')& \stackrel {F(\bar C_{T,T'})}{
 \longrightarrow }&   F(T' \underline {\otimes} T)& {}&{} \\
 {}  & \nearrow \mu (T,T') & {} &  {} &{} & \nwarrow  \mu (T',T) &
{} \\
F(T) \otimes  F(T')& { \ }& { \ }&  \longrightarrow & { \ }& { \ } &
F(T') \otimes  F(T) \\
{}&{}&{}&            C_{F(T),F(T')} &{}&{} &{}
\end {array}.  $$

Finally, we can easily check that $F$ is a braided tensor functor
and $({\cal C}^{str},\bar C)$ and $({\cal C},C)$ are braided tensor
equivalent categories.
\begin{picture}(8,8)\put(0,0){\line(0,1){8}}\put(8,8){\line(0,-1){8}}\put(0,0){\line(1,0){8}}
\put(8,8){\line(-1,0){8}}\end{picture}

Therefore, we can view every braided tensor category as
 a strict braided tensor category and use braiding diagrams freely
 ( we can also see  \cite [exercise 5 P337] {Ka95}).

\begin {Example} \label {0.4} (The tensor category of vector spaces )
The most fundamental   example of a tensor category is given by the
category ${\cal C} = {\cal V}ect (k)$ of
 vector
spaces over  field $k$. ${\cal V}ect(k)$  is equipped with tensor
structure for which $\otimes $ is the tensor product of the vector
spaces over $k$, the unit object $I$ is the ground field $k$ itself,
and the associativity constraint and unit constraint are the natural
isomorphisms
$$a_{U, V, W}((u \otimes v) \otimes w)= u \otimes (v \otimes w)  \hbox { \ \ and \ \ }
l_V(1 \otimes v) = v = r_V (v \otimes 1)$$ for any vector space $U,
V, W$ and $u \in U, v\in V, w\in W.$

Furthermore, the most fundamental example of a braided tensor
category is given by the tensor category  ${\cal V}ect (k)$, whose
braiding is usual twist map from $U \otimes V $  to $V\otimes U$
defined by sending $a \otimes b$ to $b \otimes a$ for any $a\in U,
b\in V$.
\end {Example}

Now we define some notations. If $f$ is a morphism from $U$ to $V$
and $g$  is a morphism from $V$ to $W$, we denote the composition
$gf$ by
\[
\begin{tangle}
\object{U}\\
\vstr{50}\id\\
\O {gf}\\
\vstr{50}\id\\
\object{W}\\
\end{tangle}
\step=\step
\begin{tangle}
\object{U}\\
\O {f}\\
\O {g}\\
\object{W}\\
\end{tangle}\ \ \ .
\]

We usually  omit $I$ and  the morphism $id _I$ in any diagrams. In
particular, If $f$ is a morphism from $I$ to $V$, $g$ is a morphism
from $ V$ to $I$,we denote $f$ and  $g$ by
\[
\begin{tangle}
\Q f\\
\object{V}\\
\end{tangle}
\step \hbox { , } \step
\begin{tangle}
\object{V}\\
\QQ g\\
\end{tangle} \ \ \ \ .
\]

 If
$f$ is a morphism from $U \otimes V$ to $P$, $g$ is a morphism from
$U \otimes V$ to $I$ and $\zeta$  is a morphism from $U\otimes V$
to $V\otimes U$, we denote $f$, $g$ and $\zeta$ by
\[
\begin{tangle}
\object{U}\step[2]\object{V}\\
\tu f\\
\step\object{P}\\
\end{tangle}
\step,\step
\begin{tangle}
\object{U}\step[2]\object{V}\\
\coro g\\
\end{tangle}
\step,\step
\begin{tangle}
\object{U}\step[2]\object{V}\\
\ox \zeta\\
\object{V}\step[2]\object{U}\\
\end{tangle}
\step .
\]

   In particular, we denote the braiding
$C_{U,V}$  and its inverse $C_{U,V}^{-1}$ by
\[
\begin{tangle}
\object{U}\step[2]\object{V}\\
\x\\
\object{V}\step[2]\object{U}\\
\end{tangle}
\step,\step
\begin{tangle}
\object{V}\step[2]\object{U}\\
\xx\\
\object{U}\step[2]\object{V}\\
\end{tangle}
\step . \\
\]

$\xi $ is called  an $R$-matrix of ${\cal C}$ if
 $\xi$ is a natural
isomorphism from $\otimes $ to $\otimes \tau $ and  for any $U,V,W
\in ob {\cal C}$, the  Yang-Baxter equation of ${\cal C}$:

(YBE):
\[
\begin{tangle}
\object{U}\step[2]\object{V}\step[2]\object{W}\\
\ox \xi \step [2] \id\\
\id \step [2] \ox \xi \\
\ox \xi \step [2] \id \\
\object{W}\step [2]\object{V}\step[2]\object{U}\\
\end{tangle}
\ \ = \
\begin{tangle}
\object{U}\step[2]\object{V}\step[2]\object{W}\\
\id \step [2] \ox \xi \\
\ox \xi \step [2] \id\\
\id \step [2] \ox \xi \\
\object{W}\step [2]\object{V}\step[2]\object{U}\\
\end{tangle}\
\] holds.  In particular, the above equation is called the Yang-Baxter
equation on $V$ when $U = V=W.$

\begin {Lemma}\label {0.3}

(i) The braiding $C$ of braided tensor category  $ ({ \cal C}, C)$
is an $R$-matrix of $ { \cal C}$, i.e.
\[
\begin{tangle}
\object{U}\step[2]\object{V}\step[2]\object{W}\\
\x \step [2] \id\\
\id \step [2] \x \\
\x \step [2] \id \\
\object{W}\step [2]\object{V}\step[2]\object{U}\\
\end{tangle}
\ \ = \
\begin{tangle}
\object{U}\step[2]\object{V}\step[2]\object{W}\\
\id \step [2] \x \\
\x \step [2] \id\\
\id \step [2] \x \\
\object{W}\step [2]\object{V}\step[2]\object{U}\\
\end{tangle}\
\] holds.

(ii)
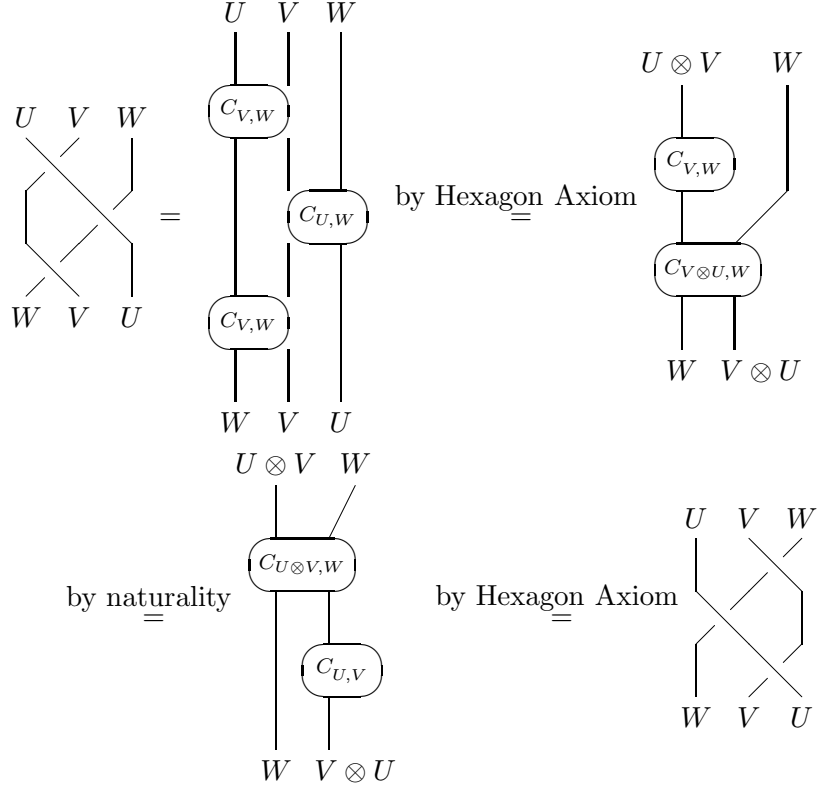
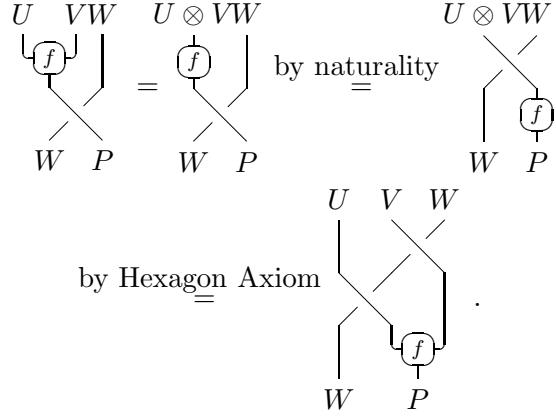
\[
\begin{tangle}
\object{U}\step[2]\object{V}\step[2]\object{W}\\
\id\step[2]\x\\
\x\step[2]\id\\
\id\step[2]\tu f\\
\object{W}\step[3]\object{P}\\
\end{tangle}
\step=\step
\begin{tangle}
\object{U}\step[2]\object{V}\step\object{W}\\
\tu f\step\id\\
\step\x\\
\step\object{W}\step[2]\object{P}\\
\end{tangle}
\step , \step
\begin{tangle}
\object{U}\step[2]\object{W}\step[2]\object{V}\\
\x\step[2]\id\\
\id\step[2]\tu f\\
\object{W}\step[3]\object{P}\\
\end{tangle}
\step=\step
\begin{tangle}
\object{U}\step[2]\object{W}\step\object{V}\\
\id\step[2]\hxx\\
\tu f\step\id\\
\step\x\\
\step\object{W}\step[2]\object{P}\\
\end{tangle} \ \ .
\]
\end {Lemma}

{\bf Proof.} (i)

\[
\begin{tangle}
\object{U}\step[2]\object{V}\step[2]\object{W}\\
\x \step [2] \id\\
\id \step [2] \x \\
\x \step [2] \id \\
\object{W}\step [2]\object{V}\step[2]\object{U}\\
\end{tangle} \ \ = \ \
\begin{tangle}

\step\object{U}\step[2]\object{V}\step[2]\object{W}\\
\step \id \step [2] \id \step [2] \id \\

\obox 3{C^{}_{V,W}}    \step [2] \id\\

 \step\id \step [2] \id \step [2] \id \\

\step\id \step [2] \obox 3 {C_{U,W}} \\

\step \id \step [2] \id \step [2] \id \\

\obox 3 {C_{V,W}} \step [2] \id \\

\step \id \step [2] \id \step [2] \id \\
\step\object{W}\step [2]\object{V}\step[2]\object{U}\\
\end{tangle} \ \ \stackrel {\mbox {by Hexagon Axiom }} {=}\begin{tangle}

\step\object{U \otimes V} \step[4]\object{W}\\
\step \id  \step [4] \id \\

\obox 3{C^{}_{V,W}}    \step [2] \id\\

 \step \id \step [3] \ne2 \step [2]\\

 \obox 4 {C_{V\otimes U,W}} \\

\step \id \step [2] \id \step [2] \\
\step\object{W}\step [3]\object{V \otimes U}\step[2]\\
\end{tangle} \ \
\]

\[\ \ \stackrel {\mbox {by naturality  }} {=}\begin{tangle}

\step\object{U \otimes V} \step[3]\object{W}\\
\step \id  \step [2] \ne1 \\

\obox 4 {C_{U\otimes V,W}}
   \step [2] \\

 \step \id \step [2] \id \step [2]\\

\step \id \step [1]\obox 3{C^{}_{U,V}}  \\

\step \id \step [2] \id \step [2] \\
\step\object{W}\step [3]\object{V \otimes U}\step[2]\\
\end{tangle} \ \
\stackrel {\mbox {by Hexagon Axiom }} {=}
\begin{tangle}
\object{U}\step[2]\object{V}\step[2]\object{W}\\
\id \step [2] \x \\
\x \step [2] \id\\
\id \step [2] \x \\
\object{W}\step [2]\object{V}\step[2]\object{U}\\
\end{tangle}\ \ \ .
\]

(ii)
 We have
\[
\begin{tangle}
\object{U}\step[2]\object{V}\step\object{W}\\
\tu f\step\id\\
\step\x\\
\step\object{W}\step[2]\object{P}\\
\end{tangle}
\step=\step
\begin{tangle}
\object{U\otimes V}\step[2]\object{W}\\
\O f\step[2]\id\\
\x\\
\object{W}\step[2]\object{P}\\
\end{tangle}
 \ \ \stackrel {\mbox {by naturality  }} {=} \step
\begin{tangle}
\object{U\otimes V}\step[2]\object{W}\\
\x\\
\id\step[2]\O f\\
\object{W}\step[2]\object{P}\\
\end{tangle} \]
\[
\stackrel {\mbox {by Hexagon Axiom }} {=}
\begin{tangle}
\object{U}\step[2]\object{V}\step[2]\object{W}\\
\id\step[2]\x\\
\x\step[2]\id\\
\id\step[2]\tu f\\
\object{W}\step[3]\object{P}\\
\end{tangle}\ \ \ .
\]
 Similarly, we can show the second equation. $\Box$

Now we give some concepts as follows: Assume  that $H, A \in ob \
{\cal C},$  and

\begin {eqnarray*}
\alpha : H \otimes A \rightarrow &A& , \hbox { \ \ \ \ }
\beta : H \otimes A \rightarrow H,    \\
\phi :  A \rightarrow H \otimes  &A&   , \hbox { \ \ \ \ }
\psi : H  \rightarrow H \otimes A,  \\
m _H : H \otimes H \rightarrow &H& , \hbox { \ \ \ \ }
m_A : A \otimes A \rightarrow A,    \\
\Delta _H :  H \rightarrow H \otimes  &H&   , \hbox { \ \ \ \ }
\Delta _A : A  \rightarrow A \otimes A, \\
\eta _H : I  \rightarrow &H& , \hbox { \ \ \ \ }
\eta _A : I \rightarrow A,    \\
\epsilon  _H :  H  \rightarrow   &I&   , \hbox { \ \ \ \ }
\epsilon  _A  : A  \rightarrow I.     \\
\end {eqnarray*}
are morphisms in ${\cal C}$.

$(A, m_A, \eta _A )$  is called an algebra living in ${\cal C}$, if
 the following conditions are
satisfied:
\[
\begin{tangle}
\object{A}\step\object{A}\step[2]\object{A}\\
\id\step\tu m\\
\tu m\\
\step\object{A}\\
\end{tangle}
\step=\step
\begin{tangle}
\object{A}\step[2]\object{A}\step\object{A}\\
\tu m\step\id\\
\step\tu m\\
\step[2]\object{A}\\
\end{tangle}
\step\step,\step\step
\begin{tangle}
\step[2]\object{A}\\
\Q {\eta_A}\step[2]\id\\
\tu m\\
\step\object{A}\\
\end{tangle}
\step=\step
\begin{tangle}
\object{A}\\
\id\step[2]\Q {\eta_A}\\
\tu m\\
\step\object{A}\\
\end{tangle}
\step=\step
\begin{tangle}
\object{A}\\
\id\\
\id\\
\object{A}\\
\end{tangle}
\quad .
\]
\noindent In this case, $\eta _A$ and $m_A$  are called  unit and
multiplication of $A$ respectively.

$(H, \Delta _H, \epsilon _H)$  is called a coalgebra living in
${\cal C},$ if the following conditions are satisfied:
\[
\begin{tangle}
\step\object{H}\\
\td \Delta\\
\id\step\td \Delta\\
\object{H}\step\object{H}\step[2]\object{H}\\
\end{tangle}
\step=\step
\begin{tangle}
\step[2]\object{H}\\
\step\td \Delta\\
\td \Delta\step\id\\
\object{H}\step[2]\object{H}\step\object{H}\\
\end{tangle}
\step\step,\step\step
\begin{tangle}
\step\object{H}\\
\td \Delta\\
\QQ {\varepsilon_H}\step[2]\id\\
\step[2]\object{H}\\
\end{tangle}
\step=\step
\begin{tangle}
\step\object{H}\\
\td \Delta\\
\id\step[2]\QQ {\varepsilon_H}\\
\object{H}\\
\end{tangle}
\step=\step
\begin{tangle}
\object{H}\\
\id\\
\id\\
\object{H}\\
\end{tangle}
\quad .
\]
 In this case, $\epsilon _H$ and $\Delta _H$  are called counit
and comultiplication of $H$ respectively.

If $A$  is an algebra and $H$  is a coalgebra, then $Hom_{\cal C}
(H, A)$  becomes  an algebra under the convolution product \[
f*g=\quad
\begin{tangle}
\step\object{H}\\
\td \Delta\\
\O f\step[2]\O g\\
\tu m\\
\step\object{A}\\
\end{tangle}\ \ \ .
\]
 \noindent and its unit element $\eta
= \eta _A \epsilon _H.$ If $S$ is  the inverse of $id _H$ in $Hom
_{\cal C} (H,H)$, then $S$ is called  antipode of $H$.

If $(H, m_H, \eta _H)$ is  an algebra,  and $(H, \Delta _H, \epsilon
_H)$ is a coalgebra living in ${\cal C}$, and the following
condition is satisfied: \[
\begin{tangle}
\object{H}\step[2]\object{H}\\
\tu m\\
\td \Delta\\
\object{H}\step[2]\object{H}\\
\end{tangle}
\step=\step
\begin{tangle}
\step\object{H}\step[3]\object{H}\\
\td \Delta\step\td \Delta\\
\id\step[2]\hx\step[2]\id\\
\tu m\step\tu m\\
\step\object{H}\step[3]\object{H}\\
\end{tangle}
\step,\step
\begin{tangle}
\object{H}\step[2]\object{H}\\
\tu m\\
\step\QQ \varepsilon\\
\end{tangle}
\step=\step
\begin{tangle}
\object{H}\step[2]\object{H}\\
\id\step[2]\id\\
\QQ \varepsilon\step[2]\QQ \varepsilon\\
\end{tangle}
\step,\step
\begin{tangle}
\step\Q \eta\\
\td \Delta\\
\end{tangle}
\step=\step
\begin{tangle}
\Q \eta\step[2]\Q \eta\\
\end{tangle}
\step,\step
\begin{tangle}
\Q \eta\\
\QQ \varepsilon\\
\end{tangle}
\step=\step
\begin{tangle}\object{I}\\
\id\\
\id\\
\object{I}
\end{tangle}
\step.\step
\]

\noindent then $H$ is called a bialgebra living in ${\cal C}$. If
$H$ is a bialgebra and there is  an inverse $S$ of $id _H$  under
convolution product in $Hom _{\cal C} (H, H)$, then $H$  is called a
Hopf algebra living in ${\cal C},$ or a braided Hopf algebra.

 For $V \in$  object ${\cal C}$, if there exists
an object $U$ and morphisms :
$$d_V : U\otimes V \rightarrow I   \hbox { \ \ \ and \ \ \ }
b_V : I \rightarrow V \otimes U$$ in ${\cal C}$ such that $$ (d_V
\otimes id _U) (id _U \otimes b_V) =id_ U \hbox {\ \ \ and \ \ \ }
(id_V \otimes d _V)(b_V \otimes id _V)= id_ V,$$ then $U$  is called
a left duality of $V,$ written as $V^{ *}.$ In this case, $d_V$ and
$b_V$ are called evaluation morphism and coevaluation morphism of
$V$ respectively.

If $U$ and $V$ have left duality $U^*$ and $V^*$, respectively, then
 $V^* \otimes U^*$ is  a left duality of
$U\otimes V$. Its  evaluation and coevaluation are
$$ d_{U\otimes V}=d_V (id _{V^*}\otimes d_U \otimes id _V ) ,
\ \ \ \ b_{U\otimes V}=(id _{U}\otimes b_V \otimes id_{U^*})b_U,
$$ respectively.

The multiplication, comultiplication, evaluation morphism and
coevaluation  are  usually denoted by
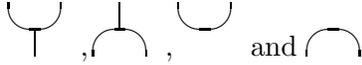
\[\begin{tangle}\cu  \   \ ,
 \ \ \cd \ \ , \ \ \ev \ \ \hbox {and } \ \  \ \ \coev \end {tangle}\] in short, respectively.

Let us first define the transpose $f^*: V^* \rightarrow U^*$ of
 morphism $f : U \rightarrow V$ by
$f^* = (d \otimes id _{U^*}) (id _{V^*} \otimes f \otimes id _{U^*})
(id _{V^*} \otimes b).$

\begin {Lemma}\label {2.1.9'}

 Let$\step\xi=\begin{tangle}
\step\object{P}\step[4]\object{Q}\\
\td f\step[2]\td g\\
\id\step[2]\x\step[2]\id\\
\tu u\step[2]\tu v\\
\step \object{Y}\step[4]\object{Z}
\end{tangle}$ be a morphism from $P\otimes Q$ to $Y\otimes Z$ and
all objects concerned have left duals . Then
$\step\xi^*=\begin{tangle}
\step\object{Z^*}\step[4.5]\object{Y^*}\\
\td {v^*}\step[2]\td {u^*}\\
\id\step[2]\x\step[2]\id\\
\tu {g^*}\step[2]\tu {f^*}\\
\step \object{Q^*}\step[4]\object{P^*}
\end{tangle}$
\end {Lemma}

{\bf Proof. } We see that

\[\begin{tangle}
\step\object{Z^*}\step[4.5]\object{Y^*}\\
\td {v^*}\step[2]\td {u^*}\\
\id\step[2]\x\step[2]\id\\
\tu {g^*}\step[2]\tu {f^*}\\
\step \object{Q^*}\step[4]\object{P^*}
\end{tangle}
\step=\step
\begin{tangle}
\object{Z^*}\step[9]\object{Y^*}\\
\id\step[2]\step\coev\step[2]\step[2]\id\step[2]\step\coev\\
\id\step[2]\ne2\coev\d\step[2]\step\id\step[2]\ne2\coev\d\step[2]\\
\id\step\tu v \step[2]\id\step\nw2\step[2]\id\step\tu
u\step[2]\id\step\id\step[2]\step\coev\\
\ev\step[2]\step\id\step[2]\step\nw2\ev\step[2]\ne2\step\id\step[2]\td
f\step\id\\
\step[2]\step[2]\step\id\step[2]\step[2]\step\x\step[2]\step\ev\step\dd\step\id\\
\step[5]\id\step[2]\step[2]\ne2\step\coev\step[-1]\nw3\step[3]\step\dd\step[2]\id\\
\step[2]\step[2]\step\d\step[2]\id\step[2]\td
g\step\id\step[2]\ev\step[2]\step\id\\
\step[6]\nw2\step\ev\step\ne2\step\id\step[7]\id\\
\step[8]\ev\step[3]\id\step[7]\id\\
 \step[13] \object{Q^*}\step[7]\object{P^*}
\end{tangle}
\]

\[
\step=\step
\begin{tangle}
\object{Z^*}\step[8]\object{Y^*}\\
\id \step[2]\step\coev\step[2]\step\id\\
\id\step[2]\ne2\coev\d\step[2]\id\step[2]\step[2]\step\coev\\
\id\step\tu v \step\dd\step\x\step[2]\step[2]\td f\step\id\\
\ev\step\dd\step\dd\step[2]\nw2\step[2]\step\id\step[2]\id\step\id\\
\step[2]\dd\step[2]\id\step[2]\step[2]\step\x\step[2]\id\step\id\\
\step[2]\id\step[2]\step\id\step[2]\step\coev\d\step\ev\step\id\\
\step[2]\id\step[2]\step\id\step[2]\ne2\coev\d\d\step[2]\step\id\\
\step[2]\id\step[2]\step\id\step\tu
u\step[2]\id\step\id\step\d\step[2]\id\\
\step[2] \id\step[2]\step\ev\step[2]\ne2\step\ev\step[2]\id\\
\step[2]\id\step[2]\step[2]\step\ne4\step[-1]\coev\step[2]\step[2]\step[2]\id\\
\step[2]\id\step[2]\id\step[2]\td
g\step\id\step[2]\step[2]\step[2]\id\\
\step[2]\d\step\ev\step\ne2\step\id\step[2]\step[2]\step[2]\id\\
\step[2]\step\Ev\step[2]\step[2]\step\id\step[2]\step[2]\step[2]\id\\
\step[9]\object{Q^*}\step[6]\object{P^*}
\end{tangle}
\step=\step
\begin{tangle}
\object{Z^*}\step[5]\object{Y^*}\\
\id\step[2]\step[2]\step\id\\
\id\step\coev\step[2]\id\\
\id\step\id\step[2]\xx\step[2]\step[3]\step\coev\\
\id\step\id\step[2]\id\step[2]\nw3\step[2]\step[2]\td f\step\id\\
\id\step\id\step[2]\id\step[2]\step[2]\step\xx\step[2]\id\step\id\\
\id\step\id\step[2]\id\step[2]\step[3]\id \step[2]\ev\step\id\\

\id\step\id\step[2]\id\step[2]\step[2]\ne2\step\coev \step[2]\id\\
\id\step\id\step[2]\d\step[2]\id\step[2]\td g
\step\id\step[2]\id\\
\d\d\step[2]\d\step\tu u\step\dd\step\id\step[2]\id\\
\step\d\nw2\step[2]\ev\step\ne3\step[2]\id\step[2]\id\\
\step[2]\d\step\tu v\step[2]\step[2]\step\id\step[2]\id\\
\step\step[2]\ev\step[2]\step[2]\step[2]\id\step[2]\id\\
\step[11]\object{Q^*}\step[2]\object{P^*}
\end{tangle}
\]

\[
\step=\step
\begin{tangle}
\object{Z^*}\step[2]\object{Y^*}\\
\id\step[2]\id\step[2]\step[2]\step\coev\\
\id\step[2]\id\step[2]\step[2]\td f\step\d\\
\id\step[2]\id\step[2]\step[2]\id\step[2]\id\step[2]\nw2\\
\id\step[2]\id\step[2]\step[2]\xx\step[2]\step[2]\nw2\\
\id\step[2]\id\step[2]\step\ne2\step[2]\id\step[2]\step\coev\step\id\\
\id\step[2]\xx\step[2]\step[2]\id\step[2]\td
g\step\id\step\id\\
\id\step[2]\id\step[2]\nw3\step[2]\step\tu
u\step\ne2\step\id\step\id\\
\id\step[2]\d\step[2]\step[2]\ev\ne2\step[2]\step\id\step\id\\
\d\step[2]\d\step[2]\step\ne2\step[2]\step[2]\step\id\step\id\\
\step\d\step[2]\tu
v\step[2]\step[2]\step[2]\step\id\step\id\\
\step[2]\Ev\step[2]\step[2]\step[6]\id\step\id\\
\step[13]\object{Q^*}\step[1]\object{P^*}
\end{tangle}\step=\step
\begin{tangle}
\object{Z^*}\step[1.5]\object{Y^*}\\
\id\step\id\step[6]\coev\\
\id\step\id\step[5]\ne4\coev\d\\
\id\step\id\step\td f\step[2]\td g\step\id\step\id\\
\id\step\id\step\id\step[2]\x\step[2]\id\step\id\step\id\\
\id\step\id\step\tu u\step[2]\tu v\step\id\step\id\\
\d\ev\step[3]\ne3\step[2]\id\step\id\\
\step\Ev\step[7]\id\step\id\\
\step[9]\object{Q^*}\step\hstep\object{P^*}
\end{tangle}\step=\step\xi^* \  .  \ \ \ \Box
\]

\begin {Proposition} \label {12.1.4}
(cf. \cite [Lemma 2.3, Pro. 2.4] {Ma95a}) (i) If $H$ is a Hopf
algebra living in ${\cal C}$, then $S m = m C_{H,H} (S \otimes S)$
and $\Delta S = C_{H,H} (S \otimes S)\Delta .$

(ii) If $(H, \Delta _H, \epsilon _H, m_H, \eta _H )$   is a
bialgebra or a Hopf algebra living in ${\cal C}$  and has a left
duality $H^*$ in ${\cal C}$, then $(H^*, (\Delta _H)^*, (\epsilon
_H)^*, (m_H)^*, (\eta _H)^* )$  is a bialgebra or  a Hopf algebra in
${\cal C}$.

\end {Proposition}
{\bf Proof.} (i)

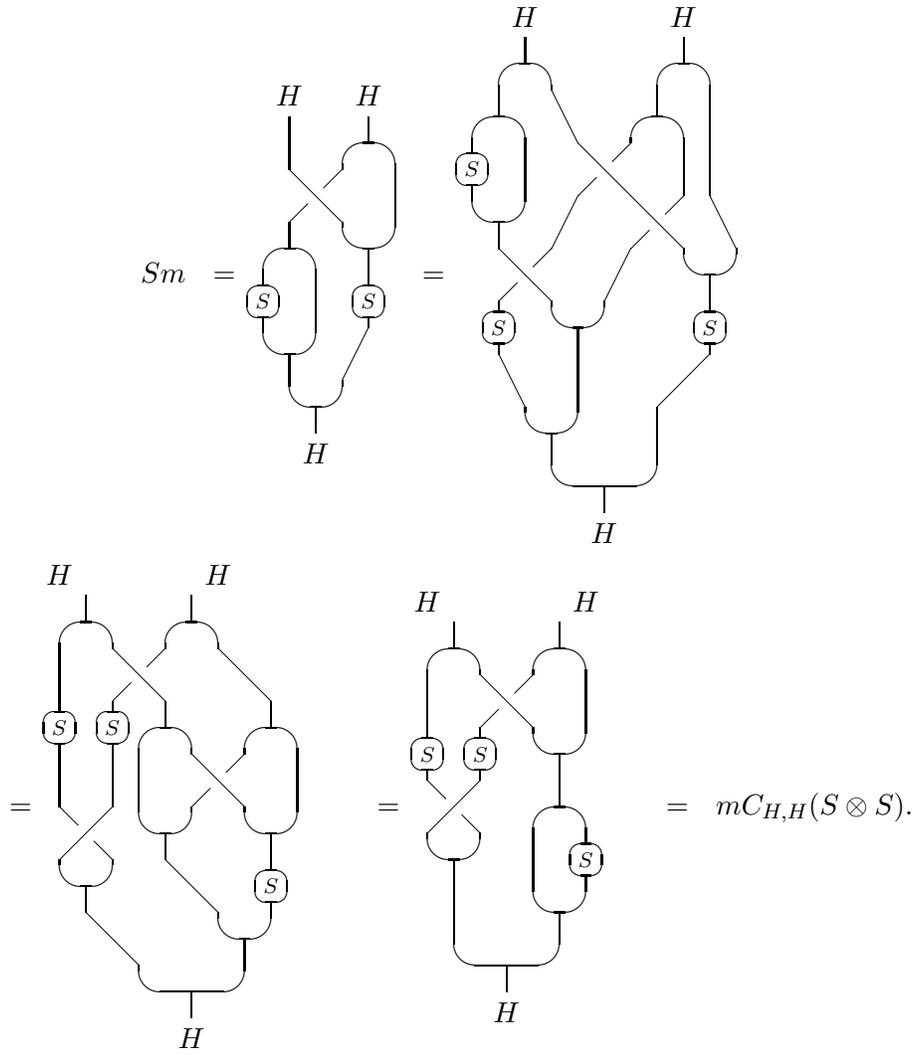
\[Sm \ \ = \ \
\begin{tangle}
\step \object{H}\step[3]\object{H} \\

\step\id\step[2]\cd\\

\step\x \step[2]\id\\

\cd \step [1] \cu  \\

 \O S \step[2]\id\step [2] \O { S }\\

 \cu  \step\ne1 \\

\step[1]\cu \\

\step [2]\object{H}\\

\end{tangle}
\ \ = \ \
\begin{tangle}
\step [2]\object{H}\step[6]\object{H} \\

\step[1]\cd\step[2]\step[2]\cd\step[2]\\

\cd\step[1]\nw1 \step[2]\cd\step[1]\id \\

\O S \step[2]\id \step[2]\x \step [2]\id\step[1]\id \\

\cu\step[1]\ne1 \step[2]\x\step[1] \nw1 \\

\step \x\step[2]\ne1 \step[2]\cu  \\

\step \O S \step[2]\cu  \step[4]\O S\\

\step \nw1  \step [2]\id  \step [4]\ne2\\

\step [2]\cu  \step [3]\id \\
\step[3]\Cu \\

\step [5]\object{H}\\

\end{tangle}
\]

\[ \ \ = \ \
\begin{tangle}
\object{H}\step[6]\object{H} \\

\cd\step[2]\cd\step[2]\\

\id \step[2]\x\step[2]\nw2 \\

\O S \step [2]\O S \step [1] \cd \step [2]\cd \step [2]\\

\id \step [2] \id  \step [1]\id \step [2]\x \step [2] \id \\

\xx \step [1]\cu \step [2] \cu \\

\cu \step [2]\nw2  \step [3] \O S\\

\step \nw2 \step [4] \cu \\

 \step [3] \Cu \step [2]\\

\step [5]\object{H}\\

\end{tangle}
 \ \ = \ \
\begin{tangle}
\object{H}\step[6]\object{H} \\

\cd\step[2]\cd\step[2]\\

\id \step[2]\x\step[2]\id \\

\O S \step [2]\O S \step [2] \cu \\

\xx \step [2]\cd \step [2] \\

\cu \step [2]\id  \step [2] \O S\\

\step \id  \step [3] \cu \\

 \step [1] \Cu \step [2]\\

\step [3]\object{H}\\

\end{tangle} \ \ = \ \ m C_{H, H} (S \otimes S).
\] Similarly, we can show  $ \Delta S = C_{H, H} (S \otimes S) \Delta.$

(ii) It is clear that $H^*$  is not only an algebra but also a
coalgebra if we set  $\Delta_{H^*} = (m_H)^*, \epsilon _{H*} = (\eta
_H)^*$ and $m_{H^*} = (\Delta _H)^* .$

 Now we show that $ \step \begin{tangle}
\object{H^*}\step[2]\object{H^*}\\
\tu {m_{H^*}}\\
\td {\Delta_{H^*}}\\
 \object{H^*}\step[2]\object{H^*}
\end{tangle}=
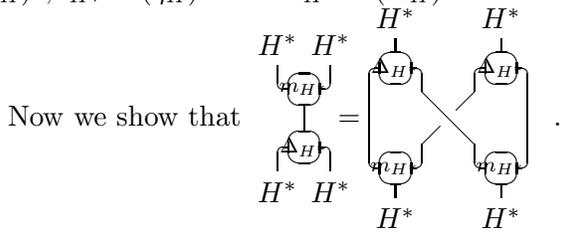
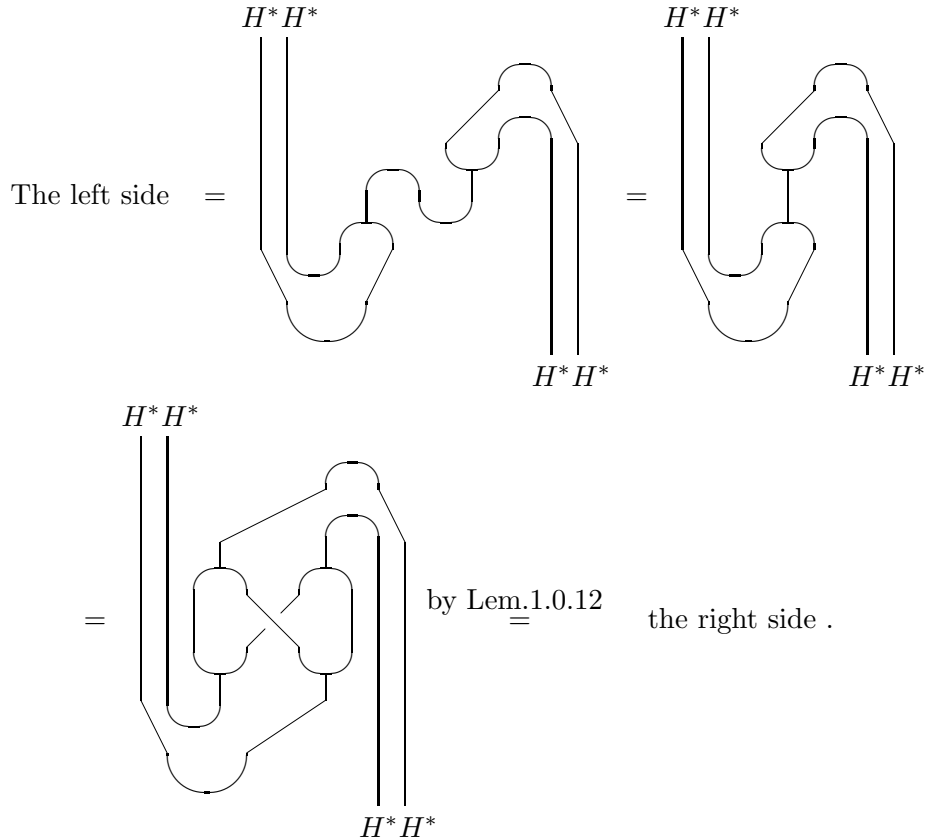
\begin{tangle}
\step\object{H^*}\step[4]\object{H^*}\\
\td {\Delta_{H^*}}\step[2]\td {\Delta_{H^*}}\\
\id\step[2]\x\step[2]\id\\
\tu {m_{H^*}}\step[2]\tu {m_{H^*}}\\
\step \object{H^*}\step[4]\object{H^*}
\end{tangle}\step .$

\[
\hbox {The left side} \step=\step
\begin{tangle}
\object{H^*}\step[1.5]\object{H^*}\\
\id\step\id\step[8]\coev\\
\id\step\id\step[7]\ne2\coev\d\\
\id\step\id\step[3]\coev\step\cu\step[2]\id\step\id\\
\id\step\id\step[2]\cd\step\ev\step[3]\id\step\id\\
\d\ev\step\dd\step[6]\id\step\id\\
\step\Ev\step[9]\id\step\id\\
\step[11]\object{H^*}\step[1.5]\object{H^*}
\end{tangle}\step=\step
\begin{tangle}
\object{H^*}\step[1.5]\object{H^*}\\
\id\step\id\step[4]\coev\\
 \id\step\id\step[3]\ne2\coev\d\\
  \id\step\id\step[2]\cu\step[2]\id\step\id\\
\id\step\id\step[2]\cd\step[2]\id\step\id\\
\d\ev\step\dd\step[2]\id\step\id\\
\step\Ev\step[5]\id\step\id\\
\step[7]\object{H^*}\step\hstep\object{H^*}
\end{tangle}
\]
\[=\step
\begin{tangle}
\object{H^*}\step[1.5]\object{H^*}\\
\id\step\id\step[6]\coev\\
\id\step\id\step[5]\ne4\coev\d\\
\id\step\id\step\cd\step[2]\cd\step\id\step\id\\
\id\step\id\step\id\step[2]\x\step[2]\id\step\id\step\id\\
\id\step\id\step\cu\step[2]\cu\step\id\step\id\\
\d\ev\step[3]\ne3\step[2]\id\step\id\\
\step\Ev\step[7]\id\step\id\\
\step[9]\object{H^*}\step\hstep\object{H^*}
\end{tangle}\stackrel{ \hbox {by Lem.\ref {2.1.9'} }}{=} \step\hbox {the right side .}
\]

 Thus $H^*$ is a bialgebra living in ${\cal C}.$
 If $H$  is  a Hopf algebra, then $H^*$  is also a Hopf algebra with
 antipode $(S_H)^*$. \begin{picture}(8,8)\put(0,0){\line(0,1){8}}\put(8,8){\line(0,-1){8}}\put(0,0){\line(1,0){8}}\put(8,8){\line(-1,0){8}}\end{picture}

If $H$  is an algebra   and the following conditions are satisfied:
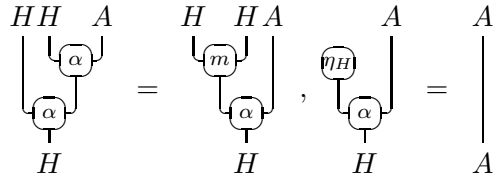
\[
\begin{tangle}
\object{H}\step \object{H} \step[2]\object{A}  \\
\id\step\tu \alpha\\
\tu \alpha\\
\step\object{H}\\
\end{tangle}
\step=\step
\begin{tangle}
\object{H}\step[2] \object{H} \step\object{A} \\
\tu {m} \step\id\\
\step\tu \alpha\\
\step[2]\object{H}\\
\end{tangle}
\step,\step
\begin{tangle}
\step[2]\object{A}\\
\Q {\eta_H}\step[2]\id \\
\tu \alpha\\
\step\object{H}\\
\end{tangle}
\step=\step
\begin{tangle}
\object{A}\\
\id \\
\id \\
\object{A}
\end{tangle}
\]

 \noindent then $(A , \alpha )$
is called a braided $H$-module or an $H$-module living in ${\cal
C}.$

If $H$  is a coalgebra   and the following conditions are satisfied:
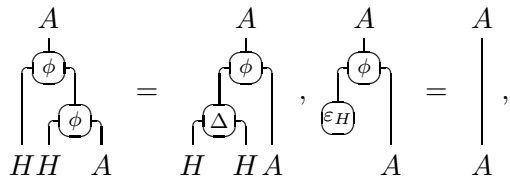
\[
\begin{tangle}
\step\object{A}\\
\td \phi\\
\id\step\td \phi\\
\object{H}\step\object{H}\step[2]\object{A}\\
\end{tangle}
\step=\step
\begin{tangle}
\step[2]\object{A}\\
\step\td \phi\\
\td \Delta\step\id\\
\object{H}\step[2]\object{H}\step\object{A}\\
\end{tangle}
\step,\step
\begin{tangle}
\step\object{A}\\
\td \phi\\
\QQ {\varepsilon_H}\step[2]\id\\
\step[2]\object{A}\\
\end{tangle}
\step=\step
\begin{tangle}
\object{A}\\
\id\\
\id\\
\object{A}\\
\end{tangle}\ \ ,
\] \noindent then $(A , \phi )$  is
called a braided $H$-comodule or an $H$-comodule living in ${\cal
C}$
 Let $_H {\mathcal M(C)}$ and $^H {\mathcal M(C)}$  denote the
 categories
 of $H$-modules in ${\mathcal C}$ and $H$-comodules in ${\mathcal
 C}$, called a braided $H$-module category and a braided $H$- comodule category, respectively.
  Let $_H {\mathcal M(C)}= {} _H {\mathcal M}$ and $^H {\mathcal M(C)}= {}^H {\mathcal M}$
  when ${\cal C} = {\cal V} ect (k)$ in Example \ref {0.4}.

Similarly, we can define the concepts of algebra morphism, coalgebra
morphism and $H$-module morphism.  An algebra  $(A, \alpha )$  is
called an $H$-module algebra, if $(A, \alpha )$ is an $H$-module and
 its unit $\eta _A $  and multiplication
$m_A$ are $H$-module morphisms. Similarly, we can define the
concepts of $H$-module coalgebra, $H$-comodule algebra and
$H$-comdule coalgebra.

Assume  that ${\cal C}$ is a concrete braided tensor category. If
$H$ is a braided bialgebra and $M$ is an $H$-module in ${\cal C}$,
then
$$M^H := \{ x \in M  \mid    h x = \epsilon(h)x,\forall  h\in H \}$$
is called the invariants of $M$. If $M$ is a regular $H$-module
(i.e. the module operation is $m$ ), then $M^H$  is called the
integral of $H$ and is written as $\int _H^l.$

A bialgebra $(H,m, \eta, \Delta,\epsilon  )$ with
convolution-invertible $R$ in $Hom _{\cal C} (I, H\otimes H)$ is
called a quasitriangular bialgebra
 living in
braided tensor category  ${\cal C}$ if there is a morphism $\bar
\Delta : H \rightarrow H \otimes H$ such that $(H, \bar \Delta ,
\epsilon )$ is a coalgebra and

(QT1):
\[
\begin{tangle}
\ro R\\
\id\step\td \Delta\\
\id\step\id\step[2]\id\\
\object{H}\step\object{H}\step[2]\object{H}
\end{tangle}
\step=\step
\begin{tangle}
\Ro R\\
\id\step\ro R\step\id\\
\hcu\step[2]\id\step\id\\
\step[0.5]\object{H}\step\object{H}\step[2]\object{H}
\end{tangle}
\]

(QT2):
\[
\begin{tangle}
\step\ro R\\
\td {\bar{\Delta}}\step\id\\
\object{H}\step[2]\object{H}\step\object{A}
\end{tangle}
\step=\step
\begin{tangle}
\Ro R\\
\id\step\ro R\step\id\\
\id\step\id\step[2]\hcu\\
\step[0.5]\object{H}\step\object{H}\step[2]\object{H}
\end{tangle}
\]

(QT3):
\[
\begin{tangle}
\step[4]\object{H}\\
\ro R\step\cd\\
\id\step[2]\hx\step[2]\id\\
\cu\step\cu\\
\step\object{H}\step[3]\object{H}
\end{tangle}
\step=\step
\begin{tangle}
\step\object{H}\\
\td {\bar{\Delta}}\step\ro R\\
\id\step[2]\hx\step[2]\id\\
\cu\step\cu\\
\step\object{H}\step[3]\object{H}
\end{tangle} \ \ .
\]
 In this case, we also say that $(H,R, \bar \Delta )$ is a braided quasitriangular
bialgebra in ${\mathcal C}$.

A bialgebra $(H, m, \eta, \Delta , \epsilon )$ with
convolution-invertible $r$ in $Hom _{\cal C} ( H\otimes H, I)$ is
called a coquasitriangular bialgebra living in braided tensor
category ${\cal C}$ if there exists a morphism $ \bar m : H \otimes
H \rightarrow H $ such that $(H, \bar m, \eta )$ is an algebra and

 (CQT1):
\[
\begin{tangle}
\object{H}\step\object{H}\step[2]\object{H}\\
\id\step\id\step[2]\id\\
\id\step\cu\\
\coro r\\
\end{tangle}
\step=\step
\begin{tangle}
\step[0.5]\object{H}\step[2.5]\object{H}\step\object{H}\\
\hcd\step[2]\id\step\id\\
\id\step\coro r \step\id\\
\coRo r
\end{tangle}
\]

(CQT2):
\[
\begin{tangle}
\object{H}\step[2]\object{H}\step[1]\object{H}\\

\tu {\bar m}\step [1]\id\step\\

\step \coro r\\
\end{tangle}
\step=\step
\begin{tangle}

\object{H}\step\object{H}\step [2]\object{H}\\
\id\step\id\step [2]\hcd\\
\id\step\coro r \step\id\\
\coRo r
\end{tangle}
\]

(CQT3):
\[
\begin{tangle}
\step\object{H}\step[3]\object{H}\\
\cd\step\cd\\
\id\step[2]\hx\step[2]\id\\
\coro r\step\cu\\
\step[4]\object{H}

\end{tangle}
\step=\step
\begin{tangle}
\step\object{H}\step[3]\object{H}\\
\cd\step\cd\\
\id\step[2]\hx \step[2]\id\\
\tu {\bar m}\step\coro r\\
\step\object{H}
\end{tangle}
\]
\noindent In this case, we also say that $(H,r, \bar m)$ is a
braided coquasitriangular bialgebra in ${\mathcal C}$.

A braided Hopf algebra $H$  is called a braided group if there
exists a braided quasitriangular structure $(H, R, \bar \Delta )$
with $\bar \Delta = \Delta _H$ or there exists a coquasitriangular
structure $(H, r, \bar m )$  with $m = \bar m.$

Let  $ m^{op}$ and $\Delta ^{cop}$ denote  $m _HC_{H,H}^{-1}  $ and
$C_{H,H}^{-1}\Delta $ respectively.

In particular, we  say that  $(H,R)$  is quasitriangular if $(H, R,
\Delta ^{cop})$ is quasitriangular. Dually, we  say that  $(H,r)$ is
quasitriangular if $(H, r, m ^{op})$ is coquasitriangular.

Let $(H, R,\bar \Delta )$ be a quasitriangular bialgebra
 living in braided tensor
category ${\cal C}$.
 $(V, \alpha )$  is said to be symmetric  with respect to $(H, R, \bar
  \Delta  )$,
if $(V, \alpha )$ is a left $H$-module and  satisfies the following
condition: \[(SWRT):  \ \ \ \ \ \ \ \
\begin{tangle}
\step\object{H}\step[3]\object{V} \\
\td {\overline{\Delta}}\step\dd\\
\id\step[2]\hx\\
\tu\alpha\step\id\\
\step\x\\
 \step\object{H}\step[2]\object{V}
 \end{tangle}=
\begin{tangle}
\step\object{H}\step[3]\object{V} \\
\cd\step[2]\id\\
\id\step[2]\tu \alpha\\
\object{H}\step[3]\object{V}
 \end{tangle}\ \ \ .
 \]

Define
$$O(H, \bar \Delta ) = \{ V \in {}_H {\cal M} \mid (V, \alpha ) \
\hbox {  is symmetric with respect to } \ (H, R, \bar \Delta  )
\}.$$

Let $(H, r, \bar m)$ be a coquasitriangular bialgebra
 living in braided tensor category.
$(V, \phi )$  is said to be symmetric  with respect to $(H, r, \bar
m )$ if $(V, \phi )$ is a left $H$-comodule and  satisfies the
following condition:

\[(SWRT): \ \ \ \ \ \ \ \
\begin{tangle}
 \step\object{H}\step[2]\object{V}\\
 \step\x\\
 \td\phi\step\id\\
 \id\step[2]\hx\\
 \tu {\overline{m}}\step\d\\
\step\object{H}\step[3]\object{V} \\
 \end{tangle}=
\begin{tangle}
\object{H}\step[3]\object{V}\\
\id\step[2]\td \phi\\
\cu\step[2]\id\\
\step\object{H}\step[3]\object{V} \\
 \end{tangle}
 \step[2] .
\]
 Define
$$O(H, \bar m ) = \{ V \in {}^H {\cal M} \mid (V, \phi )
\hbox { is symmetric with respect to } (H, r, \bar m ) \}.$$

\begin {Proposition} \label {0.5}
(i) If   $(H,R, \bar \Delta )$ is a quasitriangular bialgebra living
in braided tensor category $({\cal C},  \otimes, I,  a, l, r, C  )$,
then  $(O(H,\bar \Delta ), \otimes, I,  a, l, r,  C^R )$ is a
braided tensor category, where

\[
C^{R}_{V,W}=\step
\begin{tangle}
\step[2]\step[2]\object{V}\step[2]\object{W}\\
\ro R\step[2]\id\step[2]\id\\
\id\step[2]\x\step[2]\id\\
\tu {\alpha_{V}}\step[2]\tu {\alpha_{W}}\\
\step\d\step[2]\dd\\
\step[2]\x\\
 \step[2]\object{W}\step[2]\object{V}
 \end{tangle}
 \step[2]\hbox {and} \step[2]
 \alpha_{V\otimes W}=
\begin{tangle}
\step\object{H}\step[3]\object{V}\step[2]\object{W}\\
\cd\step[2]\id\step[2]\id\\
\id\step[2]\x\step[2]\id\\
\tu {\alpha_{V}}\step[2]\tu {\alpha_{W}}\\
 \step\object{V}\step[4]\object{W}
 \end{tangle}
\] for any $(V,\alpha
_V), (W,\alpha _W) \in O(H, \bar \Delta )$

(ii) If $(H, r, \bar m )$ is a coquasitriangular bialgebra living in
braided tensor category $({\cal C},  \otimes, I,  a, l, r  ),$ then
$(O(H, \bar m), \otimes, I,  a, l, r,  C^r )$ is a braided tensor
category, where

\[
C^{r}_{V,W}=\step
\begin{tangle}
 \step[2]\object{V}\step[2]\object{W}\\
 \step[2]\x\\
 \step\dd\step[2]\d\\
 \td {\phi_{W}}\step[2]\td {\phi_{V}}\\
 \id\step[2]\x\step[2]\id\\
 \coro r\step[2]\id\step[2]\id\\
\step[2]\step[2]\object{W}\step[2]\object{V}\\
\end{tangle}
  \step[2]\hbox {and} \step[2]
 \phi_{V\otimes W}=
\begin{tangle}
 \step\object{V}\step[4]\object{W}\\
 \td {\phi_{V}}\step[2]\td {\phi_{W}}\\
 \id\step[2]\x\step[2]\id\\
 \cu\step[2]\id\step[2]\id\\
\step\object{H}\step[3]\object{V}\step[2]\object{W}\\
 \end{tangle}
\]

 for any $(V,
\phi _V), (W, \phi _W) \in  O(H, \bar m).$

          \end {Proposition}

{\bf Proof.}  We only show part (i) since part (ii) can be shown
similarly. We first show that $O(H, \bar \Delta )$ is a tensor
category. That is, we need show that $V\otimes W \in O(H, \bar
\Delta )$  for any $V, W  \in O(H, \bar \Delta ). $
In fact,
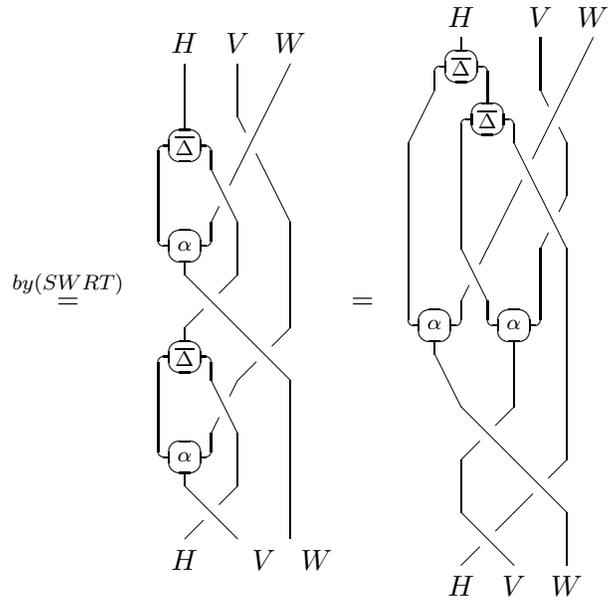
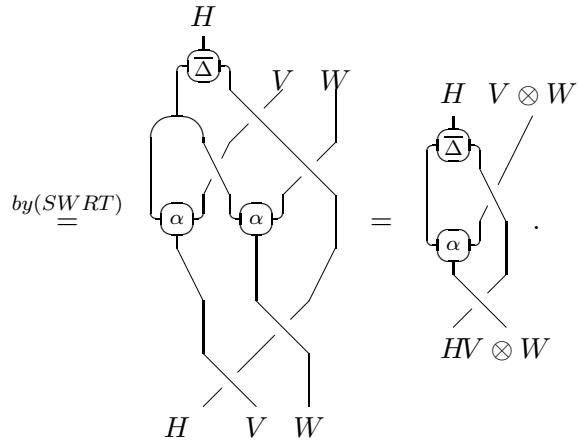
\[
\begin{tangle}
\step\object{H}\step[3]\object{V\otimes W} \\
\cd\step[2]\id\\
\id\step[2]\tu \alpha\\
\object{H}\step[3]\object{V\otimes W}
 \end{tangle}
 \step=\step
\begin{tangle}
\step[2]\object{H} \\
\step\cd\step[4]\object{V}\step[3]\object{W} \\
\cd\step\nw2 \step[2]\dd\step[2]\dd\\
\id\step[2]\d\step[2]\hx\step[2]\dd\\
\id\step[2]\step \tu \alpha\step\tu\alpha\\
\object{H}\step[4]\object{V}\step[3]\object{W} \\
 \end{tangle}
  \step=\step
\begin{tangle}
\step[2]\object{H} \\
\step\cd\step[4]\object{V}\step[2]\object{W} \\
\cd\step\nw2 \step[3]\hxx\\
\id\step[2]\d\step[2]\tu\alpha\step\id\\
\id\step[2]\step \d\step[2]\x\\
\id\step[4]\tu\alpha\step[2]\id\\
\object{H}\step[4]\object{V}\step[3]\object{W} \\
 \end{tangle}
 \]

\[
\stackrel{by (SWRT)}{=} \step
 \begin{tangle}
\step\object{H}\step[2]\object{V}\step[2]\object{W} \\
\step\id \step[2]\id\step\dd\\
\td {\overline{\Delta}}\step\hxx\\
\id\step[2]\hx\step\d\\
\tu\alpha\step\id\step[2]\id\\
\step\x\step[2]\id\\
\td {\overline{\Delta}}\step\x\\
\id\step[2]\hx\step[2]\id\\
\tu\alpha\step\id\step[2]\id\\
\step\x\step[2]\id\\
\step\object{H}\step[3]\object{V}\step[2]\object{W} \\
 \end{tangle}\step=\step
 \begin{tangle}
\step[2]\object{H}\step[3]\object{V}\step[2]\object{W} \\
\step\td {\overline{\Delta}}\step[2]\id\step\dd\\
\dd\step\td {\overline{\Delta}}\step\hxx\\
\id\step[2]\id\step[2]\hx\step\id\\
\id\step[2]\id\step\dd\step\hx\\
\id\step[2]\hx\step[2]\id\step\id\\
\tu \alpha\step\tu \alpha\step\id\\
\step\d\step[2]\id\step[2]\id\\
\step[2]\x\step[2]\id\\
\step[2]\id\step[2]\x\\
\step[2]\x\step[2]\id\\
\step[2]\object{H}\step[2]\object{V}\step[2]\object{W} \\
 \end{tangle}
\]
\[\ \ \  \stackrel{by (SWRT)}{=} \ \
 \begin{tangle}
\step[2]\object{H}\\
\step\td {\overline{\Delta}}\step[2]\object{V}\step[2]\object{W} \\
\cd\step\x\step[2]\id\\
\id\step[2]\hx\step[2]\x\\
\tu \alpha\step\tu\alpha\step[2]\id\\
\step\d\step[2]\id\step[2]\dd\\
\step[2]\id\step[2]\x\\
\step[2]\x\step[2]\id\\
 \step\object{H}\step[3]\object{V}\step[2]\object{W}
 \end{tangle}
 \step=\step
 \begin{tangle}
\step\object{H}\step[3]\object{V\otimes W} \\
\td {\overline{\Delta}}\step\dd\\
\id\step[2]\hx\\
\tu\alpha\step\id\\
\step\x\\
 \step\object{H}\step[2]\object{V\otimes W}
 \end{tangle}.
\]Therefore $O(H,\bar \Delta )$  is a tensor  subcategory of ${\cal
C}.$ Next we show that $C^{R}$ is a braiding in $O(H, \bar \Delta
).$

It is clear that the inverse of $C^R_{V,W}$ is
 $\step \begin{tangle}
\step[4]\object{W}\step[2]\object{V} \\
\ro { \bar R}\step[2]\xx\\
\id\step[2]\x\step[2]\id\\
\tu\alpha\step[2]\tu\alpha\\
\step\object{V}\step[4]\object{W} \\
 \end{tangle}\step .$

We see that
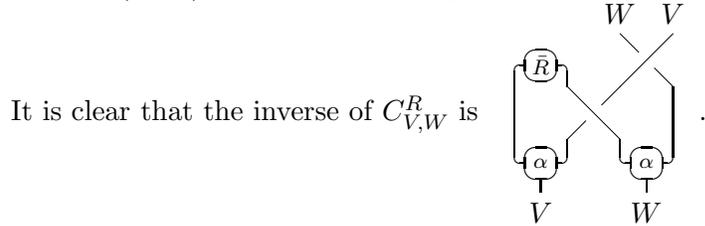
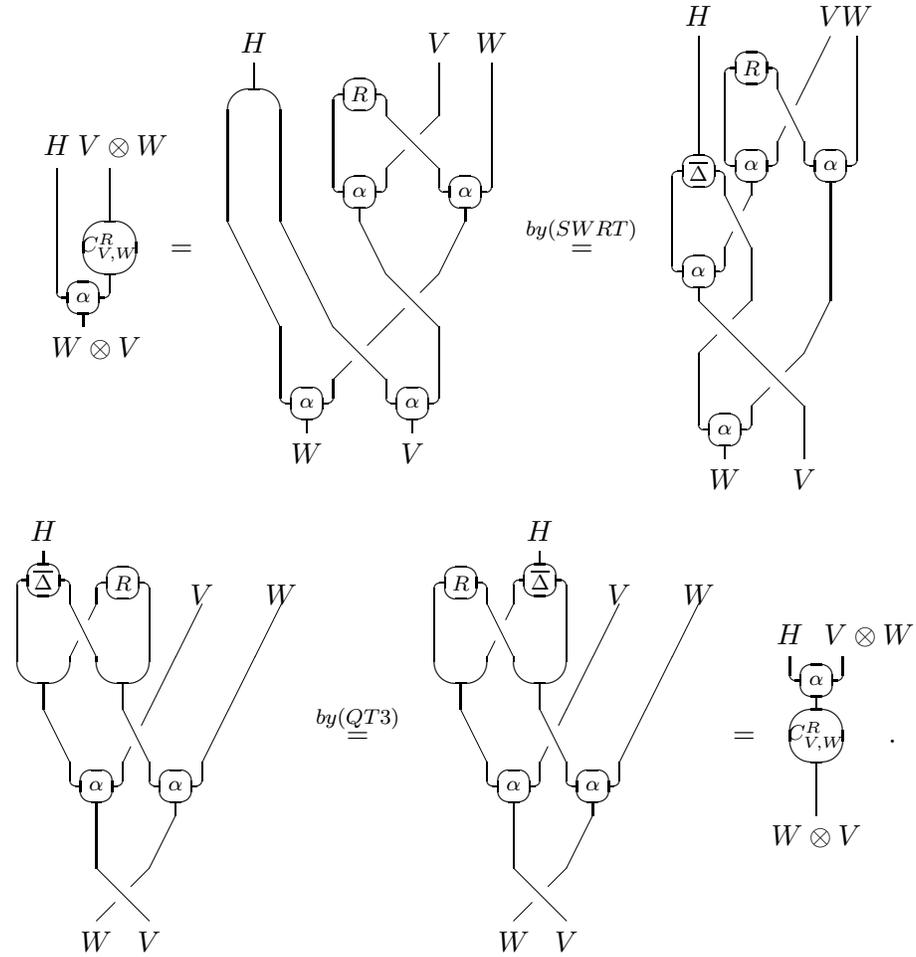
\[
\begin{tangle}
\object{H}\step[2.5]\object{V\otimes W} \\
\id\step[2]\id\\
\id\step\obox 2{C^{R}_{V,W}} \\
\tu \alpha\\
\step[1.5]\object{W\otimes V} \\
 \end{tangle}\step=\step
  \begin{tangle}
 \step\object{H}\step[7]\object{V}\step[2]\object{W} \\
\cd\step[2]\ro R \step[2]\id\step[2]\id\\
\id\step[2]\id\step[2]\id\step[2]\x\step[2]\id\\
\id\step[2]\id\step[2]\tu\alpha\step[2]\tu \alpha\\
\d\step\d\step[2]\d\step[2]\dd\\
\step\d\step\d\step[2]\x\\
\step[2]\id\step[2]\x\step[2]\id\\
\step[2]\tu\alpha\step[2]\tu\alpha\\
\step[3]\object{W}\step[4]\object{V} \\
 \end{tangle}\step\stackrel{by (SWRT)}{=}\step
  \begin{tangle}
 \step\object{H}\step[5]\object{V}\step\object{W} \\
 \step\id\step\ro R\step\dd\step\id\\
 \step\id\step\id\step[2]\hx\step[2]\id\\
 \td {\overline{\Delta}}\tu \alpha\step\tu \alpha\\
 \id\step[2]\hx\step[2]\step\id\\
 \tu\alpha\step\id\step[2]\step\id\\
 \step\x\step[2]\dd\\
 \step\id\step[2]\x\\
 \step\tu \alpha\step[2]\id\\

\step[2]\object{W}\step[3]\object{V}
  \end{tangle}
\]
\[
  \begin{tangle}
 \step\object{H}\\
 \td {\overline{\Delta}}\step\ro R\step[2]\object{V}\step[3]\object{W} \\
 \id\step[2]\hx\step[2]\id\step\dd\step[2]\dd\\
 \cu\step\cu\dd\step[2]\dd\\
 \step\d\step[2]\hx\step[2]\dd\\
 \step[2]\tu\alpha\step\tu\alpha\\
\step[3]\id \step[2]\ne1\\
\step[3]\x\\
 \step[3]\object{W}\step[2]\object{V}
  \end{tangle}\step\stackrel{by (QT3)}{=}\step
    \begin{tangle}
 \step[4]\object{H}\\
\ro R \step \td {\overline{\Delta}}\step[2]\object{V}\step[3]\object{W} \\
 \id\step[2]\hx\step[2]\id\step\dd\step[2]\dd\\
 \cu\step\cu\dd\step[2]\dd\\
 \step\d\step[2]\hx\step[2]\dd\\
 \step[2]\tu\alpha\step\tu\alpha\\
\step[3]\id \step[2]\ne1\\
\step[3]\x\\
 \step[3]\object{W}\step[2]\object{V}
  \end{tangle}\step=\step
  \begin{tangle}
\object{H}\step[3]\object{V\otimes W} \\
\tu \alpha\\
\obox 2{C^{R}_{V,W}} \\
\step\id\\
\step[1]\object{W\otimes V} \\
 \end{tangle} \ \ .
\] Thus $C^R_{V, W}$ is an $H$-module morphism.

We also see that

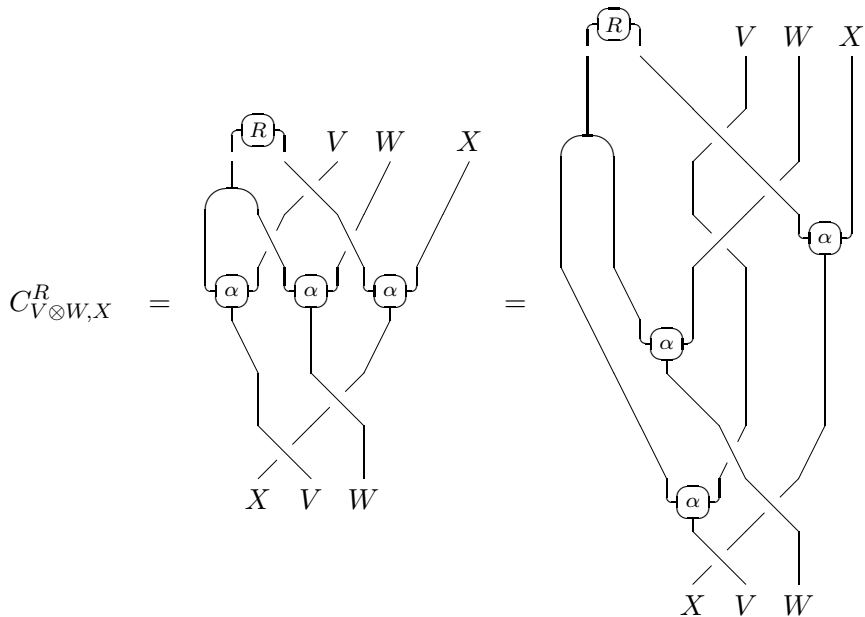
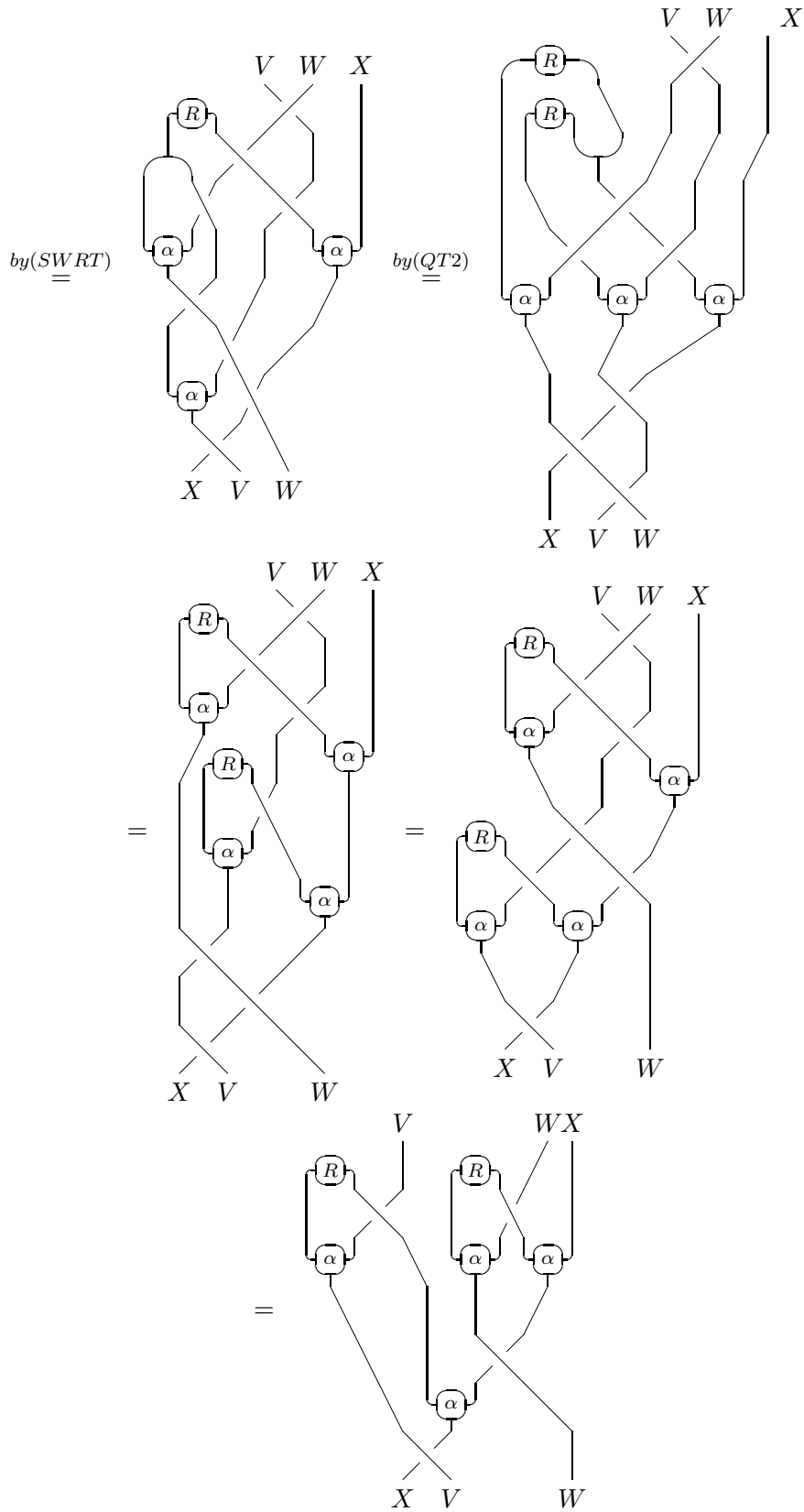
\[
C^{R}_{V\otimes W ,X}\step=\step
   \begin{tangle}
\step\ro R  \step[2]\object{V}\step[2]\object{W}\step[3]\object{X} \\
\cd\step\x\step\dd\step[2]\dd\\
\id \step[2]\hx \step[2]\hx\step[2]\dd\\
\tu\alpha\step\tu\alpha\step\tu\alpha\\
\step\d\step[2]\id\step[2]\dd\\
\step[2]\id\step[2]\x\\
\step[2]\x\step[2]\id\\
 \step[2]\object{X} \step[2]\object{V}\step[2]\object{W}
  \end{tangle}\step=\step
   \begin{tangle}
\step\ro R  \step[4]\object{V}\step[2]\object{W}\step[2]\object{X} \\
\step\id\step[2]\nw2\step[2]\step\id\step[2]\id\step[2]\id\\

\cd\step[3]\x\step[2]\id\step[2]\id\\

\id\step[2]\id\step[3]\id\step[2]\x\step[2]\id\\

\id\step[2]\id\step[3]\xx\step[2]\tu \alpha \step[2]\\

\d \step\d \step[2]\id\step[2]\id \step[3]\id\\

\step\d\step\tu\alpha\step[2]\id\step[3] \id \\
\step[2]\d\step\nw2\step[2]\id\step[3]\id\\
\step[3]\d\step[2]\hx\step[2]\dd\\
\step[4]\tu\alpha\step\x\\
\step[5]\x\step[2]\id\\
 \step[5]\object{X} \step[2]\object{V}\step[2]\object{W}
  \end{tangle}
\]

\[
\stackrel{by (SWRT)}{=}\step
   \begin{tangle}
\step[5]\object{V}\step[2]\object{W}\step[2]\object{X} \\
\step\ro R\step[2]\xx\step[2]\id\\
\cd\step\x\step[2]\id\step[2]\id\\
\id\step[2]\hx\step[2]\x\step[2]\id\\
\tu\alpha\step\id\step[2]\id\step[2]\tu\alpha\\
\step\x\step\dd\step[2]\dd\\
\step\id\step[2]\hx\step[2]\ne2\\
\step\tu\alpha\step\hx\\
\step[2]\x\step\d\\
 \step[2]\object{X} \step[2]\object{V}\step[2]\object{W}
  \end{tangle}
\step\stackrel{by (QT2)}{=}\step
\begin{tangle}
 \step[7]\object{V}\step[2]\object{W}\step[3]\object{X} \\
\Ro R\step[3]\xx\step[2]\id\\

\id \step\ro R\step[1]\nw1 \step[2]\id\step[2]\id\step[2]\id\\

\id \step\id \step[2]\cu \step[1]\ne1\step[1]\ne1\step [1] \ne1\\

\id \step\nw1 \step[2]\xx \step[2]\id\step[2]\id\\

\id  \step[2]\xx \step[2]\xx\step[2]\id\\

\tu\alpha \step[2]\tu\alpha \step[2]\tu\alpha \step[2]\\

\step\d\step[2]\ne1\step[3]\ne3\\
\step[2]\id\step[2]\x\\
\step[2]\x\step[2]\id\\

\step[2]\id \step[2]\x\\

 \step[2]\object{X} \step[2]\object{V}\step[2]\object{W}
  \end{tangle}
  \]

\[
\step=\step
   \begin{tangle}
 \step[4]\object{V}\step[2]\object{W}\step[2]\object{X} \\
\ro R \step[2]\xx\step[2]\id\\
\id \step[2]\x \step[2]\id\step[2]\id\\
\tu\alpha\step[2]\x\step[2]\id\\
\dd\ro R \step\id\step[2]\tu\alpha\\
\id\step\id\step[2]\hx\step[2]\step\id\\
\id\step[1]\tu\alpha\step\d\step[2]\id\\
\id\step[2]\id\step[2]\step\tu\alpha\\
\x\step[3]\ne2\\
\id\step[2]\x\\
\x\step[2]\nw2\\
\object{X} \step[2]\object{V}\step[4]\object{W}
  \end{tangle}
  \step=\step
   \begin{tangle}
 \step[6]\object{V}\step[2]\object{W}\step[2]\object{X} \\
\step[2]\ro R \step[2]\xx\step[2]\id\\
\step[2]\id \step[2]\x \step[2]\id\step[2]\id\\
\step[2]\tu\alpha\step[2]\x\step[2]\id\\
\step[2]\step\d\step[2]\id\step[2]\tu\alpha\\
\ro R\step[2]\x\step[2]\dd\\
\id\step[2]\x\step[2]\x\\
\tu\alpha\step[2]\tu\alpha\step[2]\id\\
\step\d\step[2]\dd\step[2]\step\id\\
\step[2]\x\step[2]\step[2]\id\\
\step[2]\object{X} \step[2]\object{V}\step[4]\object{W}
  \end{tangle}
  \]
\[
  \step=\step
   \begin{tangle}
 \step[4]\object{V}\step[6]\object{W}\step[1]\object{X} \\
\ro R \step[2]\id\step[2]\ro R\step\dd\step\id\\
\id \step[2]\x \step[2]\id\step[2]\hx\step[2]\id\\
\tu\alpha\step[2]\d\step\tu\alpha\step\tu\alpha\\
\step\d\step[2]\step\id\step[2]\id\step[2]\dd\\
\step[2]\d\step[2]\id\step[2]\x\\
\step[2]\step\d\step\tu\alpha\step[2]\nw2\\
\step[2]\step[2]\x\step[2]\step[3]\id\\
 \step[4]\object{X}
\step[2]\object{V}\step[5]\object{W}
  \end{tangle}
  \]
$= (C_{V,X}^R \otimes id _W)(id _V \otimes C_{W,X}^R)$. Thus
$$C_{V\otimes W, X}^R = (C_{V,X}^R \otimes id _W)(id _V \otimes C_{W,X}^R).$$
Similarly, we have
$$C_{V, W\otimes X}= (id _W\otimes C_{V,X}^R)(C_{V,W}^R \otimes id _X).$$

 Therefore $O(H,\bar
\Delta )$  is a braided tensor category.
\begin{picture}(8,8)\put(0,0){\line(0,1){8}}\put(8,8){\line(0,-1){8}}\put(0,0){\line(1,0){8}}\put(8,8){\line(-1,0){8}}\end{picture}

It is clear that $O(H, \Delta _H^{cop} ) =\ {} _H \!{\cal M (C)}$
and $O(H, m_H ^{op} ) = \ {}^H \!{\cal M (C)}$ when $({\cal C}, C)$
is a symmetric braided tensor category.

\begin {Corollary} \label {5.1.3''} Let $({\cal C}, \otimes, I,  a, l, r ,          C)$ be a
symmetric braided tensor category.

 (i)  If   $(H,R)$ is a quasitriangular bialgebra  in ${\cal C}$, then  $({}_H {\cal M}, \otimes, I,  a, l, r,
C^R )$ is a braided tensor category.

(ii) If  $(H,r)$ is a coquasitriangular bialgebra in ${\cal C}$,
then  $({}^H {\cal M}, \otimes, I,  a, l, r,  C^r )$ is a braided
tensor category.
\end {Corollary}

Let $H$ be a bialgebra algebra in braided tensor category
$({\mathcal C}, C)$, and $(M, \alpha, \phi )$ a left  $H$-module and
left $H$ -comodule in $({\mathcal C}, C)$. If

\[
(YD):\step[2]
\begin{tangle}
\step[1]\object{H}\step[3]\object{A}\\
\cd\step[2]\id\\
\id\step[2]\x\\
\tu \alpha\step[2]\id\\
\td \phi\step[2]\id\\
\id\step[2]\x\\
\cu \step[2]\id\\
\step[1]\object{H}\step[3]\object{A}
\end{tangle}
\step=\step
\begin{tangle}
\step[1]\object{H}\step[4]\object{A}\\
\cd\step[2]\td \phi\\
\id\step[2]\x\step[2]\id\\
\cu \step[2]\tu \alpha\\
\step[1]\object{H}\step[4]\object{A}
\end{tangle}\ \ \ .
\]
then $(M, \alpha, \phi )$ is called  a Yetter-Drinfeld
  $H$-module in ${\mathcal C}$, written as YD $H$-module in short.
When $H$ is a Hopf algebra algebra in braided tensor category
$({\mathcal C}, C)$, then the condition above is equivalent to

(YD): \[ \phi \alpha =  \begin{tangle}
\step [3]\object{H}\step[7]\object{M} \\

\step \Cd \step [4] \td \phi \\

\cd \step [3] \O S \step [3] \ne2 \step [2] \id \\

\id \step [2] \nw1  \step [2]  \x \step [4] \id  \\

\id \step [2]  \step \x  \step[2] \id \step [4] \id\\

\id \step [2] \ne1  \step [2] \x \step[3]   \ne2\\

\cu  \step[2] \ne2  \step[2] \tu \alpha  \\

\step \cu  \step[4]  \step \id  \\
\step [2]\object{H}\step[6]\object{M}
  \end{tangle} \ \ \ .
  \]
  Let $^H_H {\mathcal YD (C)}$ denote the category of all  Yetter-Drinfeld
  $H$-modules in ${\mathcal C}$. If $({\mathcal C}, C) = {\mathcal V}ect(k)$, we write
  $^H_H {\mathcal YD (C)} = ^H_H {\mathcal YD}$, called Yetter-Drinfeld
  category.
  It follows from  \cite {RT93} and \cite [Theorem 4.1.1]{BD98} that
  $^H_H {\mathcal YD}({\mathcal C})$ is a braided tensor
  category with $^{ YD}C _{U,V} = (\alpha _V \otimes id_U) (id _H \otimes C_{U, V})
  (\phi _U\otimes id _V)$ for any
  two Yetter-Drinfeld modules  $(U, \phi _U, \alpha _U)$ and $(V, \phi _V, \alpha
  _V)$  when $H$ has an invertible antipode. In this case,
$^{YD}C _{U,V}^{-1} = (id _V \otimes \alpha _U)( C_{H,V}^{-1}
\otimes id _U) (S^{-1} \otimes C_{U, V}^{-1}) (C_{U,H}^{-1} \otimes
id _V)
  ( id _U \otimes \phi _V).$ Algebras, coalgebras and Hopf algebras
  and so on
  in $^H_H {\mathcal YD}$ are called  Yetter-Drinfeld ones or YD ones in
  short.

The structure of a Nichols algebra appeared first in the paper
\cite{Ni78} and was rediscovered later by several authors.

Let $V$ be a YD $H$-module in ${\mathcal V}ect (k)$. A graded
braided Hopf algebra $R=\bigoplus_{n\geq 0}R(n)$ in $^H_H{\cal YD}$
is called a {\it Nichols algebra} of $V$ if the following are
satisfied:\\
$\bullet$  $k\cong R(0)$,\\
$\bullet$  $V\cong R(1)$ in $^H_H{\cal YD}$,\\
$\bullet$  $R(1)=P(R)$, the vector space of the primitive elements
of $R$,\\
$\bullet$  $R$ is generated as an algebra by $R(1)$.

Remark \footnote { For any YD $H$-module $V$, it follows from
\cite[Proposition 2.2]{AS02} that there exists a unique Nichols
algebra ${\cal B}(V)$ of $V$ up to graded Hopf algebra isomorphism
in $^H_H{\cal YD}$. The structure of ${\cal B}(V)$ can be described
as follows:

Let $V$ be a YD $H$-module. Then the tensor algebra
$T(V)=\bigoplus_{n\geq 0}T(V)(n)$ of the vector space $V$ admits a
natural structure of a YD $H$-module, since $^H_H{\cal YD}$ is a
braided monoidal category. It is then an algebra in $^H_H{\cal YD}$.
Let $T(V)\underline{\otimes}T(V)$ be the tensor product algebra of
$T(V)$ with itself in $^H_H{\cal YD}$. Then there is a unique
algebra map $\Delta: T(V)\rightarrow T(V)\underline{\otimes}T(V)$
such that $\Delta(v)=v\otimes 1+1\otimes v$ for all $v\in V$. With
this structure, $T(V)$ is a graded braided Hopf algebra in
$^H_H{\cal YD}$. Let $\tilde I(V)$
be the sum of all YD $H$-submodules $I$ of $T(V)$ such that\\
$\bullet$ $I$ is a graded ideal generated by homogeneous elements
of degree $\geq 2$,\\
$\bullet$ $I$ is also a coideal, i.e., $\Delta(I)\subset I\otimes
T(V)+T(V)\otimes I$.\\
Then ${\cal B}(V)=T(V)/{\tilde I}(V)$. }

\begin {Example} \label {0.6}
Let $V$  be a finite-dimensional vector space over field $k$ and
$V^* = Hom _k (V, k)$. Define maps $b_V: k \rightarrow V\otimes V^*$
and $d_V: V ^* \otimes V \rightarrow k $  by
$$b_V(1) = \sum _iv_i \otimes v^i   \hbox { \ \ and \ \ }
\sum _{i, j } d_V (v^i \otimes v_j ) = v^i (v_j)$$ where $\{ v_i
\mid i = 1, 2, \cdots , n \}$  is any basis of $V$ and $\{v^i \mid i
= 1, 2, \cdots , n \}$ is the dual basis in $V^*$. It is clear that
$b_V$  and $d_V$ are  evaluation and coevaluation of $V$ in ${\cal
V}ect(k),$ respectively.

\end {Example}

\begin {Proposition} \label {0.7}
(i) Let $({\cal C}, C)$  be a braided tensor category and $H$  be an
algebra living in ${\cal C}. $ If $(V, \alpha )$ is a left
$H$-module and has a left duality $V^*$  in ${\cal C}, C)$, then
$V^*$ is a left duality of $V$ in ${}_H{\cal M (C)}$.  Here the
module operation  $\alpha _ {V^*}$ of
              $V^*$ is defined as follows:
\[ \alpha _{V^*} =
     \begin{tangle}
\object{H}\step[2]\object{V^*}\\
\S\step\dd\\
\hx\step\step\coev\\
\id\step\tu \alpha\step[2]\id\\
\ev \step[2]\step\id\\
\step[4]\step\object{V^*}
  \end{tangle}
  \]

(ii) Let $({\cal C}, C)$  be a braided tensor category and $H$  be a
coalgebra living in $({\cal C},C). $ If $(V, \psi )$ is a right
$H$-comodule and has a left duality $V^*$  in $({\cal C}, C)$, then
$V^*$ is a left duality of $V$ in ${\cal M(C)}^{H}$.
                      Here the comodule operation  $\psi  _ {V^*}$ of
                        $V^*$ is defined as follows:
                   \[
                   \psi_{V^*}=
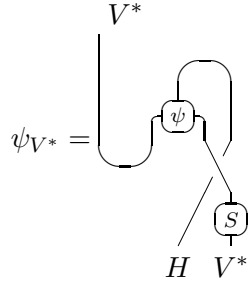
   \begin{tangle}
\step\object{V^*}\\
\id \step[2]\step\coev\\
\id\step[2]\td \psi\step\id\\
\ev\step\step\hx\\
\step[3]\dd\step\S\\
\step[3]\object{H}\step[2]\object{V^*}\\
  \end{tangle}
  \]

\end {Proposition}

{\bf Proof .}  (i)  We easily check that $(V^*, \alpha _{V^*})$ is a
$H$-module. Now we show that  $d_V$  and $b_V$  are $H$-module
homomorphisms. We see that
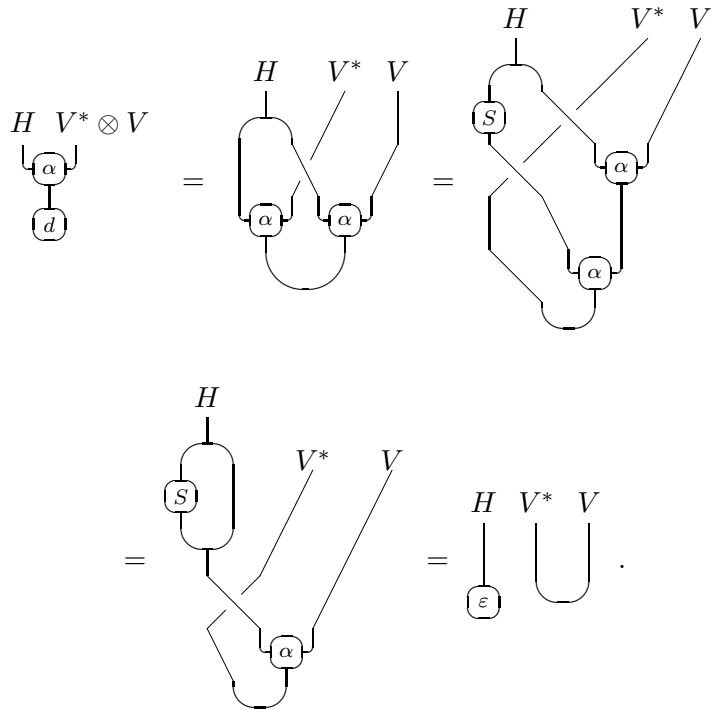
\[
   \begin{tangle}
\object{H}\step[3]\object{V^*\otimes V}\\
\tu \alpha \\
\step\QQ d
  \end{tangle}
  \step [2] \ \ =\step
   \begin{tangle}
\step\object{H}\step[3]\object{V^*}\step[2]\object{V}\\
\cd\step\dd\step[2]\id\\
\id\step[2]\hx\step[2]\dd\\
\tu \alpha\step\tu\alpha\\
\step\Ev
  \end{tangle}
  \step=\step
   \begin{tangle}
\step\object{H}\step[5]\object{V^*}\step[2]\object{V}\\
\cd\step[3]\ne2\step[1]\dd\\
\S\step[2]\x\step[2]\dd\\
\x\step[2]\tu\alpha\\
\id\step[2]\d\step[2]\id\\
\nw2\step[2]\tu\alpha\\
\step[2]\ev
     \end{tangle}
     \]
     \[
  \step=\step
   \begin{tangle}
\step\object{H}\\
\cd\step[2]\step[1]\object{V^*}\step[3]\object{V}\\
\S\step[2]\id\step[2]\dd\step[2]\dd\\
\cu\step\dd\step[2]\dd\\
\step\x\step[2]\dd\\
\step\d\step\tu\alpha\\
\step[2]\ev
     \end{tangle}
       \step=\step
   \begin{tangle}
\object{H}\step[2]\object{V^*}\step[2]\object{V}\\
\id\step[2]\id\step[2]\id\\
\QQ \varepsilon \step[2]\ev
 \end{tangle}\ \ \ .
\] Thus $d_V$
 is an $H$-module homomorphism.

 We also see that
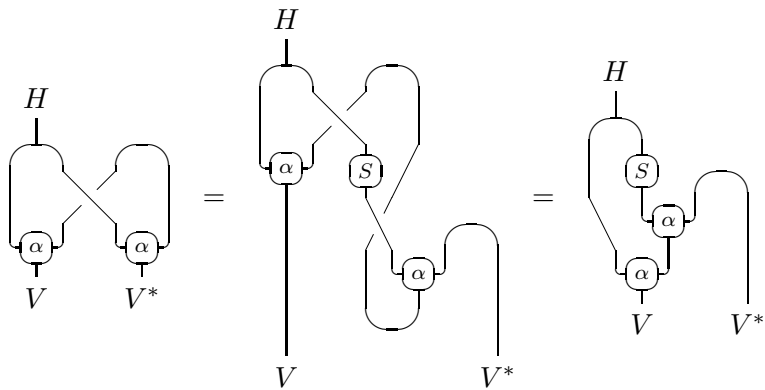
\[
   \begin{tangle}
\step\object{H}\\
\cd\step[2]\coev\\
\id\step[2]\x\step[2]\id\\
\tu \alpha\step[2]\tu\alpha\\
\step\object{V}\step[4]\object{V^*}\\
     \end{tangle}
       \step=\step
\begin{tangle}
\step\object{H}\\
\cd\step[2]\coev\\
\id\step[2]\x\step[2]\id\\
\tu\alpha\step[2]\S\step\dd\\
\step\id\step[2]\step\hx\step[2]\coev\\
\step\id\step[2]\step\id\step\tu \alpha\step[2]\id\\
\step\id\step[2]\step\ev\step[2]\step\id\\
\step\object{V}\step[8]\object{V^*}\\
     \end{tangle}
     \step=\step
\begin{tangle}
\step\object{H}\\
\cd\\
\id\step[2]\S\step[2]\coev\\
\d\step\tu\alpha\step[2]\id\\
\step\tu \alpha\step[2]\step\id\\
\step[2]\object{V}\step[4]\object{V^*}\\
 \end{tangle}
 \]
 \[
     \step=\step
\begin{tangle}
\step\object{H}\\
\cd\\
\id\step[2]\S\\
\cu\step\coev\\
\step\tu\alpha\step[2]\id\\
\step[2]\object{V}\step[3]\object{V^*}\\
 \end{tangle}
     \step=\step
     \begin{tangle}
\object{H}\\
\id \step[2]\coev\\
\QQ \varepsilon \step[2]\object{V}\step[2]\object{V^*}\\
  \end{tangle} \ \ .
\] Thus $b_V$ is an $H$-module homomorphism.

         (ii)  is the dual case of part (i). \begin{picture}(8,8)\put(0,0){\line(0,1){8}}\put(8,8){\line(0,-1){8}}\put(0,0){\line(1,0){8}}\put(8,8){\line(-1,0){8}}\end{picture}

  Every algebra (coalgebra, Hopf algebra, quasitriangular Hopf algebra )
    living in
  ${\cal V}ect (k)$  is called an algebra ( coalgebra, Hopf algebra,
  quasitriangular Hopf algebra )  or an ordinary algebra ( ordinary coalgebra, ordinary Hopf algebra,
  ordinary quasitriangular Hopf algebra ). We usually use $\phi $ and
$\delta ^-$ to denote  left comodule operations; $\psi $ and $\delta
^+$ to denote  right comodule operations; $\alpha  $ and $\alpha ^-$
to denote  left module operations; $\beta   $ and $\alpha ^+$ to
denote right module operations. We use The Sweedler's and
Montgomery's sigma notations ( see \cite {Sw69a} \cite {Mo93}) for
coalgebras and comodules are $\Delta (x) = \sum x_{1}\otimes x_{2}$,
$\phi (x)= \sum x^{(-1)} \otimes x^{(0)}$, $\psi (x)= \sum x^{(0)}
\otimes x^{(1)}$;  or $\Delta (x) = \sum x_{(1)}\otimes x_{(2)}$,
$\phi (x)= \sum x_{(-1)} \otimes x_{(0)}$, $\psi (x)= \sum x_{(0)}
\otimes x_{(1)}$.

In fact, if $H$  is an ordinary algebra and $V$  is a
finite-dimensional left $H$-module, then $V$ has a left duality
$H^*$ in category ${}_H{\cal M}.$

\chapter {
The Maschke's Theorem for Braided Hopf Algebras}\label {c2}

In this chapter, we obtain the fundamental theorem of Hopf modules
living in braided tensor   categories. In Yetter-Drinfeld category,
we obtain that the antipode $S$  of $H$ is invertible and the
integral $\int _H^l$  of $H$  is a one-dimensional space for any
finite-dimensional braided  Hopf algebra $H$. We give Maschke's
Theorem for YD Hopf algebras.

Supersymmetry has attracted a great deal of interest from physicists
and mathematicians (see \cite {Ka77}, \cite {Sc79}, \cite {MR94} ).
It finds numerous applications in particle physics and statistical
mechanics ( see \cite {CNS75} ) . Of central importance in
supersymmetric theories is the  $Z_2$-graded structure which
permeates them. Superalgebras and super-Hopf algebras are  most
naturally very important.

Braided tensor categories were introduced by Joyal and  Street in
1986 \cite {JS86}. It is a generalization of super case. Majid,
Joyal, Street and Lyubasheko have obtained many interesting
conclusions in braided tensor categories, for example, the braided
reconstruction theorem, transmutation and bosonisation, integral,
q-Fourier transform, q-Mikowski space, random walk and so on  (see
\cite {Ma93a}, \cite {Ma95a}, \cite {Ma95b}, \cite {MR94} , \cite
{Ly95} ).

In many applications one needs to know the integrals  of  braided
Hopf algebras and whether a braided Hopf algebra is semisimple or
finite-dimensional. A systematic study of braided Hopf algebras,
especially  of those braided Hopf algebras that play a role in
physics (as, e.g. braided line algebra ${\bf C} [x]$ in \cite
{Ma93a} and anyonic line algebra ${\bf C} [\xi] \cong {\bf C} \{x\}
/<x^n>$ in \cite {MR94} )  is very useful.

V.Lyubasheko in\cite {Ly95} introduced the integral for braided Hopf
algebras and gave a formula about integral. S. Majid in \cite {MR94}
obtained the integral of anyonic line algebra. In this chapter, we
obtain that the antipode $S$  of $H$ is invertible and the integral
$\int _H^l$  of $H$  is a one-dimensional space for any
finite-dimensional (super) anyonic  or q-statistical Hopf algebra
$H$. In particular, if $H$  is a finite-dimensional (super) anyonic
Hopf algebra, then $\int _H^l$ is a homogeneous space. We show that
a finite-dimensional (super) anyonic or q-statistical Hopf algebra
$H$  is semisimple iff $\epsilon (\int _H^l) \not=0.$ Furthermore,
if $H$ is a semisimple (super) anyonic Hopf algebra, then $H$ is
finite-dimensional.

We begin by recalling some definitions and properties. Attempts to
find solutions of the Yang-Baxter equation in a systematic way have
led to the theory of quantum groups, in which the concept of
(co)quasitriangular Hopf algebra is introduced .

       A pair $(H,r)$  is called coquasitriangular Hopf algebra if
 $H$  is a Hopf algebra over field $k$ and there exists a convolution-
 invertible $k$-bilinear  form $r$: $H \otimes H \rightarrow k$
 such that the following conditions hold:  for any $a, b, c \in H,$

(CQT1)  $r(a, bc) = \sum r(a_1, c) r(a_2, b);$

(CQT2)  $r(ab, c) = \sum r(a, c_1) r(b, c_2);$

 (CTQ3)  $ \sum r(a_1 , b_1) a_2 b_2 = \sum b_1 a_1 r(a_2, b_2).$

Furthermore, if $r^{-1} (a, b) = r(b, a)$  for any $a, b \in H,$
then $(H, r)$ is called    cotriangular.

For coquasitriangular Hopf algebra $(H,r),$  we can define a
braiding $C^r$ in the category ${}^ H {\cal M}$  of $H$-comodules as
follows:

%$  \left C_{U,V} ^r  \begin {array} {ll} H \otimes R & \longrightarrow R \\
% u \otimes v & \mapsto   r(v^{(-1)}, u^{(-1)})u^{(0)} v^{(0)}
% \end {array} \right. $

$$   C_{U,V} ^r  : U \otimes V  \longrightarrow  V \otimes U  $$
sending $ u \otimes v $  to  $\sum  r(v^{(-1)}, u^{(-1)})\otimes
v^{(0)} u^{(0)} $ \ \ \
   for any $u \in U, v\in V,$ where
  $(U, \phi _U)$ and $(V,\phi _V)$  are $H$-comodule and
  $\phi _V (v) = \sum v ^{(-1)} \otimes v^{(0)}$ (see \cite [Proposition
  VIII.5.2 ] {Ka95}).
  This braided tensor category is called  one determined by coquasitriangular
  structure. Dually, we have the concept of quasitriangular Hopf algebra.
  We also have that $({}_H{\cal M} , C^R)$  is a braided tensor category,
  which is called one determined by quasitriangular structure.

For example, let $H={\bf C Z}_n$ and $r(a,b) = (e ^{\frac {2 \pi
i}{n}})^{ab}$  for any $a, b \in {\bf Z}_n$, where ${\bf C}$ is the
complex field. It is clear that $({\bf CZ}_n, r )$  is a
coquasitriangular Hopf algebra. Thus $ ({} ^{{\bf CZ} _n} {\cal M} ,
C^r)$ is a braided tensor category, usually  written as ${\cal
C}_n.$  Every algebra or Hopf algebra living   in ${\cal C}_n $ is
called an anyonic algebra or anyonic Hopf algebra (see  \cite
[Example 9.2.4] {Ma95b}).
 Every algebra or Hopf algebra living   in ${\cal C}_2 $
is called a superalgebra or super-Hopf algebra.

If $H = {\bf CZ} $ and $r(a,b) = q ^{ab}$ for any $a, b \in {\bf
Z},$ then $ ({}^{\bf CZ} {\cal M}, C^r)$ is a braided tensor
category, where $0\not=q \in {\bf C}$.
  Every algebra or Hopf algebra living   in it
is called  a q-statistical algebra or q-statistical  Hopf algebra
(see \cite [Definition 10.1.9] {Ma95b}).

For example,  let $(H,r)$ be a cocommutative cotriangular Hopf
algebra. If $(L, \phi )$ is a left $H$-comodule with an $H$-comodule
homomorphism $[\ ,\ ]$ from $L\otimes L$  to $L$ in $({}^H{\cal M},
C^r)$ such that the following conditions hold:

(i) $[\ ,\ ] (id + C_{L,L} ^r) =0;$

(ii)  $[\ ,\ ] ([\ ,\ ]\otimes id ) (id + (id _L \otimes C_{L,L} ^r)
(C^r_{L,L} \otimes id _L) + (C^r_{L,L} \otimes id _L)(id _L \otimes
C^r_{L,L} ))=0,$

then $(L,[\ ,\ ])$  is called a Lie $H$-algebra or a Lie algebra
living in     $({}^H {\cal M}, C^r_{L,L})$ (see \cite [Definition
1.5] {BFM96}).  If $H = {\bf CZ}_2$ and $r(a,b) = (-1)^{ab}$ for any
$a,b \in {\bf Z}_2$, then every Lie $H$-algebra is called a Lie
superalgebra (see \cite {Ka77}). If $H = kG$ and $r$  is a
bicharacter of abelian group $G$ with $r^{-1}(a,b)= r(b,a)$ for any
$a, b \in G$, then every Lie $H$-algebra is called a Lie colour
algebra (see \cite {Mo93}) or $\epsilon $ Lie algebra (see \cite
[Definition 2.1] {Sc79}). If $(L, [\ ,\ ])$ is a Lie colour algebra,
then the universal enveloping  algebra $U(L)$ of $L$  is an algebra
living in $({}^H {\cal M}, C^r)$ \ \ \ (see \cite {Sc79}). Thus the
universal enveloping algebras of Lie superalgebras are
superalgebras.

\section {The integrals in braided tensor categories} \label {s1}

 In this section, we obtain the fundamental theorem of Hopf
modules living in braided tensor   categories.

An object $E$ in category ${\cal C}$ is called a final object if for
every  object  $V$ in ${\cal C}$, there exists a unique morphism $f$
from $V$ to $E$.  For two morphisms $f_1 , f_2 : A \rightarrow B$ in
${\cal C}$,  set $target (A, f_1, f_2) =: \{ (D, g) \mid g$ is a
morphism from $D$ to $A$ such that $ f_1 g = f_2g  \}$ and $source
(B, f_1, f_2) =: \{ (D, g) \mid g$  is a morphism from $B$ to $D$
such that $ gf_1 = gf_2 \}$, the final object of full subcategory
$target (A, f_1, f_2)$
 and the initial object of full subcategory $source  (B, f_1, f_2)$
are called equivaliser and coequivaliser of $f_1$ and $f_2$, written
as $equivaliser (f_1. f_2)$ and
 $coequivaliser (f_1, f_2)$, respectively.
If  $(M, \phi)$ is
  a left comodule over bialgebra $H$ in ${\cal C}$,  then  $equivaliser (\phi ,  \eta _H \otimes
id _M)$ is called  the coinvariants of $M$, written $(M ^{coH}, i).$
That is, for every $H$-comodule $N$ in ${\cal M}$ and a morphism $f$
from $N$ to $M$ in $^H{\cal M(C)}$, then there exists a unique a
morphism $\bar f$ from $N$ to $M^{co \ H}$ in $^H{\cal M(C)}$ such
that $f = i \bar f $.

\begin {Definition} \label {1.1.2}
Let $H$ be a bialgebra, $(M, \alpha)$  be a left $H$-module and $(M,
\phi)$ a left $H$-comodule. If $\alpha$ is an $H$-comodule morphism
and $\phi$ is an $H$-module morphism, then $(M,\alpha , \phi)$ is
called an $H$-Hopf module.
\end {Definition}
 We easily know that $(M,\alpha , \phi)$ is an $H$-Hopf module iff
 $\step\begin{tangle}\step[1]\object{H\otimes M}\\
\morph {\alpha }\\
\td \phi\\
\object{H}\step[2]\object{M}
\end{tangle}\step=\step
\begin{tangle}
\step\object{H}\step[4]\object{M}\\
\cd\step[2]\td\phi\\
\id\step[2]\x\step[2]\id\\
\cu\step[2]\tu\alpha\\
\step \object{H}\step[4]\object{M}
\end{tangle}\ \ \ .$

\begin {Theorem} \label {1.1.3}   (The fundamental theorem of Hopf modules)
Let $({\mathcal C}, C)$ be a braided tensor category
 If  $(M,
\alpha , \phi)$ is an $H$-Hopf module in braided tensor category
$({\cal C},C)$ and there exists the  equivaliser $(\phi ,  \eta _H
\otimes id _M)$, then
             $$H \otimes M^{coH} \cong M   \hbox { \ \ \ ( as } H\hbox {-Hopf
             modules.})$$
\end {Theorem}
{\bf Proof. } Let $u=\step\begin{tangle}\step \object{M}\\
\td \phi\\
\S \step[2]\id\\
\tu \alpha\\
\step \object{M}\end{tangle}\step$. We have that $ \phi u=\step\begin{tangle}\step \object{M}\\
\td \phi\\
\S \step[2]\id\\
\tu \alpha\\
\td \phi\\
\object{H} \step[2] \object{M}\end{tangle}\step=\step
\begin{tangle}\step[2]\object{M}\\
\step\td \phi\\
\step\S \step[2]\d\\
\cd\step\td \phi\\
\id \step[2]\hx\step[2]\id\\
\cu\step\tu\alpha\\
\step \object{H} \step[3] \object{M}\end{tangle}\step \stackrel
 {\hbox { by
Pro.\ref {12.1.4}(i)} }{=} \step \begin{tangle}\step[2]\object{M}\\
\step\td \phi\\
\cd \step\nw2\\
\S\step[2]\S\step[2]\td \phi\\
\d\step\id\step[2]\id\step[2]\id\\
\step\hx\step\dd\step[2]\id\\
\step\id\step\hx\step[2]\dd\\
\step\hcu\step\tu\alpha\\
 \step\hstep \object{H} \step[2.5]
\object{M}\end{tangle}\step =\step\begin{tangle}\Q {\eta_{H}}\\
\id\\\object{H} \end{tangle}\step[2]\begin{tangle}\step \object{M}\\
\td \phi\\
\S \step[2]\id\\
\tu \alpha\\
\step \object{M}\end{tangle}$ \ \ . Thus there exists an
$H$-comodule homomorphism $\bar u$ from $M$ to $M^{coH}$ such that
$u = i \bar u.$.

Let $\xi=\begin{tangle}\object{H}\step[2]\object{M^{coH}}\\
\id \step [1] \morph i\\
\tu \alpha\\
\step \object{M}\end{tangle}\step$ and  $\step\psi =\begin{tangle}\step\object{M}\\
\td \phi\\
\id\step\morph {\bar u}\\
\object{H}\step[2]\object{M}\end{tangle}$. We see that $\xi\psi=\begin{tangle}\step \object{M}\\
\td \phi\\
\id\step\morph u\\
\tu \alpha\\
\step \object{M}\end{tangle}=\begin{tangle}\step \object{M}\\
\td \phi\\
\id\step\td\phi\\
\id\step\S \step[2]\id\\
\id\step\tu \alpha\\
\tu \alpha\\
\step \object{M}\end{tangle}=\begin{tangle} \object{M}\\
\id\\
\object{M}\end{tangle}\step$ and
\[
\step\psi\xi\hstep=\hstep\begin{tangle}\object{H}\step[2] \object{M^{coH}}\\
\tu \alpha\\
\td\phi\\
\id \step\td\phi\\
\id\step\S \step[2]\id\\
\id\step\tu \alpha\\
\object{M}\step[3] \object{M^{coH}}\end{tangle} \hstep=\hstep
\begin{tangle}
\step\object{H}\step[4]\object{M}\\
\cd\step[2]\td \phi\\
\id\step[2]\x\step[2]\nw2\\
\cu\step\cd\step[2]\td\phi\\
\step\id\step[2]\id\step[2]\x\step[2]\id\\
\step\id\step[2]\cu\step[2]\tu\alpha\\
\step\id\step[3]\S\step[3]\ne2\\
\step\object{H}\step[3]\tu \alpha\\
\step[5]\object{M^{coH}}
\end{tangle}\hstep=\hstep
\begin{tangle}
\step\object{H}\step[5]\object{M^{coH}}\\
\cd\step[4]\id\\
\id\step\cd\step[2]\td \phi\\
\id\step\id\step[2]\x\step[2]\id\\
\id\step\cu\step[2]\tu \alpha\\
\id\step[2]\S\step[2]\step\ne2\\
\id\step[2]\tu\alpha\\
\object{H}\step[4]\object{M^{coH}}
\end{tangle}\hstep=\hstep
\begin{tangle}
\step\object{H}\step[4]\object{M^{coH}}\\
\cd\step[3]\id\\
\id\step\cd\step[2]\id\\
\id\step[1]\S\step[2]\tu\alpha\\
\id\step\id\step[2]\dd\\
\id\step[1]\tu\alpha\\
\object{H}\step[3]\object{M^{coH}}
\end{tangle}\hstep=\hstep
\begin{tangle}
\object{H}\step[2]\object{M^{coH}}\\
\id\step[2]\id\\
\object{H}\step[2]\object{M^{coH}}
\end{tangle}\step.
\]
Thus $\xi$ is an isomorphism.

Since $\step\begin{tangle}\step[2]\object{H\otimes M^{coH}}\\
\morph {\xi}\\
\td \phi\\
\object{H}\step[3]\object{M^{coH}}
\end{tangle}\step=\step
\begin{tangle}
\step\object{H}\step[5]\object{M^{coH}}\\
\cd\step[2]\td\phi\\
\id\step[2]\x\step[2]\id\\
\cu\step[2]\tu\alpha\\
\step \object{H}\step[5]\object{M^{coH}}
\end{tangle}\step=\step\begin{tangle}\step[2]\object{H\otimes M^{coH}}\\
\td \phi\\
\id\step\morph {\xi}\\
\object{H}\step[3]\object{M^{coH}}
\end{tangle}\step[2]$ and

\[
\begin{tangle}
\step[2]\object{H}\\
\step\cd\step[2]\object{H}\step[3]\object{M^{coH}}\\
\cd\step\id\step\dd\step[2]\dd\\
\id\step[2]\S\step\hx\step[2]\dd\\
\id\step[2]\hx\step\tu\alpha\\
\cu\step\id\step[2]\id\\
\step\cu\step\dd\\
\step[2]\tu \alpha\\
\step[3]\object{M}
\end{tangle}\step=\step
\begin{tangle}
\step[2]\object{H}\\
\step\cd\step[3]\object{H}\step[3]\object{M^{coH}}\\
\dd\step\cd\step\dd\step[2]\id\\
\id\step[2]\S\step[2]\hx\step[2]\step\id\\
\id\step[2]\x\step\d\step[2]\id\\
\cu\step[2]\cu\step\dd\\
\step\d\step[2]\dd\step\dd\\
 \step[2]\cu\step\dd\\
 \step[2]\step[1]\tu\alpha\\
 \step[4]\object{M}
\end{tangle}\step=\step
\begin{tangle}
\object{H}\step[2]\object{H}\step[2]\object{M^{coH}}\\
\cu\step[1]\id\\
\step\tu\alpha\\
\step[2]\object{M}
\end{tangle}\step=\step
\begin{tangle}
\object{H}\step[4]\object{H\otimes M^{coH}}\\
\d\step[1]\morph \xi\\
\step\tu\alpha\\
\step[2]\object{M}
\end{tangle}\ \ \ .
\]
Consequently, $\xi$ is an $H$-Hopf module isomorphism. $\Box$

\begin {Lemma} \label {1.1.4}
If $H$ is a Hopf algebra with left duality, then $(H^*, \alpha ,
\phi )$ is a left $H$-Hopf module, where the operations of module
and comodule are
 defined
as follows:\[
\begin{tangle}
\object{H}\step[2]\object{H^*}\\
\tu\alpha\\
\step[1]\object{H^*}
\end{tangle}\step=\step
\begin{tangle}
\object{H}\step[2]\object{H^*}\\
\tu\rightharpoondown\\
\step[1]\object{H^*}
\end{tangle}\step=\step
\begin{tangle}
\object{H}\step[3]\object{H^*}\\
\S\step[2]\cd\\
\x\step[2]\id\\
\id\step[2]\x\\
\id\step[2]\ev\\
\object{H^*}
\end{tangle}\ \ \ '\ \
\begin{tangle}
\step[1.5]\object{H^*}\\
\td\phi\\
\object{H}\step[2]\object{H^*}\\
\end{tangle}\step=\step
\begin{tangle}
\step[2]\step[2]\object{H^*}\\
\coev\step[2]\id\\
\id\step\step\cu\\
 \object{H}\step[3]\object{H^*}
\end{tangle}\ \ \ .
\]

\end {Lemma}

{\bf Proof.} We see that

\[
\begin{tangle}
\object{H}\step[2]\object{H}\step\hstep\object{H^*}\\
\cu\step\id\\
\step\tu\alpha\\
\step[2]\object{H^*}
\end{tangle}\step=\step
\begin{tangle}
\object{H}\step[2]\object{H}\step[2]\object{H^*}\\
\cu\step[2]\id\\
\step\S\step[2]\cd\\
\step\x\step[2]\id\\
\step\id\step[2]\x\\
\step\id\step[2]\ev\\
\step\object{H^*}
\end{tangle}
\step\stackrel{ \hbox {by Pro.\ref {12.1.4}(i)} }{=}\step
\begin{tangle}
\object{H}\step[2]\object{H}\step[2]\object{H^*}\\
\S\step[2]\S\step[2]\id\\
\x \step [2 ]\id\\

\cu\step\cd\\
\step\x\step[2]\id\\
\step\id\step[2]\x\\
\step\id\step[2]\ev\\
\step\object{H^*}
\end{tangle}
\]

\[
\step\stackrel{ \hbox {by Pro.\ref {12.1.4}(ii)}}{=}\step
\begin{tangle}
\object{H}\step[2]\object{H}\step[3]\object{H^*}\\
\S\step[2]\S\step[2]\cd\\
\id\step[2]\x\step[2]\id\\
\x\step[2]\x\\
\id\step[2]\x\step[2]\id\\
\id\step\cd\step\x\\
\id\step\d\step\hev\step\ne2\\
\id\step[2]\ev\\
 \object{H^*}
\end{tangle}
\step=\step
\begin{tangle}
\object{H}\step[2]\object{H}\step[3]\object{H^*}\\
\S\step[2]\S\step[2]\cd\\
\id\step[2]\x \step\cd\\
\x\step[2]\hx\step[2]\id\\
\id\step[2]\x\step\x\\
\id\step[2]\ev\step\ev\\
 \object{H^*}
\end{tangle}
\step=\step
\begin{tangle}
\object{H}\step[1]\object{H}\step[2]\object{H^*}\\
\id\step\tu\alpha\\
\tu\alpha\\
 \step\object{H^*}
\end{tangle}
\]
\[
\begin{tangle}
\Q {\eta_{H}}\step[2]\object{H^*}\\
\tu\alpha\\
\step[1]\object{H^*}
\end{tangle}\step=\step
\begin{tangle}
\step[3]\object{H^*}\\
\Q {\eta_{H}}\step[2]\cd\\
\x\step[2]\id\\
\id\step[2]\x\\
\id\step[2]\ev\\
\object{H^*}
\end{tangle}
\step=\step
\begin{tangle}
\object{H^*}\\
\id\\
 \object{H^*}
\end{tangle}\ \ \ .
\]
Thus $(H^*, \alpha )$ is an $H$-module.

We see that
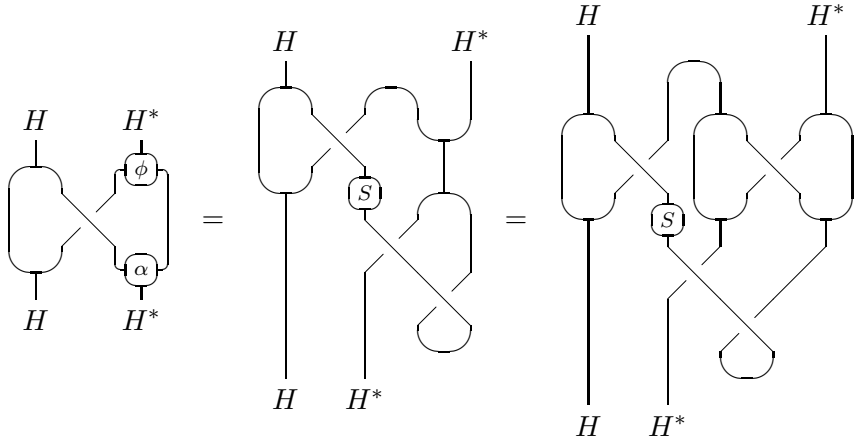
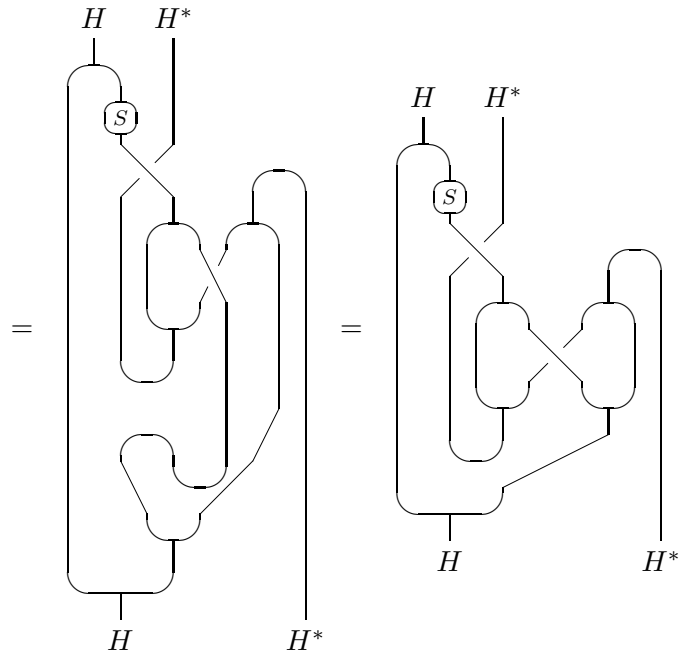
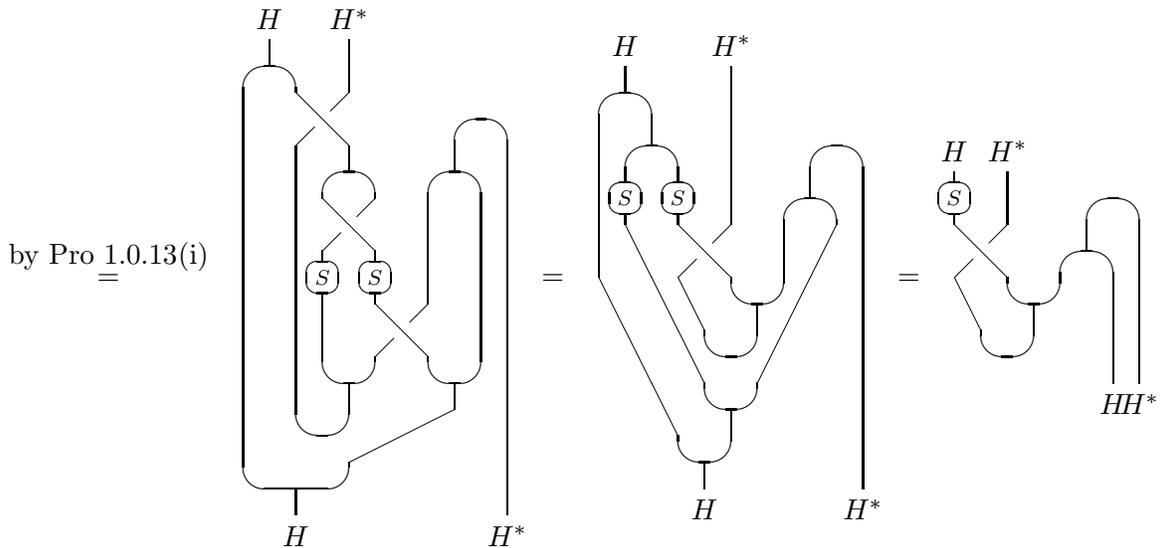
\[
\begin{tangle}
\step\object{H}\step[4]\object{H^*}\\
\cd\step[2]\td\phi\\
\id\step[2]\x\step[2]\id\\
\cu\step[2]\tu\alpha\\
\step \object{H}\step[4]\object{H^*}
\end{tangle}\step=\step
\begin{tangle}
\step\object{H}\step[7]\object{H^*}\\
\cd\step[2]\coev\step[2]\id\\
\id\step[2]\x\step[2]\cu\\
\cu\step[2]\S\step[2]\cd\\
\step\id\step[2]\step\x\step[2]\id\\
\step\id\step[2]\step\id\step[2]\x\\
\step\id\step[2]\step\id\step[2]\ev\\
 \step\object{H}\step[3]\object{H^*}\\
\end{tangle}
\step=\step
\begin{tangle}
\step\object{H}\step[9]\object{H^*}\\
\step\id\step[2]\step[1]\coev\step[2]\step[2]\id\\
\cd\step\step\id\step\cd\step[2]\cd\\
\id\step[2]\x\step\id\step[2]\x\step[2]\id\\
\cu\step[2]\S\step\cu\step[2]\cu\\
\step\id\step[2]\step\x\step[2]\sw2\\
\step\id\step[2]\step\id\step[2]\x\\
\step\id\step[2]\step\id\step[2]\ev\\
 \step\object{H}\step[3]\object{H^*}\\
\end{tangle}
\]

\[
\step =\step
\begin{tangle}
\step\object{H}\step[3]\object{H^*}\\
\cd \step[2]\id\\
\id\step\step\S\step[2]\id\\
\id \step[2]\x \step[2] \step[1]\coev\\
\id \step[2]\id\step\cd \step[1]\cd\step\id\\
\id \step[2]\id\step \id\step[2]\hx \step[2]\id\step[1]\id\\
\id \step[2]\id\step\cu \step[1]\id\step[2]\id\step[1]\id\\
\id \step[2]\ev \step[2]\id \step[2]\id\step\id\\
\id\step[2]\coev\step[2]\id\step\dd\step\id\\
\id\step[2]\d\step\ev\ne2\step[2]\id\\
\id\step[2]\step\cu\step[2]\step[2]\id\\
\Cu\step[2]\step[2]\step\id\\
 \step[2]\object{H}\step[7]\object{H^*}\\
\end{tangle}
\step =\step
\begin{tangle}
\step\object{H}\step[3]\object{H^*}\\
\cd \step[2]\id\\
\id\step\step\S\step[2]\id\\
\id \step[2]\x \step[2] \step[2]\coev\\
\id \step[2]\id\step\cd \step[2]\cd\step\id\\
\id \step[2]\id\step \id\step[2]\x \step[2]\id\step\id\\
\id \step[2]\id\step\cu \step[2]\cu\step\id\\
\id \step[2]\ev \step[2]\step\ne4 \step[2]\id\\
\Cu \step[2] \step[2] \step[2]\id\\
 \step[2]\object{H}\step[8]\object{H^*}\\
\end{tangle}
\]

\[
\step\stackrel{\hbox {by Pro \ref {12.1.4}(i)}}{=}\step
\begin{tangle}
\step\object{H}\step[3]\object{H^*}\\
\cd \step[2]\id\\
\id \step[2]\x \step[2] \step[2]\coev\\
\id \step[2]\id\step\cd \step[2]\cd\step\id\\
\id \step[2]\id\step\x \step[2]\id \step [2]\id\step\id\\

\id \step[2]\id\step\S \step[2]\S \step[2]\id \step[2]\id \step\id\\
\id \step[2]\id\step \id\step[2]\x \step[2]\id\step\id\\
\id \step[2]\id\step\cu \step[2]\cu\step\id\\
\id \step[2]\ev \step[2]\step\ne4 \step[2]\id\\
\Cu \step[2] \step[2] \step[2]\id\\
 \step[2]\object{H}\step[8]\object{H^*}\\
\end{tangle}
\step=\step
\begin{tangle}
\step\object{H}\step[4]\object{H^*}\\
\cd\step[2]\step\id\\
\id\step\cd\step[2]\id\step[2]\step\coev\\
\id\step\S\step[2]\S\step[2]\id\step[2]\cd\step\id\\
\id\step\d\step\x\step[2]\id\step\dd\step\id\\
\d\step\d\d\step\cu\dd\step[2]\id\\
\step\d\step\d\ev\dd\step[2]\step\id\\
\step[2]\d\step\cu\step[2]\step[2]\id\\
\step[2]\step\cu\step[2]\step[2]\step\id\\
 \step[4]\object{H}\step[6]\object{H^*}\\
\end{tangle}
\step=\step
\begin{tangle}
\object{H}\step[2]\object{H^*}\\
\S\step[2]\id\step[2]\step\coev\\
\x\step[2]\cd\step\id\\
\d\step\cu\step[2]\id\step\id\\
\step\ev\step[2]\step\id\step\id\\
 \step[6]\object{H}\step[1]\object{H^*}\\
\end{tangle}
\]
and
\[
\begin{tangle}
\object{H}\step[2]\object{H^*}\\
\tu\alpha\\
\td\phi\\
\object{H}\step[2]\object{H^*}\\
\end{tangle}
\step=\step
\begin{tangle}
\step[2]\step[2]\object{H}\step[3]\object{H^*}\\
\step[2]\step[2]\S\step[2]\cd\\
\coev\step[2]\x\step[2]\id\\
\id\step[2]\cu\step[2]\x\\
\id\step[2]\step\id\step[2]\step\ev\\
\object{H}\step[3]\object{H^*}\\
\end{tangle}
\step=\step
\begin{tangle}
\object{H}\step[3]\object{H^*}\\
\S\step[2]\step\id\step[2]\step\coev\\
\id\step[2]\step\id\step[2]\ne2\coev\d\\
\nw2\step[2]\id\step\cu\step[2]\id\step\id\\
\step[2]\nw2\ev\step[2]\ne2\step\id\\
\coev\step[2]\x\step[2]\dd\\
\id\step[2]\cu\step[2]\x\\
\id\step[2]\step\id\step[2]\step\ev\\
\object{H}\step[3]\object{H^*}\\
\end{tangle}
\]
\[
\step=\step
\begin{tangle}
\object{H}\step[2]\object{H^*}\\
\S\step[2]\id\\
\x\step[2]\coev\\
\d\step\cu\step[2]\id\\
\step\ev\step[2]\step\id\\
\step[2]\coev\step[2]\id\\
\step[2]\id\step[2]\cu\\
\step[2]\object{H}\step[3]\object{H^*}\\
\end{tangle}
\step=\step
\begin{tangle}
\object{H}\step[2]\object{H^*}\\
\S\step[2]\id\step\step[2]\coev\\
\x\step[2]\cd\step\id\\
\d\step\cu\step[2]\id\step\id\\
\step\ev\step\step[2]\id\step\id\\
\step[6]\object{H}\step[1]\object{H^*}\\
\end{tangle}\ \ \ .
\]

Thus
\[
\begin{tangle}
\object{H}\step[2]\object{H^*}\\
\tu\alpha\\
\td\phi\\
\object{H}\step[2]\object{H^*}\\
\end{tangle}
\step=\step
\begin{tangle}
\step\object{H}\step[4]\object{H^*}\\
\cd\step[2]\td\phi\\
\id\step[2]\x\step[2]\id\\
\cu\step[2]\tu\alpha\\
\step \object{H}\step[4]\object{H^*}
\end{tangle}\ \ \ .
\]

Consequently, $(H^*,\alpha , \phi)$ is an $H$-Hopf module. $\Box$

\section {The Integrals of Braided  Hopf algebras in
Yetter-Drinfeld Categories }\label {s2}

In this section, using the fundamental theorem of  Hopf  modules in
 braided tensor  categories,
 we give the relation between the integrals and semisimplicity
 of  Braided  Hopf algebras in
Yetter-Drinfeld Categories.

Assume that $H$ is a YD Hopf algebra and $(M, \phi)$ is a $H$-
comodule in $^B_B {\cal YD} $. Obviously, $$\{ x \in M \mid \phi (x)
=\eta _H \otimes   x \}= ker (\phi - \eta _H \otimes id _M).$$ It is
a straightforward verification that $(M ^{coH}, i)$ is the
$equivaliser (\phi ,  \eta _H \otimes id _M)$ with
$$M^{coH} := \{ x \in M \mid \phi (x) =\eta _H \otimes   x \}$$ and
inclusion map $i$. Since $\phi - \eta _H \otimes id _M$ is a
morphism in $^B_B {\cal YD}$, $ker (\phi - \eta _H \otimes id _M)$
is object in $^B_B {\cal YD}$. Let $\int _H ^l := \{ x \in H \mid hx
= \epsilon (h)x , \hbox { for any } h\in H\}$, called  the integral
space of $H.$

If $H$ is a finite-dimensional Hopf algebra living in a
Yetter-Drinfeld category, then $H$ has a left duality $H^*$ ( see
\cite [Example 1 P347]{Ka95} or Lemma \ref {0.7}) and  $\int
_{H^*}^l = (H^*)^{co H},$ where the $H$-comodule operation of $H^*$
is defined in Lemma \ref {1.1.4}.

\begin {Theorem} \label {1.2.1} (The Maschke's theorem) Assume $B$
is an ordinary Hopf algebra over field $k$ with invertible antipode.
If $H$ is a finite-dimensional Hopf algebra living in
Yetter-Drinfeld category $^B_B {\cal YD}$, then

(1)         $H \otimes \int _{H^*}^l  \cong H^*$   \ \ \ (as
$H$-Hopf    modules in  $^B_B {\cal YD})$;

(2) $\int _H^l$  is a one dimensional space;

(3) the antipode  $S$ of $H$ is bijective;

(4)  $H$ is  semisimple  iff $\epsilon (\int _H^l)\not=0;$

(5) Every semisimple Hopf algebra $H$ living in ${}_B^B{\cal YD}$ is
finite-dimensional.

\end {Theorem}

{\bf Proof.} (1) It follows from the above discussion.

(2)   It follows from part (1).

(3) According to the proof of Lemma \ref {1.1.3} and Lemma \ref
{1.1.4}, we have that
$$\xi = (id _{H^*} \otimes d_H) (id _{H^*} \otimes C _{H,H^*})
(C_{H,H^*} \otimes id _{ H^*}) (id _H \otimes \Delta) (S \otimes id
_ {\int _{H^*} ^l})$$ is an isomorphism from $H \otimes \int
_{H^*}^l$  to $H^*$,  which implies that $S$ is injective.
Considering $H$ is a finite-dimensional vector space, we obtain that
$S$ is bijective.

(4) If $H$ is semisimple then there is a left ideal $I$ such that
$$H = I\oplus ker \epsilon. $$
For any $y\in I, h\in H$, we see that
\begin {eqnarray*}
 hy &=& ((h- \epsilon (h) 1_H ) + \epsilon (h)1_H )y \\
&=& (h- \epsilon (h) 1_H )y + \epsilon (h) y  \\
  &=& \epsilon (h) y    \hbox { \ \ \ \ \ since }   (h-\epsilon (h) 1_H) y \in
  (ker \epsilon )I =0 \ .
\end {eqnarray*}
Thus $y \in \int _H^l$, and so  $I \subseteq \int _H^l,$ which
implies $\epsilon (\int _H^l)\not=0.$

Conversely, if $\epsilon (\int _H^l) \not=0,$ let $z \in \int _H^l$
with $\epsilon (z) =1.$ Say $M$ is a left $H$-module and $N$ is an
$H$-submodule of $M$.
 Assume that $\xi$  is a $k$-linear projection from $M$ to $N$.
 We define $$\mu (m) = \sum z_1 \cdot \xi (S(z_2)\cdot m)$$
 for any $m \in M.$ It is sufficient to show that $\mu$ is an $H$-module projection
 from $M$ to $N$. Obviously, $\mu$  is a $k$-linear projection. Now we only need to
 show that it is an $H$-module map.
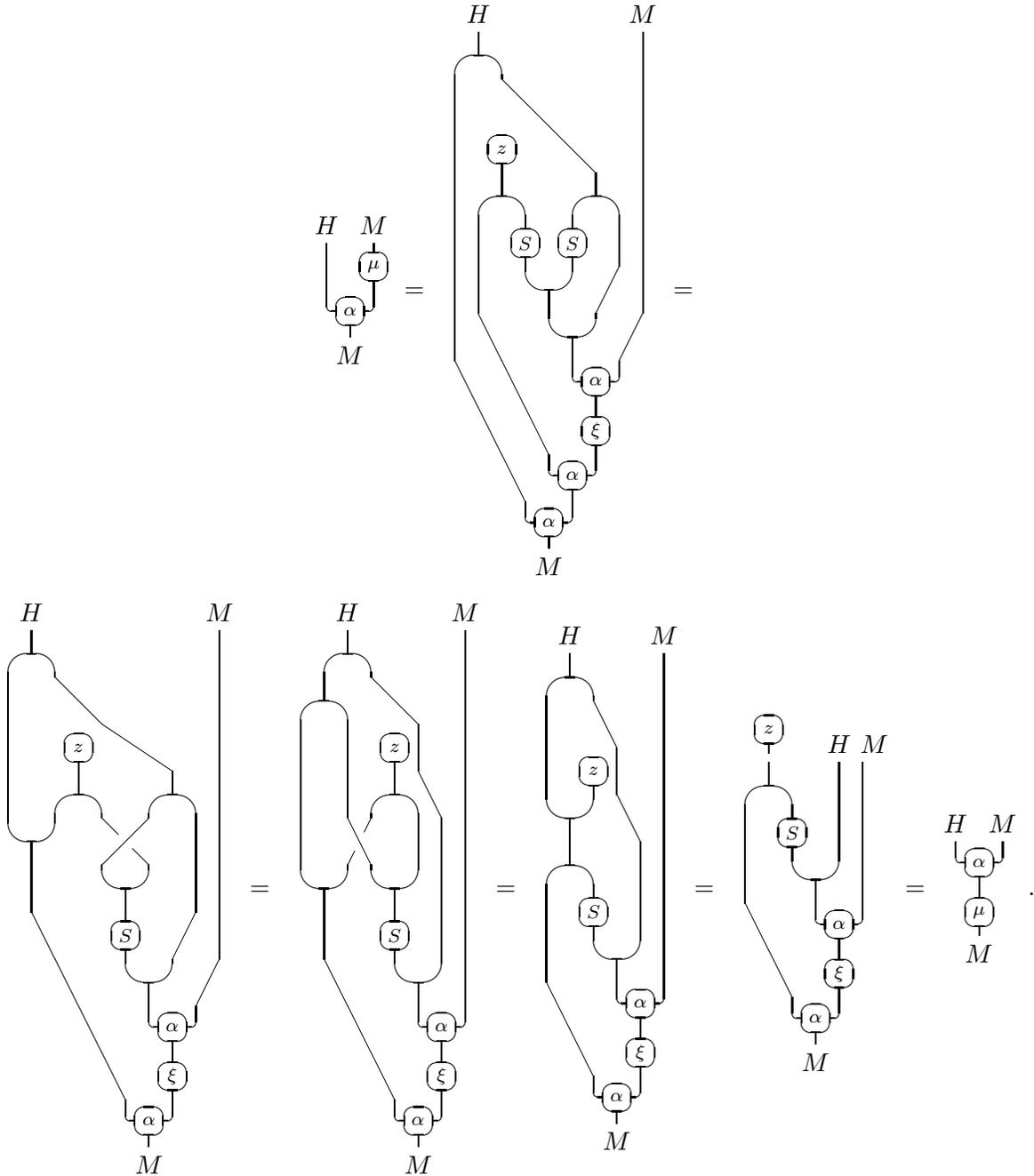
\[
\begin{tangle}
\object{H}\step[2]\object{M}\\
\id\step[2]\O {\mu}\\
\tu {\alpha}\\
\step\object{M}\\
\end{tangle}
 \step=\step
\begin{tangle}
\step\object{H}\step[7]\object{M}\\
\cd\step[6]\id\\
\id\step[2]\nw2\step[5]\id\\
\id\step[2]\Q z\step[2]\nw2\step[3]\id\\
\id\step\cd\step[2]\cd\step\id\\
\id\step\id\step[2]\S\step[2]\S\step[2]\id\step\id\\
\id\step\id\step[2]\cu\step\dd\step\id\\
\id\step\d\step[2]\cu\step\dd\\
\d\step\d\step[2]\tu {\alpha}\\
\step\d\step\d\step[2]\O {\xi}\\
\step[2]\d\step\tu {\alpha}\\
\step[3]\tu {\alpha}\\
\step[4]\object{M}\\
\end{tangle}
\step=\step
\]
\[
\begin{tangle}
\step\object{H}\step[8]\object{M}\\
\cd\step[7]\id\\
\id\step[2]\nw2\step[6]\id\\
\id\step[3]\Q z\step[1]\nw3\step[4]\id\\
\id\step[2]\cd\step[2]\cd\step\id\\
\cu\step[2]\xx\step[2]\id\step\id\\
\step\id\step[3]\cu\step[2]\id\step\id\\
\step\d\step[3]\S\step[2]\dd\step\id\\
\step[2]\d\step[2]\cu\step\dd\\
\step[3]\d\step[2]\tu {\alpha}\\
\step[4]\d\step[2]\O {\xi}\\
\step[5]\tu {\alpha}\\
\step[6]\object{M}\\
\end{tangle}
\step=\step
\begin{tangle}
\step[2]\object{H}\step[5]\object{M}\\
\step\cd\step[4]\id\\
\cd\step\nw2\step[3]\id\\
\id\step[2]\id\step[2]\Q z\step \id\step[2]\id\\
\id\step[2]\id\step\cd\d\step\id\\
\id\step[2]\hx\step[2]\id\step\id\step\id\\
\cu\step\cu\step\id\step\id\\
\step\id\step[3]\S\step[2]\id\step\id\\
\step\d\step[2]\cu\step\id\\
\step[2]\d\step[2]\tu {\alpha}\\
\step[3]\d\step[2]\O {\xi}\\
\step[4]\tu {\alpha}\\
\step[5]\object{M}\\
\end{tangle}
\step=\step
\begin{tangle}
\step\object{H}\step[4]\object{M}\\
\cd\step[3]\id\\
\id\step[2]\nw1 \step[2]\id\\

\id\step[2]\Q z\step \id\step[2]\id\\
\cu\step\d\step\id\\
\cd\step[2]\id\step\id\\
\id\step[2]\S\step[2]\id\step\id\\
\id\step[2]\cu\step\id\\
\d\step[2]\tu {\alpha}\\
\step\d\step[2]\O {\xi}\\
\step[2]\tu {\alpha}\\
\step[3]\object{M}\\
\end{tangle}
\step=\step
\begin{tangle}
\step\Q z\step[3]\object{H}\step[1.5]\object{M}\\
\cd\step[2]\id\step\id\\
\id\step[2]\S\step[2]\id\step\id\\
\id\step[2]\cu\step\id\\
\d\step[2]\tu {\alpha}\\
\step\d\step[2]\O {\xi}\\
\step[2]\tu {\alpha}\\
\step[3]\object{M}\\
\end{tangle}
\step=\step
\begin{tangle}
\object{H}\step[2]\object{M}\\
\tu {\alpha}\\
\step\O {\mu}\\
\step\object{M}\\
\end{tangle}\ \ \ .
\]
Thus, $\mu$ is an $H$-module morphism. $\Box$

(5)  The smash  product  $H \# B$     is an ordinary semisimple Hopf
algebra by bosonisation
 \cite [Theorem 4.1]{Ma94a} and  by \cite [Theorem 7.4.2] {Mo93}.
Thus   $H \# B$ is finite-dimensional, which implies  that $H$ is
finite dimensional.
\begin{picture}(8,8)\put(0,0){\line(0,1){8}}\put(8,8){\line(0,-1){8}}\put(0,0){\line(1,0){8}}\put(8,8){\line(-1,0){8}}\end{picture}

 \begin {Lemma} \label {1.2.2}
 Let $B$  be an ordinary finite-dimensional Hopf algebra with left
 duality
 $B^*$. If $(H,\phi)$   is a left $B$-comodule algebra and $B$-comodule
 coalgebra, then $(H,\alpha )$ is a left ${B^*} ^{ \ cop}$-module algebra and
 $B^*$-module coalgebra, where
\[
 \begin{tangle}

 \object{B^*}\step [2] \object{H}\\
\tu \alpha \\

\step  \object{H}
\end{tangle} \ \ = \ \
\begin{tangle}

 \object{B^*}\step [3] \object{H}\\
\id \step [2]\td {\phi} \\

\ev \step [2]\id\\
 \step [4]\object{H}
\end{tangle} \ \ .
\]
              \end {Lemma}
{\bf Proof.}   It is trivial.
\begin{picture}(8,8)\put(0,0){\line(0,1){8}}\put(8,8){\line(0,-1){8}}\put(0,0){\line(1,0){8}}\put(8,8){\line(-1,0){8}}\end{picture}

 \begin {Proposition} \label {1.2.3} Let $char \ k=0$.

(i) If $H$ is not only a $kG$-module algebra but also a $kG$-module
coalgebra, then the action of $G$ on $\int _H^l$ is stable;

(ii) If $H$  is a Hopf algebra living in $({\cal M} ^{kG}, C^r)$
and $G$  is a commutative finite group with a bicharacter $r$, then
$\int _H^l$  is a homogeneous space of $G$-graded space $H$.

In particular,
 $\int _H^l $   is an (super) anyonic space for  every
  finite-dimensional (super) anyonic Hopf algebra $H$.

 \end {Proposition}

{\bf Proof.} (i)  For any $ g \in G, x \in H, t\in \int _H^l,$ we
see that
\begin {eqnarray*}
x(g\cdot t) &=& (g \cdot y)( g \cdot t)  \hbox {\ \ \ since } H
\hbox { \ \ is a } kG \hbox {-module algebra, where } y = g ^{-1}
\cdot x  \\
&=& g \cdot (yt) \\
&=& \epsilon (y)(g \cdot t) \\
&=& \epsilon (x)(g \cdot t)    \hbox { \ \ since }  H \hbox { is a }
kG
\hbox {-module coalgebra. }\\
\end {eqnarray*}
Therefore action of $G$ on $\int _H^l$ is stable.

(ii) Since $H$  is a left $kG$-comodule algebra and $kG$-comodule
coalgebra, we have that $H$ is a left $(kG)^*$-module algebra and
$(kG)^*$-module coalgebra by Lemma \ref {1.2.2}. However, $(kG)^* =
k \hat G$ (as Hopf algebra)  by \cite [Corollary 1.5.6] {Ma95b},
where $\hat G$  denotes the character group of $G$. Thus $H$ is a
left $k\hat G$-module algebra and a $k\hat G$-module coalgebra,
which implies that the action of $\hat G$ on $\int _H^l$ is stable
by part (1). That is, $\int _H^l$  is a sub-module of
$(kG)^*$-module $H$, which implies that $\int_H^l$  is a homogeneous
space of $G$-graded space $H$. $\Box$

 \begin {Example} \label {1.2.4}   Let $H$  denote  the anyonic line algebra,
 i.e.  $H = {\bf C} \{x \} / < x^n >$ , where  ${\bf C} \{x\}$
 is a free algebra
  over complex
 field {\bf C}  and $<x^n>$  is an ideal generalized by $x^n$ of
  ${\bf C} \{x\}$.
Set ${\bf C} \{x\} / <x^n> = {\bf C} [\xi]$  with $\xi ^n=0$.
 Its comultiplication and counit and antipode are
 $$ \Delta (\xi) = \xi \otimes 1 + 1 \otimes \xi \hbox { \ \ \ and \ \ \ }
 \epsilon (\xi ) =0 \hbox { \ \ \ and \ \ \ } S(\xi ) =- \xi$$
 respectively.
 It is straightforward to check that  $H$  is an anyonic Hopf algebra
 (see \cite {MR94}). Obviously, the dimension of  $H$  is  $n$.
 It is clear that $\int _H^l = {\bf C} \xi ^{n-1}$. Since $\epsilon
 (\int _{H}^l) =0$,
 $H$ is not semisimple.

 \end {Example}

 \begin {Example} \label {1.2.5}
 (see \cite [ P510 ] {Ma95b} )  Let $H= {\bf C} [x]$  denote
  the braided  line algebra. It is just the usual algebra ${\bf C} [x]$
  of polynomials in $x$ over complex field   ${\bf C}$, but we regard it
  as a q-statistical Hopf algebra with

 $$ \Delta (x) = x \otimes 1 + 1 \otimes x ,  \ \ \
 \epsilon (x ) =0 , \ \ \  S(x ) = -x , \ \ \ \mid x ^n \mid = n$$
 and $$C^r (x^n , x^m) =q^{nm} (x^m \otimes x^n).$$
If $y = \sum _0^n a_i x^i \in \int_H^l,$ then $a_i=0 $  for $i =1,2,
\cdots , n$ since $xy = \epsilon (x)y =0.$ Thus  $\int _H^l =0$. It
follows from the proof of Theorem \ref {1.2.1} (5) that
 $H$ is not semisimple.

 \end {Example}
 By the way, we  give some properties of group graded (co)algebras.

An algebra $A$ is called a $G$-graded algebra if $A= \sum _{g\in G}
\oplus A_g$, $A_g A_h \subseteq A_{gh}$  and $1_A \in A_e$, where
$e$ is the unity element of $G$. A coalgebra $D$ is called a
$G$-graded coalgebra if $D= \sum _{g\in G} \oplus D_g$, $\Delta
(D_g) \subseteq \sum _{x, y \in G, xy =g} D_{x} \otimes D_y$  and
$\epsilon (D_h) = 0 $ for any $h, g \in G$ and $h\not= e.$

\begin {Proposition} \label {1.2.6}

(1) $A$ is a right $kG$-comodule algebra iff $A$ is a $G$-graded
algebra;

(2) $D$ is a right $kG$-comodule coalgebra iff $D$ is a $G$-graded
coalgebra.
\end {Proposition}

{\bf Proof.} (1)  It was obtained in \cite [Proposition 1.3]
{CM84b}.

 (2)  If $(D,\psi)$ is a $kG$-comodule coalgebra, then $D=\sum _{g\in G}
 \oplus D_g$  is a $G$-graded space by \cite [Lemma 4.8] {BM85}.
 For any homogeneous $0\not=d \in D_g,$ let

$$\Delta (d)= \sum _{i = 1, 2, \cdots m} (d_{x_i} \otimes d_{y_i})$$
where $0 \not=d_{x_i}\in D_{x_i}, 0\not= d_{y_i} \in D_{y_i}, x_i ,
y_i \in G$
 for $i= 1, 2, \cdots , m.$
 Since $\Delta $ is a $kG$-comodule homomorphism we have that

$$ \sum _{i = 1, 2, \cdots m} (d_{x_i} \otimes d_{y_i} \otimes g) =
 \sum _{i = 1, 2, \cdots m} (d_{x_i} \otimes d_{y_i} \otimes x_iy_i).$$

 Thus $$x_{i}y_i = g$$ for $i= 1, 2, \cdots ,m,$
 which implies that $$\Delta (D_g) \subseteq
 \sum _{x, y \in G, xy =g} D_x \otimes D_y.$$

Since $\epsilon $ is a $kG$-comodule homomorphism, we have that
$$\epsilon (d)g= \epsilon (d)e$$
for any $g\in G, d\in D_g,$ which implies that $\epsilon (D_g)=0$
for any $e\not=g \in G.$ Consequently, $D$ is a $G$-graded
coalgebra.

Conversely, if $D= \sum _{g\in G}\oplus D_g$ is a $G$-graded
coalgebra, then $(D,\psi)$ is a right $kG$-comodule by \cite [Lemma
4.8] {BM85}. It is sufficient to show that $\Delta $ and $\epsilon $
are $kG$-comodule maps. For any $d\in D_g,$ we assume that
$$\Delta (d)= \sum _{x, y \in G, xy=g} d_x \otimes d_y$$
where $d_x \in D_x, d_y \in D_y.$ We see that
\begin {eqnarray*}
(\Delta \otimes id_D) \psi (d) &=& (\Delta \otimes id _D)(d \otimes g) \\
&=& \sum _{x, y \in G, xy =g} d_x \otimes d_y \otimes g \\
&=& (id _{D \otimes D} \otimes m_{kG}) (id _D \otimes \tau \otimes
id _ {kG})(\psi \otimes \psi) \Delta (d)
\end {eqnarray*}
 where  $\tau $ is an ordinary twist map.

Thus $\Delta $ is a $kG$-comodule map. Since $$(\epsilon \otimes id
_{kG} )\psi (d) = \epsilon (d)g = \epsilon (d)e =\epsilon
(d)1_{kG}$$  for any $d \in D_g,$ we have that $\epsilon $  is a
$kG$-comodule map.
 \begin{picture}(8,8)\put(0,0){\line(0,1){8}}\put(8,8){\line(0,-1){8}}\put(0,0){\line(1,0){8}}\put(8,8){\line(-1,0){8}}\end{picture}

\chapter {  The Double Bicrossproducts in Braided Tensor Categories }\label {c3}

In this chapter, we construct   the double bicrossproduct
 $D =  A_\alpha ^\phi \bowtie _\beta ^\psi H  $ of two bialgebras $A$  and $H$
  in a braided tensor category and give   the necessary  and
 sufficient conditions for  $D$ to be a  bialgebra. We study
 the universal property of
 double bicrossproduct.

Throughout this chapter, we assume  that   $H$ and $A$ are
 two bialgebras in  braided tensor categories and
\begin {eqnarray*}
\alpha : H \otimes A \rightarrow &A& , \hbox { \ \ \ \ }
\beta : H \otimes A \rightarrow H,    \\
\phi :  A \rightarrow H \otimes  &A&   , \hbox { \ \ \ \ } \psi : H
\rightarrow H \otimes A
\end {eqnarray*}
                        such that $(A, \alpha )$ is a left $H$-module coalgebra,
 $(H, \beta )$ is a right $A$-module coalgebra,
 $(A, \phi )$ is a left $H$-comodule algebra, and
 $(H, \psi )$ is a right $A$-comodule algebra.

We define   the multiplication $m_D$ , unit $\eta _D$,
comultiplication  $\Delta _D$ and counit $\epsilon _D$  in $A
\otimes H$  as follows:

\[
\begin{tangle}
\Delta_{D}
\end{tangle}
 \ =\enspace
\begin{tangles}{clr}
\step[1]\object{A}\step[6]\object{H}\\
\cd\step[4]\cd\\
\id\step[1]\td \phi\step[2]\td \psi\step[1]\id\\
\id\step[1]\id\step[2]\x\step[2]\id\step[1]\id\\
\id\step[1]\cu \step[2]\cu \step[1]\id\\
\object{A}\step[2]\object{H}\step[4]\object{A}\step[2]\object{H}
\end{tangles}
\step\,\step\\
\begin{tangle}
m_{D}
\end{tangle}
\\ \ =\enspace
\begin{tangles}{clr}
\object{A}\step[2]\object{H}\step[4]\object{A}\step[2]\object{H}\\
\id\step[1]\cd\step[2]\cd\step[1]\id\\
\id\step[1]\id\step[2]\x\step[2]\id\step[1]\id\\
\id\step[1]\tu \alpha\step[2]\tu \beta\step[1]\id\\
\cu \step[4]\cu \\
\step[1]\object{A}\step[6]\object{H}
\end{tangles}
\]
and $\epsilon _D = \epsilon _A \otimes \epsilon _H$ , $ \eta _D =
\eta _A \otimes \eta _H.$ We denote   $(A \otimes H, m_D, \eta _D,
\Delta _D, \epsilon _D) $
  by                      $$ A ^{\phi} _{\alpha}  {\bowtie }
^{\psi}_{\beta} H ,$$ which is
 called the
double bicrossproduct of $A$ and $H$. We  denote it by  $A \stackrel
{b} {\bowtie }  H$ in short.

If one of $\alpha , \beta , \phi , \psi $  is trivial, we omit the
trivial one
 in                      the notation                  $ A ^{\phi} _{\alpha}  {\bowtie }
^{\psi}_{\beta} H .$  For example, when $\phi $  and $\psi$ are
trivial,  we   denote $ A_\alpha ^\phi \bowtie _\beta ^\psi H $  by
 $A _\alpha {\bowtie } _\beta H$. When
 $\alpha $ and $\beta $ are trivial, we denote
  $ A_\alpha ^\phi \bowtie _\beta ^\psi H $  by
$ A ^\phi \bowtie  ^\psi H  $. We call $ A_\alpha  \bowtie _\beta H
$  a double cross product and denote it by  $A  {\bowtie }  H$ in
short. We call
  $ A ^\phi \bowtie  ^\psi H $  a double cross coproduct
and denote it by  $A \stackrel {c} {\bowtie }  H$ in short.
  We call
$ A_\alpha ^\phi \bowtie H $
 ( $ A \bowtie _\beta ^\psi H $) and
       $ A  _{\alpha}  {\bowtie }
^{\psi} H $ (       $ A ^{\phi}   {\bowtie } _{\beta} H $ ) a
biproduct and  a bicrossproduct respectively. Note that S. Majid
uses different notations.

We give some notation as follows.

\[
(M1):\step[2]
\begin{tangle}
\object{H}\\
\id\step[2]\Q {\eta_{A}}\\
\tu \alpha\\
\step[1]\object{A}
\end{tangle}
\;=\enspace
\begin{tangles}{clr}
\object{H}\\
\QQ {\epsilon_{H}}\\
\Q {\eta_{A}}\\
\object{A}
\end{tangles}
\;\step[1],\step[1]\enspace
\begin{tangle}
\step[2]\object{A}\\
\Q { \eta_{H}}\step[2]\id\\
\tu \beta\\
\step[1]\object{H}
\end{tangle}
\;=\enspace
\begin{tangles}{clr}
\object{A}\\
\QQ {\epsilon _{A}}\\
\Q {\eta_{H}}\\
\object{H}
\end{tangles} \ \ \  ;
\]

\[
(M2):\step[2]
\begin{tangle}
\object{H}\step[1]\object{A}\step[2]\object{A}\\
\id\step[1]\cu \\
\tu \alpha\\
\step[1]\object{A}
\end{tangle}
\;=\enspace
\begin{tangles}{clr}
\step[1]\object{H}\step[4]\object{A}\step[2]\object{A}\\
\cd\step[2]\cd\step[1]\id\\
\id\step[2]\x\step[2]\id\step[1]\id\\
\tu \alpha\step[2]\tu \beta\step[1]\id\\
\step[1]\nw3\step[3]\tu \alpha\\
\step[3]\cu \\
\step[4]\object{A}
\end{tangles} \ \ ;
\]

\[
(M3):\step[2]
\begin{tangle}
\object{H}\step[2]\object{H}\step[2]\object{A}\\
\cu \step[1]\id\\
\step[1]\tu \beta\\
\step[2]\object{H}
\end{tangle}
\;=\enspace
\begin{tangles}{clr}
\step[1]\object{H}\step[2]\object{H}\step[4]\object{A}\\
\step[1]\id\step[1]\cd\step[2]\cd\\
\step[1]\id\step[1]\id\step[2]\x\step[2]\id\\
\step[1]\id\step[1]\tu \alpha\step[2]\tu \beta\\
\tu \beta\step[3]\ne3\\
\cu \step[2]\\
\object{H}\step[2]
\end{tangles} \ \ ;
\]

\[
(M4):\step[2]
\begin{tangle}
\step[1]\object{H}\step[4]\object{A}\\
\cd\step[2]\cd\\
\id\step[2]\x\step[2]\id\\
\tu \beta\step[2]\tu \alpha\\
\step[1]\object{H}\step[4]\object{A}
\end{tangle}
\;=\enspace
\begin{tangles}{clr}
\step[1]\object{H}\step[4]\object{A}\\
\cd\step[2]\cd\\
\id\step[2]\x\step[2]\id\\
\tu \alpha\step[2]\tu \beta\\
\step[1]\nw1\step[2]\ne1\step[1]\\
\step[2]\x\step[2]\\
\step[2]\object{H}\step[2]\object{A}\step[2]
\end{tangles} \ \ ;
\]

\[
(CM1):\step[2]
\begin{tangle}
\step[1]\object{A}\\
\td \phi\\
\id\step[2]\QQ {\epsilon _{A}}\\
\object{H}\\
\end{tangle}
\step=\step
\begin{tangle}
\object{A}\\
\QQ {\epsilon _{A}}\\
\Q {\eta_{H}}\\
\object{H}\\
\end{tangle}
\step,\step
\begin{tangle}
\step[1]\object{H}\\
\td \psi\\
\QQ {\epsilon _{H}}\step[2]\id\\
\step[2]\object{A}
\end{tangle}
\step=\step
\begin{tangle}
\object{H}\\
\QQ {\epsilon _{H}}\\
\Q { \eta_{A}}\\
\object{A}
\end{tangle} \ \ \ ;
\]

\[
(CM2):\step[2]
\begin{tangle}
\step[1]\object{A}\\
\td \phi\\
\id\step[1]\cd\\
\object{H}\step[1]\object{A}\step[2]\object{A}
\end{tangle}
\;=\enspace
\begin{tangles}{clr}
\step[3]\object{A}\\
\step[3]\cd\\
\step[3]\ne3\step[1]\td \phi\\
\td \phi\step[2]\td \psi\step[1]\id\\
\id\step[2]\x\step[2]\id\step[1]\id\\
\cu \step[2]\cu \step[1]\id\\
\step[1]\object{H}\step[4]\object{A}\step[2]\object{A}
\end{tangles} \ \ \ ;
\]

\[
(CM3):\step[2]
\begin{tangle}
\step[2]\object{H}\\
\step[1]\td \psi\\
\cd\step[1]\id\\
\id\step[2]\id\step[1]\id\\
\object{H}\step[2]\object{H}\step[1]\object{A}
\end{tangle}
\step=\step
\begin{tangle}
\step[2]\object{H}\\
\step[1]\cd\\
\td \psi\step[1]\nw3\\
\id\step[1]\td \phi\step[2]\td \psi\\
\id\step[1]\id\step[2]\x\step[2]\id\\
\id\step[1]\cu \step[2]\cu \\
\object{H}\step[2]\object{H}\step[4]\object{A}
\end{tangle}\ \ \ ;
\]

\[
(CM4):\step[2]
\begin{tangle}
\step[1]\object{H}\step[4]\object{A}\\
\td \psi\step[2]\td \phi\\
\id\step[2]\x\step[2]\id\\
\cu \step[2]\cu \\
\step[1]\object{H}\step[4]\object{A}
\end{tangle}
\step=\step
\begin{tangle}
\step[1]\object{H}\step[2]\object{A}\\
\step[1]\x\\
\step[1]\id\step[2]\nw2\\
\td \phi\step[2]\td \psi\\
\id\step[2]\x\step[2]\id\\
\cu \step[2]\cu \\
\step[1]\object{H}\step[4]\object{A}
\end{tangle}\ \ \ ;
\]

\[
(B1):\step[2]
\begin{tangle}
\object{A}\step[2]\object{A}\\
\cu \\
\cd\\
\object{A}\step[2]\object{A}\\
\end{tangle}
\step=\step
\begin{tangle}
\step[1]\object{A}\step[5]\object{A}\\
\cd\step[4]\id\\
\id\step[1]\td \phi\step[2]\cd\\
\id\step[1]\id\step[2]\x\step[2]\id\\
\id\step[1]\tu \alpha \step[2]\cu \\
\cu \step[4]\id\\
\step[1]\object{A}\step[5]\object{A}
\end{tangle} \ \ \ ;
\]

\[
(B2):\step[2]
\begin{tangle}
\object{H}\step[2]\object{A}\\
\tu \alpha\\
\cd\\
\object{A}\step[2]\object{A}
\end{tangle}
\step=\step
\begin{tangle}
\step[2]\object{H}\step[3]\object{H}\\
\step[1]\cd\step[1]\cd\\
\td \psi\step[1]\hx\step[2]\id\\
\id\step[2]\hx\step[1]\tu \alpha\\
\tu \alpha\step[1]\cu \\
\step[1]\object{A}\step[3]\object{A}
\end{tangle}\ \ \ ;
\]

\[
(B3):\step[2]
\begin{tangle}
\object{H}\step[2]\object{H}\\
\cu \\
\cd\\
\object{H}\step[2]\object{H}
\end{tangle}
\step=\step
\begin{tangle}
\step[1]\object{H}\step[5]\object{H}\\
\cd\step[3]\cd\\
\id\step[2]\id\step[2]\td \psi\step[1]\id\\
\id\step[2]\x\step[2]\id\step[1]\id\\
\id\step[2]\id\step[2]\tu \beta\step[1]\id\\
\cu \step[3]\cu \\
\step[1]\object{H}\step[5]\object{H}
\end{tangle}\ \ \ ;
\]

\[
(B4):\step[2]
\begin{tangle}
\object{H}\step[2]\object{A}\\
\tu \beta\\
\cd\\
\object{H}\step[2]\object{H}\\
\end{tangle}
\step=\step
\begin{tangle}
\step[1]\object{H}\step[4]\object{A}\\

\cd\step[2]\cd\\
\id\step[2]\x\step[1]\td \phi\\
\tu \beta\step[2]\hx\step[2]\id\\

\step[1]\id\step[2]\ne1\step[1]\id\step[2]\id\\
\step[1]\cu \step[2]\tu \beta\\
\step[2]\object{H}\step[4]\object{H}
\end{tangle}\ \ \ ;
\]

\[
(B5):\step[2]
\begin{tangle}
\step[1]\object{A}\step[5]\object{H}\step[3]\object{A}\step[5]\object{H}\\

\td \phi\step[3]\cd\step[1]\cd\step[3]\td \psi\\
\id\step[2]\id\step[2]\td \psi\step[1]\hx\step[1]\td \phi\step[2]\id\step[2]\id\\

\id\step[2]\x\step[2]\hx\step[1]\hx\step[2]\x\step[2]\id\\

\cu \step[2]\x\step[1]\hx\step[1]\x\step[2]\cu \\

\step[1]\id\step[2]\ne1\step[2]\hx\step[1]\hx\step[2]\nw1\step[2]\id\\

\step[1]\tu \beta\step[2]\ne1\step[1]\hx\step[1]\nw1\step[2]\tu \alpha\\
\step[2]\nw2\step[2]\cu \step[1]\cu \step[2]\ne2\\
\step[4]\cu \step[3]\cu \\
\step[5]\object{H}\step[5]\object{A}
\end{tangle}
\step=\step
\begin{tangle}
\step[1]\object{A}\step[4]\object{H}\step[4]\object{A}\step[4]\object{H}\\
\td \phi\step[2]\cd\step[2]\cd\step[2]\td \psi\\
\id\step[2]\id\step[2]\id\step[2]\x\step[2]\id\step[2]\id\step[2]\id\\
\id\step[2]\id\step[2]\tu \alpha\step[2]\tu \beta\step[2]\id\step[2]\id\\
\id\step[2]\id\step[2]\td \phi\step[2]\td \psi\step[2]\id\step[2]\id\\
\id\step[2]\x\step[2]\x\step[2]\x\step[2]\id\\
\cu \step[2]\x\step[2]\x\step[2]\cu \\
\step[1]\id\step[3]\id\step[2]\x\step[2]\id\step[3]\id\\
\step[1]\nw2\step[2]\cu \step[2]\cu \step[2]\ne2\\
\step[3]\cu \step[4]\cu \\
\step[4]\object{H}\step[6]\object{A}
\end{tangle}\ \ \ ;
\]

\[
(B6):\step[2]
\begin{tangle}
\step[1]\object{H}\step[3]\object{A}\\
\cd\step[2]\id\\
\id\step[2]\x\\
\tu \alpha\step[2]\id\\
\td \phi \step[2]\id\\
\id\step[2]\x\\
\cu \step[2]\id\\
\step[1]\object{H}\step[3]\object{A}
\end{tangle}
\step=\step
\begin{tangle}
\step[1]\object{H}\step[4]\object{A}\\
\cd\step[2]\td \phi\\
\id\step[2]\x\step[2]\id\\
\cu \step[2]\tu \alpha\\
\step[1]\object{H}\step[4]\object{A}
\end{tangle}\ \ \ ;
\]

\[
(B6'):\step[2]
\begin{tangle}
\object{H}\step[3]\object{A}\\
\id\step[2]\cd\\
\x\step[2]\id\\
\id\step[2]\tu \beta\\
\id\step[2]\td \psi\\
\x\step[2]\id\\
\id\step[2]\cu \\
\object{H}\step[3]\object{A}
\end{tangle}
\step=\step
\begin{tangle}
\step[1]\object{H}\step[4]\object{A}\\
\td \psi\step[2]\cd\\
\id\step[2]\x\step[2]\id\\
\tu \beta\step[2]\cu \\
\step[1]\object{H}\step[4]\object{A}
\end{tangle}\ \ \ ;
\]

\[
(CB2):\step[2]
\begin{tangle}
\object{A}\step[2]\object{A}\\
\cu \\
\td \phi\\
\object{H}\step[2]\object{A}
\end{tangle}
\step=\step
\begin{tangle}
\step[1]\object{A}\step[4]\object{A}\\
\td \phi\step[2]\cd\\
\id\step[2]\x\step[1]\td \phi\\
\tu \beta\step[2]\X\step[2]\id\\
\step[1]\id\step[2]\ne1\step[1]\cu \\
\step[1]\cu \step[3]\id\\
\step[2]\object{H}\step[4]\object{A}
\end{tangle}\ \ \ ;
\]

\[
(CB4):\step[2]
\begin{tangle}
\object{H}\step[2]\object{H}\\
\cu \\
\td \psi\\
\object{H}\step[2]\object{A}
\end{tangle}
\step=\step
\begin{tangle}
\step[2]\object{H}\step[4]\object{H}\\
\step[1]\cd\step[2]\td \psi\\
\td \psi\step[1]\x\step[2]\id\\

\id\step[2]\hx\step[2]\tu \alpha\\
\cu \step[1]\nw1\step[2]\id\\
\step[1]\id\step[3]\cu \\
\step[1]\object{H}\step[4]\object{A}
\end{tangle}\ \ \ ;
\]

\[
(B7):\step[2]
\begin{tangle}
\step[1]\object{H}\step[3]\object{A}\\
\cd\step[2]\id\\
\id\step[2]\x\\
\tu \alpha\step[1]\td \psi\\
\step[1]\x\step[2]\id\\
\step[1]\id\step[2]\cu \\
\step[1]\object{H}\step[3]\object{A}
\end{tangle}
\step=\step
\begin{tangle}
\step[2]\object{H}\step[3]\object{A}\\
\step[1]\cd\step[2]\id\\
\td \psi\step[1]\tu \alpha\\
\id\step[2]\cu \\
\object{H}\step[3]\object{A}
\end{tangle}\ \ \ ;
\]

\[
(B7'):\step[2]
\begin{tangle}
\step[1]\object{H}\step[3]\object{A}\\
\step[1]\id\step[2]\cd\\
\step[1]\x\step[2]\id\\
\td \phi\step[1]\tu \beta\\
\id\step[2]\x\\
\cu \step[2]\id\\
\step[1]\object{H}\step[3]\object{A}
\end{tangle}
\step=\step
\begin{tangle}
\object{H}\step[3]\object{A}\\
\id\step[2]\cd\\
\tu \beta\step[1]\td \phi\\
\step[1]\cu \step[2]\id\\
\step[2]\object{H}\step[3]\object{A}
\end{tangle}\ \ \ .
\]

Let \[
\begin{tangle}
\step \object{H}\\
\td {i_H}\\

\object{A}\step[2]\object{H}
\end{tangle}
\ \ = \ \
\begin{tangle}
\step [2]\object{H}\\
\Q {\eta _A} \step [2] \id \\

\object{A}\step[2]\object{H}
\end{tangle}
\ \ ; \ \
\begin{tangle}
\step\object{A}\\
\td {i_A}\\

\object{A}\step[2]\object{H}
\end{tangle}
\ \ = \ \
\begin{tangle}
\object{H}\\
 \id  \step [2]\Q {\eta _H}\\

\object{A}\step[2]\object{H}
\end{tangle} \ \ \ ; \ \ \
\begin{tangle}
\object{A}\step[2]\object{H}\\
\tu {\pi_H}\\
\step\object{H}

\end{tangle}
\ \ = \ \
\begin{tangle}
\object{A}\step[2]\object{H}\\
\QQ {\epsilon  _A} \step [2] \id \\

\step [2]\object{H}
\end{tangle} \ \ ; \ \
\begin{tangle}
\object{A}\step[2]\object{H}\\
\tu {\pi_A}\\
\step\object{A}

\end{tangle}
\ \ = \ \
\begin{tangle}
\object{A}\step[2]\object{H}\\
\id  \step [2]\QQ {\epsilon  _A} \\

\object{A}
\end{tangle} \ \ . \ \
\]

 $\pi_A$ and $\pi _H$
 are called trivial actions; $i_A$ and $i_H$ are called trivial coactions.

  \section {Double bicrossproduct }\label {s3}

In this section, we give the sufficient and necessary conditions for
double bicrossproduct $A _\alpha ^\phi \bowtie _\beta ^\psi H $ to
be a bialgebra.

\begin {Lemma}\label {2.1.1}

(i)  If $A$  and $H$ are bialgebras and $H$ acts on $A$  trivially,
then $A$ is an $H$-module algebra and an $H$-module coalgebra.

(ii)   If $A$  and $H$ are bialgebras and $H$ coacts on $A$
trivially, then $A$ is an $H$-comodule algebra and an $H$-comodule
coalgebra.

(iii)  $C_{I,I} =id_I, $  \ \ \ \ $C_{I, I \otimes I} =
a^{-1}_{I,I,I},$ \ \ \ \ $C_{I \otimes I , I} = a _{I, I, I}$ ;

(iv)  If we define:
$$\Delta _I = r_I^{-1}, \ \ m_I = r_I,  \ \
\epsilon _I = id _I = \eta _I =S \ ,$$ then $I$ is a Hopf algebra.

(v)  If $H$ is a bialgebra, then $I$ is an $H$-module algebra,
$H$-module coalgebra, $H$-comodule algebra, $H$-comodule coalgebra.
Meantime, $H$ is
 an $I$-module algebra, $I$-module
coalgebra, $I$-comodule algebra, $I$-comodule coalgebra.

(vi) If $B$ is a coalgebra and $A$ is an algebra, then
 $Hom_ {{\cal C }} (B, A)$  is a monoid, i.e.
 $Hom_ {{\cal C }} (B, A)$  is a semigroup with unit.

(vii)  If $(B1)$ and $(B2)$ hold, then

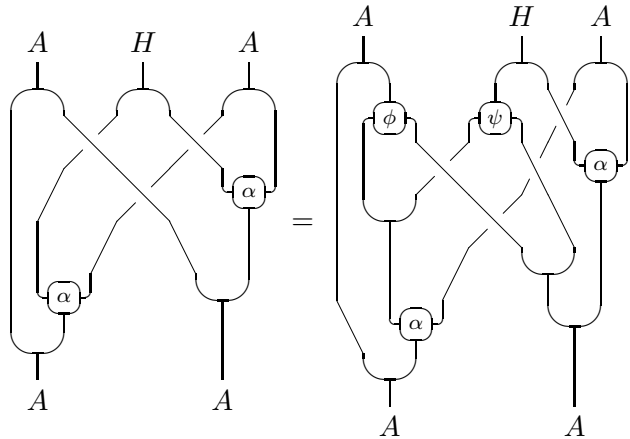
\[
\begin{tangle}
\step[1]\object{A}\step[4]\object{H}\step[4]\object{A}\\
\cd\step[2]\cd\step[2]\cd\\
\id\step[2]\x\step[2]\x\step[2]\id\\
\id\step[1]\ne1\step[2]\x\step[2]\tu \alpha\\
\id\step[1]\id\step[2]\ne1\step[2]\nw1\step[2]\id\\
\id\step[1]\tu \alpha\step[4]\cu\\
\cu \step[6]\id\\
\step[1]\object{A}\step[7]\object{A}
\end{tangle}
\;=\enspace
\begin{tangles}{lcr}
\step[1]\object{A}\step[6]\object{H}\step[3]\object{A}\\
\cd\step[4]\cd\step[1]\cd\\
\id\step[1]\td \phi\step[2]\td \psi\step[1] \hx \step[2]\id\\
\id\step[1]\id\step[2]\x\step[2]\hx\step[1]\tu \alpha\\
\id\step[1]\cu \step[2]\x\step[1]\nw1\step[1]\id\\
\id\step[2]\id\step[2]\ne1\step[2]\cu \step[1]\id\\
\nw1\step[1]\tu \alpha\step[4]\cu \\
\step[1]\cu \step[6]\id\\
\step[2]\object{A}\step[7]\object{A}
\end{tangles}\ \ \   .
\]

(viii)  If $(B3)$ and $(B4)$ hold, then

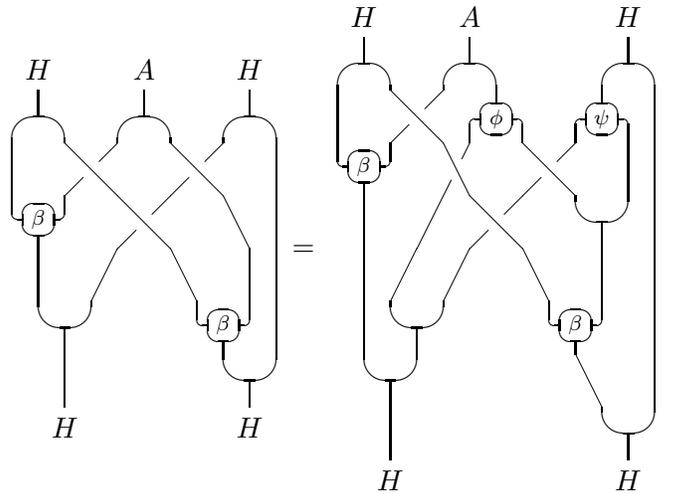
\[
\begin{tangle}
\step[1]\object{H}\step[4]\object{A}\step[4]\object{H}\\
\cd\step[2]\cd\step[2]\cd\\
\id\step[2]\x\step[2]\x\step[2]\id\\
\tu \beta\step[2]\x\step[2]\nw1\step[1]\id\\
\step[1]\id\step[2]\ne1\step[2]\nw1\step[2]\id\step[1]\id\\
\step[1]\cu \step[4]\tu \beta\step[1]\id\\
\step[2]\id\step[6]\cu \\
\step[2]\object{H}\step[7]\object{H}
\end{tangle}
\;=\enspace
\begin{tangles}{lcr}
\step[1]\object{H}\step[4]\object{A}\step[6]\object{H}\\
\cd\step[2]\cd\step[4]\cd\\
\id\step[2]\x\step[1]\td \phi\step[2]\td \psi\step[1]\id\\
\tu \beta\step[2]\hx\step[2]\x\step[2]\id\step[1]\id\\
\step[1]\id\step[2]\ne1\step[1]\x\step[2]\cu\step[1]\id\\
\step[1]\id\step[1]\ne1\step[1]\ne1\step[2]\nw1\step[2]\id\step[2]\id\\
\step[1]\id\step[1]\cu \step[4]\tu \beta\step[2]\id\\
\step[1]\cu \step[6]\nw1\step[2]\id\\
\step[2]\id\step[8]\cu \\
\step[2]\object{H}\step[9]\object{H}
\end{tangles}\ \ \   .
\]

 \end {Lemma}
 {\bf Proof.} The proof of parts (i)--(vi)  is trivial. The proof of
 part (vii) and part (viii)  can be obtained by turning the proof \cite [Lemma 1 and 2 ]
 {ZWJ98}  into braided diagrams.
 \begin{picture}(8,8)\put(0,0){\line(0,1){8}}\put(8,8){\line(0,-1){8}}\put(0,0){\line(1,0){8}}\put(8,8){\line(-1,0){8}}\end{picture}

\begin {Lemma} \label {2.1.2}      $D = A {\bowtie } H $  is an algebra iff
$(M1)$--$(M3)$ hold.
\end {Lemma}
{\bf Proof.}  If (M1)- (M3) hold, we shall show that $D$ is an
algebra.
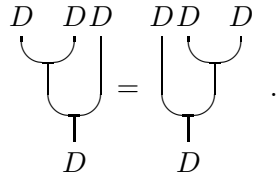
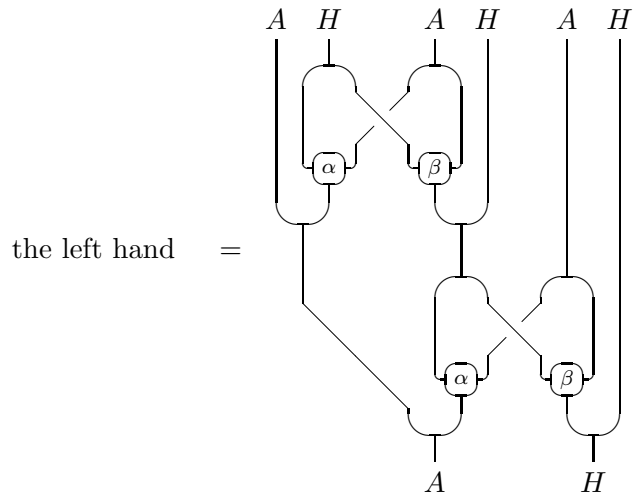
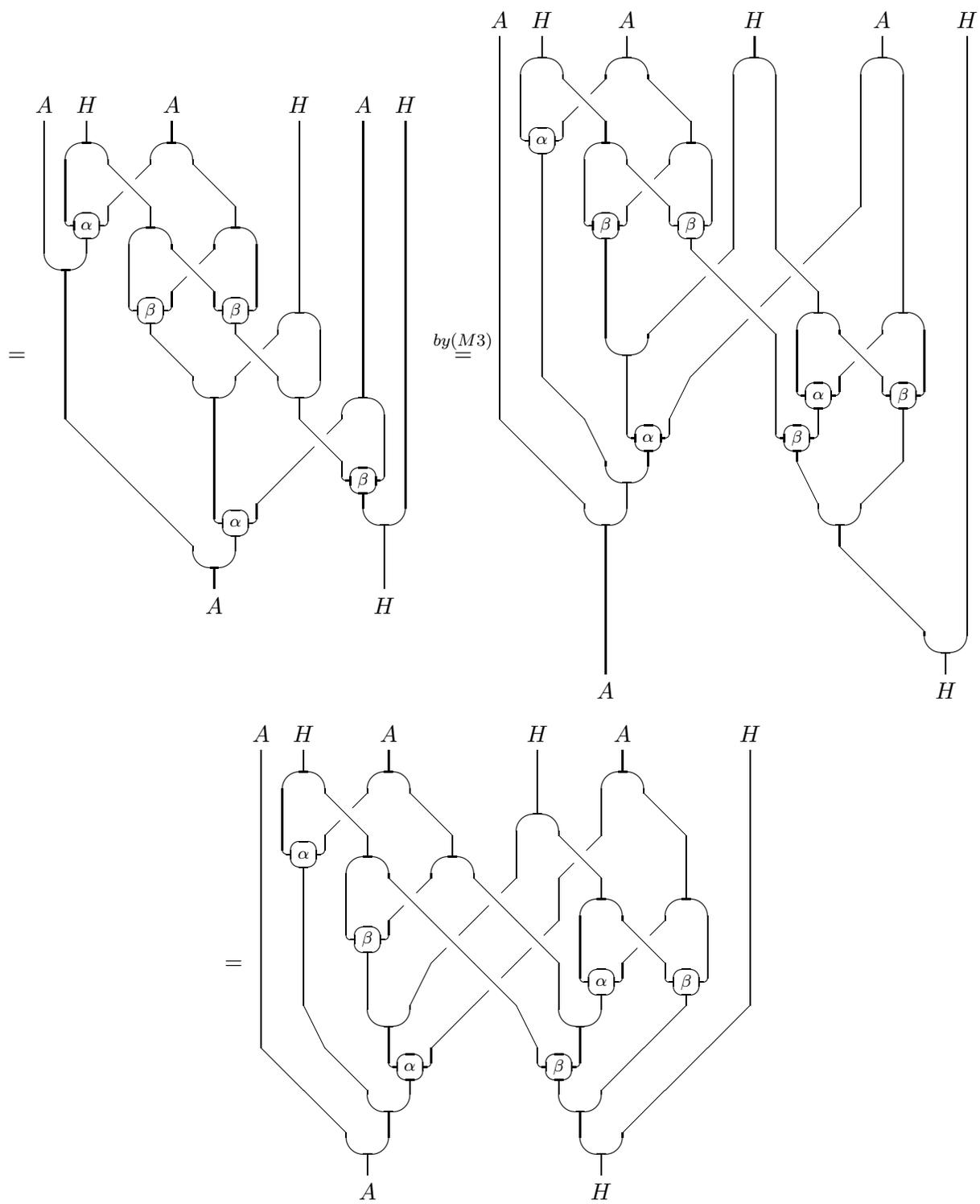
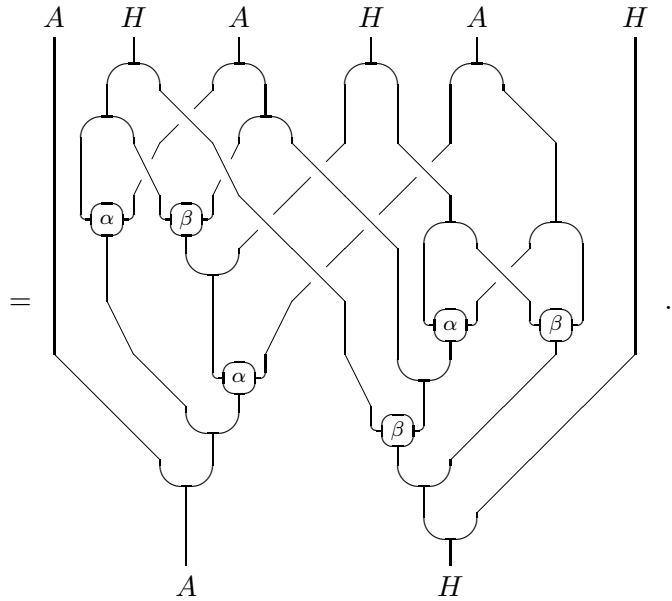
\[
\begin{tangle}
\object{D}\step[2]\object{D}\step[1]\object{D}\\
\cu \step[1]\id\\
\step[1]\cu \\
\step[2]\object{D}
\end{tangle}
\;=\enspace
\begin{tangles}{lcr}
\object{D}\step[1]\object{D}\step[2]\object{D}\\
\id\step[1]\cu \\
\cu \\
\step[1]\object{D}
\end{tangles}\ \ \   .
\]

See that
\[
\hbox {the left hand } \step=\step
\begin{tangle}
\object{A}\step[2]\object{H}\step[4]\object{A}\step[2]\object{H}\step[3]\object{A}\step[2]\object{H}\\
\id\step[1]\cd\step[2]\cd\step[1]\id\step[3]\id\step[2]\id\\
\id\step[1]\id\step[2]\x\step[2]\id\step[1]\id\step[3]\id\step[2]\id\\
\id\step[1]\tu \alpha\step[2]\tu \beta\step[1]\id\step[3]\id\step[2]\id\\
\cu \step[4]\cu \step[3]\id\step[2]\id\\
\step[1]\id\step[5]\cd\step[2]\cd\step[1]\id\\
\step[1]\nw2\step[4]\id\step[2]\x\step[2]\id\step[1]\id\\
\step[3]\nw2\step[2]\tu \alpha\step[2]\tu \beta\step[1]\id\\
\step[5]\cu \step[4]\cu \\
\step[6]\object{A}\step[6]\object{H}
\end{tangle}
\]

\[
\;=\enspace
\begin{tangle}
\object{A}\step[2]\object{H}\step[4]\object{A}\step[6]\object{H}\step[3]\object{A}\step[2]\object{H}\\
\id\step[1]\cd\step[2]\cd\step[5]\id\step[3]\id\step[2]\id\\
\id\step[1]\id\step[2]\x\step[2]\nw2\step[4]\id\step[3]\id\step[2]\id\\
\id\step[1]\tu \alpha\step[1]\cd\step[2]\cd\step[2]\id\step[3]\id\step[2]\id\\
\cu \step[2]\id\step[2]\x\step[2]\id\step[2]\id\step[3]\id\step[2]\id\\
\step[1]\id\step[3]\tu \beta\step[2]\tu \beta\step[1]\cd\step[2]\id\step[2]\id\\
\step[1]\id\step[4]\nw2\step[3]\x\step[2]\id\step[2]\id\step[2]\id\\
\step[1]\id\step[6]\cu \step[2]\cu \step[1]\cd\step[1]\id\\
\step[1]\nw2\step[6]\id\step[4]\x\step[2]\id\step[1]\id\\
\step[3]\nw2\step[4]\id\step[3]\ne2\step[2]\tu \beta\step[1]\id\\
\step[5]\nw2\step[2]\tu \alpha\step[5]\cu \\
\step[7]\cu \step[7]\id\\
\step[8]\object{A}\step[8]\object{H}
\end{tangle}
\step\stackrel{by (M3)}{=}
\begin{tangle}
\object{A}\step[2]\object{H}\step[4]\object{A}\step[6]\object{H}\step[6]\object{A}\step[4]\object{H}\\
\id\step[1]\cd\step[2]\cd\step[4]\cd\step[4]\cd\step[3]\id\\
\id\step[1]\id\step[2]\x\step[2]\nw2\step[3]\id\step[2]\id\step[4]\id\step[2]\id\step[3]\id\\
\id\step[1]\tu \alpha\step[1]\cd\step[2]\cd\step[1]\id\step[2]\id\step[4]\id\step[2]\id\step[3]\id\\
\id\step[2]\id\step[2]\id\step[2]\x\step[2]\id\step[1]\id\step[2]\id\step[4]\id\step[2]\id\step[3]\id\\
\id\step[2]\id\step[2]\tu \beta\step[2]\tu \beta\step[1]\id\step[2]\id\step[3]\ne2\step[2]\id\step[3]\id\\
\id\step[2]\id\step[3]\id\step[4]\x\step[2]\x\step[4]\id\step[3]\id\\
\id\step[2]\id\step[3]\id\step[3]\ne2\step[2]\x\step[1]\cd\step[2]\cd\step[2]\id\\
\id\step[2]\id\step[3]\cu \step[3]\ne2\step[2]\id\step[1]\id\step[2]\x\step[2]\id\step[2]\id\\
\id\step[2]\nw2\step[3]\id\step[2]\ne1\step[4]\id\step[1]\tu \alpha\step[2]\tu \beta\step[2]\id\\
\nw2\step[3]\nw1\step[1]\tu \alpha\step[5]\tu \beta\step[4]\id\step[3]\id\\
\step[2]\nw2\step[2]\cu \step[7]\nw1\step[3]\ne2\step[3]\id\\
\step[4]\cu \step[9]\cu \step[5]\id\\
\step[5]\id\step[11]\nw2\step[5]\id\\
\step[5]\id\step[13]\nw2\step[3]\id\\
\step[5]\id\step[15]\cu \\
\step[5]\object{A}\step[16]\object{H}
\end{tangle}
\]

\[
\;=\enspace
\begin{tangle}
\object{A}\step[2]\object{H}\step[4]\object{A}\step[7]\object{H}\step[4]\object{A}\step[6]\object{H}\\
\id\step[1]\cd\step[2]\cd\step[6]\id\step[3]\cd\step[5]\id\\
\id\step[1]\id\step[2]\x\step[2]\nw2\step[4]\cd\step[2]\id\step[2]\nw2\step[4]\id\\
\id\step[1]\tu \alpha\step[1]\cd\step[2]\cd\step[2]\id\step[2]\x\step[4]\id\step[3]\id\\
\id\step[2]\id\step[2]\id\step[2]\x\step[2]\x\step[2]\id\step[1]\cd\step[2]\cd\step[2]\id\\
\id\step[2]\id\step[2]\tu \beta\step[2]\x\step[2]\x\step[1]\id\step[2]\x\step[2]\id\step[2]\id\\
\id\step[2]\id\step[3]\id\step[2]\ne1\step[2]\x\step[2]\id\step[1]\tu \alpha\step[2]\tu \beta\step[2]\id\\
\id\step[2]\nw1\step[2]\cu \step[2]\ne2\step[2]\nw1\step[1]\cu \step[3]\ne2\step[2]\ne2\\
\nw2\step[2]\nw2\step[2]\tu \alpha\step[5]\tu \beta\step[2]\ne2\step[2]\ne2\\
\step[2]\nw2\step[2]\cu \step[7]\cu \step[2]\ne2\\
\step[4]\cu \step[9]\cu  \\
\step[5]\object{A}\step[11]\object{H}
\end{tangle}
\]

\[
\;=\enspace
\begin{tangle}
\object{A}\step[3]\object{H}\step[4]\object{A}\step[5]\object{H}\step[4]\object{A}\step[6]\object{H}\\
\id\step[2]\cd\step[2]\cd\step[3]\cd\step[2]\cd\step[5]\id\\
\id\step[1]\cd\step[1]\x\step[1]\cd\step[2]\id\step[2]\id\step[2]\id\step[2]\nw2\step[4]\id\\
\id\step[1]\id\step[2]\hx\step[2]\hx\step[2]\x\step[2]\x\step[4]\id\step[3]\id\\
\id\step[1]\tu \alpha\step[1]\tu \beta\step[1]\x\step[2]\x\step[1]\cd\step[2]\cd\step[2]\id\\
\id\step[2]\id\step[3]\cu \step[2]\x\step[2]\id\step[1]\id\step[2]\x\step[2]\id\step[2]\id\\
\id\step[2]\nw1\step[3]\id\step[2]\ne1\step[2]\id\step[2]\id\step[1]\tu \alpha\step[2]\tu \beta\step[2]\id\\
\nw2\step[2]\nw2\step[2]\tu \alpha\step[3]\nw1\step[1]\cu \step[3]\ne2\step[2]\ne2\\
\step[2]\nw2\step[2]\cu \step[5]\tu \beta\step[2]\ne2\step[2]\ne2\\
\step[4]\cu \step[7]\cu \step[2]\ne2\\
\step[5]\id\step[9]\cu \\
\step[5]\object{A}\step[10]\object{H}
\end{tangle}\ \ \ .
\]
Similarly, we can check that the right hand is equal to the diagram
above. Therefore, the equation holds.

Next we see that
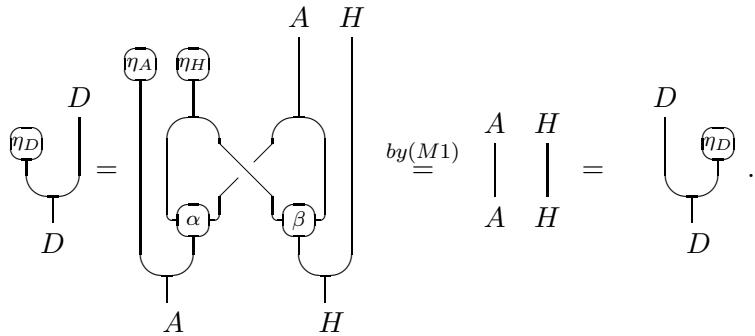
\[
\begin{tangle}
\step[2]\object{D}\\
\Q {\eta_{D}}\step[2]\id\\
\cu \\
\step[1]\object{D}
\end{tangle}
\;=\enspace
\begin{tangles}{clr}
\step[6]\object{A}\step[2]\object{H}\\
\Q { \eta_{A}}\step[2]\Q {\eta_{H}}\step[4]\id\step[2]\id\\
\id\step\cd\step[2]\cd\step[1]\id\\
\id\step[1]\id\step[2]\x\step[2]\id\step[1]\id\\
\id\step[1]\tu \alpha\step[2]\tu \beta\step[1]\id\\
\cu \step[4]\cu \\
\step[0.5]\object{A}\step[6]\object{H}
\end{tangles}
\step\stackrel{by (M1)}{=}\step
\begin{tangle}
\object{A}\step[2]\object{H}\\
\id\step[2]\id\\
\object{A}\step[2]\object{H}
\end{tangle}
\step=\step
\begin{tangles}{clr}
\object{D}\step[1]\\
\step\id\step[2]\Q {\eta_{D}}\\
\step \cu \\
\step[1.5]\object{D}
\end{tangles}\ \ \ .
\]
Thus $D$ is an algebra.

Conversely, if $D$ is an algebra, we show that $(M1),(M2)$ and
$(M3)$ hold. By
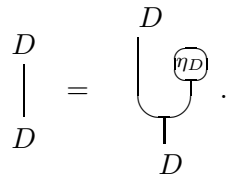
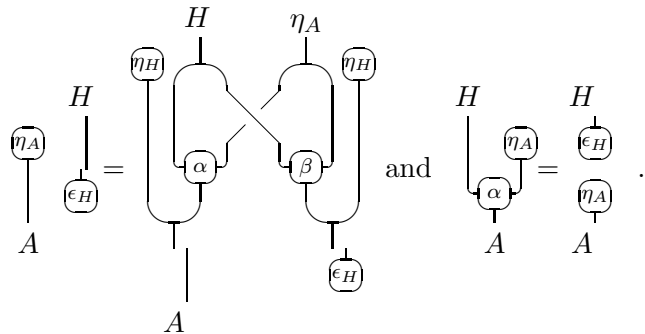
\[
\begin{tangle}
\object{D}\step[1]\\
\id\\
\object{D}
\end{tangle}
\;=\enspace
\begin{tangles}{clr}
\object{D}\\
\step\id\step[2]\Q {\eta_{D}}\\
\step \cu \\
\step[1.5]\object{D}
\end{tangles}\ \ \   .
\]
We have that

\[
\begin{tangle}
\step[2]\object{H}\\
\Q {\eta_{A}}\step[2.2]\id\\
\id\step[2]\QQ {\epsilon_{H}}\\
\object{A}
\end{tangle}
\;=\enspace
\begin{tangles}{clr}
\step[0.5]\object{H}\step[3.5]\eta_{A}\\
\Q {\eta_{H}}\step \cd\step[2]\cd\step\Q {\eta_{H}}\\
\id\step\id\step[2]\x\step[2]\id\step\id\\
\id\step[1]\tu \alpha\step[2]\tu \beta\step\id\\
\cu \step[4]\cu \\
\step[1]\id\step[6]\QQ {\epsilon_{H}}\\
\step[1]\object{A}\step[7]
\end{tangles}
\step\hbox{and }\step
\begin{tangle}
\object{H}\\
\id\step[2]\Q {\eta_{A}}\\
\tu \alpha\\
\step\object{A}
\end{tangle}
\;=\enspace
\begin{tangles}{clr}
\object{H}\step[1]\\
\QQ {\epsilon_{H}}\\
\Q {\eta_{A}}\\
\object{A}\step[1]
\end{tangles}\ \ \   .
\]

Similarly, we have that
\[
\begin{tangle}
\step[2]\object{A}\\
\Q {\eta_{H}}\step[2]\id\\
\tu \beta\\
\step[1]\object{H}
\end{tangle}
\;=\enspace
\begin{tangles}{clr}
\object{A}\\
\QQ {\epsilon_{A}}\\
\Q { \eta_{H}}\\
\object{H}
\end{tangles}\ \ \   .
\] Thus $(M1)$ holds.

By associative law, we have that
\[
\begin{tangle}
\step[2]\object{H}\step[3]\object{A}\step[6]\object{A}\\
\Q { \eta_{A}}\step[2]\id\step[3]\id\step[2]\Q {\eta_{H}}\step[4]\id\step[2]\Q {\eta_{H}}\\
\id\step[2]\id\step[3]\id\step\cd\step[2]\cd\step\id\\
\id\step[2]\id\step[3]\id\step\id\step[2]\x\step[2]\id\step\id\\
\id\step[2]\id\step[3]\id\step\tu \alpha\step[2]\tu \beta\step\id\\
\id\step[2]\id\step[3]\cu \step[4]\cu \\
\id\step\cd\step[2]\cd\step[5]\id\\
\id\step\id\step[2]\x\step[2]\id\step[4]\ne2\\
\id\step\tu \alpha\step[2]\tu \beta\step[2]\ne2\\
\cu \step[4]\cu \\
\step\id\step[6]\QQ {\epsilon_{H}}\\
\step[1]\object{A}
\end{tangle}
\;=\enspace
\begin{tangles}{clr}
\step[0.5]\object{H}\step[4]\object{A}\step[5]\object{A}\\
\Q {\eta_{A}}\step[1]\cd\step[2]\cd\step\Q {\eta_{H}}\step[3]\id\step[2]\Q {\eta_{H}}\\
\id\step\id\step[2]\x\step[2]\id\step\id\step[3]\id\step[2]\id\\
\id\step\tu \alpha\step[2]\tu \beta\step\id\step[3]\id\step[2]\id\\
\cu \step[4]\cu \step[3]\id\step[2]\id\\
\step\id\step[5]\cd\step[2]\cd \step\id\\
\step[1]\nw2\step[4]\id\step[2]\x\step[2]\id\step\id\\
\step[3]\nw2\step[2]\tu \alpha\step[2]\tu \beta\step\id\\
\step[5]\cu \step[4]\cu \\
\step[6]\id\step[6]\QQ {\epsilon_{H}}\\
\object{A}\step[1]
\end{tangles}\ \ \   .
\]
Thus (2) holds. Similarly, we can show that $(M3)$ holds. $\Box$

\begin {Lemma} \label {2.1.3}      $D =  A ^\phi \bowtie  ^\psi H
 $  is a coalgebra iff
$(CM1)$--$(CM3)$ hold.
\end {Lemma}
{\bf Proof.}  It follows by turning the diagram in proof of Theorem
\ref {2.1.2}
 upside down.                          \begin{picture}(8,8)\put(0,0){\line(0,1){8}}\put(8,8){\line(0,-1){8}}\put(0,0){\line(1,0){8}}\put(8,8){\line(-1,0){8}}\end{picture}

\begin {Theorem} \label {2.1.4}   Assume that    $H$ and $A$ are
 two bialgebras in  braided tensor category ${\cal C}$ and
$(A, \alpha )$ is a left $H$-module coalgebra,
 $(H, \beta )$ is a right $A$-module coalgebra,
 $(A, \phi )$ is a left $H$-comodule algebra, and
 $(H, \psi )$ is a right $A$-comodule algebra in ${\cal C}$. Then
     $D =  A_\alpha ^\phi \bowtie _\beta ^\psi H
 $  is a
bialgebra in ${\cal C}$ iff $(M1)$--$(M3)$, $(CM1)$--$(CM3)$,
$(B1)$--$(B5)$ hold.
\end {Theorem}

{\bf Proof.} It is sufficient to show that (B1)-(B5) hold iff

\[
\begin{tangle}
\object{D}\step[2]\object{D}\\
\cu \\
\cd\\
\object{D}\step[2]\object{D}
\end{tangle}
\;=\enspace
\begin{tangles}{clr}
\step\object{D}\step[4]\object{D}\\
\cd\step[2]\cd\\
\id\step[2]\x\step[2]\id\\
\cu \step[2]\cu \\
\step\object{D}\step[4]\object{D}
\end{tangles}\ \ \   . \ \ \ \
\cdots\cdots(1)
\]

 If $(B1)-(B5)$ hold, see that
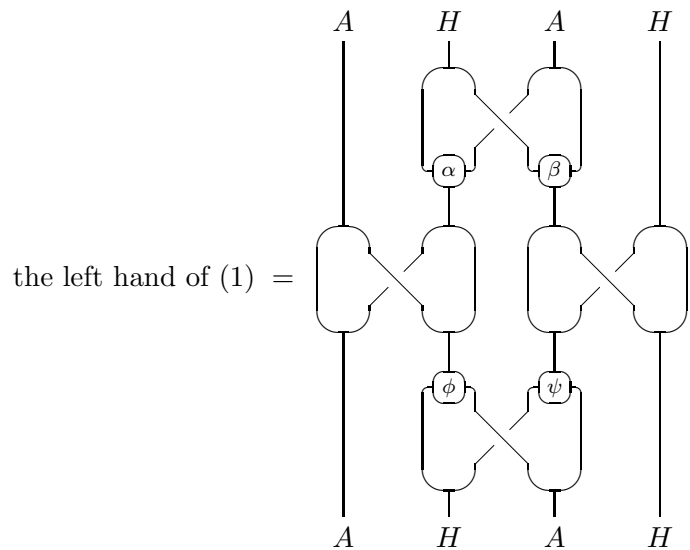
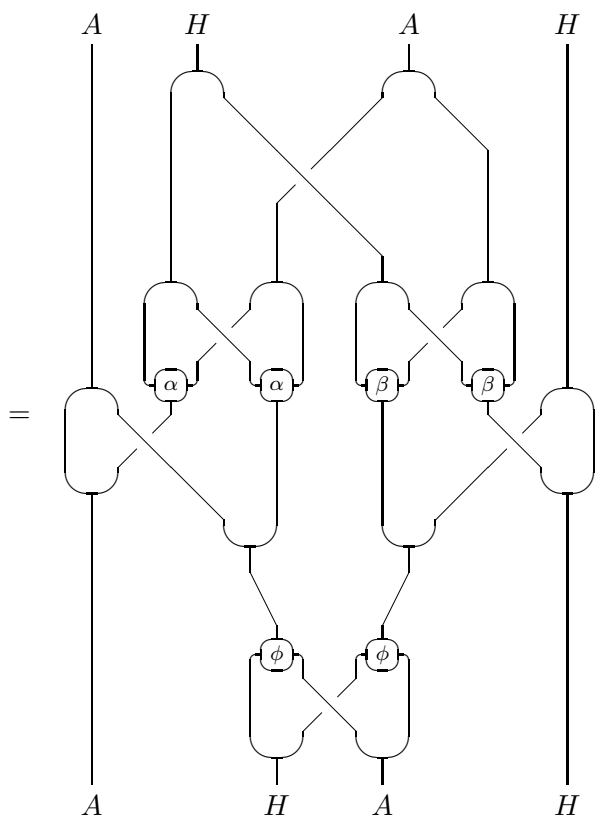
\[
\hbox {the left hand of (1)} \;=\enspace
\begin{tangle}
\step\object{A}\step[4]\object{H}\step[4]\object{A}\step[4]\object{H}\\
\step\id\step[3]\cd\step[2]\cd\step[3]\id\\
\step\id\step[3]\id\step[2]\x\step[2]\id\step[3]\id\\
\step\id\step[3]\tu \alpha\step[2]\tu \beta\step[3]\id\\
\cd\step[2]\cd\step[2]\cd\step[2]\cd\\
\id\step[2]\x\step[2]\id\step[2]\id\step[2]\x\step[2]\id\\
\cu \step[2]\cu \step[2]\cu \step[2]\cu \\
\step\id\step[3]\td \phi\step[2]\td \psi\step[3]\id\\
\step\id\step[3]\id\step[2]\x\step[2]\id\step[3]\id\\
\step\id\step[3]\cu \step[2]\cu \step[3]\id\\
\step\object{A}\step[4]\object{H}\step[4]\object{A}\step[4]\object{H}
\end{tangle}
\]

\[
\step=\step
\begin{tangle}
\step\object{A}\step[4]\object{H}\step[8]\object{A}\step[6]\object{H}\\
\step\id\step[3]\cd\step[6]\cd\step[5]\id\\
\step\id\step[3]\id\step[2]\nw2\step[4]\ne2\step[2]\nw2\step[4]\id\\
\step\id\step[3]\id\step[4]\x\step[6]\id\step[3]\id\\
\step\id\step[3]\id\step[4]\id\step[2]\nw2\step[5]\id\step[3]\id\\
\step\id\step[2]\cd\step[2]\cd\step[2]\cd\step[2]\cd\step[2]\id\\
\step\id\step[2]\id\step[2]\x\step[2]\id\step[2]\id\step[2]\x\step[2]\id\step[2]\id\\
\cd\step \tu \alpha\step[2]\tu \alpha\step[2]\tu \beta\step[2]\tu \beta\step[1]\cd\\
\id\step[2]\x\step[4]\id\step[4]\id\step[4]\x\step[2]\id\\
\cu \step[2]\nw2\step[3]\id\step[4]\id\step[3]\ne2\step[2]\cu \\
\step\id\step[5]\cu \step[4]\cu \step[5]\id\\
\step\id\step[6]\nw1\step[4]\ne1\step[6]\id\\
\step\id\step[6]\td \phi\step[2]\td \phi\step[6]\id\\
\step\id\step[6]\id\step[2]\x\step[2]\id\step[6]\id\\
\step\id\step[6]\cu \step[2]\cu  \step[6]\id\\
\step\object{A}\step[7]\object{H}\step[4]\object{A}\step[7]\object{H}
\end{tangle}
\]

\[
\step\stackrel{ \hbox {by Lemma }\ref {2.1.1}(vii)(viii)}{=}\step
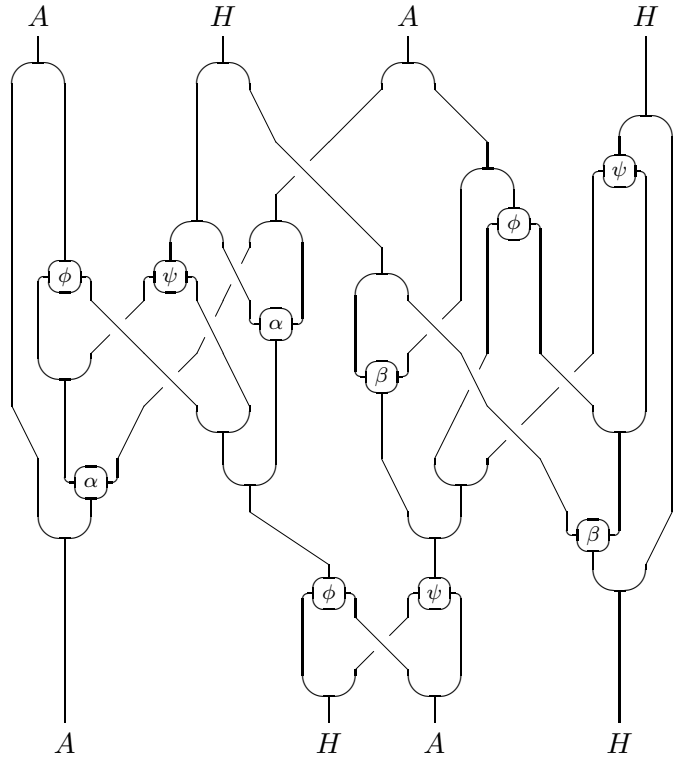
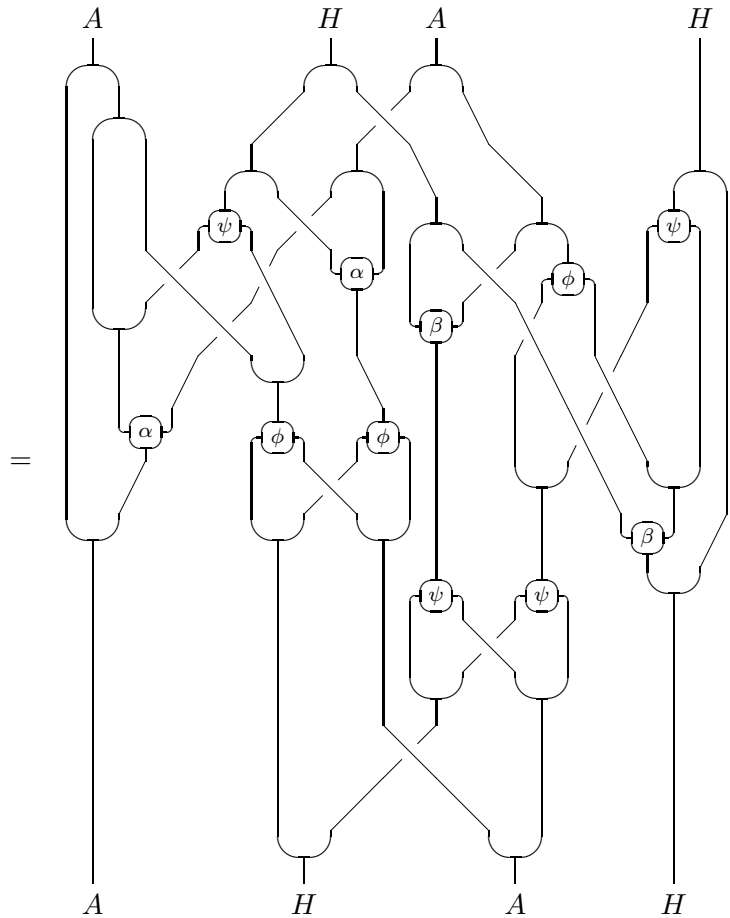
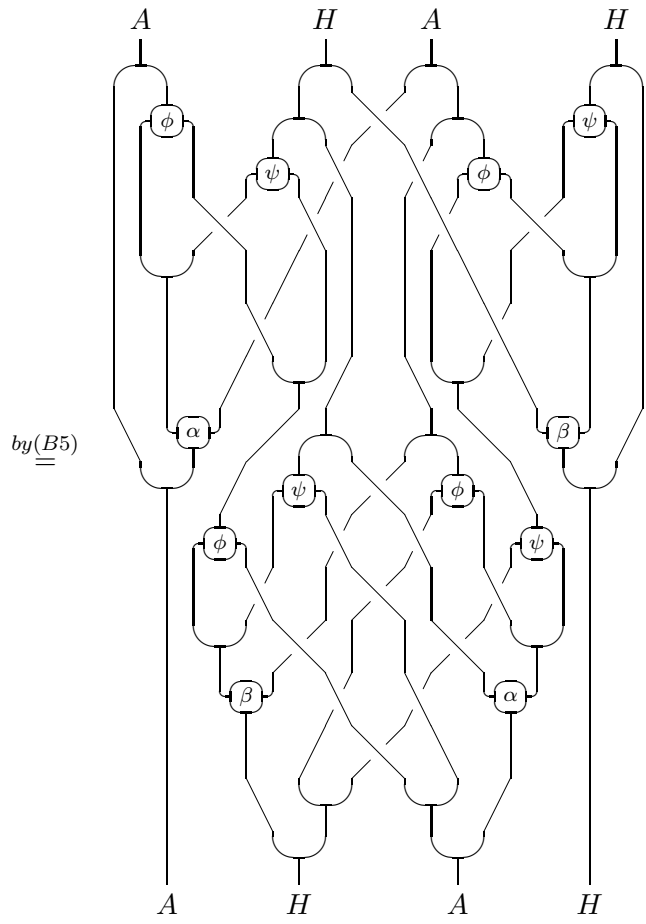
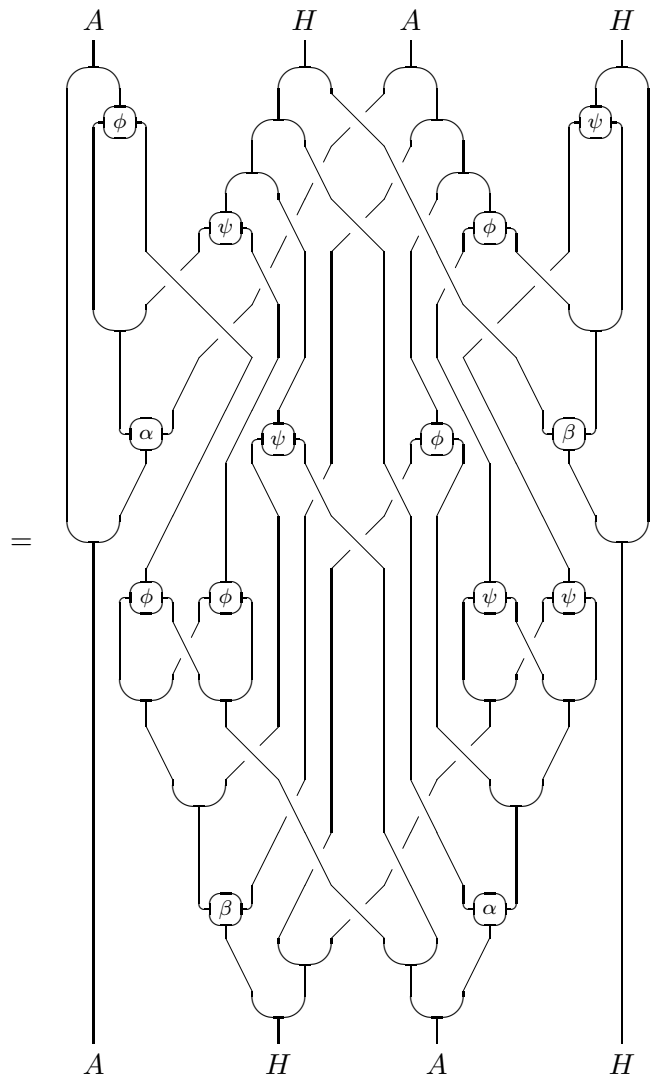
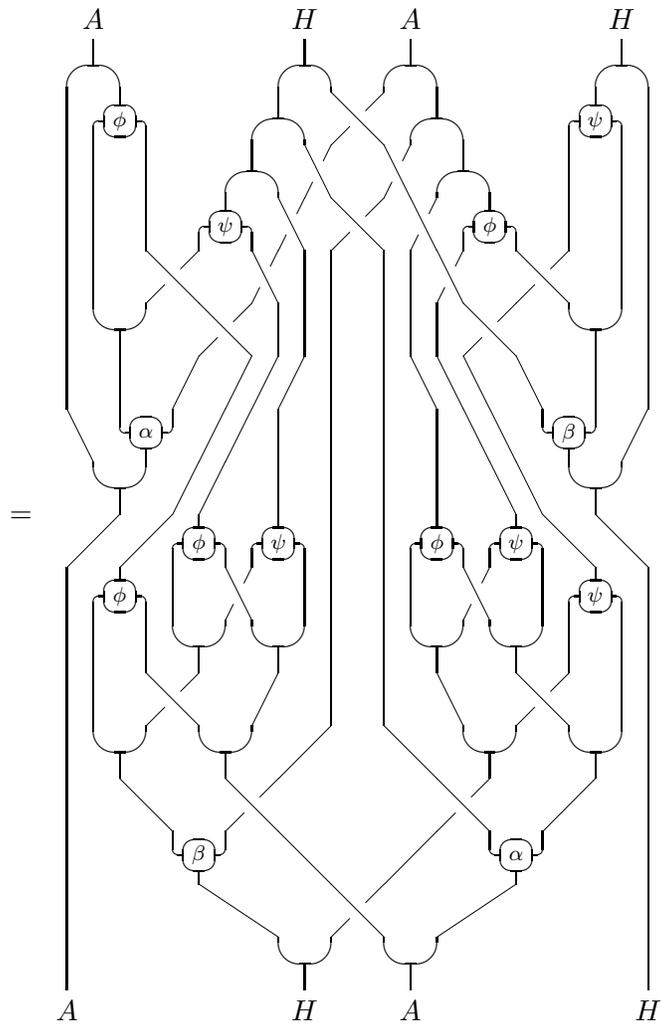
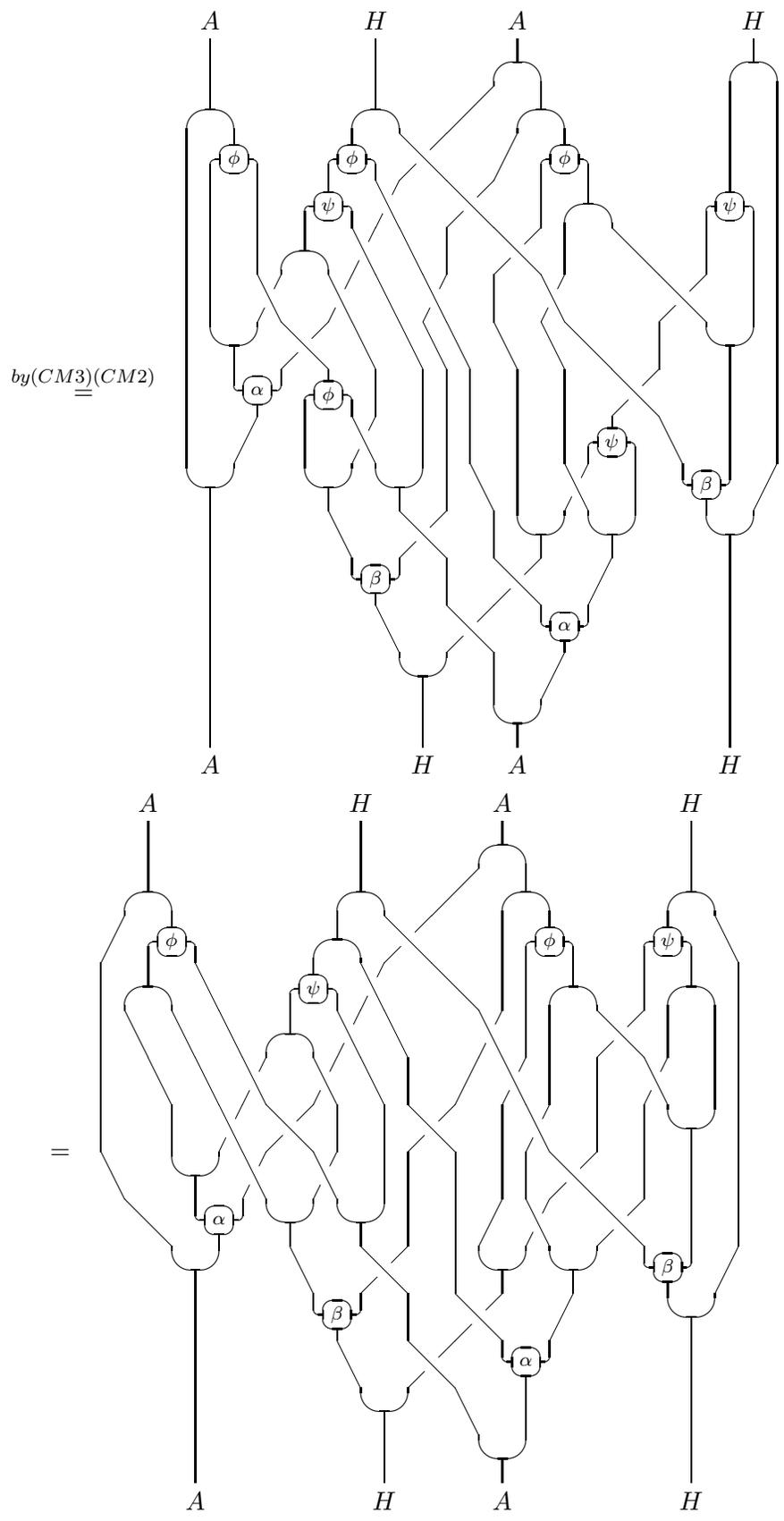
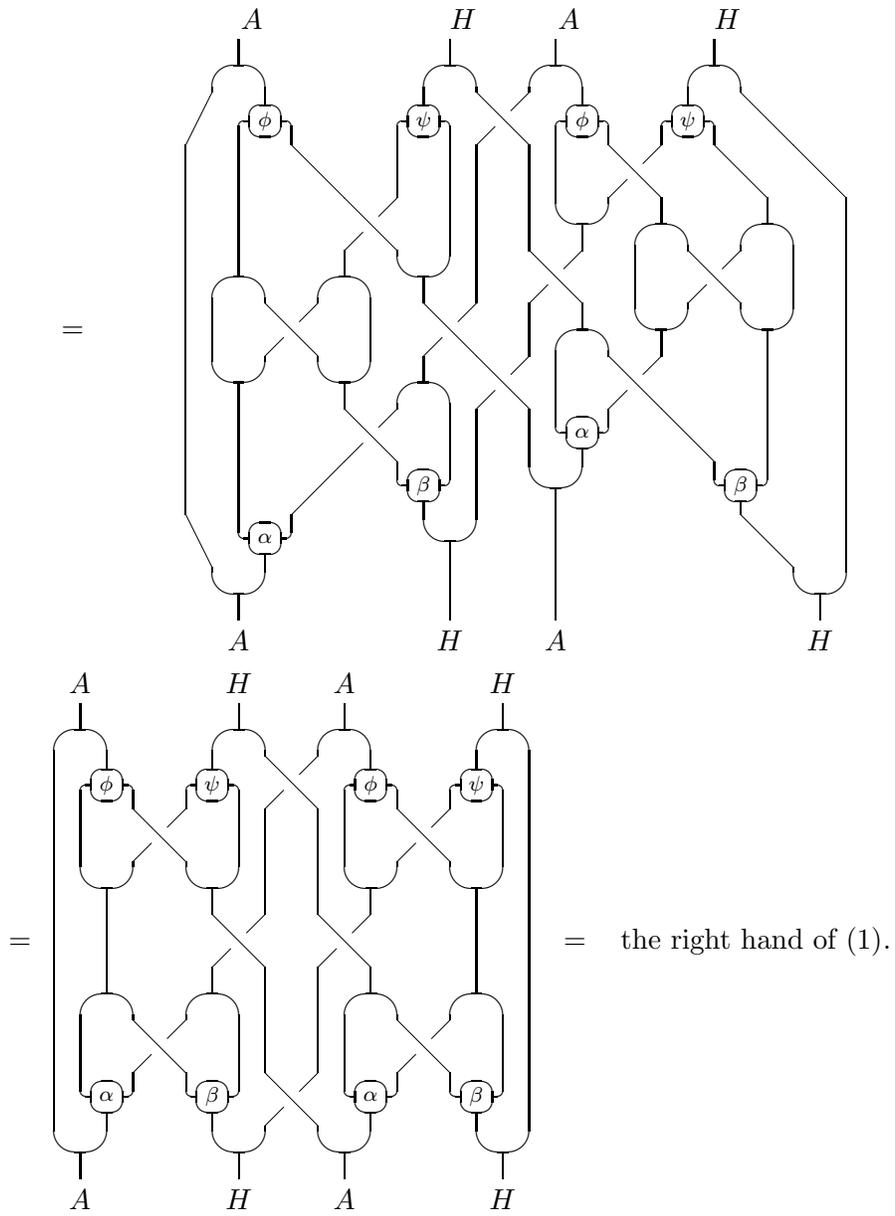
\begin{tangle}
\step[1]\object{A}\step[7]\object{H}\step[7]\object{A}\step[9]\object{H}\\
\cd\step[5]\cd\step[5]\cd\step[8]\id\\
\id\step[2]\id\step[5]\id\step[2]\nw1\step[3]\ne2\step[2]\nw2\step[6]\cd\\
\id\step[2]\id\step[5]\id\step[3]\x\step[5]\cd\step[3]\td \psi\step\id\\
\id\step[2]\id\step[4]\cd\step[1]\cd\step[1]\nw2\step[4]\id\step\td \phi\step[2]\id\step[2]\id\step\id\\
\id\step[1]\td \phi\step[2]\td \psi\step[1]\hx\step[2]\id\step[2]\cd\step[2]\id\step\id\step[2]\id\step[2]\id\step[2]\id\step\id\\
\id\step[1]\id\step[2]\x\step[2]\hx\step \tu \alpha\step[2]\id\step[2]\x\step \id\step[2]\id\step[2]\id\step[2]\id\step\id\\
\id\step\cu \step[2]\x\step \nw1\step \id\step[3]\tu \beta\step[2]\hx\step[2]\x\step[2]\id\step\id\\
\nw1\step[1]\id\step[2]\ne1\step[2]\cu \step\id\step[4]\id\step[2]\ne1\step\x\step[2]\cu \step\id\\
\step[1]\id\step[1]\tu \alpha\step[4]\cu \step[4]\nw1\step[1]\cu \step[2]\nw1\step[2]\id\step[2]\id\\
\step[1]\cu \step[6]\nw3\step[5]\cu \step[4]\tu \beta\step\ne1\\
\step[2]\id\step[9]\td \phi\step[2]\td \psi\step[5]\cu \\
\step[2]\id\step[9]\id\step[2]\x\step[2]\id\step[6]\id\\
\step[2]\id\step[9]\cu \step[2]\cu \step[6]\id\\
\step[2]\object{A}\step[10]\object{H}\step[4]\object{A}\step[7]\object{H}
\end{tangle}
\]

\[
\step=\step
\begin{tangle}
\step[1]\object{A}\step[9]\object{H}\step[4]\object{A}\step[10]\object{H}\\
\cd\step[7]\cd\step[2]\cd\step[9]\id\\
\id\step\cd\step[5]\ne2\step[2]\x\step[2]\nw1\step[8]\id\\
\id\step\id\step[2]\id\step[3]\cd\step[2]\cd\step\nw1\step[2]\nw2\step[6]\cd\\
\id\step\id\step[2]\id\step[2]\td \psi\step\x\step[2]\id\step\cd\step[2]\cd\step[3]\td \psi\step\id\\
\id\step\id\step[2]\x\step[2]\hx\step[2]\tu \alpha\step\id\step[2]\x\step\td \phi\step[2]\id\step[2]\id\step\id\\
\id\step\cu \step[2]\x\step\nw1\step[2]\id\step[2]\tu \beta\step[2] \hx\step[2]\id\step\ne1\step[2]\id\step\id\\
\id\step[2]\id\step[2]\ne1\step[2] \cu \step[2] \nw1\step[2]\id\step[3]\id\step\nw1\step\hx\step[3]\id\step\id\\
\id\step[2]\tu \alpha\step[3]\td \phi\step[2]\td \phi\step\id\step[3]\id\step[2]\hx\step\nw1\step[2]\id\step\id\\
\id\step[2]\ne1\step[4]\id\step[2]\x\step[2]\id\step\id\step[3]\cu \step\nw1\step\cu \step\id\\
\cu \step[5]\cu \step[2]\cu \step\id\step[4]\id\step[3]\tu \beta\step\ne1\\
\step[1]\id\step[7]\id\step[4]\id\step\td \psi\step[2]\td \psi\step[3]\cu \\
\step[1]\id\step[7]\id\step[4]\id\step\id\step[2]\x\step[2]\id\step[4]\id\\
\step[1]\id\step[7]\id\step[4]\id\step\cu \step[2]\cu \step[4]\id\\
\step[1]\id\step[7]\id\step[4]\x\step[4]\id\step[5]\id\\
\step[1]\id\step[7]\id\step[3]\ne2\step[2]\nw2\step[3]\id\step[5]\id\\
\step[1]\id\step[7]\cu \step[6]\cu \step[5]\id\\
\step[1]\object{A}\step[8]\object{H}\step[8]\object{A}\step[6]\object{H}
\end{tangle}
\]

\[\step\stackrel{by(B5)}{=}\step
\begin{tangle}
\step[1]\object{A}\step[7]\object{H}\step[4]\object{A}\step[7]\object{H}\\
\cd\step[5]\cd\step[2]\cd\step[5]\cd\\
\id\step[1]\td \phi\step[3]\cd\step[1]\x\step[1]\cd\step[3]\td \psi\step[1]\id\\
\id\step[1]\id\step[2]\id\step[2]\td \psi\step[1]\hx\step[2]\hx\step[1]\td \phi\step[2]\id\step[2]\id\step[1]\id\\
\id\step[1]\id\step[2]\x\step[2]\hx\step[1]\id\step[2]\id\step[1]\hx\step[2]\x\step[2]\id\step[1]\id\\
\id\step[1]\cu \step[2]\id\step[1]\ne1\step[1]\id\step[1]\id\step[2]\id\step[1]\id\step[1]\nw1\step[1]\id\step[2]\cu \step[1]\id\\
\id\step[2]\id\step[3]\hx\step[2]\id\step[1]\id\step[2]\id\step[1]\id\step[2]\hx\step[3]\id\step[2]\id\\
\id\step[2]\id\step[2]\ne1\step[1]\cu \ne1\step[2]\nw1\cu \step[1]\nw1\step[2]\id\step[2]\id\\
\nw1\step[1]\tu \alpha\step[2]\ne2\cd\step[2]\cd\nw2\step[2]\tu \beta\step[1]\ne1\\
\step[1]\cu \step[1]\ne1\step[1]\td \psi\step[1]\x\step[1]\td \phi\step[1]\nw1\step[1]\cu \\
\step[2]\id\step[1]\td \phi\step[1]\id\step[2]\hx\step[2]\hx\step[2]\id\step[1]\td \psi\step[1]\id\\
\step[2]\id\step[1]\id\step[2]\hx\step[2]\id\step[1]\x\step[1]\id\step[2]\hx\step[2]\id\step[1]\id\\
\step[2]\id\step[1]\cu \step[1]\x\step[1]\id\step[2]\id\step[1]\x\step[1]\cu \step[1]\id\\
\step[2]\id\step[2]\tu \beta\step[2]\hx\step[2]\hx\step[2]\tu \alpha\step[2]\id\\
\step[2]\id\step[3]\id\step[2]\ne1\step[1]\x\step[1]\nw1\step[2]\id\step[3]\id\\
\step[2]\id\step[3]\nw1\step[1]\cu \step[2]\cu \step[1]\ne1\step[3]\id\\
\step[2]\id\step[4]\cu \step[4]\cu \step[4]\id\\
\step[2]\object{A}\step[5]\object{H}\step[6]\object{A}\step[5]\object{H}
\end{tangle}
\]
\[\step=\step
\begin{tangle}
\step[1]\object{A}\step[8]\object{H}\step[4]\object{A}\step[8]\object{H}\\
\cd\step[6]\cd\step[2]\cd\step[6]\cd\\
\id\step[1]\td \phi\step[4]\cd\step[1]\x\step[1]\cd\step[4]\td \psi\step[1]\id\\
\id\step[1]\id\step[2]\id\step[3]\cd\step[1]\hx\step[2]\hx\step[1]\cd\step[3]\id\step[2]\id\step[1]\id\\
\id\step[1]\id\step[2]\id\step[2]\td \psi\step[1]\hx\step[1]\x\step[1]\hx\step[1]\td \phi\step[2]\id\step[2]\id\step[1]\id\\
\id\step[1]\id\step[2]\x\step[2]\hx\step[1]\id\step[1]\id\step[2]\id\step[1]\id\step[1]\hx\step[2]\x\step[2]\id\step[1]\id\\
\id\step[1]\cu \step[2]\x\step[1]\id\step[1]\id\step[1]\id\step[2]\id\step[1]\id\step[1]\id\step[1]\x\step[2]\cu \step[1]\id\\
\id\step[2]\id\step[2]\ne1\step[1]\ne1\ne1\ne1\step[1]\id\step[2]\id\step[1]\nw1\nw1\nw1\step[1]\nw1\step[2]\id\step[2]\id\\
\id\step[2]\tu \alpha\step[1]\ne1\ne1\td \psi\step[1]\id\step[2]\id\step[1]\td \phi\nw1\nw1\step[1]\tu \beta\step[2]\id\\
\id\step[2]\ne1\step[1]\ne1\step[1]\id\step[1]\nw1\step[1]\hx\step[2]\hx\step[1]\ne1\step[1]\id\step[1]\nw1\step[1]\nw1\step[2]\id\\
\cu \step[1]\ne1\step[2]\id\step[2]\id\step[1]\id\step[1]\x\step[1]\id\step[1]\id\step[2]\id\step[2]\nw1\step[1]\cu \\
\step[1]\id\step[1]\td \phi\step[1]\td \phi\step[1]\id\step[1]\id\step[1]\id\step[2]\id\step[1]\id\step[1]\id\step[1]\td \psi\step[1]\td \psi\step[1]\id\\
\step[1]\id\step[1]\id\step[2]\hx\step[2]\id\step[1]\id\step[1]\id\step[1]\id\step[2]\id\step[1]\id\step[1]\id\step[1]\id\step[2]\hx\step[2]\id\step[1]\id\\
\step[1]\id\step[1]\cu \step[1]\cu \step[1]\id\step[1]\id\step[1]\id\step[2]\id\step[1]\id\step[1]\id\step[1]\cu \step[1]\cu \step[1]\id\\
\step[1]\id\step[2]\nw1\step[2]\x\step[1]\id\step[1]\id\step[2]\id\step[1]\id\step[1]\x\step[2]\ne1\step[2]\id\\
\step[1]\id\step[3]\cu \step[2]\hx\step[1]\id\step[2]\id\step[1]\hx\step[2]\cu \step[3]\id\\
\step[1]\id\step[4]\id\step[2]\ne1\step[1]\hx\step[2]\hx\step[1]\nw1\step[2]\id\step[4]\id\\
\step[1]\id\step[4]\tu \beta\step[1]\ne1\step[1]\x\step[1]\nw1\step[1]\tu \alpha\step[4]\id\\
\step[1]\id\step[5]\nw1\step[1]\cu \step[2]\cu \step[1]\ne1\step[5]\id\\
\step[1]\id\step[6]\cu \step[4]\cu \step[6]\id\\
\step[1]\object{A}\step[7]\object{H}\step[6]\object{A}\step[7]\object{H}
\end{tangle}\
\]

\[
\step=\step
\begin{tangle}
\step[1]\object{A}\step[8]\object{H}\step[4]\object{A}\step[8]\object{H}\\
\cd\step[6]\cd\step[2]\cd\step[6]\cd\\
\id\step[1]\td \phi\step[4]\cd\step[1]\x\step[1]\cd\step[4]\td \psi\step[1]\id\\
\id\step[1]\id\step[2]\id\step[3]\cd\step[1]\hx\step[2]\hx\step[1]\cd\step[3]\id\step[2]\id\step[1]\id\\
\id\step[1]\id\step[2]\id\step[2]\td \psi\step[1]\hx\step[1]\x\step[1]\hx\step[1]\td \phi\step[2]\id\step[2]\id\step[1]\id\\
\id\step[1]\id\step[2]\x\step[2]\hx\step[1]\id\step[1]\id\step[2]\id\step[1]\id\step[1]\hx\step[2]\x\step[2]\id\step[1]\id\\
\id\step[1]\cu \step[2]\x\step[1]\id\step[1]\id\step[1]\id\step[2]\id\step[1]\id\step[1]\id\step[1]\x\step[2]\cu \step[1]\id\\
\id\step[2]\id\step[2]\ne1\step[1]\ne1\ne1\ne1\step[1]\id\step[2]\id\step[1]\nw1\nw1\nw1\step[1]\nw1\step[2]\id\step[2]\id\\
\nw1\step[1]\tu \alpha\step[1]\ne1\ne1\step[1]\id\step[2]\id\step[2]\id\step[2]\id\step[1]\nw1\nw1\step[1]\tu \beta\step[1]\ne1\\
\step[1]\cu \step[1]\ne1\ne1\step[2]\id\step[2]\id\step[2]\id\step[2]\id\step[2]\nw1\nw1\step[1]\cu \\
\step[1]\ne2\step[1]\ne2\td \phi\step[1]\td \psi\step[1]\id\step[2]\id\step[1]\td \phi\step[1]\td \psi\nw2\step[1]\nw2\\
\id\step[1]\td \phi\step[1]\id\step[2]\hx\step[2]\id\step[1]\id\step[2]\id\step[1]\id\step[2]\hx\step[2]\id\step[1]\td \psi\step[1]\id\\
\id\step[1]\id\step[2]\id\step[1]\cu \step[1]\cu \step[1]\id\step[2]\id\step[1]\cu \step[1]\cu \step[1]\id\step[2]\id\step[1]\id\\
\id\step[1]\id\step[2]\x\step[2]\ne1\step[2]\id\step[2]\id\step[2]\nw1\step[2]\x\step[2]\id\step[1]\id\\
\id\step[1]\cu \step[2]\cu \step[2]\ne2\step[2]\nw2\step[2]\cu \step[2]\cu \step[1]\id\\
\id\step[2]\nw2\step[3]\x\step[6]\x\step[3]\ne2\step[2]\id\\
\id\step[4]\tu \beta\step[2]\nw2\step[4]\ne2\step[2]\tu \alpha\step[4]\id\\
\id\step[5]\nw3\step[4]\x\step[4]\ne3\step[5]\id\\
\id\step[8]\cu \step[2]\cu \step[8]\id\\
\object{A}\step[9]\object{H}\step[4]\object{A}\step[9]\object{H}
\end{tangle}
\]

\[
\step\stackrel{by (CM3)(CM2)}{=}\step
\begin{tangle}
\step[1]\object{A}\step[7]\object{H}\step[6]\object{A}\step[10]\object{H}\\
\step[1]\id\step[7]\id\step[5]\cd\step[8]\cd\\
\cd\step[5]\cd\step[3]\ne2\step[1]\cd\step[7]\id\step[2]\id\\
\id\step[1]\td \phi\step[3]\td \phi\step[1]\x\step[2]\ne1\step[1]\td \phi\step[6]\id\step[2]\id\\
\id\step[1]\id\step[2]\id\step[2]\td \psi\step[1]\hx\step[2]\x\step[2]\id\step[1]\cd\step[4]\td \psi\step[1]\id\\
\id\step[1]\id\step[2]\id\step[1]\cd\step[1]\hx\step[1]\nw1\step[1]\id\step[2]\x\step[1]\id\step[2]\nw2\step[3]\id\step[2]\id\step[1]\id\\
\id\step[1]\id\step[2]\hx\step[2]\hx\step[1]\nw1\step[1]\hx\step[2] \id\step[2]\hx\step[4]\x\step[2]\id\step[1]\id\\
\id\step[1]\cu \step[1]\x\step[1]\nw1\step[1]\nw1\nw1\nw1\step[1]\nw1\step[1]\nw1\nw2\step[3]\id\step[2]\cu \step[1]\id\\
\id\step[2]\tu \alpha\step[1]\td \phi\step[1]\id\step[2]\id\step[1]\id\step[1]\id\step[2]\id\step[2]\id\step[2]\x\step[3]\id\step[2]\id\\
\id\step[2]\ne1\step[2]\id\step[2]\hx\step[2]\id\step[1]\id\step[1]\id\step[2]\id\step[2]\id\step[1]\td \psi\step[1]\nw1\step[2]\id\step[2]\id\\
\cu \step[3]\cu \step[1]\cu \step[1]\id\step[1]\nw1\step[1]\id \step[2]\hx\step[2]\id\step[2]\tu \beta\step[1]\ne1\\
\step[1]\id\step[5]\nw1\step[2]\x\step[2]\id\step[1]\cu \step[1]\cu \step[3]\cu \\
\step[1]\id\step[6]\tu \beta\step[2]\id\step[2]\x\step[2]\ne1\step[5]\id\\
\step[1]\id\step[7]\nw1\step[2]\x\step[2]\tu \alpha\step[6]\id\\
\step[1]\id\step[8]\cu \step[2]\id\step[2]\ne1\step[7]\id\\
\step[1]\id\step[9]\id\step[3]\cu \step[8]\id\\
\step[1]\object{A}\step[9]\object{H}\step[4]\object{A}\step[9]\object{H}
\end{tangle}
\]

\[
\step=\step
\begin{tangle}
\step[2]\object{A}\step[9]\object{H}\step[6]\object{A}\step[8]\object{H}\\
\step[2]\id\step[9]\id\step[5]\cd\step[7]\id\\
\step[1]\cd\step[7]\cd\step[3]\ne2\step[1]\cd\step[5]\cd\\
\ne1\step[1]\td \phi\step[5]\cd\step[1]\x\step[3]\id\step[1]\td \phi\step[3]\td \psi\step[1]\nw1\\
\id\step[1]\cd\step[1]\nw1\step[3]\td \psi\step[1]\hx\step[2]\nw2\step[2]\id\step[1]\id\step[1]\cd\step[2]\id\step[1]\cd\step[1]\id\\
\id\step[1]\nw1\step[1]\nw1\step[1]\nw1\step[1]\cd\step[1]\hx\step[1]\nw1\step[3]\hx\step[1]\id\step[1]\id\step[2]\x\step[1]\id\step[2]\id\step[1]\id\\
\id\step[2]\nw1\step[1]\nw1\step[1]\hx\step[2]\hx\step[1]\nw1\step[1]\id\step[2]\ne1\step[1]\hx\step[1]\id\step[2]\id\step[2]\hx\step[2]\id\step[1]\id\\
\id\step[3]\id\step[2]\hx\step[1]\x\step[1]\id\step[2]\id\step[1]\x\step[2]\id\step[1]\hx\step[2]\id\step[2]\id\step[1]\cu \step[1]\id\\
\nw1\step[2]\cu \step[1]\hx\step[2]\hx\step[2]\id\step[1]\id\step[2]\id\step[2]\id\step[1]\id\step[1]\x\step[2]\id\step[2]\id\step[2]\id\\
\step[1]\nw2\step[2]\tu \alpha\step[1]\cu \step[1]\cu \step[1]\id\step[2]\id\step[1]\ne1\step[1]\hx\step[2]\x\step[2]\id\step[2]\id\\
\step[3]\cu \step[3]\nw1\step[2]\x\step[2]\id\step[1]\cu \step[1]\cu \step[2]\tu \beta\step[1]\ne1\\
\step[4]\id\step[5]\tu \beta\step[2]\id\step[2]\x\step[2]\ne1\step[4]\cu \\
\step[4]\id\step[6]\nw1\step[2]\x\step[2]\tu \alpha\step[6]\id\\
\step[4]\id\step[7]\cu \step[2]\nw1\step[2]\id\step[7]\id\\
\step[4]\id\step[8]\id\step[4]\cu \step[7]\id\\
\step[4]\object{A}\step[8]\object{H}\step[5]\object{A}\step[8]\object{H}
\end{tangle}
\]

\[
\begin{tangle}
\end{tangle}
\;=\enspace
\begin{tangles}{clr}
\step[1]\object{A}\step[8]\object{H}\step[4]\object{A}\step[6]\object{H}\\
\cd\step[6]\cd\step[2]\cd\step[4]\cd\\
\ne1\step[1]\td \phi\step[4]\td \psi\step[1]\x\step[1]\td \phi\step[2]\td \psi\step[1]\nw2\\
\step\step[1]\id\step[2]\id\step[2]\nw2\step[3]\id\step[2]\id\step[1]\id\step[2]\id\step[1]\id\step[2]\x\step[2]\nw2\step[2]\nw2\\
\step\step[1]\step\id\step[2]\id\step[4]\x\step[2]\id\step[1]\id\step[2]\id\step[1]\cu \step[1]\cd\step[2]\cd\step[2]\id\\
\step\step[1]\step\id\step[1]\cd\step[2]\cd\step[1]\cu \step[1]\id\step[2]\x\step[2]\id\step[2]\x\step[2]\id\step[2]\id\\
\step\step[1]\step\id\step[1]\id\step[2]\x\step[2]\id\step[2]\x\step[2]\id\step[1]\cd\step[1]\cu \step[2]\cu \step[2]\id\\
\step\step[1]\step\id\step[1]\cu \step[2]\cu \step[1]\cd\step[1]\x\step[1]\id\step[2]\x\step[4]\id\step[3]\id\\
\step\step[1]\step\id\step[2]\id\step[4]\x\step[2]\id\step[1]\id\step[2]\id\step[1]\tu \alpha\step[2]\nw2\step[3]\id\step[3]\id\\
\step\step[1]\step\id\step[2]\id\step[3]\ne2\step[2]\tu \beta\step[1]\id\step[2]\cu \step[5]\tu \beta\step[3]\id\\
\step\step[1]\step\nw1\step[1]\tu \alpha\step[5]\cu \step[3]\id\step[7]\nw2\step[3]\id\\
\step\step[1]\step\step[1]\cu \step[7]\id\step[4]\id\step[9]\cu \\
\step\step[1]\step[2]\object{A}\step[8]\object{H}\step[4]\object{A}\step[10]\object{H}
\end{tangles}
\]

\[
\;=\enspace
\begin{tangle}
\step[1]\object{A}\step[6]\object{H}\step[4]\object{A}\step[6]\object{H}\\
\cd\step[4]\cd\step[2]\cd\step[4]\cd\\
\id\step[1]\td \phi\step[2]\td \psi\step[1]\x\step[1]\td \phi\step[2]\td \psi\step[1]\id\\
\id\step[1]\id\step[2]\x\step[2]\id\step[1]\id\step[2]\id\step[1]\id\step[2]\x\step[2]\id\step[1]\id\\
\id\step[1]\cu \step[2]\cu \step[1]\id\step[2]\id\step[1]\cu \step[2]\cu \step[1]\id\\
\id\step[2]\id\step[4]\x\step[2]\x\step[4]\id\step[2]\id\\
\id\step[1]\cd\step[2]\cd\step[1]\id\step[2]\id\step[1]\cd\step[2]\cd\step[1]\id\\
\id\step[1]\id\step[2]\x\step[2]\id\step[1]\id\step[2]\id\step[1]\id\step[2]\x\step[2]\id\step[1]\id\\
\id\step[1]\tu \alpha\step[2]\tu \beta\step[1]\x\step[1]\tu \alpha\step[2]\tu \beta\step[1]\id\\
\cu \step[4]\cu \step[2]\cu \step[4]\cu \\
\step[1]\object{A}\step[6]\object{H}\step[4]\object{A}\step[6]\object{H}
\end{tangle}
\step=\step \hbox {the right hand of (1)}.
\]

Conversely, if (1) holds, we have that
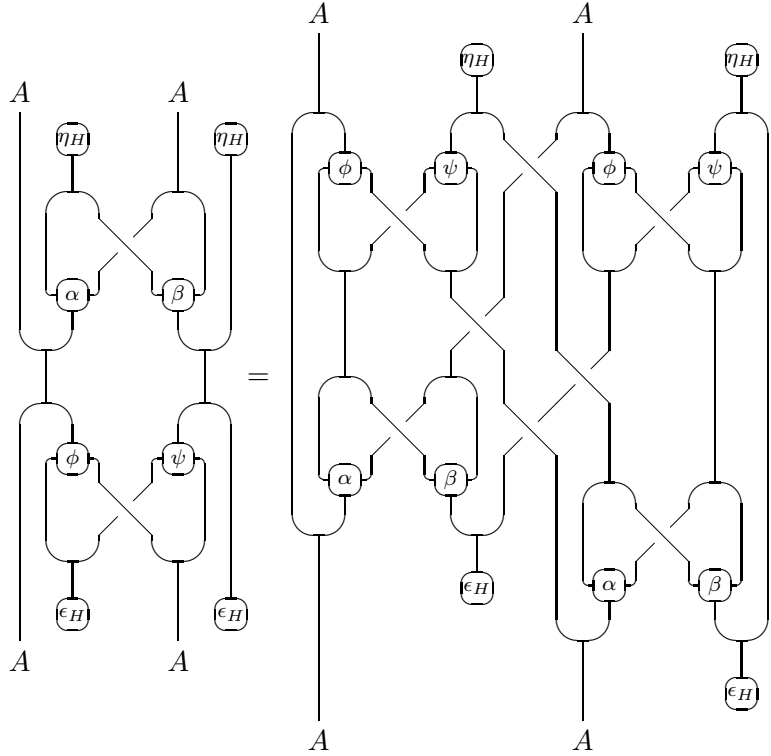
\[
\begin{tangle}
\object{A}\step[6]\object{A}\\
\id\step[2]\Q { \eta_{H}}\step[4]\id\step[2]\Q {\eta_{H}}\\
\id\step[1]\cd\step[2]\cd\step[1]\id\\
\id\step[1]\id\step[2]\x\step[2]\id\step[1]\id\\
\id\step[1]\tu \alpha\step[2]\tu \beta\step[1]\id\\
\cu \step[4]\cu \\
\cd\step[4]\cd\\
\id\step[1]\td \phi\step[2]\td \psi\step[1]\id\\
\id\step[1]\id\step[2]\x\step[2]\id\step[1]\id\\
\id\step[1]\cu \step[2]\cu \step[1]\id\\
\id\step[2]\QQ {\epsilon_{H}}\step[4]\id\step[2]\QQ {\epsilon_{H}}\\
\object{A}\step[6]\object{A}
\end{tangle}
\;=\enspace
\begin{tangles}{clr}
\object{A}\step[10]\object{A}\step[6]\\
\id\step[6]\Q {\eta_{H}}\step[4]\id\step[6]\Q {\eta_{H}}\\
\cd\step[4]\cd\step[2]\cd\step[4]\cd\\
\id\step[1]\td \phi\step[2]\td \psi\step[1]\x\step[1]\td \phi\step[2]\td \psi\step[1]\id\\
\id\step[1]\id\step[2]\x\step[2]\id\step[1]\id\step[2]\id\step[1]\id\step[2]\x\step[2]\id\step[1]\id\\
\id\step[1]\cu \step[2]\cu \step[1]\id\step[2]\id\step[1]\cu \step[2]\cu \step[1]\id\\
\id\step[2]\id\step[4]\x\step[2]\id\step[2]\id\step[4]\id\step[2]\id\\
\id\step[1]\cd\step[2]\cd\step[1]\id\step[2]\x\step[4]\id\step[2]\id\\
\id\step[1]\id\step[2]\x\step[2]\id\step[1]\x\step[2]\id\step[4]\id\step[2]\id\\
\id\step[1]\tu \alpha\step[2]\tu \beta\step[1]\id\step[2]\id\step[1]\cd\step[2]\cd\step[1]\id\\
\cu \step[4]\cu \step[2]\id\step[1]\id\step[2]\x\step[2]\id\step[1]\id\\
\step[1]\id\step[6]\QQ {\epsilon_{H}}\step[3]\id\step[1]\tu \alpha\step[2]\tu \beta\step[1]\id\\
\step[1]\id\step[9]\cu \step[4]\cu \\
\id\step[10]\id\step[6]\QQ {\epsilon_{H}}\\
\step[1]\object{A}\step[10]\object{A}\step[7]
\end{tangles}
\]
and obtain $(B1)$ by computation.

Similarly, we can obtain $(B2)$  by  applying  $\eta _A  {\otimes}
id  {\otimes} id  {\otimes}  \eta _H $  to the top of  relation (1)
and
 $ id  {\otimes} \epsilon_H  {\otimes}   id
  {\otimes} \epsilon _H$ to the bottom of relation (1) . We can also obtain $(B3)$  and $(B4)$  by
 relation (1), $\eta _A  {\otimes}
id  {\otimes} \eta _A  {\otimes}  id , \epsilon _A  {\otimes} id
{\otimes}   \epsilon_A
  {\otimes} id, \eta _A  {\otimes}
id  {\otimes} id  {\otimes}  \eta _H $ and
 $ \epsilon _A  {\otimes} id  {\otimes}  \epsilon _A
   {\otimes} id$ respectively. By applying
 $\epsilon _A  {\otimes}  id  {\otimes} id
  {\otimes} \epsilon _A$  to relation (1),  we obtain $(B5).$
\begin{picture}(8,8)\put(0,0){\line(0,1){8}}\put(8,8){\line(0,-1){8}}\put(0,0){\line(1,0){8}}\put(8,8){\line(-1,0){8}}\end{picture}

\begin {Theorem} \label {2.1.5}   Let
$D =  A_\alpha ^\phi \bowtie _\beta ^\psi H $  be a bialgebra.
  If $A$ and $H$ are Hopf algebras with
antipodes $S_A$ and $S_H$ respectively, then $D$ is a Hopf algebra
with an antipode

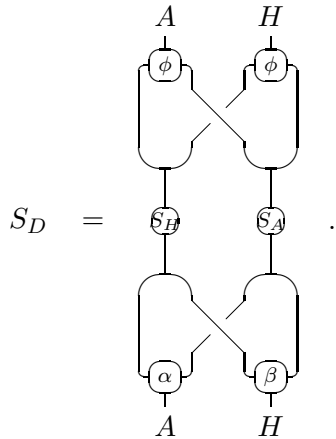
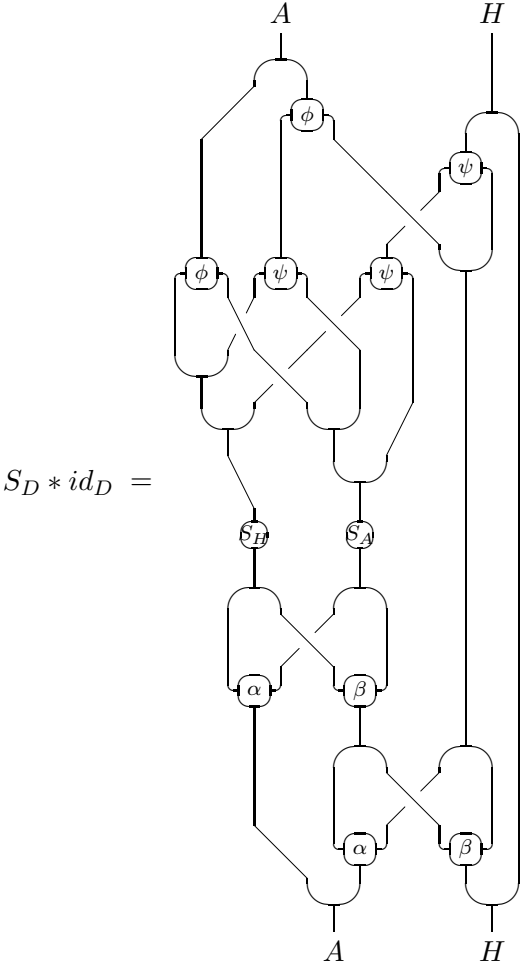
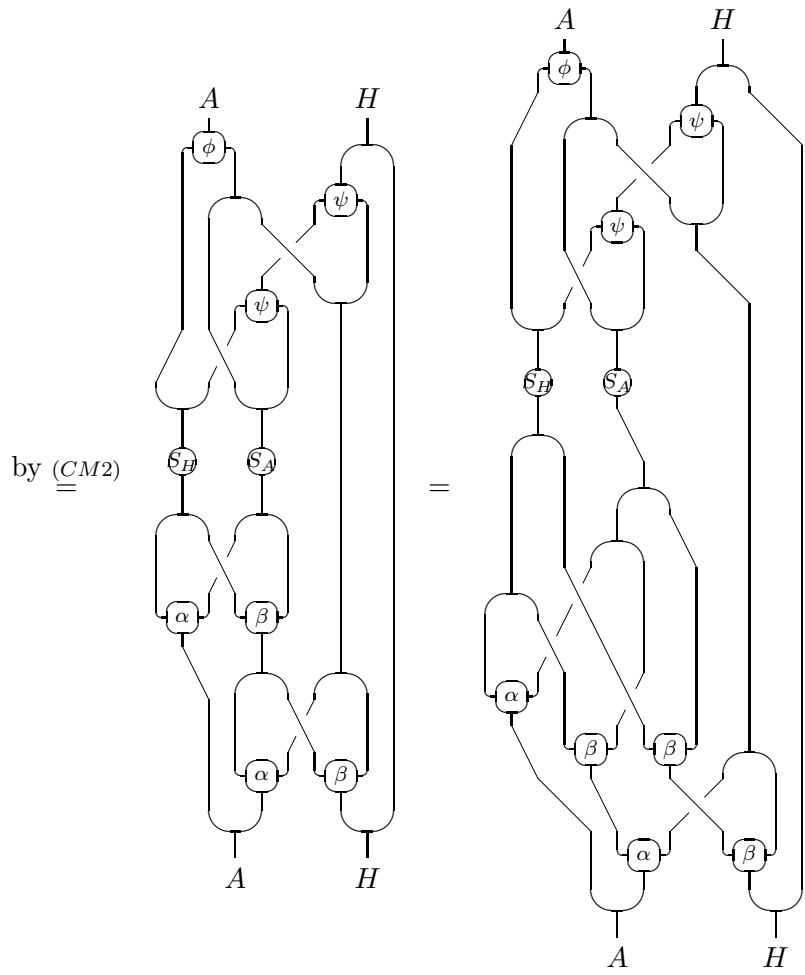
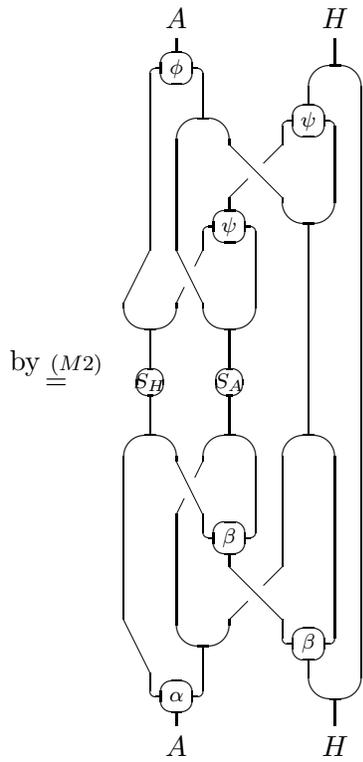
\[
S_{D} \step=\step
\begin{tangle}
\step[1]\object{A}\step[4]\object{H}\\
\td \phi\step[2]\td \phi\\
\id\step[2]\x\step[2]\id\\
\cu \step[2]\cu \\
\morph { S_{H}}\step[2]\morph {S_{A}}\\
\cd\step[2]\cd\\
\id\step[2]\x\step[2]\id\\
\tu \alpha\step[2]\tu \beta\\
\step[1]\object{A}\step[4]\object{H}
\end{tangle}\ \ \ .
\]

\end {Theorem}{\bf Proof.}
\[
\begin{tangle}
S_{D}*id_{D}
\end{tangle}
\;=\enspace
\begin{tangles}{clr}
\step[3]\object{A}\step[8]\object{H}\\
\step[2]\cd\step[7]\id\\
\step[2]\ne2\step[1]\td \phi\step[5]\cd\\
\step[1]\id\step[3]\id\step[2]\nw2\step[3]\td \psi\step[1]\id\\
\step[1]\id\step[3]\id\step[4]\x\step[2]\id\step[1]\id\\
\td \phi\step[1]\td \psi\step[2]\td \psi\step[1]\cu \step[1]\id\\
\id\step[2]\hx\step[2]\x\step[2]\id\step[2]\id\step[2]\id\\
\cu \step[1]\x\step[2]\id\step[2]\id\step[2]\id\step[2]\id\\
\step[1]\cu \step[2]\cu \step[1]\ne1\step[2]\id\step[2]\id\\
\step[2]\nw1\step[3]\cu \step[3]\id\step[2]\id\\
\step[2]\morph {S_{H}}\step[2]\morph {S_{A}}\step[3]\id\step[2]\id\\
\step[2]\cd\step[2]\cd\step[3]\id\step[2]\id\\
\step[2]\id\step[2]\x\step[2]\id\step[3]\id\step[2]\id\\
\step[2]\tu \alpha\step[2]\tu \beta\step[3]\id\step[2]\id\\
\step[3]\id\step[3]\cd\step[2]\cd\step[1]\id\\
\step[3]\id\step[3]\id\step[2]\x\step[2]\id\step[1]\id\\
\step[3]\nw2\step[2]\tu \alpha\step[2]\tu \beta\step[1]\id\\
\step[5]\cu \step[4]\cu \\
\step[5]\object{A}\step[6]\object{H}
\end{tangles}
\]

\[
\stackrel{ \hbox {by }  (CM2)}{=}\step
\begin{tangle}
\step[2]\object{A}\step[6]\object{H}\\
\step[1]\td \phi\step[4]\cd\\
\step[1]\id\step[1]\cd\step[2]\td \psi\step[1]\id\\
\step[1]\id\step[1]\id\step[2]\x\step[2]\id\step[1]\id\\
\step[1]\id\step[1]\id\step[1]\td \psi\step[1]\cu \step[1]\id\\
\ne1\step[1]\hx\step[2]\id\step[2]\id\step[2]\id\\
\cu \step[1]\cu \step[2]\id\step[2]\id\\
\morph {S_{H}}\step[1]\morph {S_{A}}\step[2]\id\step[2]\id\\
\cd\step[1]\cd\step[2]\id\step[2]\id\\
\id\step[2]\hx\step[2]\id\step[2]\id\step[2]\id\\
\tu \alpha\step[1]\tu \beta\step[2]\id\step[2]\id\\
\step[1]\nw1\step[1]\cd\step[1]\cd\step[1]\id\\
\step[2]\id\step[1]\id\step[2]\hx\step[2]\id\step[1]\id\\
\step[2]\id\step[1]\tu \alpha\step[1]\tu \beta\step[1]\id\\
\step[2]\cu \step[3]\cu \\
\step[3]\object{A}\step[5]\object{H}
\end{tangle}
\step = \step
\begin{tangle}
\step[3]\object{A}\step[6]\object{H}\\
\step[2]\td \phi\step[4]\cd\\
\step[1]\ne1\step[1]\cd\step[2]\td \psi\step[1]\nw2\\
\step[1]\id\step[2]\id\step[2]\x\step[2]\id\step[3]\id\\
\step[1]\id\step[2]\id\step[1]\td \psi\step[1]\cu \step[3]\id\\
\step[1]\id\step[2]\hx\step[2]\id\step[2]\nw2\step[3]\id\\
\step[1]\cu \step[1]\cu \step[4]\id\step[2]\id\\
\step[1]\morph {S_{H}}\step \morph {S_{A}}\step[4]\id\step[2]\id\\
\step[1]\cd\step[2]\nw1\step[4]\id\step[2]\id\\
\step[1]\id\step[2]\id\step[2]\cd\step[3]\id\step[2]\id\\
\step[1]\id\step[2]\id\step[1]\cd\step[1]\nw1\step[2]\id\step[2]\id\\
\cd\step[1]\hx\step[2]\id\step[2]\id\step[2]\id\step[2]\id\\
\id\step[2]\hx\step[1]\nw1\step[1]\id\step[2]\id\step[2]\id\step[2]\id\\
\tu \alpha\step[1]\id\step[2]\hx\step[2]\id\step[2]\id\step[2]\id\\
\step[1]\nw1\step[1]\tu \beta\step[1]\tu \beta\step[1]\cd\step[1]\id\\
\step[2]\nw2\step[1]\nw1\step[2]\x\step[2]\id\step[1]\id\\
\step[4]\id\step[1]\tu \alpha\step[2]\tu \beta\step[1]\id\\
\step[4]\cu \step[4]\cu \\
\step[5]\object{A}\step[6]\object{H}
\end{tangle}
\]

\[  \ \ \stackrel {\hbox {by } (M2)} { = } \ \
\begin{tangle}
\step[2]\object{A}\step[6]\object{H}\\
\step[1]\td \phi\step[4]\cd\\
\step[1]\id\step[1]\cd\step[2]\td \psi\step[1]\id\\
\step[1]\id\step[1]\id\step[2]\x\step[2]\id\step[1]\id\\
\step[1]\id\step[1]\id\step[1]\td \psi\step[1]\cu \step[1]\id\\
\ne1\step[1]\hx\step[2]\id\step[2]\id\step[2]\id\\
\cu \step[1]\cu \step[2]\id\step[2]\id\\
\morph {S_{H}}\step \morph {S_{A}}\step[2]\id\step[2]\id\\
\cd\step[1]\cd\step[1]\cd\step[1]\id\\
\id\step[2]\hx\step[2]\id\step[1]\id\step[2]\id\step[1]\id\\
\id\step[2]\id\step[1]\tu \beta\step[1]\id\step[2]\id\step[1]\id\\
\id\step[2]\id\step[2]\x\step[2]\id\step[1]\id\\
\nw1\step[1]\cu \step[2]\tu \beta\step[1]\id\\
\step[1]\tu \alpha\step[4]\cu \\
\step[2]\object{A}\step[6]\object{H}
\end{tangle}
\step\stackrel{\hbox {by Proposition } \ref {12.1.4}(i)}{=}\step
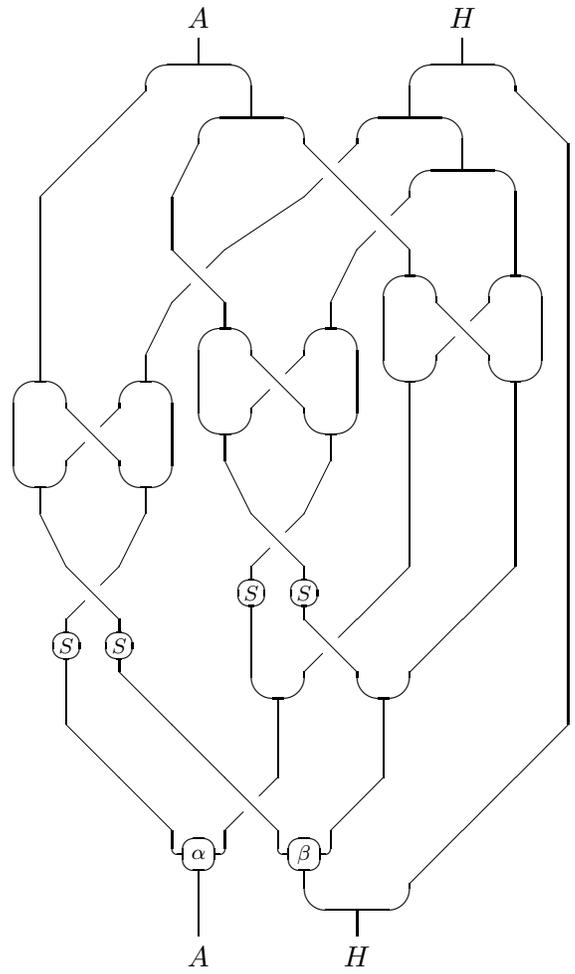
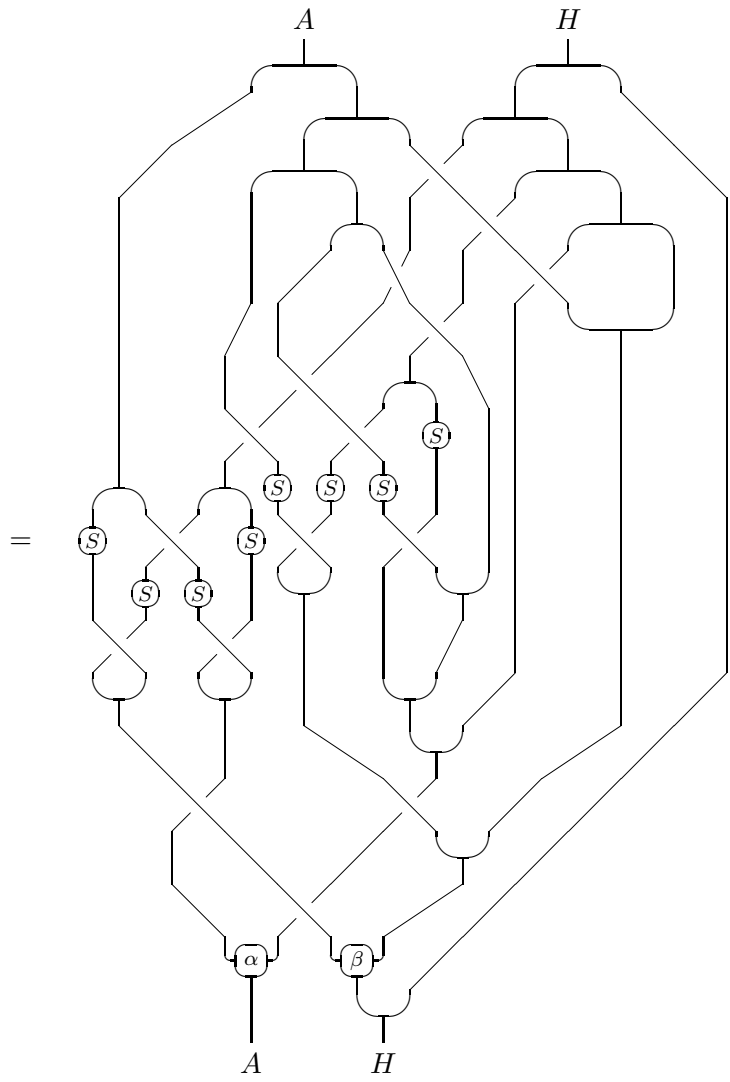
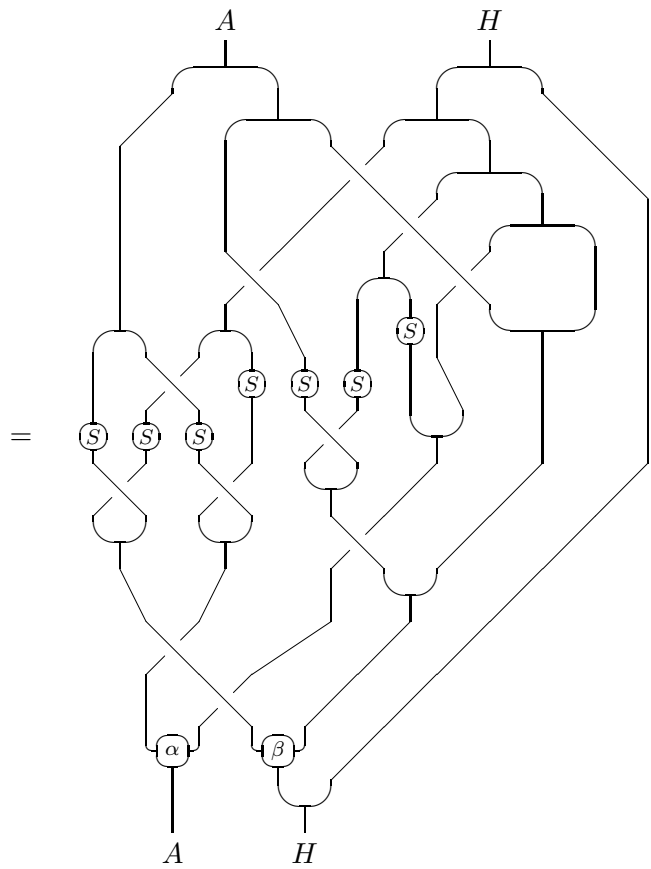
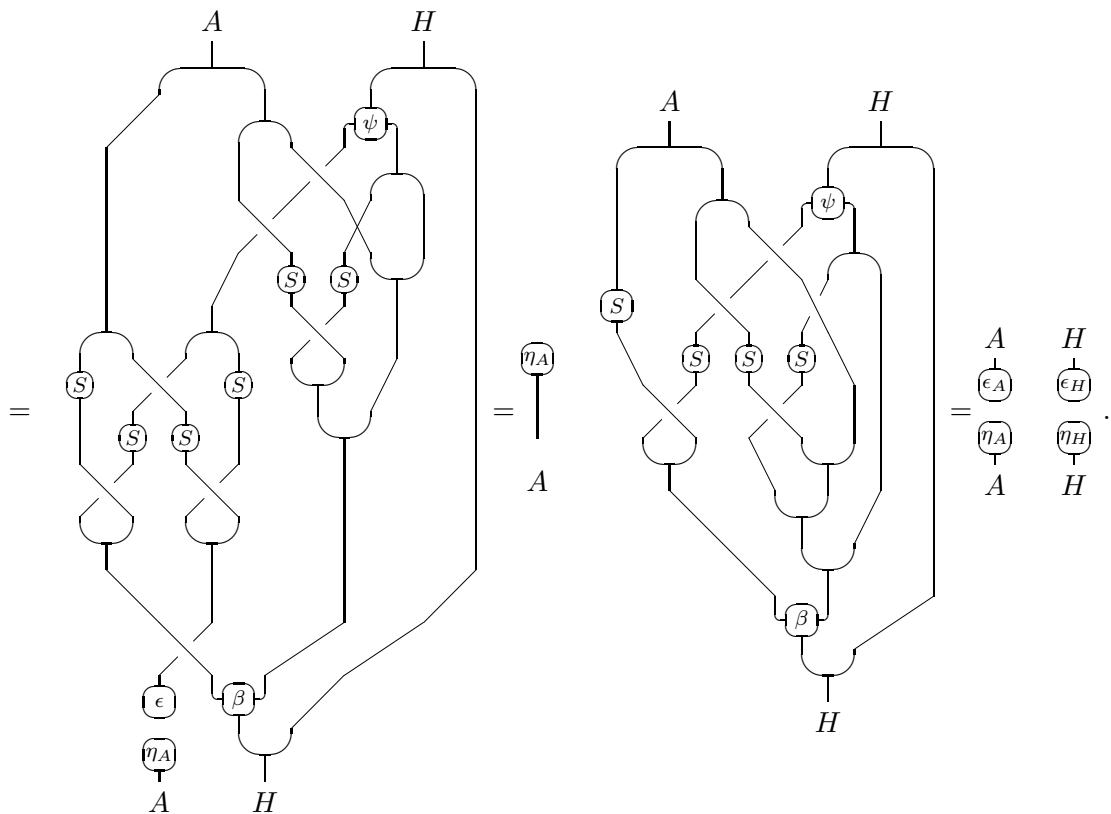
\begin{tangle}
\step[7]\object{A}\step[10]\object{H}\\
\step[5]\Cd\step[6]\Cd\\
\step[4]\ne2\step[2]\Cd\step[2]\Cd\step[2]\nw2\\
\step[2]\ne2\step[3]\ne1\step[4]\x\step[2]\Cd\step[2]\id\\
\step[1]\id\step[5]\id\step[4]\ne3\step[2]\x\step[4]\id\step[2]\id\\
\step[1]\id\step[5]\x\step[4]\ne1\step[1]\cd\step[2]\cd\step[1]\id\\
\step[1]\id\step[4]\ne1\step[1]\cd\step[2]\cd\step[1]\id\step[2]\x\step[2]\id\step[1]\id\\
\cd\step[2]\cd\step[1]\id\step[2]\x\step[2]\id\step[1]\cu \step[2]\cu \step[1]\id\\
\id\step[2]\x\step[2]\id\step[1]\cu \step[2]\cu \step[2]\id\step[4]\id\step[2]\id\\
\cu \step[2]\cu \step[2]\nw1\step[2]\ne1\step[3]\id\step[4]\id\step[2]\id\\
\step[1]\nw1\step[2]\ne1\step[4]\x\step[4]\id\step[4]\id\step[2]\id\\
\step[2]\x\step[4]\morph S\morph S\step[2]\ne2\step[3]\ne2\step[2]\id\\
\step[1]\morph S\morph S\step[4]\id\step[2]\x\step[3]\ne2\step[4]\id\\
\step[2]\id\step[2]\nw2\step[4]\cu \step[2]\cu \step[6]\id\\
\step[2]\nw2\step[3]\nw2\step[3]\id\step[4]\id\step[6]\ne2\\
\step[4]\nw2\step[3]\x\step[3]\ne2\step[4]\ne2\\
\step[6]\tu \alpha\step[2]\tu \beta\step[4]\ne2\\
\step[7]\id\step[4]\Cu\\
\step[7]\object{A}\step[6]\object{H}
\end{tangle}
\]

\[
\step=\step
\begin{tangle}
\step[9]\object{A}\step[10]\object{H}\\
\step[7]\Cd\step[6]\Cd\\
\step[6]\ne3\step[2]\Cd\step[2]\Cd\step[2]\nw2\\
\step[3]\ne2\step[3]\Cd\step[2]\x\step[2]\Cd\step[2]\nw2\\
\step[2]\id\step[5]\id\step[3]\cd\step[1]\id\step[2]\x\step[2]\Cd\step[2]\id\\
\step[2]\id\step[5]\id\step[2]\ne2\step[2]\hx\step[2]\id\step[2]\x\step[4]\id\step[2]\id\\
\step[2]\id\step[4]\ne1\step[1]\id\step[3]\ne2\step[1]\x\step[2]\id\step[2]\Cu\step[2]\id\\
\step[2]\id\step[4]\id\step[2]\x\step[2]\cd\step[1]\nw1\step[1]\id\step[4]\id\step[4]\id\\
\step[2]\id\step[4]\x\step[2]\x\step[1]\morph S\step[1]\id\step[1]\id\step[4]\id\step[4]\id\\
\step[1]\cd\step[2]\cd\morph S\morph S\morph S\step[1]\id\step[2]\id\step[1]\id\step[4]\id\step[4]\id\\
\morph S\step[1]\x\step[1]\morph S\x\step[2]\x\step[2]\id\step[1]\id\step[4]\id\step[4]\id\\
\step[1]\id\step[1]\morph S\morph S\step[1]\id\step[1]\cu \step[2]\id\step[2]\cu \step[1]\id\step[4]\id\step[4]\id\\
\step[1]\x\step[2]\x\step[2]\id\step[3]\id\step[2]\ne1\step[2]\id\step[4]\id\step[4]\id\\
\step[1]\cu \step[2]\cu \step[2]\id\step[3]\cu \step[2]\ne2\step[4]\id\step[3]\ne2\\
\step[1]\step[1]\nw2\step[3]\id\step[3]\nw3\step[3]\cu \step[5]\ne3\step[1]\ne2\\
\step[1]\step[3]\x\step[6]\x\step[3]\ne2\step[2]\ne2\\
\step[1]\step[3]\id\step[2]\nw2\step[4]\ne2\step[2]\cu \step[2]\ne2\\
\step[1]\step[3]\nw2\step[3]\x\step[4]\ne3\step[1]\ne2\\
\step[1]\step[5]\tu \alpha\step[2]\tu \beta\step[2]\ne2\\
\step[1]\step[6]\id\step[4]\cu \\
\step[1]\step[6]\object{A}\step[5]\object{H}\\
\end{tangle}
\]
\[
\step=\step
\begin{tangle}
\step[1]\step[5]\object{A}\step[10]\object{H}\\
\step[1]\step[3]\Cd\step[6]\Cd\\
\step[1]\step[2]\ne2\step[2]\Cd\step[2]\Cd\step[2]\nw2\\
\step[1]\step[1]\id\step[4]\id\step[4]\x\step[2]\Cd\step[2]\nw2\\
\step[1]\step[1]\id\step[4]\id\step[3]\ne2\step[2]\x\step[2]\Cd\step[2]\id\\
\step[1]\step[1]\id\step[4]\x\step[3]\cd\step[1]\x\step[4]\id\step[2]\id\\
\step[1]\cd\step[2]\cd\step[1]\nw1\step[2]\id\step[1]\morph S \id\step[2]\Cu\step[2]\id\\
\step[1]\id\step[2]\x\step[1]\morph S\morph S\morph S\step[1]\id\step[1]\nw1\step[3]\id\step[4]\id\\
\morph S\morph S\morph S\step[1]\id\step[2]\x\step[2]\cu \step[3]\id\step[4]\id\\
\step[1]\x\step[2]\x\step[2]\cu \step[2]\ne2\step[3]\ne2\step[3]\ne2\\
\step[1]\cu \step[2]\cu \step[3]\x\step[3]\ne2\step[3]\ne2\\
\step[1]\step[1]\nw1\step[2]\ne1\step[4]\id\step[2]\cu \step[3]\ne2\\
\step[1]\step[2]\x\step[4]\ne3\step[2]\ne2\step[2]\ne2\\
\step[1]\step[2]\id\step[2]\x\step[3]\ne2\step[2]\ne2\\
\step[1]\step[2]\tu \alpha\step[2]\tu \beta\step[2]\ne2\\
\step[1]\step[3]\id\step[4]\cu \\
\step[1]\step[3]\object{A}\step[5]\object{H}
\end{tangle}
\]

\[
\;=\enspace
\begin{tangle}
\step[6]\object{A}\step[8]\object{H}\\
\step[4]\Cd\step[4]\Cd\\
\step[3]\ne2\step[3]\cd\step[2]\td \psi\step[3]\id\\
\step[2]\id\step[5]\id\step[2]\x\step[1]\cd\step[2]\id\\
\step[2]\id\step[5]\x\step[2]\X\step[2]\id\step[2]\id\\
\step[2]\id\step[4]\ne1\step[1]\morph S\morph S\cu \step[2]\id\\
\step[1]\cd\step[2]\cd\step[2]\x\step[2]\id\step[3]\id\\
\morph S\step[1]\x\step[1]\morph S\step[1]\cu \step[1]\ne1\step[3]\id\\
\step[1]\id\step[1]\morph S\morph S\step[1]\id\step[3]\cu \step[4]\id\\
\step[1]\x\step[2]\x\step[4]\id\step[5]\id\\
\step[1]\cu \step[2]\cu \step[4]\id\step[5]\id\\
\step[2]\nw2\step[3]\id\step[5]\id\step[4]\ne2\\
\step[4]\x\step[4]\ne3\step[2]\ne3\\
\step[4]\QQ { \epsilon}\step[2]\tu \beta\step[2]\ne2\\
\step[4]\Q {\eta_{A}}\step[3]\cu \\
\step[4]\object{A}\step[4]\object{H}
\end{tangle}
\;=\enspace
\begin{tangles}{clr}
\step[3]\object{A}\step[8]\object{H}\\
\step[3]\Cd\step[4]\Cd\\
\step[3]\id\step[3]\cd\step[2]\td \psi\step[3]\id\\
\step[3]\id\step[3]\id\step[2]\x\step[1]\cd\step[2]\id\\
\step[3]\O S\step[3]\x\step[2]\hx\step[2]\id\step[2]\id\\
\Q {\eta_{A}}\step[3]\nw1\step[1]\morph S\morph S\morph S\nw1\step[1]\id\step[2]\id\\
\id\step[4]\x\step[2]\x\step[2]\id\step[1]\id\step[2]\id\\
\object{A}\step[4]\cu \step[2]\nw1\step[1]\cu \step[1]\id\step[2]\id\\
\step[5]\nw2\step[3]\cu \step[1]\ne1\step[2]\id\\
\step[7]\nw2\step[2]\cu \step[3]\id\\
\step[9]\tu \beta\step[3]\ne3\\
\step[7]\cu \\
\step[7]\object{H}
\end{tangles}
\;=\enspace
\begin{tangle}
\object{A}\step[3]\object{H}\\
\QQ {\epsilon_{A}}\step[3]\QQ {\epsilon_{H}}\\
\Q {\eta_{A}}\step[3]\Q {\eta_{H}}\\
\object{A}\step[3]\object{H}
\end{tangle}\ \ \ .
\]
Similarly, we can show that  $id _D * S_D = \epsilon _D \eta _D.$
$\Box$

\begin {Proposition} \label {2.1.6}
 Let $D = A \bowtie H$ be a bialgebra ,  and let $A$ and $H$ be  Hopf algebras  with antipodes $S_A $ and $S_H$ respectively.
 Then

\[
\begin{tangle}
\step[2]\object{D}\\
\step[2]\id\\
\obox 3 { (S_D) ^{2}} \\
\step[2]\id\\
\step[2]\object{D}
\end{tangle}
\step=\step
\begin{tangle}
\step[1]\object{A}\step[5]\object{H}\\
\step[1]\id \step[5]\id \\
\obox 3 { (S_D) ^{2}}\step[0.6]\obox 3 { (S_D) ^{2}}\\
\step[1]\nw2\step[3]\id\\
\step[3]\x\\
\step[3]\x\\
\step[3]\object{A}\step[2]\object{H}
\end{tangle}\ \ \ .
\]

\end {Proposition}

{\bf Proof.}

\[
\begin{tangle}
\step[1]\object{D}\\
\step[1]\id\\
\obox3 { (S_D)^2}\\
\step[1]\id\\
\step[1]\object{D}
\end{tangle}
\step\ \ = \ \ \step
\begin{tangle}
\object{A\otimes \eta_{H}}\step[6]\object{\eta_{H}\otimes H}\\
\step [2]\tu {m_D}\\
\step[2]\obox 3 {(S_D)^2}\\
\step[3]\id\\
\step[3]\object{D}
\end{tangle}
\step\ \ =\step
\begin{tangle}
\step \object{A\otimes \eta_{H}}\step[6]\object{\eta_{H}\otimes H}\\

\step [1]\id\step [5]\id\\

\obox 3 {(S_D)^2}\step[1.5]\obox 3 {(S_D)^2}\\
\step[1]\id\step[4]\ne3\\
\step[1]\x\\
\step[1]\x\\
\step[1]\cu _{D}\\
\step[2]\object{D}
\end{tangle}
\step=\step
\begin{tangle}
\step[1]\object{A}\step[5]\object{H}\\

\step [1]\id\step [5]\id\\
\obox 3 {(S_A)^2}\step[1.5]\obox 3 {(S_H)^2}\\
\step[1]\id\step[4]\ne3\\
\step[1]\x\\
\step[1]\x\\
\step[1]\object{A}\step[2]\object{H}
\end{tangle}\ \ \ \Box .
\]

We can easily check Lemma \ref {2.1.8}, Corollary \ref {2.1.9} and
Corollary \ref {2.1.10} by Theorem \ref {2.1.4} and by simple
computation.

\begin {Lemma}\label {2.1.8} (i) If $\phi$  and $\psi$  are trivial, then
$(B1)$--$(B5)$  hold iff  $(M4)$  holds;

 (ii) If $\alpha$  and $\beta $  are trivial, then
$(B1)$--$(B5)$  hold iff  $(CM4)$  holds.

\end {Lemma}

\begin {Corollary}\label {2.1.9}
(i) $A \bowtie H$ is a bialgebra iff $(M1)$--$(M4)$  hold;

 (ii)  $ A ^\phi \bowtie  ^\psi H $ is a bialgebra iff
$(CM1)$--$(CM4)$  hold.
\end {Corollary}

In fact, \cite  {Ma94b} already  contains the essence of Corollary
\ref {2.1.9}.

\begin {Corollary}\label {2.1.10}
(i) If $\beta $  and $\psi$  are trivial, then biproduct  $A^\phi
_\alpha \bowtie H $ is a bialgebra iff $(A, \alpha )$ is an
$H$-module algebra and $(A, \phi )$ is an $H$- comodule coalgebra,
(B1) and (B6)   hold;

 (ii) If $\alpha $ and $\phi $ are trivial, then
biproduct $ A \bowtie _\beta ^\psi H $ is a bialgebra iff $(H, \beta
)$ is an $A$-module algebra and $(H, \psi )$ is an $A$- comodule
coalgebra, (B3) and (B6')  hold.
\end {Corollary}

\begin {Corollary}\label {2.1.11}  (i) If $\beta $ and $\phi $  are trivial
, then bicrossproduct  $A _\alpha \bowtie ^\psi H $ is a bialgebra
iff $(A, \alpha )$ is an $H$-module algebra and $(H, \psi )$ is an
$A$- comodule coalgebra,(B2), (B7) and (CB4)  hold;

 (ii) If $\alpha $ and $\psi $  are trivial
then bicrossproduct $A ^\phi \bowtie _\beta H $ is a bialgebra iff
$(H, \beta )$ is an $A$-module algebra and $(A, \phi )$ is an $H$-
comodule coalgebra, (B4), (B7') and (CB2)  hold.

\end {Corollary}
{\bf Proof.} (i) The necessity is clear. For the sufficiency,  we
only show that (B5) holds since the others can similarly be shown.

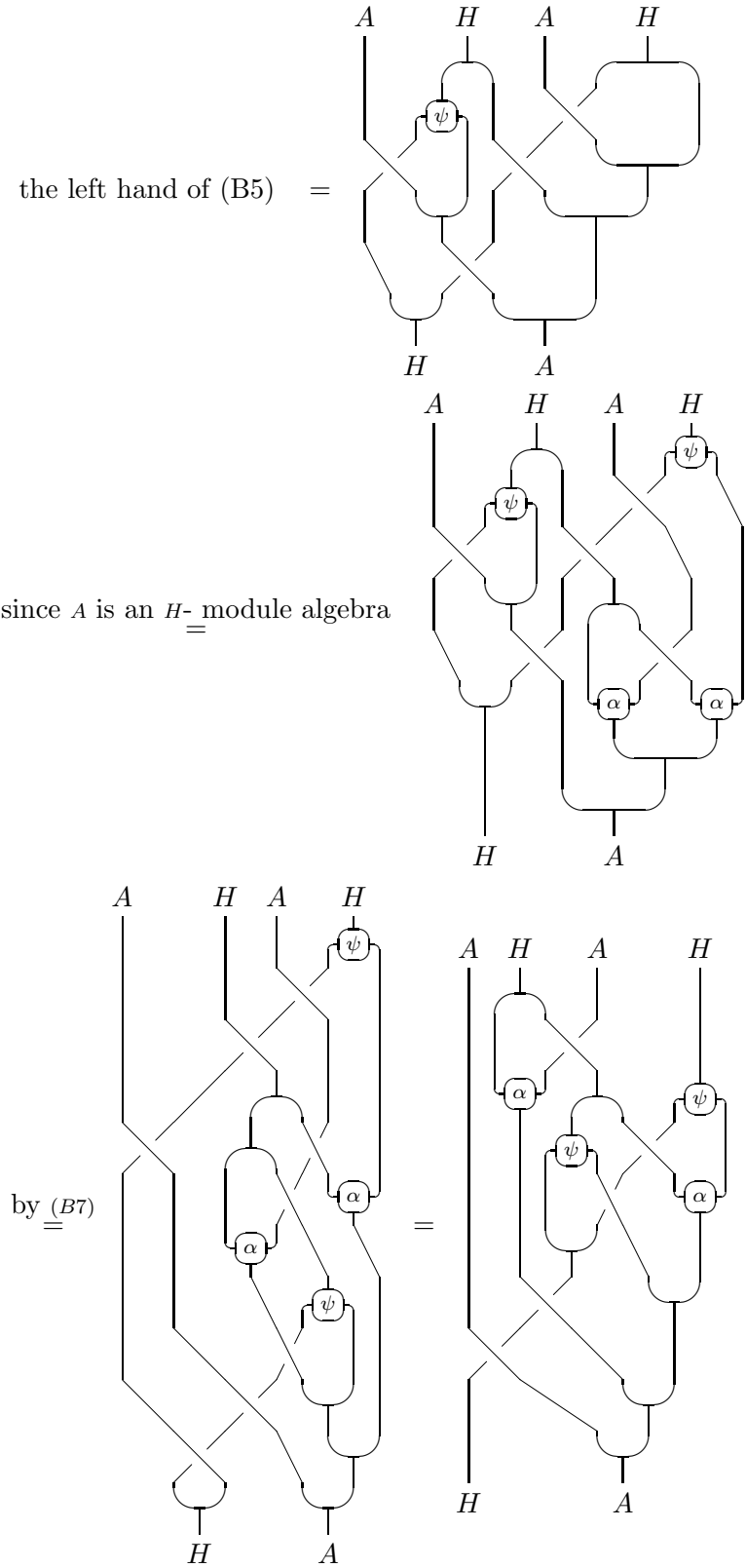
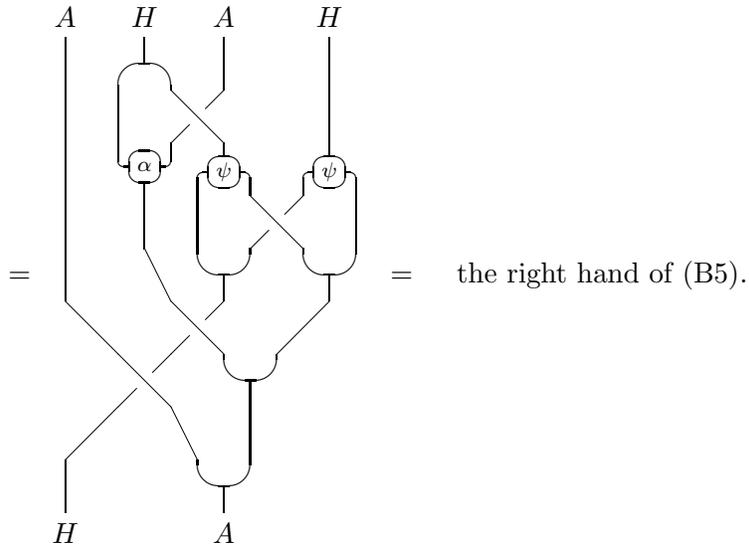
\[
\hbox { the left hand of (B5)} \step=\step
\begin{tangle}
\object{A}\step[4]\object{H}\step[3]\object{A}\step[4]\object{H}\\
\id\step[3]\cd\step[2]\id\step[2]\Cd\\
\id\step[2]\td \psi\step[1]\id\step[2]\x\step[4]\id\\
\x\step[2]\id\step[1]\x\step[2]\Cu\\
\id\step[2]\cu \step[1]\id\step[2]\Cu\\
\nw1\step[2]\x\step[4]\id\\
\step[1]\cu \step[2]\Cu\\
\step[2]\object{H}\step[5]\object{A}
\end{tangle}
\]
\[\step[1] \stackrel {\hbox { since } A \hbox { is an } H \hbox {- module algebra }}{=}\step[1]
\begin{tangle}
\object{A}\step[4]\object{H}\step[3]\object{A}\step[3]\object{H}\\
\id\step[3]\cd\step[2]\id\step[2]\td \psi\\
\id\step[2]\td \psi\step[1]\id\step[2]\x\step[2]\nw1\\
\x\step[2]\id\step[1]\x\step[2]\nw1\step[2]\id\\
\id\step[2]\cu \step[1]\id\step[1]\cd \step[2]\id\step[2]\id\\
\nw1\step[2]\x\step[1]\id\step[2]\x\step[2]\id\\
\step[1]\cu \step[2]\id\step[1]\tu \alpha\step[2]\tu \alpha\\
\step[2]\id\step[3]\id\step[2]\Cu\\
\step[2]\id\step[3]\Cu\\
\step[2]\object{H}\step[5]\object{A}
\end{tangle}
\]

\[ \ \ \
\stackrel {\hbox {by } (B7)  }{=} \ \
\begin{tangle}
\object{A}\step[4]\object{H}\step[2]\object{A}\step[3]\object{H}\\
\id\step[4]\id\step[2]\id\step[2]\td \psi\\
\id\step[4]\id\step[2]\x\step[2]\id\\
\id\step[4]\x\step[2]\id\step[2]\id\\
\id\step[3]\ne2\step[1]\cd\step[1]\id\step[2]\id\\
\x\step[2]\cd\step[1]\hx\step[2]\id\\
\id\step[2]\id\step[2]\id\step[2]\hx\step[1]\tu \alpha\\
\id\step[2]\id\step[2]\tu \alpha\step[1]\nw1\step[1]\nw1\\
\id\step[2]\id\step[3]\nw1\step[1]\td \psi\step[1]\id\\
\id\step[2]\nw2\step[3]\hx\step[2]\id\step[1]\id\\
\nw2\step[3]\x\step[1]\cu \step[1]\id\\
\step[2]\x\step[2]\nw1\step[1]\cu \\
\step[2]\cu \step[3]\cu \\
\step[3]\object{H}\step[5]\object{A}
\end{tangle}
\step[1] = \step[1]
\begin{tangle}
\object{A}\step[2]\object{H}\step[3]\object{A}\step[4]\object{H}\\
\id\step[1]\cd\step[2]\id\step[4]\id\\
\id\step[1]\id\step[2]\x\step[4]\id\\
\id\step[1]\tu \alpha\step[1]\cd\step[2]\td \psi \\
\id\step[2]\id\step[1]\td \psi\step[1]\x\step[2]\id\\
\id\step[2]\id\step[1]\id\step[2]\hx\step[2]\tu \alpha\\
\id\step[2]\id\step[1]\cu \step[1]\nw1\step[2]\id\\
\id\step[2]\x\step[3]\cu \\
\x\step[2]\nw2\step[3]\id\\
\id\step[2]\nw3\step[3]\cu \\
\id\step[5]\cu \\
\object{H}\step[6]\object{A}
\end{tangle}
\]

\[
\step[1]=\step[1]
\begin{tangle}
\object{A}\step[3]\object{H}\step[3]\object{A}\step[4]\object{H}\\
\id\step[2]\cd\step[2]\id\step[4]\id\\
\id\step[2]\id\step[2]\x\step[4]\id\\
\id\step[2]\tu \alpha\step[1]\td \psi\step[2]\td \psi\\
\id\step[3]\id\step[2]\id\step[2]\x\step[2]\id\\
\id\step[3]\nw1\step[1]\cu \step[2]\cu \\
\nw2\step[3]\x\step[3]\ne2\\
\step[2]\x\step[2]\cu \\
\step[1]\ne2\step[2]\nw1\step[2]\id\\
\id\step[5]\cu \\
\object{H}\step[6]\object{A}
\end{tangle}
\step[1]=\step[1] \hbox { the right hand of (B5)}.
\]
Thus (B5) holds.

(ii) It is a dual case of Part (i). $\Box$

\begin {Lemma}\label {2.1.10'}  Let $H$ be a Hopf algebra,
and
 $ A$ and $H$ have left duality $A^*$ and $H^*$ respectively.

(i)  If $(A, \alpha ) $ is
 a left  $H$-module coalgebra, then
 $(A^*, \alpha ^*)$  is a right $H^*$-comodule  algebra;

(ii)  If $(A, \phi ) $ is
 a left  $H$-comodule algebra, then
 $(A^*, \phi ^*)$  is a right $H^*$-module  coalgebra.
\end {Lemma}
{\bf Proof.}  We only show part (i), since part (ii) is the dual
case of part (i).
\[
\begin{tangle}
\step[2]\object{A^*}\\
\step\td {\alpha^*}\\
\td {\alpha^*}\step\id\\
\object{A^*}\step[2]\object{H^*}\step[1]\object{H^*}\\
\end{tangle}\step=\step
\begin{tangle}
\step[1]\object{A^*}\\
\td {\alpha^*}\\
\id\step\td {\Delta^*}\\
\object{A^*}\step[1]\object{H^*}\step[2]\object{H^*}\\
\end{tangle}\step\step[2]\hbox {.....(1)}
\]

\[
\hbox {The left side of (1)}\step=\step
\begin{tangle}
\object{A^*}\\
\id\step[2]\step\coev\\
\id\step[2]\ne2\coev\d\\
\id\step\tu \alpha\step[2]\id\step\d\\
\ev\step[2]\ne2\step[-1]\coev\step\d\\
\step[2]\ne2\ne2\coev\d\step\id\\
\step\id\step\tu \alpha\step[2]\id\step\id\step\id\\
\step\ev\step[2]\step\id\step\id\step\id\\
\step[6]\object{A^*}\step[1]\object{H^*}\step[1]\object{H^*}\\
\end{tangle}
\step=\step
\begin{tangle}
\object{A^*}\\
\id\step[2]\step\coev\\
\id\step[2]\ne2\step\coev\se2\\
\id\step\id\step[2]\ne2\coev\d\d\\
\id\step\id\step\tu \alpha\step[2]\id\step\id\step\id\\
\id\step\tu\alpha\step[2]\step\id\step\id\step\id\\
\ev\step[2]\step[2]\id\step\id\step\id\\
\step[6]\object{A^*}\step[1]\object{H^*}\step[1]\object{H^*}\\
\end{tangle}
\]
 and

\[
\hbox {The right side of (1)}\step=\step
\begin{tangle}
\object{A^*}\\
\id\step[2]\step\coev\step[2]\step[2]\coev\\
\id\step[2]\ne2\coev\d\step[2]\ne2\coev\d\\
\id\step\tu
\alpha\step[2]\id\step\id\step\tu\alpha\step[2]\id\step\id\\
\ev\step[2]\step\id\step\ev\step[2]\step\id\step\id\\
\step[5]\object{A^*}\step[6]\object{H^*}\step[1]\object{H^*}\\
\end{tangle}
\step=\step
\begin{tangle}
\object{A^*}\\
\id\step[2]\step[2]\coev\\
\id\step[2]\step\ne2\coev\d\\
\id\step[2]\tu \alpha\step[2]\id\step\id\\
\id\step[2]\ne2\coev\step\id\step\id\\
\id\step\tu\alpha\step[2]\id\step\id\step\id\\
\ev\step\step[2]\id\step\id\step\id\\
\step[5]\object{A^*}\step[1]\object{H^*}\step[1]\object{H^*}\\
\end{tangle}\]
$ \ \ \stackrel {\hbox {since } (A, \alpha ) \hbox {is an } H \hbox
{-module }} {=} \hbox {the left hand of (1)}. \ \ $ Thus (1) holds.

It is clear that
\[
\begin{tangle}
\step\object{A^*}\\
\td {\alpha^*}\\
\id\step[2]\QQ {\epsilon}\\
\object{A^*}\\
\end{tangle}\step=\step
\begin{tangle}
\object{A^*}\\
\id\\
 \object{A^*}\\
\end{tangle}\ \ \ .
\] Thus $(A^*, \alpha ^*)$ is a right $H^*$-comodule.

Now we show that $ \step \begin{tangle}
\object{A^*}\step[2]\object{A^*}\\
\tu {\Delta^*}\\
\td {\alpha^*}\\
 \object{A^*}\step[2]\object{H^*}
\end{tangle}=
\begin{tangle}
\step\object{A^*}\step[4]\object{A^*}\\
\td {\alpha^*}\step[2]\td {\alpha^*}\\
\id\step[2]\x\step[2]\id\\
\tu {\Delta^*}\step[2]\tu {\Delta^*}\\
\step \object{A^*}\step[4]\object{H^*}
\end{tangle}\step\hbox {......(2)}$
\medskip

\[
\hbox {The left side of (2)} \step=\step
\begin{tangle}
\object{A^*}\step[1.5]\object{A^*}\\
\id\step\id\step[8]\coev\\
\id\step\id\step[7]\ne2\coev\d\\
\id\step\id\step[3]\coev\step\tu\alpha\step[2]\id\step\id\\
\id\step\id\step[2]\cd\step\ev\step[3]\id\step\id\\
\d\ev\step\dd\step[6]\id\step\id\\
\step\Ev\step[9]\id\step\id\\
\step[11]\object{A^*}\step[1.5]\object{H^*}
\end{tangle}\step=\step
\begin{tangle}
\object{A^*}\step[1.5]\object{A^*}\\
\id\step\id\step[4]\coev\\
 \id\step\id\step[3]\ne2\coev\d\\
  \id\step\id\step[2]\tu\alpha\step[2]\id\step\id\\
\id\step\id\step[2]\cd\step[2]\id\step\id\\
\d\ev\step\dd\step[2]\id\step\id\\
\step\Ev\step[5]\id\step\id\\
\step[7]\object{A^*}\step\hstep\object{H^*}
\end{tangle}
\]
\[=\step
\begin{tangle}
\object{A^*}\step[1.5]\object{A^*}\\
\id\step\id\step[6]\coev\\
\id\step\id\step[5]\ne4\coev\d\\
\id\step\id\step\cd\step[2]\cd\step\id\step\id\\
\id\step\id\step\id\step[2]\x\step[2]\id\step\id\step\id\\
\id\step\id\step\tu\alpha\step[2]\tu\alpha\step\id\step\id\\
\d\ev\step[3]\ne3\step[2]\id\step\id\\
\step\Ev\step[7]\id\step\id\\
\step[9]\object{A^*}\step\hstep\object{H^*}
\end{tangle}\stackrel{ \hbox {by Lem.\ref {2.1.9'}}}{=} \step\hbox {the right side of (2) .}
\] Thus (2) holds. Obviously, $\eta _{A^*}$ is a comodule morphism.

Consequently, $(A^*, \alpha ^*)$ is a right $H^*$-comodule algebra.
$\Box$

\begin {Proposition} \label {2.1.11'}

 Let $A$ and $H$ have left dualities in braided tensor category  $({\cal C }, C)$.
  If $D = A ^{\varphi} _{\alpha}  {\bowtie }
^{\psi}_{\beta} H $  is a    bialgebra or    Hopf algebra, then
             $$ (A ^{\varphi} _{\alpha}  {\bowtie }
^{\psi}_{\beta} H)^* \cong { H^*} _{\psi ^*}^{\beta ^*} {\bowtie }
_{\varphi^* }^{\alpha  ^*} A^* \hbox { \ \ \ \ \ \  as bialgebras or
Hopf algebras }.$$

\end {Proposition}

{\bf Proof. }               By \cite [Proposition 2.4] {Ma95a} or
Proposition \ref {12.1.4} (ii), $D^*, H^*$ and $A^*$ are Hopf
algebras. By Lemma \ref {2.1.10'}, we have that
 $(A^*, \alpha ^* )$ is a right $H^*$-comodule algebra,
 $(H^*, \beta ^*)$ is a left $A^*$-comodule algebra,
 $(A^*, \phi ^* )$ is a right $H^*$-module coalgebra and
 $(H^*, \psi ^* )$ is a left $H^*$-module coalgebra.
Let $\bar D = { H^*} _{\psi ^*}^{\beta ^*} {\bowtie } _{\varphi^*
}^{\alpha  ^*} A^* .$ Obviously, the unit and counit of $\bar D$
 are the same of $D^*$.

 Now we show that their multiplication are the same.

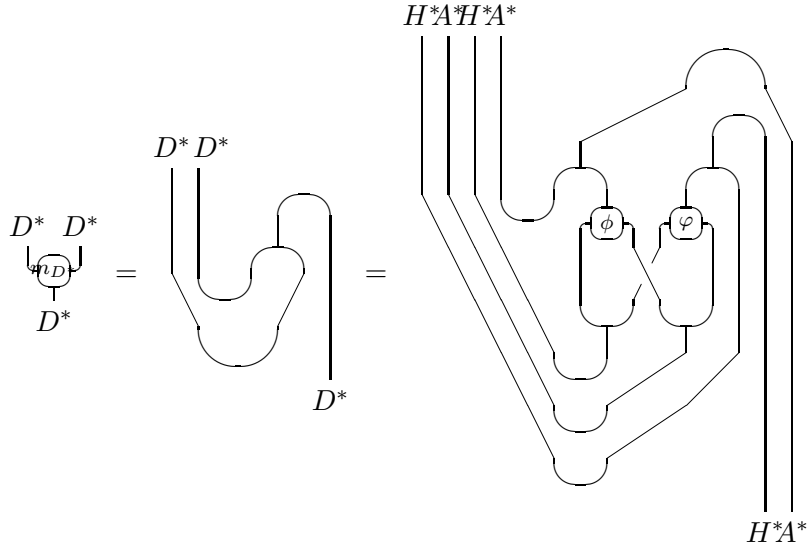
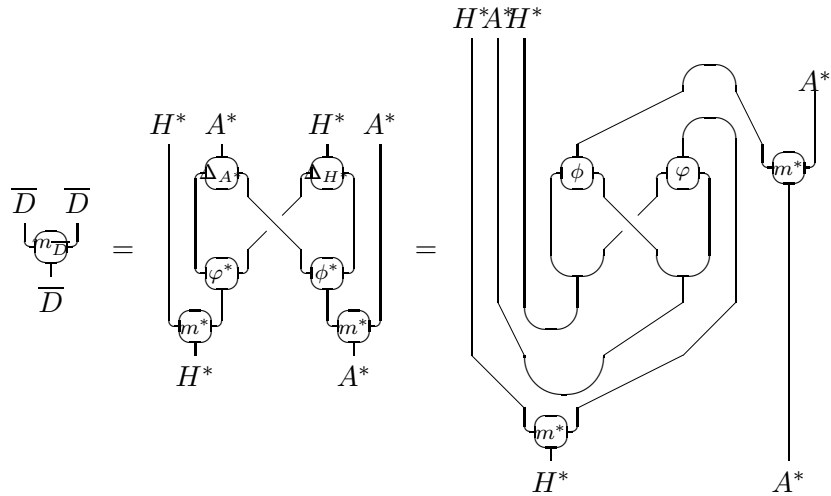
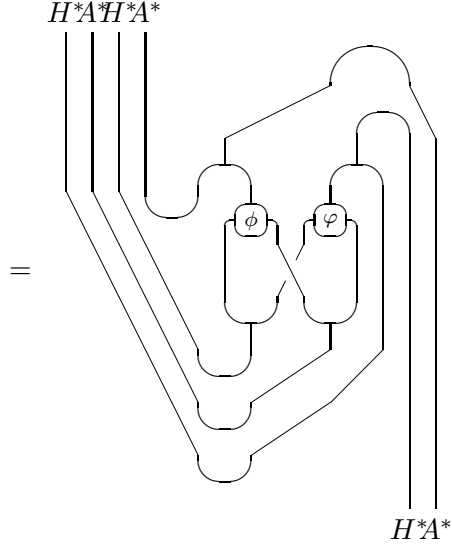
\[
\begin{tangle}
\object{D^*}\step[2]\object{D^*}\\
\tu {m_{D^*}}\\
\step[1]\object{D^*}\\
\end{tangle}
\step=\step
\begin{tangle}
\object{D^*}\step[1.5]\object{D^*}\\
 \id\step\id\step[3]\coev\\
\id\step\id\step[2]\cd\step[1]\id\\
\d\ev\step\dd\step[1]\id\\
\step\Ev\step[4]\id\\
\step[6]\object{D^*}
\end{tangle}
\step=\step
 \begin{tangle}
\object{H^*}\step[1]\object{A^*}\step[1]\object{H^*}\step[1]\object{A^*}\\
\id\step\id\step\id\step\id\step[7]\Coev\\
\id\step\id\step\id\step\id\step[6]\ne4\step\coev\d\\\
\id\step\id\step\id\step\id\step\step\cd\step[3]\cd\step\id\step[1]\id\\
\d\d\d\ev\step\td\phi \step\td\varphi\step\id\step[1]\id\step[1]\id\\
\step\d\d\d\step[2]\id\step[2]\hx\step[2]\id\step\id\step\id\step\id\\
\step\step\d\d\d\step\cu\step\cu\step\id\step\id\step\id\\
\step\step\step\d\d\ev\step\step\ne3\step\ne2\step\id\step\id\\
\step[4]\d \ev\step[2]\ne3\step[2]\step\id\step\id\\
\step[5]\ev\step[5]\step\id\step\id\\
 \step[13]\object{H^*}\step[1]\object{A^*}
\end{tangle}
\]

and
\[
\begin{tangle}
\object{\overline{D}}\step[2]\object{\overline{D}}\\
\tu {m_{\overline{D}}}\\
\step[1]\object{\overline{D}}\\
\end{tangle}
\step=\step
\begin{tangle}
\object{H^*}\step\step\object{A^*}\Step\Step\object{H^*}\step\step\object{A^*}\\
 \id \step \td {\Delta_{A^*}} \step [2] \td {\Delta_{H^*}}\step [1]\id  \\
  \id \step \id \step [2] \x \step [2]\id \step \id \\
   \id \step [1]\tu {\varphi^*} \step [2]\tu {\phi^*} \step [1] \id\\
    \tu {m^*}\step[4] \tu {m^*} \\
\step\object{H^*}\Step\Step\Step\object{A^*}\\
 \end{tangle}
\step=\step
 \begin{tangle}
\object{H^*}\step[1]\object{A^*}\step[1]\object{H^*}\\
\id\step\id\step\id\step[6]\coev\step[3]\object{A^*}\\
\id\step\id\step\id\step[5]\ne4\coev\d\step[2]\id\\
\id\step\id\step\id\step\td\phi\step[2]\td\varphi\step\id\step[1]\tu {m^*}\\
\id\step\id\step\id\step\id\step[2]\x\step[2]\id\step\id\step[2]\id\\
\id\step\id\step\id\step\cu\step[2]\cu\step\id\step[2]\id\\
\id\step\d\ev\step[3]\ne3\step\ne2\step[2]\id\\
\nw2\step\Ev\step[4]\ne4\step[2]\step[2]\id\\
\step[2]\tu {m^*}\step[6]\step[2]\id\\
 \step[3]\object{H^*}\step[9]\object{A^*}
\end{tangle}
\]

\[
\step=\step
 \begin{tangle}
\object{H^*}\step[1]\object{A^*}\step[1]\object{H^*}\step[1]\object{A^*}\\
\id\step\id\step\id\step\id\step[7]\Coev\\
\id\step\id\step\id\step\id\step[6]\ne4\step\coev\d\\\
\id\step\id\step\id\step\id\step\step\cd\step[3]\cd\step\id\step[1]\id\\
\d\d\d\ev\step\td\phi \step\td\varphi\step\id\step[1]\id\step[1]\id\\
\step\d\d\d\step[2]\id\step[2]\hx\step[2]\id\step\id\step\id\step\id\\
\step\step\d\d\d\step\cu\step\cu\step\id\step\id\step\id\\
\step\step\step\d\d\ev\step\step\ne3\step\ne2\step\id\step\id\\
\step[4]\d \ev\step[2]\ne3\step[2]\step\id\step\id\\
\step[5]\ev\step[5]\step\id\step\id\\
 \step[13]\object{H^*}\step[1]\object{A^*}
\end{tangle}\ \ \ .
\]
 Thus the two multiplications are the same. Dually, we have that their
 comultiplications are the same. Similarly, we can show that their antipodes
 are the same. Consequently, $D^*$  and $\bar D$  are Hopf algebra isomorphic.
 \begin{picture}(8,8)\put(0,0){\line(0,1){8}}\put(8,8){\line(0,-1){8}}\put(0,0){\line(1,0){8}}\put(8,8){\line(-1,0){8}}\end{picture}

We can immediately get the following by Proposition \ref {2.1.11'}

\begin {Corollary} \label {2.1.12}
Let $A$ and $H$ have left dualities in braided tensor category
$({\cal C }, C)$.

 (i)  If $D = A  _{\alpha}  {\bowtie }
_{\beta} H $  is a  bialgebra or      Hopf algebra, then
             $$ (A  _{\alpha}  {\bowtie }
_{\beta} H)^* \cong { H^*} ^{\beta ^*} {\bowtie }^{\alpha  ^*} A^*
\hbox { \ \ \ \ \ \
  as bialgebras or Hopf algebras };$$

(ii)
  If $D = A ^{\varphi}   {\bowtie }
^{\psi} H $  is a  bialgebra or      Hopf algebra, then
             $$ (A ^{\varphi}  {\bowtie }
^{\psi} H)^* \cong { H^*} _{\psi ^*} {\bowtie } _{\varphi^* } A^*
\hbox { \ \ \ \ \ \  as bialgebras or  Hopf algebras }.$$

\end {Corollary}

Now, we give an example of double cross coproduct in a strictly
braided tensor category.
\begin {Example} \label {2.1.13}  In this example, we work in
 the braided tensor category ${\cal V}ect ({\bf C})$ of
vector spaces over complex field ${\bf C}$, where  its tensor is the
usual tensor,  and its braiding is  the usual twist map. Let ${\bf
Z}_n$ be a finite cycle group of order $n,$  generated by element
$g$, and let ${\bf C Z}_n $  be its group algebra. We have that $$R
=  \frac {1 }{n}  \sum ^{n-1}_{a, b =0}
 e ^{ - \frac {2\pi i ab }{n} } g ^a
\otimes g^b$$
  is a quasitriangular structure of ${\bf C Z}_n$  and
  its module category     ${}_{{\bf CZ}_n} {\cal M}$
  becomes a   braided tensor category .
  Furthermore, the category  is strictly braided for  $n > 2$
( see \cite [ Example 2.1.6 and 9.2.5] {Ma95b} ).  Assume $n>2$ in
the following.
 By \cite [Theorem 2.2  and 2.5] {Ch98}, there is a quasitriangular structure
 $\bar R$ of ${\bf CZ}_n \bowtie ^R {\bf CZ}_n$ and $\bar R$ is not triangular.
 Let $D =  {\bf CZ}_n \bowtie ^R {\bf CZ}_n$.
  Thus  $({}_D {\cal M}, C^R)$   is a strictly braided tensor category by
 \cite [Proposition XIII.1.4] {Ka95}.
  It follows from \cite [Example 9.4.9]{Ma95b}
 (transmutation) that there is a Hopf algebra $\underline D$ living
 in the
 strictly braided  tensor category ${}_D {\cal M}$, where
 $\underline D$ is  given by the same algebra, unit and counit as $D,$
 and comultiplication
 $$\underline \Delta (b) = \sum b_1 S_D(\bar R ^{(2)} )\otimes
 (\bar R^{(1)} \rhd b_2)  $$
for any $b \in \underline D;$ notice that  ``$\rhd $'' is the
quantum action. Since $\underline D$  is  commutative we have that
 $\bar R^{(1)} \cdot b_2 = \epsilon (R ^{(1)}) b_2 $ and
$$\underline \Delta (b) = \sum b_1 S_D(\bar R ^{(2)} )\otimes
 (\bar R^{(1)} \rhd b_2) = \sum b_1 \otimes b_2 = \Delta (b). $$
 That is, $\underline D = D$  as Hopf algebra.   Therefore  there
 are non-trivial $\phi$ and non-trivial $ \psi $ such that
 ${\bf CZ}_n \bowtie ^R {\bf CZ}_n = {\bf CZ}_n {}^{\phi } \bowtie  ^ {\psi}
 {\bf CZ}_n$ . Furthermore  ${\bf CZ}_n {}^{\phi } \bowtie  ^ {\psi}
 {\bf CZ}_n$   is a double cross coproduct in the strictly braided tensor category
 ${}_D {\cal M}$  by \cite
 [Lemma 1.3] {Ch98}.
\end {Example}

\section {The universal property of double bicrossproduct}\label {s4}
In this section, we give the universal property of double
bicrossproduct. Throughout this section,  $H$ and $A$ are two
bialgebras in  braided tensor categories,
  $(A, \alpha )$ is a left $H$-module coalgebra,
 $(H, \beta )$ is a right $A$-module coalgebra,
 $(A, \phi )$ is a left $H$-comodule algebra
 and $(H, \psi )$ is a right $H$-comodule algebra.

\begin {Lemma} \label {2.3.1}            Let $B$ and $E$ be  bialgebras,
 $p$ be a coalgebra morphism from $B$ to $D$, $j$ be an algebra morphism
 from $D$ to $E$, where $D = A ^{\phi } _\alpha \bowtie ^\psi _\beta H$
is a bialgebra. Then

 (i)   $(\pi _A \otimes \pi _H) \Delta _D =id _D$, where $\pi _A$ and $\pi _H $
 are trivial actions;

(ii)  there are two coalgebra morphisms $ p_A$  from
 $B$ to $A$ and $p_H$ from  $B$ to $H$ such that
$$p = (p_A \otimes p_H) \Delta _B.$$
Furthermore, if there are two coalgebra morphisms $f$ from $B$ to
$A$ and $g$ from $B$ to $H$ such that
$$p = (f \otimes g) \Delta _B,$$
then $f= p_A = \pi _A p$  and $g=p_H =\pi _H p.$

(iii) If two coalgebra morphisms  $f$ and  $g$ from $B$ to $D$
satisfy
 $$\pi _A f= \pi g, { \ \ \ } \pi _H f = \pi _H g,$$
 then $f = g$ ;

(iv)    $  m_D (i _A \otimes i _H)  =id _D$;

 (v)  there are two algebra morphisms $ j_A$  from
 $A$ to $E$ and $j_H$ from   $H$ to $E$ such that
$$j = m_E (j_A \otimes j_H) .$$
Furthermore, if there are two algebra morphisms $f$ from  $A$ to $E$
and $g$ from  $H$ to $E$ such that
$$j = m_E (f \otimes g),$$
then $f= j_A= j i_A$  and $g=j_H =j i_H.$

(vi)  If two algebra morphisms $f$ and $g$   from $D$  to $E$
satisfy      $$ f i_A= g i_A, { \ \ \ } f i_H  =g i_H,$$
     then $f =g$ \ \ .
     \end {Lemma}

{\bf Proof.} (i) 
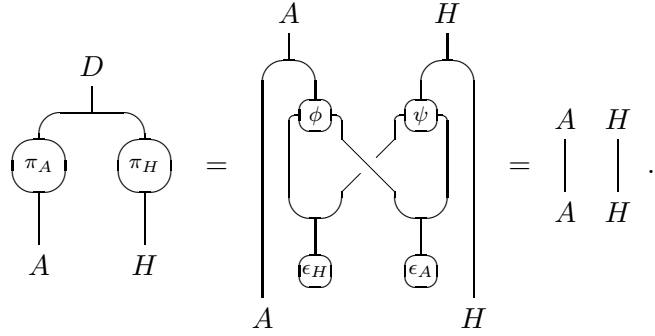
\[
\begin{tangle}
\step[3]\object{D}\\
\step\Cd\\
\obox 2{\pi_{A}}\step[2]\obox 2{\pi_{H}}\\
\step\id\step[2]\step[2]\id\\
\step\object{A}\step[2]\step[2]\object{H}\\
\end{tangle}
\step=\step
\begin{tangle}
\step\object{A}\Step\Step\Step\object{H}\\
 \cd \step[4] \cd \\
 \id \step [1]\td \phi \step [2]\td \psi \step [1] \id\\
 \id \step \id \step [2] \x \step [2]\id \step \id \\
 \id \step \cu \step [2] \cu \step [1]\id  \\
 \id \step [2] \QQ {\epsilon_{H}}\step [2] \step [2]\QQ {\epsilon_{A}} \step [2]\id  \\
\object{A}\step\step\Step\Step\step\step\object{H}
 \end{tangle}
 \step=\step
 \begin{tangle}
\object{A}\Step\object{H}\\
\id\Step\id\\
\object{A}\Step\object{H}\\
 \end{tangle}\ \ \ .
\]

(ii) Let $p_A = \pi _A p$  and $p_H= \pi _H p, $  then $p_A$ and
$p_H$ are coalgebra morphisms and
$$ p= (\pi _A \otimes \pi _H) \Delta _D p = (\pi _A \otimes \pi _H)
 (p \otimes p) \Delta _B = (p_A \otimes p_H) \Delta _B $$  by part (i).
 On the other hand, if           $$p = (f \otimes g) \Delta _B,$$
 then $p_A = \pi _A p = (id _A \otimes \epsilon _H) (f \otimes g) \Delta _B
 =f $  and $p_H = g.$

(iii) We see that $$f = (\pi _A \otimes \pi _H) \Delta _D f = (\pi
_A \otimes \pi _H) (f \otimes f) \Delta _B = (\pi _A \otimes \pi _H)
(g \otimes g) \Delta _B = (\pi _A \otimes \pi _H)  \Delta _B g = g.
$$ Therefore  $f=g.$

(iv)--(vi)   can dually be shown.
\begin{picture}(8,8)\put(0,0){\line(0,1){8}}\put(8,8){\line(0,-1){8}}\put(0,0){\line(1,0){8}}\put(8,8){\line(-1,0){8}}\end{picture}

Assume that $p_A, p_H,$ $j_A$ and $j_H$ are morphisms from $B$ to
$A,$ from $B$ to $H$, from $A$ to $B$ and from $H$ to $B$,
respectively. Set
\[
\hbox {[u1]:}\step[4]
\begin{tangle}
\step[3]\object{B}\\
\step\Cd\\
\obox 2{p_{H}}\step[2]\obox 2{p_{A}}\\
\step\id\step[2]\step[2]\id\\
\step\object{H}\step[2]\step[2]\object{A}\\
\end{tangle}
\step=\step
\begin{tangle}
\step[3]\object{B}\\
\step\Cd\\
\obox 2{p_{A}}\step[2]\obox 2{p_{H}}\\
\step\id\step[2]\step[2]\id\\
\td \phi \step [2]\td \psi\\
\id \step [2] \x \step [2]\id\\
\cu \step [2] \cu \\
\step\object{H}\step[2]\step[2]\object{A}\\
\end{tangle}\step[4];
\]

\[
\hbox {[u2]:}\step[2]
\begin{tangle}
\object{B}\step[2]\object{B}\\
\cu\\
\morph {p_{A}}\\
\step\object{A}\\
\end{tangle}
\step=\step
\begin{tangle}
\step\object{B}\step\step[2]\object{B}\\
\cd\step[2]\id\\
\step[-1]\morph {p_{A}}\morph {p_{H}}\morph {p_{A}}\\
\d\step\tu\alpha\\
\step\cu\\
\step[2]\object{A}\\
\end{tangle}\step
\hbox {and}\step
\begin{tangle}
\object{B}\step[2]\object{B}\\
\cu\\
\morph {p_{H}}\\
\step\object{H}\\
\end{tangle}
\step=\step
\begin{tangle}
\object{B}\step\step[2]\object{B}\\
\id\step[2]\cd\\
\step[-1]\morph {p_{H}}\morph {p_{A}}\morph {p_{H}}\\
\tu\beta\step\dd\\
\step\cu\\
\step[2]\object{H}\\
\end{tangle}\step[4];
\]

\[
\hbox {[u3]:}\step[2]
\begin{tangle}
\Q {\eta_{B}}\\
\step[-1]\morph {p_{H}}\\
\object{H}\\
\end{tangle}
\hstep=\hstep
\begin{tangle}
\Q {\eta_{H}}\\
\id\\
\object{H}\\
\end{tangle}
\step, \step
\begin{tangle}
\Q {\eta_{B}}\\
\step[-1]\morph {p_{A}}\\
\object{A}\\
\end{tangle}
\hstep=\hstep
\begin{tangle}
\Q {\eta_{A}}\\
\id\\
\object{A}\\
\end{tangle}\step ;
\step[2]\hbox {[u4]:}\step
\begin{tangle}
\step\object{H}\step[2]\step[2]\object{A}\\
\step\id\step[2]\step[2]\id\\
\obox 2{j_{H}}\step[2]\obox 2{j_{A}}\\
\step\Cu\\
\step[3]\object{B}\\
\end{tangle}
\step=\step
\begin{tangle}
\step\object{H}\step[2]\step[2]\object{A}\\
\cd \step [2] \cd \\
\id \step [2] \x \step [2]\id\\
\tu \alpha \step [2]\tu \beta\\
\step\id\step[2]\step[2]\id\\
\obox 2{j_{A}}\step[2]\obox 2{j_{H}}\\
\step\Cu\\
\step[3]\object{B}\\
\end{tangle};
\]

\[
\hbox {[u5]:}\step[2]
\begin{tangle}
\step\object{A}\\
\morph {j_{A}}\\
\cd\\
\object{B}\step[2]\object{B}\\
\end{tangle}
\step=\step
\begin{tangle}
\step[2]\object{A}\\
\step\cd\\
\dd\step\td\phi\\
\step[-1]\morph {j_{A}}\morph {j_{H}}\morph {j_{A}}\\
\cu\step[2]\id\\
\step\object{B}\step\step[2]\object{B}\\
\end{tangle}\step,
\step
\begin{tangle}
\step\object{H}\\
\morph {j_{H}}\\
\cd\\
\object{B}\step[2]\object{B}\\
\end{tangle}
\step=\step
\begin{tangle}
\step[2]\object{H}\\
\step\cd\\
\td\psi\step\d\\
\step[-1]\morph {j_{H}}\morph {j_{A}}\morph {j_{H}}\\
\id\step[2]\cu\\
\object{B}\step\step[2]\object{B}\\
\end{tangle}\step[4];
\]

\[
\hbox {[u6]:}\step[4]
\begin{tangle}
\object{H}\\
\step[-1]\morph {j_{H}}\\
\QQ {\epsilon_{B}}\\
\end{tangle}
\hstep=\hstep
\begin{tangle}
\object{H}\\
\id\\
\QQ {\epsilon_{H}}\\
\end{tangle}
\step, \step
\begin{tangle}
\object{A}\\
\step[-1]\morph {j_{A}}\\
\QQ {\epsilon_{B}}\\
\end{tangle}
\hstep=\hstep
\begin{tangle}
\object{A}\\
\id\\
\QQ {\epsilon_{A}}\\
\end{tangle}\step[4].
\]

\begin {Lemma} \label {2.3.2} Let  $D = A ^\phi _\alpha
{\bowtie } ^\psi _\beta H $, $B$ be a bialgebra and $p = (p_A
\otimes p_H) \Delta _B$, where $p_A$  and $p_H$ are  coalgebra
morphisms from $B$ to $A$  and from $B$  to $H$ respectively. Then

  (i) $p$ is a coalgebra morphism iff $(u1)$  holds;

(ii)  $p$ is an algebra morphism iff $(u2)$  and $(u3)$ hold;

(iii)      $p$ is a bialgebra morphism iff $(u1)$-- $(u3)$  hold.
\end {Lemma}

{\bf Proof.} (i) We see that

\[
\begin{tangle}
\step\object{B}\\
\morph p\\
\cd\\
\object{D}\step[2]\object{D}\\
\end{tangle}
\step=\step
\begin{tangle}
\step[4]\object{B}\\
\step[2]\Cd\\
\step\dd\step[4]\d\\
\morph {p_{A}}\step[4]\morph {p_{H}}\\
 \cd \step[4] \cd \\
 \id \step [1]\td \phi \step [2]\td \psi \step [1] \id\\
 \id \step \id \step [2] \x \step [2]\id \step \id \\
 \id \step \cu \step [2] \cu \step [1]\id  \\
\object{A}\step\step\object{H}\Step\Step\object{A}\step\step\object{H}\\
 \end{tangle}
 \step=\step
 \begin{tangle}
\step[4]\object{B}\\
\step[3]\cd\\
\step[2]\cd\step[1]\d\\
\step\dd\Cd\d\\
\dd\morph {p_{A}}\step[2]\morph {p_{H}}\d\\
 \id \step [1]\td \phi \step [2]\td \psi \step [1] \id\\
\step[-1]\morph {p_{A}}  \id \step [2] \x \step [2]\id  \morph {p_{H}} \\
 \id \step \cu \step [2] \cu \step [1]\id  \\
\object{A}\step\step\object{H}\Step\Step\object{A}\step\step\object{H}\\
 \end{tangle}
 \step\hbox {.....(1)}
\]

\[
\begin{tangle}
\step\object{B}\\
\cd\\
\step[-1]\morph p\morph p\\
\object{D}\step[2]\object{D}\\
\end{tangle}
\step=\step
\begin{tangle}
\step[3]\object{B}\\
\step[2]\cd\\
\step\dd\step[2]\d\\
 \cd \step[2] \cd \\
\step[-1]\morph {p_{A}}\morph {p_{H}}\morph {p_{A}}\morph
{p_{H}}\\
\object{A}\step\step\object{H}\Step\object{A}\step\step\object{H}\\
 \end{tangle}
 \step=\step
 \begin{tangle}
\step[4]\object{B}\\
\step[3]\cd\\
\step[2]\cd\step[1]\d\\
\step\dd\step\cd\step\d\\
\morph {p_{A}}\morph {p_{H}}\morph
{p_{A}}\morph {p_{H}}\\
\step\object{A}\step\step\object{H}\Step\object{A}\step\step\object{H}\\
 \end{tangle}
 \step\hbox {.....(2)}
\]
If $p$ is a coalgebra morphism, then by relation (1) and (2), we
have that
\[
 \begin{tangle}
\step[4]\object{B}\\
\step[3]\cd\\
\step[2]\cd\step[1]\d\\
\step\dd\Cd\d\\
\dd\morph {p_{A}}\step[2]\morph {p_{H}}\d\\
\step[-1]\morph {p_{A}}  \td \phi \step [2]\td \psi \morph {p_{H}}\\
 \QQ {\epsilon_{A}}\step\id \step [2] \x \step [2] \id\step\QQ {\epsilon_{A}}   \\
 \step \cu \step [2] \cu  \\
\step\step\object{H}\Step\Step\object{A}\step\step\\
 \end{tangle}
 \step=\step
 \begin{tangle}
\step[4]\object{B}\\
\step[3]\cd\\
\step[2]\cd\step[1]\d\\
\step\dd\step\cd\step\d\\
\morph {p_{A}}\morph {p_{H}}\morph
{p_{A}}\morph {p_{H}}\\
\step[1]\QQ {\epsilon_{A}}\step[2]\id\step[2]\id\step[2]\QQ
{\epsilon_{H}}\\
\step\step\step\object{H}\Step\object{A}\step\step\\
 \end{tangle} \step [4].
\]
i.e.

\[
 \begin{tangle}
\step[4]\object{B}\\
\step[2]\Cd\\
\step\morph {p_{A}}\step[2]\morph {p_{H}}\\
 \step \td \phi \step [2]\td \psi \step\\
\step\id \step [2] \x \step [2] \id\step   \\
 \step \cu \step [2] \cu  \\
\step\step\object{H}\Step\Step\object{A}\step\step\\
 \end{tangle}
 \step=\step
 \begin{tangle}
\step[1]\object{B}\\
\cd\\
\step[-1]\morph {p_{H}}\morph {p_{A}}\\
\id\step[2]\id\\
\object{H}\Step\object{A}\\
 \end{tangle} \step [4].
\]
Thus (u1) holds.

 Conversely, if $\hbox {(u1 )}$ holds, by relation (1) and (2), we have that
  \[\step\begin{tangle}
\step\object{B}\\
\morph p\\
\cd\\
\object{D}\step[2]\object{D}
\end{tangle}
\step=\step\begin{tangle}
\step\object{B}\\
\cd\\
\step[-1]\morph p\morph p\\
\object{D}\step[2]\object{D}
\end{tangle}\step \step [2] and \step [2]
\step\begin{tangle}
\step\object{B}\\
\morph p\\
\step\QQ {\epsilon_{D}}
\end{tangle}
\step=\step\begin{tangle}
\object{B}\\
\id\\ \QQ {\epsilon_{B}}
\end{tangle}\step.\]  Thus $p$ is a coalgebra morphism.

(ii) We first see that
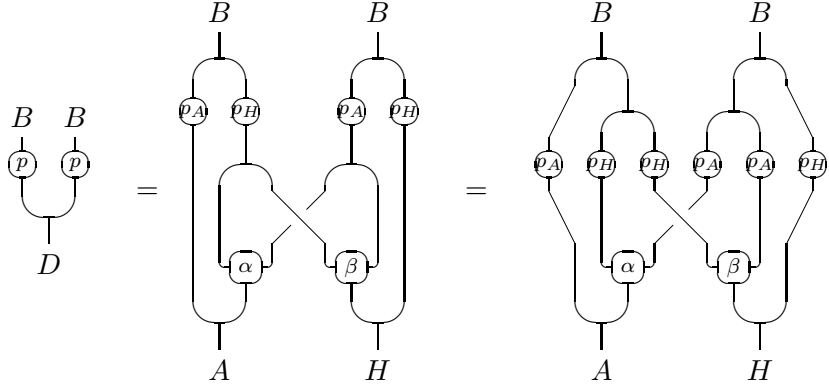
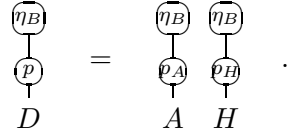
\[
\begin{tangle}
\object{B}\step[2]\object{B}\\
\cu\\
\morph p\\
\step\object{D}\\
\end{tangle}
\step=\step\begin{tangle}
\object{B}\step[2]\object{B}\\
\cu\\
\cd\\
\step[-1]\morph {p_{A}}\morph {p_{H}}\\
\object{A}\step[2]\object{H}\\
\end{tangle}
\]
and
\[
\begin{tangle}
\object{B}\step[2]\object{B}\\
\step[-1]\morph p\morph p\\
\cu\\
\step\object{D}\\
\end{tangle}
\step=\step
\begin{tangle}
\step\object{B}\step[6]\object{B}\\
\cd\step[2]\step[2]\cd\\
\step[-1]\morph {p_{A}}\morph {p_{H}}\step[2]\morph
{p_{A}}\morph {p_{H}}\\
\id\step\cd\step[2]\cd\step\id\\
\id\step\id\step[2]\x\step[2]\id\step\id\\
\id\step\tu \alpha\step[2]\tu\beta\step\id\\
\cu\step[2]\step[2]\cu\\
\step\object{A}\step[6]\object{H}\\
\end{tangle}
\step=\step \begin{tangle}
\step\step\step\object{B}\step[6]\object{B}\\
\step\step\cd\step[2]\step[2]\cd\\
\step\ne1\step\cd\step[2]\cd\step\nw1\\

\morph {p_{A}}\morph {p_{H}}\morph {p_{H}}  \morph {p_{A}}\morph
{p_{A}}\morph {p_{H}}\\
\step\nw1 \step\id\step[2]\x\step[2]\id\step\ne1\\
\step\step\id\step\tu \alpha\step[2]\tu\beta\step\id\\
\step\step\cu\step[2]\step[2]\cu\\
\step\step\step\object{A}\step[6]\object{H}\\
\end{tangle}
\]
and
\[
\begin{tangle}
\step\Q {\eta_{B}}\\
\morph p\\
\step\object{D}
\end{tangle}
\step=\step
\begin{tangle}
\step\Q {\eta_{B}}\\
\morph {p_{A}}\\
\step\object{A}
\end{tangle}
\begin{tangle}
\step\Q {\eta_{B}}\\
\morph {p_{H}}\\
\step\object{H}
\end{tangle}\ \ \ .
\]

If $p$ is an algebra morphism, by the above proof, we have that

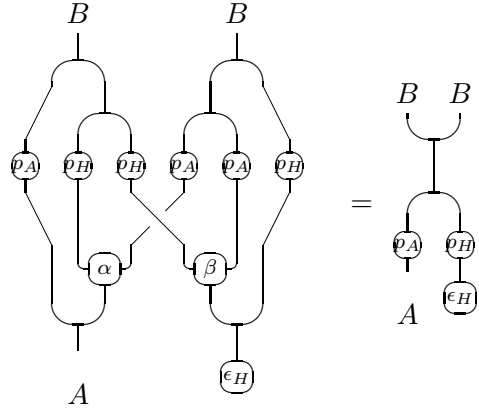
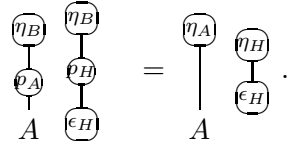
\[\begin{tangle}
\step\step\step\object{B}\step[6]\object{B}\\
\step\step\cd\step[2]\step[2]\cd\\
\step\ne1\step\cd\step[2]\cd\step\nw1\\

\morph {p_{A}}\morph {p_{H}}\morph {p_{H}}  \morph {p_{A}}\morph
{p_{A}}\morph {p_{H}}\\
\step\nw1 \step\id\step[2]\x\step[2]\id\step\ne1\\
\step\step\id\step\tu \alpha\step[2]\tu\beta\step\id\\
\step\step\cu\step[2]\step[2]\cu\\
\step\step\step\object{A}\step[6]\QQ {\epsilon _H}\\
\end{tangle}
\step=\step\begin{tangle}
\object{B}\step[2]\object{B}\\
\cu\\
\cd\\
\step[-1]\morph {p_{A}}\morph {p_{H}}\\
\object{A}\step[2]\QQ {\epsilon _{H}}\\
\end{tangle}
\]
and
\[
\begin{tangle}
\step\Q {\eta_{B}}\\
\morph {p_{A}}\\
\step\object{A}
\end{tangle}
\begin{tangle}
\step\Q {\eta_{B}}\\
\morph {p_{H}}\\
\step\QQ {\epsilon _{H}}
\end{tangle}
\step=\step
\begin{tangle}
\Q {\eta_{A}}\\
\id\\
\object{A}
\end{tangle}
\step
\begin{tangle}
\step\Q {\eta_{H}}\\
\step\QQ {\epsilon _{H}}
\end{tangle}\ \ \ .
\]

In this way, we can obtain (u2)  and (u3).

Conversely, if (u1) and (u2) hold, we can similarly show that $p$ is
an algebra morphism.

(iii) It follows from part (i) and part (ii).
\begin{picture}(8,8)\put(0,0){\line(0,1){8}}\put(8,8){\line(0,-1){8}}\put(0,0){\line(1,0){8}}\put(8,8){\line(-1,0){8}}\end{picture}

 Dually, we have
\begin {Lemma} \label {2.3.3'}  Let
$B$ be a bialgebra and $j = m_B(j_A \otimes j_H) $, where $j_A$  and
$j_H$ are algebra morphisms from $A$ to $B$  and from $H$  to $B$
respectively. Then

  (i) $j$ is an algebra morphism iff $(u4)$  holds;

(ii)  $j$ is a coalgebra morphism iff $(u5)$  and $(u6)$ hold;

(iii)      $j$ is a bialgebra morphism iff $(u4)$-- $(u6)$  hold .

\end {Lemma}

\begin {Definition} \label {2.3.3}
(UT1): Let $B$ be a bialgebra and   two coalgebra morphisms
$$ p_A : B \longrightarrow A , { \ \ \ } p_H : B \longrightarrow H $$
 satisfy conditions (u1)--(u3). $(B, p_A, p_H)$  is said
to satisfy condition (UT1),
 if for any bialgebra $E$ and any two coalgebra morphisms
 $$q_A : E \rightarrow A, q_H : E \rightarrow H \ , $$ which satisfy conditions
 (u1)--(u3),
there exists a unique bialgebra morphism $q'$  from $E$  to $B$ such
that
$$p_A q' = q_A, { \ \ \ } p_Hq' = q_H,$$
i.e. the diagram

\[
\xymatrix{& A&\\
E\ar[ur]^{q_{A}}\ar[dr]_{q_{H}}\ar@{.>}[rr]^{q'}& &B\ar[ul]_{p_{A}}\ar[dl]^{p_{H}} \\
&H & }
\]
 commutes.

(UT2): Let $B$ and   $D = A ^\phi _\alpha {\bowtie } ^\psi _\beta H$
be  bialgebras and $p$ a bialgebra morphism from $B$  to $D$.
$(B,p)$  is said  to   satisfy condition (UT2) if for any bialgebra
$E$ and any bialgebra morphism $q$ from $E$ to $D$, there exists a
unique bialgebra morphism $q'$ from $E$ to $B$  such that the
diagram
\[
\xymatrix{&D&\\
E\ar[ur]^{q}\ar@{.>}[rr]^{q'}& &B\ar[ul]_{p} }
\] commutes.

(UT3): Let $B$ be coalgebra and   two coalgebra morphisms
$$ p_A : B \longrightarrow A , { \ \ \ } p_H : B \longrightarrow H $$
 satisfy condition (u1). $(B, p_A, p_H)$  is said
to satisfy condition (UT3),
 if for any coalgebra $E$ and  any two coalgebra morphisms
 $$q_A : E \rightarrow A, q_H : E \rightarrow H\ , $$ which satisfy condition
 (u1),
there exists a unique coalgebra morphism $q'$  from $E$  to $B$ such
that
$$p_A q' = q_A, { \ \ \ } p_Hq' = q_H,$$
i.e. the diagram \[
\xymatrix{& A&\\
E\ar[ur]^{q_{A}}\ar[dr]_{q_{H}}\ar@{.>}[rr]^{q'}& &B\ar[ul]_{p_{A}}\ar[dl]^{p_{H}} \\
&H & }
\]
commutes.

(UT4): Let $B$ and   $D = A \stackrel {c} {\bowtie } H$   be
coalgebras and $p$ a coalgebra morphism from $B$  to $D$. $(B,p)$ is
said  to   satisfy condition (UT4) if for any coalgebra $E$ and any
coalgebra morphism $q$ from $E$ to $D$, there exists a unique
coalgebra morphism $q'$ from $E$ to $B$  such that the diagram
\[
\xymatrix{&D&\\
E\ar[ur]^{q}\ar@{.>}[rr]^{q'}& &B\ar[ul]_{p} }
\] commutes.

(UT1'): Let $B$ be a bialgebra and   two algebra morphisms
$$ j_A : A \longrightarrow B , { \ \ \ } j_H : H \longrightarrow B $$
 satisfy conditions (u4)--(u6). $(B, j_A, j_H)$  is said
to satisfy condition (UT1'),
 if for any bialgebra $E$ and any two algebra morphisms
 $$q_A : A \rightarrow E, q_H : H \rightarrow E \ , $$ which satisfy conditions
 (u4)--(u6),
there exists a unique bialgebra morphism $q'$  from $B$  to $E$ such
that
$$ q'j_A = q_A, { \ \ \ } q'j_H = q_H,$$
i.e. the diagram \[
\xymatrix{& A\ar[dl]_{q_{A}}\ar[dr]^{j_{A}}&\\
E& &B\ar@{.<}[ll]_{q'} \\
&H \ar[ul]^{q_{H}}\ar[ur]_{j_{H}}& }
\]
commutes.

(UT2'): Let $B$ and   $D = A ^\phi _\alpha {\bowtie }^\psi _\beta H$
be  bialgebras and $j$ a bialgebra morphism from $D$  to $B$.
$(B,j)$  is said  to   satisfy condition (UT2') if for any bialgebra
$E$ and any bialgebra morphism $q$ from $D$ to $E$, there exists a
unique bialgebra morphism $q'$ from $B$ to $E$  such that the
diagram $q'j=q$.

(UT3'): Let $B$ be algebra and   two algebra morphisms
$$ j_A : A \longrightarrow B , { \ \ \ } j_H : H \longrightarrow B $$
 satisfy condition (u1). $(B, j_A, j_H)$  is said
to satisfy condition (UT3'),
 if for any algebra $E$ and any two algebra morphisms
 $$q_A : A \rightarrow E, q_H : H \rightarrow E\ , $$ which satisfy condition
 (u4),
there exists a unique algebra morphism $q'$  from $B$  to $E$ such
that
$$q'j_A = q_A, { \ \ \ } q'j_H = q_H.$$

(UT4'): Let $B$ and   $D = A {\bowtie } H$   be  algebras and $j$ an
algebra morphism from $D$  to $B$. $(B,j)$  is said  to   satisfy
condition (UT4'), if for any algebra $E$ and any algebra morphism
$q$ from $D$ to $E$, there exists a unique algebra morphism $q'$
from $B$ to $E$  such that the diagram $q=q'j.$
\end {Definition}

 \begin {Theorem} \label {2.3.4}
 Let $B$ and
 $D = A ^\phi _\alpha
{\bowtie }^\psi _\beta  H$ be bialgebras.

 (I) (i)  $(D, \pi _A, \pi _H)$  satisfies condition (UT1);

(ii)  $(D, id _D)$ satisfies condition (UT2);

 (iii)  $(D, i _A, i _H)$  satisfies condition (UT1');

(iv)  $(D, id _D)$ satisfies condition (UT2').

(II) The following statements are equivalent;

(i) $(B, p_A, p_H)$  satisfies condition (UT1);

(ii) There exists a (necessarily unique )  bialgebra isomorphism
$p'$  from $B$ to $D$ such that
$$\pi _A p' = p_A, { \ \ \ }  \pi _H p' = p_H.$$

(iii)             $p$  is a bialgebra isomorphism from $B$ to $D$
and $p = (p_A \otimes p_H) \Delta _B;$

(iv)  $(B, p)$ satisfies condition (UT2) and
 $p = (p_A \otimes p_H) \Delta _B.$

(III)        The following statements are equivalent.

(i) $(B, j_A, j_H)$  satisfies condition (UT1');

(ii) There exists a (necessarily unique )  bialgebra isomorphism
$j'$  from $D$ to $B$ such that
$$j' i _A  = j_A, { \ \ \ }  \j' i_H  = j_H;$$

(iii)             $j$  is a bialgebra isomorphism from $D$ to $B$
and $j = m_B(j_A \otimes j_H) ;$

(iv)  $(B, j)$ satisfies condition (UT2') and
 $j = m_B(j_A \otimes j_H). $

\end {Theorem}

{\bf Proof.}  (I) (i) For any bialgebra  $E$ and two coalgebra
morphisms
$$q_A : E \longrightarrow A, { \ \ \ } q_H : E  \longrightarrow H \ , $$
 which satisfy conditions
(u1)--(u3), set $$ q' =  (q_A \otimes q_H)\Delta _E \ .$$ It is
clear that $$\pi _A q'= q_A, { \ \ \ } \pi _H q' = q_H$$  The
uniqueness of $q'$ can be obtained  by  Lemma \ref {2.3.1} (iii).
Otherwise, $id _D = (\pi _A \otimes \pi_H) \Delta _D$  by Lemma \ref
{2.3.1} (i) and $\pi _A$ and $\pi _H$ satisfy conditions (u1)--(u3)
by Lemma \ref {2.3.2} (iii). Thus
 $(D, \pi _A, \pi _H )$  satisfies condition (TU1).

(ii) For any bialgebra $E$ and any bialgebra morphism $q$  from $E$
to $D$, let $q' =q$. Obviously, $id _D q'=q$ and we have the
uniqueness of $q'$. Thus $(D,id_D)$ satisfies condition (UT2).

(iii) and (iv) can dually be checked .

(II) (i) $\Rightarrow $  (ii)
 Since   $(D, \pi _A, \pi _H )$ and
 $(B, p_A, p_H)$ satisfy condition (UT1),
we have that there exist two bialgebra morphisms
 $$f: B \rightarrow D  \hbox { \ \ and \ \ }  g: D \rightarrow B $$
such that
$$p_Ag = \pi _A, p_Hg = \pi _H, \pi _A f = p_A, \pi _H f = p_H.$$
Now we see that
$$\pi _A fg = \pi _A, \pi _H fg = \pi _H, p_A gf = p_A , p_H gf = p_H.$$
Therefore, by uniqueness,  $$fg = id _D \hbox { \ \ and \ \ }gf =
id_B.$$ That is,  $f $  is a bialgebra isomorphism from $B$ to $D$.

(ii) $\Rightarrow $ (iii)  Let $p = (p_A \otimes p_H)\Delta _B.$
Since
$$\pi _A p' = p_A, { \ \ \ }  \pi _H p' = p_H, { \ \ \ }\pi _A p = p_A, { \ \ \ }
 \pi _H p = p_H,$$
considering Lemma \ref {2.3.1} (iii), we have that $$p =p'.$$
Therefore $p$ is a bialgebra isomorphism.

(iii) $\Rightarrow $ (iv)
            For any bialgebra $E$ and any bialgebra morphism $q$  from
$E$ to $D$,  let $$q'=p^{-1} q.$$ Obviously, $pq'= q$  and we have
the uniqueness of $q'.$ Thus $(B,p)$ satisfies condition (UT2).

(iv) $\Rightarrow $ (i) For any bialgebra $E$ and two coalgebra
morphisms
 $$q_A : E \rightarrow A, q_H : E \rightarrow H \ , $$ which satisfy conditions
 (u1)--(u3), set $q = (q_A \otimes q_H) \Delta _E.$  By Lemma 3.2 (iii),
 we have that $q$  is a bialgebra morphism from $E$  to $D$. Thus
 there exists a unique bialgebra morphism $q'$  from $E$  to $B$
 such that $$pq'=q$$
 since $(B,p)$ satisfies condition (UT2).
It is clear that $$p_A q' = q_A, { \ \ \ } p_Hq' = q_H.$$ Otherwise,
$p_A$ and $p_H$ satisfy conditions (u1)--(u3) by Lemma \ref {2.3.2}
(iii) Therefore, $(B, p_A, p_H)$  satisfies condition (UT1).

(III) It can dually be obtained.

We can easily find the following two theorems from the proof of the
previous theorem.

 \begin {Theorem} \label {2.3.5}
 Let $B$ and
 $D = A \stackrel {c}
{\bowtie } H$ be coalgebras.

 (I) (i)  $(D, \pi _A, \pi _H)$  satisfies condition (UT3);

(ii)  $(D, id _D)$ satisfies condition (UT4);

(II) The following statements are equivalent

(i) $(B, p_A, p_H)$  satisfies condition (UT3).

(ii) There exists a (necessarily unique )  coalgebra isomorphism
$p'$  from $B$ to $D$ such that
$$\pi _A p' = p_A, { \ \ \ }  \pi _H p' = p_H;$$

(iii)             $p$  is a coalgebra isomorphism from $B$ to $D$
and $p = (p_A \otimes p_H) \Delta _B.$

(iv)  $(B, p)$ satisfies condition (UT4) and
 $p = (p_A \otimes p_H) \Delta _B.$
 \end {Theorem}

 \begin {Theorem} \label {2.3.6}
 Let $B$ and
 $D = A
{\bowtie } H$ be algebras.

(I)  (i)  $(D, i _A, i _H)$  satisfies condition (UT3');

(ii)  $(D, id _D)$ satisfies condition (UT4');

(II)        The following statements are equivalent.

(i) $(B, j_A, j_H)$  satisfies condition (UT3').

(ii) There exists an (necessarily unique )  algebra isomorphism $j'$
from $D$ to $B$ such that
$$j' i _A  = j_A, { \ \ \ }  \j' i_H  = j_H;$$

(iii)             $j$  is an algebra isomorphism from $D$ to $B$ and
$j = m_B(j_A \otimes j_H) .$

(iv)  $(B, j)$ satisfies condition (UT4') and
 $j = m_B(j_A \otimes j_H) $.

\end {Theorem}

           \chapter{   The Bicrossproducts in Braided Tensor Categories }\label {c4}

In the category of usual vector spaces with usual twist braiding, S.
Majid introduced bicrossproducts  and gave the necessary and
sufficient conditions for them to become bialgebras in \cite
[Theorem 2.9] {Ma94a}. D.E. Radford  introduced biproducts  and gave
the necessary and sufficient conditions for them to become
bialgebras in \cite [Theorem 1] {Ra85}. In braided tensor category,
S. Majid introduced the smash product and showed that if $(A, \alpha
)$ is an $H$-module bialgebra and $H$ is cocommutative with respect
to $(H, \alpha ),$ then  the smash product  $A _\alpha \# H$ is a
bialgebra in \cite [Theorem 2.4] {Ma94b}.

In this Chapter, we construct the bicrossproducts and biproducts in
braided tensor category  $({\cal C}, C)$ and give the necessary and
sufficient conditions for them to be bialgebras. We obtain that  for
any left module $(A, \alpha )$  of bialgebra $H$ the smash product
$A \#H$  is a bialgebra iff $A$  is an $H$-module bialgebra and $H$
is cocommutative with respect to $(A, \alpha ).$

Now we give some concepts as follows: assume  that $H, A \in ob \
{\cal C},$  and  that
\begin {eqnarray*}
\alpha : H \otimes A \rightarrow &A& , \hbox { \ \ \ \ }
\beta : H \otimes A \rightarrow H,    \\
\phi :  A \rightarrow H \otimes  &A&   , \hbox { \ \ \ \ }
\psi : H  \rightarrow H \otimes A,  \\
m _H : H \otimes H \rightarrow &H& , \hbox { \ \ \ \ }
m_A : A \otimes A \rightarrow A,    \\
\Delta _H :  H \rightarrow H \otimes  &H&   , \hbox { \ \ \ \ }
\Delta _A : A  \rightarrow A \otimes A, \\
\eta _H : I  \rightarrow &H& , \hbox { \ \ \ \ }
\eta _A : I \rightarrow A,    \\
\epsilon  _H :  H  \rightarrow   &I&   , \hbox { \ \ \ \ }
\epsilon  _A  : A  \rightarrow I.     \\
\sigma  : H \otimes H \rightarrow &A& , \hbox { \ \ \ \ }
\mu : A \otimes A \rightarrow H,    \\
P :  A \rightarrow H \otimes  &H&   , \hbox { \ \ \ \ } Q : H
\rightarrow A \otimes A
\end {eqnarray*}
are morphisms in ${\cal C}$.

$(H,\alpha )$   is said to act weakly on $A $ if the following
conditions are satisfied:
 \[
(WA):
\begin{tangle}
\object{H}\step\object{A}\step[2]\object{A}\\
\id\step\tu m\\
\tu \alpha\\
\step\object{A}\\
\end{tangle}
\step=\step
\begin{tangle}
\step\object{H}\step[3]\object{A}\step[2]\object{A}\\
\cd\step[2]\id\step[2]\id\\
\id\step[2]\x\step[2]\id\\
\tu \alpha\step[2]\tu \alpha\\
\step\Cu\\
\step[3]\object{A}\\
\end{tangle}
\step,\step
\begin{tangle}
\object{H}\\
\id\step[2]\Q {\eta_A}\\
\tu \alpha\\
\step\object{A}\\
\end{tangle}
\step=\step
\begin{tangle}
\object{H}\\
\QQ \epsilon \\
\Q {\eta_A}\\
\object{A}\\
\end{tangle}
\step,\step
\begin{tangle}
\step[2]\object{A}\\
\Q {\eta_H}\step[2]\id\\
\tu \alpha\\
\step\object{A}\\
\end{tangle}
\step=\step
\begin{tangle}
\object{A}\\
\id\\
\id\\
\object{A}\\
\end{tangle}
\step.\step
\]

$(A,\beta )$   is said to act weakly on $H $  if the following
conditions are satisfied:

\[
(WA):
\begin{tangle}
\object{H}\step[2]\object{H}\step\object{A}\\
\tu m\step\id\\
\step\tu \beta\\
\step[2]\object{H}\\
\end{tangle}
\step=\step
\begin{tangle}
\object{H}\step[2]\object{H}\step[3]\object{A}\\
\id\step[2]\id\step[2]\cd\\
\id\step[2]\x\step[2]\id\\
\tu \beta\step[2]\tu \beta\\
\step\Cu\\
\step[3]\object{H}\\
\end{tangle}
\step,\step
\begin{tangle}
\step[2]\object{A}\\
\Q {\eta_H}\step[2]\id\\
\tu \beta\\
\step\object{H}\\
\end{tangle}
\step=\step
\begin{tangle}
\object{A}\\
\QQ \epsilon \\
\Q {\eta_H}\\
\object{H}\\
\end{tangle}
\step,\step
\begin{tangle}
\object{H}\\
\id\step[2]\Q {\eta_A}\\
\tu \beta\\
\step\object{H}\\
\end{tangle}
\step=\step
\begin{tangle}
\object{A}\\
\id\\
\id\\
\object{H}\\
\end{tangle}
\step.\step
\]

$(H,\phi  )$   is said to coact weakly on $A $  if the following
conditions are satisfied:

\[
(WCA):
\begin{tangle}
\step\object{A}\\
\td \phi\\
\id\step\cd\\
\object{H}\step\object{A}\step[2]\object{A}\\
\end{tangle}
\step=\step
\begin{tangle}
\step[3]\object{A}\\
\step\Cd\\
\td \phi\step[2]\td \phi\\
\id\step[2]\x\step[2]\id\\
\cu\step[2]\id\step[2]\id\\
\step\object{H}\step[3]\object{A}\step[2]\object{A}\\
\end{tangle}
\step,\step
\begin{tangle}
\step\object{A}\\
\td \phi\\
\id\step[2]\QQ \epsilon \\
\object{H}\\
\end{tangle}
\step=\step
\begin{tangle}
\object{A}\\
\QQ \epsilon \\
\Q {\eta_H}\\
\object{H}\\
\end{tangle}
\step,\step
\begin{tangle}
\step\object{A}\\
\td \phi\\
\QQ \epsilon \step[2]\id\\
\step[2]\object{A}\\
\end{tangle}
\step=\step
\begin{tangle}
\object{A}\\
\id\\
\id\\
\object{A}\\
\end{tangle}
\step.\step
\]

$(A,\psi )$   is said to coact weakly on $H $  if the following
conditions are satisfied:

\[
(WCA):
\begin{tangle}
\step[2]\object{H}\\
\step\td \psi\\
\cd\step\id\\
\object{H}\step[2]\object{H}\step\object{A}\\
\end{tangle}
\step=\step
\begin{tangle}
\step[3]\object{H}\\
\step\Cd\\
\td \psi\step[2]\td \psi\\
\id\step[2]\x\step[2]\id\\
\id\step[2]\id\step[2]\cu\\
\object{H}\step[2]\object{H}\step[3]\object{A}\\
\end{tangle}
\step,\step
\begin{tangle}
\step\object{H}\\
\td \psi\\
\QQ \epsilon \step[2]\id\\
\step[2]\object{A}\\
\end{tangle}
\step=\step
\begin{tangle}
\object{H}\\
\QQ \epsilon \\
\Q {\eta_A}\\
\object{A}\\
\end{tangle}
\step,\step
\begin{tangle}
\step\object{H}\\
\td \psi\\
\id\step[2]\QQ \epsilon \\
\object{H}\\
\end{tangle}
\step=\step
\begin{tangle}
\object{H}\\
\id\\
\id\\
\object{H}\\
\end{tangle}
\step.\step
\]

$\sigma $   is called a 2-cocycle from $H\otimes H$  to $A $  if the
following conditions are satisfied:
\[
(2\hbox {-}COC):
\begin{tangle}
\step\object{H}\step[3]\object{H}\step[3]\object{H}\\
\cd\step\cd\step\cd\\
\id\step[2]\id\step\id\step[2]\hx\step[2]\id\\
\id\step[2]\id\step\cu\step\cu\\
\id\step[2]\x\step[2]\dd\\
\tu {\alpha}\step[2]\tu {\sigma}\\
\step\Cu\\
\step[3]\object{A}\\
\end{tangle}
\step=\step
\begin{tangle}
\step\object{H}\step[3]\object{H}\step[2]\object{H}\\
\cd\step\cd\step\id\\
\id\step[2]\hx\step[2]\id\step\id\\
\tu {\sigma}\step\cu\step\id\\
\step\id\step[3]\tu {\sigma}\\
\step\Cu\\
\step[3]\object{A}\\
\end{tangle}
\step,\step
\]
\[
\begin{tangle}
\object{H}\\
\id\step[2]\Q {\eta_H}\\
\tu {\sigma}\\
\step\object{A}\\
\end{tangle}
\step=\step
\begin{tangle}
\step[2]\object{H}\\
\Q {\eta_H}\step[2]\id\\
\tu {\sigma}\\
\step\object{A}\\
\end{tangle}
\step=\step
\begin{tangle}
\object{H}\\
\QQ {\epsilon }\\
\Q {\eta_A}\\
\object{A}\\
\end{tangle}
\step.\step
\]

$\mu$   is called a 2-cocycle from $A \otimes A$ to $H$  if the
following conditions are satisfied:
\[
(2\hbox {-}COC):
\begin{tangle}
\object{A}\step[2]\object{A}\step[3]\object{A}\\
\id\step\cd\step\cd\\
\id\step\id\step[2]\hx\step[2]\id\\
\id\step\cu\step\tu {\mu}\\
\tu {\mu}\step[3]\id\\
\step\Cu\\
\step[3]\object{H}\\
\end{tangle}
\step=\step
\begin{tangle}
\step\object{A}\step[3]\object{A}\step[3]\object{A}\\
\cd\step\cd\step\cd\\
\id\step[2]\hx\step[2]\id\step\id\step[2]\id\\
\cu\step\tu {\mu}\step\id\step[2]\id\\
\step\se1\step[2]\x\step[2]\id\\
\step[2]\tu {\mu}\step[2]\tu {\beta}\\
\step[3]\Cu\\
\step[5]\object{H}\\
\end{tangle}
\step,\step
\]
\[
\begin{tangle}
\object{A}\\
\id\step[2]\Q {\eta_A}\\
\tu {\mu}\\
\step\object{H}\\
\end{tangle}
\step=\step
\begin{tangle}
\step[2]\object{A}\\
\Q {\eta_A}\step[2]\id\\
\tu {\mu}\\
\step\object{H}\\
\end{tangle}
\step=\step
\begin{tangle}
\object{A}\\
\QQ {\epsilon }\\
\Q {\eta_H}\\
\object{H}\\
\end{tangle}
\step.\step
\]

$P $   is called a 2-cycle from $A$ to $H\otimes H$    if the
following conditions are satisfied:

\[
(2\hbox {-}C):
\begin{tangle}
\step[3]\object{A}\\
\step\Cd\\
\td {\phi}\step[2]\td {P}\\
\id\step[2]\x\step[2]\se1\\
\cu\step\td {P}\step\cd\\
\step\id\step[2]\id\step[2]\hx\step[2]\id\\
\step\id\step[2]\cu\step\cu\\
\step\object{H}\step[3]\object{H}\step[3]\object{H}\\
\end{tangle}
\step=\step
\begin{tangle}
\step[3]\object{A}\\
\step\Cd\\
\step\id\step[3]\td {P}\\
\td {P}\step\cd\step\id\\
\id\step[2]\hx\step[2]\id\step\id\\
\cu\step\cu\step\id\\
\step\object{H}\step[3]\object{H}\step[2]\object{H}\\
\end{tangle}
\step,\step
\]
\[
\begin{tangle}
\step\object{A}\\
\td {P}\\
\QQ {\epsilon }\step[2]\id\\
\step[2]\object{H}\\
\end{tangle}
\step=\step
\begin{tangle}
\step\object{A}\\
\td {P}\\
\id\step[2]\QQ {\epsilon }\\
\object{H}\\
\end{tangle}
\step=\step
\begin{tangle}
\object{A}\\
\QQ {\epsilon }\\
\Q {\eta_H}\\
\object{H}\\
\end{tangle}
\step.\step
\]

$Q$   is called a 2-cycle from $H$ to $A \otimes A$   if the
following conditions are satisfied:

\[
(2\hbox {-}C):
\begin{tangle}
\step[3]\object{H}\\
\step\Cd\\
\td {\alpha}\step[3]\id\\
\id\step\cd\step\td {\alpha}\\
\id\step\id\step[2]\hx\step[2]\id\\
\id\step\cu\step\cu\\
\object{A}\step[2]\object{A}\step[3]\object{A}\\
\end{tangle}
\step=\step
\begin{tangle}
\step[5]\object{H}\\
\step[3]\Cd\\
\step[2]\td {Q}\step[2]\td {\psi}\\
\step\ne1\step[2]\x\step[2]\id\\
\cd\step\td {Q}\step\cu\\
\id\step[2]\hx\step[2]\id\step[2]\id\\
\cu\step\cu\step[2]\id\\
\step\object{A}\step[3]\object{A}\step[3]\object{A}\\
\end{tangle}
\step,\step
\]
\[
\begin{tangle}
\step\object{H}\\
\td {Q}\\
\QQ {\epsilon }\step[2]\id\\
\step[2]\object{A}\\
\end{tangle}
\step=\step
\begin{tangle}
\step\object{H}\\
\td {Q}\\
\id\step[2]\QQ {\epsilon }\\
\object{A}\\
\end{tangle}
\step=\step
\begin{tangle}
\object{H}\\
\QQ {\epsilon }\\
\Q {\eta_H}\\
\object{H}\\
\end{tangle}
\step.\step
\]

$(A, \alpha )$   is called a twisted $H$-module   if the following
conditions are satisfied:
\[
(TM):
\begin{tangle}
\step\object{H}\step[3]\object{H}\step[3]\object{A}\\
\cd\step\cd\step[2]\id\\
\id\step[2]\id\step\id\step[2]\x\\
\id\step[2]\id\step\tu {\alpha}\step[2]\id\\
\id\step[2]\x\step[2]\ne1\\
\tu {\alpha}\step[2]\tu {\sigma}\\
\step\Cu\\
\step[3]\object{A}\\
\end{tangle}
\step=\step
\begin{tangle}
\step\object{H}\step[3]\object{H}\step[2]\object{A}\\
\cd\step\cd\step\id\\
\id\step[2]\hx\step[2]\id\step\id\\
\tu {\sigma}\step\cu\step\id\\
\step\id\step[3]\tu {\alpha}\\
\step\Cu\\
\step[3]\object{A}\\
\end{tangle}
\step,\step
\]
\[
\begin{tangle}
\step[2]\object{A}\\
\Q {\eta_H}\step[2]\id\\
\tu {\alpha}\\
\step\object{A}\\
\end{tangle}
\step=\step
\begin{tangle}
\object{A}\\
\id\\
\id\\
\object{A}\\
\end{tangle}
\step.\step
\]

$(H, \beta )$   is called a twisted $A$-module   if the following
conditions are satisfied:
\[
(TM):
\begin{tangle}
\object{H}\step[2]\object{A}\step[3]\object{A}\\
\id\step\cd\step\cd\\
\id\step\id\step[2]\hx\step[2]\id\\
\id\step\cu\step\tu {\mu}\\
\tu {\beta}\step[3]\id\\
\step\Cu\\
\step[3]\object{H}\\
\end{tangle}
\step=\step
\begin{tangle}
\object{H}\step[3]\object{A}\step[3]\object{A}\\
\id\step[2]\cd\step\cd\\
\x\step[2]\id\step\id\step[2]\id\\
\id\step[2]\tu {\beta}\step\id\step[2]\id\\
\se1\step[2]\x\step[2]\id\\
\step\tu {\mu}\step[2]\tu {\beta}\\
\step[2]\Cu\\
\step[4]\object{H}\\
\end{tangle}
\step,\step
\]
\[
\begin{tangle}
\object{H}\\
\id\step[2]\Q {\eta_A}\\
\tu {\beta}\\
\step\object{H}\\
\end{tangle}
\step=\step
\begin{tangle}
\object{H}\\
\id\\
\id\\
\object{H}\\
\end{tangle}
\step.\step
\]

$(A, \phi )$   is called a twisted $H$-comodule   if the following
conditions are satisfied:
\[
(TCM):
\begin{tangle}
\step[3]\object{A}\\
\step\Cd\\
\td {\phi}\step[2]\td {P}\\
\id\step[2]\x\step[2]\id\\
\cu\step\td {\phi}\step\id\\
\step\id\step[2]\id\step[2]\hx\\
\step\id\step[2]\cu\step\id\\
\step\object{H}\step[3]\object{H}\step[2]\object{A}\\
\end{tangle}
\step=\step
\begin{tangle}
\step[3]\object{A}\\
\step\Cd\\
\step\id\step[3]\td {\phi}\\
\td {P}\step\cd\step\id\\
\id\step[2]\hx\step[2]\id\step\id\\
\cu\step\cu\step\id\\
\step\object{H}\step[3]\object{H}\step[2]\object{A}\\
\end{tangle}
\step,\step
\]
\[
\begin{tangle}
\step\object{A}\\
\td {\phi}\\
\QQ {\epsilon _H}\step[2]\id\\
\step[2]\object{A}\\
\end{tangle}
\step=\step
\begin{tangle}
\object{A}\\
\id\\
\id\\
\object{A}\\
\end{tangle}
\step.\step
\]

$(H, \psi )$   is called a twisted $H$-comodule   if the following
conditions are satisfied:
\[
(TCM):
\begin{tangle}
\step[3]\object{H}\\
\step\Cd\\
\td {\psi}\step[3]\id\\
\id\step\cd\step\td {Q}\\
\id\step\id\step[2]\hx\step[2]\id\\
\id\step\cu\step\cu\\
\object{H}\step[2]\object{A}\step[3]\object{A}\\
\end{tangle}
\step=\step
\begin{tangle}
\step[3]\object{H}\\
\step\Cd\\
\td {Q}\step[2]\td {\psi}\\
\id\step[2]\x\step[2]\id\\
\id\step\td {\psi}\step\id\step[2]\id\\
\hx\step[2]\id\step\id\step[2]\id\\
\id\step\cu\step\cu\\
\object{H}\step[2]\object{A}\step[3]\object{A}\\
\end{tangle}
\step,\step
\]
\[
\begin{tangle}
\step\object{H}\\
\cd\\
\id\step[2]\QQ {\epsilon _A}\\
\object{H}\\
\end{tangle}
\step=\step
\begin{tangle}
\object{H}\\
\id\\
\id\\
\object{H}\\
\end{tangle}
\step.\step
\]

$H$   is said to be cocommutative with respect to $(A, \alpha )$ if
the following conditions are satisfied:
\[
(CC):
\begin{tangle}
\step\object{H}\step[2]\object{A}\\
\cd\step\id\\
\id\step[2]\hx\\
\tu {\alpha}\step\id\\
\step\x\\
\step\object{H}\step[2]\object{A}\\
\end{tangle}
\step=\step
\begin{tangle}
\step\object{H}\step[3]\object{A}\\
\cd\step[2]\id\\
\id\step[2]\tu {\alpha}\\
\object{H}\step[3]\object{A}\\
\end{tangle}
\step.\step
\]

$A$   is said to be cocommutative with respect to $(H, \beta )$ if
the following conditions are satisfied:
  \[
(CC):
\begin{tangle}
\object{H}\step[3]\object{A}\\
\id\step[2]\cd\\
\tu {\beta}\step[2]\id\\
\step\object{H}\step[3]\object{A}\\
\end{tangle}
\step=\step
\begin{tangle}
\object{H}\step[2]\object{A}\\
\id\step\cd\\
\hx\step[2]\id\\
\id\step\tu {\beta}\\
\x\\
\object{H}\step[2]\object{A}\\
\end{tangle}
\step.\step
\]

$H$   is said to be commutative with respect to $(A, \phi )$ if the
following conditions are satisfied:

\[
(C):
\begin{tangle}
\step\object{H}\step[2]\object{A}\\
\step\id\step[2]\id\\
\step\x\\
\td {\phi}\step\id\\
\id\step[2]\hx\\
\cu\step\id\\
\step\object{H}\step[2]\object{A}\\
\end{tangle}
\step=\step
\begin{tangle}
\object{H}\step[3]\object{A}\\
\id\step[2]\td {\phi}\\
\cu\step[2]\id\\
\step\object{H}\step[3]\object{A}\\
\end{tangle}
\step.\step
\]

$A$   is called commutative with respect to $(H, \psi )$ if the
following conditions are satisfied:

\[
(C):
\begin{tangle}
\object{H}\step[2]\object{A}\\
\id \step[2]\id      \step    \\
\x   \step  \\
\id\step   \td {\psi}  \\
\hx\step[2]   \id  \\
\id\step     \cu\\
\object{H}\step[2]\object{A}     \\
\end{tangle}
\step=\step
\begin{tangle}
\step \object{H}\step[3]\object{A}     \\
\td {\psi}\step[2]    \id \\
\id\step[2]   \cu  \\
\object{H}\step[3]\object{A}     \\
\end{tangle}
\step.\step
\]

If $(H, \alpha )$ weakly acts on $A$ such that $(A, \alpha)$ becomes
a left $H$-module, then  we say that $(H, \alpha )$
 acts on $A$. Dually, if $(H, \phi )$ weakly coacts on $A$
such that $(A, \phi )$ become a left $H$-comodule, then  we say that
$(H, \phi )$
 coacts on $A$.
 An action of $H$ on $A$ is called inner if  there exists a morphism
 $u: H \rightarrow A$  with a convolution inverse $u^{-1}$ such that
\[
\begin{tangle}
\object{H}\step[2]\object{A}\\
\id\step[2]\id\\
\id\step[2]\id\\
\tu {\alpha}\\
\step\id\\
\step\object{A}\\
\end{tangle}
\step=\step
\begin{tangle}
\step\object{H}\step[3]\object{A}\\
\cd\step[2]\id\\
\id\step[2]\x\\
\O {u}\step[2]\id\step[1]\morph {\bar u }\\
\se1\step\cu\\
\step\cu\\
\step[2]\object{A}\\
\end{tangle}
\step.\step
\]

Let \[
\begin{tangle}
\object{H}\step[2]\object{H}\\
\id\step[2]\id\\
\id\step[2]\id\\
\tu {ad}\\
\step\id\\
\step\object{H}\\
\end{tangle}
\step=\step
\begin{tangle}
\step\object{H}\step[3]\object{H}\\
\cd\step[2]\id\\
\id\step[2]\x\\
\id\step[2]\id\step[2]\S\\
\se1\step\cu\\
\step\cu\\
\step[2]\object{H}\\
\end{tangle}
\step.\step
\]
It is clear that $(H, ad)$ acts on itself. This action $ad$  is
called a (left) adjoint action. Dually we can define an (left)
adjoint coaction as follows:
\[
\begin{tangle}
\step\object{H}\\
\step\id\\
\td {ad}\\
\id\step[2]\id\\
\id\step[2]\id\\
\object{H}\step[2]\object{H}\\
\end{tangle}
\step=\step
\begin{tangle}
\step\object{H}\\
\cd\\
\id\step\cd\\
\id\step\id\step[2]\S\\
\id\step\x\\
\hcu\step[2]\id\\
\step[0.5]\object{H}\step[2.5]\object{H}\\
\end{tangle}
\step.\step
\]

Bialgebra $A$ is called an $H$-module bialgebra if $A$  is an
$H$-module algebra and $H$-module coalgebra . Bialgebra $A$ is
called an $H$-comodule bialgebra if $A$ is an $H$-comodule algebra
and an $H$-comodule coalgebra.

We define the relations (BB1)--(BB11)  and (BB1')--(BB11')  as
follows:

\[
(BB1):\step
\begin{tangle}
\step\object{H}\step[3]\object{A}\\
\step\id\step[2]\cd\\
\step\x\step[2]\id\\
\step\id\step[2]\tu {\beta}\\
\td {P}\step\cd\\
\id\step[2]\hx\step[2]\id\\
\cu\step\cu\\
\step\object{H}\step[3]\object{H}\\
\end{tangle}
\step=\step
\begin{tangle}
\step\object{H}\step[4]\object{A}\\
\step\id\step[3]\cd\\
\cd\step\ne1\step\cd\\
\id\step[2]\hx\step\td {\phi}\step\id\\
\tu {\beta}\step\hx\step[2]\id\step\id\\
\step\cu\step\tu {\beta}\td {P}\\
\step[2]\se1\step[2]\hx\step[2]\id\\
\step[3]\cu\step\cu\\
\step[4]\object{H}\step[3]\object{H}\\
\end{tangle}
\step,\step
\begin{tangle}
\object{H}\step[2]\object{A}\\
\tu {\beta}\\
\step\QQ {\epsilon }\\
\end{tangle}
\step=\step
\begin{tangle}
\object{H}\step[2]\object{A}\\
\id\step[2]\id\\
\QQ {\epsilon }\step[2]\QQ {\epsilon }\\
\end{tangle}
\step;\step
\]\\

\[
(BB2):\step
\begin{tangle}
\step\object{A}\step[3]\object{A}\\
\cd\step\cd\\
\id\step[2]\hx\step[2]\se1\\
\tu {\mu}\td {\phi}\step\cd\\
\step\id\step\id\step[2]\hx\step\td {\phi}\\
\step\id\step\tu {\beta}\step\hx\step[2]\id\\
\step\cu\step\ne1\step\cu\\
\step[2]\cu\step[3]\id\\
\step[3]\object{H}\step[4]\object{A}\\
\end{tangle}
\step=\step
\begin{tangle}
\step\object{A}\step[3]\object{A}\\
\cd\step\cd\\
\id\step[2]\hx\step[2]\id\\
\cu\step\tu {\mu}\\
\td {\phi}\step[2]\id\\
\id\step[2]\x\\
\cu\step[2]\id\\
\step\object{H}\step[3]\object{A}\\
\end{tangle}
\step,\step
\begin{tangle}
\step\Q {\eta_A}\\
\td {\phi}\\
\object{H}\step[2]\object{A}\\
\end{tangle}
\step=\step
\begin{tangle}
\Q {\eta_H}\step[2]\Q {\eta_A}\\
\id\step[2]\id\\
\object{H}\step[2]\object{A}\\
\end{tangle}
\step;\step
\]\\

\[
(BB3):\step
\begin{tangle}
\step\object{H}\step[3]\object{A}\\
\step\id\step[2]\cd\\
\step\x\step[2]\id\\
\td {\phi}\step\tu {\beta}\\
\id\step[2]\x\\
\cu\step[2]\id\\
\step\object{H}\step[3]\object{A}\\
\end{tangle}
\step=\step
\begin{tangle}
\object{H}\step[3]\object{A}\\
\id\step[2]\cd\\
\tu {\beta}\step\td {\phi}\\
\step\cu\step[2]\id\\
\step[2]\object{H}\step[3]\object{A}\\
\end{tangle}
\step;\step
\]\\

\[
(BB4):\step
\begin{tangle}
\step\object{A}\step[3]\object{A}\\
\cd\step\cd\\
\id\step[2]\hx\step[2]\id\\
\cu\step\tu \mu\\
\td P\step\cd\\
\id\step[2]\hx\step[2]\id\\
\cu\step\cu\\
\step\object{H}\step[3]\object{H}\\
\end{tangle}
\step=\step
\begin{tangle}
\step[3]\object{A}\step[13]\object{A}\\
\step[2]\cd\step[11]\cd\\
\step\ne1\step\cd\step[9]\cd\step\se1\\
\ne1\step\ne1\step\td {P}\step[7]\cd\step\se1\step[2]\se2\\
\id\step\td {\phi}\step\id\step[3]\se2\step[4]\cd\step\se1\step[2]\se2\step[2]\se2\\
\id\step\id\step[2]\hx\step[5]\se2\step\cd\step\id\step\td {\phi}\step\td {\phi}\step\td {P}\\
\id\step\id\step[2]\id\step\se1\step[5]\hx\step[2]\id\step\id\step\id\step[2]\hx\step[2]\hx\step[2]\id\\
\id\step\id\step[2]\se1\step[2]\se2\step[2]\ne1\step\x\step\id\step\id\step[2]\id\step\x\step\se1\step\id\\
\id\step\se1\step[2]\se1\step[2]\x\step[2]\id\step[2]\hx\step\id\step[2]\id\step\id\step[2]\id\step[2]\id\step\id\\
\id\step[2]\se1\step[2]\x\step[2]\x\step[2]\id\step\hx\step[2]\id\step\id\step[2]\id\step[2]\id\step\id\\
\se1\step[2]\x\step[2]\x\step[2]\x\step\id\step\x\step\id\step[2]\id\step[2]\id\step\id\\
\step\tu {\mu}\step[2]\tu {\beta}\step[2]\tu {\beta}\step[2]\hx\step\id\step[2]\hx\step[2]\id\step[2]\id\step\id\\
\step[2]\id\step[4]\id\step[4]\id\step[2]\ne1\step\hx\step\ne1\step\x\step[2]\id\step\id\\
\step[2]\id\step[4]\id\step[4]\id\step[2]\id\step\ne1\step\hx\step[2]\id\step[2]\tu
{\beta}\step\id\\
\step[2]\se1\step[3]\id\step[4]\id\step[2]\id\step\cu\step\tu
{\mu}\step[3]\cu\\
\step[3]\se1\step[2]\se1\step[3]\se1\step\cu\step[3]\se1\step[4]\id\\
\step[4]\se1\step[2]\se1\step[3]\cu\step[5]\Cu\\
\step[5]\se1\step[2]\Cu\step[8]\id\\
\step[6]\Cu\step[10]\id\\
\step[8]\object{H}\step[12]\object{H}\\
\end{tangle}
\step;\step
\]\\

\[
(BB5):\step
\begin{tangle}
\step\Q {\eta_A}\\
\td {P}\\
\object{H}\step[2]\object{H}\\
\end{tangle}
\step=\step
\begin{tangle}
\Q {\eta_H}\step[2]\Q {\eta_H}\\
\id\step[2]\id\\
\object{H}\step[2]\object{H}\\
\end{tangle}
\step,\step \step[5]
\begin{tangle}
\object{A}\step[2]\object{A}\\
\tu {\mu}\\
\step\QQ {\epsilon }\\
\end{tangle}
\step=\step
\begin{tangle}
\object{A}\step[2]\object{A}\\
\id\step[2]\id\\
\QQ {\epsilon }\step[2]\QQ {\epsilon }\\
\end{tangle}
\step.\step
\]

\[
(BB8):\step[2]
\begin{tangle}
\step[1]\object{H}\step[4]\object{A}\\
\cd\step[2]\td P\\
\id\step[2]\x\step[2]\id\\
\cu \step[2]\cu \\
\step[1]\object{H}\step[4]\object{H}
\end{tangle}
\step=\step
\begin{tangle}
\step[1]\object{H}\step[2]\object{A}\\
\step[1]\x\\
\step[1]\id\step[2]\nw2\\
\td P\step[2]\cd\\
\id\step[2]\x\step[2]\id\\
\cu \step[2]\cu \\
\step[1]\object{H}\step[4]\object{H}
\end{tangle} \ \ ;
\]

\[
(BB9):\step[2]
\begin{tangle}
\step[1]\object{H}\step[4]\object{A}\\
\td Q\step[2]\cd\\
\id\step[2]\x\step[2]\id\\
\cu \step[2]\cu \\
\step[1]\object{A}\step[4]\object{A}
\end{tangle}
\step=\step
\begin{tangle}
\step[1]\object{H}\step[2]\object{A}\\
\step[1]\x\\
\step[1]\id\step[2]\nw2\\
\cd\step[2]\td Q\\
\id\step[2]\x\step[2]\id\\
\cu \step[2]\cu \\
\step[1]\object{A}\step[4]\object{A}
\end{tangle} \ \ ;
\]

\[
(BB1'):\step
\begin{tangle}
 \step[4]\object{H}\step[3]\object{H}         \step                      \\
\step[3]\cd\step   \cd                        \\

 \step[2]\ne1\step[2]\hx \step[2]                      \id     \\

\step[1]\cd\step \td {\psi }              \tu {\sigma}            \\

\td {\psi }\step \hx \step[2] \id\step   \id           \step  \\

\id  \step[2] \hx  \step \tu {\alpha} \step\id        \step  \\

\cu   \step \nw1\step    \cu    \step        \\
 \step[1]\id\step[3]  \cu  \step[2]                        \\
  \step[1]\object{H} \step[4]  \object{A}       \step[3]                 \\

\end{tangle}
\step=\step
\begin{tangle}
\step[1]\object{H}\step[3] \object{H}            \step         \\
\cd\step        \cd          \\
\id \step[2]\hx \step[2]     \id              \\
\tu {\sigma} \step   \cu                  \\
\step[1]\id\step[2]  \td {\psi }                      \\
\step[1]\x\step[2]        \id             \\
\step[1]\id \step[2]    \cu                    \\
\step[1]\object{H}\step[3] \object{A}   \step               \\
\end{tangle}
\step,\step
\begin{tangle}

\step\Q {\eta_H}           \step[3]                  \\
\td {\psi }           \step[3]                 \\
\object{H}\step[2]\object{A}           \step[3]                 \\
\end{tangle}
\step=\step
\begin{tangle}
\Q {\eta_H}\step[2]\Q {\eta_A}           \step[3]                   \\
\id\step[2]\id           \step[3]             \\
\object{H}\step[2]\object{A}                          \\
\end{tangle}
\step;\step
\]

\[
(BB2'):\step
\begin{tangle}

\step  \object{H}  \step[3]  \object{A}\\

\cd\step[2]\id\\

\id\step[2]\x\step\\

 \tu {\alpha}\step[2]\id  \step       \\

\cd  \step  \td {Q}                      \\

\id   \step[2]    \hx  \step[2]           \id       \\

 \cu   \step \cu       \\

 \step \object{A}  \step[3] \object{A}  \step       \\
\end{tangle}
\step=\step
\begin{tangle}
\step[3]\object{H}\step[4]\object{A}\step[3]     \\

\step[2] \cd \step[3]\id \step \\

\step[1]\cd\step\nw1\step \cd   \\

\step[1]\id\step \td {\psi }\step    \hx \step[2]  \id  \\

\step[1]\id\step \id \step[2]\hx  \step  \tu {\alpha}\\

\td {Q}\tu {\alpha}\step\cu \step  \\

\id\step[2]\hx\step[2]\ne1   \step[2] \\

\cu\step\cu  \step[3]  \\

\step\object{A}\step[3]\object{A}   \step[4] \\

\end{tangle}
\step,\step
\begin{tangle}

\object{H}\step[2]\object{A}\\
\tu {\alpha}\\
\step\QQ {\epsilon }\\
\end{tangle}
\step=\step
\begin{tangle}
\object{H}\step[2]\object{A}\\
\id\step[2]\id\\
\QQ {\epsilon }\step[2]\QQ {\epsilon }\\
\end{tangle}
\step;\step
\]\step[3]           \\

\[
(BB3'):\step
\begin{tangle}
\step\object{H}\step[3]\object{A}                   \\
\cd  \step[2]\id       \step        \\
\id\step[2]\x     \step               \\
\tu {\alpha}\step   \td {\psi }              \\
\step\x\step[2]    \id              \\
\step\id \step[2]       \cu           \\
\step\object{A}\step[3]\object{H}        \step            \\
\end{tangle}
\step=\step
\begin{tangle}
\step\step\object{H}\step[3]\object{A}                    \\
\step\cd\step[2]     \id              \\
\td {\psi }\step    \tu {\alpha}           \\
\id \step[2]\cu        \step         \\
\object{H}\step[3]\object{A}                   \\
\end{tangle}
\step;\step
\]

\[
(BB4'):\step
\begin{tangle}
\step\object{H}\step[3]\object{H}              \\
\cd\step        \cd         \\
\id\step[2]\hx\step[2]          \id          \\
\tu \sigma \step\cu                  \\
\cd\step\cd                     \\
\id\step[2]\hx\step[2]\id                   \\
\cu\step\tu Q                     \\
\step\object{A}\step[3]\object{A}                   \\
\end{tangle}
\step=\step
\begin{tangle}
 \step[7]  \object{H}\step[13]\object{H}    \step[3]                      \\
 \step[6]\cd\step[11]\cd              \step[2]              \\
 \step[5]\ne1\step \cd \step[9]\cd  \step \nw1         \step               \\
\step[4]\ne2\step[1]\ne1\step \cd\step[7] \td {Q}\step\nw1 \step          \nw1                \\

\step[2]  \ne2\step[2]\ne2\step[1]\ne1\step \cd \step[5]\ne2\step[2]\id \step \td {\psi }\step            \id          \\

\td {Q}\step\td {\psi }\step\td {\psi
}\step\id\step\cd\step[2]\ne2\step[4]\hx\step[2]    \id\step\id
           \\

\id
\step[2]\hx\step[2]\hx\step[2]\id\step\id\step\id\step[2]\hx\step[5]\ne1\step\id\step[2]\id
\step\id
 \\

\id\step\ne1\step\x\step\id\step[2]\id\step\id\step\x\step\nw1\step[3]\ne2\step[1]\ne1\step[2]\id
\step\id
    \\

\id\step\id\step[2]\id\step[2]\id\step\id\step[2]\id\step\hx\step[2]\id\step[2]
\x\step[2]\ne1\step[2]\ne1\step
 \id
      \\

\id\step\id\step[2]\id\step[2]\id\step\id\step[2]\hx\step\id\step[2]\x\step[2]
\x\step[2]\ne1\step[2] \id
       \\

\id\step\id\step[2]\id\step[2]\id\step\x\step\id\step\x\step[2]\x\step[2]
\x\step[2]\ne1
\\

\id\step\id\step[2]\id\step[2]\hx\step[2]\id \step\hx\step[2]\tu
{\alpha}\step[2]\tu {\alpha}\step[2]\tu {\sigma}   \step
            \\

\id\step\id\step[2]\x\step\nw1\step\hx\step\nw1\step[2]\id\step[4]\id\step[4]\id
\step[2]
  \\

\id\step\tu
{\alpha}\step[2]\id\step[2]\hx\step\nw1\step\id\step[2]\id\step[4]\id\step[4]
\id\step[2]
     \\

\cu\step[3]\tu
{\sigma}\step\cu\step\id\step[2]\id\step[4]\id\step[3]\ne1
\step[2]    \\

 \step[1]\id\step[4]\ne1\step[3]\cu \step\ne1 \step[3]\ne1\step[2] \ne1         \step[3]
                \\

 \step[1]\Cu\step[5]\cu\step[3]\ne1 \step[2] \ne1                    \step[4]      \\

 \step[3]\id\step[8]\Cu\step[2]
\ne1 \step[5]                  \\

\step[3]\id   \step[10] \Cu        \step[6]                \\

 \step[3] \object{A}\step[12]\object{A}    \step[8]                       \\

\end{tangle}
\step;\step
\]

\[
(BB5'):\step
\begin{tangle}
\step\Q {\eta_H}                  \\
\td {P}                 \\
\object{A}\step[2]\object{A}                  \\
\end{tangle}
\step=\step
\begin{tangle}
\Q {\eta_A}\step[2]\Q {\eta_A}                 \\
\id\step[2]\id                \\
\object{A}\step[2]\object{A}                 \\

\end{tangle}
\step,\step \step[5]
\begin{tangle}

\object{H}\step[2]\object{H}                  \\
\tu {\sigma}                  \\
\step\QQ {\epsilon }                  \\
\end{tangle}
\step=\step
\begin{tangle}
\object{H}\step[2]\object{H}                  \\
\id\step[2]\id                  \\
\QQ {\epsilon }\step[2]\QQ {\epsilon }                 \\
\end{tangle}
\step.\step
\]

\[
(BB8'):\step[2]
\begin{tangle}
\step[1]\object{H}\step[4]\object{H}\\
\cd\step[2]\cd\\
\id\step[2]\x\step[2]\id\\
\cu \step[2]\tu \sigma \\
\step[1]\object{H}\step[4]\object{A}
\end{tangle}
\;=\enspace
\begin{tangles}{clr}
\step[1]\object{H}\step[4]\object{H}\\
\cd\step[2]\cd\\
\id\step[2]\x\step[2]\id\\
\tu \sigma \step[2]\cu\\
\step[1]\nw1\step[2]\ne1\step[1]\\
\step[2]\x\step[2]\\
\step[2]\object{H}\step[2]\object{A}\step[2]
\end{tangles} \ \ ;
\]

\[
(BB9'):\step[2]
\begin{tangle}
\step[1]\object{A}\step[4]\object{A}\\
\cd\step[2]\cd\\
\id\step[2]\x\step[2]\id\\
\tu \mu  \step[2]\cu \\
\step[1]\object{H}\step[4]\object{A}
\end{tangle}
\;=\enspace
\begin{tangles}{clr}
\step[1]\object{A}\step[4]\object{A}\\
\cd\step[2]\cd\\
\id\step[2]\x\step[2]\id\\
\cu \step[2]\tu \mu\\
\step[1]\nw1\step[2]\ne1\step[1]\\
\step[2]\x\step[2]\\
\step[2]\object{H}\step[2]\object{A}\step[2]
\end{tangles} \ \ ;
\]

For convenience,  $(BB6')$ and $(B1) $ be the same; $(BB7')$ and
$(B6) $ be the same; $(BB10')$ be the same as $(CB4) $ and $ \psi
\eta_ H = \eta _H \otimes \eta _A$; $(BB11')$ be the same as $(B2) $
and $\epsilon \alpha = \epsilon \otimes \epsilon $. Let $(BB6)$ and
$(B3) $ are the same; $(BB7)$ and $(B6') $ are the same;  $(BB10)$
is the same as $(B4) $ and $\epsilon \beta = \epsilon \otimes
\epsilon $; $(BB11)$ is the same as $(CB2) $ and $ \phi \eta _A =
\eta _H \otimes \eta _A$. Here $(B_1), (B_6)$ and so on are defined
in Chapter \ref {c3}.

\[
m_{A_{\alpha,\sigma}\# H}\step=\step
\begin{tangle}
\object{A}\step[2]\object{H}\step[3]\object{A}\step[3]\object{H}\\
\id\step\cd\step[2]\id\step[3]\id\\
\id\step\id\step[2]\x\step[3]\id\\
\id\step\tu {\alpha}\step\cd\step\cd\\
\cu\step[2]\id\step[2]\hx\step[2]\id\\
\step\id\step[3]\tu {\sigma}\step\cu\\
\step\Cu\step[3]\id\\
\step[3]\object{A}\step[5]\object{H}\\
\end{tangle} \ \ ;
\]
\[
m_{A\#_{\beta,\mu} H}\step=\step
\begin{tangle}
\step\object{A}\step[3]\object{H}\step[3]\object{A}\step[2]\object{H}\\
\step\id\step[3]\id\step[2]\cd\step\id\\
\step\id\step[3]\x\step[2]\id\step\id\\
\cd\step\cd\step\tu {\beta}\step\id\\
\id\step[2]\hx\step[2]\id\step[2]\cu\\
\cu\step\tu {\mu}\step[3]\id\\
\step\id\step[3]\Cu\\
\step\object{A}\step[5]\object{H}\\
\end{tangle} \ \ ;
\]
\[
\Delta_{A^{\phi,P}\# H}\step=\step
\begin{tangle}
\step[3]\object{A}\step[5]\object{H}\\
\step\Cd\step[3]\id\\
\cd\step[2]\td {P}\step\cd\\
\id\step\td {\phi}\step\id\step[2]\hx\step[2]\id\\
\id\step\id\step[2]\id\step\cu\step\cu\\
\id\step\id\step[2]\x\step[3]\id\\
\id\step\cu\step[2]\id\step[3]\id\\
\object{A}\step[2]\object{H}\step[3]\object{A}\step[3]\object{H}\\
\end{tangle} \ \ ;
\]
\[
\Delta_{A\#^{\psi,Q} H}\step=\step
\begin{tangle}
\step\object{A}\step[5]\object{H}\\
\step\id\step[3]\Cd\\
\cd\step\td {Q}\step[2]\cd\\
\id\step[2]\hx\step[2]\id\step\td {\psi}\step\id\\
\cu\step\cu\step\id\step[2]\id\step\id\\
\step\id\step[3]\x\step[2]\id\step\id\\
\step\id\step[3]\id\step[2]\cu\step\id\\
\step\object{A}\step[3]\object{H}\step[3]\object{A}\step[2]\object{H}\\
\end{tangle}\ \ ;
\]
$\eta _{A _{\alpha , \sigma } \# H} = \eta _A \otimes \eta _H= \eta
_{A \# _{\beta , \mu } H }$   and $\epsilon _{A ^{\phi , P} \# H} =
\epsilon _A \otimes \epsilon _H = \epsilon _{A \# ^{\psi , Q} H}.$
$(A _{\alpha , \sigma  } \# H,  $ \ $ m_{A_ {\alpha , \sigma } \# H
}, \eta _{A _{\alpha , \sigma } \# H})$  \ \   and  \ \   $(A  \#
_{\beta , \mu } H,  m_{A \# _{\beta , \mu } H },$ \ $ \eta _{A  \#
_{\beta , \mu} H})$ are called crossed products of $A$  and $H$. $(
A ^{\phi , P} \# H ,$ \ $ \Delta _ {A ^{\phi , P} \# H} , \epsilon _
 {A ^{\phi , P} \# H}) $ and  $( A \# ^{\psi , Q} H , \Delta _
 {A \# ^{\psi , Q} H}, \epsilon _{A \# ^{\psi , Q} H})$  are
 called crossed coproducts of $A$  and $H$.
 If $\sigma $ and $\mu$  or $P$ and $Q$  are trivial, then they are
 called the smash products or smash coproducts, written  as
 $A_\alpha \# H$, $A \# _\beta H$, $A ^\phi \# H$  and $A \# ^\psi H .$
Every smash product is called a semi-direct product in physics.
 In particular, the smash product $A _\alpha \# H$  is often  written as
 $A \# H$.
$(A _{\alpha , \sigma  } \#  ^{\psi , Q} H ,$ \ $ \Delta _ {A  \#
^{\psi , Q} H }, $ \ $ \epsilon _{A \# ^{\psi , Q} H} ) , m _{
A_{\alpha ,\sigma }  \# H } , \eta _{ A_{\alpha ,\sigma }  \# H } )$
  and
  $(A ^{\phi , P } \# _{\beta , \mu} H ,
  \Delta _ {A ^{ \phi , P} \#  H } , \epsilon _{A ^{ \phi , P} \#  H } , $
 \ \  $ m_{A  \# _{\beta , \mu} H } , \eta _ {A \# _{\beta , \mu} H } $
 are called bicrossproducts.
$(A ^\phi _\alpha \# H,$ \ $ \Delta _{A ^\phi  \# H } ,  \epsilon _
{A ^\phi  \# H}, $ \ $m_{ A  _\alpha \# H}, \eta _{ A  _\alpha \# H
} )$ and $(A \# ^ \psi _\beta H , \Delta _{A \# ^ \psi H } ,
\epsilon _ {A \# ^ \psi  H }, m_{ A \#  _\beta H}, \eta _{ A \#
_\beta H } )$
     are called biproducts.

Note that S. Majid uses different notations.

\section {Crossed Products  }\label {s5}

In this section, we give the necessary and  sufficient  conditions
for crossed (co)products to become (co)algebras .

\begin {Theorem}\label {3.1.1}
If $A$  is an algebra and $H$ is a bialgebra with $\alpha (\eta _H
\otimes id _A) = id _A$ and $\alpha (id  _H \otimes \eta _A) =
\epsilon _H \eta _A$, then $D=A _{\alpha , \sigma } \# H$ is an
algebra with unity element $\eta _D = \eta _A \otimes \eta _H$
 iff $(H, \alpha )$  acts weakly on $A$ and
$(A, \alpha )$ is a twisted $H$-module with 2-cocycle $\sigma .$

 \end {Theorem}
 {\bf Proof.}
 It is clear that       $\eta _D = \eta _A \otimes \eta _H$
 iff  $\sigma (id _H \otimes \eta _H) = \sigma (\eta _H \otimes id _H)
 =\eta _A \epsilon _H.$

 If $D$  is an algebra, then the associative law holds, i.e.
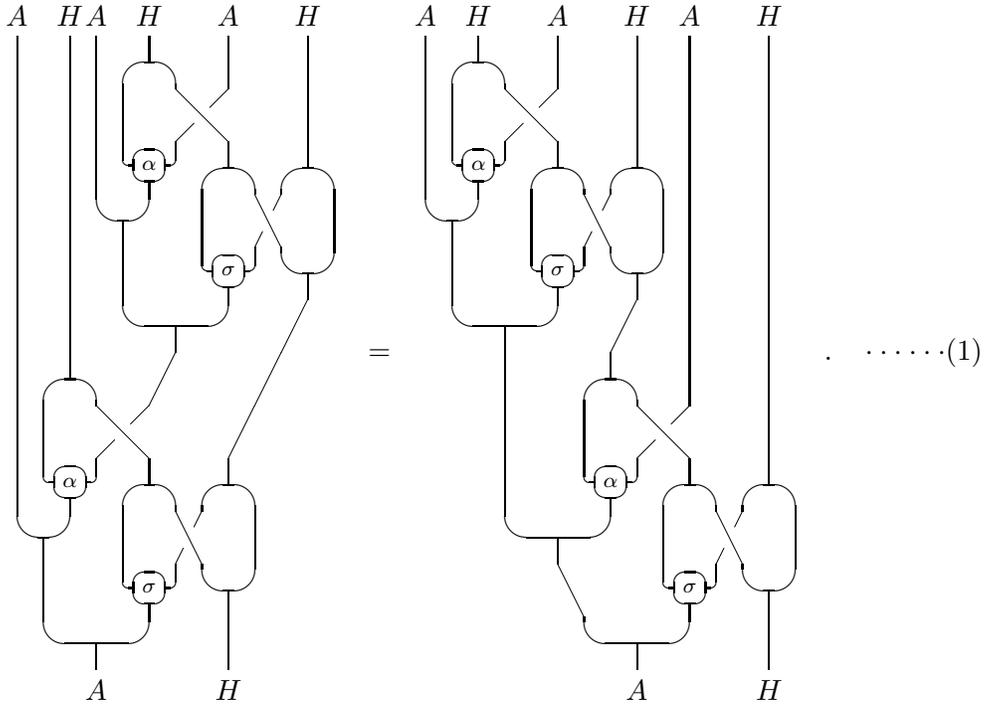
\[
\begin{tangle}
\object{A}\step[2]\object{H}\step\object{A}\step[2]\object{H}\step[3]\object{A}\step[3]\object{H}\\
\id\step[2]\id\step\id\step\cd\step[2]\id\step[3]\id\\
\id\step[2]\id\step\id\step\id\step[2]\x\step[3]\id\\
\id\step[2]\id\step\id\step\tu {\alpha}\step\cd\step\cd\\
\id\step[2]\id\step\cu\step[2]\id\step[2]\hx\step[2]\id\\
\id\step[2]\id\step[2]\id\step[3]\tu {\sigma}\step\cu\\
\id\step[2]\id\step[2]\Cu\step[2]\sw1\\
\id\step\cd\step[2]\sw1\step[3]\sw1\\
\id\step\id\step[2]\x\step[3]\sw1\\
\id\step\tu {\alpha}\step\cd\step\cd\\
\cu\step[2]\id\step[2]\hx\step[2]\id\\
\step\id\step[3]\tu {\sigma}\step\cu\\
\step\Cu\step[3]\id\\
\step[3]\object{A}\step[5]\object{H}\\
\end{tangle}
\step=\step
\begin{tangle}
\object{A}\step[2]\object{H}\step[3]\object{A}\step[3]\object{H}\step[2]\object{A}\step[3]\object{H}\\
\id\step\cd\step[2]\id\step[3]\id\step[2]\id\step[3]\id\\
\id\step\id\step[2]\x\step[3]\id\step[2]\id\step[3]\id\\
\id\step\tu {\alpha}\step\cd\step\cd\step\id\step[3]\id\\
\cu\step[2]\id\step[2]\hx\step[2]\id\step\id\step[3]\id\\
\step\id\step[3]\tu {\sigma}\step\cu\step\id\step[3]\id\\
\step\Cu\step[2]\sw1\step[2]\id\step[3]\id\\
\step[3]\id\step[3]\cd\step[2]\id\step[3]\id\\
\step[3]\id\step[3]\id\step[2]\x\step[3]\id\\
\step[3]\id\step[3]\tu {\alpha}\step\cd\step\cd\\
\step[3]\Cu\step[2]\id\step[2]\hx\step[2]\id\\
\step[5]\se1\step[3]\tu {\sigma}\step\cu\\
\step[6]\Cu\step[3]\id\\
\step[8]\object{A}\step[5]\object{H}\\
\end{tangle}\ \ \ .
\step\cdot\cdot\cdot\cdot\cdot\cdot(1)
\]

By relation (1), we have that
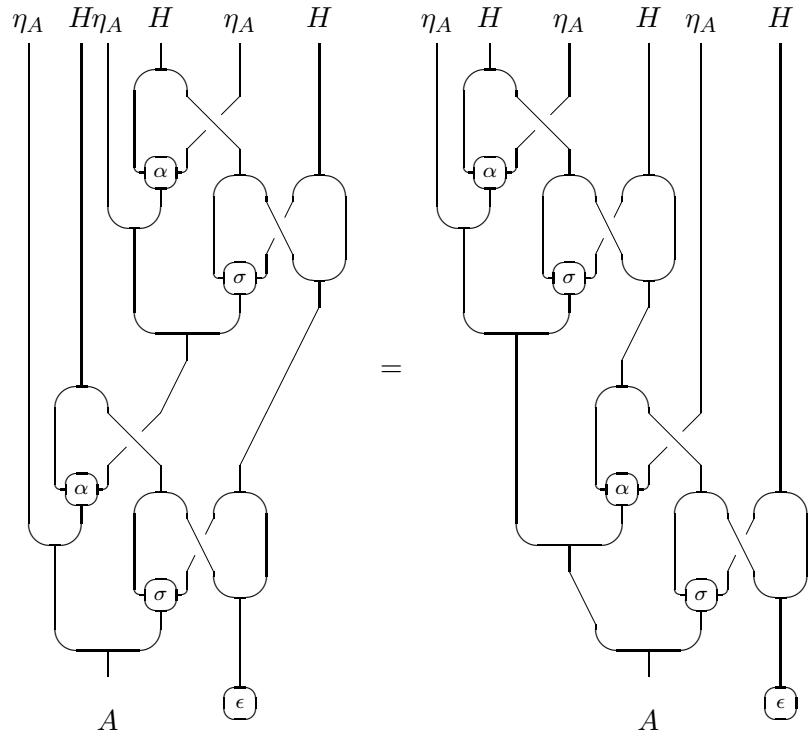
\[
\begin{tangle}
\object{\eta_A}\step[2]\object{H}\step\object{\eta_A}\step[2]\object{H}\step[3]\object{\eta_A}\step[3]\object{H}\\
\id\step[2]\id\step\id\step\cd\step[2]\id\step[3]\id\\
\id\step[2]\id\step\id\step\id\step[2]\x\step[3]\id\\
\id\step[2]\id\step\id\step\tu {\alpha}\step\cd\step\cd\\
\id\step[2]\id\step\cu\step[2]\id\step[2]\hx\step[2]\id\\
\id\step[2]\id\step[2]\id\step[3]\tu {\sigma}\step\cu\\
\id\step[2]\id\step[2]\Cu\step[2]\sw1\\
\id\step\cd\step[2]\sw1\step[3]\sw1\\
\id\step\id\step[2]\x\step[3]\sw1\\
\id\step\tu {\alpha}\step\cd\step\cd\\
\cu\step[2]\id\step[2]\hx\step[2]\id\\
\step\id\step[3]\tu {\sigma}\step\cu\\
\step\Cu\step[3]\id\\
\step[3]\object{A}\step[5]\QQ {\epsilon }\\
\end{tangle}
\step=\step
\begin{tangle}
\object{\eta_A}\step[2]\object{H}\step[3]\object{\eta_A}\step[3]\object{H}\step[2]\object{\eta_A}\step[3]\object{H}\\
\id\step\cd\step[2]\id\step[3]\id\step[2]\id\step[3]\id\\
\id\step\id\step[2]\x\step[3]\id\step[2]\id\step[3]\id\\
\id\step\tu {\alpha}\step\cd\step\cd\step\id\step[3]\id\\
\cu\step[2]\id\step[2]\hx\step[2]\id\step\id\step[3]\id\\
\step\id\step[3]\tu {\sigma}\step\cu\step\id\step[3]\id\\
\step\Cu\step[2]\sw1\step[2]\id\step[3]\id\\
\step[3]\id\step[3]\cd\step[2]\id\step[3]\id\\
\step[3]\id\step[3]\id\step[2]\x\step[3]\id\\
\step[3]\id\step[3]\tu {\alpha}\step\cd\step\cd\\
\step[3]\Cu\step[2]\id\step[2]\hx\step[2]\id\\
\step[5]\se1\step[3]\tu {\sigma}\step\cu\\
\step[6]\Cu\step[3]\id\\
\step[8]\object{A}\step[5]\QQ {\epsilon }\\
\end{tangle}
\]
and $\sigma$is a 2-cocycle from $H\otimes H$ to $A$.

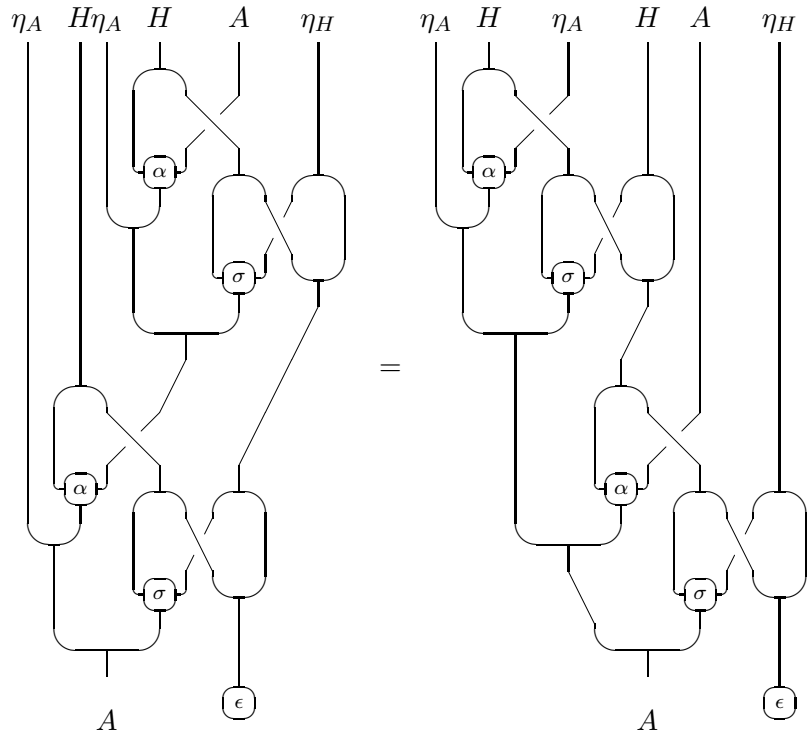
\[
\begin{tangle}
\object{\eta_A}\step[2]\object{H}\step\object{\eta_A}\step[2]\object{H}\step[3]\object{A}\step[3]\object{\eta_H}\\
\id\step[2]\id\step\id\step\cd\step[2]\id\step[3]\id\\
\id\step[2]\id\step\id\step\id\step[2]\x\step[3]\id\\
\id\step[2]\id\step\id\step\tu {\alpha}\step\cd\step\cd\\
\id\step[2]\id\step\cu\step[2]\id\step[2]\hx\step[2]\id\\
\id\step[2]\id\step[2]\id\step[3]\tu {\sigma}\step\cu\\
\id\step[2]\id\step[2]\Cu\step[2]\sw1\\
\id\step\cd\step[2]\sw1\step[3]\sw1\\
\id\step\id\step[2]\x\step[3]\sw1\\
\id\step\tu {\alpha}\step\cd\step\cd\\
\cu\step[2]\id\step[2]\hx\step[2]\id\\
\step\id\step[3]\tu {\sigma}\step\cu\\
\step\Cu\step[3]\id\\
\step[3]\object{A}\step[5]\QQ {\epsilon }\\
\end{tangle}
\step=\step
\begin{tangle}
\object{\eta_A}\step[2]\object{H}\step[3]\object{\eta_A}\step[3]\object{H}\step[2]\object{A}\step[3]\object{\eta_H}\\
\id\step\cd\step[2]\id\step[3]\id\step[2]\id\step[3]\id\\
\id\step\id\step[2]\x\step[3]\id\step[2]\id\step[3]\id\\
\id\step\tu {\alpha}\step\cd\step\cd\step\id\step[3]\id\\
\cu\step[2]\id\step[2]\hx\step[2]\id\step\id\step[3]\id\\
\step\id\step[3]\tu {\sigma}\step\cu\step\id\step[3]\id\\
\step\Cu\step[2]\sw1\step[2]\id\step[3]\id\\
\step[3]\id\step[3]\cd\step[2]\id\step[3]\id\\
\step[3]\id\step[3]\id\step[2]\x\step[3]\id\\
\step[3]\id\step[3]\tu {\alpha}\step\cd\step\cd\\
\step[3]\Cu\step[2]\id\step[2]\hx\step[2]\id\\
\step[5]\se1\step[3]\tu {\sigma}\step\cu\\
\step[6]\Cu\step[3]\id\\
\step[8]\object{A}\step[5]\QQ {\epsilon }\\
\end{tangle}
\]
and  $(A,\alpha)$ is a twisted $H$ -module.

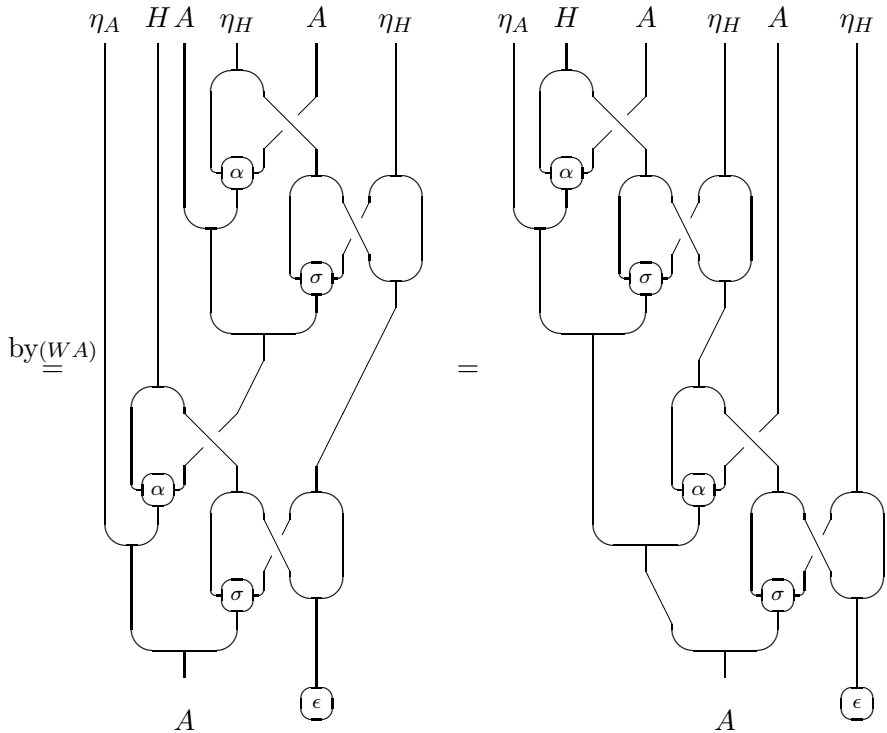
\[
\stackrel{ \hbox {by} (WA)}{=}
\begin{tangle}
\object{\eta_A}\step[2]\object{H}\step\object{A}\step[2]\object{\eta_H}\step[3]\object{A}\step[3]\object{\eta_H}\\
\id\step[2]\id\step\id\step\cd\step[2]\id\step[3]\id\\
\id\step[2]\id\step\id\step\id\step[2]\x\step[3]\id\\
\id\step[2]\id\step\id\step\tu {\alpha}\step\cd\step\cd\\
\id\step[2]\id\step\cu\step[2]\id\step[2]\hx\step[2]\id\\
\id\step[2]\id\step[2]\id\step[3]\tu {\sigma}\step\cu\\
\id\step[2]\id\step[2]\Cu\step[2]\sw1\\
\id\step\cd\step[2]\sw1\step[3]\sw1\\
\id\step\id\step[2]\x\step[3]\sw1\\
\id\step\tu {\alpha}\step\cd\step\cd\\
\cu\step[2]\id\step[2]\hx\step[2]\id\\
\step\id\step[3]\tu {\sigma}\step\cu\\
\step\Cu\step[3]\id\\
\step[3]\object{A}\step[5]\QQ {\epsilon }\\
\end{tangle}
\step=\step
\begin{tangle}
\object{\eta_A}\step[2]\object{H}\step[3]\object{A}\step[3]\object{\eta_H}\step[2]\object{A}\step[3]\object{\eta_H}\\
\id\step\cd\step[2]\id\step[3]\id\step[2]\id\step[3]\id\\
\id\step\id\step[2]\x\step[3]\id\step[2]\id\step[3]\id\\
\id\step\tu {\alpha}\step\cd\step\cd\step\id\step[3]\id\\
\cu\step[2]\id\step[2]\hx\step[2]\id\step\id\step[3]\id\\
\step\id\step[3]\tu {\sigma}\step\cu\step\id\step[3]\id\\
\step\Cu\step[2]\sw1\step[2]\id\step[3]\id\\
\step[3]\id\step[3]\cd\step[2]\id\step[3]\id\\
\step[3]\id\step[3]\id\step[2]\x\step[3]\id\\
\step[3]\id\step[3]\tu {\alpha}\step\cd\step\cd\\
\step[3]\Cu\step[2]\id\step[2]\hx\step[2]\id\\
\step[5]\se1\step[3]\tu {\sigma}\step\cu\\
\step[6]\Cu\step[3]\id\\
\step[8]\object{A}\step[5]\QQ {\epsilon }\\
\end{tangle}
\]
 and $(H,\alpha)$ acts weakly on $A$.

Conversely, if $(H,\alpha)$acts weakly on $A$ and $(A,\alpha)$ is a
twisted $H$-module with 2-cocycle $\sigma$， we show that the
associative law holds. That is, we show that relation (1) holds.

$\hbox {the left hand of } (1)\stackrel{\hbox {by }
(WA)}{=}$\quad\quad\quad
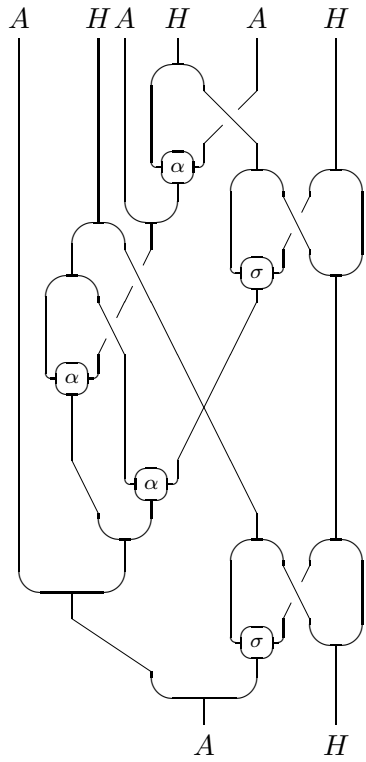
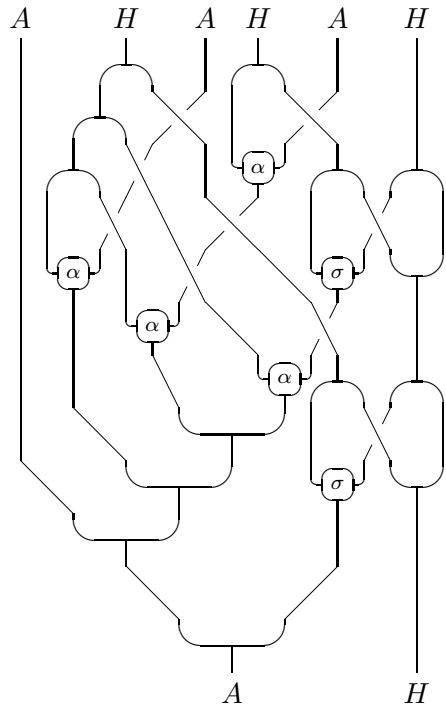
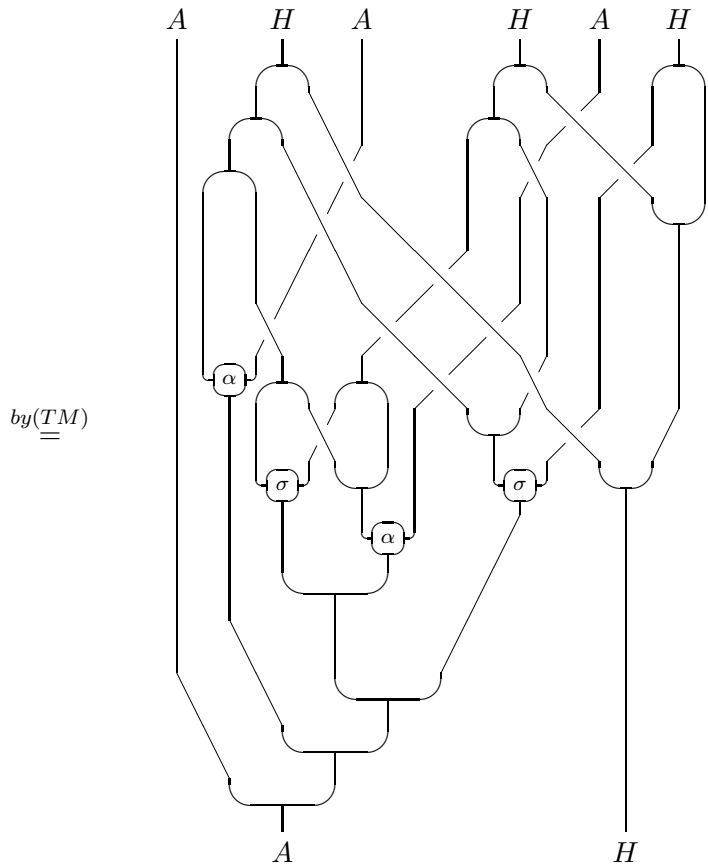
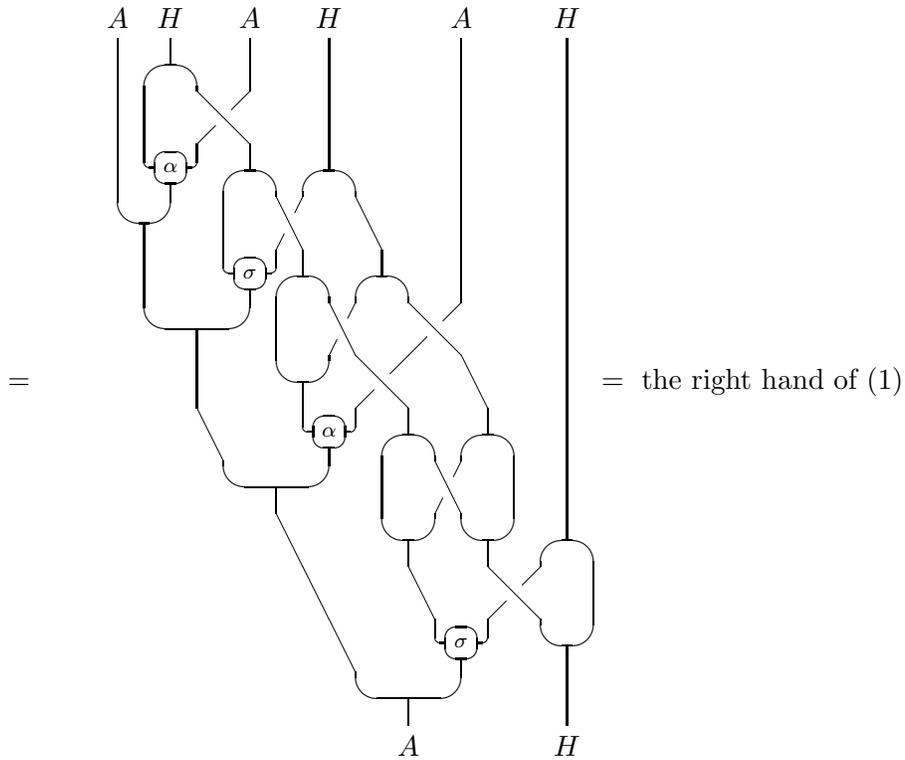
\[
\begin{tangle}
\object{A}\step[3]\object{H}\step\object{A}\step[2]\object{H}\step[3]\object{A}\step[3]\object{H}\\
\id\step[3]\id\step\id\step\cd\step[2]\id\step[3]\id\\
\id\step[3]\id\step\id\step\id\step[2]\x\step[3]\id\\
\id\step[3]\id\step\id\step\tu {\alpha}\step\cd\step\cd\\
\id\step[2]\cd\cu\step[2]\id\step[2]\hx\step[2]\id\\
\id\step\cd\step\hx\step[3]\tu {\sigma}\step\cu\\
\id\step\id\step[2]\hx\step\se1\step[2]\sw1\step[3]\id\\
\id\step\tu {\alpha}\step\id\step[2]\se1\sw1\step[4]\id\\
\id\step[2]\id\step[2]\id\step[2]\sw1\se1\step[4]\id\\
\id\step[2]\se1\step\tu {\alpha}\step[2]\se1\step[3]\id\\
\id\step[3]\cu\step[3]\cd\step\cd\\
\Cu\step[4]\id\step[2]\hx\step[2]\id\\
\step[4]\se3\step[3]\tu {\sigma}\step\cu\\
\step[5]\Cu\step[3]\id\\
\step[7]\object{A}\step[5]\object{H}\\
\end{tangle}
\]

\[
\stackrel{\hbox {by } (WA)}{=}\step[3]
\begin{tangle}
\object{A}\step[4]\object{H}\step[3]\object{A}\step[2]\object{H}\step[3]\object{A}\step[3]\object{H}\\
\id\step[3]\cd\step[2]\id\step\cd\step[2]\id\step[3]\id\\
\id\step[2]\cd\step\x\step\id\step[2]\x\step[3]\id\\
\id\step\cd\step\hx\step[2]\id\step\tu {\alpha}\step\cd\step\cd\\
\id\step\id\step[2]\hx\step\se1\step\x\step[2]\id\step[2]\hx\step[2]\id\\
\id\step\tu {\alpha}\step\id\step[2]\hx\step[3]\se2\tu
{\sigma}\step\cu\\
\id\step[2]\id\step[2]\tu
{\alpha}\step[2]\se2\step[2]\hx\step[3]\id\\
\id\step[2]\id\step[3]\se1\step[3]\tu {\alpha}\cd\step\cd\\
\id\step[3]\se2\step[2]\Cu\step\id\step[2]\hx\step[2]\id\\
\step\se2\step[2]\Cu\step[3]\tu {\sigma}\step\cu\\
\step[2]\Cu\step[6]\id\step[3]\id\\
\step[5]\se2\step[4]\sw2\step[4]\id\\
\step[6]\Cu\step[5]\id\\
\step[8]\object{A}\step[7]\object{H}\\
\end{tangle}
\]

\[
=\step[3]
\begin{tangle}
\object{A}\step[5]\object{H}\step[3]\object{A}\step[2]\object{H}
\step[3]\object{A}\step[4]\object{H}\\
\id\step[4]\cd\step[2]\id\step\cd\step[2]\id\step[4]\id\\
\id\step[3]\cd\step\se1\step\id\step\id\step[2]\x\step[4]\id\\
\id\step[2]\sw1\step\cd\step\hx\step\tu
{\alpha}\step\cd\step[2]\cd\\
\id\step\cd\step\id\step[2]\hx\step\x\step[2]\id\step[2]\x\step[2]\id\\
\vstr{50}\id\step\id\step[2]\id\step\id\step[2]\id\step\id\step\id\step[2]\id\step[2]\id\step[2]\se1\step\id\step[2]\id\\
\id\step\id\step[2]\id\step\x\step\hx\step[2]\id\step\cd\step\cd\cu\\
\id\step\id\step[2]\hx\step[2]\hx\step\se1\step\id\step\id\step[2]\hx\step[2]\id\step\id\\
\id\step\tu {\alpha}\step\tu
{\alpha}\step\id\step[2]\id\step\id\step\tu {\sigma}\step\cu\step\id\\
\id\step[2]\id\step[3]\id\step[2]\id\step[2]\id\step\x\step[2]\ne1\step[2]\id\\
\se1\step\se1\step[2]\id\step[2]\id\step[2]\hx\step[2]\x\step[2]\ne1\\
\step\se1\step\se1\step\se1\step\tu {\alpha}\step\tu
{\sigma}\step[2]\cu\\
\step[2]\se1\step\se1\step\se1\step\se1\step[2]\id\step[4]\id\\
\step[3]\se1\step\se1\step\se1\step\cu\step[4]\id\\
\step[4]\se1\step\se1\step\cu\step[5]\id\\
\step[5]\se1\step\cu\step[6]\id\\
\step[6]\cu\step[7]\id\\
\step[7]\object{A}\step[8]\object{H}\\
\end{tangle}
\]

\[
\stackrel{ \hbox {by (2-}COC)}{=}\step[3]
\begin{tangle}
\object{A}\step[4]\object{H}\step[3]\object{A}\step[3]\object{H}\step[3]\object{A}\step[3]\object{H}\\
\id\step[3]\cd\step[2]\id\step[2]\cd\step[2]\id\step[3]\id\\
\id\step[2]\cd\step\x\step[2]\id\step[2]\x\step[3]\id\\
\id\step\cd\step\hx\step[3]\se2\tu {\alpha}\step\cd\step\cd\\
\id\step\id\step[2]\hx\step\se1\step[3]\hx\step[2]\id\step[2]\hx\step[2]\id\\
\id\step\tu {\alpha}\step\se1\step\se1\step\ne1\step\se1\step\id\step[2]\id\step\cu\\
\id\step[2]\id\step[3]\id\step[2]\hx\step[3]\hx\step[2]\id\step[2]\id\\
\id\step[2]\id\step[3]\tu{\alpha}\cd\step\cd\se1\step\id\step[2]\id\\
\id\step[2]\id\step[4]\id\step\id\step[2]\hx\step[2]\id\step\hx\step[2]\id\\
\id\step[2]\id\step[4]\id\step\tu {\sigma}\step\cu\step\id\step\cu\\
\id\step[2]\id\step[4]\id\step[2]\id\step[3]\tu {\sigma}\step[2]\id\\
\id\step[2]\se1\step[3]\id\step[2]\Cu\step[3]\id\\
\se1\step[2]\se1\step[2]\Cu\step[5]\id\\
\step\se1\step[2]\Cu\step[7]\id\\
\step[2]\Cu\step[9]\id\\
\step[4]\object{A}\step[11]\object{H}\\
\end{tangle}
\]

\[
\stackrel{by (TM)}{=}\step[3]
\begin{tangle}
\object{A}\step[4]\object{H}\step[3]\object{A}\step[6]\object{H}\step[3]\object{A}\step[3]\object{H}\\
\id\step[3]\cd\step[2]\id\step[5]\cd\step[2]\id\step[2]\cd\\
\id\step[2]\cd\step\se1\step\id\step[4]\cd\step\x\step[2]\id\step[2]\id\\
\id\step\cd\step\se1\step\hx\step[4]\id\step[2]\hx\step[2]\x\step[2]\id\\
\id\step\id\step[2]\id\step[2]\hx\step[2]\se2\step[2]\id\step[2]\id\step\id\step[2]\id\step[2]\cu\\
\id\step\id\step[2]\id\step\ne1\step\se1\step[2]\x\step[2]\id\step\id\step[2]\id\step[3]\id\\
\id\step\id\step[2]\hx\step[3]\x\step[2]\x\step\id\step[2]\id\step[3]\id\\
\id\step\tu
{\alpha}\cd\step\cd\step\x\step[2]\hx\step[2]\id\step[3]\id\\
\id\step[2]\id\step\id\step[2]\hx\step[2]\id\step\id\step[2]\cu\step\x\step[2]\ne1\\
\id\step[2]\id\step\tu {\sigma}\step\cu\step\id\step[3]\tu
{\sigma}\step[2]\cu\\
\id\step[2]\id\step[2]\id\step[3]\tu
{\alpha}\step[3]\ne1\step[4]\id\\
\id\step[2]\id\step[2]\Cu\step[3]\ne1\step[5]\id\\
\id\step[2]\se1\step[3]\id\step[4]\ne1\step[6]\id\\
\se1\step[2]\se1\step[2]\Cu\step[7]\id\\
\step\se1\step[2]\Cu\step[9]\id\\
\step[2]\Cu\step[11]\id\\
\step[4]\object{A}\step[13]\object{H}
\end{tangle}
\]

\[
=\step[3]
\begin{tangle}
\object{A}\step[2]\object{H}\step[3]\object{A}\step[3]\object{H}\step[5]\object{A}\step[4]\object{H}\\
\id\step\cd\step[2]\id\step[3]\id\step[5]\id\step[4]\id\\
\id\step\id\step[2]\x\step[3]\id\step[5]\id\step[4]\id\\
\id\step\tu {\alpha}\step\cd\step\cd\step[4]\id\step[4]\id\\
\cu\step[2]\id\step[2]\hx\step[2]\se1\step[3]\id\step[4]\id\\
\step\id\step[3]\tu {\sigma}\cd\step\cd\step[2]\id\step[4]\id\\
\step\Cu\step\id\step[2]\hx\step[2]\x\step[4]\id\\
\step[3]\id\step[3]\cu\step\x\step[2]\se1\step[3]\id\\
\step[3]\se1\step[3]\tu {\alpha}\step\cd\step\cd\step[2]\id\\
\step[4]\Cu\step[2]\id\step[2]\hx\step[2]\id\step[2]\id\\
\step[6]\se1\step[3]\cu\step\cu\step\cd\\
\step[7]\se1\step[3]\se1\step[2]\x\step[2]\id\\
\step[8]\se1\step[3]\tu {\sigma}\step[2]\cu\\
\step[9]\Cu\step[4]\id\\
\step[11]\object{A}\step[6]\object{H}\\
\end{tangle} = \hbox { the right hand of }(1)
\]

 Thus  $      A _{\alpha , \sigma } \# H$ is an
algebra.
\begin{picture}(8,8)\put(0,0){\line(0,1){8}}\put(8,8){\line(0,-1){8}}\put(0,0){\line(1,0){8}}\put(8,8){\line(-1,0){8}}\end{picture}

Dually, we have the followings:

Theorem \ref {3.1.1}' \ If $H$  is an algebra and $A$ is a bialgebra
with $\beta (\eta _H \otimes id _A) = \epsilon _A\eta_{H}$ and
$\beta (id _H \otimes \eta _A) = id_A$,  then $A \# _{\beta , \mu
}H$is an algebra iff $(A, \beta )$  acts weakly on $H$ and $(H,
\beta )$ is a twisted $A$-module with 2-cocycle $\mu .$

Theorem \ref {3.1.1}" If $A$  is a coalgebra and $H$ is a bialgebra
with $(\epsilon _H \otimes id _A)\phi = id _A$ and $ (id  _H \otimes
\epsilon _A)\phi = \eta_H \epsilon _A$, then $A ^{\phi , P } \# H$
is a coalgebra iff $(H, \phi )$  coacts weakly on $A$ and $(A, \phi
)$ is a twisted $H$-comodule with 2-cycle $P.$

Theorem \ref {3.1.1}" \ If $H$  is a coalgebra and $A$ is a
bialgebra with $(\epsilon _H \otimes id _A)\varphi = \eta
_A\epsilon_{H}$ and $(id _H \otimes \epsilon _A) = id_H$, then $A\#
^{\psi , Q}H$
  is a coalgebra iff $(A, \psi )$ coacts weakly on $H$ and
$(H, \psi )$ is a twisted $A$-comodule with 2-cycle $Q .$

\begin {Corollary}\label {3.1.2}
(i) If $A$  is an algebra and $H$ is a bialgebra, then $A _{\alpha }
\# H$   is an algebra iff $(A, \alpha )$ is an $H$-module algebra;

(ii) If $H$  is an algebra and $A$ is a bialgebra, then $A \#
_{\beta  }H$
  is an algebra iff
$(H, \beta )$ is an $A$-module algebra;

(iii) If $A$  is a coalgebra and $H$ is a bialgebra, then $A ^{\phi
} \# H$   is a coalgebra iff $(A, \phi )$ is an $H$-comodule
coalgebra;

(iv) If $H$  is a coalgebra and $A$ is a bialgebra, then $A\# ^{\psi
}H$
  is a coalgebra iff
$(H, \psi )$ is an $A$-comodule coalgebra.

 \end {Corollary}

\section {Bicrossproducts }\label {s6}
In this section, we give the necessary and sufficient conditions for
bicrossproducts to become bialgebras.

\begin {Lemma}\label {3.2.1.1}
Let  $A$  and $H$ be  two algebras in ${\cal C}$. Let $D= A \otimes
H$ is an object of ${\cal C}$. Assume that $(D, m_D, \eta _D)$ and
$(D, \Delta _D, \epsilon _D)$ are an algebra  and a coalgebra in
${\cal C}$ respectively such that

\[
\begin{tangle}
\object{A\otimes H}\step[4]\object{\eta_A\otimes H}\\
\Cu\\
\step[2]\object{D}\\
\end{tangle}
\step=\step
\begin{tangle}
\object{A}\step\object{H}\step[2]\object{H}\\
\id\step\cu\\
\object{A}\step[2]\object{H}\\
\end{tangle}
\step,\step
\] \ \
 or \ \
\[
\begin{tangle}
\object{A\otimes \eta_H}\step[4]\object{A\otimes H}\\
\Cu\\
\step[2]\object{D}\\
\end{tangle}
\step=\step
\begin{tangle}
\object{A}\step[2]\object{A}\step\object{H}\\
\cu\step\id\\
\step\object{A}\step[2]\object{H}\\
\end{tangle}
\step.\step
\]\\

If
\[
\begin{tangle}
\object{\eta_A\otimes H}\step[4]\object{\eta_A\otimes H}\\
\Cu\\
\Cd\\
\object{D}\step[4]\object{D}\\
\end{tangle}
\step=\step
\begin{tangle}
\step\object{\eta_A\otimes H}\step[4]\object{\eta_A\otimes H}\\
\cd\step[2]\cd\\
\id\step[2]\x\step[2]\id\\
\cu\step[2]\cu\\
\step\object{D}\step[4]\object{D}\\
\end{tangle}
\step,\step
\]

\[
\begin{tangle}
\object{\eta_A\otimes H}\step[4]\object{A\otimes \eta_H}\\
\Cu\\
\Cd\\
\object{D}\step[4]\object{D}\\
\end{tangle}
\step=\step
\begin{tangle}
\step\object{\eta_A\otimes H}\step[4]\object{A\otimes \eta_H}\\
\cd\step[2]\cd\\
\id\step[2]\x\step[2]\id\\
\cu\step[2]\cu\\
\step\object{D}\step[4]\object{D}\\
\end{tangle}
\step,\step
\]

\[
\begin{tangle}
\object{A\otimes \eta_H}\step[4]\object{\eta_A\otimes H}\\
\Cu\\
\Cd\\
\object{D}\step[4]\object{D}\\
\end{tangle}
\step=\step
\begin{tangle}
\step\object{A\otimes \eta_H}\step[4]\object{\eta_A\otimes H}\\
\cd\step[2]\cd\\
\id\step[2]\x\step[2]\id\\
\cu\step[2]\cu\\
\step\object{D}\step[4]\object{D}\\
\end{tangle}
\step,\step
\]

\[
\begin{tangle}
\object{A\otimes \eta_H}\step[4]\object{A\otimes \eta_H}\\
\Cu\\
\Cd\\
\object{D}\step[4]\object{D}\\
\end{tangle}
\step=\step
\begin{tangle}
\step\object{A\otimes \eta_H}\step[4]\object{A\otimes \eta_H}\\
\cd\step[2]\cd\\
\id\step[2]\x\step[2]\id\\
\cu\step[2]\cu\\
\step\object{D}\step[4]\object{D}\\
\end{tangle}
\step,\step
\]
then
\[
\begin{tangle}
\object{D}\step[4]\object{D}\\
\Cu\\
\Cd\\
\object{D}\step[4]\object{D}\\
\end{tangle}
\step=\step
\begin{tangle}
\step\object{D}\step[4]\object{D}\\
\cd\step[2]\cd\\
\id\step[2]\x\step[2]\id\\
\cu\step[2]\cu\\
\step\object{D}\step[4]\object{D}\\
\end{tangle}
\step.\step
\]\\

 \end {Lemma}

{\bf Proof.} (i) Assume
\[
\begin{tangle}
\object{A\otimes H}\step[4]\object{\eta_A\otimes H}\\
\Cu\\
\step[2]\object{D}\\
\end{tangle}
\step=\step
\begin{tangle}
\object{A}\step\object{H}\step[2]\object{H}\\
\id\step\cu\\
\object{A}\step[2]\object{H}\\
\end{tangle}
\step.\step
\]

We first show that
\[
\begin{tangle}
\object{A\otimes H}\step[4]\object{\eta_A\otimes H}\\
\Cu\\
\Cd\\
\object{D}\step[4]\object{D}\\
\end{tangle}
\step=\step
\begin{tangle}
\step\object{A\otimes H}\step[4]\object{\eta_A\otimes H}\\
\cd\step[2]\cd\\
\id\step[2]\x\step[2]\id\\
\cu\step[2]\cu\\
\step\object{D}\step[4]\object{D}\\
\end{tangle}
\step[2]\cdot\cdot\cdot\cdot\cdot\cdot(1)
\]

\[ \hbox { the left hand of } (1)=
\begin{tangle}
\object{A\otimes \eta_H}\step[3]\object{A\otimes H}\step[4]\object{\eta_A\otimes H}\\
\se1\step[2]\Cu\\
\step\Cu\\
\step\Cd\\
\step\object{D}\step[4]\object{D}\\
\end{tangle}
\step\stackrel{ \hbox {by\ assumption }}{=}\step
\begin{tangle}
\object{A\otimes \eta_H}\step[4]\object{\eta_A}\step\object{H}\step[2]\object{H}\\
\se1\step[3]\id\step\cu\\
\step\id\step[3] \tu {\otimes } \\
\step\id\step[4]\id\\
\step\Cu\\
\step\Cd\\
\step\object{D}\step[4]\object{D}\\
\end{tangle}
\]

\[
\begin{tangle}
\object{A\otimes \eta_H}\step[3]\object{\eta_A}\step\object{H}\step[2]\object{H}\\
\step\id\step[2]\id\step\cu\\
\step\id\step[2]\tu \otimes \\
\cd\step\cd\\
\id\step[2]\hx\step[2]\id\\
\cu\step\cu\\
\step\object{D}\step[3]\object{D}\\
\end{tangle}
\step=\step
\begin{tangle}
\step[2]\object{A\otimes \eta_H}\step[4]\object{\eta_A\otimes H}\step[5]\object{\eta_A\otimes H}\\
\step[2]\id\step[2]\Cd\step\Cd\\
\step[2]\id\step[2]\id\step[4]\hx\step[4]\id\\
\Cd\Cu\step\Cu\\
\id\step[4]\x\step[4]\ne1\\
\Cu\step[2]\Cu\\
\step[2]\object{D}\step[6]\object{D}\\
\end{tangle}
\]

\[
\step=\step
\begin{tangle}
\step[2]\object{A\otimes \eta_H}\step[5]\object{\eta_A\otimes H}\step[4]\object{\eta_A\otimes H}\\
\Cd\step\Cd\step[2]\id\\
\id\step[4]\hx\step[4]\id\step[2]\id\\
\Cu\step\Cu\Cd\\
\step[2]\se1\step[4]\x\step[4]\id\\
\step[3]\Cu\step[2]\Cu\\
\step[5]\object{D}\step[6]\object{D}\\
\end{tangle}
\step=\step
\begin{tangle}
\object{A\otimes \eta_H}\step[4]\object{\eta_A\otimes H}\step[4]\object{\eta_A\otimes H}\\
\Cu\step[3]\ne1\\
\Cd\step\Cd\\
\id\step[4]\hx\step[4]\id\\
\Cu\step\Cu\\
\step[2]\object{D}\step[5]\object{D}\\
\end{tangle}
= \hbox { the right hand of } (1).\] Thus (1) holds. Next we show
that
\[
\begin{tangle}
\object{A\otimes \eta_H}\step[4]\object{A\otimes H}\\
\Cu\\
\Cd\\
\object{D}\step[4]\object{D}\\
\end{tangle}
\step=\step
\begin{tangle}
\step\object{A\otimes \eta_H}\step[4]\object{A\otimes H}\\
\cd\step[2]\cd\\
\id\step[2]\x\step[2]\id\\
\cu\step[2]\cu\\
\step\object{D}\step[4]\object{D}\\
\end{tangle}
\step\cdot\cdot\cdot\cdot\cdot\cdot(2)
\]

\[\hbox { the left hand of } (2)=
\begin{tangle}
\object{A\otimes \eta_H}\step[4]\object{A\otimes \eta_H}\step[4]\object{\eta_A\otimes H}\\
\step\se2\step[2]\Cu\\
\step[2]\Cu\\
\step[2]\Cd\\
\step[2]\object{D}\step[4]\object{D}\\
\end{tangle}
\step=\step
\begin{tangle}
\object{A\otimes \eta_H}\step[4]\object{A\otimes \eta_H}\step[4]\object{\eta_A\otimes H}\\
\Cu\step[3]\ne2\\
\step[2]\Cu\\
\step[2]\Cd\\
\step[2]\object{D}\step[4]\object{D}\\
\end{tangle}
\]

\[
\step\stackrel{by\ (1)}{=}\step
\begin{tangle}
\object{A\otimes \eta_H}\step[4]\object{A\otimes \eta_H}\step[4]\object{\eta_A\otimes H}\\
\Cu\step[3]\ne1\\
\Cd\step\Cd\\
\id\step[4]\hx\step[4]\id\\
\Cu\step\Cu\\
\step[2]\object{D}\step[5]\object{D}\\
\end{tangle}
=\step
\begin{tangle}
\step[2]\object{A\otimes \eta_H}\step[5]\object{A\otimes \eta_H}\step[4]\object{\eta_A\otimes H}\\
\Cd\step\Cd\step[2]\id\\
\id\step[4]\hx\step[4]\id\step[2]\id\\
\Cu\step\Cu\Cd\\
\step[2]\se1\step[4]\x\step[4]\id\\
\step[3]\Cu\step[2]\Cu\\
\step[5]\object{D}\step[6]\object{D}\\
\end{tangle}
\]

\[
\step=\step
\begin{tangle}
\step[2]\object{A\otimes \eta_H}\step[4]\object{A\otimes \eta_H}\step[5]\object{\eta_A\otimes H}\\
\step[2]\id\step[2]\Cd\step\Cd\\
\step[2]\id\step[2]\id\step[4]\hx\step[4]\id\\
\Cd\Cu\step\Cu\\
\id\step[4]\x\step[4]\ne1\\
\Cu\step[2]\Cu\\
\step[2]\object{D}\step[6]\object{D}\\
\end{tangle}
\step=\step
\begin{tangle}
\step[2]\object{A\otimes \eta_H}\step[4]\object{A\otimes \eta_H}\step[4]\object{\eta_A\otimes H}\\
\step[2]\id\step[3]\Cu\\
\Cd\step\Cd\\
\id\step[4]\hx\step[4]\id\\
\Cu\step\Cu\\
\step[2]\object{D}\step[5]\object{D}\\
\end{tangle} \ \ \  =\hbox {the right hand of  }(2).
\]
Thus relation (2) holds.

Finally, we show that
\[
\begin{tangle}
\object{D}\step[2]\object{D}\\
\cu\\
\cd\\
\object{D}\step[2]\object{D}\\
\end{tangle}
\step=\step
\begin{tangle}
\step\object{D}\step[3]\object{D}\\
\cd\step\cd\\
\id\step[2]\hx\step[2]\id\\
\cu\step\cu\\
\step\object{D}\step[3]\object{D}\\
\end{tangle}
\step[3]\cdot\cdot\cdot\cdot\cdot\cdot(3)
\]\\

\[ \hbox {the left hand of } (3)=
\begin{tangle}
\object{A\otimes \eta_H}\step[4]\object{\eta_A\otimes H}\step[4]\object{A\otimes \eta_H}\step[4]\object{\eta_A\otimes H}\\
\Cu\step[4]\Cu\\
\step[3]\se2\step[5]\ne2\\
\step[4]\Cu\\
\step[4]\Cd\\
\step[4]\object{D}\step[4]\object{D}\\
\end{tangle}
\]

\[
\step[3]\stackrel{by\ (1)}{=}\step
\begin{tangle}
\object{A\otimes \eta_H}\step[4]\object{\eta_A\otimes H}\step[4]\object{A\otimes \eta_H}\step[5]\object{\eta_A\otimes H}\\
\Cu\step[3]\ne2\step[4]\ne2\\
\step[2]\Cu\step[4]\ne2\\
\step[2]\Cd\step\Cd\\
\step[2]\id\step[4]\hx\step[4]\id\\
\step[2]\Cu\step\Cu\\
\step[4]\object{D}\step[5]\object{D}\\
\end{tangle}
\]

\[
\step[3]\stackrel{by\ (2)}{=}\step
\begin{tangle}
\step[2]\object{A\otimes \eta_H}\step[4]\object{\eta_A\otimes H}\step[4]\object{A\otimes \eta_H}\step[4]\object{\eta_A\otimes H}\\
\step[2]\id\step[4]\Cu\step[4]\id\\
\Cd\step[2]\Cd\step[4]\id\\
\id\step[4]\x\step[4]\id\step[4]\id\\
\Cu\step[2]\Cu\step[4]\id\\
\step[3]\se2\step[5]\se2\step[2]\Cd\\
\step[5]\se2\step[4]\x\step[4]\id\\
\step[6]\Cu\step[2]\Cu\\
\step[8]\object{D}\step[6]\object{D}\\
\end{tangle}
\]

\[
\step[3]=\step
\begin{tangle}
\step[2]\object{A\otimes \eta_H}\step[4]\object{\eta_A\otimes H}\step[5]\object{A\otimes \eta_H}\step[5]\object{\eta_A\otimes H}\\
\step[2]\id\step[2]\Cd\step\Cd\step[3]\id\\
\step[2]\id\step[2]\id\step[4]\hx\step[4]\id\step[3]\id\\
\Cd\Cu\step\Cu\step[2]\ne2\\
\id\step[4]\x\step[4]\ne1\step[2]\ne2\\
\Cu\step[2]\Cu\Cd\\
\step[3]\se2\step[4]\x\step[4]\id\\
\step[4]\Cu\step[2]\Cu\\
\step[6]\object{D}\step[6]\object{D}\\
\end{tangle}
\]

\[
\step[3]=\step
\begin{tangle}
\step[2]\object{A\otimes \eta_H}\step[5]\object{\eta_A\otimes H}\step[4]\object{A\otimes \eta_H}\step[4]\object{\eta_A\otimes H}\\
\Cd\step\Cd\step[2]\id\step[4]\id\\
\id\step[4]\hx\step[4]\id\step[2]\id\step[4]\id\\
\Cu\step\Cu\Cd\step[2]\id\\
\step[2]\se1\step[4]\x\step[4]\id\step[2]\id\\
\step[3]\Cu\step[2]\Cu\Cd\\
\step[6]\se2\step[4]\x\step[4]\id\\
\step[7]\Cu\step[2]\Cu\\
\step[9]\object{D}\step[6]\object{D}\\
\end{tangle}
\]

\[
\step[3]=\step
\begin{tangle}
\object{A\otimes \eta_H}\step[4]\object{\eta_A\otimes H}\step[4]\object{A\otimes \eta_H}\step[4]\object{\eta_A\otimes H}\\
\Cu\step[4]\id\step[4]\id\\
\Cd\step[2]\Cd\step[2]\id\\
\id\step[4]\x\step[4]\id\step[2]\id\\
\Cu\step[2]\Cu\Cd\\
\step[3]\se2\step[4]\x\step[4]\id\\
\step[4]\Cu\step[2]\Cu\\
\step[6]\object{D}\step[6]\object{D}\\
\end{tangle}
\]

\[
\step[3]=\step
\begin{tangle}
\step[2]\object{A\otimes H}\step[4]\object{A\otimes \eta_H}\step[5]\object{\eta_A\otimes H}\\
\step[2]\id\step[2]\Cd\step\Cd\\
\step[2]\id\step[2]\id\step[4]\hx\step[4]\id\\
\Cd\Cu\step\Cu\\
\id\step[4]\x\step[4]\ne1\\
\Cu\step[2]\Cu\\
\step[2]\object{D}\step[6]\object{D}\\
\end{tangle}
=\hbox { the right hand of  }(3).\] This relation (3) holds.

 (ii) We can show
that  relation (3) holds similarly
\[
\begin{tangle}
\object{A\otimes \eta_H}\step[4]\object{A\otimes H}\\
\Cu\\
\step[2]\object{D}\\
\end{tangle}
\step=\step
\begin{tangle}
\object{A}\step[2]\object{A}\step\object{H}\\
\cu\step\id\\
\step\object{A}\step[2]\object{H}\\
\end{tangle}
\step.\step  \ \Box
\]

\begin {Theorem}\label {3.2.1} Let  $A$  and $H$ be
bialgebras with $(\epsilon _H \otimes id _A)\phi = id _A$ and $ (id
_H \otimes \epsilon _A)\phi = \eta_H \epsilon _A$ and $\beta (\eta
_H \otimes id _A) = \epsilon _A\eta_{H}$ and $\beta (id _H \otimes
\eta _A) = id_A$ and $\epsilon_H \otimes \epsilon _A = \epsilon
_A\beta$ and $\phi\eta _A) = \eta_H\otimes\eta_{A}$. Then
 $ D = A ^{\phi , P }  \# _{\beta , \mu } H$   is a bialgebra iff
$D$ is an algebra and a coalgebra,
  and relations $(BB1)$--$(BB5)$ hold.

 \end {Theorem}
{\bf Proof.} (See the proof of \cite [Theorem 2.9] {Ma94a}) Let  $ D
= A ^{\phi , P }  \# _{\beta , \mu } H$   be a bialgebra. We have
that
\[
\begin{tangle}
\object{A\otimes H}\step[4]\object{A\otimes H}\\
\Cu\\
\Cd\\
\object{D}\step[4]\object{D}\\
\end{tangle}
\step=\step
\begin{tangle}
\step[2]\object{A\otimes H}\step[5]\object{A\otimes H}\\
\Cd\step\Cd\\
\id\step[4]\hx\step[4]\id\\
\Cu\step\Cu\\
\step[2]\object{D}\step[5]\object{D}\\
\end{tangle}
\step[3]\cdot\cdot\cdot\cdot\cdot\cdot\step(1)
\]
By (1), we have that
\[
\begin{tangle}
\footnotesize\object{A\otimes \eta_H}\step[4]\object{A\otimes \eta_H}\\
\Cu\\
\Cd\\
\object{\epsilon \otimes \epsilon }\step[4]\object{\epsilon \otimes \epsilon }\\
\end{tangle}
\step=\step
\begin{tangle}
\step[2]\object{A\otimes \eta_H}\step[5]\object{A\otimes \eta_H}\\
\Cd\step\Cd\\
\id\step[4]\hx\step[4]\id\\
\Cu\step\Cu\\
\step[2]\object{\epsilon \otimes \epsilon }\step[5]\object{\epsilon \otimes \epsilon }\\
\end{tangle}
\step,\step
\begin{tangle}
\object{A}\step[2]\object{A}\\
\tu {\mu}\\
\step\QQ {\epsilon }\\
\end{tangle}
\step=\step
\begin{tangle}
\object{A}\step[2]\object{A}\\
\id\step[2]\id\\
\QQ {\epsilon }\step[2]\QQ {\epsilon }\\
\end{tangle}
\step.\step
\]

\[
\begin{tangle}
\object{\eta_{A}\otimes \eta_{H}}\step[4]\object{\eta_{A}\otimes \eta_{H}}\\
\Cu\\
\Cd\\
\object{\epsilon \otimes id_{H}}\step[4]\object{\epsilon \otimes id_{H}}\\
\end{tangle}
\step=\step
\begin{tangle}
\step[2]\object{\eta_{A}\otimes \eta_{H}}\step[5]\object{\eta_{A}\otimes \eta_{H}}\\
\Cd\step\Cd\\
\id\step[4]\hx\step[4]\id\\
\Cu\step\Cu\\
\step[2]\object{\epsilon \otimes id_{H}}\step[5]\object{\epsilon \otimes id_{H}}\\
\end{tangle}
\step \ and \ \step
\begin{tangle}
\step\object{\eta_{A}}\\
\td {P}\\
\object{H}\step[2]\object{H}\\
\end{tangle}
\step=\step
\begin{tangle}
\object{\eta_{H}}\step[2]\object{\eta_{H}}\\
\id\step[2]\id\\
\object{H}\step[2]\object{H}\\
\end{tangle}
\step.\step
\]

Similarly, we can obtain (BB1)-(BB5) by relation (1).

Conversely, if $D$ is an algebra and coalgebra ,
  and relations $(BB1)$--$(BB5)$ hold, we show that
 the conditions in Lemma \ref {3.2.1.1} hold.
We first show that
\[
\begin{tangle}
\object{A\otimes \eta_{H}}\step[4]\object{A\otimes \eta_{H}}\\
\Cu\\
\Cd\\
\object{D}\step[4]\object{D}\\
\end{tangle}
\step=\step
\begin{tangle}
\step[2]\object{A\otimes \eta_{H}}\step[5]\object{A\otimes \eta_{H}}\\
\Cd\step\Cd\\
\id\step[4]\hx\step[4]\id\\
\Cu\step\Cu\\
\step[2]\object{D}\step[5]\object{D}\\
\end{tangle}
\step.\step
\]

\[\hbox { the left hand }
\step\stackrel{by\ (BB4),(WA),(WCA)}{=}\step
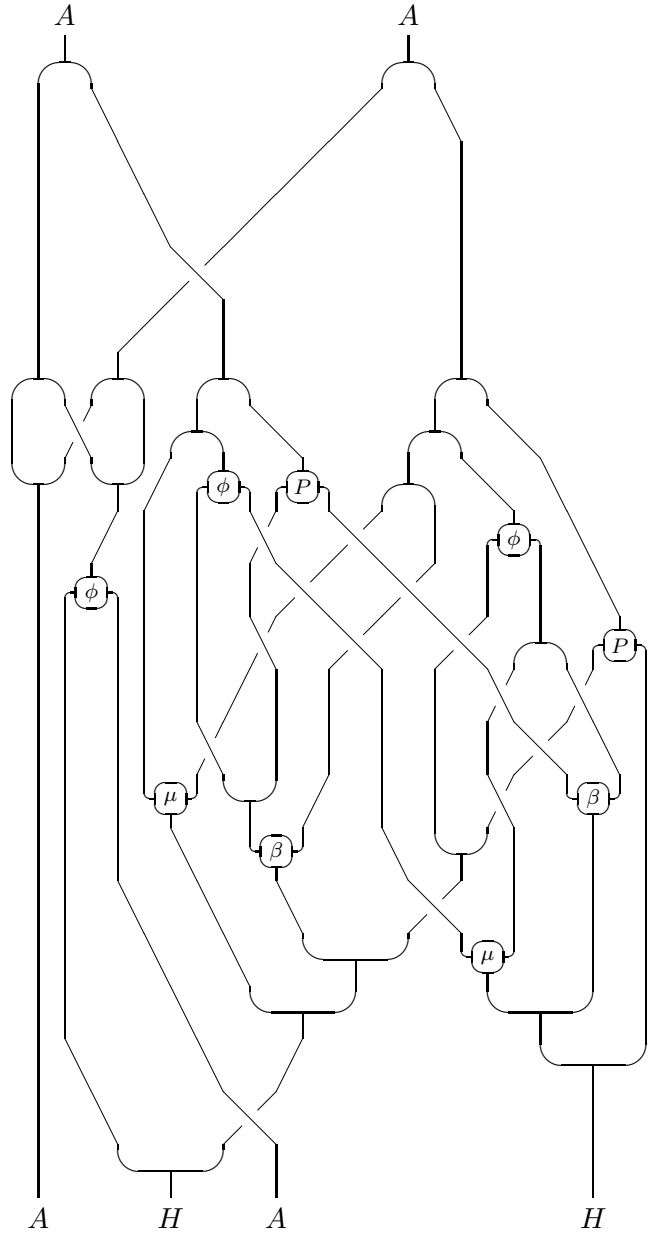
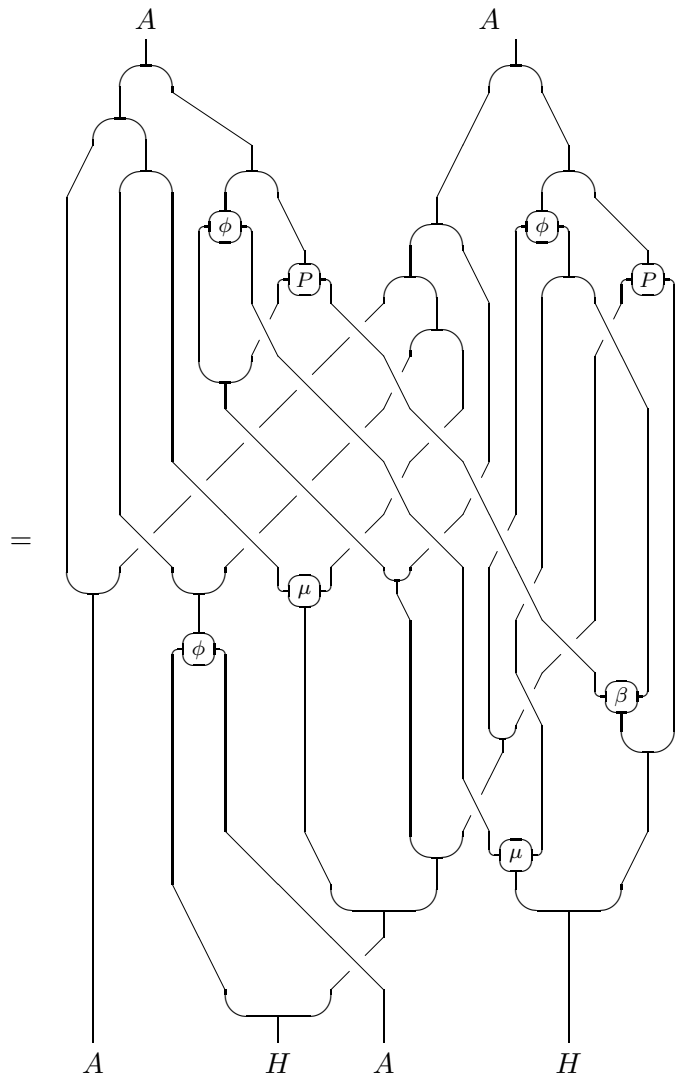
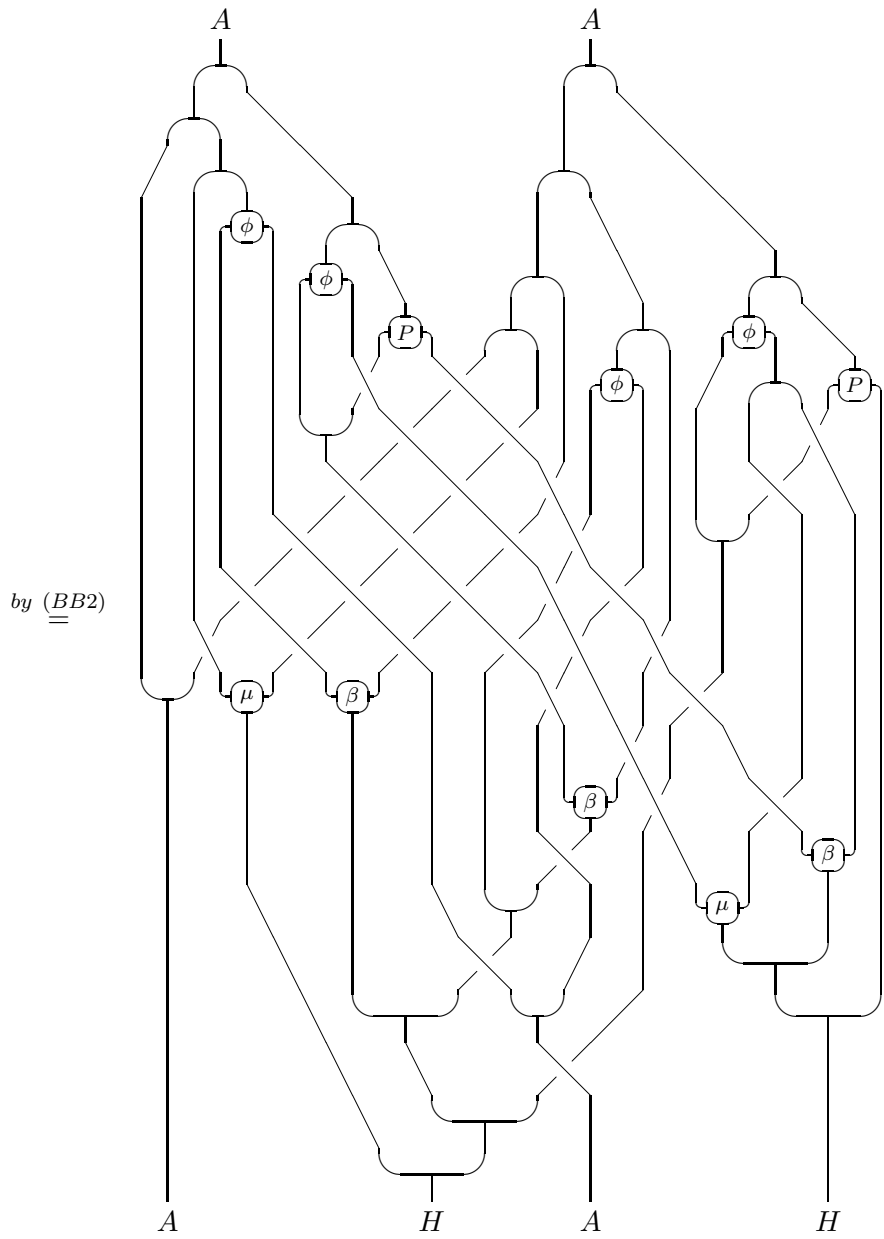
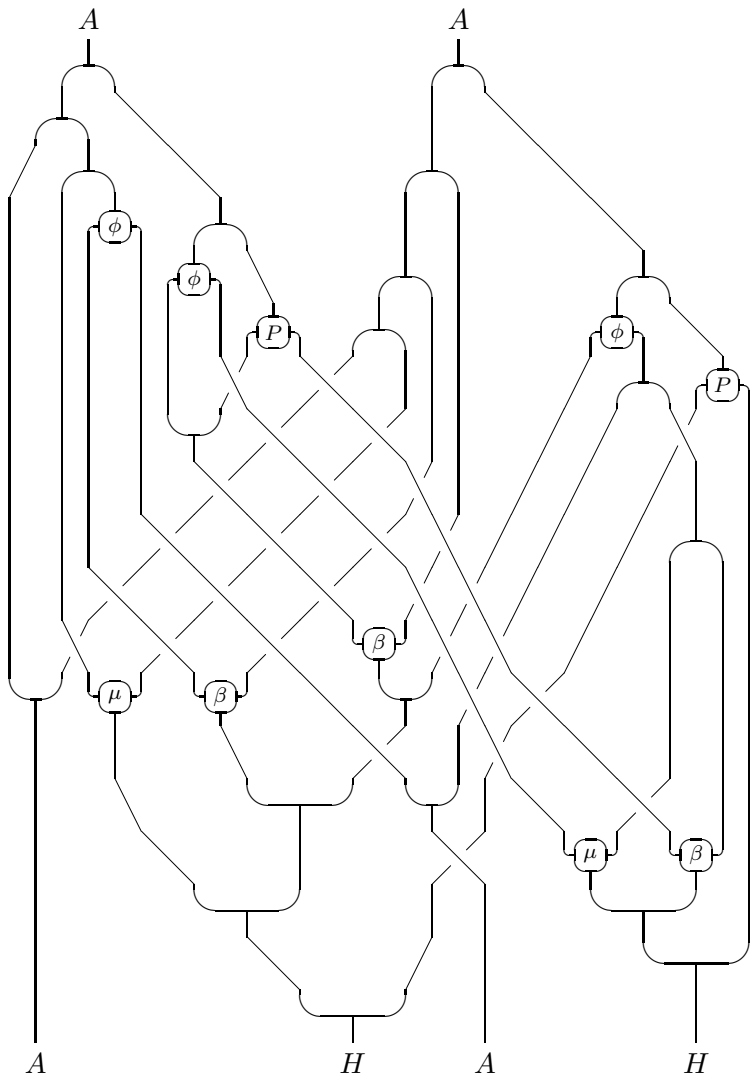
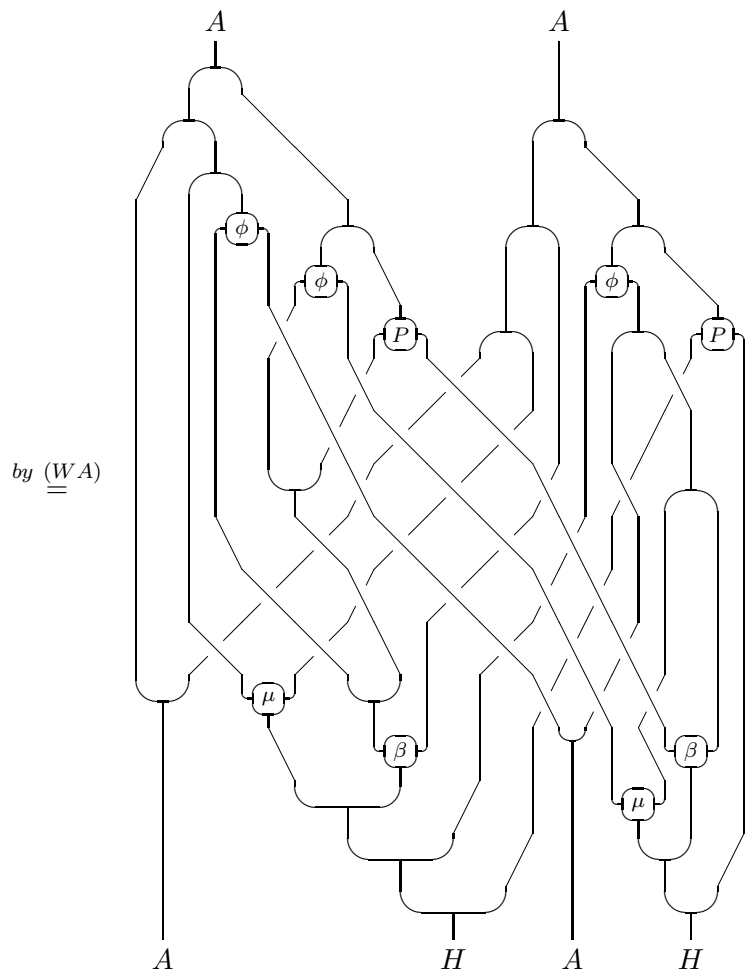
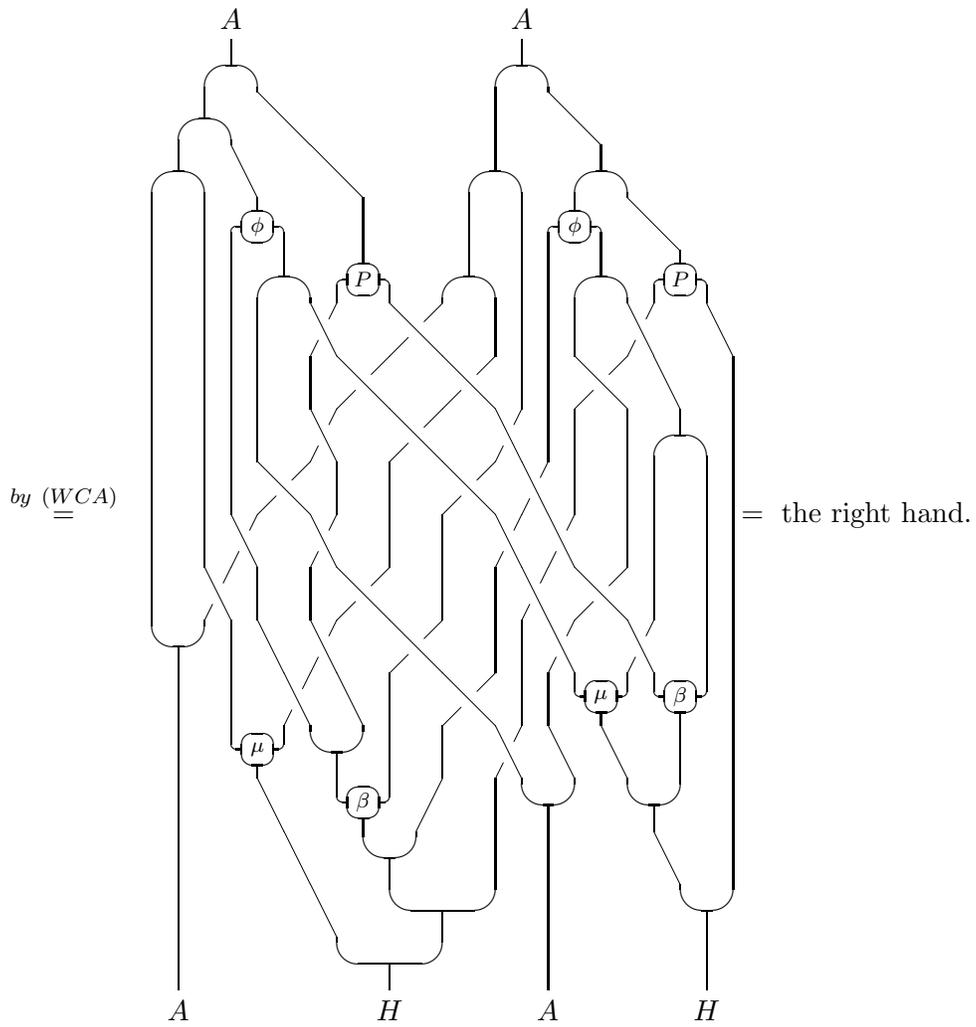
\begin{tangle}
\step[2]\object{A}\step[13]\object{A}\\
\step\cd\step[11]\cd\\
\step\id\step[2]\se1\step[9]\ne2\step[2]\se1\\
\step\id\step[3]\se1\step[6]\ne2\step[5]\id\\
\step\id\step[4]\se1\step[3]\ne2\step[7]\id\\
\step\id\step[5]\x\step[9]\id\\
\step\id\step[4]\ne2\step[2]\id\step[9]\id\\
\cd\step\cd\step[2]\cd\step[7]\cd\\
\id\step[2]\hx\step[2]\id\step\cd\step[2]\se2\step[4]\cd\step[2]\se2\\
\cu\step\cu\ne1\step\td {\phi}\step\td
{P}\step[2]\cd\step[2]\se2\step\se1\\
\step\id\step[2]\ne1\step\id\step[2]\id\step[2]\hx\step[2]\x\step[2]\id\step[2]\td
{\phi}\step\se1\\
\step\id\step\td
{\phi}\step\id\step[2]\id\step[2]\id\step\x\step[2]\x\step[2]\id\step[2]\id\step[2]\se1\\
\step\id\step\id\step[2]\id\step\id\step[2]\id\step[2]\hx\step[2]\x\step[2]\x\step\cd\step\td
{P}\\
\step\id\step\id\step[2]\id\step\id\step[2]\id\step\ne1\step\id\step[2]\id\step[2]\id\step[2]\id\step[2]\hx\step[2]\hx\step[2]\id\\
\step\id\step\id\step[2]\id\step\id\step[2]\hx\step[2]\id\step[2]\id\step[2]\id\step[2]\id\step[2]\id\step\x\step\se1\step\id\\
\step\id\step\id\step[2]\id\step\tu
{\mu}\step\cu\step\ne1\step[2]\id\step[2]\id\step[2]\hx\step[2]\tu
{\beta}\step\id\\
\step\id\step\id\step[2]\id\step[2]\se1\step[2]\tu
{\beta}\step[3]\se1\step\cu\step\id\step[3]\id\step[2]\id\\
\step\id\step\id\step[2]\se1\step[2]\se1\step[2]\se1\step[4]\x\step[2]\id\step[3]\id\step[2]\id\\
\step\id\step\id\step[3]\se1\step[2]\se1\step[2]\Cu\step[2]\tu
{\mu}\step[3]\id\step[2]\id\\
\step\id\step\id\step[4]\se1\step[2]\Cu\step[5]\Cu\step[2]\id\\
\step\id\step\se1\step[4]\se1\step[2]\ne1\step[9]\Cu\\
\step\id\step[2]\se1\step[4]\x\step[12]\id\\
\step\id\step[3]\Cu\step[2]\id\step[12]\id\\
\step\object{A}\step[5]\object{H}\step[4]\object{A}\step[12]\object{H}\\
\end{tangle}
\]

\[
\step=\step
\begin{tangle}
\step[3]\object{A}\step[13]\object{A}\\
\step[2]\cd\step[12]\cd\\
\step\cd\step[3]\se3\step[8]\ne1\step[2]\se1\\
\ne1\step\cd\step[2]\cd\step[6]\ne1\step[3]\cd\\
\id\step[2]\id\step[2]\id\step\td
{\phi}\step\se1\step[4]\cd\step[2]\td {\phi}\step[2]\se2\\
\id\step[2]\id\step[2]\id\step\id\step[2]\id\step\td
{P}\step[2]\cd\step\se1\step\id\step\cd\step\td {P}\\
\id\step[2]\id\step[2]\id\step\id\step[2]\hx\step[2]\x\step\cd\step\id\step\id\step\id\step[2]\hx\step[2]\id\\
\id\step[2]\id\step[2]\id\step\cu\step\x\step[2]\hx\step[2]\id\step\id\step\id\step\id\step[2]\id\step\se1\step\id\\
\id\step[2]\id\step[2]\id\step[2]\x\step[2]\x\step\x\step\id\step\id\step\id\step[2]\id\step[2]\id\step\id\\
\id\step[2]\id\step[2]\x\step[2]\x\step[2]\hx\step[2]\hx\step\id\step\id\step[2]\id\step[2]\id\step\id\\
\id\step[2]\x\step[2]\x\step[2]\x\step\x\step\hx\step\id\step[2]\id\step[2]\id\step\id\\
\cu\step[2]\cu\step[2]\tu {\mu}\step[2]\hdcu\step[2]\id\step\id\step\hx\step[2]\id\step[2]\id\step\id\\
\step\id\step[3]\td
{\phi}\step[3]\id\step[4]\id\step[2]\id\step\id\step\id\step\x\step[2]\id\step\id\\
\step\id\step[3]\id\step[2]\id\step[3]\id\step[4]\id\step[2]\id\step\id\step\hx\step[2]\tu
{\beta}\step\id\\
\step\id\step[3]\id\step[2]\id\step[3]\id\step[4]\id\step[2]\id\step\hddcu\step\id\step[3]\cu\\
\step\id\step[3]\id\step[2]\id\step[3]\id\step[4]\id\step[2]\hx\step[2]\id\step[4]\id\\
\step\id\step[3]\id\step[3]\se2\step\se1\step[3]\cu\step\tu
{\mu}\step[3]\ne1\\
\step\id\step[3]\se1\step[4]\se2\Cu\step[3]\Cu\\
\step\id\step[4]\se1\step[4]\x\step[7]\id\\
\step\id\step[5]\Cu\step[2]\id\step[7]\id\\
\step\object{A}\step[7]\object{H}\step[4]\object{A}\step[7]\object{H}\\
\end{tangle}
\]

\[
\step\stackrel{by\ (BB2)}{=}\step
\begin{tangle}
\step[3]\object{A}\step[14]\object{A}\\
\step[2]\cd\step[12]\cd\\
\step\cd\step[2]\se2\step[10]\id\step[3]\se2\\
\ne1\step\cd\step[3]\se2\step[7]\cd\step[4]\se2\\
\id\step[2]\id\step\td
{\phi}\step[2]\cd\step[6]\id\step[2]\se1\step[5]\se2\\
\id\step[2]\id\step\id\step[2]\id\step\td
{\phi}\step\se1\step[4]\cd\step[2]\se1\step[4]\cd\\
\id\step[2]\id\step\id\step[2]\id\step\id\step[2]\id\step\td
{P}\step[2]\cd\step\id\step[2]\cd\step[2]\td {\phi}\step[2]\se2\\
\id\step[2]\id\step\id\step[2]\id\step\id\step[2]\hx\step[2]\x\step[2]\id\step\id\step\td
{\phi}\step\id\step\ne1\step\cd\step\td {P}\\
\id\step[2]\id\step\id\step[2]\id\step\cu\step\x\step[2]\x\step\id\step\id\step[2]\id\step\id\step\id\step[2]\id\step[2]\hx\step[2]\id\\
\id\step[2]\id\step\id\step[2]\id\step[2]\x\step[2]\x\step[2]\hx\step\id\step[2]\id\step\id\step\id\step[2]\x\step\se1\step\id\\
\id\step[2]\id\step\id\step[2]\x\step[2]\x\step[2]\x\step\hx\step[2]\id\step\id\step\cu\step[2]\id\step[2]\id\step\id\\
\id\step[2]\id\step\x\step[2]\x\step[2]\x\step[2]\hx\step\x\step\id\step[2]\id\step[3]\id\step[2]\id\step\id\\
\id\step[2]\hx\step[2]\x\step[2]\x\step[2]\x\step\hx\step[2]\hx\step[2]\id\step[3]\id\step[2]\id\step\id\\

\cu\step\tu \mu\step[2]\tu
{\beta}\step[2]\id\step[2]\id\step[2]\hx\step\se1\step\id\step\x\step[3]\id\step[2]\id\step\id\\
\step\id\step[3]\id\step[4]\id\step[3]\id\step[2]\id\step[2]\id\step\id\step[2]\hx\step\id\step[2]\se1\step[2]\id\step[2]\id\step\id\\
\step\id\step[3]\id\step[4]\id\step[3]\id\step[2]\id\step[2]\id\step\tu
{\beta}\step\hx\step[3]\x\step[2]\id\step\id\\
\step\id\step[3]\id\step[4]\id\step[3]\id\step[2]\id\step[2]\x\step[2]\id\step\se1\step[2]\id\step[2]\tu
{\beta}\step\id\\
\step\id\step[3]\se1\step[3]\id\step[3]\se1\step\cu\step[2]\id\step[2]\id\step[2]\tu
{\mu}\step[3]\id\step[2]\id\\
\step\id\step[4]\se1\step[2]\id\step[4]\x\step[2]\ne1\step[2]\id\step[3]\Cu\step[2]\id\\
\step\id\step[5]\se1\step\Cu\step[2]\cu\step[2]\ne2\step[5]\Cu\\
\step\id\step[6]\se1\step[2]\se1\step[4]\x\step[9]\id\\
\step\id\step[7]\se1\step[2]\Cu\step[2]\id\step[9]\id\\
\step\id\step[8]\Cu\step[4]\id\step[9]\id\\
\step\object{A}\step[10]\object{H}\step[6]\object{A}\step[9]\object{H}\\
\end{tangle}
\]

\[
\step\stackrel{by\ (BB3)}{=}\step
\begin{tangle}
\step[3]\object{A}\step[14]\object{A}\\
\step[2]\cd\step[12]\cd\\
\step\cd\step[2]\se2\step[10]\id\step[3]\se2\\
\ne1\step\cd\step[3]\se2\step[7]\cd\step[4]\se2\\
\id\step[2]\id\step\td {\phi}\step[2]\cd\step[6]\id\step[2]\id\step[6]\se2\\
\id\step[2]\id\step\id\step[2]\id\step\td {\phi}\step\se1\step[4]\cd\step\id\step[6]\cd\\
\id\step[2]\id\step\id\step[2]\id\step\id\step[2]\id\step\td {P}\step[2]\cd\step\id\step\id\step[5]\td {\phi}\step[2]\se2\\
\id\step[2]\id\step\id\step[2]\id\step\id\step[2]\hx\step[2]\x\step[2]\id\step\id\step\id\step[4]\ne1\step\cd\step\td {P}\\
\id\step[2]\id\step\id\step[2]\id\step\cu\step\x\step[2]\x\step\id\step\id\step[3]\ne1\step\ne1\step[2]\hx\step[2]\id\\
\id\step[2]\id\step\id\step[2]\id\step[2]\x\step[2]\x\step[2]\hx\step\id\step[2]\ne1\step\ne1\step[2]\ne1\step\id\step[2]\id\\
\id\step[2]\id\step\id\step[2]\x\step[2]\x\step[2]\x\step\hx\step\ne1\step\ne1\step[2]\ne1\step\cd\step\id\\
\id\step[2]\id\step\x\step[2]\x\step[2]\x\step[2]\hx\step\hx\step\ne1\step[2]\ne1\step[2]\id\step[2]\id\step\id\\
\id\step[2]\hx\step[2]\x\step[2]\x\step[2]\tu {\beta}\step\hx\step\hx\step[2]\ne1\step[3]\id\step[2]\id\step\id\\
\cu\step\tu {\mu}\step[2]\tu {\beta}\step[3]\se2\step\cu\step\hx\step\x\step[4]\id\step[2]\id\step\id\\
\step\id\step[3]\id\step[4]\se1\step[4]\x\step[2]\id\step\hx\step[3]\se2\step[2]\id\step[2]\id\step\id\\
\step\id\step[3]\se1\step[4]\Cu\step[2]\cu\step\id\step[2]\se2\step[2]\x\step[2]\id\step\id\\
\step\id\step[5]\se2\step[4]\id\step[5]\x\step[3]\tu {\mu}\step[2]\tu {\beta}\step\id\\
\step\id\step[6]\Cu\step[5]\id\step[2]\id\step[4]\Cu\step[2]\id\\
\step\id\step[9]\se2\step[4]\ne1\step[2]\id\step[6]\Cu\\
\step\id\step[10]\Cu\step[3]\id\step[8]\id\\
\step\object{A}\step[12]\object{H}\step[5]\object{A}\step[8]\object{H}\\
\end{tangle}
\]

\[
\step\stackrel{by\ (WA)}{=}\step
\begin{tangle}
\step[3]\object{A}\step[13]\object{A}\\
\step[2]\cd\step[12]\id\\
\step\cd\step[2]\se2\step[9]\cd\\
\ne1\step\cd\step[3]\se2\step[7]\id\step[3]\se2\\
\id\step[2]\id\step\td {\phi}\step[2]\cd\step[5]\cd\step[2]\cd\\
\id\step[2]\id\step\id\step[2]\id\step\td {\phi}\step\se1\step[4]\id\step[2]\id\step\td {\phi}\step[2]\se2\\
\id\step[2]\id\step\id\step[2]\hx\step[2]\id\step\td {P}\step[2]\cd\step\id\step\id\step\cd\step\td {P}\\
\id\step[2]\id\step\id\step[2]\id\step\se1\step\hx\step[2]\x\step[2]\id\step\id\step\id\step\id\step[2]\hx\step[2]\id\\
\id\step[2]\id\step\id\step[2]\id\step[2]\hx\step\x\step[2]\x\step\id\step\id\step\id\step\ne1\step\id\step[2]\id\\
\id\step[2]\id\step\id\step[2]\cu\step\hx\step[2]\x\step[2]\hx\step\id\step\hx\step\cd\step\id\\
\id\step[2]\id\step\se1\step[2]\x\step\x\step[2]\x\step\hx\step\id\step\id\step\id\step[2]\id\step\id\\
\id\step[2]\id\step[2]\x\step[2]\hx\step[2]\x\step[2]\hx\step\hx\step\id\step\id\step[2]\id\step\id\\
\id\step[2]\x\step[2]\x\step\se1\step\id\step[2]\x\step\hx\step\hx\step\id\step[2]\id\step\id\\
\cu\step[2]\tu {\mu}\step[2]\cu\step\id\step[2]\id\step[2]\hx\step\hx\step\hx\step[2]\id\step\id\\
\step\id\step[4]\se1\step[3]\tu {\beta}\step[2]\id\step[2]\id\step\hcu\step\id\step\se1\tu{\beta}\step\id\\
\step\id\step[5]\Cu\step[2]\ne1\step[2]\id\step[1.5]\id\step[1.5]\tu {\mu}\step\id\step[2]\id\\
\step\id\step[7]\Cu\step[2]\ne1\step[1.5]\id\step[2.5]\cu\step\ne1\\
\step\id\step[9]\Cu\step[2.5]\id\step[3.5]\cu\\
\step\object{A}\step[11]\object{H}\step[4.5]\object{A}\step[4.5]\object{H}\\
\end{tangle}
\]

\[
\step\stackrel{by\ (WCA)}{=}\step
\begin{tangle}
\step[3]\object{A}\step[11]\object{A}\\
\step[2]\cd\step[9]\cd\\
\step\cd\step[2]\se2\step[7]\id\step[3]\se2\\
\cd\step\se1\step[3]\se2\step[4]\cd\step[2]\cd\\
\id\step[2]\id\step\td {\phi}\step[3]\id\step[4]\id\step[2]\id\step\td {\phi}\step[2]\se2\\
\id\step[2]\id\step\id\step\cd\step\td {P}\step[2]\cd\step\id\step\id\step\cd\step\td {P}\\
\id\step[2]\id\step\id\step\id\step[2]\hx\step[2]\x\step[2]\id\step\id\step\id\step\id\step[2]\hx\step[2]\se1\\
\id\step[2]\id\step\id\step\id\step[2]\id\step\x\step[2]\x\step\id\step\id\step\x\step\se1\step[2]\id\\
\id\step[2]\id\step\id\step\id\step[2]\hx\step[2]\x\step[2]\hx\step\id\step\id\step[2]\id\step\cd\step\id\\
\id\step[2]\id\step\id\step\x\step\id\step[2]\id\step[2]\x\step\hx\step\id\step[2]\id\step\id\step[2]\id\step\id\\
\id\step[2]\id\step\hx\step[2]\hx\step[2]\id\step[2]\id\step[2]\hx\step\hx\step[2]\id\step\id\step[2]\id\step\id\\
\id\step[2]\hx\step\id\step[2]\id\step\x\step[2]\id\step[2]\id\step\hx\step\x\step\id\step[2]\id\step\id\\
\cu\step\id\step\se1\step\hx\step[2]\x\step[2]\id\step\id\step\hx\step[2]\hx\step[2]\id\step\id\\
\step\id\step[2]\id\step[2]\hx\step\se1\step\id\step[2]\x\step\id\step\id\step\tu {\mu}\step\tu {\beta}\step\id\\
\step\id\step[2]\tu {\mu}\step\cu\step\id\step[2]\id\step[2]\hx\step\se1\step\se1\step[2]\id\step[2]\id\\
\step\id\step[3]\se1\step[2]\tu {\beta}\step\ne1\step[2]\id\step\cu\step[2]\cu\step[2]\id\\
\step\id\step[4]\se1\step[2]\cu\step[3]\id\step[2]\id\step[4]\se1\step[2]\id\\
\step\id\step[5]\se1\step[2]\Cu\step[2]\id\step[5]\cu\\
\step\id\step[6]\Cu\step[4]\id\step[6]\id\\
\step\object{A}\step[8]\object{H}\step[6]\object{A}\step[6]\object{H}\\
\end{tangle}
=\hbox { the right hand}. \]

Similarly, we have that
\[
\begin{tangle}
\object{\eta_A\otimes H}\step[4]\object{\eta_A\otimes H}\\
\Cu\\
\Cd\\
\object{D}\step[4]\object{D}\\
\end{tangle}
\step=\step
\begin{tangle}
\step\object{\eta_A\otimes H}\step[4]\object{\eta_A\otimes H}\\
\cd\step[2]\cd\\
\id\step[2]\x\step[2]\id\\
\cu\step[2]\cu\\
\step\object{D}\step[4]\object{D}\\
\end{tangle}
\step,\step
\]

\[
\begin{tangle}
\object{\eta_A\otimes H}\step[4]\object{A\otimes \eta_H}\\
\Cu\\
\Cd\\
\object{D}\step[4]\object{D}\\
\end{tangle}
\step=\step
\begin{tangle}
\step\object{\eta_A\otimes H}\step[4]\object{A\otimes \eta_H}\\
\cd\step[2]\cd\\
\id\step[2]\x\step[2]\id\\
\cu\step[2]\cu\\
\step\object{D}\step[4]\object{D}\\
\end{tangle}
\step,\step
\]

\[
\begin{tangle}
\object{A\otimes \eta_H}\step[4]\object{\eta_A\otimes H}\\
\Cu\\
\Cd\\
\object{D}\step[4]\object{D}\\
\end{tangle}
\step=\step
\begin{tangle}
\step\object{A\otimes \eta_H}\step[4]\object{\eta_A\otimes H}\\
\cd\step[2]\cd\\
\id\step[2]\x\step[2]\id\\
\cu\step[2]\cu\\
\step\object{D}\step[4]\object{D}\\
\end{tangle}
\step.\step
\]
Obviously,
\[
\begin{tangle}
\object{A\otimes \eta_H}\step[4]\object{A\otimes H}\\
\Cu\\
\step[2]\object{D}\\
\end{tangle}
\step=\step
\begin{tangle}
\object{A}\step[2]\object{A}\step\object{H}\\
\cu\step\id\\
\step\object{A}\step[2]\object{H}\\
\end{tangle}
\step.\step
\]
By Lemma \ref {3.2.1.1}, $D$  is a bialgebra.
\begin{picture}(8,8)\put(0,0){\line(0,1){8}}\put(8,8){\line(0,-1){8}}\put(0,0){\line(1,0){8}}\put(8,8){\line(-1,0){8}}\end{picture}

Dually we have Theorem \ref {3.2.1}'\ Let  $A$  and $H$ be
bialgebras with $\alpha(id_H \otimes \eta _A) = \epsilon
_H\epsilon_{A}$ and $\alpha (\eta _H \otimes id _A) = id _A$ and $
(\epsilon _H \otimes id _A)\varphi = \epsilon _H\eta_{A}$ and $ (id
_H \otimes \epsilon _A)\varphi = id_H$ and $\epsilon_H \otimes
\epsilon _A) = \epsilon _A\alpha$ and $\varphi\eta _H=
\eta_H\otimes\eta_{A}$. Then
 $ D =
 A _{\alpha , \sigma }  \# ^{\psi , Q} H$   is a bialgebra iff
$D$ is an algebra and coalgebra ,
  and relations $(BB1')$--$(BB5')$ hold.

\begin {Corollary}\label {3.2.2}   Let $A$ and $H$  be bialgebras. Then

(i)   $A _{\alpha  }  \# ^{\psi } H$   is a bialgebra iff $(A,
\alpha )$ is an $H$-module algebra, $(H, \psi )$ is an $A$-comodule
coalgebra,  and (BB3'), (BB10'), (BB11') hold;

(ii) $A ^{\phi  }  \# _{\beta } H$   is a bialgebra iff $(A, \phi )$
is an $H$-comodule coalgebra, $(H, \beta )$ is an $A$-module
algebra,  and (BB3), (BB10), (BB11) hold.

 \end {Corollary}

\begin {Corollary}\label {3.2.3}   Let $A$ and $H$  be bialgebras. Then

(i) $A _{\alpha , \sigma }  \#   H$   is a bialgebra iff $A _{\alpha
, \sigma }  \#   H$ is an algebra, $\sigma $ and $\alpha $  are
coalgebra morphisms, $H$ is cocommutative with respect to $(A,
\alpha )$, and (BB8')  holds;

(ii)   $A  \#_{\beta , \mu }   H$   is a bialgebra iff
 $A  \#_{\beta , \mu }   H$   is an algebra,
$\mu $ and $\beta $  are coalgebra morphisms, $A$ is cocommutative
with respect to $(H, \beta )$, and (BB9')  holds;

(iii)       $A ^{\phi , P }  \#   H$   is a bialgebra iff $A ^{\phi
, P }  \#   H$   is an algebra, $P $ and $\phi $  are coalgebra
morphisms, $H$ is commutative with respect to $(A, \phi )$, and
(BB8)  holds;

(iv) $A   \# ^{\psi , Q} H$   is a bialgebra iff
 $A   \# ^{\psi , Q} H$   is an algebra,
$Q $ and $\psi $  are coalgebra morphisms, $H$ is commutative with
respect to $(A, \psi  )$, and (BB9)  holds;

 \end {Corollary}

\section {Biproducts }\label {s7}
In this section, we give the necessary and sufficient conditions for
biproducts to become bialgebras.

\begin {Lemma}\label {3.3.1.1}Let  $H$  be a bialgebra,  $(A,\alpha )$  an  $H$-module algebra
and $(A,\phi )$ is an $H$-comodule coalgebra. Set $D = A ^\phi
_\sigma \# H.$ Then (i) Relation (1) and relation (2) are
equivalent.
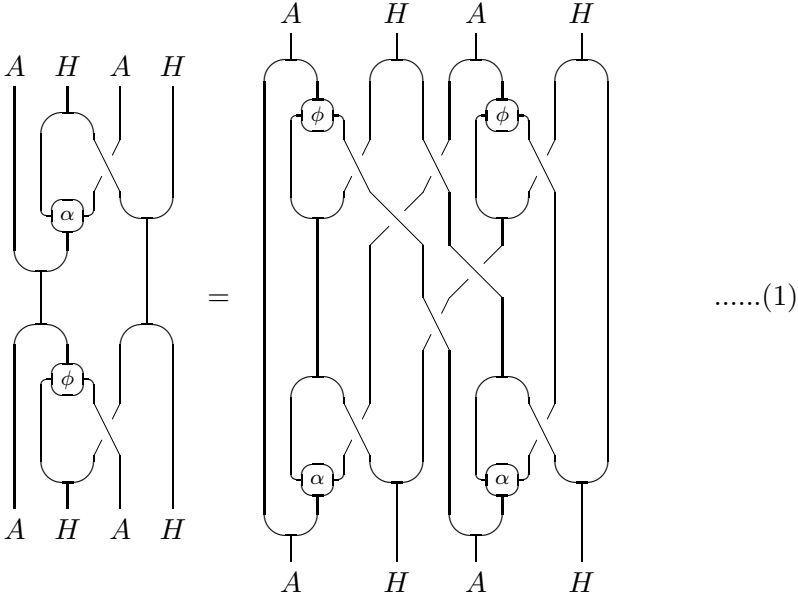
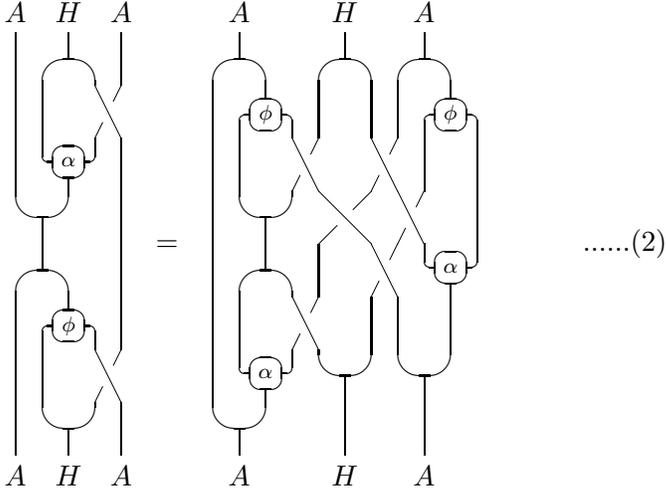
\[
\begin{tangle}
\object{A}\step\step\object{H}\step\step\object{A}\Step\object{H}\\
\id\step\cd\step\id\Step\id\\
\id\step\id\Step\hx\Step\id\\
\id\step\tu \alpha \step\cu\\
\cu\Step\step\id\\
\cd\Step\cd\\
\id\step\td \phi \step\id\Step\id\\
\id\step\id\Step\hx\Step\id\\
\id\step\cu\step\id\Step\id\\
\object{A}\Step\object{H}\Step\object{A}\Step\object{H}\\
\end{tangle}
\step=\step
\begin{tangle}
\step\object{A}\Step\Step\object{H}\Step\step\object{A}\Step\Step\object{H}\\
\cd\Step\cd\step\cd\Step\cd\\
\id\step\td \phi\step\id\Step\id\step\id\step\td
\phi\step\id\Step\id\\
\id\step\id\Step\hx\Step\hx\step\id\Step\hx\Step\id\\
\id\step\cu\step\x\step\id\step\cu\step\id\Step\id\\
\id\Step\id\Step\id\Step\id\step\x\Step\id\Step\id\\
\id\Step\id\Step\id\Step\hx\Step\id\Step\id\Step\id\\
\id\step\cd\step\id\Step\id\step\id\step\cd\step\id\Step\id\\
\id\step\id\Step\hx\Step\id\step\id\step\id\Step\hx\Step\id\\
\id\step\tu \alpha\step\cu\step\id\step\tu \alpha \step\cu\\
\cu\Step\step\id\Step\cu\Step\step\id\\
\step\object{A}\Step\Step\object{H}\step\Step\object{A}\Step\Step\object{H}\\
\end{tangle}
\Step\Step ......(1)
\]

\[
\begin{tangle}
\object{A}\Step\object{H}\Step\object{A}\\
\id\step\cd\step\id\\
\id\step\id\Step\hx\\
\id\step\tu \alpha\step\id\\
\cu\Step\id\\
\cd\Step\id\\
\id\step\td \phi\step\id\\
\id\step\id\Step\hx\\
\id\step\cu\step\id\\
\object{A}\Step\object{H}\Step\object{A}\\
\end{tangle}
\step=\step
\begin{tangle}
\step\object{A}\Step\Step\object{H}\step\Step\object{A}\\
\cd\Step\cd\step\cd\\
\id\step\td \phi\step\id\Step\id\step\id\step\td \phi\\
\id\step\id\Step\hx\Step\hx\step\id\Step\id\\
\id\step\cu\step\x\step\hx\Step\id\\
\id\step\cd\step\id\Step\hx\step\tu \alpha\\
\id\step\id\Step\hx\Step\id\step\id\Step\id\\
\id\step\tu \alpha\step\cu\step\cu\\
\cu\Step\step\id\Step\step\id\\
\step\object{A}\Step\Step\object{H}\step\Step\object{A}\\
\end{tangle}
\Step\Step ......(2)
\]

(ii) Relation (3) is equivalent to (BB7') and relation (4)
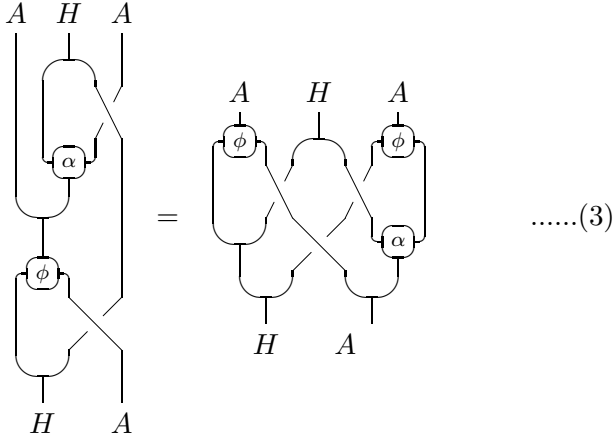
\[
\begin{tangle}
\object{A}\Step\object{H}\Step\object{A}\\
\id\step\cd\step\id\\
\id\step\id\Step\hx\\
\id\step\tu \alpha\step\id\\
\cu\Step\id\\
\td \phi\Step\id\\
\id\Step\x\\
\cu\Step\id\\
\step\object{H}\step\Step\object{A}\\
\end{tangle}
\step=\step
\begin{tangle}
\step\object{A}\step\Step\object{H}\step\Step\object{A}\\
\td \phi\step\cd\step\td \phi\\
\id\Step\hx\Step\hx\Step\id\\
\cu\step\x\step\tu \alpha\\
\step\cu\Step\cu\\
\Step\object{H}\step\Step\object{A}\\
\end{tangle}
\Step\Step ......(3)
\]

\[
\begin{tangle}
\object{A}\Step\object{A}\\
\cu\\
\td \phi\\
\object{H}\Step\object{A}\\
\end{tangle}
\step=\step
\begin{tangle}
\step\object{A}\step\Step\object{A}\\
\td \phi\step\td \phi\\
\id\Step\hx\Step\id\\
\cu\step\cu\\
\step\object{H}\step\Step\object{A}\\
\end{tangle}
\Step\Step ......(4)
\]

(iii)  Relation (2) implies relation (5) (6).

\[
\begin{tangle}
\object{A}\Step\object{A}\\
\cu\\
\cd\\
\id\step\td \phi\\
\object{A}\step\object{H}\Step\object{A}\\
\end{tangle}
\step=\step
\begin{tangle}
\Step\object{A}\Step\Step\step\object{A}\\
\Cd\Step\cd\\
\id\Step\step\td \phi\step\id\step\td \phi\\
\id\Step\dd\Step\hx\step\id\Step\id\\
\id\step\cd\step\dd\step\hx\Step\id\\
\id\step\id\Step\hx\Step\id\step\cu\\
\id\step\tu \alpha\step\cu\Step\id\\
\cu\step\Step\id\step\Step\id\\
\step\object{A}\Step\Step\object{H}\step\Step\object{A}\\
\end{tangle}
\Step\Step ......(5)
\]

\[
\begin{tangle}
\object{A}\step\object{H}\Step\object{A}\\
\id\step\tu\alpha\\
\cu\\
\cd\\
\object{A}\Step\object{A}\\
\end{tangle}
\step=\step
\begin{tangle}
\step\object{A}\Step\Step\object{H}\step\Step\object{A}\\
\cd\Step\cd\step\cd\\
\id\step\td \phi\step\id\Step\id\step\id\Step\id\\
\id\step\id\Step\hx\Step\hx\Step\id\\
\id\step\cu\step\x\step\tu \alpha\\
\d\step\tu \alpha\Step\cu\\
\step\cu\Step\Step\id\\
\Step\object{A}\step\Step\Step\object{A}\\
\end{tangle}
\Step\Step ......(6)
\]
(iv) Relation (6) is equivalent to (BB6') and relation (7)
\[
\begin{tangle}
\object{H}\Step\object{A}\\
\tu \alpha\\
\cd\\
\object{A}\Step\object{A}\\
\end{tangle}
\step=\step
\begin{tangle}
\step\object{H}\step\Step\object{A}\\
\cd\step\cd\\
\id\Step\hx\Step\id\\
\tu \alpha\step\tu \alpha\\
\step\object{A}\step\Step\object{A}\\
\end{tangle}
\Step\Step ......(7)
\]
(v)  Relation (2)  implies (3) when $\epsilon _A (\alpha ) =
\epsilon _A \otimes \epsilon _H;$ (vi) $\epsilon _D$ is an algebra
morphism iff $\epsilon _A$  is an algebra morphism and
\[
\begin{tangle}
\object{H}\Step\object{A}\\
\tu \alpha\\
\step\QQ \epsilon \\
\end{tangle}
\step=\step
\begin{tangle}
\object{H}\Step\object{A}\\
\id\Step\id\\
\QQ \epsilon \Step\QQ \epsilon \\
\end{tangle}\ \ \ .
\]

(vii)
\[
\begin{tangle}
\step\object{\eta_D}\\
\cd\\
\object{D}\Step\object{D}\\
\end{tangle}
\step=\step
\begin{tangle}
\object{\eta_A}\Step\object{\eta_H}\Step\object{\eta_A}\Step\object{\eta_H}\\
\id\Step\id\Step\id\Step\id\\
\object{A}\Step\object{H}\Step\object{A}\Step\object{H}\\
\end{tangle}
\Step\Step \hbox { iff }\Step\Step
\begin{tangle}
\step\object{\eta_A}\\
\td \phi\\
\object{H}\Step\object{A}\\
\end{tangle}
\step=\step
\begin{tangle}
\object{\eta_H}\Step\object{\eta_A}\\
\id\Step\id\\
\object{H}\Step\object{A}\\
\end{tangle}
\Step\Step and\Step\Step
\begin{tangle}
\step\object{\eta_A}\\
\cd\\
\object{A}\Step\object{A}\\
\end{tangle}
\step=\step
\begin{tangle}
\object{\eta_A}\Step\object{\eta_A}\\
\id\Step\id\\
\object{A}\Step\object{A}\\
\end{tangle}\ \ \ .
\]
\end {Lemma}
{\bf Proof.} (i) It is clear that (1) implies (2). Conversely, if
relation (2)  holds, we have that
\[
\hbox{the left side of (1)}\ =\
\begin{tangle}
\object{A}\Step\object{H}\Step\object{A}\Step\step\object{H}\\
\id\step\cd\step\id\step\Step\id\\
\id\step\id\Step\hx\step\Step\id\\
\id\step\tu \alpha\cd\step\cd\\
\cu\step\id\Step\hx\Step\id\\
\cd\step\cu\step\cu\\
\id\step\td \phi\step\id\step\Step\id\\
\id\step\id\Step\hx\Step\step\id\\
\id\step\cu\step\id\step\Step\id\\
\object{A}\Step\object{H}\Step\object{A}\Step\step\object{H}\\
\end{tangle}
\step=\step
\begin{tangle}
\object{A}\step\Step\object{H}\Step\object{A}\Step\object{H}\\
\id\Step\cd\step\id\step\cd\\
\id\step\cd\step\hx\step\id\Step\id\\
\id\step\id\Step\hx\step\hx\Step\id\\
\id\step\tu \alpha\step\id\step\id\step\cu\\
\cu\Step\id\step\id\Step\id\\
\cd\Step\id\step\id\Step\id\\
\id\step\td \phi\step\id\step\id\Step\id\\
\id\step\id\Step\hx\step\id\Step\id\\
\id\step\cu\step\hx\Step\id\\
\id\Step\cu\step\id\Step\id\\
\object{A}\step\Step\object{H}\Step\object{A}\Step\object{H}\\
\end{tangle}
\ \stackrel{by\ (2)}{=}\
\begin{tangle}
\step\object{A}\step\Step\Step\object{H}\step\Step\object{A}\Step\Step\object{H}\\
\cd\step\Step\cd\step\cd\Step\cd\\
\id\step\td \phi\step\cd\step\hx\step\td \phi\step\id\Step\id\\
\id\step\id\Step\hx\Step\hx\step\hx\Step\id\step\id\Step\id\\
\id\step\cu\step\x\step\hx\step\d\step\id\step\id\Step\id\\
\id\Step\id\Step\id\Step\hx\step\id\Step\hx\step\id\Step\id\\
\id\step\cd\step\id\Step\id\step\id\step\tu\alpha\step\hx\Step\id\\
\id\step\id\Step\hx\Step\id\step\id\Step\x\step\cu\\
\id\step\tu \alpha\step\cu\step\x\Step\id\Step\id\\
\cu\step\Step\cu\Step\cu\Step\id\\
\step\object{A}\step\Step\Step\object{H}\step\Step\object{A}\Step\Step\object{H}\\
\end{tangle}
\]
\ =\ \hbox{the right side of (1).}

(ii) Assume relation (3)  holds. Applying $(id _A \otimes \eta _H
 \otimes id _A)$ on relation (3), we get relation (4).
                                  Applying $( \eta _A \otimes id _H
 \otimes id _A)$ on relation (3), we get (BB7').

Conversely, assume relation (4) and (BB7') hold. We see that
\[
\hbox{the left side of (3)}\ \stackrel{by\ (4)}{=}\
\begin{tangle}
\step\object{A}\Step\step\object{H}\Step\object{A}\\
\step\id\Step\cd\step\id\\
\step\id\Step\id\Step\hx\\
\step\id\Step\tu \alpha\step\id\\
\td \phi\step\td \phi\step\id\\
\id\Step\hx\Step\id\step\id\\
\cu\step\cu\step\id\\
\step\d\Step\x\\
\Step\cu\Step\id\\
\step\Step\object{H}\step\Step\object{A}\\
\end{tangle}
\step=\step
\begin{tangle}
\step\object{A}\step\Step\object{H}\Step\object{A}\\
\td \phi\step\cd\step\id\\
\id\Step\id\step\id\Step\hx\\
\id\Step\id\step\tu \alpha\step\id\\
\id\Step\id\step\td \phi\step\id\\
\id\Step\hx\Step\id\step\id\\
\id\Step\id\step\cu\step\id\\
\id\Step\id\Step\x\\
\d\step\cu\Step\id\\
\step\cu\step\Step\id\\
\Step\object{H}\Step\Step\object{A}\\
\end{tangle}
\step\stackrel{(BB7')}{=}\step
\begin{tangle}
\step\object{A}\step\Step\object{H}\step\Step\object{A}\\
\td \phi\step\cd\step\td \phi\\
\id\Step\id\step\id\Step\hx\Step\id\\
\id\Step\id\step\cu\step\tu \alpha\\
\id\Step\x\Step\dd\\
\cu\Step\cu\\
\step\object{H}\Step\Step\object{A}\\
\end{tangle}\]
\ =\ \hbox{the right side of (3).}

(iii)  Applying $id _A \otimes \eta _H \otimes id _A$ on relation
(2), we can get relation (5).
       Applying $id _A \otimes \epsilon  _H \otimes id _A$ on relation (2),
we can get relation (6).

(iv) Assume relation (6) holds.
       Applying $\eta _A \otimes id _H \otimes id _A$ on relation (6),
we can get relation (7).
       Applying $id _A \otimes \eta _H \otimes id _A$ on relation (6),
we can get relation (BB6'). Conversely, assume (BB6') and relation
(7) hold. We see that

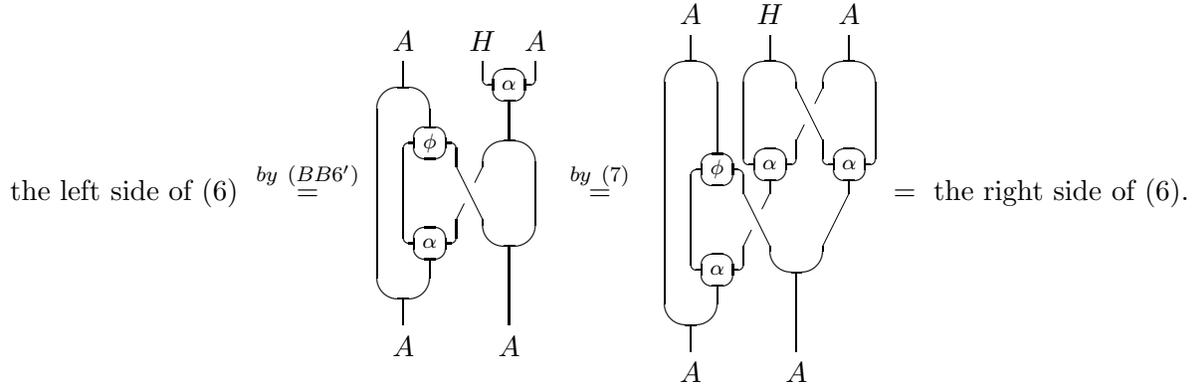
\[
\hbox{the left side of (6)}\ \stackrel{by\ (BB6')}{=}\
\begin{tangle}
\step\object{A}\step\Step\object{H}\Step\object{A}\\
\cd\Step\tu \alpha\\
\id\step\td \phi\step\cd\\
\id\step\id\Step\hx\Step\id\\
\id\step\tu \alpha\step\cu\\
\cu\Step\step\id\\
\step\object{A}\Step\Step\object{A}\\
\end{tangle}
\step\stackrel{by\ (7)}{=}\step
\begin{tangle}
\step\object{A}\step\Step\object{H}\step\Step\object{A}\\
\cd\step\cd\step\cd\\
\id\Step\id\step\id\Step\hx\Step\id\\
\id\step\td \phi\tu \alpha\step\tu \alpha\\
\id\step\id\Step\hx\Step\dd\\
\id\step\tu \alpha\step\cu\\
\cu\Step\step\id\\
\step\object{A}\Step\Step\object{A}\\
\end{tangle}
\ =\ \hbox{the right side of (6).}
\]
(v),(vi) and (vii) are clear.

\begin {Theorem}\label {3.3.1}
Let  $H$  be a bialgebra, $(A, \alpha )$   an  $H$-module algebra
and $(A, \phi )$ an $H$-comodule coalgebra. Then  the following
conditions are equivalent.

(i) $A ^\phi _\alpha \# H$  is a bialgebra.

(ii) $(A, \alpha )$  is an $H$-module coalgebra, $(A, \phi )$  is an
$H$- comodule algebra, $\Delta _A(\eta _A) = \eta _A \otimes \eta
_A,  \epsilon _A $  is an algebra morphism, and (BB6')--(BB7')
hold.

(iii)  $\epsilon _A$  and $\phi$  are algebra morphisms, $\alpha $
is a coalgebra morphism, $\Delta _A(\eta _A) = \eta _A \otimes \eta
_A,$ and (BB6')--(BB7')  hold.
\end {Theorem}
{\bf Proof.} $(ii) \Leftrightarrow  (iii)$ is clear. $(i)
\Rightarrow (ii)$  follows from Lemma \ref {3.3.1.1},
 so it remains to show that
$(ii) \Rightarrow (i).$  By Lemma \ref {3.3.1.1}, we only need show
that relation (2) holds. We see that

\[
\hbox{the left side of (2)}\ \stackrel{by\ (6)}{=}\
\begin{tangle}
\step\object{A}\step\Step\Step\object{H}\Step\Step\object{A}\\
\cd\step\Step\cd\step\Step\id\\
\id\step\td \phi\step\cd\step\d\Step\id\\
\id\step\id\Step\hx\Step\id\Step\x\\
\id\step\cu\step\id\Step\id\step\cd\step\id\\
\id\Step\id\Step\id\Step\hx\Step\id\step\id\\
\id\Step\id\Step\x\step\tu \alpha\step\id\\
\d\step\tu \alpha\Step\cu\Step\id\\
\step\cu\step\Step\td \phi\Step\id\\
\Step\id\Step\Step\id\Step\x\\
\Step\id\Step\Step\cu\Step\id\\
\Step\object{A}\step\Step\Step\object{H}\step\Step\object{A}\\
\end{tangle}
\step\stackrel{by\ (3)}{=}\step
\begin{tangle}
\step\object{A}\Step\step\Step\object{H}\Step\Step\object{A}\\
\cd\Step\Cd\step\cd\\
\id\step\td \phi\step\id\Step\Step\hx\step\td \phi\\
\id\step\id\Step\hx\Step\step\ne2\dd\step\id\Step\id\\
\id\step\cu\step\x\step\cd\step\id\Step\id\\
\d\step\tu \alpha\step\td \phi\d\step\hx\Step\id\\
\step\cu\Step\id\Step\hx\step\id\step\tu \alpha\\
\Step\id\step\Step\cu\step\hx\Step\id\\
\Step\id\Step\Step\cu\step\cu\\
\Step\object{A}\Step\Step\step\object{H}\Step\step\object{A}\\
\end{tangle}\]
\ =\ \hbox{the right side of (2).}  $\Box$

 Dually we have

\begin {Corollary}\label {3.3.2}
(i) If $A$ is an algebra and  coalgebra, and $H$ is a bialgbra, then
$A _\alpha \# H$ is a bialgebra iff $(A, \alpha )$ is an $H$-module
bialgebra and $H$ is cocommutative with respect to $(A, \alpha );$

\end {Corollary}

\begin {Theorem}\label {3.3.8}
Let  $H$  be a Hopf algebra, $(A, \alpha )$  an  $H$-module algebra
and $(A, \phi )$  an $H$-comodule coalgebra with an antipode $S_A$.
Then $D= A^\phi _\alpha \# H$ has an antipode
\[
S_D\ =\
\begin{tangle}
\step\object{A}\step[3]\object{H}\\
\td \phi\Step\id\\
\id\step[2]\x\\
\cu\Step\id\\
\obox 2{{S_H}}\step\obox 2{{S_A}}\\
\cd\Step\id\\
\id\step[2]\x\\
\tu \alpha\Step\id\\
\step\object{A}\step[3]\object{H}\\
\end{tangle}\ \ \ .
\]
\end {Theorem}
{\bf Proof.} We see that
\[
S_D*id_D\ =\
\begin{tangle}
\Step\object{A}\Step\Step\object{H}\\
\step\cd\Step\cd\\
\step\id\step\td \phi\step\id\Step\id\\
\step\id\step\id\Step\hx\Step\id\\
\step\id\step\cu\step\id\Step\id\\
\td \phi\step\id\Step\id\Step\id\\
\id\Step\hx\Step\id\Step\id\\
\cu\step\id\Step\id\Step\id\\
\step\O S\Step\O S\Step\id\Step\id\\
\cd\step\id\Step\id\Step\id\\
\id\Step\hx\Step\id\Step\id\\
\tu \alpha\cd\step\id\Step\id\\
\step\id\step\id\Step\hx\Step\id\\
\step\id\step\tu \alpha\step\cu\\
\step\cu\step\Step\id\\
\Step\object{A}\step\Step\object{H}\\
\end{tangle}
\step\stackrel {\hbox {since (A,$\phi$) is
 an H-comodule coalgebra}}{=}\step
\begin{tangle}
\Step\object{A}\Step\Step\step\object{H}\\
\step\td \phi\Step\step\cd\\
\step\id\step\cd\Step\id\Step\id\\
\step\id\step\d\step\x\Step\id\\
\step\id\Step\hx\Step\id\Step\id\\
\step\cu\step\id\Step\id\Step\id\\
\step\step\O S\Step\O S\Step\id\Step\id\\
\step\cd\step\id\Step\id\Step\id\\
\cd\step\hx\Step\id\Step\id\\
\id\Step\hx\step\d\step\id\Step\id\\
\tu \alpha\step\id\Step\hx\Step\id\\
\step\d\step\tu \alpha\step\cu\\
\Step\cu\Step\step\id\\
\step\Step\object{A}\Step\Step\object{H}\\
\end{tangle}\]
\[
\step\stackrel{\hbox {since (A,$\alpha$) is an H-module
algebra}}{=}\step
\begin{tangle}
\step\object{A}\Step\Step\object{H}\\
\td \phi\Step\cd\\
\id\step\cd\step\id\Step\id\\
\id\step\d\step\hx\Step\id\\
\id\Step\hx\step\id\Step\id\\
\cu\step\id\step\id\Step\id\\
\step\O S\Step\O S\step\id\Step\id\\
\cd\step\id\step\id\Step\id\\
\id\Step\hx\step\id\Step\id\\
\id\step\dd\step\hx\Step\id\\
\id\step\cu\step\cu\\
\tu \alpha\Step\step\id\\
\step\object{A}\Step\Step\object{H}\\
\end{tangle}
\step=\step
\begin{tangle}
\object{A}\Step\object{H}\\
\id\Step\id\\
\QQ \epsilon \Step\QQ \epsilon \\
\Q \eta\Step\Q \eta\\
\id\Step\id\\
\object{A}\Step\object{H}\\
\end{tangle}
\]

and
\[
id_D*S_D\step=\step
\begin{tangle}
\step\object{A}\Step\Step\object{H}\\
\cd\Step\cd\\
\id\step\td \phi\step\id\Step\id\\
\id\step\id\Step\hx\Step\id\\
\id\step\cu\td \phi\step\id\\
\id\Step\id\step\id\Step\hx\\
\id\Step\id\step\cu\step\id\\
\id\Step\id\step\step\O S\Step\O S\\
\id\Step\id\step\cd\step\id\\
\id\Step\id\step\id\Step\hx\\
\id\Step\id\step\tu \alpha\step\id\\
\id\step\cd\step\id\Step\id\\
\id\step\id\Step\hx\Step\id\\
\id\step\tu \alpha\step\cu\\
\cu\step\Step\id\\
\step\object{A}\Step\Step\object{H}\\
\end{tangle}
\step=\step
\begin{tangle}
\step\object{A}\Step\Step\Step\object{H}\\
\cd\step\Step\Step\id\\
\id\step\td \phi\Step\Step\id\\
\id\step\d\step\nw2\Step\step\id\\
\id\Step\id\step\Step\x\\
\id\step\cd\step\cd\step\id\\
\id\step\id\Step\hx\Step\id\step\id\\
\id\step\cu\step\cu\step\id\\
\id\Step\id\Step\step\O S\Step\O S\\
\id\step\cd\step\cd\step\id\\
\id\step\id\Step\hx\Step\id\step\id\\
\id\step\cu\step\cu\step\id\\
\nw2\step\d\Step\xx\\
\Step\id\step\tu\alpha\Step\id\\
\Step\cu\step\Step\id\\
\Step\step\object{A}\Step\Step\object{H}\\
\end{tangle}
\step=\step
\begin{tangle}
\step\object{A}\step\Step\object{H}\\
\cd\Step\id\\
\id\step\td \phi\step\id\\
\id\step\id\Step\hx\\
\id\step\cu\step\id\\
\id\step\cd\step\O S\\
\id\step\id\Step\O S\step\id\\
\id\step\cu\step\id\\
\id\step\cd\step\id\\
\id\step\id\Step\hx\\
\id\step\tu \alpha\step\id\\
\cu\Step\id\\
\step\object{A}\step\Step\object{H}\\
\end{tangle}
\step=\step
\begin{tangle}
\object{A}\Step\object{H}\\
\id\Step\id\\
\QQ \epsilon \Step\QQ \epsilon \\
\Q \eta\Step\Q \eta\\
\id\Step\id\\
\object{A}\Step\object{H}\\
\end{tangle}\ \ \ .
\]
Thus $S_D$ is an antipode of $D$.
\begin{picture}(8,8)
\put(0,0){\line(0,1){8}}\put(8,8){\line(0,-1){8}}
\put(0,0){\line(1,0){8}}\put(8,8){\line(-1,0){8}}
\end{picture}\\

\begin {Lemma}\label {3.3.3'} Let ${\cal C}$ be a symmetric braided
tensor category and $H$ a bialgebra in ${\cal C}$.

(i) If $(M, \alpha, \phi )$ is a left $H$-module and left
$H$-comodule in ${\cal C}$, then (YD) and (BB7') are equivalent.

(ii) Let $(H, R)$  be a quasitriangular Hopf algebra in ${\cal C}$
and $(A, \alpha )$  a Hopf algebra in $({}_H{  \cal M (C)}, C^R)$.
Then $(A,  \phi, \alpha )$  is a Yetter-Drinfeld  $H$- module in
$^H_H{\cal YD( C)}$, where $\phi = (id_H \otimes \alpha ) (C_{H,H}
\otimes id _A ) (R \otimes id _A). $

(iii) Let  $(H, r)$  be a coquasitriangular Hopf algebra in ${\cal
C}$ and $(A, \phi )$  a Hopf algebra in $({}_H{\cal M( C)}, C^r)$.
Then $(A,  \phi, \alpha )$  is a Yetter-Drinfeld  $H$- module in
$^H_H{\cal YD( C)}$, where $\alpha  = ( r \otimes id_A ) (C_{H,H}
\otimes id _A ) (id _H \otimes \phi ). $
\end {Lemma}

Consequently, by Lemma \ref {3.3.3'}, Theorem \ref {3.3.1} and
Theorem \ref {3.3.8} we have following three corollaries

\begin {Corollary}\label {3.3.3}  Let $H$  be a braided
Hopf algebra in symmetric braided tensor category ${\cal C}$  and
$(A, \alpha , \phi )$ be a bialgebra or Hopf algebra in
$({}^H_H{\cal YD(C)}, ^{YD}C)$. Then $A ^\phi _\alpha \# H$
is a bialgebra or Hopf algebra in ${\cal C}$.\\
 \end {Corollary}

 We can easily get the well-known bosonisation theorem in \cite
[Theorem 9.4.12]{Ma95b}.

\begin {Corollary}\label {3.3.3''}  Let $(H, R)$  be a quasitriangular
Hopf algebra in symmetric ${\cal C}$  and $(A, \alpha )$ be a Hopf
algebra in $({}_H{\cal M(C)}, C^R)$. Then $A ^\phi _\alpha \# H$  is
a Hopf algebra in ${\cal C}$, where $\phi = (id_H \otimes \alpha )
(C_{H,H} \otimes id
_A ) (R \otimes id _A). $ \\

 \end {Corollary}

\begin {Corollary}\label {3.3.3'''} Let $(H, r)$  be a coquasitriangular Hopf algebra in
symmetric ${\cal C}$ and $(A, \phi )$  a Hopf algebra in $({}_H{\cal
M( C)}, C^r)$. Then $A ^\phi _\alpha \# H$  is a Hopf algebra in
${\cal C}$, where $\alpha  = ( r \otimes id_A ) (C_{H,H} \otimes id
_A ) (id _H \otimes \phi ). $

 \end {Corollary}

 In particular, these corollaries hold for ordinary (co)quasitriangular
Hopf algebras.

 We now get   the Lagrange's theorem  for braided Hopf algebras.
\begin {Theorem}\label {3.3.4}  Let $H$ be a finite-dimensional
  ordinary  Hopf algebra,
 If $(A, \alpha )$ is a finite-dimensional
Hopf algebra in  ${}_H ^H{\cal YD }$ and $B$ is a subHopf algbra of
$A$ in ${}_{H}^H {\cal YD}$, then
$$dim { \ } B \mid dim { \ } A .$$
That is, the dimension of $B$  divides the dimension of $A$.
\end {Theorem}
{\bf Proof.}  By Corolloary \ref {3.3.3} and \cite [ Corollary
3.2.1] {Mo93},
$$dim { \ } (B ^\phi _\alpha  \# H) \mid dim { \ } (A ^\phi _\alpha \# H).$$ Consequently,
$$dim { \ } B \mid dim { \ } A \hbox { \ \ .   \ \ \ \  } \Box $$

\begin {Definition} \label {3.3.6}
A coalgebra  $A$  is called super cocommutative or quantum
cocommutative if  $\Delta = C_{A,A} \Delta; $ an algebra  $A$  is
called super commutative or quantum commutative
 if  $m_A  = m_A C_{A,A}. $
\end {Definition}
Remark : In the category ${\cal V}$ect (k), a coalgebra  $A$  is
called  cocommutative  if  $\Delta = C_{A,A} \Delta; $ an algebra
$A$  is called  commutative
 if  $m_A  = m_A C_{A,A}. $ Here braiding  $C$  is usual twist map.

\begin {Proposition} \label {3.3.7}
(i) If $H$ is a quantum commutative, or quantum cocommutative Hopf
algebra, then $S^2 = id_H$.

(ii) If $H$ is a Hopf algebra, then the following conditions are
equivalent:

 (1) $ m C_{H,H}^{-1} (id _H \otimes S) \Delta _H = \eta _H\epsilon
 _H;$

 (2) $ m C_{H,H}^{-1} ( S\otimes id _H) \Delta _H = \eta _H\epsilon
 _H;$

 (3) $S^2 = id _H$;

\end {Proposition}

{\bf Proof.} (ii) We first show that (2) implies (3). If (2) holds,
we see that
\begin {eqnarray*}
  m  (id _H \otimes S) (S \otimes S)\Delta _H &=&
  m C_{H,H}^{-1} ( S\otimes id _H) \Delta _H S
  \ \  \ \ ( \hbox {by Pro. \ref {12.1.4}})\\
&= &\eta _H\epsilon
 _H \ \ \ \ \ ( \hbox {by (2) })
\end {eqnarray*} and
\begin {eqnarray*}
  m  (S\otimes id _H  ) (S \otimes S)\Delta _H &=&
  S m C_{H,H}^{-1} ( S\otimes id _H) \Delta _H
  \ \  \ \ ( \hbox {by Pro. \ref {12.1.4}})\\
&= &\eta _H\epsilon
 _H \ \ \ \ \ ( \hbox {by (2) }).
\end {eqnarray*} Thus $S^2$ is the  convolution inverse of $S$,
which implies $S^2 = id _H.$ Therefore (2) implies (3). Similarly,
(1) implies (3).

We next show that (3) implies (1). If (3) holds, we see that
\begin {eqnarray*}
\eta _H\epsilon &=&   m   (id _H \otimes S)\Delta _H \\
&=& m  (S\otimes id _H  ) (S \otimes S)\Delta _H \\
&=&   m  ( S\otimes id _H) C_{H,H}^{-1}\Delta _H S
  \ \  \ \ ( \hbox {by Pro. \ref {12.1.4}})\\
\end {eqnarray*}
and
\begin {eqnarray*}
\eta _H\epsilon &=&   m   (id _H \otimes S)\Delta _H \\
&=& m  (S\otimes id _H  ) (S \otimes S)\Delta _H \\
&=&   m  ( S\otimes id _H) C_{H,H}^{-1}\Delta _H S
  \ \  \ \ ( \hbox {by Pro. \ref {12.1.4}})\\
&=&\eta _H\epsilon S^{-1}\\
m C_{H,H}^{-1} (id _H \otimes S) \Delta _H.
\end {eqnarray*}

Thus (1) holds. Similarly (3) implies (2).

(i) It follows from (ii) $\Box$

 In this subsection, we construct a braided Hopf algebra in a
 Yetter-Drinfeld category by means of a graded Hopf algebra.

\begin {Lemma} \label {4.1'}
(see \cite [p.1530] {Ni78},\cite {AS98a} and \cite {Ra85}) Assume
that $H$ and $\Lambda$ are two ordinary Hopf algebras> If there
exist two bialgebra homomorphisms $\pi _0 : \ \Lambda \rightarrow H$
and $\iota _0 : \ H \rightarrow \Lambda $ with $\pi_0\iota _0= id
_H$, then

(i) $(\Lambda, \delta ^+, \alpha ^+)$ is a right $H$-Hopf module
with
 $\delta ^+ =: (id \otimes \pi _0)
\Delta $ and  $\alpha ^+ =: m  (id \otimes \iota_0)$

(ii) $A\# H \stackrel {\alpha ^+} {\cong } \Lambda $  \ \ ( as
algebras) with $A=:\Lambda^{co H}$.

(iii) $\Lambda^{coH} =Im \ \omega $ with $\omega =: id _\Lambda
* (\iota _0 \circ S_H \circ \pi_H) : \Lambda \rightarrow \Lambda $.

\end{Lemma}
{\bf Proof.} (i) It is easy to check directly that $(\Lambda, \delta
^+, \alpha ^+)$ is a right $H$-Hopf module.

(ii) By the fundamental theorem of Hopf modules, we have $
\Lambda^{co H } \otimes H \cong \Lambda$ as $H$-Hopf modules with
isomorphism $\alpha ^+$. We can also check that $ \Lambda^{co H}$ is
an $H$-module algebra with adjoint action  $ad $. Moreover, $\alpha
^+$ preserves the algebra operations of $\Lambda^{co H} \# H$.

(iii) Obviously, for any $x\in A$, $\omega (x) =x$, which implies $A
\subseteq Im \ (\omega ) $. On the other hand, for any $x\in \Lambda
$, see that
\begin {eqnarray*} \delta ^+ \circ \omega (x) &=&
\sum (id \otimes \pi_0) \Delta (x_1 (\iota _0 \circ S \circ \phi _0 (x_2)))\\
&=&\sum x_1\iota _0 \circ S \circ \phi _0 (x_4) \otimes \phi _0
(x_2) S \circ \phi _0 (x_3) \\
&=& \omega (x) \otimes 1.
\end {eqnarray*}
Thus $A = Im \ (\omega).$ $\Box$

\begin {Theorem} \label {4.1}
(see \cite [p.1530] {Ni78},\cite {AS98a} and \cite {Ra85})Assume
that $H$ and $\Lambda$ are two ordinary Hopf algebras> If there
exist two bialgebra homomorphisms $\pi _0 : \ \Lambda \rightarrow H$
and $\iota _0 : \ H \rightarrow \Lambda $ with $\pi_0\iota _0= id
_H$, then
 $(A, \alpha, \phi )$ is a braided Hopf algebra in
$^H_H {\cal YD}$
  and  $A   \ ^\phi _\alpha \#  H \stackrel {\alpha^+} { \cong } \Lambda $ \ \ ( as Hopf algebras),
where $\alpha (h \otimes a) = \sum \iota _0(h_1) a \iota _0
(S(h_2))$, $\phi (a)= (\pi_0 \otimes id ) \Delta _{\Lambda }(a),$
for any $h\in H, a\in A.$
\end{Theorem}

{\bf Proof.} Define $\Delta _A = (\omega \otimes id ) \Delta $ and
$\epsilon
 _A= \epsilon _\Lambda $. By the proof of \cite [Theorem
3]{Ra85}, $(A, m  , \eta , \Delta _A, \epsilon )$ become a braided
Hopf algebra in $^H_H {\cal YD}$ under module operation $ \alpha$
and comodule operation $\phi$.

 Set $\Phi = \alpha ^ +.$  For any $x\in A, h\in H$, see that
\begin {eqnarray*} \Delta \circ \Phi  (x \#  h) &=&
\sum x_1 \iota _0 (h_1) \otimes x_2 \iota _0 (h_2)\\
&=& \sum \alpha^ +( x_1, h_1) \otimes  \alpha ^+( x_2, h_2) \\
&=& (\Phi \otimes \Phi ) \Delta ( x \# h).
\end {eqnarray*}
Thus $\Phi $ is a coalgebra map. By Lemma \ref {4.1'} (ii), we have
$\Phi$ is a Hopf algebra isomorphism. $\Box$

Assume that  $\Lambda  =\oplus _{i=0} ^{\infty} \Lambda _i$ is a
graded Hopf algebra with $H = \Lambda _0$. Let $\pi _0 : \Lambda
\rightarrow \Lambda _0$ and $  \iota _0 : \Lambda _0 \rightarrow
\Lambda $ denote the canonical projection and injection. It is clear
that $\pi_0$ and $\iota_0$ are two bialgebra homomorphisms.
Consequently, Theorem above holds. $\Lambda^{co H}$ is called the
diagram of $\Lambda$, written $diag (\Lambda ).$

\chapter{The Factorization Theorem of Bialgebras in Braided
Tensor Categories }\label {c5}

It is well-known  that  the factorization of domain  plays an
important role in ring theory. S. Majid in \cite [Theorem 7.2.3]
{Ma95b} studied the factorization  of Hopf algebra and showed that
$H \cong A \bowtie H$ for two subbialgebras $A$ and $B$ when
multiplication $m_H$ is bijective.

S. Majid found a method to turn an ordinary Hopf algebra $H$ into a
braided Hopf algebra $\underline H$. This method is called the
transmutation, and $\underline H$ is called the braided group
analogue of $H$ (see \cite {Ma95a} and \cite {Ma95b} ).
 Huixiang Chen in
\cite {Ch98} showed that the double cross coproduct $A  \bowtie ^R
H$ of two quasitriangular Hopf algebras is a quasitriangular Hopf
algebra. One needs to know the relation among the braided group
analogues  of the double cross coproduct $A \bowtie ^R H$, $A$ and
$H$.

In this chapter we generalize  Majid's factorization theorem
  into  braided tensor categories.
We show that the braided group analogue of double cross (co)product
is double cross (co)product of braided group analogues. We give  the
factorization theorem  of Hopf algebras  and the relation between
Hopf algebras and their factors in braided tensor categories.

\section{the braided reconstruction theorem}\label {s8}
In this section, we introduce the method ( transmutation ) turning
 (co)quasi- triangular Hopf algebra into a
braided Hopf algebra, which is due to S. Majid \cite {Ma95a}.

Let ${\cal C}$ be a tensor category,  ${\cal D}$  a
 braided tensor category and
$(F, \mu _0 , \mu) $   a tensor functor from ${\cal C}$  to ${\cal
D}$ with $\mu_0 = id _I$. Let $Nat (G,T)$ denote all  the natural
transformations from functor $G$ to functor $ T$. Assume  that there
is  an object  $B$ of ${\cal D}$ and a natural transformation
$\alpha $ in $Nat (B\otimes F, F).$ Here $(B \otimes F) (X) = B
\otimes F(X)$  for any object $X $ in ${\cal D}$. Let
\[
\begin{tangle}
\theta_V(g)_X
\end{tangle}
\step=\step
\begin{tangle}
\object{V}\step[5]\object{F(X)}\\
\O g\step[4]\dd\\
\d\step[2]\dd\\
\obox 4{\alpha_X}\\
\Step\id\\
\Step\object{F(X)}\\
\end{tangle}
\step[3],\step[2]
\begin{tangle}
\theta_V^{(2)}(h)_{X_1\otimes X_2}
\end{tangle}
\step=\step
\begin{tangle}
\step\object{V}\step[4]\object{F(X_1)}\step[6]\object{F(X_2)}\\
\td h\Step\dd\step[5]\ne2\\
\d\step\x\step[4]\ne2\\
\obox 3{\alpha_{X_1}}\step\obox 3{\alpha_{X_2}}\\
\step[1.5]\id\step[4]\id\\
\step\object{F(X_1)}\step[6]\object{F(X_2)}\\
\end{tangle}
\step[3],\step[2]
\]
\[
\begin{tangle}
\theta_V^{(3)}(p)_{X_1\otimes X_2\otimes X_3}
\end{tangle}
 \step=\step
\begin{tangle}
\step[3]\object{V}\step[4]\object{F(X_1)}\step[6]\object{F(X_2)}\step[6]\object{F(X_3)}\\
\step[3]\O p\step[3]\dd\step[5]\dd\step[4]\dd\\
\Step\dd\id\d\Step\id\step[5]\dd\step[4]\dd\\
\step\dd\step\id\step\x\step[4]\dd\step[4]\dd\\
\dd\Step\hx\Step\d\Step\dd\step[4]\dd\\
\d\step\dd\step\d\step[2]\x\step[4]\ne2\\
\obox 3{\alpha_{X_1}}\step\obox 3{\alpha_{X_2}}\Step\obox 4{\alpha_{X_2}}\\
\step[1.5]\id\step[4]\id\step[5]\id\\\\
\object{F(X_1)}\step[5.5]\object{F(X_2)}\step[6]\object{F(X_3)}\\
\end{tangle}
\]
for any object $X, X_1, X_2, X_3\in ob{\cal C}, g\in Hom_{\cal
D}(V,B),h\in Hom_{\cal D}(V,B\otimes B), p\in Hom_{\cal
D}(V,B\otimes B\otimes B).$

Let $V \in ob  {\cal D},$ $g, g'\in Hom_{\cal D}(V,B).$  $\theta _V$
is called injective if  $\theta _V (g) _{X} = \theta (g')_X$ for any
$X \in ob {\cal D}$ implies $g =g'$. Similarly, we can obtain the
definitions  about  $\theta _V^{(2)}$ and $\theta _V^{(3)}$.

\begin {Proposition} \label {12.1.2} (see \cite [Proposition 3.8]
{{Ma95a}} \cite [Theorem 9.4.6, Proposition 9.4.7] {Ma95b}
 If $\theta _V$, $\theta _V^{(2)}$ and  $\theta _V^{(3)}$ are injective
 for any
object $V =I,  B, B\otimes B, B\otimes B\otimes B$, then

(i)  $B$  is a bialgebra living in ${\cal D}$,
 called the braided  bialgebra determined by
 braided reconstruction;

(ii) Furthermore, if ${\cal C}$  is
  a rigid tensor category ( i.e. every object has a
left duality ),   then $B$  is Hopf algebra,
 called the braided Hopf algebra determined by
 braided reconstruction.

(iii) Furthermore, if ${\cal D}$  and ${\cal C }$ is  braided tensor
categories, then $B$ is  quasitriangular.
\end {Proposition}

{\bf Proof.} (i) We define the multiplication, unit ,
comultiplication and counit of $B$  as follows:

\[
\begin{tangle}
\theta_{B\otimes B}(m_B)_X
\end{tangle}
\step=\step
\begin{tangle}
\object{B}\Step\object{B}\step[4]\object{F(X)}\\
\id\Step\d\step\dd\\
\id\Step\obox 3{\alpha_X}\\
\d\step[1.5]\dd\\
\obox 3{\alpha_X}\\
\step[1.5]\id\\
\step[1.5]\object{F(X)}\\
\end{tangle}
\step[3],\step[2]
\begin{tangle}
\theta_I(\eta_B)
\end{tangle}
\step=\step[4]
\begin{tangle}
\object{F(X)}\\
\id\\
\object{F(X)}\\
\end{tangle}
\step[3],\step[2]
\]\[
\begin{tangle}
\theta_B^{(2)}(\Delta_B)_{X\otimes Y}
\end{tangle}
\step=\step
\begin{tangle}
\object{B}\step[4]\object{F(X)}\step[6]\object{F(Y)}\\
\id\step[4]\d\step[3]\ne2\\
\d\step[4]\tu \mu\\
\step\d\step[3.5]\ne2\\
\step\obox 5{\alpha_{X\otimes Y}}\\
\step[3.5]\id\\
\step[2]\obox 4{\mu^{-1}}\\
\step[2]\dd\step\d\\
\step[1.5]\object{F(X)}\step[6]\object{F(Y)}
\end{tangle}
\step[3],\step[2]
\begin{tangle}
\epsilon _B
\end{tangle}
\step=\step
\begin{tangle}
\object{B}\step[4]\object{F(I)}\\
\d\step\dd\\
\obox 3{\alpha_I}\\
\end{tangle}
\]
for any object $X , Y, Z \in ob {\cal C}.$\\

We only show that the coassociative law holds, the others can be
shown similarly. That is,  we need show that

\[
\begin{tangle}
\Step\object{B}\\
\step\cd\\
\cd\step\id\\
\object{B}\Step\object{B}\step[1.5]\object{B}\\
\end{tangle}
\step=\step
\begin{tangle}
\step\object{B}\\
\cd\\
\id\step\cd\\
\object{B}\step[1.5]\object{B}\step[1.5]\object{B}\\
\end{tangle}
\step[5]......(1)
\]

We see that
\[
\theta_B^{(3)}(
\begin{tangle}
\Step\object{B}\\
\step\cd\\
\cd\step\id\\
\object{B}\Step\object{B}\step[1.5]\object{B}\\
\end{tangle})_{X\otimes Y\otimes Z}
\step=\step
\begin{tangle}
\step[3]\object{B}\step[4]\object{F(X)}\step[5]\object{F(Y)}\step[5]\object{F(Z)}\\

\Step\cd\Step\dd\step[3]\dd\step[4]\dd\\

\step\cd\step\x\step[3]\dd\step[4]\dd\\

\step\id\Step\hx\Step\d\step\dd\step[4]\dd\\

\obox 3{\alpha_X}\step\d\step[2]\hx\step[4]\dd\\
\step[1.5]\id\step[2.5]\obox 3{\alpha_Y}\step\d\step[2]\dd\\
\step[1.5]\id\step[4]\id\step[2.5]\obox 3{\alpha_Z}\\
\step[1.5]\id\step[4]\id\step[4]\id\\
\step\object{F(X)}\step[5]\object{F(Y)}\step[5]\object{F(Z)}\\
\end{tangle}\]\[
\step[2]=\step
\begin{tangle}
\step\object{B}\object{F(X)}\step[5]\object{F(Y)}\step[5]\object{F(Z)}\\
\cd\step[3]\d\Step\dd\step[5]\dd\\
\id\Step\d\step[3]\tu \mu\step[5]\dd\\
\id\step[3]\d\Step\dd\step[5]\dd\\
\id\step[4]\x\step[5]\dd\\
\d\step[2]\dd\Step\d\step[3]\dd\\
\obox 5{\alpha_{X\otimes Y}}\step[2]\id\Step\dd\\
\step[2.5]\id\step[3.5]\obox 3{\alpha_Z}\\
\step\obox 3{\mu^{-1}}\step[3.5]\id\\
\step\dd\step\d\step[3.5]\id\\
\object{F(X)}\step[5]\object{F(Y)}\step[5]\object{F(Z)}\\
\end{tangle}
\step=\step
\begin{tangle}
\object{B}\step[4]\object{F(X)}\step[5]\object{F(Y)}\step[5]\object{F(Z)}\\
\d\step[3]\d\Step\dd\step[4]\ne2\\
\step\d\step[3]\tu \mu\step[3]\ne3\\
\Step\d\step[3]\tu \mu\\
\step[3]\d\Step\dd\\
\step\obox 9{\alpha_{(X\otimes Y)\otimes Z}}\\
\step[5.5]\id\\
\step[4]\obox 3{\mu^{-1}}\\
\step[4]\dd\step\d\\
\step[1.5]\obox 3{\mu^{-1}}\step[2.5]\d\\
\Step\dd\step\d\step[3]\d\\
\object{F(X)}\step[5]\object{F(Y)}\step[5]\object{F(Z)}
\end{tangle}
\]
\[\ \ \
\stackrel {\hbox {since } (F, \mu) \hbox { is tensor functor} } {=}
\ \
\begin{tangle}
\object{B}\step[4]\object{F(X)}\step[5]\object{F(Y)}\step[5]\object{F(Z)}\\
\d\step[3]\nw2\step[3]\d\Step\dd\\
\step\d\step[4]\nw2\step[2]\tu \mu\\
\Step\d\step[5]\tu \mu\\
\step[3]\d\step[4]\dd\\
\step[1.5]\obox 9{\alpha_{(X\otimes Y)\otimes Z}}\\
\step[6]\id\\
\step[4.5]\obox 3{\mu^{-1}}\\
\step[5]\dd\step\d\\
\step[4]\dd\step[1.5]\obox 3{\mu^{-1}}\\
\step[3]\dd\step[2]\dd\step[2]\d\\
\step[2]\object{F(X)}\step[5]\object{F(Y)}\step[5]\object{F(Z)}
\end{tangle}
\step=\step \theta_B^{(3)}(
\begin{tangle}
\step\object{B}\\
\cd\\
\id\step\cd\\
\object{B}\step[1.5]\object{B}\step[1.5]\object{B}\\
\end{tangle})_{X\otimes Y\otimes Z}.
\]
Thus relation (1) holds since $\theta_B^{(3)}$ is injective.\\

(ii)We define
\[
\begin{tangle}
\theta_B(S_B)_X
\end{tangle}
\step=\step
\begin{tangle}
\object{B}\step[10]\object{F(X)}\\
\id\step\obox 5{b_{F(X)}}\step[4]\id\\
\x\step[2]\dd\step[5]\id\\
\id\step\obox 4{\alpha_{X^*}}\step[4]\dd\\
\id\step[3]\d\step[4]\ne2\\
\id\step[3.5]\obox 5{d_{F(X)}}\\
\object{F(X)}\\
\end{tangle}
\]
for any $X \in ob {\cal C}.$ We can show that $S_B$  is the antipode
of $B$ .In fact\\
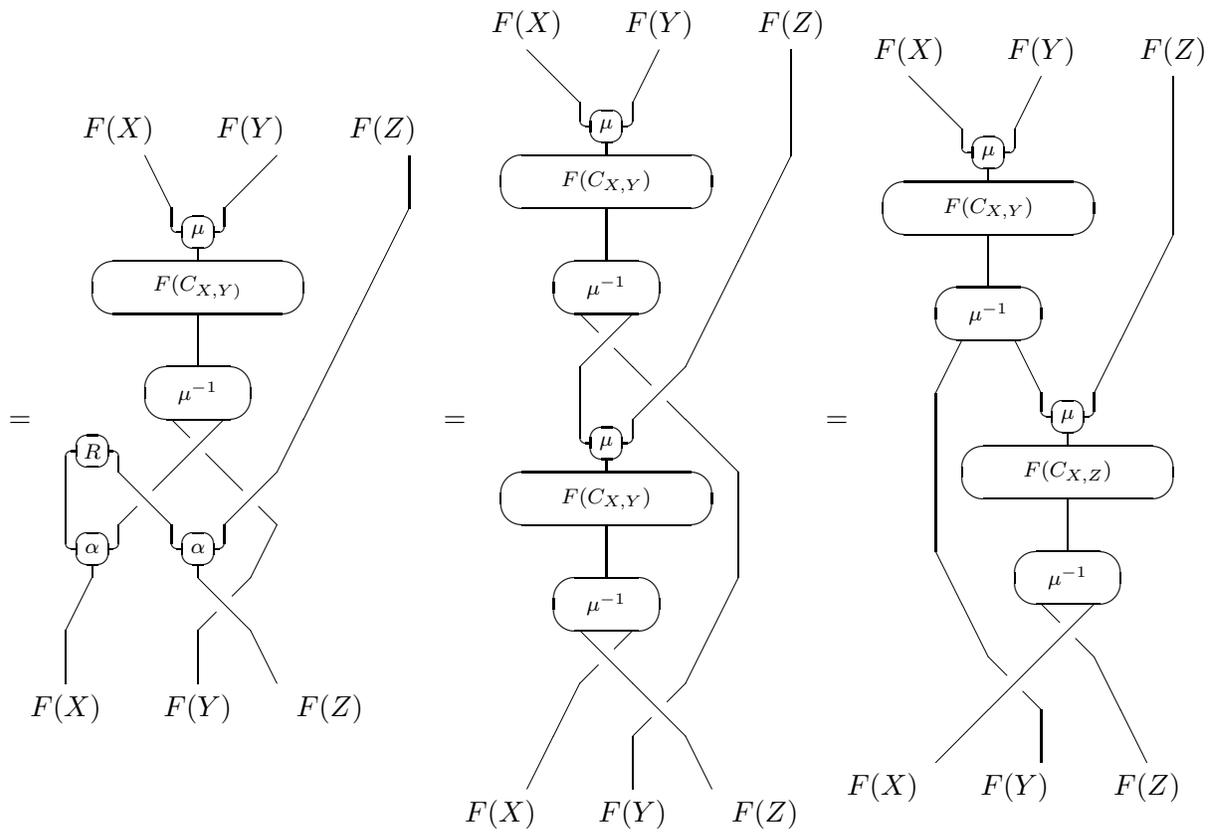
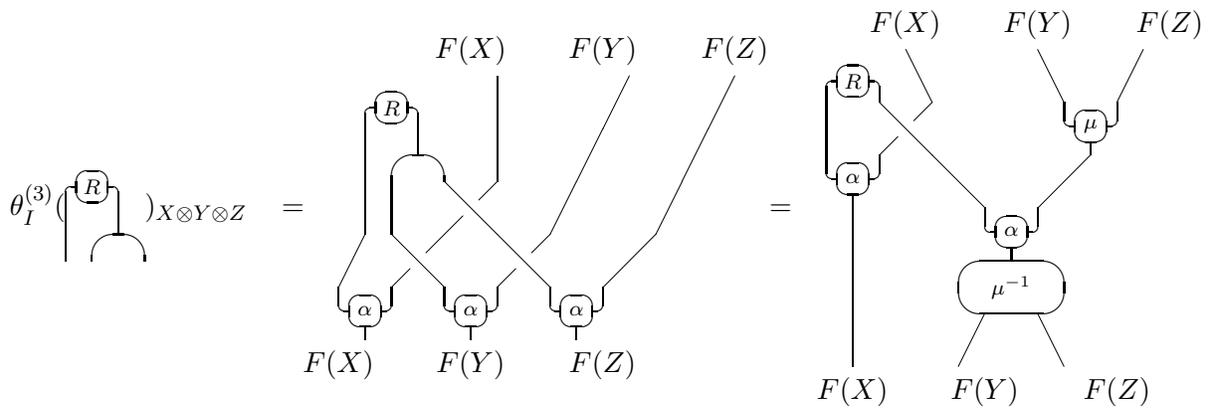
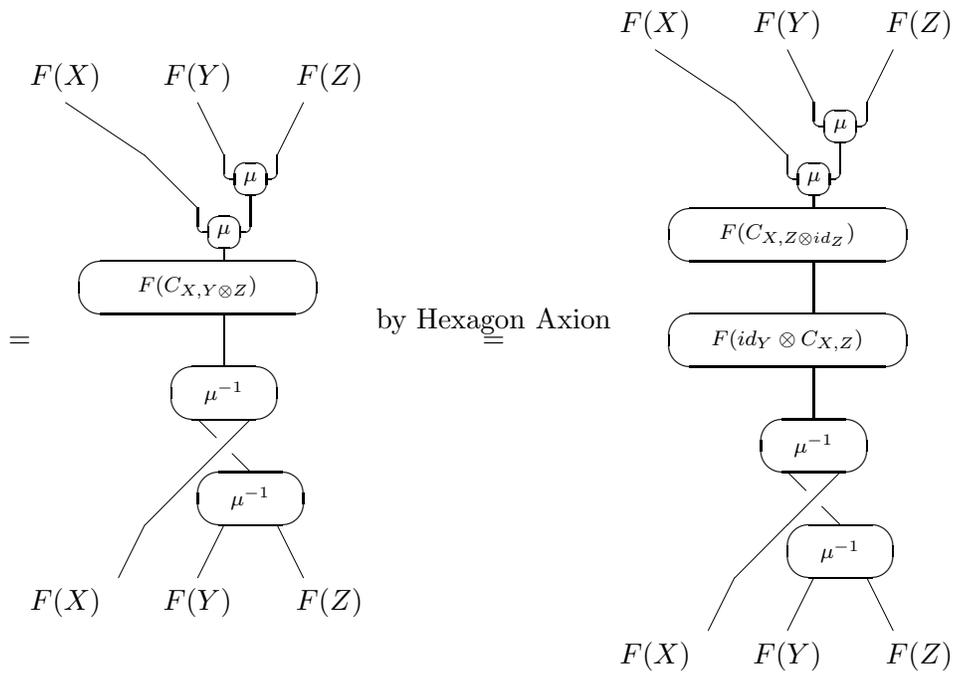
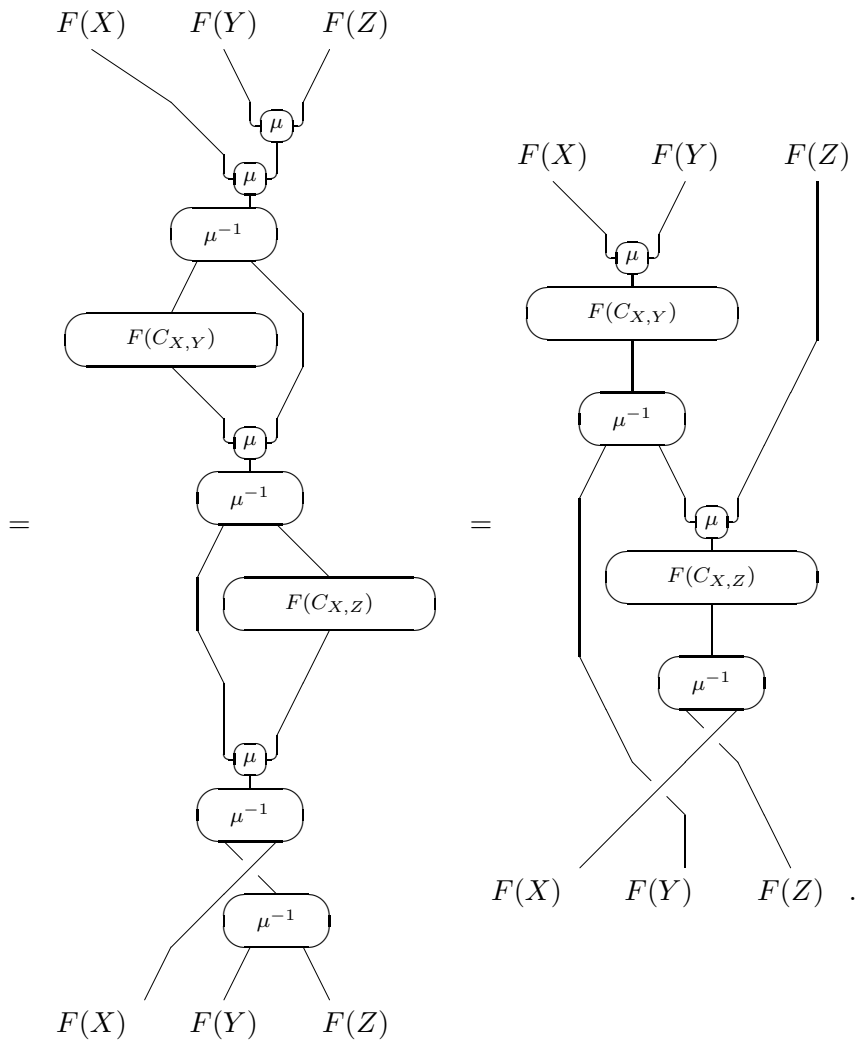
\[
\begin{tangle}
\theta_B(S_B\ast I)_X
\end{tangle}
\step=\step
\begin{tangle}
\step\object{B}\step[6]\object{F(X)}\\
\cd\step[3]\dd \\
\S\step[2]\id\step[2]\dd\\
\cu\step\dd\\
\step\tu {\alpha_{X}}\\
\step[2]\id \\
\step[2]\object{F(X)}
     \end{tangle}
\step=\step
\begin{tangle}
\step\object{B}\step[6]\object{F(X)}\\
\cd\step[3]\dd \\
\S\step[2]\id\step[2]\dd\\
\d\step\tu{\alpha_{X}}\\
\step\tu {\alpha_{X}}\\
\step[2]\id \\
\step[2]\object{F(X)}
     \end{tangle}
\step=\step
\begin{tangle}
\step\object{B}\step[9]\object{F(X)}\\
\step\id\step[3]\Coev\step[4]\id\\
\step\d\step[2]\id\step[3]\id\step[2]\id\\
\step[2]\x\step[3]\id\step[2]\id\\
\step\dd\step\cd\step\dd\step[2]\id\\
\step\id\step[2]\id\step[2]\hx\step[2]\dd\\
\step\id\step[2]\tu {\alpha_{X*}}\step\tu{\alpha_{X}}\\
\step\id\step[3]\Ev \\
\step\object{F(X)}
\end{tangle}\]\[
\step =\step
\begin{tangle}
\step\object{B}\step[11]\object{F(X)}\\
\step\id\step[4]\Coev\step[6]\id\\
\step\d\step[3]\id\step[3]\id\step[4]\id\\
\step[2]\d\step[2]\id\step[3]\id\step[4]\id\\
\step[3]\x\step[3]\id\step[4]\id\\
\step[2]\dd\step\cd\step[2]\id\step[4]\id\\
\step[2]\id\step[2]\id\step[2]\x\step[4]\id\\
\step[2]\id\step[2]\id\step[2]\id\step[2]\d\step[2]\dd\\
\step[2]\id\step[2]\id\step[2]\id\step[3]\tu{\alpha_{X}}\\
\step[2]\id\step[2]\S\step[2]\id\step[3]\dd\\
\step[2]\id\step[2]\x\step[2]\dd\\
\step[2]\id\step[2]\id\step[2]\tu {\alpha_{X}}\\
\step[2]\id\step[2]\Ev\\
\step[2]\object{F(X)}
\end{tangle}
\step=\step
\begin{tangle}
\step\object{B}\step[8]\object{F(X)}\\
\step\id\step[3]\coev\step[3]\id\\
\step\d\step[2]\id\step[2]\id\step[3]\id\\
\step[2]\x\step[2]\id\step[3]\id\\
\step\dd\step\cd\step\id\step[3]\id\\
\step\id\step[2]\S\step[2]\id\step\id\step[3]\id\\
\step\id\step[2]\cu\step\id\step[3]\id\\
\step\id\step[3]\x\step[2]\dd\\
\step\id\step[3]\id\step[2]\tu{\alpha_{X}}\\
\step\id\step[3]\Ev\\
\step\object{F(X)}
\end{tangle}
\step=\step
\begin{tangle}
\object{B}\step[3]\object{F(X)}\\
\id\step[3]\id\\
\QQ \epsilon\step[3]\id\\
\step[3]\object{F(X)}
\end{tangle}\]
Similarly, we get that
\[
\begin{tangle}
\theta_B(I\ast S_{B})_X
\end{tangle}
\step=\step
\begin{tangle}
\object{B}\step[3]\object{F(X)}\\
\id\step[3]\id\\
\QQ \epsilon\step[3]\id\\
\step[3]\object{F(X)}
\end{tangle}\]

Obviously,\[
\begin{tangle}
\theta_B(I\ast S_{B})_X
\end{tangle}
\step=\step
\begin{tangle}
\object{B}\step[5]\object{F(X)}\\
\id\step[4]\id\\
\QQ \epsilon\step[3]\dd\\
\Q\eta\step[2]\dd\\
\tu {\alpha_{X}}\\
\step\object{F(X)}
\end{tangle}\]
We showed that $S_{B}$ is the antipode of $B$.

(iii) We define
\[
\begin{tangle}
\theta_I^{(2)}(R)_{X\otimes Y}
\end{tangle}
\step=\step
\begin{tangle}
\step\object{F(X)}\step[6]\object{F(Y)}\\
\step[2]\d\step[2]\dd\\
\step[3]\tu \mu\\
\obox 8{F(C_{X,Y})}\\
\step[4]\id\\
\step[2]\obox 4{\mu^{-1}}\\
\step[3]\xx\\
\Step\dd\Step\d\\
\step\object{F(X)}\step[6]\object{F(Y)}\\
\end{tangle}
\step[3]and\step[3]
\begin{tangle}
\theta_B^{(2)}(\bar{\Delta})_{X\otimes Y}
\end{tangle}
\step=\step
\begin{tangle}
\step\object{B}\step[4]\object{F(X)}\step[6]\object{F(Y)}\\
\cd\step[3]\d\step[2]\dd\\
\d\step\nw2\step[3]\x\\
\step\d\step[2]\x\Step\id\\
\step[2]\tu \alpha\step[2]\tu \alpha\\
\step[3]\d\step[2]\dd\\
\step[4]\xx\\
\step[3]\ne2\Step\nw2\\
\step\object{F(X)}\step[6]\object{F(Y)}\\
\end{tangle}\ \ \ .
\]
 We easily check that $(B,\bar{\Delta},\epsilon )$ is a
 coalgebra by modifying the proof of part (i).\\

We see that
\[
\theta_I^{(3)}(
\begin{tangle}
\Ro R\\
\id\step[4]\d\\
\id\Step\ro R\step\id\\
\cu\Step\id\step\id
\end{tangle})_{X\otimes Y\otimes Z}
\step=\step
\begin{tangle}
\step[7]\object{F(X)}\step[5]\object{F(Y)}\step[5]\object{F(Y)}\\
\Ro R\step[3]\id\step[4]\dd\step[4]\dd\\
\id\step[4]\d\Step\id\step[3]\dd\step[4]\dd\\
\id\Step\ro R\step\x\Step\dd\step[4]\ne2\\
\cu\Step\hx\Step\x\step[3]\ne2\\
\step\id\Step\dd\step\tu \alpha\Step\tu \alpha\\
\step\tu \alpha\step[3]\id\step[4]\id\\
\step\object{F(X)}\step[5]\object{F(Y)}\step[5]\object{F(Z)}\\
\end{tangle}
\]\[
\step=\step
\begin{tangle}
\step[7]\object{F(X)}\step[5]\object{F(Y)}\step[5]\object{F(Y)}\\
\Ro R\step[3]\id\step[4]\dd\step[4]\dd\\
\id\step[4]\d\Step\id\step[3]\dd\step[4]\dd\\
\id\Step\ro R\step\x\Step\dd\step[4]\ne2\\
\id\Step\id\Step\hx\Step\x\step[3]\ne2\\
\d\step\tu \alpha\step\tu \alpha\Step\tu \alpha\\
\step\tu \alpha\step[3]\id\step[4]\id\\
\step\object{F(X)}\step[5]\object{F(Y)}\step[5]\object{F(Z)}\\
\end{tangle}
\step=\step
\begin{tangle}
\step[2]\object{F(X)}\step[5]\object{F(Y)}\step[5]\object{F(Z)}\\
\step[3]\d\step[3]\ne2\step[5]\id\\
\step[4]\tu \mu\step[7]\id\\
\step\obox 8{F(C_{X,Y)}}\step[4]\id\\
\step[5]\id\step[7]\dd\\
\step[3]\obox 4{\mu^{-1}}\step[4]\dd\\
\ro R\Step\xx\step[4]\dd\\
\id\Step\x\Step\id\step[3]\dd\\
\tu \alpha\Step\x\Step\dd\\
\dd\step[3]\id\Step\tu \alpha\\
\id\step[4]\id\step[3]\d\\
\object{F(X)}\step[5]\object{F(Y)}\step[5]\object{F(Z)}\\
\end{tangle}
\]\[
\step=\step
\begin{tangle}
\step[2]\object{F(X)}\step[5]\object{F(Y)}\step[5]\object{F(Z)}\\
\step[3]\d\step[3]\ne2\step[5]\id\\
\step[4]\tu \mu\step[6]\dd\\
\step\obox 8{F(C_{X,Y)}}\step[2]\dd\\
\step[5]\id\step[5]\dd\\
\step[3]\obox 4{\mu^{-1}}\step[2]\dd\\
\ro R\Step\xx\step[2]\dd\\
\id\Step\x\Step\xx\\
\tu \alpha\Step\tu \alpha\step\dd\\
\dd\step[4]\x\\
\id\step[5]\id\step[2]\d\\
\object{F(X)}\step[5]\object{F(Y)}\step[5]\object{F(Z)}\\
\end{tangle}
\step=\step
\begin{tangle}
\step\object{F(X)}\step[5]\object{F(Y)}\step[5]\object{F(Z)}\\
\step\nw2\step[3]\dd\step[5]\id\\
\step[3]\tu \mu\step[6]\id\\
\obox 8{F(C_{X,Y})}\Step\dd\\
\step[4]\id\step[5]\dd\\
\Step\obox 4{\mu^{-1}}\step[2]\dd\\
\step[3]\xx\Step\dd\\
\step[3]\id\Step\xx\\
\step[3]\tu \mu\Step\nw2\\
\obox 8{F(C_{X,Y})}\step\id\\
\step[4]\id\step[5]\id\\
\Step\obox 4{\mu^{-1}}\step[2]\dd\\
\step[3]\x\Step\dd\\
\step[2]\dd\Step\x\\
\step\dd\step[3]\id\Step\d\\
\object{F(X)}\step[5]\object{F(Y)}\step[5]\object{F(Z)}\\
\end{tangle}
\step=\step
\begin{tangle}
\step\object{F(X)}\step[5]\object{F(Y)}\step[5]\object{F(Z)}\\
\step\nw2\step[3]\dd\step[5]\id\\
\step[3]\tu \mu\step[6]\id\\
\obox 8{F(C_{X,Y})}\step[3]\id\\
\step[4]\id\step[6]\dd\\
\Step\obox 4{\mu^{-1}}\step[3]\dd\\
\step[2]\dd\Step\d\Step\dd\\
\step[2]\id\step[4]\tu \mu\\
\step[2]\id\step\obox 8{F(C_{X,Z})}\\
\step[2]\id\step[5]\id\\
\step[2]\d\Step\obox 4{\mu^{-1}}\\
\step[3]\d\Step\xx\\
\step[4]\xx\Step\d\\
\step[3]\ne2\Step\id\step[3]\d\\
\object{F(X)}\step[5]\object{F(Y)}\step[5]\object{F(Z)}\\
\end{tangle}
\]

and
\[
\theta_I^{(3)}(
\begin{tangle}
\ro R\\
\id\step\cd\\
\end{tangle})_{X\otimes Y\otimes Z}
\step=\step
\begin{tangle}
\step[5]\object{F(X)}\step[5]\object{F(Y)}\step[5]\object{F(Z)}\\
\step\ro R\step[3]\id\step[4]\dd\step[3]\dd\\
\step\id\step\cd\Step\id\step[3]\dd\step[3]\dd\\
\step\id\step\id\Step\x\Step\dd\step[3]\dd\\
\dd\step\x\Step\x\step[3]\ne2\\
\tu \alpha\Step\tu \alpha\Step\tu \alpha\\
\object{F(X)}\step[5]\object{F(Y)}\step[5]\object{F(Z)}\\
\end{tangle}
\step=\step
\begin{tangle}
\step[3]\object{F(X)}\step[5]\object{F(Y)}\step[5]\object{F(Z)}\\
\ro R\step\d\step[4]\d\Step\dd\\
\id\Step\x\step[5]\tu \mu\\
\tu \alpha\step[2]\nw2\step[4]\ne2\\
\step\id\step[5]\tu \alpha\\
\step\id\step[4]\obox 4{\mu^{-1}}\\
\step\id\step[4]\dd\Step\d\\
\step\object{F(X)}\step[5]\object{F(Y)}\step[5]\object{F(Z)}\\
\end{tangle}
\]
\[
\step=\step
\begin{tangle}
\object{F(X)}\step[5]\object{F(Y)}\step[5]\object{F(Z)}\\
\nw3\step[4]\d\Step\dd\\
\step[3]\nw2\step[2]\tu \mu\\
\step[5]\tu \mu\\
\step[0.5]\obox 9{F(C_{X,Y\otimes Z})}\\
\step[6]\id\\
\step[4]\obox 4{\mu^{-1}}\\
\step[5]\xx\\
\step[4]\ne2\obox 4{\mu^{-1}}\\
\step[2]\dd\Step\dd\Step\d\\
\object{F(X)}\step[5]\object{F(Y)}\step[5]\object{F(Z)}\\
\end{tangle} \ \ \
\stackrel { \hbox { by Hexagon Axion } } {= }\step
\begin{tangle}
\object{F(X)}\step[5]\object{F(Y)}\step[5]\object{F(Z)}\\
\nw3\step[4]\d\Step\dd\\
\step[3]\nw2\step[2]\tu \mu\\
\step[5]\tu \mu\\
\step[0.5]\obox 9{F(C_{X,Z\otimes id_Z})}\\
\step[6]\id\\
\step[0.5]\obox 9{F(id_Y\otimes C_{X,Z})}\\
\step[6]\id\\
\step[4]\obox 4{\mu^{-1}}\\
\step[5]\xx\\
\step[4]\ne2\obox 4{\mu^{-1}}\\
\step[2]\dd\Step\dd\Step\d\\
\object{F(X)}\step[5]\object{F(Y)}\step[5]\object{F(Z)}\\
\end{tangle}
\]
\[
\step=\step
\begin{tangle}
\step\object{F(X)}\step[5]\object{F(Y)}\step[5]\object{F(Z)}\\
\step\nw3\step[4]\d\Step\dd\\
\step[4]\nw2\step[2]\tu \mu\\
\step[6]\tu \mu\\
\step[4]\obox 4{\mu^{-1}}\\
\step[4]\dd\step[2]\nw2\\
\obox 8{F(C_{X,Y})}\step\id\\
\step[4]\nw2\step[3]\dd\\
\step[6]\tu \mu\\
\step[5]\obox 4{\mu^{-1}}\\
\step[5]\dd\step[2]\nw2\\
\step[5]\id\step\obox 8{F(C_{X,Z})}\\
\step[5]\d\step[3]\dd\\
\step[6]\id\Step\dd\\
\step[6]\tu \mu\\
\step[5]\obox 4{\mu^{-1}}\\
\step[6]\xx\\
\step[5]\ne2\obox 4{\mu^{-1}}\\
\step[3]\dd\Step\dd\Step\d\\
\step\object{F(X)}\step[5]\object{F(Y)}\step[5]\object{F(Z)}\\
\end{tangle}
\step=\step
\begin{tangle}
\step\object{F(X)}\step[5]\object{F(Y)}\step[5]\object{F(Z)}\\
\step\nw2\step[3]\dd\step[5]\id\\
\step[3]\tu \mu\step[6]\id\\
\obox 8{F(C_{X,Y})}\step[3]\id\\
\step[4]\id\step[6]\dd\\
\Step\obox 4{\mu^{-1}}\step[3]\dd\\
\step[2]\dd\Step\d\Step\dd\\
\step[2]\id\step[4]\tu \mu\\
\step[2]\id\step\obox 8{F(C_{X,Z})}\\
\step[2]\id\step[5]\id\\
\step[2]\d\Step\obox 4{\mu^{-1}}\\
\step[3]\d\Step\xx\\
\step[4]\xx\Step\d\\
\step[3]\ne2\Step\id\step[3]\d\\
\object{F(X)}\step[5]\object{F(Y)}\step[5]\object{F(Z)} \ \ \ \ \ \ .\\
\end{tangle}\ \ \ .
\]

Thus
\[
\theta_I^{(3)}(
\begin{tangle}
\Ro R\\
\id\step[4]\d\\
\id\Step\ro R\step\id\\
\cu\Step\id\step\id
\end{tangle})_{X\otimes Y\otimes Z}
\step=\step \theta_I^{(3)}(
\begin{tangle}
\ro R\\
\id\step\cd\\
\end{tangle})_{X\otimes Y\otimes Z}
\step[3] and \step[3]
\begin{tangle}
\ro R\\
\id\step\cd\\
\object{H}\step[1.5]\object{H}\step[1.5]\object{H}\\
\end{tangle}
\step=\step
\begin{tangle}
\Ro R\\
\id\step[4]\d\\
\id\Step\ro R\step\id\\
\cu\Step\id\step\id\\
\step\object{H}\step[3]\object{H}\step[1.5]\object{H}\\
\end{tangle}\ \ \ .
\]

Similarly, we can show that (QT2)and (QT3) hold. Therefore,
$(B,R,\bar{\Delta})$ is a quasitriangular bialgebra. $\Box$

In the subsection,
 the diagram
 \[
 \begin{tangle}
 \object{U}\Step\object{V}\\
 \XX \\
 \object{V}\Step\object{U}
 \end{tangle}
 \]
 always denotes the ordinary twisted map:
 $x\otimes y \longrightarrow y\otimes x$.

\begin {Theorem} \label {12.1.3}
Let $H$  be an ordinary bialgebra and $(H_1, R)$   an ordinary
quasitriangular Hopf algebra over field $k$. Let  $f $  be
 a  bialgebra homomorphism from $H_1 $  to $H$.
 Then there exists  a   bialgebra
 $B$, written as  $B(H_1 , f, H)$,  living in $({}_{H_1} {\cal M}, C^R).$
  Here
  $B(H_1, f, H) = H$ as algebra,
  its counit is $\epsilon _H$,  and its comultiplication and antipode
are

\[
\begin{tangle}
\step\object{B}\\
\cd\\
\object{H}\Step\object{H}\\
\end{tangle}
\step=\step
\begin{tangle}
\step[3]\object{B}\\
\step\Cd\\
\dd\step\ro R\step\d\\
\id\Step\O f\step\obox 2{fS}\step\id\\
\id\Step\id\Step\id\Step\id\\
\id\Step\XX\step\dd\\
\cu\step[1.5]\obox 2{ad}\\
\step\id\step[3.5]\id\\
\step\object{H}\step[3]\object{H}\\
\end{tangle}\step[2] and \step[3]
\begin{tangle}
\step\object{B}\\
\step\id\\
\obox 2{S_B}\\
\step\id\\
\step\object{B}\\
\end{tangle}
\step=\step
\begin{tangle}
\step[3]\object{B}\\
\ro R\step\id\\
\O f\Step\O f\step\id\\
\XX\dd\\
\id\step\obox 2{ad}\\
\id\Step\id\\
\id\Step\O S\\
\cu\\
\step\object{B}\\
\end{tangle}
\step[2] respectively.
\]
Then $B$  is a braided Hopf algebra if $H$
  is a finite-dimensional
  Hopf algebra. Furthermore, if $H$  is an ordinary quasitriangular bialgebra, then
  $B$ is a braided quasitriangular bialgebra.
  In particular, when $H=H_1$ and $f= id _H$,  $B(H_1, f, H)$ is  a
  braided group, called
  the braided group analogue of $H$ and written as $\underline H$.

\end {Theorem}
{\bf Proof.} (i) Set ${\cal C} = {}_{H}{\cal M}$  and ${\cal D} =
({}_{H_1} {\cal M}, C^R ) $. Let $F$  be the  functor by pull-back
along $f$. That is, for any $(X, \alpha _X ) \in {}_H {\cal M}$, we
obtain an $H_1$-module $(X, \alpha _X')$ with $\alpha _X'  = \alpha
_X (f \otimes id _X),$ written as $(X, \alpha _X') = F(X).$ For any
morphism $g \in Hom _{\cal C} (U, V)$ , define $F(g) = g$. $B$ is a
left  $B$-module by adjoint action.
 Let $_BB$ denote the left regular $B$-module. It is clear that $\alpha
 $ is a natural transformation from $B\otimes F$ to $F$.
We first show that  $\theta _V$  is  injective for any
 $V\in {}_{H_1} {\cal M }$.
If $ \theta _V (g) = \theta _V (h)$, i.e.

\[
\begin{tangle}
\step[0.5]\object{V}\step[4]\object{F(_BB)}\\
\step[0.5]\O g\step[3]\dd\\
\step[0.5]\d\step\dd\\
\obox 4{\alpha_{_BB}}\\
\step[2]\id\\
\step[2]\object{F(_BB)}\\
\end{tangle}
\step[3]=\step[3]
\begin{tangle}
\step[0.5]\object{V}\step[4]\object{F(_BB)}\\
\step[0.5]\O h\step[3]\dd\\
\step[0.5]\d\step\dd\\
\obox 4{\alpha_{_BB}}\\
\step[2]\id\\
\step[2]\object{F(_BB)}\\
\end{tangle}
\step[2],
 \hbox {which implies that} \step[2]
\begin{tangle}
\object{V}\\
\id\\
\O g\\
\id\\
\object{B}\\
\end{tangle}
\step=\step
\begin{tangle}
\object{V}\\
\id\\
\O h\\
\id\\
\object{B}\\
\end{tangle},
\]
where $g$  and $h$  are $H_1$-module homomorphisms from
  $V$ to $B.$  Thus $\theta _V$  is injective.
Similarly, we can show that $\theta _V^{(2)}$ and $\theta _V^{(3)}$
are injective.

By Proposition \ref {12.1.2}, $B$ is a quasitriangular bialgebra
living in $({}_{H_1} {\cal M}, C^R)$ determined by braided
reconstruction.

We first show that  the comultiplication of $B$   is the same as
stated. Now we see that

\[
\begin{tangle}
\theta_B^{(2)}(\Delta_B)_{X\otimes Y}
\end{tangle}
\step=\step
\begin{tangle}
\object{B}\step[4]\object{F(X)}\step[6]\object{F(Y)}\\
\id\step[4]\d\step[3]\ne2\\
\d\step[4]\tu \mu\\
\step\d\step[3.5]\ne2\\
\step\obox 5{\alpha_{X\otimes Y}}\\
\step[3.5]\id\\
\step[2]\obox 4{\mu^{-1}}\\
\step[2]\dd\step\d\\
\step[1.5]\object{F(X)}\step[6]\object{F(Y)}
\end{tangle}
\step[5]......(1)
\]

\[\hbox {The left side of (1)}
\step=\step
\begin{tangle}
\step[1.5]\object{B}\step[4]\object{F(X)}\step[5]\object{F(Y)}\\
\step[1.5]\id\step[3.5]\id\step[5]\id\\
\obox 3{\Delta_B}\Step\id\step[5]\id\\
\dd\step\d\Step\id\step[4]\dd\\
\id\step[3]\ox {C^R} \step[3]\dd\\

\d\step\dd\Step\d\step\dd\\
\obox 3{\alpha_X}\Step\obox 3{\alpha_Y}\\
\step[1.5]\id\step[5]\id\\
\step\object{F(X)}\step[6]\object{F(X)}\\
\end{tangle}
\step=\step
\begin{tangle}
\step[3]\object{B}\step[7]\object{F(X)}\step[5]\object{F(Y)}\\
\step[3]\id\step[7]\id\step[5]\id\\
\step[1.5]\obox 3{\Delta_H}\step[5.5]\id\step[5]\id\\
\step\dd\ro R\d\step[5]\id\step[5]\id\\
\dd\step\O f\Step\O S\step\d\step[4]\id\step[5]\id\\
\id\Step\id\Step\O f\step[2]\d\step[3]\id\step[5]\id\\
\id\Step\XX\step[3]\d\step[2]\id\step[5]\id\\
\cu\ro R\d\step[3]\id\step[2]\id\step[5]\id\\
\dd\step\O f\Step\O f\step\d\step[2]\id\step[2]\id\step[5]\id\\
\id\step[2]\XX\Step\id\Step\id\step[2]\id\step[5]\id\\
\cu\Step\cu\Step\id\step[2]\id\step[5]\id\\
\step\d\step[3]\d\step\dd\step[2]\id\step[5]\id\\
\step[2]\d\step[2]\obox 3{ad}\step[2]\id\step[5]\id\\
\step[3]\d\step[3]\d\step[2]\id\step[4]\dd\\

\step[4]\d\step[3]\XX \step[3]\dd\\
\step[5]\d\step\dd\Step\d\step\dd\\
\step[5]\obox 3{\alpha_X}\step[2]\obox 3{\alpha_Y}\\
\step[6.5]\id\step[5]\id\\
\step[6]\object{F(X)}\step[5]\object{F(Y)}\\
\end{tangle}
\]\[
\step=\step
\begin{tangle}
\step[3]\object{B}\step[5]\object{F(X)}\step[5]\object{F(Y)}\\
\step[3]\id\step[5]\id\step[5]\id\\
\step[1.5]\obox 3{\Delta_H}\step[3.5]\id\step[5]\id\\
\step\dd\ro R\d\step[3]\id\step[5]\id\\
\dd\step\O f\step\cd\d\Step\id\step[5]\id\\
\id\Step\id\step\O S\Step\id\step\id\Step\id\step[5]\id\\
\id\Step\id\step\cu\step\id\Step\id\step[5]\id\\
\id\Step\id\Step\O f\Step\id\Step\id\step[5]\id\\
\id\Step\XX\step\dd\Step\id\step[4]\dd\\
\cu\step\obox 3{ad}\step[2]\id\step[3]\dd\\
\step\d\step[3]\d\Step\id\step[2]\dd\\

\step[2]\d\step[3]\XX\step\dd\\
\step[3]\d\step\dd\step\obox 3{\alpha_Y}\\
\step[3]\obox 3{\alpha_X}\step[2.5]\id\\
\step[4.5]\id\step[4]\id\\
\step[4]\object{F(X)}\step[5]\object{F(Y)}\\
\end{tangle}
\step=\step
\begin{tangle}
\step[1.5]\object{B}\step[4]\object{F(X)}\step[5]\object{F(Y)}\\
\step[1.5]\id\step[3.5]\id\step[5]\id\\
\obox 3{\Delta_H}\Step\id\step[5]\id\\
\dd\step\d\Step\id\step[4]\dd\\
\id\step[3]\XX\step[3]\dd\\
\d\step\dd\Step\d\step\dd\\
\obox 3{\alpha_X}\Step\obox 3{\alpha_Y}\\
\step[1.5]\id\step[5]\id\\
\step\object{F(X)}\step[6]\object{F(Y)}\\
\end{tangle}
\step=\step\hbox {the right side of (1).}
\]
Thus relation (1) holds.\\

Furthermore,  if $H$  is a finite-dimensional ordinary Hopf algebra,
set       ${\cal C}  = \{ M \in  {}_{H}{\cal M} \mid M \hbox {is
finite -dimensional }\}$ and the others are the same as part (i). It
is clear that ${\cal C}$ is rigid. Thus $B$  has a braided-antipode
by Proposition \ref {12.1.2} (ii).

We see that
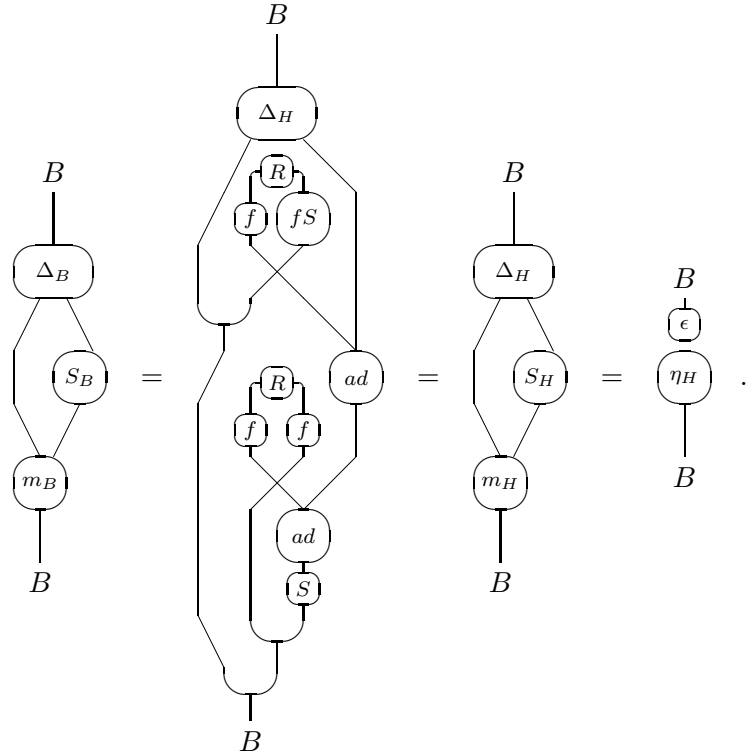
\[
\begin{tangle}
\step[1.5]\object{B}\\
\step[1.5]\id\\
\obox 3{\Delta_B}\\
\dd\step\d\\
\id\step[1.5]\obox 2{S_B}\\
\d\step[0.5]\dd\\
\obox 2{m_B}\\
\step\id\\
\step\object{B}\\
\end{tangle}
\step=\step
\begin{tangle}
\step[3]\object{B}\\
\step[3]\id\\
\step[1.5]\obox 3{\Delta_H}\\
\step\dd\ro R\nw2\\
\dd\step\O f\step\obox 2{fS}\step\id\\
\id\Step\XX\step[2]\id\\
\cu\step[2]\nw2\step\id\\
\dd\step\ro R\step\obox 2{ad}\\
\id\Step\O f\step[2]\O f\Step\id\\
\id\Step\XX\step\ne2\\
\id\Step\id\step\obox 2{ad}\\
\id\Step\id\Step\O S\\
\d\step\cu\\
\step\cu\\
\step[2]\object{B}\\
\end{tangle}
\step=\step
\begin{tangle}
\step[1.5]\object{B}\\
\step[1.5]\id\\
\obox 3{\Delta_H}\\
\dd\step\d\\
\id\step[1.5]\obox 2{S_H}\\
\d\step[0.5]\dd\\
\obox 2{m_H}\\
\step\id\\
\step\object{B}\\
\end{tangle}
\step=\step
\begin{tangle}
\step\object{B}\\
\step\QQ \epsilon \\
\obox 2{\eta_H}\\
\step\id\\
\step\object{B}\\
\end{tangle}\ \ \ .
\]

Thus the antipode of $B$  as stated.

Assume  $H_1 = H$. Our functor $F$  is  the identity Functor form
${}_H {\cal M}$ to ${}_H {\cal M}.$  It follows from the proof of
Proposition \ref {12.1.2} (iii) that the quasitriangular structure
of $\underline H$  is trivial and $\bar \Delta = \Delta $. That is,
$\underline H$ is a braided group. $\Box$

Dually we have

\begin {Theorem} \label {12.1.3}
Let $H$  be an ordinary bialgebra and $(H_1, r)$  be an ordinary
coquasitriangular Hopf algebra over field $k$. Let  $f $  be
 a Hopf algebra homomorphism from $H $  to $H_1$.
 Then there exists a braided bialgebra
 $B$, written  as $B(H, f, H_1)$,  living in $({}^{B} {\cal M}, C^r).$
  Here
  $B(H_1, f, H) = H$ as coalgebra,
   its unit is $\eta  _H$,  and its multiplication and antipode are
  respectively,
 
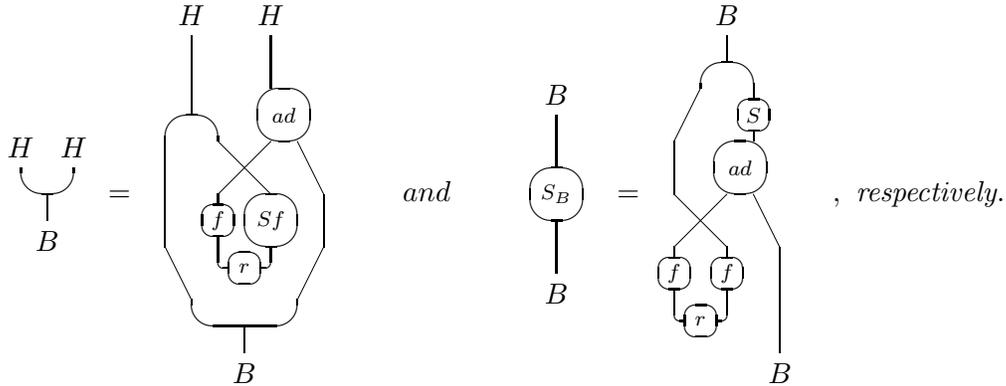
\[
\begin{tangle}
\object{H}\Step\object{H}\\
\cu\\
\step\object{B}\\
\end{tangle}
\step=\step
\begin{tangle}
\step\object{H}\step[3]\object{H}\\
\step\id\step[3]\id\\
\cd\step[1.5]\obox 2{ad}\\
\id\Step\XX\step\d\\
\id\Step\O f\step\obox 2{Sf}\step\id\\
\d\step\coro r\step\dd\\
\step\Cu\\
\step[3]\object{B}\\
\end{tangle}
\step[3] and\step[3]
\begin{tangle}
\step\object{B}\\
\step\id\\
\obox 2{S_B}\\
\step\id\\
\step\object{B}\\
\end{tangle}
\step=\step
\begin{tangle}
\step[2]\object{B}\\
\step\cd\\
\dd\Step\O S\\
\id\step[1.5]\obox 2{ad}\\
\XX\step\d\\
\O f\Step\O f\Step\id\\
\coro r\Step\id\\
\step[4]\object{B}\\
\end{tangle}
\step[2] ,\ \hbox {respectively}.
\]
  \noindent Then  $B$ is a braided Hopf algebra
  if $H$ is a finite-dimensional
  Hopf algebra. Furthermore, if $H$  is an ordinary coquasitriangular
   bialgebra, then
  $B$ is a braided coquasitriangular bialgebra.
  In particular, when $H=H_1$ and $f= id _H$,  $B(H, f, H_1)$ is
  a braided group, called
  the braided group analogue of $H$ and written as $\underline H$.

\end {Theorem}

\section {Braided group analogues of double cross products }\label {s9}
In this section we show that
 the braided group analogue of double cross
(co)product is double cross (co)product of braided group analogues.

 \begin {Theorem} \label {4.1.3}  (Factorisation theorem) (See
\cite [Theorem 7.2.3]{Ma95b}) Let $X$ , $A$ and  $H$  be bialgebras
or  Hopf  algebras. Assume that $j_A$ and $j_H$  are bialgebra or
Hopf algebra morphisms from $A$ to $X$ and $H$  to $X$ respectively.
If $   \xi =:  m_X (j_A \otimes j_H)$ is an isomorphism  from $ A
\otimes H$ onto $X$ as objects in ${\cal C}$, then there exist
morphisms
 $$ \alpha  : H\otimes  A \rightarrow A \hbox { \ \ \ and \ \ \ }
 \beta  :  H \otimes A \rightarrow H $$
 such that $ A {}_ \alpha \bowtie _\beta H $  becomes  a bialgebra or
 Hopf algebra
  and  $\xi $   is a bialgebra or Hopf algebra isomorphism from $ A {}_\alpha \bowtie _\beta H$  onto
  $X.$
\end {Theorem}
\textbf{Proof}. Set
\[  \zeta  =:
\begin{tangle}
\step\object{H}\Step\object{A}\\
 \step \O {j_H} \step[2] \O {j_A} \\
\step \cu \\
\step \td {\bar \xi}\\
\step \object{A}\step\step\object{H}
 \end{tangle}
 \Step,\Step \alpha =:
 \begin{tangle}
\object{H}\step\step\object{A}\\
 \ox \zeta \\
 \id \step[2] \QQ \epsilon \\
 \object{A}
 \end{tangle}\step \hbox { \ and } \step
 \beta =:
 \begin{tangle}
\step \object{H}\step\step\object{A}\\
 \step \ox \zeta \\
\step \QQ \epsilon \step[2] \id  \\
 \step [3]\object{H} \ \ \
 \end{tangle}
 \ \ .\]
We see
\[
\begin{tangle}
\object{H}\Step\object{H}\step\step\object{A}\\
\id\step\step\ox \zeta \\
        \ox \zeta\step\step\id \\
              \O {j_{A}}\step\step\cu \\
       \id\step\step\step\O {j_{H}}\\
       \id\step\step\dd\\
       \cu\\
       \step\object{X}\end{tangle}
\;=\enspace
\begin{tangle}
\object{H}\Step\object{H}\step\step\object{A}\\
\id\step\step\ox \zeta \\
        \ox \zeta \step\step\id \\
        \O{j_{A}}\Step\O{j_{H}}\Step\O{j_{H}}\\
       \id\Step\cu \\
       \id\Step\dd \\
       \cu \\
       \step\object{X}\end{tangle}
\;=\enspace
\begin{tangle}
\object{H}\Step\object{H}\step\step\object{A}\\
\id\step\step\ox \zeta \\
        \ox\zeta \step\step\id \\
        \O{j_{A}}\Step\O{j_{H}}\Step\O{j_{H}}\\
       \cu\Step\id \\
       \step\d\Step\id \\
       \Step\cu \\
       \Step\step\object{X}\end{tangle}
       \;=\enspace
\begin{tangle}
\object{H}\Step\object{H}\step\step\object{A}\\
\id\step\step\ox \zeta \\
  \O{j_{H}}\Step\O{j_{A}}\Step\O{j_{H}}\\
       \cu\Step\id \\
       \step\d\Step\id \\
       \Step\cu\\
       \Step\step\object{X}\end{tangle}
  \;=\enspace
\begin{tangle}
\object{H}\Step\object{H}\step\step\object{A}\\
\id\step\step\id\Step\id\\
  \O{j_{H}}\Step\O{j_{H}}\Step\O{j_{A}}\\
       \cu\Step\id \\
       \step\d\Step\id \\
       \Step\cu \\
       \Step\step\object{X}\end{tangle}
\ =\
\begin{tangle}
\object{H}\Step\object{H}\step\step\object{A}\\
\cu\Step\id\\
\step\id\Step\dd\\
\step \ox \zeta \\
\step\O{j_{A}}\Step\O{j_{H}}\\
\step\cu \\
\step\step\object{X} \ \ \end{tangle} \ . \] Thus
\[
\begin{tangle}
\object{H}\Step\object{H}\step\step\object{A}\\
\id\step\step\ox \zeta \\
\ox \zeta \Step\id\\
\id\Step\cu\\
\object{A}\Step\step\object{H}
\end{tangle}
\;=\enspace
\begin{tangle}
\object{H}\Step\object{H}\step\step\object{A}\\
\cu \Step\id\\
\step\id\Step\dd\\
\step \ox \zeta \\
\step\object{A}\Step\object{H}\end{tangle}\ \ . \Step\Step \cdots
\cdots  (1)
\]
Similarly we have
\[
\begin{tangle}
\object{H}\Step\object{A}\step\step\object{A}\\
\ox \zeta \Step\id\\
\id\Step\ox \zeta \\
\cu \Step\id\\
\step\object{A}\Step\object{H}
\end{tangle}
\;=\enspace
\begin{tangle}
\object{H}\Step\object{A}\step\step\object{A}\\
\id\Step\cu\\
\d\Step\id\\
\step \ox \zeta \\
\step\object{A}\Step\object{H}\end{tangle} \ \ . \Step\Step \cdots
\cdots  (2)
\]
We also have
\begin {eqnarray*} \zeta (\eta \otimes id ) = id \otimes \eta
\hbox { \ and \ } \zeta (id \otimes \eta  ) = \eta\otimes id \  . \
\ \ \ \ \ \ \ \cdots \cdots (3)
\end {eqnarray*}

 It is clear that $\zeta$ is a coalgebra morphism from
$H\otimes A$ to $A\otimes H$, since $j_{A}$,  $ j_{H}$ and $m_{X}$
all are coalgebra homorphisms. Thus we have
\[
\begin{tangle}
\step\object{H}\Step\Step\object{A}\\
\cd \Step\cd \\
\id\Step\x\Step\id\\
\ox \zeta \Step\ox \zeta \\
\object{A}\Step\object{H}\Step\object{A}\Step\object{H}
\end{tangle}
\;=\enspace
\begin{tangle}
\Step\object{H}\Step\object{A}\\
\Step\ox \zeta \\
\step\dd\Step\d\\
\cd\Step\cd \\
\id\Step\x\Step\id \\
\object{A}\Step\object{H}\Step\object{A}\Step\object{H}
\end{tangle}
\step \hbox { and \ } (\epsilon \otimes \epsilon )\zeta = (\epsilon
\otimes \epsilon ) . \ \ \ \ \ \ \ \ \ \ \cdots \cdots (4)
\]
We now show that $(A,\alpha)$ is an $H$-module coalgebra:\\
\[
\begin{tangle}
\object{H}\Step\object{H}\step\object{A}\\
\cu \step\id\\
\step\tu \alpha\\
\Step\object{A}
\end{tangle} \ \
= \begin{tangle}
\object{H}\Step\object{H}\step\step\object{A}\\
\cu \Step\id\\
\step\id\Step\dd\\
\step \ox \zeta \\
\step\object{A}\Step\object{\HH\obox 1\epsilon}
\end{tangle}
 \Step \stackrel { \hbox{ by } (1)}{= } \Step
\begin{tangle}
\object{H}\Step\object{H}\step\step\object{A}\\
\id\step\step\ox \zeta \\
\ox \zeta \Step\id\\
\id\Step\cu\\
\object{A}\Step\step\object{\HH\obox 1\epsilon}
\end{tangle}
 \;=\enspace
\begin{tangle}
\object{H}\step\object{H}\Step\object{A}\\
\id\step\tu \alpha\\
\tu \alpha\\
\step\object{A}
\end{tangle}
\]
and \ \  $\alpha (\eta \otimes id_A ) = (id_A \otimes \epsilon )
\zeta (\eta \otimes id_A ) \stackrel {\hbox {by} (3)}{=} id_A .$

We see that  $\epsilon \circ \alpha  = (\epsilon \otimes \epsilon
)\zeta \stackrel {\hbox {by} (4)}{=}\epsilon \otimes \epsilon $
 \ \  and \[
\begin{tangle}
\step\object{H}\Step\Step\object{A}\\
\cd \Step\cd\\
\id\Step\x\Step\id\\
\cu\Step\cu\\
\step\object{A}\Step\Step\object{A}
 \end{tangle}
 \;=\enspace
\begin{tangle}
\step\object{H}\Step\Step\object{A}\\
\cd \Step\cd \\
\id\Step\x\Step\id\\
\ox \zeta \Step\ox \zeta \\
\object{A}\Step \object{\HH\obox
1\epsilon}\Step\object{A}\Step\object{\HH\obox 1\epsilon}\Step
\end{tangle}
\ \ \stackrel {\hbox {by} (4)}{= } \ \
\begin{tangle}
\Step\object{H}\Step\object{A}\\
\Step\ox \zeta \\
\step\dd\Step\d\\
\cd \Step\cd \\
\id\Step\ox \zeta \Step\id\\
\object{A}\Step\object{\HH\obox
1\epsilon}\Step\object{A}\Step\object{\HH\obox 1\epsilon}
\end{tangle}
\;=\enspace
\begin{tangle}
\object{H}\Step\object{A}\\
\cu\\
\cd\\
\object{A}\Step\object{A}
\end{tangle} \ \ .
\]Thus $(A,\alpha)$ is an $H$-module
coalgebra. Similarly, we can show that $(H,\beta)$ is an $A$-module
coalgebra.

  Now we show that conditions $(M1)$--$(M4)$ in [12,p37] hold. By (3), we easily
know that$(M1)$ holds. Next we show that$(M2)$ holds.
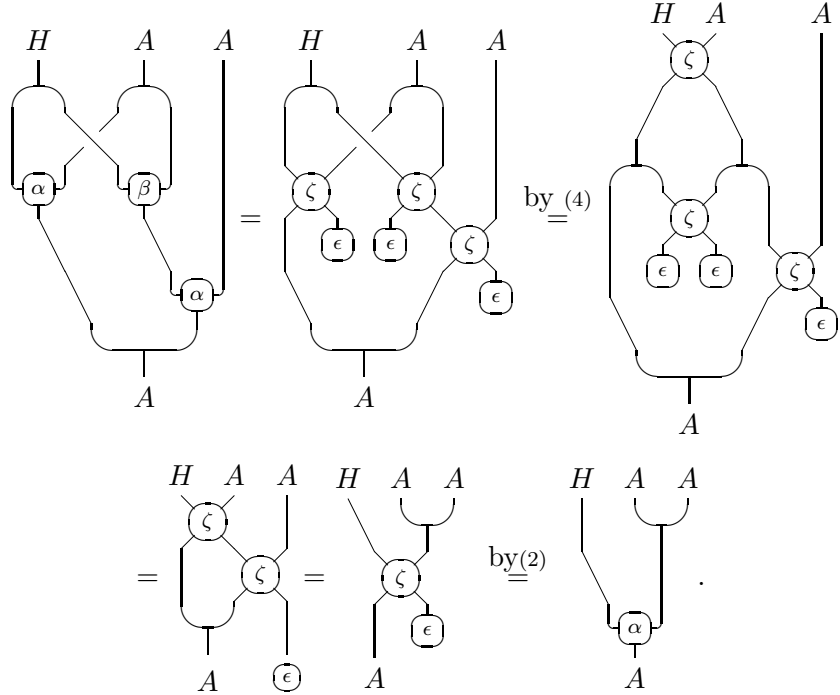
\[
\begin{tangle}
\step\object{H}\Step\Step\object{A}\Step\step\object{A}\\
\cd \Step\cd \Step\id\\
\id\Step\x\Step\id\Step\id\\
\tu \alpha \Step\tu \beta \Step\id\\
\step\d\Step \step \d\Step\id\\
\step\step\d\step\step\step\tu \alpha\\\
\step\step\step\Cu\\
\Step\Step\step\object{A}
\end{tangle}
\;=\enspace
\begin{tangle}
\step\object{H}\Step\Step\object{A}\Step\step\object{A}\\
\cd \Step\cd \Step\id\\
\id\Step\x\Step\id\Step\id\\
\ox \zeta \Step\ox \zeta \Step\id\\
\id\Step\QQ \epsilon \Step\QQ \epsilon \Step\ox \zeta  \\
\d\Step\Step\dd \Step \QQ \epsilon \\
\step\Cu\\
\Step\step\object{A}\Step
\end{tangle}
\ \ \stackrel {\hbox {by } (4)} {= } \
\begin{tangle}
\Step\object{H}\Step\object{A}\Step\Step\object{A}\\
\Step\ox \zeta \Step\Step\id\\
\step\dd\Step\d\Step\step\id\\
\cd \Step\cd \Step\id\\
\id\Step\ox \zeta \Step\id\Step\id\\
\id\Step\QQ \epsilon\Step\QQ \epsilon\Step\ox \zeta \\
\d\Step\Step\dd\Step\QQ \epsilon\\
\step\Cu\\
\Step\step\object{A}
\end{tangle}
\]
\[
\;=\enspace
\begin{tangle}
\object{H}\Step\object{A}\step\step\object{A}\\
\ox \zeta \Step\id\\
\id\Step\ox \zeta \\
\cu \Step\id\\
\step\object{A}\Step\step\object{\HH\obox 1\epsilon}
\end{tangle}
\;=\enspace
\begin{tangle}
\object{H}\Step\object{A}\step\step\object{A}\\
\nw1 \step\cu\\
\step \ox \zeta \\
\step \id \step [2] \QQ \epsilon \\
 \step\object{A}\end{tangle} \step \stackrel
{\hbox {by} (2)}{=} \step
\begin{tangle}
\object{H}\Step\object{A}\step\step\object{A}\\
\id\Step\cu\\
\d\Step\id\\
\step\tu \alpha\\
\Step\object{A}
\end{tangle} \ .
\]
 Thus $(M2)$ holds. Similarly, we can get the proofs of
$(M3)$ and$(M4)$. Consequently,   $ A {}_\alpha \bowtie _\beta H$ is
a bialgebra or Hopf algebra by [12,  Corollary 1.8]. It suffices to
show that $\zeta$ is a bialgebra morphism from$ A {} _\alpha \bowtie
_\beta H$ to $X$. Let $ D=A {}_\alpha \bowtie _\beta H$. Since
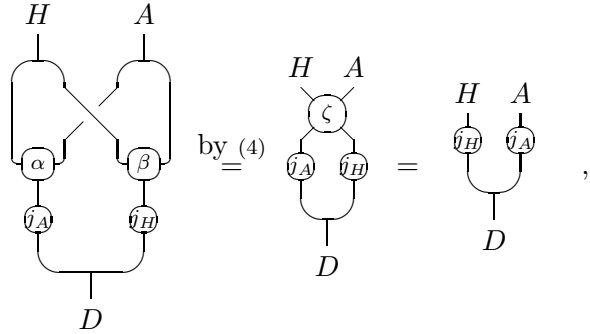
\[
\begin{tangle}
\step\object{H}\Step\Step\object{A}\\
\cd\Step\cd\\
\id\Step\x\Step\id\\
\tu \alpha \Step\tu \beta\\
\morph {j_{A}} \Step\morph {j_{H}}\\
\step\Cu\\
 \step\Step\object{D}
\end{tangle}
\ \ \stackrel {\hbox {by } (4)} {= }
\begin{tangle}
\step\object{H}\Step\object{A}\\
\step\ox \zeta \\
\morph {j_{A}} \morph {j_{H}}\\
\step\cu\\
\Step\object{D}
\end{tangle}
\;=\enspace
\begin{tangle}
\step\object{H}\Step\object{A}\\
\morph {j_{H}} \morph {j_{A}}\\
\step\cu\\
\Step\object{D}
\end{tangle}
\step\;,  \enspace
\]
 we have that $\xi$ is a bialgebra morphism from $ A{} _\alpha \bowtie _\beta H$ to $X$ by \cite
 [Lemma 2.5] {ZC99}. \begin{picture}(5,5)
\put(0,0){\line(0,1){5}}\put(5,5){\line(0,-1){5}}
\put(0,0){\line(1,0){5}}\put(5,5){\line(-1,0){5}}
\end{picture}\\

\begin {Theorem} \label {4.1.4}
(Co-factorisation theorem) Let $X$ ,   $A$ and  $H$  be bialgebras
or  Hopf  algebras. Assume that $p_A$ and $p_H$  are bialgebra or
Hopf algebra morphisms from $X$ to $A$ and $X$  to $H$,
respectively. If  $ \xi =   (p_A \otimes p_H)\Delta _X$ is an
isomorphism  from $X$ onto $ A \otimes H$ as objects in ${\cal C}$,
then there exist morphisms:
 $$ \phi  :  A \rightarrow H\otimes A \hbox { \ \ \ and \ \ \ }
 \psi  :  H \rightarrow H \otimes A$$
 such that $ A {}^\phi \bowtie ^\psi H $  becomes  a bialgebra or Hopf
 algebra   and  $\xi $   is a bialgebra or Hopf algebra isomorphism from
  $X$ to $ A {}^\phi \bowtie ^\psi H$.
  \end {Theorem}

  {\bf Proof.}  Set \\
\[
\begin{tangle}
\object{A}\Step\object{H}\\
\ox \zeta \\
\object{H}\Step\object{A}
\end{tangle}
\;=\enspace
\begin{tangle}
\step\object{A}\Step\object{H}\\

\step\tu {\overline{\xi}}\\
\step\td {\Delta_{X}} \\
\morph {p_{H}} \morph {p_{A}}\\
\step\object{H}\Step\object{A}
\end{tangle}
\Step\Step
%%%%%%%%%%%%%%%%%%%%%%%%%%%%%%%%%%%%%%%%%%%%%%%%
\begin{tangle}
\step\object{A}\\
\td \phi \\
\id\step\step\id\\
\object{H}\Step\object{A}
\end{tangle}
\;=\enspace
\begin{tangle}
\object{A}\Step\object{\HH\obox 1{\eta_{H}}}\\
\ox \zeta \\
\object{H}\Step\object{A}\\
\end{tangle}
%%%%%%%%%%%%%%%%%%%%%%%%%%%%%
\Step\Step
\begin{tangle}
\step\object{H}\\
\td \psi \\
\id\step\step\id\\
\object{H}\Step\object{A}
\end{tangle}
\;=\enspace
\begin{tangle}
\object{\HH\obox 1{\eta_{A}}}\Step\object{H}\\
\ox \zeta \\
\object{H}\Step\object{A}\\
\end{tangle}\ \ \ .
\] We can complete the proof by turning upside down the diagrams in
 the proof of the preceding theorem. \begin{picture}(5,5)
\put(0,0){\line(0,1){5}}\put(5,5){\line(0,-1){5}}
\put(0,0){\line(1,0){5}}\put(5,5){\line(-1,0){5}}
\end{picture}\\

 In the subsection,   we always consider Hopf algebras over field $k$ and
 the diagram
 \[
 \begin{tangle}
 \object{U}\Step\object{V}\\
 \XX \\
 \object{V}\Step\object{U}
 \end{tangle}
 \]
 always denotes the ordinary twisted map:
 $x\otimes y \longrightarrow y\otimes x$.
 Our diagrams only denote homomorphisms between vector spaces,   so
two diagrams can have the additive operation.

  We now investigate the relation among braided group analogues of
  quasitriangular Hopf algebras $A$ and $H$ and their double cross
  coproduct $D = A\bowtie^{R}H$. Let us recall transmutation. We
  denote the braided group analogue of any ordinary quasitriangular Hopf
  algebras $H$ by $\underline{H}$. $R$ is called a weak $R$-matix of  $A\otimes H$ if $R$ is invertible
under convolution and
\[
\begin{tangle}
\step \ro R  \\
 \cd \step \id\\
\object{A}\Step\object{A}\step\object{H}\\
\end{tangle}
\;= \enspace
\begin{tangle}
\ro R  \Step \ro R  \\
\id \Step \XX \Step\id\\
\id \Step \id \Step \cu\\
\object{A}\Step\object{A}\Step\step\object{H}
\end{tangle} \hbox { \ \ and \ \ }
\begin{tangle}
\ro R  \\
\id \step \cd\\
\object{A}\step\object{H}\Step\object{H}
\end{tangle}
\;= \enspace
\begin{tangle}
\step \Ro R\\
\dd \step \ro R \step  \id\\
\cu \Step \id \step \id\\
\step\object{A}\Step\step\object{H}\step\object{H}
\end{tangle}\ \ \ .
\]

Let $(A,  P)$ and $(H,  Q)$ be ordinary finite-dimensional
quasitriangular Hopf algebras over field $k$. Let $R$ be a weak
$R$-matrix of $A\otimes H$. For any $U,  V \in CW(A\otimes H)=:$$\{U
\in A\otimes H \ | \ U$  is a weak $R$-matrix and in the center of $
A\otimes H\}$,
 $$R_{D}=:\sum R'P'U' \otimes Q'(R^{-1})''V'' \otimes P'' (R^{-1})' V'
 \otimes R''Q''U''$$
is a quasitriangular structure of $D$ and every quasitriangular
structure of $D$ is of this form (\cite [Theorem 2.9 ] {Ch98}),
where $R = \sum R'\otimes R''$,   etc.

\begin {Lemma} \label {1.3}
Under the above discussion,   then

(i) $\pi _A : D \rightarrow A$ and $\pi _H : D \rightarrow H$ are
bialgebra  or Hopf algebra homomorphisms,   respectively. Here $\pi
_A $ and $\pi _H$ are trivial action,   that is,   $\pi _A (h
\otimes a ) = \epsilon (h)a$ for any $a\in A,   h \in H.$

(ii) $B(D,   \pi _A,   A) = \underline A$ and $B(D,   \pi _H,   H) =
\underline H.$

(iii) $\pi _{\underline A} : \underline D \rightarrow \underline A$
and $\pi _{\underline H }: \underline D \rightarrow \underline H$
are bialgebra  or Hopf algebra homomorphisms,   respectively.

  \end {Lemma}
  {\bf Proof.} (i) It is clear.

  (ii) It is enough to show $\Delta _B = \Delta _{\underline
  {A}}$ since $B= \underline A$ as algebras,   where $B =: B(D,   \pi_ A,
  A).$  See\\
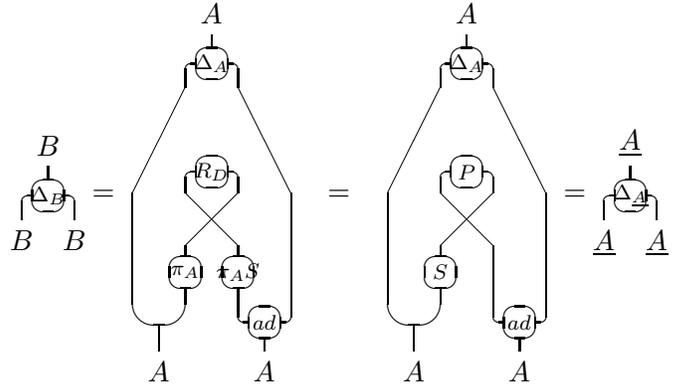
\[
\begin{tangle}
\step \object{B}\\
\td {\Delta_{B}} \\
\object{B}\Step\object{B}\\
\end{tangle}
\ = \
\begin{tangle}
\step\Step\object{A}\\

\step [2] \td {\Delta_{A}} \\

\step \ne1 \step [2] \nw1\\
 \dd \step \ro {R_{D}}\step \d\\
\id\step\step\XX  \step\step\id \\
\id \Step \O {\pi_{A}} \Step\O {\pi_{A}S}\Step\id\\
\cu \step  \step \tu {ad}\\
\step\object{A}\Step\Step\object{A}
 \end{tangle}
\ \ \ = \ \ \
\begin{tangle}
\step\Step\object{A}\\

\step [2] \td {\Delta_{A}} \\

\step \ne1 \step [2] \nw1\\
 \dd \step \ro {P} \step \d\\
\id\step\step\XX  \step\step\id \\
\id\Step \S\Step\id\Step\id\\
\cu \step  \step \tu {ad}\\
\step\object{A}\Step\Step\object{A}
 \end{tangle}
\ = \
 \begin{tangle}
\step \object{\underline{A}}\\
\td {\Delta_{\underline{A}}} \\
\object{\underline{A}}\Step\object{\underline{A}}\\
\end{tangle} \ \ .
\]
  Thus $\Delta _B = \Delta _{\underline A}$.
Similarly, we have $B(D,   \pi _H,   H)= \underline H.$

 (iii) See\\
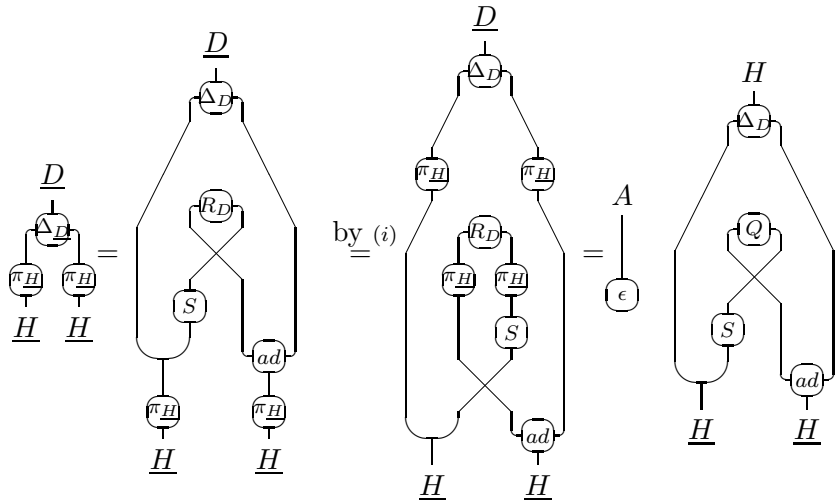
  \[
 \begin{tangle}
\step \object{\underline{D}}\\
\td {\Delta_{\underline{D}}} \\
\O {\pi_{\underline{H}}}\Step\O {\pi_{\underline{H}}}\\
\object{\underline{H}}\Step\object{\underline{H}}\\
\end{tangle}
\ = \
\begin{tangle}
\step\Step\object{\underline{D}}\\

\step [2] \td {\Delta_{D}} \\

\step \ne1 \step [2] \nw1\\

 \dd \step \ro {R_{D}}\step \d\\
\id\step\step\XX  \step\step\id \\
\id\Step \S\Step\id\Step\id\\
\cu \step  \step \tu {ad}\\
\step\O {\pi_{\underline{H}}}\Step\Step\O {\pi_{\underline{H}}}\\
\step\object{\underline{H}}\Step\Step\object{\underline{H}}
 \end{tangle}
\ \ \stackrel {\hbox { by }(i)}{ = }
\begin{tangle}
\Step \step\object{\underline{D}}\\

\step [2] \td {\Delta_{D}} \\

\step \ne1 \step [2] \nw1\\
\step\O {\pi_{\underline{H}}}\Step\step\step\O {\pi_{\underline{H}}}\\
\dd \step  \ro {R_{D}} \step \d\\
\id \Step \O {\pi_{\underline{H}}}\Step\O {\pi_{\underline{H}}}\Step\id\\
\id\Step\id\Step\S\Step\id\\
\id\step\step\XX  \step\step\id \\
\cu \step  \step \tu {ad}\\
\step\object{\underline{H}}\Step\Step\object{\underline{H}}
 \end{tangle}
   \ = \
\begin{tangle}
\object{A}\\
\id\\
\QQ \epsilon
 \end{tangle}
 \Step
 \begin{tangle}
\step\Step\object{H}\\

\step [2]\td {\Delta_{D}} \\

\step \ne1 \step [2] \nw1\\
 \dd \step \ro {Q}\step \d\\
\id\step\step\XX \step\step\id \\
\id\Step \S\Step\id\Step\id\\
\cu \step  \step \tu {ad}\\
\step\object{\underline{H}}\Step\Step\object{\underline{H}}
 \end{tangle}
\]
 and $\epsilon \circ \pi _{\underline H}= \epsilon.$ Thus $\pi
_{\underline H}
 $ is a coalgebra homomorphism.

 Since the multiplications in $\underline D$ and $\underline H$
 are the same as in $D$ and $H$, respectively,   we have that $\pi _{\underline H}
 $ is an algebra homomorphism by  (i). Similarly,   we can show
 that $\pi _{\underline A}$ is a bialgebra homomorphism.  \begin{picture}(5,5)
\put(0,0){\line(0,1){5}}\put(5,5){\line(0,-1){5}}
\put(0,0){\line(1,0){5}}\put(5,5){\line(-1,0){5}}
\end{picture}\\

 We now investigate the relation among
 braided group analogues of  quasitriangular Hopf algebras $A$  and $H$
and their double cross coproduct $ D= A  \bowtie ^R H.$

 \begin {Theorem}  \label {4.1.5}
  Under the above discussion,   let
  $\xi = (\pi _{\underline A}
  \otimes \pi _{\underline H}) \Delta _{\underline D}.$
 Then
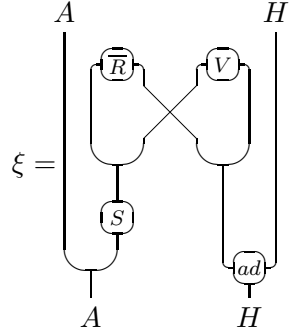
\[ \xi =
 \begin{tangle}
\object{A}\Step\Step\Step\Step \object{H}\\
\id\step\ro {\overline R} \Step \ro V
\step\id\\
\id\step\id\step\step\XX\step\step\id\step\id\\
\id\step\cu\step\step\cu\step\id\\
\id\step\step\S\step\step\step\step\id\step\step\id\\
\cu\step\step\step\step\tu {ad}\\
\step\object{A}\Step\Step\Step\object{H}
 \end{tangle}
\]
 and  $\xi $  is surjective,
where $ad$ denotes the left adjoint action of $H$.

(ii) Furthermore,   if $A$ and $H$  are finite-dimensional,   then
$\xi$  is a bijective map from $\underline D$ onto $\underline A
\otimes \underline H$. That is,
 in braided tensor category $(_{D} {\cal M},   C^{R_D}),  $
 there exist morphisms $\phi $  and $\psi $
 such that $$ \underline D \cong
 \underline A{} ^\phi \bowtie ^\psi \underline H
\hbox { \ \ \ ( as Hopf algebras )}$$ and the isomorphism is $(\pi
_{\underline A} \otimes \pi _{\underline H} ) \Delta _{\underline
D}$.

(iii) If $H$ is commutative or $V=R$,     then $\xi = id
_{\underline D} $.

\end {Theorem}
{\bf Proof.} (i)
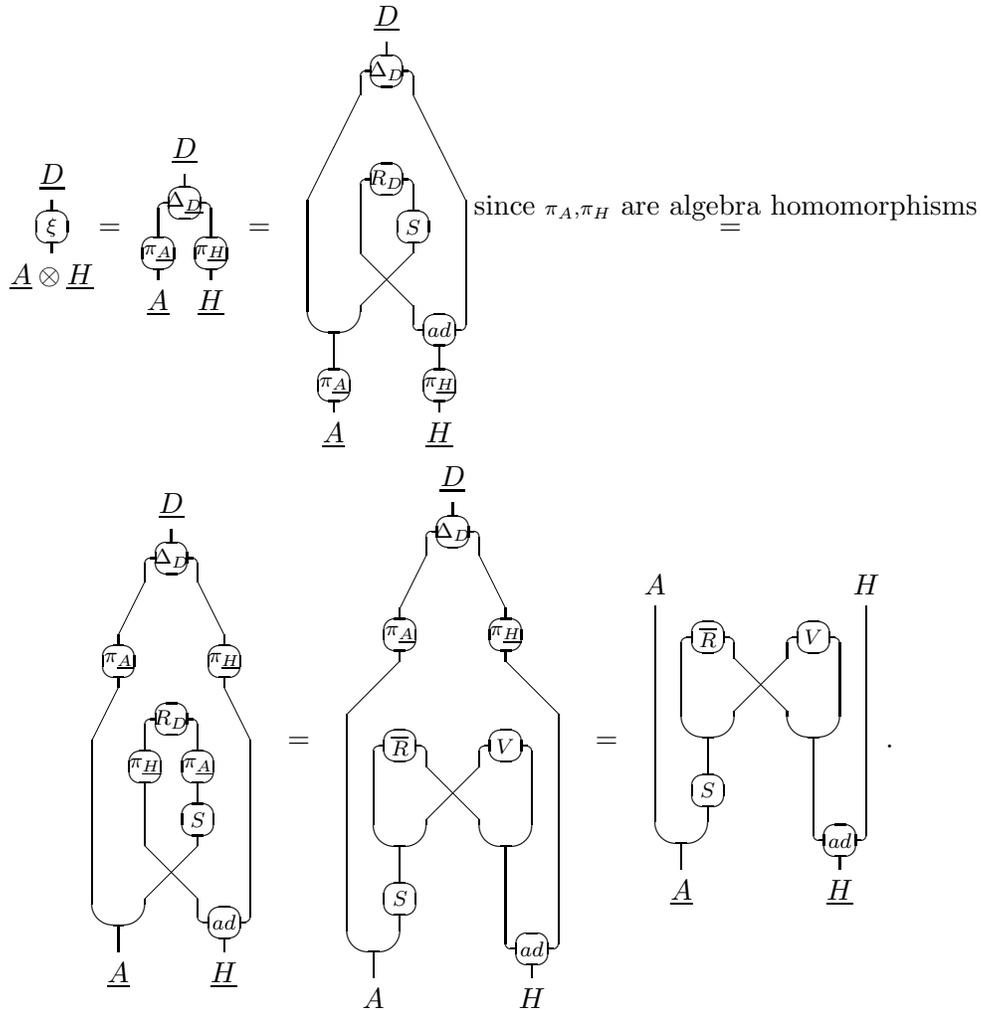
\[
\step\step
\begin{tangle}
\object{\underline{D}}\\
\O \xi\\
\object{\underline{A}\otimes \underline{H}}
\end{tangle}
\\ \ \ \ \ = \ \ \
\begin{tangle}
\step \object{\underline{D}}\\
\td {\Delta_{\underline{D}}} \\
\O {\pi_{\underline{A}}}\Step\O {\pi_{\underline{H}}}\\
\object{\underline{A}}\Step\object{\underline{H}}\\
\end{tangle}
\ \ \ = \ \ \
\begin{tangle}
\step\Step\object{\underline{D}}\\

\step [2] \td {\Delta_{D}} \\

\step \ne1 \step [2] \nw1\\

 \dd \step\ro {R_{D}}\step \d\\
\id\Step \id\Step\S\Step\id\\
\id\step\step\XX \step\step\id \\
\cu \step  \step \tu {ad}\\
\step\O {\pi_{\underline{A}}}\Step\Step\O {\pi_{\underline{H}}}\\
\step\object{\underline{A}}\Step\Step\object{\underline{H}}
 \end{tangle}
\stackrel { \hbox {since } \pi _A, \pi _H \hbox { are algebra
homomorphisms }}{=}
\]\[
\begin{tangle}
\Step \step\object{\underline{D}}\\

\step [2] \td {\Delta_{D}} \\

\step \ne1 \step [2] \nw1\\

\step\O {\pi_{\underline{A}}}\Step\step\step\O {\pi_{\underline{H}}}\\
\dd  \step [1] \ro {R_{D}} \step \d\\
\id \Step \O {\pi_{\underline{H}}}\Step\O {\pi_{\underline{A}}}\Step\id\\
\id\Step \id\Step\S\Step\id\\
\id\step\step\XX \step\step\id \\
\cu \step  \step \tu {ad}\\
\step\object{\underline{A}}\Step\Step\object{\underline{H}}
 \end{tangle}
 \ \ \ = \ \ \
 \begin{tangle}
\Step\Step \object{\underline{D}}\\

\step\step [2] \td {\Delta_{D}} \\

\step\step \ne1 \step [2] \nw1\\
\step\step\O {\pi_{\underline{A}}}\Step\step\step\O {\pi_{\underline{H}}}\\
\step\ne2 \Step\Step\nw2\\
 \id\step\ro {\overline R}
\Step
\ro {V} \step\id\\
\id\step\id\step\step\XX\step\step\id\step\id\\
\id\step\cu\step\step\cu\step\id\\
\id\step\step\S\step\step\step\step\id\step\step\id\\
\cu\step\step\step\step\tu {ad}\\
\step\object{\underline{A}}\Step\Step\Step\object{\underline{H}}
 \end{tangle}
\ \ \ = \ \ \
 \begin{tangle}
\object{A}\Step\Step\Step\Step \object{H}\\
\id\step\ro {\overline{R}} \Step \ro {V}
\step\id\\
\id\step\id\step\step\XX\step\step\id\step\id\\
\id\step\cu\step\step\cu\step\id\\
\id\step\step\S\step\step\step\step\id\step\step\id\\
\cu\step\step\step\step\tu {ad}\\
\step\object{\underline{A}}\Step\Step\Step\object{\underline{H}}
 \end{tangle} \ \ .
\]

(ii) By the proof of  (i),   $\xi$ is bijective and\\
 \[ \bar \xi =
 \begin{tangle}
\object{A}\Step\Step\Step\Step \object{H}\\
\id\step\ro {\overline{V}} \Step \ro R
\step\id\\
\id\step\id\step\step\XX\step\step\id\step\id\\
\id\step\cu\step\step\cu\step\id\\
\id\step\step\S\step\step\step\step\id\step\step\id\\
\cu\step\step\step\step\tu {ad}\\
\step\object{A}\Step\Step\Step\object{H}
 \end{tangle} \ \ .
\]
Applying the co-factorization theorem \ref {4.1.4},   we complete
the proof of (ii).

(iii) follows from (i). \begin{picture}(5,5)
\put(0,0){\line(0,1){5}}\put(5,5){\line(0,-1){5}}
\put(0,0){\line(1,0){5}}\put(5,5){\line(-1,0){5}}
\end{picture}\\

\begin {Remark} \label{r1}  { \  }\end {Remark}

 Under the assumption of Theorem \ref {4.1.5},
if we set \[
 \begin{tangle}
\step\object{A}\\
\td { \phi^{'}}\\
\object{H}\Step\object{A}
 \end{tangle}
 \ \ \ =: \ \ \
  \begin{tangle}
  \Step\Step\object{A}\\
  \ro R\Step\id\Step \ro {\overline{R}}\\
  \XX\Step\id\Step\XX\\
  \id\Step\d\step\XX\Step\id\\
  \id\Step\sw2\nw2\step\step\cu\\
 \cu\Step\step\cu\\
\step\object{H}\Step\Step\object{A}
   \end{tangle}
     \Step \hbox {\ \hbox {and} \ }\Step
    \begin{tangle}
\step\object{H}\\
\td {\psi^{'}}\\
\object{H}\Step\object{A}
 \end{tangle}
 \ \ \ =: \ \ \
  \begin{tangle}
  \Step\Step\object{H}\\
  \ro R\Step\id\Step\ro {\overline{R}}\\
  \XX\Step\id\Step\XX\\
  \id\Step\XX\step\dd\Step\id\\
  \cu\step\sw2\nw2\step\step\step\id\\
\step \cu\Step\step\cu\\
\step\step\object{H}\Step\Step\step\object{A}
   \end{tangle} \ \  ,
\]
 then $A{}^{\phi '} \bowtie ^{\psi '} H= A
\bowtie ^R H$   by \cite [Lemma 1.3] {Ch98}.  We now see the
relation among $\phi,   \psi,   R, \phi'$ and $\psi '$:
\[ \phi \ \  \stackrel {\hbox {by proof of Th \ref {4.1.4} }}{=}\  \zeta (id \otimes \eta
)\ \  = \ \
 \begin{tangle}
\object{A}\Step\object{\HH\obox 1\eta}\\
\tu {\overline{\xi}}\\
  \td {\Delta_{\underline{D}}}\\
\O {\pi_{\underline{H}}}\Step\O {\pi_{\underline{A}}}\\
 \object{H}\Step\object{A}
\end{tangle} \ \ \ \  \stackrel {\hbox {by Th \ref {4.1.5} (iii)}} {=}\ \ \
 \]
 \[\begin{tangle}
\step\object{A}\Step\Step\Step\Step\object{\HH\obox 1\eta}\\
\step\id\step\ro {\overline{V}}\step\step
\ro R\step\id\\
\step\id\step\id\step\step\XX\step\step\id\step\id\\
\step\id\step\cu\step\step\cu\step\id\\
\step\id\step\step\S\step\step\step\step\id\step\step\id\\
\step\cu\Step\Step\tu {ad}\\
\step\cd\Step\Step\cd\\
\step\id\step\td {\phi^{'}}\Step
\td {\psi^{'}}\step\id\\
\step\QQ \epsilon \step\id\Step\XX\Step\id \step\QQ \epsilon\\
\step\step\cu\Step\cu\\
\Step\ne2\Step\step\step\nw2\\
\step\id\step  \ro R\Step
 \ro U \step\id\\
\step\id\step\id\Step\XX\Step\id\step\id\\
\step\id\step\cu\Step\cu\step\id\\
\step\id\Step\id\Step\Step\S\Step\id\\
\step\id\Step\nw2\step\step\ne2\Step\id\\
\step\id\Step\sw2\step\step\se2\Step\id\\
\step\cu\Step\Step\tu {ad}\\
\step\step\object{H}\Step\Step\Step\object{A}
\end{tangle}
  \ = \
 \begin{tangle}
 \Step \step\Step\object{A}\\
 \step [4]\td {\phi ^{'}}\\

\step[3] \ne1 \step [2] \nw1\\
 \Step\ne2\Step\step\step\nw2\\
\step\id\step  \ro R\Step
 \ro U\step\id\\
\step\id\step\id\Step\XX\Step\id\step\id\\
\step\id\step\cu\Step\cu\step\id\\
\step\id\Step\id\Step\Step\S\Step\id\\
\step\id\Step\nw2\step\step\ne2\Step\id\\
\step\id\Step\sw2\step\step\se2\Step\id\\
\step\cu\Step\Step\tu {ad}\\
\step\step\object{H}\Step\Step\Step\object{A}
\end{tangle}
 \ = \
 \begin{tangle}
 \Step\Step\Step\object{A}\\
 \ro R\Step\Step\id\Step\ro {\overline{R}}\\
 \XX\Step\Step\id \step[2] \XX\\

 \id\Step\nw2\Step\step[1]\XX\Step\id\\

 \id\Step\Step\XX\Step\id\step[2]\id\\

 \id\Step \Step\id\Step\nw1\step [1]\cu\\

 \id\Step \Step\id\step [3]\cu \\

 \id\Step\Step\id\Step\Step\id\Step\\
 \Cu\Step\Step \nw2\\
 \Step\id\step  \ro {R}\Step
 \ro U\step\id\\
\Step\id\step\id\Step\XX\Step\id\step\id\\
\Step\id\step\cu\Step\cu\step\id\\
\Step\id\Step\id\Step\Step\S\Step\id\\
\Step\id\Step\nw2\step\step\ne2\Step\id\\
\Step\id\Step\sw2\step\step\se2\Step\id\\
\Step\cu\Step\Step\tu {ad}\\
\Step\step\object{H}\Step\Step\Step\object{A}
 \end{tangle}
\]
 $ \hbox {and \ \ \ } \psi \ \  \stackrel { \hbox { by proof of Th \ref {4.1.4} }}{=}
  \ \ \  \zeta (\eta  \otimes
id  ) = $
 \[\begin{tangle}
\object{\HH\obox 1\eta}\Step\object{H}\\
\tu {\overline{\xi}}\\
 \td {\Delta _{\underline D}} \\
\O {\pi_{\underline{H}}}\Step\O {\pi_{\underline{A}}}\\
 \object{H}\Step\object{A}
\end{tangle} \ \ \  \stackrel { \hbox{ by Th \ref {4.1.5} (iii) } }{=} \ \ \
 \begin{tangle}
\step\object{\HH\obox 1\eta}\Step\Step\Step\Step\object{H}\\
\step\id\step\ro {\overline{V}}\step\step
\ro R\step\id\\
\step\id\step\id\step\step\XX\step\step\id\step\id\\
\step\id\step\cu\step\step\cu\step\id\\
\step\id\step\step\S\step\step\step\step\id\step\step\id\\
\step\cu\Step\Step\tu {ad}\\
\step\cd\Step\Step\cd\\
\step\id\step\td {\phi^{'}}\Step
\td {\psi^{'}}\step\id\\
\step\QQ \epsilon \step\id\Step\XX\Step\id \step\QQ \epsilon\\
\step\step\cu\Step\cu\\
\Step\ne2\Step\step\step\nw2\\
\step\id\step  \ro R\Step
 \ro U\step\id\\
\step\id\step\id\Step\XX\Step\id\step\id\\
\step\id\step\cu\Step\cu\step\id\\
\step\id\Step\id\Step\Step\S\Step\id\\
\step\id\Step\nw2\step\step\ne2\Step\id\\
\step\id\Step\sw2\step\step\se2\Step\id\\
\step\cu\Step\Step\tu {ad}\\
\step\step\object{H}\Step\Step\Step\object{A}
\end{tangle}
\]\[
 \ \ \ = \ \ \
  \begin{tangle}
  \step\Step\Step\Step\Step\object{H}\\
\step\step\ro {\overline{V}}\step\step
\ro {R}\step\id\\
\Step\id\step\step\XX\step\step\id\step\id\\
\step\step\cu\step\step\cu\step\id\\
\step\step\step\S\step\step\step\step\id\step\step\id\\
\Step\step\id\Step\Step\tu {ad}\\
\Step\step\id\Step\step\step\dd\\
\Step\td {\phi^{'}}\Step
\td {\psi^{'}}\\
\step\step \id\Step\XX\Step\id \\
\step\step\cu\step\step\cu\\
\Step\ne2\Step\step\step\nw2\\
\step\id\step  \ro R\Step
 \ro U\step\id\\
\step\id\step\id\Step\XX\Step\id\step\id\\
\step\id\step\cu\Step\cu\step\id\\
\step\id\Step\id\Step\Step\S\Step\id\\
\step\id\Step\nw2\step\step\ne2\Step\id\\
\step\id\Step\sw2\step\step\se2\Step\id\\
\step\cu\Step\Step\tu {ad}\\
\step\step\object{H}\Step\Step\Step\object{A}
 \end{tangle}
 \ \ \ = \ \ \
  \begin{tangle}
  \Step\step\Step\Step\Step\Step\object{H}\\
\Step\step\step\ro {\overline{V}}\step\step
\ro R\step\id\\
\Step\Step\id\step\step\XX\step\step\id\step\id\\
\Step\step\step\cu\step\step\cu\step\id\\
\Step\step\step\step\S\step\step\step\step\id\step\step\id\\
\Step\Step\step\id\Step\Step\tu {ad}\\
 \Step\Step\step\d \Step\Step\id\\
 \ro {\overline{R}}\Step\Step\id\Step\Step\id\Step\Step
 \ro {R}\\
 \XX\Step\Step\nw2\Step\ne2\Step\Step\XX\\
 \id\Step\nw2\Step\Step\ne2
 \nw2\Step\Step\ne2\Step\id\\
 \id\Step\Step\XX\Step\Step\XX\Step\Step\id\\
\Cu\step\step\nw2\Step\ne2\Step\Cu\\
 \Step\id\Step\Step\sw2\step\step\se2\Step\Step\id\\
 \Step\id\Step\Step\id\Step\Step\id\Step\Step\id\\
 \Step\Cu\Step\Step\Cu\\
 \Step\Step\id\step  \ro R\Step
 \ro U\step\id\\
\Step\Step\id\step\id\Step\XX\Step\id\step\id\\
\Step\Step\id\step\cu\Step\cu\step\id\\
\Step\Step\id\Step\id\Step\Step\S\Step\id\\
\Step\Step\id\Step\nw2\step\step\ne2\Step\id\\
\Step\Step\id\Step\sw2\step\step\se2\Step\id\\
\Step\Step\cu\Step\Step\tu {ad}\\
\Step\Step\step\object{H}\Step\Step\Step\object{A}
 \end{tangle}\ \ \ .
\]

Furthermore,   $\psi = \psi ' $ and $\phi = \phi '$  \ \ when $H$ is
commutative or $R = V= U^{-1}.$ In this case,   $\underline A{}
^{\phi }\bowtie ^{\psi}\underline H= \underline
 {A {}^{ \phi '} \bowtie ^ {\psi '} H}
  = \underline {A \bowtie ^R H}$
  \ \  as Hopf algebras living
  in braided tensor category $( {}_D {\cal M},   C^{R_D}).$

\begin {Example} \label {4.1.7}
(cf. \cite [Example 2.11] {Ch98})  Let  $A $ be Sweedler's four
dimensional Hopf algebra $H_4$ and $H= k{\bf Z}_2$ with $char k
\not= 2 $. $A$ is generated by $g$ and $x$ with relations

$$g^2 =1, \ \ x^2 = 0, \  \ xg =-gx .$$
$H$ is generated by $a$ with $a^2 =1.$
 Let

               $P = \frac {1} {2} (1 \otimes 1 + 1 \otimes g + g \otimes 1
- g \otimes g) + \frac {\mu} {2} (x \otimes x + x \otimes gx + gx
\otimes gx - gx \otimes x);$

$Q = \frac {1} {2} (1 \otimes 1 + 1 \otimes a + a \otimes 1 - a
\otimes a); $

$R = \frac {1} {2} (1 \otimes 1 + 1 \otimes a + g \otimes 1 - g
\otimes a).$

It is clear that $P$ and $Q$ are  quasitriangular structures of $A$
and $H$ respectively. $R$  is a weak $R$-matrix of $A \otimes H$
with $R= R^{-1}.$  Let $D = A \bowtie ^R H.$ Thus, for any
quasitriangular structure  $R_D$ of $D$, the braided group analogue
$\underline D$ of $D$ is a double cross coproduct of braided  group
analogues $\underline A$  and $\underline H$. That is, in braided
tensor category $({}_D{\cal M}, C^{R_D}),$
$$\underline {A \bowtie ^R H} = \underline  A ^\phi \bowtie ^\psi
\underline H \hbox { \ \ \ (as Hopf algebras )}.$$ Here $\psi $ is
trivial coaction and
 $\phi (z) =
\sum R^{(1)} z \bar R^{(1)} \otimes R^{(2)} \bar R^{(2)}$   for any
$z \in A$ with $R = \bar R = \sum R^{(1)} \otimes R^{(2)}.$

\end {Example}

\section {The factorization  of   ordinary bialgebras }\label {s10}

Throughout  this section, we work in  the braided tensor category
 of  vector spaces over
field $k$ with ordinary twist braiding. In this section, we give the
relation between bialgebra or Hopf algebra $H$ and its factors.

\begin {Lemma} \label {4.3.1}  (cf.  \cite  {Mo93}, \cite  {Ra94} ).
If $H$ is a finite dimensional Hopf algebra  with $char { \  \ } k
=0 ,$ then the following conditions are equivalent.

(i) $H$ is semisimple.

(ii) $H$ is cosemisimple.

(iii) $S_H^2 = id _H$.

(iv) $tr(S_H^2) \not=0.$
\end {Lemma}
{\bf Proof.} (i) and (ii) are equivalent  by \cite [Theorem 2.5.2]
{Mo93}.

(i) and (iii) are equivalent  by \cite [Proposition 2 (c)] {Ra94}.

(iii) and (iv) are equivalent, since  (iii) implies (iv) and (iv)
implies (iii) by \cite [Theorem 2.5.2] {Mo93}.
\begin{picture}(8,8)\put(0,0){\line(0,1){8}}\put(8,8){\line(0,-1){8}}\put(0,0){\line(1,0){8}}\put(8,8){\line(-1,0){8}}\end{picture}

\begin {Lemma} \label {4.3.2}  (cf.  \cite  {Mo93}, \cite  {Ra94} )
If $H$ is a finite dimensional Hopf algebra
 $char { \ \ } k > (dim { \ \ } H )^2 $,
then the following conditions are equivalent.

(i) $H$ is semisimple and cosemisimple.

(ii) $S_H^2 = id _H$.

(iii) $tr(S_H^2) \not=0.$
\end {Lemma}
{\bf Proof.} (i) and (iii) are equivalent  by \cite [Proposition 2
(c)] {Ra94}.

 (i) implies (ii)  by   \cite [Theorem 2.5.3] {Mo93}.

Obviously, (ii) implies (iii).
\begin{picture}(8,8)\put(0,0){\line(0,1){8}}\put(8,8){\line(0,-1){8}}\put(0,0){\line(1,0){8}}\put(8,8){\line(-1,0){8}}\end{picture}

When $H$ is a finite-dimensional Hopf algebra, let $\Lambda _H ^l$
and $\Lambda _H^r$  denote a non-zero left integral and non-zero
right integral of $H$ respectively.

\begin {Lemma} \label {4.3.3}  Let $A$, $H$  and
$D = A {\bowtie} H$ all be finite dimensional Hopf algebra. Then

(i)  there are $u\in H, v\in A$ such that
 $$\Lambda _D {}^l =
 \Lambda _A^l \otimes u, \hbox { \ \ \ }
\Lambda _D^r = v \otimes \Lambda _H^r ; $$

(ii) If  $ A { \bowtie } H  $ is semisimple, then $A$ and $H$ are
semisimple;

(iii) If $A { \bowtie } H$ is unimodular, then
             $ \Lambda_ D =
 \Lambda _A ^l \otimes \Lambda _H^r,$
 and              $A  { \bowtie } H$  is semisimple iff
 $A$ and $H$ are semisimple.
 \end {Lemma}
 {\bf Proof .}  (i)
 Let $a  ^{(1)},  a  ^{(2)}, \cdots, a  ^{(n)}$  and $ h ^{(1)},  h  ^{(2)},
 \cdots ,  h^{(m)}$ be the basis of $A$  and  $H$ respectively.
Assuming $$ \Lambda _D^l = \sum k_{ij} (a ^{(i)} \otimes h^{(j)}) $$
where $k_{ij} \in k$, we have that

$$a \Lambda _D^l = \epsilon (a) \Lambda _D^l $$ and
$$   \sum \epsilon (a) k_{ij} (a ^{(i)} \otimes h^{(j)})=
 \sum k_{ij} (a a ^{(i)} \otimes h^{(j)}) ,$$   for any $a \in A.$

 Let $x_j = \sum _i k_{ij}a  ^{(i)}$. Considering $\{h^{(j)} \}$ is a base of $H$,
 we get $ x_j $ is a left integral of $A$  and there exists  $k_j \in k$ such that
 $x_j = k_j \Lambda _A^l$  for $j = 1, 2, \cdots m$.

Thus $$\Lambda _D^l = \sum _j k_j (\Lambda  _A \otimes h^{(j)}) =
\Lambda _A^l \otimes  u, $$

where $u = \sum _j k_j h^{(j)}.$

Similarly, we have  that $\Lambda _D^r = v \otimes \Lambda _H^r$.

(ii) and (iii)  follow from part (i) .
\begin{picture}(8,8)\put(0,0){\line(0,1){8}}\put(8,8){\line(0,-1){8}}\put(0,0){\line(1,0){8}}\put(8,8){\line(-1,0){8}}\end{picture}

\begin {Proposition}  \label {4.3.4}

If $A, H$ and $D=A \bowtie H$ are finite dimensional  Hopf algebras,
then

(i) $A  \bowtie H$ is  semisimple and cosemisimple  iff $A$ and $H$
are semisimple and cosemisimple;

(ii) If $A$ and $H$ are involutory and character $char k $ of $k$
does  not divides $dim D$,
 then  $A  \bowtie H$, $A$ and  $H$ are  semisimple and cosemisimple.

(iii)
  $A  \bowtie H$ is   cosemisimple  iff $A$ and $H$ are
 cosemisimple;

(iv) If $char k=0$, then $D$ is semisimple iff $D$ is cosemisimple
iff
                        $A$ and $H$ are semisimple iff $A$ and $H$
                        are  cosemisimple iff  $S_D^2 =id$
                        iff $S_A^2 = id$ and $S_H^2 = id.$

\end {Proposition}

{\bf Proof .} (i) It follows from  \cite [Proposition 2] {Ra94}.

(ii)  It follows from \cite [Theorem 4.3] {La71} and part (i);

(iii) If $A \bowtie H$ is  cosemisimple, then $(A \bowtie H) ^*$  is
semisimple. Thus $A$ and $H$  are cosemisimple.

Conversely, if $A$ and $H$ are cosemisimple, then $(A\bowtie H)^* $
is semisimple. Thus
 $A \bowtie H$ is cosemisimple.

(iv) It follows from  Lemma \ref {4.3.1}.
\begin{picture}(8,8)\put(0,0){\line(0,1){8}}\put(8,8){\line(0,-1){8}}\put(0,0){\line(1,0){8}}\put(8,8){\line(-1,0){8}}\end{picture}

\chapter {Braided  L
ie Algebras}\label {c4'}

The theory of Lie superalgebras has been developed systematically,
which includes  the representation theory and classifications of
simple Lie superalgebras and their varieties \cite {Ka77} \cite
{BMZP92}.  In many physical applications or in pure mathematical
interest, one has to consider not only ${\bf Z}_2$- or ${\bf Z}$-
grading but also $G$-grading of Lie algebras, where $G$ is an
abelian group equipped with a skew symmetric bilinear form given by
a 2-cocycle. Lie algebras in symmetric and more general categories
were discussed in \cite {Gu86} and \cite {GRR95}. A sophisticated
multilinear version of the Lie bracket was considered in \cite
{Kh99} \cite {Pa98}. Various generalized Lie algebras have already
appeared  under different names, e.g. Lie color algebras, $\epsilon
$ Lie algebras \cite {Sc79}, quantum and braided Lie algebras,
generalized Lie algebras \cite {BFM96} and $H$-Lie algebras \cite
{BFM01}.

In  \cite {Ma94c}, Majid introduced braided Lie algebras from
geometrical point of view, which have attracted attention in
mathematics and mathematical physics (see e.g. \cite {Ma95b} and
references therein).

In this chapter, braided m-Lie algebras induced by multiplication
are introduced, which generalize Lie algebras, Lie color algebras
and quantum Lie algebras. The necessary and sufficient conditions
for the braided m-Lie  algebras to be strict Jacobi braided Lie
algebras are given. Two classes of braided m-Lie algebras are given,
which are generalized matrix braided m-Lie algebras and braided
m-Lie subalgebras of $End _F M$, where $M$ is a Yetter-Drinfeld
module over $B$ with dim $B< \infty $ . In particular, generalized
classical braided m-Lie algebras $sl_{q, f}( GM_G(A),  F)$  and
$osp_{q, t} (GM_G(A), M, F)$  of generalized matrix algebra
$GM_G(A)$ are constructed  and their connection with  special
generalized matrix Lie superalgebra $sl_{s, f}( GM_{{\bf Z}_2}(A^s),
F)$  and orthosymplectic generalized matrix Lie super algebra
$osp_{s, t} (GM_{{\bf Z}_2}(A^s), M^s, F)$  are established. The
relationship between representations of braided m-Lie algebras and
their  associated algebras are established.

Throughout, $ F$ is a field, $G$ is an additive group and $r$ is a
bicharacter of $G$; $\mid $$ x $$\mid $ denotes the degree of $x$
and  $({\cal C}, C)$ is a braided tensor category with braiding $C$.
We  write $W \otimes f $ for $id _W \otimes f$ and $f \otimes W$ for
$f \otimes id _W$. Algebras discussed here may not have unity
element.

\section {Braided  m-Lie Algebras}
In this section we introduce  braided m-Lie  algebras and (strict)
Jacobi braided Lie algebras.  We give the necessary and sufficient
conditions for the braided m-Lie  algebras to be strict Jacobi
braided Lie algebras.

\begin {Definition} \label {4'.1.1}
Let  $(L, [\ \ ])$ be an object in the braided tensor category $
({\cal C }, C)$ with morphism $[\ \ ] : L \otimes L \rightarrow L$.
If there exists an algebra $(A, m)$ in  $ ({\cal C }, C)$ and
monomorphism $\phi : L \rightarrow A$ such that   $\phi [\  \ ] = m
(\phi \otimes \phi ) - m (\phi \otimes \phi )  C_{L, L},$   then
$(L, [ \ \  ])$ is called a braided m-Lie algebra in $ ({\cal C },
C)$ induced by multiplication of $A$ through $\phi $. Algebra $(A,
m) $ is called an algebra  associated to $(L, [ \ \ ])$.

\end {Definition}
A Lie  algebra is a braided m-Lie algebra in the category of
ordinary vector spaces, a Lie color algebra is a braided  m-Lie
algebra
 in symmetric braided tensor category $ ({\cal M} ^{FG}, C^r)$ since the
 canonical map $\sigma: L \rightarrow U(L)$ is injective (see
\cite [Proposition 4.1]{Sc79}), a  quantum Lie  algebra is a braided
m-Lie algebra in the Yetter-Drinfeld category $ (^B_B{\cal YD}, C)$
by \cite [Definition 2.1 and Lemma 2.2]{GM03}), and a ``good"
braided Lie  algebra is a braided  m-Lie algebra
 in the Yetter-Drinfeld category $ (^B_B{\cal YD}, C)$ by
 \cite [Definition 3.6 and Lemma 3.7]{GM03}).
For  a cotriangular Hopf algebra $(H, r)$, the $(H,r)$-Lie algebra
defined in \cite [4.1] {BFM01} is a
 braided m-Lie  algebra in the braided
tensor category $({}^H{\cal M}, C^r)$.  Therefore, the braided m-Lie
algebras generalize most known generalized Lie algebras.

For an algebra $(A, m)$ in $({\cal C}, C)$, obviously $L = A$ is  a
braided m-Lie algebra under operation $ [\ \ ] = m  - m   C_{L, L}$,
which  is induced by $A$ through $id _A$. This  braided  m-Lie
algebra is written as $A^-$.

If $V$ is an object in ${\cal C}$ and $C_{V,V} = C_{V, V}^{-1},$
then we say that the braiding is symmetric on $V$.

\begin {Example} \label {4'.1.2} If $H$ is a  braided Hopf algebra in the
Yetter-Drinfeld module category $(^B_B {\cal YD}, C)$ with $B = FG$
and $C(x\otimes y) = r(\mid $$ y $$\mid , \mid $$ x $ $\mid )y
\otimes x $ for any homogeneous elements $x, y \in H$, then $P(H)=:
\{ x \in H \mid x \hbox { is a primitive element }\} $ is a braided
m-Lie  algebra iff the braiding $C$ is symmetric on $P(H)$.
\end {Example}

Indeed, it is easy to check  that $P(H)$ is the Yetter-Drinfeld
module. By simple computation we have $\Delta ([x, y ]) = [x, y]
\otimes 1 + 1 \otimes [x, y] + (1 - r(\mid $$x$$ \mid , \mid $$y$$
\mid )r(\mid $$y$$ \mid , \mid $$x$ $ \mid )) x \otimes y$ for any
homogeneous elements $x, y \in P(H).$
 Thus $[x, y] \in P(H)$ iff $r(\mid $$x$$ \mid , \mid $$y$$ \mid )
 r(\mid $$y$$ \mid , \mid $$x$$ \mid )=1,$ as asserted.

\begin {Theorem} \label {4'.1.3} Let  $(L, [\ \ ])$ be a braided  m-Lie
 algebra in $({\cal C}, C)$.

 (i)  $(L, [\  \ ])$ satisfies the braided anti-symmetry (or quantum
 anti-symmetry):

 (BAS): $[\  \ ] = - [\  \ ] C_{L, L}$

\noindent  if and only if
 $mC_{L,L}= mC_{L,L}^{-1}.$

 (ii) If the braided anti-symmetry holds, then braided m-Lie  algebra
 $(L, [\  \  ], m)$ satisfies the (left) braided primitive  Jacobi identity:

  (BJI):\ \ $[\ \ ] (L\otimes [\ \ ]
)+[\ \ ] (L\otimes [\ \ ] ) (L  \otimes C_{L,L}^{-1}) (C_{L,L}
\otimes L) + [\ \ ] (L \otimes [\ \ ]
)(C_{L,L}^{-1} \otimes L)(L \otimes C_{L,L} ))=0$, \\
 and the right  braided primitive Jacobi identity:

 (BJI'):\ \ $[\ \ ] ([\  \ ]\otimes L
)+[\ \ ] ([\  \ ]\otimes L ) (L \otimes C_{L,L}{}^{-1}) (C_{L,L}
\otimes L) + [\ \ ] ([\  \ ]\otimes L)(C_{L,L}{}^{-1} \otimes L) (L
\otimes C_{L,L} ))=0.$
 \end {Theorem}

{\bf Proof .} (i) Assume that $(L, [\  \ ])$ satisfies the braided
anti-symmetry. Since $[\  \ ] = - [\  \ ]C$ we have $m- m C = mCC -
mC$ and $m =m CC,$ which implies $mC= m C^{-1}.$ The necessity is
clear.

(ii) \begin {eqnarray*}
 \hbox { l.h.s. of  (BJI)} &=& m(L\otimes m)- m(L \otimes m)(L \otimes C) \\
& &- m C(L \otimes m) + m C(L \otimes m) (L \otimes C)\\
 & &+ m (L \otimes m) (L\otimes C^{-1}) (C \otimes L) \\
& &+ m (L \otimes m)(L \otimes C)(L \otimes C^{-1})(C \otimes L)\\
& &- m C (L \otimes m)(L \otimes C^{-1})(C \otimes L) \\
& &+ m C(L \otimes m)(L \otimes C)(L \otimes C^{-1})(C \otimes L)\\
& &+ m (L \otimes m) (C^{-1} \otimes L )(L \otimes C)\\
& &- m(L \otimes m)(L \otimes C)(C^{-1} \otimes L)(L \otimes C)\\
& &- mC(L \otimes m) (C^{-1} \otimes L)(L \otimes C)\\
& &+ m C(L \otimes m)(L \otimes C)(C^{-1} \otimes L)(L \otimes C) .
\end {eqnarray*}

We first check that $-$(6th term) =  12th term and 4th term =
$-$(10th term).
 Indeed,\begin {eqnarray*}
\hbox { 12th term } &=& m C^{-1}(L \otimes m)
(L \otimes C^{-1})(C^{-1} \otimes L)(L \otimes C)\\
&=&m(m \otimes L)(L \otimes C^{-1})(C^{-1} \otimes L)(L \otimes
C^{-1})
(C^{-1} \otimes L)(L \otimes C)\\
&=&m (L \otimes m)(L \otimes C)(C^{-1} \otimes L)
(L \otimes C^{-1})(C^{-1} \otimes L)(L \otimes C)\\
&=&m (L \otimes m)(C^{-1} \otimes L)(L \otimes C^{-1})
(C \otimes L)(C^{-1} \otimes L)(L \otimes C)\\
&=& - \hbox {(6th term )}.
\end {eqnarray*}
\begin {eqnarray*}
\hbox { 4th term }&=& m (m \otimes L)(C^{-1} \otimes L)
(L \otimes C^{-1})(C^{-1} \otimes L)\\
&=&m(L \otimes m)(L \otimes C^{-1})(C^{-1} \otimes L)(L \otimes C) \\
&=& - \hbox {(10th term )}.
\end {eqnarray*}
We can similarly show that 1st term = $-$ (11th term), $-$(2nd term)
= 8th term, $-$ (3rd term) = 5th term,
 $-$(7th term) = 9th term.
Consequently, $(BJI)$ holds. We can similarly show that $(BJI)'$
holds.  $\Box$

Readers can prove the above with the help of   braiding diagrams.

\begin {Definition} \label {4'.1.4}  Let $[\  \ ]$  be a morphism from
$L\otimes L$  to $L$ in ${\cal C}$. If $(BJI)$ holds, then  $(L, [\
\ ])$  is called a  (left )  Jacobi braided Lie  algebra. If both of
(BAS) and (BJI) hold then $(L, [\  \ ])$  is called a  (left )
strict Jacobi braided  Lie  algebra.
\end {Definition}
Dually, we can define  right Jacobi braided  Lie  algebras and right
strict Jacobi braided Lie  algebras. Left (strict) Jacobi braided
Lie  algebras  are called  (strict) Jacobi braided   Lie  algebras
in general.

By Theorem \ref {4'.1.3} we have
 \begin {Corollary} \label  {4'.1.5} If  $(L , [\  \ ] )$ is a braided m-Lie
algebra, then the following conditions are equivalent:

(i)  $(L, [\  \ ])$  is a left (right) strict Jacobi braided Lie
algebra.

(ii) $m C_{L,L} = m C_{L,L}^{-1}.$

(iii)  $[\  \ ] C_{L,L} = [\ \ ] C_{L,L}^{-1}.$

 \end {Corollary}

Furthermore, if $L$ is a space graded by $G$ with bicharacter $r$,
then the braided primitive Jacobi identity (BJI) becomes:
$$r (\mid c \mid , \mid a \mid) [a, [b, c]] +
 r (\mid b \mid , \mid a \mid) [b, [c, a]]+ r (\mid c \mid , \mid b \mid)
 [c, [a, b]] =0   \ \ \ \ \ \ \ \  \ \ \ \ \ \ \ \ (*) $$
for any homogeneous elements $a, b, c \in L.$  That is, $L$ is a
Jacobi braided Lie  algebra if and only if (*) holds. For
convenience, we let $J(a, b, c)$ denote the left hand side of $(*)$.

We now recall the (left ) braided Lie algebra defined by Majid \cite
[Definition 4.1]{Ma94c}.  A (left) braided Lie algebra
  in ${\cal C}$ is a coalgebra $(L, \Delta , \epsilon  )$ in ${\cal C}$,
 equipped with a morphism  $[\  \ ]: L \otimes L \rightarrow L$
 satisfying the axioms:

(L1) $([\  \ ])(L \otimes [\  \ ])= [\  \ ]([\  \ ] \otimes [\  \
])(L \otimes C \otimes L)(\Delta \otimes L \otimes L)  $

(L2)  $C ([\  \ ] \otimes L) (L \otimes C) (\Delta \otimes L) = (L
\otimes [\  \ ])(\Delta \otimes L)$

(L3) $[\  \ ]$ is a coalgebra morphism in $({\cal C}, C).$

\noindent Axiom (L1) is called the left braided Jacobi identity.

There exist braided m-Lie  algebras which are neither  braided Lie
algebras nor Jacobi braided Lie  algebras as is seen from the
following example.
\begin {Example} \label {4'.1.6} (see \cite  {Ma95b} )
Let $L = F\{x\} / <x^n>$ be an algebra in $(^{F {\bf Z}_n} {\cal M},
C^r)$ with a primitive nth root $q$ of 1 and $r (k, m) = q ^{km}$
for any $k, m \in {\bf Z}_n,$ where $n$ is a natural number.

(i) If  $n > 3$,  then  $(L , [\  \ ])$ is  a braided m-Lie algebra
but is neither a Jacobi braided Lie  algebra nor a braided Lie
algebra of Majid.

(ii) If  $n = 3$,  then  $(L , [\  \ ])$ is a braided m-Lie algebra
and a Jacobi braided Lie  algebra  but is neither a strict Jacobi
braided Lie algebra nor a braided Lie algebra of Majid.
\end {Example}

{\bf Proof.} (i) Since $J(x, x, x) =3qx^3 ( 1 -q - q^2 + q^3) \not=0
$ we have that $(L, [\  \ ] )$ is not a Jacobi braided Lie algebra.
If $(L, [\  \ ], \Delta , \epsilon )$ is a braided Lie algebra,
since $[\  \  ]$ is a coalgebra homomorphism, we have  that
\begin {eqnarray*}\label {e1}
\hbox {[1  1]} &=& 0, \hbox { implying } \epsilon (1)=0\\
\hbox {[ } x , x\hbox{ ]} &=& x^2 (1 - q),  \hbox { implying }
\epsilon (x^2) (1-q)= \epsilon (x)^2\\
\hbox {[}x , x^2\hbox {]} &=& x^3 (1 - q^2),  \hbox { implying }
\epsilon (x^3) (1-q^2)= \epsilon (x) \epsilon (x^2)\\
&\cdots & \\
\hbox {[}x , x^{n-1}\hbox {]} &=& x^n (1 - q^{n-1}),  \hbox {
implying }
0 = \epsilon (x^n) (1-q^{n-1})= \epsilon (x) \epsilon (x^{n-1}).\\
\end {eqnarray*}
Thus $\epsilon (x^m ) = 0$ for $m =0, 1, 2, \cdots , n-1$, which
contradicts the fact that $(L, \Delta , \epsilon )$ is a coalgebra.

(ii) Since  $J(x^k,x^l,x^m)=0$ for  any $k, l , m \in {\bf Z}_3$, we
have that   $(L , [\  \ ])$ is a Jacobi braided Lie  algebra. It
follows from $[x, x] = x^2(1-q) \not= - q [x, x]$ that  $(L , [\  \
])$ is not a strict Jacobi braided Lie  algebra. $\Box$

Note that  $L$ in the above example may never be a braided Lie
algebra in the vector space category $_F {\cal M}$ with the ordinary
flip  $\tau $ (i.e. $\tau (x\otimes y) = y \otimes x$) as braiding,
although an extension of $L$ may become a braided Lie algebra in $_F
{\cal M}$. Furthermore, for the algebra $L = F\{ x \} $  in $(^{F
{\bf Z}} {\cal M}, C^r)$ with $q ^2\not=1$ and $r (k, m) = q ^{km}$
for any $k, m \in {\bf Z}$, the above conclusion holds.

\begin {Definition} \label {4'.1.7}
Let $(L, [\ \  ])$ be a braided  m-Lie algebra in $({\cal C}, C)$.
If $M$ is an object and  there exists a morphism $\alpha : L \otimes
M \rightarrow M$ such that $\alpha ([\ \ ] \otimes M) = \alpha (L
\otimes \alpha ) - \alpha (L \otimes \alpha ) (C \otimes M),$
 then $(M, \alpha )$ is called an $L$-module.
\end {Definition}

\section {Generalized Matrix Braided m-Lie Algebras}
As examples of the braided m-Lie algebras, we  introduce the
concepts of generalized matrix algebras (see \cite {Zh93}) and
generalized matrix braided m-Lie algebras.  We construct generalized
classical braided m-Lie algebras $sl_{q, f}( GM_G(A), F)$ and
$osp_{q, t} (GM_G(A), M, F)$ of generalized matrix algebra
$GM_G(A)$. We show  how generalized matrix  Lie color algebras are
related to Lie superalgebras for any abelian group $G$. That is,
 we establish  the relationship
between generalized matrix  Lie color algebras and Lie
superalgebras.

 Let $I$ be a set. For any $ i, j, l, k \in I,$ we choose a vector space $A_{ij}$
 over field $F$ and an  $F$-linear map
$\mu_{ijl}$ from  $A_{ij}\otimes A_{jl}$ into $A_{il}$ (written $\mu
_{ijl} (x, y)=xy)$ such that  $x(yz)=(xy)z$ for any $x\in A_{ij}$ ,
$y\in A_{jl} , z\in A_{lk}$.  Let $A$ be the external direct sum of
$\{ A_{ij} \mid i, j\in  I \}$. We define the multiplication in $A$
as
$$xy = \{ \sum _k x_{ik}y_{kj} \}$$
 for any $x=\{x_{ij}\},  y=\{y_{ij}\}\in A$ .
  It is easy to check that  $A$ is an algebra (possibly without unit ).
  We call $A$ a generalized matrix algebra,  or a gm algebra in short,
written as $A=\sum \{A_{ij} \mid i,  j\in  I\}$ or $GM_I(A)$. Every
element in $A$ is called a generalized matrix. We can easily define
gm ideals and gm subalgebras.  It is easy to define upper triangular
generalized matrices, strictly upper triangular generalized matrices
and diagonal generalized matrices under some total order $ \prec $
of $I$.

\begin {Proposition} \label {4'.2.1}
Let  $A=\sum \{A_{ij} \mid i, j\in I \}$ be  a gm algebra and  $G$
 an abelian group with $G=I$.
 Then $A$ is an algebra graded by $G$ with
$A_g = \sum _{i = j +g} A_{ij}$ for any $g\in G$. In this case, the
gradation is called a generalized matrix gradation, or gm gradation
in short.
\end {Proposition}
{\bf Proof.} For any $g, h \in G$, see that
\begin{eqnarray*}
A_g A_h &=& (\sum _{i = j +g} A_{ij})( \sum _{s = t +h}
A_{st})\\
 &\subseteq & \sum _{i = t+h +g}A_{i,t+h}A_{t+h,t}\\
 &\subseteq & A_{g +h}. \ \ \ \ \ \
 \end{eqnarray*}
Thus $A=\sum\{A_{ij} \mid i, j\in I \}= \sum _{g\in G} A_g$ is a
$G$-graded algebra. $\Box$

If ${\cal  C }$ is a small preadditive category  and
$A_{ij}=Hom_{\cal C}(j, i)$ is a vector space over $F$ for any $i,
j\in  I =$ objects in $ \ {\cal C}$, then we may easily show that
$\sum\{A_{ij} | i, j\in I\}$ is a generalized matrix algebra.
 Furthermore, if $V = \oplus_{g\in G} V_g$ is a
graded vector space over field $F$ with $A_{gh}=Hom_{F}(V_h, V_g)$
for any $g, h \in G$, then we call braided m-Lie  algebra
$A=\sum\{A_{ij} \mid i, j\in G \}$ the general linear braided m-Lie
algebra, written as $gl (\{V_g\}, F)$. If dim $V_g = n_g < \infty $
for any $g\in G$, then  $gl (\{V_g\}, F)$ is written as $gl
(\{n_g\}, F)$. Its braided m-Lie  subalgebras are called linear
braided m-Lie  algebras. In fact, $gl (\{n_g\}, F) = \{ f \in End _F
V \mid ker f \hbox { is finite codimensional } \}.$ When $G$ is
finite, we may  view $gl (\{n_g\}, F)$ as a block matrix algebra
over $F.$ When $G=0$, we denote $gl (\{n_g\}, F)$ by $gl (n, F)$,
which  is the ordinary ungraded  general linear Lie algebra.

Assume that $D$ is a directed graph ($D$ is possibly  an infinite
directed graph and also possibly not a simple graph ).  Let $I$
denote the vertex set of $D$, $x_{ij}$ an arrow from $i$ to $j$ and
$x=(x_{i_1i_2}, x_{i_2i_3}, \cdots , x_{i_{n-1}i_{n}})$ a path from
$i_1$ to $i_n$  via arrows $x_{i_1i_2},x_{i_2i_3}, \cdots ,
x_{i_{n-1}i_{n}}$. For two paths $x=(x_{i_1i_2},x_{i_2i_3}, \cdots ,
x_{i_{n-1}i_{n}})$ and $y=(y_{j_1j_2},y_{j_2j_3}, \cdots ,
y_{j_{m-1}j_{m}})$ of $D$ with $i_n=j_1$, we define the
multiplication of $x$ and $y$ as  $$xy = (x_{i_1i_2}, x_{i_2i_3},
\cdots , x_{i_{n-1}i_{n}}, y_{j_1j_2}, y_{j_2j_3}, \cdots ,
y_{j_{m-1}j_{m}}) .$$ \noindent Let $A_{ij}$ denote the vector space
over field $F$ with basis being all paths from $i$ to $j$, where $i,
j\in I$. Notice that we view every vertex $i$ of $D$ as a path from
$i$ to $i$, written $e_{ii}$ and $e_{ii}x_{ij} = x_{ij}e_{jj} = x
_{ij}$.   We can naturally define a linear  map from $A_{ij} \otimes
A_{jk}$ to $A_{ik}$ as $x\otimes y = xy $ for any two pathes $x \in
A_{ij}, y \in A_{jk}$. We may easily show that $\sum \{A_{ij}\mid i,
j \in  I \}$ is a generalized matrix algebra, which is  called a
path algebra,  written as $A(D)$ (see, \cite [Chapter 3]{ARS95}).

\begin {Example} \label {4'.2.3}
There are   finite-dimensional braided  m-Lie  algebras in braided
tensor category $(^{F {\bf Z}_3} {\cal M}, C^r)$ with a primitive
3th root $q$ of 1 and $r (k, m) = q ^{km}$ for any $k, m \in {\bf
Z}_3$. Indeed, for any natural number $n$,  we can construct a
generalized matrix braided  m-Lie algebra $A = \sum \{A_{ij} \mid i,
j \in {\bf Z}_3\}$ such that $dim \ A =n.$

(i) Let $A_{ij} = 0$ when $i\neq j$ but $A_{11}=F$. Thus $dim \ A =
1$.

(ii) Let $A_{ij} = 0$ when $i\neq j$ but $A_{11}=A_{22}=F$. Thus
$dim \ A = 2$.

(iii) Let $A_{ij} = 0$ when $i\neq j$ but $A_{11}= A_{22}=A_{33}=F$.
Thus $dim \ A = 3$.

(iv) Let $D$ be a directed graph with vertex set ${\bf Z}_3$ and
only one arrow from 1 to 2. Set $A = A(D).$ It is clear $dim (
A_{ij}) = 0$ when $i\neq j$ but $ dim(A_{11})=
dim(A_{22})=dim(A_{33})=dim ( A_{12})=1$. Thus $dim \  A = 4$.

(v) Let $D$ be a directed graph with vertex set ${\bf Z}_3$ and only
two arrows: one  from 1 to 2 and other one from 1 to 3. Set
$A=A(D)$. It is clear  $dim (A_{ij}) = 0 $ but $ dim (A_{11})= dim (
A_{22})=dim ( A_{33})=dim ( A_{12})= dim (A_{13})=1$ Thus $dim \ A =
5$.

vi) Let $D$ be a directed graph with vertex set ${\bf Z}_3$ and only
$n+2$ arrows: one  from 1 to 2,  one from 2 to 3 and the others from
1 to 3. Set $A= A(D).$ It is clear $dim (A_{ij}) = 0 $ but $ dim
(A_{11})= dim ( A_{22})=dim ( A_{33})=dim ( A_{12})= dim (A_{23})=1$
and $dim (A_{13}) = n +1$.  Thus $dim \ A = n+6 $ for $n =0, 1,
\cdots $.
\end {Example}

Let $A$ be a braided m-Lie algebra, $G$ be an abelian group with a
bicharacter $r $ and $W$ be a vector space over $F$.
\begin{Definition}\label{quantum-trace}
If $f$ is an $F$-linear map from $GM_G(A)$ to $W$ and satisfies the
following

(i)  $f (a) = \sum _{g \in G} r (g, g) f (a_{gg})$

(ii) $f (a_{ij} b_{ji}) =  f (b_{ji} a_{ij})$  for any $a, b \in A$
and $i, j \in G,$ then $f$ is called a generalized quantum trace
function from gm  algebra $GM_G(A)$ to $W$, written $tr_{q, f}$.
\end{Definition}
Set
$$sl_{q, f}(GM_G(A), F) = \{ a \in GM_G(A)  \mid  tr_{q,f} (a)=0 \}.$$

By computation,
$$tr_{q,f} [a, b] = r(u, u) ^{-1} r(g, g)\sum _{g \in G} (1-  r(u, u)^2 r
(u, g)r (g, u)) tr_{q,f} (b _{g, g +u} a_{g+u, g}) $$ for any
homogeneous elements $ a, b \in A, u \in G. $ Thus we get
\begin {Lemma}\label {4'.2.4}
$sl_{q, f}(GM_G(A), F)$ is a braided m-Lie algebra  with $[A,A]
\subseteq sl_{q, f}(GM_G(A), F)$ iff
$$\sum _{g \in G} (1-  r(u, u)^2 r (u, g)r (g, u)) tr_{q, f}
(b _{g, g +u} a_{g+u, g}) =0$$ \noindent for any homogeneous
elements $ a, b \in A, u \in G $ with $\mid $$a$$ \mid = u$ and
$\mid $$b$$ \mid = -u.$ If, in addition,  $tr_{q, f} (A_{ij}A_{ji})
\not=0$ for any $i, j \in G$, then
 $sl_{q, f}(GM_G(A), F)$ is a braided m-Lie algebra with
 $[A,A] \subseteq sl_{q, f}(GM_G(A), F)$
 iff $r$ is a skew symmetric bicharacter.
 \end {Lemma}

If $f(a_{gg}) = a_{gg}$ for any $g \in G,$ then from definition
\ref{quantum-trace} $f$ is a generalized quantum trace  from
$GM_G(A)$  to $GM_G(A).$ If $GM_G(A) = M_n(F)$ is the full matrix
algebra over $F$ with $G = {\bf Z}_n$, $r (i, j)=1$ and $f(a_{gg}) =
a_{gg} $ for any $i , j , g \in {\bf Z}_n$, then $f$ is the ordinary
trace function. If $GM_G(A) = gl (\{ n_g\}, F)$ and $f (a_{gg})= tr
(a_{gg})$ (i.e. $f (a_{gg})$ is the ordinary trace of the matrix
$a_{gg}$ )
 for any $g \in G$,
then $f$ is the quantum trace of the graded matrix algebra $gl (\{
n_g\}, F)$. In this case, $sl_{q, f} (gl (\{n_g\}, F))$ is simply
written as $sl _q (\{n_g\}, F).$

Let $G$ be an abelian group with a  bicharacter $r $ .  If $t$ is an
$F$-linear map from $GM_G(A)$ to $GM_G(A)$ such that  $t (A_{ij})
\subseteq A_{ji}$, $(t(a))_{ij} = t(a _{ji})$ and $t(ab)= t(b)t(a)$
for any $i, j \in G, a, b \in GM_G(A)$,
 then $t$ is called a generalized transpose on $GM_G(A).$
Given   $ 0\not=M \in GM_G(A) $, for any $u\in G,$ let $ osp _{q, t}
(GM_G(A), M, F )_u $ $= \{ a \in GM_G(A) _u \mid t(a_{u+g,
g})M_{u+g,h} = -  r(g,u) M_{g, u+h } a_{u+h, h}
    \hbox { for any } g, h \in G\}$
and   $osp _{q, t}(GM_G(A), M, F) = \oplus  _{u\in G} osp (GM_G(A),
M, F)_u.$

\begin {Lemma}\label {4'.2.5}
$osp_{q, t} (GM_G(A), M, F)$ is a braided m-Lie  algebra iff
$$\sum _{g \in G} (1-  r (u, v)r (v, u)) M_{g, h +u +v} a _{h +u +v, v+h}
b _{v +h, h} =0$$ \noindent for any homogeneous elements $ a, b \in
A, u, v \in G $ with $\mid $$a$$ \mid = u$ and $\mid $$b$$ \mid =
v.$ If, in addition,  for any $i, j,  k\in G$,  $a_{ij} \not=0$
implies $a_{ij}A_{jk} \not=0$, then $osp_{q, t} (GM_G(A), M, F))$ is
a braided m-Lie algebra  iff $r$ is a skew symmetric bicharacter.
\end {Lemma}
{\bf Proof.}  Obviously,  $osp_{q, t} (GM_G(A), M, F)_u$ is a
subspace.  It remains to check that $osp _{q, t} (GM_G(A), M, F)$ is
closed under bracket operation. For any $a \in  osp_{q, t} (GM_G(A),
M, F)_u$, $b \in osp _{q, t} (GM_G(A), M, F)_v$ and  $ u, v, g, h
\in G,$ set $w = u+v.$ See that
\begin {eqnarray*}
t([a, b] _{w+g, g}) M _ {w+g, h}
&=& t((ab - r (v, u)ba) _{w + g, g})M _ {w+g, h} \\
&=&t(b_{v + g, g})t(a_ {w+g, g+v})M _ {w+g, h} \\
 &{\ \ \ }& -
r (v, u)t(a_{u + g, g})t(b_ {w+g, g+u})M _ {w+g, h}\\
&=& r (g +v, u)r(g, v)M _{g, h +w}b _{h+w, u+h}a _{u+h, h} \\
&{\ \ \ }&- r(v,u)r (g +u, v)r(g, u)M _{g, h +w}a _{h+w, v+h}b _{v+h, h}, \\
-  r(g,w) M_{g, w+h } [ a, b]_{w+h, h}
 &=& - r(g, w) M_{g, w +h} a _{w+h, v+h}b _{v+h, h}\\
& &- r(g, w) r (v , u) M _{g, w+h} b _{w+h, u +h} a _{u +h, h}.
\end {eqnarray*}
Thus $[a, b] \in osp _{q, t} (GM_G(A), M, F)_w $ iff $(r (u, v)r(v,
u) -1) M_{g, w+h}
 a_{w+h, v+h} b _{v+h, h} =0$ for any $g,  h \in G.$
 $\Box$

We now  consider   $sl_{q, f}( GM_G(A),  F)$  and $osp_{q, t}
(GM_G(A), M, F)$ when the  bicharacter $r$ is skew symmetric. In
this case, they become  Lie color algebras in $(^{FG} {\cal M},
C^r)$,  called special gm Lie color algebra and ortho-symplectic gm
Lie color algebra, respectively.

It is well-known that a Lie color algebra $(^{FG}{\cal M}, C^r)$
with finitely generated $G$ is related to a Lie super algebra by
\cite [Theorem 2] {Sc79}.  We now show  how the above gm Lie color
algebras  are related to  Lie superalgebras for any abelian group
$G$.

Let $G$ be an abelian group with a skew symmetric bicharacter $r $.
Set $G_ {\bar 0} = \{g \in G \mid r (g, g )=1 \}$ and $G_ {\bar 1} =
\{g \in G \mid r (g, g )=-1 \}$. We define a new bicharacter: $r_0
(g, h) = -1$ for $g, h\in G_{\bar 1}$ and $r_0 (g, h) =1$ otherwise.
It is clear that $r_0$ is a bicharacter too.
 Obviously, $(r_0)_0 = r_0$ for any skew symmetric bicharacter $r$ on $G$.
 For convenience,
let $L^s $ denote the Lie superalgebra $ L= L_{\bar 0} \oplus
L_{\bar 1}$ for Lie color algebra $L$ in $(^{FG}{\cal M}, C^r)$,
 where $L_{\bar 0} = \oplus _{ i \in G_{\bar 0}} L_i$ and
 $L_{\bar 1} = \oplus _{ i \in G_{\bar 1}} L_i$ with
 $[ x, y ] = xy - r_0(\mid $$y$$ \mid  , \mid $$x$$ \mid )yx $ for any
$x \in L_g, y\in L_h, g, h \in G$ (see \cite [Page 718] {Sc79}).

Let $A = \sum \{A_{ij} \mid i, j \in G\} = GM_G(A)$ be a generalized
matrix algebra. Set $B_{\bar i , \bar j} = \sum _{g \in G_{\bar i},
h \in G_{\bar j}} A_{gh}$ for any $\bar i , \bar j \in {\bf Z}_2$
and $B = \sum \{B_{\bar i , \bar j } \mid  \bar i , \bar j \in {\bf
Z}_2\}$. We denote  the generalized matrix algebra $GM_{{\bf
Z}_2}(B)$ by $ GM_{{\bf Z}_2}(A^s) $.  For $a\in GM_G(A),$  if
$b_{\bar i , \bar j} = \sum _{g \in G_{\bar i}, h\in G_{\bar j} }
a_{gh}$ for any $\bar i , \bar j \in {\bf Z}_2$, then we denote the
element $b$ by $a^s$. When $G= {\bf Z}_2$ with $r (g,h) = (-1)
^{gh}$ for any $g, h \in {\bf Z}_2,$ we denote $tr_{q, f}$ by  $
tr_{s, f}$ and $osp _{q, t}$ by  $ osp _{s, t}.$ We have
\begin {Theorem}\label {4'.2.6}
% Let  $sl(\{n_g\}, F)=: \{ A \in gl (\{n_g\}, F) \mid tr_q(A)=0\}.$

(i)  $sl_{q, f}( GM_G(A),  F)^s  = sl _{s, f} (GM_{{\bf Z}_2}(A^s),
F).$

(ii) $osp_{q,t}(GM_G(A), M, F) ^s = osp _{s, t} (GM_{{\bf
Z}_2}(A^s), M^s, F).$
 \end {Theorem}
{\bf Proof.} (i) For any $a \in  GM_G(A)$, see that
 \begin{eqnarray*}
tr _{q, f} (a) &=& \sum _{g \in G_{\bar 0}} r_0 (g, g) tr_{q,f}
(a_{gg}) +
 \sum _{g \in G_{\bar 1}} r_0 (g, g) tr_{q,f} (a_{gg})\\
&=& tr_{s, f} (a).
\end {eqnarray*}
This completes the proof of (i).

(ii) For any $a \in (osp _{q, t} (GM_G(A), M, F)^s )_{\bar 0}$ with
$a \in osp _{q, t} (GM_G(A), M, F)_u$ and $u \in G_{\bar 0},$ let
$\bar i , \bar j \in {\bf Z}_{2}$ and  we see that
\begin{eqnarray*}
t(a _{ \bar i , \bar  i}) M_{\bar  i , \bar  j}
&=& \sum_{g, h \in G_ {\bar i}, k \in G_{\bar j}} t(a _{gh}) M _{gk} \\
&=& - \sum_{ h \in G_ {\bar i}, k \in G_{\bar j}} r_0
(h, u) M _{h, k +u} a_{k+u, k} \\
&=& - M_{\bar i, \bar j } a_{\bar j , \bar j }.
\end{eqnarray*}
\noindent This shows $a \in osp _{s, t} (GM_{{\bf Z}_2}(A^s), M^s,
F)$. We can similarly prove the others.  $\Box$

In fact, the relations \ \ \ (i) and (ii)\ \ \  above define a
$G$-grading of Lie superalgebras \ \ \ \ \ \ \ \ \ \ $sl _{s, f}
(GM_{{\bf Z}_2}(A^s), F)$ and $ osp _{s, t} (GM_{{\bf Z}_2}(A^s),
M^s, F)$, respectively.

We may apply the above results to $gl (\{n _g\}, F)$ with the
ordinary quantum trace $tr_q (a) = \sum _{g\in G} r(g,g) tr(a_{gg})$
and the ordinary transpose  $t (a) = a'.$ In this case, $sl _{q, f}$
and $osp _{q, f}$ are denoted by $sl _q $ and $osp _q,$
respectively. We have
\begin {Corollary}\label {4'.2.7}
(i)  $(sl_{q}( gl (\{V_g\}, F))^s = sl _{s} (V_{\bar 0}, V_{\bar 1},
F).$

(ii) $ (osp_q gl(\{V_g\}, F))^s = osp _s (V _ {\bar 0 }, V_{\bar 1},
M^s, F).$
 \end {Corollary}

\section {Braided m-Lie Algebras in the Yetter-Drinfeld Category and
Their Representations } In this section, we  give another class of
braided m-Lie algebras. We shall show that representations of an
algebra $A$ associated to a braided m-Lie algebra $L$ are also
representations of $L$. Furthermore, we show that if $(M, \psi)$ is
a faithful representation of $L$, then
 representations of $End _F M$ are also representations of  $L$.

The category is the Yetter-Drinfeld category $(^B_B{\cal YD}, C)$,
where $B$ is a finite dimensional Hopf algebra and $C$ is a braiding
with $C(x,y) = \sum (x_{(-1)}\cdot y)\otimes x_{(0)}$ for any $x\in
M, y\in N$.

We use the Sweedler's notation for  coproducts and comodules, i.e.
$$\Delta (a) = \sum _a a_1 \otimes a_2  \hbox { \ \ \ and  \ \ \ \ }
\psi (x) = \sum _x x_{(-1)} \otimes x_{(0)}$$ when $a\in H$ a
coalgebra and $x\in M$ a left $H$-comodule.

\begin {Lemma} \label {4'.3.1}

(i) If $(V, \alpha_V, \phi _V)$ and $(W, \alpha_W, \phi _W)$ are two
Yetter-Drinfeld modules over $B$ with dim $B<\infty $, then
 $Hom _F (V, W)$ is a Yetter-Drinfeld module under the following module
 operation and comodule operation:
 $(b \cdot f)(x) = \sum b_1\cdot f(S(b_2)\cdot x )$ and
$\phi (f) = (S^{-1} \otimes \hat \alpha ) (b_B \otimes f )$,
 where  $ \hat \alpha $ is defined by
$(b^* \cdot f)(x) = <b^*, x_{(-1)} S(f (x_{(0)})_{(-1)})>_{ev} (f
(x_{(0)}))_{(0)}$ for any $x\in V, f \in Hom _F (V,W), b^* \in B^*.$
Here $b_B$ denotes a coevoluation and $<, >_{ev}$ an evoluation of
$B.$

(ii) If  $(M, \alpha_M, \phi _M)$ is  a Yetter-Drinfeld modules over
$B$, $End _F M$ is an algebra in $(^B_B {\cal YD}, C)$
 \end {Lemma}

{\bf Proof.}
 (i)  It is clear  that
$\sum f _{(-1)} f _{(0)}(x) = \sum (f (x_{(0)}))_{(-1)}
S^{-1}(x_{(-1)}) \otimes (f (x_{(0)}))_{(0)} $ for any $x\in V, f
\in Hom _F (V, W), b \in B.$ Using this, we can show that $Hom _F
(V, W)$ is a $B$-comodule. Similarly, we can show that $Hom _F (V,
W)$ is a $B$-module. We now show
\begin {eqnarray} \label {e4'.3.1}
\sum (b \cdot f )_{(-1)} \otimes (b \cdot f )_{(0)} = \sum b_1 f
_{(-1)}S(b_3) \otimes b_2 \cdot f _{(0)}
\end {eqnarray}
for any $f \in Hom _F(V, W), b \in B$. For any $x\in V,$  see that
\begin {eqnarray*}
 \sum (b \cdot f )_{(-1)} \otimes (b \cdot f )_{(0)}(x)&=&
 \sum b_1 (f (S(b_4))\cdot x_{(0)})_{(-1)} S(b_3)S^{-1}(x_{(-1)})\\
&{\ }& \otimes b_2 \cdot (f (S(b_4))\cdot x_{(0})_{(0)}, \\
b_1 f _{(-1)}S(b_3) \otimes (b_2 \cdot f _{(0)})(x)&=& b_1 f_{(-1)}
S(b_4)
\otimes  b_2 \cdot f_{(0)} ((S(b_3)\cdot x))\\
&=& \sum b _1 (f (S(b_4) \cdot x _{(0)}))_{(-1)}
S(b_3)S^{-1}(x_{(-1)})
b_5S(b_6) \\
&{\ }& \otimes b_2 \cdot (f (S(b_4) \cdot x _{(0)}))_{(0)} \\
&=& \sum b_1 (f (S(b_4))\cdot x_{(0)})_{(-1)} S(b_3)S^{-1}(x_{(-1)})\\
&{\ }& \otimes b_2 \cdot (f (S(b_4))\cdot x_{(0)})_{(0)}.
\end {eqnarray*}
Thus (\ref {e4'.3.1}) holds and $Hom _F (V, W)$ is a Yetter-Drinfeld
module.

(ii) Let $E = End _F M$ and $m$ denote the multiplication of $E$.
 Now we show that $m$ is a homomorphism of  $B$-comodules.
It is sufficient to show
\begin {eqnarray} \label {e4'.3.11}
\sum _{fg}(fg )_{(-1)} \otimes (fg )_{(0)} = \sum _{f, g} f_{(-1)}
g_{(-1)} \otimes f_{(0)}g_{(0)}.
\end {eqnarray}
for any $f, g \in E$. Indeed, for any $x\in M,$  see that
 \begin {eqnarray*}
\sum _{fg}(fg )_{(-1)} \otimes (fg )_{(0)}(x)  &=& \sum _{x } (fg
(x_{0}))_{(-1)} S^{-1}(x_{(-1)})
\otimes (fg (x_{0}))_{(0)},\\
\sum _{f, g} f_{(-1)} g_{(-1)} \otimes f_{(0)}g_{(0)} (x)
 &=& \sum _{f, x} f_{(-1)} (g (x_{(0)}))_{(-1)} S^{-1} (x_{(-1)})
 \otimes f_{(0)} ((g (x_{(0)}))_{(0)}) \\
 &=& \sum _{f, x} (f((g (x_{(0)}))_{(0)}))_{(-1)} \\
&{\ }& S^{-1}((g (x_{(0)}))_{(0)}))_{(-1)}{}_2)
(g (x_{(0)}))_{(-1)} {}_1 S^{-1} (x_{(-1)}) \\
&{\ }& \otimes  (f((g (x_{(0)}))_{(0)}))_{(0)} \\
&=&   (fg(x_{(0)}))_{(-1)} S^{-1}(x_{(-1)})\otimes
(fg(x_{(0)}))_{(0)}.
\end {eqnarray*}
Thus (\ref {e4'.3.11}) holds. Similarly, we can show that $m$ is a
homomorphism of  $B$-modules.
 $\Box$

\begin {Example} \label {4'.3.2}
Let  $(M, \alpha_M, \phi _M)$ be  a Yetter-Drinfeld modules over $B$
and $L$ a subobject of  $E = End _F M$. If $L$ is  closed under
operation  $[\ \ ] = m - m C_{L, L},$ then $(L, [\  \ ] )$ is a
braided  m-Lie  subalgebra of $E^-$.
\end {Example}

By Lemma \ref {4'.3.1}, we may define representations of braided
m-Lie algebras in the Yetter-Drinfeld category.

\begin {Definition} \label {4'.3.3}
Let $(L, [\ \  ])$ be a braided  m-Lie algebra in $(^B_B{\cal YD},
C)$. If $M$ is an object in $(^B_B{\cal YD}, C)$ with  morphism
$\psi : L \rightarrow End _F M$ such that $\psi $ is a homomorphism
of braided m-Lie algebras, i.e.
 $\psi [\  \ ] = m (\psi \otimes \psi ) -
m (\psi \otimes \psi )  C_{L, L},$ then $(M, \psi)$ is called a
representation of $(L, [\ \  ])$.
\end {Definition}
Obviously, $(M, \psi )$ is a  representation of $L$  iff  $(M,
\alpha )$ is an $L$-module (see definition \ref {4'.1.7}), where the
relation between two operations  is $\alpha (a,  x)  = \psi (a)(x)$
for any $a\in L, x\in M$.

\begin {Proposition} \label {4'.3.4}
Let $(L, [\ \  ])$ be an object in  $(^B_B{\cal YD}, C)$ with a
morphism $[\ \ ] : L \otimes L \rightarrow L$. If $M$ is an object
in $(^B_B{\cal YD}, C)$ with monomorphism
 $\psi : L \rightarrow End _F M$
such that $\psi [\  \ ] = m (\psi \otimes \psi ) - m (\psi \otimes
\psi )  C_{L, L},$ then  $L$ is a braided m-Lie algebra and $(M,
\psi)$ is a faithful representation of $(L, [\ \  ])$.
\end {Proposition}
{\bf Proof.} By Lemma \ref {4'.3.1} $E= End _F M$ is an algebra in
$(^B_B {\cal YD}, C)$. By Definition \ref {4'.1.1}, $(L, [\ \ ])$ is
a braided m-Lie algebra. $\Box$

\begin {Proposition} \label {4'.3.5}
 Let $(L, [ \ \ ])$ be  a braided m-Lie algebra in $(_B^B {\cal YD}, C)$
 induced by multiplication of $A$ through $\phi .$

(i) If $( M , \psi )$  is a representation of algebra $A$, then $(M,
\psi \phi)$ is a representation of $L.$

(ii) If $( M , \psi )$  is a representation of braided m-Lie algebra
$A{}^-$, then $(M, \psi \phi)$ is a representation of $L.$

(iii) If   $( M , \psi )$  is  a faithful representation of $L$ and
$(N, \varphi )$ is a representation of algebra $End_F M$, then
 $(N, \varphi \psi  )$   is a representation of  $L$.
\end {Proposition}

{\bf Proof.} (i) and (ii) follow from  Definition \ref {4'.3.3}.

 (iii) By Lemma \ref {4'.3.1}, $E = End_F M$ is
an algebra in $(^B_B{\cal YD}, C)$. Thus $L$ is a braided m-Lie
algebra induced by $E$ through $\psi$. By Proposition \ref {4'.3.5},
we complete the proof. $\Box$

Let $(L, [\ \ ])$ be  a braided  m-Lie subalgebra of the path
algebra $(F(D, \rho ))^-$ with relations,  then a  representation
$(V, f )$ of $D$ with $f_\sigma =0$ for any $\sigma \in \rho$ is
also a representation of $L$ (see \cite [Proposition
II.1.7]{ARS95}).

\part {Hopf Algebras Living in Symmetric Tensor Categories}

\chapter{  The Quasitriangular Hopf Algebras
 in Symmetric Braided Tensor Categories }\label {c6}
In this chapter, we study the structures of Hopf algebras living in
  a symmetric braided tensor category  $({\cal C},C)$. We obtain that
$(H,R)$ is a quasitriangular bialgebra  living in ${\cal C}$ iff
$({}_H {\cal M}, C^R)$ is a braided tensor category. We show that
the antipode of (co)quasitriangular Hopf algebra living in $ {\cal
C}$ is invertible. Next we also prove that
 $S^2$ is inner   when almost cocommutative Hopf
 algebra $(H,R)$ living in ${\cal C}$ has an invertible antipode $S$.
    In particular, we structure the Drinfeld (co)double
 $D(H)$ in symmetric braided tensor category with left duality and we prove
 that it
 is (co)quasitriangular.

        We give some basic concepts as follows:

 Let  $H$ and $A$ be two bialgebras in a symmetric  braided tensor categories
 $({\cal C},C)$, and for any  $U, V, W \in ob {\cal C},$
 assume that

\begin {eqnarray*}
 R :  I \rightarrow H \otimes  H                            &,& \hbox { \ \ \ \ }
r : H \otimes H \rightarrow I , \hbox { \ \ \ \ }  \\
\tau : H \otimes A \rightarrow I &,& \hbox { \ \ \ \ }
\sigma : H \otimes H \rightarrow I,    \\
P :  I \rightarrow H \otimes  A   &,& \hbox { \ \ \ \ } Q : I
\rightarrow H \otimes H.\\
\end {eqnarray*}
where  $r$  and $R$ are invertible under convolution, and they are
morphisms in ${\cal C}.$

A bialgebra $(H, m, \eta, \Delta, \epsilon )$ with
convolution-invertible $R$ in $Hom _{\cal C} (I, H\otimes H)$ is
called a quasitriangular bialgebra
 living in
braided tensor category  ${\cal C}$ if the following conditions
hold\\
(QT1):
\[
\begin{tangle}
\ro R\\
\id\step\cd\\
\id\step\id\step[2]\id\\
\object{H}\step\object{H}\step[2]\object{H}
\end{tangle}
\step=\step
\begin{tangle}
\Ro R\\
\id\step\ro R\step\id\\
\hcu\step[2]\id\step\id\\
\step[0.5]\object{H}\step[2]\object{H}\step[2]\object{H}
\end{tangle}\ \ ;
\]

(QT2):
\[
\begin{tangle}
\step\ro R\\
\cd\step\id\\
\object{H}\step[2]\object{H}\step\object{H}
\end{tangle}
\step=\step
\begin{tangle}
\ro R\step\ro R\\
\id\step[2]\X\step[2]\id\\
\id\step[2]\id\step\cu\\
\object{H}\step[2]\object{H}\step[2]\object{H}
\end{tangle}\ \ ;
\]

(QT3):
\[
\begin{tangle}
\step[4]\object{H}\\
\ro R\step\cd\\
\id\step[2]\X\step[2]\id\\
\cu\step\cu\\
\step\object{H}\step[3]\object{H}
\end{tangle}
\step=\step
\begin{tangle}
\step\object{H}\\

\cd \step [2] \\

\XX\step [2]\\

\id \step [2] \id \step\ro R\\
\id\step[2]\X\step[2]\id\\
\cu\step\cu\\
\step\object{H}\step[3]\object{H}
\end{tangle}\ \ \ \ .
\]

A bialgebra $(H, m, \eta, \Delta, \epsilon )$ with
convolution-invertible $r$ in $Hom _{\cal C} ( H\otimes H, I)$ is
called a coquasitriangular bialgebra
 living in
braided tensor category  ${\cal C}$ if the following conditions
hold\\
 (CQT1):
\[
\begin{tangle}
\object{H}\step\object{H}\step[2]\object{H}\\
\id\step\id\step[2]\id\\
\id\step\cu\\
\coro r\\
\end{tangle}
\step=\step
\begin{tangle}
\step[0.5]\object{H}\step[2.5]\object{H}\step\object{H}\\
\hcd\step[2]\id\step\id\\
\id\step\coro r \step\id\\
\coRo r
\end{tangle}\ \ ;
\]

(CQT2):
\[
\begin{tangle}
\object{H}\step[2]\object{H}\step[1]\object{H}\\
\cu \step[1]\id\\

\step \coro r\\
\end{tangle}
\step=\step
\begin{tangle}
\object{H}\step[2]\object{H}\step[2]\object{H}\\
\id\step[2]\id\step\cd\\
\id\step[2]\X\step[2]\id\\
\coro  r \step\coro  r
\end{tangle}\ \ ;
\]
(CQT3):
\[
\begin{tangle}
\step\object{H}\step[3]\object{H}\\
\cd\step\cd\\
\id\step[2]\X\step[2]\id\\
\coro r\step\cu\\
\step[4]\object{H}

\end{tangle}
\step=\step
\begin{tangle}
\step\object{H}\step[3]\object{H}\\
\cd\step\cd\\
\id\step[2]\X\step[2]\id\\
\id \step[2]\id \step [1] \coro r\\
\XX \step [2]\\
\cu\\
 \step\object{H}
\end{tangle} \ \ \ .
\]
$\tau $ is called a skew pairing on $H\otimes A$ if
 the following
conditions are satisfied:

 (SP1)：
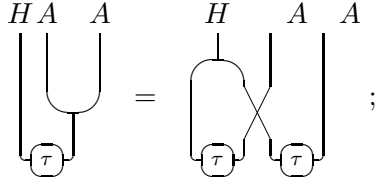
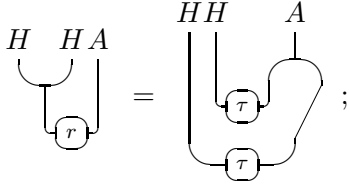
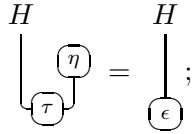
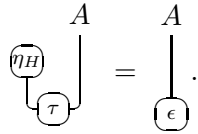
\[
\begin{tangle}
\object{H}\step\object{A}\step[2]\object{A}\\
\id\step\id\step[2]\id\\
\id\step\cu\\
\coro \tau\\
\end{tangle}
\step=\step
\begin{tangle}
\step\object{H}\step[3]\object{A}\step[2]\object{A}\\
\cd\step\id\step[2]\id\\
\id\step[2]\X\step[2]\id\\
\coro \tau\step\coro\tau
\end{tangle}\ \ ;
\]
(SP2):
\[
\begin{tangle}
\object{H}\step[2]\object{H}\step[1]\object{A}\\
\cu \step [1]\id \\

\step \coro r\\
\end{tangle}
\step=\step
\begin{tangle}
\object{H}\step[1]\object{H}\step[3]\object{A}\\
\id \step[1]\id\step[2]\cd\\
\id\step\coro \tau \step\ne1\\
\coRo \tau
\end{tangle}\ \ ;
\]
(SP3)：
\[
\begin{tangle}
\object{H}\\
\id\step[2]\Q \eta\\
\coro \tau
\end{tangle}
\step=\step
\begin{tangle}
\object{H}\\
\id\\
\QQ \epsilon
\end{tangle}\ \ ;
\]
(SP4)：
\[
\begin{tangle}
\step[2]\object{A}\\
\Q {\eta_H}\step[2]\id\\
\coro \tau
\end{tangle}
\step=\step
\begin{tangle}
\object{A}\\
\id\\
\QQ \epsilon
\end{tangle}\ \ .
\]

$P $ is called a skew copairing     of $H\otimes A$ if the following
conditions are satisfied:

(CSP1):
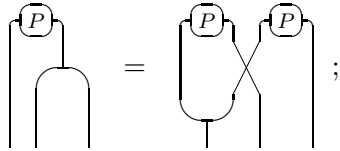
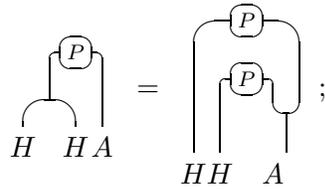
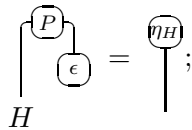
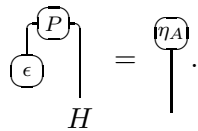
\[
\begin{tangle}
\ro P\\
\id\step\cd\\
\id\step\id\step[2]\id
\end{tangle}
\step=\step
\begin{tangle}
\ro P\step\ro P\\
\id\step[2]\X\step[2]\id\\
\cu\step\id\step[2]\id
\end{tangle}\ \ ;
\]

(CSP2):
\[
\begin{tangle}
\step\ro P\\
\cd\step\id\\
\object{H}\step[2]\object{H}\step\object{A}
\end{tangle}
\step=\step
\begin{tangle}
\Ro P\\
\id\step\ro P\step\id\\
\id\step\id\step[2]\hcu\\
\object{H}\step\object{H}\step[2]\object{A}
\end{tangle}\ \ ;
\]

(CSP3):
\[
\begin{tangle}
\ro P\\
\id\step[2]\QQ \epsilon \\
\object{H}
\end{tangle}
\step=\step
\begin{tangle}
\Q {\eta_H}\\
\id\\
\end{tangle}\ \ ;
\]
(CSP 4):
\[
\begin{tangle}
\ro P\\
\QQ \epsilon \step[2]\id\\
\step[2]\object{H}
\end{tangle}
\step=\step
\begin{tangle}
\Q {\eta_A}\\
\id
\end{tangle}\ \ .
\]

$\tau $ is called a pairing     on $H\otimes A$ if the following
conditions are satisfied:

(P1):
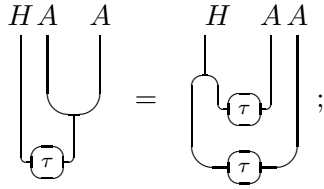
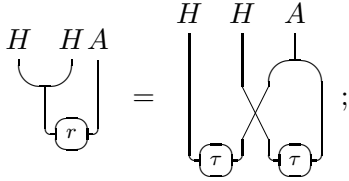
\[
\begin{tangle}
\object{H}\step\object{A}\step[2]\object{A}\\
\id\step\id\step[2]\id\\
\id\step\cu\\
\coro \tau\\
\end{tangle}
\step=\step
\begin{tangle}
\step\object{H}\step[2]\object{A}\step\object{A}\\
\hcd\step[2]\id\step\id\\
\id\step\coro \tau\step\id\\
\coRo \tau
\end{tangle}\ \ ;
\]

(P2):
\[
\begin{tangle}
\object{H}\step[2]\object{H}\step[1]\object{A}\\
\cu \step[1]\id\\

\step \coro r\\
\end{tangle}
\step=\step
\begin{tangle}
\object{H}\step[2]\object{H}\step[2]\object{A}\\
\id\step[2]\id\step\cd\\
\id\step[2]\X\step[2]\id\\
\coro \tau\step\coro \tau
\end{tangle}\ \ ;
\]

(P3) and (SP3) are the same, (P4) and (SP4) are the same.

 $P $ is called a copairing     of $H\otimes A$ if the following
conditions are satisfied: (CP1):
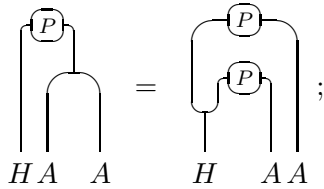
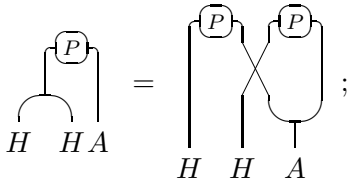
\[
\begin{tangle}
\ro P\\
\id\step\cd\\
\id\step\id\step[2]\id\\
\object{H}\step\object{A}\step[2]\object{A}
\end{tangle}
\step=\step
\begin{tangle}
\Ro P\\
\id\step\ro P\step\id\\
\hcu\step[2]\id\step\id\\
\step[0.5]\object{H}\step[2.5]\object{A}\step\object{A}
\end{tangle}\ \ ;
\]

(CP2):\ \ \
\[
\begin{tangle}
\step\ro P\\
\cd\step\id\\
\object{H}\step[2]\object{H}\step\object{A}
\end{tangle}
\step=\step
\begin{tangle}
\ro P\step\ro P\\
\id\step[2]\X\step[2]\id\\
\id\step[2]\id\step\cu\\
\object{H}\step[2]\object{H}\step[2]\object{A}
\end{tangle}\ \ ;
\]

(CP3) and (CSP4) are the same. (CP4) ane (CSP4) are the same.

$r$  is called almost commutative if $r$ satisfies 
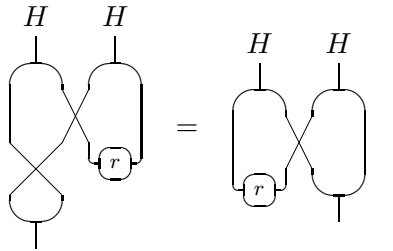
\[
\begin{tangle}
\step\object{H}\step[3]\object{H}\\
\cd\step\cd\\
\id\step[2]\X\step[2]\id\\
\XX\step\coro r\\
\cu
\end{tangle}
\step=\step
\begin{tangle}
\step\object{H}\step[3]\object{H}\\
\cd\step\cd\\
\id\step[2]\X\step[2]\id\\
\coro r\step\cu
\end{tangle}\ \ \ .
\]

$R$  is called almost cocommutative if $R$ satisfies
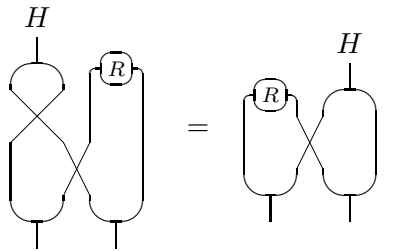
\[
\begin{tangle}
\step\object{H}\\
\cd\step\ro R\\
\XX\step\id\step[2]\id\\
\id\step[2]\X\step[2]\id\\
\cu\step\cu\
\end{tangle}
\step=\step
\begin{tangle}
\step[4]\object{H}\\
\ro R\step\cd\\
\id\step[2]\X\step[2]\id\\
\cu\step\cu
\end{tangle}\ \ \ .
\]

$r $ is called an universal co-$R$-matrix on $H \otimes H$ if $r$ is
invertible and  satisfies $(P1), (P2)$ and $(AC);$

 $R $ is called an universal $R$-matrix     of $H\otimes H$ if
   $R$  is invertible and  satisfies $(CP1), (CP2)$ and $(ACO);$

$\sigma  $ is called a 2-cocycle  on $H\otimes H$ if the following
conditions are satisfied:
\[ \hbox {2-}COC:
\begin{tangle}
\step\object{H}\step[3]\object{H}\step[2]\object{H}\\
\cd\step\cd\step\id\\
\id\step[2]\X\step[2]\id\step\id\\
\coro \sigma\step\cu\step\id\\
\step[4]\coro \sigma
\end{tangle}
\step=\step
\begin{tangle}
\object{H}\step[3]\object{H}\step[3]\object{H}\\
\id\step\cd\step\cd\\
\id\step\id\step[2]\X\step[2]\id\\
\d\coro \sigma\step\cu\\
\step\coRo \sigma
\end{tangle}\ \ \ .
\]

$\sigma  $ is called an anti-2-cocycle  on $H\otimes H$ if the
following conditions are satisfied:
\[A\hbox {2-}COC:
\begin{tangle}
\step\object{H}\step[3]\object{H}\step[2]\object{H}\\
\cd\step\cd\step\id\\
\id\step[2]\X\step[2]\id\step\id\\

\cu\step\coro {\sigma} \ne1 \\

\step[1]\coRo { \sigma }
\end{tangle}
\step=\step
\begin{tangle}
\object{H}\step[3]\object{H}\step[3]\object{H}\\
\id\step\cd\step\cd\\
\id\step\id\step[2]\X\step[2]\id\\
\id \step \cu\step\coro \sigma\\
\coro \sigma
\end{tangle}\ \ \ .
\]

$Q  $ is called a 2-cycle  of $H\otimes H$ if the following
conditions are satisfied:

\[(2\hbox {-}C):
\begin{tangle}
\step[4]\ro Q\\
\ro Q\step\cd\step\id\\
\id\step[2]\X\step[2]\id\step\id\\
\cu\step\cu\step\id\\
\step\object{H}\step[3]\object{H}\step[2]\object{H}
\end{tangle}
\step=\step
\begin{tangle}
\step \Ro Q\\

\ne1\ro Q\step\cd\\
\id\step\id\step[2]\X\step[2]\id\\
\id\step\cu\step\cu\\
\object{H}\step[2]\object{H}\step[3]\object{H}
\end{tangle}\ \ \ .
\]
$Q  $ is called an anti-2-cycle  of $H\otimes H$ if the following
conditions are satisfied:
\[A\hbox {2-}C:
\begin{tangle}
\step[1]\Ro Q\\

\cd\step\ro Q\nw1\\
\id\step[2]\X\step[2]\id\step\id\\

\cu\step\cu\step\id\\
\step\object{H}\step[3]\object{H}\step[2]\object{H}
\end{tangle}
\step=\step
\begin{tangle}
\ro Q\\
\id\step\cd\step\ro Q\\
\id\step\id\step[2]\X\step[2]\id\\
\id\step\cu\step\cu\\
\object{H}\step[2]\object{H}\step[3]\object{H}
\end{tangle}\ \ \ .
\]

Let $\tau ^{-1}$ and $P^{-1}$  denote the inverse of $\tau$ and $P$
under the convolution respectively. It is easy to check the
following:

(i) $\tau$ is a skew pairing iff $\tau ^{-1}$ is a pairing;
    $P$ is a skew copairing iff $P ^{-1}$ is a copairing;

(ii)  $\sigma$  is a 2-cocycle iff $\sigma ^{-1}$ is an
anti-2-cocycle;
     $Q$  is a 2-cycle iff $Q ^{-1}$ is an anti-2-cycle;

(iii) $\tau ^{-1} = \tau (S \otimes id_A ) (= \tau (id _H \otimes
S^{-1}))$ if $H$ is a Hopf algebra ( or $A$  is a Hopf algebra with
invertible antipode) and $\tau $ is a pairing;

(iv) $\tau ^{-1} =  \tau (id _H \otimes S (= \tau (S^{-1} \otimes
id_A ))$ if $A$ is a Hopf algebra ( or  $H$  is a Hopf algebra with
invertible antipode  )  and $\tau $ is a  skew pairing;

(v)  $P^{-1} =  (S \otimes id _A)P ( = (id_H \otimes S^{-1}) P ) $
if $H$ is a Hopf algebra ( or $A$ is a Hopf algebra with invertible
antipode)  and $P$  is a  copairing;

(vi) $P^{-1} = (id_H \otimes S)P (= S^{-1} \otimes id_A)P$ if $A$ is
a Hopf algebra ( or  $H$  is a Hopf algebra with invertible antipode
)  and $P $ is a  skew copairing.

In this chapter, we assume that the braided tensor category $({\cal
C},C)$
        is symmetric and all the algebras, coalgebras, bialgebras
        and Hopf algebras are
in the category $({\cal C},C),$ unless otherwise  stated, (see \cite
[section 10.5] {Mo93} or \cite [Definition 2.2] {Ma95a}). We omit
all of the proofs of the dual  conclusions since they can be
obtained by turning the diagrams in the proofs of the conclusions
upside down. Obviously, the category ${\cal V}ect(k)$   of ordinary
vector spaces is a symmetric
 braided tensor category, in which the braiding is the ordinary transposition.
 We can easily get the corresponding conclusions in any symmetric braided
 categories for most of the conclusions in ${\cal V}ect (k)$  by turning the ordinary
 transposition into a braiding. In this case, we shall omit the proofs of
 these conclusions.

\section {Yang-Baxter equation}\label {s11}

In this section,  we give  the relations among quantum Yang-Baxter
equation, Yang-Baxter equation, (co)quasitriangular Hopf algebra and
the braided tensor category.

\begin {Proposition} \label {5.1.1}

(i)  If $(H,R)$ is a quasitriangular bialgebra, then $R$ satisfies
the quantum Yang- Baxter equation:

\[ (QYBE): \ \
\begin{tangle}
\ro R \step [2]\ro R \step [2]\ro R \step [2]\\
\id \step[2] \XX \step [2] \XX \step [2] \id\\
\cu \step [2]\cu \step [2]\cu\\
 \step\object{H}\step[4]\object{H}\step[4]\object{H}
\end{tangle}
\step=\step
\begin{tangle}
\ro R \step [2]\ro R \step [2]\ro R \step [2]\\

\id \step[2] \XX \step [2] \XX \step [2] \id\\
 \XX \step [2] \XX \step [2] \XX\\

\id \step[2] \XX \step [2] \XX \step [2] \id\\
\cu \step [2]\cu \step [2]\cu\\
  \step\object{H}\step[4]\object{H}\step[4]\object{H}
\end{tangle} \ \ .
\]

(ii)  If $(H,r)$ is a coquasitriangular bialgebra, then $r$
satisfies the quantum co-Yang- Baxter equation:
\[ (CQYBE): \ \
\begin{tangle}
\step\object{H}\step[4]\object{H}\step[4]\object{H}\\

\cd \step [2]\cd\step [2]\cd\step [2]\\
\id \step[2] \XX \step [2] \XX \step [2] \id\\
\coro r \step [2]\coro r  \step [2]\coro r \\

\end{tangle}
\step=\step
\begin{tangle}
 \step\object{H}\step[4]\object{H}\step[4]\object{H}\\
 \cd\step [2]\cd\step [2]\cd\step [2]\\

\id \step[2] \XX \step [2] \XX \step [2] \id\\
 \XX \step [2] \XX \step [2] \XX\\

\id \step[2] \XX \step [2] \XX \step [2] \id\\
\coro r  \step [2]\coro r  \step [2]\coro r \\

\end{tangle} \ \ .
\]

          \end {Proposition}
{\bf Proof.}  (i) By turning the proofs of   \cite [Theorem VIII
2.4]{Ka95}
 into braid diagrams,
we obtain the proofs.

(ii) It is the dual case of part (i).
 \begin{picture}(8,8)\put(0,0){\line(0,1){8}}\put(8,8){\line(0,-1){8}}\put(0,0){\line(1,0){8}}\put(8,8){\line(-1,0){8}}\end{picture}
\begin {Proposition} \label {5.1.2}
(i) If $(H, m, \eta )$  is an algebra with two algebraic morphisms $
\Delta : H \longrightarrow  H \otimes H $  and $\epsilon  : H
\longrightarrow I $, then
  $H$ is a bialgebra iff          ${}_H {\cal M}$  is a
tensor category;

(ii) If $(H, \Delta, \epsilon )$  is a coalgebra with two
coalgebraic morphisms $ m : H \otimes H \longrightarrow  H  $  and
$\eta  : I \longrightarrow H $, then
  $H$ is a bialgebra iff          ${}^H {\cal M}$  is a
tensor category.
\end {Proposition}
{\bf Proof.} (i) By turning the proofs of
 \cite [Proposition XI.3.1]{Ka95} into braid diagrams, we obtain the
 proofs. We also can see  Corollary \ref {5.1.3''}.

         (ii) It is a dual case of part (i).\begin{picture}(8,8)\put(0,0){\line(0,1){8}}\put(8,8){\line(0,-1){8}}\put(0,0){\line(1,0){8}}\put(8,8){\line(-1,0){8}}\end{picture}

\begin {Theorem} \label {5.1.3}
(i)    $(H,R)$ is a quasitriangular bialgebra iff $({}_H {\cal M},
C^R )$ is a braided tensor category.

(ii) $(H,r)$ is a coquasitriangular bialgebra iff $({}^H {\cal M},
C^r )$ is a braided tensor category.
\end {Theorem}
{\bf Proof.} (i) The necessity can be shown by means of
 turning the proofs of
  \cite [Proposition VIII. 3.1, Proposition XIII. 1.4 ]{Ka95}  into braid diagrams.

The sufficiency. For left regular $H$-module $V$, using
 $(id _V \otimes C ^R_{V,V})(C^R_{V,V} \otimes id_V)
= C^R_{V\otimes V,V},$ we have that
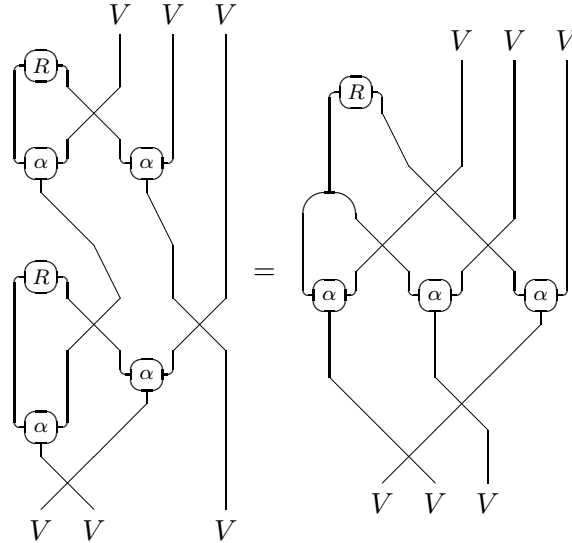
\[
\begin{tangle}
\step[4]\object{V}\step[2]\object{V}\step[2]\object{V}\\

\ro R \step [2] \id\step [2] \id\step [2] \id\\

\id \step [2] \XX\step [2]\id\step [2] \id\\

\tu \alpha \step [2]\tu \alpha \step [2] \id \\

\step [1] \nw2 \step [3] \nw1 \step [2] \id \\

\ro R \step [1] \nw1\step [2] \id\step [2] \id\\

\id \step [2] \XX \step [2] \XX \\

\id \step [2] \id \step [2]\tu \alpha \step [2] \id \\

\tu \alpha  \step [2] \ne2   \step [3] \id \\

\step \XX\step [2] \step [3]\id   \\
\step\object{V}\step[2]\object{V}\step[5]\object{V}\\

\end{tangle} \ \ = \ \
\begin{tangle}
\step[6]\object{V}\step[2]\object{V}\step[2]\object{V}\\

\step \ro R \step [3] \id\step [2] \id\step [2] \id\\

\step \id \step [2] \nw1 \step [2] \id\step [2] \id\step [2] \id\\

\cd \step [2] \XX\step [2]\id\step [2] \id\\

 \id  \step [2] \XX \step [2] \XX \step [2] \id \\

\tu \alpha \step [2]\tu \alpha \step [2]\tu \alpha \step [2]\\

 \step\id\step [4] \id \step [3] \ne2 \\

 \step\nw2\step [3] \XX \\

\step [3]\XX\step [2] \id\\

\step[3]\object{V}\step[2]\object{V}\step[2]\object{V}\\

\end{tangle} \ \ \ .
\]Applying $\eta _H \otimes \eta _H \otimes \eta _H$ on the above, we
obtain
 $(CP2).$ Similarly, we can get $(CP1)$ and $(ACO).$

(ii) It is the dual case of part (i).
\begin{picture}(8,8)\put(0,0){\line(0,1){8}}\put(8,8){\line(0,-1){8}}\put(0,0){\line(1,0){8}}\put(8,8){\line(-1,0){8}}\end{picture}

\section {The antipode of quasitriangular Hopf algebra}\label {s12}

In this section, we show that the antipode of (co)quasitriangular
Hopf algebra
 is invertible,
and $S^2$ is inner or coinner when Hopf algebra $(H,R)$ with
invertible antipode $S$ is almost cocommutative or $(H,r)$ with
invertible antipode $S$ is almost commutative.

\begin {Theorem} \label {5.2.1} Let $(H,R)$ be an almost cocommutative
 Hopf algebra with invertible antipode $S$
and define
$$u = mC_{H,H}(id_H \otimes S) R.$$

Then

(i)  $u^{-1} = m C_{H,H} (id _H \otimes S^{-1}) R^{-1};$

 (ii)  $S^2 = m (m \otimes id _H )( u  \otimes id_H  \otimes u^{-1}).$

\end {Theorem}
{\bf Proof.} By turning the proofs of
 \cite [Proposition 10.1.4]{Mo93} into braid diagrams, we obtain the proof.
 \begin{picture}(8,8)\put(0,0){\line(0,1){8}}\put(8,8){\line(0,-1){8}}\put(0,0){\line(1,0){8}}\put(8,8){\line(-1,0){8}}\end{picture}

\begin {Theorem} \label {5.2.2} Let $(H,r)$ be an almost commutative
 Hopf algebra with invertible antipode  $S$
and define
$$ \mu = r C_{H,H} (id _H \otimes S) \Delta _H \hbox { \ \ } . $$

Then

(i)  $\mu ^{-1} = r^{-1} C_{H,H} (id \otimes S^{-1}) \Delta _H  { \
\ } ;$

 (ii)  $S^2 = (\mu \otimes id _H \otimes \mu ^{-1} )(\Delta \otimes
 id _H) \Delta _H  { \ \ } .$

\end {Theorem}
{\bf Proof.} It is a dual case of Theorem \ref
{5.2.1}.\begin{picture}(8,8)\put(0,0){\line(0,1){8}}\put(8,8){\line(0,-1){8}}\put(0,0){\line(1,0){8}}\put(8,8){\line(-1,0){8}}\end{picture}

\begin {Theorem} \label {5.2.3} Let $(H,r)$ be a coquasitriangular Hopf algebra
and define
$$\lambda = r(id _H \otimes S) \Delta _H  \hbox { \ \ }. $$

Then

  (i)  $\lambda ^{-1} = r^{-1} (S \otimes id_H) \Delta _H  \hbox { \ \ } ;$

 (ii)  $S^2 = (\lambda ^{-1} \otimes id _H \otimes \lambda )(\Delta \otimes
 id _H) \Delta _H \hbox { \ \ };$

(iii)  $S^{-2} = (\lambda  \otimes id _H \otimes \lambda
^{-1})(\Delta \otimes
 id _H) \Delta _H \hbox { \ \ }.$

\end {Theorem}
{\bf Proof.}  By turning the proofs of    \cite [Theorem 1.3] {Do93}
 into braid diagrams,  we obtain the proofs.\begin{picture}(8,8)\put(0,0){\line(0,1){8}}\put(8,8){\line(0,-1){8}}\put(0,0){\line(1,0){8}}\put(8,8){\line(-1,0){8}}\end{picture}

\begin {Theorem} \label {5.2.4} Let $(H,R)$ be a quasitriangular Hopf algebra
and define
$$v = m(id _H \otimes S) R \hbox { \ \ }.$$

Then

(i)  $v^{-1} = m (S \otimes id_H) R^{-1} \hbox { \ \ };$

 (ii)  $S^2 = m (m \otimes id _H )( v^{-1}  \otimes id_H  \otimes v)
 \hbox { \ \ };$

 (iii)  $S^{-2} = m (m \otimes id _H )( v  \otimes id_H  \otimes v^{-1}) \hbox
 { \ \ }.$

\end {Theorem}
{\bf Proof.} It is a dual case of Theorem \ref {5.2.3}.
\begin{picture}(8,8)\put(0,0){\line(0,1){8}}\put(8,8){\line(0,-1){8}}\put(0,0){\line(1,0){8}}\put(8,8){\line(-1,0){8}}\end{picture}

\section {Drinfeld double}\label {s13}
In this section, we first construct $H^{\sigma}$  by a 2-cocycle
$\sigma$, as in \cite [Theorem 1.6] {Do93}. Next, we show that
Drinfeld (co)double is a (co)quasitriangular Hopf algebra.

\begin {Proposition}\label {5.3.1}
Let $H$ be a bialgebra with a 2-cocycle $\sigma$ on $H \otimes H.$
Define $H^{\sigma}$ as follows: $H^{\sigma}=H $   as coalgebra and
the multiplication

\[
\begin{tangle}
m_{H^\sigma}
\end{tangle}
\step=\step
\begin{tangle}
\step\object{H^\sigma}\step[5]\object{H^\sigma}\\
\cd\step[3]\cd\\
\id\step[2]\id\step[2]\dd\step[2]\id\\
\id\step[2]\XX\step[2]\cd\\
\coro \sigma\step\cd\step\id\step[2]\id\\
\step[3]\id\step[2]\X\step[2]\id\\
\step[3]\cu\step\coro {\bar \sigma } \\
\end{tangle} \ \ .
\]

(i)   $H^{\sigma}$ is a bialgebra;

(ii) If  $H$ is a Hopf algebra, then  $H^\sigma$ is also a Hopf
algebra with antipode

\[
\begin{tangle}
S_{H^\sigma}
\end{tangle}
\step=\step
\begin{tangle}
\step\object{H^\sigma}\\
\cd\\
\id\step[2]\d\\
\id\step[2]\cd\\
\id\step[2]\S\step\cd\\
\coro \sigma\step\S\step\cd\\
\step[4]\S\step[2]\id\\
\step[4]\coro {\bar \sigma }
\end{tangle}\ \ \ .
\]
(iii) If $H$ is a commutative algebra, then   $(H^{\sigma}, \hat
\sigma)$ is a symmetric coquasitriangular bialgebra, where
\[
\begin{tangle}
\hat{\sigma}
\end{tangle}
\step=\step
\begin{tangle}
\step\object{H^\sigma}\step[4]\object{H^\sigma}\\
\cd\step[2]\cd\\
\id\step[2]\XX\step[2]\id\\
\XX\step[2]\coro  {\bar \sigma }\\
\coro \sigma
\end{tangle}\ \ \ .
\]

(iv) If  $H$ is a commutative Hopf
 algebra, then  $S_{H^{\sigma}}$ is invertible.

          \end {Proposition}

{\bf Proof.}
 (i), (ii) and (iv)  can be obtained by turning the proofs  of
  \cite [Theorem 1.6]{Do93} into braid diagrams.
   We now show part (iii).

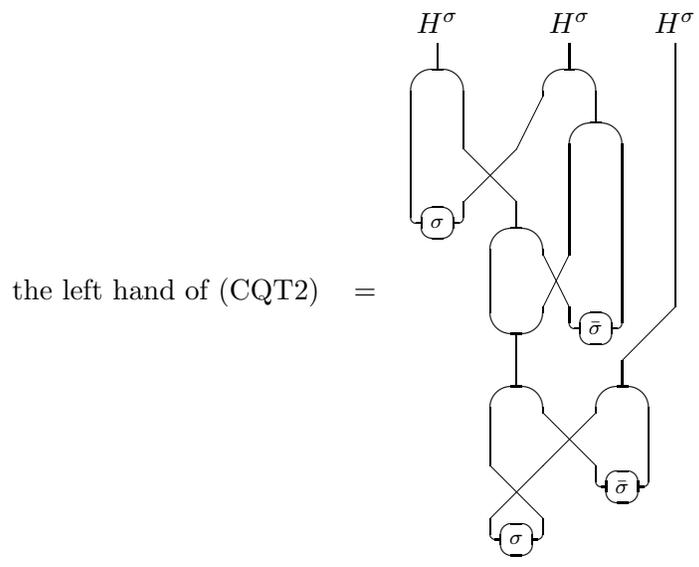
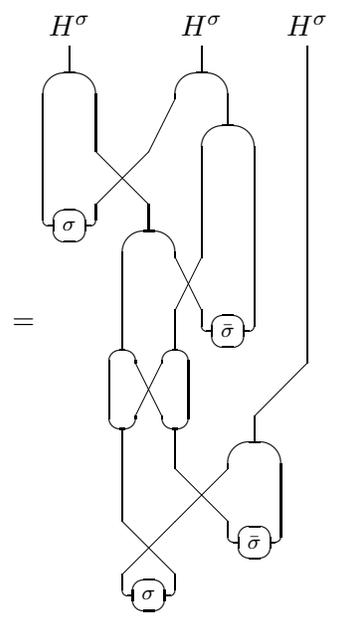
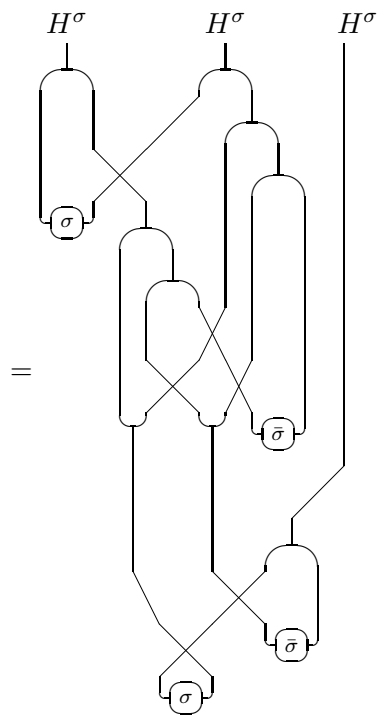
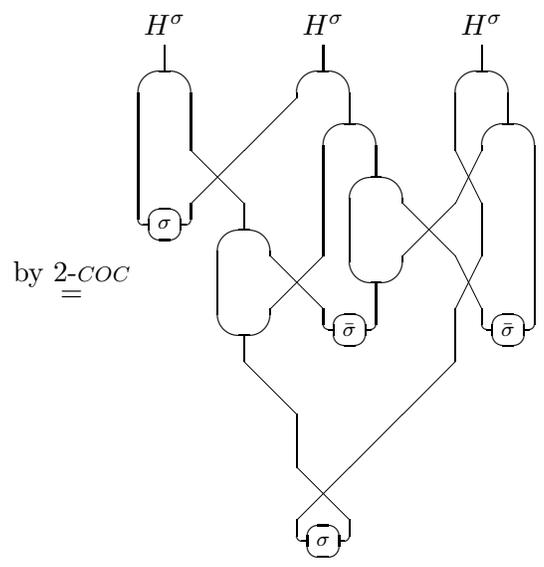
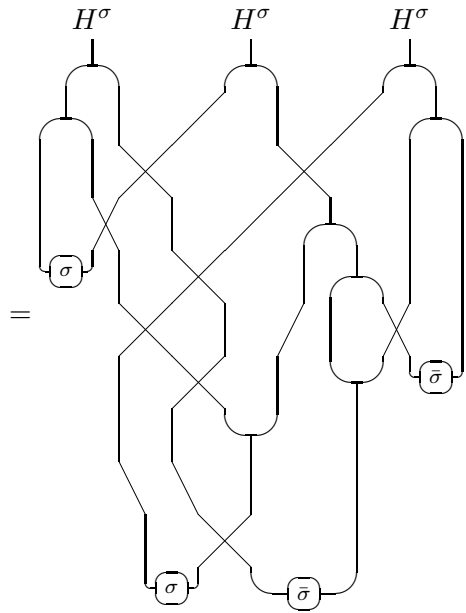
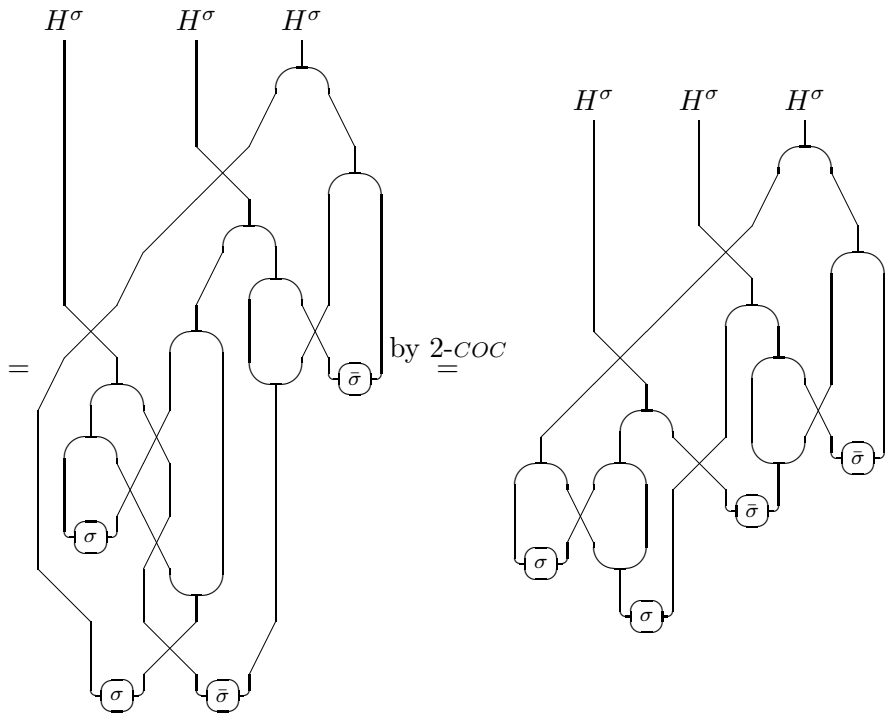
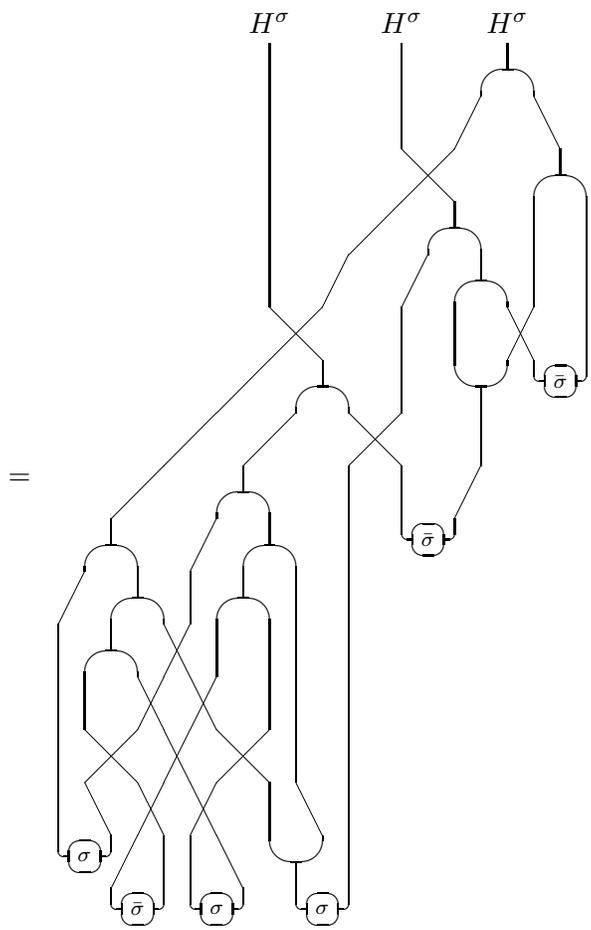
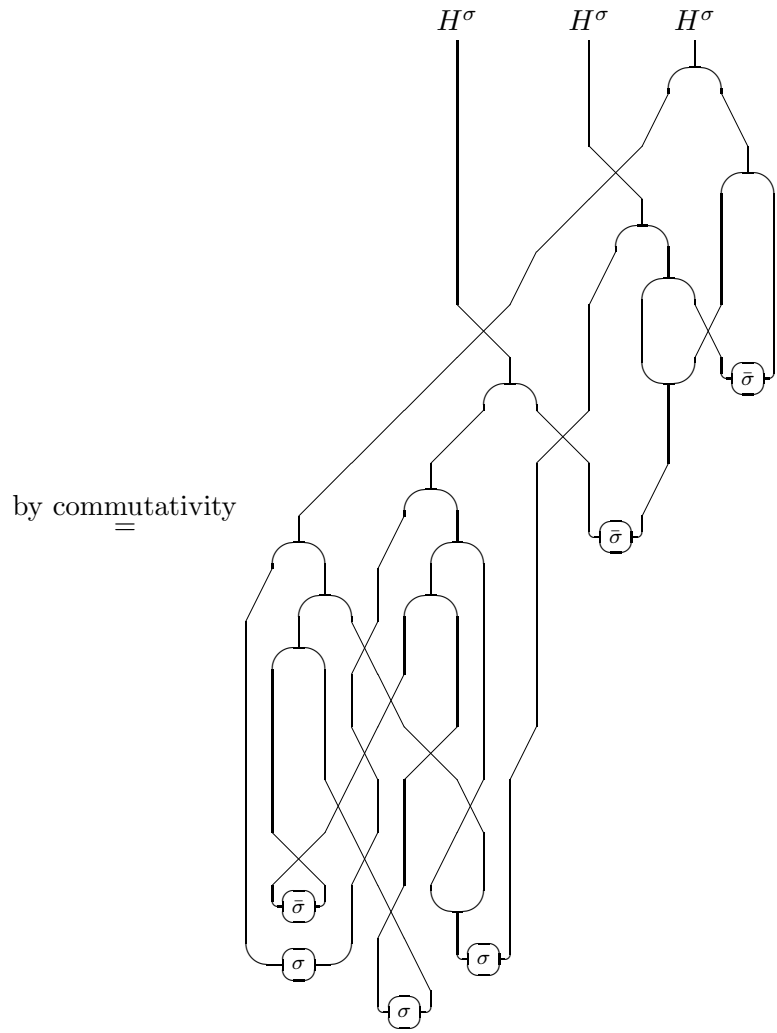
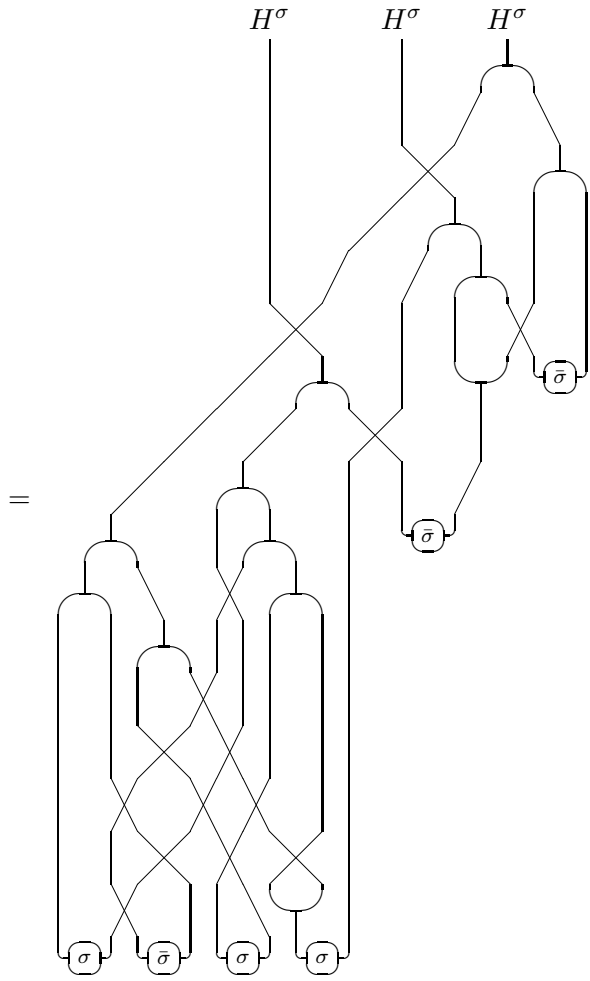
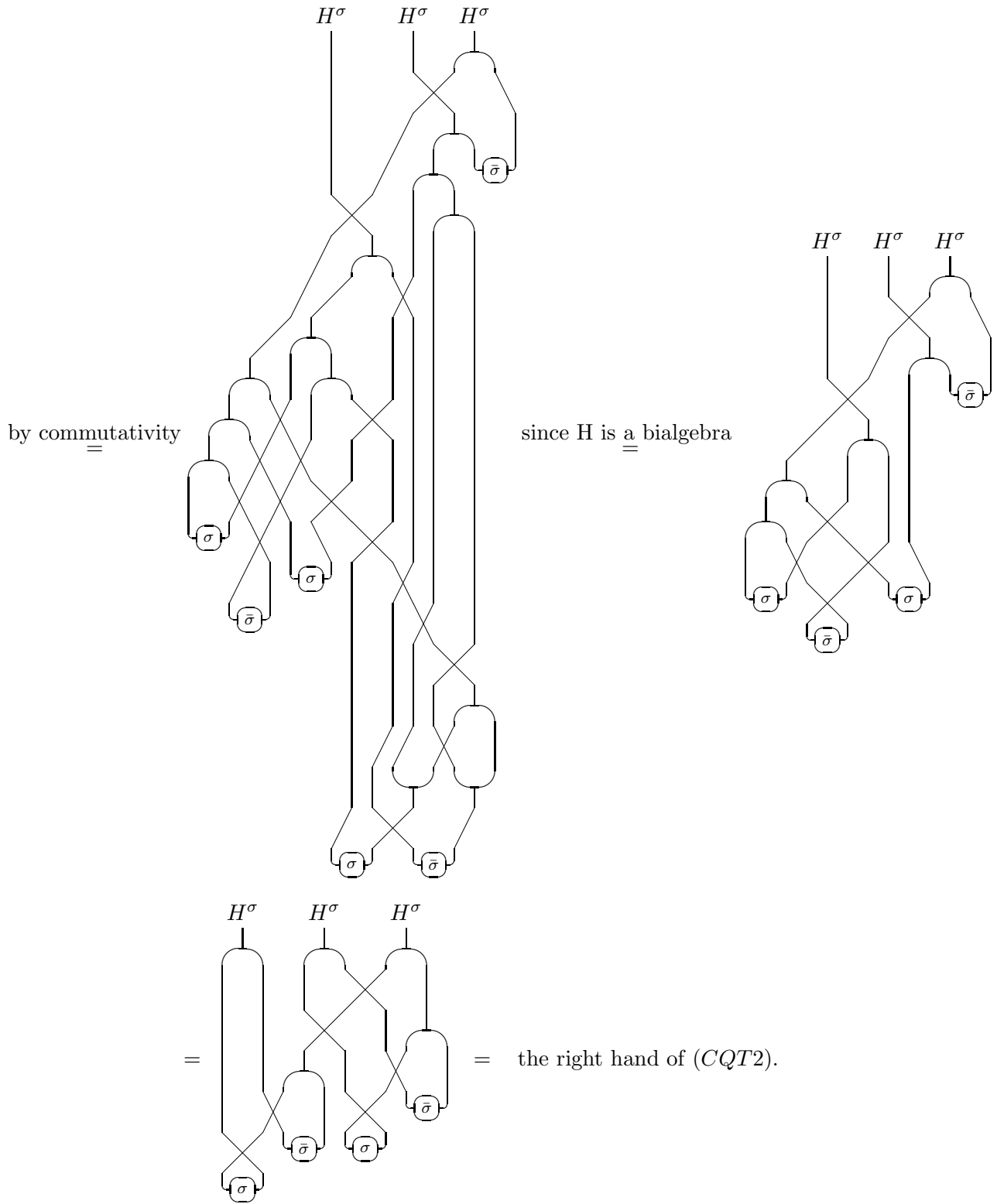
\[
\begin{tangle}
\hbox  { the left hand of (CQT2)}
\end{tangle}
\step=\step
\begin{tangle}
\step\object{H^\sigma}\step[5]\object{H^\sigma}\step[4]\object{H^\sigma}\\
\cd\step[3]\cd\step[3]\id\\
\id\step[2]\id\step[2]\dd\step\cd\step[2]\id\\
\id\step[2]\XX\step[2]\id\step[2]\id\step[2]\id\\
\coro \sigma\step\cd\step\id\step[2]\id\step[2]\id\\
\step[3]\id\step[2]\X\step[2]\id\step[2]\id\\
\step[3]\cu\step\coro {\bar  \sigma } \step  \ne2\\
\step[3]\cd\step[2]\cd\\
\step[3]\id\step[2]\XX\step[2]\id\\
\step[3]\XX\step[2]\coro {\bar  \sigma }\\
\step[3]\coro \sigma
\end{tangle}
\]

\[
=
\begin{tangle}
\step\object{H^\sigma}\step[5]\object{H^\sigma}\step[4]\object{H^\sigma}\\
\cd\step[3]\cd\step[3]\id\\
\id\step[2]\id\step[2]\dd\step\cd\step[2]\id\\
\id\step[2]\XX\step[2]\id\step[2]\id\step[2]\id\\
\coro \sigma\step\cd\step\id\step[2]\id\step[2]\id\\
\step[3]\id\step[2]\X\step[2]\id\step[2]\id\\
\step[2.5]\hcd\step\hcd\step[0.5]\coro  { \bar \sigma } \step\step\id\\
\step[2.5]\id\step\X\step\id\step[3.5]\ne2\\
\step[2.5]\hcu\step\hcu\step[1.5]\cd\\
\step[3]\id\step[2]\XX\step[2]\id\\
\step[3]\XX\step[2]\coro {\bar \sigma }\\
\step[3]\coro \sigma
\end{tangle}
\]

\[
=
\begin{tangle}
\step\object{H^\sigma}\step[6]\object{H^\sigma}\step[5]\object{H^\sigma}\\
\cd\step[4]\cd\step[3.5]\id\\
\id\step[2]\id\step[3]\ne2\step\cd\step[2.5]\id\\
\id\step[2]\XX\step[3]\id\step\cd\step[1.5]\id\\
\coro \sigma\step\cd\step[2]\id\step\id\step[2]\id\step[1.5]\id\\
\step[3]\id\step\cd\step\id\step\id\step[2]\id\step[1.5]\id\\
\step[3]\id\step\id\step[2]\X\step\id\step[2]\id\step[1.5]\id\\
\step[3]\id\step\XX\step\X\step[2]\id\step[1.5]\id\\
\step[3]\hcu\step[2]\hcu\step\coro  { \bar \sigma }\step\step[0.5]\id\\
\step[3.5]\id\step[3]\id\step[4]\ne2\\
\step[3.5]\id\step[3]\id\step[2]\cd\\
\step[3.5]\d\step[2]\XX\step[2]\id\\
\step[4.5]\XX\step[2]\coro {\bar  \sigma }\\
\step[4.5]\coro \sigma
\end{tangle}
\]
\[
\stackrel{ \hbox {by 2-}COC}{=}
\begin{tangle}
\step\object{H^\sigma}\step[6]\object{H^\sigma}\step[6]\object{H^\sigma}\\
\cd\step[4]\cd\step[4]\cd\\
\id\step[2]\id\step[3]\ne2\step\cd\step[3]\id\step\cd\\
\id\step[2]\XX\step[3]\id\step\cd\step[2]\X\step[2]\id\\
\coro \sigma\step\cd\step[2]\id\step\id\step[2]\XX\step\id\step[2]\id\\
\step[3]\id\step[2]\XX\step\cu\step[2]\X\step[2]\id\\
\step[3]\cu\step[2]\coro  { \bar \sigma }\step \step[2]\id\step\coro {\bar  \sigma }\\
\step[4]\nw2\step[6]\ne2\\
\step[6]\id\step[3]\ne2\\
\step[6]\XX\\
\step[6]\coro \sigma
\end{tangle}
\]

\[
 =
 \begin{tangle}
\step[2]\object{H^\sigma}\step[6]\object{H^\sigma}\step[6]\object{H^\sigma}\\
\step\cd\step[4]\cd\step[4]\cd\\
\cd\step\id\step[3]\ne2\step[2]\id\step[3]\ne2\step\cd\\
\id\step[2]\id\step\XX\step[4]\XX\step[3]\id\step[2]\id\\
\id\step[2]\X\step[2]\id\step[3]\ne2\step\cd\step[2]\id\step[2]\id\\
\coro \sigma\step\id\step[2]\XX\step[3]\id\step\cd\step\id\step[2]\id\\
\step[3]\XX\step[2]\id\step[2]\dd\step\id\step[2]\X\step[2]\id\\
\step[3]\id\step[2]\XX\step[2]\id\step[2]\cu\step\coro {\bar  \sigma }\\
\step[3]\id\step[2]\id\step[2]\cu\step[3]\id\\
\step[3]\d\step\d\step[2]\id\step[4]\id\\
\step[4]\id\step[2]\XX\step[4]\id\\
\step[4]\coro \sigma\step[2]\coRo {\bar \sigma } \\
\end{tangle}
\]

\[
=
 \begin{tangle}
\step\object{H^\sigma}\step[5]\object{H^\sigma}\step[4]\object{H^\sigma}\\
\step\id\step[5]\id\step[3]\cd\\
\step\id\step[5]\id\step[2]\dd\step[2]\d\\
\step\id\step[5]\XX\step[3]\cd\\
\step\id\step[4]\ne2\step\cd\step[2]\id\step[2]\id\\
\step\id\step[2]\dd\step[2]\dd\step\cd\step\id\step[2]\id\\
\step\XX\step[2]\cd\step\id\step[2]\X\step[2]\id\\
\dd\step\cd\step\id\step[2]\id\step\cu\step\coro { \bar  \sigma }\\
\id\step\cd\step\X\step[2]\id\step[2]\id\\
\id\step\id\step[2]\X\step\id\step[2]\id\step[2]\id\\
\id\step\coro \sigma\step\X\step[2]\id\step[2]\id\\
\nw2\step[3]\id\step\cu\step[2]\id\\
\step[2]\id\step[2]\XX\step[2]\dd\\
\step[2]\coro \sigma\step[2]\coro {\bar \sigma }
\end{tangle}
\stackrel {\hbox {by 2-}COC }{=}
 \begin{tangle}
\step[3]\object{H^\sigma}\step[4]\object{H^\sigma}\step[4]\object{H^\sigma}\\
\step[3]\id\step[4]\id\step[3]\cd\\
\step[3]\id\step[4]\id\step[2]\dd\step[2]\d\\
\step[3]\id\step[4]\XX\step[3]\cd\\
\step[3]\id\step[3]\ne2\step\cd\step[2]\id\step[2]\id\\
\step[3]\XX\step[3]\id\step\cd\step\id\step[2]\id\\
\step[2]\ne2\step\cd\step[2]\id\step\id\step[2]\X\step[2]\id\\
\cd\step\cd\step\XX\step\cu\step\coro {\bar \sigma }\\
\id\step[2]\X\step[2]\id\step\id\step[2]\coro {\bar \sigma }\\
\coro \sigma\step\cu\step\id\\
\step[4]\coro \sigma\\
\end{tangle}
\]

\[ \ \
 = \ \
 \begin{tangle}
 \step[8]\object{H^\sigma}\step[5]\object{H^\sigma}\step[4]\object{H^\sigma}\\
\step[8]\id\step[5]\id\step[3]\cd\\
\step[8]\id\step[5]\id\step[2]\dd\step[2]\d\\
\step[8]\id\step[5]\XX\step[3]\cd\\
\step[8]\id\step[4]\ne2\step\cd\step[2]\id\step[2]\id\\
\step[8]\id\step[2]\dd\step[2]\dd\step\cd\step\id\step[2]\id\\
\step[8]\XX\step[3]\id\step[2]\id\step[2]\X\step[2]\id\\
\step[7]\ne2\step\cd\step[2]\id\step[2]\cu\step\coro {\bar \sigma }\\
\step[5]\ne2\step[2]\ne2\step[2]\XX\step[3]\id\\
\step[3]\ne2\step[2]\cd\step[3]\id\step[2]\id\step[2]\dd\\
\step\cd\step[2]\dd\step\cd\step[2]\id\step[2]\coro {\bar \sigma }\\
\dd\step\cd\step\id\step\cd\step\id\step[2]\id\\
\id\step\cd\step\X\step\id\step[2]\id\step\id\step[2]\id\\
\id\step\id\step[2]\X\step\X\step[2]\id\step\id\step[2]\id\\
\id\step\XX\step\X\step\XX\step\id\step[2]\id\\
\id\step\d\step\X\step\X\step[2]\id\step\d\step\id\\
\coro \sigma\dd\step\id\step\id\step\d\step\cu\step\id\\
\step[2]\coro {\bar \sigma }\step \coro \sigma\step[2]\coro \sigma
\end{tangle}
\]
\[ \ \  \stackrel {\hbox  {by commutativity}}{=}
\begin{tangle}
\step[8]\object{H^\sigma}\step[5]\object{H^\sigma}\step[4]\object{H^\sigma}\\
\step[8]\id\step[5]\id\step[3]\cd\\
\step[8]\id\step[5]\id\step[2]\dd\step[2]\d\\
\step[8]\id\step[5]\XX\step[3]\cd\\
\step[8]\id\step[4]\ne2\step\cd\step[2]\id\step[2]\id\\
\step[8]\id\step[2]\dd\step[2]\dd\step\cd\step\id\step[2]\id\\
\step[8]\XX\step[3]\id\step[2]\id\step[2]\X\step[2]\id\\
\step[7]\ne2\step\cd\step[2]\id\step[2]\cu\step\coro {\bar \sigma }\\
\step[5]\ne2\step[2]\ne2\step[2]\XX\step[3]\id\\
\step[3]\ne2\step[2]\cd\step[3]\id\step[2]\id\step[2]\dd\\
\step\cd\step[2]\dd\step\cd\step[2]\id\step[2]\coro {\bar \sigma }\\
\dd\step\cd\step\id\step\cd\step\id\step[2]\id\\
\id\step\cd\step\X\step\id\step[2]\id\step\id\step[2]\id\\
\id\step\id\step[2]\id\step\id\step\X\step[2]\id\step\id\step[2]\id\\
\id\step\id\step[2]\id\step\X\step\XX\step\id\step\dd\\
\id\step\id\step[2]\X\step\id\step\id\step[2]\X\step\id\\
\id\step\XX\step\X\step\id\step\dd\step\id\step\id\\
\id\step\coro {\bar \sigma }\step \id\step\X\step\cu\step\id\\
\coRo \sigma\step\id\step\d\step\coro \sigma\\
\step[5]\coro \sigma
\end{tangle}
\]

\[ \ \
= \ \
 \begin{tangle}
\step[8]\object{H^\sigma}\step[5]\object{H^\sigma}\step[4]\object{H^\sigma}\\
\step[8]\id\step[5]\id\step[3]\cd\\
\step[8]\id\step[5]\id\step[2]\dd\step[2]\d\\
\step[8]\id\step[5]\XX\step[3]\cd\\
\step[8]\id\step[4]\ne2\step\cd\step[2]\id\step[2]\id\\
\step[8]\id\step[2]\dd\step[2]\dd\step\cd\step\id\step[2]\id\\
\step[8]\XX\step[3]\id\step[2]\id\step[2]\X\step[2]\id\\
\step[7]\ne2\step\cd\step[2]\id\step[2]\cu\step\coro {\bar \sigma }\\
\step[5]\ne2\step[2]\ne2\step[2]\XX\step[3]\id\\
\step[3]\ne2\step[2]\cd\step[3]\id\step[2]\id\step[2]\dd\\
\step\cd\step[3]\id\step\cd\step[2]\id\step[2]\coro {\bar \sigma }\\
\cd\step\d\step[2]\X\step\cd\step\id\\
\id\step[2]\id\step\cd\step\id\step\id\step\id\step[2]\id\step\id\\
\id\step[2]\id\step\id\step[2]\X\step\id\step\id\step[2]\id\step\id\\
\id\step[2]\id\step\XX\step\X\step\id\step[2]\id\step\id\\
\id\step[2]\X\step[2]\X\step\X\step[2]\id\step\id\\
\id\step[2]\id\step\XX\step\X\step\XX\step\id\\
\id\step[2]\X\step[2]\id\step\id\step\d\cu\step\id\\
\coro \sigma\step\coro {\bar \sigma }\step \coro \sigma\step\coro
\sigma
\end{tangle}
\]

\[ \ \
= \ \
\begin{tangle}
\step[6]\object{H^\sigma}\step[5]\object{H^\sigma}\step[4]\object{H^\sigma}\\
\step[6]\id\step[5]\id\step[3]\cd\\
\step[6]\id\step[5]\id\step[2]\dd\step[2]\d\\
\step[6]\id\step[5]\XX\step[3]\cd\\
\step[6]\id\step[4]\ne2\step\cd\step[2]\id\step[2]\id\\
\step[6]\id\step[2]\dd\step[2]\dd\step\cd\step\id\step[2]\id\\
\step[6]\XX\step[3]\id\step[2]\id\step[2]\X\step[2]\id\\
\step[5]\ne2\step\cd\step[2]\id\step[2]\cu\step\coro {\bar \sigma }\\
\step[3]\ne2\step[2]\ne2\step[2]\XX\step[3]\id\\
\step\cd\step\cd\step[3]\id\step[2]\id\step[2]\dd\\
\cd\step\X\step\cd\step[2]\id\step[2]\coro {\bar \sigma }\\
\id\step[2]\X\step\X\step[2]\id\step[2]\id\\
\coro \sigma\step\X\step\XX\step[2]\id\\
\step[2]\dd\step\id\step\d\step\nw2\step\nw2\\
\step[2]\coro {\bar \sigma } \step \cd\step\cd\step\id\\
\step[5]\id\step[2]\X\step[2]\id\step\id\\
\step[5]\coro \sigma\step\cu\step\id\\
\step[9]\coro \sigma
\end{tangle}
\stackrel {\hbox  {by 2-COC}}{=}
\begin{tangle}
\step[6]\object{H^\sigma}\step[5]\object{H^\sigma}\step[4]\object{H^\sigma}\\
\step[6]\id\step[5]\id\step[3]\cd\\
\step[6]\id\step[5]\id\step[2]\dd\step[2]\d\\
\step[6]\id\step[5]\XX\step[3]\cd\\
\step[6]\id\step[4]\ne2\step\cd\step[2]\id\step[2]\id\\
\step[6]\id\step[2]\dd\step[2]\dd\step\cd\step\id\step[2]\id\\
\step[6]\XX\step[2]\cd\step\id\step[2]\X\step[2]\id\\
\step[5]\ne2\step\cd\step\id\step[2]\id\step\cu\step\coro {\bar \sigma }\\
\step[3]\ne2\step[2]\ne2\step[2]\X\step[2]\id\step[2]\id\\
\step\cd\step\cd\step[3]\id\step\XX\step[2]\id\\
\cd\step\X\step\cd\step[2]\id\step\id\step[2]\coro {\bar \sigma }\\
\id\step[2]\X\step\X\step[2]\id\step[2]\id\step\id\\
\coro \sigma\step\X\step\XX\step[2]\id\step\id\\
\step[2]\dd\step\id\step\id\step\cd\step\id\step\id\\
\step[2]\coro {\bar \sigma }\step \id\step\id\step[2]\X\step\d\\
\step[5]\id\step\coro \sigma\step\cu\\
\step[5]\id\step[4]\dd\\
\step[5]\coRo \sigma\\
\end{tangle}
\]
\[
\stackrel {\hbox  {by commutativity}}{=}
\begin{tangle}
\step[7]\object{H^\sigma}\step[4]\object{H^\sigma}\step[3]\object{H^\sigma}\\
\step[7]\id\step[4]\id\step[2]\cd\\
\step[7]\id\step[4]\XX\step[2]\d\\
\step[7]\id\step[3]\dd\step\cd\step[2]\id\\
\step[7]\id\step[2]\dd\step\cd\step\coro {\bar \sigma }\\
\step[7]\XX\step[2]\id\step\cd\\
\step[6]\dd\step\cd\step\id\step\id\step[2]\id\\
\step[5]\dd\step\ne2\step[2]\X\step\id\step[2]\id\\
\step[4]\ne2\cd\step[3]\id\step\id\step\id\step[2]\id\\
\step[2]\cd\step\id\step\cd\step[2]\id\step\id\step\id\step[2]\id\\
\step\cd\step\X\step\id\step[2]\XX\step\id\step\id\step[2]\id\\
\cd\step\X\step\X\step[2]\id\step[2]\id\step\id\step\id\step[2]\id\\
\id\step[2]\X\step\X\step\XX\step[2]\id\step\id\step\id\step[2]\id\\
\coro \sigma\step\X\step\id\step\d\step\XX\step\id\step\id\step[2]\id\\
\step[2]\dd\step\id\step\coro \sigma\step\id\step[2]\X\step\id\step[2]\id\\
\step[2]\coro {\bar \sigma }\step \step[3]\id\step[2]\id\step\X\step[2]\id\\
\step[8]\id\step[2]\id\step\id\step\XX\\
\step[8]\id\step[2]\id\step\id\step\id\step\cd\\
\step[8]\id\step\dd\dd\step\X\step[2]\id\\
\step[8]\id\step\id\step\cu\step\cu\\
\step[7]\dd\step\XX\step[2]\dd\\
\step[7]\coro \sigma\step[2]\coro {\bar \sigma }
\end{tangle}
\stackrel {\hbox  {since H is a bialgebra }}{=}
\begin{tangle}
\step[4]\object{H^\sigma}\step[3]\object{H^\sigma}\step[3]\object{H^\sigma}\\
\step[4]\id\step[3]\id\step[2]\cd\\
\step[4]\id\step[3]\XX\step[2]\d\\
\step[4]\id\step[2]\dd\step\cd\step[2]\id\\
\step[4]\XX\step[2]\id\step[2]\coro {\bar \sigma }\\
\step[3]\ne2\step\cd\step\id\\
\step\cd\step[2]\id\step[2]\id\step\id\\
\cd\step\XX\step[2]\id\step\id\\
\id\step[2]\X\step[2]\XX\step\d\\
\coro \sigma\step\XX\step[2]\coro \sigma\\
\step[3]\coro {\bar \sigma }\\
\end{tangle}
\]
\[ \ \
= \ \
\begin{tangle}
\step\object{H^\sigma}\step[4]\object{H^\sigma}\step[4]\object{H^\sigma}\\
\cd\step[2]\cd\step[2]\cd\\
\id\step[2]\id\step[2]\id\step[2]\XX\step[2]\id\\
\id\step[2]\id\step[2]\XX\step[2]\id\step\cd\\
\id\step[2]\id\step\cd\step\id\step[2]\X\step[2]\id\\
\id\step[2]\X\step[2]\id\step\XX\step\coro {\bar \sigma }\\
\XX\step\coro {\bar \sigma }\step \coro \sigma\\
\coro \sigma\\
\end{tangle}
\step=\step
\begin{tangle}
\hbox {the right hand of }(CQT2).
\end{tangle}
\]

Similarly, we can check that (CQT1) and (CQT3) hold. Obviously,
$\hat \sigma $ is invertible. Thus $(H^\sigma , \hat \sigma )$ is
coquasitriangular. $\Box$

\begin {Proposition}\label {5.3.2}
Let $H$ be a bialgebra with a 2-cycle $Q$ of $H \otimes H.$ Define
$H^{Q}$ as follows: $H^{Q}=H $   as algebra and the comultiplication

\[
\begin{tangle}
\bigtriangleup_{H^Q}
\end{tangle}
\step=\step
\begin{tangle}
\step[4]\object{H^Q}\\
\ro Q\step\cd\step[2]\ro {\bar Q}\\
\id\step[2]\id\step\id\step[2]\XX\step[2]\id\\
\id\step[2]\id\step\cu\step[2]\id\step[2]\id\\
\id\step[2]\XX\step[3]\cu\\
\cu\step[2]\Cu\\
\end{tangle}\ \ \ .
\]
Then

(i)   $H^{Q}$ is a bialgebra;

(ii) Moreover, when $H$ is a Hopf algebra,  $H^Q$ is also a Hopf
algebra with antipode

 \[
\begin{tangle}
S_{H^Q}
\end{tangle}
\step=\step
\begin{tangle}
\step[4]\object{H^Q}\step[3]\ro {\bar Q }\\
\ro Q\step[2]\id\step[3]\S\step[2]\id\\
\id\step[2]\S\step[2]\S\step[3]\cu\\
\id\step[2]\id\step[2]\Cu\\
\id\step[2]\Cu\\
\Cu
\end{tangle}\ \ \ .
\]

(iii)  When $H$ is a cocommutative coalgebra,  $(H^{Q}, \hat Q)$ is
a symmetric quasitriangular bialgebra, where

\[
\begin{tangle}
\hat{Q}
\end{tangle}
\step=\step
\begin{tangle}
\ro Q\step[2]\ro {\bar Q}\\
\XX\step[2]\id\step[2]\id\\
\id\step[2]\XX\step[2]\id\\
\cu\step[2]\cu\\
\step\object{H}\step[5]\object{H}\\
\end{tangle}\ \ \ .
\]

(iv) When  $H$ is a cocommutative Hopf algebra,  $S_{H^{Q}}$ is
invertible.

\end {Proposition}

{\bf Proof.}  It is a dual case of Proposition \ref {5.3.1}.
\begin{picture}(8,8)\put(0,0){\line(0,1){8}}\put(8,8){\line(0,-1){8}}\put(0,0){\line(1,0){8}}\put(8,8){\line(-1,0){8}}\end{picture}
      \begin {Proposition} \label {5.3.3}
If $H$ and $A$ are  bialgebras and $\tau$ is a skew pairing on $(H
\otimes A),$ then

(i) $[\tau]$  is a 2-cocycle on $(A \otimes H) \otimes (A \otimes
H)$, where

\[
\begin{tangle}
[\tau]
\end{tangle}
\step=\step
\begin{tangle}
\object{A}\step[4]\object{H}\step[4]\object{A}\step[4]\object{H}\\
\id\step[4]\id\step[4]\id\step[4]\id\\
\id\step[4]\id\step[4]\id\step[4]\id\\
\QQ\epsilon \step[4]\coRo\tau\step[4]\QQ\epsilon \\
\end{tangle}\ \ \ .
\]

(ii) $\tau$  is invertible under convolution iff  $[\tau]$  is
invertible. Furthermore, $[\tau] ^{-1} = [\tau ^{-1}];$

(iii) When $\tau$ is invertible,  $(A \otimes H)^{[\tau]}$ is a
bialgebra. In this case,
    \[
\begin{tangle}
m_{A\bowtie _\tau H}
\end{tangle}
\step=\step
\begin{tangle}
\object{A}\step[4]\object{H}\step[4]\object{A}\step[4]\object{H}\\
\id\step [3]\cd\step[2]\cd\step[3]\id\\
\id\step[2]\cd\step\XX\step\cd\step[2]\id\\
\id\step[2]\id\step[2]\X\step[2]\X\step[2]\id\step[2]\id\\
\id\step[2]\coro{\tau}\step\XX\step\coro{\bar \tau}\step[2]\id\\
\id\step[4]\dd\step[2]\d\step[4]\id\\
\Cu\step[4]\Cu\\
\end{tangle}
\]

and we  write
 $A \bowtie _{\tau} H$ for $(A \otimes H)^{[\tau]}$ .

\end {Proposition}
{\bf Proof.}  The proofs of part (i) and (ii) are easily gotten by
turning  proofs of \cite [Proposition 1.5] {Do93} into braid
diagrams.

 (iii) It follows from Proposition \ref {5.3.1}. \begin{picture}(8,8)\put(0,0){\line(0,1){8}}\put(8,8){\line(0,-1){8}}\put(0,0){\line(1,0){8}}\put(8,8){\line(-1,0){8}}\end{picture}

\begin {Remark} \label{r2}  { \  }\end {Remark}
 Let $\tau$  be an invertible skew pairing
on $H \otimes A$
 and define:
\[
\begin{tangle}
\object{H}\step[2]\object{A}\\
\id\step[2]\id\\
\tu {\alpha }\\
\step\object{A}
\end{tangle}
\;=\enspace
\begin{tangles}{lcr}
\step\object{H}\step[3]\object{A}\\
\cd\step\cd\\
\id\step[2]\X\step\cd\\
\coro \tau\step\X\step[2]\id\\
\step[4]\coro {\bar{\tau}}\\
\end{tangles}
\]

and
\[
\begin{tangle}
\object{H}\step[2]\object{A}\\
\id\step[2]\id\\
\tu \beta\\
\step\object{H}
\end{tangle}
\;=\enspace
\begin{tangles}{lcr}
\step[2]\object{H}\step[3]\object{A}\\
\step\cd\step\cd\\
\cd\step\X\step[2]\id\\
\id\step[2]\X\step\coro {\bar{\tau}}\\
\coro \tau
\end{tangles} \ \ .
\]
  It is a straightforward verification to check
 that $A _\alpha \bowtie _\beta H = A \bowtie _\tau H.$

In fact, we can also show  that $A \bowtie _\tau H$ is a bialgebra.
(see Theorem \ref {2.1.4}).

\begin {Proposition} \label {5.3.3'}
Let $\tau$  be an invertible skew pairing on $H \otimes A$
 and define:

\[
\begin{tangle}
\object{H}\step[2]\object{A}\\
\id\step[2]\id\\
\tu \alpha \\
\step\object{A}
\end{tangle}
\;=\enspace
\begin{tangles}{lcr}
\step\object{H}\step[3]\object{A}\\
\cd\step\cd\\
\id\step[2]\X\step\cd\\
\coro \tau\step\X\step[2]\id\\
\step[4]\coro {\bar{\tau}}\\
\end{tangles}
\]

and
\[
\begin{tangle}
\object{H}\step[2]\object{A}\\
\id\step[2]\id\\
\tu \beta \\
\step\object{H}
\end{tangle}
\;=\enspace
\begin{tangles}{lcr}
\step[2]\object{H}\step[3]\object{A}\\
\step\cd\step\cd\\
\cd\step\X\step[2]\id\\
\id\step[2]\X\step\coro {\bar{\tau}}\\
\coro \tau
\end{tangles} \ \ .
\] Then  $A _\alpha \bowtie _\beta H$ is a bialgebra.

\end {Proposition}
{\bf Proof.}

\[
\begin{tangle}
\object{H}\\
\id\step[2]\Q {\eta _A}\\
\tu \alpha
\end{tangle}
\;=\enspace
\begin{tangles}{lcr}
\step\object{H}\\
\step\id\step[3]\Q {\eta_A}\\
\cd\step\cd\\
\id\step[2]\X\step\cd\\
\coro \tau\step\X\step[2]\id\\
\step[4]\coro {\bar \tau }
\end{tangles}
\]

\[
=
\begin{tangle}
\step\object{H}\\
\cd\step[2]\Q {\eta_A}\step[2]\Q {\eta_A}\step[2]\Q {\eta_A}\\
\id\step[2]\XX\step[2]\id\step[2]\id\\
\coro \tau\step[2]\XX\step[2]\id\\
\step[6]\coro {\bar{\tau}}
\end{tangle}
\;=\enspace
\begin{tangles}{lcr}
\object{H}\\
\QQ {\epsilon  _H}\\
\Q {\eta _A}
\end{tangles} \ \ .
\]
Similarly,
\[
\begin{tangle}
\step[2]\object{A}\\
\Q {\eta _A}\step[2]\id\\
\tu {\beta }
\end{tangle}
\;=\enspace
\begin{tangles}{lcr}
\object{A}\\
\QQ {\epsilon  _A}\\
\Q {\eta _H}
\end{tangles} .
\]
Hhus (M1) holds. We next show that (M2) holds.
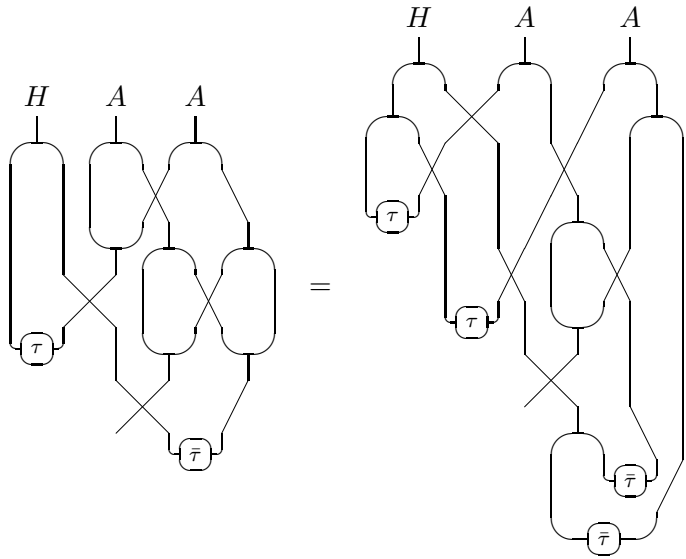
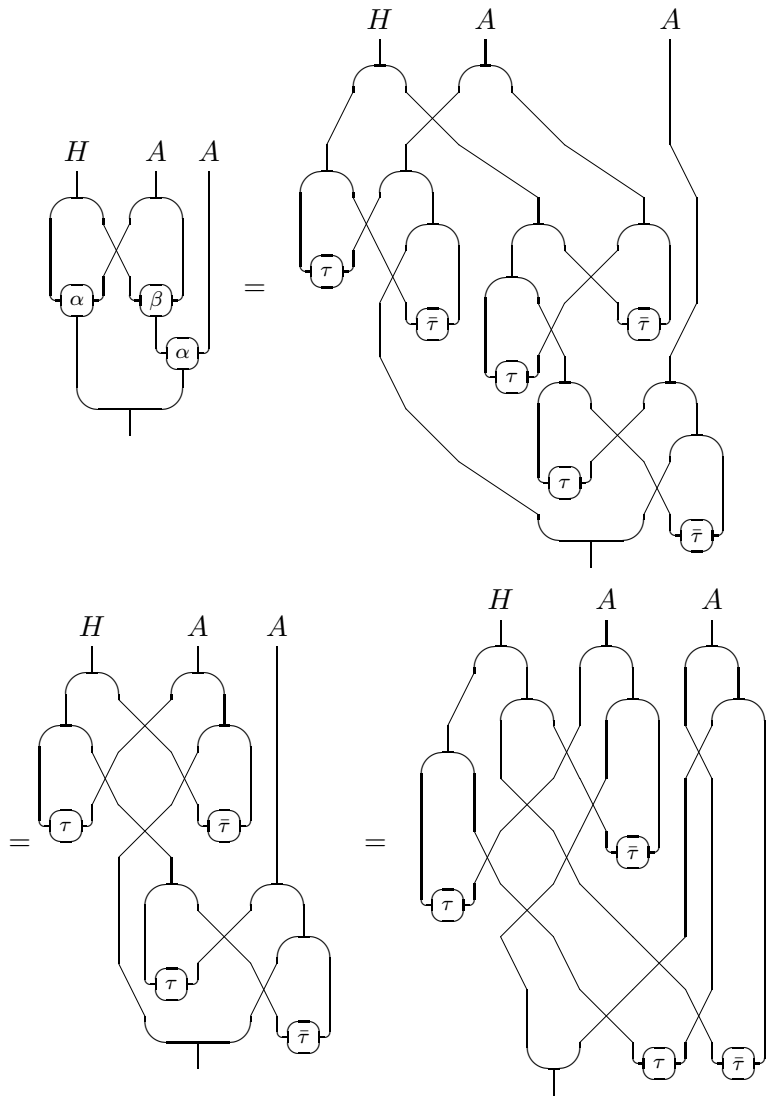
\[
\begin{tangle}
\object{H}\step\object{A}\step[2]\object{A}\\
\id\step\cu\\
\tu \alpha \\
\end{tangle}
\;=\enspace
\begin{tangles}{lcr}
\step\object{H}\step[2]\object{A}\object{A}\\
\cd\step\cu\\
\id\step[2]\id\step\cd\\
\id\step[2]\X\step\cd\\
\coro \tau\step\X\step[2]\id\\
\step[4]\coro {\bar \tau } \\
\end{tangles}
\;=\enspace
\]

\[
\begin{tangle}
\step\object{H}\step[3]\object{A}\step[3]\object{A}\\
\cd\step\cd\step\cd\\
\id\step[2]\id\step\id\step[2]\X\step[2]\d\\
\id\step[2]\id\step\cu\cd\step\cd\\
\id\step[2]\XX\step\id\step[2]\X\step[2]\id\\
\coro \tau\step[2]\id\step\cu\step\cu\\
\step[4]\XX\step[2]\dd\\
\step[6]\coro {\bar{\tau}}
\end{tangle}
\step=\step
\begin{tangle}
\step[2]\object{H}\step[4]\object{A}\step[4]\object{A}\\
\step\cd\step[2]\cd\step[2]\cd\\
\cd\step\XX\step[2]\id\step\dd\step\cd\\
\id\step[2]\X\step[2]\id\step[2]\X\step[2]\id\step[2]\id\\
\coro \tau\step\id\step[2]\id\step\dd\cd\step\id\step[2]\id\\
\step[3]\id\step[2]\X\step\id\step[2]\X\step[2]\id\\
\step[3]\coro \tau\step\id\step\cu\step\id\step[2]\id\\
\step[6]\XX\step[2]\id\step[2]\id\\
\step[7]\cd\step\d\step\id\\
\step[7]\id\step[2]\coro {\bar{\tau}}\dd\\
\step[7]\coRo {\bar{\tau}}
\end{tangle}
\] \ \
and \ \
\[
\begin{tangle}
\step\object{H}\step[3]\object{A}\step[2]\object{A}\\
\cd\step\cd\step\id\\
\id\step[2]\X\step[2]\id\step\id\\
\tu \alpha \step\tu \beta \step\id\\
\step\id\step[3]\tu \alpha \\
\step\Cu
\end{tangle}
\step=\step
\begin{tangle}
\step[3]\object{H}\step[4]\object{A}\step[7]\object{A}\\
\step[2]\cd\step[2]\cd\step[6]\id\\
\step\dd\step[2]\XX\step[2]\nw2\step[5]\id\\
\cd\step\cd\step\nw3\step[3]\nw3\step[3]\nw1\\

\id\step[2]\X\step\cd\step[2]\cd\step[2]\cd\step\id\\

\coro \tau\step\X\step[2]\id \step \cd\step\XX\step[2]\id\step\id\\

\step[3]\id\step\coro {\bar{\tau}}\step \id\step[2]\X\step[2]\coro {\bar{\tau}}\dd\\

\step[3]\d\step[3]\coro \tau\cd\step[2]\cd\\
\step[4]\nw2\step[4]\id\step[2]\XX\step\cd\\
\step[6]\nw3\step[2]\coro \tau\step[2]\X\step[2]\id\\
\step[9]\Cu\step\coro {\bar \tau }
\end{tangle}
\]

 \[
 =
\begin{tangle}
\step[2]\object{H}\step[4]\object{A}\step[3]\object{A}\\
\step\cd\step[2]\cd\step[2]\id\\
\cd\step\XX\step\cd\step\id\\
\id\step[2]\X\step[2]\X\step[2]\id\step\id\\
\coro \tau\step\XX\step\coro {\bar{\tau}}\step\id\\
\step[3]\id\step\cd\step[2]\cd\\
\step[3]\id\step\id\step[2]\XX\step\cd\\
\step[3]\d\coro \tau\step[2]\X\step[2]\id\\
\step[4]\Cu\step\coro {\bar{\tau}}
\end{tangle}
\step=\step
\begin{tangle}
\step[3]\object{H}\step[4]\object{A}\step[4]\object{A}\\
\step[2]\cd\step[2]\cd\step[2]\cd\\
\step\dd\step\cd\step\id\step\cd\step\id\step\cd\\
\cd\step\id\step[2]\X\step\id\step[2]\id\step\X\step[2]\id\\
\id\step[2]\id\step\XX\step\X\step[2]\id\step\id\step\id\step[2]\id\\
\id\step[2]\X\step[2]\X\step\coro {\bar{\tau}}\step\id\step\id\step[2]\id\\
\coro \tau\step\XX\step\nw2\step[3]\id\step\id\step[2]\id\\
\step[3]\d\step\d\step[2]\XX\step\id\step[2]\id\\
\step[4]\id\step[2]\XX\step[2]\X\step[2]\id\\
\step[4]\cu\step[2]\coro \tau\step\coro {\bar{\tau}}\\
\end{tangle}
\] Thus (M2) holds. Similarly, we can show that (M3) holds. Now we
show that (M4) holds.

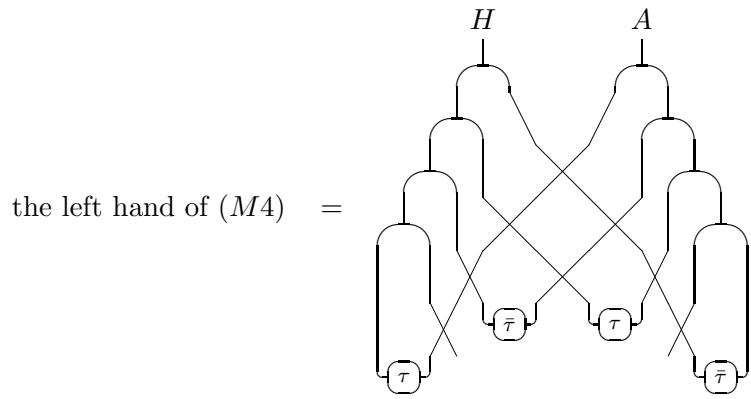
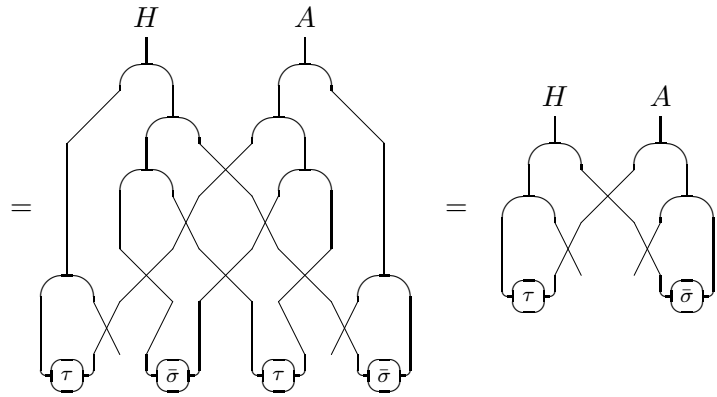
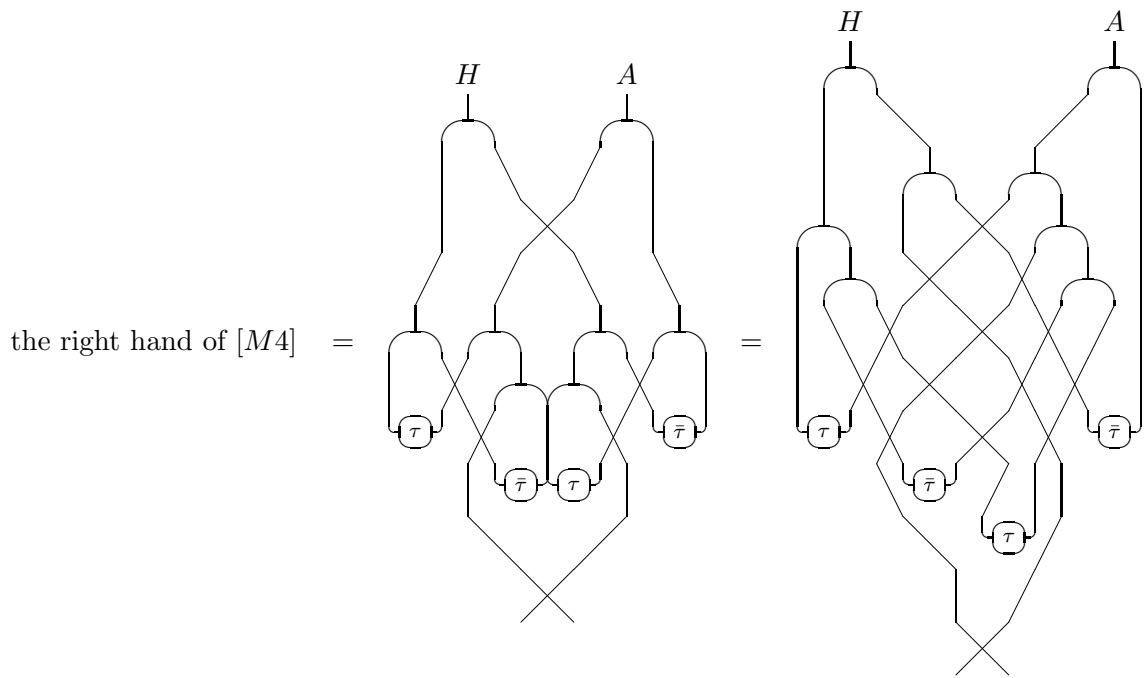
\[ \hbox { the left hand of }(M4)
\step=\step
\begin{tangle}
\step[4]\object{H}\step[6]\object{A}\\
\step[3]\cd\step[4]\cd\\
\step[2]\cd\step\d\step[2]\dd\step\cd\\
\step\cd\step\id\step[2]\XX\step[2]\id\step\cd\\
\cd\step\id\step\XX\step[2]\XX\step\id\step\cd\\
\id\step[2]\id\step\X\step[2]\XX\step[2]\X\step\id\step[2]\id\\
\id\step[2]\X\step\coro {\bar{\tau}}\step[2]\coro \tau\step\X\step[2]\id\\
\coro \tau\step[10]\coro {\bar{\tau}}
\end{tangle}
\]

\[
=
\begin{tangle}
\step[4]\object{H}\step[6]\object{A}\\
\step[3]\cd\step[4]\cd\\
\step[2]\ne2\step\cd\step[2]\cd\step\nw2\\
\step\id\step[2]\cd\step\XX\step\cd\step[2]\id\\
\step\id\step[2]\id\step[2]\X\step[2]\X\step[2]\id\step[2]\id\\
\cd\step\XX\step\XX\step\XX\step\cd\\
\id\step[2]\X\step\dd\step\id\step[2]\id\step\d\step\X\step[2]\id\\
\coro \tau\step[2]\coro {\bar{\sigma}}\step[2]\coro
\tau\step[2]\coro {\bar{\sigma}}
\end{tangle}
\step=\step
\begin{tangle}
\step[2]\object{H}\step[4]\object{A}\\
\step\cd\step[2]\cd\\
\cd\step\XX\step\cd\\
\id\step[2]\X\step[2]\X\step[2]\id\\
\coro \tau\step[4]\coro {\bar{\sigma}}
\end{tangle}
\] and

\[
\begin{tangle}
 \hbox { the right hand of } [M4]
\end{tangle}
\step=\step
\begin{tangle}
\step[3]\object{H}\step[6]\object{A}\\
\step[2]\cd\step[4]\cd\\
\step[2]\id\step[2]\d\step[2]\dd\step[2]\id\\
\step[2]\id\step[3]\XX\step[3]\id\\
\step\dd\step[2]\dd\step[2]\d\step[2]\d\\
\cd\step\cd\step[2]\cd\step\cd\\
\id\step[2]\X\step\cd\cd\step\X\step[2]\id\\
\coro \tau\step\X\step[2]\id\step[2]\X\step\coro {\bar{\tau}}\\
\step[3]\id\step\coro {\bar{\tau}}\coro \tau\step\id\\
\step[3]\nw2\step[4]\ne2\\
\step[5]\XX
\end{tangle}
\step=\step
\begin{tangle}
\step[2]\object{H}\step[10]\object{A}\\
\step\cd\step[8]\cd\\
\step\id\step[2]\nw2\step[6]\ne2\step[2]\id\\
\step\id\step[3]\cd\step[2]\cd\step[3]\id\\
\cd\step[2]\id\step[2]\XX\step\cd\step[2]\id\\
\id\step\cd\step\XX\step[2]\X\step\cd\step\id\\
\id\step\d\step\X\step[2]\XX\step\X\step\dd\step\id\\
\id\step[2]\X\step\XX\step[2]\X\step\X\step[2]\id\\
\coro \tau\step\X\step[2]\XX\step\X\step\coro {\bar{\tau}}\\
\step[3]\d\coro {\bar{\tau}}\step\dd\step\id\step\id\\
\step[4]\nw2\step[2]\coro \tau\dd\\
\step[6]\id\step[2]\dd\\
\step[6]\XX
\end{tangle}
\]

\[
=
\begin{tangle}
\step\object{H}\step[6]\object{A}\\
\cd\step[4]\cd\\
\id\step\cd\step[2]\cd\step\id\\
\id\step\d\step\XX\step\dd\step\id\\
\id\step[2]\X\step[2]\X\step[2]\id\\
\coro \tau\step[4]\coro {\bar{\tau}}
\end{tangle}
\step=\step
\begin{tangle}
 \hbox {the left hand of (M4)}.
\end{tangle}\ \ \
\] Thus (M4) holds. We finally show that $A$ is an $H$-module
coalgebra. We see that
 \[
\begin{tangle}
\object{H}\step\object{A}\step[2]\object{A}\\
\id\step\tu \alpha \\
\tu \alpha
\end{tangle}
\step=\step
\begin{tangle}
\step\object{H}\step[2]\object{H}\step[3]\object{A}\\
\step\id\step\cd\step\cd\\
\step\id\step\id\step[2]\X\step\cd\\
\step\id\step\coro \tau\step\X\step[2]\id\\
\cd\step[2]\cd\coro {\bar{\tau}}\\
\id\step[2]\XX\step\cd\\
\coro \tau\step[2]\X\step[2]\id\\
\step[5]\coro{ \bar \tau }
\end{tangle}
\]
\[
=
\begin{tangle}
\object{H}\step[2]\object{H}\step[4]\object{A}\\
\id\step\cd\step[2]\cd\\
\id\step\id\step[2]\XX\step\cd\\
\d\coro \tau\step[2]\X\step\cd\\
\cd\step[2]\dd\step\X\step\cd\\
\id\step[2]\XX\step\dd\step\X\step[2]\id\\
\coro \tau\step[2]\X\step[2]\id\step\coro {\bar{\tau}}\\
\step[5]\coro {\bar{\tau}}
\end{tangle}
\stackrel {by (P2)}  {=}
\begin{tangle}
\object{H}\step[2]\object{H}\step[4]\object{A}\\
\id\step\cd\step[2]\cd\\
\id\step\id\step[2]\XX\step\cd\\
\d\coro \tau\step[2]\X\step\cd\\
\cd\step[2]\dd\step\X\step\dd\\
\id\step[2]\XX\step\dd\step\id\step\id\\
\coro \tau\step[2]\X\step[2]\id\step\id\\
\step[5]\cu\step\id\\
\step[6]\coro {\bar{\tau}}
\end{tangle}
\]

\[
=
\begin{tangle}
\object{H}\step[2]\object{H}\step[4]\object{A}\\
\id\step\cd\step[2]\cd\\
\id\step\id\step[2]\id\step\cd\step\id\\
\id\step\id\step[2]\X\step[2]\id\step\id\\
\d\coro \tau\step\XX\step\id\\
\cd\step[2]\id\step[2]\id\step\id\\
\id\step[2]\XX\step[2]\id\step\id\\
\coro \tau\step[2]\cu\cd\\
\step[5]\X\step[2]\id\\
\step[6]\coro {\bar{\tau}}
\end{tangle}
\stackrel {by (Sp2)}  {=}
\begin{tangle}
\step\object{H}\step[3]\object{H}\step[2]\object{A}\\
\cd\step\cd\step\id\\
\id\step[2]\X\step[2]\id\step\id\\
\cu\step\cu\cd\\
\step\d\step[2]\X\step\cd\\
\step[2]\coro \tau\step\X\step[2]\id\\
\step[6]\coro {\bar{\tau}}
\end{tangle}
\]

\[
 =
\begin{tangle}
\object{H}\step[2]\object{H}\step[2]\object{A}\\
\cu\step[2]\id\\
\cd\step\cd\\
\id\step[2]\X\step\cd\\
\coro \tau\step\X\step[2]\id\\
\step[4]\coro {\bar{\tau}}
\end{tangle}
\step=\step
\begin{tangle}
\object{H}\step[2]\object{H}\step\object{A}\\
\cu\step\id\\
\step \tu \alpha
\end{tangle}
\]
and
\[
\begin{tangle}
\step[2]\object{A}\\
\Q {\eta_H}\step[2]\id\\
\tu \alpha
\end{tangle}
\step=\step
\begin{tangle}
\step[4]\object{A}\\
\step\Q {\eta_H}\step[3]\id\\
\cd\step\cd\\
\id\step[2]\X\step\cd\\
\coro \tau\step\X\step[2]\id\\
\step[4]\coro {\bar{\tau}}
\end{tangle}
\step=\step
\begin{tangle}
\step\object{A}\\
\cd\\
\QQ {\epsilon _A}\step\cd\\
\step[3]\QQ {\epsilon _A}
\end{tangle}
\step=\step
\begin{tangle}
\object{A}\\
\id\\
\id
\end{tangle}\ \ \ .
\]
Therefore, $A$ an $H$-module. We also see that
 \[
\begin{tangle}
\step\object{H}\step[3]\object{A}\\
\cd\step\cd\\
\id\step[2]\X\step[2]\id\\
\tu \alpha \step\tu \beta\\
\end{tangle}
\step=\step
\begin{tangle}
\step[2]\object{H}\step[8]\object{A}\\
\step\cd\step[6]\cd\\
\step\id\step[2]\nw2\step[4]\ne2\step[2]\id\\
\step\id\step[4]\XX\step[4]\id\\
\step\id\step[3]\dd\step[2]\d\step[3]\id\\
\cd\step\cd\step[2]\cd\step\cd\\
\id\step[2]\X\step\cd\step\id\step[2]\X\step\cd\\
\coro \tau\step\X\step[2]\id\step\coro \tau\step\X\step[2]\id\\
\step[4]\coro {\bar{\tau}}\step[5]\coro {\bar{\tau}}
\end{tangle}
 \step=\step
\begin{tangle}
\step\object{H}\step[3]\object{A}\\
\cd\step\cd\\
\id\step[2]\X\step\cd\\
\coro \tau\step\X\step[2]\id\\
\step[2]\cd\coro {\bar{\tau}}
\end{tangle}
\step=\step
\begin{tangle}
\object{H}\step[2]\object{A}\\
\tu \alpha \\
\cd
\end{tangle}
\]

and
\[
\begin{tangle}
\object{H}\step[2]\object{A}\\
\tu \alpha \\
\step\QQ {\epsilon _A}
\end{tangle}
 \step=\step
\begin{tangle}
\step\object{H}\step[3]\object{A}\\
\cd\step\cd\\
\id\step[2]\X\step\cd\\
\coro \tau\step\X\step[2]\id\\
\step[3]\QQ {\epsilon _A}\step\coro {\bar{\tau}}
\end{tangle}
 \step=\step
\begin{tangle}
\object{H}\step[2]\object{A}\\
\QQ {\epsilon _A}\step[2]\QQ {\epsilon _A}
\end{tangle}\ \ \ .
\]

Thus $A$ is an $H$-module coalgebra.  Similarly, we can show that
$H$ is an $A$-module coalgebra.   $\Box$

\begin {Proposition} \label {5.3.4}
If $H$ and $A$ are  bialgebras and $P$ is a skew copairing of $H
\otimes A,$ then

(i) $[P]$  is a 2-cycle of $A \otimes H$, where
  \[
\begin{tangle}
 $[P]$
\end{tangle}
\step=\step
\begin{tangle}
\Q {\eta_A}\step\step\ro P\step[2]\Q {\eta_H}\\
\id\step[2]\id\step[2]\id\step[2]\id
\end{tangle}\ \ \ .
\]

(ii) $P$  is invertible under convolution iff $[P]$  is invertible.
Furthermore $[P]^{-1} = [P^{-1}]$;

(iii) Moreover, when $P$  is invertible, $(A \otimes H)^{[P]}$  is a
bialgebra. in this case, \[
\begin{tangle}
\Delta_{A\bowtie ^PH}
\end{tangle}
\step=\step
\begin{tangle}
\step\object{A}\step[8]\object{H}\\
\cd\step[6]\cd\\
\id\step[2]\nw2\step[4]\ne2\step[2]\id\\
\id\step\ro P\step\XX\step\ro {\bar{P}}\step\id\\
\id\step\id\step[2]\X\step[2]\X\step[2]\id\step\id\\
\id\step\cu\step\XX\step\cu\step\id\\
\id\step[2]\cu\step[2]\cu\step[2]\id
\end{tangle}
\]
and we write
 $ A \bowtie ^{P} H$ for
 $(A \otimes H)^{[P]} \hbox { \ \ }.$
\end {Proposition}

{\bf Proof.} It is a dual case of Proposition \ref {5.3.3}.
\begin{picture}(8,8)\put(0,0){\line(0,1){8}}\put(8,8){\line(0,-1){8}}\put(0,0){\line(1,0){8}}\put(8,8){\line(-1,0){8}}\end{picture}

\begin {Lemma} \label {5.3.5}
Let $H$ and $A$ be  bialgebras. Let  $\tau$ be an invertible skew
pairing on $(H \otimes A)$ and $P$ be an invertibe copairing of $(H
\otimes A).$ Set $D = A \bowtie _{\tau} H$. Then $(D, [P])$  is
almost cocommutative  iff

 (AC1)：
 \[
\begin{tangle}
\step\object{H}\step[2]\object{A}\step[2]\object{H}\\
\cd\step\id\step[2]\id\\
\id\step[2]\X\step[2]\id\\
\coro \tau\step\cu
\end{tangle}
\step=\step
\begin{tangle}
\step\object{H}\step[5]\object{A}\step[2]\object{H}\\
\cd\step\ro P\step\XX\step\ro {\bar{P}}\step\\
\id \step [2] \id\step\id\step[2]\X\step[2]\X\step[2]\id\step\\
\id\step [2] \id\step\cu\step\XX\step\cu\step\\
\nw1\step [1] \nw1\step\cu\step[2]\cu\step[2]\\

\step \id\step [2] \XX\step[3]\ne2 \step [2]\\

\step \XX\step [2] \coro \tau \\

\step \cu \step [2]\\
\end{tangle};
\]
 (AC2)：
 \[
\begin{tangle}
\object{A}\step[2]\object{H}\step[2]\object{A}\\
\id\step[2]\id\step\cd\\
\id\step[2]\X\step[2]\id\\
\XX\step\coro \tau\\
\cu
\end{tangle}
\step=\step
\begin{tangle}
\step[3]\object{A}\step[2]\object{H}\step[4]\object{A}\\
\ro P\step\XX\step\ro {\bar{P}}\step\id\\
\id\step[2]\X\step[2]\X\step[2]\id\step\id\\
\cu\step\XX\step\cu\step\id\\
\step\cu\step[2]\cu\step\cd\\
\step[2]\id\step[4]\XX\step[2]\id\\
\step[2]\coRo \tau\step[2]\cu
\end{tangle}\ \ \ .
\]

\end {Lemma}
{\bf Proof.} It can be obtained by straightforward
verification$\Box$

\begin {Lemma} \label {5.3.6}
Let $H$ and $A$ be  bialgebras. Let  $\tau$ be an invertible skew
pairing on $(H \otimes A)$ and $P$ be an invertibe copairing of $(H
\otimes A).$ Set $D = A \bowtie _{\tau} H$. Then $(D, [P])$  is
almost cocommutative  iff

(ACO1)
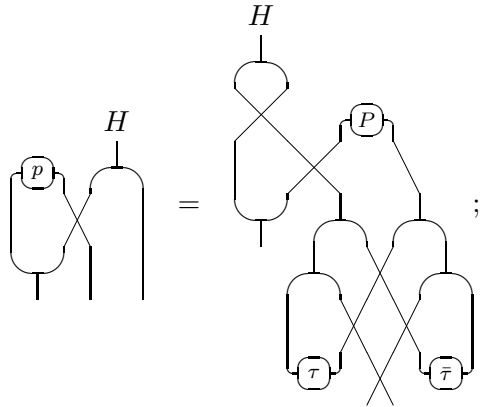
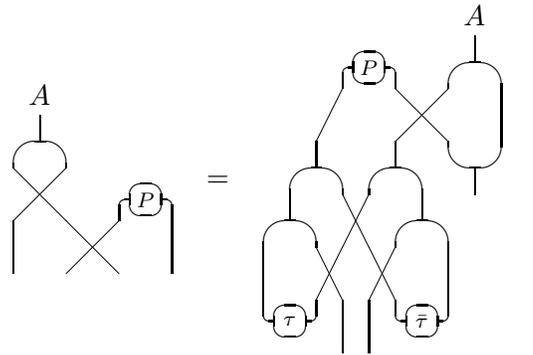
\[：
\begin{tangle}
\step[4]\object{H}\\
\ro p\step\cd\\
\id\step[2]\X\step[2]\id\\
\cu\step\id\step[2]\id\\
\end{tangle}
\step=\step
\begin{tangle}
\step\object{H}\\
\cd\\
\XX\step[2]\ro P\\
\id\step[2]\XX\step[2]\d\\
\cu\step\cd\step\cd\\
\step[2]\cd\step\X\step\cd\\
\step[2]\id\step[2]\X\step\X\step[2]\id\\
\step[2]\coro \tau\step\X\step\coro {\bar{\tau}}
\end{tangle};
\]

(ACO2)：
\[
\begin{tangle}
\step\object{A}\\
\cd\\
\XX\step[2]\ro P\\
\id\step[2]\XX\step[2]\id
\end{tangle}
\step=\step
\begin{tangle}
\step[8]\object{A}\\
\step[3]\ro P\step[2]\cd\\
\step[2]\dd\step[2]\XX\step[2]\id\\
\step\cd\step\cd\step\cu\\
\cd\step\X\step\cd\\
\id\step[2]\X\step\X\step[2]\id\\
\coro \tau\step\id\step\id\step\coro {\bar{\tau}}
\end{tangle}\ \ \ .
\]
\end {Lemma}
{\bf Proof.} It is the dual case of Lemma \ref {5.3.5}.
\begin{picture}(8,8)\put(0,0){\line(0,1){8}}\put(8,8){\line(0,-1){8}}\put(0,0){\line(1,0){8}}\put(8,8){\line(-1,0){8}}\end{picture}

\begin {Proposition} \label {5.3.7}
Let $A$ be a Hopf algebra and $H$ be a Hopf algebra with invertible
antipode. Let  $\tau$ be a pairing on $(H \otimes A)$ and $P$ be a
skew copairing of $(H \otimes A).$ Set $D = A \bowtie ^P H$. If

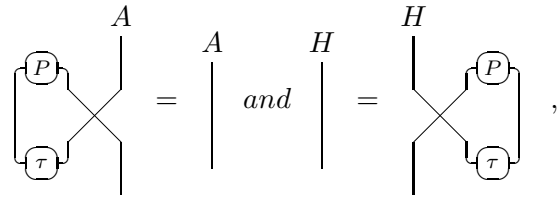
\[
\begin{tangle}
\step[4]\object{A}\\
\ro P\step[2]\id\\
\id\step[2]\XX\\
\coro \tau\step[2]\id
\end{tangle}
\step=\step
\begin{tangle}
\object{A}\\
\id\\
\id
\end{tangle} \ \ \ and \ \ \
\begin{tangle}
\object{H}\\
\id\\
\id
\end{tangle}
\step=\step
\begin{tangle}
\object{H}\\
\id\step[2]\ro P\\
\XX\step[2]\id\\
\id\step[2]\coro \tau
\end{tangle} \ \ \ ,
\] then $(D, [\tau])$  is almost commutative.
\end {Proposition}
{\bf Proof.}
\[
\begin{tangle}
 \hbox  { the right hand of (AC1)}
\end{tangle}
\stackrel {by (P1)}  {=}
\begin{tangle}
\step\object{H}\step[6]\object{A}\step[2]\object{H}\\
\step\id\step[2]\ro P\step[2]\XX\step\ro {\bar{P}}\\
\step\id\step[2]\id\step[2]\XX\step[2]\X\step[2]\id\\
\step\id\step[2]\cu\step[2]\XX\step\cu\\
\cd\step[2]\id\step[2]\dd\step[2]\id\step[2]\id\\
\id\step\cd\step\cu\step[3]\id\step[2]\id\\
\id\step\id\step[2]\XX\step[4]\id\step[2]\id\\
\id\step\XX\step[2]\coRo \tau\step[2]\id\\
\X\step[2]\nw2\step[6]\ne2\\
\id\step\d\step[3]\coRo \tau\\
\cu
\end{tangle}
\stackrel { \hbox { by assumption}} {=}
\begin{tangle}
\step\object{H}\step[4]\object{A}\step[2]\object{H}\\
\cd\step[3]\id\step[2]\id\\
\id\step\cd\step[2]\XX\step[2]\ro {\bar{P}}\\
\id\step\id\step[2]\cu\step[2]\XX\step[2]\id\\
\id\step\id\step[3]\id\step[2]\dd\step[2]\cu\\
\id\step\nw2\step[2]\cu\step[4]\id\\
\d\step[2]\XX\step[4]\dd\\
\step\XX\step[2]\coRo \tau\\
\step\cu
\end{tangle}
\]
\[
\stackrel {by (P1)}  {=}
\begin{tangle}
\step\object{H}\step[4]\object{A}\step[2]\object{H}\\
\cd\step[3]\id\step[2]\id\\
\id\step\cd\step[2]\XX\step[2]\ro {\bar{P}}\\
\id\step\id\step[2]\cu\step[2]\XX\step[2]\id\\
\id\step\id\step[3]\id\step[2]\dd\step[2]\id\step[2]\id\\
\id\step\nw2\step[2]\cu\step[3]\id\step[2]\id\\
\d\step[2]\XX\step[4]\id\step[2]\id\\
\step\XX\step\cd\step[2]\dd\step[2]\id\\
\step\cu\step\nw2\step\coro \tau\step[2]\dd\\
\step[6]\coRo \tau
\end{tangle}
\]

\[
=
\begin{tangle}
\step\object{H}\step[4]\object{A}\step[2]\object{H}\\
\cd\step[3]\id\step[2]\id\step[2]\ro P\\
\id\step\cd\step[2]\XX\step[2]\O {\bar S}\step[2]\id\\
\id\step\id\step[2]\cu\step[2]\XX\step[2]\id\\
\id\step\id\step[3]\id\step[2]\dd\step[2]\id\step[2]\id\\
\id\step\nw2\step[2]\cu\step[3]\id\step[2]\id\\
\d\step[2]\XX\step[4]\id\step[2]\id\\
\step\XX\step\cd\step[2]\dd\step[2]\id\\
\step\cu\step\nw2\step\coro\tau\step[2]\dd\\
\step[6]\coRo \tau
\end{tangle}
\stackrel {\hbox { by assumption}} {=}
\begin{tangle}

\step\object{H}\step[5]\object{A}\step[2]\object{H}\\
\cd\step[4]\id\step[2]\id\\
\id\step[2]\d\step[3]\id\step[2]\id\\
\id\step[2]\cd\step[2]\id\step[2]\id\\
\id\step\cd\step\XX\step[2]\id\\
\id\step\id\step[2]\id\step\d\step\cu\\
\id\step\nw2\step\coro \tau\step\ne2\\
\id\step[3]\XX\\
\id\step[3]\id\step[2]\O {\bar S}\\
\nw2\step[2]\cu\\
\step [2]\XX\\
\step [2]\cu
\end{tangle}
\]

\[
=
\begin{tangle}
\step[2]\object{H}\step[5]\object{A}\step[2]\object{H}\\
\step\cd\step[4]\id\step[2]\id\\
\cd\step\d\step[3]\id\step[2]\id\\
\id\step[2]\id\step\cd\step[2]\id\step[2]\id\\
\id\step[2]\id\step\id\step[2]\XX\step[2]\id\\
\id\step[2]\id\step\coro \tau\step[2]\cu\\
\id\step[2]\nw2\step[4]\ne2\\
\nw2\step[3]\XX\\
\step[2]\XX\step[2]\O {\bar S}\\
\step\dd\step[2]\XX\\
\step\id\step[3]\cu\\
\step\Cu
\end{tangle}
\step=\step
\begin{tangle}
\step\object{H}\step[3]\object{A}\step[2]\object{H}\\
\cd\step[2]\id\step[2]\id\\
\id\step[2]\XX\step[2]\id\\
\coro \tau\step[2]\cu
\end{tangle}
\step=\step
\begin{tangle}
\hbox  {the left hand of (AC1)}.
\end{tangle}\ \ \ .
\] Thus (AC1) holds. (AC2) can similarly be shown.

 Consequently, $(B,
[\tau ])$  is almost commutative by Lemma \ref {5.3.5}.
\begin{picture}(8,8)\put(0,0){\line(0,1){8}}\put(8,8){\line(0,-1){8}}\put(0,0){\line(1,0){8}}\put(8,8){\line(-1,0){8}}\end{picture}

\begin {Proposition}  \label {5.3.8}.
Let $A$ be a Hopf algebra and $H$ be a Hopf algebra with invertible
antipode.
 Let  $\tau$ be a skew pairing on
$(H \otimes A)$ and $P$ be a copairing of $(H \otimes A).$ Set $D =
A \bowtie _{\tau} H$. If
 \[
\begin{tangle}
\step[4]\object{A}\\
\ro P\step[2]\id\\
\id\step[2]\XX\\
\coro \tau
\end{tangle}
\step=\step
\begin{tangle}
\object{A}\\
\id\\
\id
\end{tangle}
 \ \ \ and \ \ \
\begin{tangle}
\object{H}\\
\id\\
\id
\end{tangle}
\step=\step
\begin{tangle}
\object{H}\\
\id\step[2]\ro P\\
\XX\step[2]\id\\
\step[2]\coro \tau
\end{tangle} \ \ \ ,
\]
 then $(D, [P])$  is almost cocommutative.
\end {Proposition}
{\bf Proof.} It is a dual case of Proposition \ref {5.3.7}.
\begin{picture}(8,8)\put(0,0){\line(0,1){8}}\put(8,8){\line(0,-1){8}}\put(0,0){\line(1,0){8}}\put(8,8){\line(-1,0){8}}\end{picture}

Let $\Delta ^{cop} = C_{H,H}\Delta$  and $m^{op} = m C_{H,H}.$ We
denote $(H, \Delta ^{cop}, \epsilon)$  by $H^{cop}$ and  $(H,
m^{op}, \eta)$  by $H^{op}$.
\begin {Theorem} \label {5.3.9}
Let $({\cal C},C) $ be a symmetric braided tensor category with left
dual and $H$ be a Hopf algebra with  antipode $S$. Let $ev_H$  and
$coev_H$ denote the evaluation and coevaluation of $H$ respectively.
Set $A= (H^*)^{op}, \tau = ev_H C_{H,H}$ and $P=coev_H.$

Then

(i)  $\tau $ is a skew pairing on $H \otimes A$;

(ii) $P$ is a  copairing of $ H\otimes A;$

(iii) Moreover, when $S$ is invertible, $(D(H), [P])$ is a
quasitriangular
 Hopf algebra, where
$D(H) =A \bowtie _{\tau} H,$ called Drinfeld double of $H$.

\end {Theorem}
{\bf Proof.} Using \cite [Proposition 2.4]{Ma95a}, we can get  the
proofs of part (i) and part (ii).

(iii) It follows from Proposition \ref {5.3.1}, Proposition \ref
{5.3.3} and Proposition \ref {5.3.8}.
\begin{picture}(8,8)\put(0,0){\line(0,1){8}}\put(8,8){\line(0,-1){8}}\put(0,0){\line(1,0){8}}\put(8,8){\line(-1,0){8}}\end{picture}

\begin {Theorem} \label {5.3.10}
Let $({\cal C},C) $ be a symmetric braided tensor category with left
duality and $H$ be a Hopf algebra with  antipode $S$. Let $ev_H$
and $coev_H$ denote the evaluation and coevaluation of $H$
respectively. Set $A= (H^*)^{cop}, \tau = ev_H C_{H,H}$ and
$P=coev_H .$

Then

(i)  $\tau $ is a  pairing on $H \otimes A$;

(ii) $P$ is a skew copairing of $ H\otimes A;$

(iii) Moreover, when $S$ is invertible, $(CD(H), [\tau])$ is a
coquasitriangular Hopf algebra, where $CD(H) =A \bowtie ^P H,$
called Drinfeld codouble of $H$.

\end {Theorem}
{\bf Proof.} Using \cite [Proposition 2.4]{Ma95a}, we can get the
proofs of  part (i) and part (ii).

(iii) It follows from Proposition \ref {5.3.2}, Proposition \ref
{5.3.4} and Proposition \ref {5.3.7}.
\begin{picture}(8,8)\put(0,0){\line(0,1){8}}\put(8,8){\line(0,-1){8}}\put(0,0){\line(1,0){8}}\put(8,8){\line(-1,0){8}}\end{picture}

 Let $({\cal C},C) $ be a symmetric braided tensor category with left duality.
 Let $A$ and $H$ be bialgebras. We easily check that
$\tau ^*$ is a (skew copairing ) copairing of $A^* \otimes H^*$  if
$\tau$ is a skew pairing (pairing) of $H\otimes A$, and $P^*$ is a
skew pairing (pairing ) on $A^* \otimes H^*$ if $P$ is a copairing
(skew pairing) of $H \otimes A$. Therefore we can easily check the
following.

\begin {Proposition} \label {5.3.11}
Let $({\cal C},C) $ be a symmetric braided tensor category with left
duality.

(i)  If         $(H,r)$ is a coquasitriangular bialgebra (Hopf
algebra), then $(H^*, r^*C_{H^*,H^*})$ is  a quasitriangular
bialgebra (Hopf algebra);

(ii) If $(H,R)$  is a  quasitriangular bialgebra (Hopf algebra),
then $(H^*, C_{H^*,H^*}R^*)$ is a coquasitriangular bialgebra (Hopf
algebra).
\end {Proposition}

Now, we see the  relation between Drinfeld  double in this chapter
and the one defined by Drinfeld in     ${\cal V}ect_f(k)$
 of ordinary finite-dimensional vector spaces over
field $k$ with ordinary twist map. Obviously, ${\cal V}ect_f(k)$  is
a symmetric braided tensor category with a left duality (see \cite
[Example 1, P347] {Ka95}). For any $H\in {\cal V}ect_f(k)$, it has
been known that   $ Hom_k(H,k)$ has  two kinds of bialgebra
 structures.
On one hand, it has a bialgebra structure under convolution.
 In this case, we write  $H^{\hat *}$  for $Hom_k(H,k)$. On the other hand,
 $Hom_k(H,k)$  is  a duality of $H$ in braided tensor category ${\cal V}ect_f(k)$.
  Therefore it is a bialgebra. In this case, we write $H^*$
 for $Hom_k(H,k).$
 The two bialgebras have the following
  relation:
$$H ^* = (H^{op \ \ cop})^{\hat *}  \hbox { \ \ \ \ \ as bialgebra. }$$
Therefore, the definition of  Drinfeld double in Theorem \ref
{5.3.9} is the same

as what Drinfeld did since $(H^{* })^{op} = (H^{\hat *})^{cop}$.

Considering   Theorem \ref {5.3.10} and \cite [Proposition 10.3.14]
{Mo93}, we obtain $(D(H^{\hat *}))^{\hat *} = CD(H)$   if we
identify
 $H^{\hat * \hat *}$  as $H$.

\begin {Remark} \label{r3}  { \  }\end {Remark} (It is possible that ${\cal C}$ is not
symmetric) If $U$ and $V$ have left duality $U^*$ and $V^*$,
respectively, then $U^*\otimes V^*$ and $V^* \otimes U^*$ both are
the left dualities of $U\otimes V$. Their evaluations and
coevaluations are $$ d_{U\otimes V}=(d_{U} \otimes d_V) (id _{U^*}
\otimes C_{V^*, U}\otimes  id _V), \ \ \ \ b_{U\otimes V}=(id
_U\otimes (C_{V, U^*})^{-1}\otimes id_ {V^*})(b_U \otimes b_V);$$\
$$ d_{U\otimes V}=d_V (id _{V^*}\otimes d_U \otimes id _V ) ,
\ \ \ \ b_{U\otimes V}=(id _{U}\otimes b_V \otimes id_{U^*})b_U,
$$ respectively.  $H^{\hat * }$ is the left duality of braided bialgebra $H$ under the
first;
 $H^{ * }$ is the left duality of  braided bialgebra $H$ under the second.

If both $V^*$ and $V^{*'}$ are  left dualities of $V$, then $
V^*\cong V^{*'}$ in ${\cal C}$. In fact, assume $d_V : V^*\otimes
V\rightarrow I $ and $b_V : I \rightarrow  V\otimes V^*$ are an
evaluation and a coevaluation of $V$, respectively;  $d_V ':
V^{*'}\otimes V\rightarrow I $ and $b_V ' : I \rightarrow  V\otimes
V^{*'}$ are another evaluation and another coevaluation  of $V$,
respectively. Let $\phi = (d_V \otimes id _{V^*{}'})(id _{V^*}
\otimes b_V'): V^* \rightarrow V^*{}'$ and $\psi = (d_V ' \otimes id
_{V^*})(id _{V^*{}'} \otimes b_V ) : \ V^*{}' \rightarrow V^*{}.$ It
is clear that $\phi \psi = id _{V^*{}'}$ and $\psi \phi = id
_{V^*}$. Thus $V^* \cong V^*{}'.$

\begin {Remark} \label{r4}  { \  }\end {Remark}
If $H$ has a left duality $H^*$ in braided tensor category ${\cal
C}$ (It is possible that ${\cal C}$ is not symmetric ), then the
 following conditions are equivalent:

(i) $H$ is of a symmetric evaluation, i.e. $$(id _U \otimes d _H)
(C_{H^*, U} \otimes id_ H)  = (d_H \otimes id _U) (id _{H^*} \otimes
C_{U, H})$$ for $U = H, \ H^*$.;

(ii) $C_{U,V} = C_{V,U}^{-1}$ for $U, V = H, \ H^*$; \ \ \ \ \ \ \ \
\ \

(iii)  $C_{H,H} = C_{H,H}^{-1}$;

(iv)  $C_{H^*, H^*} = C_{H^*, H^*}^{-1}$.

\begin {Remark} \label{r5}  { \  }\end {Remark}
Theorem \ref {5.3.9} can be stated as follows:

 Let $H$ be
a finite  Hopf algebra (i.e. $H$ has a left duality $H^*$) with
$C_{H,H}= C_{H,H}^{-1}$ (It is possible that ${\cal C}$ is not
symmetric ). Set $A= (H^*)^{op}$ and $ \tau = d_H C_{H,A}$. Then
$(D(H), [b])$ is a quasi-triangular
 Hopf algebra in ${\cal C}$ with  $[b]= \eta _A \otimes b \otimes \eta
_H$.

{\bf Proof.}
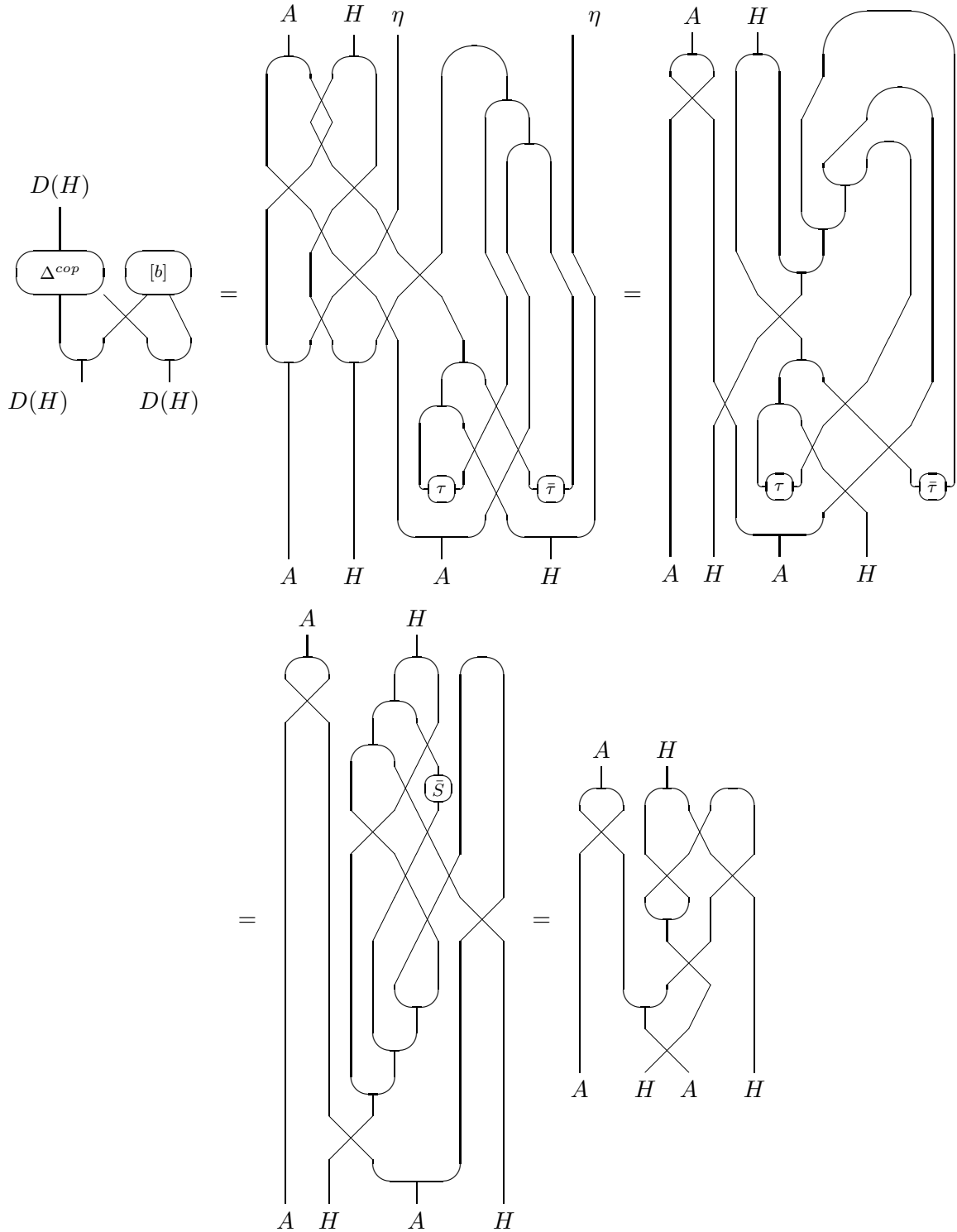
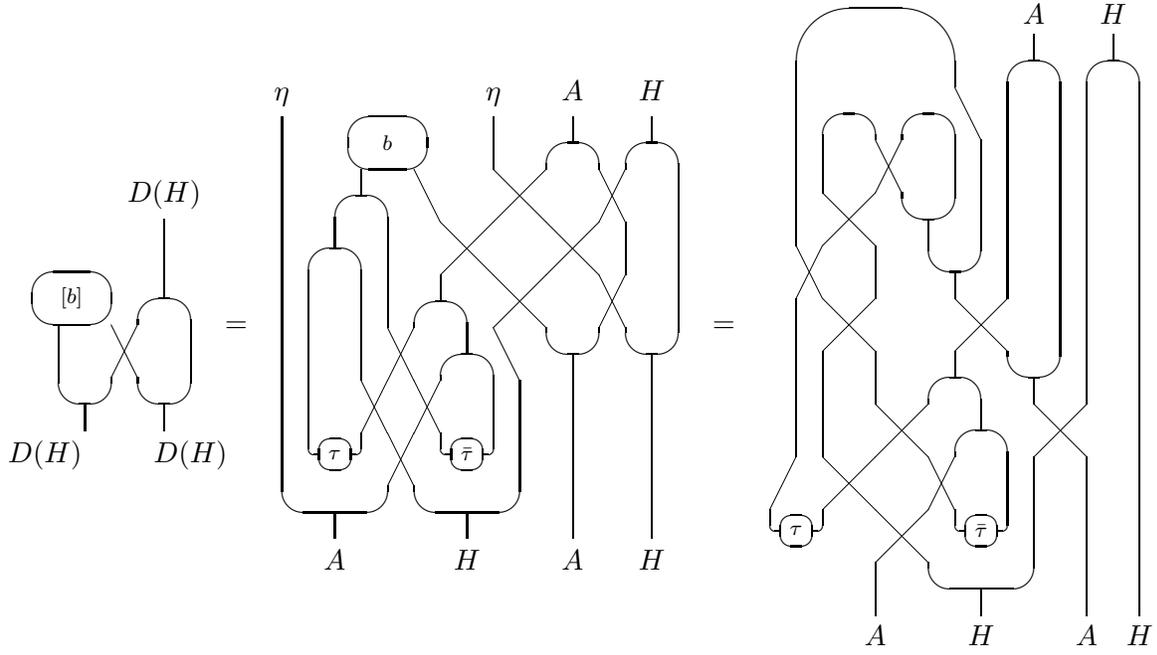
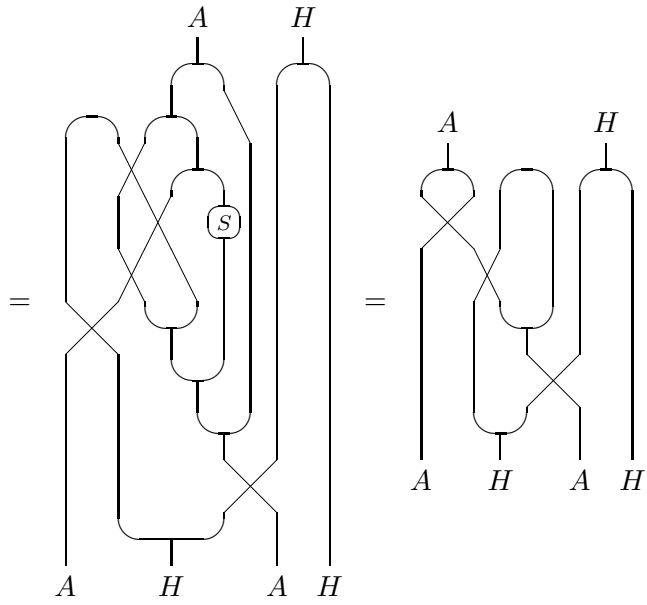
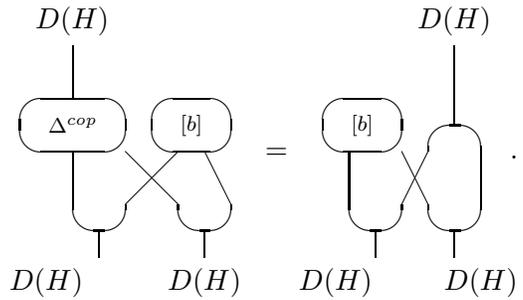
\[
\begin{tangle}
\step[2]\object{D(H)}\\
\step[2]\id\\
\obox 4{\Delta^{cop}}\step\obox 3{[b]}\\
\Step\id\step[2]\XX\step\d\\
\Step\cu\Step\cu\\
\step\object{D(H)}\step[6]\object{D(H)}\\
\end{tangle}
\step=\step
\begin{tangle}
\step\object{A}\step[3]\object{H}\Step\object{\eta}\step[9]\object{\eta}\\
\cd\step\cd\step\id\Step\Coev\step[5]\id\\
\id\Step\X\Step\id\step\id\Step\id\step[2]\cd\Step\id\\
\id\Step\X\Step\id\step\id\Step\id\Step\id\step\cd\step\id\\
\XX\step\XX\step\id\Step\id\Step\id\step\id\Step\id\step\id\\
\id\Step\X\Step\X\Step\id\Step\id\step\id\Step\id\step\id\\
\id\Step\id\step\XX\step\XX\Step\d\d\step\d\d\\
\id\Step\X\Step\X\Step\d\Step\id\step\id\Step\id\step\id\\
\cu\step\cu\step\id\step[2]\cd\step\id\step\id\Step\id\step\id\\
\step\id\step[3]\id\Step\id\step\cd\step\X\step\id\Step\id\step\id\\
\step\id\step[3]\id\Step\id\step\id\Step\X\step\X\Step\id\step\id\\
\step\id\step[3]\id\Step\id\step\coro \tau \step\X\step\coro {\bar \tau }\step\id\\
\step\id\step[3]\id\Step\Cu\step\Cu\\
\step\object{A}\step[3]\object{H}\step[4]\object{A}\step[5]\object{H}\\
\end{tangle}
\step=\step
\begin{tangle}
\step\object{A}\step[3]\object{H}\\
\cd\step\cd\Step\COEV\\
\XX\step\id\Step\id\step\dd\Step\Coev\step[3]\id\\
\id\Step\id\step\id\Step\id\step\id\Step\ne2\coev\step\id\step\id\\
\id\Step\id\step\id\Step\id\step\id\step\cu\Step\id\step\id\step\id\\
\id\Step\id\step\id\Step\id\step\cu\step[3]\id\step\id\step\id\\
\id\Step\id\step\d\step\cu\step[4]\id\step\id\step\id\\
\id\Step\id\Step\XX\step[4]\dd\step\id\step\id\\
\id\Step\id\step\dd\step\cd\Step\dd\Step\id\step\id\\
\id\Step\X\step\cd\step\XX\Step\dd\step\id\\
\id\Step\id\step\id\step\id\Step\X\Step\XX\Step\id\\
\id\Step\id\step\id\step\coro \tau \step\XX\Step\coro {\bar \tau }\\
\id\Step\id\step\Cu\Step\id\\
\object{A}\Step\object{H}\step[3]\object{A}\step[4]\object{H}\\
\end{tangle} \]
\[
\step=\step
\begin{tangle}
\step\object{A}\step[5]\object{H}\\
\cd\step[3]\cd\step\coev\\
\XX\Step\cd\step\id\step\id\Step\id\\
\id\Step\id\step\cd\step\X\step\id\Step\id\\
\id\Step\id\step\id\Step\X\step\O {\bar{S}}\step\id\Step\id\\
\id\Step\id\step\XX\step\X\step\id\Step\id\\
\id\Step\id\step\id\Step\X\step\X\Step\id\\
\id\Step\id\step\id\step\dd\step\X\step\XX\\
\id\Step\id\step\id\step\id\step\dd\step\id\step\id\Step\id\\
\id\Step\id\step\id\step\id\step\cu\step\id\Step\id\\
\id\Step\id\step\id\step\cu\Step\id\Step\id\\
\id\Step\id\step\cu\step[3]\id\Step\id\\
\id\Step\XX\step[4]\id\Step\id\\
\id\Step\id\Step\Cu\Step\id\\
\object{A}\Step\object{H}\step[4]\object{A}\step[4]\object{H}\\
\end{tangle}
\step=\step
\begin{tangle}
\step\object{A}\step[3]\object{H}\\
\cd\step\cd\step\coev\\
\XX\step\id\Step\X\Step\id\\
\id\Step\id\step\XX\step\XX\\
\id\Step\id\step\cu\step\id\Step\id\\
\id\Step\id\Step\XX\Step\id\\
\id\Step\cu\step\dd\Step\id\\
\id\step[3]\XX\step[3]\id\\
\object{A}\step[3]\object{H}\Step\object{A}\step[3]\object{H}\\
\end{tangle}
\]

and
\[
\begin{tangle}
\step[5]\object{D(H)}\\
\step[5]\id\\
\obox 3{[b]}\step\cd\\
\step\id\step[2]\X\Step\id\\
\step\cu\step\cu\\
\step[0.5]\object{D(H)}\step[5.5]\object{D(H)}\\
\end{tangle}
\step=\step
\begin{tangle}
\object{\eta}\step[8]\object{\eta}\step[3]\object{A}\step[3]\object{H}\\
\id\step[2.5]\obox 3{b}\step[2.5]\id\Step\cd\step\cd\\
\id\Step\cd\step\d\Step\XX\Step\X\Step\id\\
\id\step\cd\step\id\Step\XX\Step\XX\step\id\Step\id\\
\id\step\id\Step\id\step\id\step\cd\step\XX\Step\X\Step\id\\
\id\step\id\Step\id\step\X\step\cd\d\step\cu\step\cu\\
\id\step\id\Step\X\step\X\Step\id\step\id\Step\id\step[3]\id\\
\id\step\coro \tau \step\X\step\coro {\bar \tau }\step\id\Step\id\step[3]\id\\
\Cu\step\Cu\Step\id\step[3]\id\\
\Step\object{A}\step[5]\object{H}\step[4]\object{A}\step[3]\object{H}\\
\end{tangle}
\step=\step
\begin{tangle}
\step \step[9]\object{A}\step[3]\object{H}\\
\step \COEV\step[6]\cd\step\cd\\
\step \id\step\coev\step\coev\d\step\id\Step\id\step\id\Step\id\\
\step \id\step\id\Step\X\Step\id\step\id\step\id\Step\id\step\id\Step\id\\
\step \id\step\XX\step\cu\step\id\step\id\Step\id\step\id\Step\id\\
\step \X\Step\id\Step\cu\step\id\Step\id\step\id\Step\id\\
\step \id\step\XX\step[3]\XX\Step\id\step\id\Step\id\\
\step \id\step\id\Step\id\Step\cd\step\cu\step\id\Step\id\\
\step \id\step\id\Step\XX\step\cd\step\XX\Step\id\\
\dd\step\XX\Step\X\Step\id\step\id\Step\id\Step\id\\
\coro \tau \Step\XX\step\coro {\bar \tau }\step\id\Step\id\Step\id\\
\step \step[3]\id\Step\Cu\Step\id\Step\id\\
\step \step[3]\object{A}\step[4]\object{H}\step[4]\object{A}\Step\object{H}\\
\end{tangle}\]
\[
\step=\step
\begin{tangle}
\step[5]\object{A}\step[4]\object{H}\\
\step[4]\cd\Step\cd\\
\coev\step\cd\step\d\step\id\Step\id\\
\id\Step\X\step\cd\step\id\step\id\Step\id\\
\id\Step\id\step\X\Step\O S\step\id\step\id\Step\id\\
\id\Step\X\step\d\step\id\step\id\step\id\Step\id\\
\XX\step\cu\step\id\step\id\step\id\Step\id\\
\id\Step\id\Step\cu\step\id\step\id\Step\id\\
\id\Step\id\step[3]\cu\step\id\Step\id\\
\id\Step\id\step[4]\XX\Step\id\\
\id\Step\Cu\Step\id\Step\id\\
\object{A}\step[4]\object{H}\step[4]\object{A}\Step\object{H}\\
\end{tangle}
\step=\step
\begin{tangle}
\step\object{A}\step[6]\object{H}\\
\cd\step\coev\step\cd\\
\XX\step\id\Step\id\step\id\Step\id\\
\id\Step\X\Step\id\step\id\Step\id\\
\id\Step\id\step\cu\step\id\Step\id\\
\id\Step\id\Step\XX\Step\id\\
\id\Step\cu\Step\id\Step\id\\
\object{A}\step[3]\object{H}\step[3]\object{A}\step[2]\object{H}\\
\end{tangle}\ \ \ .
\]

Thus
\[
\begin{tangle}
\step[2]\object{D(H)}\\
\step[2]\id\\
\obox 4{\Delta^{cop}}\step\obox 3{[b]}\\
\Step\id\step[2]\XX\step\d\\
\Step\cu\Step\cu\\
\step\object{D(H)}\step[6]\object{D(H)}\\
\end{tangle}
\step=\step
\begin{tangle}
\step[5]\object{D(H)}\\
\step[5]\id\\
\obox 3{[b]}\step\cd\\
\step\id\step[2]\X\Step\id\\
\step\cu\step\cu\\
\step[0.5]\object{D(H)}\step[5.5]\object{D(H)}\\
\end{tangle}\ \ \ .
\]
We complete the proof.\ \
\begin{picture}(5,5)
\put(0,0){\line(0,1){5}}\put(5,5){\line(0,-1){5}}
\put(0,0){\line(1,0){5}}\put(5,5){\line(-1,0){5}}
\end{picture}\\

          \section {(Co) quasitriangular Hopf algebras of double
          cross (co)products}\label {s14}
In this section, we study the relation among (co)quasitriangular
structures of the double cross coproduct $A \bowtie ^R H$ ( double
cross product $A \bowtie _r H$ ) and $A$ and $H$. We omit all of the
proof in this section  since we  can  obtain them
 by turning the proof
in \cite {Ch98} into diagrams.

Assume that
\begin {eqnarray*}
&r& : A \otimes H \rightarrow I  \hbox { \ \ and \ \ }
R :  I \rightarrow A \otimes  H,\\
\end {eqnarray*}
 are morphisms in ${\cal C}.$

$R $ is called a weak $R$-matrix  of $A\otimes H$ if $R$ is an
invertible copairing of $A \otimes H$. $r $ is called a weak
$r$-comatrix  on $A\otimes H$ if $r$ is an invertible pairing on $A
\otimes H$. Obviously,  $R $ is a weak $R$-matrix  of $A\otimes H$
iff $ C_{A, H} R $ is an invertible skew copairing of $H \otimes A$.
$r C_{A, H} $ is  a weak $r$-comatrix  on $H\otimes A$ iff $r$ is an
invertible copairing on $A \otimes H$.

Therefore, if  $r $  is a weak $r$-comatrix on $A\otimes H$ and let
$\tau = r C_{A,H}$, then $A \bowtie _\tau H$ is a bialgebra by
Theorem \ref  {5.3.3}. For convenience, we also denote $A \bowtie
_\tau H $  by $A \bowtie _r H.$

\begin {Lemma}\label {5.4.1}
           Let $R $  be a weak $R$-matrix of $A\otimes H$, and define

\[ \phi =
\begin{tangle}

\step [4] \object{A}\\

\ro R \step [2] \id \step [2] \ro  {\bar R }\\

\XX \step [2] \id \step [2] \XX\step [2] \\

\id \step [2] \id \step[2] \XX \step [2]\id  \\

\id  \step [2] \XX \step [2] \cu \\

\cu \step [2] \nw1 \step [2] \id \\

\step \id \step [2]  \step [2] \cu \\

\step \object{H}\step [5] \object{A}\\

\end{tangle}
 \ \ \hbox { and } \ \
\psi = \begin{tangle}

\step [4] \object{H}\\

\ro R \step [2] \id \step [2] \ro  {\bar R }\\

\XX \step [2] \id \step [2] \XX\step [2] \\

\id \step [2] \XX\step[2] \id \step [2]\id  \\

\cu \step [2] \XX \step [2] \id \\

\step \id \step [2] \ne1 \step [2] \cu \\

\step \cu \step  \step [2]  \step [1] \id \\

\step[2] \object{H}\step [5] \object{A}\\

\end{tangle} \ \ .
\]
 Then $A {}^\phi
\bowtie ^\psi H$  is a bialgebra. In this case, we denote   $A
{}^\phi \bowtie ^\psi H$ by $A  \bowtie ^R H$. Furthermore,
                                $A {}^\phi \bowtie ^\psi H$
                                is   a Hopf algebra
                                if $A$ and $H$  are  Hopf algebras.

\end {Lemma}

Considering the universal property of double bicrossproducts in
Lemma \ref {2.3.2}, we have

\begin {Lemma}\label {5.4.5} Let $D = A \bowtie ^R H$ be a bialgebra.
 Assume that $p_A : B \longrightarrow A$
and $p_H : B \longrightarrow H$ are bialgebra morphisms. Then $p =
(p_A \otimes p_H) \Delta _B$  is a bialgebra morphism from $B
\rightarrow D$  iff (UT):
 \[
 \begin{tangle}

\step[2] \step [3] \object{B}\\

\ro R\step [2] \cd  \\

 \id \step [2] \id \step [2]\O {p_A}  \step [2] \O {p_H}  \\

\id \step [2] \XX  \step [2] \id \step [2]\\

\cu \step [2] \cu   \\

\step [1] \object{A}\step [4] \object{H}\\
\end{tangle} \ \ = \ \
 \begin{tangle}

\step[1] \object{B}\\

 \cd\step [2]\ro R  \\

 \O {p_H}  \step [2] \O {p_A} \step [2] \id \step [2] \id \step [2]\\

\XX \step [2] \id  \step [2] \id \step [2]\\

\id \step [2] \XX  \step [2] \id \step [2]\\

\cu \step [2] \cu   \\

\step [1] \object{A}\step [4] \object{H}\\
\end{tangle} \ \ .
 \]

\end {Lemma}

\begin {Corollary}\label {5.4.6}
Let $R $ be a weak $R$-matrix of $H \otimes H.$  Then $(H,R)$ is  a
quasitriangular bialgebra iff the comultiplication $\Delta$ of $H$
is a bialgebra morphism from $H$ to $H\otimes H.$
                        \end {Corollary}

\begin {Theorem}\label {5.4.7}
$D= A \bowtie ^RH$ has the following universal property with respect
 to $(\pi _A,
\pi_H)$.
 Assume that $p_A : B \longrightarrow A$
and $p_H : B \longrightarrow H$ are bialgebra morphisms. If $p_A $
and $p_H$   satisfy $(UT)$ defined in Lemma \ref {5.4.5}, then there
exists a unique bialgebra morphism $p$ from $B$ to $D$ such that
$\pi _A p = p_A$ and  $\pi _H p = p_H.$

\end {Theorem}

\begin {Theorem}\label {5.5.1}
If $A {}^\phi \bowtie ^ \psi H$  admits a quasitriangular structure,
then $A$ and $H$  admit  quasitriangular structures, and there
exists a weak $R$-matrix $R$ of $A \otimes H$ such that $A {}^\phi
\bowtie ^ \psi H = A \bowtie ^R H$ and

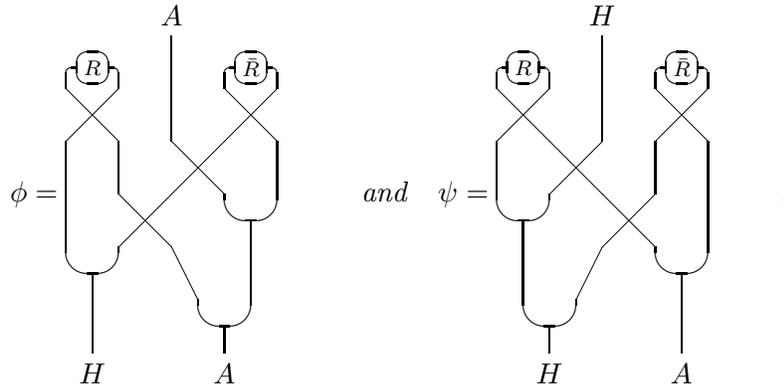
\[ \phi =
\begin{tangle}

\step [4] \object{A}\\

\ro R \step [2] \id \step [2] \ro  {\bar R }\\

\XX \step [2] \id \step [2] \XX\step [2] \\

\id \step [2] \id \step[2] \XX \step [2]\id  \\

\id  \step [2] \XX \step [2] \cu \\

\cu \step [2] \nw1 \step [2] \id \\

\step \id \step [2]  \step [2] \cu \\

\step \object{H}\step [5] \object{A}\\

\end{tangle}
 \ \ \hbox { and } \ \
\psi = \begin{tangle}

\step [4] \object{H}\\

\ro R \step [2] \id \step [2] \ro  {\bar R }\\

\XX \step [2] \id \step [2] \XX\step [2] \\

\id \step [2] \XX\step[2] \id \step [2]\id  \\

\cu \step [2] \XX \step [2] \id \\

\step \id \step [2] \ne1 \step [2] \cu \\

\step \cu \step  \step [2]  \step [1] \id \\

\step[2] \object{H}\step [5] \object{A}\\

\end{tangle} \ \ .
\]

\end {Theorem}
If \begin {eqnarray*}
  &U& : I \rightarrow A \otimes H , \hbox { \ \ \ \ }
  V : I \rightarrow A \otimes H ,  \\
  &P& : I \rightarrow A \otimes A , \hbox { \ \ \ \ }
  Q : I \rightarrow H \otimes H ,    \\
  \end {eqnarray*}
 are morphisms in ${\cal C},$
 we define

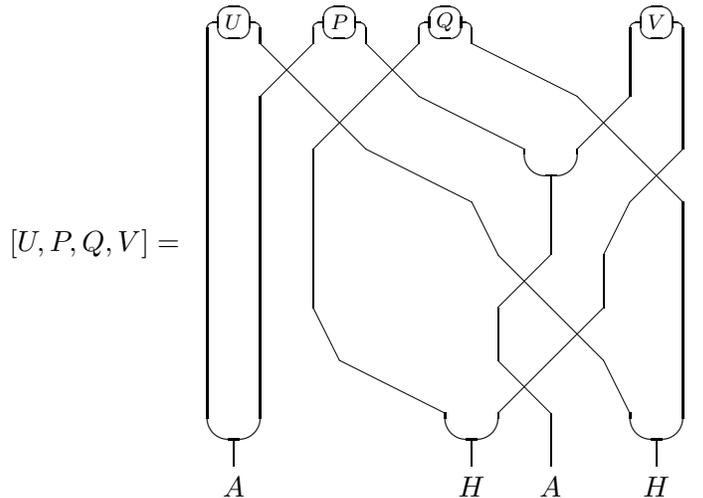
\[ [U,P,Q,V] = \ \
 \begin{tangle}

 \ro U\step [2]\ro P  \step [2]\ro Q\step [6]\ro V \\

\id  \step [2] \XX\step [2] \XX\step [2] \nw4 \step [5] \id \step [2] \id \\

\id  \step [2] \id \step [2] \XX\step [2]\nw4  \step [5]  \XX  \step [2] \id \\

\id  \step [2] \id \step [2] \id \step [2] \nw4 \step [5] \cu \step [2]  \XX \\

\id  \step [2] \id \step [2] \id \step [6] \nw1\step [2] \id \step [2] \ne1\step [2]\id   \\

\id  \step [2] \id \step [2] \id \step [7] \XX \step [2] \id  \step [3]  \id\\

\id  \step [2] \id   \step [2] \nw1 \step [5] \step \id   \step [2] \XX  \step [3] \id  \\

\id  \step [2] \id \step [3] \nw4  \step [5] \XX \step [2]\nw1 \step [2] \id \\

\cu  \step [2] \step [3]  \step [2] \cu \step [2] \id  \step [3] \cu\\

 \step \object{A}\step [9] \object{H} \step [3] \object{A} \step [4] \object{H}

\end{tangle} \ \ .
 \]

From now on, unless otherwise stated, $R$ is a fixed weak $R$-matrix
of $A \otimes H.$

\begin {Theorem}\label {5.5.2}
If $(A,P)$ and $(H,Q)$  are quasitriangular bialgebras,
then \\
$(A\bowtie ^R H, [R, P,Q,R^{-1}])$  is a quasitriangular bialgebra.

\end {Theorem}

\begin {Corollary}\label {5.5.3}
 $A {}^\phi \bowtie ^ \psi H$  admits a quasitriangular structure iff
$A$ and $H$  admit a quasitriangular structure, and there exists a
weak $R$-matrix $R$ of $A \otimes H$ such that $A {}^\phi \bowtie ^
\psi H = A \bowtie ^R H$.
\end {Corollary}

\begin {Corollary}\label {5.5.4} Let $R$  be a weak
$R$-matrix of $H\otimes H.$ Then $(H \bowtie ^R H,  [R,R,R,R^{-1}])$
is a quasitriangular bialgebra iff $(H,R)$  is a quasitriangular
bialgebra.
\end {Corollary}

\begin {Lemma}\label {5.5.5}
Let $R$  be a weak $R$-matrix of $H\otimes H$. Then the following
conditions are equivalent.

(i) $(H,R)$ is a triangular bialgebra.

(ii)  $(A \bowtie ^R H, [R, R, R, R^{-1}])$  is a triangular
bialgebra.

(iii)  $\Delta _H: H \longrightarrow H\bowtie ^R H$ is a bialgebra
morphism and $[R,R,R,R^{-1}] = (\Delta _H\otimes \Delta _H) (R)$,
i.e. $\Delta _H$ is a quasitriangular morphism.
\end {Lemma}

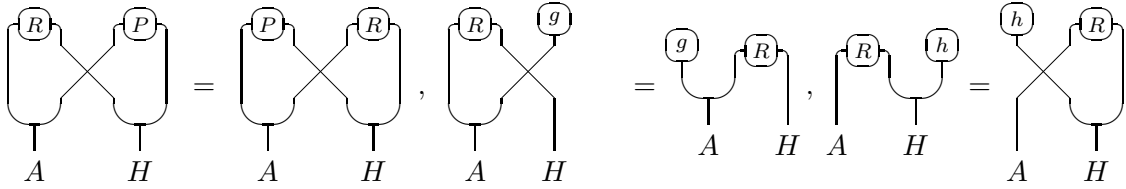
\[
 \begin{tangle}

 \ro R  \step [2]\ro P  \\

\id  \step [2] \XX\step [2] \id \\

\cu \step [2] \cu \\
\step \object{A}\step [4] \object{H}
\end{tangle} \ \ = \ \
\begin{tangle}

 \ro P \step [2]\ro R  \\

\id  \step [2] \XX\step [2] \id \\

\cu \step [2] \cu \\
 \step \object{A}\step [4] \object{H}
\end{tangle} \ \ , \ \
\begin{tangle}

 \ro R \step [2]\Q g \\

\id  \step [2] \XX\step [2]\\

\cu \step [2]\id \\
\step \object{A}\step [3] \object{H}
\end{tangle} \ \ = \ \
\begin{tangle}

 \Q g  \step [2]\ro R\\

\cu\step [2] \id\\
\step  \object{A}\step [3] \object{H}
\end{tangle} \ \ , \ \
\begin{tangle}

\ro R  \step [2] \Q h\\

\id\step [2] \cu \\
 \object{A}\step [3] \object{H}
\end{tangle} \ \ = \ \ \begin{tangle}

\Q h   \step [2]\ro R\\

\XX \step [2] \id \step [2]\\

\id \step [2] \cu\\
 \object{A}\step [3] \object{H}
\end{tangle} \ \  \ \
\] for any morphisms $P :
I\rightarrow A\otimes H , \ \ \ g : I \rightarrow A,  \ \ \ h : I
\rightarrow H. $

Let $$QT(H) = \{P \mid P \hbox { \  a quasitriangular structure of }
H\}.$$
\begin {Lemma}\label {5.5.7}
 $\Phi $   is a bijective map
from $CW(A,H) \times CW(A,A) \times CW(H,H) \times CW(A,H)$    onto
$CW(A \bowtie ^RH, A \bowtie ^RH),$ where $\Phi (U,P,Q,V) =
[U,P,Q,V].$
\end {Lemma}

\begin {Lemma}\label {5.5.8}
Let $R$ and $\bar R$ be two weak $R$-matrices  of $A\otimes H$. Then
$A \bowtie ^R H = A \bowtie ^{\bar R}H$ iff there exists $U \in
CW(A,H)$  such that $\bar R = RU.$
\end {Lemma}

\begin {Theorem}\label {5.5.9}
Assume $(A,P)$ and $(H,Q)$  are quasitriangular bialgebras and $U, V
\in CW(A,H)$. \noindent Then $(A\bowtie ^R H, [R, P, Q, R^{-1}][U,
\eta , \eta , V])$ is a quasitriangular bialgebra. Conversely, if
$(A \bowtie ^R H, {\cal R})$ is a quasitriangular bialgebra, then
there exist $P, Q, U, V$ such that $(A,P)$  and $(H,Q)$ are
quasitriangular bialgebras with ${\cal R} =  [R,P,Q,R^{-1}][U, \eta
, \eta , V]$ and
 $U, V \in CW(A,H).$

\end {Theorem}

\begin {Proposition}\label {5.5.10}
There exists a bijective map from
 $QT(A) \times QT(H) \times CW(A, H) \times
CW(A,H)$ onto $QT(A \bowtie ^R H)$ by sending $(P,Q,U,V)$ to $
[R,P,Q,R^{-1}][U, \eta _H, \eta _A, V]$.

\end {Proposition}

 Dually, we can get the followings about double cross coproducts.

{Lemma} \ref {5.4.1}${}^{\circ}$  \
           Let $r $  be a weak $r$-comatrix on $A\otimes H$, and define

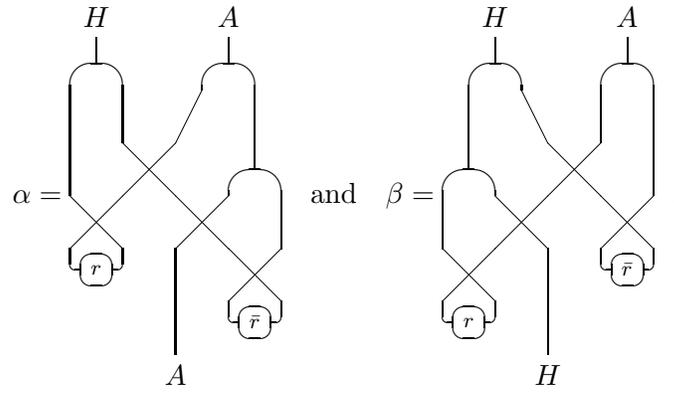
\[ \alpha =
\begin{tangle}

\step  \object{H}\step [5] \object{A}\\
 \cd \step [3]  \cd \\

\id \step [2] \id \step [2] \ne1   \step [2] \id \\

\id \step [2] \XX  \step [2] \cd \\

\XX\step [2] \XX \step[2] \id  \\

\coro r  \step [2]\id \step [2] \XX \\

\step [4] \id \step [2] \coro  {\bar r} \\

\step [4] \object{A}\\
\end{tangle}
 \ \ \hbox { and } \ \
\beta =
\begin{tangle}

\step[2]  \object{H}\step [5] \object{A}\\
\step \cd \step [3]  \cd \\

\step \id \step [2]\nw1  \step [2] \id   \step [2] \id \\

\cd \step [2] \XX  \step [2] \id \\

\id \step [2] \XX \step[2] \XX  \\

 \XX \step [2]\id \step [2] \coro {\bar r} \\

\coro  { r}\step [2]   \id\step [2] \\

\step [4] \object{H}\\
\end{tangle} \ \ .
 \]
Then $A {}_\alpha \bowtie _\beta H$  is a bialgebra. In this case,
we denote   $A {}_\alpha \bowtie _\beta H$ by $A  \bowtie _r H$.
Furthermore,
                                $A {}_\alpha  \bowtie _\beta H$
                                is   a Hopf algebra
                                if $A$ and $H$  is  Hopf algebras.

{Lemma} \ref {5.4.5}${}^{\circ}$ \  Let $D = A \bowtie _r H$ be a
bialgebra.
 Assume that $j_A : A \longrightarrow B$
and $j_H : H \longrightarrow B$ are bialgebra morphisms. Let $j =
m_B(j_A \otimes j_H) $  be a bialgebra morphism  from $D$ to $B$ iff

(CUT):

\[ \begin{tangle}

\step \object{A}\step [4] \object{H}\\

\cd    \step [2]\cd\\

\id \step [2]\XX \step [2] \id \\

\coro r \step [2]  \O {j_A} \step [2] \O {j_H} \\

 \step [4]\cu \\

\step [5]\object{B}
\end{tangle} \ \ = \ \
\begin{tangle}

\step \object{A}\step [4] \object{H}\\

\cd    \step [2]\cd\\

\id \step [2]\XX \step [2] \id \\

\XX  \step [2]\id  \step [2] \id \\

  \O {j_H} \step [2] \O {j_A}\step [2]\coro r  \\

 \cu \step [4]\\

\step [1]\object{B}
\end{tangle} \ \ .
\]

{Corollary} \ref {5.4.6}${}^{\circ}$ \ Let $r $ be a weak
$r$-comatrix on $H \otimes H.$  Then $(H,r)$ is  a coquasitriangular
bialgebra iff the multiplication $m _H$ of $H$ is a bialgebra
morphism from $H\otimes H$ to $H$.

{Theorem} \ref
 {5.4.7}${} ^{\circ}$ \
$D= A \bowtie _rH$ has the following universal property with respect
to $(i _A, i_H)$.
 Assume that $j_A : A \longrightarrow B$
and $j_H : H \longrightarrow B$ are bialgebra morphisms. If $j_A $
and $j_H$   satisfy $(CUT)$ defined in Lemma \ref {5.4.5}, then
there exists a unique bialgebra morphism $j$ from $D$ to $B$ such
that $ji _A  = j_A$ and  $ji _H  = j_H.$

{Theorem}\ref {5.5.1}${}^{\circ}$  \ If $A {}_\alpha  \bowtie _\beta
H$  admits a coquasitriangular structure, then $A$ and $H$  admit
coquasitriangular structures, and there exists a weak $r$-comatrix
$r$ on $A \otimes H$ such that $A {}_\alpha  \bowtie _\beta H = A
\bowtie _r H$ and

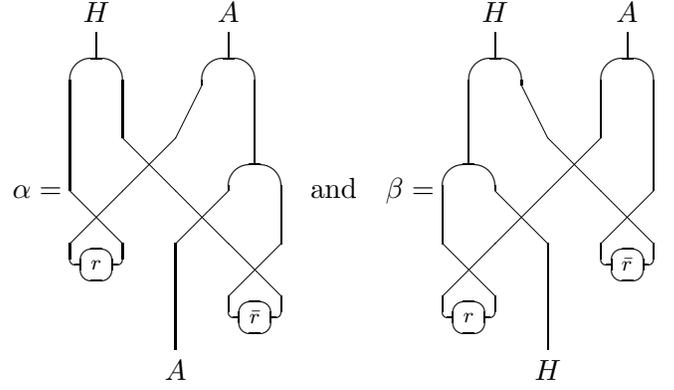
\[ \alpha =
\begin{tangle}

\step  \object{H}\step [5] \object{A}\\
 \cd \step [3]  \cd \\

\id \step [2] \id \step [2] \ne1   \step [2] \id \\

\id \step [2] \XX  \step [2] \cd \\

\XX\step [2] \XX \step[2] \id  \\

\coro r  \step [2]\id \step [2] \XX \\

\step [4] \id \step [2] \coro  {\bar r} \\

\step [4] \object{A}\\
\end{tangle}
 \ \ \hbox { and } \ \
\beta =
\begin{tangle}

\step[2]  \object{H}\step [5] \object{A}\\
\step \cd \step [3]  \cd \\

\step \id \step [2]\nw1  \step [2] \id   \step [2] \id \\

\cd \step [2] \XX  \step [2] \id \\

\id \step [2] \XX \step[2] \XX  \\

 \XX \step [2]\id \step [2] \coro {\bar r} \\

\coro  { r}\step [2]   \id\step [2] \\

\step [4] \object{H}\\
\end{tangle} \ \ .
 \]

If
\begin {eqnarray*}
u : A \otimes H \rightarrow I , \hbox { \ \ \ \ }
v : A \otimes H \rightarrow I , \\
p : A \otimes A \rightarrow I , \hbox { \ \ \ \ } q : H \otimes H
\rightarrow I ,
\end {eqnarray*}
are morphisms, we set

\[ [u, p, q, v] = \ \
 \begin{tangle}
\step \object{A}\step [4] \object{H} \step [4] \object{A} \step [4]
\object{H}\\

\cd   \step [2] \cd \step [2] \cd \step [2] \cd \\

\id  \step [2] \id \step [2]  \id \step [2]\XX  \step [2]  \XX  \step [2] \id \\

\id  \step [2] \id \step [2]  \XX  \step [2]\XX  \step [2]  \XX  \step [2]  \\

\id  \step [2] \coro p \step [2]  \XX  \step [2]\XX  \step [2]  \id  \\

\id  \step [4] \step [2]  \id   \step [2]\coro q  \step [2]  \XX  \step [2] \\

\id  \step [5] \ne2 \step [6]  \coro v\\

\coRo u  \step [4]
\end{tangle} \ \ .
 \]

{Theorem} \ref {5.5.2}${} ^{\circ}$ \ Assume $(A,p)$ and $(H,q)$ are
coquasitriangular bialgebras. Then $(A\bowtie _r H, [r,p,q,r^{-1}])$
is a coquasitriangular bialgebra.

 {Corollary} \ref
{5.5.3}${} ^{\circ} $ \
 $A {}_\alpha  \bowtie _\beta H$  admits a coquasitriangular structure iff
$A$ and $H$  admit a coquasitriangular structure, and there exists a
weak $r$-comatrix $r$ of $A \otimes H$ such that $A {}_\alpha
\bowtie _\beta  H = A \bowtie _r H$.

 {Corollary} \ref
{5.5.4}${} ^{\circ}$ \  Let $r$  be a weak $r$-comatrix on $H\otimes
H$. Then $(H \bowtie _r H, [r,r,r,r^{-1}] )$  is a coquasitriangular
bialgebra iff $(H,r)$  is a coquasitriangular bialgebra.

 {Lemma} \ref
{5.5.5}${}^{\circ}$ \  Let $r$  be a weak $r$-comatrix on $H\otimes
H$. Then  the following are equivalent.

(i) $(H,r)$ is a cotriangular bialgebra.

(ii)  $(A \bowtie _r H, [r,r,r,r^{-1}])$  is a cotriangular
bialgebra.

(iii)  $m _H: H\bowtie _r H \longrightarrow  H$ is a bialgebra
morphism and $ [r,r,r,r^{-1}] (m_H\otimes m_H)= r$, i.e. $m_H$ is a
coquasitriangular morphism.

Let $CCW(A \otimes H)$  denote the set of all weak $r$-comatrix $r$
on $A \otimes H$ satisfying the following conditions:

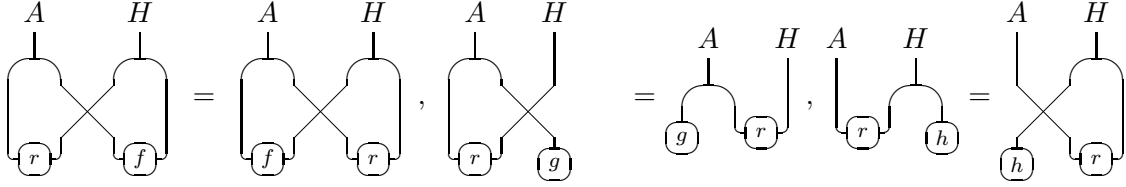
\[
 \begin{tangle}

\step \object{A}\step [4] \object{H}\\
\cd \step [2] \cd \\

\id  \step [2] \XX\step [2] \id \\

\coro r  \step [2]\coro f  \\

\end{tangle} \ \ = \ \
\begin{tangle}

\step \object{A}\step [4] \object{H}\\

\cd \step [2] \cd \\

\id  \step [2] \XX\step [2] \id \\

 \coro f \step [2]\coro r  \\

\end{tangle} \ \ , \ \
\begin{tangle}

\step \object{A}\step [3] \object{H}\\

\cd \step [2]\id \\
\id  \step [2] \XX\step [2]\\

 \coro r \step [2]\QQ g \\

\end{tangle} \ \ = \ \
\begin{tangle}

\step \object{A}\step [3] \object{H}\\

\cd\step [2] \id\\
 \QQ g  \step [2]\coro r\\
\end{tangle} \ \ , \ \
\begin{tangle}

\object{A}\step [3] \object{H}\\

\id\step [2] \cd \\
\coro r  \step [2] \QQ h\\
\end{tangle} \ \ = \ \ \begin{tangle}

\object{A}\step [3] \object{H}\\

\id \step [2] \cd\\
\XX \step [2] \id \step [2]\\

\QQ h   \step [2]\coro r\\

\end{tangle} \ \  \ \
\] for any morphisms $f :
A\otimes H \rightarrow I, \ \ \ g : A \rightarrow I,  \ \ \ h : H
\rightarrow I $. Let $$CQT(H) = \{ p \mid p \hbox {\ \ is a
coquasitriangular structure of }
 H\}.$$
From now on, unless otherwise stated, $r$ is a fixed weak
$r$-comatrix on $A \otimes H.$

{Lemma} \ref {5.5.7}${} ^{\circ}$ \
 $\phi $   is a bijective map
from $CCW(A,H) \times CCW(A,A) \times CCW(H,H) \times CCW(A,H)$ to
$CCW(A \bowtie _rH, A \bowtie _rH),$  where $\phi (u,p,q,v)
=[u,p,q,v].$

{Lemma} \ref {5.5.8}${}^{\circ}$ \ Let $r$ and $\bar r$ be two weak
$r$-comatrices  on $A\otimes H$. Then $A \bowtie _r H = A \bowtie
_{\bar r}H$ iff there exists $u \in CCW(A,H)$  such that $\bar r =
r*u.$

{Theorem} \ref {5.5.9}${} ^{\circ}$ \ Assume $(A,p)$ and $(H,q)$ are
coquasitriangular bialgebras and $u, v \in CCW(A,H)$.
\noindent Then
$(A\bowtie _r H, [r,p,q,r^{-1}]*[u, \epsilon , \epsilon , v])$ is a
coquasitriangular bialgebra. Conversely, if $(A \bowtie _r H, \bar
r)$ is a coquasitriangular bialgebra, then there exist $p, q, u, v$
such that $(A,p)$  and $(H,q)$ are coquasitriangular bialgebras with
$\bar r =  [r,p,q,r^{-1}]*[u, \epsilon , \epsilon , v]$ and  $u, v
\in CCW(A,H)$.

{Proposition} \ref {5.5.10}${}^{\circ}$ \ There exists a bijective
map from
 $CQT(A) \times CQT(H) \times CCW(A, H) \times
CCW(A,H)$ to $CQT(A \bowtie _r H)$ by sending $(s,t,u,v)$ to
$[r,s,t,r^{-1}]
*[u, \epsilon , \epsilon , v]$.

  \chapter {The Duality Theorem  for Braided Hopf Algebras}\label {c7}

 In this chapter, we give the duality theorem
 for  Hopf
 algebras living in a symmetric braided tensor category. Using the result,
 we show that

   $$
(R \# H)\# H^{\hat *}   \cong M_n(R)  \hbox { \ \ \  as algebras }
$$
when $H$ is a finite-dimensional Hopf algebra living in symmetric
braided
 category determined by (co)triangular structures. In particular,
 the duality theorem holds for any
  finite-dimensional
super-Hopf algebra $H$

Supersymmetry has attracted a great deal of interest from physicists
and mathematicians (see \cite {CNS75}, \cite {MR94} ). It finds
numerous applications in particle physics and statistical mechanics.
Of central importance in supersymmetric theories is the $Z_2$-graded
structure which permeates them. Superalgebras and super-Hopf
algebras are  most naturally
 very important.

               It is well-known that in the work of $C^*$-algebras and
               von Neumann algebras for abelian groups the crossed product
 $R$ crossed by $G$ crossed by $G ^ { \hat  { \  }  } $  is isomorphic to $R$ tensor the
 compact operators. This result does have implications for symmetries in
 quantum theory. One may see for example  Pedersen's book on $C^*$
 -algebras \cite {Pe79}.
 This offers some clues as to why duality theorem might be useful in mathematical
 physics. Its generalization to Hopf-von Neumann algebras was known again
, see for example the book of Stratila on modular theory
 \cite {St81} or the papers on Hopf-von Neumann algebras (or the book on
 Kac algebras) by Enock and Schwartz. Blattner and Montgomery
strip off the functional analysis and  duplicate the result at the
level of Hopf algebras \cite {BM85}. But if one looks more
 deeply into the literature one finds indeed its origins in mathematical
 physics. However, since most algebras in physics are not finite-dimensional
 one does need to say some words about what one thinks the
 operator-algebra version should be. What plays the role of compact operators
 in the super case ?

          In this chapter, we show that the subalgebras $H \bar \otimes
          H^{ \hat *}$  of $End _k H$ plays the role of compact operators
          in the super case.
We obtain the duality theory  in symmetric tensor category. Using
this theory, we study super-Hopf algebras.
 If super(quantum)group $H$ acts on a superalgebra $R$, then the duality of
 the super(quantum)group acts on the crossed  product. We show
 that the repeated  crossed  product is isomorphic to $M_n(R)$
 when $dim { \ } H =n $

\section {   The duality theory for Hopf algebras living in
symmetric tensor categories   }\label {s15}

                            In this section we obtain the duality theory
 for Hopf algebras living in symmetric tensor category ${\cal
 C}$.

 Throughout this section $H$  is a Hopf algebra living in {\cal C}
 with a left duality $H^{ *}$ and
 $(R, \alpha )$  is a left $H$-module algebra in ${\cal C}.$
 Let $H^{\hat * }= (H^ {  * })^{ op \ cop } $  and
 $R\# H$  be an algebra.

It follows from \cite {Ma93a}  that

\begin {Lemma} \label {6.1.2}
If $(R, \alpha )$ is an $H$-module algebra, let $\phi = (id \otimes
\alpha )(b' \otimes id ):$ \ $R \rightarrow  H^{*\ op } \otimes R,$
where $b' = C_{H, H^{* \ op }}b$. Then
$$\phi m =
(m\otimes m )(id \otimes C_{R, H^{*\ op}}\otimes id)(\phi \otimes
\phi ). $$

\end {Lemma}

\begin {Lemma} \label {6.1.3}  If $H$  is a Hopf algebra living in
symmetric tensor category  with a left duality  $H^{ *}$ , then

(1) $(H, \rightharpoonup )$  is a left  $H^{ \hat *}$-module algebra
under the module operation  $\rightharpoonup \  = (id \otimes d)
(C_{H^{\hat *}, H} \otimes id  ) (id \otimes \Delta ) $.

(2)  $(H^{\hat *}, \rightharpoonup )$  is a left $H$-module algebra
under the module operation  $\rightharpoonup \ = (id \otimes d)(id
\otimes C_{H, H^{\hat *}}) (C_{ H,H^{\hat *}} \otimes id  ) (id
\otimes \Delta )$.

(3)  $(H, \leftharpoonup )$  is a right $H^{ \hat *}$-module algebra
under the module operation  $\leftharpoonup \ = ( d\otimes id)
(C_{H,H^{\hat *} } \otimes id  )(id \otimes C_{H, H^{\hat *}})
(\Delta \otimes id )$.

(4)  $(H^{\hat *}, \leftharpoonup )$  is a right $H$-module algebra
under the module operation  $\leftharpoonup \  = ( d\otimes id) (id
\otimes C_{ H^{\hat *},H}) (\Delta \otimes id )$.

\end {Lemma}

\begin {Lemma}  \label {6.1.31} Let $H$ be a Hopf algebra living in symmetric
tensor category with a left duality $H^*.$  If we define
multiplication and unity on $H\otimes H^{\hat *}$  as follows:
\[
\begin{tangle}
\object{H \otimes H^{\hat *}}\step[6]\object{H \otimes H^{\hat *}}\\

\Cu\step[3]\\

\step[2]\object{H \otimes H^{\hat *}}\\
 \end{tangle} \ \ = \ \
\begin{tangle}
\object{H }\step[2]\object{H ^{\hat *}}  \step [2]\object{H } \step [2] \object { H^{\hat *}}\\

\id \step [2] \ev \step [2] \id \\
\object{H }\step [6]\object{H }

 \end{tangle} \ \
\]
and  $\eta _{H \otimes H^{\hat *}} = b_H,$
          then $H \otimes H^{\hat *}$   is  an algebra, written as
          $H \bar \otimes H^{\hat *}.$
          \end {Lemma}

          In fact, if ${\cal C}$  is a symmetric Yetter-Drinfeld category,  then $H\bar \otimes
          H^{\hat *}$  can be viewed as a subalgebra of $End _k(H)$.

\begin {Lemma} \label {6.1.4}
$\lambda $ is an algebraic morphism and $\rho $ is an anti-algebraic
morphism, where
 $\lambda = (m \otimes d \otimes id ) (id \otimes C_{H^{\hat
* }, H}\otimes id \otimes id )(id \otimes id \otimes \Delta \otimes id )(id \otimes id \otimes
b)$  and $\rho  =  ( d \otimes m\otimes id ) (id \otimes id \otimes
C_{ H, H}\otimes id  )(id \otimes C_{H, H}\otimes id \otimes id )(id
\otimes id \otimes \Delta \otimes id )(id \otimes id \otimes b)$.

\end {Lemma}
Proof.
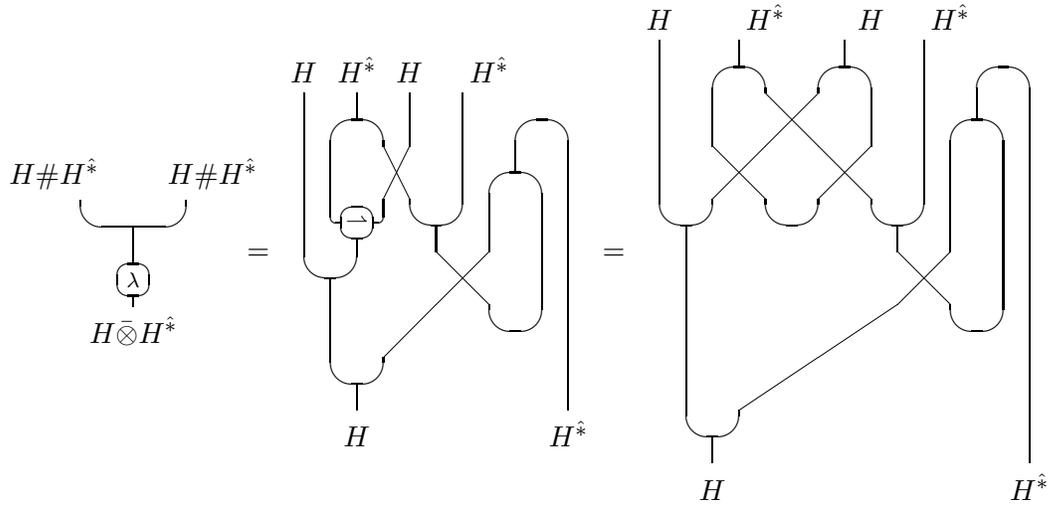
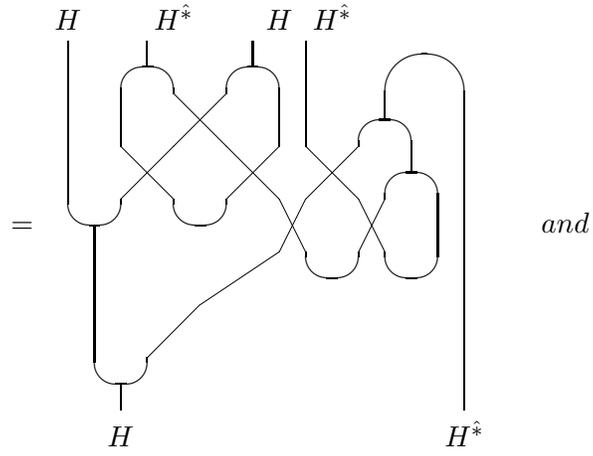
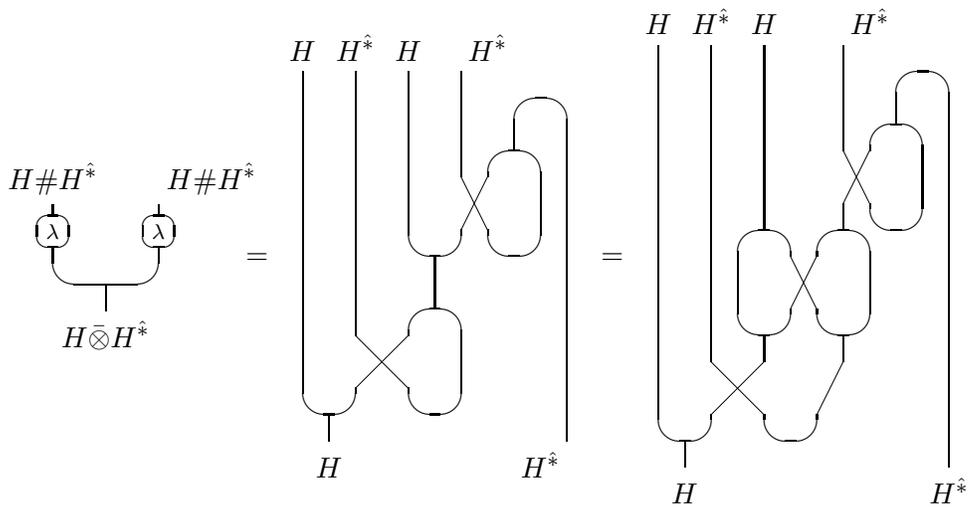
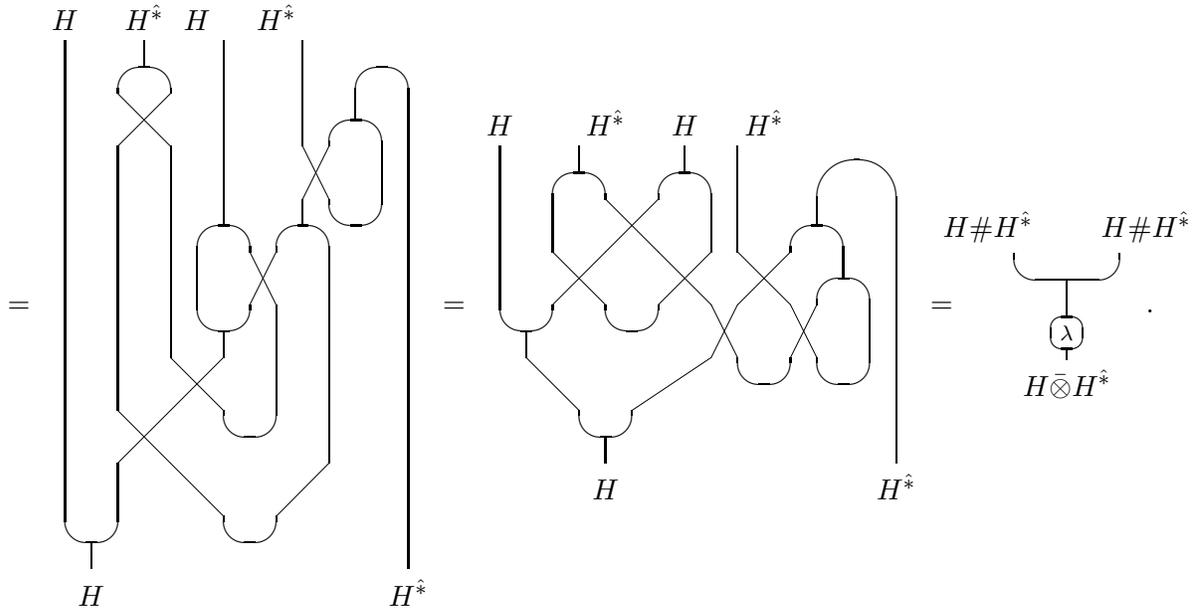
\[
\begin{tangle}
\object{H\# {H^{\hat{*}}}}\step[6]\object{H\#
{H^{\hat{*}}}}\\
\step\Cu\\
\step\Step\O \lambda\\
\step[3]\object{H\bar{\otimes}{H^{\hat{*}}}}
\end{tangle}
\step=\step
\begin{tangle}
\object{H}\step[2]\object{H^{\hat{*}}}\Step\object{H}\step[3]\object{H^{\hat{*}}}\\
\id\step\cd\step\id\Step\id\Step\coev\\
\id\step\id\Step\X\Step\id\step\cd\step\id\\
\id\step\tu \rightharpoonup\step\cu\step\id\Step\id\step\id\\
\cu\step[3]\XX\Step\id\step\id\\
\step\id\step[3]\ne2\Step\ev\step\id\\
\step\cu\step[7]\id\\
\Step\object{H}\step[8]\object{H^{\hat{*}}}\\
\end{tangle}
\step=\step
\begin{tangle}
\object{H}\step[4]\object{H^{\hat{*}}}\step[4]\object{H}\step[3]\object{H^{\hat{*}}}\\
\id\Step\cd\Step\cd\Step\id\Step\coev\\
\id\Step\id\Step\XX\Step\id\Step\id\step\cd\step\id\\
\id\Step\XX\Step\XX\Step\id\step\id\Step\id\step\id\\
\cu\Step\ev\Step\cu\step\id\Step\id\step\id\\
\step\id\step[8]\XX\Step\id\step\id\\
\step\id\step[7]\ne3\Step\ev\step\id\\
\step\id\step[4]\ne3\step[8]\id\\
\step\cu\step[11]\id\\
\Step\object{H}\step[12]\object{H^{\hat{*}}}\\
\end{tangle} \]
\[
\step=\step
\begin{tangle}
\object{H}\step[4]\object{H^{\hat{*}}}\step[4]\object{H}\step[2]\object{H^{\hat{*}}}\\
\id\Step\cd\Step\cd\step\id\step[3]\Coev\\
\id\Step\id\Step\XX\Step\id\step\id\Step\cd\step[2]\id\\
\id\Step\XX\Step\XX\step\XX\step\cd\step\id\\
\cu\Step\ev\Step\X\Step\X\Step\id\step\id\\
\step\id\step[6]\ne3\step\ev\step\ev\step\id\\
\step\id\step[3]\ne2\step[10]\id\\
\step\cu\step[12]\id\\
\Step\object{H}\step[13]\object{H^{\hat{*}}}\\
\end{tangle} \ \ \ \ \ \ \ \
and \]
\[
\begin{tangle}
\object{H\# {H^{\hat{*}}}}\step[6]\object{H\#
{H^{\hat{*}}}}\\
\O \lambda\step[4]\O \lambda\\
\Cu\\
\step[2]\object{H\bar{\otimes}{H^{\hat{*}}}}
\end{tangle}
\step=\step
\begin{tangle}
\object{H}\Step\object{H^{\hat{*}}}\Step\object{H}\step[3]\object{H^{\hat{*}}}\\
\id\Step\id\Step\id\Step\id\Step\coev\\
\id\Step\id\Step\id\Step\id\step\cd\step\id\\
\id\Step\id\Step\id\Step\X\Step\id\step\id\\
\id\Step\id\Step\cu\step\ev\step\id\\
\id\Step\id\Step\cd\step[4]\id\\
\id\Step\XX\Step\id\step[4]\id\\
\cu\Step\ev\step[4]\id\\
\step\object{H}\step[8]\object{H^{\hat{*}}}\\
\end{tangle}
\step=\step
\begin{tangle}
\object{H}\Step\object{H^{\hat{*}}}\Step\object{H}\step[4]\object{H^{\hat{*}}}\\
\id\Step\id\Step\id\step[3]\id\Step\coev\\
\id\Step\id\Step\id\step[3]\id\step\cd\step\id\\
\id\Step\id\Step\id\step[3]\X\Step\id\step\id\\
\id\Step\id\step\cd\step\cd\ev\step\id\\
\id\Step\id\step\id\Step\X\Step\id\step[3]\id\\
\id\Step\id\step\cu\step\cu\step[3]\id\\
\id\Step\XX\Step\dd\step[4]\id\\
\cu\Step\ev\step[5]\id\\
\step\object{H}\step[10]\object{H^{\hat{*}}}\\
\end{tangle} \]
\[
\step=\step
\begin{tangle}
\object{H}\step[3]\object{H^{\hat{*}}}\Step\object{H}\step[3]\object{H^{\hat{*}}}\\
\id\Step\cd\Step\id\step[3]\id\Step\coev\\
\id\Step\XX\Step\id\step[3]\id\step\cd\step\id\\
\id\Step\id\Step\id\Step\id\step[3]\X\Step\id\step\id\\
\id\Step\id\Step\id\step\cd\step\cd\ev\step\id\\
\id\Step\id\Step\id\step\id\Step\X\Step\id\step[3]\id\\
\id\Step\id\Step\id\step\cu\step\id\Step\id\step[3]\id\\
\id\Step\id\Step\XX\Step\id\Step\id\step[3]\id\\
\id\Step\XX\Step\ev\Step\id\step[3]\id\\
\id\Step\id\Step\nw2\step[4]\ne2\step[3]\id\\
\cu\step[4]\ev\step[5]\id\\
\step\object{H}\step[12]\object{H^{\hat{*}}}\\
\end{tangle}
\step=\step
\begin{tangle}
\object{H}\step[4]\object{H^{\hat{*}}}\step[3]\object{H}\step[3]\object{H^{\hat{*}}}\\
\id\Step\cd\Step\cd\step\id\step[3]\Coev\\
\id\Step\id\Step\XX\Step\id\step\id\Step\cd\step[2]\id\\
\id\Step\XX\Step\XX\step\XX\step\cd\step\id\\
\cu\Step\ev\Step\X\Step\X\Step\id\step\id\\
\step\nw2\step[5]\ne3\step\ev\step\ev\step\id\\
\step[3]\cu\step[10]\id\\
\step[4]\object{H}\step[11]\object{H^{\hat{*}}}\\
\end{tangle}
\step=\step
\begin{tangle}
\object{H\# {H^{\hat{*}}}}\step[6]\object{H\#
{H^{\hat{*}}}}\\
\step\Cu\\
\step\Step\O \lambda\\
\step[3]\object{H\bar{\otimes}{H^{\hat{*}}}}\\
\end{tangle} .
\]
Thus $\lambda $ is an algebraic morphism.

Similarly, we can show that $\rho $ is an anti-algebraic morphism.\
\ \ \ \begin{picture}(5,5)
\put(0,0){\line(0,1){5}}\put(5,5){\line(0,-1){5}}
\put(0,0){\line(1,0){5}}\put(5,5){\line(-1,0){5}}
\end{picture}\\

\begin {Lemma} \label {6.1.5}
$\lambda $  and $\rho $ satisfy the following:

(1)
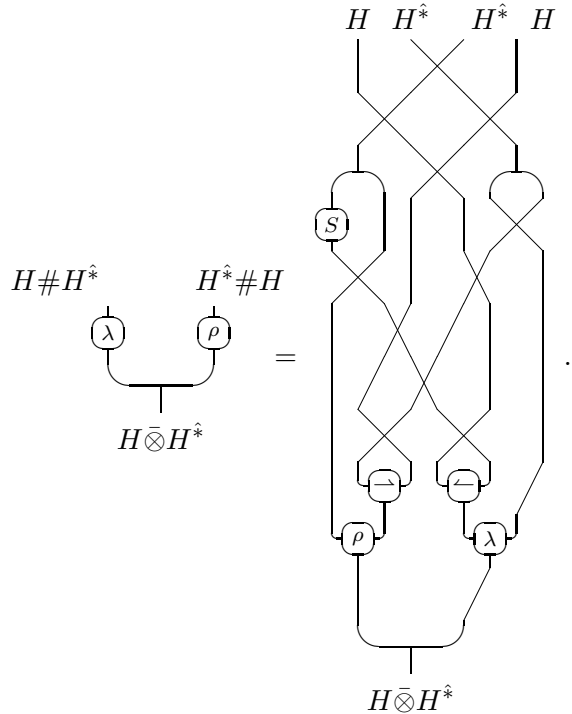
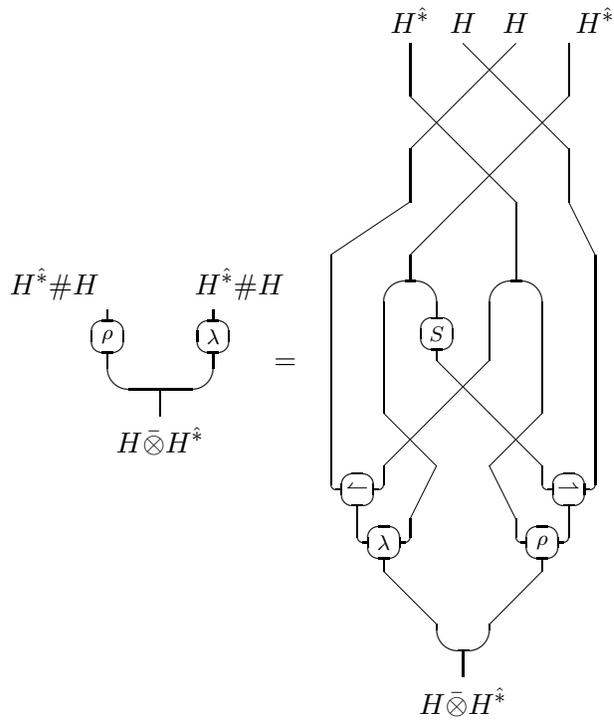
\[
\begin{tangle}
\object{H\# {{H^{\hat{*}}}}}\step[7]\object{{H^{\hat{*}}}\#
{H}}\\
\Step\O \lambda\step[4]\O \rho\\
\Step\Cu\\
\step[4]\object{H\bar{\otimes}{H^{\hat{*}}}}\\
\end{tangle}
\step=\step
\begin{tangle}
\step\object{H}\step[2]\object{H^{\hat{*}}}\step[3]\object{H^{\hat{*}}}\Step\object{H}\\
\step\id\Step\XX\Step\id\\
\step\XX\Step\XX\\
\cd\step\XX\step\cd\\
\O S\Step\id\step\id\Step\id\step\XX\\
\XX\step\id\Step\X\Step\id\\
\id\Step\X\step\dd\step\id\Step\id\\
\id\step\dd\step\X\Step\id\Step\id\\
\id\step\XX\step\XX\Step\id\\
\id\step\tu \rightharpoonup\step\tu \leftharpoonup\step\dd\\
\tu \rho\step[3]\tu \lambda\\
\step\id\step[4]\dd\\
\step\Cu\\
\step[3]\object{H\bar{\otimes}{H^{\hat{*}}}}\\
\end{tangle}
\ \ . \ \
\]

(2)

 \[
\begin{tangle}
\object{H^{\hat*}\#H }\step[7]\object{H^{\hat* }\#H}\\
\Step\O \rho\step[4]\O \lambda\\
\Step\Cu\\
\step[4]\object{H\bar{\otimes}{H^{\hat{*}}}}\\
\end{tangle}
\step=\step
\begin{tangle}
\step[3]\object{H^{\hat{*}}}\step[2]\object{H}\Step\object{H}\step[3]\object{{H^{\hat{*}}}}\\
\step[3]\id\Step\XX\Step\id\\
\step[3]\XX\Step\XX\\
\step[3]\id\Step\XX\Step\id\\

\step [2]\ne3 \step [1] \ne2 \step [2] \id \step [2] \nw1\\

\id\Step\cd \step[2]\cd\step[2]\id\\

\id\step[2]\id  \step[2]\O S\step[2]\id\step[2]\id\step[2]\id\\

\id\step[2]\id  \step[2]\XX\step[2]\id\step[2]\id\step[2]\\

\id\step[2]\XX  \step[2]\XX\step[2]\id\\

\tu \leftharpoonup\step\ne1\step[2] \nw1 \step [1]\tu \rightharpoonup\\

\step\tu \lambda\step[4]\tu \rho\\

\step[2]\nw2\step[4]\ne2\\

\step[4]\cu\\
\step[5]\object{H\bar{\otimes}{H^{\hat{*}}}}\\
\end{tangle}\ \ \ .
\]

\end {Lemma}

{\bf Proof.} We show (1) by following five steps. It is easy to
check the following (i) and (ii).

(i)
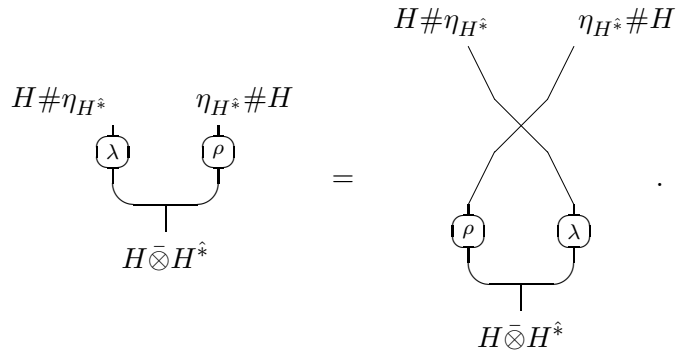
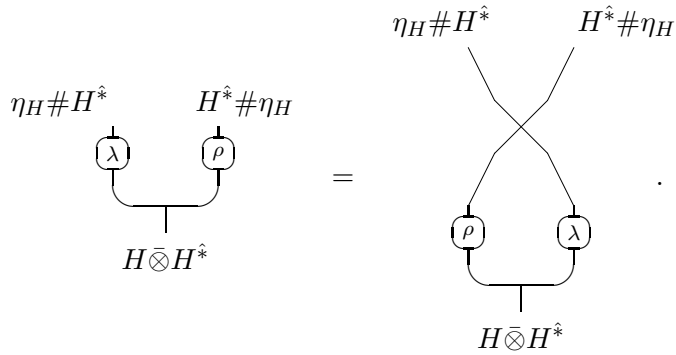
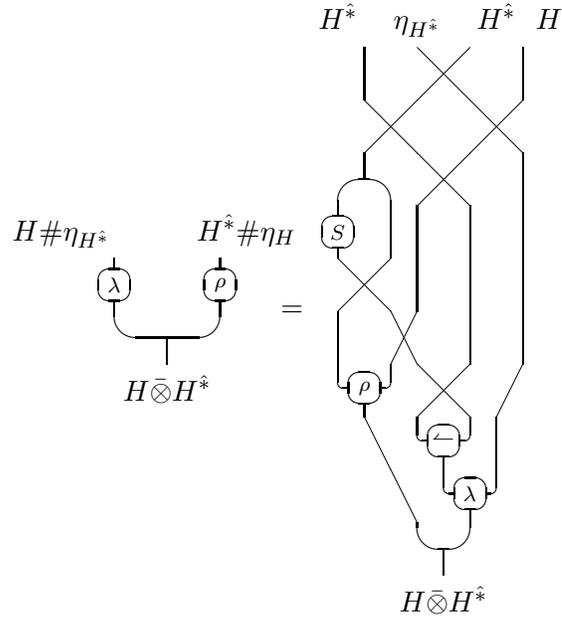
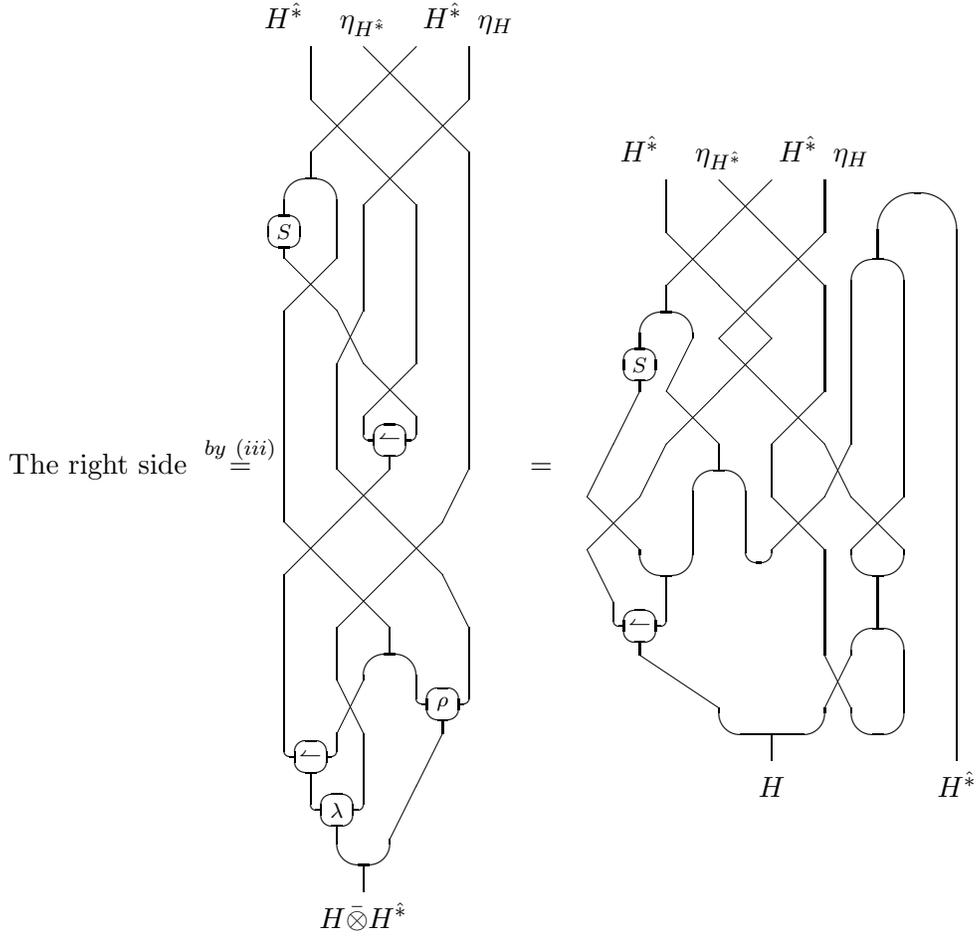
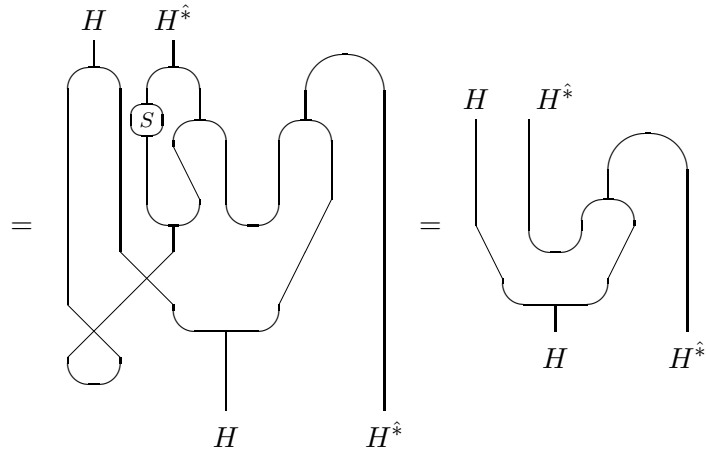
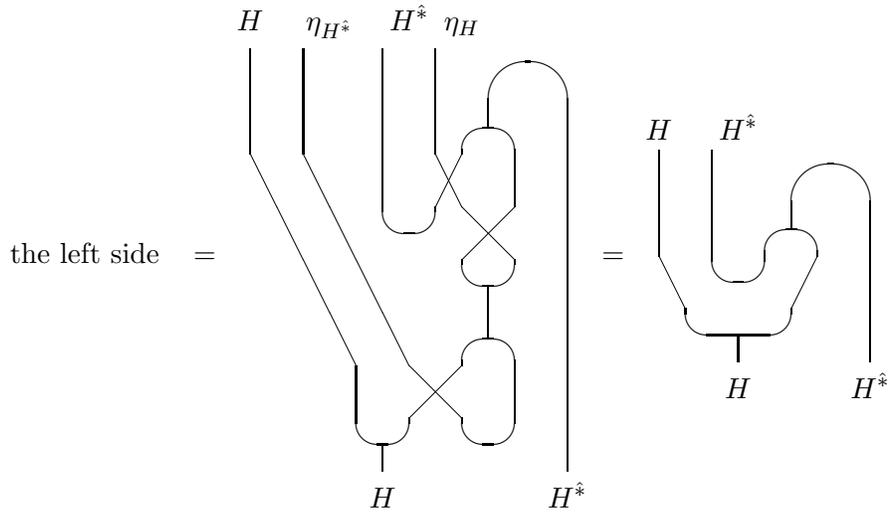
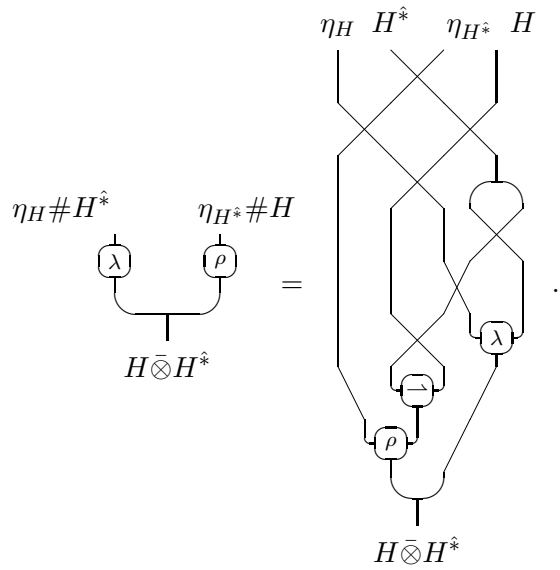
\[
\begin{tangle}
\object{H\#
{\eta_{H^{\hat{*}}}}}\step[7]\object{{\eta_{H^{\hat{*}}}}\#
H}\\
\Step\O \lambda\step[4]\O \rho\\
\Step\Cu\\
\step[4]\object{H\bar{\otimes}{H^{\hat{*}}}}\\
\end{tangle}
\step[3]=\step[3]
\begin{tangle}
\object{H\#
{\eta_{H^{\hat{*}}}}}\step[7]\object{{\eta_{H^{\hat{*}}}}\#
H}\\
\step\d\Step\dd\\
\Step\XX\\
\step\dd\Step\d\\
\step\O \rho\step[4]\O \lambda\\
\step\Cu\\
\step[3]\object{H\bar{\otimes}{H^{\hat{*}}}}\\
\end{tangle}\ \ \ .
\]

(ii)
\[
\begin{tangle}
\object{\eta_H\# {H^{\hat{*}}}}\step[7]\object{H^{\hat{*}}\#
\eta_H}\\
\Step\O \lambda\step[4]\O \rho\\
\Step\Cu\\
\step[4]\object{H\bar{\otimes}{H^{\hat{*}}}}\\
\end{tangle}
\step[3]=\step[3]
\begin{tangle}
\object{\eta_H\# {H^{\hat{*}}}}\step[7]\object{H^{\hat{*}}\#
\eta_H}\\
\step\d\Step\dd\\
\Step\XX\\
\step\dd\Step\d\\
\step\O \rho\step[4]\O \lambda\\
\step\Cu\\
\step[3]\object{H\bar{\otimes}{H^{\hat{*}}}}\\
\end{tangle}\ \ \ .
\]

(iii)
\[
\begin{tangle}
\object{{H^{\hat{*}}}\#\eta_H}\step[7]\object{{H^{\hat{*}}}\#\eta_H
}\\
\Step\O \rho\step[4]\O \lambda\\
\Step\Cu\\
\step[4]\object{H\bar{\otimes}{H^{\hat{*}}}}\\
\end{tangle}
\step=\step
\begin{tangle}
\object{H^{\hat{*}}}\step[3]\object{\eta_H}\Step\object{H}\step[3]\object{\eta_{H^{\hat{*}}}}\\
\id\Step\XX\Step\id\\
\XX\Step\XX\\
\id\Step\XX\Step\id\\
\id\Step\id\step\cd\step\id\\
\id\Step\X\step\dd\step\id\\
\tu \leftharpoonup\step\id\step\tu \rho\\
\step\tu \lambda\step\dd\\
\Step\cu\\
\step[3]\object{H\bar{\otimes}{H^{\hat{*}}}}\\
\end{tangle}\ \ \ .
\]

\[
\hbox {In fact,the right side} \step=\step
\begin{tangle}
\object{H^{\hat{*}}}\step[3]\object{\eta_H}\step[2]\object{H}\step[3]\object{\eta_{H^{\hat{*}}}}\\
\step\id\Step\XX\Step\id\step[3]\coev\\
\step\XX\Step\XX\Step\cd\step\id\\
\cd\step\XX\Step\XX\Step\id\step\id\\
\id\Step\id\step\id\step\cd\step\d\step\XX\step\id\\
\id\Step\id\step\X\Step\ev\step\cu\step\id\\
\id\Step\X\step\nw3\step[4]\cd\step\id\\
\XX\step\nw2\step[3]\XX\Step\id\step\id\\
\ev\step[3]\cu\Step\ev\step\id\\
\step[6]\object{H}\step[6]\object{H^{\hat{*}}}\\
\end{tangle} \]
\[
\step=\step
\begin{tangle}
\step\object{H^{\hat{*}}}\step[4]\object{H}\\
\cd\step\cd\Step\Coev\\
\id\Step\X\Step\id\step\cd\Step\id\\
\ev\step\id\Step\X\Step\id\Step\id\\
\step[3]\ev\step\cu\Step\id\\
\step[7]\object{H}\step[3]\object{H^{\hat{*}}}\\
\end{tangle}
\step=\step \hbox{the left side.}
\]
Thus (iii) holds.

(iv)
\[
\begin{tangle}
\object{H\# {\eta_{H^{\hat{*}}}}}\step[7]\object{{H^{\hat{*}}}\#
{\eta_H}}\\
\Step\O \lambda\step[4]\O \rho\\
\Step\Cu\\
\step[4]\object{H\bar{\otimes}{H^{\hat{*}}}}\\
\end{tangle}
\step=\step
\begin{tangle}
\object{H^{\hat{*}}}\step[3]\object{\eta_{H^{\hat{*}}}}\step[3]\object{H^{\hat{*}}}\step[2]\object{H}\\
\step\id\Step\XX\Step\id\\
\step\XX\Step\XX\\
\cd\step\XX\Step\id\\
\O S\Step\id\step\id\Step\id\Step\id\\
\XX\step\id\Step\id\Step\id\\
\id\Step\X\Step\id\Step\id\\
\tu \rho\step\XX\step\dd\\
\step\d\step\tu \leftharpoonup\step\id\\
\Step\d\step\tu \lambda\\
\step[3]\cu\\
\step[4]\object{H\bar{\otimes}{H^{\hat{*}}}}\\
\end{tangle}\ \ \ .
\]

\[
\hbox{The right side }\stackrel{by\ (iii)}{=}
\begin{tangle}
\object{H^{\hat{*}}}\step[3]\object{\eta_{H^{\hat{*}}}}\step[3]\object{H^{\hat{*}}}\step[2]\object{\eta_{H}}\\
\step\id\Step\XX\Step\id\\
\step\XX\Step\XX\\
\cd\step\XX\Step\id\\
\O S\Step\id\step\id\Step\id\Step\id\\
\XX\step\id\Step\id\Step\id\\
\id\Step\X\Step\id\Step\id\\
\id\Step\id\step\XX\Step\id\\
\id\Step\id\step\tu \leftharpoonup\Step\id\\
\id\Step\XX\Step\dd\\
\XX\Step\XX\\
\id\Step\XX\Step\d\\
\id\Step\id\step\cd\Step\id\\
\id\Step\X\Step\tu \rho\\
\tu \leftharpoonup\step\id\Step\dd\\
\step\tu \lambda\step\dd\\
\Step\cu\\
\step[3]\object{H\bar{\otimes}{H^{\hat{*}}}}\\
\end{tangle}
\step=\step
\begin{tangle}
\Step\object{H^{\hat{*}}}\step[3]\object{\eta_{H^{\hat{*}}}}\step[3]\object{H^{\hat{*}}}\step[2]\object{\eta_{H}}\\
\step[3]\id\Step\XX\Step\id\Step\Coev\\
\step[3]\XX\Step\XX\step\cd\Step\id\\
\Step\cd\step\XX\Step\id\step\id\Step\id\Step\id\\
\Step\O S\step\dd\step\XX\Step\id\step\id\Step\id\Step\id\\
\step\dd\step\XX\Step\XX\step\id\Step\id\Step\id\\
\dd\step\dd\step\cd\step\id\Step\X\Step\id\Step\id\\
\XX\Step\id\Step\id\step\XX\step\XX\Step\id\\
\d\step\cu\Step\hev\Step\id\step\cu\Step\id\\
\step\tu \leftharpoonup\step[6]\id\step\cd\Step\id\\
\Step\nw3\step[6]\X\Step\id\Step\id\\
\step[5]\Cu\step\ev\Step\id\\
\step[7]\object{H}\step[7]\object{H^{\hat{*}}}\\
\end{tangle}\]
\[
\step=\step
\begin{tangle}
\step\object{H}\step[3]\object{H^{\hat{*}}}\\
\cd\step\cd\step[4]\Coev\\
\id\Step\id\step\O S\step\cd\Step\cd\Step\id\\
\id\Step\id\step\id\step\d\step\id\Step\id\Step\id\Step\id\\
\id\Step\id\step\cu\step\ev\step\dd\Step\id\\
\id\Step\XX\step[4]\dd\step[3]\id\\
\XX\Step\Cu\step[4]\id\\
\ev\step[4]\id\step[6]\id\\
\step[6]\object{H}\step[6]\object{H^{\hat{*}}}\\
\end{tangle}
\step=\step
\begin{tangle}
\object{H}\step[3]\object{H^{\hat{*}}}\\
\id\Step\id\step[3]\Coev\\
\id\Step\id\Step\cd\Step\id\\
\d\step\ev\step\dd\Step\id\\
\step\Cu\step[3]\id\\
\step[3]\object{H}\step[5]\object{H^{\hat{*}}}\\
\end{tangle}
\]

and
\[
\hbox{the left side}\step=\step
\begin{tangle}
\object{H}\step[3]\object{\eta_{H^{\hat{*}}}}\step[3]\object{H^{\hat{*}}}\step[2]\object{\eta_H}\\
\id\Step\id\step[3]\id\Step\id\Step\Coev\\
\id\Step\id\step[3]\id\Step\id\step\cd\Step\id\\
\d\step\d\step[2]\id\Step\X\Step\id\Step\id\\
\step\d\step\d\step\ev\step\XX\Step\id\\
\Step\d\step\d\step[3]\cu\Step\id\\
\step[3]\d\step\d\Step\cd\Step\id\\
\step[4]\id\Step\XX\Step\id\Step\id\\
\step[4]\cu\Step\ev\Step\id\\
\step[5]\object{H}\step[7]\object{H^{\hat{*}}}\\
\end{tangle}
\step=\step
\begin{tangle}
\object{H}\step[3]\object{H^{\hat{*}}}\\
\id\Step\id\step[3]\Coev\\
\id\Step\id\Step\cd\Step\id\\
\d\step\ev\step\dd\Step\id\\
\step\Cu\step[3]\id\\
\step[3]\object{H}\step[5]\object{H^{\hat{*}}}\\
\end{tangle}\ \ .
\]

Thus (iv) holds.

(v)
\[
\begin{tangle}
\object{\eta_H\# {H^{\hat{*}}}}\step[7]\object{\eta_{H^{\hat{*}}}\#
H}\\
\Step\O \lambda\step[4]\O \rho\\
\Step\Cu\\
\step[4]\object{H\bar{\otimes}{H^{\hat{*}}}}\\
\end{tangle}
\step=\step
\begin{tangle}
\object{\eta_H}\Step\object{H^{\hat{*}}}\step[3]\object{\eta_{H^{\hat{*}}}}\step[2]\object{H}\\
\id\Step\XX\Step\id\\
\XX\Step\XX\\
\id\Step\XX\step\cd\\
\id\Step\id\Step\id\step\XX\\
\id\Step\id\Step\X\Step\id\\
\id\Step\XX\step\tu \lambda\\
\d\step\tu \rightharpoonup\step\dd\\
\step\tu \rho\step\dd\\
\Step\cu\\
\step[3]\object{H\bar{\otimes}{H^{\hat{*}}}}\\
\end{tangle}\ \ \ .
\]

\[
\hbox{In fact, the right side}\step=\step
\begin{tangle}
\object{\eta_H}\Step\object{H^{\hat{*}}}\step[3]\object{\eta_{H^{\hat{*}}}}\step[2]\object{H}\\
\id\Step\XX\Step\id\step[3]\Coev\\
\XX\Step\XX\Step\cd\Step\id\\
\id\Step\XX\step\cd\step\id\Step\id\Step\id\\
\id\step\cd\step\id\step\XX\step\id\Step\id\Step\id\\
\id\step\id\Step\id\step\X\Step\X\Step\id\Step\id\\
\id\step\d\step\X\step\cu\step\ev\Step\id\\
\d\step\d\hev\step\cd\step[5]\id\\
\step\id\Step\XX\Step\id\step[5]\id\\
\step\ev\Step\XX\step[5]\id\\
\step[5]\cu\step[5]\id\\
\step[6]\object{H}\step[6]\object{H^{\hat{*}}}\\
\end{tangle}
\step=\step
\begin{tangle}
\step\object{H^{\hat{*}}}\step[3]\object{H}\\
\cd\step\cd\Step\Coev\\
\XX\step\XX\step\cd\Step\id\\
\id\Step\X\Step\id\step\id\Step\id\Step\id\\
\ev\step\XX\step\id\Step\id\Step\id\\
\step[3]\id\Step\X\Step\id\Step\id\\
\step[3]\XX\step\ev\step[2]\id\\
\step[3]\cu\step[5]\id\\
\step[4]\object{H}\step[6]\object{H^{\hat{*}}}\\
\end{tangle}
\]

and
\[
\hbox{the left side}\step=\step
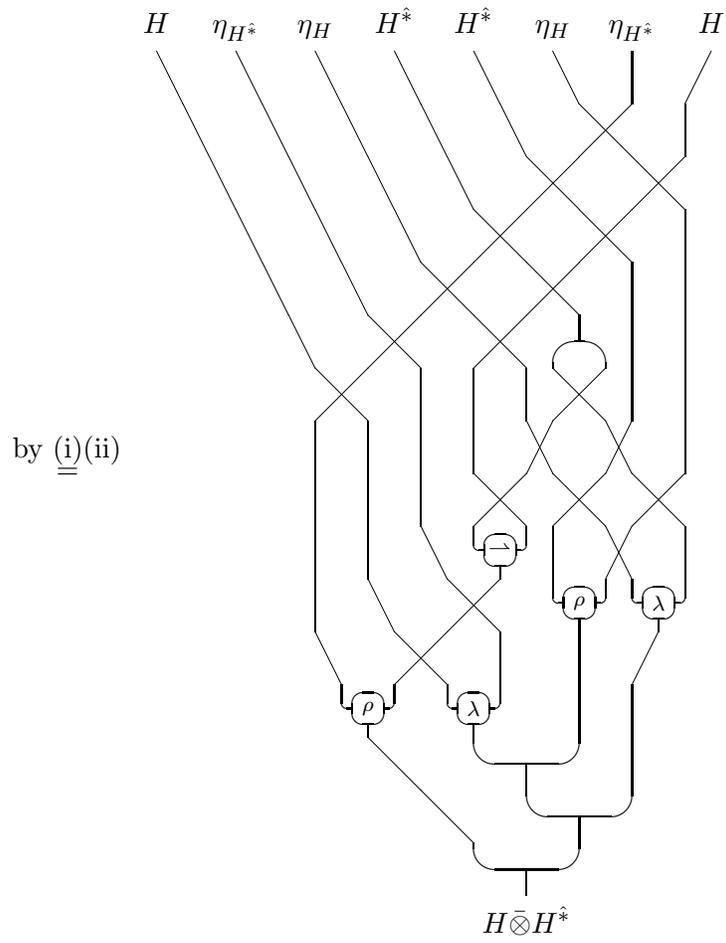
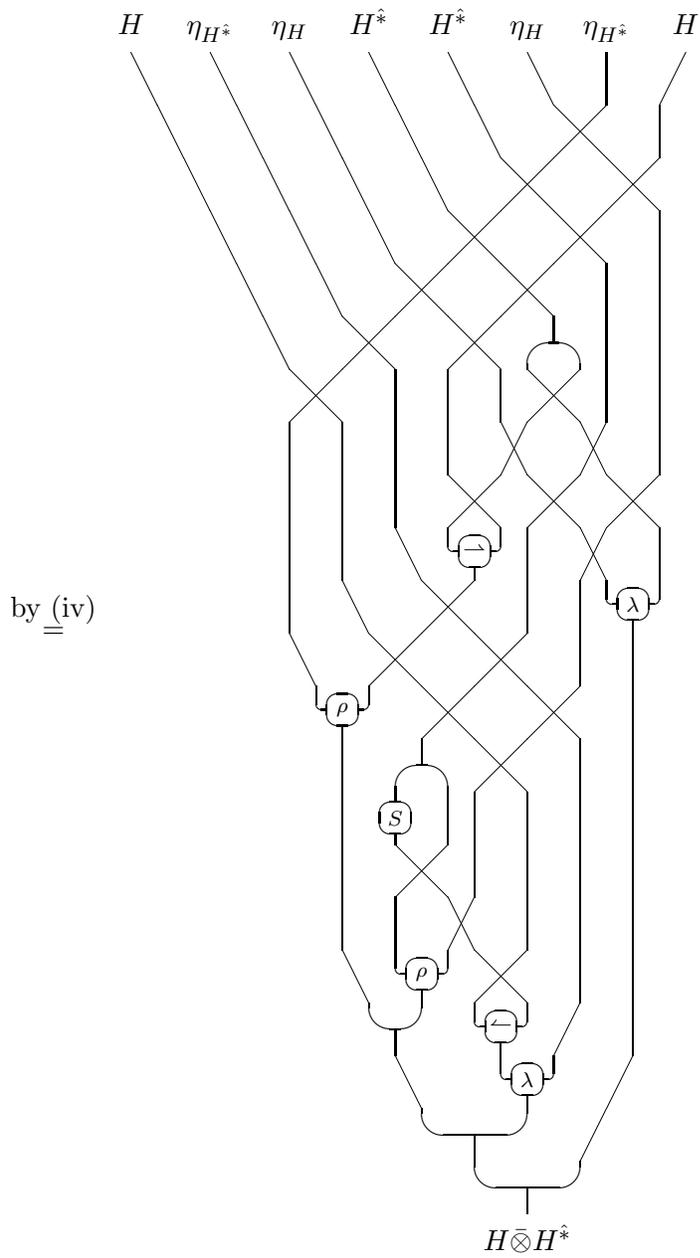
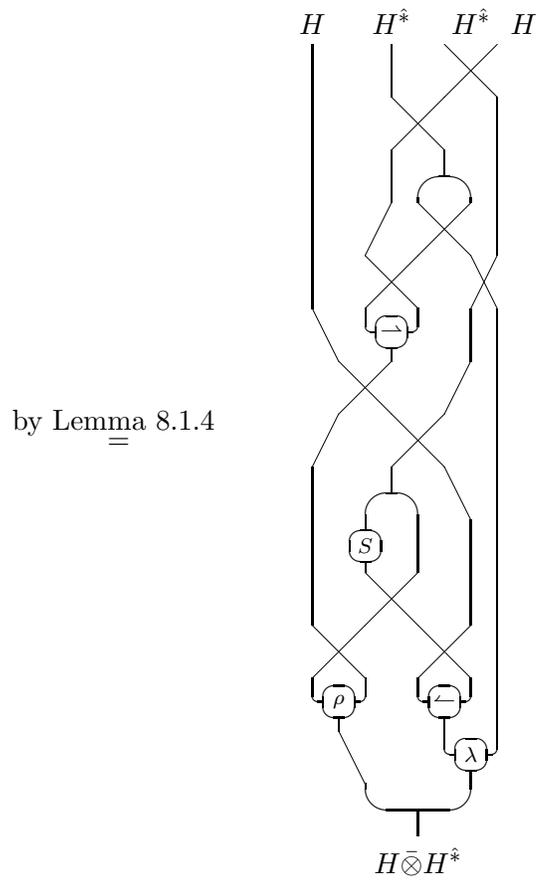
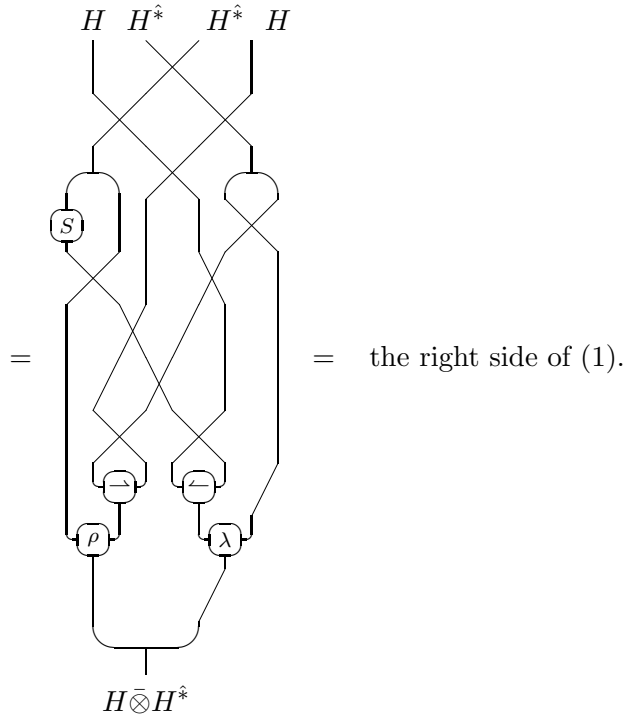
\begin{tangle}
\object{H^{\hat{*}}}\step[3]\object{H}\\
\id\Step\id\Step\coev\\
\id\Step\XX\Step\id\\
\id\Step\cu\Step\id\\
\id\Step\cd\Step\id\\
\XX\Step\id\Step\id\\
\id\Step\ev\Step\id\\
\object{H}\step[6]\object{H^{\hat{*}}}\\
\end{tangle}
\step=\step
\begin{tangle}
\object{H^{\hat{*}}}\step[3]\object{H}\\
\id\Step\id\Step\Coev\\
\id\Step\XX\step[3]\id\\
\id\Step\id\Step\d\step[2]\id\\
\id\step\cd\step\cd\step\id\\
\id\step\id\Step\X\Step\id\step\id\\
\id\step\cu\step\cu\step\id\\
\XX\Step\dd\step[2]\id\\
\id\Step\ev\step[3]\id\\
\object{H}\step[7]\object{H^{\hat{*}}}\\
\end{tangle}\]
\[
\step=\step
\begin{tangle}
\step\object{H^{\hat{*}}}\step[3]\object{H}\\
\cd\step[2]\id\Step\Coev\\
\id\Step\id\step[2]\XX\step[3]\id\\
\id\step[2]\id\Step\id\Step\d\Step\id\\
\id\step[2]\id\step\cd\step\cd\step\id\\
\id\Step\id\step\id\Step\X\Step\id\step\id\\
\id\Step\id\step\cu\step\id\Step\id\step\id\\
\id\Step\XX\Step\id\Step\id\step\id\\
\XX\Step\XX\Step\id\step\id\\
\id\Step\ev\Step\ev\step\id\\
\object{H}\step[9]\object{H^{\hat{*}}}\\
\end{tangle}
\step=\step
\begin{tangle}
\step\object{H^{\hat{*}}}\step[3]\object{H}\\
\cd\step\cd\Step\Coev\\
\id\Step\id\step\id\Step\id\step\cd\Step\id\\
\id\Step\id\step\id\Step\X\Step\id\Step\id\\
\id\Step\id\step\XX\step\id\Step\id\Step\id\\
\id\Step\id\step\cu\step\id\Step\id\Step\id\\
\id\Step\XX\Step\id\step\dd\Step\id\\
\XX\Step\ev\dd\step[3]\id\\
\id\Step\d\step[3]\id\step[4]\id\\
\id\step[3]\Ev\step[6]\id\\
\object{H}\step[10]\object{H^{\hat{*}}}\\
\end{tangle}
\step=\step\hbox{the right side.}
\]
Thus (v) holds.

Now we show that the relation (1) holds.
\[
\hbox{the left side of (1)}\stackrel{\hbox{by Lemma \ref {6.1.4}
}}{\step=\step[5]}
\begin{tangle}
\object{H\# {\eta_{H^{\hat{*}}}}}\step[8]\object{{\eta_H}\#
H^{\hat{*}}}\step[9]\object{H^{\hat{*}}\#\eta_H}\step[6]\object{\eta_{H^{\hat{*}}}\# H}\\
\O \lambda\step[6]\Step\O \lambda\step[11]\XX \\
\nw3\step[7]\nw3\step[10]\O \rho\step[2]\O \rho\\
\step[3]\nw3\step[7]\nw3\step[7]\cu\\
\step[6]\nw3\step[7]\d\step[4]\dd\\
\step[9]\nw3\step[5]\Cu\\
\step[12]\d\step[4]\id\\
\step[13]\Cu\\
\step[15]\object{H\bar{\otimes}{H^{\hat{*}}}}\\
\end{tangle} \]
\[
\step\stackrel{\hbox{by (v)}}{=}\step
\begin{tangle}
\object{H\#
{\eta_{H^{\hat{*}}}}}\step[5]\object{\eta_H}\step[3]\object{H^{\hat{*}}}\step[5]
\object{H^{\hat{*}}\#\eta_H}\step[5]\object{\eta_{H^{\hat{*}}}}\step[3]\object{H}\\
\O \lambda\step[5]\nw3\step[2]\nw3\step[4]\d\Step\dd\Step\dd\\
\d\step[7]\nw2\Step\d\Step\XX\Step\dd\\
\step\d\step[8]\id\Step\XX\Step\XX\\
\Step\d\step[7]\XX\Step\XX\Step\O \rho\\
\step[3]\d\step[6]\id\Step\XX\step\cd\step\id\\
\step[4]\d\step[5]\id\Step\id\Step\id\step\XX\step\id\\
\step[5]\d\step[4]\id\Step\id\Step\X\Step\id\step\id\\
\step[6]\d\step[3]\id\Step\XX\step\tu \lambda\step\id\\
\step[7]\d\Step\d\step\tu \rightharpoonup\Step\cu\\
\step[8]\d\Step\tu \rho\step[3]\dd\\
\step[9]\d\Step\Cu\\
\step[10]\Cu\\
\step[12]\object{H\bar{\otimes}{H^{\hat{*}}}}\\
\end{tangle}
\]

\[
\step\stackrel{\hbox{by (i)(ii)}}{=}\step
\begin{tangle}
\object{H}\step[3]\object{\eta_{H^{\hat{*}}}}\step[3]\object{\eta_H}\step[3]\object{H^{\hat{*}}}
\step[3]\object{H^{\hat{*}}}\step[3]\object{\eta_H}\step[3]\object{\eta_{H^{\hat{*}}}}\step[3]\object{H}\\
\d\Step\d\Step\d\Step\d\Step\d\Step\d\Step\id\Step\dd\\
\step\d\Step\d\Step\d\Step\d\Step\d\Step\XX\Step\id\\
\Step\d\Step\d\Step\d\Step\d\Step\XX\Step\XX\\
\step[3]\d\Step\d\Step\d\Step\XX\Step\XX\Step\id\\
\step[4]\d\Step\d\Step\XX\Step\XX\Step\id\Step\id\\
\step[5]\d\Step\XX\Step\XX\step\cd\step\id\Step\id\\
\step[6]\XX\Step\id\Step\id\Step\id\step\XX\step\id\Step\id\\
\step[6]\id\Step\id\Step\id\Step\id\Step\X\Step\X\Step\id\\
\step[6]\id\Step\id\Step\id\Step\XX\step\XX\step\XX\\
\step[6]\id\Step\id\Step\d\step\tu
\rightharpoonup\step\id\Step\X\Step\id\\
\step[6]\id\Step\d\Step\XX\Step\tu \rho\step\tu \lambda\\
\step[6]\d\Step\XX\Step\id\step[3]\id\Step\dd\\
\step[7]\tu \rho\Step\tu \lambda\step[3]\id\Step\id\\
\step[8]\nw2\step[3]\Cu\Step\id\\
\step[10]\nw2\step[3]\Cu\\
\step[12]\Cu\\
\step[14]\object{H\bar{\otimes}{H^{\hat{*}}}}\\
\end{tangle}
\]
\[
\step\stackrel{\hbox{by (iv)}}{=}\step
\begin{tangle}
\object{H}\step[3]\object{\eta_{H^{\hat{*}}}}\step[3]\object{\eta_H}\step[3]\object{H^{\hat{*}}}
\step[3]\object{H^{\hat{*}}}\step[3]\object{\eta_H}\step[3]\object{\eta_{H^{\hat{*}}}}\step[3]\object{H}\\
\d\Step\d\Step\d\Step\d\Step\d\Step\d\Step\id\Step\dd\\
\step\d\Step\d\Step\d\Step\d\Step\d\Step\XX\Step\id\\
\Step\d\Step\d\Step\d\Step\d\Step\XX\Step\XX\\
\step[3]\d\Step\d\Step\d\Step\XX\Step\XX\Step\id\\
\step[4]\d\Step\d\Step\XX\Step\XX\Step\id\Step\id\\
\step[5]\d\Step\XX\Step\XX\step\cd\step\id\Step\id\\
\step[6]\XX\Step\id\Step\id\Step\id\step\XX\step\id\Step\id\\
\step[6]\id\Step\id\Step\id\Step\id\Step\X\Step\X\Step\id\\
\step[6]\id\Step\id\Step\id\Step\XX\step\XX\step\XX\\
\step[6]\id\Step\id\Step\d\step\tu
\rightharpoonup\step\id\Step\X\Step\id\\
\step[6]\id\Step\d\Step\XX\Step\id\Step\id\step\tu \lambda\\
\step[6]\d\Step\XX\Step\XX\Step\id\Step\id\\
\step[7]\tu \rho\Step\XX\Step\XX\Step\id\\
\step[8]\id\Step\cd\step\XX\Step\id\Step\id\\
\step[8]\id\Step\O S\Step\id\step\id\Step\id\Step\id\Step\id\\
\step[8]\id\Step\XX\step\id\Step\id\Step\id\Step\id\\
\step[8]\id\Step\id\Step\X\Step\id\Step\id\Step\id\\
\step[8]\d\step\tu \rho\step\XX\Step\id\Step\id\\
\step[9]\cu\Step\tu \leftharpoonup\step\dd\Step\id\\
\step[10]\d\step[3]\tu \lambda\step[2]\dd\\
\step[11]\Cu\step[2]\dd\\
\step[13]\Cu\\
\step[15]\object{H\bar{\otimes}{H^{\hat{*}}}}\\
\end{tangle}
\step\stackrel{\hbox{by Lemma \ref {6.1.4}  }}{=}\step[3]
\begin{tangle}
\object{H}\step[3]\object{H^{\hat{*}}}\step[3]\object{H^{\hat{*}}}\Step\object{H}\\
\id\step[3]\id\Step\XX\\
\id\step[3]\XX\Step\id\\
\id\step[3]\id\step\cd\step\id\\
\id\Step\dd\step\XX\step\id\\
\id\Step\XX\Step\X\\
\d\step\tu \rightharpoonup\Step\id\step\id\\
\step\XX\Step\dd\step\id\\
\dd\Step\XX\Step\id\\
\id\Step\cd\step\d\step\id\\
\id\Step\O S\Step\id\Step\id\step\id\\
\id\Step\XX\Step\id\step\id\\
\XX\Step\XX\step\id\\
\tu \rho\Step\tu \leftharpoonup\step\id\\
\step\d\step[3]\tu \lambda\\
\Step\Cu\\
\step[4]\object{H\bar{\otimes}{H^{\hat{*}}}}\\
\end{tangle}
\]
\[
\step=\step
\begin{tangle}
\step\object{H}\step[2]\object{H^{\hat{*}}}\step[3]\object{H^{\hat{*}}}\Step\object{H}\\
\step\id\Step\XX\Step\id\\
\step\XX\Step\XX\\
\cd\step\XX\step\cd\\
\O S\Step\id\step\id\Step\id\step\XX\\
\XX\step\id\Step\X\Step\id\\
\id\Step\X\step\dd\step\id\Step\id\\
\id\step\dd\step\X\Step\id\Step\id\\
\id\step\XX\step\XX\Step\id\\
\id\step\tu \rightharpoonup\step\tu \leftharpoonup\step\dd\\
\tu \rho\step[3]\tu \lambda\\
\step\id\step[4]\dd\\
\step\Cu\\
\step[3]\object{H\bar{\otimes}{H^{\hat{*}}}}\\
\end{tangle}
\step=\step
\begin{tangle}
\hbox{the right side of (1). }
\end{tangle}
\]

(2) We can similarly show relation (2). \ \ \ \ $\Box$

\begin {Lemma} \label {6.1.6}
If the antipode of $H$ is  invertible, then  $\lambda $ is
invertible with inverse
\[  \lambda ^{-1}\step=\step
\begin{tangle}
\object{H}\step[2]\object{H^{\hat{*}}} \\

\id\step[2]\id\step [3] \coev \Step\\

\id\step[2]\id\step [2] \cd \step \id \\

\id\step[2]\id\step [2] \O {\bar S }\step \step\id \step\id\\

\id\step[2]\XX\step [2] \id  \step\id \\

\cu \step[2]\ev  \step \id \\

\object{H}\step[2]\object{H^{\hat{*}}}

\end{tangle}\ \ .
\]
\end {Lemma}

\begin {Lemma} \label {1.7}
$R \# H$ becomes    an  $H^{\hat  *}$-module algebra under the
 module operation  $\rightharpoonup ' \ =  (id_{H^{\hat *}} \otimes
\rightharpoonup )(C_{H^{\hat *}, R} \otimes id_H ).$
\end {Lemma}

          \begin {Theorem} \label {6.1.8}
Let $H$ be a Hopf algebra living in a symmetric  tensor category
${\cal C}$. If $H$ has a left duality $H^*$ and the antipode $S$ of
$H$ is invertible, then
  $$   (R \# H)\# H^{\hat *}   \cong R \otimes (H \bar
  \otimes H^{\hat *} ) \hbox { \ \ \  as algebras }$$
where  $H \bar  \otimes H^{\hat *}$  is defined in Lemma \ref
{6.1.31}.

\end {Theorem}

(i)
\[
\step Let\step[3]
\begin{tangle}
\object{H^{\hat{*}}}\\
\id\\
\O w\\
\id\\
\object{H\# H^{\hat{*}}}\\
\end{tangle}
\step=\step
\begin{tangle}
\object{H^{\hat{*}}}\step[3]\object{\eta_H}\\
\step[0.5]\O { \bar S}\Step\id\\
\step[0.5]\tu \rho\\
\obox 3{\lambda^{-1}}\\
\step[1.5]\id\\
\step[1.5]\object{H\# H^{\hat{*}}}\\
\end{tangle}
\Step and\Step
\begin{tangle}
\step\object{R}\\
\td \phi\\
\object{H^{\hat{*}}}\Step\object{R}\\
\end{tangle}
\step=\step
\begin{tangle}
\step[4]\object{R}\\
\ro {b'} \Step\id\\
\id\Step\tu \alpha\\
\object{H^{\hat{*}}}\step[3]\object{R}\\
\end{tangle}
\step=\step
\begin{tangle}
\step[4]\object{R}\\
\coev\Step\id\\
\XX\Step\id\\
\id\Step\tu \alpha\\
\object{H^{\hat{*}}}\step[3]\object{R}\\
\end{tangle}\ \ \ ,\]
where  $b' = C_{H, H^*}b_H . $

We define

\[
\begin{tangle}
\object{R\otimes(H\# H^{\hat{*}})}\\
\id\\
\O \Psi\\
\id\\
\object{(R\# H)\# H^{\hat{*}}}\\
\end{tangle}
\step[3]=\step[3]
\begin{tangle}
\step\object{R}\step[4]\object{H}\Step\object{H^{\hat{*}}}\\
\td \phi\step[3]\id\Step\id\\
\O S\Step\id\step[3]\id\Step\id\\
\x\step[3]\id\Step\id\\
\id\step\td w\Step\id\Step\id\\
\id\step\id\step\cd\step\id\Step\id\\
\id\step\id\step\id\Step\X\Step\id\\
\id\step\id\step\tu \rightharpoonup\step\cu\\
\id\step\cu\step[3]\id\\
\object{R}\Step\object{H}\step[4]\object{H^{\hat{*}}}\\
\end{tangle}
\step[3]and\step[5]
\begin{tangle}
\object{(R\# H)\# H^{\hat{*}}}\\
\id\\
\O \Phi\\
\id\\
\object{R\otimes(H\# H^{\hat{*}})}\\
\end{tangle}
\step[3]=\step[3]
\begin{tangle}
\step\object{R}\step[4]\object{H}\Step\object{H^{\hat{*}}}\\
\td \phi\step[3]\id\Step\id\\
\x\step[3]\id\Step\id\\
\id\step\td w\Step\id\Step\id\\
\id\step\id\step\cd\step\id\Step\id\\
\id\step\id\step\id\Step\X\Step\id\\
\id\step\id\step\tu \rightharpoonup\step\cu\\
\id\step\cu\step[3]\id\\
\object{R}\Step\object{H}\step[4]\object{H^{\hat{*}}}\\
\end{tangle}\ \ \ .
\]

We see that
\[
\Psi\Phi\step=\step
\begin{tangle}
\step[7]\object{R}\step[5]\object{H\# H^{\hat{*}}}\\
\ro {b'} \step\ro {b'}\Step\id\step[4]\dd\\
\O S\Step\id\step\d\step\tu \alpha\step[3]\dd\\
\O w\Step\id\Step\x\step[3]\dd\\
\d\step\tu \alpha\Step\O w\Step\dd\\
\step\x\Step\dd\step\dd\\
\step\id\Step\cu\step\dd\\
\step\id\step[3]\cu\\
\step\object{R}\step[4]\object{H\# H^{\hat{*}}}\\
\end{tangle}
\step \stackrel{\hbox{ since } w \hbox{ is algebraic }}{=}\step
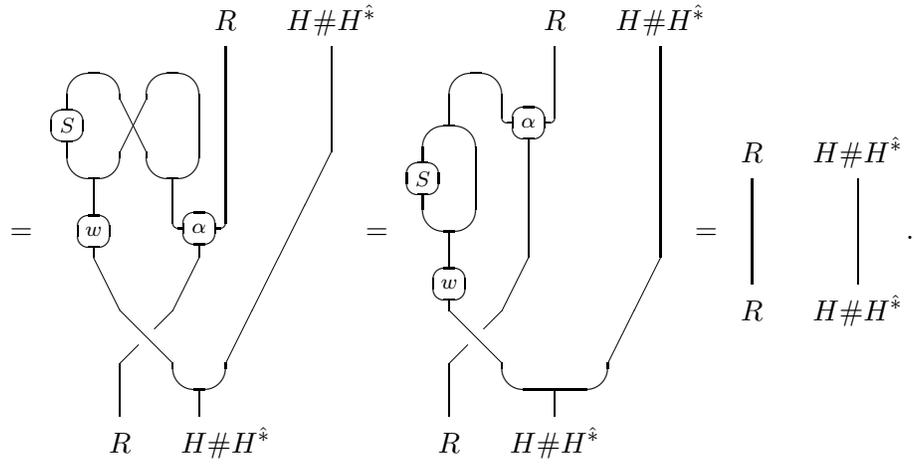
\begin{tangle}
\step[7]\object{R}\step[4]\object{H\# H^{\hat{*}}}\\
\step[3]\coev\Step\id\step[4]\id\\
\coev\step\d\step\tu \alpha\step[4]\id\\
\O S\Step\id\Step\x\step[4]\dd\\
\d\step\tu \alpha\Step\id\step[3]\dd\\
\step\x\Step\dd\step[2]\dd\\
\step\id\Step\cu\step[2]\dd\\
\step\id\step[3]\O w\Step\dd\\
\step\id\step[3]\cu\\
\step\object{R}\step[5]\object{H\# H^{\hat{*}}}\\
\end{tangle} \]
\[
\step=\step
\begin{tangle}
\step[6]\object{R}\step[4]\object{H\# H^{\hat{*}}}\\
\coev\step\coev\step\id\step[4]\id\\
\O S\Step\X\Step\id\step\id\step[4]\id\\
\cu\step\cu\step\id\step[3]\dd\\
\step\O w\step[3]\tu \alpha\Step\dd\\
\step\d\Step\dd\Step\dd\\
\Step\x\Step\dd\\
\Step\id\step[2]\cu\\
\Step\object{R}\step[4]\object{H\# H^{\hat{*}}}\\
\end{tangle}
\step=\step
\begin{tangle}
\step[5]\object{R}\step[4]\object{H\# H^{\hat{*}}}\\
\step\coev\Step\id\step[4]\id\\
\cd\step\tu \alpha\step[4]\id\\
\O S\Step\id\Step\id\step[5]\id\\
\cu\Step\id\step[5]\id\\
\step\O w\Step\dd\step[4]\dd\\
\step\x\step[4]\dd\\
\step\id\Step\Cu\\
\step\object{R}\step[4]\object{H\# H^{\hat{*}}}\\
\end{tangle}
\step=\step
\begin{tangle}
\object{R}\step[4]\object{H\# H^{\hat{*}}}\\
\id\step[4]\id\\
\id\step[4]\id\\
\object{R}\step[4]\object{H\# H^{\hat{*}}}\\
\end{tangle} \ \ \ \ \ .
\]
Similarly, we have  $\Phi \Psi = id $. Thus $\Phi$ is invertible.

Now we show that $\Phi$ is algebraic.
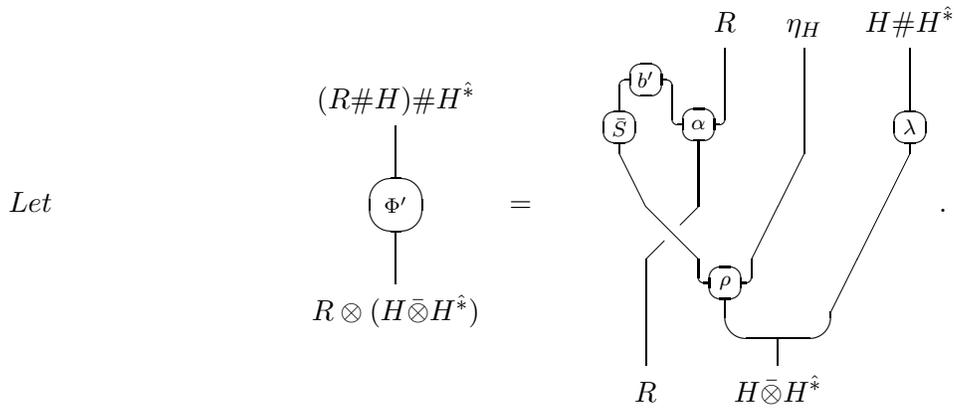
\[
Let\step[12]
\begin{tangle}
\step\object{(R\# H)\# H^{\hat{*}}}\\
\step\id\\
\obox 2{\Phi'}\\
\step\id\\
\step\object{R\otimes(H\bar{\otimes}H^{\hat{*}})}\\
\end{tangle}
\step[3]=\step[3]
\begin{tangle}
\step[4]\object{R}\step[3]\object{\eta_H}\step[4]\object{H\# H^{\hat{*}}}\\
\ro {b'}\Step\id\step[3]\id\step[4]\id\\
\O {\bar{S}}\Step\tu \alpha\step[3]\id\step[4]\O \lambda\\
\d\Step\id\step[3]\dd\step[3]\dd\\
\step\x\Step\dd\step[3]\dd\\
\step\id\Step\tu \rho\step[3]\dd\\
\step\id\step[3]\Cu\\
\step\object{R}\step[5]\object{H\bar{\otimes}H^{\hat{*}}}\\
\end{tangle}\ \ \ .
\]

It is clear that $\Phi = (id \otimes \lambda ^{-1})\Phi '$.
Consequently, we only need show that $\Phi '$ is alegebraic.

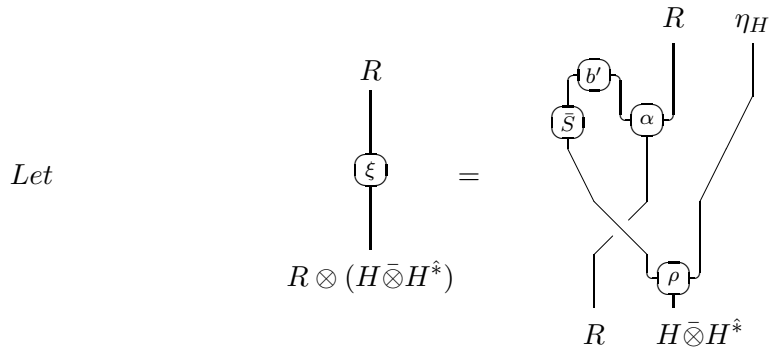
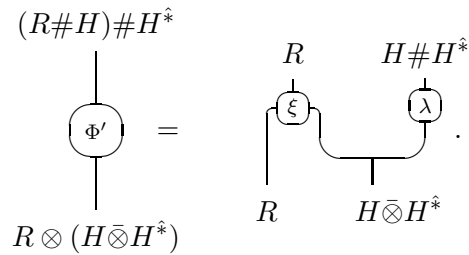
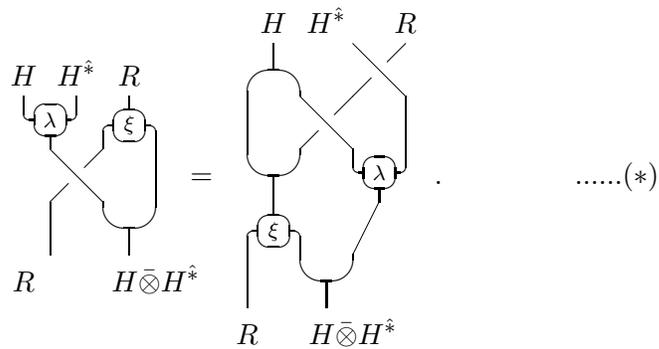
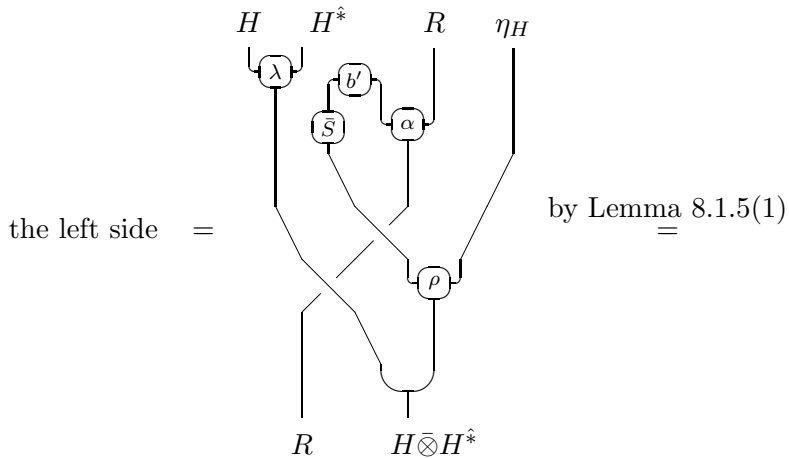
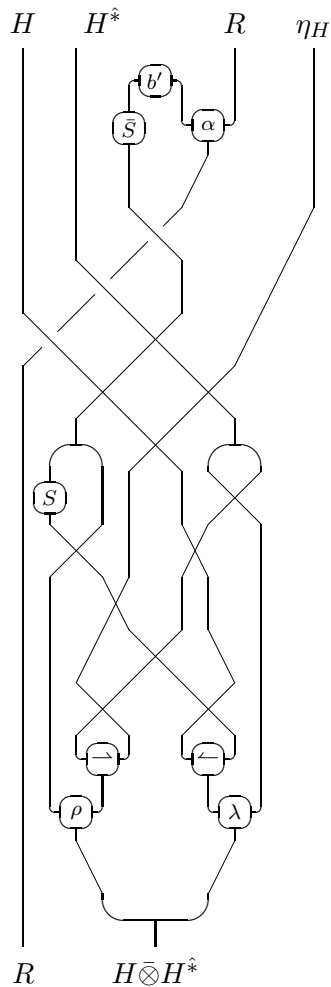
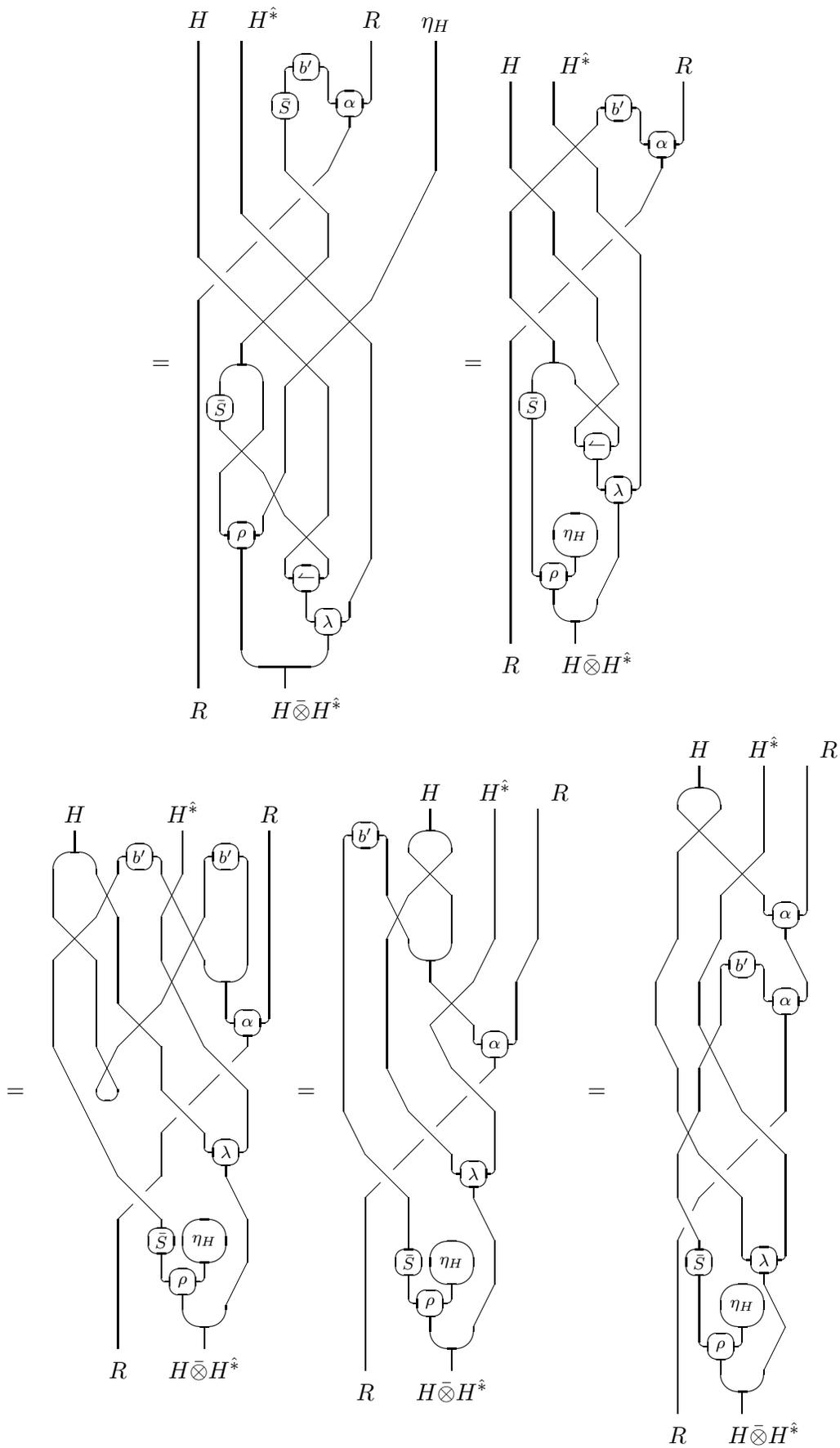
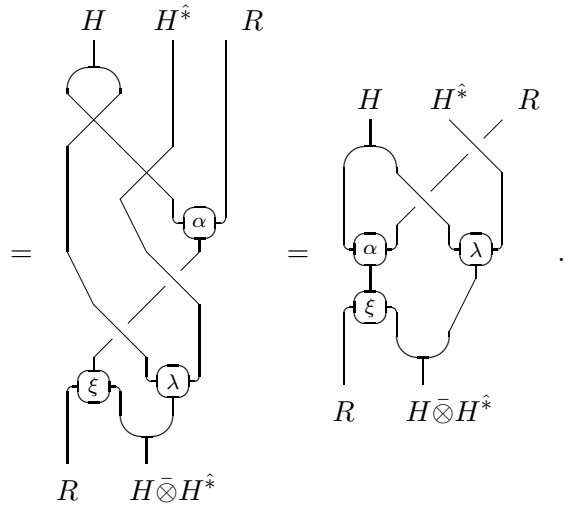
\[
Let\step[12]
\begin{tangle}
\object{R}\\
\id\\
\O \xi\\
\id\\
\object{R\otimes(H\bar{\otimes}H^{\hat{*}})}\\
\end{tangle}
\step[3]=\step[3]
\begin{tangle}
\step[4]\object{R}\step[3]\object{\eta_H}\\
\ro {b'}\Step\id\step[3]\id\\
\O {\bar{S}}\Step\tu \alpha\Step\dd\\
\d\Step\id\Step\dd\\
\step\x\step[2]\id\\
\step\id\Step\tu \rho\\
\step\object{R}\step[4]\object{H\bar{\otimes}H^{\hat{*}}}\\
\end{tangle}\ \ \ .
\]

We have that

\[
\begin{tangle}
\step\object{(R\# H)\# H^{\hat{*}}}\\
\step\id\\
\obox 2{\Phi'}\\
\step\id\\
\step\object{R\otimes(H\bar{\otimes}H^{\hat{*}})}\\
\end{tangle}
\step=\step[3]
\begin{tangle}
\step\object{R}\step[5]\object{H\# H^{\hat{*}}}\\
\td \xi\step[4]\O \lambda\\
\id\Step\Cu\\
\object{R}\step[5]\object{H\bar{\otimes} H^{\hat{*}}}\\
\end{tangle}\ \ \ .
\]

We claim that

\[
\begin{tangle}
\object{H}\step[2]\object{H^{\hat{*}}}\Step\object{R}\\
\tu \lambda\step\td \xi\\
\step\x\Step\id\\
\step\id\Step\cu\\
\object{R}\step[5]\object{H\bar{\otimes} H^{\hat{*}}}\\
\end{tangle}
\step=\step
\begin{tangle}
\step\object{H}\step[2]\object{H^{\hat{*}}}\step[3]\object{R}\\
\cd\Step\x\\
\id\Step\x\Step\id\\
\cu\Step\tu \lambda\\
\td \xi\Step\dd\\
\id\Step\cu\\
\object{R}\step[4]\object{H\bar{\otimes} H^{\hat{*}}}\\
\end{tangle}\ \ \ .
\step[5] ......(*)
\]

\[
\hbox{the left side}\step=\step
\begin{tangle}
\object{H}\step[3]\object{H^{\hat{*}}}\step[4]\object{R}\step[3]\object{\eta_H}\\
\tu \lambda\step\ro {b'}\Step\id\step[3]\id\\
\step\id\Step\O {\bar{S}}\Step\tu \alpha\step[3]\id\\
\step\id\Step\d\Step\id\step[3]\dd\\
\step\d\Step\x\Step\dd\\
\Step\x\Step\tu \rho\\
\Step\id\Step\d\Step\id\\
\Step\id\step[3]\cu\\
\Step\object{R}\step[5]\object{H\bar{\otimes} H^{\hat{*}}}\\
\end{tangle}
\step\stackrel{\hbox{by Lemma \ref
{6.1.5}(1)}}{\step[3]=\step[3]}\step
\begin{tangle}
\object{H}\step[3]\object{H^{\hat{*}}}\step[5]\object{R}\step[3]\object{\eta_H}\\
\id\Step\id\Step\ro {b'}\Step\id\step[3]\id\\
\id\Step\id\Step\O {\bar{S}}\Step\tu \alpha\step[3]\id\\
\id\Step\id\Step\id\Step\dd\step[4]\id\\
\id\Step\id\Step\x\step[4]\dd\\
\id\Step\x\Step\id\step[3]\dd\\
\x\Step\XX\step[2]\dd\\
\id\Step\XX\Step\XX\\
\id\step\cd\step\XX\step\cd\\
\id\step\O S\Step\id\step\id\Step\id\step\XX\\
\id\step\XX\step\id\Step\X\Step\id\\
\id\step\id\Step\X\Step\id\step\id\Step\id\\
\id\step\id\step\dd\step\XX\step\d\step\id\\
\id\step\id\step\XX\Step\XX\step\id\\
\id\step\id\step\tu \rightharpoonup\Step\tu
\leftharpoonup\step\id\\
\id\step\tu \rho\step[4]\tu \lambda\\
\id\Step\d\step[4]\dd\\
\id\step[3]\Cu\\
\object{R}\step[5]\object{H\bar{\otimes} H^{\hat{*}}}\\
\end{tangle} \]
\[
\step=\step
\begin{tangle}
\object{H}\step[3]\object{H^{\hat{*}}}\step[5]\object{R}\step[3]\object{\eta_H}\\
\id\Step\id\Step\ro {b'}\Step\id\step[3]\id\\
\id\Step\id\Step\O {\bar{S}}\Step\tu \alpha\step[3]\id\\
\id\Step\id\Step\id\Step\dd\step[4]\id\\
\id\Step\id\Step\x\step[4]\dd\\
\id\Step\x\Step\id\step[3]\dd\\
\x\Step\XX\step[2]\dd\\
\id\Step\XX\Step\XX\\
\id\step\cd\step\XX\Step\id\\
\id\step\O {\bar S}\Step\id\step\id\Step\id\Step\id\\
\id\step\XX\step\id\Step\id\Step\id\\
\id\step\id\Step\X\Step\id\Step\id\\
\id\step\tu \rho\step\XX\Step\id\\
\id\Step\id\Step\tu \leftharpoonup\step\dd\\
\id\Step\id\step[3]\tu \lambda\\
\id\Step\Cu\\
\object{R}\step[5]\object{H\bar{\otimes} H^{\hat{*}}}\\
\end{tangle}
\step=\step
\begin{tangle}
\object{H}\step[3]\object{H^{\hat{*}}}\step[5]\object{R}\\
\id\Step\id\Step\ro {b'}\Step\id\\
\id\Step\XX\Step\tu \alpha\\
\XX\Step\id\Step\dd\\
\id\Step\id\Step\x\\
\id\Step\x\Step\id\\
\x\Step\id\Step\id\\
\id\step\cd\step\d\step\id\\
\id\step\O {\bar S}\Step\XX\step\id\\
\id\step\id\Step\tu \leftharpoonup\step\id\\
\id\step\id\step[3]\tu \lambda\\
\id\step\id\step\obox 2{\eta_H}\step\id\\
\id\step\tu \rho\step\dd\\
\id\Step\cu\\
\object{R}\step[4]\object{H\bar{\otimes} H^{\hat{*}}}\\
\end{tangle}
\]
\[
\step=\step
\begin{tangle}
\step\object{H}\step[5]\object{H^{\hat{*}}}\step[4]\object{R}\\
\cd\step\ro {b'}\step\id\step\ro {b'}\step\id\\
\id\Step\X\Step\X\step\id\Step\id\step\id\\
\XX\step\id\Step\id\step\X\Step\id\step\id\\
\id\Step\id\step\id\Step\X\step\cu\step\id\\
\id\Step\id\step\XX\step\d\step\tu \alpha\\
\d\step\X\Step\id\Step\x\\
\step\d\hev\Step\x\Step\id\\
\Step\d\Step\id\Step\tu \lambda\\
\step[3]\x\step[3]\d\\
\step[3]\id\Step\O {\bar{S}}\step\obox 2{\eta_H}\step\id\\
\step[3]\id\Step\tu \rho\step\dd\\
\step[3]\id\step[3]\cu\\
\step[3]\object{R}\step[4]\object{H\bar{\otimes} H^{\hat{*}}}\\
\end{tangle}
\step=\step
\begin{tangle}
\step[4]\object{H}\step[3]\object{H^{\hat{*}}}\step[3]\object{R}\\
\ro {b'}\step\cd\step[2]\id\Step\id\\
\id\Step\id\step\XX\Step\id\Step\id\\
\id\Step\X\Step\id\Step\id\Step\id\\
\id\Step\id\step\cu\step\dd\step\dd\\
\id\Step\id\step [2]\XX\step [2]\id\\
\id\Step\id\step [2]\d\step\tu \alpha\\
\id\Step\d\Step\x\\
\d\Step\x\Step\id\\
\step\x\Step\tu \lambda\\
\step\id\Step\id\step[3]\d\\
\step\id\Step\O {\bar{S}}\step\obox 2{\eta_H}\step\id\\
\step\id\Step\tu \rho\step\dd\\
\step\id\step[3]\cu\\
\step\object{R}\step[4]\object{H\bar{\otimes} H^{\hat{*}}}\\
\end{tangle}
\step=\step
\begin{tangle}
\step[3]\object{H}\step[3]\object{H^{\hat{*}}}\step[3]\object{R}\\
\Step\cd\Step\id\Step\id\\
\Step\XX\Step\id\Step\id\\
\step [2]\id \step[2]\XX  \Step\id\\
\step [2]\id \step[2]\id \step [2]\tu \alpha\\
\step \dd \step \dd\ro {b'}\step [1]\d\\
\step \id \step[2]\id \step \id\step [2]\tu \alpha\\
\step \d \step \X\step[3]\id\\
\step[2]\id \step[1]\id \step \d \step [2] \id \\
\step[2]\X\step[2]\x \\
\step [2] \id \step [1] \x \step [2]\id \\
\step[2]\hx\step[2]\id\Step\id\\
\step [2]\id \step\O {\bar S}\step \step \tu \lambda\\
\step [2]\id\step\id\step\obox 2{\eta_H}\d\\
\step [2]\id\step\tu \rho\step\dd\\
\step [2]\id\step[2]\cu\\
\step [2]\object{R}\step[4]\object{H\bar{\otimes} H^{\hat{*}}}\\
\end{tangle}\]
\[
\step=\step
\begin{tangle}
\step\object{H}\step[3]\object{H^{\hat{*}}}\step[3]\object{R}\\
\cd\Step\id\Step\id\\
\XX\Step\id\Step\id\\
\id\Step\XX\Step\id\\
\id\Step\d\step\tu \alpha\\
\d\Step\x\\
\step\x\Step\id\\
\td \xi\step\tu \lambda\\
\id\Step\cu\\
\object{R}\step[4]\object{H\bar{\otimes} H^{\hat{*}}}\\
\end{tangle}
\step=\step
\begin{tangle}
\step\object{H}\step[3]\object{H^{\hat{*}}}\step[3]\object{R}\\
\cd\Step\x\\
\id\Step\x\Step\id\\
\tu \alpha\Step\tu \lambda\\
\td \xi\Step\dd\\
\id\Step\cu\\
\object{R}\step[4]\object{H\bar{\otimes} H^{\hat{*}}}\\
\end{tangle}\ \ \ .
\]
Thus relation (*) holds.

Next, we check that $\xi$ is algebraic. We see that
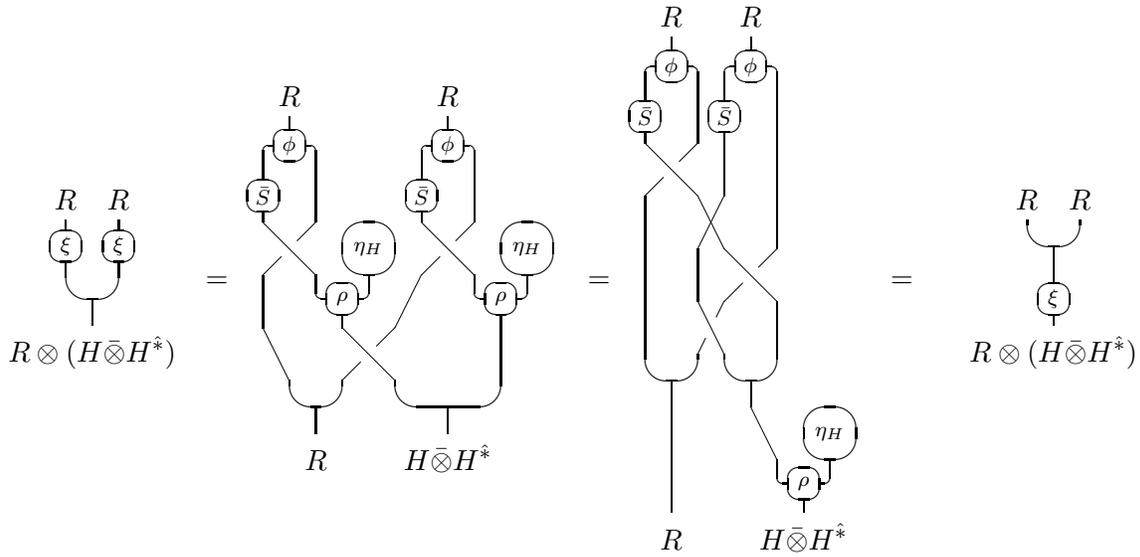
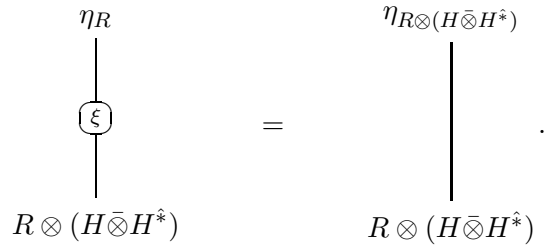
\[
\begin{tangle}
\object{R}\Step\object{R}\\
\O \xi\Step\O \xi\\
\cu\\
\step\object{R\otimes(H\bar{\otimes}H^{\hat{*}})}\\
\end{tangle}
\step[3]=\step
\begin{tangle}
\step\object{R}\step[6]\object{R}\\
\td \phi\step[4]\td \phi\\
\O {\bar{S}}\Step\id\step[4]\O {\bar{S}}\Step\id\\
\x\step\obox 2{\eta_H}\step\x\step\obox 2{\eta_H}\\
\id\Step\tu \rho\step\dd\Step\tu \rho\\
\d\Step\x\step[4]\id\\
\step\cu\Step\Cu\\
\Step\object{R}\step[5]\object{H\bar{\otimes} H^{\hat{*}}}\\
\end{tangle}
\step=\step
\begin{tangle}
\step\object{R}\step[3]\object{R}\\
\td \phi\step\td \phi\\
\O {\bar{S}}\Step\id\step\O {\bar{S}}\Step\id\\
\x\step\id\Step\id\\
\id\Step\X\Step\id\\
\id\Step\id\step\x\\
\id\Step\hx\Step\id\\
\cu\step\cu\\\
\step\id\step[3]\d\step\obox 2{\eta_H}\\
\step\id\step[4]\tu \rho\\
\step\object{R}\step[5]\object{H\bar{\otimes} H^{\hat{*}}}\\
\end{tangle}
\step=\step[4]
\begin{tangle}
\object{R}\Step\object{R}\\
\cu\\
\step\O \xi\\
\step\object{R\otimes(H\bar{\otimes}H^{\hat{*}})}\\
\end{tangle}
\]
and obviously
\[
\begin{tangle}
\object{\eta_R}\\
\id\\
\O \xi\\
\id\\
\object{R\otimes(H\bar{\otimes}H^{\hat{*}})}\\
\end{tangle}
\step[6]=\step[6]
\begin{tangle}
\object{\eta_{R\otimes(H\bar{\otimes}H^{\hat{*}})}}\\
\id\\
\id\\
\id\\
\object{R\otimes(H\bar{\otimes}H^{\hat{*}})}\\
\end{tangle}\ \ \ \ \ \ \ \ \ .
\]
Thus $\xi$ is algebraic.

Now we show that $\Phi'$ is algebraic.
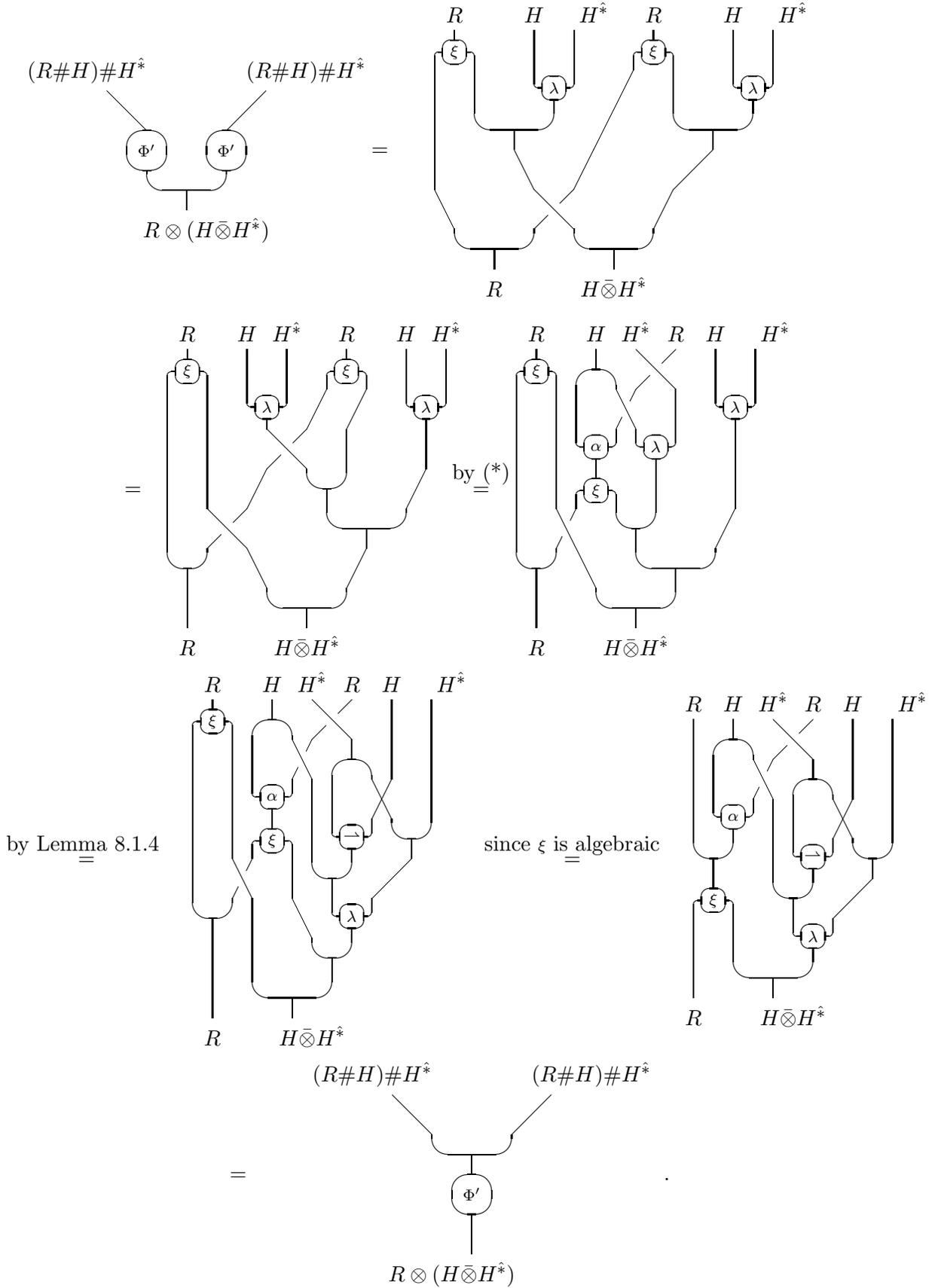
\[
\begin{tangle}
\object{(R\# H)\# H^{\hat{*}}}\step[11]\object{(R\# H)\# H^{\hat{*}}}\\
\step\nw2\step[6]\ne2\\
\step[2]\obox 2{\Phi'}\step[2]\obox 2{\Phi'}\\
\step[3]\Cu\\
\step[6]\object{R\otimes(H\bar{\otimes}H^{\hat{*}})}\\
\end{tangle}
\step[3]=\Step
\begin{tangle}
\step\object{R}\step[4]\object{H}\step[3]\object{H^{\hat{*}}}
\step[3]\object{R}\step[4]\object{H}\step[3]\object{H^{\hat{*}}}\\
\td \xi\step[3]\id\Step\id\step[3]\td \xi\step[3]\id\Step\id\\
\id\Step\id\step[3]\tu \lambda\Step\dd\Step\id\step[3]\tu \lambda\\
\id\Step\Cu\Step\dd\step[3]\Cu\\
\id\step[4]\d\Step\dd\step[5]\ne2\\
\d\step[4]\x\step[4]\dd\\
\step\Cu\Step\Cu\\
\step[3]\object{R}\step[6]\object{H\bar{\otimes}H^{\hat{*}}}\\
\end{tangle}\]
\[
\step=\step
\begin{tangle}
\step\object{R}\step[3]\object{H}\step[2]\object{H^{\hat{*}}}
\step[3]\object{R}\step[3]\object{H}\step[2]\object{H^{\hat{*}}}\\
\td \xi\step[2]\id\Step\id\step[2]\td \xi\step[2]\id\Step\id\\
\id\Step\id\Step\tu \lambda\step\dd\step\dd\Step\tu \lambda\\
\id\Step\id\step[3]\x\Step\id\step[4]\id\\
\id\Step\id\Step\dd\Step\cu\step[3]\dd\\
\id\Step\x\step[4]\Cu\\
\cu\Step\d\step[4]\dd\\
\step\id\step[4]\Cu\\
\step\object{R}\step[6]\object{H\bar{\otimes}H^{\hat{*}}}\\
\end{tangle}
\stackrel{\hbox{by (*)}}{=}
\begin{tangle}
\step\object{R}\step[3]\object{H}\step[2]\object{H^{\hat{*}}}
\step[2]\object{R}\step[2]\object{H}\step[3]\object{H^{\hat{*}}}\\
\td \xi\step\cd\step\x\Step\id\Step\id\\
\id\Step\id\step\id\Step\hx\Step\id\Step\tu \lambda\\
\id\Step\id\step\tu \alpha\step\tu \lambda\step[3]\id\\
\id\Step\id\step\td \xi\Step\id\step[4]\id\\
\id\Step\hx\Step\cu\step[3]\dd\\
\cu\step\d\Step\Cu\\
\step\id\step[3]\Cu\\
\step\object{R}\step[5]\object{H\bar{\otimes}H^{\hat{*}}}\\
\end{tangle}
\]
\[
\step\stackrel{\hbox{by Lemma  \ref  {6.1.4}   }}{=}\step
\begin{tangle}
\step\object{R}\step[3]\object{H}\step[2]\object{H^{\hat{*}}}
\step[2]\object{R}\step[2]\object{H}\step[3]\object{H^{\hat{*}}}\\
\td \xi\step\cd\step\x\Step\id\Step\id\\
\id\Step\id\step\id\Step\hx\step\cd\step\id\Step\id\\
\id\Step\id\step\tu \alpha\step\id\step\id\Step\X\Step\id\\
\id\Step\id\step\td \xi\step\id\step\tu \rightharpoonup\step\cu\\
\id\Step\hx\Step\id\step\cu\step[2]\ne2\\
\cu\step\id\Step\d\step\tu \lambda\\
\step\id\Step\id\step[3]\cu\\
\step\id\Step\Cu\\
\step\object{R}\step[5]\object{H\bar{\otimes}H^{\hat{*}}}\\
\end{tangle}
\step\stackrel{\hbox{ since } \xi \hbox{ is algebraic }}{=}\step
\begin{tangle}
\object{R}\step[2]\object{H}\step[2]\object{H^{\hat{*}}}
\step[2]\object{R}\step[2]\object{H}\step[3]\object{H^{\hat{*}}}\\
\id\step\cd\step\x\Step\id\Step\id\\
\id\step\id\Step\hx\step\cd\step\id\Step\id\\
\id\step\tu \alpha\step\id\step\id\Step\X\Step\id\\
\cu\Step\id\step\tu \rightharpoonup\step\cu\\
\td \xi\Step\cu\Step\ne2\\
\id\Step\id\step[3]\tu \lambda\\
\id\Step\Cu\\
\object{R}\step[5]\object{H\bar{\otimes}H^{\hat{*}}}\\
\end{tangle} \]
\[
\step=\step[6]
\begin{tangle}
\object{(R\# H)\# H^{\hat{*}}}\step[11]\object{(R\# H)\# H^{\hat{*}}}\\
\step\nw2\step[6]\ne2\\
\step[3]\Cu\\
\step[4]\obox 2{\Phi'}\\
\step[5]\id\\
\step[4]\object{R\otimes(H\bar{\otimes}H^{\hat{*}})}\\
\end{tangle}\ \ \ \ \ \ \ \ \ \ .
\]

It is clear that
\[
\begin{tangle}
\object{\eta_{(R\# H)\# H^{\hat{*}}}}\\
\step\id\\
\obox 2{\Phi'}\\
\step\id\\
\step\object{R\otimes(H\bar{\otimes}H^{\hat{*}})}\\
\end{tangle}
\step[6]=\step[6]
\begin{tangle}
\object{\eta_{R\otimes(H\bar{\otimes}H^{\hat{*}})}}\\
\id\\
\id\\
\id\\
\object{R\otimes(H\bar{\otimes}H^{\hat{*}})}\\
\end{tangle}\ \ \ \ \ \ \ \ \ \ \ \ \ .\] Thus $\Phi '$ is algebraic. \ \
\begin{picture}(5,5)
\put(0,0){\line(0,1){5}}\put(5,5){\line(0,-1){5}}
\put(0,0){\line(1,0){5}}\put(5,5){\line(-1,0){5}}
\end{picture}\\

\section { The duality theory for Yetter-Drinfeld  Hopf algebras  }\label {s16}

In this section,
 we give the duality theory  for  Yetter-Drinfeld  Hopf algebra $H$.

A braided Hopf algebra in  Yetter -Drinfeld categories is  called a
Yetter-Drinfeld  Hopf algebra.

In this section,  all Yetter -Drinfeld categories are symmetric.

If $H$ is a finite-dimensional Hopf algebra living in a symmetric
tensor category ${\cal C}$ determined by (co)triangular structure,
then $H \bar \otimes H^{\hat *} = End_k (H)$ since $ dim (H \bar
\otimes H^{\hat *}) = dim (End _kH)$ and $H \bar \otimes H^{\hat *}$
is a subalgebra of $End _k H.$ We can easily check that $$R \otimes
End_k H \cong M_n (R) \hbox { \ \ \  as algebras }
$$
where $M_n(R)  = \{x \mid x \hbox { is a } n \times n\hbox {-matrix
over } k \}$

Considering  Theorem \ref {6.1.8}, we have :

\begin {Theorem} \label {6.2.2}  (Duality theorem)
Let $H$ be a finite-dimensional Yetter-Drinfeld Hopf algebra  Then

 $$(R \# H)\# H^{\hat *}   \cong M_n(R) \hbox { \ \ \  (as algebras). }$$

\end {Theorem}
    Some interesting braided group crossed products are the  Weyl algebras,
       which have been studied,  for example in \cite {Ma93b}. Now we give an example
       in super case.

       \begin {Example} \label {6.2.3}
       (see \cite [P502]{Ma95b}) Let $R$ denote the Gassmann plane
       $C^{0 \mid 0 }$. That is, $R= {\bf C} (\theta )$ with odd coordinate
       $\theta $ and relation $\theta ^2 =0.$  The comultiplication
       and antipode are
       $$\Delta (\theta ) = \theta \otimes 1 + 1 \otimes \theta
       \hbox { \ \ \
        and \ \ \ } S(\theta ) = -\theta .$$
        $R$  is a super-Hopf algebra with $dim \ R = 2.$
        Let $H=R$  and
       $\alpha $ be  the adjoint action.
        Thus  $(R, \alpha )$  is a left $H$-module algebra in super case.
        Since $H$ is finite-dimensional we have that

 $$(R \# H)\# H^{\hat *}   \cong M_n(R) \hbox { \ \ \  (as superalgebras). }$$
              by Theorem \ref {6.2.2}.  \begin{picture}(8,8)\put(0,0){\line(0,1){8}}\put(8,8){\line(0,-1){8}}\put(0,0){\line(1,0){8}}\put(8,8){\line(-1,0){8}}\end{picture}

       \end {Example}

       In fact, we can obtain duality theory for crossed coproducts  by
       only
       turning the above diagrams upside down.

\part {Hopf Algebras Living in the  Category of Vector spaces}

 \chapter { The Relation Between Decomposition of
Comodules and Coalgebras}\label {c8}

The decomposition of coalgebras and comodules is an important
subject in study of Hopf algebras.
  T. Shudo and H. Miyamito \cite{SM78} showed that
  $C$ can  be decomposed
into a direct sum  of its indecomposable subcoalgebras of $C$.
 Y.H. Xu \cite {XF92}  showed that the decomposition  was unique.
He also  showed that
  $M$ can uniquely  be decomposed
into a direct sum  of the weak-closed indecomposable subcomodules of
$M$(we call the decomposition
 the weak-closed indecomposable
decomposition ) in  \cite{XSF94}. In this chapter, we   give
 the relation between the
   two decomposition.
   We show that if $M$  is a full, $W$-relational
   hereditary $C$-comodule, then the following conclusions hold:

(1) $M$ is indecomposable iff $C$ is indecomposable;

(2) $M$ is relative-irreducible iff $C$ is irreducible;

(3) $M$ can be decomposed into a direct sum of
 its weak-closed  relative-irreducible subcomodules iff
$C$ can be decomposed into a direct sum of its irreducible
subcoalgebras.     \\
We also obtain the relation between  coradical  of $C$- comodule $M$
and radical of algebra  $C(M)^*$

Let $k$ be a field, $M$ be a $C$-comodule, $N$ be a subcomodule of
$M$, $E$ be a subcoalgebras of  $C$ and $P$ be an ideal of  $C^*$.
As in  \cite{XSF94},  we define:

$E^{ \bot C^* } = \{ f \in C^* \mid  f(E) = 0 \}$.

$P^{ \bot C } = \{ c \in C \mid  P(c) = 0 \}$.

$N^{ \bot C^* } = \{ f \in C^* \mid  f \cdot N = 0 \}$.

$P^{ \bot M } = \{ x \in M \mid  P \cdot x = 0 \}$.

Let $<N>$ denote  $N^{ \bot C^*  \bot M } $.
 $<N>$ is called the closure of  $N$. If  $<N> = N$,
then  $N$ is called closed. If  $N  = C^*x$, then $<N>$ is denoted
by $<x>$. If  $<x> \subseteq N $ for any $x \in N$, then $N$ is
called weak-closed. It is clear that any closed subcomodule is
weak-closed. If  $\rho (N) \subseteq N \otimes E  $, then $N$ is
called an $E$-subcomodule of $M$. Let
$$M_E = \sum \{ N \mid N \hbox {  is  a  subcomodule  of } M
\hbox { and } \rho (N) \subseteq N \otimes E \}. $$ We call  $M_E$ a
component of $M$ over $E$. If $M_E$ is some component of $M$ and
$M_F \subseteq M_E$ always implies  $M_F = M_E$ for any non-zero
component $M_F$,  then  $M_E$ is called a minimal component of $M$.

Let $\{m_{ \lambda }  \mid  \lambda \in \Lambda \}$ be the basis of
$M$ and $C(M)$ denote the subspace of $C$ spanned by

$W(M) = \{ c \in C \mid  {~there~ exists~ an~} m \in M
 {~with~}
\rho(m) = \sum m_{ \lambda } \otimes c_{\lambda } {~such ~that~ }
 c_{\lambda _0 } = c  {~for ~some~}   \lambda _0 \in \Lambda  \}$.\\
E. Abe  in \cite[P129]{Ab80} checked that  $C(M)$ is a subcoalgebra
of  $C$. It is easy to know that if $E$ is subcoalgebra of $C$ and $
\rho (M) \subseteq M \otimes E$, then $C(M) \subseteq E$. If $C(M) =
C$, then $M$ is called a full $C$-comodule. If $D$ is a simple
subcoalgbra  of   $C$ and  $M_D \not= 0$ or  $D=0$, then $D$ is
called faithful to  $M$. If  every simple subcoalgbra of $C$ is
faithful to  $M$, then  $M$ is called a component faithfulness
$C$-comodule.

     Let     $X$ and $Y$  be subspaces of coalgebra $C$.
     Define $X \wedge Y$  to be the kernel of the composite
     $$C \stackrel {\Delta } {\rightarrow } C \otimes C \rightarrow
     C/X \otimes C/Y. $$
     $X \wedge Y$  is called a wedge of $X$  and $Y$.

\section{The relation between the
decomposition  of comodules and coalgebras}\label {s17}
\begin{Lemma} \label{7.1.1}
 Let $N$ and $L$ be subcomodules of  $M$.
Let $D$ and $E$ be subcoalgebras of  $C$. Then

(1)  $N^{ \bot C^* } = (C(N))^{ \bot C^* }$;

(2)  $ M_E = M_{C(M_E)}$;

(3) $N$ is closed iff there exists a subcoalgebra  $E$ such that $N
=  M_E $;

(4)  $M_D \not= 0$ iff $D^{\bot C^* \bot M} \not= 0 $;

(5)  If $D \cap C(M) = 0$, then $M_D=0$;

(6)  If $D$ is a simple subcoalgebra of  $C$, then
$C(M_D) =   \left \{  \begin{array}{l l} D & { \rm if~} M_D \not= 0\\
0 & {\rm if~ }  M_D = 0
\end{array} \right.$;

(7) If $D\cap E=0$, then $M_D\cap M_E = 0$;

(8)  If  $D$ and $E$ are simple subcoalgebras and  $M_D \cap M_E =
0$ with  $M_D \not= 0$ or  $M_E \not= 0$,
 then  $D \cap E = 0$.

(9)  If $M_E$ is the minimal component of $M$ and  $0 \not= N
\subseteq M_E$,
 then  $C(M_E) = C(N)$.

\end{Lemma}
{\bf Proof}.
 (1)
Let $\{m_{\lambda} \mid \lambda\in\Lambda\}$ be a basis of $N$.
 For any  $f\in N^{\bot C^*}$, if $c\in W(N)$, then there exists
$m\in N$ with $\rho(m)= \sum_{\lambda\in\Lambda}m_{\lambda}\otimes
c_{\lambda}$ such that
 $c_{\lambda_0}=c$ for some $\lambda_0 \in\Lambda$.
  Since $f\cdot m=\sum m_{\lambda}
f(c_{\lambda})=0$, $f(c_{\lambda})=0$ for any $\lambda\in\Lambda$.
 Obviously, $f(c)=f(c_{\lambda _0})=0$.
Considering  that $C(N)$ is the space spanned by $W(N)$, we have
$f\in C(N)^{\bot C^*}$. Conversely, if $f\in C(N)^{\bot C^*}$, then
$f\cdot m=\sum m_{\lambda}f(c_{\lambda})=0$ for any $m\in N$, i.e.
$f\in N^{\bot C^*}$. This shows that $N^{\bot C^*}=C(N)^{\bot C^*}$.

(2) Since $M_E$ is an $E$-subcomodule of $M$, $C(M_E)\subseteq E$
and $M_{ C(M_E)}\subseteq M_E$. Since $M_E$ is a
$C(M_E)$-subcomodule, $M_E \subseteq M_{C(M_E)}$.

(3) If $N$ is closed, then $N^{\bot C^*\bot M}=N$. Let $E=C(N)$.
Obviously, $N\subseteq M_E$. By Lemma \ref{7.1.1}(1),
 $N^{ \bot C^{\star}}\cdot M_E=E^{\bot C^{\star}}
\cdot M_E=0$. Thus
 $M_E\subseteq N^{\bot C^{\star}\bot M}$ and $M_E\subseteq N$.
This shows that $M_E=N$. Conversely, if $M_E=N$,  obviously, $(M_E)
\subseteq (M_E)^{ \bot C^*  \bot M }.$ Thus it is sufficient to show
that
    $$<M_E> \subseteq  M_E.$$
    Let $L= (M_E)^{\bot C^* \bot M} = <M_E>.$ We see that
  $C(L)^{ \bot C^* } = L^{ \bot C^* }
= (M_E)^{ \bot C^* \bot M \bot C^* } = (M_E)^{ \bot C^* }= C(M_E)^{
\bot C^* }$. Thus $C(L) = C(M_E)$ and  $ L \subseteq M_{C(L)} = M
_{C(M_E)} = M_E $.

(4)  If  $D^{ \bot C^* \bot M } \not= 0$, let $\{m_{\lambda} \mid
\lambda \in\Lambda \}$ be a basis of $D^{\bot C^* \bot M}$ and $\rho
(x) = \sum m_ \lambda \otimes d_ \lambda $ for any $x \in D^{\bot
C^* \bot M } $.
 Since   $D^{ \bot C^* } \cdot x  =
\sum_{\lambda } m_\lambda  D^{ \bot C^* } (d_ \lambda) = 0$,
$d_\lambda
 \in
D^{ \bot C^* \bot C } = D$ for any $\lambda \in \Lambda$, which
implies
 that
$D^{ \bot C^* \bot M}$  is a $D$-subcomodule. Therefore $ 0 \not=
D^{ \bot C^* \bot M} \subseteq M_D$, i.e. $M_D \not= 0$.
 Conversely, if  $M_D\not= 0$,  then we have that
 \begin {eqnarray*}
 0 \not= M_D &\subseteq & (M_D)^{ \bot C^* \bot M} \\
 &  = & C(M_D)^{ \bot C^* \bot M} \hbox { \ \ (by part (1)) }  \\
& \subseteq & D^{ \bot C^* \bot M}.
\end {eqnarray*}

(5)  Since $M_D \subseteq M$, $C(M_D) \subseteq C(M)$. Obviously,
$C(M_D) \subseteq D$. Thus $C(M_D) \subseteq C(M) \cap D = 0$, which
implies  that   $M_D = 0$.

(6)  If  $M_D = 0$, then $C(M_D) = 0$. If $M_D \not= 0$,
 then $0 \not= C(M_D) \subseteq D$. Since   $D$ is a simple
subcoalgebra, $C(M_D) = D$.

(7)  If $x \in M_D \cap M_E$, then
 $\rho (x) \in (M_D \otimes D) \cap (M_E \otimes E)$,
and
 \begin {eqnarray} \rho (x) =   \sum_{ i = 1 }^{n}  x_i \otimes d_i =
 \sum_{ j = 1 }^{m}  y_j \otimes e_j,
 \label {e.1.1}
 \end {eqnarray}
where $x_i, \cdots , x_n$ is linearly independent and   $d_i \in D$
and $ e_j \in E$. Let $ f_l \in M^*$ with  $ f_l(x_i) = \delta
_{li}$ for  $ i, l = 1, 2, \cdots n$ . Let  $f_l \otimes id$ act on
equation ({\ref{e.1.1}).
 We have that
$d_l = \sum_{ j = 1 }^{m} f_l(y_j)e_j \in E,$  which implies
  $d_l \in D \cap E = 0$ and  $d_l = 0$  for $ l = 1, \cdots , n$.
 Therefore  $M_D \cap M_E = 0$.

(8)  If $D \cap E \not= 0$, then $D =E$. Thus $M_D = M_E$ and   $M_D
\cap M_E = M_E = M_D = 0$.
 We get a contradiction. Therefore  $D \cap E = 0$.

(9) Obviously, $C(N) \subseteq C(M_E)$. Conversely,since  $0 \not= N
 \subseteq M_{C(N)} \subseteq M_{C(M_E)} = M_E$  and
$M_E$ is a minimal component, $M_{C(N)} = M_E$. By the definition of
component, $C(M_E) \subseteq C(N)$. Thus $C(N) = C(M_E)$.
 \begin{picture}(8,8)\put(0,0){\line(0,1){8}}\put(8,8){\line(0,-1){8}}\put(0,0){\line(1,0){8}}\put(8,8){\line(-1,0){8}}\end{picture}

\begin{Proposition} \label{7.1.2}
 If $E$ is a subcoalgebra of $C$,then
the following conditions are equivalent.

(1)  $M_E$ is a minimal component of $M$.

(2)  $C(M_E)$ is a simple subcoalgebra of $C$.

(3)  $M_E$ is a minimal closed subcomodule of $M$.
\end{Proposition}
{\bf Proof}.  It is easy to check that $M_E = 0$ iff  $C(M_E)= 0$.
Thus  (1), (2) and (3) are equivalent when $C(M_E) =0$,. We now
assume that $M_E \not= 0$.

(1)  $\Longrightarrow$ (2) Since   $ 0 \not= M_E$, there exists a
non-zero finite dimensional simple subcomodule $N$  of $M$ such that
 $N \subseteq M_E$. By  \cite[Lemma 1.1]{XSF94},
$N$ is  a simple $C^*$-submodule of $M$. Since $M_E$ is a minimal
component, $C(N) = C(M_E)$  by Lemma \ref {7.1.1} (9). Let $D = C(N)
= C(M_E)$. By Lemma \ref{7.1.1}(1), $(0:N)_{C^*} = N^{ \bot C^*} =
C(N)^{ \bot C^*} = D^{ \bot C^*}$. Thus $N$ is a faithful simple
$C^*/D^{ \bot C^*}$-module, and so  $C^*/D^{ \bot C^*}$ is a simple
algebra.It is clear that $D$ is a simple subcoalgebra of $C.$

$(2 ) \Longrightarrow (3)$  If $0\not= N \subseteq M_E$ and  $N$ is
a closed subcomodule of  $M$, then by Lemma \ref {7.1.1}(3)
 there exists a subcoalgebra $F$ of  $C$ such that  $N = M_F$.
Since $0 \not= C(N) = C(M_F) \subseteq C(M_E)$ and   $C(M_E)$ is
simple, $C(M_F) = C(M_E)$. By Lemma \ref{7.1.1}(2), $ N = M_F =
M_{C(M_F)} = M_{C(M_E)} = M_E$, which implies  that $M_E$ is a
minimal closed subcomodule.

$(3) \Longrightarrow (1)$, If  $0\not= M_F \subseteq M_E$,
  then $M_F$ is a closed subcomodule by Lemma \ref{7.1.1}(3) and
   $M_E = M_F$, i.e.
$M_E$ is a minimal component.

This completes the proof.
 \begin{picture}(8,8)\put(0,0){\line(0,1){8}}\put(8,8){\line(0,-1){8}}\put(0,0){\line(1,0){8}}\put(8,8){\line(-1,0){8}}\end{picture}

Let

${ \bf \cal C}_0 =  \{ D  \mid  D$ is a simple subcoalgebra of  $C$
\};

${ \bf \cal C}(M)_0 =  \{ D  \mid  D$ is a simple subcoalgebra of
$C$ and  $D \subseteq C(M)$ \};

${ \bf \cal M}_0 =  \{ N  \mid  N$ is a minimal closed subcomodule
of  $M$ \};

${ \bf \cal C}(M)_1 =  \{ D  \mid  D$ is a faithful simple
subcoalgebra of  $C$ to $M$ \};

$C_0 = \sum  \{ D  \mid  D$ is  a simple subcoalgebra of $C$ \};

$M_0 = \sum \{ N  \mid  N$ is a minimal closed subcomodule of  $M$
\};

$C_0$ and  $M_0$  are called the coradical of coalgebra  $C$ and the
coradical of comodule  $M$ respectively. If $M_0 =M$, then M is
called cosemisimple.  By  Lemma \ref {7.1.1}(5),

${ \bf \cal C}(M)_1 \subseteq { \bf \cal C}(M)_0$

\begin {Theorem}  \label {7.1.3}
$\psi  \left \{  \begin{array}{l} { \bf \cal C}(M)_1  \longrightarrow { \bf \cal M}_0\\
 D \longmapsto M_D
\end{array} \right.$
is bijective.

\end{Theorem}
{\bf Proof}.  By Proposition \ref{7.1.2}, $\psi$ is a map. Let $D$
and $E \in$ { \bf \cal C}$(M)_1$ and $\psi (D) = \psi (E)$, i.e.
$M_D = M_E $. By Lemma \ref {7.1.1} (6), we have that $D =C(D_D)=
C(M_E)=E.$ If $N \in { \bf \cal M}_0$, then $N = M_{C(N)}$ and
$C(N)$ is a simple subcoalgebra by Lemma \ref{7.1.1}(3) and
Proposition \ref{7.1.2}. Thus $\psi (C(N)) = N,$  which implies that
$ \psi$ is surjective.
\begin{picture}(8,8)\put(0,0){\line(0,1){8}}\put(8,8){\line(0,-1){8}}\put(0,0){\line(1,0){8}}\put(8,8){\line(-1,0){8}}\end{picture}

 In \cite{XSF94} and \cite{XF92}, Xu defined the equivalence relation
for coalgebra and for comodule as follows:

\begin{Definition} \label{7.1.4}
   We say that
 $D \sim E$ for $D$ and  $E \in { \bf \cal C}_0$ iff for any pair of
subclasses ${ \bf \cal C}_1$  and   ${ \bf \cal C}_2$
 of  ${ \bf \cal C}_0$ with
 $D \in    { \bf \cal C}_1$ and $ E \in { \bf \cal C}_2 $ such that
${ \bf \cal C}_1 \cup { \bf \cal C}_2  = { \bf \cal C}_0$ and ${ \bf
\cal C}_1 \cap { \bf \cal C}_2 = \emptyset $, there exist elements
$D_1 \in { \bf \cal C}_1$ and $ E_1 \in { \bf \cal C}_2 $ such that
 $D_1 \wedge E_1 \not= E_1 \wedge D_1 $. Let $[D]$ denote the
equivalence class which contains  $D$.

   We say that
 $N \sim L$ for $N$ and  $L \in { \bf \cal M}_0$ iff for any pair of
subclasses ${ \bf \cal M}_1$  and   ${ \bf \cal M}_2$
 of  ${ \bf \cal M}_0$ with
 $N \in    { \bf \cal M}_1$ and $ L \in { \bf \cal M}_2 $ such that
${ \bf \cal M}_1 \cup { \bf \cal M}_2  = { \bf \cal M}_0$ and ${ \bf
\cal M}_1 \cap { \bf \cal M}_2 = \emptyset $, there exist elements
$N_1 \in { \bf \cal M}_1$              and $ L_1 \in { \bf \cal M}_2
$ such that
 $N_1 \wedge L_1 \not= L_1 \wedge N_1 $. Let $[N]$ denote the
equivalence class which contains  $N$.
\end {Definition}

\begin{Definition} \label{7.1.5}

 If $D \wedge E = E \wedge D$ for  any simple subcoalgebras $D$ and $E$
 of $C$,
 then $C$ is called $\pi$-commutative.
  If $N \wedge L = L \wedge N$ for any minimal closed subcomodules
$N$ and $L$ of $M$, then $M$ is called  $\pi$-commutative.
\end {Definition}

Obviously, every cocommutative coalgebra is $\pi$-commutative.
  By \cite[Theorem 3.8 and Theorem 4.18]{XSF94},
$M$ is  $\pi$-commutative iff    $M$ can   be decomposed into a
direct sum  of the weak-closed relative-irreducible
 subcomodules of $M$ iff every equivalence class of $M$ contains only one element.
  By \cite{XF92}, $C$ is  $\pi$-commutative iff
   $C$ can   be decomposed
into a direct sum  of irreducible
 subcoalgebras  of $C$ iff equivalence every class of $C$ contains only one
 element.

\begin{Lemma} \label{7.1.6}          Let
 $D$, $E$ and $F$ be subcoalgebras of  $C$. $N$, $L$ and $T$
 be subcomodules of  $M$. Then

(1)  $M_D \wedge M_F = M_{C(M_D) \wedge C(M_F)} \subseteq M_{D
\wedge F}$;

(2)  If  $D$ and  $E$ are faithful simple subcoalgebras of  $C$
 to $M$, then
$ M_D \wedge M_E =  M_{D \wedge F}$;

(3)   $ M_{D + E} \supseteq  M_D + M_E$;

(4)  If $F = \sum $\{ $D_{\alpha} \mid  \alpha \in \Omega $\} and
$\{ D_{\alpha} \mid  \alpha \in \Omega \}  \subseteq  { \bf \cal
C}_0$, then $M_F = \sum \{ M_{D_{\alpha}} \mid \alpha \in \Omega
\}$.In particular, $M_{C_0} = M_0$.

(5) $  (N + L) \wedge T \supseteq N \wedge T + L \wedge T$;

(6)  $(D + E) \wedge F \supseteq D \wedge F + E \wedge F$.
\end{Lemma}
{\bf Proof}. (1)     We see that
\begin {eqnarray*}
 M_D \wedge M_E &=&
\rho ^{-1}(M \otimes (M_D)^{ \bot C^* \bot C } \wedge
 (M_E)^{ \bot C^* \bot C }) \hbox { \ ( by \cite[Proposition 2.2(1)]{XSF94}) }\\
&=&  \rho ^{-1}(M \otimes C(M_D) \wedge C(M_E)) \hbox { \ (by Lemma
\ref{7.1.1}(1)).}
\end {eqnarray*}
By the definition of component, subcomodule $M_D \wedge M_E
\subseteq M_{C(M_D) \wedge C(M_E)}$. It follows  from the equation
above  that
 $M_{C(M_D) \wedge C(M_E) } \subseteq M_D \wedge M_E$.
Thus
$$M_D \wedge M_E = M_{C(M_D) \wedge C(M_E)}$$
and
 $$M_{C(M_D) \wedge C(M_E) } \subseteq M_{D \wedge E}.$$

(2)  Since $D$ and $E$ are faithful simple subcoalgebras of  $C$ to
 $M$, $C(M_D) = D$ and $C(M_E) = E$  by  Lemma \ref {7.1.1}(6). By  Lemma \ref{7.1.6}(1),
 $M_D \wedge M_E =  M_{D \wedge E}$.

(3) It is trivial.

(4) Obviously $ M_F \supseteq \sum \{ M_{D_{ \alpha }} \mid \alpha
\in \Omega \} $. Conversely, let $N = M_F$. Obviously $N$ is  an
$F$-comodule. For any $x \in N$, let $ L = C^*x$. it is clear
 that  $L$ is a finite dimensional comodule over $F$.
By \cite[Lemma 14.0.1] {Sw69a}, $L$ is a completely reducible module
over
 $F^*$. Thus  $L$ can be decomposed into a direct sum of simple
$F^*$-submodules:
$$L = L_1 \oplus L_2 \oplus  \cdots  \oplus L_n,$$
where  $L_i$ is a simple $F^*$-submodule. By \cite [Proposition
1.16] {XSF94},
 $<L_i> =
(L_i)^{ \bot F^* \bot N}$ is a minimal closed  $F$-subcomodule of
$N$.
 By Theorem \ref{7.1.3}, there exists a simple subcoalgebra
$D_{\alpha _i}$ of $F$ such that
 $<L_i> =
N_{ D_{\alpha_i}}$. Obviously, $N_{ D_{\alpha_i}} \subseteq M_{
D_{\alpha _i}}$. Thus $L \subseteq \sum_{ i = 1}^{n} <L_i> = \sum_{i
= 1}^{ n} N_{D_{\alpha _i}} \subseteq \sum \{ M_{D_ {\alpha}}  \mid
\alpha  \in \Omega \}$. Therefore
 $M_F =  \sum \{ M_{ D_{\alpha}}  \mid
\alpha  \in \Omega \}$. If $C_0 = F = \sum \{ D  \mid  D \in { \bf
\cal C}_0 \}$, then $M_{C_0} = M_F = \sum \{ M_D  \mid  D  \in { \bf
\cal C}_0  \}=M_0$ by Theorem \ref{7.1.3},

(5) and (6)  are
trivial.\begin{picture}(8,8)\put(0,0){\line(0,1){8}}\put(8,8){\line(0,-1){8}}\put(0,0){\line(1,0){8}}\put(8,8){\line(-1,0){8}}\end{picture}

\begin{Lemma} \label{7.1.7} Let $N$ be a subcomodule
of  $M$, and let $D$, $E$ and $F$ be simple subcoalgebras of $C$.
Then

(1) $D \sim 0$ iff $D = 0$;  {~}{~} $M_D \sim 0$ iff $M_D = 0$;

We called the equivalence class which contains zero a zero
equivalence class.

(2)  If  $D$ and $E$ are faithful to  $M$ and
  $M_D \sim M_E$, then  $D \sim E$;

(3)  $[M_D] \subseteq \{ M_E \mid  E \in [D] \}$;

(4)  If  $D$ is faithful to  $M$, then
   $[M_D] \subseteq \{ M_E \mid  E \in [D]$ and  $E$  is faithful to  $M$  \}.
   \end{Lemma}
 {\bf Proof}.
(1) If $ D \sim 0$ and $ D \not=  0$, let  ${ \bf \cal C}_1 = \{ 0
\}$ and ${ \bf \cal C}_2 = $\{ $F \mid F \not= 0$, $F \in { \bf \cal
C}_0 $\}. Thus ${ \bf \cal C}_1 \cup { \bf \cal C}_2 =  { \bf \cal
C}_0$ and
 ${ \bf \cal C}_1 \cap { \bf \cal C}_2 = \emptyset $ with
$0 \in  { \bf \cal C}_1$ and $D \in { \bf \cal C}_2$. But for any
$D_1 \in { \bf \cal C}_1$ and $E_1 \in { \bf \cal C}_2$, since  $D_1
= 0$,  $D_1 \wedge E_1 = E_1 \wedge D_1 = E_1$. We get a
contradiction. Thus  $D = 0$. Conversely, if  $D = 0$, obviously $ D
\sim 0$. Similarly, we can show that   $M_D \sim 0$ iff $M_D = 0$.

(2) For any pair of subclasses
 ${ \bf \cal C}_1$ and
 ${ \bf \cal C}_2$ of ${\bf \cal C}_0$, if
${ \bf \cal C}_1 \cap { \bf \cal C}_2 = \emptyset $ and   ${ \bf
\cal C}_1  \cup { \bf \cal C}_2 = { \bf \cal C}_0$ with
 $ D \in { \bf \cal C}_1$ and
$E \in { \bf \cal C}_2$, let
 ${ \bf \cal M}_1  = $\{ $M_F \mid F \in { \bf \cal C}_1$  and  $F$
 is faithful to  $M$  \} and
 ${ \bf \cal M}_2  =$ \{ $M_F \mid F \in { \bf \cal C}_2 $ and  $F$
 is faithful to  $M$ \}.
By Theorem \ref{7.1.3},
 ${ \bf \cal M}_0 = { \bf \cal M}_1 \cup { \bf \cal M}_2$ and
 ${ \bf \cal M}_1 \cap { \bf \cal M}_2 =
\emptyset $ . Obviously, $M_D \in  { \bf \cal M}_1$ and $M_E \in {
\bf \cal M}_2$. Since  $ M_D \sim M_E$,there exist  $M_{D_1} \in {
\bf \cal M}_1$ and
 $M_{E_1} \in { \bf \cal M}_2$ such that
 $M_{D_1} \wedge  M_{E_1} \not=   M_{E_1} \wedge  M_{D_1}$, where
$D_1 \in { \cal C}_1$ and  $E_1 \in { \cal C}_2$. By Lemma
\ref{7.1.6}(2), $M_{ D_1 \wedge E_1} \not=  M_{E_1 \wedge D_1}$.
Thus $D_1 \wedge E_1 \not= E_1 \wedge D_1$. Obviously
 $D {~and~}  D_1 \in { \bf \cal C}_1$. Meantime
$E {~and~} E_1 \in  { \bf \cal C}_2$. By Definition \ref{7.1.4}, $D
\sim E$.

(3)  Obviously, $M_F \sim M_D$ for any $M_F \in [M_D]$. If $M_D
\not= 0$,
  then $M_F \not= 0$ by Lemma \ref {7.1.7} (1),. By Lemma \ref {7.1.7}(2), $ F \sim D$.
Thus  $M_F \in \{ M_E \mid E \in [D] \}$. If $M_D = 0$, by Lemma
\ref{7.1.7}(1), $M_F = 0 = M_D \in \{ M_E \mid E \in [D] \}$.

(4)  If $M_F \in [M_D]$, then $M_D \sim M_F$. If $D = 0$, then $M_D
= 0$. By Lemma \ref {7.1.7}(1), $ M_F = 0$. Thus $M_F = M_D = 0 \in$
\{$ M_E  \mid E \in [D]$ and $E$ is faithful to  $M$ \}. If $ D
\not= 0$, we have that
 $M_D \not= 0$ since  $D$ is faithful to  $M$.
By  Lemma \ref {7.1.7}(1), $M_F \not= 0$. By Lemma \ref{7.1.7}(3),
$M_F\in $ \{ $M_E  \mid E \in [D]$ and  $E$ is faithful to $M$ \}.
\begin{picture}(8,8)\put(0,0){\line(0,1){8}}\put(8,8){\line(0,-1){8}}\put(0,0){\line(1,0){8}}\put(8,8){\line(-1,0){8}}\end{picture}

\begin{Theorem} \label{7.1.8}
Let $\{ { \bf \cal E}_{\alpha} \mid \alpha \in \bar{\Omega} \}$ be
all of the equivalence classes of   $C$. and $E_\alpha = \sum \{ E
\mid E \in { \bf \cal E}_\alpha \}$. Then

(1)  For any  $\alpha \in \bar{ \Omega }$, there exists a set $I_
\alpha$ and subclasses
 ${ \bf \cal E}(\alpha , i ) \subseteq { \bf \cal E}_{ \alpha }$
such that $\cup \{ { \bf \cal E}( \alpha,i) \mid i \in I_{ \alpha }
\} = { \bf \cal E}_{ \alpha }$ and $\{ M_{ { \bf \cal E}( \alpha,i)}
\mid   \alpha \in \bar \Omega ,
i  \in I_ \alpha  \}$ \\
 is the set of the equivalence classes of  $M$
(they are distinct except for  zero equivalence class), where $M_{ {
\bf \cal E}( \alpha ,i)}$ denotes $ \{ M_D \mid  D \in { \bf \cal
E}( \alpha,i) \}$.

(2)  If  $M$ is a component faithfulness  $C$-comodule, then
 $ \{ M_{ { \bf \cal E}(\alpha , i)}  \mid
 \alpha \in \bar{ \Omega } , i \in I_{ \alpha } \}$
is the set  of the distinct equivalence  classes of  $M$.

(3) $$ M = \sum_{ \alpha \in \bar{ \Omega }}  \oplus M_{(E_{ \alpha
})^{ (\infty)}} =  \sum_{ \alpha \in \bar{ \Omega }} \sum_{ i \in
I_{ \alpha } }
 \oplus (M_{E( \alpha,i)})^{ (\infty)}
= \sum_{ \alpha \in \bar{ \Omega }}  \oplus (M_{E_{ \alpha }})^{
(\infty)}
$$
and for any
 $\alpha \in \bar{ \Omega }$,

$$ M_{(E_{ \alpha })^{ (\infty)}}
= \sum_{  i \in I_{ \alpha } }
 (M_{E( \alpha,i)})^{ (\infty)}
= ( M_{E_{ \alpha }})^{ (\infty)} $$ where $E( \alpha,i) = \sum \{ E
\mid  E  \in { \bf \cal E}( \alpha,i) \}$.

\end{Theorem}
 {\bf Proof}. (1) By Theorem \ref{7.1.3} and  Lemma \ref{7.1.7}(3) we can
immediately get part (1).

(2)  If  $M$ is a component faithfulness $C$-comodule, then ${ \bf
\cal C}(M)_1  =  { \bf \cal C}(M)_0 =  { \bf \cal C}_0  $. It
follows from Theorem \ref {7.1.3} and part (1)  that $ \{ M_{ { \bf
\cal E}(\alpha , i)}  \mid  \alpha \in \bar{ \Omega }, i \in I_{
\alpha } \}$ consists of  all  the distinct equivalence classes of
$M$.

(3) We see that

 $M = M_C =
M_{ \sum \{ (E_{ \alpha})^{ (\infty)} \mid  \alpha \in \bar{ \Omega
} \}}$
 (by  \cite{XF92})

$\supseteq \sum_{ \alpha \in \bar{ \Omega }} M_{{(E_{\alpha})}^{
(\infty)}}$ (by  Lemma \ref{7.1.6}(3))

$\supseteq \sum _{ \alpha \in \bar{ \Omega }} \sum _{ n = 0 }^{
\infty} M_{ \wedge ^{n+1}E_{ \alpha}}$ ( by Lemma \ref{7.1.6}(3))

$\supseteq \sum _{ \alpha \in \bar{ \Omega }} \sum _{ n = 0 }^{
\infty} {\wedge}^{n+1} M_{E_{ \alpha}}$ (by  Lemma \ref{7.1.6}(1))

$= \sum_{ \alpha \in \bar{ \Omega }} \sum _{ n = 0 } ^{ \infty}
 {\wedge} ^{n+1} M_{ \sum_{ i \in I_{ \alpha }}E( \alpha,i)}$  (by Theorem
\ref{7.1.8}(1) )

$\supseteq \sum _{ \alpha \in \bar{ \Omega }} \sum _{ n = 0 } ^{
\infty} {\wedge} ^{n+1} \sum _{ i \in I_{ \alpha }} M_{E( \alpha,
i)}$ (by Lemma \ref{7.1.6}(3))

$\supseteq \sum _{ \alpha \in \bar{ \Omega } }\sum _{ n = 0 } ^{
\infty}
 \sum_{ i \in I_{ \alpha }} {\wedge} ^{n+1} M_{E(\alpha,i)}$
( by  Lemma \ref{7.1.6}(5)(6))

$=  \sum _{ \alpha \in \bar{ \Omega }} \sum_{ i \in I_{ \alpha}}
 \sum _{ n = 0 } ^{ \infty} \wedge^{n+1}M_{E( \alpha,i)}$

$= \sum _{ \alpha \in \bar{ \Omega }} \sum_{ i \in I_{ \alpha }}
(M_{ E(\alpha,i)})^{ (\infty)}$

$=M$  (by \cite[(4.10) in Theorem4.15]{XSF94} and Lemma
\ref{7.1.6}(4)
 and part (1)).                                           \\
Thus
\begin {eqnarray}
 M = \sum_{ \alpha \in \bar{ \Omega }} M_{(E_{\alpha })^{ (\infty)}}
= \sum_{ \alpha \in \bar{ \Omega }} \sum_{ i \in I_{ \alpha }}
(M_{E( \alpha,i)})^{ (\infty)} \label {e.8.4}
\end {eqnarray}
and

$$M = \sum _{ \alpha \in \bar{ \Omega }} \sum_{ i \in I_{ \alpha }}\oplus
(M_{E( \alpha,i)})^{ (\infty)}$$ by \cite[Theorem 4.15]{XSF94} and
part (1). We see that
\begin {eqnarray*}
M_{E(\alpha,i)} \wedge M_{E(\alpha,i)}
&\subseteq & M_{(E_{\alpha})^{ (\infty)}} \wedge M_{(E_{\alpha})^{ (\infty)}} \\
&\subseteq & M_{(E_{\alpha})^{ (\infty)} \wedge (E_{\alpha})^{
(\infty)}}
\hbox { \ (by Lemma \ref{7.1.6}(1)) }  \\
&=& M_{(E_{\alpha})^{ (\infty)}} \hbox { (by \cite[ Proposition
2.1.1 ]  {HR74})},
\end {eqnarray*}
Thus $$M_{(E_{\alpha})^{ (\infty)}} \supseteq (M_{E( \alpha,i)})^{
(\infty)}$$  for any  $i  \in I_{\alpha }$ and   for any  $\alpha
\in \bar{ \Omega }$,  and
 \begin {eqnarray}
 M_{(E_{ \alpha})^{ (\infty)}} \supseteq
\sum_{ i \in I_{ \alpha }}(M_{E( \alpha,i)})^{ (\infty)} \label
{e.8.5}
\end {eqnarray}
If  $M_{(E_{\alpha})^{(\infty)}} \cap \sum_{\beta \in \bar{ \Omega
}, \beta \not= \alpha } M_{(E_{ \beta})^{ (\infty)}} \not= 0$, then
there exists a non-zero simple subcomodule $C^*x \subseteq M_{(E_{
\alpha})^{ (\infty)}} \cap
 \sum_{ \beta \in \bar{ \Omega },
\beta \not= \alpha } M_{(E_{ \beta})^{ (\infty)}}$.
 By
\cite[Proposition 1.16]{XSF94}, $<x>$ is a minimal closed
subcomodule of $M$. By Lemma \ref{7.1.1}(3),
   $M_{(E_{ \alpha})^{ (\infty)}}$
is a closed subcomodule of $M$. By \cite[Lemma 3.3]{XSF94},
 there exists  $\gamma \in  \bar{ \Omega }$ with
$ \gamma \not= \alpha$ such that   $C^*x \subseteq M_{(E_{
\gamma})^{(\infty)}}$.
 Thus  $M_{(E_{ \gamma})^{(\infty)}} \cap M_{(E_{ \alpha })^{(\infty)}}
\not= 0$. By Lemma \ref {7.1.1} (7), $ (E_{ \gamma})^{(\infty)} \cap
(E_{ \alpha })^{(\infty)} \not= 0,$ which  contradicts \cite{XF92}.
Thus for any $\alpha \in \bar{ \Omega }$, $M_{(E_{ \alpha})^{
(\infty)}} \cap \sum_{ \beta \in \bar{ \Omega }, \beta \not= \alpha
} M_{(E_{ \alpha})^{ (\infty)}}= 0,$ which implies  that
$$M = \sum_{ \alpha \in \bar{ \Omega }} \oplus
 M_{(E_{ \alpha})^{ (\infty)}}.$$

It follows from equations ( \ref {e.8.4}) and (\ref {e.8.5} ) that
\begin {eqnarray}
                M_{(E_{ \alpha})^{ (\infty)}} &=&
\sum_{ i \in I_{ \alpha }}(M_{E( \alpha,i)})^{ (\infty)} \label
{e.8.6}
\end {eqnarray}
We see that
\begin {eqnarray*}
M_{E_{\alpha}} \wedge M_{E_{\alpha}} & \subseteq & M_{E_{\alpha}
\wedge E_{\alpha}}
 \hbox { \ (by Lemma \ref {7.1.6} (1)) } \\
 & \subseteq & M_{(E_{\alpha})^{ (\infty)}}.
 \end {eqnarray*}
 Thus
$$M_{(E_{ \alpha})^{(\infty)}} \supseteq (M_{E_{ \alpha}})^{(\infty)}
\supseteq \sum_{ i \in I_{ \alpha }}(M_{E( \alpha,i)})^{ (\infty)}$$
and   $$M_{(E_{ \alpha})^{(\infty)}}= (M_{E_{ \alpha}})^{(\infty)}$$
by relation ( \ref {e.8.6} ).
 This completes the proof. \begin{picture}(8,8)\put(0,0){\line(0,1){8}}\put(8,8){\line(0,-1){8}}\put(0,0){\line(1,0){8}}\put(8,8){\line(-1,0){8}}\end{picture}

\begin {Definition} \label{7.1.9}  If  $M$
 can   be decomposed
into a direct sum  of two non-zero  weak-closed subcomodules, then
$M$ is called decomposable.      If $N$ is a subcomodule of $M$ and
$N$ contains exactly one non-zero minimal closed submodule, then $N$
is said to be relative-irreducible.

\end {Definition}

\begin{Corollary} \label{7.1.10}   $C$ is a coalgebra.

(1)  If $C$ is   $\pi$-commutative, then every  $C$-comodule
 $M$ is also $\pi$-commutative;

(2)  If  $C$ can be decomposed into a direct sum of its irreducible
subcoalgebras, then every $C$-comodule
 $M$ can also be decomposed into a direct sum of its relative-irreducible
subcomodules;

(3)  If  $C$ is decomposable, then every component faithfulness
 $C$-comodule  $M$ is decomposable;

(4) If  $C$ is irreducible,  then every non-zero  $C$-comodule  $M$
is relative-irreducible;

(5) $C$ is irreducible iff every component faithfulness $C$-comodule
$M$ is relative-irreducible.

\end{Corollary}
{\bf Proof}. (1) For any pair of non-zero closed subcomodules $N$
and $L$ of $M$, by Theorem \ref{7.1.3}, there exist  faithful simple
subcoalgebras $D$ and $E$ of $C$ to $M$  such that $N=M_D$ and $L
=M_E$.
 By Lemma \ref{7.1.6}(2),
$ N \wedge L = M_D \wedge M_E = M_{D \wedge E} = M_{E \wedge D } =
M_E \wedge M_D = L \wedge N$. Thus $M$ is  $\pi$-commutative.

(2) Since  $C$ can be decomposed  into a direct sum of its
irreducible
 subcoalgebras,   every equivalence class of $C$
contains only one element by \cite{XF92}. By Theorem \ref{7.1.8}(1),
every equivalence class of $M$ also contains only one element. Thus
it follows from
 \cite [Theorem 4.18]{XSF94}  that
 $M$ can be decomposed into a direct sum of its relative-irreducible
 subcomodules.

(3)  If  $C$ is decomposable, then there are at least two non-zero
equivalence classes in  $C$. By   Theorem \ref{7.1.8}(2),
 there are at least two non-zero equivalence classes in $M$.
Thus  $M$ is decomposable.

(4)  If  $C$ is irreducible, then  there is only one non-zero simple
subcoalgebra of $C$ and there is at most one non-zero minimal closed
subcomodule in $M$ by Theorem \ref {7.1.3}. Considering  $M \not=
0$, we have that $M$ is relative-irreducible.

(5)  If  $C$ is irreducible, then every component faithfulness
$C$-comodule  $M$ is relative-irreducible
 by  Corollary \ref{7.1.10}(4).
Conversely, let $M = C$  be the regular  $C$-comodule. If  $D$ is a
non-zero simple subcoalgebra of  $C$, then
 $ 0 \not= D \subseteq M_D$. Thus $M$ is a component faithfulness
$C$-comodule. Since  $M$
 is a relative-irreducible $C$-comodule by assumption,
  there is only one non-zero
minimal closed subcomodule in  $M$ and so
 there is also only  one non-zero simple subcoalgebra in $C$
 by Theorem \ref {7.1.3}.
Thus  $C$ is irreducible.
\begin{picture}(8,8)\put(0,0){\line(0,1){8}}\put(8,8){\line(0,-1){8}}\put(0,0){\line(1,0){8}}\put(8,8){\line(-1,0){8}}\end{picture}

\begin{Lemma} \label{7.1.11} Let $N$ be a
$C$-subcomodule of $M$ and $\emptyset \not= L \subseteq M$. Then

(1)  $C^* \cdot L = C(M)^* \cdot L$; {~}{~} $N^{ \bot C^* \bot M } =
N^{ \bot C(M)^* \bot M }$;

(2)  $N$ is a (weak-) closed  $C$-subcomodule iff  $N$ is a
(weak-)closed $C(M)$-subcomodule;

(3)  $N$ is a minimal closed  $C$-subcomodule iff $N$ is a minimal
closed $C(M)$-subcomodule;

(4) $N$ is an indecomposable  $C$-subcomodule iff $N$ is an
indecomposable  $C(M)$-subcomodule;

(5)  $N$  is a relative-irreducible  $C$-subcomodule iff $N$ is a
relative-irreducible  $C(M)$-subcomodule.
\end{Lemma}
{\bf Proof}. (1) Let $C=C(M) \oplus V,$ where $V$ is a subspace of
$C$. If  $f \in V^*$, then  $ f \cdot L = 0$. Thus
  $C^* \cdot L  = (C(M)^* + V^* ) \cdot L = C(M)^* \cdot L$.
We now show the second equation. Obviously,
 $$N^{ \bot C^* \bot M } \subseteq
N^{ \bot C(M)^* \bot M }.$$ Conversely,  for any $x \in N^{ \bot
C(M)^* \bot M }$ and $f \in N^{ \bot C^* }$, there exist  $f_1 \in
C(M)^*$ and  $f_2 \in V^*$ such that
 $f = f_1 + f_2$.     Obviously,
  $f \cdot x = f_1 \cdot x = 0.$  Thus
$x \in N^{ \bot C^* \bot M }$ and $N^{ \bot C(M)^* \bot M }
\subseteq N^{ \bot C^* \bot M }$.
 Therefore
 $$N^{ \bot C^* \bot M } = N^{ \bot C(M)^* \bot M }.$$

(2)  If  $N$ is a weak-closed  $C$-subcomodule, then $(C^* \cdot
x)^{ \bot C^* \bot M } \subseteq N $ for any $x \in N$. We see that
\begin {eqnarray*}  (C^* \cdot x)^{ \bot C^* \bot M }
&=& (C(M)^* \cdot x)^{ \bot C^* \bot M } \hbox { ( \ by  Lemma \ref{7.1.11}(1))  } \\
&=& (C(M)^* \cdot x)^{ \bot C(M)^* \bot M } \hbox { \ ( \ by Lemma
\ref{7.1.11}(1)) }
\end {eqnarray*}
Thus
 $(C(M)^* \cdot x)^{ \bot C(M)^* \bot M } \subseteq N$ and so
 $N$ is a weak-closed
 $C(M)$-subcomodule.
Conversely, if  $N$ is a weak-closed  $C(M)$-subcomodule. Similarly,
we can show that  $N$ is a weak-closed  $C$-subcomodule. We now show
the second assertion.  If $N$ is a closed $C$-subcomodule, then $N^{
\bot C^* \bot M } =  N$ and
 $N= N^{ \bot C^* \bot M } = N^{ \bot C(M)^* \bot M } $ by part (1),
 which implies that
 $N$ is a closed  $C(M)$-subcomodule. Conversely,
if  $N$ is a closed  $C(M)$-subcomodule,
 similarly, we can  show that $N$ is a closed $C$-subcomodule.

Similarly the others can  be proved.
\begin{picture}(8,8)\put(0,0){\line(0,1){8}}\put(8,8){\line(0,-1){8}}\put(0,0){\line(1,0){8}}\put(8,8){\line(-1,0){8}}\end{picture}

\begin{Proposition} \label{7.1.12}  Let
every simple subcoalgebra in $C(M)$ be faithful to  $M$.

(1)  If $M$ is an indecomposable $C$-comodule, then $C(M)$ is also
an indecomposable subcoalgebra.

(2)  $M$ is a  relative-irreducible  $C$-comodule iff $C(M)$ is an
irreducible subcoalgebra.
\end{Proposition}
 {\bf Proof}.
(1)  If  $M$ is an indecomposable  $C$-comodule, then $M$ is an
indecomposable  $C(M)$-comodule by Lemma \ref{7.1.11} and $M$ is a
component faithfulness $C(M)$-comodule. By Corollary
\ref{7.1.10}(3),  $C(M)$ is indecomposable.

(2)  If  $M$ is a relative-irreducible $C$-comodule, then
  $M$ is a relative-irreducible  $C(M)$-comodule
by Lemma \ref{7.1.11}(5)  and
  $C(M)$ is irreducible. Conversely,
if  $C(M)$ is irreducible, then $M$ is a relative-irreducible
$C(M)$-comodule  by Corollary \ref{7.1.10}(4) and so $M$ is a
relative-irreducible $C$-comodule by Lemma \ref{7.1.11}(5).
\begin{picture}(8,8)\put(0,0){\line(0,1){8}}\put(8,8){\line(0,-1){8}}\put(0,0){\line(1,0){8}}\put(8,8){\line(-1,0){8}}\end{picture}

\begin{Definition} \label{7.1.13}
If
 $D \sim E$ implies  $M_D \sim M_E$ for any simple subcoalgebras
$D$ and $E$ in $C(M)$, then $M$ is called a $W$-relational
hereditary
 $C$-comodule.
\end {Definition}

 If $M$ is a  $W$-relational hereditary $C$-comodule, then
 ${ \cal C}(M)_1 = {\cal C}(M)_0$.
In fact, by Lemma \ref {7.1.1}(5) , ${ \cal C}(M)_1 \subseteq  {\cal
C}(M)_0$. If ${ \cal C}(M)_1 \not= {\cal C}(M)_0$, then there exists
$ 0 \not= D \in { \cal C}(M)_0$ such that  $M_D = 0$. Since  $M_D
\sim M_0$ and  $M$ is  $W$-relational hereditary,  we have that $D
\sim 0$ and $D = 0$ by Lemma \ref {7.1.7}(1). We get a
contradiction. Thus
 ${ \cal C}(M)_1 = {\cal C}(M)_0$.

Obviously, every $\pi$-commutative comodule is  $W$-relational
hereditary. If $C$ is $\pi$-commutative, then $M$ is
$\pi$-commutative by Corollary \ref {7.1.10}(1) and  every
$C$-comodule $M$ is $W$-relational hereditary.
 Furthermore, $M$ is also a component faithfulness $C(M)$-comodule.

\begin{Proposition} \label{7.1.14}
Let  the notation be the same as in Theorem \ref{7.1.8}. Then the
following conditions are equivalent.

(1)  $M$ is  $W$-relational hereditary.

(2)  For any  $\alpha \in  \bar { \Omega }$,
 there is at most one non-zero equivalence class in
 $\{ M_{ { \bf \cal E}( \alpha, i)} \mid i \in I_{ \alpha } \},$
and   ${ \cal C}(M)_1 = {\cal C}(M)_0.$

(3)    For any   $ D \in { \bf \cal C}(M)_0$,

 $[M_D] = $\{ $M_F \mid F \in [D]$  and  $ F \subseteq C(M)$ \}

(4)  For any  $D$ and  $E$ $\in { \bf \cal C}(M)_0$ ,
 $M_D \sim M_E$ iff $D \sim E$.

(5)  For any   $\alpha \in  \bar { \Omega }$, $M_{(E_{\alpha})^{
(\infty) }}$ is indecomposable, and  ${ \cal C}(M)_1 = {\cal
C}(M)_0.$

(6) For any  $ \alpha \in  \bar { \Omega } $, $(M_{E_{\alpha}})^{
(\infty) }$
 is indecomposable, and
  ${ \cal C}(M)_1 = {\cal C}(M)_0;$

(7) $ M = \sum_{ \alpha \in \bar{ \Omega }}  \oplus (M_{E_{ \alpha
}})^{ (\infty)}$ is its weak-closed indecomposable decomposition,
and
 ${ \cal C}(M)_1 = {\cal C}(M)_0;$

(8)  $ M = \sum_{ \alpha \in \bar{ \Omega }}  \oplus M_{(E_{ \alpha
})^{ (\infty)}}$ is its weak-closed indecomposable decomposition,
and
 ${ \cal C}(M)_1 = {\cal C}(M)_0.$
 \end{Proposition}
{\bf Proof}. We prove it along with the following lines:
 $(1) \Longrightarrow  (4)
 \Longrightarrow  (3) \Longrightarrow  (1)$
$(3) \Longleftrightarrow  (2) \Longleftrightarrow  (5)
 \Longleftrightarrow  (6)  \Longleftrightarrow  (7)  \Longleftrightarrow  (8)$.

$(1) \Longrightarrow  (4)$ It follows from the discussion above that
$${ \cal C}(M)_1 = {\cal C}(M)_0.$$
If  $D$ and  $E \in {  \cal C}(M)_0$
 and $M_D \sim M_E $, then   $ D \sim E  $ by  Lemma \ref {7.1.7} (2). If
$D \sim E$ and  $D$ and $E \in {  \cal C}(M)_0$, then $M_ D \sim
M_E$.

$(4) \Longrightarrow (3)$  Considering  Lemma \ref{7.1.7}(4) and
 ${ \cal C}(M)_1 = {\cal C}(M)_0$,
 we only need to show that
$$ \{ M_F \mid F \in [D] \hbox { \  and \ } F \subseteq C(M) \}
 \subseteq [M_D].$$ For any $F \in [D]$ with $F \subseteq C(M)$,
$  F \sim D$ and  $M_F \sim M_D$ by part (4), which implies that
$M_F  \in  [M_D]$.

$(3) \Longrightarrow (1)$  It is trivial.

$(5) \Longleftrightarrow (6)
 \Longleftrightarrow (7)
 \Longleftrightarrow (8)$
 It follows from  Theorem \ref{7.1.8}(3).

$(5) \Longleftrightarrow (2)$ By Theorem \ref{7.1.8} (3),
$$M_{(E_{\alpha})^{ (\infty )}}= \sum_{ i \in I_{ \alpha}}
(M_{E(\alpha,i)})^{ (\infty) }.$$ Thus part (2) and part (5) are
equivalent.

$(2) \Longrightarrow  (3)$ It follows from  Lemma \ref {7.1.7} (4).

 $(3) \Longrightarrow  (2)$ If there are two non-zero
equivalence classes  $M_{ { \bf \cal E}( \alpha, i_1)}  \not= 0 $
and $ M_{ { \bf \cal E}( \alpha, i_2)} \not= 0$ in $\{ M_{ { \bf
\cal E}( \alpha, i) } \mid i \in I_{ \alpha } \}$, then there exist
 $D_1 \in   { \bf \cal E}( \alpha, i_1)$ and
$D_2 \in  { \bf \cal E}( \alpha, i_2)$ such that $M_{D_1} \not= 0$
and $M_{D_2} \not=0$. Let  $D = D_1$. Since  $M_{D_1}$ and $M_{D_2}
\in$
 \{$M_F \mid F \in [D]$  and $F \subseteq C(M)$ \} and
 $M_{D_2} \not\in [M_D]
 = M_{ { \bf \cal E}( \alpha, i_1)}$, this contradicts part (3).
Thus there is at most one non-zero equivalence class in
 $\{ M_{ { \bf \cal E}( \alpha, i)} \mid i \in I_{ \alpha } \}$.
\begin{picture}(8,8)\put(0,0){\line(0,1){8}}\put(8,8){\line(0,-1){8}}\put(0,0){\line(1,0){8}}\put(8,8){\line(-1,0){8}}\end{picture}

\begin{Proposition} \label{7.1.15} If
  $M$ is a full,   $W$-relational hereditary $C$-comodule,
  then

(1) $M$ is indecomposable iff $C$ is indecomposable;

(2) $M$ is relative-irreducible iff $C$ is irreducible.

(3) $M$ can be decomposed into a direct sum of
 its weak-closed  relative-irreducible subcomodules iff
$C$ can be decomposed into a direct sum of its irreducible
subcoalgebras.

(4) $M$ is  $\pi$-commutative iff $C$ is  $\pi$-commutative.

\end{Proposition}
 {\bf Proof}. Since  $M$ is a full $C$-comodule,  $C(M) =C$.
Since  $M$ is  $W$-relational hereditary,  ${ \cal C}(M)_1 = {\cal
C}(M)_0 = { \cal C}_0$. Thus  $M$ is a component faithfulness
$C$-comodule.

(1) If  $M$ is indecomposable, then $C$ is indecomposable by
Proposition \ref {7.1.12}(1). Conversely, if  $C$ is indecomposable,
then there is at most one non-zero equivalence class in  $C$. By
Proposition \ref {7.1.14}(2), there is at most one non-zero
equivalence class in  $M$. Thus
 $M$ is indecomposable.

(4) If  $M$  is  $\pi$-commutative, then there is only one element
in every  equivalence class of $M$. By Proposition \ref {7.1.14}(4)
and Theorem \ref{7.1.3}, there is only one element in every
equivalence class of   $C$.
 Thus  $C$ is  $\pi$-commutative by \cite {XF92}.
Conversely, if $C$  is $\pi$-commutative, then $M$ is
$\pi$-commutative by Corollary \ref {7.1.10}(1).

 (2) It follows from the above discussion
 and Proposition \ref {7.1.12}.

$(3) \Longleftrightarrow (4)$
 By  \cite[Theorem3.8 and Theorem 4.18]{XSF94} and \cite{XF92},
 it is easy to check that (3) and (4) are equivalent. \begin{picture}(8,8)\put(0,0){\line(0,1){8}}\put(8,8){\line(0,-1){8}}\put(0,0){\line(1,0){8}}\put(8,8){\line(-1,0){8}}\end{picture}

\begin {Proposition} \label {7.1.16}
 If  $M_D
 \wedge M_E =M_E \wedge M_D$ implies  $D \wedge E = E \wedge D$
for any simple coalgebras $D$ and $E$ in $C(M)$, then  $M$ is
$W$-relational hereditary.
\end {Proposition}
{\bf Proof}.  Let $ D \sim E$. For any  pair of subclasses ${ \bf
\cal M}_1$ and  ${ \bf \cal M}_2$ of ${ \bf \cal M}_0$ with
 ${ \bf \cal M}_1 \cap { \bf \cal M}_2
= \emptyset $ and   ${ \bf \cal M}_1  \cup { \bf \cal M}_2 = { \bf
\cal M}_0$, if  $ M_D \in { \bf \cal M}_1$ and $M_E \in { \bf \cal
M}_2$, let
 ${ \bf \cal C}_1  =
\{ F \in { \bf \cal C}_0  \mid M_F \in { \bf \cal M}_1  \}$ and
 ${ \bf \cal C}_2  = \{ F \in { \bf \cal C}_0  \mid M_F \in { \bf \cal M}_2  \}$.
By  Theorem \ref{7.1.3},
 $ { \bf \cal C}_0 = { \bf \cal C}_1 \cup { \bf \cal C}_2$ and
 ${ \bf \cal C}_1 \cap { \bf \cal C}_2 =
\emptyset $. Obviously, $D \in  { \bf \cal C}_1 $ and $ E \in { \bf
\cal C}_2$. By Definition \ref {7.1.4}, there exist $D_1 \in { \bf
\cal C}_1$ and $E_1 \in { \bf \cal C}_2 $ such that
 $D_1 \wedge  E_1 \not=   E_1 \wedge  D_1$. By the assumption condition,
we have that
 $M_ {D_1} \wedge M_{E_1} \not=  M_{E_1} \wedge M_{D_1}$. Thus
$M_D \sim M_E$, i.e. $M$ is  $W$-relational hereditary.
\begin{picture}(8,8)\put(0,0){\line(0,1){8}}\put(8,8){\line(0,-1){8}}\put(0,0){\line(1,0){8}}\put(8,8){\line(-1,0){8}}\end{picture}

\begin{Proposition} \label{7.1.17}
Let  $M=C$ as a right $C$-comodule. Let $N=D$ and $L=E$ with $N$ and
$L$ as subcomodules of $M$ with $D$ and $E$ as right coideals of
$C$. Let $X$ be an ideal of $C^*$ and $F$ be a subcoalgebra of $C$.
Then:

(1) $ X^{ \bot C } = X^{ \bot M}$;{~}{~}
  $ N^{ \bot C^* \bot M}  = C(N)$;     {~}{~}
$C(M_F) = F$;

(2)  $N$ is a closed subcomodule of  $M$ iff $D$ is a subcoalgebra
of $C$;

(3)   $N$ is a closed subcomodule iff  $N$ is a weak-closed
subcomodule of $M$;

(4) $N$ is a minimal closed subcomodule of $M$ iff  $D$  is  a
simple subcoalgebra of $C$;

(5) When  $N$ and  $L$ are closed subcomodules, $N \wedge _M L  = D
\wedge _C E$, where  $ \wedge _M$ and  $\wedge _C $ denote  wedge in
comodule $M$ and  in coalgebra $C$ respectively;

(6)   $M$ is a full and  $W$-relational hereditary $C$-comodule and
a component faithfulness $C$-comodule.

(7) The weak-closed indecomposable decomposition of $M$
 as a $C$-comodule and
the indecomposable decomposition of $C$ as coalgebra are the same.

\end{Proposition}
{\bf Proof}. (1)  If $ x \in  X^{ \bot M}$, then $  f \cdot x = 0$
for any  $ f \in X$.
 Let  $\rho (x) = \sum_{ i = 1 } ^{ n} x_i \otimes c_i $ and
$ x_1$, $\cdots $, $x_n$ be linearly independent. Thus  $ f(c_i) =
0$ and  $i = 1$, $ \cdots$ ,$ n$.
 Since
$x = \sum _{ i = 1 }^{n} \epsilon (x_i)c_i$,
 $f(x) = 0$, which implies  $x \in X^{ \bot C}$.

Conversely, if  $x \in X^{ \bot C}$, we have that $\rho (x) \in
  X^{ \bot C } \otimes  X^{ \bot C}$
since  $X^ {\bot C}$  is subcoalgebra of $C$. Thus  $f \cdot x = 0$
for any  $f \in X$, which implies $ x \in X^{ \bot M}$. Thus $ X^{
\bot M}  =  X^{ \bot C }$. By Lemma \ref{7.1.1}(1),  $ N^{ \bot C^*
\bot M}  = C(N)^{ \bot C^* \bot M}$. By part (1),  $ C(N)^{ \bot C^*
\bot M}
 = C(N)^{ \bot C^* \bot C}$.
Thus $ N^{ \bot C^* \bot M}  =C(N)^{ \bot C^* \bot C}  = C(N)$.

Finally, we show that  $C(M_F) = F$. Obviously,  $ C(M_F) \subseteq
F$. If we view $F$ as a $C$-subcomodule of $M$, then
  $F \subseteq C(M_F)$. Thus  $F = C(M_F)$.

(2)  If  $N$ is a closed subcomodule of  $M$, then $  N^{ \bot C^*
\bot M}  = N$. By Proposition \ref {7.1.17}(1),
  $ N^{ \bot C^* \bot M}  = C(N) = N$. Thus
 $\rho (N) = \Delta (D)\subseteq
N \otimes N = D \otimes D$, which implies  that $D$ is a
subcoalgebra of $C$. Conversely, if  $D$ is a subcoalgebra of $C$,
then
 $\rho (N) \subseteq N \otimes N$.
Thus  $C(N) \subseteq D = N$. By Proposition \ref{7.1.17}(1),
   $N^{ \bot C^* \bot M}  = C(N) \subseteq N$. Thus
   $N^{ \bot C^* \bot M}  = N$, i.e. $N$ is closed.

(3) If  $N$ is a closed subcomodule, then  $N$ is weak-closed.
Conversely, if  $N$ is weak-closed, then $ <x> \subseteq N$
 for  $x \in N$. Since  $<x>$
is a closed subcomodule of  $M$,
 $ <x>$ is
subcoalgebra of $C$ if let $<x>$ with structure of coalgebra $C$.
This shows that $ \Delta (x) \in <x> \otimes <x> \subseteq D \otimes
D$. Thus  $D$ is a subcoalgebra of $C$. By Proposition \ref
{7.1.17}(2), $N$ is a closed subcomodule of $M$.

(4) It follows from  part (2).

(5) We only need to show that
$$  ( N^{ \bot C^* } L^{ \bot C^*})^{ \bot M }=
(  D^{ \bot C^* }E^{ \bot C^*})^{\bot C}.$$ Since  $N$ and  $L$ are
closed subcomodules,
  $N^{ \bot C^* } = C(N)^{ \bot C^*} = D^{ \bot C^* }$ and
$L ^{ \bot C^* }  = C(L)^{ \bot C^* } = E^{\bot C^*}$ by Proposition
\ref {7.1.17}(1) and Lemma \ref{7.1.1}(1).
 Thus we only need to show that

$(D^{ \bot C^*} E^{ \bot C^*})^{ \bot M}  =
(D^{ \bot C^* } E^{ \bot C^*})^{ \bot C }$.  \\
The above formula follows from Proposition \ref {7.1.17}(1).

(6)  By the proof of Corollary \ref {7.1.10}(5), we know that $M$ is
a component faithfulness $C$-comodule. Let $\{m_{ \lambda }  \mid
\lambda \in \Lambda \}$ be a basis of $M$. For any   $c \in C$,
$\Delta(c) =
        \sum m_{ \lambda } \otimes c_{\lambda }$,
by the definition of  $C(M)$, $ c_{\lambda } \in C(M)$. Since   $c=
\sum \epsilon(m_{ \lambda }) c_{\lambda } \in C(M)$, $ C \subseteq
C(M)$, i.e. $ C = C(M)$. Consequently, it follows from part (4)(5)
that  $M$ is  $W$-relational hereditary.

(7)  Since  $M$ is  $W$-relational hereditary,
 $ M = \sum_{ \alpha \in \bar{ \Omega }}
 \oplus M_{(E_{ \alpha })^{ (\infty)}}$
is a weak-closed indecomposable decomposition  of $M$ by Proposition
\ref{7.1.14}. By part (1) (3), $ C(M_{(E_{ \alpha })^{ (\infty)}})
=(E_{ \alpha })^{ (\infty)}  =
 (M_{(E_{ \alpha })^{ (\infty)}})^{ \bot C^* \bot M}
=  M_{(E_{ \alpha })^{ (\infty)}}$. Thus   $ M = \sum_{ \alpha \in
\bar{ \Omega }} \oplus M_{(E_{ \alpha })^{ (\infty)}} = \sum_{
\alpha \in \bar{ \Omega }}  \oplus (E_{ \alpha })^{ (\infty)}$. By
\cite{XF92},
 $ \sum_{ \alpha \in \bar{ \Omega }}  \oplus (E_{ \alpha })^{ (\infty)}
= C$ is a indecomposable decomposition  of $C$. Thus the weak-closed
indecomposable decomposition  of $M$
 as a $C$-comodule and
the indecomposable decomposition  of $C$ as coalgebra are the same.
        This completes the proof. \begin{picture}(8,8)\put(0,0){\line(0,1){8}}\put(8,8){\line(0,-1){8}}\put(0,0){\line(1,0){8}}\put(8,8){\line(-1,0){8}}\end{picture}

By Lemma \ref{7.1.6}(4) and Proposition \ref {7.1.17}, $C$ is
cosemisimple iff every $C$-comodule $M$ is cosemisimple.

\section { The coradicals of comodules}\label {s18}
\begin{Proposition} \label{7.2.1}
   Let  $M$ be a  $C$-comodule,  $J$ denote the
Jacobson radical of $C^*$ and   $r_j(C(M)^*)$ denote
 the Jacobson radical of $C(M)^*.$

(1)   $$M_0 = (r_j(C(M)^*))^{ \bot M } = Soc_{C^*}M \hbox { \ \ and
\ \ } r_j(C(M)^*)= M_0^{ \bot C^*};$$

(2) If we view    $C$  as  a right  $C$-comodule, then
$$C_0 = Soc_{C^*}C = J^{ \bot C }
= \sum \{ D \mid D \hbox { \ is a minimal right coideal of \ }  C
\}.$$

\end{Proposition}
 {\bf Proof}.  (1)  We first show that
$$ M_0^{ \bot C^*}
\subseteq J$$ when  $M$ is a full  $C$-comodule. We only need to
show that  $ \epsilon - f$  is invertible in
 $C^*$ for any  $f \in (M_0)^{\bot C^*} $.
Let  $I = (M_0)^{ \bot C^*}$. For any natural number $n$, $f^{n+1}
\cdot (M_0)^{(n)} =0$ since
 $f^{n+1} \in I^{n+1}$,
 where $(M_0)^{(n)} = \wedge ^{n+1}M_0$.
 Thus
\begin {eqnarray}
f^{n+1}(C((M_0)^{(n)})) &=& 0 \label {e.18.1}.
\end {eqnarray}
 by  Lemma \ref {7.1.1}(1) and

\begin {eqnarray}      M=(M_0)^{ (\infty)}
&=& \cup \{ (M_0)^{(n)} \mid n = 0, 1, \cdots \}
 \label {e.18.2}
 \end {eqnarray} by \cite[Theorem 4.7]{XSF94}.
We now show that
\begin {eqnarray}
C  &=& \cup \{ C((M_0)^{(n)}) \mid n = 0, 1, \cdots \} \label
{e.18.3}
\end {eqnarray}
Since  $(M_0)^{(n)} \subseteq (M_0)^{(n+1)},$  we have that there
exists a basis  $ \{ m_\lambda \mid  \lambda \in \Lambda  \}$ of $M$
such that for every given natural number  $n$ there exists a subset
of $ \{ m_\lambda \mid  \lambda \in \Lambda  \},$ which is a basis
of $(M_0)^{(n)}$. For any $c \in W(M)$, there exists  $m \in M$ with
$\rho(m) = \sum m_{ \lambda } \otimes c_{\lambda }$ such that
   $c_{\lambda _0 } = c $ for some $\lambda_0 \in \Lambda $.
  By equation (\ref {e.18.2}), there exists a natural number $n$ such that
   $m \in (M_0)^{(n)},$ which implies that
$c \in C((M_0)^{(n)})$ and $C(M) \subseteq  \cup _{0}^{\infty}
C((M_0)^{(n)})$. Thus
$$C= \cup \{ C((M_0)^{(n)}) \mid n = 0, 1, \cdots \}.$$
Let
$$g_n = \epsilon + f + \cdots + f^n , n = 1, 2, \cdots.$$

For any $c \in C$, there exists a natural number $n$
 such that
$c \in  C((M_0)^{(n)})$ by relation (\ref {e.18.3}). We define
\begin {eqnarray}
g(c) = g_n(c)\ . \label {e.18.4}
\end {eqnarray}
Considering relation (\ref{e.18.1}), we have  that $g$ is
well-defined. Thus  $g \in C^*$.

We next show that  $g$ is an inverse of  $ \epsilon -f $ in $C^*$.
For any $c \in C^*$,  there exists a natural number  $n$ such that
$c \in C((M_0)^{(n)})$ by relation (\ref {e.18.3}). Thus  $\Delta
(c) \in C((M_0)^{(n)}) \otimes C((M_0)^{(n)})$. We see that
\begin {eqnarray*}
g*( \epsilon -f)(c) &=& \sum g(c_1)(\epsilon -f)(c_2) \\
 &=& \sum g_n(c_1)(\epsilon -f)(c_2) \hbox { (by relation (\ref {e.18.4} )) } \\
 &=& (g_n *(\epsilon -f))(c)   \\
 &=& (\epsilon -f^{n+1})(c)  \\
 &=& \epsilon (c) \hbox { ( by relation  (\ref {e.18.1})) }.
 \end {eqnarray*}
Thus  $g*(\epsilon -f) = \epsilon $. Similarly, $(\epsilon -f)*g =
\epsilon $. Thus  $\epsilon -f$ has an inverse in  $C^*$.

        We next show that  $$ M_0^{ \bot C^*} = J$$ when $C(M)=C$.
It follows from Lemma \ref {7.1.1}(1) that

\begin {eqnarray*} (M_0)^{ \bot C^* }
&=& \cap \{ N^{\bot C^*} \mid N \hbox { \  is
a minimal closed subcomodule of \ } M \} \\
&=& \cap \{ C(N)^{ \bot C^*} \mid N \hbox { is a minimal closed
subcomodule
 of } M \}       \\
& \supseteq & J  \hbox { (  \ by \cite [Proposition 2.1.4] {HR74}
and
  Proposition \ref {7.1.2}) }
\end {eqnarray*}
Thus $M_0^{ \bot C^*} = J$

We now show that
 $$r_j(C(M)^*) = M_0 ^ {\bot C^*} \hbox { \ \ \ and \ \ \ }
 (r_j(C(M)^*))^{\bot M} = M_0$$
 for any $C$-comodule $M$.
If $M$ is  a $C$-comodule,  then $M$ is full $C(M)$-comodule and so
$$(r_j(C(M)^*) = M_0^{\bot C^*}.$$ By  \cite [Proposition 4.6]
{XSF94}, $M_0$ is closed. Thus   $((r_j(C(M)^*))^{\bot M} = M_0$.

 Finally, we  show that   $M_0 = Soc_{C^*}M$ for any $C$-comodule $M$. By \cite
[Proposition 1.16]{XSF94}, $Soc_{C^*}M \subseteq M_0$. Conversely,
if  $x \in M_0$, let   $N = C^*x$. By Lemma \ref{7.1.6}(4),
 $M_{C_0} = M_0$. Thus
 $N$ is a finite dimensional  $C_0$-comodule. By  \cite [Theorem14.0.1]{Sw69a},
$N = N_1 \oplus \cdots \oplus N_n$ and
 $N_i$ is a simple  $C_0^*$-submodule. By  \cite [Lemma 1.1]{XSF94},
$N_i$ is a  $C_0$-subcomodule. Thus  $N_i$ is a $C$-comodule. By
\cite[Lemma 1.1]{XSF94}, $N_i$ is a  $C^*$-submodule. If  $L$ is  a
non-zero $C^*$-submodule of $M$ and $L \subseteq N_i$, then  $L$ is
also a  $(C_0)^*$-submodule. Thus
 $L = N_i$. This shows     that
 $N_i$ is also a simple  $C^*$-submodule. Thus  $N_i
\subseteq Soc_{C^*}M$ and  $N \subseteq Soc_{C^*}M$.
 It follows from the above proof that $M_0 = Soc_{C^*}M$.

(2) It follows from  Proposition \ref{7.1.17} and part (1).

\chapter { The Homological Dimensions of Crossed Products}\label {c9}

Throughout this chapter, $k$ is a field, $R$  is an algebra over
$k$, and  $H$ is a Hopf algebra over $k$ . We say that $R \#_\sigma
H$  is the crossed product of $R$ and $H$ if $R\#_\sigma H$ becomes
an algebra over $k$ by multiplication:
$$(a \# h) (b \# g) =
\sum _{h, g } a (h_1 \cdot b) \sigma (h_2, g_1) \# h_3g_2 $$ for any
$a, b \in R, h, g \in H,$ where $\sigma $ is a k-linear map from
$H\otimes H$ to $R$ and
 $\Delta (h) = \sum h_1\otimes h_2$
\ \ \ (  see \cite [Definition 7.1.1] {Mo93}.)

Let $lpd ({}_RM)$, $lid({}_RM)$ \ \  and \ \  $lfd ({}_RM)$ denote
the left projective dimension, the left injective dimension and the
left flat dimension of left $R$-module $M$ respectively. Let $lgD
(R)$ \ \  and \ \  $wD (R)$ denote the left global  dimension and
the weak dimension of algebra $R$ respectively.

Crossed products are very  important algebraic structures. The
relation between homological dimensions of algebra $R$  and crossed
product $R \#_\sigma H$ is often  studied. J.C. Mconnell and J.C.
Robson in \cite [Theorem 7.5.6] {MR87} obtained that
$$rgD (R) = rgD (R*G)$$ for
any finite group $G$ with $\mid G \mid ^{-1} \in k$, where $R*G$ is
the skew group ring. It is clear that every skew group ring $R*G$ is
crossed product $R \# _\sigma kG$ with trivial $\sigma$. Zhong Yi in
\cite {Yi84}
 obtained            that
the global dimension of crossed product $R*G $ is finite when the
global dimension of $R$ is finite  and some other conditions hold.

In this chapter, we obtain that the global dimensions of $R$ and
 the crossed product
$R \# _\sigma H$ are the same; meantime,
 their weak dimensions are also the same,  when
 $H$ is a finite-dimensional semisimple and cosemisimple Hopf algebra.

                \section {The homological dimensions of modules over
                crossed products}\label {s19}

In this section, we give a relation between homological dimensions
of modules over $R$  and  $R\#_\sigma H$.

If $M$ is a left (right) $R\#_\sigma H$-module, then $M$ is also a
left (right ) $R$-module since we can view $R$ as a subalgebra of
$R\#_\sigma H$.

\begin {Lemma} \label {8.1.1}
Let $R$  be a subalgebra of algebra $A$.

(i) If $M$ is a free $A$-module and $A$ is a free $R$-module, then
$M$  is   a free $R$-module;

(ii) If $P$ is a projective left $R\#_\sigma H$-module, then $P$ is
a projective left $R$-module;

(iii) If $P$  is a projective right $R\#_\sigma H$-module and $H$ is
a Hopf algebra with invertible antipode, then $P$ is a projective
right $R$ -module;

(iv) If $${\cal P}_M : \hbox { \ \ \ \ } \hbox { \ \ \ \ } \cdots
P_n \stackrel {d_n} {\rightarrow} P_{n-1} \cdots \rightarrow P_0
\stackrel {d_0}{\rightarrow}  M \rightarrow 0 $$
 is  a projective resolution of left $R\#_\sigma H$-module $M$,
 then ${\cal P}_M$ is a projective resolution of left $R$-module $M$;

 (v) If
 $${\cal P}_M: \hbox { \ \ \ \ } \cdots  P_n \stackrel {d_n} {\rightarrow} P_{n-1} \cdots
\rightarrow P_0 \stackrel {d_0}{\rightarrow}  M \rightarrow 0 $$
 is a projective resolution of right  $R\#_\sigma H$-module $M$ and
 $H$ is a Hopf algebra with invertible antipode, then
  ${\cal P}_M$  is a projective resolution of right $R$-module $M$ .

\end {Lemma}

{\bf Proof.}  (i) It   is  obvious.

(ii)  Since $P$ is a projective $R\#_\sigma H$-module, we have that
there exists a free $R \#_\sigma H$-module $F$ such that $P$ is a
summand of $F$.
 It is clear that $R\# _\sigma H \cong R \otimes H$  as left $R$-modules,
 which implies that  $R\#_\sigma H$ is a free $R$-module.
Thus it follows from part (i)  that $F$ is a free $R$-module and $P$
is a summand of $F$ as $R$-module. Consequently, $P$ is a projective
$R$-module.

(iii)  By \cite [Corollary 7.2.11] {Mo93}, $R\#_\sigma H \cong H
\otimes R$ as right $R$-modules. Thus $R \#_\sigma H$  is a free
right $R$-module. Using the  method  in the proof of part (i), we
have that $P$ is a projective right $R$-module.

(iv) and (v) can be obtained by part (ii) and (iii).
\begin{picture}(8,8)\put(0,0){\line(0,1){8}}\put(8,8){\line(0,-1){8}}\put(0,0){\line(1,0){8}}\put(8,8){\line(-1,0){8}}\end{picture}

\begin {Lemma} \label {8.1.2.1}
(i) Let $R$ be a subalgebra of an algebra $A$. If $M$ is a  flat
right ( left ) $A$-module and $A$ is a flat right ( left )
$R$-module, then $M$  is   a flat right ( left ) $R$-module;

(ii) If $F$ is a flat left $R\#_\sigma H$-module, then $F$ is a flat
left $R$-module;

(iii) If $F$  is a flat right $R\#_\sigma H$-module and $H$ is a
Hopf algebra with invertible antipode, then $F$ is a flat right $R$
-module;

(iv) If $${\cal F}_M: \hbox { \ \ \ \ } \cdots  F_n \stackrel {d_n}
{\rightarrow} F_{n-1} \cdots \rightarrow F_0 \stackrel
{d_0}{\rightarrow}  M \rightarrow 0 $$
 is  a flat resolution of left $R\#_\sigma H$-module $M$,
 then ${\cal F}_M$ is a flat resolution of left $R$-module $M$;

 (v) If
 $${\cal F}_M: \hbox { \ \ \ \ } \cdots  F_n \stackrel {d_n} {\rightarrow} F_{n-1} \cdots
\rightarrow F_0 \stackrel {d_0}{\rightarrow}  M \rightarrow 0 $$
 is a flat resolution of right  $R\#_\sigma H$-module $M$ and
 $H$ is a Hopf algebra with invertible antipode, then
  ${\cal F}_M$  is a flat resolution of $M$.
  \end {Lemma}
{\bf Proof.} (i)  We only show part (i) in the case which $M$ is a
right $A$-module and $A$ is a right $R$-module; the other cases
 can    be shown similarly .
Let $$0 \rightarrow X \stackrel {f}{\rightarrow} Y$$ be an exact
left $R\#_\sigma H$-module sequence. By umptions,
    $$0 \rightarrow A  \otimes _R X \stackrel { A \otimes f}{\rightarrow}
    A \otimes _R Y$$
    and
    $$0 \rightarrow M \otimes _A (A \otimes _R X )\stackrel {M \otimes
    (A\otimes f)}
    {\rightarrow} M \otimes _A (A \otimes _R Y)$$
    are exact sequences.
    Obviously, $$M \otimes _A (A \otimes _R X) \cong M \otimes _RX
    \hbox { \ \  \ and }  M \otimes _A (A \otimes _R Y) \cong M \otimes _RY $$
    as additive groups.
    Thus
    $$0 \rightarrow M  \otimes _R X \stackrel {M\otimes f}
    {\rightarrow} M  \otimes _R Y$$
    is an exact  sequence, which implies
    $M$ is a flat $R$-module.

    (ii)--(v)  are immediate consequence of part (i).
    \begin{picture}(8,8)\put(0,0){\line(0,1){8}}\put(8,8){\line(0,-1){8}}\put(0,0){\line(1,0){8}}\put(8,8){\line(-1,0){8}}\end{picture}

              The following is an immediate consequence of Lemma \ref {8.1.1} and
             \ref {8.1.2.1}.

\begin {Proposition} \label {8.1.2}
(i) If $M$ is a left $R\#_\sigma H$-module, then $$lpd ({}_RM) \leq
lpd({}_{R \#_\sigma H}M);$$

(ii) If $M$ is a right $R\#_\sigma H$-module and $H$ is a Hopf
algebra with invertible antipode, then
 $$rpd (M_R) \leq rpd({M}_{R \#_\sigma H});$$

(iii) If $M$ is a left $R\#_\sigma H$-module, then $$lfd ({}_RM)
\leq lfd({}_{R \#_\sigma H}M);$$

(iv) If $M$ is a right $R\#_\sigma H$-module and $H$ is a Hopf
algebra with invertible antipode, then
 $$rfd (M_R) \leq rfd({M}_{R \#_\sigma H}).$$

\end {Proposition}

\begin {Lemma} \label {8.1.4}
 Let $H$ be a finite-dimensional semisimple  Hopf algebra,
 and let $M$ and $N$ be left $R\#_\sigma H$-modules.
  If
 $ f$  is an $R $-module homomorphism from $M$ to $N$,
  and
  $$\bar f (m) = \sum \gamma ^{-1}(t_1)f(\gamma (t_2)m)$$
 for any $m \in M,$ then
 $\bar f$  is an $R \#_\sigma H$-module homomorphism from $M$ to $N$,
 where  $t \in \int _H^r$ with $\epsilon (t)=1$,
 $\gamma $  is a map from $H$ to $R \#_\sigma H$ sending $h$ to $1 \# h$
 and  the convolution-inverse
$$\gamma ^{-1} (h) = \sum \sigma ^{-1}(S(h_2), h_3) \#S(h_1)$$
for any $h \in H.$

\end {Lemma}
{\bf Proof.} (See the proof of \cite [Theorem 7.4.2] {Mo93})
 For any $a\in R, h \in H , m \in M$, we see that
\begin {eqnarray*}
\bar f (am) &=&\sum \gamma ^{-1}(t_1)f(\gamma (t_2)am ) \\
&=&\sum \gamma ^{-1}(t_1)f((t_2 \cdot a) \gamma (t_3)m ) \\
&=&\sum \gamma ^{-1}(t_1)(t_2 \cdot a)f(\gamma (t_3)m ) \\
 &=&\sum a \gamma ^{-1}(t_1)f(\gamma (t_2)m ) \\
 &=&a \bar f (m)
 \end {eqnarray*}
and
\begin {eqnarray*}
 \bar f (\gamma (h)m) &=&\sum \gamma ^{-1}(t_1)f(\gamma (t_2)
 \gamma (h)m ) \\
&=&\sum \gamma ^{-1}(t_1)f(\sigma (t_2, h_1) \gamma (t_3h_2)m )
\hbox {\ \ \ by \cite [Definition 7.1.1] {Mo93}} \\
&=&\sum \gamma ^{-1}(t_1)\sigma (t_2 , h_1)f(\gamma (t_3h_2)m ) \\
&=& \sum \gamma (h_1)\gamma ^{-1}(t_1h_2) f(\gamma (t_2h_3)m)\\
&=&\sum  \gamma (h) \gamma ^{-1}(t_1)f(\gamma (t_2)m )
 \hbox { \ \ \ since } \sum h_1 \otimes t_1h_2 \otimes t_2 h_3 =
 \sum h \otimes t_1 \otimes t_2 \\
  &=&\gamma (h) \bar f (m)
 \end {eqnarray*}
   Thus $\bar f$ is an $R\#_\sigma H$-module homomorphism.
  \begin{picture}(8,8)\put(0,0){\line(0,1){8}}\put(8,8){\line(0,-1){8}}\put(0,0){\line(1,0){8}}\put(8,8){\line(-1,0){8}}\end{picture}

  In fact, we can obtain a functor by Lemma \ref {8.1.4}. Let ${}_{R\#_\sigma
  H} \overline {\cal M}$ denote the full subcategory of ${}_R {\cal M},$
  its objects are all the left $R\#_\sigma H$-modules and its morphisms
from $M$  to $N$  are all the $R$-module homomorphisms from $M$ to
$N$. For any  $M, N \in ob
   {} _{R\#_\sigma H} \overline {\cal M}$  and $R$-module homomorphism
   $f$ from $M$  to $N$, we define that
   $$F:  {}_{R\#_\sigma H}\overline {\cal M} \longrightarrow
   {}_ {R\#_\sigma H} {\cal M}$$
   such that $$F(M) = M \hbox   {  \ \ \ \ and \ \ \ \  } F(f) = \bar f,$$
   where $\bar f$ is defined in Lemma \ref {8.1.4}. It is clear that
   $F$  is a functor.

\begin {Lemma} \label {8.1.5}
 Let $H$ be a finite-dimensional semisimple  Hopf algebra,
 and let $M$ and $N$ be right $R\#_\sigma H$-modules.
  If
 $ f$  is an $R$-module homomorphism from $M$ to $N$,
  and
  $$\bar f (m) = \sum f(m \gamma ^{-1}(t_1))\gamma (t_2)$$
 for any $m \in M,$ then
 $\bar f$  is an $R \#_\sigma H$-module homomorphism from $M$ to $N$,
 where  $t \in \int _H^r$ with $\epsilon (t)=1,$
 $\gamma $  is a map from $H$ to $R \#_\sigma H$ sending $h$ to $1 \# h.$

\end {Lemma}
{\bf Proof.} (See the proof of \cite [Theorem 7.4.2] {Mo93}.)
 For any $a\in R$ and $ h \in H , m \in M$, we see that
\begin {eqnarray*}
\bar f (ma) &=&\sum f(ma \gamma ^{-1}(t_1))\gamma (t_2)  \\
&=&\sum f(m \gamma ^{-1}(t_1)(t_2 \cdot a)) \gamma (t_3) \\
&=&\sum f(m \gamma ^{-1}(t_1))(t_2 \cdot a) \gamma (t_3) \\
 &=&\sum f(m \gamma ^{-1}(t_1))\gamma (t_2) a \\
 &=& \bar f (m) a
 \end {eqnarray*}
and
\begin {eqnarray*}
 \bar f (m \gamma (h)) &=&\sum f(m \gamma (h) \gamma ^{-1}(t_1))
 \gamma (t_2)  \\
&=&\sum f(m \gamma (h_1) \gamma ^{-1} (t_1 h_2)) \gamma (t_2h_3 )
 \hbox { \ \ \ since } \sum h_1 \otimes t_1h_2 \otimes t_2 h_3 =
 \sum h \otimes t_1 \otimes t_2 \\
&=&\sum f(m \gamma ^{-1}(t_1) \sigma (t_2, h_1)) \gamma (t_3h_2 )
\hbox {\ \ \ by \cite [Definition 7.1.1] {Mo93}} \\
&=&\sum f(m \gamma ^{-1}(t_1))\sigma (t_2 , h_1)\gamma (t_3h_2 )
  \\
 &=&\sum f(m \gamma ^{-1}(t_1)) \gamma (t_2)\gamma (h) \\
  &=& \bar f(m)) \gamma (h)
 \end {eqnarray*}
 Thus $\bar f$  is an $R\#_\sigma H$-module homomorphism.
  \begin{picture}(8,8)\put(0,0){\line(0,1){8}}\put(8,8){\line(0,-1){8}}\put(0,0){\line(1,0){8}}\put(8,8){\line(-1,0){8}}\end{picture}

 \begin {Proposition} \label {8.1.6}  Let
 $H$ be a finite-dimensional semisimple Hopf algebra.

(i) If $P$  is a left (right) $R\#_\sigma H$-modules and a
projective left (right) $R$-module, then $P$ is a projective left
(right) $R \#_\sigma H$-module;

(ii) If $E$  is a left (right) $R\#_\sigma H$-modules and an
injective left (right) $R$-module, then $E$ is an injective left
(right) $R \#_\sigma H$-module;

(iii) If $F$  is a left (right) $R\#_\sigma H$-modules and a flat
left (right) $R$-module, then $F$ is a flat left (right) $R
\#_\sigma H$-module.

\end {Proposition}
{\bf Proof.} (i)  Let $$X \stackrel {f} {\rightarrow} Y \rightarrow
0$$ be an exact sequence of left (right ) $R\#_\sigma H$-modules and
$g$  be a $R\#_\sigma H$-module homomorphism from $P$ to $Y$. Since
$P$ is a projective left (right) $R$-module, we have that there
exists an $R$-module homomorphism $\varphi $ from $P$ to $X$  such
that  $$f \varphi = g.$$ By Lemma \ref {8.1.4} and \ref {8.1.5},
there exists an $R\#_\sigma H$-module homomorphism  $\bar \varphi $
from $P$  to $X$ such that  $$f \bar \varphi = g.$$ Thus $P$ is a
projective left (right )  $R\#_\sigma H$-module.

Similarly, we can obtain the proof of part (ii).

(iii) Since $F$  is a flat left (right ) $R$-module, we have that
characteristic module $ Hom _{\cal Z} (F, {\cal Q}/{\cal Z})$ of $F$
is an injective left (right )  $R$-module by \cite [Theorem 2.3.6]
{To98}. Obviously, $Hom_{\cal Z} (F, {\cal Q}/{\cal Z})$ is a left
(right ) $R\#_\sigma H$-module. By part (ii),
 $Hom_{\cal Z} (F, {\cal Q}/{\cal Z})$ is an injective left ( right ) $R \#_\sigma H$
 -module.  Thus $F$ is a flat left (right ) $R\#_\sigma H$-module.
 \begin{picture}(8,8)\put(0,0){\line(0,1){8}}\put(8,8){\line(0,-1){8}}\put(0,0){\line(1,0){8}}\put(8,8){\line(-1,0){8}}\end{picture}

\begin {Proposition} \label {8.1.7}
Let $H$ be a finite-dimensional semisimple Hopf algebra. Then
 for  left (right ) $R \#_\sigma H $ -modules $M$ and $N$,
 $$Ext_{R\#_\sigma H}^n (M, N) \subseteq Ext_R ^n (M,N), $$
 where $n$ is any natural number.

\end {Proposition}
{\bf Proof.} We view the $Ext^n(M,N)$ as the equivalent classes of
$n$- extensions of $M$ and $N$ (see \cite [Definition 3.3.7]
{To98}). We only prove  this result for the case $n=1$. For  other
cases, we can
 prove them similarly.
We denote the equivalent classes in $Ext_{R\#_\sigma H}^1 (M,N)$ and
$Ext_R^1(M,N)$ by  [E] and [F]'  respectively, where  $E$ is an
extension of $R\#_\sigma H$-modules $M$ and $N$,
 and $F$ is an extension  of $R$-modules $M$ and $N$.
We define a map $$\Psi :Ext _{R\#_\sigma H}^1 (M,N) \rightarrow
Ext_R^1(M,N), \hbox { \ \ \ by sending  \ } [E] \hbox { \ \ to }
[E]' .$$ Obviously, $\Psi$ is a map. Now we show that $\Psi $  is
injective. Let
 $$0 \rightarrow M \stackrel {f}{\rightarrow}  E
 \stackrel {g}{\rightarrow } N \rightarrow 0
\hbox {  \ \ \ and \ } 0 \rightarrow M \stackrel {f'} {\rightarrow}
E'
  \stackrel {g'}{\rightarrow}  M \rightarrow 0 $$
  are two extensions of $R\#_\sigma H$-modules $M$ and $N$, and they
  are equivalent in $Ext_R^1 (M,N)$. Thus there exists an
  $R$-module homomorphism $\varphi $ from $E$ to $E'$ such that
$$\varphi f =f' \hbox { \ \ \ and \ \ } \varphi g = g'.$$
  By lemma \ref {8.1.4} , there exists an $R\#_\sigma H$-module homomorphism
  $\bar \varphi $  from $E$  to $E'$ such that
   $$\bar \varphi f =f' \hbox { \ \ \ and \ \ } \bar \varphi g = g'.$$
Thus $E$ and $E'$  are equivalent in $Ext _{R\#_\sigma H}^1 (M,N),$
which implies that $\Psi$ is  injective.
\begin{picture}(8,8)\put(0,0){\line(0,1){8}}\put(8,8){\line(0,-1){8}}\put(0,0){\line(1,0){8}}\put(8,8){\line(-1,0){8}}\end{picture}

\begin {Lemma} \label {8.1.7}
 For any $M \in {\cal M}_{R \#_\sigma H} $  and $N \in {}_{R \#_\sigma H}
 {\cal M}$, there exists an additive group  homomorphism
 $$\xi:   M \otimes _{R} N  \rightarrow M \otimes _{R \#_\sigma H}N$$
    by sending  $(m \otimes n)$  to $m \otimes n$,
    where $m \in M, n \in N$.

\end {Lemma}
{\bf Proof.} It is trivial.
\begin{picture}(8,8)\put(0,0){\line(0,1){8}}\put(8,8){\line(0,-1){8}}\put(0,0){\line(1,0){8}}\put(8,8){\line(-1,0){8}}\end{picture}

 \begin {Proposition} \label {8.1.8}
 If $M$   is a right $R\#_\sigma H$-modules and
$N$ is a left    $R \#_\sigma H$-module, then
     there exists an additive group homomorphism
     $$\xi_* : {\ } Tor ^R_n (M,N)  \longrightarrow
     Tor ^{R\#_\sigma H}_n (M,N)$$
     such that   $\xi_*([z_n] ) = [\xi(z_n)]$,
     where $\xi$ is the same as in Lemma \ref {8.1.7}.

\end {Proposition}
{\bf Proof.}
 Let
$${\cal P}_M: \hbox { \ \ \ \ } \cdots  P_n \stackrel {d_n} {\rightarrow} P_{n-1} \cdots
\rightarrow P_0 \stackrel {d_0}{\rightarrow}  M \rightarrow 0 $$
 be a projective resolution of right  $R\#_\sigma H$-module $M$,  and
 set
 $$T = -   {\otimes }_{R\#_\sigma H} N \hbox { \ \ \  and  \ \ }
 T^R = - {\otimes }_ {R} N .$$
We have that
  $$ T{\cal P}_{\hat M}: \hbox { \ \ \ \ } \cdots  T(P_n) \stackrel {Td_n} {\rightarrow}
  T(P_{n-1})
  \cdots
\rightarrow T(P_1) \stackrel {Td_1}{\rightarrow}  TP_0 \rightarrow 0
$$ and
$$T^R{\cal P}_{ \hat  M}:\hbox  { \ \ \ } \cdots T^R( P_n) \stackrel {T^Rd_n} {\rightarrow}
T^R( P_{n-1}) \cdots \rightarrow T^R(P_1) \stackrel
{T^Rd_1}{\rightarrow} T^R (P_0) \rightarrow 0 $$ are complexes .
     Thus $\xi$ is a complex homomorphism from
     $T^R {\cal P}_{\hat M}$ to $T{\cal P}_  {\hat M}$,
     which implies that
     $\xi_*$  is an additive group homomorphism.
     \begin{picture}(8,8)\put(0,0){\line(0,1){8}}\put(8,8){\line(0,-1){8}}\put(0,0){\line(1,0){8}}\put(8,8){\line(-1,0){8}}\end{picture}

\section {The global dimensions and weak dimensions of crossed products}\label {s20}

In this section we give the relation between homological dimensions
of $R$
 and  $R \#_\sigma H$.

\begin {Lemma} \label {8.2.2}
If $R$ and $R'$  are Morita equivalent rings, then

(i) rgD(R) = rgD(R');

(ii) lgD(R) = lgD(R');

(iii) wD(R)= wD(R').
\end {Lemma}

{\bf Proof.}  It is an immediate consequence of \cite [Proposition
21.6, Exercise  22.12] {AF74}.
\begin{picture}(8,8)\put(0,0){\line(0,1){8}}\put(8,8){\line(0,-1){8}}\put(0,0){\line(1,0){8}}\put(8,8){\line(-1,0){8}}\end{picture}

\begin {Theorem} \label {8.2.4}
If $H$  is a finite-dimensional semisimple Hopf algebra, then

(i)  $rgD(R\#_\sigma H)\leq rgD(R) ;$

(ii)  $ lgD(R\#_\sigma H) \leq lgD(R);$

(iii) $ wD(R\# _\sigma H) \leq wD(R)$.

\end {Theorem}

{\bf Proof.} (i) When $lgD(R)$ is infinite, obviously  part (i)
holds. Now we assume $lgD(R) =n.$ For any left $R\#_\sigma H$-module
$M,$ and
 a projective resolution of left  $R\#_\sigma H$-module $M$ :
  $${\cal P}_M: \hbox { \ \ \ \ } \hbox { \ \ \ \ } \cdots  P_n \stackrel {d_n} {\rightarrow} P_{n-1} \cdots
\rightarrow P_0 \stackrel {d_0}{\rightarrow}  M \rightarrow 0 ,$$ we
have  that ${\cal P} _M$ is also  a projective resolution of left
$R$-module $M$ by Lemma \ref {8.1.1}. Let $K_n = ker { \ } d_n$  be
syzygy $n$ of ${\cal P}_M$. Since $ lgD(R) = n$, $Ext _R^{n+1} (M,
N)=0$ for any left $R$-module $N$  by \cite [Corollary 3.3.6]
{To98}. Thus $Ext _R^1 (K_n , N)=0,$  which implies $K_n$ is a
projective $R$-module. By Lemma \ref {8.1.6} (i), $K_n$ is a
projective $R\#_\sigma H$-module and $Ext_{R\#_\sigma H}^{n+1}(M,N)=
0$  for any $R\#_\sigma H$-module $N$. Consequently,
$$lgD(R\#_\sigma H) \leq n = lgD(R) \hbox { \ \ \ by \cite
[Corollary 3.3.6] {To98} }.$$

We complete the proof of part (i).

We can  show part (ii)  and part (iii) similarly.
\begin{picture}(8,8)\put(0,0){\line(0,1){8}}\put(8,8){\line(0,-1){8}}\put(0,0){\line(1,0){8}}\put(8,8){\line(-1,0){8}}\end{picture}

\begin {Theorem} \label {8.2.5}  Let
$H$  be a finite-dimensional semisimple and cosemisimple Hopf
algebra. Then

(i) $rgD(R) = rgD(R\#_\sigma H);$

(ii) $rgD(R) = rgD(R\#_\sigma H);$

(iii) $wD(R)= wD(R\# _\sigma H)$.

\end {Theorem}

{\bf Proof.} (i)  By the duality theorem (see \cite [Corollary
9.4.17] {Mo93}), we have that $(R \#_\sigma H)\# H^*$  and $R$  are
Morita equivalent  algebras. Thus $lgD(R) =  lgD ((R \#_\sigma
H)\#H^*)$ by Lemma \ref {8.2.2} (i). Considering Theorem \ref
{8.2.4} (i), we have that
$$ lgD ((R \#_\sigma H)\#H^*) \leq lgD (R \#_\sigma H) \leq lgD (R) .$$
Consequently,
$$lgD (R) = lgD (R\#_\sigma H).$$

Similarly, we can prove (ii) and (iii) .
\begin{picture}(8,8)\put(0,0){\line(0,1){8}}\put(8,8){\line(0,-1){8}}\put(0,0){\line(1,0){8}}\put(8,8){\line(-1,0){8}}\end{picture}

\begin {Corollary} \label {8.2.6}  Let
$H$  be a finite-dimensional semisimple  Hopf algebra.

(i) If $R $ is a left (right ) semi-hereditary, then  so is
$R\#_\sigma H;$

(ii) If $R $  is von Neumann regular,  then  so is $R\#_\sigma H.$

\end {Corollary}

{\bf Proof.} (i)  It follows from Theorem \ref {8.2.4} and \cite
[Theorem 2.2.9] {To98}          .

(ii) It follows from Theorem \ref {8.2.4}  and \cite [Theorem
3.4.13] {To98}.
     \begin{picture}(8,8)\put(0,0){\line(0,1){8}}\put(8,8){\line(0,-1){8}}\put(0,0){\line(1,0){8}}\put(8,8){\line(-1,0){8}}\end{picture}

By the way, part (ii) of Corollary \ref {8.2.6} gives one case about
the semiprime question in \cite [Question 7.4.9] {Mo93}. That is, if
$H$ is a finite-dimensional  semisimple Hopf algebra and $R$ is a
von Neumann regular  algebra (notice that every von Neumann regular
algebra is  semiprime ),  then    $R\#_\sigma H$ is semiprime.

\begin {Corollary} \label {8.2.8}  Let
$H$  be a finite-dimensional semisimple  and cosemisimple Hopf
algebra. Then

 (i) $R$  is  semisimple artinian iff $R\#_\sigma H$  \ \ is semisimple artinian;

(ii) $R $ is left (right ) semi-hereditary iff  $R\#_\sigma H$ \ \
is left (right ) semi-hereditary;

(iii)  $R $  is von Neumann regular iff $R\#_\sigma H $ \ \ is von
Neumann regular.
\end {Corollary}

{\bf Proof.} (i) It follows from Theorem \ref {8.2.5}  and \cite
[Theorem 2.2.9] {To98}.

(ii)  It follows from Theorem \ref {8.2.5} and \cite [Theorem 2.2.9]
{To98}.

(iii) It follows from Theorem \ref {8.2.5}  and \cite [Theorem
3.4.13] {To98}.
\begin{picture}(8,8)\put(0,0){\line(0,1){8}}\put(8,8){\line(0,-1){8}}\put(0,0){\line(1,0){8}}\put(8,8){\line(-1,0){8}}\end{picture}

If $H$ is commutative or cocommutative, then $S^2 = id_H$  by \cite
{Sw69a}. Consequently, by \cite [Proposition 2 (c)] {Ra94},
 $H$  is semisimple and cosemisimple iff
the character $char k $ of $k$ does  not divides $dim H$.
Considering Theorem \ref {8.2.5} and Corollary \ref {8.2.8}, we have

\begin {Corollary} \label {8.2.9}
Let $H$ be a finite-dimensional commutative or cocommutative Hopf
algebra. If the character $char k $ of $k$ does  not divides $dim
H$, then

(i) $rgD(R) = rgD(R\#_\sigma H);$

(ii) $rgD(R) = rgD(R\#_\sigma H);$

(iii) $wD(R)= wD(R\# _\sigma H);$

 (iv) $R$  is  semisimple artinian iff $R\#_\sigma H$  \ \  is semisimple artinian;

(v) $R $ is left (right ) semi-hereditary iff  $R\#_\sigma H$ \ \ is
left (right ) semi-hereditary;

(vi)  $R $  is von Neumann regular iff $R\#_\sigma H $  \ \ is von
Neumann regular.

\end {Corollary}

Since group algebra $kG$  is a cocommutative Hopf algebra, we have
that
     $$rgD(R) = rgD(R*G).$$
     Thus Corollary \ref {8.2.9} implies  \cite
[Theorem 7.5.6] {MR87}.

\chapter   { The Radicals of Hopf Module Algebras }\label {c10}
Remark \footnote {This chapter can be omitted }

In this chapter, the characterization of $H$-prime radical is given
in many ways. Meantime, the relations between the radical of smash
product $R \# H$ and the $H$-radical of Hopf module algebra $R$ are
obtained.

In this chapter, let $k$ be a commutative associative ring with
unit, $H$ be an algebra with unit and comultiplication
$\bigtriangleup$ ( i.e. $\Delta $ is a linear map: $H \rightarrow H
\otimes H$),
  $R$ be an algebra
over $k$ ($R$ may be without unit) and $R$ be an $H$-module algebra.

We define some necessary concept as follows.

If   there exists a linear map $
\left \{ \begin {array} {ll} H \otimes R & \longrightarrow R \\
 h \otimes r & \mapsto   h \cdot r \end {array} \right. $ such that
     $$h\cdot rs = \sum (h_1\cdot r)(h_2 \cdot s) \hbox { \ \ and \ }
     1_H \cdot r = r $$
     for all $r, s \in R, h\in H,$  then we say that
     $H$ weakly acts on $R.$
  For any ideal $I$ of $R$, set
  $$(I:H):= \{ x \in R \mid h\cdot x \in I \hbox { for all } h \in H  \}. $$
 $I$ is called an $H$-ideal if $h\cdot I
\subseteq I$ for any $h \in H$. Let $I_H$ denote the maximal
$H$-ideal of $R$ in $I$. It is clear that $I_H= (I:H).$ An
$H$-module algebra $R$ is called an $H$-simple module algebra if $R$
has not any non-trivial $H$-ideals and $R^2 \not=0.$
 $R$ is said to be  $H$-semiprime if there are no
non-zero nilpotent $H$-ideals in $R$.  $R$ is said to be
  $H$-prime
if $IJ =0$ implies $I=0$ or $J=0$ for any $H$-ideals $I$ and $J$ of
$R$. An $H$-ideal $I$ is called an   $H$-(semi)prime ideal of $R$ if
$R/I$ is  $H$-(semi)prime.
   $ \{ a_{n} \}$  is called an
  $H$-$m$-sequence  in  $R$ with beginning $a$
  if there exist $h_n, h_n' \in H$ such that $a_1 = a \in R$ and
  $a_{n+1} =
  (h_{n}.a_{n})b_{n}(h_{n}'.a_{n})$   for any natural number $n$.
If every $H$-$m$-sequence $\{ a_{n} \}$
  with  $a_{7.1.1} = a$, there exists a natural number $k$ such that $a_{k}= 0,$
  then $a$ is called an $H$-$m$-nilpotent element.
  Set
  $$W_{H}(R) = \{ a \in R \mid a \hbox { \ is an \ } H \hbox{-}m
  \hbox {-nilpotent element} \}. $$
 $R$ is called an $H$-module algebra
 if the following
conditions hold:

   (i)  $R$ is a unital left $H$-module(i.e. $R$ is a left $H$-module and
   $1_H \cdot a = a$
   for any $a \in R$);

   (ii)  $h\cdot ab = \sum (h_1\cdot a)(h_2 \cdot b)$ for any $a, b \in R$,
      $h\in H$, where
$\Delta (h) = \sum h_1 \otimes h_2$.  \\
$H$-module algebra is sometimes called a Hopf module.

If $R$ is an $H$-module algebra with a unit $1_R$, then
 $$h \cdot 1_R = \sum_h (h_1 \cdot 1_R)(h_2S(h_3)\cdot 1_R)$$
 $$ = \sum _h h_1 \cdot
 (1_R (S(h_2)\cdot 1_R)) = \sum _h h_1S(h_2) \cdot 1_R= \epsilon (h)1_R, $$
 i.e.         \ \ \ \ $ \hbox { \ \ \ \ \ }   h \cdot 1_R = \epsilon (h)1_R$ \\
for any $h \in H.$

An $H$-module algebra $R$ is called a unital $H$-module algebra if
$R$ has a unit $1_R$ such that $h \cdot 1_R = \epsilon (h)1_R$ for
any $h \in H$. Therefore, every $H$-module algebra with unit is a
unital $H$-module algebra. A left $R$-module $M$ is called an
$R$-$H$-module if $M$ is also a left unital $H$-module
 with $h  (am)= \sum (h_1 \cdot a)(h_2m)$  for all
$h \in H, a \in R, m \in M$. An $R$-$H$-module $M$ is called an
$R$-$H$- irreducible module if there are no non-trivial
$R$-$H$-submodules in $M$ and $RM \not=0$. An algebra homomorphism
$\psi: R \rightarrow R'$ is called an $H$-homomorphism  if $\psi (h
\cdot a) = h \cdot \psi (a)$ for any $h \in H, a \in R.$ Let $r_b,
r_j, r_l, r_{bm}$  denote the Baer radical, the Jacobson radical,
the locally nilpotent radical,
 the Brown-MacCoy radical
  of algebras respectively.   Let
   $I \lhd_H R$ denote that $I$ is an $H$-ideal of $R.$

\section{ The $H$-special radicals for $H$-module algebras}\label {s21}
J.R. Fisher \cite{Fi75} built up the general theory of $H$-radicals
for $H$-module
 algebras. We can easily give the definitions of the $H$-upper radical and
 the $H$-lower radical for $H$-module algebras as in \cite{Sz82}.
  In this section, we obtain
 some   properties of $H$-special radicals for $H$-module algebras.

\begin{Lemma}\label{9.1.1}
(1)   If $R$ is an $H$-module algebra and $E$ is a non-empty subset
of $R$, then
   $(E) = H\cdot E + R(H \cdot E) + (H\cdot E)R + R(H \cdot E)R$,
    where $(E)$ denotes the $H$-ideal generated by $E$ in  $R$.

(2) If $B$ is an $H$-ideal of $R$ and $C$ is an $H$-ideal of $B$,
then $(C)^3 \subseteq C,$  where $(C)$ denotes the $H$-ideal
generated by $C$ in $R$.
\end {Lemma}
{ \bf Proof. }  It is trivial.
\begin{picture}(8,8)\put(0,0){\line(0,1){8}}\put(8,8){\line(0,-1){8}}\put(0,0){\line(1,0){8}}\put(8,8){\line(-1,0){8}}\end{picture}

  \begin{Proposition}
  \label {9.1.2}   (1)  $R$ is $H$-semiprime iff $(H\cdot a)R(H\cdot a) = 0$ always implies
    $a = 0$ for any $a\in R$.

 (2) $R$ is $H$-prime iff $(H\cdot a)R(H\cdot b) = 0$ always implies
 $a = 0$  or $b = 0$ for any $a, b \in R$.
\end{Proposition}
{\bf Proof.}
 If $R$ is an $H$-prime module algebra and
$(H\cdot a)R(H\cdot b) = 0$ for $a, b \in R$, then $(a)^2 (b)^2 =
0$, where $(a)$ and $(b)$ are the $H$-ideals generated by  $a$ and
$b$  in $R$ respectively. Since $R$ is $H$-prime, $(a) = 0$ or $(b)=
0$. Conversely, if $B$ and $C$ are  $H$-ideals of $R$ and $ BC = 0$,
then $(H \cdot a)R(H\cdot b) = 0$  and $a = 0$ or $b = 0$
 for any $ a \in B, b \in C,$  which implies that $B = 0 $ or $C = 0$,
  i.e. $R$ is an $H$-prime module algebra.

Similarly, part (1) holds.
\begin{picture}(8,8)\put(0,0){\line(0,1){8}}\put(8,8){\line(0,-1){8}}\put(0,0){\line(1,0){8}}\put(8,8){\line(-1,0){8}}\end{picture}

 \begin{Proposition}\label {9.1.3}
  If $ I \lhd_{H} R$ and $I$ is an $H$-semiprime module algebra,
  then    \\
     (1)  $I \cap I^{*} = 0$;
     (2)  $I_{r} = I_{l} = I^{*}$;
     (3)  $I^{*} \lhd_H R$,
     where $I_{r} = \{ a \in R \mid  I(H\cdot a) = 0 \}$,
     $I_{l} = \{ a \in R \mid  (H\cdot a)I = 0 \}$,
     $I^{*} = \{ a \in R \mid  (H\cdot a)I = 0 = I(H\cdot a) \}.$
 \end{Proposition}
 {\bf Proof .} For any $x \in I^* \cap I$, we have that  $I(H\cdot x) = 0$
 and  $(H\cdot x)I(H\cdot x) = 0$. Since $I$ is an $H$-semiprime
 module algebra,  $x =0$,
 i.e. $I \cap I^* =0$.

 To show $I^* = I_r$,
 we only need to show that $(H \cdot x)I = 0$ for any $x \in I_r$.
 For any $y \in I, h \in H$, let $z = (h \cdot x)y$.
 It is clear that   $(H\cdot z)I(H\cdot z) = 0$.
 Since $I$ is an $H$-semiprime module algebra, $z =0$, i.e.
 $(H\cdot x)I = 0$. Thus $I^* = I_r$. Similarly, we can show that
 $I_l = I^*$.

 Obviously, $I^*$ is an ideal of $R$. For any $x \in I^*, h \in H$,
 we have $(H\cdot(h \cdot x)) I = 0 $. Thus $h\cdot x \in I^*$  by
 part (2), i.e. $I^*$ is an $H$-ideal of $R$.
 \begin{picture}(8,8)\put(0,0){\line(0,1){8}}\put(8,8){\line(0,-1){8}}\put(0,0){\line(1,0){8}}\put(8,8){\line(-1,0){8}}\end{picture}

 \begin{Definition}\label {9.1.4} ${\cal K}$ is called an $H$-(weakly )special class if

     (S1)  ${\cal K}$ consists of $H$-(semiprime)prime module algebras.

     (S2)  For any $R \in {\cal K}$, if $0 \not= I \lhd_{H} R$  then
     $I \in {\cal K}$.

     (S3)  If $R$ is an $H$-module algebra and  $B\lhd_{H} R$ with
      $B \in {\cal K}$,
  then    $ R/B^{*}\in{\cal K}$,

where  $B^{*} = \{ a \in R\mid  (H\cdot a)B = 0 = B(H\cdot a)\}$.
    \end{Definition}

      It is clear that (S3) may be replaced by one of the following conditions:

     (S3')   If $B$ is an essential $H$-ideal of $R$(i.e. $B \cap I \not=0$
     for any non-zero $H$-ideal $I$ of $R$) and $B \in {\cal K}$, then
     $R\in {\cal K}$.

(S3")  If there exists an $H$-ideal $B$ of $R$ with $B^{*} = 0$ and
$B\in{\cal K}$,
 then $R\in{\cal K}$.

 It is easy to check that if ${ \cal K}$ is an $H$-special class, then
 ${ \cal K}$ is an $H$-weakly special class.

 \begin{Theorem}\label {9.1.5}   If ${\cal K}$ is an $H$-weakly special class,
  then
   $r^{{\cal K}}(R) = \cap \{ I \lhd_{H} R \mid R/I \in {\cal K} \}$,
   where    $r^{{\cal K}}$ denotes the $H$-upper radical determined by $ {\cal K} $.
   \end{Theorem}
{\bf Proof.} If $I$ is a non-zero $H$-ideal of $R$ and $I \in {\cal
K}$,
 then $R/I^* \in {\cal K}$  by (S3) in Definition \ref{9.1.4}
 and  $I \not\subseteq I^* $ by Proposition  \ref {9.1.3}.
 Consequently, it follows from  \cite [Proposition 5]{Fi75} that
 $$r^{\cal K} (R) = \cap \{ I \mid I \hbox { \ is an \ } H \hbox {-ideal of } R
 \hbox { and } R/I \in {\cal K} \} \hbox { \ . \ \ }  \Box $$

    \begin{Definition}\label {9.1.6} If $r$ is a hereditary
    $H$-radical(i.e. if $R$ is
    an $r$-$H$-module algebra and   $B$ is an $H$-ideal of $R$, then so is $B$ )
    and any nilpotent $H$-module algebra is an $r$-$H$-module algebra,
   then $r$ is called a supernilpotent $H$-radical.
   \end{Definition}
\begin{Proposition}\label {9.1.7}  $r$ is a supernilpotent $H$-radical, then $r$
  is $H$-strongly hereditary,
  i.e. $r(I) = r(R) \cap I$   for any $I \lhd _{H} R$.
\end{Proposition}
{\bf Proof.} It follows from \cite [Proposition 4]{Fi75} .
\begin{picture}(8,8)\put(0,0){\line(0,1){8}}\put(8,8){\line(0,-1){8}}\put(0,0){\line(1,0){8}}\put(8,8){\line(-1,0){8}}\end{picture}

\begin{Theorem}\label {9.1.8}  If ${\cal K}$ is an $H$-weakly special class,
then     $r^{{\cal K}}$ is a supernilpotent $H$-radical.
\end{Theorem}
{\bf Proof.}     Let $r= r^{\cal K}.$
      Since every non-zero $H$-homomorphic image $R'$ of a nilpotent
$H$-module algebra $R$ is  nilpotent and is not $H$-semiprime, we
have that $R$ is an $r$-$H$-module algebra by Theorem  \ref {9.1.5}.
It remains to show that any $H$-ideal $I$ of $r$-$H$-module
 algebra $R$ is an $r$-$H$-ideal.
If $I$ is not an $r$-$H$-module algebra,
 then there exists an $H$-ideal $J$ of $I$ such that
 $ 0 \not= I/J \in {\cal K}$. By (S3), $(R/J)/(I/J)^* \in {\cal K}$.
 Let
 $Q = \{x \in R \mid (H \cdot x)I
 \subseteq J$ and $I(H \cdot x) \subseteq J$ \}.
It is clear that
 $J$ and $Q$ are $H$-ideals of $R$ and  $Q/J= (I/J)^*$.
  Since $R/Q \cong (R/J)/(Q/J) = (R/J)/(I/J)^*$ and $R/Q$ is an $r$-$H$-
  module algebra, we have $(R/J)/(I/J)^*$ is an $r$-$H$-module algebra.
  Thus $R/Q = 0$ and $I^2 \subseteq J$, which  contradicts that
$I/J$ is a non-zero $H$-semiprime module algebra. Thus $I$ is an
$r$-$H$-ideal.
     \begin{picture}(8,8)\put(0,0){\line(0,1){8}}\put(8,8){\line(0,-1){8}}\put(0,0){\line(1,0){8}}\put(8,8){\line(-1,0){8}}\end{picture}

\begin{Proposition}\label {9.1.9}  $R$ is $H$-semiprime iff for any  $0\not= a \in R,$
  there exists an $H$-$m$-sequence
  $\{ a_{n} \}$ in $R$ with $a_{7.1.1} = a$ such that $a_{n} \not= 0$ for all $n$.
  \end{Proposition}
{ \bf Proof.} If $R$ is $H$-semiprime, then for any $0 \not= a \in
R$, there exist
 $b_1 \in R$ , $h_1$ and $h_1' \in  H$  such that
 $ 0 \not= a_2 = (h_1 \cdot a_1) b_1(h_1' \cdot a_1)  \in
 (H \cdot a_1)R(H \cdot a_1)$ by Proposition \ref {9.1.2}, where $a_1 = a$.
 Similarly, for $0 \not= a_2 \in R$, there exist
 $b_2 \in R$ and $h_2$ and $h_2' \in  H$  such that
 $ 0 \not= a_3 = (h_2 \cdot a_2) b_2(h_2' \cdot a_2)  \in
 (H \cdot a_2)R(H \cdot a_2),$  which implies that
 there exists an $H$-$m$-sequence $\{ a_n \}$  such that
  $a_n \not= 0$ for any natural number $n$.  Conversely, it is trivial.
\begin{picture}(8,8)\put(0,0){\line(0,1){8}}\put(8,8){\line(0,-1){8}}\put(0,0){\line(1,0){8}}\put(8,8){\line(-1,0){8}}\end{picture}

\section{$H$-Baer radical}\label {s22}

In this section, we give the characterization of $H$-Baer
radical($H$-prime radical) in many ways.

\begin{Theorem}\label {9.2.1} We define a property $r_{Hb}$ for
$H$-module algebras as follows: $R$ is an $r_{Hb}$-$H$-module
algebra iff every non-zero
   $H$-homomorphic image of $R$ contains a non-zero nilpotent $H$-ideal;
   then $r_{Hb}$ is an $H$-radical property.
   \end{Theorem}
  {\bf Proof.}   It is clear that every $H$-homomorphic image of
$r_{Hb}$-$H$-module algebra is an $r_{Hb}$-$H$-module algebra. If
every non-zero $H$-homomorphic image  $B$ of $H$-module algebra $R$
contains a non-zero  $r_{Hb}$-$H$-ideal $I$, then $I$ contains a
non-zero nilpotent $H$-ideal $J$. It is clear that $(J)$ is a
non-zero nilpotent $H$-ideal of $B$, where (J) denotes the $H$-ideal
generated by $J$ in $B$. Thus $R$ is an $r_{Hb}$-$H$-module algebra.
Consequently, $r_{Hb}$ is an $H$-radical property.
  \begin{picture}(8,8)\put(0,0){\line(0,1){8}}\put(8,8){\line(0,-1){8}}\put(0,0){\line(1,0){8}}\put(8,8){\line(-1,0){8}}\end{picture}

$r_{Hb}$    is called $H$-prime radical or $H$-Baer radical.

  \begin{Theorem}\label {9.2.2}  Let
  $${\cal E} = \{ R \mid R \hbox { is a nilpotent } H\hbox {-module
   algebra } \},$$ then
   $r_{\cal E} = r_{Hb}$,  where    $r_{\cal E}$   denotes the $H$-lower radical determined
 by ${\cal E}$.
  \end{Theorem}
{\bf Proof.}
 If $R$  is an $r_{Hb}$-$H$-module algebra, then every non-zero
$H$-homomorphic image $B$ of $R$ contains a non-zero nilpotent
$H$-ideal $I$. By the definition of the lower $H$-radical, $I$ is an
$r_{\cal E}$-$H$-module algebra. Consequently, $R$ is an $r_{{\cal
E}}$-$H$-module algebra.
 Conversely, since every nilpotent $H$-module algebra is an
 $r_{Hb}$-$H$-module algebra, $r_{\cal E} \leq r_{Hb}$.
\begin{picture}(8,8)\put(0,0){\line(0,1){8}}\put(8,8){\line(0,-1){8}}\put(0,0){\line(1,0){8}}\put(8,8){\line(-1,0){8}}\end{picture}

   \begin{Proposition}\label {9.2.3} R is $H$-semiprime if and only if $r_{Hb}(R) = 0$.
    \end{Proposition}
 {\bf Proof.}
 If $R$ is $H$-semiprime with $r_{Hb}(R) \not= 0$, then there exists
 a non-zero nilpotent $H$-ideal $I$ of $r_{Hb}(R)$.
 It is clear that $H$-ideal $(I)$,
 which the $H$-ideal  generated by $I$ in $R$,
 is  a non-zero nilpotent $H$-ideal of $R$. This contradicts that
 $R$ is  $H$-semiprime. Thus $r_{Hb}(R) =0$.
 Conversely, if $R$ is an $H$-module algebra with $r_{Hb}(R)= 0$ and
 there exists a non-zero nilpotent $H$-ideal $I$ of $R$,
 then $I \subseteq r_{Hb}(R)$.
 We get a contradiction. Thus $R$ is $H$-semiprime if
 $r_{Hb}(R)=0$.
 \begin{picture}(8,8)\put(0,0){\line(0,1){8}}\put(8,8){\line(0,-1){8}}\put(0,0){\line(1,0){8}}\put(8,8){\line(-1,0){8}}\end{picture}

 \begin{Theorem}\label {9.2.4}  If  ${\cal K} = \{ R \mid\ R$ is an $H$-prime module algebra\},
 then $\cal K$ is an $H$-special class and
      $r_{Hb} = r^{{\cal K}}.$
       \end{Theorem}
{\bf Proof.}       Obviously, $(S1)$ holds.
 If $I$ is a non-zero $H$-ideal of an $H$-prime module
algebra $R$  and $BC = 0 $ for
 $H$-ideals $B$ and $C$  of $I$, then $(B)^3(C)^3 =0$
 where $(B)$ and $(C)$ denote the $H$-ideals generated by $B$ and $C$ in $R$
    respectively. Since $R$ is
 $H$-prime, $(B)= 0$ or $(C) =0$, i.e. $B =0$ or $C = 0$.
Consequently, $(S2)$ holds.
 Now we shows that (S3) holds. Let $B$ be an $H$-prime module algebra
 and be an $H$-ideal of $R$. If $ JI \subseteq B^*$ for $H$-ideals
 $I$ and $J$ of $R$, then $(BJ)(IB) = 0$, where
 $B^* = \{ x \in R \mid (H \cdot x)B = 0 = B(H \cdot x) \}$.
 Since $B$ is an $H$-prime module algebra, $BJ=0$ or $IB = 0$. Considering $I$ and $J$
 are $H$-ideals, we have that $B(H \cdot J)= 0$ or
 $(H \cdot I)B=0$. By Proposition \ref {9.1.3}, $J \subseteq B^*$
 or $I \subseteq B^*$, which implies that $R/B^* $  is an $H$-prime
 module algebra. Consequently, $(S3)$  holds and so ${\cal K}$ is an $H$-special
 class.

Next we show that $r_{Hb} = r^{\cal K}$. By Proposition \ref
{9.1.5}, $r^{\cal K}(R) = \cap \{ I \mid I$ is an $H$-ideal of $R$
and $R/I \in {\cal K}\}$. If $R$ is a nilpotent $H$-module algebra,
then $R$ is an $r^{\cal K}$-$H$-module algebra. It follows from
Theorem \ref {9.2.2} that $r_{Hb} \leq r^{\cal K}$. Conversely, if
$r_{Hb}(R)=0$, then $R$ is an $H$-semiprime module algebra
 by Proposition \ref {9.2.3}.
 For any $0 \not= a \in R$, there exist
 $b_1 \in R$ , $h_1, h_1' \in  H$  such that
 $ 0 \not= a_2 = (h_1 \cdot a_1) b_1(h_1' \cdot a_1)  \in
 (H \cdot a_1)R(H \cdot a_1)$, where $a_1 = a$.
 Similarly, for $0 \not= a_2 \in R$, there exist
 $b_2 \in R$ and $h_2, h_2' \in  H$  such that
 $ 0 \not= a_3 = (h_2 \cdot a_2) b_2(h_2' \cdot a_2)  \in
 (H \cdot a_2)R(H \cdot a_2)$. Thus
 there exists an $H$-$m$-sequence $\{ a_n \}$  such that
  $a_n \not= 0$ for any natural number $n$. Let
  $${\cal F} = \{ I
  \mid I \hbox { is  an } H \hbox {-ideal  of } R \hbox { and }
  I \cap \{ a_1, a_2, \cdots \}
  = \emptyset \}.$$
 By Zorn's Lemma, there exists a
maximal element $P$  in ${\cal F}$. If $I$ and $J$ are $H$-ideals of
$R$ and $I \not\subseteq P$ and
 $J \not\subseteq P$, then there exist natural numbers
 $n$  and $m$ such that  $a_n \in I$  and
 $a_m \in J$. Since $0 \not= a_{n +m + 1}
 = (h_{n+m}\cdot a_{n+m})b_{n+m}(h'_{n+m}\cdot a_{n+m}) \in IJ$,
 which implies that
 $IJ \not\subseteq P$ and so $P$ is an $H$-prime ideal of $R$. Obviously, $a \not\in P,$
 which implies that $a \not\in r^{\cal K}(R)$
 and $r^{\cal K}(R)= 0$. Consequently,
  $r^{\cal K} = r_{Hb}$.  \begin{picture}(8,8)\put(0,0){\line(0,1){8}}\put(8,8){\line(0,-1){8}}\put(0,0){\line(1,0){8}}\put(8,8){\line(-1,0){8}}\end{picture}

   \begin{Theorem}\label {9.2.5}  $r_{Hb}(R) = W_{H}(R)$.
     \end {Theorem}
{\bf Proof.} If $0 \not= a \not\in W_H(R)$, then there exists an
$H$-prime ideal $P$ such that $a \not\in P$  by the proof of Thoerem
\ref {9.2.4}. Thus $a \not\in r_{Hb}(R),$ which implies that
$r_{Hb}(R) \subseteq W_H(R)$.
          Conversely, for any $x \in W_H(R)$,
let $\bar R = R/r_{Hb}(R)$. Since $r_{Hb}(\bar R)=0$, $\bar R $ is
an
 $H$-semiprime module algebra by Proposition \ref {9.2.3}. By
 the proof of Theorem \ref {9.2.4}, $W_H(\bar R)=0$.
 For   an $H$-$m$-sequence $\{ \bar a_n \}$
 with $\bar a_1 = \bar x$ in  $\bar R$,  there exist
 $\overline b_n \in \overline R$ and $h_n, h_n' \in H$ such that
 $$ \overline a_{n+1} =
 (h_n \cdot \overline a_n) \overline b_n (h_n' \cdot \overline a_n)$$
                             for any natural number $n$.
Thus there exists
 $ a_n' \in  R$  such that  $a_1'=x $  and
 $  a_{n+1}' = (h_n \cdot  a_n')  b_n (h_n' \cdot  a_n')$ {~}{~}
                             for any natural number $n$.
 Since $ \{ a_n'\}$ is an $H$-$m$-sequence with $a_1' = x$ in $R$,
 there exists a natural number $k$ such that $a_k' =0$. It is easy
 to show that $\overline a_n= \overline a_n'$
  for any natural number $n$ by induction. Thus $ \bar a_k = 0$
  and $\overline x \in W_H(\overline R)$.
  Considering $W_H(\overline R)=0$, we have $x \in r_{Hb}(R)$ and
  $W_H(R) \subseteq r_{Hb}(R)$. Therefore
  $W_H(R) = r_{Hb}(R)$. \begin{picture}(8,8)\put(0,0){\line(0,1){8}}\put(8,8){\line(0,-1){8}}\put(0,0){\line(1,0){8}}\put(8,8){\line(-1,0){8}}\end{picture}

 \begin{Definition}\label {9.2.6}  We define an $H$-ideal  $N_{\alpha}$ in $H$-module
 algebra  $R$ for every ordinal number $\alpha$ as follows:

 (i)  $N_{0} = 0$.

 Let us assume that $N_{\alpha}$ is already defined for $\alpha
\prec\beta$.

  (ii)  If $\beta = \alpha + 1, N_{\beta}/N_{\alpha}$ is the sum
   of all nilpotent  $H$-ideals  of $R/N_{\alpha}$

  (iii)  If $\beta$ is a limit ordinal number, $N_{\beta} =
  {\sum}_{\alpha\prec\beta} N_{\alpha}$.

  By set theory, there exists an ordinal number $\tau$ such that
  $N_{\tau} =  N_{\tau + 1}$.
   \end{Definition}

\begin{Theorem} \label {9.2.7} $N_{\tau} = r_{Hb}(R) =
\cap \{ I \mid I \hbox { is
   an } H \hbox {-semiprime ideal of } R \}$.
   \end{Theorem}
{\bf Proof.}
  Let $ D =\cap \{ I \mid I \hbox { is
   an } H \hbox {-semiprime ideal of } R \}.$
 Since  $R/N_{\tau}$ has not any non-zero nilpotent $H$-ideal, we have  that
$r_{Hb}(R) \subseteq N_{\tau}$ by Proposition \ref {9.2.3}.
Obviously, $D \subseteq r_{Hb}(R).$
 Using transfinite  induction, we can show that
 $N_{\alpha } \subseteq I$ for every $H$-semiprime ideal $I$ of $R$ and every
 ordinal number $\alpha $
 (see the proof of \cite [Theorem 3.7] {ZC91} ).
 Thus $ N_{\tau} \subseteq D,$  which completes the proof.
                        \begin{picture}(8,8)\put(0,0){\line(0,1){8}}\put(8,8){\line(0,-1){8}}\put(0,0){\line(1,0){8}}\put(8,8){\line(-1,0){8}}\end{picture}

   \begin{Definition}
   \label {9.2.8} Let $ \emptyset \not= L \subseteq H$.
   An    $H$-$m$-sequence $\{ a_n \}$ in
    $R$ is called an $L$-$m$-sequence with beginning $a$ if
     $a_{7.1.1} = a$  and $a_{n+1} = (h_{n}.a_{n})b_{n}(h_{n}'.a_{n})$
    such that  $h_{n}, h_{n}'\in L$  for all $n$.
  For every $L$-$m$-sequence $\{a_{n} \}$ with $a_{7.1.1} = a$, there exists a
  natural number $k$ such that $a_{k} = 0$, then $R$ is called an $L$-$m$-nilpotent
  element,  written as  $W_{L}(R) = \{ a \in R \mid a$ is an $L$-$m$-nilpotent
  element\}.
     \end{Definition}

Similarly, we have
  \begin{Proposition}\label {9.2.9}  If $L \subseteq H$ and $H = kL$, then

(i)  $R$ is $H$-semiprime iff $(L.a)R(L.a) = 0$ always implies
  $ a = 0 $ for any  $a\in R$.

 (ii)    $R$ is $H$-prime iff $(L.a)R(L.b) = 0$ always implies
$ a = 0$  or $b = 0$       for any     $a, b \in R$.

(iii) $R$ is $H$-semiprime if and only if for any $0 \not= a \in R$,
there exists an $L$-$m$-sequence $\{a_{n} \}$  with $a_1 = a$ such
that $a_{n}\not= 0 $ for all $n$.

(iv)  $W_{H}(R) = W_{L}(R)$.
\end{Proposition}

\section { The $H$-module theoretical characterization of $H$-special
radicals }\label {s23}
 If $V$ is an algebra over $k$ with unit and $x \otimes 1_V=0$  always implies
 that $x=0$ for any right $k$-module $M$ and for any $x \in M$,
 then $V$ is called a  faithful algebra to tensor. For example,
 if $k$ is a field,   then
 $V$ is  faithful to tensor for any algebra $V$ with unit.

In this section, we need to add the following condition:
  $H$ is  faithful  to tensor.

  We shall characterize
  $H$-Baer radical  $r_{Hb}$, $H$-locally nil radical $r_{Hl}$, $H$-Jacobson
 radical $r_{Hj}$ and $H$-Brown-McCoy radical $r_{Hbm}$ by $R$-$H$-modules.

 We can view  every $H$-module algebra $R$ as a sub-algebra of
 $R \# H$ since $H$ is   faithful to tensor.
 By computation, we have that $$ h\cdot a =
 \sum (1 \# h_1) a (1 \# S(h_2)) $$
 for any $h\in H, a \in R,$ where $S$ is the antipode of $H.$

  \begin{Definition}\label {9.3.1}
 An $R$-$H$-module $M$ is called an $R$-$H$-prime module
if for $M$ the following conditions are fulfilled:

(i)   $RM \not= 0$;

(ii)  If $x$ is an element of $M$ and $I$ is an $H$-ideal of $R$,
then $I(Hx) = 0$ always implies $x = 0$ or $I \subseteq  (0:M)_{R},$
where $(0:M)_R = \{ a \in R \mid  aM=0 \}$.
\end{Definition}
    \begin{Definition} \label {9.3.2}
   We associate to every $H$-module algebra $R$ a class ${\cal M}_{R}$ of
    $R$-$H$-modules. Then the class ${\cal M} = \cup {\cal M}_{R}$
is  called an $H$-special class of modules if the following
conditions are fulfilled:

(M1)  If $M \in {\cal M}_{R}$, then $M$ is an $R$-$H$-prime module.

(M2)  If $I$ is an $H$-ideal of $R$ and $M \in {\cal M}_{I}$, then
$IM \in {\cal M}_{R}$.

(M3)  If $M \in {\cal M}_{R}$ and $I$ is an $H$-ideal of $R$ with
$IM \not= 0$, then $M \in {\cal M}_{I}$.

(M4)  Let $I$ be an $H$-ideal of $R$ and $\bar{R} = R/I$. If
 $M \in {\cal M}_{R}$
and $I \subseteq (0:M)_R$, then $M \in {\cal M}_{\bar {R}}$.
Conversely, if $M \in {\cal M}_{\bar {R}}$, then $M \in {\cal
M}_{R}$.
\end{Definition}

  Let ${\cal M}(R)$ denote $\cap \{ (0:M)_R \mid M \in {\cal M}_{R} \},$
  or $R$ when ${\cal M}_R = \emptyset $.
  \begin {Lemma} \label {9.3.3}
  (1)  If $M$ is an $R$-$H$-module, then $M$ is an $R\#H$-module.
In this case, $(0:M)_{R\#H}\cap R = (0:M)_{R}$ and $(0:M)_{R}$ is an
$H$-ideal of $R$;

 (2)  $R$ is a non-zero $H$-prime module algebra iff there exists a faithful
  $R$-$H$-prime
 module $M$;

(3) Let $I$ be an $H$-ideal of $R$ and $\bar R= R/I$. If $M$ is an
$R$-$H$-(resp. prime, irreducible)module and $I \subseteq (0:M)_R $,
then $M$ is an $\overline R$-$H$-(resp. prime, irreducible)module
(defined by $h \cdot (a + I) = h \cdot a$ and $(a+I) x = ax )$.
Conversely, if $M$ is an $\bar R$-$H$-(resp. prime
irreducible)module, then
 $M$ is an $R$-$H$-(resp. prime, irreducible)module(defined by $h \cdot a = h \cdot (a +I)$
 and $ax = (a+I)x$). In the both cases, it is always true that
$R/(0:M)_R \cong \overline R/(0:M)_{\overline R}$;

 (4)  $I$ is an $H$-prime ideal of $R$ with $I \not= R$ iff there exists an $R$-$H$-prime
 module $M$ such that $I=(0:M)_R$;

(5) If $I$ is an $H$-ideal of $R$ and $M$ is an $I$-$H$-prime
module, then $IM$ is an $R$-$H$-prime module with $ (0:M)_I =
(0:IM)_R\cap I$;

(6) If $M$ is an $R$-$H$-prime module and $I$ is an $H$-ideal of $R$
with $IM \not=0$, then $M$ is an $I$-$H$-prime module;

(7) If $R$ is an $H$-semiprime module algebra with  one side unit,
then $R$ has a unit.
 \end{Lemma}
 {\bf Proof}. (1)  Obviously, $(0:M)_R = (0:M)_{R\#H} \cap R$. For any $h \in H,
  a \in (0:M)_R$,
 we see that $(h\cdot a)M = \sum (1\# h_1) a ( 1 \#S(h_2))M
 \subseteq \sum (1 \# h_1)aM=0$ for any $h\in H, a \in R$. Thus $h\cdot a \in (0:M)_R$,
 which implies $(0:M)_R$ is an $H$-ideal of $R$.

 (2) If $R$ is an $H$-prime module algebra, view  $M = R$ as an
$R$-$H$-module. Obviously, $M$ is faithful. If $I(H\cdot x) = 0$
 for $0 \not= x \in M $
 and an $H$-ideal $I$ of $R$, then
$ I(x)= 0$  and $I = 0$, where (x) denotes the $H$-ideal generated
by $x$  in $R$. Consequently,  $M$ is a faithful $R$-$H$-prime
module. Conversely, let $M$ be a faithful $R$-$H$-prime module. If
$IJ = 0$ for two $H$-ideals  $I$ and $J$ of $R$ with $J \not= 0$,
then $JM \not=0$ and there exists $0 \not= x \in JM$ such that $I(H
x) = 0$. Since $M$ is a faithful  $R$-$H$-prime module, $I = 0$.
Consequently,  $R$ is $H$-prime.

(3)  If $M$ is an $R$-$H$-module, then it is clear that $M$ is a
(left)$\overline R$-module and $h  (\overline ax) = h  ( a x) =
  \sum (h_1 \cdot  a )(h_2 x)
  = \sum \overline{ (h_1 \cdot a) }(h_2  x)
  = \sum (h_1 \cdot \overline a)(h_2  x)$ {~}{~}
  for any $h \in H$,  $a \in R$ and $x \in M$. Thus
  $M$ is an $\overline R$-$H$-module. Conversely, if $M$ is an $\overline R$-$H$-module,
  then $M$ is an (left) $R$-module and
 $$h  (ax)
= h (\overline a x)  =
  \sum (h_1 \cdot \overline a )(h_2 x)
  = \sum \overline{h_1 \cdot a}(h_2  x)
  = \sum (h_1 \cdot  a)(h_2  x)$$
    {~}{~} for any $h \in H$, $a \in R$
  and $x \in M$.  This shows that
  $M$ is an $R$-$H$-module.

  Let  $M$ be an $R$-$H$-prime module and $I$ be an $H$-ideal of $R$
  with $I \subseteq (0:M)_R$. If $\overline J(Hx) =0$ for $0 \not= x \in M$
   and an $H$-ideal $J$ of $R$, then $J(Hx)=0$  and $J \subseteq (0:M)_R$.
   This shows that $\bar J \subseteq (0:M)_{ \overline R}$.
   Thus $M$ is an $R$-$H$-prime module. Similarly,
   we can show the other assert.

(4)  If $I$ is an $H$-prime ideal of $R$ with $R\not=I$,
 then $\overline R = R/I$ is an $H$-prime
module algebra. By Part (2), there exists a faithful $\overline
R$-$H$-prime module $M$. By part (3), $M$ is an $R$-$H$-prime module
with $(0:M)_R = I$. Conversely, if there exists a $R$-$H$-prime $M$
with $I = (0:M)_R$, then $M$ is a faithful $\overline R$-$H$-prime
module by part (3) and $I$ is an $H$-prime ideal of $R$ by part (2).

(5)  First, we show that $IM$ is an $R$-module. We define
 \begin {eqnarray}  \label {e2.3.3.1} a (\sum_{i}a_i x_i )= \sum_{i} (a a_i)x_i
\end {eqnarray}
for any $a \in R$ and $\sum_i a_ix_i \in IM$, where $a_i \in I$ and
$x_i \in M$. If $\sum_{i}a_i x_i = \sum_{i}  a_i' x_i'$ with $a_i$,
$a_i'  \in  R$, $x_i$,  $x_i' \in M$, let $y = \sum_{i}(aa_i )x_i -
\sum_{i}  (aa_i') x_i'$. For any $b \in I$ and $h \in H$, we see
that
\begin {eqnarray*}
b(h y) &=& \sum_{i} b\{h  [(aa_i) x_i -  (aa_i') x_i']\} \\
 &=& \sum_{i}\sum_{(h)} b\{[(h_1 \cdot (a a_i)]( h_2  x_i) -
 [h_1 \cdot (a a_i')] (h_2 x_i') \}  \\
  &=& \sum_{(h)} \sum_i \{b[(h_1 \cdot a)(h_2 \cdot a_i)](h_3 x_i)) -
 b[(h_1 \cdot a)(h_2 \cdot a_i')] (h_3  x_i')\}  \\
  &=& \sum_{(h)} \sum_{i} b(h_1 \cdot a)[h_2  (a_i x_i) -
 h_2 (a_i'x_i')] \\
 &=& \sum_{(h)}  b(h_1 \cdot a)h_2  \sum_i [a_i x_i -
a_i'x_i'] = 0.
 \end {eqnarray*}
 Thus  $I(H y)=0$.
 Since $M$ is an $I$-$H$-prime module and $IM \not= 0$,
 we have that $y=0$. Thus
 this definition in (\ref {e2.3.3.1}) is well-defined.
 It is easy to check that $IM$ is an $R$-module.
 We see that
\begin {eqnarray*}
h (a \sum_i a_ix_i) &=& \sum_i h[(aa_i)x_i]   \\
 &=& \sum_i \sum_h [h_1 \cdot (aa_i)][h_2 x_i]  \\
 &=&  \sum_i \sum_h [(h_1 \cdot a)(h_2 \cdot a_i)](h_3 x) \\
 &=& \sum_h (h_1 \cdot a) \sum_i (h_2 \cdot a_i)(h_3 x_i) \\
  &=& \sum_h (h_1 \cdot a) [ h_2 \sum_i (a_ix_i)]
\end {eqnarray*}
for any $h \in H$ and $\sum_i a_ix_i \in IM.$ Thus  $IM$ is an
$R$-$H$-module.

 Next, we show that $(0:M)_I = (0:IM)_R \cap I$.
 If $a \in (0:M)_I$, then $aM =0$  and $aIM=0$,
 i.e. $a \in (0:IM)_A \cap I$. Conversely, if $a \in (0:IM)_R \cap I$,
 then $aIM=0$. By part (1), $(0:IM)_R$ is an $H$-ideal of $R$.
 Thus $(H\cdot a)IM =0$ and $(H \cdot a)I \subseteq (0:M)_I$.
 Since $(0:M)_I$ is an $H$-prime ideal of $I$ by part (4),
 $a \in (0:M)_I$. Consequently,  $(0:M)_I = (0:IM)_R \cap I$.

Finally, we show that $IM$ is an $R$-$H$-prime module. If $RIM=0$,
then $RI \subseteq (0:M)_R$ and  $I \subseteq (0:M)_R,$ which
contradicts that $M$ is an $I$-$H$-prime module. Thus $RIM \not=0$.
If $J(Hx) = 0$ for $ 0\not= x \in IM $  and an $H$-ideal $J$  of
$R$, then $JI(Hx) \subseteq J(Hx) = 0$. Since $M$ is an
$I$-$H$-prime module, $JI \subseteq (0:M)_I$ and $J(IM) =0$.
Consequently,  $IM$ is an $R$-$H$-prime module.

(6)  Obviously, $M$ is an $I$-$H$-module. If $J(Hx) =0$ for $0 \not=
x \in M$ and an $H$-ideal $J$ of $I$, then $(J)^3(H x) = 0$  and
$(J)^3 \subseteq (0:M)_R$, where (J) denotes the $H$-ideal generated
by $J$ in $R$.
 Since $(0:M)_R$  is an $H$-prime ideal of $R$, $(J) \subseteq (0:M)_R$
  and $J \subseteq (0:M)_I$. Consequently,  $M$ is an $I$-$H$-prime module.

   (7) We can assume that $u$ is a right unit of $R$.
   We see that
   $$(h \cdot (au-a))b = \sum (1 \# h_1 ) (au-a)(1 \# S(h_2)) b =0$$
   for any $a, b \in R, h \in H.$  Therefore
   $(H\cdot (au-a)) R=0$ and $au =a$, which implies that $R$ has a unit.
     \begin{picture}(8,8)\put(0,0){\line(0,1){8}}\put(8,8){\line(0,-1){8}}\put(0,0){\line(1,0){8}}\put(8,8){\line(-1,0){8}}\end{picture}

 \begin{Theorem} \label {9.3.4}
 (1)  If ${\cal M}$ is an $H$-special class of modules and ${\cal K}$
= \{ $R \mid$ there exists a faithful $R$-$H$-module $M \in {\cal
M}_{R}$\}, then ${\cal K}$ is an $H$-special class and $r^{{\cal
K}}(R) = {\cal M}(R)$.

(2)  If ${\cal K}$ is an $H$-special class and ${\cal M}_{R}$ = \{
$M \mid M$ is an $R$-$H$-prime module and $R/(0:M)_R \in {\cal
K}$\}, then ${\cal M} = \cup {\cal M}_{R}$ is an $H$-special class
of modules and $r^{{\cal K}}(R) = {{\cal M}}(R)$.
 \end{Theorem}
 { \bf Proof.}
 (1) By Lemma \ref {9.3.3}(2), $(S1)$ is satisfied. If $I$ is a non-zero
$H$-ideal of $R$  and $R \in {\cal K}$, then there exists a faithful
$R$-$H$-prime module  $M \in {\cal M}_R$. Since $M$ is faithful, $IM
\not=0$ and $M \in {\cal M}_I$ with $(0:M)_I = (0:M)_R\cap I=0$ by
$(M3)$. Thus $I \in {\cal K}$ and $(S2)$ is satisfied. Now we show
that $(S3)$ holds. If $I$ is an $H$-ideal of $R$ with $I \in {\cal
K}$, then there exists a faithful $I$-$H$-prime module $M\in {\cal
M}_I$. By $(M2)$ and Lemma \ref {9.3.3}(5),
 $IM\in {\cal M}_R$ and $0=(0:M)_I=(0:IM)_R \cap I$.
Thus $(0:IM)_R \subseteq I^* $. Obviously, $I^* \subseteq (0:IM)_R$.
Thus $I^* = (0:IM)_R$.
 Using $(M4)$, we have that $IM \in {\cal M}_{\overline R}$
 and $IM$ is a faithful $\overline R$-$H$-module with $\overline R = R/I^*$.
 Thus $R/I^* \in {\cal K}$. Therefore ${\cal K}$ is an $H$-special class.

It is clear that
\begin {eqnarray*}
\{ I \mid I \hbox { is an } H \hbox {-ideal of } R \hbox { and } R/I
\in {\cal K} \}  &=& \{(0:M)_R \mid M \in {\cal M}_R  \}.
  \end {eqnarray*}
  Thus $r^{\cal K}(R) = {\cal M}(R)$.

 (2) It is clear that $(M1)$ is satisfied. If $I$ is an $H$-ideal of
$R$ with $M \in {\cal M}_I$, then $M$ is an $I$-$H$-prime module
with $I/(0:M)_I\in {\cal K}$. By Lemma \ref {9.3.3}(5), $IM$ is an
$R$-$H$-prime module with $(0:M)_I=(0:IM)_R \cap I$. It is clear
that
$$(0:IM)_R = \{ a \in R \mid (H \cdot a )I \subseteq (0:M)_I
\hbox { and } I(H \cdot a)  \subseteq (0:M)_I \} $$ and
$$(0:IM)_R/(0:M)_I= (I/(0:M)_I)^*.$$
Thus $R/(0:IM)_R \cong (R/(0:M)_I)/((0:IM)_R/(0:M)_I) =
(R/(0:M)_I)/(I/(0:M)_I)^* \in {\cal K},$ which implies that $IM \in
{\cal M}_R$ and $(M2)$ holds.
 Let $M \in {\cal M}_R$ and $I$ be an $H$-ideal of $R$
 with $IM \not=0$. By Lemma \ref {9.3.3}(6), $M$ is an $I$-$H$-prime
 module and
 $I/(0:M)_I = I/((0:M)_R \cap I) \cong (I + (0:M)_R)/(0:M)_R$.
 Since $R/(0:M)_R \in {\cal K}$, $I/(0:M)_I \in {\cal K}$ and
 $M \in {\cal M}_I$.
 Thus $(M3)$ holds.
 It follows from Lemma \ref{9.3.3}(3) that $(M4)$ holds.

It is clear that
 $$\{I \mid I \hbox { is an } H \hbox{-ideal of } R \hbox { and }
 0 \not= R/I \in {\cal K} \}
 = \{ (0:M)_R \mid M \in {\cal M}_R \}.$$
 Thus  $r^{\cal K}(R) = {\cal M}(R)$.
$\Box $

  \begin{Theorem} \label {9.3.5}
  Let  ${\cal M}_{R}$ =\{ $M \mid M$
  is an  $R$-$H$-prime module\} for any $H$-module algebra $R$
   and ${\cal M} = \cup {\cal M}_{R}$.
   Then ${\cal M}$ is an $H$-special class of modules and
  ${\cal M}(R) = r_{Hb}(R)$.
  \end{Theorem}
{\bf Proof.} It follows from Lemma \ref {9.3.3}(3)(5)(6) that
  ${\cal M}$ is an $H$-special class of modules. By Lemma \ref {9.3.3}(2),
    $$\{ R \mid R  \hbox { is an } H \hbox {-prime module algebra with }
    R \not=0 \} = $$
    $$\{   R \mid \hbox { there exists a faithful } R \hbox {-}H
   \hbox {-prime module }\}. $$ \\
   Thus  $r_{Hb}(R)= {\cal M}(R)$ by Theorem \ref {9.2.4}(1).
\begin{picture}(8,8)\put(0,0){\line(0,1){8}}\put(8,8){\line(0,-1){8}}\put(0,0){\line(1,0){8}}\put(8,8){\line(-1,0){8}}\end{picture}

  \begin{Theorem} \label {9.3.6}
  Let  ${\cal M}_{R}$ =\{ $M \mid M$
  is an  $R$-$H$-irreducible module\} for any $H$-module algebra $R$
   and ${\cal M} = \cup {\cal M}_{R}.$
  Then ${\cal M}$ is an $H$-special class of modules and
  ${\cal M}(R) = r_{Hj}(R),$            where $r_{Hj}$ is the $H$-Jacobson
   radical of $R$ defined in \cite {Fi75}.
    \end{Theorem}

{\bf Proof.}                       If $M$ is an $R$-$H$-irreducible
module and $J (Hx) = 0$  for $0 \not= x \in M$ and an $H$-ideal $J$
of $R$, let  $N= \{ m \in M \mid J(H  m) = 0 \}$. Since $J(h (am)) =
J(\sum_h (h_1 \cdot a)(h_2  m)) = 0$, $am \in N$ for any $m \in N, h
\in H, a \in R,$  we have
  that $N$ is an $R$-submodule of $M$.
Obviously, $N$ is an $H$-submodule of $M$. Thus $N$ is an
$R$-$H$-submodule of $M$. Since $N \not= 0$, we have that $N = M$
and $JM= 0$, i.e. $J \subseteq (0:M)_R$. Thus $M$ is an
$R$-$H$-prime  module and (M1) is satisfied. If $M$ is an
$I$-$H$-irreducible module and $I$ is an $H$-ideal, then $IM$ is an
$R$-$H$-module. If $N$ is an $R$-$H$-submodule of $IM$, then $N$ is
also an $I$-$H$-submodule of $M,$ which implies that $N=0$ or $N=M$.
Thus $(M2)$ is satisfied. If $M$ is an $R$-$H$-irreducible module
and $I$ is an $H$-ideal of $R$ with $IM\not=0$, then $IM=M$. If $N$
is an non-zero $I$-$H$-submodule of $M$, then $IN$ is an
$R$-$H$-submodule of $M$ by Lemma \ref{9.3.3}(5) and $IN = 0 $  or
$IN=M$. If $IN=0$, then $I \subseteq (0:M)_R$ by the above proof and
$IM=0$. We get a contradiction. If $IN=M$, then $N=M$. Thus $M$ is
an $I$-$H$-irreducible module and $(M3)$ is satisfied.

It follows from Lemma \ref {9.3.3}(3) that (M4) holds. By Theorem
\ref {9.3.4}(1), ${\cal M}(R)= r_{Hj}(R)$.
\begin{picture}(8,8)\put(0,0){\line(0,1){8}}\put(8,8){\line(0,-1){8}}\put(0,0){\line(1,0){8}}\put(8,8){\line(-1,0){8}}\end{picture}

 J.R. Fisher
  \cite[Proposition 2]{Fi75} constructed
  an $H$-radical $r_{H}$ by a common hereditary radical $r$ for algebras,
  i.e. $r_H(R) = (r(R):H) =
  \{ a \in R \mid h\cdot a \in r(R)  \hbox { for any }
  h \in H \}.$  Thus we can get
  $H$-radicals $r_{bH}, r_{lH}, r_{jH}, r_{bmH}$.

  \begin{Definition} \label {9.3.7}
  An $R$-$H$-module $M$ is called an $R$-$H$-$BM$-module, if for
    $M$ the following conditions are fulfilled:

    (i)  $RM \not=0$;

    (ii)  If $I$ is an $H$-ideal of $R$ and $I \not\subseteq (0:M)_R$, then there exists
    an element $u \in I$ such that $m = um$ for all $m  \in M$.
\end{Definition}

\begin{Theorem} \label {9.3.8}
  Let ${\cal M}_{R}$ = \{ $M \mid M$ is an $R$-$H$-$BM$-module\} for
  every $H$-module algebra $R$ and ${\cal M} = \cup {\cal M}_{R}$. Then ${\cal M}$ is an
  $H$-special class of  modules.
\end{Theorem}
{\bf Proof.} It is clear that {\cal M} satisfies $(M_1)$  and
$(M_4)$. To prove $(M_2)$ we exhibit: if $I \lhd _H R$  and $M \in
{\cal M}_I$, then $M$ is an $I$-$H$-prime module and $IM$ is an
$R$-$H$-prime module. If $J $ is an $H$-ideal of $R$ with $J
\not\subseteq (0:M)_R, $  then $JI$ is an $H$-ideal of $I$ with $JI
\not\subseteq (0:M)_I.$  Thus there exists an element $u \in JI
\subseteq J$ such that
 $um = m $ for every $m \in M$. Hence
$IM \in {\cal M}_R.$

To prove $(M_3)$, we exhibit: if $M \in {\cal M}_R$ and $I$ is an
$H$-ideal of $R$ with $IM \not=0.$ If $J$ is an $H$-ideal of $I$
with $J \not\subseteq (0:M)_I,$ then $ (J) \not\subseteq (0:M)_R,$
where $ (J)$  is the $H$-ideal generated by $J$ in $R.$  Thus  there
exists  an elements $u\in  (J)$ such that  $um=m$  for every $m \in
M.$ Moreover, $$m=um=uum=uuum = u^3 m $$ and $u^3 \in J.$  Thus $M
\in {\cal M}_I.$
\begin{picture}(8,8)\put(0,0){\line(0,1){8}}\put(8,8){\line(0,-1){8}}\put(0,0){\line(1,0){8}}\put(8,8){\line(-1,0){8}}\end{picture}

\begin{Proposition} \label {9.3.9}
  If $M$ is an $R$-$H$-$BM$-module, then
  $R/(0:M)_R$ is an  $H$-simple module algebra with unit.
\end {Proposition}

{\bf Proof. }  Let $I$ be any $H$- ideal of $R$  with $I
\not\subseteq (0:M)_R.$ Since $M$ is an $R$-$H$-$BM$-module, there
exists an element $u \in I$ such that $uam =am$  for every $m\in M,
a \in R.$ It follows that $a -ua \in (0:M)_R,$ whence $R=I+(0:M)_R.$
Thus $(0:M)_R$ is a maximal $H$-ideal of $R.$ Therefore $R/(0:M)_R$
is an $H$-simple module algebra.

Next we shall show that $R/(0:M)_R$ has a unit. Now $R \not\subseteq
(0:M)_R,$ since $RM \not=0.$  By the above proof, there exists an
element $u \in R$ such that $a-ua \in (0:M)_R$ for any $a \in R.$
Hence $R/(0:M)_R$ has a left unit. Furthermore,  by Lemma \ref
{9.3.7} (7) it has a unity element.
\begin{picture}(8,8)\put(0,0){\line(0,1){8}}\put(8,8){\line(0,-1){8}}\put(0,0){\line(1,0){8}}\put(8,8){\line(-1,0){8}}\end{picture}

\begin{Proposition} \label {9.3.10}
  If $R$ is an $H$-simple-module algebra with unit, then
  there exists a faithful $R$-$H$-$BM$-module.
  \end {Proposition}
{\bf Proof. }  Let $M= R.$ It is clear that $M$ is a faithful
$R$-$H$-$BM$- module.
\begin{picture}(8,8)\put(0,0){\line(0,1){8}}\put(8,8){\line(0,-1){8}}\put(0,0){\line(1,0){8}}\put(8,8){\line(-1,0){8}}\end{picture}

\begin{Theorem} \label {9.3.11}
  Let ${\cal M}_{R}$ = \{ $M \mid M$ is an $R$-$H$-$BM$-module\} for
  every $H$-module algebra $R$ and ${\cal M} = \cup {\cal M}_{R}$. Then
$r_{Hbm}(R) = {\cal M}(R)$,
  where $r_{Hbm}$ denotes the $H$-upper radical determined by $\{R \mid R$ is
  an $H$-simple     module algebra with unit \}.
\end{Theorem}
 {\bf Proof.} By Theorem \ref {9.3.8}, ${\cal M}$ is an $H$-special class
 of modules. Let $${\cal K} = \{   R \mid \hbox { there exists a faithful }
 R  \hbox {-}H \hbox {-}BM\hbox {-module }   \}.$$  By Theorem \ref {9.3.4}(1),
   ${\cal K}$ is an $H$-special class  and  $r ^{\cal K}(R) = {\cal M}(R).$
  Using Proposition \ref  {9.3.9} and \ref {9.3.10}, we have that

$$ {\cal K} = \{  R \mid R \hbox { is an $H$-simple module algebra with
unit }\}.$$ Therefore ${\cal M}(R) = r_{Hbm}(R).$
\begin{picture}(8,8)\put(0,0){\line(0,1){8}}\put(8,8){\line(0,-1){8}}\put(0,0){\line(1,0){8}}\put(8,8){\line(-1,0){8}}\end{picture}

  Assume that
$H$ is a finite-dimensional semisimple Hopf
  algebra with
 $t \in \int_{H}^{l} $ and $\epsilon(t) = 1$. Let
 $$G_{t}(a) =
  \{z \mid z = x + (t.a)x + \sum(x_{i}(t.a)y_{i} + x_{i}y_{i})
  \hbox { \ for all \ } x_{i}, y_{i}, x \in R \}.$$
   $R$ is called an $r_{gt}$-$H$-module algebra, if $a \in G_{t}(a)$ for all
    $a \in R$.

\begin{Theorem}\label {9.3.12}
$r_{gr}$ is an $H$-radical property of $H$-module algebra and
$r_{gt} =r_{Hbm}$.
\end{Theorem}

{\bf Proof.}  It is clear that any $H$-homomorphic image of
$r_{gt}$-$H$- module algebra is an $r_{gt}$-$H$-module algebra. Let
$$N= \sum \{I \lhd _H \mid I \hbox { is an }  r_ {gt} \hbox {-} H
\hbox {-ideal of } R\}.$$ Now we show that $N$ is an
$r_{gt}$-$H$-ideal of $R$. In fact, we only need to show  that $I_1
+I_2$ is an $r_{gt}$-$H$-ideal  for any two $r_{gt}$-$H$-ideals
$I_1$ and $I_2$. For any $a \in I_1, b \in I_2$, there exist $x,
x_i, y_i \in R$ such that
$$a = x + (t \cdot a)x + \sum_i (x_i (t \cdot a)y_i + x_iy_i).$$
Let $$c = x + (t \cdot (a +b))x +   \sum x_i (t \cdot (a +b))y_i +
x_iy_i  \hbox { \ \ \ } \in G_t(a +b).$$ Obviously, $$ a +b -c = b -
(t \cdot b)x - \sum x_i (t\cdot b)y _i \hbox { \ \  \ } \in I_2.$$
Thus there exist $w, u_j, v_j \in R$ such that
$$a+b-c= w + (t \cdot (a+b-c))w + \sum _j (u_j(t \cdot (a+b-c))v_j + u_jv_j).$$
Let $d= (t\cdot (a +b))w +w + \sum _j (u_j(t \cdot (a +b))v_j
+u_jv_j)$ and $e = c - \sum _j u_j(t \cdot c)v_j - (t \cdot c)w.$ By
computation, we have that $$a +b =d +e.$$ Since $c \in G_t (a +b)$
and $d \in G_t(a+b)$, we get that $e \in G_(a +b)$ and $a +b\in
G_t(a+b),$ which implies that $I_1 +I_2$
 is an $r_{gt}$-$H$-ideal.

Let $\bar R = R/N$ and $\bar B $ be an $r_{gt}$-$H$-ideal of $\bar
R.$ For any $a \in B$, there exist $ x , x _i , y_i \in R$ such that
$$\bar a = \bar x+ (t \cdot \bar a) \bar x + \sum (\bar x_i (t \cdot
\bar a ) \bar y_i + \bar x_i \bar  y_i)$$ and  $$x + (t \cdot  a)  x
+ \sum ( x_i (t \cdot  a )  y_i +  x_i   y_i) -a \in N. $$ Let $$c =
x + (t \cdot a)x + \sum (x_i (t \cdot a) y _i + x_i y_i) \in
G_t(a).$$ Thus there exist $w, u_j, v_j \in R$ such that
$$a -c = (t\cdot (a-c))w +w + \sum (u_j (t \cdot (a-c))v_j + u_jv_j) $$ and
$$a = (t\cdot a)w +w + \sum u_j (t \cdot a)v_j + u_jv_j +c
-(t \cdot c)w - \sum u_j(t \cdot c)v_j \hbox { \ \  \ }\in G_t(a)
,$$ which implies that $B$ is an $r_{gt}$-$H$-ideal and $\bar B=0.$
Therefore $r_{gt}$ is an $H$ -radical property.
\begin{picture}(8,8)\put(0,0){\line(0,1){8}}\put(8,8){\line(0,-1){8}}\put(0,0){\line(1,0){8}}\put(8,8){\line(-1,0){8}}\end{picture}

\begin{Proposition}\label {9.3.13}
If $R$ is an $H$-simple module algebra, then $r_{gr}(R) =0$ iff $R$
has a unit.
\end{Proposition}
{\bf Proof.} If $R$ is an $H$-simple module algebra with unit $1$,
then $-1 \not\in G_t(-1)$ since
$$x + (t \cdot (-1))x + \sum (x_i (t \cdot (-1))y_i + x_iy_i )=0 $$
for any $x, x_i, y_i \in R$. Thus $R$ is $r_{gt}$-$H$-semisimple.
Conversely, if $r_{gt}(R)=0,$ then there exists $0 \not=a \not\in
G_t(a)$ and $G_t(a)=0,$ which implies that $ax +x=0$  for any $x \in
R.$  It follows from Lemma {9.3.3} (7) that $R$ has a unit.
\begin{picture}(8,8)\put(0,0){\line(0,1){8}}\put(8,8){\line(0,-1){8}}\put(0,0){\line(1,0){8}}\put(8,8){\line(-1,0){8}}\end{picture}

\begin{Theorem}\label {9.3.14}
 $r_{gt} =r_{Hbm}$.
 \end{Theorem}
{\bf Proof.}  By Proposition \ref {9.3.13}, $r_{gt}(R) \subseteq
r_{Hbm}(R) $ for any $H$-module algebra $R.$  It  remains
 to show that if $a \not\in r_{gt}(R)$ then
$a \not\in r_{Hbm}(R) $. Obviously, there exists $b \in (a)$ such
that $b \not\in G_t(b),$  where $(a)$  denotes the $H$-ideal
generated by $a$ in $R.$
 Let $${\cal E} = \{ I \lhd _H R \mid G_t(b) \subseteq I , b \not\in I \}.$$
 By Zorn's Lemma, there exists a maximal element $P$ in ${ \cal E }$. $P$ is a
 maximal $H$-ideal of $R$, for, if $Q$  is an $H$-ideal of $R$  with
 $P \subseteq Q$ and $P\not= Q,$  then $b \in Q$   and
 $x = -bx + (bx +x) \in Q$  for any $x \in R.$ Consequently,
 $R/P$ is an $H$-simple module algebra with $r_{gt}(R/P)=0$.
 It follows from Proposition \ref {9.3.13}  that
 $R/P$  is an $H$-simple module algebra with unit and $r_{Hbm}(R) \subseteq P.$
 Therefore  $b \not\in r_{Hbm}(R)$  and so $a \not\in r_{Hbm}(R).$
 \begin{picture}(8,8)\put(0,0){\line(0,1){8}}\put(8,8){\line(0,-1){8}}\put(0,0){\line(1,0){8}}\put(8,8){\line(-1,0){8}}\end{picture}

                                 \begin{Definition} \label {9.3.15}
    Let $I$ be an $H$-ideal of $H$-module algebra $R$, $N$ be an $R$-$H$-submodule
    of $R$-$H$-module $M$. $(N,I)$ are said to have ``L-condition'',
    if for any finite subset $F \subseteq I$, there exists a positive integer $k$
    such that $F^{k}N = 0$.
 \end{Definition}
  \begin{Definition} \label {9.3.16}
  An $R$-$H$-module $M$ is called an $R$-$H$-$L$-module, if for
   $ M$ the following conditions are fulfilled:

    (i)  $RM \not= 0$.

    (ii)  For every non-zero $R$-$H$-submodule $N$ of $M$ and every $H$-ideal $I$ of
    $R$, if $(N,I)$ has ``$L$-condition'', then $I \subseteq (0:M)_R$.
 \end{Definition}

    \begin{Proposition} \label {9.3.17}
If $M$ is an $R$-$H$-$L$-module, then $R/(0:M)_R$ is an
$r_{lH}$-$H$-semisimple and $H$-prime module algebra.
     \end {Proposition}
{\bf Proof.}  If $M$ is an $R$-$H$-$L$-module, let $\bar R =
R/(0:M)_R.$ Obviously, $\bar R $ is $H$-prime.
 If $\bar B$ is an   $r_{lH}$-$H$-ideal of $\bar R$, then
 $(M, B)$  has "$L$-condition" in $R$-$H$-module $M$, since for any
 finite subset $F$ of $B$, there exists a natural number $n$ such that
 $F^n \subseteq (0:M)_R$ and $F^nM=0.$  Consequently,
 $B \subseteq (0:M)_R$  and $\bar R$ is  $r_{lH}$-semisimple.
 \begin{picture}(8,8)\put(0,0){\line(0,1){8}}\put(8,8){\line(0,-1){8}}\put(0,0){\line(1,0){8}}\put(8,8){\line(-1,0){8}}\end{picture}

    \begin{Proposition} \label {9.3.18}
    $R$ is a non-zero $r_{lH}$-$H$-semisimple and $H$-prime module algebra iff
    there exists a faithful $R$-$H$-$L$-module.
    \end {Proposition}
 {\bf Proof.}  If    $R$ is a non-zero $r_{lH}$-$H$-semisimple and
 $H$-prime module algebra, let $M = R$. Since $R$ is an $H$-prime module
 algebra, $(0:M)_R=0.$
 If $(N,B)$ has "$L$-condition" for non-zero $R$-$H$-submodule of $M$ and
 $H$-ideal $B$, then, for any finite subset $F$ of $B$, there exists an natural number
 $n$, such that  $F^n N =0$ and $F^n (NR)=0$, which implies that
 $F^n =0$ and $B$ is an $r_{lH}$-$H$-ideal, i.e.
 $B=0 \subseteq (0:M)_R$. Consequently,  $M$ is a faithful $R$-$H$-$L$-
 module.

 Conversely, if $M$ is a faithful $R$-$H$-$L$-module, then $R$ is an $H$-prime
 module algebra. If $I$ is an $r_{lH}$-$H$-ideal of $R$, then
 $(M,I)$ has ``$L$-condition", which implies $I=0 $  and
 $R$ is an $r_{lH}$-$H$-semisimple module algebra.
\begin{picture}(8,8)\put(0,0){\line(0,1){8}}\put(8,8){\line(0,-1){8}}\put(0,0){\line(1,0){8}}\put(8,8){\line(-1,0){8}}\end{picture}

 \begin{Theorem} \label {9.3.19}
 Let ${\cal M}_{R}$ =
    \{ $M \mid M$ is an $R$-$H$-$L$-module\} for any $H$-module algebra
    $R$ and ${\cal M} = \cup {\cal M}_{R}$. Then
    ${\cal M}$ is an $H$-special class of modules
   and      ${\cal M} (R) = r_{Hl}(R),$
   where ${\cal K} = \{ R \mid R \hbox { is an } H \hbox{-prime module
 algebra with }  r_{lH}(R)=0 \}$  and
  $r_{Hl}= r^{\cal K}$.
  \end {Theorem}

{\bf Proof.}          Obviously, $(M1)$  holds. To show  that $(M2)$
holds, we only need to show that if $I$ is an $H$-ideal of $R$ and
$M \in {\cal M}_I,$ then $IM \in {\cal M}_R.$ By Lemma \ref
{9.3.3}(5), $IM$ is an $R$-$H$-prime module. If $(N,B)$  has  the
"$L$-condition"  for non-zero  $R$-$H$-submodule $N$ of $IM$ and
$H$-ideal $B$ of $R$, i.e. for any finite subset $F$  of $B$, there
exists a natural number $n$ such that $F^n N=0$, then $(N, BI)$ has
"$L$-condition" in $I$-$H$-module $M$. Thus $BI \subseteq (0:M)_I =
(0:IM)_R \cap I$. Considering  $(0:IM)_R$  is an $H$-prime ideal of
$R$, we have that $B \subseteq (0:IM)_R$ or $I \subseteq (0:IM)_R.$
If $I \subseteq (0:IM)_R$, then $I^2 \subseteq (0:M)_I$ and $I
\subseteq (0:M)_I,$ which contradicts $IM \not=0.$ Therefore $B
\subseteq (0:IM)_R$ and so $IM$   is an $R$-$H$-$L$- module.

To show that $(M3)$ holds, we only need to show that if $M \in {\cal
M}_R$ and $I \lhd _H R$ with $IM \not=0$, then $M\in {\cal M}_I.$ By
Lemma \ref {9.3.3}(6), $M$ is an $I$-$H$-prime module. If $(N,B)$
has  the "$L$-condition"  for non-zero  $I$-$H$-submodule $N$ of $M$
and $H$-ideal $B$ of $I,$ then
 $IN$ is an $R$-$H$-prime module and $(IN, (B))$ has "$L$-condition"
in $R$-$H$-module $M,$ since for any finite subset $F$  of $(B)$,
$F^3 \subseteq B$ and  there exists a natural number $n$ such that
$F^{3n}I N\subseteq F^{3n}N=0$, where $(B)$ is the $H$-ideal
generated by $B$ in $R$. Therefore, $(B) \subseteq (0:M)_R$ and $B
\subseteq (0:M)_I,$  which implies $M \in {\cal M}_I$.

Finally, we show that $(M4)$ holds.  Let $I \lhd _H R$ and $\bar R =
R/I$. If $M \in {\cal M}_R$ and $I \subseteq (0:M)_R,$ then $M$ is
an $\bar R$-$H$- prime module. If $(N, \bar B)$ has "$L$-condition"
for $H$-ideal $\bar B$ of $\bar R$ and $\bar R$-$H$-submodule $N$ of
$M$, then subset $F\subseteq B$ and there exists a natural number
$n$ such that $F^nN = (\bar F)^nN =0.$   Consequently, $M\in {\cal
M}_{\bar R}$. Conversely, if $M \in {\cal M}_{\bar R}$, we can
similarly show that $M \in {\cal M}_R.$

The second claim follows from Proposition \ref {9.3.18}  and Theorem
\ref {9.3.4}(1).
\begin{picture}(8,8)\put(0,0){\line(0,1){8}}\put(8,8){\line(0,-1){8}}\put(0,0){\line(1,0){8}}\put(8,8){\line(-1,0){8}}\end{picture}

\begin {Theorem}\label {9.3.20}
$r_{Hl}= r_{lH}$.
\end {Theorem}

{\bf Proof.} Obviously, $r_{lH} \le r_{Hl}.$ It remains to show that
$r_{Hl}(R) \not= R$ if $r_{lH}(R) \not=R.$ There exists a finite
subset $F$ of $R$ such that $F^n \not=0$    for  any natural number
$n$. Let $${\cal F} = \{ I \mid I \hbox { is an } H \hbox {-ideal of
} R \hbox { with } F^n \not\subseteq I \hbox { for any natural
number } n \}.$$ By Zorn's lemma, there exists a maximal element $P$
in ${\cal F}$. It is clear that $P$ is an $H$-prime ideal of $R.$
Now we show that $r_{lH}(R/P)=0$. If $0\not=B/P$ is an $H$-ideal of
$R/P$, then there exists a natural number $m$ such that $F^m
\subseteq B$. Since $(F^m + P)^n \not=0+P$ for any natural number
$n$, we have that $B/P$ is not locally nilpotent and
$r_{lH}(R/P)=0.$  Consequently, $r_{Hl}(R) \not=R$.
\begin{picture}(8,8)\put(0,0){\line(0,1){8}}\put(8,8){\line(0,-1){8}}\put(0,0){\line(1,0){8}}\put(8,8){\line(-1,0){8}}\end{picture}

In fact, all of the results hold in braided tensor categories
determined by (co)quasitriangular structure.

\chapter { The $H$-Radicals of Twisted $H$-Module Algebras }\label {c11}
Remark \footnote {This chapter can be omitted }

       J.R. Fisher \cite{Fi75} built up the general theory of $H$-radicals for
$H$-module algebras. He studied $H$-Jacobson radical and  obtained
\begin {eqnarray}
  r_{j}(R\#H) \cap R  = r_{Hj}(R)
 \label {e (1)}
 \end {eqnarray}
\noindent for any irreducible Hopf algebra $H$(\cite [Theorem
4]{Fi75}).
 J.R. Fisher \cite{Fi75} asked when is
 \begin {eqnarray}
 r_{j}(R\#H) =r_{Hj}(R) \#H
\label {e (2)}
\end {eqnarray}
\noindent and asked if
\begin {eqnarray}
 r_{j}(R\#H) \subseteq (r_{j}(R):H) \#H
\label {e (3)}
\end {eqnarray}
R.J. Blattner, M. Cohen and S. Montgomery in \cite {BCM86} asked
whether  $R\#_{\sigma}H$ is semiprime with a finite-dimensional
 semisimple Hopf algebra $H$ when $R$ is semiprime,
 which is called the semiprime problem.

   If $H$ is a finite-dimensional semisimple Hopf algebra and $R$ is semiprime,
    then $R\#_{\sigma}H$ is semiprime
in  the following five cases:

(i)  $k$ is a perfect field and $H$ is cocommutative;

(ii)  $H$ is irreducible cocommutative;

(iii)  The weak action of $H$ on $R$ is inner;

 (iv) $H=(kG)^*$, where $G$ is a finite group;

(v) $H$ is cocommutative.

 Part (i) (ii) are  due to W. Chin
 \cite[Theorem 2, Corollary 1]{Ch92}.
Part (iii) is due to B.J. Blattner and S. Montgomery \cite [Theorem
2.7]{BM89}. Part (iv) is due to M. Cohen and S. Montgomery
\cite[Theorem 2.9] {CM84b}. Part (v) is due to S. Montgomery and
H.J. Schneider \cite [Corollary 7.13]{MS95}.

         If $H= (kG)^*$, then  relation (\ref {e (2)}) holds, due to
         M. Cohen and
 S. Montgomery \cite[Theorem 4.1] {CM84b}

In  this chapter   we obtain      the relation between $H$-radical
of $H$-module algebra $R$ and radical of $R \# H$. We give
 some sufficient conditions for  (\ref {e (2)})
and (\ref {e (3)}) and  the formulae, which are similar to (\ref {e
(1)}),
 (\ref {e (2)})  and (\ref {e (3)}) for $H$-prime radical respectively.
 We show that (\ref {e (1)}) holds for any
 Hopf algebra $H.$
 Using radical theory and the conclusions in \cite {MS95},
 we also obtain
that    if $H$ is a finite-dimensional semisimple, cosemisimle and
either commutative
   or cocommutative Hopf algebra, then $R$ is $H$-semiprime iff
      $R$ is semiprime iff
$R\#_{\sigma}H$ is semiprime.

In this chapter, unless otherwise stated, let $k$ be a field, $R$ be
an algebra with unit over $k$,
 $H$ be a Hopf algebra over $k$ and $H^*$ denote the dual space of $H$.

 $R$ is called a twisted $H$-module algebra if
 the following
conditions are satisfied:

   (i) $H$ weakly acts on  $R$;

(ii) $R$ is a twisted $H$-module, that is,
 there exists a linear map
 $\sigma \in Hom_k (H \otimes H, R) $ such that
$h \cdot (k \cdot r) = \sum \sigma (h_1,k_1)(h_2k_2 \cdot r) \sigma
^{-1}(h_3,k_3)$ for all $h, k \in H$ and $r \in R$.

  It is clear that if $\sigma $  is trivial, then
   twisted $H$-module algebra $R$ is an $H$-module algebra.

Set
$$ Spec(R) = \{ I \mid I \hbox { \ is a prime ideal of \ } R \};   $$
 $$H\hbox {-}Spec(R) =
\{ I \mid I \hbox { \ is an \ } H\hbox{-prime ideal of } \ R \}.$$

 \section {The Baer radical of twisted  H-module algebras }\label {s24}

In this section, let $k$ be a commutative associative ring with
unit, $H$ be an algebra
 with unit and comultiplication $\bigtriangleup$,
  $R$ be an algebra
over $k$ ($R$ may be without unit) and
 $R$ be  a twisted  $H$-module algebra.

\begin {Definition} \label {10.1.1}
$r_{Hb}(R) := \cap \{I \mid I$ is  an $H$-semiprime ideal of $R$ \};

$r_{bH}(R):= (r_b(R):H)$

$r_{Hb}(R)$ is called the $H$-Baer radical ( or $H$-prime radical )
of twisted $H$-module algebra $R$.

\end {Definition}

\begin {Lemma} \label {10.1.3}
(1) If $E$ is a non-empty subset of $R$,
  then  $(E) = (H \cdot E) + R(H \cdot E) + (H \cdot E)R + R(H \cdot E)R$,
   where $(E)$ denotes the $H$-ideal generated by $E$ in  $R$;

(2) If $I$ is a nilpotent $H$-ideal of $R$, then $I \subseteq
r_{Hb}(R)$.
\end {Lemma}

  {\bf Proof}.
(1) It is trivial.

(2)  If $I$ is a nilpotent $H$-ideal and $P$ is an $H$-semiprime
ideal,
                      then  $(I +P)/P$ is nilpotent simply because
$ (I +P)/P \cong I/(I \cap P) $  \ \ \ \ (as \ algebras) . Thus  $I
\subseteq P$ and  $I \subseteq r_{Hb}(R)$.
\begin{picture}(8,8)\put(0,0){\line(0,1){8}}\put(8,8){\line(0,-1){8}}\put(0,0){\line(1,0){8}}\put(8,8){\line(-1,0){8}}\end{picture}

\begin {Proposition}\label {10.1.4}

(1) $r_{Hb}(R)=0$ iff $R$ is $H$-semiprime;

(2) $r_{Hb}(R/r_{Hb}(R))=0$;

(3)  $R$ is $H$-semiprime iff $(H\cdot a)R(H\cdot a) = 0$ always
implies
    $a = 0$ for any $a\in R$;

  $R$ is $H$-prime iff $(H\cdot a)R(H\cdot b) = 0$ always implies
 $a = 0$  or $b = 0$ for any $a$, $b \in R$;

(4) If $R$ is $H$-semiprime, then $W_H(R)=0$.
\end {Proposition}

  {\bf Proof}.
(1) If $r_{Hb}(R)=0$, then $R$ is $H$-semiprime by Lemma \ref
{10.1.3} (2). Conversely, if $R$ is $H$-semiprime, then $0$ is an
$H$-semiprime ideal and so $r_{Hb}(R)=0$ by Definition \ref
{10.1.1}.

(2) If $B/r_{Hb}(R) $ is a nilpotent $H$-ideal of $R/r_{Hb}(R)$,
then $B^k \subseteq r_{Hb}(R)$  for some natural number $k$ and so
$B \subseteq r_{Hb}(R)$, which implies that $R/r_{Hb}(R)$  is
$H$-semiprime. Thus $r_{Hb}(R/r_{Hb}(R))=0$ by part (1).

(3) If $R$ is  $H$-prime and $(H\cdot a)R(H\cdot b) = 0$ for $a$ and
$b \in R$, then $(a)^2 (b)^2 = 0$ by Lemma \ref {10.1.3} (1), where
$(a)$ and $(b)$ are the $H$-ideals generated by  $a$ and $b$  in $R$
respectively. Since $R$ is $H$-prime, $a = 0$ or $b= 0$. Conversely,
if both $B$ and $C$ are  $H$-ideals of $R$ and $ BC = 0$, then $(H
\cdot a)R(H\cdot b) = 0$  and $a = 0$ or $b = 0$
 for any $ a \in B$ and $b \in C$, which implies  that $B = 0 $ or $C = 0$.
 Thus $R$ is an $H$-prime.
Similarly, the other assertion holds.

(4)    For any $0 \not= a \in R$, there exist
 $b_1 \in R$ and  $h_1, h_1' \in  H$  such that
 $ 0 \not= a_2 = (h_1 \cdot a_1) b_1(h_1' \cdot a_1)  \in
 (H \cdot a_1)R(H \cdot a_1)$ by part (3), where $a_1 = a$.
 Similarly, for $0 \not= a_2 \in R$, there exist
 $b_2 \in A$, $h_2$, $h_2' \in  H$  such that
 $ 0 \not= a_3 = (h_2 \cdot a_2) b_2(h_2' \cdot a_2)  \in
 (H \cdot a_2)R(H \cdot a_2)$, which implies that
 there exists an $H$-$m$-sequence $\{ a_n \}$  such that
  $a_n \not= 0$ for any natural number $n$.
   Thus $W_{H}(R)=0$.  \begin{picture}(8,8)\put(0,0){\line(0,1){8}}\put(8,8){\line(0,-1){8}}\put(0,0){\line(1,0){8}}\put(8,8){\line(-1,0){8}}\end{picture}

\begin {Theorem} \label {10.1.5}
$r_{Hb}(R)= W_H(R) = \cap  \{ I \mid I { \ is \ an \ } H{
\hbox{-}prime \ ideal \ of \  } R    \} $.
\end {Theorem}

{\bf Proof.} Let  $D = \cap \{ I \mid I $ is an $H$-prime ideal of
$R$ \}. Obviously, $r_{Hb}(R) \subseteq D$.

    If $0 \not= a \not\in W_H(R)$, then there exists
    an $m$-sequence $\{ a_i \}$  in $R$ with
    $a_1=a$ and $a_{n+1} = (h_n \cdot a_n )b_n (h'_n \cdot a_n) \not= 0$
    for $n = 1, 2, \cdots $.
Let  ${\cal F} = \{ I
  \mid I$ is an $H$-ideal of $R$ and $I \cap \{ a_1, a_2, \cdots \}
  = \emptyset \}$.
 By Zorn's Lemma, there exists a
maximal $P$  in ${\cal F}$. If both  $I$ and $J$ are $H$-ideals of
$R$ with  $I \not\subseteq P$ and
 $J \not\subseteq P$
 such that $IJ \subseteq P$, then there exist natural numbers
 $n$  and $m$ such that  $a_n \in I +P$  and
 $a_m \in J+P$. Since $ a_{n +m + 1}
 = (h_{n+m}\cdot a_{n+m})b_{n+m}(h'_{n+m}\cdot a_{n+m}) \in (I+P)(J+P)
 \subseteq P$, we get a contradiction. Thus
 $P$ is an $H$-prime ideal of $R$. Obviously, $a \not\in P$, which implies
 that
  $a \not\in D$.
Therefore $D \subseteq W_H(R)$.

    For any $x \in W_H(R)$,
let $\bar R = R/r_{Hb}(R)$. It follows from
 Proposition \ref {10.1.4} (1) (2) (4) that   $W_H(\bar R)=0$.
 For  any  $H$-$m$-sequence $\{ \bar a_n \}$
 with $\bar a_1 = \bar x$ in  $\bar R$,  there exist
 $\overline b_n \in \overline R$ and $h_n, h_n' \in H$ such that
 $ \overline a_{n+1} =
 (h_n \cdot \overline a_n) \overline b_n (h_n' \cdot \overline a_n)$
                           {~}{~}  for any natural number $n$.
Let   $a_1'=x $  and
 $  a_{n+1}' = (h_n \cdot  a_n')  b_n (h_n' \cdot  a_n')$ {~}{~}
                             for any natural number $n$.
 Since $ \{ a_n'\}$ is an $H$-$m$-sequence with $a_1' = x$ in $R$,
 there exists a natural number $k$ such that $a_k' =0$. It is
 clear that $\overline a_n= \overline {a_n'}$
  for any natural number $n$ by induction. Thus $ \bar a_k = 0$ and
  $\overline x \in W_H(\overline R)$.
  Considering $W_H(\overline R)=0$, we have $x \in r_{Hb}(R),$
  which implies that
  $W_H(R) \subseteq r_{Hb}(R)$. Therefore
  $W_H(R) = r_{Hb}(R)=D$. \begin{picture}(8,8)\put(0,0){\line(0,1){8}}\put(8,8){\line(0,-1){8}}\put(0,0){\line(1,0){8}}\put(8,8){\line(-1,0){8}}\end{picture}

  \section{ The Baer and Jacobson radicals of crossed products }\label {s25}

By \cite [Lemma 7.1.2] {Mo93}, if $R\#_\sigma H$  is crossed product
defined in \cite [Definition 7.1.1] {Mo93}, then $R$ is a twisted
$H$-module algebra.

Let $R$ be an algebra and $M_{m \times n }(R)$ be the algebra of $m
\times n$ matrices with entries in $R$. For  $ i = 1, 2, \cdots m $
and  $ j = 1, 2, \cdots n $.  Let $(e_{ij})_{m \times n}$ denote the
matrix in $M_{m \times n}(R),$ where $(i,j)$-entry is $1_R$  and the
others are zero. Set

 ${ \cal I}(R) = \{ I \mid I \hbox { is an ideal of } R \}$.

$r_{Hj}(R) := r_j(R\#_{\sigma} H) \cap R$;

$r_{jH}(R):=(r_{j}(R):H)$.

\begin {Lemma}\label{10.2.2}
Let $M_R$ be a free $R$-module with finite rank  and $R'= End
(M_R)$. Then there exists a  unique bijective  map
 $$\Phi: {\cal I}(R) \longrightarrow {\cal I}(R')$$
such that \ \ \  $\Phi (I)M=MI$ and

 (1) $\Phi $ is a map  preserving containments, finite products,
and infinite intersections;

(2) $ \Phi (I) \cong   M_{n \times n}(I)$ for any ideal $I$ of $R$;

(3) $I$ is a (semi)prime ideal of $R$ iff $\Phi(I)$ is
 a (semi)prime ideal of $R'$;

(4) $r_b(R')= \Phi(r_b(R))$;

(5) $r_j(R')= \Phi(r_j(R))$.
\end {Lemma}

{\bf Proof} Since $M_R$ is a free $R$-module with rank $n$, we can
assume $M=M_{n \times 1}(R)$. Thus $R' = End (M_R)= M_{n \times
n}(R)$ and the module operation of $M$ over $R$ becomes the matrix
operation. Set $M'= M_{1 \times n}(R)$. Obviously, $M'M=R$. Since
$(e_{i1})_{n \times 1} (e_{1j})_{1 \times n} = (e_{ij})_{n \times
n}$  for $i, j = 1, 2, \cdots n$,  $MM' = M_{n \times n}(R) =R'$.
Define
$$\Phi (I) = MIM'$$
\noindent for any ideal $I$ of $R$. By simple computation, we have
that $\Phi (I)$ is an ideal of $R'$ and $\Phi (I)M = MI$. If $J$ is
an ideal of $R'$ such that $JM=MI$, then $JM= \Phi (I)M$ and $J=
\Phi(I)$, which implies $\Phi $ is unique. In order to show that
$\Phi$ is a bijection from ${ \cal I}(R)$ onto ${ \cal I}(R')$, we
define  a map  $\Psi$ from ${\cal I}(R')$ to ${ \cal I}(R)$ sending
$I'$ to $M'I'M$ for any ideal $I'$ of $R'$.  Since $ \Phi \Psi (I')=
MM'I'MM' = I'$ and $\Psi \Phi (I)= M'MIM'M = I$ for any ideal $I'$
of $R'$ and ideal $I$ of
 $R$,  we have that $\Phi$ is bijective.

(1) Obviously $\Phi$ preserves containments. We see that
 $$\Phi (IJ) = MIJM'= (MIM')(MJM') = \Phi (I) \Phi(J)$$
  for any ideals $I$
and $J$ of $R$. Thus $\Phi$ preserves finite products.
 To show that $\Phi$ preserves infinite
intersections, we first show that
\begin {eqnarray} M (\cap \{ I_\alpha \mid \alpha \in \Omega \}) =
 \cap \{ MI_\alpha \mid \alpha \in \Omega \}  \label {e1.2 (1)}
 \end {eqnarray}
for any  $\{I_\alpha \mid \alpha \in \Omega \} \subseteq {\cal
I}(R)$.
 Obviously, the right side of  relation (\ref {e1.2 (1)})   contains the left side of
  relation (\ref {e1.2 (1)}).
 Let \{$u_1, u_2, \cdots, u_n $ \}  be a basis of $M$ over $R$.
 For any $x \in  \cap \{ MI_{\alpha} \mid \alpha \in \Omega \}$, any
 $\alpha, \alpha' \in \Omega $,
  there exist
 $r_i \in I_\alpha$ and $r_i' \in I_{\alpha'}$ such that
 $x = \sum u_ir_i= \sum u_ir_i'$.
 Since $\{u_i\}$ is a basis, $r_i = r_i'$,  which implies  $x \in
 M (\cap \{ I_\alpha \mid \alpha \in \Omega \})$.
 Thus the  relation (\ref {e1.2 (1)})
 holds.  It follows from  relation (\ref {e1.2 (1)})  that
 $$\Phi(\cap \{ I_\alpha \mid \alpha \in \Omega \})M =
 \cap \{ \Phi(I_\alpha)M \mid \alpha
 \in \Omega \}$$
Since $  \Phi (\cap  \{ I_\alpha  \mid \alpha \in \Omega \})M =
  \cap  \{ \Phi (I_\alpha )M \mid \alpha \in \Omega \}
\supseteq
 (\cap \{ \Phi(I_\alpha) \mid \alpha
 \in \Omega \})M$, we have that
   $$\Phi(\cap \{ I_\alpha \mid \alpha \in \Omega \})
 \supseteq \cap \{ \Phi(I_\alpha) \mid \alpha \in \Omega \}.$$
Obviously,  $$\Phi(\cap \{ I_\alpha \mid \alpha \in \Omega \})
 \subseteq \cap \{ \Phi(I_\alpha) \mid \alpha \in \Omega \}.$$
 Thus           $$\Phi(\cap \{ I_\alpha \mid \alpha \in \Omega \})
 = \cap \{ \Phi(I_\alpha) \mid \alpha \in \Omega \}.$$

(2) Obviously, $ \Phi (I) = MIM' = M_{n \times 1}(R) I M_{1 \times
n}(R) \subseteq M_{n \times n}(I)$. Since $ a(e_{ij})_{n \times n} =
(e_{i 1})_{n \times 1}a (e_{1j})_{1 \times n} \in MIM'$ for all $a
\in I$ and $i, j = 1, 2, \cdots n$,
$$ \Phi (I) = MIM' =
M_{n \times 1}(R) I M_{1 \times n}(R) \supseteq M_{n \times n}(I).$$
Thus part (2) holds.

(3) Since bijection $\Phi$ preserves  products,  part (3) holds.

(4) We see that
\begin {eqnarray*} \Phi (r_b(R)) &=&
\Phi (\cap \{ I \mid I \hbox { is a prime ideal of }  R \} )  \\
&=& \cap \{ \Phi (I) \mid I \hbox { is a prime ideal of }  R \}
\hbox { by part (1) } \\
&=& \cap \{ \Phi (I) \mid \Phi (I) \hbox { is a prime ideal of }  R'
\}
\hbox { by part (3) }   \\
&=& \cap \{ I' \mid I' \hbox { is a prime ideal of }  R' \}
\hbox { since  }  \Phi \hbox { is surjective } \\
&=& r_b(R')
\end {eqnarray*}

(5) We see that
\begin {eqnarray*} \Phi (r_j(R)) &=&
M_{n \times n}(r_j(R)) \hbox { by part (2) }     \\
&=& r_j(M_{n \times n}(R)) \hbox { by \cite [Theorem 30.1] {Sz82} } \\
&=& r_j(R')  { \ . \ \ \ \ \ } \Box
\end {eqnarray*}

     Let $H$ be a finite-dimensional   Hopf algebra and $A = R\#_{\sigma} H$.
Then $A$ is a free right $R$-module  with finite rank by \cite
[Proposition 7.2.11]{Mo93} and $End (A_R) \cong (R \#_{\sigma}
H)\#H^*$ by \cite [Corollary 9.4.17]{Mo93}. By part (a) in the proof
of \cite [Theorem 7.2] {MS95}, it follows that $\Phi$ in Lemma 1.2
is the same as  in \cite [Theorem 7.2] {MS95}.

 \begin {Lemma} \label {10.2.3} Let $H$ be a finite-dimensional
      Hopf algebra and $A = R\#_{\sigma} H$. Then

(1) If $P$ is an $H^*$-ideal of $A$, then $P=(P\cap R) \#_\sigma H$.

(2) $ \Phi (I) = (I \#_{\sigma}H)\#H^*$ for  every  $H$-ideal $I$ of
$R$;

(3)         \begin {eqnarray}
   \{ P \mid P \hbox { is an } H \hbox {-ideal of } A \#H^* \}
&=& \{ (I \#_{\sigma}H)\#H^* \mid I \hbox { is an \ } H\hbox {-ideal
of } R \}
 \label {e1.3 (1)}
 \end {eqnarray}
 \begin {eqnarray}
 \{ P \mid P \hbox { is an } H^*\hbox {-ideal of } A  \}
 &=& \{ I \#_{\sigma }H \mid I \hbox { is an } H\hbox {-ideal of } R \}
\label {e1.3 (2)}
\end {eqnarray}
$ \{ P \mid P \hbox { is an } H\hbox{-prime ideal of } A \#H^* \}$
           \begin {eqnarray}
&=& \{ (I \#_{\sigma}H)\#H^* \mid I \hbox { is an } H\hbox{-prime
ideal of } R \} \label {e1.3 (3)}
\end {eqnarray}

(4)    $H$-Spec$(R) = \{ (I:H) \mid I  \in$ Spec $(R) \}$;

 (5) \begin {eqnarray}    (\cap  \{ I_{\alpha } \mid  \alpha \in \Omega \} : H)
 &=&
\cap \{ (I_{\alpha } : H) \mid \alpha \in \Omega \} \label {e1.3
(4)}
\end {eqnarray}
 where $I_{\alpha}$  is an ideal
of $R$ for all $\alpha \in \Omega ; $

 (6)  \begin {eqnarray}   ( \cap  \{ I_{\alpha } \mid \alpha \in
 \Omega \} )\#_{\alpha}H &=&
\cap \{ (I_{\alpha } \#_{\sigma } H) \mid \alpha \in \Omega \}
\label {e1.3 (5)}
 \end {eqnarray}
 where $I_{\alpha}$  is an $H$-ideal
of $R$ for all $\alpha \in \Omega $;

(7) $\Phi(r_b(R))= r_b(A\#H^*)$;

(8)  $\Phi(r_j(R))= r_j(A\#H^*)$;

(9) $\Phi(r_{Hb}(R))= r_{Hb}(A\#H^*) = (r_{Hb}(R) \#_\sigma H) \#
H^*$.
\end {Lemma}

{\bf Proof} (1) By  \cite [Corollary 8.3.11] {Mo93}, we have that
 $P = (P \cap R)A = (P\cap R)\#_\sigma H$.

(2) By the part (b) in the proof of \cite [Theorem 7.2]{MS95}, it
follows
 that
 $$\Phi(I)= (I \#_{\sigma}H)\#H^*$$
 for every $H$-ideal $I$ of $R$.

(3)  Obviously, the left side of  relation (\ref {e1.3 (1)})
contains
 the right side of  relation (\ref {e1.3 (1)}). If $P$ is an $H$-ideal of
 $A\#H^*$, then
 $P= (P \cap A) \# H^* = (((P \cap A) \cap R) \#_\sigma H) \# H^*$ by
 part (1), which implies that the right side of  relation (\ref {e1.3 (1)}) contains
 the left side. Thus  relation (\ref {e1.3 (1)}) holds.
 Similarly,  relation (\ref {e1.3 (2)}) holds.
 Now, we show that  relation (\ref {e1.3 (3)}) holds. If $P$
 is an $H$-prime ideal of $A \# H^*$, there exists an $H$-ideal $I$ of $R$
 such that $P = (I \#_\sigma H) \#H^*$ by  relation (\ref {e1.3 (1)}). For any
 $H$-ideals $J$ and $J'$ of $R$ with $JJ' \subseteq I$,
 since $\Phi (JJ') = \Phi (J) \Phi(J') \subseteq \Phi (I) =P$  by Lemma \ref {10.2.2} (1), we have that $\Phi (J) \subseteq \Phi(I)$
or $\Phi(J')\subseteq \Phi(I)$, which implies that $J \subseteq I$
or $J' \subseteq I$ by Lemma \ref {10.2.2}. Thus $I$ is an $H$-prime
ideal of $R$. Conversely, if $I$ is an $H$-prime ideal of $R$ and
$P= (I \#_\sigma H)\#H^*$, we claim that $P$ is an $H$-prime of $A
\#H^*$. For any $H$-ideals $Q$ and $Q'$ of $A \#H^*$ with  $QQ'
\subseteq P$, there exist two $H$-ideals $J$ and $J'$ of $R$ such
that $(J\#_\sigma H)\#H^* = Q$ and $(J'\#_\sigma H)\#H^* = Q'$ by
 relation (\ref {e1.3 (1)}).
Since $\Phi(JJ') = \Phi(J) \Phi(J') = QQ'\subseteq P = \Phi(I) $,
$JJ' \subseteq I$, which implies $J \subseteq I$, or $J' \subseteq
I$, and so $Q \subseteq P$ or $Q' \subseteq P$. Thus $P$ is an
$H$-prime ideal of $A\#H^*$. Consequently,  relation (\ref {e1.3
(3)}) holds.

 (4) It follows from \cite [Lemma 7.3 (1) (2)] {MS95}.

(5) Obviously,  the right side  of relation (\ref {e1.3 (4)})
contains the left side. Conversely, if $x \in \cap \{(I_\alpha:H)
\mid \alpha \in \Omega \}$, then $x \in (I_{\sigma} :H)$ and $h\cdot
x \in I_\alpha $ for all $\alpha \in \Omega, h \in H$, which implies
that $h \cdot x \in \cap \{ I_\alpha \mid \alpha \in \Omega \} $ and
$x \in (\cap \{I_\alpha \mid \alpha \in \Omega \}:H)$. Thus relation
(\ref {e1.3 (4)}) holds.

(6) Let $\{ h^{(1)}, \cdots, h^{(n)} \}$  be a basis of $H$.
Obviously,
 the right side of
 relation (\ref {e1.3 (5)}) contains
 the left side of  relation (\ref {e1.3 (5)}). Conversely,
 for   $ u \in \cap \{ (I_{\alpha } \#_{\sigma } H)
 \mid \alpha \in \Omega \}$ and   $\alpha, \alpha' \in \Omega, $
there exist   $r_{i} \in I_\alpha$ and $r_i' \in I_{\alpha'}$
 such that $u = \sum_{i=1}^{n} r_i \# h^{(i)} =
  \sum_{i=1}^{n} r_i' \# h^{(i)}$.
   Since  $\{ h^{(1)}, \cdots, h^{(n)} \}$  is linearly independent,
   we have that $r_i= r_i'$, which implies that
 $ u \in (\cap  \{ I_{\alpha } \mid \alpha \in \Omega \} ) \#_{\alpha} H$.
Thus  relation (\ref {e1.3 (5)}) holds.

(7) and (8) follow from Lemma \ref{10.2.2}(4)(5).

(9)  We see  that
\begin {eqnarray*}
r_{Hb}(A \#H^*) &=& \cap \{ P \mid P \hbox { is an } H\hbox {-prime
ideal of }
 A \#H^*      \} \hbox { \ \  by Theorem \ref {10.1.5} }\\
 &=& \cap \{ (I \#_\sigma H)\#H^*  \mid I \hbox { is an }
 H\hbox {-prime ideal of }  R       \}  \hbox { \ \ \
 by relation (\ref {e1.3 (3)}) } \\
  &=& ( \cap \{ I \#_\sigma H \mid I \hbox { is an } H\hbox {-prime ideal of }
   R   \} ) \# H^*    \hbox { \ \ \ by  part (6) }  \\
 &=& ((\cap \{ I \mid I \hbox { is an }
 H\hbox {-prime ideal of } R \}) \#_\sigma H) \# H^*
  \hbox { \ \ \ by part (6) }  \\
 &=& (r_{Hb}(R) \#_\sigma H) \# H^*  \hbox { \ \ by Theorem \ref {10.1.5} }\\
&=& \Phi (r_{Hb}(R))    \hbox { \ \ \ by part (2) }.
\end {eqnarray*}

  (10) If $H$ is cosemisimple, then $H$ is semisimple  by \cite [Theorem 2.5.2]
{Mo93}. Conversely, if $H$ is semisimple, then $H^*$ is
cosemisimple. By \cite [Theorem 2.5.2]{Mo93}, $H^*$ is semisimple.
Thus $H$ is cosemisimple.
\begin{picture}(8,8)\put(0,0){\line(0,1){8}}\put(8,8){\line(0,-1){8}}\put(0,0){\line(1,0){8}}\put(8,8){\line(-1,0){8}}\end{picture}

\begin {Proposition}\label {10.2.4}
  (1)  $r_{Hb}(R) \subseteq  r_{b}(R\#_{\sigma}H) \cap R \subseteq r_{bH}(R)$;

  (2)  $ r_{Hb}(R)\#_{\sigma}H \subseteq r_{b}(R \#_{\sigma} H)$.
  \end {Proposition}

  {\bf Proof}.
   (1) If $P$ is a prime ideal of $R \#_{\sigma} H$, then $P \cap R$ is
    an $H$-prime    ideal of $R$ by \cite [Lemma 1.6]{Ch91}.
   Thus  $r_b(R \#_{\sigma} H) \cap R= \cap \{ P\cap R \mid
   P$ is prime ideal of $R \#_{\sigma} H \}  \supseteq r_{Hb}(R)$.
For any $a \in r_b(R \#_\sigma H) \cap R$ and any $m$-sequence $ \{
a_i \} $ in $R$ with $a_1 = a $, it is easy to check that $\{ a_i
\}$ is also an $m$-sequence in $R \#_\sigma H$. Thus $a_n = 0 $  for
some natural $n$, which implies $a \in r_b(R)$. Thus $r_b(R\#_\sigma
H) \cap R \subseteq r_{bH}(R)$ by \cite [Lemma 1.6]{BM89}

(2) We see that
\begin {eqnarray*}
    r_{Hb}(R) \#_{\sigma} H &=&    (r_{Hb}(R) \#_{\sigma} 1)(1 \#_{\sigma} H)\\
    &\subseteq &  r_b(R \#_{\sigma} H) ( 1 \#_\sigma H)
    \hbox { \ \  by  part  (1)}\\
    &\subseteq & r_b(R \#_{\sigma} H). { \ \ \ \ }  \Box
    \end {eqnarray*}

 \begin {Proposition} \label {10.2.5} Let $H$ be finite-dimensional
      Hopf algebra and $A = R\#_{\sigma} H$.  Then

(1) $r_{H^*b}(R\#_{\sigma}H) =  r_{Hb}(R)\#_{\sigma}H$;

(2) $r_{Hb}(R)= r_{bH}(R)= r_b(R\#_{\sigma}H) \cap R$.
\end {Proposition}

{\bf Proof}

(1) We see that
\begin {eqnarray*}
r_{H^*b}(R \#_\sigma H) &=& \cap \{ P \mid P
\hbox { is an } H^*\hbox {-prime ideal of } A \} \\
&=& \cap \{ I \#_\sigma H \mid I \hbox { is an } H \hbox {-prime
ideal of } R \}  \hbox
{ \ \ \ ( by  \cite [Lemma 7.3 (4)] {MS95} ) }\\
&=& ( \cap \{ I \mid I \hbox { is an } H \hbox {-prime ideal of } R
\}) \#_\sigma H
\hbox { \ \ \ ( by Lemma \ref {10.2.3} (6)) }\\
&=& r_{Hb}(R) \#_\sigma H.
\end {eqnarray*}

 (2) We see that
\begin {eqnarray*} r_{Hb}(R) &=& \cap \{ P \mid P { \ is \ an \ }
H \hbox {-prime  ideal  of  } R \} \\
   &=& \cap \{(I:H) \mid I \in { \ Spec} (R) \}
   \hbox { \ by  Lemma \ref {10.2.3} \  part (4) }\\
   &=& ( \cap \{ I \mid I \in {Spec}(R) \}:H  )
 \hbox   { \ by \ Lemma \ref {10.2.3} \  part (5) }  \\
   &=& (r_b(R):H) \\
   &=& r_{bH}(R).
   \end{eqnarray*}

 Thus it follows from Proposition \ref {10.2.4}(1)
   that $r_{Hb}(R)= r_b(R\#_{\sigma}H)\cap R = r_{bH}(R)$.
\begin{picture}(8,8)\put(0,0){\line(0,1){8}}\put(8,8){\line(0,-1){8}}\put(0,0){\line(1,0){8}}\put(8,8){\line(-1,0){8}}\end{picture}

      \begin {Theorem} \label {10.2.6}.
Let $H$ be a finite-dimensional Hopf algebra and  the weak  action
of $H$ be  inner. Then

(1) $r_{Hb}(R)= r_b(R)= r_{bH}(R)$;

 Moreover, if $H$ is semisimple, then

 (2) $r_{b}(R\#_{\sigma} H)= r_{Hb}(R)\#_{\sigma} H$.
 \end {Theorem}

{\bf Proof} (1) Since the weak action is inner, every ideal of $R$
is an $H$-ideal, which implies that
 $r_{Hb}(R)= r_b(R) = r_{bH}(R)$ by Proposition \ref {10.2.5} (2).

(2) Considering  Proposition \ref {10.2.4}(2), it suffices to show
 $r_{b}(R\#_{\sigma} H)\subseteq  r_{Hb}(R)\#_{\sigma} H.$
It is clear that
\begin {eqnarray} (R\#_{\sigma} H)/ (r_{Hb}(R)\#_{\sigma} H)
\cong (R/r_{Hb}(R))\#_{\sigma} H  \hbox { \ \ \ \ (  \ as \ algebras
)}. \label {e1.6 (1)}
\end {eqnarray}
It follows by \cite [Theorem 7.4.7]{Mo93} that $
(R/r_{Hb}(R))\#_{\sigma} H$
 is semiprime.
Therefore $$r_{b}(R \#_{\sigma} H) \subseteq r_{Hb}(R)\#_{\sigma} H.
\ \ \Box$$

   \begin {Theorem}\label {10.2.7} Let $H$ be  a finite-dimensional, semisimple
   and either commutative or cocommutative Hopf
algebra and let $A = R\#_{\sigma} H$. Then

(1) ${r_{b}(R\#_{\sigma}H) =  r_{Hb}(R)\#_{\sigma}H};$

(2) $R$ is $H$ semiprime iff $R \#_\sigma H$ is semiprime.

Moreover, if $H$ is cosemisimple, or char $k$ does not divide dim
$H$, then  both part (3) and part (4) hold:

(3) $r_{Hb}(R)= r_{bH}(R)= r_b(R)$;

(4) $R$ is $H$-semiprime iff $R$ is semiprime iff $R\#_{\sigma}H$ is
semiprime.

\end {Theorem}
{\bf Proof} (1) Considering Proposition \ref {10.2.4}(2), it
suffices to show
 $$r_{b}(R\#_{\sigma} H)\subseteq  r_{Hb}(R)\#_{\sigma} H.$$
It follows by \cite [Theorem 7.12 (3)]{MS95} that $
(R/r_{Hb}(R))\#_{\sigma} H$
 is semiprime. Using relation (\ref {e1.6 (1)}), we have
 that $r_{b}(R \#_{\sigma} H) \subseteq r_{Hb}(R)\#_{\sigma} H$.

(2)  It follows from part (1) and Proposition \ref {10.1.4} (1).

(3)  By \cite [Theorem 4.3 (1)] {La71}, we have that $H$ is
semisimple and cosemisimple.

 We see  that
 \begin{eqnarray*}
\Phi(r_b(R)) &=& r_{b}(A \# H^*) \hbox { \ \ \ by \ Lemma \ref {10.2.3} \ (7)}\\
 &=& r_{H^*b}(A) \# H^*    \hbox { \ \ \  by \  part \  (1) } \\
 &=& (r_{Hb}(R) \#_{\sigma} H) \# H^* \hbox { \ \ \ by \ Proposition \ref {10.2.5} (1) } \\
 &=& \Phi (r_{Hb}(R))  \hbox { \ \ \ by \ Lemma \ \ref {10.2.3} (2). }
 \end {eqnarray*}
 Thus $r_b(R)= r_{Hb}(R).$

(4) It immediately follows from  part (2) and part (3).
\begin{picture}(8,8)\put(0,0){\line(0,1){8}}\put(8,8){\line(0,-1){8}}\put(0,0){\line(1,0){8}}\put(8,8){\line(-1,0){8}}\end{picture}

We now provide an example to show that the Baer radical $r_b(R)$
 of $R$ is not $H$-stable  when $H$ is not cosemisimple.

Example: Let $k$ be a field of characteristic $p > 0$ and $R =
k[x]/(x^p)$. Then we can define a derivation $d$ on $R$ by sending
$x$  to $x+1$. Then $d^2(x)=d(x+1)=d(x)$ and then, by induction,
$d^p(x)=d(x)$. It follows that $d^p=d$ on all of $R$.Thus $H=u(kd)$,
the restricted enveloping  algebra, is semisimple by \cite [Theorem
2.3.3]{Mo93}.
 clearly $H$ acts on $R$, but $H$
does not stabilize the Baer radical of $R$ which is the principal
ideal generated by $x$. Note also that $H$ is commutative  and
cocommutative.

\begin {Proposition}\label {10.2.8}
If $R$ is an $H$-module algebra, then
$$r_{Hj}(R) = \cap \{(0:M)_R \mid M { \ is \ an \ irreducible \  }
R\hbox {-}H\hbox {-} { module} \}.$$ That is, $r_{Hj}(R)$ is the
$H$-Jacobson radical of the $H$-module algebra $R$ defined in \cite
{Fi75}.
\end {Proposition}

{\bf Proof.}  It is easy to show that $M$ is an irreducible
$R$-$H$-module iff $M$ is an irreducible $R\#H$-module by \cite
[Lemma 1]{Fi75}. Thus
\begin {eqnarray*}
r_{Hj}(R) &=& r_{j}(R\#H)\cap R \hbox { \ by \ definition \ref {10.1.1}} \\
&=& (\cap \{ (0:M)_{R\#H} \mid M \hbox { \ is \ an \ irreducible \ }
R\# H \hbox {-module} \} ) \cap R \\
&=& \cap \{ (0:M)_R \mid M  \hbox { \ is \ an \ irreducible \ }
R\hbox {-}H \hbox {-module} \}. \Box
\end {eqnarray*}

\begin {Proposition} \label {10.2.9}

  (1)  $ r_{j}(R\#_{\sigma}H ) \cap R =r_{Hj}(R) \subseteq r_{jH}(R); $

  (2)  $r_{Hj}(R)\#_{\sigma} H \subseteq r_{j}(R \#_\sigma H)$.
  \end {Proposition}

  {\bf Proof}.
  (1)  For any $a \in r_j(R \#_{\sigma} H) \cap R$, there exists
$u = \sum_{i} a_i \# h_i \in R \#_{\sigma} H$ such that
$$ a + u + au = 0. $$
Let $(id \otimes \epsilon )$ act on the above equation. We get that
$a + \sum a_i\epsilon (h_i) + a (\sum a_i \epsilon (h_i))=0$, which
implies that $a$ is a right quasi-regular element in $R$. Thus
$r_j(R \#_{\sigma} H) \cap R \subseteq r_{jH}(R)$.

   (2)  It is similar to the proof of Proposition \ref {10.2.4} (2).
   \begin{picture}(8,8)\put(0,0){\line(0,1){8}}\put(8,8){\line(0,-1){8}}\put(0,0){\line(1,0){8}}\put(8,8){\line(-1,0){8}}\end{picture}

 \begin {Proposition} \label {10.2.10} Let $H$ be a finite-dimensional
      Hopf algebra and $A = R\#_{\sigma} H$. Then

  (1) $r_{jH}(R) \#_{\sigma} H = r_{H^*j}(R \#_\sigma H);$

(2)  $r_{Hj}(R)=r_{jH}(R);$

  (3) $r_{Hj}(A \# H^*) = (r_{Hj}(R) \#_\sigma H) \# H^*$.
  \end {Proposition}

{\bf Proof} (1)         We see that
\begin{eqnarray*}
(r_{jH}(R)\#_{\sigma} H)\#H^* &=& \Phi(r_{jH}(R))\\
 &=& (\Phi(r_j(R) )\cap A) \#H^*  \hbox { \ \  by\   \cite [Theorem 7.2]{MS95} }\\
&=&(r_j(A\#H^*)\cap A) \# H^* \hbox {  \  \ by\  Lemma \ref {10.2.3} (8) } \\
&=& r_{H^*j}(A)\#H^* \hbox { \ by \ Definition \ref{10.1.1} }.
\end {eqnarray*}
Thus $r_{H^*j}(A)= r_{jH}(R)\#_{\sigma} H.$

(2) We see that
\begin {eqnarray*}
 r_{Hj}(R) &=& r_j(A) \cap R \\
& \supseteq & r_{H^*j}(A) \cap R \hbox  { \ by \ Proposition \ref {10.2.9} (1) }\\
&=&  r_{jH}(R)  \hbox { \ by \ part  \  (1) }.
\end {eqnarray*}
It follows by Proposition \ref{10.2.2} (1) that $r_{Hj}(R)=
r_{jH}(R)$.

(3) It immediately follows from part (1) (2).
\begin{picture}(8,8)\put(0,0){\line(0,1){8}}\put(8,8){\line(0,-1){8}}\put(0,0){\line(1,0){8}}\put(8,8){\line(-1,0){8}}\end{picture}

By Proposition \ref {10.2.9} and \ref {10.2.10}, it is clear that if
$H$ is a finite-dimensional Hopf algebra, then
 relation (\ref {e (2)})  holds iff  relation (\ref {e (3)}) holds.

 \begin {Theorem} \label {10.2.11}
 Let $H$ be a finite-dimensional  Hopf algebra and
  the weak action of $H$ be inner.
Then

(1) $r_{Hj}(R)= r_j(R)= r_{jH}(R).$

Moreover, if $H$ is semisimple, then

(2) $r_{j}(R\#_{\sigma} H)= r_{Hj}(R)\#_{\sigma} H$.
\end {Theorem}

{\bf Proof} (1) Since the weak action is inner, every ideal of $R$
is an $H$-ideal and $r_j(R) = r_{jH}(R)$. It follows from
Proposition \ref {10.2.3}(2) that $r_{Hj}(R)= r_{jH}(R)= r_j(R)$.

(2) Considering Proposition \ref {10.2.9}(2),
 it suffices to show
 $$r_{j}(R\#_{\sigma} H)\subseteq  r_{Hj}(R)\#_{\sigma} H.$$
It is clear that  $$(R\#_{\sigma} H)/ (r_{Hj}(R)\#_{\sigma} H) \cong
(R/r_{Hj}(R))\#_{\sigma} H  \hbox { \ \ \ (as algebras). }$$  It
follows by \cite [Corollary 7.4.3]{Mo93}  and part (1) that $
(R/r_{Hj}(R))\#_{\sigma} H$  is semiprimitive. Therefore $r_{j}(R
\#_{\sigma} H) \subseteq r_{Hj}(R)\#_{\sigma} H$.
\begin{picture}(8,8)\put(0,0){\line(0,1){8}}\put(8,8){\line(0,-1){8}}\put(0,0){\line(1,0){8}}\put(8,8){\line(-1,0){8}}\end{picture}

\begin {Theorem}\label {10.2.12} Let $H$ be  a finite-dimensional, semisimple Hopf
algebra, let $k$ be an algebraically closed field and let $A =
R\#_{\sigma} H$. Assume  $H$ is cosemisimple or char $k$ does not
divide dim $H$.

 (1) If $H$ is cocommutative, then
 $$r_{Hj}(R)= r_{jH}(R) = r_j(R);$$

(2) If $H$ is commutative, then
$$r_{j}(R \#_{\sigma}H)= r_{Hj}(R) \#_{\sigma}H.$$
\end {Theorem}

{\bf Proof}     By \cite [Theorem 4.3 (1)] {La71} , we have that $H$
is semisimple and cosemisimple.

(1)  If $g \in G(H)$, then the weak action of $g$ on $R$  is an
algebraic homomorphism, which implies that $g \cdot r_j(R) \subseteq
r_j(R)$. Let $H_0$ be the coradical of $H, H_1 = H_0 \wedge H_0,
H_{i+1} = H_0 \wedge H_i$ for $i = 1, \cdots, n$, where $n$ is the
dimension $dimH$ of $H$. It is clear that $H_0=kG$ with $G =G(H)$ by
\cite [Theorem 8.0.1 (c)]{Sw69a} and $H= \cup H_i$. It is easy to
show that if $k > i $, then  $$H_i \cdot (r_j(R))^k\subseteq
r_j(R)$$ by induction for $i$. Thus
$$H \cdot (r_j(R))^{dimH+1} \subseteq r_j(R),$$
which implies that $(r_j(R))^{dimH +1}\subseteq r_{jH}(R)$.

We see that
\begin{eqnarray*}
 r_{j}(R/r_{jH}(R)) &=& r_j(R)/r_{jH}(R)\\
 &=& r_b(R/r_{jH}(R)) \hbox { \ since \ }  r_j(R)/r_{jH}(R) \hbox { \ is
 \ nilpotent } \\
 &=& r_{bH}(R/r_{jH}(R)) \hbox { \ by \ Theorem \ \ref {10.2.7} (3)} \\
 &\subseteq &  r_{jH}(R/r_{jH}(R)) \\
 &=& 0 \ \ .
\end {eqnarray*}

Thus       $r_{j}(R) \subseteq r_{jH}(R)$, which implies that
  $r_j(R)=r_{jH}(R)$.

(2) It immediately follows from part (1) and Proposition
 \ref {10.2.10}(1) (2).
\begin{picture}(8,8)\put(0,0){\line(0,1){8}}\put(8,8){\line(0,-1){8}}\put(0,0){\line(1,0){8}}\put(8,8){\line(-1,0){8}}\end{picture}

\section {The general theory of $H$-radicals for twisted $H$-module algebras }\label {s26}

In this section we give the general theory of $H$-radicals for
twisted $H$-module algebras.

\begin {Definition}\label {10.3.1}
Let $r$ be a property of $H$-ideals of twisted $H$-module algebras.
An $H$- ideal $I$ of twisted $H$-module algebra $R$ is called   an
$r$-$H$-ideal of $R$ if it is of the $r$-property. A twisted
$H$-module algebra $R$  is called   an $r$-twisted $H$-module
algebra if it is $r$-$H$-ideal of itself. A property $r$ of
$H$-ideals of twisted $H$-module algebras   is called an
 $H$-radical property if the following conditions are satisfied:

(R1) Every twisted $H$-homomorphic image of $r$-twisted $H$-module
algebra is an $r$ twisted $H$-module algebra;

(R2) Every twisted $H$-module algebra $R$ has  the maximal
$r$-$H$-ideal $r(R)$;

(R3)  $R/r(R)$ has not any non-zero $r$-$H$-ideal.

\end {Definition}
We call $r(R)$ the $H$-radical of $R$.

\begin {Proposition}\label {10.3.2}
Let $r$ be an ordinary hereditary radical property for rings. An
$H$-ideal $I$ of twisted $H$-module algebra $R$ is called  an
$r_H$-$H$-ideal  of $R$ if $I$ is an $r$-ideal of ring $R$. Then
$r_H$ is an $H$-radical property for twisted $H$-module algebras.
\end {Proposition}

{\bf Proof.} (R1). If $(R,\sigma )$ is an $r_H$-twisted $H$-module
algebra
 and $(R, \sigma )  \stackrel {f} { \sim }  (R', \sigma ')$, then $r(R') =R'$
 by ring theory. Consequently, $R'$  is an $r_H$-twisted $H$-module algebra.

(R2). For any  twisted $H$-module algebra $R$, $r(R)$  is the
maximal $r$-ideal of $R$ by ring theory. It is clear that $r(R)_H$
is the maximal $r$-$H$-ideal,  which is an $r_H$-$H$-ideal of $R$.
Consequently, $r_H(R) = r(R)_H$ is the maximal $r_H$-$H$-ideal of
$R$.

(R3). If $I/r_H(R)$  is an $r_H$-$H$-ideal of $R/r_H(R)$, then $I$
is an $r$-ideal of algebra $R$ by ring theory. Consequently, $I
\subseteq r(R)$ and $I \subseteq r_H(R).$
\begin{picture}(8,8)\put(0,0){\line(0,1){8}}\put(8,8){\line(0,-1){8}}\put(0,0){\line(1,0){8}}\put(8,8){\line(-1,0){8}}\end{picture}

\begin {Proposition}\label {10.3.3}
$r_{Hb}$ is an $H$-radical property.
\end {Proposition}

{\bf Proof.} (R1). Let $(R,\sigma )$ be an $r_{Hb}$-twisted
$H$-module algebra
 and $(R, \sigma )  \stackrel {f} { \sim }  (R' \sigma ')$.
    For any $x' \in R'$
and  any  $H$-$m$-sequence $\{  a_n' \}$    in $R'$
 with $ a_1' =  x' $,  there exist
 $ b_n' \in  R'$ and $h_n, h_n' \in H$ such that
 $  a_{n+1}' =
 (h_n \cdot  a_n')  b_n' (h_n' \cdot  a_n')$
                           {~}{~}  for any natural number $n$.
Let   $a_1, b_i \in R $  such that $f(a_1) =x'$ and   $f(b_i)=b_i'$
for $i = 1, 2, \cdots .$  Set
 $  a_{n+1} = (h_n \cdot  a_n)  b_n (h_n' \cdot  a_n)$ {~}{~}
                             for any natural number $n$.
 Since $ \{ a_n \}$ is an $H$-$m$-sequence in $R$,
 there exists a natural number $k$ such that $a_k =0$. It is
 clear that $f( a_n)= a_n'$
  for any natural number $n$ by induction. Thus $  a_k' = 0,$
  which implies that $x' $  is an $H$-$m$-nilpotent element.
Consequently, $R'$ is an $r_{Hb}$-twisted $H$-module algebra.

(R2). By \cite [Theorem 1.5] {Zh98a}, $r_{Hb}(R)= W_H(R) = \{ a \mid
a$ is an $H$-$m$-nilpotent element in $R \}$. Thus $r_{Hb}(R)$ is
the maximal $r_{Hb}$-$H$-ideal of $R$.

(R3).  It immediately follows from  \cite [Proposition 1.4] {Zh98a}.
\begin{picture}(8,8)\put(0,0){\line(0,1){8}}\put(8,8){\line(0,-1){8}}\put(0,0){\line(1,0){8}}\put(8,8){\line(-1,0){8}}\end{picture}

\section {The relations among  radical of $R$  , radical of
  $R \# _\sigma H$, and $H$-radical of $R$ }\label {s27}

In this section we give the relation among the Jacobson radical $r_j
(R)$ of $R$  ,the Jacobson
 radical $r_j(R \# _\sigma H)$ of
  $R \# _\sigma H$, and $H$-Jacobson radical $r_{Hj}(R)$ of $R$.

In this section, let $k$ be a field, $R$ an algebra with unit, $H$ a
Hopf algebra over $k$ and $R \# _\sigma H$ an algebra with unit. Let
$r$ be a hereditary radical property for rings which satisfies
$$r(M_{n \times n} (R)) = M_{n \times n} (r(R))$$ for any twisted
$H$-module algebra $R$.

 Example.  $r_j$, $r_{bm}$ and $r_n$ satisfy the above conditions
by \cite {Sz82}. Using \cite [Lemma 2.1 (2)]{Zh98a},we can easily
prove that $r_b$ and $r_l$   also satisfy the above conditions.

\begin {Definition} \label {10.4.1}
$\bar r_H (R) := r(R \# _\sigma H) \cap R$  and $r_H (R) :=
(r(R):H)$.
\end {Definition}

     If $H$ is a finite-dimensional   Hopf algebra and $M = R\#_{\sigma} H$,
     then $M$ is a free right $R$-module  with finite rank by
\cite [Proposition 7.2.11]{Mo93} and $End (M_R) \cong (R \#_{\sigma}
H)\#H^*$ by \cite [Corollary 9.4.17]{Mo93}. It follows from  part
(a) in the proof of \cite [Theorem 7.2] {MS95} that
 there exists a  unique bijective  map
 $$\Phi: {\cal I}(R) \longrightarrow {\cal I}(R')$$
such that \ \ \  $\Phi (I)M=MI,$ where $R' = (R \#_\sigma H)\# H^*$
and

 ${ \cal I}(R) = \{ I \mid I \hbox { is an ideal of } R \}$.

\begin {Lemma}\label {10.4.2}  If $H$ is a finite-dimensional Hopf algebra, then

  $ \Phi (r(R)) = r((R \# _\sigma H) \# H^*).$

\end {Lemma}
{\bf Proof.}  It is similar to the proof of \cite [lemma 2.1 (5)]
{Zh98a}.
\begin{picture}(8,8)\put(0,0){\line(0,1){8}}\put(8,8){\line(0,-1){8}}\put(0,0){\line(1,0){8}}\put(8,8){\line(-1,0){8}}\end{picture}

\begin {Proposition}\label {10.4.3}
$\bar r_H (R) \# _\sigma H \subseteq r_{H^*}(R \# _\sigma H)
                           \subseteq r(R \# _\sigma H). $

\end {Proposition}
{\bf Proof.}
 We see that
\begin {eqnarray*}
   \bar r_{H}(R) \#_{\sigma} H &=&    (\bar r_{H}(R) \#_{\sigma} 1)(1 \#_{\sigma} H)\\
    &\subseteq &  r(R \#_{\sigma} H) ( 1 \#_\sigma H) \\
    &\subseteq & r(R \#_{\sigma} H). { \ \ \ \  }
\end {eqnarray*}
Thus      $\bar r_H (R) \# _\sigma H \subseteq r_{H^*}(R \# _\sigma
H)$ since     $\bar r_H(R) \# _\sigma H $  is an $H^*$-ideal of $R
\# _\sigma H $.
\begin{picture}(8,8)\put(0,0){\line(0,1){8}}\put(8,8){\line(0,-1){8}}\put(0,0){\line(1,0){8}}\put(8,8){\line(-1,0){8}}\end{picture}

\begin {Proposition}\label {10.4.4}
If $H$ is a finite-dimensional
      Hopf algebra,  then

(1) $r_H (R) \# _\sigma H =  \bar r_{H^*}(R \# _\sigma H)$;

Furthermore, if               $\bar r_H \le r_H $, then

(2)    $\bar r_H = r_H$ and $r_H (R) \# _\sigma H \subseteq r(R \#
_\sigma H)$;

(3)    $R\# _\sigma H$  is $r$-semisimple for any $r_H$-semisimple
$R$ iff
   $$r(R \#_\sigma H) = r_H(R) \# _\sigma H.$$

\end {Proposition}
{\bf Proof.} Let $A = R \# _\sigma H.$

(1)         We see that
\begin{eqnarray*}
(r_{H}(R)\#_{\sigma} H)\#H^* &=& \Phi(r_{H}(R))\\
 &=& (\Phi(r(R) )\cap A) \#H^*  \hbox { \ \  by\   \cite [Theorem 7.2]{MS95} }\\
&=&(r(A\#H^*)\cap A) \# H^* \hbox {  \  \ \ by  Lemma \ref {10.4.2}  } \\
&=& \bar r_{H^*}(A)\#H^* \hbox { \ by \ Definition \ref{10.4.1} }.
\end {eqnarray*}
Thus $ \bar r_{H^*}(A)= r_{H}(R)\#_{\sigma} H.$

(2) We see that
\begin {eqnarray*}
\bar  r_{H}(R) &=& r(A) \cap R \\
& \supseteq & \bar r_{H^*}(A) \cap R \hbox  { \ by assumption }\\
&=&  r_{H}(R)  \hbox { \ by \ part  \  (1) }.
\end {eqnarray*}
Thus  $\bar r_{H}(R)= r_{H}(R)$ by assumption.

(3) Sufficiency is obvious. Now we show the necessity.
            Since $$r((R \# _\sigma H)/(r_H (R) \#_\sigma H))
            \cong r(R/r_H(R) \#_{\sigma '}H) =0 ,$$
            we have $r (R \#_\sigma H) \subseteq r_H(R) \#_\sigma H.$
Considering part (2), we have

   $$r(R \#_\sigma H) = r_H(R) \# _\sigma H. { \ \ \ } \Box  $$

\begin {Corollary}\label {10.4.5}
Let $r$ denote $r_b, r_l,  r_j, r_{bm}$ and $r_n.$ Then

(1) $\bar r_H \le r_H;$

Furthermore, if $H$ is a finite-dimensional Hopf algebra, then

(2) $\bar r_H = r_H $;

(3) $R \# _\sigma H $ is $r$- semisimple for any $r_{H}$-semisimple
$R$  iff $r(R \# _\sigma H)= r_{H}(R)\#_\sigma H$;

(4)  $R \# _\sigma H $ is $r_j$- semisimple for any
$r_{Hj}$-semisimple $R$  iff $r_j(R \# _\sigma H)=
r_{Hj}(R)\#_\sigma H.$

\end {Corollary}
{\bf Proof.}  (1) When $r= r_b$  or $r = r_j$, it has been proved in
\cite [Proposition 2.3 (1) and 3.2 (1)] {Zh98a} and in the preceding
sections. The others can similarly be proved.

(2) It follows from Proposition \ref {10.4.4} (2).

(3) and (4)  follow from part (1) and Proposition \ref {10.4.4} (3).
\begin{picture}(8,8)\put(0,0){\line(0,1){8}}\put(8,8){\line(0,-1){8}}\put(0,0){\line(1,0){8}}\put(8,8){\line(-1,0){8}}\end{picture}
\begin {Proposition}\label {10.4.6}
 If $H = kG$ or  the weak action of $H$ on $R$ is inner,
then

(1). $r_H (R) = r(R)$;

(2)  If, in  addition,   $H$ is a finite-dimensional
 Hopf algebra and $\bar r _H \le r _H$, then
 $r_H(R) = \bar r_H(R) = r(R).$

\end {Proposition}

{\bf Proof.} (1)  It is trivial.

(2) It immediately follows from part (1) and Proposition \ref
{10.4.1} (1) (2).
\begin{picture}(8,8)\put(0,0){\line(0,1){8}}\put(8,8){\line(0,-1){8}}\put(0,0){\line(1,0){8}}\put(8,8){\line(-1,0){8}}\end{picture}

\begin {Theorem}\label {10.4.8}
Let $G$ be a finite group and $\mid G \mid ^{-1} \in k$. If $H= kG$
or $H= (kG)^*$, then

(1)  $r_j (R) = r_{Hj}(R)= r_{jH}(R)$;

  (2)   $r_{j}(R\# _\sigma  H)= r_{Hj}(R)\# _\sigma H.$
                  \end {Theorem}
{\bf Proof.} (1) Let $H =kG$.  We  can easily check $r_j(R) =
r_{jH}(R)$  using the method similar to  the proof of
 \cite [Proposition 4.6]  {Zh97b}. By  \cite [Proposition 3.3 (2) ] {Zh98a},
 $r_{Hj}(R) =r_{jH}(R)$. Now,
we only need to show that
  $$r_j (R) = r_{H^*j}(R). $$
 We see that
\begin {eqnarray*}
  r_{j}((R\# _\sigma  H^*) \# H) &=&
   r_{H^*j}((R\# _\sigma H^*) \# H)   \hbox { \ \ by \cite [Theorem 4.4 (3)]
   {CM84b} } \\
   &=&  r_{Hj}(R\# _\sigma  H^*) \# H  \hbox { \ \ by \cite [ Proposition
   3.3 (1)] {Zh98a}}  \\
   &=& (r_{H^*j}(R)\# _\sigma H^* ) \# H  \hbox { \ \ by \cite [ Proposition
   3.3 (1)] {Zh98a}}.
\end {eqnarray*}
On the one hand, by       \cite [ Lemma 2.2 (8)] {Zh98a}, $\Phi (r_j
(R)) =  r_{j}((R\# _\sigma  H^*) \# H).$ On the other hand, we have
that $ \Phi (r_{H^*j} (R)) =  ( r_{H^*j}(R )\# _\sigma H^*) \# H $ \
\
 by \cite [Lemma 2.2 (2)]    {Zh98a}.
 Consequently, $r_j (R) = r_{H^*j}(R)$.

(2) It immediately follows from part (1) and \cite [Proposition 3.3
(1) (2)] {Zh98a}.
  \begin{picture}(8,8)\put(0,0){\line(0,1){8}}\put(8,8){\line(0,-1){8}}\put(0,0){\line(1,0){8}}\put(8,8){\line(-1,0){8}}\end{picture}

\begin {Corollary}\label {10.4.9}
Let $H$ be a semisimple and cosemisimple Hopf algebra over
algebraically closed field $k$. If $H$ is commutative or
cocommutative, then
  $$r_j (R) = r_{Hj}(R)= r_{jH}(R)  \hbox { \ \ \ \ and  \ \ \ }
  r_{j}(R\# _\sigma  H)= r_{Hj}(R)\# _\sigma H.$$
                  \end {Corollary}
{\bf Proof.} It immediately follows from Theorem \ref {10.4.8} and
\cite [Lemma 8.0.1 (c)] {Sw69a}.
\begin{picture}(8,8)\put(0,0){\line(0,1){8}}\put(8,8){\line(0,-1){8}}\put(0,0){\line(1,0){8}}\put(8,8){\line(-1,0){8}}\end{picture}

    We give an example to show that conditions in Corollary \ref {10.4.9} can not
be omitted.

\begin {Example}\label {10.4.10} (see \cite [Example P20] {Zh98a})
Let $k$ be a field of characteristic $p > 0 $, $R = k[x]/(x^p)$. We
can define a derivation on $R$ by sending $x$ to $x +1$. Set
$H=u(kd)$, the restricted enveloping algebra, and $A = R\#H.$ Then

(1) $r_b (A \# H^*) \not= r_{H^*b}(A) \# H^*$;

(2) $r_j (A \# H^*) \not= r_{H^*j}(A) \# H^*$;

(3) $r_j (A \# H^*) \not\subseteq  r_{jH^*}(A) \# H^*$.

\end {Example}
{\bf Proof.}        (1) By \cite [Example P20]  {Zh98a}, we have
$r_b(R) \not=0$ and $r_{bH}(R)=0.$ Since $\Phi (r_b(R) ) = r_b(A
\#H^*)\not=0$  and $\Phi (r_{bH}(R))= r_{bH^*}(A) \#H^*=0,$ we have
that  part (1) holds.

(3) We see that $r_j(A \# H^*) = \Phi (r_j(R))$  and $r_{Hj}(A)
\#H^*= \Phi (r_{Hj}(R)).$  Since $R$ is commutative, $r_j(R)=
r_b(R).$  Thus $r_{Hj}(R)= r_{jH}(R)= r_{bH}(R)=0$ and $r_j(R) =
r_b(R) \not=0$, which implies
 $r_j(A\#H^*) \not\subseteq r_{jH^*}(A)\#H^*$.

(2) It follows from part (3).
 \begin{picture}(8,8)\put(0,0){\line(0,1){8}}\put(8,8){\line(0,-1){8}}\put(0,0){\line(1,0){8}}\put(8,8){\line(-1,0){8}}\end{picture}

This example also answer the question J.R. Fisher asked in \cite
{Fi75} :
$$\hbox {Is  \ \ \ } r_j (R \# H) \subseteq  r_{jH}(R) \# H \hbox { \ \ \  ?} $$

If $F$ is an extension field of $k$, we write  $R^F$  for $R \otimes
_k F$ (see \cite [P49 ]{MS95}) .

\begin {Lemma}\label {10.4.11}  If $F$ is an extension field  of $k$, then

(1)  $H$ is a semisimple  Hopf algebra over $k$  iff $H^F$ is a
semisimple Hopf algebra over $F$;

(2)  Furthermore, if $H$ is a finite-dimensional Hopf algebra, then
  $H$ is a cosemisimple  Hopf algebra over $k$  iff $H^F$ is a cosemisimple
Hopf algebra over $F$.

\end {Lemma}
{\bf Proof.} (1) It is clear that $\int _H^l \otimes F = \int
_{H^F}^l.$   Thus
 $H$ is a semisimple  Hopf algebra over $k$  iff $H^F$ is a semisimple
Hopf algebra over $F$.

(2)  $(H \otimes F)^* = H^* \otimes F $ since $H^* \otimes F
\subseteq  (H \otimes F)^*$  and $dim_F (H \otimes F ) = dim_F (H^*
\otimes F)= dim _k H$. Thus we can obtain part (2) by  Part (1).
\begin{picture}(8,8)\put(0,0){\line(0,1){8}}\put(8,8){\line(0,-1){8}}\put(0,0){\line(1,0){8}}\put(8,8){\line(-1,0){8}}\end{picture}

By the way, if $H$ is a semisimple Hopf algebra, then $H$ is a
separable algebra by Lemma \ref {10.4.11} (see \cite [P284] {Pa77}).

\begin {Proposition}\label {10.4.12}
Let $F$ be an algebraic closure of $k$, $R$ an algebra over $k$ and
$$r(R \otimes _k F) = r(R) \otimes _k F.$$
 If $H$ is a finite-dimensional  Hopf algebra  with cocommutative
 coradical over $k$ ,
then
     $$r(R) ^{dim H} \subseteq r_H(R).$$

\end {Proposition}
{\bf Proof.} It is clear that $H^F$  is   a finite-dimensional  Hopf
algebra
 over $F$ and $dim H = dim  H^F =n$. Let $H^F_0$  be the coradical
 of  $H^F$,  $H^F_1 = H^F_0 \wedge H^F_0,
H^F_{i+1} = H^F_0 \wedge H^F_i$ for $i = 1, \cdots, n-1$. Notice
$H^F_0 \subseteq H_0 \otimes F $. Thus $H^F _0 $  is cocommutative.
It is clear that $H^F_0=kG$  by \cite [Lemma 8.0.1 (c)]{Sw69a} and
$H^F= \cup H^F_i$. It is easy to show that if $k > i $, then
$$H^F_i \cdot (r(R ^F))^k\subseteq r(R ^F)$$ by induction for $i$.
Thus
$$H^F \cdot (r(R ^F))^{ dim H} \subseteq r(R^F),$$
which implies that $(r(R^F))^{dim H }\subseteq r(R^F)_{H^F}$. By
assumption, we have that $(r(R) \otimes F)^{dim H }\subseteq (r(R)
\otimes F)_{H^F}$. It is clear that  $(I \otimes F)_{H^F} = I_H
\otimes F$  for any ideal $I$ of $R$. Consequently, $(r(R))^{dim H
}\subseteq r(R)_H$.
\begin{picture}(8,8)\put(0,0){\line(0,1){8}}\put(8,8){\line(0,-1){8}}\put(0,0){\line(1,0){8}}\put(8,8){\line(-1,0){8}}\end{picture}

\begin {Theorem}\label {10.4.13}
Let $H$ be a  semisimple, cosemisimple and either commutative or
cocommutative Hopf algebra
 over $k$.
If there exists an algebraic closure $F$ of $k$ such that
$$r_j(R \otimes  F) = r_j(R) \otimes F \hbox { \ \ and \ \ }
r_j((R \# _\sigma H) \otimes  F) = r_j(R \# _\sigma H) \otimes F,$$
then

(1) $r_j (R) = r_{Hj}(R)= r_{jH}(R); $

 (2)
  $r_{j}(R\# _\sigma  H)= r_{Hj}(R)\# _\sigma H.$
\end {Theorem}
{\bf Proof.}
 (1). By Lemma \ref {10.4.11},
$H^F$ is semisimple and cosemisimple. Considering Corollary \ref
{10.4.9}, we have that
   $r_j (R^F) = r_{H^Fj}(R^F)= r_{jH^F}(R^F).$
On the one hand, by assumption, $r_j(R^F) = r_j(R) \otimes F$. On
the other hand, $r_{jH^F} (R^F) = (r_j (R ) \otimes F)_{H^F} =
r_{jH}(R) \otimes F$. Thus $r_j(R) = r_{jH}(R)$.

(2).  It immediately  follows  from part (1).
\begin{picture}(8,8)\put(0,0){\line(0,1){8}}\put(8,8){\line(0,-1){8}}\put(0,0){\line(1,0){8}}\put(8,8){\line(-1,0){8}}\end{picture}

Considering Theorem \ref {10.4.13} and \cite [Theorem 7.2.13]
{Pa77}, we have

\begin {Corollary}\label {10.4.14}
Let $H$ be a semisimple,  cosemisimple and either commutative or
cocommutative Hopf
 algebra over $k$. If there exists an algebraic closure $F$ of $k$ such that
$F/k$ is separable and algebraic, then

  (1) $r_j (R) = r_{Hj}(R)= r_{jH}(R)$;

  (2)   $r_{j}(R\# _\sigma  H)= r_{Hj}(R)\# _\sigma H.$
                  \end {Corollary}

\begin {Lemma}\label {10.4.15}
(Szasz \cite {Sz82})   $$r_j (R) = r_k(R)$$ holds in the following
three cases:

(1) Every element in $R$ is  algebraic over $k$ (\cite [Proposition
31.2] {Sz82});

(2) The cardinality of $k$  is strictly greater  than the dimension
of $R$ and $k$ is infinite (\cite [Theorem 31.4] {Sz82});

(3) $k$ is uncountable and $R$ is  finitely generated (\cite
[Proposition 31.5] {Sz82}).
\end {Lemma}

\begin {Proposition}\label {10.4.16}
Let $F$ be an extension of  $k$. Then $r(R) \otimes  F \subseteq r(R
\otimes  F),$ where $r$ denotes $r_b, r_k, r_l, r_n  .$

\end {Proposition}
{\bf Proof.} When $r =r_n$, for any $x \otimes a  \in r_n (R)
\otimes F$ with $a \not=0$, there exists $y \in R$ such that $x
=xyx$. Thus $x \otimes a = (x \otimes a)(y \otimes a^{-1}) (x
\otimes a)$, which implies $r_n(R) \otimes F \subseteq r_n(R \otimes
F)$.

Similarly, we can obtain the others.
\begin{picture}(8,8)\put(0,0){\line(0,1){8}}\put(8,8){\line(0,-1){8}}\put(0,0){\line(1,0){8}}\put(8,8){\line(-1,0){8}}\end{picture}

\begin {Corollary}\label {10.4.17}
Let $H$ be a semisimple, cosemisimple and commutative or
cocommutative Hopf algebra. If there exists an algebraic closure $F$
of $k$ such that $F/k$ is a pure transcendental extension and one of
the following three conditions holds:

(i) every element in $R\#_\sigma H$ is  algebraic over $k$;

(ii) the cardinality of $k$  is strictly greater  than the dimension
of $R$ and $k$ is infinite;

(iii) $k$ is uncountable and $R$ is  finitely generated;

then

(1) $r_j (R) = r_{Hj}(R)= r_{jH}(R)$;

(2)    $r_{j}(R\# _\sigma  H)= r_{Hj}(R)\# _\sigma H$;

(3) $r_j(R)= r_k(R) $  and $r_j(R \#_\sigma H) = r_k(R\#_\sigma H).$
                  \end {Corollary}
{\bf Proof.} First, we have that part (3) holds by Lemma \ref
{10.4.15}. We next see that
\begin {eqnarray*}
r_j(R \otimes F) &\subseteq & r_j(R) \otimes F  \hbox { \ \  \cite
[Theorem
7.3.4] {Pa77}} \\
&=&  r_k(R) \otimes F  \hbox { \ \  part (3)} \\
&\subseteq & r_k(R \otimes F ) \hbox { \ \ proposition \ref {10.4.16} } \\
 &\subseteq & r_j(R \otimes F).
\end {eqnarray*}
Thus $r_j(R \otimes F) = r_j (R) \otimes F.$ Similarly, we can show
that
 $r_j((R \#_\sigma H) \otimes F) = r_j (R\#_\sigma H) \otimes F.$

Finally, using Theorem \ref {10.4.13}, we complete the proof.
\begin{picture}(8,8)\put(0,0){\line(0,1){8}}\put(8,8){\line(0,-1){8}}\put(0,0){\line(1,0){8}}\put(8,8){\line(-1,0){8}}\end{picture}

\section {The $H$-Von Neumann regular radical}\label {s28}

In this section, we construct the $H$-von Neumann regular radical
for $H$-module algebras and show that it is an $H$-radical property.

\begin {Definition} \label {10.5.1}
Let $a \in R$. If $a\in (H \cdot a) R (H\cdot a)$, then $a$ is
called an $H$-von Neumann regular element, or an $H$-regular element
in short. If every element of $R$ is an $H$-regular, then $R$ is
called an $H$-regular module algebra, written as $r_{Hn}$-$H$-module
algebra.  $I$ is an $H$-ideal
 of $R$ and every element in $I$ is $H$-regular, then $I$ is called an $H$-
 regular ideal.
\end {Definition}

\begin {Lemma}\label {10.5.2} If $I$ is an $H$-ideal of $R$ and $a \in I$,
then $a$ is $H$-regular in $I$ iff $a$ is $H$-regular in $H$.
\end {Lemma}
{\bf Proof.} The necessity is clear.

Sufficiency: If $a \in (H \cdot a)R (H\cdot a),$ then there exist
$h_i, h_i' \in H, b_i \in R,$  such that
$$a = \sum (h_i \cdot a) b_i (h_i' \cdot a).$$
We see that
\begin {eqnarray*}
a &=& \sum _{i, j} [ h_i \cdot (( h_j \cdot a) b_j (h_j' \cdot
a))]b_i
(h_i' \cdot a) \\
&=&      \sum _{i, j} [ ((h_i)_1 \cdot ( h_j \cdot a)) ((h_i)_2
\cdot b_j)
      ((h_i)_3 \cdot (h_j' \cdot a))]b_i
(h_i' \cdot a)    \\
&\in & (H\cdot a) I (H\cdot a).
\end {eqnarray*}
Thus $a$ is an $H$-regular in $I$.
\begin{picture}(8,8)\put(0,0){\line(0,1){8}}\put(8,8){\line(0,-1){8}}\put(0,0){\line(1,0){8}}\put(8,8){\line(-1,0){8}}\end{picture}

\begin {Lemma}\label {10.5.3}
If  $x - \sum _i (h_i \cdot x) b_i (h_i' \cdot x)$ is $H$-regular,
then $x$ is $H$-regular, where $x, b_i \in R, h_i , h_i' \in H.$
\end {Lemma}
{\bf Proof.} Since   $x - \sum _i (h_i \cdot x) b_i (h_i' \cdot x)$
is $H$-regular, there exist $g_i, g_i' \in H, c_i \in R $ such that
$$x - \sum _i (h_i \cdot x) b_i (h_i' \cdot x)  =
\sum _j ( g_j \cdot (x - \sum _i (h_i \cdot x) b_i (h_i' \cdot x)))
c_j ( g_j' \cdot ( x - \sum _i (h_i \cdot x) b_i (h_i' \cdot x))).
$$ Consequently, $x \in (H\cdot x)R (H \cdot x).$
\begin{picture}(8,8)\put(0,0){\line(0,1){8}}\put(8,8){\line(0,-1){8}}\put(0,0){\line(1,0){8}}\put(8,8){\line(-1,0){8}}\end{picture}

\begin {Definition}\label {10.5.4}
$$r_{Hn}(R):=  \{ a \in R \mid \hbox  { \ the \ }   H\hbox {-ideal } (a)
\hbox { generated by } a \hbox { is } H \hbox {-regular } \}.$$
\end {Definition}

\begin {Theorem}\label {10.5.5}
$r_{Hn}(R)$ is an $H$-ideal of $R$.
\end {Theorem}
{\bf Proof.} We first show that $R r_{Hn}(R) \subseteq r_{Hn}(R).$
 For any $a \in r_{Hn}(R), x \in R,$ we have that
 $(xa) $  is $H$-regular since $(xa) \subseteq (a).$
 We next show that $a-b \in r_{Hn}(R)$  for any $a , b \in r_{Hn}(R).$
 For any $x \in (a-b),$ since $(a-b) \subseteq (a) + (b)$, we have that
 $x = u -v$ and $u \in (a), v\in (b).$
 Say $u = \sum _i (h_i \cdot u) c_i (h_i' \cdot u)$  and $h_i, h_i' \in H,
  c_i \in R.$
We see that
\begin {eqnarray*}
&x& - \sum _i (h_i \cdot x)c_i (h_i' \cdot x)\\
 &=& (u-v) - \sum _i  (h_i \cdot (u-v))c_i (h_i' \cdot (u-v))   \\
&=& -v - \sum _i [ -  (h_i \cdot u)c_i (h_i' \cdot v) -   (h_i \cdot
v)c_i (h_i' \cdot u)
+  (h_i \cdot v)c_i (h_i' \cdot v)] \\
& \in & (v) .
\end {eqnarray*}
Thus  $x - \sum _i (h_i \cdot x)c_i (h_i' \cdot x) $ is $H$-regular
and $x$ is $H$-regular by Lemma \ref {10.5.3}. Therefore $a-b \in
r_{Hn}(R).$ Obviously, $r_{Hn}(R)$ is $H$-stable. Consequently,
$r_{Hn}(R)$ is an $H$-ideal of $R$.
\begin{picture}(8,8)\put(0,0){\line(0,1){8}}\put(8,8){\line(0,-1){8}}\put(0,0){\line(1,0){8}}\put(8,8){\line(-1,0){8}}\end{picture}

\begin {Theorem}\label {10.5.5}
$r_{Hn}(R/r_{Hn}(R)) = 0.$

\end {Theorem}
{\bf Proof.}   Let $\bar R = r/r_{Hn}(R)$ and  $\bar b = b +
r_{Hn}(R) \in r_{Hn}(R/r_{Hn}(R)).$ It is sufficient to show that $b
\in r_{Hn}(R).$  For any $a \in (b)$, it is clear that  $\bar a \in
(\bar b)$.  Thus there exist $h_i, h_i' \in H, \bar c_i \in \bar R$
such that
$$\bar a = \sum _i (h_i \cdot \bar a) \bar c_i (h_i' \cdot \bar a)
= \sum _i \overline { (h_i \cdot  a)  c_i (h_i' \cdot  a)}.$$ Thus $
a-  \sum _i  (h_i \cdot  a)  c_i (h_i' \cdot  a) \in r_{Hn}(R),$
which implies that $ a $  is $H$-regular. Consequently,   $b \in
r_{Hn}(R).$ Namely, $\bar b =0$ and $r_{Hn}(R)=0.$
\begin{picture}(8,8)\put(0,0){\line(0,1){8}}\put(8,8){\line(0,-1){8}}\put(0,0){\line(1,0){8}}\put(8,8){\line(-1,0){8}}\end{picture}

\begin {Corollary}\label {10.5.7.1}
$r_{Hn}$ is an $H$-radical property for $H$-module algebras and
$r_{nH} \le r_{Hn}$.
\end {Corollary}
{\bf Proof.} (R1) If $R \stackrel {f} {\sim }  R'$  and $R$  is an
$r_{Hn}$ -$H$-module algebra, then for any $f(a) \in R'$, $f (a) \in
(H \cdot f(a)) R' (H \cdot f(a)).$  Thus $R'$ is also an $r_{Hn}$-
$H$-module algebra.

(R2) If $I$ is an $r_{Hn}$-$H$-ideal of $R$  and $ r_{Hn}(R)
\subseteq I$ then, for any $a\in I,$  $(a)$  is $H$-regular since
$(a) \subseteq I.$  Thus $I \subseteq r_{Hn}(R).$

  (R3) It follows from Theorem \ref {10.5.6}.

  Consequently $r_{Hn}$ is an $H$-radical property for $H$-module algebras.
  It is straightforward  to check  $r_{nH} \le r_{Hn}.$
 \begin{picture}(8,8)\put(0,0){\line(0,1){8}}\put(8,8){\line(0,-1){8}}\put(0,0){\line(1,0){8}}\put(8,8){\line(-1,0){8}}\end{picture}

  $r_{Hn}$ is called the  $H$-von Neumann regular radical.

\begin {Theorem}\label {10.5.6}
If $I$ is an $H$-ideal of $R$, then $r_{Hn}(I) = r_{Hn}(R) \cap I $.
Namely, $r_{Hn}$  is a strongly hereditary $H$-radical property.

\end {Theorem}
{\bf Proof.} By Lemma \ref {10.5.2}, $r_{Hn} (R) \cap I \subseteq
r_{Hn}(I).$ Now, it is sufficient to show that $(x)_I = (x)_R$ for
any $x \in r_{Hn}(I)$, where $(x)_I$ and $(x)_R$ denote the
$H$-ideals generated by $x$ in $I$ and $R$ respectively. Let $x =
\sum (h_i \cdot x) b_i (h_i ' \cdot x)$ , where $h_i, h_i' \in H,
b_i \in I.$  We see that
\begin {eqnarray*}
R(H\cdot x) &=& R( H \cdot ( \sum (h_i \cdot x) b_i (h_i ' \cdot x)) \\
 &\subseteq &R (H\cdot x) I (H \cdot x)  \\
 &\subseteq & I (H \cdot x).
\end {eqnarray*}
Similarly, $$(H\cdot x)R \subseteq (H\cdot x)I.$$ Thus $(x)_I=
(x)_R.$
\begin{picture}(8,8)\put(0,0){\line(0,1){8}}\put(8,8){\line(0,-1){8}}\put(0,0){\line(1,0){8}}\put(8,8){\line(-1,0){8}}\end{picture}

A graded algebra $R$ of type $G$ is said to be Gr-regular if for
every homogeneous $a \in R_g$ there exists $b\in R$ such that $a=
aba$
 \ \ \ ( see \cite {NO82} P258 ). Now, we give the relations between
 Gr-regularity and  $H$-regularity.

\begin {Theorem}\label {10.5.7}
              If $G$ is a finite group, $R$ is a graded algebra of type $G$,
               and $H=(kG)^*$, then
$R$ is  Gr-regular iff  $R$ is  $H$-regular.

\end {Theorem}
{\bf Proof.} Let $\{ p_g \mid g\in G \}$ be the dual base of base
$\{ g \mid g\in G \}$. If $R$ is Gr-regular for any $a\in R$, then
$a = \sum _{g \in G}  a_g$ with $a_g \in R_g$. Since $R$ is
Gr-regular, there exist $b_{g^{-1}} \in R_{g^{-1}}$ such that $a_g =
a_g b_{g^{-1}} a_g $  and
$$ a = \sum _{g \in G} a_g = \sum _{g \in G } a_g b_{g^{-1}} a_g
= \sum _{g \in G}( p_g \cdot a ) b_{g^{-1}} (p_g \cdot a).$$
Consequently, $R$ is $H$-regular.

Conversely, if $R$ is $H$-regular, then for any $a \in R_g,$ there
exists $b_{x, y}\in R $  such that
$$a = \sum _{x,y \in G} (p_x \cdot a) c_{x,y} (p_y \cdot a).$$
Considering $a \in R_g$, we have that $a = a b_{g,g} a.$ Thus $R$ is
Gr-regular
.\begin{picture}(8,8)\put(0,0){\line(0,1){8}}\put(8,8){\line(0,-1){8}}\put(0,0){\line(1,0){8}}\put(8,8){\line(-1,0){8}}\end{picture}

\section { About J.R. Fisher's question }\label {s29}
 In this section,
 we answer the question  J.R. Fisher asked in \cite {Fi75}. Namely,
 we give a necessary and sufficient condition for validity of relation (2) .

Throughout this section, let $k$ be a commutative ring with unit,
$R$ an $H$- module algebra and $H$ a Hopf algebra over $k.$

\begin {Theorem}\label {10.6.1}
Let  ${\cal K}$        be an ordinary special class of rings and
closed with respect to isomorphism. Set $r= r^{\cal K}$ and $\bar
r_H (R) = r(R \# H) \cap R$ for any $H$-module algebra $R$. Then
$\bar r_H$ is an $H$-radical property of $H$-module algebras.
Furthermore, it is an $H$-special radical.

\end {Theorem}

{\bf Proof.} Let $\bar {\cal M}_R = \{ M \mid M $ is an $R$-prime
module and $R/(0:M)_R \in {\cal K } \}$    for any ring $R$ and
$\bar {\cal M} = \cup \bar {\cal M}_R.$ Set ${\cal M}_R =\{ M \mid M
\in \bar {\cal M }_{R \# H} \}$ for any $H$- module algebra $R$ and
${\cal M}= \cup {\cal M}_R.$  It is straightforward to check that
$\bar {\cal M}$  satisfies the conditions of \cite [Proposition 4.3]
{Zh97b}. Thus ${\cal M}$ is an $H$-special module by      \cite
[Proposition 4.3] {Zh97b}.    It is clear that ${\cal M}(R) = \bar {
\cal M} (R \#H) \cap R = r(R \#H) \cap R$  for any $H$-module
algebra $R$. Thus $\bar r_H$ is an $H$-special radical by \cite
[Theorem 3.1] {Zh97b}.
\begin{picture}(8,8)\put(0,0){\line(0,1){8}}\put(8,8){\line(0,-1){8}}\put(0,0){\line(1,0){8}}\put(8,8){\line(-1,0){8}}\end{picture}

Using the Theorem \ref {10.6.1}, we have that
 $\bar r_{bH}, \bar r_{lH}, \bar r_{kH}, \bar r_{jH}, \bar r_{bmH}$
are all $H$-special radicals.

\begin {Proposition}\label {10.6.2}
Let ${\cal K}$ be a  special class of rings and closed with respect
to isomorphism.  Set $r= r^{\cal K}.$ Then

(1) $\bar r_H (R) \# H \subseteq r(R \#H);$

(2) $\bar r_H (R) \# H = r(R \#H)$  iff there exists an $H$-ideal
$I$  of $R$ such that $r(R \#H) = I \#H$;

(3)  $R$ is an $ \bar r_H$-$H$-module algebra iff $r(R \#H) =R\# H$;

(4)  $I$ is an $ \bar r_H$-$H$-ideal of $R$ iff $r(I \#H)= I \#H$;

(5) $r(\bar r _H (R) \# H) = \bar r _H(R) \# H$ ;

(6) $r(R \# H) = \bar r _H(R) \# H$  iff
 $r(\bar r _H (R) \# H) =  r (R \# H).$

\end {Proposition}

{\bf Proof.} (1)  It is similar to the proof of Proposition \ref
{10.4.3}.

(2) It is a straightforward verification.

(3) If $R$ is an $\bar r_H$-module algebra, then $R \# H \subseteq
r(R \# H)$
 by part (1). Thus $R\# H = r(R \# H)$. The sufficiency is obvious.

(4), (5) and (6) immediately follow from part (3) .
\begin{picture}(8,8)\put(0,0){\line(0,1){8}}\put(8,8){\line(0,-1){8}}\put(0,0){\line(1,0){8}}\put(8,8){\line(-1,0){8}}\end{picture}

\begin {Theorem}\label {10.6.3}  If $R$ is an algebra over field $k$ with
unit and $H$ is a Hopf algebra over field $k$, then

(1)  $\bar r_{jH} (R) = r_{Hj}(R)$ and  $r_j(r_{Hj}(R)\# H) =
r_{Hj}(R) \#H;$

(2)   $r_j(R\# H) = r_{Hj}(R) \#H$  iff $r_j(r_{Hj}(R)\# H) =
r_{j}(R \#H)$    iff $r_j(r_j(R \# H)\cap R \# H) = r_j(R\# H);$

(3) Furthermore, if   $H$ is finite-dimensional, then
    $r_j(R\# H) = r_{Hj}(R) \#H$  { \ \ \ \ } iff \\
$r_j(r_{jH}(R)\# H) = r_{j}(R \#H).$

\end {Theorem}

{\bf Proof.}  (1) By \cite [Proposition 3.1 ]  {Zh98a}, we have
$\bar r_{jH} (R) = r_{Hj}(R)$. Consequently, $r_j(r_{Hj}(R)  \# H)
=r_{Hj}(R) \#H$ by Proposition \ref {10.6.2} (5).

(2) It immediately  follows from part (1) and Proposition \ref
{10.6.2} (6).

(3) It can easily be proved by part (2) and \cite [Proposition 3.3
(2)] {Zh98a}.
\begin{picture}(8,8)\put(0,0){\line(0,1){8}}\put(8,8){\line(0,-1){8}}\put(0,0){\line(1,0){8}}\put(8,8){\line(-1,0){8}}\end{picture}

The theorem answers the question J.R. Fisher asked in \cite {Fi75} :
When is $r_j(R \#H) = r_{Hj}(R)\#H$ ?

\begin {Proposition}\label {10.6.4}  If $R$ is an algebra over field $k$ with
unit and $H$ is a finite-dimensional Hopf algebra over field $k$,
then

(1)  $\bar r_{bH} (R) = r_{Hb}(R)=r_{bH}(R)$ and  $r_b(r_{Hb}(R)\#
H) = r_{Hb}(R) \#H;$

(2)   $r_b(R\# H) = r_{Hb}(R) \#H$  iff $r_b(r_{Hb}(R)\# H) =
r_{b}(R \#H)$  iff $r_b(r_{bH}(R)\# H) = r_{b}(R \#H)$  iff
$r_b(r_{b}(R \#H) \cap R \# H) = r_{b}(R \#H).$

\end {Proposition}

{\bf Proof.}  (1) By \cite [Proposition 2.4 ]  {Zh98a}, we have
$\bar r_{bH} (R) = r_{Hb}(R)$.  Thus $r_b(r_{Hb}(R)\# H) = r_{Hb}(R
)\#H$  by Proposition \ref {10.6.2} (5).

(2) It follows from part (1) and Proposition \ref {10.6.2} (6) .
\begin{picture}(8,8)\put(0,0){\line(0,1){8}}\put(8,8){\line(0,-1){8}}\put(0,0){\line(1,0){8}}\put(8,8){\line(-1,0){8}}\end{picture}

In fact,  if $H$ is commutative or cocommutative, then $S^2 = id_H$
 by \cite [Proposition 4.0.1] {Sw69a},
and  $H$  is semisimple and cosemisimple iff the character $char k $
of $k$ does  not divides $dim H$   \ \ ( see \cite [Proposition 2
(c)] {Ra94} ). It is clear that if $H$ is a finite-dimensional
commutative or cocommutative Hopf algebra and
 the character $char k $ of $k$ does  not
divides $dim H$, then  $H$ is a finite-dimensional semisimple,
cosemisimple, commutative or cocommutative Hopf algebra.
Consequently, the conditions in Corollary \ref {10.4.9}, Theorem
\ref {10.4.13}, Corollary  \ref{10.4.14} and \ref {10.4.17}
  can be simplified

\chapter{Classical Yang-Baxter Equation And Low Dimensional Triangular
Lie Algebras }\label {c12}

The concept and structures of Lie coalgebras were  introduced and
studied by   W. Michaelis in \cite {Mi80} \cite {Mi85}.
 V.G.Drinfel'd and   A.A.Belavin in    \cite {Dr86}
\cite {BD82}  introduced the notion  of triangular, coboundary  Lie
bialgebra $L$ associated to a solution $r \in L \otimes L$ of the
CYBE and gave a classification of  solutions of CYBE with parameter
for  simple Lie algebras. W. Michaelis in \cite {Mi94} obtained the
structure of a triangular, coboundary Lie bialgebra on any Lie
algebra containing linearly independent elements $a$ and $b$
satisfying $[a,b]= \alpha b$ for some non-zero $\alpha \in k$ by
setting $r= a\otimes b - b \otimes a $.

The Yang-Baxter equation first came up in a paper by Yang as
factorization condition of the scattering S-matrix in the many-body
problem in one dimension and in the work of Baxter on exactly
solvable models in statistical mechanics. It has been playing an
important role in mathematics and physics ( see \cite {BD82} , \cite
{YG89} ). Attempts to find solutions of The Yang-Baxter equation in
a systematic way have led to the theory of quantum groups. The
Yang-Baxter equation is of many forms. The classical Yang-Baxter
equation is one of these.

 In many applications one needs to know the solutions of classical Yang-Baxter
 equation and know if a Lie algebra is a coboundary Lie bialgebra or
 a triangular Lie bialgebra.  A systematic study of low dimensional Lie
 algebras,
  specially, of those Lie algebras that play a role in physics
 (as e.g. $sl (2, {\bf C})$, or the Heisenberg algebra), is very useful.

 In this chapter, we obtain all  solutions of constant classical Yang-Baxter
 equation (CYBE)  in Lie algebra $L$ with dim $L \le 3$
  and  give
 the sufficient and  necessary conditions for  $(L, \hbox {[ \ ]},
 \Delta _r, r)$ to be a coboundary
 (or triangular ) Lie bialgebra.
We find the strongly symmetric elements in $L\otimes L $ and show
they  are all
 solutions of CYBE in $L$ with $dim \ \ L \le 3$.
  Using these conclusions, we study the
 Lie algebra  $sl(2)$.

In order to make the chapter somewhat self-contained, we begin by
recalling the definition of Lie bialgebra.

Let $k$ be a field with char $k \not= 2 $ and $L$ be a Lie algebra
over $k$. If $r = \sum a_i \otimes b_i  \in L \otimes L $ and $x\in
L$, we define
 \begin{eqnarray*}
\left[ r^{12}, r^{13} \right] &:=& \sum_{i, j} [a_i, a_j] \otimes
b_i
\otimes b_j \\
\left[ r^{12}, r^{23} \right] &:=& \sum_{i, j} a_i\otimes [b_i,a_j]
 \otimes  b_j \\
\left[ r^{13}, r^{23} \right] &:=& \sum_{i, j} a_i\otimes a_j \otimes [b_i, b_j] \\
       x \cdot r &:=& \sum _{i} [x,a_i] \otimes b_i + a_i \otimes [x,b_i] \\
      \Delta _r (x) &:=& x \cdot r
\end{eqnarray*}

and call $$     [r^{12}, r^{13}]+  [r^{12}, r^{23}]+ [r^{13},
r^{23}]=0$$ the classical Yang-Baxter equation (CYBE).

Let $$ \tau: L \otimes L \longrightarrow L \otimes L  $$ denote the
natural twist map ( defined by $x \otimes y \mapsto y \otimes x$ ),
 and let $$
 \xi : L \otimes  L \otimes L \longrightarrow L \otimes  L \otimes L  $$
 the map defined by
 $ x \otimes y \otimes z \mapsto  y \otimes z \otimes x$
 for any $x, y, z \in L$.

A vector space $L$ is called a Lie coalgebra, if there exists a
linear map $$\Delta : L \longrightarrow L \otimes L $$ such that

(i)  $Im \Delta \subseteq Im (1- \tau )$ and

(ii) $(1 + \xi + \xi ^2)(1 \otimes \Delta ) \Delta =0$

A vector space $(L, [\hbox { \ }], \Delta )$  is called a Lie
bialgebra if

 (i)  $(L,[\hbox { \ }])$   is a Lie algebra;

(ii)  $(L, \Delta ) $  is a Lie coalgebra; and

(iii) for all $x, y \in L,$
$$ \Delta [x,y]= x \cdot \Delta (y) - y \cdot \Delta (x), $$
 where, for all $x, a_i, b_i  \in L$,

$(L,[\hbox { \ }], \Delta, r )$   is called a coboundary Lie
 bialgebra, if  $(L,[\hbox { \ }], \Delta )$  is  a  Lie
bialgebra and  $r \in Im (1-\tau ) \subseteq L \otimes L $ such that
$$\Delta (x) = x \cdot r$$  for all $x\in L$.
A coboundary Lie bialgebra $(L,  [\hbox { \ }], \Delta, r )$ is
called triangular, if $r$ is a solution of CYBE.

Jacobson gave  a classification of Lie algebras with their dimension
$dim L \le 3$ in \cite [ P 11-14] {Ja62}. We now write their whole
operations as follows:

(I)  If $L$ is an abelian Lie algebra, then its operation is
trivial;

(II) If  dim $ L $= dim $L'=3$, then there exists a base
 $\{ e_1, e_2, e_3 \}$ and $\alpha, \beta \in k$ with $\alpha \beta \not=0$
such that $[e_1, e_2]= e_3, [e_2, e_3]= \alpha  e_1,  [e_3,
e_1]=\beta e_2$;

(III) If  dim $ L =3$ and $L' \subseteq $ the center of $L$ with dim
$L'=1$, then there exist a base
 $\{ e_1, e_2, e_3 \}$ and $\alpha, \beta \in k$ with $\alpha= \beta =0$
such that $[e_1, e_2]= e_3, [e_2, e_3]= \alpha  e_1,  [e_3,
e_1]=\beta e_2;$

(IV) If  $k$ is an algebraically closed field and dim $ L =3$
 with dim $L'=2$, then there exists a base
 $\{ e_1, e_2, e_3 \}$ and $\beta, \delta \in k$ with $\delta \not=0$
such that $[e_1, e_2]= 0, [e_1, e_3]= e_1 + \beta  e_2,  [e_2,
e_3]=\delta e_2;$

(V) If  dim $ L =3$ and $L' \not\subseteq$ the  center of $L$,
 then there exists a base
 $\{ e_1, e_2, e_3 \}$ and $\beta, \delta \in k$ with $\beta= \delta=0$
such that $[e_1, e_2]= 0, [e_1, e_3]=  e_1 + \beta e_2,  [e_2,
e_3]=\delta e_2;$

(VI) If  dim $ L=2$ with dim $L'=1$, then there exists a base
 $\{ e,f \}$
such that $[e,f]= e.$

Remark: In Case (IV), the condition which $k$ is an algebraically
closed field is required. In other word, A lie algebra $L$ with dim
$ L =3$ and  dim $L'=2$, it is possible that there is not the basis
in (IV). We shall consider the general case  in Section 6--9.

\section {The solutions of CYBE }\label {s30}
In this section, we find the general solution of CYBE for Lie
algebra $L$ with $dim \ \ L \le 3.$

\begin {Definition} \label {11.1.1}
  Let  $\{ e_1, e_2, \cdots, e_n \}$ be a base of vector space $V$
and     $r = \sum_{i,j= 1 }^{n} k_{ij}(e_{i} \otimes e_{j}) \in V
\otimes V,$ where   $k_{ij} \in k,$ for  $i, j =1, 2, \cdots, n.$

(i) If  $$k_{ij}=k_{ji}  \hbox { \ \ \ } k_{ij} k_{lm} =
k_{il}k_{jm}$$ for $i, j, l, m = 1, 2, \cdots, n,$   then  $r$ is
called strongly symmetric to the base
 $\{ e_1, e_2, \cdots, e_n \}$;

(ii) If   $k_{ij} =- k_{ji}$, for $i,j = 1, 2, \cdots, n$,
 then  $r$ is called
 skew symmetric.

(iii) Let $dim { \ } V =3$,
 $k_{11}=x, k_{22}=y,  k_{33}=z, k_{12}=p,  k_{21} =q, k_{13}=s,
 k_{31} =t, k_{23}=u,  k_{32}=v, \alpha, \beta \in k.$
If $p=-q,s=-t, u=-v,  x = \alpha z, y= \beta z$ and
$$\alpha \beta z^2 +   \beta s^2 + \alpha u^2 + p^2 =0,$$
then  $r$ is called
 $\alpha , \beta $-skew symmetric to the
base  $\{ e_1, e_2, e_3 \}.$
\end {Definition}
Obviously, skew symmetry does not depend on the particular choice of
bases of $V$ . We need to know if strong  symmetry  depends on the
particular choice of bases of $V$.

\begin {Lemma}
\label {11.1.2} (I)  If $V$ is a finite-dimensional vector space,
then the strong symmetry  does not depend on the particular choice
of bases of $V;$

(II)  Let  $\{ e_1, e_2, e_3 \}$ be a base of vector space $V,$
  $ p,q, s,t, u, v, x, y, z, \alpha, \beta \in k $  and
    $r=  p(e_1 \otimes e_2)+q (e_2 \otimes e_1)
+ s (e_1 \otimes e_3)+t (e_3 \otimes e_1) +u (e_2 \otimes e_3) +v
(e_3 \otimes e_2) +x (e_1 \otimes e_1) +y (e_2 \otimes e_2) + z (e_3
\otimes e_3).$

Then   $r$ is strongly symmetric iff
 $$ xy =p^2, xz = s^2,  yz=u^2, xu=sp, ys=pu, zp=su,
p=q, s=t, u=v$$ iff  $$xy = p^2, xz = s^2,  yz=u^2, xu=sp, p=q, s=t,
u=v$$ iff
 $$ xy =p^2, xz = s^2,  yz=u^2, ys=up,
p=q, s=t, u=v$$ iff
 $$ xy =p^2, xz = s^2,  yz=u^2, zp=su,
p=q, s=t, u=v$$ iff

 $r=   z^{-1}su  (e_1 \otimes e_2)+ z^{-1}su   (e_2 \otimes e_1)
+ s (e_1 \otimes e_3)+ s (e_3 \otimes e_1) +u (e_2 \otimes e_3) +u
(e_3 \otimes e_2) +z^{-1}s^2  (e_1 \otimes e_1) + z^{-1} u^2 (e_2
\otimes e_2) + z (e_3 \otimes e_3)$
                 with $z \not= 0;$   or

  $r= p(e_1 \otimes e_2)+ p (e_2 \otimes e_1)
+x (e_1 \otimes e_1) +x^{-1}p^2 (e_2 \otimes e_2)$ with $x \not= 0;$
or

  $r= y (e_2 \otimes e_2)$.
  \end {Lemma}

{ \bf Proof. } (I)  Let $\{  e_1, e_2, \cdots, e_n \}$  and  $\{
e_1', e_2', \cdots, e_n' \}$ are two bases of $V$  and
    $r = \sum_{i,j= 1 }^{n} k_{ij}(e_{i} \otimes e_{j})$  be strongly symmetric to
the base   $\{  e_1, e_2, \cdots, e_n \}.$ It is sufficient to show
that  $r$ is strongly symmetric to the base $\{  e_1', e_2', \cdots,
e_n' \}.$ Obviously, there exists $q_{ij} \in k$ such that  $e_{i} =
\sum _{s} e_s'q_{si}  $    for $i = 1, 2, \cdots, n.$ By
computation,  we have that
$$r = \sum_{ s, t }( \sum _{i,j} k_{ij} q_{si}q_{tj} )(e_{s} \otimes e_{t}).$$
If we set   $k_{st}'= \sum _{i,j} k_{ij} q_{si}q_{tj}$, then for $l,
m, u, v =1, 2, \cdots , n,  $   we have that
 \begin {eqnarray*}
 k_{lm}'k_{uv}'
  &=&  \sum_{i,j,s,t} k_{ij} k_{st} q_{li}q_{mj} q_{us}q_{vt}   \\
  &=&  \sum_{i,j,s,t} k_{is} k_{jt} q_{li}q_{mj} q_{us}q_{vt} { \ \ }
 (  \hbox { \ \  by  }   k_{ij}k_{st} = k_{is}k_{jt}  ) \\
&=&  k_{lu}'k_{mv}'
     \end {eqnarray*}
 and
 \begin {eqnarray*}
 k_{lm}'&=&k_{ml}',
 \end {eqnarray*}
which implies that
 $r$ is strongly symmetric to
the base $\{  e_1', e_2', \cdots, e_n' \}.$

(II)  We only show that  if
 $$ xy =p^2, xz = s^2,  yz=u^2, xu=sp,
p=q, s=t, u=v,$$  then
$$ ys=pu, zp=su.$$
If $p \not= 0,$  then $ysp=yxu=up^2$  and $usp=uxu = xyz = zp^2,$
which implies $ys=up$  and $us=zp.$ If $p=0,$ then $x=0 $ or $y=0,$
which implies $s=0$ or $u=0.$
                      Consequently,  $ys=up$  and $us=zp.$
\begin{picture}(8,8)\put(0,0){\line(0,1){8}}\put(8,8){\line(0,-1){8}}\put(0,0){\line(1,0){8}}\put(8,8){\line(-1,0){8}}\end{picture}

\begin {Proposition}\label {11.1.3}
Let $L$ be a Lie algebra with $dim { \ } L =2$. Then $r$ is a
solution of CYBE iff $r$ is strongly symmetric or skew symmetric.
\end {Proposition}

{\bf Proof.}   It is not hard since the tensor $r$ has only 4
coefficients.
 \begin{picture}(8,8)\put(0,0){\line(0,1){8}}\put(8,8){\line(0,-1){8}}\put(0,0){\line(1,0){8}}\put(8,8){\line(-1,0){8}}\end{picture}

\begin {Proposition}
\label {11.1.4} Let $L$ be a Lie algebra with a basis $\{ e_1, e_2,
e_3 \}$ such that $[e_1,e_2]= e_3, [e_2,e_3]=\alpha e_1,
[e_3,e_1]=\beta e_2$, where $\alpha, \beta \in k$. Let $ p, q, s, t,
u, v, x, y, z \in k$.

(I)  If $r$ is strongly symmetric or $\alpha, \beta $-skew symmetric
to basis $\{ e_1, e_2, e_3 \},$  then $r$ is a solution of CYBE;

(II)  If  $\alpha \not= 0 , \beta \not=0$, then $r$ is a   solution
of CYBE in $L$ iff   $r$ is strongly symmetric or $\alpha, \beta
$-skew symmetric;

(III) If  $\alpha = \beta =0,$ then $r$ is a solution of CYBE in $L$
iff

 $r= p (e_1 \otimes e_2) +p (e_2 \otimes e_1)+
s (e_1 \otimes e_3)+ t (e_3 \otimes e_1) + u (e_2 \otimes e_3)+ v
(e_3 \otimes e_2) + x (e_1 \otimes e_1) +        y (e_2 \otimes e_2)
 +z (e_3 \otimes e_3)$

with  $p \not=0, p^2 = xy, xu = sp, xv =tp, tu =vs,$ or

 $r= s (e_1 \otimes e_3)+ t (e_3 \otimes e_1) + u (e_2 \otimes e_3)+
 v (e_3 \otimes e_2) +
x (e_1 \otimes e_1) + y (e_2 \otimes e_2)
 +z (e_3 \otimes e_3)$

with $xy =xu=xv=ys=yt=0$ and $tu=vs.$

\end {Proposition}

{ \bf Proof .} Let   $r = \sum_{i,j= 1 }^{3} k_{ij}(e_{i} \otimes
e_{j}) \in L \otimes L$ and  $k_{ij} \in k$ with  $i, j =1, 2, 3.$
It is clear that
\begin {eqnarray*}
\left[ r^{12}, r^{13} \right] &=& \sum _{i,j=1}^3 \sum _{s,t =1}^3
k_{ij}k_{st}
 [e_i,e_s] \otimes e_j \otimes e_t  \\
\left[r^{12}, r^{23} \right] &=& \sum _{i,j=1}^3 \sum _{s,t =1}^3
k_{ij}k_{st}
e_i \otimes[ e_j, e_s ]\otimes e_t \\
\left[r^{13}, r^{23} \right] &=& \sum _{i,j=1}^3 \sum _{s,t =1}^3
k_{ij}k_{st} e_i \otimes e_s \otimes[e_j, e_t].
\end {eqnarray*}
 By computation, for  $i, j, n =1, 2, 3,$   we have that the coefficient of
  $e_j \otimes e_i \otimes e_i$ in
$[r^{12}, r^{13}] $ is zero and
 $e_i \otimes e_i \otimes e_j$
in  $[r^{13}, r^{23}]$ is zero.

We now see the coefficient
 of  $e_i \otimes e_j \otimes e_n$

   in  $[r^{12}, r^{13}] + [r^{12}, r^{23}] +[r^{13}, r^{23}]$.

(1) $e_1 \otimes e_1 \otimes e_1 ( \alpha k_{12}  k_{31}- \alpha
k_{13} k_{21});$

(2) $e_2 \otimes e_2 \otimes e_2 (
 - \beta k_{21}  k_{32}+\beta k_{23} k_{12} );$

(3) $e_3 \otimes e_3 \otimes e_3 ( k_{31}  k_{23}- k_{32} k_{13});$

(4) $e_1 \otimes e_2 \otimes e_3 ( \alpha k_{22}k_{33} - \alpha
k_{32} k_{23} -  \beta k_{11}  k_{33} +\beta k_{13} k_{13} + k_{11}
k_{22} -k_{12} k_{21} );$

(5) $e_2 \otimes e_3 \otimes e_1 ( \beta k_{33}k_{11} - \beta k_{13}
k_{31} - k_{22}  k_{11}+ k_{21} k_{21} + \alpha k_{22}  k_{33}
-\alpha k_{23} k_{32} );$

(6) $e_3 \otimes e_1 \otimes e_2 ( k_{11}k_{22} - k_{21} k_{12} -
\alpha k_{33}  k_{22}+ \alpha k_{32} k_{32} + \beta k_{33}  k_{11}
-\beta k_{31} k_{13} );$

(7) $e_1 \otimes e_3 \otimes e_2 ( -\alpha k_{33}k_{22} + \alpha
k_{23} k_{32} + k_{11}  k_{22}- k_{12} k_{12} -\beta k_{11}  k_{33}
+ \beta k_{13} k_{31} );$

(8) $e_3 \otimes e_2 \otimes e_1 ( -k_{22}k_{11} + k_{12} k_{21} +
\beta k_{33}  k_{11}- \beta k_{31} k_{31} - \alpha k_{33}  k_{22} +
\alpha k_{32} k_{23} );$

(9) $e_2 \otimes e_1 \otimes e_3 ( -\beta k_{11}k_{33} + \beta
k_{31} k_{13}  + \alpha k_{22}  k_{33} - \alpha k_{23} k_{23} -
k_{22}  k_{11} +k_{21} k_{12} );$

(10) $e_1 \otimes e_1 \otimes e_2 ( - \alpha k_{31}k_{22} + \alpha
k_{21} k_{32} + \alpha k_{12}  k_{32} - \alpha k_{13} k_{22});$

(11) $e_2 \otimes e_1 \otimes e_1 (
 -\alpha k_{23}  k_{21}+
\alpha k_{22} k_{31} +\alpha k_{22}  k_{13} - \alpha k_{23} k_{12}
);$

(12) $e_1 \otimes e_1 \otimes e_3 ( \alpha k_{21}k_{33} - \alpha
k_{31} k_{23} + \alpha k_{12}  k_{33}
 - \alpha k_{13} k_{23} );$

(13) $e_3 \otimes e_1 \otimes e_1 (
 \alpha k_{32}  k_{31}- \alpha k_{33} k_{21}
+ \alpha k_{32}  k_{13} -\alpha k_{33} k_{12} );$

(14) $e_2 \otimes e_2 \otimes e_1 ( - \beta k_{12}k_{31} + \beta
k_{32} k_{11} + \beta k_{23}  k_{11} -\beta  k_{21} k_{31});$

(15) $e_1 \otimes e_2 \otimes e_2 (
 -\beta  k_{11}  k_{32}+
\beta k_{13} k_{12} -\beta k_{11}  k_{23} + \beta k_{13} k_{21} );$

(16) $e_2 \otimes e_2 \otimes e_3 ( \beta k_{32}k_{13} -\beta k_{12}
k_{33} -\beta k_{21}  k_{33} +\beta k_{23} k_{13});$

(17) $e_3 \otimes e_2 \otimes e_2 ( -\beta k_{31}  k_{32}+ \beta
k_{33} k_{12} + \beta k_{33}  k_{21} - \beta k_{31} k_{23} );$

(18) $e_3 \otimes e_3 \otimes e_1 ( k_{13}k_{21} - k_{23} k_{11} +
k_{31}  k_{21}- k_{32} k_{11} );$

(19) $e_3 \otimes e_3 \otimes e_2 ( -k_{23}k_{12} + k_{13} k_{22} -
k_{32}  k_{12}+ k_{31} k_{22});$

(20) $e_2 \otimes e_3 \otimes e_3 (
 k_{21}  k_{23} - k_{22} k_{13}
+ k_{21}  k_{32} -k_{22} k_{31} );$

(21) $e_1 \otimes e_3 \otimes e_3 ( - k_{12}  k_{13}+ k_{11} k_{23}
- k_{12}  k_{31} +k_{11} k_{32} );$

(22) $e_1 \otimes e_3 \otimes e_1 ( \alpha k_{23}k_{31} - \alpha
k_{33} k_{21} - k_{12}  k_{11}+ k_{11} k_{21} +\alpha k_{12}  k_{33}
-\alpha k_{13} k_{32} );$

(23) $e_1 \otimes e_2 \otimes e_1 ( -\alpha k_{32}k_{21} +\alpha
k_{22} k_{31} - \beta k_{11}  k_{31} + \beta k_{13} k_{11} -\alpha
k_{13}  k_{22} +\alpha k_{12} k_{23} );$

(24) $e_2 \otimes e_1 \otimes e_2 ( \beta k_{31}k_{12} -\beta k_{11}
k_{32} + \alpha k_{22}  k_{32} -\alpha  k_{23} k_{22} - \beta k_{21}
k_{13} +\beta k_{23} k_{11} );$

(25) $e_2 \otimes e_3 \otimes e_2 ( -\beta k_{13}k_{32} + \beta
k_{33} k_{12} + k_{21}  k_{22}- k_{22} k_{12} -\beta  k_{21}  k_{33}
+\beta k_{23} k_{31} );$

(26) $e_3 \otimes e_2 \otimes e_3 ( k_{12}k_{23} - k_{22} k_{13}
-\beta k_{31}  k_{33}+ \beta k_{33} k_{13} + k_{31}  k_{22} -k_{32}
k_{21} );$

(27) $e_3 \otimes e_1 \otimes e_3 ( -k_{21}k_{13} + k_{11} k_{23} +
\alpha k_{32}  k_{33}- \alpha k_{33} k_{23} + k_{31}  k_{12} -k_{32}
k_{11} ).$

Let $k_{11}=x, k_{22}=y,  k_{33}=z, k_{12}=p,  k_{21} =q, k_{13}=s,
 k_{31} =t, k_{23}=u,  k_{32}=v.$

It follows from (1)--(27)  that

(28) $\alpha pt= \alpha qs ;$

(29) $\beta qv= \beta pu;$

(30) $tu=vs;$

(31) $  \alpha yz - \beta xz + xy - \alpha uv + \beta s^2 - pq =0;$

(32) $ \beta zx -yx+ \alpha yz - \beta st + q^2-\alpha uv =0;$

(33) $ xy  -\alpha zy + \beta zx -  pq+ \alpha v^2 - \beta st =0;$

(34) $- \alpha zy+ xy- \beta xz+  \alpha uv - p^2+ \beta st =0;$

(35) $-xy +\beta xz- \alpha yz + pq - \beta t^2+ \alpha uv =0;$

(36) $ -\beta xz+ \alpha yz - yx+  \beta st -  \alpha u^2 + pq =0;$

(37)  $\alpha (-ty + qv +pv -sy) = 0;$

(38)  $\alpha (-uq + yt +ys -up) = 0;$

(39)  $\alpha (qz- tu +pz -su) = 0;$

(40)  $\alpha (vt- zq +vs -zp) = 0;$

(41)  $\beta (-pt + vx +ux -qt) = 0;$

(42)  $\beta (-xv + sp - xu +sq) = 0;$

(43)  $\beta (vs-pz +us -qz) = 0;$

(44)  $\beta (-tv + zp +zq -tu) = 0;$

(45)  $sq -ux + tq- vx = 0;$

(46)  $-up + sy -vp +ty = 0;$

(47)  $qu-ys + qv  -yt = 0;$

(48)  $-ps + xu-pt  +xv = 0;$

(49)  $ \alpha ut - \alpha zq -px +xq +\alpha pz- \alpha sv = 0;$

(50)  $- \alpha vq + \alpha yt  -\beta xt + \beta sx - \alpha sy
+\alpha pu = 0;$

(51)  $\beta tp - \beta xv + \alpha yv - \alpha uy - \beta qs+ \beta
ux = 0;$

(52)  $ - \beta sv + \beta zp+ qy - yp - \beta qz+ \beta ut = 0;$

(53)  $pu - ys -\beta tz + \beta zs + ty - vq = 0;$

(54)  $-qs + xu + \alpha vz -\alpha zu + tp - vx  = 0.$

It is clear that $r$ is the solution of CYBE iff relations
(28)--(54) hold.

(I)  By computation, we have that if
              $r$ is strongly symmetric or $\alpha, \beta $-skew symmetric,
 then  relations (28)--(54) hold and so  $r$ is a solution of CYBE;

(II)  Let $\alpha \not=0$  and $\beta \not=0.$
 By part (I), we only need show that if $r$ is a solution of  CYBE, then
 $r$ is strongly symmetric or $\alpha, \beta $-skew symmetric.

By computation, we have

(55)  $q^2 =  p^2$   (by  (34)$+$ (32));

(56)  $u^2 =  v^2$  ( by  (33)$+$ (36));

(57)  $t^2 =  s^2$  ( by  (31)$+$ (35) ).

(a) If  $s=t \not=0,$ then we have that
 \begin {eqnarray*}
 p&=&q  \hbox { \ \  by (28)}; \\
 u&=&v \hbox { \ \ by (30) }; \\
 (58) { \ \ \ \ \ \ } yz &=& u^2  \hbox { \ by } (35)-(32); \\
 (59) { \ \ \ \ \ \ } xz &=& s^2  \hbox { \ by } (32)+(33); \\
 (60) { \ \ \ \ \ \ } xy &=& p^2  \hbox { \ by } (36)-(31); \\
(61)  { \ \ \ \ \ \ } xu&=&sp  \hbox { \ by } (41).
 \end {eqnarray*}
Thus $r$ is strongly symmetric.

 Similarly, we can show that if $p=q\not=0,$ or  $u=v \not=0,$
then $r$  is strongly symmetric.

(b) If  $s = -t \not=0,$ then we have that
   \begin {eqnarray*}
 p&=&-q  \hbox { \ \  by (28)}; \\
 u&=&-v \hbox { \ \ by (30) }; \\
  y &=& \beta z  \hbox { \ by } (53); \\
 x &=& \alpha z  \hbox { \ by } (50) ;   \\
 \alpha \beta z^2 + \alpha u^2 + \beta s^2 + p^2 &=& 0   \hbox { \ by } (31).
\end {eqnarray*}
Then  $r$  is  $\alpha , \beta $-skew symmetric.

Similarly, we have that if $p=-q \not=0$ or  $u=-v \not=0,$  then
$r$ is $\alpha, \beta $-skew symmetric.

(c) If $s=t=u=v=p=q=0$,
  then we have that
   \begin {eqnarray*}
xz&=&0 \hbox { \ \ by \ }(32)+(33); \\
 yz&=&0 \hbox { \ \ by } (32)+(31);   \\
xy&=& 0 \hbox  { \ \ by } (31).
\end {eqnarray*}
Thus  $r$  is strongly symmetric.

(III) Let $\alpha = \beta =0.$ It is clear that the system of
equations  (28)--(54) is equivalent to the below
$$  \left    \{  \begin{array} {l}
 tu=vs, q=p, p^2=xy \\
(s+t)p -(u+v )x =0  \\
(s+t)y-(u+v )p =0  \\
(t-s)y+ (u-v)p  =0  \\
(t-s)p +(u-v )x =0  \\
\end{array} \right.$$

It is easy to check that  if $r$ is  one  case in  part (III), then
the system of equations  hold. Conversely, if $r$ is a solution of
CYBE, we shall show that $r$ is one of the two cases in  part (III).
If $p \not=0$, then  $r$  is the first case in part (III). If $p=0$,
then $r$ is the second case in part (III).
        \begin{picture}(8,8)\put(0,0){\line(0,1){8}}\put(8,8){\line(0,-1){8}}\put(0,0){\line(1,0){8}}\put(8,8){\line(-1,0){8}}\end{picture}

In particular, Proposition \ref {11.1.4} implies:

\begin {Example} \label {11.1.5}
 Let $$  sl(2) := \{ x \mid  x  \hbox{ \ is a }  2 \times 2
\hbox { \ matrix
 with trace zero over \  } k   \}$$
 and
 $$e_1 = \left ( \begin {array} {cc}
 0 & 1\\
 1&0
 \end {array}
 \right ),
 e_2 = \left ( \begin {array} {cc}
 0 & -1\\
 1&0
 \end {array}
 \right ),
 e_3 = \left ( \begin {array} {cc}
 2 & 0\\
 0&-2
 \end {array}
 \right ).$$
 Thus  $L$ is a Lie algebra (defined by $[x,y] = xy -yx$)
 with a basis $\{ e_1, e_2, e_3 \}.$ It is clear that
 $$[e_1, e_2] = e_3, [e_2, e_3] = 4e_1, [e_3, e_1] =-4 e_2.$$
Consequently,
 $r$ is a solution of CYBE iff $r$  is strongly symmetric or
4, $-4$-skew symmetric to the basis $\{ e_1, e_2, e_3 \}.$
\end {Example}

\begin {Proposition}
\label {11.1.6}  Let $L$ be a Lie algebra with a basis $\{ e_1, e_2,
e_3 \}$ such that $[e_1,e_2]=0,  [e_1,e_3]= e_1+ \beta e_2,
[e_2,e_3]=\delta e_2$, where $\beta, \delta \in k$. Let $ p, q, s,
t, u, v, x, y, z \in k$.

(I)  If  $r$ is strongly symmetric,  then $r$ is a solution of CYBE.

(II)  If  $\beta = 0 , \delta \not=0$, then $r$   is a solution of
CYBE in $L$ iff $r$ is  strongly symmetric, or

  $r=  p(e_1 \otimes e_2)+ q (e_2 \otimes e_1)
+ s (e_1 \otimes e_3)- s (e_3 \otimes e_1) +u (e_2 \otimes e_3) -u
(e_3 \otimes e_2) +x (e_1 \otimes e_1) +y (e_2 \otimes e_2) $

with $xu= xs=ys=yu= (1- \delta)us=(1+\delta)s (q +p) = (1+\delta) u
(q +p)= 0$

(III)  If   $\beta \not= 0$ and  $\delta =1,$   then $r$ is a
solution of CYBE in $L$  iff
 $r$ is strongly symmetric, or

  $r=  p(e_1 \otimes e_2)+ q (e_2 \otimes e_1)
+u (e_2 \otimes e_3) -u (e_3 \otimes e_2) +x (e_1 \otimes e_1) +y
(e_2 \otimes e_2)$

with $xu= yu=  u (q +p)= 0$

(IV)   If  $ \beta = \delta = 0,$ then $r$ is a  solution  of CYBE
in $L$ iff

  $r=  p(e_1 \otimes e_2)+ q (e_2 \otimes e_1)
+ s (e_1 \otimes e_3)+s (e_3 \otimes e_1) +u (e_2 \otimes e_3) +v
(e_3 \otimes e_2) +x (e_1 \otimes e_1)+ y (e_2 \otimes e_2) + z (e_3
\otimes e_3)$

with $ z \not=0, zp=vs, zq =us, zx =s^2,$  or

  $r=  p(e_1 \otimes e_2)+ q (e_2 \otimes e_1)
+ s (e_1 \otimes e_3)-s (e_3 \otimes e_1) +u (e_2 \otimes e_3) +v
(e_3 \otimes e_2) +x (e_1 \otimes e_1)+ y (e_2 \otimes e_2)$

with $ us=vs=xs= xu=xv=0, up=qv,  s(p+q)=0$.
\end {Proposition}

{ \bf Proof.} Let   $r = \sum_{i,j= 1 }^{3} k_{ij}(e_{i} \otimes
e_{j}) \in L \otimes L$ and  $k_{ij} \in k$ with  $i, j =1, 2, 3.$
 By computation, for all $i, j, n = 1, 2,  3,$   we have that the coefficient of
  $e_j \otimes e_i \otimes e_i$ in
$[r^{12}, r^{13}] $ is zero and
 $e_i \otimes e_i \otimes e_j$
in  $[r^{13}, r^{23}]$ is zero.

We now see the coefficient
 of  $e_i \otimes e_j \otimes e_n$

   in  $[r^{12}, r^{13}] + [r^{12}, r^{23}] +[r^{13}, r^{23}]$.

(1) $e_1 \otimes e_1 \otimes e_1 (
 -k_{13}k_{11} +  k_{11} k_{31})$;

(2) $e_2 \otimes e_2 \otimes e_2 (
 -\beta k_{23}k_{12} + \beta k_{21} k_{32}
  - \delta k_{23}k_{22} + \delta k_{22} k_{32})$;

(3) $e_3 \otimes e_3 \otimes e_3 ( 0);$

(4) $e_1 \otimes e_2 \otimes e_3 (
 -k_{32}k_{13} +  k_{12} k_{33} -  \beta k_{13}  k_{13}+  \beta k_{11} k_{33}
- \delta k_{13}  k_{23} + \delta k_{12} k_{33} );$

(5) $e_2 \otimes e_3 \otimes e_1 (
 -\beta k_{33}k_{11} + \beta k_{13} k_{31}
 - \delta k_{33}  k_{21}+ \delta k_{23} k_{31}
 -  k_{23}  k_{31} + k_{21} k_{33} );$

(6) $e_3 \otimes e_1 \otimes e_2 (
 -k_{33}k_{12} +  k_{31} k_{32} -  \beta k_{33}  k_{11}+  \beta k_{31} k_{13}
- \delta k_{33}  k_{12} + \delta k_{32} k_{13} );$

(7) $e_1 \otimes e_3 \otimes e_2 (
 -k_{33}k_{12} +  k_{13} k_{32} -  \beta k_{13}  k_{32}+  \beta k_{11} k_{33}
- \delta k_{13}  k_{32} + \delta k_{12} k_{33} );$

(8) $e_3 \otimes e_2 \otimes e_1 (
 -\beta k_{33}k_{11} + \beta k_{31} k_{31} -  \delta k_{33}  k_{21}
 +  \delta k_{32} k_{31}
- k_{33}  k_{21} +  k_{31} k_{23} );$

(9) $e_2 \otimes e_1 \otimes e_3 (
 -\beta k_{31}k_{13} + \beta  k_{11} k_{33}
 - \delta k_{31}  k_{23}+  \delta k_{21} k_{33}
-  k_{23}  k_{13} +  k_{21} k_{33} );$

(10) $e_1 \otimes e_1 \otimes e_2 (
 -k_{31}k_{12} +  k_{11} k_{32} -   k_{13}  k_{12}+   k_{11} k_{32});$

(11) $e_2 \otimes e_1 \otimes e_1 (
 -k_{23}k_{11} +  k_{21} k_{31} -   k_{23}  k_{11}+  k_{21} k_{13});$

(12) $e_1 \otimes e_1 \otimes e_3 (
 -k_{31}k_{13} +  k_{11} k_{33} -   k_{13}  k_{13}+   k_{11} k_{33}
 );$

(13) $e_3 \otimes e_1 \otimes e_1 (
 -k_{33}k_{11} +  k_{31} k_{31} -   k_{33}  k_{11}+  k_{31} k_{13}
);$

(14) $e_2 \otimes e_2 \otimes e_1 (
 -\beta k_{32}k_{11} +  \beta k_{12} k_{31} -  \delta k_{32}  k_{21}
 +  \delta k_{22} k_{31}
- \beta k_{23}  k_{11} + \beta k_{21} k_{31} - \delta k_{23}  k_{21}
+ \delta k_{22} k_{31}
 );$

(15) $e_1 \otimes e_2 \otimes e_2 (
 -\beta k_{13}k_{12} + \beta  k_{11} k_{32} -  \delta k_{13}  k_{22}
 +  \delta k_{12} k_{32}
- \beta k_{13}  k_{21} + \beta k_{11} k_{23} - \delta k_{13}  k_{22}
+ \delta k_{12} k_{23}
 );$

(16) $e_2 \otimes e_2 \otimes e_3 (
 -\beta k_{32}k_{13} +  \beta k_{12} k_{33} -  \delta k_{32}  k_{23}
 + \delta k_{22} k_{33}
- \beta k_{23}  k_{13} + \beta k_{21} k_{33} - \delta k_{23}  k_{23}
+ \delta k_{22} k_{33} );$

(17) $e_3 \otimes e_2 \otimes e_2 (
 -\beta k_{33}k_{12} + \beta  k_{31} k_{32}- \delta k_{33}k_{22}
 +\delta k_{32}k_{32}
  -  \beta k_{33} k_{21} + \beta k_{31}  k_{23}
  - \delta k_{33} k_{22}
  + \delta k_{32} k_{23}
  );$

(18) $e_3 \otimes e_3 \otimes e_1 (0);$

(19) $e_1 \otimes e_3 \otimes e_3 (0);$

(20) $e_3 \otimes e_3 \otimes e_2 (0);$

(21) $e_2 \otimes e_3 \otimes e_3 (0);$

(22) $e_1 \otimes e_3 \otimes e_1 (
 -k_{33}k_{11} +  k_{13} k_{31} -   k_{13}  k_{31}+   k_{11} k_{33})
 = e_1 \otimes e_3 \otimes e_1 (0)$;

(23) $e_1 \otimes e_2 \otimes e_1 (
 -k_{32}k_{11} +  k_{12} k_{31} -  \beta k_{13}  k_{11}+  \beta k_{11} k_{31}
- \delta k_{13}  k_{21} + \delta k_{12} k_{31} -  k_{13}  k_{21} +
k_{11} k_{23});$

(24) $e_2 \otimes e_1 \otimes e_2 (
 -\beta k_{31}k_{12} + \beta k_{11} k_{32}
 - \delta k_{31}  k_{22}+  \delta k_{21} k_{32}
- k_{23}  k_{12} +  k_{21} k_{32}
 -\beta k_{23}k_{11} + \beta k_{21} k_{13}
  - \delta k_{23}  k_{12}+  \delta k_{22} k_{13}  );$

(25) $e_2 \otimes e_3 \otimes e_2 (
 -\beta k_{33}k_{12} + \beta k_{13} k_{32}
  -  \delta k_{33}  k_{22}+  \delta k_{23} k_{32}
- \beta k_{23}  k_{31} + \beta k_{21} k_{33}
 + \delta k_{33}  k_{22}-  \delta k_{22} k_{33}
              );$

(26) $e_3 \otimes e_2 \otimes e_3 (
 -\beta k_{33}k_{13} + \beta k_{31} k_{33}
  - \delta k_{33}  k_{23}+  \delta k_{32} k_{33}
);$

(27) $e_3 \otimes e_1 \otimes e_3 (
 - k_{33}k_{13} +   k_{31} k_{33}
 );$

Let $k_{11}=x, k_{22}=y,  k_{33}=z, k_{12}=p,  k_{21} =q, k_{13}=s,
 k_{31} =t, k_{23}=u,  k_{32}=v.$

   It follows from (1)--(27) that

(28) $-sx + xt=0$;

(29) $-\beta up +\beta qv -\delta uy + \delta yv=0$;

(30) $- vs + pz -\beta s^2 + \beta xz -\delta su + \delta zp =0$;

(31) $-\beta xz + \beta st - \delta zq + \delta ut - ut   + qz =0$;

(32) $- zp + tv -\beta zx + \beta ts -\delta zp + \delta vs =0$;

(33) $- zp + sv -\beta st + \beta xz -\delta sv + \delta zp =0$;

(34) $- \beta zx + \beta tt - \delta zq + \delta vt - zq + tu =0$;

(35) $- \beta st + \beta xz  -\delta tu + \delta qz -us + qz=0 $;

(36) $- tp + xv - sp + vx =0$;

(37) $- ux + qt-ux + qs =0$;

(38) $- st + xz - s^2 +  xz =0$;

(39) $-zx +tt -zx +st =0$;

(40) $-\beta vx +\beta pt -\delta vq + \delta yt -\beta ux + \beta
qt -\delta uq + \delta yt =0$;

(41) $-\beta sp + \beta xv  -\delta  sy + \delta pv -\beta sq +
\beta xu - \delta sy + \delta pu =0$;

(42) $-\beta vs + \beta pz -\delta vu + \delta yz -\beta su + \beta
zq -\delta u^2 + \delta yz  =0$;

(43) $-\beta  zp + \beta tv - \delta zy + \delta vv  -\beta zq +
\beta tu -\delta zy +\delta vu =0$;

(44) $- vx + pt -\beta sx + \beta xt -\delta sq + \delta pt -sq + xu
=0$;

(45) $- \beta pt + \beta vx -\delta ty + \delta qv  -up + qv -\beta
ux + \beta qs - \delta up + \delta ys =0$;

(46) $-\beta zp + \beta sv -\beta ut + \beta qz =0$;

(47) $- \beta zs + \beta tz -\delta zu + \delta vz =0$;

(48) $- zs + tz =0$.

It is clear that   $r$   is a solution of CYBE iff
 (28)--(48) hold.

(I)  It is trivial.

(II)  Let  $ \beta =0$  and  $\delta \not=0.$ We only show that if
$r$ is a solution of CYBE, then $r$ is one of the two cases in part
(II). If $z\not=0$, we have that $u = v, s=t, yz = u^2, xz = s^2$
and $ux = sp $   by (47), (48), (42), (38) and (36) respectively. It
follows from (30) and (33) that $pz = us$ and from (31) and (34)
that $qz = us.$  Thus $q=p$ and $xy =p^2$ since $xyz =xu^2=spu
=zp^2. $  By Lemma \ref {11.1.2}, $r$ is strongly symmetric.

If $z=0$, we have that $u=-v$ and $s=-t$ by (43), (42), (38) and
(39). It follows from (28), (29), (36), (40), (30), (44) and (45)
that $sx =0, yu=0, xv=0, ys=0, (1-\delta)us=0, (1+\delta)s(p+q)=0$
and $(1+\delta)u(p+q)=0$  respectively. Consequently, $r$ is the
second case.

       (III) Let  $\beta \not=0,$ $\delta =1.$
We only show that if $r$  is a solution of CYBE, then  $r$ is one
case in part (III). If $z\not=0$, we have that $t=s , u=v, q=p, xz
=s^2, zp = su, sp =xu$ and $u ^2 = yz$ by (48), (47), (46), (38),
(32), (36) and (42) respectively. Since $xyz =xu^2 = spu = zp^2$,
$xy=p^2.$   Thus
 $r$ is strongly symmetric by Lemma \ref {11.1.2}.

If $z=0$, then $st =0, t=0$ and $s=0$ by (31), (39) and (38)
respectively. It follows from (42) and (43) that $u=-v$. By (36),
(45) and (29), $xu=0, u(p+q)=0$ and $uy=0$. Consequently,  $r$  is
the second case.

(IV)        Let $\beta = \delta =0.$ We only show that if  $r$ is a
solution of CYBE, then $r$ is one case in part  (IV). If $z\not=0$,
we have that $s=t, s^2=xz, zp=vs $ and $zq=us$ by (48), (38), (32)
and (34) respectively. Consequently,
 $r$ is the first case in part (IV). If $z=0$,

 then $s= -t$ by (38) and (39). It follows from (30), (31), (28), (37), (36),
 (45) and (44) that $vs=0, us=0, xs=0, xu=0, xv=0 , up = qv$ and $s(p+q)=0.$
Consequently,  $r$ is the second case in part (IV).
\begin{picture}(8,8)\put(0,0){\line(0,1){8}}\put(8,8){\line(0,-1){8}}\put(0,0){\line(1,0){8}}\put(8,8){\line(-1,0){8}}\end{picture}

\begin {Corollary}\label {11.1.7}
Let $L$ be a Lie algebra with $dim { \ } L \le 3$  and $r \in L
\otimes L.$ If $r$ is strongly symmetric, then $r$ is a solution of
CYBE in $L$.

\end{Corollary}

{\bf Proof .}   If $k$ is algebraically closed, then  $r$ is a
solution of
 CYBE by Proposition  \ref {11.1.4}, \ref {11.1.6} and
 \ref {11.1.3}.
 If $k$ is not algebraically closed, let
  $P$ be algebraically closure of $k$.
  We can construct a Lie algebra   $ L_P = P \otimes L$ over $P$,
  as in \cite [section 8] {Ja62}.
Set $ \Psi : L \longrightarrow  L_P $  by sending $x $  to
  $1 \otimes x$. It is clear that $L_P$ is a Lie algebra over
  $P$ and  $\Psi$  is  homomorphic with $ker \Psi =0$ over $k.$
  Let  $$\bar r = (\Psi \otimes \Psi )(r).$$
  Obviously, $\bar r$  is  strongly symmetric. Therefore  $\bar r$
 is a solution of CYBE in $L_P$ and so is $r$  in $L$.
\begin{picture}(8,8)\put(0,0){\line(0,1){8}}\put(8,8){\line(0,-1){8}}\put(0,0){\line(1,0){8}}\put(8,8){\line(-1,0){8}}\end{picture}

\section {Coboundary Lie bialgebras}\label {s31}

In this section, using the general solution, which are obtained in
the preceding section,  of CYBE in Lie algebra $L$ with $dim \ \ L
\le 3, $ we give the
       the sufficient and  necessary conditions for  $(L, \hbox {[ \ ]},
 \Delta _r, r)$ to be a coboundary
 (or triangular ) Lie bialgebra.

We now observe the  connection between solutions of CYBE and
triangular Lie bialgebra structures. It is clear that if
 $(L, [ \hbox { \ } ], \Delta _r, r)$   is a triangular
Lie bialgebra then  $r$ is a solution of CYBE. Conversely, if $r$ is
a solution of CYBE and $r$ is skew symmetric with $r \in L \otimes
L$ , then
               $(L, [ \hbox { \ } ], \Delta _r, r)$   is a triangular
Lie bialgebra  by \cite [Proposition 1] {Ta93} .
\begin {Theorem}
\label {11.2.1} Let $L$ be a  Lie algebra with
       a basis $\{ e_1, e_2, e_3 \}$
such that $[e_1,e_2]= e_3, [e_2,e_3]=\alpha e_1, [e_3,e_1]=\beta
e_2$, where $\alpha, \beta \in k$ and $\alpha \beta \not=0$ or
$\alpha =\beta =0$. Let $p, u, s \in k$ and
      $r=  p(e_1 \otimes e_2)- p (e_2 \otimes e_1)
+ s (e_1 \otimes e_3)- s (e_3 \otimes e_1) + u (e_2 \otimes e_3)- u
(e_3 \otimes e_2)$. Then

(I)   $(L, [\hbox { \ }], \Delta _r, r )$ is  a coboundary Lie
bialgebra iff $r$  is skew symmetric;

(II)       $(L, [\hbox { \ }],\Delta _r, r )$ is  a triangular  Lie
bialgebra iff  $$\beta s^2 + \alpha u^2 + p^2 =0. $$
\end {Theorem}

 {\bf Proof.}
(I) It is sufficient to show that
 $$ (1+ \xi + \xi ^2) (1 \otimes \Delta ) \Delta (e_i) =0$$
for $i=1, 2, 3.$ First, by computation, we have that
\begin {eqnarray*}
&{ \ }& (1 \otimes \Delta ) \Delta (e_1) \\
&=& \{ ( e_1 \otimes e_2\otimes e_3 ) (\beta ps -  \beta ps)
 +(e_3 \otimes e_1\otimes e_2)(\beta ps)
 + (e_2 \otimes e_3\otimes e_1) (- \beta sp) \} \\
&+&\{ ( e_1 \otimes e_3\otimes e_2 ) (-\beta ps + \beta ps)
 +(e_3 \otimes e_2\otimes e_1)( - \beta ps)
 + (e_2 \otimes e_1\otimes e_3) (\beta sp) \} \\
&+&\{ ( e_1 \otimes e_1\otimes e_2 ) (\alpha \beta su )
 +(e_1 \otimes e_2\otimes e_1)( - \beta \alpha su)
 + (e_2 \otimes e_1\otimes e_1) (0) \} \\
&+& \{ ( e_1 \otimes e_1\otimes e_3 ) ( - \alpha pu)
 +(e_1 \otimes e_3\otimes e_1)(\alpha  pu)
 + (e_3 \otimes e_1\otimes e_1) (0) \} \\
&+&\{   ( e_2 \otimes e_1\otimes e_2 ) (- \beta ^2 s^2)
 +(e_1 \otimes e_2\otimes e_2)(0)
 + (e_2 \otimes e_2\otimes e_1) ( \beta ^2 s^2) \} \\
&+&\{  ( e_3 \otimes e_3\otimes e_1 ) ( p^2 )
 +(e_3 \otimes e_1\otimes e_3)(- p^2)
 + (e_1 \otimes e_3\otimes e_3) (0) \}.
 \end {eqnarray*}
Thus
 $$ (1+ \xi + \xi ^2) (1 \otimes \Delta ) \Delta (e_1) =0.$$
Similarly, we have that
 $$ (1+ \xi + \xi ^2) (1 \otimes \Delta ) \Delta (e_2) =0, \ \ \
  (1+ \xi + \xi ^2) (1 \otimes \Delta ) \Delta (e_3) =0.$$
Thus  $(L, [\hbox { \ }],\Delta _r,r)$ is a coboundary Lie
bialgebra.

(II) By part (I), it is sufficient to show that $r$ is a solution of
CYBE iff
  $$\beta s^2 + \alpha u^2 + p^2 =0. $$
(a) Let $\alpha \not=0 $ and $\beta \not=0.$ Considering $r$ is skew
symmetric, by Proposition \ref {11.1.4} (II), we have
   that                          $r$ is a solution of CYBE
iff
  $$\beta s^2 + \alpha u^2 + p^2 =0. $$
(b) Let $\alpha = \beta=0.$                  Considering $r$ is skew
symmetric, by Proposition \ref {11.1.4} (III), we have
       that                       $r$ is a solution of CYBE
iff
  $$\beta s^2 + \alpha u^2 + p^2 =0. $$
\begin{picture}(8,8)\put(0,0){\line(0,1){8}}\put(8,8){\line(0,-1){8}}\put(0,0){\line(1,0){8}}\put(8,8){\line(-1,0){8}}\end{picture}

Let $(L,[ \hbox{ \ }])$  be a familiar Lie algebra, namely,
 Euclidean 3-space
under vector cross product. By Theorem \ref {11.2.1},  $(L, [ \hbox
{ \ }], \Delta _r, r) $
 is not  a triangular Lie bialgebra for any   $0 \not= r \in L \otimes L.$
In fact, \cite [Example 2.14] {Mi94} already  contains the essence
of this observation.

\begin {Example} \label {11.2.2} Under Example \ref {11.1.5}, we have the
following:

(i)   $(sl(2), [\hbox { \ }],\Delta _r, r )$ is  a coboundary Lie
bialgebra iff $r$ is skew symmetric;

(ii)       $(sl(2), [\hbox { \ }],\Delta _r, r )$ is  a triangular
Lie bialgebra iff
$$-4 s^2 +4 u^2 + p^2 =0, $$
where    $r=  p(e_1 \otimes e_2)- p (e_2 \otimes e_1) + s (e_1
\otimes e_3)- s (e_3 \otimes e_1) +u (e_2 \otimes e_3) -u (e_3
\otimes e_2)$  and $p, s, u \in k.$
\end {Example}

 \begin {Theorem} \label  {11.2.3}
Let  $L$ be a Lie algebra
 with a basis   $\{ e_1, e_2, e_3 \}$  such that
$$[e_1e_2]= 0, [e_1, e_3]= e_1 + \beta e_2, [e_2, e_3]=\delta e_2,$$
where $\delta, \beta \in k.$  Let $p, s, u \in k$ and
 $r=  p(e_1 \otimes e_2)- p (e_2 \otimes e_1)
+ s (e_1 \otimes e_3)- s (e_3 \otimes e_1) + u (e_2 \otimes e_3)- u
(e_3 \otimes e_2)$.

(I)  $(L, [\hbox  { \ }],\Delta _r, r )$ is a coboundary Lie
bialgebra iff $$ (\delta +1) ((\delta -1)u + \beta s)s =0; $$

(II) If  $\beta =0,$  then
 $(L, [\hbox  { \ }],\Delta _r, r )$ is a triangular Lie bialgebra
 iff $$(1- \delta )  us =0;$$

(III) If   $\beta \not=0$ and  $\delta =1,$ then
 $(L, [\hbox  { \ }],\Delta _r,r)$ is a triangular Lie bialgebra iff
 $$ s =0$$
iff   $(L, [\hbox  { \ }],\Delta _r, r )$ is a coboundary Lie
bialgebra.
\end {Theorem}
{ \bf Proof }.
 (I) We get by computation
 \begin {eqnarray*}
&{ \ }& ( 1 \otimes \Delta ) \Delta (e_1) \\
 &=& \{ ( e_1 \otimes e_1\otimes e_2 ) (\beta \delta ss -  \delta us)
 +(e_1 \otimes e_2\otimes e_1)(- \beta \delta ss + \delta us)
 + (e_2 \otimes e_1\otimes e_1) (0) \} \\
&+&  \{ ( e_2 \otimes e_2\otimes e_1 ) (-\beta us + uu+ \beta ^2 s^2
- \beta su)                                       \\
 &+& (e_2 \otimes e_1\otimes e_2)( \beta us -uu + \beta us - \beta ^2 s^2)
  + (e_1 \otimes e_2\otimes e_2) (0) \} \\
&{ \ }& ( 1 \otimes \Delta ) \Delta (e_2) \\
 &=& \{ ( e_1 \otimes e_1\otimes e_2 ) ( \delta ^2 ss)
 +(e_1 \otimes e_2\otimes e_1)(- \delta ^2 s^2 )
 + (e_2 \otimes e_1\otimes e_1) (0) \} \\
 &+& \{ ( e_2 \otimes e_2\otimes e_1 ) ( \delta \beta ss - \delta su)
  +(e_2 \otimes e_1\otimes e_2)( \delta  su - \delta \beta ss )
 + (e_1 \otimes e_2\otimes e_2) (0) \} \\
&{ \ }& ( 1 \otimes \Delta ) \Delta (e_3)  \\
 &=& \{ ( e_1 \otimes e_2\otimes e_3 ) (\delta su +  \beta ss)
 +(e_2 \otimes e_3\otimes e_1)(-\delta us - \beta ss) \\
 &+& (e_3 \otimes e_1\otimes e_2) ( \delta \delta us + \beta \delta ss
 + \beta ss -su) \} \\
&+&\{ ( e_1 \otimes e_3\otimes e_2 ) (-\delta su-  \beta ss)
 +(e_3 \otimes e_2\otimes e_1)( - \beta \delta  ss - \delta \delta us
 -\beta ss + su ) \\
 &+&  (e_2 \otimes e_1\otimes e_3) (\delta us + \beta ss) \} \\
&+& \{ ( e_1 \otimes e_1\otimes e_2 ) (- \delta  sp - \delta \delta
ps
+ sp + \delta sp ) \\
 &+& (e_1 \otimes e_2\otimes e_1)( \delta ps + \delta \delta ps -\delta sp
 -sp) +
 (e_2 \otimes e_1\otimes e_1) (0) \} \\
&+& \{
 ( e_1 \otimes e_1\otimes e_3 ) ( ss)
 +(e_1 \otimes e_3\otimes e_1)(-ss)
 + (e_3 \otimes e_1\otimes e_1) (0) \} \\
&+& \{   ( e_2 \otimes e_2\otimes e_1 ) (- \beta  ps + pu - \delta
\beta sp
-\delta \delta up - \delta \beta ps + \delta pu - \beta sp - \delta up)  \\
&+& (e_2 \otimes e_1\otimes e_2)( \beta  ps - pu + \delta \beta sp
+\delta \delta up + \delta \beta ps - \delta pu + \beta sp + \delta up) ) \\
&-& (e_1 \otimes e_2\otimes e_2) (0) \} \\
&+& \{ ( e_2 \otimes e_2\otimes e_3 ) ( \delta \delta uu + \delta
\beta us
+ \beta \delta us + \beta \beta ss )   \\
 &+& (e_2 \otimes e_3\otimes e_2)(   -\delta \delta uu - \delta \beta us
- \beta \delta us - \beta \beta ss )
                    +        (e_3 \otimes e_2\otimes e_2) (0) \}.
                    \end {eqnarray*}
Consequently,   $$ (1+ \xi + \xi ^2) (1 \otimes \Delta ) \Delta
(e_1) =0,$$
 $$ (1+ \xi + \xi ^2) (1 \otimes \Delta ) \Delta (e_2) =0$$
 and
 $$ (1+ \xi + \xi ^2) (1 \otimes \Delta ) \Delta (e_3) =0$$
iff
  $$ \delta ^2 us + \delta \beta s^2 + \beta s^2 - us =0.$$
This implies that
    $(L, [\hbox  { \ }],\Delta _r,r)$ is a coboundary Lie bialgebra
    iff
 $$ \delta ^2 us + \delta \beta s^2 + \beta s^2 - us =0   { \ \ . \ \ }$$

 (II) If $(L, [ \hbox { \ } ], \Delta _r, r)$  is a triangular Lie bialgebra,
  then
 $(\delta ^2 -1)us =0$ by part (I). If $(\delta +1) \not=0,$ then
 $( 1 - \delta )us =0.$  If $\delta +1 =0, $  then $(1 - \delta )us =0$
 by Proposition \ref {11.1.6} (II). Conversely,
 if $(1 - \delta )us =0$, then we have
 that $(L, [ \hbox { \ } ], \Delta _r, r)$ is a coboundary Lie bialgebra
 by part (I).
 Since $r$ is skew symmetric, we have that $r$ is a solution of CYBE by
 Proposition   \ref {11.1.6} (II) (IV). Thus
 $(L, [ \hbox { \ } ], \Delta _r, r)$ is a triangular Lie bialgebra.

(III)  Similarly, we can show that part (III) holds by part (I) and
Proposition \ref {11.1.6} (III).
\begin{picture}(8,8)\put(0,0){\line(0,1){8}}\put(8,8){\line(0,-1){8}}\put(0,0){\line(1,0){8}}\put(8,8){\line(-1,0){8}}\end{picture}

\begin {Theorem}  \label   {11.2.4}
If $L$ is a Lie algebra with $dim L =2$ and $r \in L \otimes L$,
then
  $(L, [\hbox  { \ }],\Delta _r,r)$ is a triangular  Lie bialgebra   iff
  $(L, [\hbox  { \ }],\Delta _r,r)$ is a coboundary   Lie bialgebra
  iff   $r$ is skew symmetric.
 \end {Theorem}
 { Proof.}  It is an immediate consequence of the main result of \cite {Mi94}.
  \begin{picture}(8,8)\put(0,0){\line(0,1){8}}\put(8,8){\line(0,-1){8}}\put(0,0){\line(1,0){8}}\put(8,8){\line(-1,0){8}}\end{picture}

\section {Strongly symmetric elements  }\label {s32}
In this section, we show  that all strongly symmetric  elements are
solutions of
  constant classical Yang-Baxter equation  in Lie algebra,     or               of
  quantum  Yang-Baxter  equation  in algebra. Throughout this section
  $k$  is an arbitrary  field.

\begin {Definition} \label {11.3.1}
 Let  $\{ e_i \mid i \in \Omega \}$ be a basis of vector space $V$
and     $r = \sum_{i,j \in \Omega }
 k_{ij}(e_{i} \otimes e_{j}) \in V \otimes V,$
where   $k_{ij} \in k.$

If  $$k_{ij}=k_{ji}  \hbox { \ \ and \ \ } k_{ij} k_{lm} =
k_{il}k_{jm}$$ for $i, j, l, m \in \Omega$, then  $r$ is said to be
strongly symmetric to the basis
 $\{ e_i \mid  i \in \Omega  \}$.
 \end {Definition}

Obviously,  symmetry does not depend on the particular choice of
basis of $V$ . We need to know if strong  symmetry  depends on the
particular choice of basis of $V$.

\begin {Lemma}
\label {11.3.2}
 Let  $\{ e_i \mid i \in \Omega \}$ be a basis of vector space $V$
and     $r = \sum_{i,j \in \Omega } k_{ij}(e_{i} \otimes e_{j}) \in
V \otimes V,$ where   $k_{ij} \in k.$  Then

(I) the strong symmetry  does not depend on the particular choice of
basis of $V;$

(II)  $r$  is strongly symmetric iff
   $$   k_{ij} k_{lm} = k_{il}k_{jm}$$
for  any  $i, j, l, m \in \Omega ;$

(III)  $r$ is strongly symmetric iff
   $$   k_{ij} k_{lm} = k_{st}k_{uv}$$
where $\{ s, t, u, v \} =\{i, j, l, m \} $  as sets for  any  $i, j,
l, m \in \Omega ;$

(IV)  $r$ is strongly symmetric iff

$r=0$  or

$r = \sum_{i,j \in \Omega }  k_{i_0 i}k_{i_0j} k_{i_0i_0}^{-1}
        (e_{i} \otimes e_{j}) \in V \otimes V,$
where  $k_{i_0 i_0} \not=0.$
\end {Lemma}

{ \bf Proof. } (I)  Let  $\{  e_i' \mid i \in \Omega ' \}$ be
another  basis of $V$. It is sufficient to show that  $r$ is
strongly symmetric to the basis  $\{  e_i' \mid i \in \Omega ' \}$.
Obviously, there exists $q_{ij} \in k$ such that  $e_{i} = \sum _{s}
e_s'q_{si}  $    for $i \in \Omega.$ By computation,  we have that
$$r = \sum_{ s, t  \in \Omega }( \sum _{i,j \in \Omega '} k_{ij} q_{si}q_{tj} )(e_{s}' \otimes e_{t}').$$
If set   $k_{st}'= \sum _{i,j} k_{ij} q_{si}q_{tj}$, then for $l, m,
u, v \in \Omega '  $   we have that
 \begin {eqnarray*}
 k_{lm}'k_{uv}'
  &=&  \sum_{i,j,s,t} k_{ij} k_{st} q_{li}q_{mj} q_{us}q_{vt}   \\
  &=&  \sum_{i,j,s,t} k_{is} k_{jt} q_{li}q_{mj} q_{us}q_{vt} { \ \ }
 (  \hbox { \ \  by  }   k_{ij}k_{st} = k_{is}k_{jt}  ) \\
&=&  k_{lu}'k_{mv}'
     \end {eqnarray*}
 and
 \begin {eqnarray*}
 k_{lm}'&=&k_{ml}',
 \end {eqnarray*}
which implies that
 $r$ is strongly symmetric to
the basis $\{  e_i' \mid i \in \Omega ' \}$.

(II) The necessity is obvious. Now we show the sufficiency. For any
$i, j \in \Omega  $, if $k_{ii} \not=0$, then $k_{ij} = k_{ji}$
since $ k_{ij}k_{ii} = k_{ii}k_{ji}$. If $k_{ii} = 0$, then $k_{ij}=
k_{ji} = 0$ since $k_{ij} ^2 = k_{ii}k_{jj} =0 = k_{ji}^2.$

(III) The sufficiency follows from part (II). Now we show the
necessity.
 For any $i, j , l, m
\in \Omega  $, we see that
 $k_{ij}k_{lm} = k_{ij}k_{ml} = k_{il}k_{jm} =k_{il}k_{mj} =k_{im}k_{lj}. $

 Similarly, we can get that $k_{ij}k_{lm}=k_{st}k_{uv}$  when $s=j$ or $s=l$
 or $s=m$.

(IV)  If $r$ is a case in part (IV),
 then $r$  is strongly symmetric by straightforward verification. Conversely,
 if $r$ is strongly symmetric and $r\not=0$, then there exists
 $i_0 \in \Omega$  such that  $k_{i_0i_0} \not=0.$
 Thus $k_{ij}  =k_{i_0i}k_{i_0j} k_{i_0i_0}^{-1}$ .
\begin{picture}(8,8)\put(0,0){\line(0,1){8}}\put(8,8){\line(0,-1){8}}\put(0,0){\line(1,0){8}}\put(8,8){\line(-1,0){8}}\end{picture}

If $A$ is an algebra, then we can get a Lie algebra $L(A)$ by
defining $[x,y] =xy -yx $ for any $x, y \in A$.

\begin {Theorem} \label {11.3.3}      Let $A$  be an algebra and $r \in A\otimes A$ be strongly symmetric.
Then  $r$ is a solution of CYBE in $L(A)$; meantime, it is also a
solution of  QYBE in $A$.
\end {Theorem}
{\bf Proof.}   Let $\{ e_i \mid i \in \Omega \}$  be a basis of $A$
and the multiplication of $A$ be defined as follows:
$$e_i e_j = \sum _m a_{ij}^m e_m $$
for any $i, j \in \Omega , a^m_{ij} \in k.$

  We first show that $r$ is  a solution of CYBE in Lie algebra $L(A)$.
   By computation, for all $i, j, l \in \Omega ,$   we have that the coefficient of
  $e_i \otimes e_j \otimes e_l$ in
$[r^{12}, r^{13}] $ is
\begin {eqnarray*}
&{}&\sum _{s t} a_{s,t} ^ik_{sj}k_{tl}-a_{ts} ^ik_{sj}k_{tl}   \\
&=& \sum _{s, t} a_{st}^i k_{sj}k_{tl}-
\sum _{s,t} a_{st}^i k_{tj}k_{sl}    \\
&=& 0   \hbox { \ \ \ \ \ \  by Lemma \ref {11.3.2} (III)}.
\end {eqnarray*}
Similarly,  we have that the coefficients of
  $e_i \otimes e_j \otimes e_l$ in
$[r^{12}, r^{23}] $ and
 $[r^{13}, r^{23}] $ are zero.
Thus $r$ is  a solution  of CYBE in $L(A)$.

  Next we show that $r$ is a solution of QYBE in algebra $A$.
We denote $r$ by $R$ since  the solution of QYBE is usually denoted
by $R$.

 By computation, for all $i, j, l \in \Omega ,$   we have that the coefficient of
  $e_i \otimes e_j \otimes e_l$ in
$R^{12}R^{13} R^{23} $ is
\begin {eqnarray*}
&{}&\sum _{s, t, u,v, w,m} a_{st} ^i a_{um} ^j a_{vw} ^l
k_{su}k_{tv}k_{mw}
\end {eqnarray*}

and the coefficient of
  $e_i \otimes e_j \otimes e_l$ in
$R^{23}R^{13} R^{12} $ is
\begin {eqnarray*}
&{}&\sum _{s, t, u,v, w,m} a_{tm} ^i a_{sw} ^j
a_{uv} ^l k_{su}k_{tv}k_{mw} \\
&=& \sum _{s, t, u,v, w,m} a_{st} ^i a_{um} ^j a_{vw} ^l
k_{su}k_{tv}k_{mw}  \hbox { \ \ \  by Lemma \ref {11.3.2} (III)}.
\end {eqnarray*}
 Thus $r$ is a solution of QYBE in algebra $A.$  \begin{picture}(8,8)\put(0,0){\line(0,1){8}}\put(8,8){\line(0,-1){8}}\put(0,0){\line(1,0){8}}\put(8,8){\line(-1,0){8}}\end{picture}

\begin {Corollary}  \label {11.3.4}   Let $L$ be a Lie algebra and $r\in
L\otimes L$ be strongly symmetric. Then $r$ is a solution of CYBE.
\end {Corollary}

{\bf Proof:} Let $A$ denote the universal enveloping algebra of $L$
. It is clear that  $r$ is also strongly symmetric in $A \otimes A$.
Thus

 $$     [r^{12}, r^{13}]+  [r^{12}, r^{23}]+ [r^{13}, r^{23}]=0$$
by Theorem \ref {11.3.3} . Consequently, $r$ is a solution of CYBE
in $L$.
\begin{picture}(8,8)\put(0,0){\line(0,1){8}}\put(8,8){\line(0,-1){8}}\put(0,0){\line(1,0){8}}\put(8,8){\line(-1,0){8}}\end{picture}

\section {The solutions of CYBE  in the case of $char k =2$ }\label {s33}
In this section, we find the general solution of CYBE for Lie
algebra $L$ over field $k$ of characteristic 2 with $dim \ \ L \le
3.$ Throughout  this section, the characteristic of $k$ is 2.

\begin {Proposition}\label {11.4.3}
Let $L$ be a Lie algebra with $dim { \ } L =2$. Then $r$ is a
solution of CYBE iff $r$ is  symmetric.
\end {Proposition}

{\bf Proof.}   The proof  is not hard since the tensor $r$ has only
4 coefficients.
\begin{picture}(8,8)\put(0,0){\line(0,1){8}}\put(8,8){\line(0,-1){8}}\put(0,0){\line(1,0){8}}\put(8,8){\line(-1,0){8}}\end{picture}

  If $r = x (e _1 \otimes e_1) + y (e _2 \otimes e_2) +
z (e _3 \otimes e_3) + p (e _1 \otimes e_2) +p(e _2 \otimes e_1)+ s
(e _1 \otimes e_3)+s(e _3 \otimes e_1)+ u (e _2 \otimes e_3) +u (e
_3 \otimes e_2)$  with  $\alpha y z + \beta xz + xy +
  \beta s^2 + \alpha u^2 + p^2 =0,$
then  $r$ is called
 $\alpha , \beta $-symmetric to the
basis  $\{ e_1, e_2, e_3 \}.$

\begin {Proposition}
\label {11.4.4} Let $L$ be a Lie algebra with a basis $\{ e_1, e_2,
e_3 \}$ such that $[e_1,e_2]= e_3, [e_2,e_3]=\alpha e_1,
[e_3,e_1]=\beta e_2$, where $\alpha, \beta \in k$. Let $ p, q, s, t,
u, v, x, y, z \in k$.

Then  (I) if $r$  is $\alpha , \beta $-symmetric to basis $\{ e_1,
e_2, e_3 \}$, then $r$ is a solution of CYBE in $L$;

(II) if $\alpha \beta \not=0 $, then  $r$ is a solution of the CYBE
in $L$ iff  $r$ is
 $\alpha, \beta $- symmetric
to the basis $\{ e_1, e_2, e_3 \};$

(III) if $\alpha = \beta =0$, then $r$ is a solution of the CYBE in
$L$ iff

 $r= p (e_1 \otimes e_2) +p (e_2 \otimes e_1)+
s (e_1 \otimes e_3)+ t (e_3 \otimes e_1) + u (e_2 \otimes e_3)+ v
(e_3 \otimes e_2) + x (e_1 \otimes e_1) +        y (e_2 \otimes e_2)
 +z (e_3 \otimes e_3)$

with $tu=vs, \ \  xy=p^2, \ \  (s+t)p + (u+v)x=0  \   \ \ (s+t)y +
(u+v)p =0.$

\end {Proposition}

{ \bf Proof .} Let   $r = \sum_{i,j= 1 }^{3} k_{ij}(e_{i} \otimes
e_{j}) \in L \otimes L$ and  $k_{ij} \in k,$ with  $i, j =1, 2, 3.$
It is clear that
\begin {eqnarray*}
\left[ r^{12}, r^{13} \right] &=& \sum _{i,j=1}^3 \sum _{s,t =1}^3
k_{ij}k_{st}
 [e_i,e_s] \otimes e_j \otimes e_t\ ,  \\
\left[r^{12}, r^{23} \right] &=& \sum _{i,j=1}^3 \sum _{s,t =1}^3
k_{ij}k_{st}
e_i \otimes[ e_j, e_s ]\otimes e_t \ ,\\
\left[r^{13}, r^{23} \right] &=& \sum _{i,j=1}^3 \sum _{s,t =1}^3
k_{ij}k_{st} e_i \otimes e_s \otimes[e_j, e_t].
\end {eqnarray*}
 By computation, for all $i, j, n =1, 2, 3,$   we have that the coefficient of
  $e_j \otimes e_i \otimes e_i$ in
$[r^{12}, r^{13}] $ is zero and
 $e_i \otimes e_i \otimes e_j$
in  $[r^{13}, r^{23}]$ is zero.

We now see the coefficient
 of  $e_i \otimes e_j \otimes e_n$

   in  $[r^{12}, r^{13}] + [r^{12}, r^{23}] +[r^{13}, r^{23}]$.

(1) $e_1 \otimes e_1 \otimes e_1 ( \alpha k_{12}  k_{31}- \alpha
k_{13} k_{21});$

(2) $e_2 \otimes e_2 \otimes e_2 (
 - \beta k_{21}  k_{32}+\beta k_{23} k_{12} );$

(3) $e_3 \otimes e_3 \otimes e_3 ( k_{31}  k_{23}- k_{32} k_{13});$

(4) $e_1 \otimes e_2 \otimes e_3 ( \alpha k_{22}k_{33} - \alpha
k_{32} k_{23} -  \beta k_{11}  k_{33} +\beta k_{13} k_{13} + k_{11}
k_{22} -k_{12} k_{21} );$

(5) $e_2 \otimes e_3 \otimes e_1 ( \beta k_{33}k_{11} - \beta k_{13}
k_{31} - k_{22}  k_{11}+ k_{21} k_{21} + \alpha k_{22}  k_{33}
-\alpha k_{23} k_{32} );$

(6) $e_3 \otimes e_1 \otimes e_2 ( k_{11}k_{22} - k_{21} k_{12} -
\alpha k_{33}  k_{22}+ \alpha k_{32} k_{32} + \beta k_{33}  k_{11}
-\beta k_{31} k_{13} );$

(7) $e_1 \otimes e_3 \otimes e_2 ( -\alpha k_{33}k_{22} + \alpha
k_{23} k_{32} + k_{11}  k_{22}- k_{12} k_{12} -\beta k_{11}  k_{33}
+ \beta k_{13} k_{31} );$

(8) $e_3 \otimes e_2 \otimes e_1 ( -k_{22}k_{11} + k_{12} k_{21} +
\beta k_{33}  k_{11}- \beta k_{31} k_{31} - \alpha k_{33}  k_{22} +
\alpha k_{32} k_{23} );$

(9) $e_2 \otimes e_1 \otimes e_3 ( -\beta k_{11}k_{33} + \beta
k_{31} k_{13}  + \alpha k_{22}  k_{33} - \alpha k_{23} k_{23} -
k_{22}  k_{11} +k_{21} k_{12} );$

(10) $e_1 \otimes e_1 \otimes e_2 ( - \alpha k_{31}k_{22} + \alpha
k_{21} k_{32} + \alpha k_{12}  k_{32} - \alpha k_{13} k_{22});$

(11) $e_2 \otimes e_1 \otimes e_1 (
 -\alpha k_{23}  k_{21}+
\alpha k_{22} k_{31} +\alpha k_{22}  k_{13} - \alpha k_{23} k_{12}
);$

(12) $e_1 \otimes e_1 \otimes e_3 ( \alpha k_{21}k_{33} - \alpha
k_{31} k_{23} + \alpha k_{12}  k_{33}
 - \alpha k_{13} k_{23} );$

(13) $e_3 \otimes e_1 \otimes e_1 (
 \alpha k_{32}  k_{31}- \alpha k_{33} k_{21}
+ \alpha k_{32}  k_{13} -\alpha k_{33} k_{12} );$

(14) $e_2 \otimes e_2 \otimes e_1 ( - \beta k_{12}k_{31} + \beta
k_{32} k_{11} + \beta k_{23}  k_{11} -\beta  k_{21} k_{31});$

(15) $e_1 \otimes e_2 \otimes e_2 (
 -\beta  k_{11}  k_{32}+
\beta k_{13} k_{12} -\beta k_{11}  k_{23} + \beta k_{13} k_{21} );$

(16) $e_2 \otimes e_2 \otimes e_3 ( \beta k_{32}k_{13} -\beta k_{12}
k_{33} -\beta k_{21}  k_{33} +\beta k_{23} k_{13});$

(17) $e_3 \otimes e_2 \otimes e_2 ( -\beta k_{31}  k_{32}+ \beta
k_{33} k_{12} + \beta k_{33}  k_{21} - \beta k_{31} k_{23} );$

(18) $e_3 \otimes e_3 \otimes e_1 ( k_{13}k_{21} - k_{23} k_{11} +
k_{31}  k_{21}- k_{32} k_{11} );$

(19) $e_3 \otimes e_3 \otimes e_2 ( -k_{23}k_{12} + k_{13} k_{22} -
k_{32}  k_{12}+ k_{31} k_{22});$

(20) $e_2 \otimes e_3 \otimes e_3 (
 k_{21}  k_{23} - k_{22} k_{13}
+ k_{21}  k_{32} -k_{22} k_{31} );$

(21) $e_1 \otimes e_3 \otimes e_3 ( - k_{12}  k_{13}+ k_{11} k_{23}
- k_{12}  k_{31} +k_{11} k_{32} );$

(22) $e_1 \otimes e_3 \otimes e_1 ( \alpha k_{23}k_{31} - \alpha
k_{33} k_{21} - k_{12}  k_{11}+ k_{11} k_{21} +\alpha k_{12}  k_{33}
-\alpha k_{13} k_{32} );$

(23) $e_1 \otimes e_2 \otimes e_1 ( -\alpha k_{32}k_{21} +\alpha
k_{22} k_{31} - \beta k_{11}  k_{31} + \beta k_{13} k_{11} -\alpha
k_{13}  k_{22} +\alpha k_{12} k_{23} );$

(24) $e_2 \otimes e_1 \otimes e_2 ( \beta k_{31}k_{12} -\beta k_{11}
k_{32} + \alpha k_{22}  k_{32} -\alpha  k_{23} k_{22} - \beta k_{21}
k_{13} +\beta k_{23} k_{11} );$

(25) $e_2 \otimes e_3 \otimes e_2 ( -\beta k_{13}k_{32} + \beta
k_{33} k_{12} + k_{21}  k_{22}- k_{22} k_{12} -\beta  k_{21}  k_{33}
+\beta k_{23} k_{31} );$

(26) $e_3 \otimes e_2 \otimes e_3 ( k_{12}k_{23} - k_{22} k_{13}
-\beta k_{31}  k_{33}+ \beta k_{33} k_{13} + k_{31}  k_{22} -k_{32}
k_{21} );$

(27) $e_3 \otimes e_1 \otimes e_3 ( -k_{21}k_{13} + k_{11} k_{23} +
\alpha k_{32}  k_{33}- \alpha k_{33} k_{23} + k_{31}  k_{12} -k_{32}
k_{11} ).$

Let $k_{11}=x, k_{22}=y,  k_{33}=z, k_{12}=p,  k_{21} =q, k_{13}=s,
 k_{31} =t, k_{23}=u,  k_{32}=v.$

It follows from (1)--(27)  that

(28) $\alpha pt= \alpha qs; $

(29) $\beta qv= \beta pu;$

(30) $tu=vs;$

(31) $  \alpha yz - \beta xz + xy - \alpha uv + \beta s^2 - pq =0;$

(32) $ \beta zx -yx+ \alpha yz - \beta st + q^2-\alpha uv =0;$

(33) $ xy  -\alpha zy + \beta zx -  pq+ \alpha v^2 - \beta st =0;$

(34) $- \alpha zy+ xy- \beta xz+  \alpha uv - p^2+ \beta st =0;$

(35) $-xy +\beta xz- \alpha yz + pq - \beta t^2+ \alpha uv =0;$

(36) $ -\beta xz+ \alpha yz - yx+  \beta st -  \alpha u^2 + pq =0;$

(37)  $\alpha (-ty + qv +pv -sy) = 0;$

(38)  $\alpha (-uq + yt +ys -up) = 0;$

(39)  $\alpha (qz- tu +pz -su) = 0;$

(40)  $\alpha (vt- zq +vs -zp) = 0;$

(41)  $\beta (-pt + vx +ux -qt) = 0;$

(42)  $\beta (-xv + sp - xu +sq) = 0;$

(43)  $\beta (vs-pz +us -qz) = 0;$

(44)  $\beta (-tv + zp +zq -tu) = 0;$

(45)  $sq -ux + tq- vx = 0;$

(46)  $-up + sy -vp +ty = 0;$

(47)  $qu-ys + qv  -yt = 0;$

(48)  $-ps + xu-pt  +xv = 0;$

(49)  $ \alpha ut - \alpha zq -px +xq +\alpha pz- \alpha sv = 0;$

(50)  $- \alpha vq + \alpha yt  -\beta xt + \beta sx - \alpha sy
+\alpha pu = 0;$

(51)  $\beta tp - \beta xv + \alpha yv - \alpha uy - \beta qs+ \beta
ux = 0;$

(52)  $ - \beta sv + \beta zp+ qy - yp - \beta qz+ \beta ut = 0;$

(53)  $pu - ys -\beta tz + \beta zs + ty - vq = 0;$

(54)  $-qs + xu + \alpha vz -\alpha zu + tp - vx  = 0.$

It is clear that $r$ is the solution of CYBE iff relations
(28)--(54) hold.

(55)  $q^2 =  p^2$   (by  (34)$+$ (32));

(56)  $\alpha u^2 = \alpha  v^2$  ( by  (33)$+$ (36));

(57)  $\beta t^2 = \beta  s^2$  ( by  (31)$+$ (35) );

  By computation, we have that if $r$ is  $\alpha , \beta$-symmetric then
relations (28)--(54) hold and so $r$ is a solution of CYBE. That is,
part (I) holds.

(II)  By part (I), we have that the sufficiency of part (II) holds.
For the necessity,  if $r$ is a solution of CYBE, then  $p=q, s=t$
and $u=v$ by (55)--(57),  since char  $k =2$. By relation (31), we
have that
$$\alpha y z + \beta xz +xy +   \beta s^2 + \alpha u^2 + p^2 =0.$$
Thus $r$ is  $\alpha ,\beta $-symmetric.

(III)  It is trivial .
\begin{picture}(8,8)\put(0,0){\line(0,1){8}}\put(8,8){\line(0,-1){8}}\put(0,0){\line(1,0){8}}\put(8,8){\line(-1,0){8}}\end{picture}

In particular, Proposition \ref {11.4.4} implies:

\begin {Example} \label {11.4.5}
 Let $$  sl(2) := \{ x \mid  x  \hbox{ \ is a }  2 \times 2
\hbox { \ matrix
 with trace zero over \  } k   \}$$
 and
 $$e_1 = \left ( \begin {array} {cc}
 0 & 0\\
 1&0
 \end {array}
 \right ),
 e_2 = \left ( \begin {array} {cc}
 0 & 1\\
 0&0
 \end {array}
 \right ),
 e_3 = \left ( \begin {array} {cc}
 1 & 0\\
 0&1
 \end {array}
 \right ).$$
 Thus  $L$ is a Lie algebra (defined by $[x,y] = xy -yx$)
 with a basis $\{ e_1, e_2, e_3 \}.$ It is clear that
 $$[e_1, e_2] = e_3, [e_2, e_3] = 0e_1, [e_3, e_1] =0 e_2.$$
Consequently,
 $r$ is a solution of CYBE iff

 $r= p (e_1 \otimes e_2) +p (e_2 \otimes e_1)+
s (e_1 \otimes e_3)+ t (e_3 \otimes e_1) + u (e_2 \otimes e_3)+ v
(e_3 \otimes e_2) + x (e_1 \otimes e_1) +        y (e_2 \otimes e_2)
 +z (e_3 \otimes e_3)$

with $tu=vs, \ \  xy=p^2, \ \  (s+t)p + (u+v)x=0, \ \ (s+t)y +
(u+v)p =0.$

\end {Example}

 By the way, we have that $\alpha , \beta $-symmetry does
 not depend on the particular choice of
basis of $L$ when $\alpha \beta \not=0$.

\begin {Proposition}
\label {11.4.6} Let $L$ be a Lie algebra with a basis $\{ e_1, e_2,
e_3 \}$ such that $[e_1,e_2]=0,  [e_1,e_3]= e_1+ \beta e_2,
[e_2,e_3] = \delta e_2$, where $\beta, \delta \in k$. Let $ p, q, s,
t, u, v, x, y, z \in k$.

(I) If $\beta =0,  \delta \not=0,$ then $r$   is a solution of CYBE
in $L$ iff

  $r=  p(e_1 \otimes e_2)+ q (e_2 \otimes e_1)
+ s (e_1 \otimes e_3)+ s (e_3 \otimes e_1) +u (e_2 \otimes e_3) +u
(e_3 \otimes e_2) +x (e_1 \otimes e_1) +y (e_2 \otimes e_2) +   z
(e_3 \otimes e_3) $

where $ (\delta +1 )zp =(\delta +1 ) qz =(\delta +1) us,    (\delta
+ 1)u q= (\delta +1)up$  and $(\delta +1) ps =(\delta +1) qs$.

(II) If  $\beta \not= 0 , \delta=1$, then $r$   is a solution of
CYBE in $L$ iff

  $r=  p(e_1 \otimes e_2)+ q (e_2 \otimes e_1)
+ s (e_1 \otimes e_3)+ s (e_3 \otimes e_1) +u (e_2 \otimes e_3) +u
(e_3 \otimes e_2) +x (e_1 \otimes e_1) +y (e_2 \otimes e_2)
                      +z (e_3 \otimes e_3) $

with $sp =sq, up =qu , zq =zp $ and $s^2=xz.$

(III)   If  $ \beta = \delta = 0,$ then $r$ is a  solution  of CYBE
in $L$ iff

$r=  p(e_1 \otimes e_2)+ q (e_2 \otimes e_1) + s (e_1 \otimes e_3)+
s (e_3 \otimes e_1) +u (e_2 \otimes e_3) +v (e_3 \otimes e_2) +x
(e_1 \otimes e_1) +y (e_2 \otimes e_2)
                      +z (e_3 \otimes e_3) $

with $vs =pz, us = qz, qv =pu$ and $(u+v)x = (p+q)s$.

\end {Proposition}

{ \bf Proof.} Let   $r = \sum_{i,j= 1 }^{3} k_{ij}(e_{i} \otimes
e_{j}) \in L \otimes L$ and  $k_{ij} \in k,$ with  $i, j =1, 2, 3.$
 By computation, for all $i, j, n = 1, 2, 3,$   we have that the coefficient of
  $e_j \otimes e_i \otimes e_i$ in
$[r^{12}, r^{13}] $ is zero and
 $e_i \otimes e_i \otimes e_j$
in  $[r^{13}, r^{23}]$ is zero.

We now see the coefficient
 of  $e_i \otimes e_j \otimes e_n$

   in  $[r^{12}, r^{13}] + [r^{12}, r^{23}] +[r^{13}, r^{23}]$.

(1) $e_1 \otimes e_1 \otimes e_1 (
 -k_{13}k_{11} +  k_{11} k_{31})$;

(2) $e_2 \otimes e_2 \otimes e_2 (
 -\beta k_{23}k_{12} + \beta k_{21} k_{32}
  - \delta k_{23}k_{22} + \delta k_{22} k_{32})$;

(3) $e_3 \otimes e_3 \otimes e_3 ( 0);$

(4) $e_1 \otimes e_2 \otimes e_3 (
 -k_{32}k_{13} +  k_{12} k_{33} -  \beta k_{13}  k_{13}+  \beta k_{11} k_{33}
- \delta k_{13}  k_{23} + \delta k_{12} k_{33} );$

(5) $e_2 \otimes e_3 \otimes e_1 (
 -\beta k_{33}k_{11} + \beta k_{13} k_{31}
 - \delta k_{33}  k_{21}+ \delta k_{23} k_{31}
 -  k_{23}  k_{31} + k_{21} k_{33} );$

(6) $e_3 \otimes e_1 \otimes e_2 (
 -k_{33}k_{12} +  k_{31} k_{32} -  \beta k_{33}  k_{11}+  \beta k_{31} k_{13}
- \delta k_{33}  k_{12} + \delta k_{32} k_{13} );$

(7) $e_1 \otimes e_3 \otimes e_2 (
 -k_{33}k_{12} +  k_{13} k_{32} -  \beta k_{13}  k_{32}+  \beta k_{11} k_{33}
- \delta k_{13}  k_{32} + \delta k_{12} k_{33} );$

(8) $e_3 \otimes e_2 \otimes e_1 (
 -\beta k_{33}k_{11} + \beta k_{31} k_{31} -  \delta k_{33}  k_{21}
 +  \delta k_{32} k_{31}
- k_{33}  k_{21} +  k_{31} k_{23} );$

(9) $e_2 \otimes e_1 \otimes e_3 (
 -\beta k_{31}k_{13} + \beta  k_{11} k_{33}
 - \delta k_{31}  k_{23}+  \delta k_{21} k_{33}
-  k_{23}  k_{13} +  k_{21} k_{33} );$

(10) $e_1 \otimes e_1 \otimes e_2 (
 -k_{31}k_{12} +  k_{11} k_{32} -   k_{13}  k_{12}+   k_{11} k_{32});$

(11) $e_2 \otimes e_1 \otimes e_1 (
 -k_{23}k_{11} +  k_{21} k_{31} -   k_{23}  k_{11}+  k_{21} k_{13});$

(12) $e_1 \otimes e_1 \otimes e_3 (
 -k_{31}k_{13} +  k_{11} k_{33} -   k_{13}  k_{13}+   k_{11} k_{33}
 );$

(13) $e_3 \otimes e_1 \otimes e_1 (
 -k_{33}k_{11} +  k_{31} k_{31} -   k_{33}  k_{11}+  k_{31} k_{13}
);$

(14) $e_2 \otimes e_2 \otimes e_1 (
 -\beta k_{32}k_{11} +  \beta k_{12} k_{31} -  \delta k_{32}  k_{21}
 +  \delta k_{22} k_{31}
- \beta k_{23}  k_{11} + \beta k_{21} k_{31} - \delta k_{23}  k_{21}
+ \delta k_{22} k_{31}
 );$

(15) $e_1 \otimes e_2 \otimes e_2 (
 -\beta k_{13}k_{12} + \beta  k_{11} k_{32} -  \delta k_{13}  k_{22}
 +  \delta k_{12} k_{32}
- \beta k_{13}  k_{21} + \beta k_{11} k_{23} - \delta k_{13}  k_{22}
+ \delta k_{12} k_{23}
 );$

(16) $e_2 \otimes e_2 \otimes e_3 (
 -\beta k_{32}k_{13} +  \beta k_{12} k_{33} -  \delta k_{32}  k_{23}
 + \delta k_{22} k_{33}
- \beta k_{23}  k_{13} + \beta k_{21} k_{33} - \delta k_{23}  k_{23}
+ \delta k_{22} k_{33} );$

(17) $e_3 \otimes e_2 \otimes e_2 (
 -\beta k_{33}k_{12} + \beta  k_{31} k_{32}- \delta k_{33}k_{22}
 +\delta k_{32}k_{32}
  -  \beta k_{33} k_{21} + \beta k_{31}  k_{23}
  - \delta k_{33} k_{22}
  + \delta k_{32} k_{23}
  );$

(18) $e_3 \otimes e_3 \otimes e_1 (0);$

(19) $e_1 \otimes e_3 \otimes e_3 (0);$

(20) $e_3 \otimes e_3 \otimes e_2 (0);$

(21) $e_2 \otimes e_3 \otimes e_3 (0);$

(22) $e_1 \otimes e_3 \otimes e_1 (
 -k_{33}k_{11} +  k_{13} k_{31} -   k_{13}  k_{31}+   k_{11} k_{33})
 = e_1 \otimes e_3 \otimes e_1 (0)$;

(23) $e_1 \otimes e_2 \otimes e_1 (
 -k_{32}k_{11} +  k_{12} k_{31} -  \beta k_{13}  k_{11}+  \beta k_{11} k_{31}
- \delta k_{13}  k_{21} + \delta k_{12} k_{31} -  k_{13}  k_{21} +
k_{11} k_{23});$

(24) $e_2 \otimes e_1 \otimes e_2 (
 -\beta k_{31}k_{12} + \beta k_{11} k_{32}
 - \delta k_{31}  k_{22}+  \delta k_{21} k_{32}
- k_{23}  k_{12} +  k_{21} k_{32}
 -\beta k_{23}k_{11} + \beta k_{21} k_{13}
  - \delta k_{23}  k_{12}+  \delta k_{22} k_{13}  );$

(25) $e_2 \otimes e_3 \otimes e_2 (
 -\beta k_{33}k_{12} + \beta k_{13} k_{32}
  -  \delta k_{33}  k_{22}+  \delta k_{23} k_{32}
- \beta k_{23}  k_{31} + \beta k_{21} k_{33}
 + \delta k_{33}  k_{22}-  \delta k_{22} k_{33}
              );$

(26) $e_3 \otimes e_2 \otimes e_3 (
 -\beta k_{33}k_{13} + \beta k_{31} k_{33}
  - \delta k_{33}  k_{23}+  \delta k_{32} k_{33}
);$

(27) $e_3 \otimes e_1 \otimes e_3 (
 - k_{33}k_{13} +   k_{31} k_{33}
 );$

Let $k_{11}=x, k_{22}=y,  k_{33}=z, k_{12}=p,  k_{21} =q, k_{13}=s,
 k_{31} =t, k_{23}=u,  k_{32}=v.$

   It follows from (1)--(27) that

(28) $-sx + xt=0$;

(29) $-\beta up +\beta qv -\delta uy + \delta yv=0$;

(30) $- vs + pz -\beta s^2 + \beta xz -\delta su + \delta zp =0$;

(31) $-\beta xz + \beta st - \delta zq + \delta ut - ut   + qz =0$;

(32) $- zp + tv -\beta zx + \beta ts -\delta zp + \delta vs =0$;

(33) $- zp + sv -\beta st + \beta xz -\delta sv + \delta zp =0$;

(34) $- \beta zx + \beta tt - \delta zq + \delta vt - zq + tu =0$;

(35) $- \beta st + \beta xz  -\delta tu + \delta qz -us + qz=0 $;

(36) $- tp + xv - sp + vx =0$;

(37) $- ux + qt-ux + qs =0$;

(38) $- st + xz - s^2 +  xz =0$;

(39) $-zx +tt -zx +st =0$;

(40) $-\beta vx +\beta pt -\delta vq + \delta yt -\beta ux + \beta
qt -\delta uq + \delta yt =0$;

(41) $-\beta sp + \beta xv  -\delta  sy + \delta pv -\beta sq +
\beta xu - \delta sy + \delta pu =0$;

(42) $-\beta vs + \beta pz -\delta vu + \delta yz -\beta su + \beta
zq -\delta u^2 + \delta yz  =0$;

(43) $-\beta  zp + \beta tv - \delta zy + \delta vv  -\beta zq +
\beta tu -\delta zy +\delta vu =0$;

(44) $- vx + pt -\beta sx + \beta xt -\delta sq + \delta pt -sq + xu
=0$;

(45) $- \beta pt + \beta vx -\delta ty + \delta qv  -up + qv -\beta
ux + \beta qs - \delta up + \delta ys =0$;

(46) $-\beta zp + \beta sv -\beta ut + \beta qz =0$;

(47) $- \beta zs + \beta tz -\delta zu + \delta vz =0$;

(48) $- zs + tz =0$;

It is clear that   $r$   is a solution of CYBE iff
 (28)--(48) hold.

By computation we have that

(49)  $t =s$  (by  (39) + (38) );

 (50)  $\delta u = \delta v$  (by (42) + (43) ).

(I) Let $\beta =0 $ and $\delta \not=0.$ If $r$ is the case in part
(II), then $r$ is a solution of CYBE by straightforward
verification. Conversely, if $r$ is a solution of CYBE, then $u=v$
by (50) and we have that
\begin {eqnarray*}
(1 +\delta )zp &=& (1+\delta ) us    \hbox { \ \ \ \ by }  (30)   \\
(1 +\delta )zq &=& (1+\delta ) us    \hbox { \ \ \ \ by }  (31)   \\
(1 +\delta )sp &=& (1+\delta ) qs    \hbox { \ \ \ \ by }  (44)   \\
(1 +\delta )qu &=& (1+\delta ) up    \hbox { \ \ \ \ by }  (45)
\end {eqnarray*}
Thus $r$ is the case in part (II).

 (II) Let $\beta \not=0$ and $ \delta =1$. If $r$ is the case in part (III),
  then $r$ is a solution of CYBE by straightforward verification. Conversely,
 if $r$ is a solution of CYBE then $u=v$ by (50) and we have that
\begin {eqnarray*}
up &=&qu    \hbox { \ \ \ \ by }  (29)   \\
s^2 &=&xz    \hbox { \ \ \ \ by }  (30)   \\
sp &=&  qs    \hbox { \ \ \ \ by }  (40)   \\
zp &=& zq    \hbox { \ \ \ \ by }  (46)
\end {eqnarray*}
Thus $r$ is the case in part (III).

(III) Let $\beta =\delta =0.$ It is clear that the system of
equations  (28)--(48) is equivalent to the following:
$$  \left    \{  \begin{array} {l}
s=t \\
us =qz   \\
vs =zp  \\
up =qv  \\
(p+q)s = (u+v )x
\end{array} \right.$$
                 Thus we complete the proof.  \begin{picture}(8,8)\put(0,0){\line(0,1){8}}\put(8,8){\line(0,-1){8}}\put(0,0){\line(1,0){8}}\put(8,8){\line(-1,0){8}}\end{picture}

\section {Coboundary Lie bialgebras in the  case of  $char k =2$ }\label {s34}

In this section, using the general solution, which are obtained in
the preceding section,  of CYBE in Lie algebra $L$ with $dim \ \ L
\le 3, $ we give the
       the sufficient and  necessary conditions for  $(L, \hbox {[ \ ]},
 \Delta _r, r)$ to be a coboundary
 (or triangular ) Lie bialgebra over field $k$ of characteristic 2.
Throughout this section, the characteristic of field $k$  is 2.

We now observe the  connection between solutions of CYBE and
triangular Lie bialgebra structures.

\begin {Lemma}  \label {11.5.1}
If $r \in Im (1-\tau )$, then $$x\cdot C(r) = (1+\xi +\xi ^2)
 (1 \otimes \Delta) \Delta (x)$$ for any $x \in L$.
 \end {Lemma}

{\bf Proof:} Let $r= \sum _i (a_i \otimes b_i - b_i \otimes a_i ).$

We see that
\begin {eqnarray*}
C(r) &=& \sum _{i,k}  ( [a_k , a_i] \otimes b_k \otimes b_i -[a_k ,
b_i] \otimes b_k \otimes a_i
-[b_k , a_i] \otimes a_k \otimes b_i \\
 &+&  [b_k, b_i] \otimes a_k \otimes a_i
 +  a _k \otimes [b_k, a_i] \otimes b_i
  -a _k \otimes [b_k, b_i] \otimes a_i  \\
  &-&  b _k \otimes [a_k, a_i] \otimes b_i
   -  b _k \otimes [a_k, a_i] \otimes b_i
   +  b _k \otimes [a_k, b_i] \otimes a_i \\
&+&   a _k \otimes a_i \otimes [b_k, b_i] -  a _k \otimes b_i
\otimes [b_k, a_i]
-   b _k \otimes a_i \otimes [a_k, b_i] \\
 &+& b _k \otimes b_i \otimes [a_k, a_i])
 \end {eqnarray*}

and

\begin {eqnarray*}
&{}& ( 1\otimes \Delta )\Delta (x) \\
 &=&
\sum _{i,k} ( [x,a_k] \otimes [b_k,a_i] \otimes b_i +
 [x,a_k] \otimes  a_i \otimes [b_k, b_i]  -
                 [x,a_k] \otimes [b_k,b_i] \otimes a_i \\
&-& [x,a_k] \otimes  b_i \otimes [b_k, a_i] -             [x,b_k]
\otimes [a_k,a_i] \otimes b_i
- [x,b_k] \otimes  a_i \otimes [a_k, b_i] \\
&+&                 [x,b_k] \otimes [a_k,b_i] \otimes a_i + [x,b_k]
\otimes  b_i \otimes [a_k, a_i]
+a_k \otimes [[x,b_k],a_i] \otimes b_i \\
&+& a_k \otimes  a_i \otimes [[x,b_k], b_i]
 -                    a_k \otimes [[x,b_k],b_i] \otimes a_i
 - a_k \otimes  b_i \otimes [[x,b_k], a_i]   \\
                 &-&              b_k \otimes [[x,a_k],a_i] \otimes b_i
                 -b_k \otimes  a_i \otimes [[x,a_k], b_i]  +
                 b_k \otimes [[x,a_k],b_i] \otimes a_i \\
                 &+&
 b_k \otimes  b_i \otimes [[x,a_k], a_i] ).
\end {eqnarray*}

It is easy to check that the sum of the terms whose third factor in
$(1+\xi +\xi ^2)(1+\Delta )\Delta (x)$ includes element $x$  is
equal to

\begin {eqnarray*}
 &{}& \sum _{i,k} \{
 [b_k,a_i] \otimes b_i  \otimes [x,a_k] +
 a_i \otimes [b_k, b_i]  \otimes  [x,a_k] -
                [b_k,b_i] \otimes a_i \otimes [x,a_k]  \\
                &-&  b_i \otimes [b_k, a_i]  \otimes [x,a_k]
-             [a_k,a_i] \otimes b_i \otimes  [x,b_k]
-  a_i \otimes [a_k, b_i] \otimes [x,b_k] \\
&+&              [a_k,b_i] \otimes a_i  \otimes  [x,b_k]
+ b_i \otimes [a_k, a_i]  \otimes   [x,b_k]   \\
&+& (b_i \otimes a_k  \otimes [[x,b_k],a_i]
- b_k \otimes  a_i \otimes [[x,a_k], b_i])  \\
&+& ( a_k \otimes  a_i \otimes [[x,b_k], b_i]
 -    a_i \otimes a_k\otimes [[x,b_k],b_i] ) \\
 &+&  ( -  a_k \otimes  b_i \otimes [[x,b_k], a_i]
 +  a_i   \otimes  b_k \otimes  [[x,a_k],b_i] ) \\
&+& (- b_i \otimes  b_k \otimes  [[x,a_k],a_i]
+ b_k \otimes  b_i \otimes [[x,a_k], a_i] ) \}  \\
              &=& (1 \otimes 1 \otimes L_x ) C(r)    \hbox { \ \ \ \
              by Jacobi identity, }
\end {eqnarray*}
 where    $L_x$  denotes the adjoint
              action  $L_x(y) = [x,y]$  of $L$ on $L$.

Consequently,   $$(1 +\xi + \xi ^2) (1 \otimes \Delta )\Delta (x)  =
 (L_x \otimes L_x \otimes L_x ) C(r) = x \cdot C(r)$$
for any $x \in L.$
\begin{picture}(8,8)\put(0,0){\line(0,1){8}}\put(8,8){\line(0,-1){8}}\put(0,0){\line(1,0){8}}\put(8,8){\line(-1,0){8}}\end{picture}

By Lemma \ref {11.5.1} and \cite [ Proposition 2.11 ] {Mi94}, we
have

\begin {Theorem} \label {11.5.2}
$(L, [ \hbox { \ }], \Delta _r ,r)$ is a triangular Lie bialgebra
iff
 $r$ is a solution of CYBE in $L$ and $r \in Im (1-\tau ).$
 \end {Theorem}

 Consequently, we can easily get a  triangular Lie bialgebra
 structure  by means of a solution of CYBE.

\begin {Theorem}
\label {11.5.3} Let $L$ be a  Lie algebra with
       a basis $\{ e_1, e_2, e_3 \}$
such that $[e_1,e_2]= e_3, [e_2,e_3]=\alpha e_1, [e_3,e_1]=\beta
e_2$, where $\alpha, \beta \in k$ , $\alpha \beta \not=0$ or $\alpha
= \beta =0$. Then

(I)   $(L, [\hbox { \ }], \Delta _r, r )$ is  a coboundary Lie
bialgebra iff $r \in Im(1 -\tau )$ ;

(II)       $(L, [\hbox { \ }],\Delta _r, r )$ is  a triangular  Lie
bialgebra iff  $r \in Im (1-\tau )$ and $r$ is  $\alpha , \beta $-
symmetric to the basis of  $\{ e_1, e_2, e_3 \}$.
\end {Theorem}

 {\bf Proof.}
(I) Obviously, $r\in Im (1-\tau )$ when
 $(L, [\hbox { \ }], \Delta _r, r )$ is  a coboundary Lie bialgebra. Conversely,
 if $r\in Im (1-\tau )$, let
 $r=  p(e_1 \otimes e_2)- p (e_2 \otimes e_1)
+ s (e_1 \otimes e_3)- s (e_3 \otimes e_1) + u (e_2 \otimes e_3)- u
(e_3 \otimes e_2)$. It is sufficient to show that
 $$ (1+ \xi + \xi ^2) (1 \otimes \Delta ) \Delta (e_i) =0$$
for $i=1, 2, 3$  by \cite [Proposition 2.11]{Mi94}. First, by
computation, we have that
\begin {eqnarray*}
&{ \ }& (1 \otimes \Delta ) \Delta (e_1) \\
&=& \{ ( e_1 \otimes e_2\otimes e_3 ) (\beta ps -  \beta ps)
 +(e_3 \otimes e_1\otimes e_2)(\beta ps)
 + (e_2 \otimes e_3\otimes e_1) (- \beta sp) \} \\
&+&\{ ( e_1 \otimes e_3\otimes e_2 ) (-\beta ps + \beta ps)
 +(e_3 \otimes e_2\otimes e_1)( - \beta ps)
 + (e_2 \otimes e_1\otimes e_3) (\beta sp) \} \\
&+&\{ ( e_1 \otimes e_1\otimes e_2 ) (\alpha \beta su )
 +(e_1 \otimes e_2\otimes e_1)( - \beta \alpha su)
 + (e_2 \otimes e_1\otimes e_1) (0) \} \\
&+& \{ ( e_1 \otimes e_1\otimes e_3 ) ( - \alpha pu)
 +(e_1 \otimes e_3\otimes e_1)(\alpha  pu)
 + (e_3 \otimes e_1\otimes e_1) (0) \} \\
&+&\{   ( e_2 \otimes e_1\otimes e_2 ) (- \beta ^2 s^2)
 +(e_1 \otimes e_2\otimes e_2)(0)
 + (e_2 \otimes e_2\otimes e_1) ( \beta ^2 s^2) \} \\
&+&\{  ( e_3 \otimes e_3\otimes e_1 ) ( p^2 )
 +(e_3 \otimes e_1\otimes e_3)(- p^2)
 + (e_1 \otimes e_3\otimes e_3) (0) \}.
\end {eqnarray*}
Thus
 $$ (1+ \xi + \xi ^2) (1 \otimes \Delta ) \Delta (e_1) =0.$$
Similarly, we have that
 $$ (1+ \xi + \xi ^2) (1 \otimes \Delta ) \Delta (e_2) =0,
 \ \ \  (1+ \xi + \xi ^2) (1 \otimes \Delta ) \Delta (e_3) =0.$$
Thus  $(L, [\hbox { \ }],\Delta _r,r)$ is a coboundary Lie
bialgebra.

(II) It follows from Theorem \ref {11.5.2} and Proposition \ref
{11.4.4}.
\begin{picture}(8,8)\put(0,0){\line(0,1){8}}\put(8,8){\line(0,-1){8}}\put(0,0){\line(1,0){8}}\put(8,8){\line(-1,0){8}}\end{picture}

\begin
{Example} \label {11.5.4} Under Example \ref {11.4.5}, we have the
following:

(i)   $(sl(2), [\hbox { \ }],\Delta _r, r )$ is  a coboundary Lie
bialgebra iff $r\in Im (1-\tau )$ ;

(ii)       $(sl(2), [\hbox { \ }],\Delta _r, r )$ is  a triangular
Lie bialgebra iff $r$ is  $0 , 0 $- symmetric to the basis of  $\{
e_1, e_2, e_3 \}$ iff
 $r=
 s (e_1 \otimes e_3)+ s (e_3 \otimes e_1)
+ u (e_2 \otimes e_3)+ u (e_3 \otimes e_2)$.
\end {Example}

 \begin {Theorem} \label  {11.5.5}
Let  $L$ be a Lie algebra
 with a basis   $\{ e_1, e_2, e_3 \}$  such that
$$[e_1,e_2]= 0, [e_1, e_3]= e_1 + \beta e_2, [e_2, e_3]=\delta e_2,$$
where $\delta, \beta \in k$ and $\delta =1$ when $\beta \not=0.$ Let
$p, s, u \in k$ and
 $r=  p(e_1 \otimes e_2)- p (e_2 \otimes e_1)
+ s (e_1 \otimes e_3)- s (e_3 \otimes e_1) + u (e_2 \otimes e_3)- u
(e_3 \otimes e_2)$. Then

(I)  $(L, [\hbox  { \ }],\Delta _r, r )$ is a coboundary Lie
bialgebra iff $r\in Im (1-\tau )$ and  $$ (\delta +1) ((\delta +1)u
+ \beta s)s =0; $$

(II)  $(L, [\hbox  { \ }],\Delta _r, r )$ is a triangular Lie
bialgebra
 iff $r \in Im(1-\tau )$ and $$ \beta s +(1+ \delta )  us =0.$$

\end {Theorem}
{ \bf Proof }.
 (I) We get by computation
 \begin {eqnarray*}
&{ \ }& ( 1 \otimes \Delta ) \Delta (e_1) \\
 &=& \{ ( e_1 \otimes e_1\otimes e_2 ) (\beta \delta ss -  \delta us)
 +(e_1 \otimes e_2\otimes e_1)(- \beta \delta ss + \delta us)
 + (e_2 \otimes e_1\otimes e_1) (0) \} \\
&+&  \{ ( e_2 \otimes e_2\otimes e_1 ) (-\beta us + uu+ \beta ^2 s^2
- \beta su)                                       \\
 &+& (e_2 \otimes e_1\otimes e_2)( \beta us -uu + \beta us - \beta ^2 s^2)
  + (e_1 \otimes e_2\otimes e_2) (0) \} \\
&{ \ }& ( 1 \otimes \Delta ) \Delta (e_2) \\
 &=& \{ ( e_1 \otimes e_1\otimes e_2 ) ( \delta ^2 ss)
 +(e_1 \otimes e_2\otimes e_1)(- \delta ^2 s^2 )
 + (e_2 \otimes e_1\otimes e_1) (0) \} \\
 &+& \{ ( e_2 \otimes e_2\otimes e_1 ) ( \delta \beta ss - \delta su)
  +(e_2 \otimes e_1\otimes e_2)( \delta  su - \delta \beta ss )
 + (e_1 \otimes e_2\otimes e_2) (0) \} \\
&{ \ }& ( 1 \otimes \Delta ) \Delta (e_3)  \\
 &=& \{ ( e_1 \otimes e_2\otimes e_3 ) (\delta su +  \beta ss)
 +(e_2 \otimes e_3\otimes e_1)(-\delta us - \beta ss) \\
 &+& (e_3 \otimes e_1\otimes e_2) ( \delta \delta us + \beta \delta ss
 + \beta ss -su) \} \\
&+&\{ ( e_1 \otimes e_3\otimes e_2 ) (-\delta su-  \beta ss)
 +(e_3 \otimes e_2\otimes e_1)( - \beta \delta  ss - \delta \delta us
 -\beta ss + su ) \\
 &+&  (e_2 \otimes e_1\otimes e_3) (\delta us + \beta ss) \} \\
&+& \{ ( e_1 \otimes e_1\otimes e_2 ) (- \delta  sp - \delta \delta
ps
+ sp + \delta sp ) \\
 &+& (e_1 \otimes e_2\otimes e_1)( \delta ps + \delta \delta ps -\delta sp
 -sp) +
 (e_2 \otimes e_1\otimes e_1) (0) \} \\
&+& \{
 ( e_1 \otimes e_1\otimes e_3 ) ( ss)
 +(e_1 \otimes e_3\otimes e_1)(-ss)
 + (e_3 \otimes e_1\otimes e_1) (0) \} \\
&+& \{   ( e_2 \otimes e_2\otimes e_1 ) (- \beta  ps + pu - \delta
\beta sp
-\delta \delta up - \delta \beta ps + \delta pu - \beta sp - \delta up)  \\
&+& (e_2 \otimes e_1\otimes e_2)( \beta  ps - pu + \delta \beta sp
+\delta \delta up + \delta \beta ps - \delta pu + \beta sp + \delta up) ) \\
&-& (e_1 \otimes e_2\otimes e_2) (0) \} \\
&+& \{ ( e_2 \otimes e_2\otimes e_3 ) ( \delta \delta uu + \delta
\beta us
+ \beta \delta us + \beta \beta ss )   \\
 &+& (e_2 \otimes e_3\otimes e_2)(   -\delta \delta uu - \delta \beta us
- \beta \delta us - \beta \beta ss )
                    +        (e_3 \otimes e_2\otimes e_2) (0) \}.
\end {eqnarray*}
Consequently,   $$ (1+ \xi + \xi ^2) (1 \otimes \Delta ) \Delta
(e_1) =0,$$
 $$ (1+ \xi + \xi ^2) (1 \otimes \Delta ) \Delta (e_2) =0$$
 and
 $$ (1+ \xi + \xi ^2) (1 \otimes \Delta ) \Delta (e_3) =0$$
iff
  $$ \delta ^2 us + \delta \beta s^2 + \beta s^2 - us =0.$$
This implies that
    $(L, [\hbox  { \ }],\Delta _r,r)$ is a coboundary Lie bialgebra
    iff    $r\in Im (1-\tau )$ and
 $$ (\delta +1)(( \delta +1) u + \beta s )s =0   { \ \ \ \ }.$$

(II) It follows from Theorem \ref {11.5.2}  and Proposition  \ref
{11.4.6}.
 \begin{picture}(8,8)\put(0,0){\line(0,1){8}}\put(8,8){\line(0,-1){8}}\put(0,0){\line(1,0){8}}\put(8,8){\line(-1,0){8}}\end{picture}

\begin {Theorem}  \label   {11.5.6}
If $L$ is a Lie algebra with $dim L =2$ and $r \in L \otimes L$,
then
  $(L, [\hbox  { \ }],\Delta _r,r)$ is a triangular  Lie bialgebra   iff
  $(L, [\hbox  { \ }],\Delta _r,r)$ is a coboundary   Lie bialgebra
  iff   $r\in Im (1-\tau )$.
   \end {Theorem}
 { Proof.}
It is an immediate consequence of the main result of \cite {Mi94}.
  \begin{picture}(8,8)\put(0,0){\line(0,1){8}}\put(8,8){\line(0,-1){8}}\put(0,0){\line(1,0){8}}\put(8,8){\line(-1,0){8}}\end{picture}

\section {The solutions of CYBE with char  $ k\neq 2  $ }

In this section, we find the general solution of CYBE for Lie
algebra  $ L $  with dim  $ L $ =3 and dim  $ L'= 2 $, where char $
k\neq 2 $.
 If $ k $ is not
algebraically   closed, let  $ P $  be algebraic closure of $ k $.
we can construct a Lie algebra  $ L_{P}=P \otimes L $ over $ P $, as
in \cite [Section 8]{Ja62}.

By \cite [P11--14]{Ja62}, we have that

\begin {Lemma} \label {11.6.1}
Let  $ L $  be a vector space over  $ k $. Then  $ L $  is a Lie
algebra over field $k$ with dim  $ L $ =3 and dim  $ L^{'}=2  $  iff
there is a basis  $ e_{1} $,  $ e_{2} $,  $ e_{3} $  in  $ L $ such
that  $ [e_{1}, e_{2}] $ =0,  $ [e_{1}, e_{3}]=\alpha e_{1}+\beta
e_{2} $,  $ [e_{2}, e_{3}]=\gamma e_{1}+\delta  e_{2} $, where  $
\alpha, \beta, \gamma, \delta  \in k $, and  $ \alpha \delta -\beta
\gamma \neq 0 $.
\end {Lemma}

In following three sections, we only study the Lie algebra $ L $ in
Lemma \ref {11.6.1}. Set  $ A $ = $
      \left(\begin{array}{cc}
      \alpha  & \gamma  \\
      \beta  & \delta
      \end{array}\right ).
      $
Thus  $ A $  is similar to  $
                          \left(\begin{array}{cc}
                          \lambda_{1}& 0\\
                          0 & \lambda_{2}
                          \end{array}\right )
                         $
or $
     \left(\begin{array}{cc}
     \lambda_{1}& 0\\
     1& \lambda_{1}
      \end{array}\right )
      $
in the algebraic closure  $ P $  of  $ k $. Therefore, there is an
invertible matrix  $ D $  over  $ P $  such that  $ AD=D
                                                \left(\begin{array}{cc}
                                                 \lambda_{1}& 0\\
                                                  0 & \lambda_{2}
                                                 \end{array}\right )
                                                 $, or
 $ AD=D
                                                \left(\begin{array}{cc}
                                                 \lambda_{1}& 0\\
                                                  1& \lambda_{2}
                                                 \end{array}\right )
                                                 $.
                                                Let  $ Q=
                                                                 \left(\begin{array}{cc}
                                                                 D& 0\\
                                                                 0 &
                                                                \frac {1} { \lambda_{1}}
                                                                 \end{array}\right )
                                                               $
and  $ (e_{1}^{'}, e_{2}^{'}, e_{3}^{'})=(e_{1}, e_{2}, e_{3})Q $.
By computation, we have that  $ [e_{1}^{'}, e_{2}^{'}]=0 $,
 $ [e_{1}^{'}, e_{3}^{'}]=e_{1}^{'}+\beta ^{'}e_{2}^{'} $,  $ [e_{2}^{'},
e_{3}^{'}]=\delta ^{'}e_{2}^{'} $,  where  $ \beta ^{'}=0 $  and
 $ \delta ^{'}=\frac { \lambda_{1}} {\lambda_{2}} $  when  $ A $  is similar
to  $
                                       \left(\begin{array}{cc}
                                       \lambda_{1}& 0\\
                                        0 & \lambda_{2}
                                        \end{array}\right )
                                       $;
 $ \beta ^{'}=\frac {1} {\lambda_{1}} $  and  $ \delta ^{'}=1 $  when  $ A $  is
similar to  $
              \left(\begin{array}{cc}
              \lambda_{1}& 0\\
               1& \lambda_{1}
               \end{array}\right )
               $.
Let  $ Q=(q_{ij}{})_{3\times 3} $  and  $ Q^{-1}=(\bar
q_{ij})_{3\times 3} $. If  $ r =\sum_{i, j=1} ^{3}
k_{ij}(e_{i}\otimes e_{j})=\sum_{i, j=1} ^{3}
k_{ij}^{'}(e_{i}^{'}\otimes  e_{j}^{'}) $, where  $ k_{ij}\in k $, $
k_{ij}^{'}\in  P $  for  $ i, j $ =1, 2, 3, then  $$
k_{ij}^{'}=\sum_{m, n}^{3}k_{mn}\overline{q}_{im}\overline{q}_{jn}
 \hbox { and }
  k_{ij}=\sum_{m, n}^{3}k_{mn}^{'}q_{im}q_{jn} $$  for  $ i, j=1, 2, 3 $.
Obviously,  $ k_{33}=k_{33}^{'} $.

\begin {Lemma} \label{11.6.2}
(i)  $ k_{i3}=k_{3i} $, for  $ i $ =1, 2, 3 iff
 $ k_{i3}^{'}=k_{3i}^{'} $  for $ i $ = 1, 2, 3;

(ii)  $ k_{i3}=-k_{3i} $  for  $ i $ =1, 2, 3 iff
 $ k_{i3}^{'}=- k_{3i}^{'} $  for  $ i $ =1, 2, 3;

(iii)  $ k_{ij}=k_{ji} $  for  $ i, j $ =1, 2, 3 iff
 $ k_{ij}^{'}=k_{ji}^{'} $  for  $ i, j=1, 2, 3 $;

(iv)  $ k_{ij}=-k_{ji} $  for  $ i, j $ =1, 2, 3 iff
 $ k_{ij}^{'}=-k_{ji}^{'} $  for  $ i, j $ =1, 2, 3.

\end {Lemma}
{\bf Proof } (i) If  $ k_{i3}^{'}=k_{3i}^{'} $  for  $ i $ =1, 2, 3,
we see that
\begin {eqnarray*}
  k_{i3}  &=& \sum_{m, n}^{3} k_{mn}^{'}q_{im}q_{3n} \\
        &=& \sum_{m}^{3}k_{m3}^{'}q_{im}q_{33} \ \ ( \hbox {    since }  q_{31}=q_{32}=0
        )\\
       & =& \sum_{m}^{3}k_{3m}^{'}q_{im}q_{33}  \ \ ( \hbox {    by
       assumption})\\
       & = &  \sum_{m}^{3}k_{3m}^{'}q_{33}q_{im} \\
        &=& \sum_{m, n}^{3}k_{mn}^{'}q_{3n}q_{im} \\
        &=&  k_{3i}.
        \end {eqnarray*}
Therefore,  $ k_{i3}=k_{3i} $  for  $ i $ =1, 2, 3. The others can
be proved similarly. $\Box$

\begin{Theorem} \label {11.7.1}
Let  $ L $  be a Lie algebra with a basis  $ {e_{1}, e_{2}, e_{3}} $
such that $ [e_{1}, e_{2}] $ =0,  $ [e_{1}, e_{3}]=\alpha
e_{1}+\beta e_{2} $,  $ [e_{2}, e_{3}]=\gamma  e_{1}+\delta e_{2} $,
where  $ \alpha, \beta, \gamma, \delta  \in  k $, and  $ \alpha
\delta -\beta \gamma \neq 0 $. Let $ p, q, s, t, u, v, x, y, z \in k
$. Then  $ r $  is a solution of CYBE iff  $ r $  is strongly
symmetric, or  $ r=p(e_{1}\otimes e_{2})+q(e_{2}\otimes
e_{1})+s(e_{1}\otimes e_{3})-s(e_{3}\otimes e_{1})+u(e_{2}\otimes
e_{3})-u(e_{3}\otimes e_{2})+x(e_{1}\otimes e_{1})+y(e_{2}\otimes
e_{2}) $  with  $ s(2\alpha x+\gamma (p+q))=u(2\delta y+\beta
(q+p))=u(2\alpha x+\gamma (q+p))=s(2\delta y+\beta (q+p))=(\alpha
-\delta )us+\gamma u^{2}-\beta s^{2} =s(2\gamma y+2\beta x+(\alpha
+\delta )(q+p))=u(2\gamma y+2\beta x+(\alpha +\delta )(p+q))=0 $.
\end{Theorem}

{\bf  Proof } Let  $ r =\sum_{i, j=1}^{3}k_{ij}(e_{i}\otimes
e_{j})\in  L \otimes L $, and  $ k_{ij} \in  k$, with  $ i, j $ =1,
2, 3. By computation, for all  $ i, j, n $ =1, 2, 3, we have that
the cofficient of  $ e_{j}\otimes  e_{i}\otimes  e_{i} $  in
 $ [r ^{12}, r ^{23}] $  is zero and the cofficient of  $ e_{i}\otimes  e_{i}\otimes  e_{j} $  in
 $ [r ^{13}, r ^{23}] $  is zero.

We can obtain the following equations by seeing the cofficient of
 $ e_{i} \otimes  e_{j} \otimes  e_{n} $  in  $ [r ^{12}, r
^{13}]+[r ^{12}, r ^{23}]+[r ^{13}, r ^{23}] $, as in \cite
[Proposition 2.6]{Zh98b}. To simplify notation, let
 $ k_{11}=x, k_{22}=y, k_{33}=z, k_{12}=p, k_{21}=q, k_{13}=s,
k_{31}=t, k_{23}=u, k_{32}=v $.

  (1)  $ -\alpha sx+\alpha  xt-\gamma sq+\gamma pt=0 $;

  (2)  $ -\beta up+\beta qv-\delta uy+\delta yv=0 $;

  (3) $ -\alpha vs+\alpha pz-\gamma uv+\gamma yz-\beta s^{2}+\beta xz-\delta su+\delta zp=0 $;

  (4) $ -\beta xz+\beta st-\delta zq+\delta ut-\alpha ut+\alpha qz-\gamma uv+\gamma yz=0 $;

  (5) $ -\alpha zp+\alpha tv-\gamma zy+\gamma uv-\beta zx+\beta ts-\delta zp+\delta vs=0 $;

  (6) $ -\alpha zp+\alpha sv-\gamma zy+\gamma uv-\beta st+\beta xz-\delta sv+\delta zp=0
  $;

  (7) $ -\beta zx+\beta tt-\delta zq+\delta vt-\alpha zq+\alpha tu-\gamma zy+\gamma uv=0 $;

  (8) $ -\beta st+\beta xz-\delta tu+\delta qz-\alpha us+\alpha qz-\gamma uu+\gamma yz=0 $;

  (9) $ -\alpha tp+\alpha xv-\gamma ty+\gamma qv-\alpha sp+\alpha vx-\gamma sy+\gamma pv=0 $;

  (10) $ -\alpha ux+\alpha qt-\gamma uq+\gamma yt-\alpha ux+\alpha qs-\gamma up+\gamma ys=0 $;

  (11) $ -\alpha st+\alpha xz-\gamma tu+\gamma qz-\alpha s^{2}+\alpha xz-\gamma su+ypz=0 $;

  (12) $ -\alpha zx+\alpha tt-\gamma zq+\gamma vt-\alpha zx+\alpha st-\gamma zp+\gamma vs=0 $;

  (13) $ -\beta vx+\beta pt-\delta vq+\delta yt-\beta ux+\beta qt-\delta uq+\delta yt=0 $;

  (14) $ -\beta sp-\beta xv-\delta sy+\delta pv-\beta sq+\beta xu-\delta sy+\delta pu=0 $;

  (15) $ -\beta vs+\beta pz-\delta vu+\delta yz-\beta su+\beta zq-\delta u^{2}+\delta yz=0 $;

  (16) $ -\beta zp+\beta tv-\delta zy+\delta vv-\beta zq+\beta tu-\delta zy+\delta vu=0 $;

  (17) $ -\gamma zq+\gamma ut-\gamma sv+\gamma pz=0 $;

  (18) $ -\alpha vx+\alpha pt-\gamma vq+\gamma yt-\beta sx+\beta xt-\delta sq+\delta pt-\alpha sq+\alpha xu-\gamma sy+\gamma pu=0 $;

  (19) $ -\beta pt+\beta vx-\delta ty+\delta qv-\alpha up+\alpha qv-\gamma uy+\gamma yv-\beta ux+\beta qs-\delta up+\delta ys=0 $;

  (20) $ -\beta zp+\beta sv-\delta zy+\delta uv-\beta ut+\beta qz-\delta uv+\delta yz=0 $;

  (21) $ -\beta zs+\beta tz-\delta zu+\delta vz=0 $;

  (22) $ -\alpha zs+\alpha tz-\gamma zu+\gamma vz=0 $;

It is clear that  $ r $  is a solution of CYBE iff (1)-(22) hold.

By simple computation, we have the sufficiency. Now we show the
necessity, If  $ k_{33}\neq 0 $, then  $ k_{33}^{'}\neq 0' $  and so
$ r $  is a strongly symmetric element in  $ L_{P}\otimes L_{P} $ by
\cite [The proof of Proposition 1.6 ] {Zh98b}. Thus  $ r $  is a
strongly symmetric element in  $ L\otimes L $. If  $ k_{33} $ =0,
then  $ k_{33}^{'} $ =0. By \cite [Proposition 1.6] {Zh98b}, we have
that  $ k_{i3}^{'}=-k_{3i}^{'} $  for  $ i $ =1, 2, 3, which implies
that  $ k_{i3}=-k_{3i} $  for  $ i $ =1, 2, 3 by Lemma \ref
{11.6.2}.

It immediately follows from (1)-(22) that

  (23) $ s(-2\alpha x-\gamma (q+p))=0 $;

  (24) $ u(2\delta y+\beta (q+p))=0 $;

  (25) $ \gamma u^{2}-\beta s^{2}+(\alpha -\delta )us=0 $;

  (26) $ u(2\alpha x+\gamma (q+p))=0 $;

  (27) $ s(2\delta y+\beta (q+p))=0 $;

  (28) $ 2\alpha ux-2\gamma ys-2\beta xs+(-s(\alpha +\delta )+\gamma u)(q+p)=0 $;

  (29) $ -2\beta ux+2\delta ys-2\gamma uy+(-u(\alpha +\delta )+\beta s)(q+p)=0
  $;\\
By (26), (27), (28) and (29), we have that

  (30)  $ s(2\gamma y+2\beta x+(\alpha +\delta )(q+p))=0 $;

  (31)  $ u(2\gamma y+2\beta x+(\alpha +\delta )(p+q))=0 $. $\Box$

\begin{Example} \label {11.7.2}
Let  $ L $  be a Lie algebra over real field  $ R $  with dim  $ L $
=3 and dim  $ L^{'}=2. $  If there is complex characteristic root
 $ \lambda_{1}=a+bi $  of  $ A $  and the root is not real, then  $ r\in
L\otimes L $  is a solution of CYBE iff  $ r $  is strongly
symmetric, or  $ r=p(e_{1}\otimes e_{2})+q(e_{2}\otimes
e_{1})+x(e_{1}\otimes e_{1})+y(e_{2}\otimes e_{2}) $  for any  $ p,
q, x, y \in R $.
\end{Example}

{\bf Proof.} There are two different characteristic roots:
 $ \lambda_{1}=a+bi $  and  $ \lambda_{2}=a-bi $, where  $ a, b\in  R $, Thus
 $ A $  must be similar to $
       \left(\begin{array}{cc}
       a& b\\
         -b& a
         \end{array}\right )
        $.  By Theorem \ref {11.7.1}, we can complete the proof. $\Box$

\section{The solutions of CYBE with char  $ k $ =2}

In this section, we find the general solution of CYBE for Lie
algebra  $ L $  with dim  $ L $ =3 and dim  $ L^{'}=2 $, where char
 $ k $ =2.

\begin{Theorem} \label{11.8.1}
Let  $ L $  be a Lie algebra with a basis  $ {e_{1}, e_{2}, e_{3}} $
such that  $ [e_{1}, e_{2}] $ =0,  $ [e_{1}, e_{3}]=\alpha
e_{1}+\beta e_{2} $,  $ [e_{2}, e_{3}]=\gamma e_{1}+\delta e_{2} $,
where  $ \alpha, \beta, \gamma, \delta  \in  k $, and  $ \alpha
\delta -\beta \gamma \neq 0 $. Let $ p, q, s, t, u, v, x, y, z \in k
$  and  $ A=
       \left(\begin{array}{cc}
       \alpha & \gamma \\
         \beta & \delta
         \end{array}\right )
      $.

 (I) If two characteristic roots of  $ A $  are equal and  $ A $  is
 similar to a diagonal matrix in the algebraic closure  $ P $  of
  $ k $, then  $ r $  is a solution of CYBE in  $ L $  for any  $ r \in L\otimes L $;

 (II) If the condition in Part (1)does not hold, then  $ r $  is a
 solution of CYBE in  $ L $  iff  $ r=p(e_{1}\otimes  e_{1})+p(e_{2}\otimes
e_{1})+s(e_{1}\otimes  e_{3})+s(e_{3}\otimes  e_{1})+u(e_{2}\otimes
e_{3})-+u(e_{3}\otimes  e_{2})+x(e_{1}\otimes  e_{1})+y(e_{2}\otimes
e_{2})+z(e_{3}\otimes  e_{3}) $  with $ \alpha us+\alpha pz+\gamma
u^{2}+\gamma yz+\beta s^{2}+\beta xz+\delta su+\delta zp=0 $  and
 $ z\neq 0 $; or $ r=p(e_{1}\otimes  e_{2})+q(e_{2}\otimes
e_{1})+s(e_{1}\otimes  e_{3})+s(e_{3}\otimes  e_{1})+u(e_{2}\otimes
e_{3})+u(e_{3}\otimes  e_{2})+x(e_{1}\otimes  e_{1})+y(e_{2}\otimes
e_{2})$  with  $ s\gamma (p+q)=u\beta (p+q)=u\gamma (p+q)=s\beta
(p+q)=(\alpha +\delta )us+\gamma u^{2}+\beta s^{2} $  = $ s(\alpha
+\delta )(p+q)=u(\alpha +\delta )(p+q)=0 $.

\end{Theorem}

{\bf Proof}
 We only show the necessity since the sufficiency  can easily be shown.
By the proof of Theorem  \ref {11.7.1}, there exists an invertible
matrix  $ Q $  such that  $ (e_{1}^{'}, e_{2}^{'},
e_{3}^{'})=(e_{1}, e_{2}, e_{3})Q $  and  $ [e_{1}^{'}, e_{2}^{'}]=0
$,  $ [e_{1}^{'}, e_{3}^{'}]=e_{1}^{'}+\beta ^{'}e_{2}^{'} $,
 $ [e_{2}^{'}, e_{3}^{'}]=\delta ^{'}e_{2}^{'} $. We use the notations
before Lemma \ref {11.6.2}.  By [11, Proposition2.4],
 $ k_{i3}^{'}=k_{3i}^{'} $  for  $ i $ =1, 2, 3, which implies that
 $ k_{3i}=k_{i3} $  for  $ i $ =1, 2, 3 by Lemma \ref {11.6.2}.

(I) If  $ A $  is similar to   $
                            \left(\begin{array}{cc}
                              \lambda_{1}& 0\\
                                0& \lambda_{1}
                              \end{array}\right )
                            $, then  $ \beta ^{'}=0 $  and
  $ \delta'=1 $. By [11, Proposition 2.4], we have Part (I).

 (II) Let  $ A $  be not similar to   $
                                  \left(\begin{array}{cc}
                                    \lambda_{1}& 0\\
                                     0& \lambda_{1}
                                   \end{array}\right )
                                    $.

  (a). If  $ z\neq 0 $,  then  $ k_{33}^{'}\neq 0 $ \. Thus
   $ k_{12}^{'}=k_{21}^{'} $  by \cite [Proposition 2.4]{Zh99f}, which implies
   $ k_{12}=k_{21} $. It is straightforward to check that relation
  (1)-(22) in the proof of Theorem \ref {11.7.1} hold iff
   $ \alpha us+\alpha pz+\gamma u^{2}+\gamma yz+\beta s^{2}+\beta xz+\delta su+\delta zp=0. $

  (b). If  $ z $ =0, then we can obtain that  $ r $  is the second case
  in Part (II) by using the method similar to the proof of
  Theorem \ref {11.7.1}. $\Box$

  \section{Coboundary Lie bialgebras}

In this section, using the general solution, which are obtained in
the section above, of CYBE in Lie algebra  $ L $  with dim $ L $ =3
and dim $ L^{'} $ =2, we give the sufficient and necessary
conditions which ( $ L, [\ ], \Delta_{r}, r $ ) is a coboundary (or
triangular)Lie bialgebra.

\begin{Theorem} \label {11.9.1}
Let  $ L $  be a Lie algebra with a basis  $ {e_{1}, e_{2}, e_{3}} $
such that $ [e_{1}, e_{2}] $ =0,  $ [e_{1}, e_{3}]=\alpha
e_{1}+\beta e_{2} $,  $ [e_{2}, e_{3}]=\gamma e_{1}+\delta e_{2} $,
where  $ \alpha, \beta, \gamma, \delta  \in  k $, and  $ \alpha
\delta -\beta \gamma \neq 0 $. Set
  $ A=
       \left(\begin{array}{cc}
      \alpha & \gamma \\
         \beta & \delta
         \end{array}\right )
      $. Let  $ p, u, s\in  k $,
      $ r\in  Im(1-\tau) $  and  $ r=p(e_{1}\otimes  e_{2})-p(e_{2}\otimes
     e_{1})+s(e_{1}\otimes  e_{3})-s(e_{3}\otimes  e_{1})+u(e_{2}\otimes  e_{3})-u(e_{3}\otimes
     e_{2}) $.
Then

(I) \ ( $ L,  [\ ], \Delta_{r},  r $ ) is a coboundary  Lie
bialgebra iff  $$ (s, u)
                       \left(\begin{array}{cc}
                       \beta \delta +\alpha \beta & -\beta \gamma -\alpha ^{2}\\
                         \delta ^{2}+\gamma \beta & -\delta \gamma -\gamma \alpha
                        \end{array}\right )
                                                   \left(\begin{array}{cc}
                             s\\
                             u
                             \end{array}\right )
                            =0; $$

(II) If two characteristic roots of  $ A $  are equal and  $ A $  is
similar to a diagonal matrix in the algebraic closure  $ P $  of  $
k $  with char $ k $ =2,  then ( $ L, [\ ], \Delta _ r,  r $ ) is a
triangular Lie bialgebra for any  $ r\in  Im(1-\tau ) $;

(III) If the condition in Part(II) does not hold, then( $ L, [\ ],
\Delta _r,  r $ ) is a triangular Lie bialgebra iff  $ -\beta
s^{2}+\gamma u^{2}+(\alpha -\delta )us=0 $.
\end {Theorem}

{\bf Proof } We can complete the proof as in the proof of \cite
[Theorem 3.3]{Zh98b}. $\Box$

\begin{Example} \label{11.9.2}

Under Example \ref {11.7.2}, and  $ r\in  Im(1-\tau) $, we have the
following :

(I) ( $ L, [\ ], \Delta_r,  r $ ) is a coboundary Lie bialgebra iff
 $ r=p(e_{1}\otimes  e_{2})-p(e_{2}\otimes  e_{1})+s(e_{1}\otimes
e_{3})-s(e_{3}\otimes  e_{1})+u(e_{2}\otimes  e_{3})-u(e_{3}\otimes
e_{2}) $  with  $ a(s^{2}+u^{2})=0 $;

(II) ( $ L, [\ ], \Delta_r, r $ )is a triangular Lie bialgebra iff
 $ r=p(e_{1}\otimes  e_{2})-p(e_{2}\otimes  e_{1}) $.

\end {Example}

        \newpage

\chapter {Classification of PM Quiver Hopf Algebras } \label {c14}

In noncommutative algebras quivers and associated path algebras are
intensively studied because of a remarkable interaction between
homology theory, algebraic geometry, and Lie theory. Ringel's
approach to quantum groups via quivers suggested that there may be
an interesting overlapping also between the theory of quivers and
Hopf algebras. Some years ago Cibils and Rosso   \cite {CR02, CR97}
started to study quivers admitting a (graded) Hopf algebra structure
and there are many papers (e.g. \cite {CHYZ04, OZ04} ) follow these
works. The representation theory of (graded) Hopf algebras has
applications in many branches of physics and mathematics such as the
construction of solutions to the quantum Yang-Baxter equation and
topological invariants (see e.g. \cite{Ma90a,RT90}. Quivers have
also been used in gauge theory and string theory (see e.g.
\cite{Zh05,RR04}).

This chapter can be viewed an extension of the analysis of Cibils
and Rosso. Let $G$ be a group and $kG$ be the group algebra of $G$
over a field $k$. It is well-known \cite{CR97} that the $kG$-Hopf
bimodule category $^{kG}_{kG}{\mathcal M}^{kG}_{kG}$ is equivalent
to the direct product category $\prod_{C\in{\mathcal K}(G)}{\mathcal
M}_{kZ_{u(C)}}$, where ${\mathcal K}(G)$ is the set of conjugate
classes in $G$, $u:{\mathcal K}(G)\rightarrow G$ is a map such that
$u(C)\in C$ for any $C\in {\mathcal K}(G)$, $Z_{u(C)}=\{g\in G\mid
gu(C)=u(C)g\}$ and ${\mathcal M}_{kZ_{u(C)}}$ denotes the category
of right $kZ_{u(C)}$-modules (see \cite [Proposition 3.3]{CR97}
\cite [Theorem 4.1]{CR02} and Theorem \ref{14.1.1}). Thus for any
Hopf quiver $(Q,G,r)$, the $kG$-Hopf bimodule structures on the
arrow comodule $kQ_1^c$ can be derived from the right
$kZ_{u(C)}$-module structures on $^{u(C)}\! (kQ_1^c)^1$ for all
$C\in {\mathcal K}(G)$. If the arrow comodule $kQ_1^c$ admits a
$kG$-Hopf bimodule structure, then there exist six graded Hopf
algebras: co-path Hopf algebra $kQ^c$, one-type-co-path Hopf
algebras $kG [kQ_1^c]$, one-type-path  Hopf algebras $(kG)^*
[kQ_1^a]$,  semi-path Hopf algebra $kQ^s$, semi-co-path Hopf algebra
$kQ^{sc}$ and path Hopf algebra $kQ^a$. We call these Hopf algebras
quiver Hopf algebras over $G$. If the corresponding
$kZ_{u(C)}$-modules $^{u(C)}\! (kQ_1^c)^1$ is pointed (i.e. it is
zero or a direct sum of one dimensional $kZ_{u(C)}$-modules) for all
$C\in {\mathcal K}(G)$, then $kQ_1^c$ is called a {\rm PM} $kG$-Hopf
bimodule. The six graded Hopf algebras derived from the {\rm PM}
$kG$-Hopf bimodule $kQ_1^c$ are called {\rm PM} quiver Hopf
algebras. A Yetter-Drinfeld
 $H$-module $V$  is {\it pointed} if $V=0$ or $V$ is a direct
sum of one dimensional {\rm YD} $H$-modules. If $V$ is a pointed
{\rm YD} $H$-module, then the corresponding Nichols algebra
${\mathcal B}(V)$ is called a {\rm PM} {\it Nichols algebra}. For
example, when $G$ is a finite abelian group of exponent $m$ and $k$
contains a primitive $m$-th root of 1 (e.g. \ $k$ is the complex
field ), all quiver Hopf algebras are {\rm PM} ( see Lemma \ref
{14.1.2}), all {\rm YD } $kG$-modules are
 pointed and all Nichols algebras of all {\rm YD}
$kG$-modules are {\rm PM}
 (see Lemma \ref {14.2.3}).

The aim of this chapter is to provide parametrization of {\rm PM}
quiver Hopf algebras, multiple Taft algebras and {\rm PM} Nichols
algebras. In other words, this chapter provides a kind of
classification of these Hopf algebras.
 This chapter is organized as follows. In
Section \ref {s14.1}, we examine the {\rm PM} quiver Hopf algebras
by means of ramification system with characters. In Section \ref
{s14.2}, we describe {\rm PM} Nichols algebras and multiple Taft
algebras by means of element system with characters. In Section \ref
{s14.3}, we show that the diagram of a quantum weakly commutative
multiple Taft algebra is not only a Nichols algebra but also a
quantum linear space in $^{kG}_{kG}{\cal YD}$; the diagram of a
semi-path Hopf algebra of ${\rm ESC}$   is  a quantum tensor algebra
in $^{kG}_{kG}{\cal YD}$; the quantum enveloping algebra of a
complex semisimple Lie algebra is a quotient of a semi-path Hopf
algebra.

%%\section{\bf Preliminaries}\label{s14.0}

Throughout this chapter, we work over a fixed field $k$. All
algebras, coalgebras, Hopf algebras, and so on, are defined over
$k$; dim, $\otimes$ and Hom stand for ${\rm dim}_k$, $\otimes_k$ and
${\rm Hom}_k$, respectively. Books \cite {DNR01,Mo93,Sw69a} provide
the necessary background for Hopf algebras and book \cite {ARS95}
provides a nice description of the path algebra approach.

Let $\mathbb{Z}$, $\mathbb{Z}^+$ and $\mathbb{N}$ denote sets of all
integers, all positive integers and all non-negative integers,
respectively. For sets $X$ and $Y$, we denote by $|X|$ the cardinal
number of $X$ and by $X^Y$ or $X^{\mid \! Y \! \mid }$ the Cartesian
product $ \Pi _{y \in Y} X_y$ with $X_y = X$ for any $y \in Y$. If
$X$ is finite, then $|X|$ is the number of elements in $X$. If $X =
\oplus _{i\in I} X_{(i)}$ as vector spaces, then we denote by $\iota
_i$ the natural injection from $X_{(i)}$ to $X$ and by $\pi _i$ the
corresponding projection from $X$ to $X_{(i)}$. We will use $\mu$ to
denote the multiplication  of an algebra, $\Delta$ to denote the
comultiplication of a coalgebra, $\alpha ^-$, $\alpha ^+$, $\delta
^-$ and $\delta ^+$ to denote the left module, right module, left
comodule and right comodule structure maps, respectively. The
Sweedler's sigma notations for coalgebras and comodules are $\Delta
(x) = \sum x_{(1)}\otimes x_{(2)}$, $\delta ^- (x)= \sum x_{(-1)}
\otimes x_{(0)}$, $\delta ^+ (x)= \sum x_{(0)} \otimes x_{(1)}$. Let
$G$ be a group. We  denote by $Z(G)$ the center of $G$. Let
$\widehat G$ denote  the set  of characters of all one-dimensional
representations  of $G$. It is clear that $\widehat G$ = $\{ \chi
\mid $  $\chi $ is a group homomorphism from $G$ to the
multiplicative group of all non-zero elements in $k$ \}.

A quiver $Q=(Q_0,Q_1,s,t)$ is an oriented graph, where  $Q_0$ and
$Q_1$ are the sets of vertices and arrows, respectively; $s$ and $t$
are two maps from  $Q_1$ to $Q_0$. For any arrow $a \in Q_1$, $s(a)$
and $t(a)$ are called its start vertex and end vertex, respectively,
and $a$ is called an arrow from $s(a)$ to $t(a)$. For any $n\geq 0$,
an $n$-path or a path of length $n$ in the quiver $Q$ is an ordered
sequence of arrows $p=a_na_{n-1}\cdots a_1$ with $t(a_i)=s(a_{i+1})$
for all $1\leq i\leq n-1$. Note that a 0-path is exactly a vertex
and a 1-path is exactly an arrow. In this case, we define
$s(p)=s(a_1)$, the start vertex of $p$, and $t(p)=t(a_n)$, the end
vertex of $p$. For a 0-path $x$, we have $s(x)=t(x)=x$. Let $Q_n$ be
the set of $n$-paths. Let $^yQ_n^x$ denote the set of all $n$-paths
from $x$ to $y$, $x, y\in Q_0$. That is, $^yQ_n^x=\{p\in Q_n\mid
s(p)=x, t(p)=y\}$.

A quiver $Q$ is {\it finite} if $Q_0$ and $Q_1$ are finite sets. A
quiver $Q$ is {\it locally finite} if $^yQ_1^x$ is a finite set for
any $x, y\in Q_0$.

Let $G$ be a group. Let ${\mathcal K}(G)$ denote the set of
conjugate classes in $G$. A formal sum $r=\sum_{C\in {\mathcal
K}(G)}r_CC$  of conjugate classes of $G$  with cardinal number
coefficients is called a {\it ramification} (or {\it ramification
data} ) of $G$, i.e.  for any $C\in{\mathcal K}(G)$, \  $r_C$ is a
cardinal number. In particular, a formal sum $r=\sum_{C\in {\mathcal
K}(G)}r_CC$  of conjugate classes of $G$ with non-negative integer
coefficients is a ramification of $G$.

 For any ramification $r$ and a $C \in {\cal K}(G)$, since $r_C$ is
 a cardinal number,
we can choice a set $I_C(r)$ such that its cardinal number is $r_C$
without lost generality.
 Let ${\mathcal K}_r(G):=\{C\in{\mathcal
K}(G)\mid r_C\not=0\}=\{C\in{\mathcal K}(G)\mid
I_C(r)\not=\emptyset\}$.  If there exists a ramification $r$ of $G$
such that the cardinal number of $^yQ_1^x$ is equal to $r_C$ for any
$x, y\in G$ with $x^{-1}y \in C\in {\mathcal K}(G)$, then $Q$ is
called a {\it Hopf quiver with respect to the ramification data
$r$}. In this case, there is a bijection from $I_C(r)$ to $^yQ_1^x$,
and hence we write  ${\ }^yQ_1^x=\{a_{y,x}^{(i)}\mid i\in I_C(r)\}$
for any $x, y\in G$ with $x^{-1}y \in C\in {\mathcal K}(G)$. Denote
by $ (Q, G, r)$ the Hopf quiver of $G$ with respect to $r$.

 The coset decomposition of $Z_{u(C)}$ in $G$ is
\begin {eqnarray} \label {e0.1}
G &=&\bigcup_{\theta\in\Theta_C}Z_{u(C)}g_{\theta},
\end {eqnarray}
where $\Theta_C$ is an index set. It is easy to check that
$|\Theta_C|=|C|$. We always assume that the representative element
of the coset $Z_{u(C)}$ is the identity $1$ of $G$. We claim that
$\theta=\eta$ if
$g_{\theta}^{-1}u(C)g_{\theta}=g_{\eta}^{-1}u(C)g_{\eta}$. In fact,
by the equation
$g_{\theta}^{-1}u(C)g_{\theta}=g_{\eta}^{-1}u(C)g_{\eta}$ one gets
$g_{\eta} g_{\theta}^{-1}u(C) = u(C)g_{\eta} g_{\theta}^{-1}$. Thus
$g_{\eta} g_{\theta}^{-1}\in Z_{u(C)}$, and hence $\theta= \eta$.
For any $x, y \in G$ with $x^{-1}y\in C\in {\cal K}(G)$, there
exists a unique $\theta\in \Theta _C$ such that \begin {eqnarray}
\label {e0.2} x^{-1}y = g_{\theta}^{-1}u(C)g_{\theta}.
\end {eqnarray}
Without specification, we will always assume that $x, y , \theta$
and $C$ satisfy the above relation (\ref{e0.2}). Note that $\theta$
is only determined by $x^{-1}y$. For any $h\in G$ and
$\theta\in\Theta_C$, there exist unique $h'\in Z_{u(C)}$ and
$\theta'\in\Theta_C$ such that $g_{\theta}h = h'g_{\theta'}$. Let
$\zeta_{\theta}(h)=h'$. Then we have
\begin {eqnarray} \label {e0.3} g_{\theta}h&=&\zeta_{\theta}(h)g_{\theta'}.
\end {eqnarray}
If $u(C)$ lies in the center $Z(G)$ of $G$, we have
$\zeta_{\theta}=id_G$. In particular, if G is abelian, then
$\zeta_{\theta}=id_G$ since $Z_{u(C)}=G$.

Let $H$ be a Hopf algebra. A (left-left) {\it Yetter-Drinfeld
module} $V$ over $H$ (simply, {\rm YD} $H$-module) is simultaneously
a left $H$-module and a left $H$-comodule satisfying the following
compatibility condition:
\begin{eqnarray}\label{ydm}
\sum(h\cdot v)_{(-1)}\otimes(h\cdot v)_{(0)}=\sum
h_{(1)}v_{(-1)}S(h_{(3)})\otimes h_{(2)}\cdot v_{(0)},\ \ v\in V,
h\in H.
\end{eqnarray}
We denote by $^H_H{\mathcal YD}$ the category of {\rm YD}
$H$-modules; the morphisms in this category preserve both the action
and the coaction of $H$.

The structure of a Nichols algebra appeared first in the paper
\cite{Ni78},  and N. Andruskiewitsch and H. J. Schneider used it to
classify finite-dimensional pointed Hopf algebras \cite {AS98a,
AS98b, AS02, AS00}. Its definition can be found in \cite [Definition
2.1] {AS02}.

If  $\phi: A\rightarrow A'$ is  an algebra homomorphism and  $(M,
\alpha ^-)$ is a left $A'$-module, then $M$ becomes a left
$A$-module with the $A$-action given by $a \cdot x =\phi (a) \cdot x
$ for any $a\in A$, $x\in M$, called a pullback $A$-module through
$\phi$, written as  $_{\phi}M$.  Dually, if  $\phi: C\rightarrow C'$
be a coalgebra homomorphism and  $(M, \delta ^- )$ is  a left
$C$-comodule, then $M$ is a left $C'$-comodule with the
$C'$-comodule structure given by $ {\delta'}^-:=(\phi\otimes{\rm
id})\delta^-$, called  a push-out $C'$-comodule through $\phi$,
written as  $^{\phi}M$.

Let $A$ be an algebra and $M$ be an $A$-bimodule. Then the tensor
algebra $T_A(M)$ of $M$ over $A$ is a graded algebra with
$T_A(M)_{(0)}=A$, $T_A(M)_{(1)}=M$ and $T_A(M)_{(n)}=\otimes^n_AM$
for $n>1$. That is, $T_A(M)=A\oplus(\bigoplus_{n>0}\otimes^n_AM)$.
Let $D$ be another algebra. If $h$ is an algebra map from $A$ to $D$
and $f$ is an $A$-bimodule map from $M$ to $_hD_h$, then by the
universal property of $T_A(M)$ (see \cite [Proposition 1.4.1]
{Ni78}) there is a unique algebra map $T_A(h,f): T_A(M)\rightarrow
D$ such that $T_A(h,f)\iota_0=h$ and $T_A(h,f)\iota_1=f$.  One can
easily see that $T_A ( h, f ) = h + \sum _{n>0} \mu ^{n-1}T_n (f )$,
where $T_n(f)$ is the map from $\otimes _A^n M$ to $\otimes_A^nD$
given by $T_n(f)(x_1\otimes x_2 \otimes \cdots \otimes x_n) =
f(x_1)\otimes f(x_2) \otimes \cdots \otimes f(x_n)$, i.e.,
$T_n(f)=f\otimes _A f\otimes_A\cdots\otimes_A f$. Note that $\mu$
can be viewed as a map from $D\otimes _A D$ to $D$.

Dually, let $C$ be a coalgebra and let $M$ be a $C$-bicomodule. Then
the cotensor coalgebra $T_C^c(M)$ of $M$ over $C$ is a graded
coalgebra with $T_C^c(M)_{(0)}=C$, $T_C^c(M)_{(1)}=M$ and
$T_C^c(M)_n=\Box^n_CM$ for $n>1$. That is,
$T_C^c(M)=C\oplus(\bigoplus_{n>0}\Box^n_CM)$. Let $D$ be another
coalgebra. If $h$ is a coalgebra map from $D$ to $C$ and $f$ is a
$C$-bicomodule map from $^hD^h$ to $M$ such that $f({\rm
corad}(D))=0$, then by the universal property of $T_C^c(M)$ (see
\cite [Proposition 1.4.2] {Ni78}) there is a unique coalgebra map
$T_C^c(h,f)$ from $D$ to $T_C^c(M)$ such that $\pi_0T_C^c(h,f)=h$
and $\pi_1T_C^c(h,f)=f$. It is not difficult to see that
$T_C^c(h,f)=h+\sum_{n>0}T_n^c(f)\Delta_{n-1}$, where $T_n^c(f)$ is
the map from $\Box_C^n D$ to $\Box_C^n M $ induced by
$T_n^c(f)(x_1\otimes x_2 \otimes \cdots \otimes x_n)=f(x_1)\otimes
f(x_2) \otimes \cdots \otimes f(x_n)$, i.e., $T^c_n(f)=f\otimes
f\otimes\cdots\otimes f$.

\section{\bf {\rm PM} quiver Hopf algebras }\label{s14.1}

In this section we describe {\rm PM} quiver Hopf algebras by
parameters.

We first describe the category of Hopf bimodules by categories of
modules.

Let $G$ be a group and let $(B, \delta ^-, \delta ^+)$ be a
$kG$-bicomodule. Then the $(x, y)$-isotypic component of $B$ is
$$^yB^x= \{b \in B \mid  \delta ^-(b) = y \otimes b, \delta ^+ (b)
  = b \otimes x  \},$$
where $x, y \in G$. Let $M$ be another $kG$-bicomodule and $f:
B\rightarrow M$ be a $kG$-bicomodule homomorphism. Then
$f(^yB^x)\subseteq\ ^yM^x$ for any $x, y\in G$. Denote by $^yf^x$
the restriction map $f|_{^yB^x}:\ ^yB^x\rightarrow\ ^yM^x$, $x, y\in
G$.

\begin {Theorem} \label {14.2} (See \cite [Proposition 3.3]{CR97} and \cite [Theorem 4.1]{CR02})
The category $^{kG}_{kG}{\mathcal M}^{kG}_{kG}$ of $kG$-Hopf
bimodules is equivalent to the Cartesian product category  $\prod_{C
\in {\mathcal K}(G)} {\mathcal M}_{kZ_{u(C)}}$ of categories
${\mathcal M}_{kZ_{u(C)}}$  of right $kZ_{u(C)}$-modules for all $C
\in {\mathcal K}(G)$.
\end {Theorem}

For later use, we give the mutually inverse functors between the two
categories here. The functors $W$ from $^{kG}_{kG}{\mathcal
M}^{kG}_{kG}$ to $\prod_{C \in {\mathcal K}(G)}{\mathcal
M}_{kZ_{u(C)}}$ is defined by
$$W(B)=\{^{u(C)}\! B^1\}_{C\in{\mathcal K}(G)},\
W(f)=\{^{u(C)}\! f^1\}_{C\in{\mathcal K}(G)}$$ for any object $B$
and morphism $f$ in $^{kG}_{kG}{\mathcal M}^{kG}_{kG}$, where the
right $kZ_{u(C)}$-module action on ${}^{u(C)}\!  B^1$ is given by
\begin {eqnarray}\label {th1e1}
b \lhd  h = h^{-1} \cdot b \cdot h,\ \ b\in\ ^{u(C)}\! B^1,\ h\in
Z_{u(C)}.
\end {eqnarray}
The functor $V$ from $\prod_{C\in{\mathcal K}(G)}{\mathcal
M}_{kZ_{u(C)}}$ to $^{kG}_{kG}{\mathcal M}^{kG}_{kG}$ is defined as
follows:

For $M =\{ M(C)\}_{C \in {\mathcal K}(G)}\in \prod_{C \in {\mathcal
K}(G)} {\mathcal M}_{kZ_{u(C)}}$, $V(M)$ is given by
\begin {eqnarray}\label {e1.2} \begin {array}{llll}
^yV(M)^x &=& x \otimes M(C)\otimes _{kZ_{u(C)}} g_{\theta},\\
h \cdot (x \otimes m \otimes _{kZ_{u(C)}} g_{\theta})  &=&
hx \otimes m \otimes _{kZ_{u(C)}} g_{\theta},\\
 (x \otimes m \otimes _{kZ_{u(C)}} g_{\theta}) \cdot h  &=&
xh \otimes m \otimes _{kZ_{u(C)}} g_{\theta}h=xh \otimes (m\lhd
\zeta_{\theta}(h))\otimes_{kZ_{u(C)}}g_{\theta'},
\end {array}
\end {eqnarray}
where $h, x, y \in G$ with $x^{-1}y\in C$ and the relation
(\ref{e0.2}) and (\ref {e0.3}), $m\in M(C)$. For any morphism
$f=\{f_C\}_{C\in{\mathcal K}(G)}: \{M(C)\}_{C\in{\mathcal
K}(G)}\rightarrow \{N(C)\}_{C\in{\mathcal K}(G)}$, $V(f)(x\otimes
m\otimes_{kZ_{u(C)}}g_{\theta})=x\otimes
f_C(m)\otimes_{kZ_{(u(C))}}g_{\theta}$ for any $m\in M_{C}$, $x, y
\in G$ with $x^{-1}y=g^{-1}_{\theta}u(C)g_{\theta}$. That is,
$^yV(f)^x={\rm id}\otimes f_C\otimes{\rm id}$.

An $A$-module $M$ is called {\it pointed} if $M=0$ or $M$ is a
direct sum of one dimensional $A$-modules.

\begin{Definition}\label{14.1.1}
A $kG$-Hopf bimodule $N$ is called a {\it
 $kG$-Hopf bimodule with pointed module structure ( or a  \ {\rm PM} $kG$-Hopf bimodule in
short) } if there exists an object \ \ \  \ \ \ \ $\prod_{C \in
{\mathcal K}(G)} M(C) $ in $\prod_{C \in {\mathcal K}(G)} {\mathcal
M}_{kZ_{u(C)}}$ such that $M(C)$ is a right pointed
$kZ_{u(C)}$-module for  any $C \in {\mathcal K }(G)$ and $N\cong
V(\{ M(C)\}_{C \in {\mathcal K}(G)}):=\bigoplus _{y = g_\theta ^{-1}
u(C) g_\theta, \  x, y \in G} \ x \otimes M(C)\otimes _{kZ_{u(C)}}
g_{\theta}$ as $kG$-Hopf bimodules. Here $V(\{ M(C)\}_{C \in
{\mathcal K}(G)}):=\bigoplus _{y = g_\theta ^{-1} u(C) g_\theta, \
x, y \in G}\   x \otimes M(C)\otimes _{kZ_{u(C)}} g_{\theta}$ is
defined in the proof of Theorem \ref {14.2}.
\end {Definition}

Lemma \ref {14.1.2} -- \ref {14.1.6}, Theorem \ref {14.3} and Lemma
\ref {14.2.7} are well-known.

\begin{Lemma}\label{14.1.2}  Assume that $G$ is a finite commutative group
of exponent $m$. If $k$ contains a primitive $m$-th root of 1, then
(i) every $kG$-module is a pointed module;
 (ii)  every $kG$-Hopf
bimodule is {\rm PM}.
\end {Lemma}

\begin{Lemma}\label{14.1.3} Let  $H=\bigoplus_{i\geq 0}H_{(i)}$ be  a graded Hopf algebra.
Set $B:=H_{(0)}$ and $M:=H_{(1)}$. Then $M$ is a $B$-Hopf bimodule
with the $B$-actions and $B$-coactions given by
\begin {eqnarray*}\label {e2.1.3}
\alpha^-=\pi_1\mu(\iota_0\otimes\iota_1),&
&\alpha^+=\pi_1\mu(\iota_1\otimes\iota_0), \
\delta^-=(\pi_0\otimes\pi_1)\Delta\iota_1, \
 \delta^+ =(\pi_1\otimes\pi_0)\Delta\iota_1.
\end {eqnarray*}
%%% \begin {eqnarray}\label {e2.1} \end {eqnarray}
\end {Lemma}

\begin{Lemma}\label{14.1.4}
{\rm (i)} \ Let $B$ be an algebra and $M$ be a $B$-bimodule. Then
the tensor algebra $T_B(M)$ of $M$ over $B$ admits a graded Hopf
algebra structure if and only if $B$ admits a Hopf algebra structure
and $M$ admits a $B$-Hopf bimodule structure.

 {\rm (ii)} Let $B$ be
a coalgebra and  $M$ be a $B$-bicomodule. Then the cotensor
coalgebra $T_B^c(M)$ of $M$ over $B$ admits a graded Hopf algebra
structure if and only if $B$ admits a Hopf algebra structure and $M$
 admits a $B$-Hopf bimodule structure.
\end {Lemma}

Let $B$ be a bialgebra (Hopf algebra) and $M$ be a $B$-Hopf
bimodule. Then $T_B^c(M)$ is a graded bialgebra (Hopf algebra) by
Lemma \ref{14.1.4}. Let $B[M]$ denote the subalgebra of $T_B^c(M)$
generated by $B$ and $M$. Then $B[M]$ is a bialgebra (Hopf algebra)
of type one by \cite [section 2.2, p.1533]{Ni78}.  $B[M]$ is a
graded subspace of  $T_B^c (M).$
\begin {Lemma} \label {14.1.5}
Let $B$ and $B'$ be two  Hopf algebras. Let $M$ and $M'$ be $B$-Hopf
bimodule and  $B'$-Hopf bimodule, respectively. Assume that $\phi:
B\rightarrow B'$ be a Hopf algebra map. If $\psi $ is simultaneously
a $B$-bimodule and $B'$-bicomodule map from $^{\phi}M^{\phi}$ to
$_{\phi}M'_{\phi}$,  then {\rm(i)} \ $T_B(\iota_0\phi,
\iota_1\psi):=\iota_0\phi+ \sum_{n>0}\mu^{n-1}T_n(\iota_1\psi)$ is a
graded Hopf algebra map from $T_B(M)$ to $T_{B'}(M')$.
 {\rm(ii)} $T_{B'}^c(\phi\pi_0,
\psi\pi_1):=\phi\pi_0+ \sum_{n>0}T_n^c (\psi\pi_1)\Delta_{n-1}$ is a
graded Hopf algebra map from $T_B^c(M)$ to $T_{B'}^c(M')$.
\end {Lemma}

\begin {Lemma} \label {14.1.6} Let $B$ and $B'$ be two  Hopf algebras. Let
 $M$ and $M'$ be  $B$-Hopf
bimodule and  $B'$-Hopf bimodule, respectively.   Then the following
statements are equivalent: {\rm(i)} There exists a Hopf algebra
isomorphism $\phi: B\rightarrow B'$ such that $M\cong\ _{\phi}
^{\phi ^{-1}}M'{}_{\phi} ^{\phi ^{-1}}$ as $B$-Hopf bimodules.
 {\rm(ii)} $T_{B}(M)$ and $T_{B'}(M')$
are isomorphic as graded Hopf algebras.
 {\rm(iii)} $T_{B}^c(M)$ and
$T_{B'}^c(M')$ are isomorphic as graded Hopf algebras.
 {\rm(iv)} $B[M]$ and
$B'[M']$ are isomorphic as graded Hopf algebras.
\end{Lemma}

Let $Q=(G, Q_1, s, t)$ be a quiver of a group $G$. Then $kQ_1$
becomes  a $kG$-bicomodule  under  the natural comodule structures:
\begin{eqnarray}\label{arcom}
\delta^-(a)=t(a)\otimes a,\ \ \delta^+(a)=a\otimes s(a),\ \ a\in
Q_1,
\end{eqnarray} called an {\it arrow comodule}, written as $kQ_1 ^c$. In this case,
the path coalgebra $kQ^c$ is exactly isomorphic to the cotensor
coalgebra $T^c_{kG}(kQ_1^c)$  over $kG$ in a natural way (see
\cite{CM97} and \cite{CR02}). We will regard $kQ^c=T^c_{kG}(kQ_1^c)$
in the following. Moreover, when $G$ is finite,  $kQ_1$ becomes a
$(kG)^*$-bimodule with the module structures defined by
\begin{eqnarray}\label{armod}
\mbox{\hspace{1cm}}p\cdot a:=\langle p, t(a)\rangle a,\ \ a\cdot
p:=\langle p, s(a)\rangle a,\ \ p\in(kG)^*, a\in Q_1,
\end {eqnarray}
written as $kQ_1^a$, called an {\it arrow module}. Therefore, we
have a tensor algebra $T_{(kG)^*}(kQ_1^a)$. Note that the tensor
algebra $T_{(kG)^*}(kQ_1^a)$ of $kQ_1^a$ over $(kG)^*$ is exactly
isomorphic to the path algebra $kQ^a$. We will regard
$kQ^a=T_{(kG)^*}(kQ_1^a)$ in the following.

 Assume that $Q$ is a finite  quiver on finite
group $G$. Let $\xi _{kQ_1^a}$ denote the linear map from $kQ_1^a$
to $(kQ_1^c)^*$ by sending $a$ to $a^*$ for  any $a\in Q_1$ and $\xi
_{kQ_1^c}$ denote the linear map from $kQ_1^c$ to $(kQ_1^a)^*$ by
sending $a$ to $a^*$ for  any $a\in Q_1$. Here $\{ a^* \mid a^*  \in
(kQ_1)^*\}$ is the dual basis of $\{ a \mid a \in Q_1 \}$.

\begin {Lemma} \label {14.1.7} Assume that $Q$ is a finite Hopf quiver on finite
group $G$. Then

(i) If $(M, \alpha ^-, \alpha ^+, \delta^-, \delta ^+)$ is a finite
dimensional $B$-Hopf bimodule and $B$ is a finite dimensional Hopf
algebra, then $(M^*, \delta^{-*}, \delta ^{+*}, \alpha ^{-*}, \alpha
^{+*})$ is a $B^*$-Hopf bimodule.

(ii) If $(kQ_1^c, \alpha ^-, \alpha ^+, \delta^-, \delta ^+)$ is a
$kG$- Hopf bimodule, then there exist  unique left $(kG)^*$-comodule
operation $ \delta _{kQ_1^a}^-$ and right $(kG)^*$-comodule
operation
 $\delta
_{kQ_1^a}^+$ such that $(kQ_1^a, \alpha _{kQ_1^a} ^-, \alpha
_{kQ_1^a}^+,$ \ $ \delta_{kQ_1^a}^-, $ \ $ \delta_{kQ_1^a} ^+)$
becomes a $(kG)^*$-Hopf bimodule and $\xi_{kQ_1^a}$ becomes a
$(kG)^*$-Hopf bimodule isomorphism from $(kQ_1^a, \alpha _{kQ_1^a}
^-, \alpha _{kQ_1^a}^+, \delta_{kQ_1^a}^-, \delta_{kQ_1^a} ^+)$ to
$((kQ_1^c)^*, \delta^-{}^*, \delta ^+{}^*, \alpha ^-{}^*, \alpha
^+{}^* )$.

(iii) If $(kQ_1^a, \alpha ^-, \alpha ^+, \delta^-, \delta ^+)$ is a
$(kG)^*$- Hopf bimodule, then there exist  unique left $kG$-module
operation $ \alpha _{kQ_1^c}^-$ and right $kG$-module
 operation $\alpha
_{kQ_1^c}^+$ such that $(kQ_1^c, \alpha _{kQ_1^c} ^-, \alpha
_{kQ_1^c}^+, \delta_{kQ_1^c}^-,$ $ \delta_{kQ_1^c} ^+)$ become a
$kG$-Hopf bimodule and $\xi_{kQ_1^c}$ becomes a $kG$-Hopf bimodule
isomorphism from $(kQ_1^c, \alpha _{kQ_1^c} ^-, \alpha _{kQ_1^c}^+,
\delta_{kQ_1^c}^-, \delta_{kQ_1^c} ^+)$ to $((kQ_1^a)^*,
\delta^-{}^*, \delta ^+{}^*, \alpha ^-{}^*, \alpha ^+{}^* )$.

(iv)
 $\xi_{kQ_1^a}$ is a $(kG)^*$-Hopf
 bimodule isomorphism from $(kQ_1^a, \alpha _{kQ_1^a} ^-, \alpha
 _{kQ_1^a}^+, \delta_{kQ_1^a}^-, \delta_{kQ_1^a} ^+)$ to
 $((kQ_1^c)^*, \delta_{kQ_1^c}^-{}^*, \delta_{kQ_1^c} ^+{}^*, \alpha_{kQ_1^c} ^-{}^*, \alpha _{kQ_1^c}
 ^+{}^* )$ if and only if
 $\xi_{kQ_1^c}$ becomes a $kG$-Hopf bimodule
 isomorphism from $(kQ_1^c, \alpha _{kQ_1^c} ^-, \alpha _{kQ_1^c}^+,
 \delta_{kQ_1^c}^-, \delta_{kQ_1^c} ^+)$ to $((kQ_1^a)^*,$ $
 \delta_{kQ_1^a}^-{}^*, $ $ \delta_{kQ_1^a} ^+{}^*, $ $\alpha _{kQ_1^a}^-{}^*, \alpha _{kQ_1^a}^+{}^* )$.
\end {Lemma}

{\bf Proof.} It is easy to check (i)--(iii). Now we show (iv). Let
$A:= kQ^a_1,$ $ B:= kQ_1^c.$  Let  $\sigma_ B$ denote the canonical
linear isomorphism from  $B$ to $ B^{**}$ by sending $b$ to $b^{**}$
for any $b \in B$, where  $<b^{**}, f> = <f, b>$ for any $f \in
B^*.$ If $A \stackrel {\xi _A } {\cong } B^* $ as $(kG)^*$-Hopf
bimodules, then $B \stackrel {\sigma _B} {\cong } B^{**}  \stackrel
{(\xi _A )^*} {\cong }A^*$ as $kG$-Hopf bimodules. It is easy to
check $\xi _B =(\xi _A)^* \sigma _B $. Therefore $\xi _B$ is a
$kG$-Hopf bimodule isomorphism. Conversely, if $B \stackrel {\xi _B
} {\cong } A^* $ as $kG$-Hopf bimodules, we can similarly show that
$A \cong  B^* $ as $(kG)^*$-Hopf bimodules. \ $\Box$

\begin {Theorem} \label {14.3} (see \cite [Theorem 3.3]{CR02} and \cite [Theorem 3.1] {CR97})
Let $Q$ be a quiver over  group $G$. Then the following two
statements are equivalent:

(i) $Q$ is a Hopf quiver.

(ii) Arrow comodule $kQ_1^c$ admits a $kG$-Hopf bimodule
structure.\\
Furthermore, if $Q$ is finite, then the above are equivalent to the
following:

(iii) Arrow module $kQ_1^a$ admits a $(kG)^*$-Hopf bimodule
structure.
\end {Theorem}

Assume that $Q$ is a Hopf quiver. It follows from Theorem \ref
{14.3} that there exist a left $kG$-module structure $\alpha ^-$ and
a right $kG$-module structure $\alpha ^+$ on arrow comodule
$(kQ_1^c, \delta^-, \delta ^+)$ such that  $(kQ_1^c, \alpha ^-,
\alpha ^+, \delta^-, \delta ^+)$ becomes a $kG$-Hopf bimodule,
called a $kG$-Hopf bimodule with arrow comodule, written  $(kQ_1^c,
\alpha ^-, \alpha ^+ )$ in short. We obtain three graded Hopf
algebras $T_{kG} (kQ_1^c)$, $T_{kG}^c(kQ_1^c)$ and $kG[kQ_1^c]$,
called semi-path Hopf algebra, co-path Hopf algebra and
one-type-co-path Hopf algebra, written $kQ^{s}(\alpha ^-, \alpha
^+)$, $kQ^{c}( \alpha ^-, \alpha ^+)$ and $kG[kQ_1^c, \alpha ^-,
\alpha ^+]$, respectively. Dually, when $Q$ is finite, it follows
from Theorem \ref {14.3} that there exist a left $(kG)^*$-comodule
structure $\delta ^-$ and a right $(kG)^*$-comodule structure
$\delta ^+$ on arrow module $(kQ_1^a, \alpha^-, \alpha ^+)$ such
that  $(kQ_1^a, \alpha ^-, \alpha ^+, \delta^-, \delta ^+)$ becomes
a $(kG)^*$-Hopf bimodule, called a $(kG)^*$-Hopf bimodule with arrow
module, written $(kQ_1^a, \delta^-, \delta ^+)$ in short. We obtain
three  graded Hopf algebras $T_{(kG)^*} (kQ_1^a)$,
$T_{(kG)^*}^c(kQ_1^a)$ and $(kG)^* [kQ_1^a]$, called path Hopf
algebra, semi-co-path Hopf algebra and one-type-path Hopf algebra,
written $kQ^a( \delta ^-, \delta ^+)$, $kQ^{sc}( \delta ^-, \delta
^+)$ and $(kG)^*[kQ_1^a ,\alpha _1, \alpha ^+]$, respectively.  We
call the six graded Hopf algebras the quiver Hopf algebras (over
$G$). We usually omit the (co)module operations when we write these
quiver Hopf algebras.

If $\xi _{kQ_1^a}$ or $\xi _{kQ_1^c}$ is a Hopf bimodule
isomorphism, then, by Lemma \ref {14.1.6} and Lemma \ref {14.1.7},
$T_{(kG)^*} (\iota _0,\iota_1 \xi _{kQ_1^a})$ and $T_{kG} ^c (\pi
_0, \xi _{kQ_1^c} \pi_1)$ are graded Hopf algebra isomorphisms from
$T_{(kG)^*} (kQ_1^a)$ to $T_{(kG)^*} ((kQ_1^c)^*)$ and from $T_{kG}
^c(kQ_1^c)$ to $T_{kG}^c((kQ_1^a)^*)$, respectively; $T_{(kG)^*} ^c
(\pi _0,$ $ \xi _{kQ_1^a} \pi_1)$ and $T_{kG}  (\iota _0, \iota _1
\xi _{kQ_1^c})$ are graded Hopf algebra isomorphisms from
$T_{(kG)^*} ^c(kQ_1^a)$ to $T_{(kG)^*} ^c((kQ_1^c)^*)$ and from
$T_{kG} (kQ_1^c)$ to $T_{kG}((kQ_1^a)^*)$, respectively. In this
case, $(kQ_1^a, kQ_1^c)$, $(kQ^a, kQ^c)$ and $(kQ^{s}, kQ^{sc})$ are
said to be arrow dual pairings.

If $(kQ_1^c,\alpha^-,\alpha^+,\delta^-,\delta^+)$ is a {\rm PM}
$kG$-Hopf bimodule, and $(kQ_1^c, kQ_1^a)$ is an arrow pairing, then
$(kQ_1^a, \alpha ^{-*}, \alpha ^{+*}) $ is called a {\rm PM}
$(kG)^*$-Hopf bimodule and  six quiver Hopf algebras induced by
$kQ_1^c$ and $kQ_1^a$ are called {\rm PM} quiver Hopf algebras.

Now we are going to describe the structure of all {\rm PM} $kG$-Hopf
bimodules and the corresponding graded Hopf algebras.

\begin{Definition}\label {14.1.8}
$(G, r, \overrightarrow \chi, u)$ is called a ramification system
with characters   (or {\rm RSC } in short ), if $r$ is a
ramification of $G$, $u$ is a map from ${\mathcal K}(G)$ to $G$ with
$u(C)\in C$ for any $C\in {\mathcal K}(G)$, and $\overrightarrow
\chi=\{\chi_C^{(i)} \}_ { i\in I_C(r), C\in{\mathcal K}_r(G)} \ \in
\prod _ { C\in{\mathcal K}_r(G)} (\widehat{Z_{u(C)}}) ^{r_C}$ with
$\chi_C^{(i)} \in \widehat{Z_{u(C)}} $ for any
 $ i \in I_C(r), C\in {\mathcal
K}_r(G)$.

${\rm RSC} (G, r, \overrightarrow \chi, u)$ and ${\rm RSC} (G', r',
\overrightarrow {\chi'}, u')$ are said to be {\it isomorphic} if the
following conditions are satisfied:

$\bullet$ There exists a group isomorphism $\phi: G\rightarrow G'$.

$\bullet$ For any $C\in{\mathcal K}(G)$, there exists an element
$h_C\in G$ such that $\phi(h_C^{-1}u(C)h_C)=u'(\phi(C))$.

$\bullet$ For any $C\in{\mathcal K}_r(G)$, there exists a bijective
map $\phi_C : I_C(r) \rightarrow I_{\phi(C)}(r')$ such that $\chi
'{}_{\phi(C)}^{(\phi_C(i))}(\phi(h_C^{-1}hh_C))=\chi_C^{(i)}(h)$ for
all $h\in Z_{u(C)}$ and $i\in I_C(r)$.
\end {Definition}

{\bf Remark.} Assume that $G=G'$, $r=r'$ and $u(C) = u'(C)$ for any
$C \in {\cal K}_r(G)$. If there is a permutation $\phi_C$ on
$I_C(r)$ for any $C\in{\mathcal K}_r(G)$ such that
${\chi'}_C^{(\phi_C(i))}=\chi_C^{(i)}$ for all $i\in I_C(r)$, then
obviously
 ${\rm RSC} (G, r, \overrightarrow \chi, u)\cong {\rm RSC} (G, r, \overrightarrow{\chi'},
 u)$.

\begin{Proposition}\label {14.1.10} If $N$ is a {\rm PM} $kG$-Hopf bimodule, then
there exist a Hopf quiver $(Q, G,r )$, an ${\rm RSC} (G, r,
\overrightarrow \chi, u)$ and a $kG$-Hopf bimodule $(kQ_1^c, \alpha
^-, \alpha ^+)$ with
$$ \alpha ^- (h \otimes a^{(i)}_{y,x}) :=   h\cdot a^{(i)}_{y,x}=a^{(i)}_{hy,hx},\ \
\alpha ^+  ( a^{(i)}_{y,x}\otimes h) := a^{(i)}_{y,x}\cdot
h=\chi^{(i)}_C(\zeta_{\theta}(h))a^{(i)}_{yh,xh}$$ where $x, y, h\in
G$ with $x^{-1}y=g^{-1}_{\theta}u(C)g_{\theta}$, $\zeta_{\theta}$ is
given by {\rm(\ref{e0.3})}, $C\in{\mathcal K}_r(G)$ and $i\in
I_C(r)$,  such that $N \cong (kQ_1^c, \alpha ^-, \alpha ^+)$ as
$kG$-Hopf bimodules.

\end {Proposition}
{\bf Proof.} Since $N$  is a {\rm PM} $kG$-Hopf bimodule, there
exists an object \ \ \  \ \ \ \ $\prod_{C \in {\mathcal K}(G)} M(C)
$ in $\prod_{C \in {\mathcal K}(G)} {\mathcal M}_{kZ_{u(C)}}$ such
that $M(C)$ is a pointed $kZ_{u(C)}$-module for  any $C \in
{\mathcal K }(G)$ and $N \cong V(\{ M(C)\}_{C \in {\mathcal
K}(G)})=\bigoplus _{y = g_\theta ^{-1} u(C) g_\theta, \  x, y \in
G}\   x \otimes M(C)\otimes _{kZ_{u(C)}} g_{\theta}$ as $kG$-Hopf
bimodules. Let $r = \sum _{C \in {\mathcal K} (G)}r_C C$ with $ r_C
={\rm dim}M(C)$ for any $C \in {\mathcal K} (G)$. Notice that ${\rm
dim}M(C)$ denotes the  cardinal number  of a basis  of a basis of
$M(C)$ when $M(C)$ is infinite dimensional.  Since $(M(C), \alpha
_C)$ is a pointed $kZ_{u(C)}$-module, there exist a $k$-basis
$\{x_C^{(i)}\mid i\in I_C(r)\}$ in $M(C)$ and a family of characters
$\{\chi_C^{(i)}\in\widehat{Z_{u(C)}}\mid i\in I_C(r)\}$ such that $
\alpha _C (x_C^{(i)} \otimes  h)= x_C^{(i)}\lhd
h=\chi_C^{(i)}(h)x_C^{(i)}$ for any $i\in I_C(r)$ and $h\in
Z_{u(C)}$.

We have to show that $(kQ_1^c, \alpha ^-, \alpha ^+)$ is isomorphic
to $\bigoplus _{y = g_\theta ^{-1} u(C) g_\theta, \  x, y \in G}\ x
\otimes M(C)\otimes _{kZ_{u(C)}} g_{\theta}$ as $kG$-Hopf bimodules.
Observe that there is a canonical $kG$-bicomodule isomorphism
$\varphi: kQ_1\rightarrow \bigoplus _{y = g_\theta ^{-1} u(C)
g_\theta, \  x, y \in G}\   x \otimes M(C)\otimes _{kZ_{u(C)}}
g_\theta$ given by
\begin{eqnarray}\label{indbya}
\varphi(a^{(i)}_{y,x})=x\otimes
x_C^{(i)}\otimes_{kZ_{u(C)}}g_{\theta}
\end{eqnarray}
where $x, y\in G$ with $x^{-1}y=g^{-1}_{\theta}u(C)g_{\theta}$,
$C\in{\mathcal K}_r(G)$ and $i\in I_C(r)$. Now we have

\begin{eqnarray*}
\varphi (\alpha ^- (h \otimes a^{(i)}_{y,x})) &=& \varphi ( a^{(i)}_{hy,hx})=    hx \otimes x^{(i)}_C \otimes_{kZ_{u(C)}}g_{\theta}  \ \ \\
&=&    h \cdot (x \otimes x^{(i)}_C \otimes_{kZ_{u(C)}}g_{\theta})
\ \ \ (\hbox {see }(
\ref {e1.2} )) \\
&=& h \cdot  \varphi ( a^{(i)}_{y,x})
\end{eqnarray*}
and
\begin{eqnarray*}
\varphi (\alpha ^+ (a^{(i)}_{y,x} \otimes h))
&=&\chi^{(i)}_C(\zeta_{\theta}(h))\varphi ( a^{(i)}_{yh,xh})
= \chi^{(i)}_C(\zeta_{\theta}(h))   xh \otimes x^{(i)}_C \otimes_{kZ_{u(C)}}g_{\theta'} \\
&{}&  (\hbox {since } h g_{\theta'}^{-1}=g_\theta ^{-1} \zeta
_\theta (h)
 \hbox { and }
yh = x h g_{\theta '}^{-1}u(C)g_{\theta '})\\
&=&   (x \otimes x^{(i)}_C \otimes_{kZ_{u(C)}}g_{\theta}) \cdot h\ \ \ (\hbox { see } (\ref {e1.2} )) \\
&=&  \varphi ( a^{(i)}_{y,x}) \cdot h,
\end{eqnarray*}
 where $x, y, h\in
G$ with $x^{-1}y=g^{-1}_{\theta}u(C)g_{\theta}$, $\zeta_{\theta}$ is
given by {\rm(\ref{e0.3})}, $C\in{\mathcal K}_r(G)$ and $i\in
I_C(r)$. Consequently, $\varphi$ is a $kG$-Hopf bimodule
isomorphism. $\Box$

Let $(kQ_1^c, G, r, \overrightarrow \chi, u)$ denote the $kG$-Hopf
bimodule $(kQ_1^c, \alpha ^-, \alpha ^+)$ given in Lemma \ref
{14.1.10}. Furthermore, if $(kQ_1^c, kQ_1^a)$ is an arrow dual
pairing, then we denote the $(kG)^*$-Hopf bimodule $kQ_1^a $ by
$(kQ_1^a, G, r, \overrightarrow \chi, u)$. We obtain six quiver Hopf
algebras $kQ^c (G, r, \overrightarrow \chi, u)$, $kQ^s (G, r,
\overrightarrow \chi, u)$, $kG[ kQ_1^c, G, r, \overrightarrow \chi,
u]$, $kQ^a (G, r, \overrightarrow \chi, u)$, $kQ^{sc} (G, r,
\overrightarrow \chi, u),$ $(kG)^*[ kQ_1^a, G, r, \overrightarrow
\chi,$ $ u]$, called the quiver Hopf algebras determined by ${\rm
RSC} (G,$ $ r,$ $ \overrightarrow \chi,u)$.

From Proposition \ref{14.1.10}, it seems that the right $kG$-action
on $(kQ_1^c,G, r, \overrightarrow \chi, u)$ depends on the choice of
the set $\{g_{\theta}\mid \theta\in\Theta_C\}$ of coset
representatives of $Z_{u(C)}$ in $G$ (see, Eq.(\ref{e0.1})). The
following lemma shows that $(kQ_1^c,G, r, \overrightarrow \chi, u)$
is, in fact, independent of the choice of the coset representative
set $\{g_{\theta}\mid \theta\in\Theta_C\}$, up to $kG$- Hopf
bimodule isomorphisms. For a while, we write $(kQ_1^c,G, r,
\overrightarrow \chi, u)=(kQ_1^c,G, r, \overrightarrow \chi, u,
\{g_{\theta}\})$ given before. Now let $\{h_{\theta}\in G\mid
\theta\in\Theta_C\}$ be another coset representative set of
$Z_{u(C)}$ in $G$ for any $C\in{\mathcal K}(G)$. That is,
\begin{eqnarray}\label{ncosetde}
G=\bigcup_{\theta\in\Theta_C}Z_{u(C)}h_{\theta}.
\end{eqnarray}
\begin{Lemma}\label{14.1.13}
With the above notations, $(kQ_1^c,G, r, \overrightarrow \chi, u,
\{g_{\theta}\})$ and $(kQ_1^c,G, r, \overrightarrow \chi, u, $
$\{h_{\theta}\})$ are isomorphic $kG$-Hopf bimodules.
\end{Lemma}

{\bf Proof.} We may assume \ \
$Z_{u(C)}h_{\theta}=Z_{u(C)}g_{\theta}$ \ \ for any \ $C\in{\mathcal
K}(G)$ and \ \  $\theta\in\Theta_C$. Then \ \ \
$g_{\theta}h^{-1}_{\theta}\in Z_{u(C)}$. Now let \ $x, y, h\in G$
with $x^{-1}y=g^{-1}_{\theta}u(C)g_{\theta}$. Then $x^{-1}y =$\ $
h^{-1}_{\theta}(g_{\theta}h^{-1}_{\theta})^{-1}u(C)(g_{\theta}h^{-1}_{\theta})h_{\theta}
$ \ $=h^{-1}_{\theta}u(C)h_{\theta}$ and
$h_{\theta}h=(h_{\theta}g^{-1}_{\theta})g_{\theta}h =
(h_{\theta}g^{-1}_{\theta})\zeta_{\theta}(h)g_{\theta'}=
(h_{\theta}g^{-1}_{\theta})\zeta_{\theta}(h)(g_{\theta'}h^{-1}_{\theta'})h_{\theta'}$,
where $g_{\theta}h=\zeta_{\theta}(h)g_{\theta'}$. Hence from
Proposition \ref{14.1.10} we know that the right $kG$-action on
$(kQ_1^c,G, r, \overrightarrow \chi, u, \{h_{\theta}\})$ is given by
$$\begin{array}{rcl}
a^{(i)}_{y,x}\cdot h& =&
\chi^{(i)}_C((h_{\theta}g^{-1}_{\theta})\zeta_{\theta}(h)(g_{\theta'}h^{-1}_{\theta'}))
a^{(i)}_{yh,xh}\\
&=&
\chi^{(i)}_C((g_{\theta}h^{-1}_{\theta})^{-1})\chi^{(i)}_C(\zeta_{\theta}(h))
\chi^{(i)}_C(g_{\theta'}h_{\theta'}^{-1}) a^{(i)}_{yh,xh},\ i\in
I_C(r).
\end{array}$$
However, we also have
$$(xh)^{-1}(yh)=h^{-1}(x^{-1}y)h=h^{-1}g^{-1}_{\theta}u(C)g_{\theta}h=g^{-1}_{\theta'}u(C)g_{\theta'}.$$
It follows that the $k$-linear isomorphism $f: kQ_1\rightarrow kQ_1$
given by
$$f(a^{(i)}_{y,x})=\chi^{(i)}_C(g_{\theta}h^{-1}_{\theta})a^{(i)}_{y,x}$$
for any $x,y\in G$ with $x^{-1}y=g^{-1}_{\theta}u(C)g_{\theta}$,
$C\in{\mathcal K}_r(G)$ and $i\in I_C(r)$, is a $kG$-Hopf bimodule
isomorphism from $(kQ_1^c,G, r, \overrightarrow \chi, u,
\{g_{\theta}\})$ to $(kQ_1^c,G, r, \overrightarrow \chi, u,
\{h_{\theta}\})$. \  \ $\Box$

Now we state one of our main results, which classifies the {\rm PM}
(co-)path Hopf algebras, {\rm PM} semi-(co-)path Hopf algebras and
{\rm PM} one-type-path Hopf algebras.

\begin {Theorem} \label {14.4}
  Let $(G, r, \overrightarrow{\chi}, u)$ and
$(G', r', \overrightarrow{\chi'}, u')$ are two {\rm RSC}'s. Then the
following statements are equivalent:

{\rm(i)} ${\rm RSC} (G, r , \overrightarrow{\chi}, u)$ $\cong $
${\rm RSC} (G', r', \overrightarrow{\chi'}, u')$.

{\rm(ii)} There exists a Hopf algebra isomorphism $\phi:
kG\rightarrow kG'$ such that $(kQ_1^c, G, r , \overrightarrow{\chi},
u)  \cong\ _{\phi} ^{\phi ^{-1}}
 ( ( kQ_1' {}^c, G', r', \overrightarrow{\chi'}, u') ){}_{\phi} ^{\phi ^{-1}}$ as $kG$-Hopf bimodules.

{\rm(iii)} $kQ^c(G, r, \overrightarrow \chi, u)\cong k{Q'}^c(G, r,
\overrightarrow \chi, u)$.
 {\rm(iv)}  $kQ^s(G, r, \overrightarrow \chi, u) \cong
k{Q'}^s(G', r', \overrightarrow \chi', u')$.

 {\rm(v)} \
$kG[kQ_1^c,G, r, \overrightarrow \chi, u] $ $\cong kG'[kQ_1'{}^c,G', r', \overrightarrow \chi', u']$.\\
Furthermore, if $Q$ is finite, then the above are equivalent to the
following:

{\rm(vi)} $kQ^a(G, r, \overrightarrow \chi, u) $ $ \cong k{Q'}^a(G',
r', \overrightarrow \chi', u')$. {\rm(vii)} $kQ^{sc}(G, r,
\overrightarrow \chi, u) \cong k{Q'}^{sc}(G', r', \overrightarrow
\chi', u')$. {\rm(viii)} \ $(kG)^*[kQ_1^a, G, r, \overrightarrow
\chi, u]\cong (kG')^*[kQ_1'{}^a, G', r', \overrightarrow \chi',
u']$. Notice that the isomorphisms above are ones of graded Hopf
algebras but (i) (ii).
\end {Theorem}

{\bf Proof.}  By Lemma \ref{14.1.6} and Lemma \ref {14.1.7}, we only
have to prove (i) $\Leftrightarrow$ (ii).

(i) $\Rightarrow$ (ii). Assume that ${\rm RSC} (G, r ,
\overrightarrow{\chi}, u)$ $\cong $ ${\rm RSC} (G', r',
\overrightarrow{\chi'}, u')$. Then  there exist a group isomorphism
$\phi: G\rightarrow G'$,  an element $h_C\in G$ such that
$\phi(h^{-1}_Cu(C)h_C)=u'(\phi(C))$ for any $C\in{\mathcal K}(G)$
and a bijective map $\phi_C: I_C(r)\rightarrow I_{\phi(C)}(r')$ such
that
${\chi'}_{\phi(C)}^{(\phi_C(i))}(\phi(h^{-1}_Chh_C))=\chi_C^{(i)}(h)$
for any $h\in Z_{u(C)}$, $C\in{\mathcal K}_r(G)$ and $i\in I_C(r)$.
Then $\phi(h^{-1}_CZ_{u(C)}h_C)=Z_{u'(\phi(C))}$ and $\phi:
kG\rightarrow kG'$ is a Hopf algebra isomorphism. Now let
$G=\bigcup_{\theta\in\Theta_C}Z_{u(C)}g_{\theta}$ be given as in
(\ref{e0.2}) for any $C\in{\mathcal K}(G)$, and assume that the
$kG$-Hopf bimodule $(kQ_1^c, G, r, \overrightarrow \chi, u)$ is
defined by using these coset decompositions. Then
\begin{eqnarray}\label{coset}
G'=\bigcup_{\theta\in\Theta_C}Z_{u'(\phi(C))}(\phi(h^{-1}_Cg_{\theta}h_C))
\end{eqnarray}
is a coset decomposition of $Z_{u'(\phi(C))}$ in $G'$ for any
$\phi(C)\in{\mathcal K}(G')$. By Lemma \ref{14.1.13}, we may assume
that the structure of the $kG'$-Hopf bimodule $(kQ_1'{}^c, G', r',
\overrightarrow \chi', u')$ is obtained by using these coset
decompositions (\ref{coset}). Define a $k$-linear isomorphism $\psi:
(kQ_1, G, r, \overrightarrow \chi, u) \rightarrow (kQ_1', G', r',
\overrightarrow \chi', u')$ by
$$\psi(a_{y,x}^{(i)})=\chi_C^{(i)}(\zeta_{\theta}(h^{-1}_C))a_{\phi(y),\phi(x)}^{(\phi_C(i))}$$
for any $x, y\in G$ with $x^{-1}y=g^{-1}_{\theta}u(C)g_{\theta}$,
and $i\in I_C(r)$, where $C\in{\mathcal K}_r(G)$ and
$g_{\theta}h^{-1}_C=\zeta_{\theta}(h^{-1}_C)g_{\eta}$ with
$\zeta_{\theta}(h^{-1}_C)\in Z_{u(C)}$ and $\theta, \eta\in
\Theta_C$. It is easy to see that $\psi$ is a $kG$-bicomodule
homomorphism from $(kQ_1^c, G, r, \overrightarrow \chi, u)$ to
$^{\phi^{-1}}_{\phi}(kQ_1'{}^c, G', r', \overrightarrow \chi',
u')^{\phi^{-1}}_{\phi}$. Since $(hx)^{-1}(hy)=x^{-1}y$ for any $x,
y, h\in G$, it follows from Proposition \ref{14.1.10} that $\psi$ is
also a left $kG$-module homomorphism from $(kQ_1^c, G, r,
\overrightarrow \chi, u)$ to $^{\phi^{-1}}_{\phi}(kQ_1'{}^c, G, r,
\overrightarrow \chi, u)^{\phi^{-1}}_{\phi}$.

Now let $x, y, h\in G$ with $x^{-1}y=g^{-1}_{\theta}u(C)g_{\theta}$,
$C\in{\mathcal K}_r(G)$ and $\theta\in\Theta_C$. Assume that
$g_{\theta}h^{-1}_C=\zeta_{\theta}(h^{-1}_C)g_{\eta}$,
$g_{\theta}h=\zeta_{\theta}(h)g_{\theta'}$,
$g_{\eta}(h_Chh^{-1}_C)=\zeta_{\eta}(h_Chh^{-1}_C)g_{\eta'}$ and
$g_{\theta'}h^{-1}_C=\zeta_{\theta'}(h^{-1}_C)g_{\theta''}$ with
$\zeta_{\theta}(h^{-1}_C), \zeta_{\theta}(h),
\zeta_{\eta}(h_Chh^{-1}_C), \zeta_{\theta'}(h^{-1}_C)\in Z_{u(C)}$
and $\eta, \theta', \eta', \theta''\in \Theta_C$. Then we have
$$\begin{array}{rcl}
g_{\theta}hh^{-1}_C&=&\zeta_{\theta}(h)g_{\theta'}h^{-1}_C\ = \zeta_{\theta}(h)\zeta_{\theta'}(h^{-1}_C)g_{\theta''}\\
\end{array}$$
and
$$\begin{array}{rcl}
g_{\theta}hh^{-1}_C&=&(g_{\theta}h^{-1}_C)(h_Chh^{-1}_C)=\zeta_{\theta}(h^{-1}_C)g_{\eta}(h_Chh^{-1}_C)\
=\zeta_{\theta}(h^{-1}_C)\zeta_{\eta}(h_Chh^{-1}_C)g_{\eta'}.
\end{array}$$
It follows that $\theta''=\eta'$ and
\begin{eqnarray}\label{coeff}
\zeta_{\theta}(h)\zeta_{\theta'}(h^{-1}_C)=
\zeta_{\theta}(h^{-1}_C)\zeta_{\eta}(h_Chh^{-1}_C).
\end{eqnarray}
By Proposition \ref{14.1.10}, $a^{(i)}_{y,x}\cdot
h=\chi^{(i)}_C(\zeta_{\theta}(h))a^{(i)}_{yh,xh}$ for any $i\in
I_C(r)$. Moreover, we have
$(xh)^{-1}(yh)=h^{-1}g^{-1}_{\theta}u(C)g_{\theta}h=
g^{-1}_{\theta'}u(C)g_{\theta'}$. This implies that
$$\psi(a^{(i)}_{y,x}\cdot
h)=\chi^{(i)}_C(\zeta_{\theta}(h))\psi(a^{(i)}_{yh,xh})
=\chi^{(i)}_C(\zeta_{\theta}(h))\chi^{(i)}_C(\zeta_{\theta'}(h^{-1}_C))
a^{(\phi_C(i))}_{\phi(yh),\phi(xh)}.$$ On the other hand, we have
$g_{\theta}=g_{\theta}h^{-1}_Ch_C=\zeta_{\theta}(h^{-1}_C)g_{\eta}h_C$,
and hence
$$\begin{array}{rcl}
\phi(x)^{-1}\phi(y)&=&\phi(x^{-1}y)=\phi(g^{-1}_{\theta}u(C)g_{\theta})
=\phi(h^{-1}_Cg^{-1}_{\eta}u(C)g_{\eta}h_C)\\
&=& \phi(h^{-1}_Cg^{-1}_{\eta}h_C)\phi(h^{-1}_Cu(C)h_C)\phi(h^{-1}_Cg_{\eta}h_C)\\
&=&\phi(h^{-1}_Cg_{\eta}h_C)^{-1}u'(\phi(C))\phi(h^{-1}_Cg_{\eta}h_C).
\end{array}$$
We also have
$$\begin{array}{rcl}
\phi(h^{-1}_Cg_{\eta}h_C)\phi(h)&=&\phi(h^{-1}_Cg_{\eta}h_Chh^{-1}_Ch_C)\\
&=&\phi(h^{-1}_C\zeta_{\eta}(h_Chh^{-1}_C)g_{\eta'}h_C)\\
&=&\phi(h^{-1}_C\zeta_{\eta}(h_Chh^{-1}_C)h_C)\phi(h^{-1}_Cg_{\eta'}h_C).
\end{array}$$
Thus by Proposition \ref{14.1.10} one gets
$$\begin{array}{rcl}
a^{(\phi_C(i))}_{\phi(y),\phi(x)}\cdot\phi(h)&=&
{\chi'}^{(\phi_C(i))}_{\phi(C)}(\phi(h^{-1}_C\zeta_{\eta}(h_Chh^{-1}_C)h_C))
a^{(\phi_C(i))}_{\phi(yh),\phi(xh)}\\
 &=&
\chi^{(i)}_C(\zeta_{\eta}(h_Chh^{-1}_C))
a^{(\phi_C(i))}_{\phi(yh),\phi(xh)}.\\
\end{array}$$
Now it follows from Eq.(\ref{coeff}) that
$$\begin{array}{rcl}
\psi(a^{(i)}_{y,x})\cdot\phi(h)&=&
\chi_C^{(i)}(\zeta_{\theta}(h^{-1}_C))a_{\phi(y),\phi(x)}^{(\phi_C(i))}
\cdot\phi(h)\\
&=&\chi_C^{(i)}(\zeta_{\theta}(h^{-1}_C))
\chi^{(i)}_C(\zeta_{\eta}(h_Chh^{-1}_C))
a^{(\phi_C(i))}_{\phi(yh),\phi(xh)}\\
&=&\psi(a^{(i)}_{y,x}\cdot h).
\end{array}$$
This shows that $\phi$ is a right $kG$-module homomorphism, and
hence a $kG$-Hopf bimodule isomorphism from $(kQ_1^c, G, r,
\overrightarrow \chi, u)$ to $^{\phi^{-1}}_{\phi}(kQ_1'{}^c, G', r',
\overrightarrow \chi', u')^{\phi^{-1}}_{\phi}$.

(ii) $\Rightarrow$ (i). Assume that there exist a Hopf algebra
isomorphism $\phi: kG\rightarrow kG'$ and a $kG$-Hopf bimodule
isomorphism $\psi: (kQ_1^c, G, r, \overrightarrow \chi,
u)\rightarrow\ ^{\phi^{-1}}_{\phi}(kQ_1'{}^c, G', r',
\overrightarrow \chi', u')^{\phi^{-1}}_{\phi}$. Then $\phi:
G\rightarrow G'$ is a group isomorphism. Let $C\in{\mathcal K}(G)$.
Then $\phi(u(C))$, $u'(\phi(C))\in\phi(C)\in{\mathcal K}(G')$, and
hence
$u'(\phi(C))=\phi(h_C)^{-1}\phi(u(C))\phi(h_C)=\phi(h^{-1}_Cu(C)h_C)$
for some $h_C\in G$. Since $\psi$ is a $kG'$-bicomodule isomorphism
from $^{\phi}(kQ_1^c, G, r, \overrightarrow \chi, u)^{\phi}$ to
$(kQ_1'{}^c, G', r', \overrightarrow \chi', u')$ and
$\phi(h^{-1}_Cu(C)h_C)=u'(\phi(C))$, by restriction one gets a
$k$-linear isomorphism
$$\psi_C:\ ^{h^{-1}_Cu(C)h_C}(kQ_1)^1\rightarrow\ ^{u'(\phi(C))}\! (kQ'_1)^1,\ x\mapsto \psi(x).$$
We also have a $k$-linear isomorphism
$$f_C:\ ^{u(C)}\! (kQ_1)^1\rightarrow\ ^{h^{-1}_Cu(C)h_C}(kQ_1)^1,\ x\mapsto h^{-1}_C\cdot x\cdot h_C.$$
Since $\phi(h^{-1}_Cu(C)h_C)=u'(\phi(C))$ and
$h^{-1}_CZ_{u(C)}h_C=Z_{h^{-1}_Cu(C)h_C}$, one gets
$\phi(h^{-1}_CZ_{u(C)}h_C)=Z_{u'(\phi(C))}$. Hence $\phi$ and $h_C$
induce an algebra isomorphism
$$\sigma_C:\ kZ_{u(C)}\rightarrow kZ_{u'(\phi(C))},\ h\mapsto \phi(h^{-1}_Chh_C).$$
Using the hypothesis that $\psi$ is a $kG$-bimodules homomorphism
from  $(kQ_1^c, G, r, \overrightarrow \chi, u)$   to
$_{\phi}(kQ_1'{}^c,$\ $ G', r', \overrightarrow \chi', u')_{\phi}$,
one can easily check that the composition $\psi_Cf_C$ is a right
$kZ_{u(C)}$-module isomorphism from $^{u(C)}\! (kQ_1)^1$ to
$(^{u'(\phi(C))}\! (kQ'_1)^1)_{\sigma_C}$. Since both  $^{u(C)}\!
(kQ_1)^1$   and $(^{u'(\phi(C))}\! (kQ'_1)^1)_{\sigma_C}$ \ \ are
pointed right $kZ_{u(C)}$-modules, they are semisimple
$kZ_{u(C)}$-modules for any $C\in{\mathcal K}_r(G)$. Moreover,
$ka^{(i)}_{u(C),1}$ and $ka^{(j)}_{u'(\phi(C)),1}$ are simple
submodules of $^{u(C)}\! (kQ_1)^1$ and $(^{u'(\phi(C))}\!
(kQ'_1)^1)_{\sigma_C}$, respectively, for any $i\in I_C(r)$ and
$j\in I_{\phi(C)}(r')$, where $C\in{\mathcal K}_r(G)$. Thus for any
$C\in{\mathcal K}_r(G)$, there exists a bijective map $\phi_C:
I_C(r)\rightarrow I_{\phi(C)}(r')$ such that $ka^{(i)}_{u(C),1}$ and
$(ka^{(\phi_C(i))}_{u'(\phi(C)),1})_{\sigma_C}$ are isomorphic right
$kZ_{u(C)}$-modules for any $i\in I_C(r)$, which implies
$\chi_C^{(i)}(h)=\chi'{}_{\phi(C)}^{(\phi_C(i))}(\phi(h^{-1}_Chh_C))$
for any $h\in Z_{u(C)}$ and $i\in I_C(r)$. It follows that ${\rm
RSC} (G, r , \overrightarrow{\chi}, u)$ $\cong $ ${\rm RSC} (G', r',
\overrightarrow{\chi'}, u')$. \ \ $\Box $

Up to now we have classified the {\rm PM} quiver Hopf algebras by
means of {\rm RSC}'s. In other words,  ramification systems with
characters uniquely determine  their corresponding PM quiver Hopf
algebras  up to graded Hopf algebra isomorphisms.

\begin {Example}\label{14.1.16} Assume that  $k$ is a field with char$(k)\not=2$.
Let $G=\{1, g\}\cong{\bf  Z}_2$ be the cyclic group of order $2$
with the generator $g$. Let $r$ be a ramification data of $G$ with
$r_1=m$ and $r_{\{g\}}=0$ and $(Q,G,r)$ be the corresponding Hopf
quiver, where $m$ is a positive integer. Then
$^1(Q_1)^1=\{a^{(i)}_{1,1}\mid i=1, 2, \cdots, m\}$,
$^g(Q_1)^g=\{a^{(i)}_{g,g}\mid i=1, 2, \cdots, m\}$, $^1(Q_1)^g$ and
$^g(Q_1)^1$ are two empty sets. For simplification, we write
$x_i=a^{(i)}_{1,1}$ and $y_i=a^{(i)}_{g,g}$ for any $1\leq i\leq m$.
Clearly, $Z_{u(\{1 \})}=G$ and $\widehat{G}=\{\chi_+, \chi_-\}$,
where $\chi_{\pm}(g)=\pm 1$. For any $0\leq n\leq m$, put
$\overrightarrow{\chi_n} \in (\hat G )^m $ with
$\chi_{n\{1\}}^{(i)}=\left\{\begin{array}{ll}
\chi_-,& \mbox{if } i>n;\\
\chi_+,& \mbox{otherwise}.
\end{array}\right.$
 Then $\{ {\rm RSC} (G, r, \overrightarrow{\chi _n}, u_n) \mid  n=0, 1, 2, \cdots, m \}$ are all non-isomorphic ${\rm RSC}$'s.
 Thus by Theorem \ref{14.4} we know that the path coalgebra
$kQ^c$ exactly admits $m+1$ distinct {\rm PM} co-path Hopf algebra
structures $kQ^c(G, r, \overrightarrow {\chi_n}, u_n)$, $0\leq n\leq
m$, up to graded Hopf algebra isomorphism. Now let $0\leq n\leq m$.
Then by Proposition \ref{14.1.10}, the $kG$-actions on $(kQ_1^c, G,
r, \overrightarrow {\chi_n}, u_n)$ are given by $g\cdot x_i=y_i,\ \
g\cdot y_i=x_i,\ \ 1\leq i\leq m;$$ \ \
$$x_i\cdot g=\left\{\begin{array}{rl}
-y_i,& \mbox{if }i>n,\\
y_i,& \mbox{otherwise},\\
\end{array}\right.\ \ \
 \ \ y_i\cdot g=\left\{\begin{array}{rl}
-x_i,& \mbox{if }i>n,\\
x_i,& \mbox{otherwise}.
\end{array}\right.$
Thus by \cite[p.245 or Theorem 3.8]{CR02}, the products of these
arrows $x_i, y_j$ in $kQ^c(G, r, \overrightarrow {\chi_n}, u_n)$ can
be described as follows. For any $i, j =1, 2, \cdots, m$, $x_i.x_j =
x_ix_j + x_jx_i,\ \ \  y_i.x_j=y_iy_j+y_jy_i,$$ \ \
$$x_i.y_j=\left\{\begin{array}{ccl}
-(y_iy_j + y_jy_i)&,& \mbox{if }i>n,\\
y_iy_j + y_jy_i&,& \mbox{otherwise},\\
\end{array}\right. \ \ \
\ y_i.y_j=\left\{\begin{array}{ccl}
-(x_ix_j + x_jx_i)&,& \mbox{if }i>n,\\
x_ix_j +x_jx_i&,& \mbox{otherwise},\\
\end{array}\right.$
where $x.y$ denotes the product of $x$ and $y$ in $kQ^c(G, r,
\overrightarrow {\chi_n}, u_n)$ for any $x, y\in kQ^c(G, r,
\overrightarrow {\chi_n}, u_n)$, $x_ix_j$ and $y_iy_j$ denote the
$2$-paths in the quiver $Q$ as usual for any $1\leq i, j\leq m$.
\end {Example}

\section{\bf  Multiple Taft algebras}\label{s14.2}

In this section  we discuss the {\rm PM} quiver Hopf algebras
determined by the {\rm RSC}'s with $ \cup {\mathcal K}_r(G)\subseteq
Z(G)$. We give the classification of {\rm PM} Nichols algebras and
multiple Taft algebras by means of  element system with characters
when $G$ is finite abelian group and $k$ is the complex field.

Let $r$ be a ramification data of $G$ and $(Q, G,r)$ be the
corresponding Hopf quiver. If $C$ contains only one element of $G$
for any $C\in {\mathcal K}_r(G)$, then $C=\{g\}$ for some $g\in
Z(G)$, the center of $G$. In this case, we say that the ramification
$r$ is {\it central}, and that ${\rm RSC} (G, r,
\overrightarrow{\chi}, u)$ {\it a central ramification system with
characters}, or a {\rm {\rm CRSC}} in short.  If ${\rm RSC} (G, r,
\overrightarrow{\chi}, u)$ is {\rm CRSC}, then the {\rm PM} co-path
Hopf algebra $kQ^c(G, r, \overrightarrow \chi, u)$ is called a {\it
multiple crown algebra} and $kG[kQ_1^c, G, r, \overrightarrow \chi,
u]$ is called a {\it multiple Taft algebra}.

\begin {Definition}\label{14.2.1}

$(G, \overrightarrow{g}, \overrightarrow{\chi}, J)$ is called an
{\it element system with characters} $($simply, {\rm ESC}$)$ if $G$
is a group, $J$ is a set, $\overrightarrow{ g } = \{g_i\} _{ i\in J}
\in Z(G)^J$ and $\overrightarrow{ \chi }= \{\chi_i\}_{ i \in J} \in
\widehat{G}^J $ with $ g_i \in Z(G)$ and $\chi _i \in \widehat G$.
${\rm ESC} (G, \overrightarrow{g}, \overrightarrow{\chi}, J)$ and
${\rm ESC} (G', \overrightarrow{g'}, \overrightarrow{\chi'}, J')$
are said to be isomorphic if there exist a group isomorphism $\phi:
G \rightarrow G'$ and a bijective map $\sigma: J\rightarrow J'$ such
that $\phi(g_i)=g'_{\sigma(i)}$ and $\chi'_{\sigma(i)}\phi=\chi_i$
for any $i \in J$.
\end {Definition}

$ {\rm ESC}(G,\overrightarrow{g},\overrightarrow{\chi},J)$ can be
written as ${\rm ESC}(G,g_i,\chi_i;i\in J)$ for convenience.
Throughout this chapter, let $q_{ji}:=\chi_i(g_j)$, $q_i=q_{ii}$ and
$N_i$ be the order of $q_i$ ($N_i=\infty$ when $q_i$ is not a root
of unit, or $q_i=1$) for $i,j\in J$.

Let $H$ be a Hopf algebra with a bijective antipode $S$. A {\rm YD}
$H$-module $V$ is {\it pointed} if $V=0$ or $V$ is a direct sum of
one dimensional {\rm YD} $H$-modules. If $V$ is a pointed {\rm YD}
$H$-module, then the corresponding Nichols algebra ${\mathcal B}(V)$
is called a {\rm PM} {\it Nichols algebra}.

\begin{Lemma}\label{14.2.2}
Let $(V, \alpha ^-, \delta ^- )$ be a {\rm YD} $kG$-module. Then
$(V, \alpha ^-, \delta^- )$ is  a pointed  {\rm YD} $kG$-module  if
and only if $(V, \alpha ^-)$ is a pointed $($left$)$ $kG$-module
with $\delta^-(V)\subseteq kZ(G)\otimes V$.
\end{Lemma}
{\bf Proof.} We may assume $V\not=0$. Assume that $V$ is a pointed
YD $kG$-module. Then $V$ is obviously a pointed $kG$-module. Now let
$U$ be a one dimensional YD $kG$-submodule of $V$ and let $0\not=
x\in U$. Then $\delta^-(x)=g\otimes x$ for some $g\in G$. Let $h\in
G$. Then $h\cdot x=\beta x$ for some $0\not=\beta\in k$. By
Eq.(\ref{ydm}), one gets that $\beta g\otimes x=\delta^-(\beta
x)=\delta^-(h\cdot x)=hgh^{-1}\otimes h\cdot x=\beta hgh^{-1}\otimes
x$. This implies that $hgh^{-1}=g$ for any $h\in G$, and hence $g\in
Z(G)$. It follows that $\delta^-(V)\subseteq kZ(G)\otimes V$.
Conversely, assume that $V$ is a pointed $kG$-module and
$\delta^-(V)\subseteq kZ(G)\otimes V$. For any $g\in G$, let
$$V_g=\{v\in V\mid \delta^-(v)=g\otimes v\}.$$
Then $V_g$ is a $kG$-subcomodule of $V$ and any subspace of $V_g$ is
a $kG$-subcomodule for all $g\in G$. By $\delta^-(V)\subseteq
kZ(G)\otimes V$, we know that $V=\bigoplus_{g\in Z(G)}V_g$. Now let
$g\in Z(G)$. Then it follows from Eq.(\ref{ydm}) that $V_g$ is a
$kG$-submodule of $V$. Since $V$ is a pointed $kG$-module, any
submodule of $V$ is pointed, and hence $V_g$ is a pointed
$kG$-module. Thus if $V_g\not=0$ then $V_g$ is a direct sum of some
one dimensional $kG$-submodules of $V_g$. However, any one
dimensional $kG$-submodule of $V_g$ is also a $kG$-subcomodule, and
hence a YD $kG$-submodule. Thus $V$ is a pointed YD $kG$-module.\  \
$\Box$

\begin{Lemma}\label{14.2.3} Assume that $G$ is a finite abelian group
of exponent $m$. If $k$ contains a primitive $m$-th root of 1, then
every {\rm YD} $kG$-module is  pointed  and  Nichols algebra  of
every {\rm YD} $kG$ -module is {\rm PM}.
\end{Lemma}
{\bf Proof.}  It follows from  Lemma \ref {14.1.2} and Lemma \ref
{14.2.2}. $\Box$

Let $(G, g_i, \chi_i; i\in J)$ be an ${\rm ESC}$. Let $V$ be a
$k$-vector space with ${\rm dim}(V)=|J|$. Let $\{x_i \mid i\in J\}$
be a basis of $V$ over $k$. Define a left $kG$-action and a left
$kG$-coaction on $V$ by
$$g\cdot x_i = \chi_i(g)x_i,\ \delta^-(x_i )= g_i \otimes
x_i,\ i\in J,\ g\in G.$$ Then it is easy to see that $V$ is a
pointed {\rm YD} $kG$-module and $kx_i$ is a one dimensional {\rm
YD} $kG$-submodule of $V$ for any $i\in J$. Denote by $V(G, g_i,
\chi_i; i\in J)$ the pointed {\rm YD} $kG$-module $V$. Note that
$V(G, g_i, \chi_i; i\in J)=0$ if $J$ is empty.

\begin{Proposition}\label{14.2.4}
$V$ is pointed {\rm YD} $kG$-module if and only if $V$ is isomorphic
to $V(G, g_i, \chi_i;$ $i\in J)$ for some  {\rm ESC} $(G, g_i,
\chi_i; i\in J)$.
\end{Proposition}

{\bf Proof.} If $V\cong V(G, g_i,\chi_i;i\in J)$ for some {\rm ESC}
$(G, g_i,\chi_i;i\in J)$ of $G$, then $V$ is obviously a pointed
{\rm YD} $kG$-module. Conversely, assume  that $V$ is a nonzero
pointed {\rm YD} $kG$-module. By Lemma \ref{14.2.2},
$V=\bigoplus_{g\in Z(G)}V_g$ and $V_g=\{v\in V|\delta^-(v)= g\otimes
v\}$ is a pointed {\rm YD} $kG$-submodule of $V$ for any $g\in
Z(G)$. Let $g\in Z(G)$ with $V_g\not=0$. Then $V_g$ is a nonzero
pointed $kG$-module. Hence there is a $k$-basis $\{x_i\mid i\in
J_g\}$ such that $kx_i$ is a $kG$-submodule of $V_g$ for any $i\in
J_g$. It follows that there is a character $\chi_i\in\widehat G$ for
any $i\in J_g$ such that $h\cdot x_i=\chi_i(h)x_i$ for all $h\in G$.
For any $i\in J_g$, put $g_i=g$. We may assume that these index sets
$J_g$ are disjoint, that is, $J_g\cap J_h=\emptyset$ for any $g\not=
h$ in $Z(G)$ with $V_g\not=0$ and $V_h\not=0$. Now let $J$ be the
union of all the $J_g$ with $g\in Z(G)$ and $V_g\not=0$. Then one
can see that $(G, g_i,\chi_i;i\in J)$ is an {\rm ESC} of $G$, and
that $V$ is isomorphic to $V(G, g_i,\chi_i;i\in J)$ as a {\rm YD}
$kG$-module.\ \ $\Box$

Now we give the relation between ${\rm RSC}$ and ${\rm ESC}$. Assume
that $(G, g_i, \chi_i; i\in J)$ is an {\rm ESC} of $G$. We define a
binary relation $\sim$ on $J$ by
$$i \sim j \Leftrightarrow g_i = g_j,$$
where $i, j \in J$. Clearly, this is an equivalence relation. Denote
by $J/\!\!\sim$ the quotient set of $J$ modulo $\sim$. For any $i\in
J$, let $[i]$ denote the equivalence class containing $i$. That is,
$[i]:=\{j\in J \mid j\sim i \}$. Choose a subset $\bar J \subseteq
J$ such that the assignment $i\mapsto [i]$ is a bijective map from
$\bar J$ to $J/\!\!\sim$. That is, $\bar J$ is a set of
representative elements of the equivalence class.  Then
$J=\bigcup_{i\in \bar J}[i]$ is a disjoint union. Let $r_{ \{g_i\}}
=|[i]|$ for any $i\in\bar J$. Then $r=\sum_{i\in \bar J} r_{\{g_i\}}
\{g_i\}$ is a central ramification of $G$ with $I_{\{g_i\}}(r)=[i]$
for $i\in \bar J$. Moreover, ${\mathcal K}_r(G)=\{ \{ g_i \}\mid
i\in\bar J\}$. Put $\chi_{\{g_i\}}^{(j)}:=\chi_j$ for any $i\in\bar
J$ and $j\in [i]$. We obtain an ${\rm  CRSC}$, written ${\rm
CRSC}(G, r(g_i, \chi_i; i\in J),$ \ $ \overrightarrow \chi (g_i,
\chi_i; i\in J), u) $.
 Let $(Q, G, r)$ be the corresponding Hopf quiver
with $r=r(g_i,\chi_i;i\in J)$ and denote by $(kQ_1^c, g_i, \chi_i;
i\in J)$ the $kG$-Hopf bimodule $(kQ_1^c, G, r, \overrightarrow
\chi, u)$. Denote by $kQ^c(G, g_i, \chi_i; i\in J)$ and $kG[kQ_1^c,
G, g_i, \chi_i; i\in J]$ the corresponding multiple crown algebra
$kQ^c(G, r, \overrightarrow \chi, u)$ and multiple Taft algebra
$kG[kQ_1^c, G, r, \overrightarrow \chi, u]$, respectively. We also
denote by $kQ^a(G, g_i,\chi_i;i\in J)$, $kQ^{sc}(G, g_i,\chi_i;i\in
J)$, $kQ^s(G, g_i,\chi_i;i\in J)$ and    $(kG)^*[kQ_1^a, G,  g_i,
\chi_i; i\in J]$ the corresponding path Hopf algebra $kQ^a(G, r,
\overrightarrow \chi, u)$, semi-co-path Hopf algebra $kQ^{sc}(G, r,
\overrightarrow \chi, u)$, semi-path Hopf algebra $kQ^s(G, r,
\overrightarrow \chi, u)$ and one-type- path Hopf algebra
$(kG)^*[kQ_1^a, G, r, \overrightarrow \chi, u]$, respectively.

Conversely, assume that  $(G, r, \overrightarrow \chi, u) $ is a
${\rm CRSC}$. We may assume $I_{\{g\}}(r)\cap
I_{\{h\}}(r)=\emptyset$ for any $\{g\}\not=\{ h\}$ in ${\mathcal
K}_r(G)$. Let $J:=\bigcup_{\{g\}\in{\mathcal K}_r(G)}I_{\{g\}}(r)$.
For any $i\in J$, put $g_i:=g$ and $\chi_i:=\chi_{\{g\}}^{(i)}$ if
$i\in I_{\{g\}}(r)$ with $\{g\}\in{\mathcal K}_r(G)$. We obtain an
${\rm ESC}$, written ${\rm ESC} (G, \overrightarrow{g} (r,
\overrightarrow \chi, u),$ \ $ \overrightarrow{\chi } (r,
\overrightarrow \chi, u), J).$

From now on, assume $I_{\{g\}}(r)\cap I_{\{h\}}(r)= \emptyset $ for
any $\{g\}, \{h\}\in {\cal K}_r(G)$ with $\{g\} \not=\{ h\}$.  Note
that, in two cases above, for any $i, j \in I_{\{g\}}(r)$, we have
\begin {eqnarray} \label {e2.41} g_i = g_j =g, \ \  a^{(j)} _{g_i,
1} = a^{(j)}_{g_j,1}, \ \  \chi ^{(j)} _{\{g_i\}} = \chi
^{(j)}_{\{g_j\}} = \chi _j.
\end {eqnarray} Throughout this chapter, let $E_j := a_{g_j,
1}^{(j)}$ for any $j\in J$.
\begin {Proposition} \label {14.2.5}    ${\rm CRSC}(G, r,
\overrightarrow \chi, u) $ \ $ \cong
 {\rm CRSC }(G', r',
\overrightarrow {\chi'}, u')$ \   if and only if \ ${\rm ESC} (G, $
\ $ \overrightarrow{g} (r, \overrightarrow \chi, u),$ \ $
\overrightarrow{\chi } (r, \overrightarrow \chi, u), J) $ \ $ \cong
$ \ $ {\rm ESC} (G', \overrightarrow{g'}(r', \overrightarrow {
\chi'}, $ \ $ u'),\overrightarrow{\chi ' } (r', \overrightarrow
{\chi'}, u'), J')$.
\end {Proposition}

{\bf Proof.} We use  notations above. If ${\rm CRSC}(G, r,
\overrightarrow \chi, u) \cong {\rm  CRSC} (G', r', \overrightarrow
\chi', u')$, then there exists a group isomorphism $\phi:
G\rightarrow G'$; for any $C\in{\mathcal K}_r(G)$, there exists a
bijective map $\phi_C : I_C(r) \rightarrow I_{\phi(C)}(r')$ such
that $\chi '{}_{\phi(C)}^{(\phi_C(i))}\phi =\chi_C^{(i)}$.
 Let $\sigma$ be the bijection from $J$ to $J'$ such
 that $\sigma (i) = \phi _C(i)$  for any $C \in {\cal K}_r(G)$, $i \in
 I_C(r)$. It is clear that
  $\phi(g_j)=g'_{\sigma(j)}$ and $\chi'_{\sigma(j)}\phi =\chi_j$ for any $j \in J.$  Thus ${\rm ESC} (G,
\overrightarrow{g} (r, \overrightarrow \chi, u),
\overrightarrow{\chi } (r, \overrightarrow \chi, u), J) $ \ $ \cong
$ \ $ {\rm ESC} (G', \overrightarrow{g'}(r', \overrightarrow {
\chi'}, u'), $ \ $\overrightarrow{\chi ' } (r', \overrightarrow
{\chi'}, u'), J')$.

Conversely, if ${\rm ESC} (G,$ \ $ \overrightarrow{g} (r,
\overrightarrow \chi, u),$ \ $ \overrightarrow{\chi } (r,
\overrightarrow \chi, u), J) $ \ \  $ \cong $ \ \ $ {\rm ESC} (G', $
$ \overrightarrow{g'}(r', $ $ \overrightarrow { \chi'}, u'), $ \ \ \
$\overrightarrow{\chi ' } (r', \overrightarrow {\chi'}, u'), J')$,
then there exist a group isomorphism $\phi: G \rightarrow G'$,  a
bijective map $\sigma: J\rightarrow J'$ such that
$\phi(g_j)=g'_{\sigma(j)}$ and $\chi'_{\sigma(j)}\phi =\chi_j $ for
any $j \in J.$ For any $C = \{g_i\} \in {\cal K}_r(G), j \in
I_C(r)$, we define $\phi _C = \sigma \mid _{I_C(r)}$  and  have
\begin {eqnarray*}\chi '{}_{\phi (C)}^{(\phi _C (j))} \phi &=& \chi '{}_{g'_{\sigma (i)}}^{(\sigma (j))}
\phi =\chi '{}_{\sigma (j)} \phi = \chi_ j = \chi ^{(j)}_C.\ \ \Box
\end {eqnarray*}

If  $ (G, g_i,\chi_i;i\in J)$ is  an {\rm ESC}, then $(kQ_1^c, G,
g_i, \chi_i; i\in J)$ is a $kG$-Hopf bimodule with module operations
$\alpha^-$ and $\alpha^+$. Define a new left $kG$-action on $kQ_1$
by
$$g\rhd x:=g\cdot x\cdot g^{-1},\ g\in G, x\in kQ_1,$$
where $g\cdot x=\alpha^-(g\otimes x)$ and $x\cdot
g=\alpha^+(x\otimes g)$ for any $g\in G$ and $x\in kQ_1$. With this
left $kG$-action and the original left (arrow) $kG$-coaction
$\delta^-$, $kQ_1$ is a YD $kG$-module. Let  $Q_1^1:=\{a\in Q_1 \mid
s(a)=1\}$. It is clear that $kQ_1^1$ is a YD $kG$-submodule of
$kQ_1$, denoted  by $(kQ_1^1, ad(G, g_i, \chi_i; i\in J))$.

\begin{Lemma}\label{14.2.6}
 $(kQ_1^1, ad(G, g_i,\chi_i;i\in J))$ and
$V(G, g_i,\chi^{-1}_i;i\in J)$   are isomorphic {\rm YD}
$kG$-modules.
\end{Lemma}

{\bf Proof.}   By definition, $V(G, g_i,\chi^{-1}_i;i\in J)$ has a
$k$-basis $\{x_i\mid i\in J\}$ such that $\delta^-(x_i)=g_i\otimes
x_i$ and $g\cdot x_i=\chi^{-1}_i(g)x_i$ for all $i\in J$ and $g\in
G$.  By Proposition \ref{14.1.10}, for any $j \in J, $ we have that
$g\rhd a_{g_j,1}^{(j)}=\chi_j(g^{-1})a_{g_j,1}^{(j)}$  and
$\delta^-(a_{g_j,1}^{(j)})=g_j\otimes a_{g_j,1}^{(j)}$. It follows
that there is a {\rm YD} $kG$-module isomorphism from $(kQ_1^1,
ad(G, g_i,\chi_i;i\in J))$ to $V(G, g_i,\chi^{-1}_i;i\in J)$ given
by $a_{g_j,1}^{(j)}\mapsto x_j$ for any  $j\in J$. $\Box$

%Consequently,  $(kQ_1^1, ad(g_i,\chi_i;i\in J))$ is a pointed YD
%$kG$-module and ${\cal B}(kQ_1^1, ad(g_i,\chi_i;i\in J))$ is a PM
%Nichols algebra, and a PM  Nichols algebra respectively.

\begin{Lemma}\label{14.2.7}
Let $B$ and $B'$ be two Hopf algebras with bijective antipodes. Let
$V$ be a {\rm YD} $B$-module. Assume that there is a Hopf algebra
isomorphism $\phi: B'\rightarrow B$. Then
$^{\phi^{-1}}_{\phi}{\mathcal B}(V)$ $\cong $ ${\mathcal
B}({}^{\phi^{-1}}_{\phi}V)$ as graded braided Hopf algebras in
$^{B'}_{B'} {\cal YD}.$
\end{Lemma}
{\bf Proof.} We first  show that $^{\phi^{-1}}_{\phi}{\mathcal
B}(V)$ is a Nichols algebra of $^{\phi^{-1}}_{\phi}V$ by following
steps. Let $R:= ^{\phi^{-1}}_{\phi}{\mathcal B}(V) $. (i) Obviously,
$R_0 =k,$  $R_1 = P(R)$ and  $R$ is generated by $R_1$ as algebras.
Since $C(x\otimes y ) = C'(x'\otimes y')$ for any $x' , y' \in R,$ $
\ x, y \in {\cal B}(V)$ with $x'=x, y'=y$, where $C$ and $C'$ denote
the braidings in $^B_B {\cal YD}$ and $^{B'}_{B'} {\cal YD}$,
respectively, we have that $R$ is a graded braided Hopf algebra in
$^{B'}_{B'} {\cal YD}$. (ii) $R_1 =  ^{\phi ^{-1}} _\phi  {\cal B}
(V)_1 $ $ \cong ^{\phi ^{-1}} _\phi V $ as YD $B'$-modules, since
${\cal B} (V)_1 \cong V $ as YD $B$-modules.

Consequently, $^{\phi^{-1}}_{\phi}{\mathcal B}(V)$ is a Nichols
algebra of $^{\phi^{-1}}_{\phi}V$.  By \cite [Proposition 2.2 (iv)
]{AS02}, $R \cong {\cal B} ( ^{\phi^{-1}}_{\phi}V)$ as graded
braided Hopf algebra in $^{B'}_{B'}{\cal YD}$. $\Box$

\begin{Theorem}\label{14.5} Assume that
 $(G, g_i, \chi_i; i\in J)$ and $ (G', g_i', \chi_i'; i\in J')$ are
 two {\rm ESC}'s.
Then the following statements are equivalent:

{\rm(i)}    ${\rm ESC} (G, g_i, \chi_i; i\in J) \cong \ {\rm ESC}
(G', g_i', \chi_i'; i\in J')$.

 {\rm(ii)} \ ${ \rm CRSC}(G, r(G, g_i, \chi_i; i\in J),$ \ $ \overrightarrow \chi (G, g_i, \chi_i; i\in J), u) $  \ $ \cong $\ $
{\rm  CRSC }(G', r'(G', g_i', \chi_i'; i\in J'),$ \ $
\overrightarrow {\chi'}(G', g_i', \chi_i'; i\in J'), $\ $ u')$.

 {\rm(iii)} There is a Hopf
algebra isomorphism $\phi: kG\rightarrow kG'$ such that $V(G, g_i,
\chi_i;$ $i\in J)\cong\ ^{\phi^{-1}}_{\phi}V'(G' g_i', \chi_i';$
$i\in J')$ as {\rm YD} $kG$-modules.

 {\rm(iv)} There is
a Hopf algebra isomorphism $\phi: kG\rightarrow kG'$ such that
${\mathcal B}(V(G,g_i, \chi_i;$ $i\in J))\cong\
^{\phi^{-1}}_{\phi}{\mathcal B}(V'(G', g_i', \chi_i';$ $i\in J'))$
as graded braided Hopf algebras in $^{kG}_{kG}{\mathcal YD}$.

 {\rm(v)} There is a Hopf
algebra isomorphism $\phi: kG\rightarrow kG'$ such that $(kQ_1^1,
ad(G, g_i,\chi_i;i\in J))\cong\ ^{\phi^{-1}}_{\phi}(kQ'{}_1^1,
ad(G', g'_i,\chi'_i;i\in J'))$ as {\rm YD} $kG$-modules.

{\rm(vi)}    $kG[ kQ^c_1, G, g_i, \chi_i; i\in J] \cong \ kG'[
kQ'{}^c_1,G', g_i', \chi_i'; i\in J'].$

\end{Theorem}
{\bf Proof. }  We use he notations before Proposition \ref {14.2.5}.

(i) $\Rightarrow$ (ii). There exist a group isomorphism $\phi: G
\rightarrow G'$ and   a bijective map $\sigma: J\rightarrow J'$ such
that $\phi(g_j)=g'_{\sigma(j)}$ and $\chi'_{\sigma(j)}\phi =\chi_j $
for any $j \in J$. For any $C=\{g_i\} \in {\cal K}_r (G)$ and $j \in
I_{C}(r)$, we have $g_i = g_j$ and
\begin {eqnarray*} \chi_{ C}^{(j)} = \chi_{\{
g_i\}}^{(j)}&=& \chi_j= \chi'_{\sigma(j)}\phi\\
&=& {\chi'}_{\{g'_{\sigma (i)}\}}^{(\sigma (j))}\phi =
{\chi'}_{\{\phi(g _{i})\}}^{(\phi_{g_i} (j))}\phi =
{\chi'}_{\phi(C)}^{(\phi_{C} (j))}\phi \ .
\end {eqnarray*}

(ii) $\Rightarrow$ (i). There is  a group isomorphism $\phi: G
\rightarrow G'$ and   a bijection $\phi _C: \ I_C(r)\rightarrow I
_{\phi (C)} (r')$ such that $\chi '{}_{\phi (C)}^{\phi _C(j)} \phi =
\chi _C ^{(j)}$ for any $C = \{ g_i \} \in {\cal K}_r(G)$, $j \in
I_C(r)$. Define a map $\sigma : \ J \rightarrow J'$  such that
$\sigma \mid _{I_C(r)} = \phi _C$ for any $C \in {\cal K}_r(G)$.
Thus $\sigma $ is  bijective. For any $C= \{g_i\} \in {\cal K}_r(G)$
and $j \in I_{\{g_i\}}(r)$,  we have
$\phi(g_j)=\phi(g_i)=g'_{\sigma(j)}$  and
${\chi'}_{\sigma(j)}\phi={\chi'}_{\phi_{\{g_i\}}(j)}\phi$ $ =
{\chi'}_{\phi(\{g_i\})}^{(\phi _{\{g_i\}}(j))}\phi
=\chi_{\{g_{i}\}}^{(j)}=\chi_j$. This shows that ${\rm ESC} (G,
\overrightarrow{g}, $\ $\overrightarrow{\chi }, J)$ \ $\cong $ ${\rm
ESC} (G', \overrightarrow{g'}, $\ $\overrightarrow{\chi '}, J')$.

(i) $\Leftrightarrow$ (iii). Let  $V:= V(G, g_i, \chi_i;$ $i\in J)$
and $ V':= V'(G' g_i', \chi_i';$ $i\in J')$. By definition $V$ has a
$k$-basis $\{x_i\mid i\in J\}$ such that $g\cdot x_i=\chi_i(g)x_i$
and $\delta^-(x_i)=g_i\otimes x_i$ for any $i\in J$ and $g\in G$.
Similarly, $V'$ has a $k$-basis $\{y_j\mid j\in J'\}$ such that
$h\cdot y_j=\chi'_j(h)y_j$ and $\delta^-(y_j)=g'_j\otimes y_j$ for
all $j\in J'$ and $h\in G'$.

Assume ${\rm ESC} (G, \overrightarrow{g}, $\ $\overrightarrow{\chi
}, J)$ \ $\cong $ ${\rm ESC} (G', \overrightarrow{g'}, $\
$\overrightarrow{\chi '}, J')$. Then there is a group isomorphism
$\phi: G \rightarrow G'$ and a bijective map $\sigma : J \rightarrow
J'$ such that $\phi(g_i) = g'_{\sigma (i)}$ and
$\chi'_{\sigma(i)}\phi=\chi_i$ for any $i \in J$. Hence $\phi:
kG\rightarrow kG'$ is a Hopf algebra isomorphism. Define a
$k$-linear isomorphism $\psi: V\rightarrow V'$ by
$\psi(x_i)=y_{\sigma(i)}$ for any $i\in J$. Then it is
straightforward to check that $\psi$ is a {\rm YD} $kG$-module
homomorphism from $V$ to $^{\phi^{-1}}_{\phi}V'$.

Conversely, assume that $\phi: kG\rightarrow kG'$ is a Hopf algebra
isomorphism and $\psi: V\rightarrow\ ^{\phi^{-1}}_{\phi}V'$ is a
{\rm  YD} $kG$-module isomorphism. Then $\phi: G\rightarrow G'$ is a
group isomorphism. We use the notations in the proof of Lemma
\ref{14.2.2} and the notations above. Then $\psi(V_g)=V'_{\phi(g)}$
for any $g\in G$. Since $V=\bigoplus_{i\in\bar J}V_{g_i}$ and
$V'=\bigoplus_{i\in\bar{J'}}V_{g'_i}$, there is a bijection $\tau:
\bar J \rightarrow \bar J'$ such that
$\psi(V_{g_i})=V'_{g'_{\tau(i)}}$ for any $i\in\bar J$. This shows
that $\phi(g_i)=g'_{\tau(i)}$ and $V_{g_i}\cong\
_{\phi}(V'_{g'_{\tau(i)}})$ as left $kG$-modules for any $i\in \bar
J$. However, $V_{g_i}$ is a pointed $kG$-module and $kx_j$ is its
one dimensional submodule for any $j\in[i]$. Similarly,
$V_{g'_{\tau(i)}}$ is a pointed $kG'$-module and $ky_j$ is its one
dimensional submodule for any $j\in[\tau(i)]$. Hence there is a
bijection $\phi_{\{g_i\}}: [i]\rightarrow [\tau(i)]$ for any $i\in
\bar J$ such that $kx_j$ and $_{\phi}(ky_{\phi_{\{g_i\}}(j)})$ are
isomorphic $kG$-modules for all $j\in[i]$. This implies that
$\chi_j=\chi'_{\phi_{\{g_i\}}(j)}\phi$ for all $j\in[i]$ and
$i\in\bar J$. Then the same argument as in the proof of (ii)
$\Rightarrow$ (i) shows that ${\rm ESC} (G, \overrightarrow{g}, $\
$\overrightarrow{\chi }, J)$ \ $\cong $ ${\rm ESC} (G',
\overrightarrow{g'}, $\ $\overrightarrow{\chi '}, J')$. (iii)
$\Leftrightarrow$ (iv)  It follows from Lemma \ref{14.2.7}. \ (iii)
$\Leftrightarrow$ (v) It follows from Lemma \ref {14.2.6}. \ (ii)
$\Leftrightarrow$ (vi) It follows from Theorem \ref {14.4}. $\Box$

Up to now we have classified  Nichols algebras and {\rm YD} modules
over finite abelian group and  the complex field up to isomorphisms,
which are under  means of Theorem  \ref {14.5} (iv)(iii),
respectively. In fact, we can explain these facts above by
introducing some new concepts about isomorphisms. For convenience,
if $B$ is a Hopf algebra and $M$ is a $B$-Hopf bimodule, then  we
say that $(B, M)$ is a  Hopf bimodules. For any two Hopf bimodules
$(B,M)$ and $(B', M')$, if $\phi$ is a Hopf algebra homomorphism
from $B$ to $B'$
 and $\psi$ is simultaneously a $B$-bimodule homomorphism from
$M$ to $_\phi M'{}_\phi $ and a $B'$-bicomodule homomorphism from
$^\phi M ^\phi $ to $M'$, then $(\phi, \psi)$ is called a pull-push
Hopf bimodule homomorphism. Similarly, we say that $(B, M)$ and $(B,
X)$ are a {\rm YD} module and a {\rm YD} Hopf algebra if $M$ is a
{\rm YD} $B$-module and $X$ is a braided Hopf algebra in
Yetter-Drinfeld category $^B_B {\cal YD}$, respectively.
 For any two
{\rm YD} modules $(B,M)$ and $(B', M')$, if $\phi$ is a Hopf algebra
homomorphism from $B$ to $B'$,  and $\psi$ is simultaneously a left
$B$-module homomorphism from  $M$ to $_\phi M' $ and a left
$B'$-comodule homomorphism from $^\phi M  $ to $M'$, then $(\phi,
\psi)$ is called a pull-push {\rm YD} module homomorphism. For any
two {\rm YD} Hopf algebra $(B,X)$ and $(B', X')$, if $\phi$ is a
Hopf algebra homomorphism from $B$ to $B'$,   $\psi$ is
simultaneously a left $B$-module homomorphism from  $X$ to $_\phi X'
$ and a left $B'$-comodule homomorphism from $^\phi X  $ to $X'$,
meantime, $\psi$ also is algebra and coalgebra homomorphism from $X$
to $X'$, then $(\phi, \psi)$ is called a pull-push {\rm YD} Hopf
algebra homomorphism.

Consequently, we have classified  Nichols algebras  over finite
abelian group and the complex field up to  pull-push graded {\rm YD}
Hopf algebra isomorphisms and {\rm YD} modules  over finite abelian
group and  the complex field up to pull-push  {\rm YD} module
isomorphisms, respectively. In other words,  element systems with
characters uniquely determine  their corresponding Nichols algebras
and {\rm YD} modules up to their isomorphisms.

\section    {The relation between quiver Hopf algebras and quotients  of free algebras
}\label {s14.3}

In this section we show that the diagram of a quantum weakly
commutative multiple Taft algebra is   not only a Nichols algebra
but also a quantum linear space in $^{kG}_{kG}{\cal YD}$; the
diagram of a semi-path Hopf algebra of ${\rm ESC}$   is  a quantum
tensor algebra in $^{kG}_{kG}{\cal YD}$; the quantum enveloping
algebra of a complex semisimple Lie algebra is a quotient of a
semi-path Hopf algebra.

\subsection {The structure of multiple Taft algebras and semi-path Hopf algebras}\label {s3.1}
Assume that $H=\bigoplus_{i\geq 0}H_{(i)}$ is a graded Hopf algebra
with invertible antipode $S$. Let $B=H_{(0)}$, and let $\pi_0:
H\rightarrow H_{(0)}=B$ and $\iota_0: B=H_{(0)}\rightarrow H$ denote
the canonical projection and injection. Set
$\omega:=id_H*(\iota_0\pi_0S): H\rightarrow H$. Then it is clear
that $(H, \delta ^+, \alpha ^+)$ is a right $B$-Hopf module with
$\delta^+:=({\rm id}\otimes\pi_0)\Delta $ and $\alpha^+:=\mu({\rm
id}\otimes\iota_0)$.   Let $R:=H^{co B}:=\{h\in
H\mid\delta^+(h)=h\otimes 1\}$, which is a graded subspace of $H$.
Then it is known that $R={\rm Im}(\omega)$ and $\Delta(R)\subseteq
H\otimes R$. Hence $R$ is a left coideal subalgebra of $H$, and so
$R$ is a left $H$-comodule algebra. It is well known that $R$ is a
graded braided Hopf algebra in $^B_B{\mathcal YD}$ with the same
multiplication, unit and counit as in $H$,  the comultiplication
$\Delta_R=(\omega \otimes id )\Delta$ , where the left $B$-action
$\alpha _R$ and left $B$-coaction $\delta _R$ on $R$ are given by
\begin {eqnarray}\label {e3.11}
\alpha _R (b \otimes x) = b\rightharpoonup_{ad}x=\sum
b_{(1)}xS(b_{(2)}),\ \delta_R^-(x)=\sum \pi_0(x_{(1)})\otimes
x_{(2)},\ b\in B,\ x\in R
\end {eqnarray}( see the proof of \cite [Theorem 3]{Ra85}). $R$ is
called the diagram of $H$, written $diag (H)$. Note that diagram $R$
of $H$ is dependent on  the gradation of $H$.  By \cite [Theorem
1]{Ra85}, the biproduct of $R$ and $B$ is a Hopf algebra, written $R
^{\delta _R}_{\alpha _R} \# B$, or $R \# B$ in short. The biproduct
$R ^{\delta _R}_{\alpha _R} \# B$ is also called the bosonization of
$R$. Furthermore, we have the following well known result.

\begin{Theorem}\label{14.6}
$($see \cite[p.1530]{Ni78}, \cite{AS98a} and \cite{Ra85}$)$  Under
notations  above, if $H=\bigoplus_{i \geq 0 }H _{(i)}$ is a graded
Hopf algebra, then $R $ is a graded braided Hopf algebra in
$^B_B{\mathcal YD}$ and $ R\#B\cong H$ as graded Hopf algebras,
where the isomorphism is $\alpha ^+  := \mu _H(id _H \otimes \iota
_0) $.
\end{Theorem}

{ \bf Remark}: If  $A$ be a Hopf algebra whose coradical $A_{0}$ is
a Hopf subalgebra, then it is clear that $H:=gr A$ is a graded Hopf
algebra. The diagram of $H$ with respect to gradation of $gr A$ is
called the diagram  of $A$ in \cite [Introduction ]{AS98b}.

\begin {Lemma}\label {14.3.1} (i) Assume that  $H$ and $H'$ are two graded Hopf
algebras with $B= H_{(0)}$ and $B' = H'_{(0)}.$  Then $H \cong H'$
as graded Hopf algebras if and only if there exists a Hopf algebra
isomorphism $\phi: B\rightarrow B'$ such that $ diag (H)\cong\
_{\phi} ^{\phi ^{-1}}diag (H')$ as {\rm YD} $B$-modules and as
graded  braided Hopf algebras in $^B_B {\cal YD}$.

(ii) Let $B$ and $B'$ be two  Hopf algebras. Let $M$ and $M'$ be
$B$-Hopf bimodule and  $B'$-Hopf bimodule, respectively. Then $B[M]
\cong B'[M']$ as graded Hopf algebras if and only if there exists a
Hopf algebra isomorphism $\phi: B\rightarrow B'$ such that $ diag
(B[M])\cong\ _{\phi} ^{\phi ^{-1}}diag (B'[M'])$ as {\rm YD}
$B$-modules and as graded braided Hopf algebras in $^B_B {\cal YD}$.

\end {Lemma}
{\bf Proof.} (i) Assume that $\xi$ is a graded Hopf algebra
isomorphism from $H$ to $H'$. Let $R := diag (H)$, $R' := diag
(H')$, $\phi := \xi \mid _ B$ and $\psi := \xi \mid _ R$. It is easy
to check that $\psi$ is the map required.

 Conversely, by Theorem
\ref {14.6}, $R \# B\cong H$ and $R'\# B' \cong H'$ as graded Hopf
algebras.  Let $\xi $ be a linear map from $ R \# B$ to $R' \# B'$
by sending $r \# b$ to $\psi (r) \# \phi (b)$ for any $r \in R$,
$b\in B$. Let $\nu $ be a linear map from $ R' \# B'$ to $R \# B$ by
sending $r' \# b'$ to $\psi ^{-1} (r') \# \phi ^{-1} (b')$ for any
$r' \in R'$, $b'\in B'$. Obviously, $\nu$ is the inverse of $\xi.$
Since $\psi$ is graded, so is $\xi$.

 Now we show that $\xi$ is an algebra homomorphism.
For any $r, r'\in R, b, b' \in B$, see
\begin {eqnarray*}
\xi \mu _{R\# B}( (r \# b) \otimes (r' \# b'))&=&  \psi (r (b _{(1)}
\cdot r'))\# \phi (b_{(2)}b')\\
&=&\psi (r) (\phi (b _{(1)})\cdot \psi (r' ) ) \# \phi (b_{(2)})
\phi (b) \ \\
&& (\hbox {since } \psi  \hbox { is a pullback module homomorphism }
\\
&& \hbox {
and  an algebra homomorphism. } ) \\
&=&\mu _{R'\# B'} (\xi (r \# b) \otimes \xi (r' \# b')).
\end {eqnarray*}
Similarly, we can show that $\xi$ is a coalgebra homomorphism.

(ii) It follows from (i). $\Box$

%\begin {Lemma}\label {3.9'} Assume that $A$ is an $H$-module algebra and $B$ is a subalgebra of $A$, as well as an
% $H$-submodule of A. If

%$A \# H = B\# H$, then $A = B$. \end {Lemma}

\begin{Lemma}\label{14.3.2} (i) $(kQ^c (G, r, \overrightarrow \chi, u))^{co \ kG} = span \{ \beta \mid \beta \hbox { is a path with  }
s(\beta ) =1 \}$. (ii) $(kG[kQ_1^c, $ $G, r, \overrightarrow \chi, u
])^{co \ kG}$ is the subalgebra of $kG[kQ^c_1, G, r, \overrightarrow
\chi, u]$ generated by $Q_1^1$ as algebras. (iii)  $(kG[kQ_1^c, $ $
G, r, \overrightarrow \chi, u] )^{co \ kG}$ $ \# kG \cong kG[kQ_1^c,
G, r, \overrightarrow \chi, u ]$ as graded Hopf algebra isomorphism.
(iv) $(kQ^s(G, r, \overrightarrow \chi, u))^{co \ kG}$ is the
subalgebra of $kG^s( G, r, \overrightarrow \chi, u)$ generated by
$Q_1^1$ as algebras.

\end{Lemma}

{\bf Proof.} (i)
 For a path $\beta $, see that
 \begin {eqnarray*}
 \delta ^+ (\beta ) = (id \otimes \pi _0)\Delta (\beta ) &=&  \beta \otimes s(\beta). \end {eqnarray*} This implies
$(kQ^c)^{co \ kG} = span \{ \beta \mid \beta \hbox { is a path with
} s(\beta ) =1 \}$.

(ii)  Since every path generated by arrows  in $Q_1^1$ is of start
vertex 1, this  path belongs to $(kG[kQ_1^c,G, r, \overrightarrow
\chi, u ])^{co \ kG}$. Let $R := (kG[kQ_1^c, G, r, \overrightarrow
\chi, u])^{co \ kG}$ and
 $A: = $ the subalgebra of
$kG[kQ^c_1, G, r, \overrightarrow \chi, u]$ generated by $Q_1^1$ as
algebras. Obviously, $A
 \subseteq
R$. It
 is clear that $ \alpha ^+ (R \# kG)  = \alpha ^+ (A \# kG) = kG[kQ_1^c,G, r, \overrightarrow \chi, u]$ and $\alpha ^+$ is injective. Thus $R \# kG  = A \#
 kG$ and $R = A$.

(iii) It follows from   Theorem \ref{14.6}.

(iv) We first show   $(kQ^s)^{co\ kG} =$  $span \{ \beta \mid \beta
=1,  \hbox {or }  \beta = \beta _n \otimes _{kG} \beta _{n-1}
\otimes _{kG} \cdots \otimes _{kG} \beta _1 $ \ $ \hbox {with }  $ $
\prod _{i =1}^n s( \beta _i ) = 1 \hbox { and } \beta _i \in Q_1, \
i = 1, 2, \cdots, n; n \in {\mathbb Z}^+ \}$. Indeed, obviously
right hand side $ \subseteq $ the left hand side. For any $\beta =
\beta _n \otimes _{kG} \beta _{n-1} \otimes _{kG} \cdots \otimes
_{kG} \beta _1 $\ $ \hbox {with } \beta _i \in Q_1 $, called a
monomial, define $s(\beta ) = \prod _{i =1}^n s( \beta _i )$. For
any $0 \not= u \in (kQ^s)^{co\ kG}$ with $u \not\in kG$, there exist
linearly independent monomials $u_1, u_2, \cdots, u _n$ such that
 $u = \sum
_{i=1}^n b_iu_i$ with $0\not= b_i \in k$  for $i = 1, 2, \cdots n$.
See $\delta ^+ (u)$ $= \sum _{i=1}^n b_iu_i \otimes s(u_i) $ $= u
\otimes 1$. Consequently, $s(u_i) =1$ for $i =1,2, \cdots, n$. This
implies that $u $ belongs to the right hand side.

For any $\beta = \beta _n \otimes _{kG} \beta _{n-1} \otimes _{kG}
\cdots \otimes _{kG} \beta _1 \hbox {with }   \prod _{i =1}^n s(
\beta _i ) = 1 \hbox { and } \beta _i \in Q_1, \ i = 1, 2, \cdots,
n, $ we show that $\beta $ can be written as multiplication of
arrows in $Q_1^1$ by induction. When $n=1$, it is clear. For $n >1$,
see $\beta = \beta _n \otimes _{kG} \beta _{n-1} \otimes _{kG}
\cdots \otimes _{kG} \beta _2 \cdot  s(\beta _1)\otimes _{kG} (
s(\beta _1)^{-1}\cdot \beta _1 )$. Thus  $\beta $ can be written as
multiplication of arrows in $Q_1^1$. Consequently, we complete the
proof of (iv). $\Box$

Recall that a braided algebra $A$ in braided tensor category $({\cal
C}, C)$ with braiding $C$ is said to be  braided commutative or
quantum commutative, if $ab = \mu C( a \otimes b)$ for any $a, b\in
A$. An $ {\rm ESC}(G, g_i, \chi _i; i\in J ) $ is said to be quantum
commutative if
$$\chi _i (g_j) \chi _j (g_i)=1$$ for any $i, j \in J$. An $ {\rm ESC}(G, g_i, \chi _i; i\in J ) $ is said to be quantum weakly
commutative if
$$\chi _i (g_j) \chi _j (g_i)=1$$ for any $i, j \in J$ with
$i\not=j$.

\begin {Lemma} \label {14.3.3}

(i) $ {\rm ESC}(G, g_i, \chi _i; i\in J ) $ is   quantum weakly
commutative if and only if in $diag (kG [kQ_1^c, G, g_i, \chi _i; i
\in J])$,
\begin {eqnarray}\label {14.3.10''e1} E_i\cdot E_j = \chi _j
(g_i^{-1}) E_j\cdot E_i \end {eqnarray} for any $i, j \in J$ with $i
\not= j.$

(ii) $diag (kG [kQ_1^c, G, g_i, \chi _i; i \in J])$ is quantum
commutative in $^{kG}_{kG} {\cal YD}$ if and only if  $ {\rm ESC}(G,
g_i, \chi _i; i\in J ) $ is   quantum commutative.

\end {Lemma}

{\bf Proof.} (ii)  By Theorem \ref {14.6}, $R:= diag (kG [kQ_1^c,
g_i, \chi _i; i \in J])$ is a braided Hopf algebra in $^{kG}_{kG}
{\cal YD}$. If $diag (kG [kQ_1^c, g_i, \chi _i; i \in J])$ is
quantum commutative in $^{kG}_{kG} {\cal YD}$, then
\begin {eqnarray}\label {3.10''e1}
E_i\cdot E_j = \chi _j (g_i^{-1}) E_j\cdot E_i \end {eqnarray} for
any $i, j \in J$. By \cite [Theorem 3.8]{CR02} and Proposition \ref
{14.1.10}, we have
$$\begin{array}{c}
E_i\cdot E_j= a^{(j)}_{g_ig_j,g_i}a^{(i)}_{g_i,1}+q_{ij}a^{(i)}_{g_ig_j,g_j}a^{(j)}_{g_j,1}, \\
E_j\cdot E_i = q_{ji}a^{(j)}_{g_ig_j,g_i}a^{(i)}_{g_i,1}+a^{(i)}_{g_ig_j,g_j}a^{(j)}_{g_j,1}.\\
\end{array}$$
Thus $\chi _i (g_j^{-1}) = \chi _j (g_i)$.

Conversely, if $\chi _i (g_j^{-1}) = \chi _j (g_i)$ for any $i, j
\in J, $ see
$$\begin{array}{rcl}
E_i\cdot E_j&=&a^{(i)}_{g_i,1}\cdot a^{(j)}_{g_j,1}\\
&=&a^{(j)}_{g_ig_j,g_i}a^{(i)}_{g_i,1}+q_{ij}a^{(i)}_{g_ig_j,g_j}a^{(j)}_{g_j,1}\\
&=&q_{ij}(q_{ji}a^{(j)}_{g_ig_j,g_i}a^{(i)}_{g_i,1}+a^{(i)}_{g_ig_j,g_j}a^{(j)}_{g_j,1})\\
&=&q_{ij}a^{(j)}_{g_j,1}\cdot a^{(i)}_{g_i,1}\\
&=&q_{ij}E_j\cdot E_i. \\
\end{array}$$ Since $E_i's$ generate $R$ as algebras, $R$ is quantum
commutative.

(i) It is similar to the proof of (ii).
 \  $\Box$

For any positive integers $m$ and $n$, let
$$D_n^{n+m}=\{d=(d_{n+m},d_{n+m-1},\cdots ,d_1)\mid d_i=0\mbox{ or
}1, \sum_{i=1}^{n+m}d_i=n\}.$$ Let $d\in D_n^{n+m}$ and let
$A=a_na_{n-1}\cdots a_1\in Q_n$ be an $n$-path. We define a sequence
$dA=((dA)_{n+m},\cdots,(dA)_1)$ by
$$(dA)_i=\left\{\begin{array}{ccl}
t(a_{d(i)})&,& \mbox{ if }d_i=0;\\
a_{d(i)}&,& \mbox{ if }d_i=1,\\
\end{array}\right.$$
where $1\leq i\leq n+m$ and $d(i)=\sum_{j=1}^id_j$. Such a sequence
$dA$ is called an $n+m$-thin splits of the $n$-path $A$. Note that
if $d(i)=0$ then we regard $t(a_{d(i)})=s(a_1)$, since
$s(a_{d(i)+1})=s(a_1)$ in this case.

%Now let $A$ be an $n$-path and $B$ be an $m$-path. Let $d\in
%D^{n+m}_n$ and let $\overline d\in D_m^{n+m}$ be the complement
%sequence obtained from $d$ by replacing each 0 by 1 and each 1 \ by
%0. In $kQ^c(\alpha)$, consider the element
%$$(A,B)_d=[(dA)_{n+m}\cdot(\overline{d}B)_{n+m}]\cdots
%[(dA)_1\cdot(\overline{d}B)_1]$$ which lies in the $(n+m)$-cotensor
%power $kQ_{n+m}$ of $kQ_1$ and is contained in the isotypic
%component $^{t(A)t(B)}(kQ_{n+m})^{s(A)s(B)}$. By \cite [Theorem
%3.8]{CR02}, we have
%\begin{eqnarray}\label{e4.5}
 %  A\cdot B=\sum_{d\in D^{n+m}_n}(A, B)_d.
%\end{eqnarray}

 If $0 \not= q \in
k$ and $0 \le i \le n < ord (q) $ (the order of $q$), we set
$(0)_{q}! =1$,
$$ \left ( \begin {array} {c} n\\
i
\end {array} \right )_q
 =  \frac{(n)_{q}!}{(i)_{q}!(n - i)_{q}!}, \quad \hbox {where
}(n)_{q}! = \prod_{1 \le i \le n} (i)_{q}, \quad (n)_{q} =
\frac{q^{n} - 1}{q - 1}.$$ In particular, $(n)_q = n$ when $q=1.$

\begin {Lemma}\label{14.3.4} In $kQ^c(G, r, \overrightarrow \chi, u)$, assume $\{ g\} \in
{\mathcal K}_r (G)$ and  $j\in I_{\{g\}}(r)$. Let $q:= \chi
^{(j)}_{\{g\}} (g)$.
 If $i_1, i_2, \cdots, i_m$ be non-negative integers, then
$$\begin{array}{rcl}
a^{(j)}_{g^{i_m +1},g^{i_m }}\cdot a^{(j)}_{g^{i_{m-1}
+1},g^{i_{m-1}}}\cdot\cdots\cdot a^{(j)}_{g^{i_1 +1},g^{i_1 }}
&=&q^{\beta_m}(m)_q! P^{(j)}_{g^{\alpha_m}}(g,m)\\
\end{array}$$
where $\alpha _m = i_1 + i_2 + \cdots + i_m $, $P^{(j)}_h(g,m) =$ \
$ a^{(j)}_{g^mh,g^{m-1}h}a^{(j)}_{g^{m-1}h, g^{m-2}h}\cdots
a^{(j)}_{gh,h}$, $\beta_1=0$ and $\beta_m=\sum_{j
=1}^{m-1}(i_1+i_2+\cdots+i_j )$ if $m>1$.
\end{Lemma}
{\bf Proof.} We prove the  equality by induction on $m$. For $m=1$,
it is easy to see that the equality holds. Now suppose $m>1$. We
have
$$\begin{array}{rl}
&a^{(j)}_{g^{i_m +1},g^{i_m }}\cdot a^{(j)}_{g^{i_{m-1}
+1},g^{i_{m-1}}}\cdot\cdots\cdot a^{(j)}_{g^{i_1 +1},g^{i_1 }}\\
=&a^{(j)}_{g^{i_m +1},g^{i_m }}\cdot(a^{(j)}_{g^{i_{m-1}
+1},g^{i_{m-1}}}\cdot\cdots\cdot a^{(j)}_{g^{i_1 +1},g^{i_1 }})\\
=&q^{\beta_{m-1}}(m-1)_q! a^{(j)}_{g^{i_m +1},g^{i_m }}\cdot
P^{(j)}_{g^{\alpha_{m-1}}}(g,m-1) \ \ \ ( \hbox {by inductive
assumption })\\
=&q^{\beta_{m-1}}(m-1)_q! a^{(j)}_{g^{i_m +1},g^{i_m}}\cdot
(a^{(j)}_{g^{\alpha_{m-1}+m-1},g^{\alpha_{m-1}+m-2}}\cdots
a^{(j)}_{g^{\alpha_{m-1}+1},g^{\alpha_{m-1}}}) \\
 =&q^{\beta_{m-1}}(m-1)_q! \sum_{l=1}^m[(g^{i_m+1}\cdot
a^{(j)}_{g^{\alpha_{m-1}+m-1},g^{\alpha_{m-1}+m-2}})
\cdots(g^{i_m+1}\cdot a^{(j)}_{g^{\alpha_{m-1}+l},g^{\alpha_{m-1}+l-1}})\\
&(a^{(j)}_{g^{i_m+1},g^{i_m}}\cdot g^{\alpha_{m-1}+l-1})
(g^{i_m}\cdot
a^{(j)}_{g^{\alpha_{m-1}+l-1},g^{\alpha_{m-1}+l-2}})\cdots
(g^{i_m}\cdot a^{(j)}_{g^{\alpha_{m-1}+1},g^{\alpha_{m-1}}})]\\
&\ \ \ ( \hbox {by \cite[Theorem 3.8]{CR02} })
\\
=&q^{\beta_{m-1}}(m-1)_q!
\sum_{l=1}^m[a^{(j)}_{g^{\alpha_m+m},g^{\alpha_m+m-1}}
\cdots a^{(j)}_{g^{\alpha_m+l+1},g^{\alpha_m+l}}\\
&(\chi_{\{g\}}^{(j)}(g^{\alpha_{m-1}+l-1})a^{(j)}_{g^{\alpha_m+l},g^{\alpha_m+l-1}})
a^{(j)}_{g^{\alpha_m+l-1},g^{\alpha_m+l-2}}\cdots
a^{(j)}_{g^{\alpha_m+1},g^{\alpha_m}}]    \ \ ( \hbox {by } Proposition  \ref {14.1.10}) \\
=&q^{\beta_{m-1}}(m-1)_q! \sum_{l=1}^m
q^{\alpha_{m-1}+l-1}P^{(j)}_{g^{\alpha_m }}(g,m)\\
=&q^{\beta_{m-1}+\alpha_{m-1}}(m)_q!P^{(j)}_{g^{\alpha_m }}(g,m)\\
=&q^{\beta_m}(m)_q! P^{(j)}_{g^{\alpha_m }}(g,m). \ \ \Box
\end{array}$$

\begin {Lemma} \label {14.3.5} (See \cite [Lemma 3.3]{AS98b}) Let $B$ be a Hopf algebra and $R$  a braided Hopf algebra in
${}_B^B {\cal YD}$ with a linearly independent set $ \{ x_{1} \dots,
x_{t} \}$ $\subseteq $ $P(R)$.  Assume that there exist $g_{j} \in
G(B)$ (the set of all group-like elements in $B$) and  $\chi_{j} \in
Alg(B,k)$ such that
 $$\delta (x_{j}) = g_{j}\otimes x_{j}, \ h\cdot  x_{j} = \chi_{j}(h)x_{j},
 \hbox { for all } h\in B, j =1, 2, \cdots, t . $$
Then \begin{eqnarray*}\{x_1^{m_1} x_2^{m_2}\cdots x_t^{m_t} \mid 0
\leq m_j<N_j, 1\le j \le t \}. \end{eqnarray*} is linearly
independent, where $N_i$ is the order of $q_i := \chi _i (g_i)$  \ (
$N_i = \infty $ when $q_i$ is not a root of unit, or $q_i =1$ ) for
$1 \le i \le t.$
\end {Lemma}
{\bf Proof.} By the quantum binomial formula, if $ 1 \le n_j < N_j$,
then
$$\Delta (x_{j}^{n_{j}}) = \sum_{0 \le i_{j} \le n_{j}}
 \left( \begin {array}{c} n_j\\
i_j
\end {array} \right )
_{q_{j}}x_{j}^{i_{j}}\otimes x_{j}^{n_{j}-i_{j}}.$$

 We use
the  following notation:  $${\bf n } = (n_{1}, \cdots, n_{j}, \cdots
, n_{t}), \quad  x^{\bf n}= x_{1}^{n_{1}} \cdots x_{j}^{n_{j}}
\cdots x_{t}^{n_{t}}, \quad \vert {\bf n }\vert = n_{1} + \cdots +
n_{j} + \cdots  + n_{t};$$ accordingly, ${\bf N}= (N_{1}, \cdots,
N_{t})$, ${\bf  1 }= (1,\cdots, 1)$. Also, we set
$${\bf i }\le {\bf n}\quad \hbox {if } i_{j} \le n_{j}, \, j = 1,
\cdots, t.$$ Then, for ${\bf n } <{\bf N }$, we deduce from the
quantum binomial formula that \begin {eqnarray} \label
{e3.511}\Delta (x^{{\bf n}}) = x^{{\bf n}}\otimes 1 + 1\otimes
x^{\bf  n} + \sum_{ 0 \le {\bf i} \le {\bf n} , \  0 \ne {\bf i }\ne
{\bf n}} c_{\bf {n, i}}x^{\bf i}\otimes x^{\bf  n- i},\end
{eqnarray} where $c_{\bf n,  i} \ne 0$ for all ${\bf i}$.

We shall prove   by induction on $r$ that the set
$$\{x^{\bf n} \mid  \quad \vert {\bf n} \vert \le r, \quad  {\bf n}
< {\bf N} \}$$ is linearly independent.

Let $r=1$ and let $a_{0} + \sum_{i=1}^{t} a_{i}x_{i} = 0$, with
$a_{j}\in k$, $0\le j \le t$. Applying $\epsilon$, we see that
$a_{0} = 0$; by hypothesis we conclude that the other $a_{j}$'s are
also 0.

Now let $r> 1$ and  suppose that $z = \sum_{   \vert {\bf n}\vert
\le r, {\bf n} < {\bf N} } a_{\bf  n} x^{\bf n} = 0$. Applying
$\epsilon$, we see that $a_{0} = 0$. Then
\begin {eqnarray*} 0
&=& \Delta (z) =z\otimes 1 + 1\otimes z + \sum_{ 1\le \vert{\bf n
}\vert \le r,  {\bf n} < {\bf N}} a_{\bf n}\sum_{ 0 \le {\bf i } \le
{\bf n}, \ 0 \ne {\bf i} \ne { \bf
n } }  c_{\bf n, i}x^{\bf i}\otimes x^{\bf n-  i} \\
&=& \sum_{1\le \vert {\bf n }\vert \le r,  {\bf n} < {\bf N}}\ \
\sum_{ 0 \le {\bf i} \le {\bf n}, \  0 \ne {\bf i }\ne { \bf n}}
a_{\bf n} c_{\bf n, \bf i}x^{\bf i}\otimes x^{\bf n-\bf i}. \end
{eqnarray*} Now, if $\vert{\bf n}\vert \le r$, $0 \le {\bf i} \le
{\bf n}$, and $0 \ne {\bf i }\ne {\bf n}$, then $\vert{\bf i}\vert <
r$ and $\vert{\bf n }- {\bf i}\vert < r$. By inductive hypothesis,
the elements $x^{\bf i}\otimes x^{\bf n-\bf i}$ are linearly
independent. Hence $a_{\bf n} c_{\bf n, \bf i} = 0$ and $a_{\bf n} =
0$ for all ${\bf n}$, $\vert{ \bf n }\vert \ge 1$.  Thus $a_{\bf n}
= 0$ for all ${\bf n}$. $\Box$

 Assume that $(G, g_i, \chi_i; j\in J)$ is an ${\rm ESC}.$ Let
${\cal T} (G, g_i, \chi_i; j\in J)$ be the free algebra generated by
set $ \{ x_j \mid j\in J \}$. Let ${\cal S} (G, g_i, \chi_i; j\in
J)$ be the algebra generated by set $ \{ x_j \mid j\in J \}$ with
relations
\begin {eqnarray} \label {qse1} x_{i}x_{j} = \chi_{j}(g_{i})
x_{j}x_{i}  \ \ \ \hbox { for any } i, j \in J \hbox { with } i
\not= j.
\end {eqnarray}
Let ${\cal R} (G, g_i, \chi_i; j\in J)$ be the algebra generated by
set $ \{ x_j \mid j\in J \}$  with relations
\begin {eqnarray} \label {qlse1} x_l ^{N_l}=0,\  x_{i}x_{j} = \chi_{j}(g_{i})
x_{j}x_{i}  \ \ \ \hbox { for any } i, j,l  \in J \hbox { with } N_l
< \infty,  i \not= j.
\end {eqnarray}
 Define their coalgebra operations and $kG$-(co-)module operations  as follows:
\begin {eqnarray} \label {e3.522}\Delta x_i = x_i \otimes 1 + 1 \otimes x_i, \ \ \epsilon (x_i)
=0, \ \ \delta ^-(x_{i}) = g_{i}\otimes x_{i}, \qquad h \cdot x_{i}
= \chi_{i}(h)x_{i}.\end {eqnarray} $ {\cal T} (G, g_i, \chi_i; j\in
J)$ is called a quantum tensor algebra in $^{kG}_{kG} {\cal YD} $,
${\cal S} (G, g_i, \chi_i; j\in J)$ is called a quantum symmetric
algebra in $^{kG}_{kG} {\cal YD} $ and ${\cal R} (G, g_i, \chi_i;
j\in J)$ is called a quantum linear space in $^{kG}_{kG} {\cal YD}
$. Note that when ${\rm ESC} (G, g_i, \chi _i; i \in J)$ is quantum
weakly commutative with finite $J$ and  finite $N_j$ for any $j\in
J$, the definition of quantum linear space is the same as in \cite
[Lemma 3.4]{AS98b}. Obviously, if $N_i$ is infinite for all $i \in
J$, then ${\cal S} (G, g_i, \chi_i; j\in J)$ = ${\cal R} (G, g_i,
\chi_i; j\in J)$.

\begin{Theorem}\label{14.7}
Assume that ${\rm ESC} (G, g_i, \chi _i; i \in J)$ is quantum weakly
commutative.  Let $\prec$ be  a total order of $J$. Then

{\rm(i)} \ The multiple Taft algebra $kG[kQ_1^c,  G,  g_i, \chi_i;
i\in J]$ has a $k$-basis
\begin{eqnarray*}\label{basis} \begin {array} {c} \{  g\cdot E_{\nu _1}^{m_1}\cdot
 E_{\nu _2}^{m_2}\cdot\cdots\cdot  E_{\nu _t}^{m_t} \mid 0 \le
m_j < N_j; \nu _j \prec \nu _{j+1}, \nu _j \in J,  j= 1, 2, \cdots ,
t;  t \in {\mathbb Z}^+, g\in G\}. {} \end {array} \end{eqnarray*}
Moreover, $kG[kQ_1^c, G g_i, \chi_i; i\in J]$ is finite dimensional
if and only if $\mid \! G \! \mid,$  $\mid \! J \! \mid $  and $N_j
$ are finite for any $j\in J$. In this case, ${\rm dim}_k(kG[kQ_1^c,
G, g_i, \chi_i; i\in J])=|G|N_1N_2\cdots N_t$ with $J = \{1, 2,
\cdots, t\}$.

 {\rm(ii)}\ $diag ( kG[kQ_1^c, G, g_i, \chi _i; i\in J])$ has a $k$-basis
\begin{eqnarray}\label{basis2}\{ E_{\nu _1}^{m_1}\cdot
 E_{\nu _2}^{m_2}\cdot\cdots\cdot  E_{\nu _t}^{m_t} \mid 0 \le
m_j < N_j; \nu _j \prec \nu _{j+1}, \nu _j \in J,  j= 1, 2, \cdots ,
t;  t \in {\mathbb Z}^+\}.\end{eqnarray}

{\rm(iii)} \ $diag ( kG[kQ_1^c, G, g_i, \chi _i; i\in J])$ is a
Nichols algebra in   ${}_{kG}^{kG}{\cal YD}$  and   ${\cal R} (G,
g_i, \chi_i ^{-1}; j\in J) \cong diag ( kG[kQ_1^c, G,  g_i, \chi _i;
i\in J])$ as graded braided Hopf algebras in ${}_{kG}^{kG}{\cal
YD}$,  by sending $x_j$ to $a_{g_j, 1}^{(j)}$ for any $j\in J$.

{\rm(iv)} \ ${\cal T} (G, g_i, \chi_i ^{-1}; j\in J) \cong  diag (
kQ^s( g_i, \chi _i; i\in J))$ as graded braided Hopf algebras in
${}_{kG}^{kG}{\cal YD}$ algebras, by sending $x_j$ to $a_{g_j,
1}^{(j)}$ for any $j\in J$.

{\rm(v)} \  $kQ^s ( G, g_i, \chi_i; i\in J) \cong {\cal T} (G, g_i,
\chi_i ^{-1}; j\in J)\# kG$ \ as graded Hopf algebras and  $kQ^s (G,
g_i, \chi_i; i\in J)$ has a $k$-basis
\begin{eqnarray*} \{  g\cdot E_{\nu _1}\otimes _{kG}
 E_{\nu _2}\otimes _{kG}\cdots \otimes _{kG} E_{\nu _t} \mid  \nu _j \in J,  j= 1, 2, \cdots
, t;  t \in {\mathbb Z}^+ \cup \{0\}, g\in G\}, \end{eqnarray*}
where $g\cdot E_{\nu _1}\otimes _{kG}
 E_{\nu _2}\otimes _{kG}\cdots \otimes _{kG} E_{\nu _t}=g$ when $t=0.$

 Note that (iv) and (v) still hold without quantum weakly commutative
condition.
 \end{Theorem}
{\bf Proof.} (ii)    Since $(N_j)_{q_j}!=0$, it follows from Lemma
\ref{14.3.4} that $E_j^{N_j}=0$ when $N_j <\infty$. By Lemma \ref
{14.3.3}, $E_i\cdot E_j = \chi _j (g_i ^{-1}) E_j\cdot E_i$ for any
$i, j \in J$ with $i \not= j$. Considering Lemma \ref{14.3.5}, we
complete the proof.

(iii) By Lemma \ref {14.3.4} and Eq.(\ref {14.3.10''e1}), there
exists an algebra homomorphism  $\psi $ from ${\cal R} (G, g_i,
\chi_i^{-1}; j\in J) $ to $diag ( kG[kQ_1^c, G,  g_i, \chi _i; i\in
J])$ by sending $x_j$ to $a_{g_j, 1}^{(j)}$ for any $j\in J$. By
(ii), $\psi$ is bijective. It is clear that $\psi$ is a graded
braided Hopf algebra isomorphism.

Let $R : = {\cal R}(G, g_i, \chi_i; i\in J)$.  Obviously, $R_{(1)}
\subseteq P(R)$. It is sufficient to show that any non-zero
homogeneous element $z\in R$, whose  degree $deg (z)$ is not equal
to 1, is not a primitive element. Obviously, $z$ is not a primitive
element when $deg (z)=0$. Now  $deg (z)>1$. We can assume, without
lost generality,  that there exist $\nu _1, \nu _2, \cdots, \nu _t
\in J $ such that
  $z = \sum _{\mid {\bf i} \mid = n} k_{ \bf i}
x^{\bf i}$, where $k_i \in k $, $x ^{\bf i} = x_{\nu_1}
^{i_1}x_{\nu_2} ^{i_2}\cdots x_{\nu_t} ^{i_t}$ with $ i_1 + i_2
+\cdots + i_t =n$.   It is clear
\begin {eqnarray} \label
{e611}  \Delta (z) = \sum _{\mid {\bf i} \mid =n} k_{\bf i}\Delta
(x^{{\bf i}}) = z\otimes 1 + 1\otimes z + \sum _{\mid {\bf i} \mid
=n} \ \sum_{ 0 \le {\bf j} \le {\bf i} , \  0 \ne {\bf j }\ne {\bf
i}} k _{\bf i} c_{\bf {i, j}}x^{\bf j}\otimes x^{\bf  i- j}.\end
{eqnarray} If $z$ is a primitive element, then   $\sum _{\mid {\bf
i} \mid =n} \ \sum_{ 0 \le {\bf j} \le {\bf i} , \  0 \ne {\bf j
}\ne {\bf i}} k _{\bf i} c_{\bf {i, j}}x^{\bf j}\otimes x^{\bf  i-
j}=0$. Since $c_{{\bf i}, {\bf j}} \not=0$, we have $k_{\bf i} =0$
for any ${\bf i}$ with $\mid {\bf i} \mid  =n$, hence $z=0$. We get
a contradiction. Thus $z$ is not a primitive element. This show
$R_{(1)} = P(R)$ and $R$ is a Nichols algebra.

(iv) and (v). Let $A:= {\cal T} (G, g_i, \chi_i ^{-1}; j\in J)$ and
$R = diag ( kQ^s( G, g_i, \chi _i; i\in J))$. Let $\psi$ be an
algebra homomorphism from ${\cal T} (G, g_i, \chi_i ^{-1}; j\in J)$
to $ diag ( kQ^s ( G,  g_i, \chi _i; i\in J))$ by sending $x_j$ to
$a_{g_j, 1}^{(j)}$.

It is clear that ${\cal T} (G, g_i, \chi_i^{-1}; j\in J)$ is a
$kG$-module algebra. Define a linear map $\nu $ from  $A \# kG$ to
$kQ^s$ by sending $x_j \# g$ to $a _{g_j, 1}^{(j)} \cdot g = \alpha
^+ (a _{g_j, 1}^{(j)} \otimes g)$ for any $g\in G, j \in J$. That
is, $\nu$ is the composition of
$$A \# kG \stackrel {\psi \otimes id} {\rightarrow } R \# kG
\stackrel {\alpha ^+} {\cong} kQ^s,
$$ where $\alpha ^+ = \mu _{kQ^s}( id \otimes \iota _0)$ (see
Theorem \ref {14.6}).  Define  a linear  map  $\lambda$ from $kG$ to
$A\#kG$ by sending $g$ to $1\# g$ for any $g \in G$ and another
linear map $\gamma$ from $kQ_1^c$ to $A\#kG$ by sending $a_{g_ih,
h}^{(i)}$ to $\chi _i^{-1}(h)x_i \# h $ for any $h \in G, i \in J.$
It is clear that $\gamma$ is a $kG$-bimodule homomorphism from
$kQ_1^c$ to $ _\lambda( A\#kG)_ \lambda$. Considering $kQ^s = T_{kG}
(kQ_1^c)$ and universal property of tensor algebra over $kG$, we
have that there exists an algebra homomorphism $\phi = T_{kG}
(\lambda, \gamma) $ from $kQ^s$ to $A\# kG$. Obviously, $\phi$ is
the inverse of $\nu$.
  Thus $\phi$ is bijective. It is easy to
check that $\phi$ is graded  Hopf algebra isomorphism. Obviously,
$\{  x_{\nu _1}
 x_{\nu _2} \cdots  x_{\nu _t} \# g \mid  \nu _j \in J,  j= 1, 2, \cdots
, t;  t \in {\mathbb Z}^+ \cup \{0\}, g\in G\}$ is a basis of ${\cal
T} (G, g_i, \chi_i^{-1}; j\in J) \# kG$. See
\begin {eqnarray*} && \nu ( x_{\nu _1}
 x_{\nu _2} \cdots  x_{\nu _t} \# g )\\
 &=&E_{\nu _1}\otimes _{kG}
 E_{\nu _2}\otimes _{kG}\cdots \otimes _{kG} E_{\nu _t} \cdot g \\
 &= & \chi _{\nu _1 } (g) \chi _{\nu _2 } (g) \cdots \chi _{\nu _t } (g) g\cdot E_{\nu _1}\otimes
_{kG}
 E_{\nu _2}\otimes _{kG}\cdots \otimes _{kG} E_{\nu _t}.
\end {eqnarray*}
Thus $ \{  g\cdot E_{\nu _1}\otimes _{kG}
 E_{\nu _2}\otimes _{kG}\cdots \otimes _{kG} E_{\nu _t} \mid  \nu _j \in J,  j= 1, 2, \cdots
, t;  t \in {\mathbb Z}^+ \cup \{0\}, g\in G\}$ is a basis of
$kQ^s$. It is easy to check that $\psi$ is graded braided Hopf
algebra isomorphism.

 (i) Considering (ii),  Lemma \ref {14.3.2} and Theorem \ref {14.6},
we complete the proof. \ \ $\Box$

\subsection {A characterization of multiple Taft algebras
}\label {s3.3}

In this subsection we characterize multiple Taft algebras by means
of elements in themselves.
\begin{Definition}\label{14.3.6} For a quantum weakly commutative   ${\rm ESC} (G, g_i, \chi _i; i\in
J)$,  let $A$ be the Hopf algebra to satisfy the following
conditions: {\rm(i)} $G$ is  a subgroup of $G(A)$;
 {\rm(ii)} there exists a linearly independent  subset $\{X_i \mid i\in J\}$ of $A$ such
 that $A$ is generated by set $\{X_i \mid i\in J\} \cup G$ as
 algebras;
 {\rm(iii)}  $X_j$ is $(1, g_j)$-primitive, i.e.,
$\Delta(X_j)=X_j\otimes 1+g_j\otimes X_j$, for any $j\in J;$
{\rm(iv)} $X_jg=\chi_j(g)gX_j$, for any $j\in J,$, $g\in G$;
 {\rm(v)}
$X_jX_i = \chi_j(g_i) X_iX_j$, for $i, j \in J$ with $i\not=j;$
 {\rm(vi)} $A_{(0)} \cap A_{(1)} =0$, where $A_{(0)}:=kG$ and
 $A_{(1)}$ is the vector space spanned by set $\{hX_i \mid  i \in J; h \in
 G\}$.
Furthermore, let  $J(A)$ denote   the ideal of $A$ generated by the
set
$$\{X_i^{N_i}  \mid N_i < \infty , i \in J \} $$ and
$H(G,g_i, \chi _i; i\in J )$  the quotient algebra $A/J(A)$.
\end{Definition}

\begin {Lemma} \label {14.3.5'} (See \cite [Lemma 3.3]{AS98b}  ) Let  $H$ be a Hopf algebra  with a linearly independent
set $ \{ x_{1} \dots, x_{t} \}$ and $G$ a subgroup of   $ G(H)$.
Assume that $g_{i} \in Z(G)$ and $\chi_{i} \in \hat G$ such that
$\Delta (x_i) = x_i \otimes 1+ g_i \otimes x_i$, \
 $x_ih=\chi_i(h)hx_i$,  for $i = 1, 2, \cdots,  t$, $h\in G.$ If the
intersection of $kG$ and span $\{hx_i \mid 1\le i \le t, h\in G\}$
is zero, then
\begin{eqnarray*}\{hx_1^{m_1} x_2^{m_2}\cdots x_t^{m_t} \mid 0 \leq
m_j<N_j, 1\le j \le t; h \in G \}.
\end{eqnarray*} is linearly independent, where $N_i$ is the order of
$q_i := \chi _i (g_i)$  \ ( $N_i = \infty $ when $q_i$ is not a root
of unit, or $q_i =1$ ) for $1 \le i \le t.$
\end {Lemma}

{\bf Proof.} (i) We first show  set $\{h x_i \mid i=1, 2, \cdots, t;
h\in G\}$ is linearly independent. Indeed, suppose  that $ z = \sum
_{h\in G} \sum _{i=1}^t k_{i, h} hx_i =0$ with $k_{i, h} \in k$ for
$h\in G$, $1\le i \le t$. Applying $\Delta$, we get
$$\sum _{h\in G} \sum _{i=1}^t k _{i, h} ( hg_i \otimes hx_i +  hx_i \otimes h)
=0, $$ which implies that $\sum _{i=1}^t k_{i, h} h x_i=0$ for any
$h\in H$, hence $k_{i, h} =0$ for any $h\in H$, $i =1, 2,\cdots, t$.

(ii) We use the  notation in the proof of Lemma \ref {14.3.5}. We
shall prove by induction on $r$ that the set
$$\{hx^{\bf n} \mid  \quad \vert {\bf n} \vert \le r, \quad  {\bf n}
< {\bf N}; h\in G \}$$ is linearly independent.

 By the quantum binomial formula, if $ 1 \le n_j < N_j$,
then \begin {eqnarray} \label {e3.522'} \Delta (x_{j}^{n_{j}}) =
\sum_{0 \le i_{j} \le n_{j}}
 \left( \begin {array}{c} n_j\\
i_j
\end {array} \right )
_{q_{j}}g_{j}^{i_j}x_{j}^{n_j - i_{j}}\otimes x_{j}^{i_{j}}.\end
{eqnarray}
  For ${\bf n } <{\bf N }$, we
deduce from Eq. (\ref {e3.522'}) or  \cite [Eq. (5.9)]{DNR01})
 that \begin {eqnarray} \label {e3.511'}\Delta (x^{{\bf n}}) = x^{{\bf
n}}\otimes 1 + g^{\bf n}\otimes x^{\bf  n} + \sum_{ 0 \le {\bf i}
\le {\bf n} , \  0 \ne {\bf i }\ne {\bf n}} c_{\bf {n, i}} g^{\bf
i}x^{\bf n- i}\otimes x^{\bf  i},\end {eqnarray} where $g^{\bf i}
=g_1^{i_1} g_2^{i_2} \cdots  g_t^{i_t}$ and $0 \ne c_{\bf n, i} \in
k$ for all ${\bf i}$. Obviously, the claim holds when $r=0$.

 Let $r=1$ and let $\sum _{h\in G }b_h h + \sum _{h\in
G}\sum_{i=1}^{t} b_{h, i}hx_{i} = 0$, with $b_{h, j}, b_h \in k$,
$0\le j \le t$. It follows from hypothesis that $\sum _{h\in G}b_h h
=0 $ and $ \sum _{h\in G}\sum_{i=1}^{t} b_{h, i}hx_{i} =0$, hence
$b_h=0$ and $b_{h, i}=0$ by (i) for any $h\in G,$ $1\le i\le t.$

Now let $r> 1$ and  suppose that $z =\sum _{h\in G} \sum_{   \vert
{\bf n}\vert \le r, \  {\bf n} < {\bf N}} b_{h, {\bf  n}} h x^{\bf
n} = 0$. Then
\begin {eqnarray*} 0
&=& \Delta (z) =z\otimes 1 + 1\otimes z +     \sum _{h\in G} b_{h,
{\bf 0}} h \otimes h + \sum _{h\in G}\sum_{ 1\le \vert\! {\bf n} \!
\vert \le r, {\bf n} < {\bf N}} b_{h, {\bf n}}\sum_{ 0 \le {\bf i }
\le {\bf n}, \ 0 \ne {\bf
i} \ne { \bf n }}  c_{\bf n, i}hg^{\bf i}x^{\bf n- i}\otimes x^{\bf i} \\
&=&\sum _{h\in G} \sum_{1\le \vert {\bf n }\vert \le r,   {\bf n} <
{\bf N}}\ \ \sum_{ 0 \le {\bf i} \le {\bf n}, \  0 \ne {\bf i }\ne {
\bf n}} b_{h, {\bf n}} c_{\bf n, \bf i}hg^{\bf i}x^{\bf n- i}\otimes
x^{\bf i}. \end {eqnarray*} Now, if $\vert{\bf n}\vert \le r$, $0
\le {\bf i} \le {\bf n}$, and $0 \ne {\bf i }\ne {\bf n}$, then
$\vert{\bf i}\vert < r$ and $\vert{\bf n }- {\bf i}\vert < r$. By
inductive hypothesis, the set $\{ h g ^{\bf i} x ^{\bf n-i}\otimes h
x ^{\bf i} \mid  {\bf 0}\le {\bf i} \le { \bf n}, {\bf 0}\ne {\bf i}
\ne { \bf n}, \mid \! {\bf n} \! \mid < r,  {\bf n} < {\bf N}, h\in
G \}$ $\cup \{h \otimes h \mid h\in G\}$ are linearly independent.
Hence $b _{h, {\bf o}} =0$ and $b_{h, {\bf n}} c_{\bf n, \bf i} =
0$, which implies $b_{h, {\bf n}} = 0$ for any $h\in G, $ all ${\bf
n}$. $\Box$

\begin{Proposition}\label{14.3.7} If    ${\rm ESC} (G, g_i, \chi _i; i\in
J)$ is  quantum weakly commutative, then
 $H(G;$ \ $ g_i, \chi_i; $ $ i\in J)$ and multiple Taft algebra $kG[kQ_1^c, G,  g_i,
\chi_i; $ $i\in J]$ are isomorphic as graded Hopf algebras.
\end{Proposition}

{\bf Proof.}  Let $H:= H(G; g_i, \chi_i; i\in J)$. Denote by $X_j$
the image of $X_j$ under the canonical projection $A\rightarrow
A/J(A)$. We say that
 $hX_{\nu _1}
^{i_1} X_{\nu _2} ^{i_2}\cdots X_{\nu _t} ^{i_t}$ is a monomial with
degree $n$ if $h\in G$ and $i_1 + i _2 + \cdots + i_t =n$, where
$\nu _1, \nu_2,\cdots , \nu _t \in J$, $t \in {\mathbb Z}^+$. Let
  $H_{(0)} : =
kG$ and  $H_{(n)}$ be the vector space spanned by all monomials with
degree $n$.  It follows from Lemma \ref {14.3.5'} that $H$ is a
graded Hopf algebra.
 By the method similar to the proof of Lemma
\ref {14.3.2}(ii), we can obtain that $H^{co \ kG}$ is a subalgebra
of $H$ generated by $X_i's$. By Theorem \ref {14.6}, $H^{co \ kG}$
is a braided Hopf algebra in $^{kG}_{kG} {\cal YD}$ and $H$ is the
bosonization of $H^{co \ kG}$. Considering Theorem \ref {14.7}
(iii), we only need to show that ${\cal R} (G, g_i, \chi_i ^{-1};
j\in J) \cong H^{co\ kG}$ as graded braided Hopf algebras. Let
$\prec$ be  a total order of $J$.
 Since $(N_j)_{q_j}!=0$, it follows from Lemma
\ref{14.3.4} that $X_j^{N_j}=0$ when $N_j <\infty$. Considering
Lemma \ref{14.3.5},  we have that
\begin{eqnarray*}\{ X_{\nu _1}^{m_1}
 X_{\nu _2}^{m_2}\cdots X_{\nu _t}^{m_t} \mid 0 \le
m_j < N_j; \nu _j \prec \nu _{j+1}, \nu _j \in J,  j= 1, 2, \cdots ,
t;  t \in {\mathbb Z}^+\}\end{eqnarray*} is a basis of $diag (H)$.
Consequently, we have  an algebra isomorphism    $\psi $ from ${\cal
R} (G, g_i, \chi_i ^{-1}; i \in J)$ to $H^{co \ kG}$ by sending $
x_i$ to $X_i$. It is easy to check that $\psi$ is also a graded
braided Hopf algebra isomorphism. $\Box$

 {\bf Remark:} $H(G,g_i, \chi _i; i\in J )$ just is
$H(C, n, c, c^*, 0, 0)$ in \cite [ Definition 5.6.8 and Definition
5.6.15] {DNR01} with $G = C$, $J = \{1, 2, \cdots , t\}$, $1<n_i =
N_i < \infty$, $g_i =c_i, c_i^* = \chi _i$ for $i = 1, 2, \cdots,
t.$

\subsection {The relation between semi-path Hopf algebras and quantum enveloping
algebras}\label {14.3.4} If $0\not= q \in k$ and $0 \le i \le n <
ord (q) $ (the order of $q$), we set
$$ \left [ \begin {array} {c} n\\
i
\end {array} \right ]_q
 =  \frac{[n]_{q}!}{[i]_{q}![n - i]_{q}!}, \quad \hbox {where
}[n]_{q}! = \prod_{1 \le i \le n} [i]_{q}, \quad [n]_{q} =
\frac{q^{n} - q^{-n}}{q - q^{-1}}.$$

Let $B$ be a Hopf algebra and $R$  a braided Hopf algebra in ${}_B^B
{\cal YD}$. For convenience,  we denote   $r\# 1$ by $r$  and $1 \#
b$ by  $b$ in    biproduct $R \# B$ for any $r\in R,\  b \in B$.

\begin {Lemma} \label {14.3.7'}  Let $B$ be a Hopf algebra and $R$  a braided Hopf algebra in
${}_B^B {\cal YD}$ with $x_1, x_2 \in P(R)$ and $g_1, g_2 \in
Z(G(B))$. Assume that there exist $\chi_1, \chi _2 \in Alg(B,k)$
with $\sqrt{\chi _i (g_j)}\in k $ such that
 $$\delta _R(x_{j}) = g_{j}\otimes x_{j}, \ h\cdot  x_{j} = \chi_{j}(h ^{-1})x_{j},
$$
  for all $ h\in B, j =1, 2. $

(i) If $r$ is a positive integer and
  \begin {eqnarray} \label {quantume1'}
\begin {array} {cc}
 \chi _2 (g_1) \chi_1(g_2) \chi
_1 (g_1) ^{r-1} =1, \  \chi _1 (g_1) ^ {\frac {1}{2} (r-1)} \chi _2
(g_1) =1,\ r-1 < ord (\chi _1(g_1)),
\end {array}
\end {eqnarray} then
\begin {eqnarray}\label {quantume1'''}
 \sum_{m=0}^{r}(-1)^m
\left[\begin{array}{c}
r\\
m\\
\end{array}\right]_{\sqrt{\chi _1(g_1^{-1}})}
x_1^{r-m}x_2x_1^m \end {eqnarray} = $(ad_c x_1)^{r}x_2$ is a
primitive element of $ R$ and a $(1, g_1^rg_2)$-primitive element of
biproduct $R\#B$,  where $(ad_cx_1)x_2 = x_1x_2 - \chi _2
(g_1^{-1})x_2 x_1$.

  (ii) If $ \ \sqrt{\chi _1 (g_2)}  \sqrt{\chi _2(g_1)}=1$ and  $x _i g_j = \chi _i (g_j)g_jx_i$ for $i, j =1, 2$, then
\begin {eqnarray*}  \sqrt{\chi _2 (g_1)} x_1x_2- \sqrt{\chi
_1 (g_2)}x_2 x_1
\end {eqnarray*} is a primitive element of
$ R$  and
  \begin {eqnarray*}  \sqrt{\chi _2 (g_1)}x_1x_2- \sqrt{\chi
_1 (g_2)}x_2 x_1- \beta (g_1 g_2-1)
\end {eqnarray*} is a $(1, g_1g_2)$-primitive element of biproduct
$ R\#B$ for any $\beta \in k$.

\end {Lemma}
{\bf Proof.} (i)   Let  $ z:= \sum_{m=0}^{r}(-1)^m
\left[\begin{array}{c}
r\\
m\\
\end{array}\right]_{\sqrt {\chi _1(g_1^{-1})}}
x_1^{r-m}x_2x_1^m $.   By \cite [Lemma A.1] {AS00}, $(ad_c
x_1)^{r}x_2 $ is primitive.  However,
\begin {eqnarray*} && (ad_c x_1)^{r}x_2
\\ &=& \sum_{m=0}^{r}(-1)^m \left(\begin{array}{c}
r\\
m\\
\end{array}\right)_{\chi _1(g_1^{-1})} \chi _1 (g_1^{-1})) ^{\frac {1}{2}m
(m-1)} \chi _2 (g_1^{-1}))^m x_1^{r-m}x_2x_1^m \ \
 (\hbox {by }  \cite [(A.8)] {AS00} )\\
&{=}&\sum_{m=0}^{r}(-1)^m \left[\begin{array}{c}
r\\
m\\
\end{array}\right]_{\sqrt {\chi _1(g_1^{-1})}} \chi _1 (g_1^{-1})) ^ {\frac {1}{2} m
(r-1)} \chi _2 (g_1^{-1})) ^m x_1^{r-m} x_2x_1^m   \ \ (\hbox
{by \cite [P40] {AS00}}) \\
&{=}&\sum_{m=0}^{r}(-1)^m \left[\begin{array}{c}
r\\
m\\
\end{array}\right]_{\sqrt {\chi _1(g_1^{-1})}} x_1^{r-m} x_2x_1^m   \ \ (\hbox
{by assumption  } (\ref {quantume1'})). \\
\end {eqnarray*}
Thus $z$ is primitive in $R$,
 i.e. \begin {eqnarray}
\label {quantume3} \Delta _R (z) = z \otimes 1 + 1 \otimes z. \end
{eqnarray}

Since \begin {eqnarray*}\delta_R (z) &=& \sum_{m=0}^{r}(-1)^m
\left[\begin{array}{c}
r\\
m\\
\end{array}\right]_{\sqrt {\chi _1(g_1^{-1})}}\delta_R (
x_1^{r-m}x_2x_1^m)\\
&=&\sum_{m=0}^{r}(-1)^m \left[\begin{array}{c}
r\\
m\\
\end{array}\right]_{\sqrt {\chi_1(g_1^{-1})}} g_1^ rg_2 \otimes
x_1^{r-m}x_2x_1^m\\
&=&g_1^ rg_2 \otimes z,
\end {eqnarray*}
we have

 \begin {eqnarray*}\Delta _{R\#B} (z) &=& \Delta _{R \# B} (z\# 1)\\
&=&   (id _R \otimes \mu _B \otimes id _R \otimes id _B) (id _R
\otimes id _B \otimes C_{R, B} \otimes id
_B) \\
&{}& (id _R \otimes \delta _R \otimes \Delta _B) (1\# z \otimes 1) +
z
 \otimes 1  \ \ ( \hbox {by Theorem \ref
{14.6}}) \\
&=& g_1^rg_2 \otimes z +  z \otimes 1,
\end {eqnarray*} where $C_{R, B}$ denotes the ordinary twist map
from  $R \otimes B$ to $B \otimes R$ by sending $r \otimes b$ to $b
\otimes r$ for any $r\in R, b\in B.$

(ii) See
\begin {eqnarray*}
&&\Delta _{R\#B}(\sqrt {\chi_2(g_1)} x_1x_2-\sqrt {\chi_1(g_2)} x_2
x_1 )
\\
&=& \sqrt {\chi_2(g_1)}  (g_1g_2 \otimes x_1x_2 + x_1 g_2 \otimes
x_2 +
g_1x_2 \otimes x_1 + x_1x_2 \otimes 1 )\\
&& -\sqrt {\chi _1 (g_2)}(g_1g_2 \otimes x_2x_1 + x_2g_1 \otimes x_1
+ g_2x_1 \otimes x_2 + x_2x_i
\otimes 1 )\\
&=&g_1g_2 \otimes   (    \sqrt {\chi _2 (g_1)}x_1x_2- \sqrt {\chi _1
(g_2)} x_2 x_1) + ( \sqrt {\chi _2 (g_1)}x_1x_2-\sqrt {\chi _1
(g_2)}x_2 x_1)\otimes 1 \\
&& \ \ (\hbox {by assumption of (ii)}  )
\end {eqnarray*}
 and
\begin {eqnarray*}
&&\Delta _{R\#B} ( \sqrt {\chi _2 (g_1)} x_1x_2-\sqrt {\chi _1
(g_2)} x_2 x_1 ) -\Delta _{R\#B} ( \beta (g_1g_2 -1))
\\
&=&g_1g_2   \otimes( \sqrt {\chi _2 (g_1)} x_1x_2-\sqrt {\chi _1
(g_2)} x_2 x_1  -
\beta (g_1g_2 -1)) \\
 &&
 + ( \sqrt {\chi _2 (g_1)} x_1x_2-\sqrt {\chi _1 (g_2)} x_2
x_1  - \beta (g_1g_2 -1)) \otimes 1.
\end {eqnarray*} Thus  $\sqrt{\chi _2 (g_1)}x_1x_2- \sqrt{\chi
_1 (g_2)}x_2 x_1- \beta (g_1 g_2-1)$
 is a $(1, g_1g_2)$-primitive element of
$ R\#B$. See
\begin {eqnarray*} &&\Delta _R (\sqrt {\chi_2(g_1)}
x_1x_2-\sqrt {\chi_1(g_2)} x_2 x_1 )
\\
&=& \sqrt {\chi_2(g_1)}  (1 \otimes x_1x_2 + x_1 \otimes x_2 +
   \chi _2 (g_1) ^{-1}x_2 \otimes x_1 + x_1x_2 \otimes 1 )\\
&& -\sqrt {\chi _1 (g_2)}(1 \otimes x_2x_1 + x_2 \otimes x_1 + \chi
_1 (g_2)^{-1}x_1 \otimes x_2  + x_2x_1
\otimes 1 )\\
&=&1 \otimes   (    \sqrt {\chi _2 (g_1)}x_1x_2- \sqrt {\chi _1
(g_2)} x_2 x_1) + ( \sqrt {\chi _2 (g_1)}x_1x_2-\sqrt {\chi _1
(g_2)}x_2 x_1)\otimes 1 \\
&& \ \ (\hbox {by assumption of (ii)}  ).
\end {eqnarray*} Thus  $\sqrt{\chi _2 (g_1)}x_1x_2- \sqrt{\chi
_1 (g_2)}x_2 x_1$
 is a primitive element of $R$.
$\Box$

For an ${\rm ESC}(G, \chi_i, g_i; i\in J)$, we give the follows
notations:

(FL1) $N$  is a set and $J = \cup _{s\in N} (J_s \cup J_{s}')$ is a
disjoint union.

(FL2) There exists a bijection $\sigma : J^{(1)}\rightarrow J^{(2)}$
such that $\sigma \mid _{J_u}$ is a bijection from $J_u$ to $J_u'$
for $u \in N,$  where  $J^{(1)} := \cup _{u\in N} J_u $ and $J^{(2)}
:= \cup _{u\in N} J_u'.$

(FL3) There exists a  $J^{(1)}\times J^{(1)}$-matrix $A = (a_{ij})$
with $a_{ii} =2$ and non-positive integer $a_{ij} $ for any $i, j
\in J^{(1)}$ and $i \not= j$. For any $u \in N$, there exists an
integer $d_{i}^{(u)}$ such that $d_i ^{(u)}a_{ij} = d_j
^{(u)}a_{ji}$ for any $i, j \in J_u$.

(FL4) For any $u \in N$, there exists $0 \not= q_u \in k$  such that
$\chi _i (g_j) = q_u ^{-2d_i^{(u)}a_{ij}}$, $\chi _{\sigma (i)}
(g_j) = \chi _i ^{-1}(g_j)$  and $g_{\sigma (j)} = g_j $ for $i, j
\in J_u$.

(FL5) There exists $\xi _i\in G$ such that $\chi _{\sigma (i)} (\xi
_j) = \chi_i ^{-1}(\xi _j)$, $\xi _{\sigma (i)} := \xi _i ^{-1}$ and
$g_i= g_{\sigma (i)} := \xi _i ^2$  for any $i, j \in J_u, \ u \in
N;$ there exists a positive integer $r_{ij}$ such that  $r_{\sigma
(i), \sigma (j)}= r_{ij}$ for any $i, j \in J^{(1)} $ with $i\not= j
$.

(FL6) $
 \chi _j (g_i) \chi_i(g_j) \chi
_i (g_i) ^{r_{ij}-1} =1,$ $ \chi _i(\xi_j) = \chi _j(\xi _i), $ $
\chi _i (g_i) ^ {\frac {1}{2}  (r_{ij}-1)} \chi _j (g_i) =1,$ for
any  $i, j \in J_u $, $i \not= j,$ $ u\in N.$

(FL7) $G$ is a free commutative  group generated by  generator set
$\{\xi_i \mid i \in J^{(1)}\}$.

 An ${\rm ESC}(G, g_i, \chi _i; i\in J)$ is said to be a  local
FL-matrix type   (see \cite [P.4]{AS00}) if (FL1)--(FL4) hold. An
${\rm ESC}(G, g_i, \chi _i; i\in J)$ is said to be a  local FL-type
if (FL1), (FL2), (FL5) and (FL6) hold. An ${\rm ESC}(G, g_i, \chi
_i; i\in J)$ is said to be a  local FL-free type    if (FL1), (FL2),
(FL5), (FL6) and (FL7) hold. An ${\rm ESC}(G, g_i, \chi _i; i\in J)$
is said to be a  local FL-quantum group type    if (FL1)-- (FL4) and
(FL7) hold. If $N $ only contains one element  and $J^{(1)} = \{1,
2, \cdots, n\}$, then  we delete `local' in  the terms above.

Let ${\rm ESC}(G, \chi_i, g_i; i\in J)$ be a local FL-free type.
 Let $I$ be the ideal of
$kQ^s(G, g_i,\chi_i;i\in J)$ generated by the following elements:
\begin {eqnarray}\label {quantume5}
\chi _j (\xi_i)E_iE_j-\chi _i ^{-1}(\xi_j)E_j E_i- \delta_{\sigma
(i), j}\ \frac{g_i^2-1}{ \chi _i (\xi_i) -\chi _i (\xi _i) ^{-1}} ,
\ \hbox {for }          i \in  J_u, j \in J_u', u\in N;
\end {eqnarray}
\begin {eqnarray}\label {quantume6}
\sum_{m=0}^{r_{ij}}(-1)^m \left[\begin{array}{c}
r_{ij}\\
m\\
\end{array}\right]_{\chi _i(\xi_i^{-1})}
E_i^{r_{ij}-m}E_jE_i^m, \end {eqnarray} $  \hbox {  for any } i, j
\in J_u \hbox { or } i, j \in J_{u}',   i \not= j \hbox { and }
r_{ij}-1 < ord (\chi _i(g_i)) ; u\in N.$

  Let  $U$  be the  algebra generated  by set $\{K_i,   X_i \mid  i \in
J\}$ with relations \begin {eqnarray}\label
{quantume2}\begin{array}{c} X_iX_j-X_j X_i= \delta_{ \sigma(i), j} \
\frac{K _i^2- K_i ^{-2}}{\chi
_i(\xi_i) -\chi _i (\xi_i) ^{-1}},   \hbox { \ for any  } i \in J_u, j\in J_u'; u\in N ;\\
 K_iK_{\sigma (i)}=K_{\sigma (i)}K_i=1
 \hbox { \ for any  } i \in J ^{(1)} ;\\
\chi _j (\xi_i)K_iX_j=X_jK_i, \ \ K_iK_j-K_jK_i=0, \  \hbox { for any  } i, j  \in J ;\\
 \sum_{m=0}^{r_{ij}}(-1)^m \left[\begin{array}{c}
r_{ij}\\
m\\
\end{array}\right]_{\chi _i(\xi_i^{-1})}
X_i^{r_{ij}-m}X_jX_i^m =0,  \\
\hbox { for } i, j \in J_{u}  \hbox { or } i, j \in J_{u}',   i
\not= j \hbox { and }  r_{ij}-1 < ord (\chi _i(g_i)) ; u\in N.
\end {array} \end {eqnarray}
The comultiplication, counit and antipode of $U$ are defined by
\begin {eqnarray}\label {e3.811}\begin{array}{lll}
\Delta(X_j)=X_j\otimes K_{\sigma (j)}+K_j\otimes X_j,&
S(X_j)=- \chi _j (\xi _j)X_j,&\epsilon(X_i)=0,\\
\Delta(K_i)=K_i\otimes K_i,& S(K_j)=K_{\sigma (j)},
&\varepsilon(K_i)=1, \\
\Delta(X_{\sigma (j)})=X_{\sigma (j)}\otimes K_{\sigma
(j)}+K_j\otimes X_{\sigma (j)}, & S(X_{\sigma (j)})=- \chi _{\sigma
(j)} (\xi _j)X_{\sigma (j)}
\end{array}\end {eqnarray}
for any $j \in J^{(1)}, i \in J .$

 In fact,  $K_{\sigma (j)}= K_j ^{-1}$ in
$U$ for any $j \in J^{(1)}.$

\begin{Theorem}\label{14.8} Under notation above, if ${\rm ESC}(G, g_i, \chi _i; i\in
J)$ is  a  local FL-free type, then  $kQ^s(G, g_i,\chi_i;i\in J)/ I
\cong U$ as Hopf algebras.

\end {Theorem}

{\bf Proof}. For any $i, j \in J_u, u \in N, $  see $\chi _{\sigma
(i)} (\xi _{\sigma (j)}) = \chi _i (\xi _j) = \chi _j (\xi_i) = \chi
_{\sigma (j)} (\xi _{\sigma (i)}),$ $\chi _{\sigma (i)}(\xi_j) =
\chi _i ^{-1}(\xi _j) = \chi _j ^{-1} (\xi_i) = \chi _j (\xi
_{\sigma (i)}),$ $\chi _{i} (\xi_ {\sigma (j)}) = \chi _i (\xi _j)
^{-1} = \chi _j (\xi_i) ^{-1} = \chi _{\sigma (j)} (\xi _i ).$
Therefore,
\begin {eqnarray}\label {quantume13}
\chi _i (\xi _j) = \chi _j (\xi _i)
\end {eqnarray} For any $i, j \in J_u\cup J_u', u \in N.$
Obviously, for $i, j \in J_u'$, $i \not= j,$ $u\in N$, (FL6) holds.

We show this theorem by following several steps.

(i) There is a  algebra homomorphism $\Phi$ from $kQ^s$ to $U$ such
that $\Phi(\xi _i)= K_i$, $\Phi( a_{hg_i, h}^{(i)})= \Phi (h)K_iX_i$
and $\Phi( a_{hg_{\sigma (i)}, h}^{(\sigma (i))})= \Phi
(h)K_iX_{\sigma (i)}$ for all $h \in G$ and $i \in J ^{(1)}$.
Indeed, define algebra homomorphism $\phi: kG\rightarrow U$ given by
$\phi(\xi_i)=K_i$ for $i\in J$ and a $k$-linear map $\psi:
kQ_1^c\rightarrow U$ by $\psi(a_{gg_j,g}^{(j)})=\phi(g)K_jX_j $ and
$\psi(a_{gg_{\sigma (j)},g}^{(\sigma (j))})=\phi(g)K_jX_ {\sigma
(j)} $  for any $j \in J^{(1)},$ $ g \in G.$ For any $g, h\in G,
j\in J ^{(1)}$, see
$$\begin{array}{rcl}
\psi(h\cdot
a^{(j)}_{gg_j,g})&=&\psi(a^{(j)}_{hgg_j,hg})=\phi(hg)K_jX_j
=\phi(h)\phi(g)K_jX_j=\phi(h)\psi(a^{(j)}_{gg_j,g}),\\
\end{array}$$
and
$$\begin{array}{rcl}
\psi(a^{(j)}_{gg_j,g}\cdot h)&=&\chi_j(h)\psi(a^{(j)}_{hgg_j,hg})=\chi_j(h)\phi(hg)K_jX_j\\
&=&\phi(g)K_jX_j\phi(h) \ \ (\hbox { since } X_j \phi (h) = \phi (h) \chi _j (h) X_j)\\
&=&\psi(a^{(j)}_{gg_j,g})\phi(h).
\end{array}$$
Similarly, $\psi(h\cdot a^{(\sigma (j))}_{gg_{\sigma (j)},g}) =
\phi(h)\psi(a^{(\sigma (j))}_{gg_{\sigma (j)},g}) $ and $
\psi(a^{(\sigma (j))}_{gg_{\sigma (j)},g}\cdot h)= \psi(a^{(\sigma
(j))}_{gg_{\sigma (j)},g})\phi(h).$
 This implies that $\psi$ is a  $kG$-bimodule map from
$(kQ_1^c,g_i,\chi_i;i\in J)$ to $_{\phi}U_{\phi}$. Using the
universal property of tensor algebra over $kG$, we complete the
proof.

(ii) $\Phi (I)=0$. For any $i \in J_u, j \in J_{u}'$, see that
$$\begin{array}{cl}
&\Phi(\chi _j (\xi_i)E_iE_j-\chi _i (\xi_j)^{-1}E_j E_i-
\delta_{\sigma (i), j} \ \frac{g_i^2-1}{ \chi _i (\xi _i) -\chi _i
(\xi _i)
^{-1}})\\
 =& \chi _j (\xi_i) K_i X_iK_j^{-1} X_j-\chi _i (\xi_j)^{-1}
 K_j^{-1} X_jK_i X_i-
\delta_{\sigma (i),
j} \ \frac{K_i^4-1}{ \chi _i (\xi) -\chi _i (\xi _i) ^{-1}}\\
=& K_iK_j^{-1}X_iX_j-K_j^{-1}K_iX_jX_i- \delta_{\sigma (i),
j}\ \frac{K_i^4-1}{ \chi _i (\xi) -\chi _i (\xi _i) ^{-1}}\\
\\=&K_iK_j^{-1}(X_iX_j-X_jX_i-\delta_{\sigma (i),
j}
  \ \frac {K_i^3K_j-K_i^{-1}K_j}{ \chi _i (\xi _i) -\chi _i (\xi _i)^{-1} }) \\
=&K_iK_j^{-1}(X_iX_j-X_jX_i- \delta_{\sigma (i), j} \
\frac{K_i^2-K_i^{-2}}{ \chi _i (\xi _i) -\chi _i (\xi _i)^{-1} })=0.\\
\end{array}$$

For $i, j\in J_u,$ $i\not= j,$ see that
$$\begin{array}{cl}
&\Phi(E_i^{r_{ij}-m}E_jE_i^m)\\
=&(K_iX_i)^{r_{ij}-m}K_jX_j(K_iX_i)^m\\
=& \chi _i (\xi _i) ^{ \frac {1} {2} ( (r_{ij} -m) (r_{ij} -m -1) +
m (m-1) + 2 (r_{ij} -m) m)} \chi _i (\xi _j) ^{r_{ij} -m} \chi _j
(\xi _i)
^m  X_i^{r_{ij}-m}X_jX_i^m\\
=& \chi _i (\xi _i) ^{\frac {1}{2} (r_{ij}-1)r_{ij}} \chi _i (\xi
_j) ^{r_{ij}}X_i^{r_{ij}-m}X_jX_i^m \ \  ((\hbox {by }
(FL6))\\
\end{array}$$
and
$$\Phi(\sum_{m=0}^{r_{ij}}(-1)^m \left[\begin{array}{c}
r_{ij}\\
m\\
\end{array}\right]_{\chi _i(\xi_i^{-1})}
E_i^{r_{ij}-m}E_jE_i^m)=0.$$ Similarly, the equation above holds for
$i, j \in J_u'$, $i \not= j,$ $u\in N.$

By (ii), there exists an algebra homomorphism $\bar \Phi$ from
$kQ^s/I$ to $U$ such that $\bar \Phi (x +I) = \Phi (x)$ for any
$x\in kQ^s.$ For convenience, we will still use $x$ to denote $x
+I$.

(iii) It follows from the definition of $U$ that there exists a
unique algebra map $\Psi: U\rightarrow kQ^s/I$ such that
$\Psi(K_i)=\xi_i$, $\Psi(X_i)=\xi_i^{-1}E_i$ and $\Psi(X_{\sigma
(i)})=\xi_i^{-1}E_{\sigma (i)}$  for all $i \in J^{(1)}$. It is easy
to see that $\Psi\overline{\Phi}= id$ and $\bar \Phi \Psi =id$.

(iv) We  show that $I$ is a Hopf ideal of $kQ^s.$ It follows from
Theorem \ref {14.7} that $R:=diag (kQ^s)$ is a braided Hopf algebra
with $\delta ^- (E_i) =g_i \otimes E_i $, $h \rhd E_i = \chi _i
^{-1} (h) E_i$ for any $h \in G$, $i \in J.$ By Theorem \ref {14.6},
$R \#B\cong kQ^s$ as Hopf algebras. It follows from Lemma \ref
{14.3.7'} that $I$ is a Hopf ideal of $kQ^s$.

 Obviously, $\Psi $ preserve the comultiplication and counit.
Since $kQ^s/I$ is a Hopf algebra, then  $\Psi$ is a Hopf algebra
isomorphism. \ \ $\Box$

\begin{Corollary}\label{14.3.14} The quantum enveloping algebra of a complex  semisimple Lie
algebra is isomorphic to a quotient of a semi-path Hopf algebra as
Hopf algebras.

\end {Corollary}
{\bf Proof.} Let  $k$ be the complex field and $L$ a complex
semisimple Lie algebra determined by $A = (a_{ij})_{n\times n}$. So
$A$ is a symmetrizable Cartan matrix with $d_i \in \{1, 2, 3\}$ such
that $d_ia_{ij} = d_j a_{ji}$ for any $i, j \in J ^{(1)} = \{1, 2,
\cdots, n\}$. Let $\mid \! N \! \mid =1, $  $J^{(2)} = \{n+1, n+2,
\cdots, n+n\}$, $J = J^{(1)}\cup J^{(2)}$ and $\sigma : J^{(1)}
\rightarrow J^{(2)}$ by sending $i$ to $i+n$. Let $G$ be a free
commutative group generated by  generator set $\{\xi_i \mid i \in
J^{(1)}\}$. Set $\xi _{\sigma (i)} = \xi _i^{-1}$,  $g_i= g_{\sigma
(i)} := \xi _i ^2$, $r_{ij}= r_{\sigma (i), \sigma (j)} = 1-a_{ij}$
for any $i, j \in J^{(1)}$, $i \not= j$. Define $\chi _i (\xi _j) =
q ^{-d_ia_{ij}}$ and $\chi _{\sigma (i)} (\xi _j) = \chi _i
^{-1}(\xi_j)$ for any $i, j \in J^{(1)}$, where $q$ is not a root of
1 with $0\not=q\in k$. It is easy to check (FL1)-- (FL7) hold, i.e.
${\rm ESC} (G, g_i, \chi_i; i\in J)$ is an FL-quantum group type.
 $U$ in Theorem \ref {14.8}
exactly  is the quantum enveloping algebra $U_q(L)$ of $L$ (see
\cite[p.218]{Mo93} or  \cite {Lu93}).  Therefore the conclusion
follows from Theorem \ref {14.8}. $\Box$

In fact, for any a generalized Cartan matrix $A$, we can obtain an
FL-quantum group type ${\rm ESC} (G, g_i, \chi_i; i\in J)$ as in the
proof above. By the way, if ${\rm ESC} (G, g_i, \chi_i; i\in J)$ is
a Local FL-type, then there exist the  Hopf ideals, which generated
by (\ref {quantume5}) and (\ref {quantume6}),  in co-path Hopf
algebra $kQ^c (G, g_i, \chi_i; i\in J)$ and multiple Taft algebra
$kQ^c (G, g_i, \chi_i; i\in J)$, respectively.

\section{Classification of ramification systems for the symmetric group}

Unless specified otherwise,  in the section we have the following
assumptions and notations. $n$ is a positive  integer with $n
\not=6$; $F$ is a field containing a primitive $n$-th root of 1; $G=
S_n$,  the symmetric group; Aut$G$ and Inn$G$ denote the
automorphism group and inner automorphism group,  respectively;
$S_I$ denotes the group of full permutations of set $I$; $Z(G)$
denotes the center of $G$ and $Z_x$ the centralizer of $x$; $1$
denotes the unity element of $G;$ $\widehat{G}$ denotes the set of
characters of all one-dimensional representations of $G$. $\phi_h$
denotes the inner automorphism induced by $h$ given by
$\phi_h(x)=hxh^{-1}, $ for any $x\in G$. $C_n$ is cyclic group of
order $n$; $D^B$ denotes the cartesian product $\prod \limits_ {i\in
B}D_i$,  where $D_i = D$ for any $i \in B$.

\begin{Proposition}\label{14.7.2}
{(i)} There are
$\frac{n!}{\lambda_1!\lambda_2!\cdots\lambda_n!1^{\lambda_1}2^{\lambda_2}\cdots
n^{\lambda_n}}$ permutations of type
$1^{\lambda_1}2^{\lambda_2}\cdots n^{\lambda_n}$ in $S_n$;

{(ii)} If $n\neq 6$,   then Aut$(S_n)$=Inn$(S_n)$; If $n\neq 2, 6$,
then Inn$(S_n) \cong S_n$;

{(iii)} $S_n' = A_n$,  where $A_n$ is the alternating group. In
addition,
$$ S_n/A_n\cong\left\{
\begin{array}{lll}
C_2& n\geq2\\
C_1 &n =1
\end{array}
\right. . $$

\end{Proposition}

{\bf Proof. } (i) It follows from  \cite [ Exercise  2.7.7] {Hu}.

(ii) By  \cite [ Theorem 1.12.7] {Zh82},  $\mathrm{Aut}(A_n) \cong
S_n$ when $n >3$ and  $n \not=6$.  By  \cite [Theorem 1.6.10]
{Hu74},  $A_n $ is a non-commutative simple group when  $n \not= 4$
and $n>1$. Thus it follows from  \cite [Theorem 1.12.6 ] {Zh82} that
$\mathrm{Aut} (S_n) = \mathrm{Inn} (S_n)\cong S_n$ when  $n \not= 6$
and $n>4$. Applying \cite [Example in Page 100] {Zh82}, we obtain
 $\mathrm{Aut} (S_n) = \mathrm{Inn} (S_n)\cong S_n$ when  $n=3,  4$.
Obviously,  $\mathrm{Aut} (S_n) = \mathrm{Inn} (S_n)$ when  $n=1,
2$.

 (iii) It follows from  \cite [Theorem 2.3.10] {Xu95}. $\Box$

\begin{Definition}\label{14.7.4}
A group $G$ is said  to act on a non-empty set $\Omega$,  if there
is a map  $G\times \Omega\rightarrow\Omega$,  denoted by $(g,
x)\mapsto g\circ x$,  such that for all $x\in \Omega$ and $g_1,
g_2\in G$:
$$
(g_1g_2)\circ x=g_1\circ(g_2\circ x)\qquad and\qquad   1\circ x=x,
$$ where 1 denote the unity element of $G.$

\end{Definition}

For each $x\in \Omega$,  let
$$ G_x:=\{g\circ x \mid  g\in G\}, $$
 called the orbit of $G$ on $\Omega$. For each $g\in G$,  let
 $$ F_g:=\{x\in \Omega \mid  g\circ x=x\}, $$
 called the  fixed point set of  $g$.

Burnside's lemma,  sometimes also called Burnside's counting
theorem,  which is  useful to compute the number of orbits.

%\subsection{Burnside 引理}

\begin{Theorem}\label{14.7.7} (See \cite [Theorem 2.9.3.1]{Hu})
( Burnside's Lemma ) Let $G$ be a finite group that acts on  a
finite set $\Omega$,  then the number $\mathcal {N}$ of orbits is
given by the following formula:
$$
\mathcal {N}=\frac{1}{ \mid  G \mid  }\sum\limits_{g\in G} \mid F_g
\mid  . $$

\end{Theorem}

%\subsection{特征标}

\begin{Definition}\label{14.7.8}

 Let $G$ be a  group. A character of a  one-dimensional representation of $G$ is  a homomorphism from $G$ to $F-\{0\}$,  that
is ,  a map $\chi:G\rightarrow F$ such that
$\chi(xy)=\chi(x)\chi(y)$ \ and\ $\chi(1)=1$ for any $ x, y\in G$.
The set of all characters of all one-dimensional representations of
$G$ is denoted by $\widehat{G}$. \end{Definition} For convenience,
a character of a  one-dimensional representation of $G$ is called a
character of $G$ in short. \cite[Page 36, example 4]{char} indicated
that there is a one-to-one correspondence between characters  of $G$
and characters  of $G/G'$.

\begin{Proposition}\label{14.7.9}

 Let $G$ be a finite abelian group  and $F$ contains $\mid \! G \! \mid $th primitive  root of 1,  then
there are exactly $\mid \! G \! \mid $  characters of
 $G$,  i.e. $ \mid \widehat{G}
\mid  = \mid  G \mid  . $

\end{Proposition}

\begin{Corollary}\label{14.7.10}

 Let $G$ be a finite  group  and $F$ contains $\mid \! G \! \mid $th primitive  root of 1,  then
there are exactly $\mid \! G/G' \! \mid $  characters of $G$, i.e. $
\mid \widehat{G} \mid  = \mid  G/G' \mid  . $

\end{Corollary}

\subsection{Ramification Systems and Isomorphisms}

We begin with defining  ramification and ramification system.
\begin{Definition}\label{14.7.11}
Let $G$ be a finite  group and $\mathcal {K}(G)$ the set of all
 conjugate classes of $G$. $r$ is called a ramification data or ramification of $G$ if $r=\sum\limits_{C\in \mathcal
{K}(G)}r_{_C}C$,  where  $r_C$ is a non-negative integer for any $C$
in $\mathcal {K}(G)$. For convenience,  we choice a set $I_C(r)$
such that  $r_{_C}= \mid  I_C(r) \mid  $ for any  $C$ in $\mathcal
{K}(G)$.

Let $\mathcal {K}_r(G)=\{C\in \mathcal {K}(G) \mid  r_{_C}\neq
0\}=\{C\in \mathcal {K}(G) \mid  I_C(r)\neq \emptyset\}$.

\end{Definition}

\begin{Definition}\label{14.7.12}
$(G, r, \overrightarrow{\chi}, u)$ is called a ramification system
 with characters (or RSC in short), if r is a ramification of G, u is a
 map from $\mathcal{K}(G)$ to G with u(C)$\in$ C for any
 C$\in\mathcal{K}(G)$, and $\overrightarrow{\chi}=\{\chi_C^{(i)}\}_{i\in I_C(r), C\in
\mathcal {K}_r(G)}\in\prod\limits_{C\in \mathcal
{K}_r(G)}(\widehat{Z_{u(C)}})^{r_{_C}}$.

\end{Definition}

\begin{Definition}\label{14.7.13}

 $RSC(G, r, \overrightarrow{\chi}, u)$ and $
RSC(G', r', \overrightarrow{\chi'}, u')$ are said to be
isomorphic(denoted $RSC(G, r, \overrightarrow{\chi}, u)\cong RSC(G',
r', \overrightarrow{\chi'}, u')$) if the following conditions are
satisfied:

{(i)} There exists a group isomorphism $\phi:G\rightarrow G';$

{(ii)} For any $ C\in \mathcal {K}(G), $there exists an element
$h_C\in G$ such that $\phi(h_Cu(C)h_C^{-1})=u'(\phi(C))$;

{(iii)} For any $ C\in\mathcal {K}_r(G), $there exists a bijective
map $\phi_C:I_C(r)\rightarrow I_{\phi(C)}(r')$ such that
$\chi_{\phi(C)}^{\prime(\phi_C(i))}(\phi(h_Chh_C^{-1}))=\chi_C^{(i)}(h),
$ for all $h\in Z_{u(C)}, $ for all $i\in I_C(r)$.
\end {Definition}

{\bf Remark.} Obviously,  if $u(C)=u'(C)$ for any $ C\in \mathcal
{K}_r(G), $  then $RSC(G, r, \overrightarrow{\chi}, u)\cong RSC(G,
r, \overrightarrow{\chi}, u')$. In fact,  let $\phi=id_G$ and
$\phi_C=id_{I_C(r)}$ for any $ C\in \mathcal {K}_r(G)$. It is
straightforward to check the conditions in Definition\ref{14.7.13}
hold.

\begin{Lemma}\label{14.8.4}

Given an $RSC(G, r, \overrightarrow{\chi}, u)$.If $RSC(G', r',
\overrightarrow{\chi'}, u')\cong RSC(G, r, \overrightarrow{\chi},
u)$, then there exists a group isomorphism $\phi:G\rightarrow G'$
such that

{(i)} $ r'_{\phi(C)}=r_{_C}$ for any $ C\in \mathcal {K}(G)$;

{(ii)}  $\mathcal {K}_{r'}(G')=\phi(\mathcal {K}_{r}(G))=\{\phi(C)
\mid C\in \mathcal {K}_r(G)\}.$

\end{Lemma}
{\bf Proof.} (i) By Definition \ref{14.7.13},  there exists a group
isomorphism $\phi:G\rightarrow G'$ such that $r'_{\phi(C)}=r_{_C}$
for any  $ C\in \mathcal {K}_r(G)$. If $C\notin
 \mathcal {K}_r(G)$,  i.e. $r_{_C}=0, $ then $r'_{\phi(C)}=r_{_C}=0$;

(ii) is an immediate consequence of (i).
 $\Box$

\begin{Proposition}\label{14.8.5}
Given a $RSC(G, r, \overrightarrow{\chi}, u)$,  then $ RSC(G,  r',
\overrightarrow{\chi'},  u')\cong RSC(G,  r, \overrightarrow{\chi},
u)$ if and only if there exists a group isomorphism $\phi \in \hbox
{Aut} G$ such that:

{(i)} $r'=\sum\limits_{C\in \mathcal {K}(G)}r_{\phi^{-1}(C)}C$;

{(ii)}
$\chi_{\phi(C)}^{\prime(i)}=\chi_C^{(\phi_C^{-1}(i))}(\phi\phi_{h_C})^{-1}$
and $u'(\phi(C))=\phi\phi_{h_C}(u(C))$ for any $ C\in \mathcal
{K}_r(G)$ ,  $i\in I_C(r) $,  $\phi\in \mathrm{Aut}G$ and $\phi_C\in
S_{r_{_C}}$,  $h_C\in G$.

\end{Proposition}
{\bf Proof.} If $RSC(G, r', \overrightarrow{\chi'}, u')\cong RSC(G,
r, \overrightarrow{\chi}, u) $,  then,  by Lemma \ref{14.8.4}(i),
there exists $\phi\in \mathrm{Aut}G$ such that
$$r'=\sum\limits_{C\in
\mathcal
 {K}(G)}r'_CC=\sum\limits_{C\in\mathcal
 {K}(G)}r'_{\phi(C)}\phi(C)=\sum\limits_{C\in\mathcal
 {K}(G)}r_{_C}\phi(C)=\sum\limits_{C\in\mathcal {K}(G)}r_{\phi^{-1}(C)}C.
$$
This shows Part (i).

By Definition \ref{14.7.13} (ii),  there exists $h_C\in G$ such that
$u'(\phi(C))=\phi(h_Cu(C)h_C^{-1})=\phi\phi_{h_C}(u(C))$ for any $
C\in \mathcal {K}_r(G)$; By Definition \ref{14.7.13} (iii),  there
exists $\phi_C\in S_{I_C(r)}=S_{r_{_C}}$ such that
$\chi_{\phi(C)}^{\prime(\phi_C(i))}(\phi(h_Chh_C^{-1}))
=\chi_{\phi(C)}^{\prime(\phi_C(i))}\phi\phi_{h_C}(h)=\chi_C^{(i)}(h)$
for any $h\in Z_{u(C)}$. Consequently,
$$\chi_{\phi(C)}^{\prime(\phi_C(i))}=\chi_C^{(i)}(\phi\phi_{h_C})^{-1}, $$ i.e. $$\chi_{\phi(C)}^{\prime
(i)}=\chi_C^{(\phi_C^{-1}(i))}(\phi\phi_{h_C})^{-1}, $$ where $C\in
\mathcal {K}_r(G)$. That is,  (ii) holds.

Conversely,  if (i) and (ii) hold,  let
$$r'=\sum\limits_{C\in \mathcal {K}(G)}r_{\phi^{-1}(C)}C$$
and $$I_{\phi(C)}(r')=I_C(r)$$ for any $C\in\mathcal {K}(G)$,  $
\phi\in$Aut$G$. Obviously,  $\mathcal {K}_{r'}(G)=\{\phi(C) \mid
C\in \mathcal {K}_r(G)\}.$ For any $C\in \mathcal {K}_r(G), h_C\in
G,  \phi_{_C}\in S_{I_C(r)}=S_{r_{_C}}$,  let
$$u'(\phi(C))=\phi(h_Cu(C)h_C^{-1})$$ and
$$\chi_{\phi(C)}^{\prime
(i)}=\chi_C^{(\phi_C^{-1}(i))}(\phi\phi_{h_C})^{-1}$$ for any $i\in
I_{\phi(C)}(r').$ It is easy to check $RSC(G, r',
\overrightarrow{\chi'}, u')\cong RSC(G, r, \overrightarrow{\chi},
u)$.$\Box$

\begin{Corollary}\label{14.8.6}
If $G=S_n$ with $n\neq 6$,  then $RSC(G,  r, \overrightarrow{\chi},
u)\cong RSC(G,  r', \overrightarrow{\chi'},  u') $ if and only if

(i) $r'=r$;

(ii)
$\chi_C^{\prime(i)}=\chi_C^{(\phi_C^{-1}(i))}\phi^{-1}_{g_{_C}}$ and
$u'(C)=\phi_{g_{_C}}(u(C))$ \ for any  $ C\in \mathcal {K}_r(G)$,  $
i\in I_C(r) $,   $\phi\in \mathrm{Aut}G$ and $g_{_C}\in G$.
\end{Corollary}

{\bf Proof.} By Proposition \ref{14.7.2}(ii),
$\mathrm{Aut}G=\mathrm{Inn}G$ and $\phi(C)=C$ for any $ \phi
\in$Aut$G,  C\in \mathcal {K}(G)$. There exists $g_{_C}\in G$ such
that  $\phi\phi_{h_C}=\phi_{g_{_C}}$ since $\phi\phi_{h_C}\in
\mathrm{Aut}G=\mathrm{Inn} G$ for any $h_C \in G,   \phi \in$Aut$G,
C\in \mathcal {K}(G)$.  Considering Proposition \ref{14.8.5},  we
complete the proof.$\Box$

From now on,   we suppose that $G=S_n(n\neq 6)$. For a given
ramification $r$ of $G$,  let  $ \Omega (G,  r)$ be the set of all
RSC's of G with the ramification $r$,  namely,  $ \Omega (G,  r) : =
\{ (G,  r,  \overrightarrow{ \chi },  u) \mid (G,   r,
\overrightarrow{ \chi },  u) \mbox {\ is\ an  } RSC\}$. Let ${\cal
N}(G,  r)$ be the number of isomorphic classes in $ \Omega (G, r)$.
This article is mostly devoted to investigate the formula of ${\cal
N}(G,  r)$.

\subsection{Action of Group }

Denote $RSC(G,  r,  \overrightarrow{\chi},  u)$ by
$\left(\left\{\chi^{(i)}_{C}\right\}_i,
u(C)\right)_{C\in\mathcal{K}_r(G)}$ for short. Let
$M:=\prod\limits_{C\in \mathcal {K}_r(G)}(\mbox{Aut}G\times
S_{I_C(r)})$ $=\prod\limits_{C\in \mathcal
{K}_r(G)}(\mbox{Inn}G\times S_{r_{_C}})$. Define the action of $M$
on $\Omega(G,  r)$ as follows:
\begin{eqnarray} \label {e14.8.1}
(\phi_{g_{_C}}, \phi_C)_C\circ\left(\left\{\chi^{(i)}_{C}\right\}_i,
u(C)\right)_C &=&
\left(\left\{\chi_C^{(\phi_C^{-1}(i))}\phi^{-1}_{g_{_C}}\right\}_i,
\phi_{g_{_C}}u(C) \right )_C
\end{eqnarray}
 In term of Corollary  \ref{14.8.6},  each orbit of $\Omega(G,
r)$ represents an isomorphic class of RSC's. As a result,  ${\cal
N}(G,  r)$ is equal to the number of orbits in $\Omega(G,  r )$.

We compute the the number of orbits by following steps.

\subsection{ The structure of centralizer in $S_n$ }

We recall the semidirect product and the wreath product of groups
(see \cite[Page 21]{char},  \cite [Page 268-272]{suzuki}, ).
\begin{Definition}\label{14.8.7}

Let $N$ and $K$ be groups, and assume that there exists a group
homomorphism $\alpha:K\rightarrow \mathrm{Aut}N$. The semidirect
product $N\rtimes_{\alpha} K$ of $N$ and $K$ with respect to
$\alpha$ is defined as follows:

(1) As set,  $N\rtimes_{\alpha} K$ is the Cartesian product of N and
K;

(2) The multiplication is given by
$$
(a, x)(b, y)=(a\alpha(x)(b), xy)
$$
for any a $\in$ N and x $\in$ N.
\end{Definition}

{\bf Remark.} Let $L=N\rtimes_{\alpha} K$. If we set
$$ \bar{N}=\{(a, 1) \mid a\in N\}, \bar{K}=\{(1, x) \mid x\in
K)\}, $$ then $ \bar{N}$ and $\bar{K}$ are subgroups of $L$ which
are isomorphic to $N$ and $K$,  respectively,  and will be
identified with $N$ and $K$.

\begin{Definition}\label{14.8.10}
 Let $A$ and $H$ be two  groups  and $H$ act on the set $X$. Assume that
 $B$ is the set  of all maps from $X$ to $A$. Define the multiplication of $B$ and
 the group homomorphism $\alpha $ from  $H$  to  Aut $B$ as follows:
 $$
(bb')(x)=b(x)b'(x) \hbox { \ and } \alpha (h)(b)(x)=b(h^{-1}\cdot
x)$$ for any $x\in X,  h\in H,  b,  b' \in B.$

 The  semidirect product
$B\rtimes_{\alpha} H$ is called the (general) wreath product of $A$
and $H$,  written as  $A$\:wr\:$H$.

\end{Definition}

{\bf Remark.}\cite[Page 272]{suzuki} indicated that in Definition
\ref {14.8.10} $B\cong \prod\limits_{x\in X}A_x,  $ where $A_x= A$
for any $x\in X$. Hence $W=A$\:wr\:$H=B\rtimes H=(\prod\limits_{x\in
X}A_x)\rtimes H$. In particular. the wreath product  plays  an
important  role in the structure of centralizer in $S_n$ when $H=
S_n$ and $X = \{1,  2,  \cdots,  n\}$.

Now we keep on the work in \cite[Page 295-299 ]{suzuki}.
 Let $Y_i $ be the set of all letters which
belong to those cycles of length $i$ in the independent  cycle
decomposition of $\sigma$. Clearly,  $Y_{i}\bigcap Y_{j}=\emptyset$
for $i\neq j$ and $\bigcup\limits_iY_i=\{1, 2, \cdots, n \}$.

\begin{Lemma}\label{14.8.11}
 $\rho \in Z_{\sigma }$ if and only
if  $Y_{i}$ is $\rho$-invariant,  namely $\rho(Y_i)\subset Y_i$, and
 the restriction $\rho_{i}$ of $\rho$ on $Y_{i}$ commutes with the
restriction $\sigma_{i}$ of $\sigma$ on $Y_{i}$ for  $i = 1,  2,
\cdots,  n$.
\end{Lemma}

{\bf Proof.}It follows from \cite[Page 295]{suzuki}.$\Box$

\begin{Proposition}\label{14.8.12}
       If $\sigma$ and  $ \sigma_i$ is the same as above,  then

(i) $Z_{\sigma}= \prod\limits_{i}Z_{\sigma_i};$

(ii) $Z_{\sigma_i}\cong C_i$\:wr\:$S_{\lambda_i }\cong
\left(C_i\right)^{\lambda_i}\rtimes S_{\lambda_i} ,   $ where
$Z_{\sigma_i}$ is the centralizer of $\sigma_i$ in $S_{Y_i}$;
$\left(C_i\right)^{\lambda_i}=\overbrace{C_i\times \cdots\times
C_i}^{\lambda_i\mbox{个}}$; define $Z_{\sigma_i}=1$ and $
S_{\lambda_i}=1 $ when $\lambda_i=0$.
\end{Proposition}

{\bf Proof. }(i) It follows from Lemma  \ref {14.8.11}.

 (ii) Assume $$\tau=(a_{10}a_{11}\cdots a_{1,  l-1})(a_{20}\cdots a_{2,
l-1})\cdots (a_{m0}\cdots a_{m,  l-1}) $$ and
$$
X=\{1,  2,  \cdots,  m\},  \ \  Y=\{a_{ij} \mid  1 \le i \le m,  0
\le j \le l-1 \}.$$ Obviously,  in order to prove (ii),  it suffices
to show that $Z_{\tau}\cong C_l$\:wr\:$S_m$,  where $Z_{\tau}$ is
the centralizer of $\tau$ in $S_Y$. The second index of  $a_{ij}$
ranges over $\{0, 1, \cdots, l-1\}$. It is convenient to identify
this set with $C_l=\mathbb{Z}/(l)$,  and to use notation such as
$a_{i, l+j}=a_{ij}$. We will define an isomorphism $\varphi$ from
$Z_{\tau}$ onto $C_l$\:wr\:$S_m$. For any element $\rho$ of
$Z_{\sigma}$,  we will define a permutation $\theta(\rho)\in S_m$
and a function $f(\rho)$ from $X$ into $C_l$ by
$$\rho^{-1}(a_{i0})=a_{jk}, \qquad j=\theta(\rho)^{-1}(i), \qquad
f_{\rho}(i)=g^k, $$ where $g$ is the generator of $C_l$,  $1\le i
\le m$. Let
$$\varphi(\rho)=(f_{\rho},  \theta(\rho)). $$
Then $\varphi$ maps $Z_{\sigma}$ into $C_l$\:wr\:$S_m$. First we
 show that $\varphi$ is surjective. In fact,  for any $(f,
\theta) \in C_l\ \mathrm{wr}\ S_m$ and $i\in X$,  if $f (i) = g^k, $
$\theta (j) = i$,  then there exists a permutation $\rho$ such that
$\rho ^{-1} (a_{ir})= a _{j,  k+r}$. It is easy to verify that
$\rho\in Z_{\sigma}$ and $\varphi (\rho ) = (f,  \theta)$. Next we
show that $\varphi$ is injective. If $\rho\in$ Ker$\varphi$,  then
$\theta(\rho)=1$ and $f _{\rho}(i)=g^0$,  this means that
$\rho(a_{i0})=a_{i0}$ for all $i$. Hence, $\rho(a_{ij})=a_{ij}$ for
all $i$ and $j$ and $\rho=1_Y$,  which implies that $\varphi$ is
bijective. We complete the proof if $\varphi$ is homomorphism,
which is proved below. For any $\rho, \pi\in Z_{\tau}$,  suppose
that
$$
(\rho\pi)^{-1}(a_{i0})=\pi^{-1}(a_{jk})=a_{s,  t+k},
$$
where $j=\theta(\rho)^{-1}(i), g^k=f_{\rho}(i),
s=\theta(\pi)^{-1}(j)=\theta(\pi)^{-1}\theta(\rho)^{-1}(i),
g^t=f_{\pi}(j)=f_{\pi}(\theta(\rho)^{-1}(i))$. Consequently,
$$\theta(\rho\pi)^{-1}(i)=s=(\theta(\rho)\theta(\pi))^{-1}(i),  \qquad
f_{\rho\pi}(i)=g^{t+k}=f_{\rho}(i)f_{\pi}(j)=f_{\rho}(i)f_{\pi}(\theta(\rho)^{-1}(i)).
$$

On the other hand,  in view of  Definition \ref {14.8.10},
$\theta(\rho)(f_{\pi})(i)=f_{\pi}(\theta(\rho)^{-1}(i))=f_{\rho\pi}(i)$,
whence
$$(f_{\rho},  \theta(\rho))(f_{\pi},
\theta(\pi))=(f_{\rho}\theta(\rho)(f_{\pi}),
\theta(\rho)\theta(\pi))=(f_{\rho\pi}, \theta(\rho\pi)).$$ It
follows that
$$\varphi(\rho\pi)=\varphi(\rho)\varphi(\pi), $$
which is just what we need.\ $\Box$

\begin{Corollary} \label{14.7.15}

{{ (i)}
$|Z_{\sigma}|=\lambda_1!\lambda_2!\cdots\lambda_n!1^{\lambda_1}2^{\lambda_2}\cdots
n^{\lambda_n}, $ where $\sigma$ is the same as above;

{(ii)} $\sum\limits_{h\in C}|Z_h|=|G|=n!$}  for any $ C\in \mathcal
{K}(G), $

\end{Corollary}

{\bf Proof. }  (i) It follows from \ref{14.8.12} that
$$
 \mid  Z_{\sigma} \mid  =\prod\limits_{i} \mid  Z_{\sigma_i} \mid
=\prod\limits_{i} \mid  \left(C_i\right)^{\lambda_i}\rtimes
S_{\lambda_i} \mid =\prod\limits_{i} \mid
\left(C_i\right)^{\lambda_i} \mid   \mid  S_{\lambda_i} \mid
=\prod\limits_{i}(i^{\lambda_i}\lambda_i!).
$$

(ii) Assume the type of $C$ is $1^{\lambda_1}2^{\lambda_2}\cdots
n^{\lambda_n}$. Combining (i) and Proposition \ref{14.7.2}(i),  we
have
$$
\sum\limits_{h\in C} \mid Z_h \mid
=(\lambda_1!\lambda_2!\cdots\lambda_n!1^{\lambda_1}2^{\lambda_2}\cdots
n^{\lambda_n})\frac{n!}{\lambda_1!\lambda_2!\cdots\lambda_n!1^{\lambda_1}2^{\lambda_2}\cdots
n^{\lambda_n}}=n! \ . \ \ \Box
$$

Proposition \ref{14.8.12} shows that $Z_{\sigma}$ is a direct
product of some wreath products  such as
 $C_l\:$wr\:$S_m$.
Denote by $((b_i)_i, \sigma)$ the element of
$C_l\:$wr\:$S_m=(C_l)^{m}\rtimes S_m$ where $b_i\in C_l, \sigma\in
S_m$. By Definition \ref{14.8.7} and \ref{14.8.10}, we have
\begin{eqnarray*}
((b_i)_i, \sigma)((b'_i)_i, \sigma')&=&((b_i)_i\cdot\sigma((b'_i)_i), \sigma\sigma')\nonumber\\
&=&((b_i)_i\cdot(b'_{\sigma^{-1}(i)})_i, \sigma\sigma')\nonumber\\
&=&((b_ib'_{\sigma^{-1}(i)})_i, \sigma\sigma')
\end{eqnarray*}
and
\begin{equation*}
((b_i)_i, \sigma)^{-1}=((b^{-1}_{\sigma(i)})_i, \sigma^{-1}).
\end{equation*}

\subsection { Characters of $Z_{\sigma}$}

\begin{Lemma}\label{14.8.14} If
$H=H_1\times\cdots\times H_r$,  then

(i) $H'=H'_1\times\cdots\times H'_r$;

(ii) $H/H'\cong \prod\limits_{i=1}^{r}H_i/H'_i.$
\end{Lemma}
{\bf Proof.} (i) It follows from  \cite  [ Page 77] {Zh82};

(ii) It follows from \cite[Corollary 8.11]{Hu74} and (i). \ $\Box$

According to Proposition \ref{14.8.12}(i) and Lemma
\ref{14.8.14}(ii), to investigate  the structure of
$Z_{\sigma_i}/Z'_{\sigma_i}$,  we have to study
$(C_l\:$wr\:$S_m)/(C_l\:$wr\:$S_m)'$.

\begin{Lemma}\label{14.8.15}
Let $B=(C_l)^{m}$,   $\bar{B}=\left\{(b_i)_i\in B \mid
\prod\limits_{i=1}^mb_i=1\right\} $ and $W=C_l\:$wr\:$S_m=B\rtimes
S_m$. Then $W'=\bar{B}\rtimes S'_m=\bar{B}\rtimes A_m$.
\end{Lemma}

{\bf Proof. } Obviously,  $\bar{B}\lhd B$.

For any  $((b_i)_i,  \sigma),  ((b'_i)_i,  \sigma')\in W, $ see
\begin{eqnarray*}
& &((b_i)_i,  \sigma)^{-1}((b'_i)_i,  \sigma')^{-1}((b_i)_i,  \sigma)((b'_i)_i,  \sigma')\\
& = &((b^{-1}_{\sigma(i)})_i,  \sigma^{-1})((b^{\prime-1}_{\sigma'(i)})_i,  \sigma^{\prime-1})((b_ib'_{\sigma^{-1}(i)})_i,  \sigma\sigma')\\
& = &((b^{-1}_{\sigma(i)}b^{\prime-1}_{\sigma'\sigma(i)})_i,  \sigma^{-1}\sigma^{\prime -1})((b_ib'_{\sigma^{-1}(i)})_i,  \sigma\sigma')\\
& =
&((b^{-1}_{\sigma(i)}b^{\prime-1}_{\sigma'\sigma(i)}b_{\sigma'\sigma(i)}b'_{\sigma^{-1}\sigma'\sigma(i)})_i,
\sigma^{-1}\sigma^{\prime-1}\sigma\sigma').
\end{eqnarray*}
Considering $\prod\limits_ib_{\tau(i)}=\prod\limits_ib_i$
 for any $ \tau \in S_m$,  we have
\begin{eqnarray*}
&&\prod\limits_i\left(b^{-1}_{\sigma(i)}b^{\prime-1}_{\sigma'\sigma(i)}b_{\sigma'\sigma(i)}b'_{\sigma^{-1}\sigma'\sigma(i)}\right)\\
&=&\prod\limits_ib^{-1}_i\prod\limits_ib^{\prime-1}_i\prod\limits_ib_i\prod\limits_ib'_i\\
&=&1
\end{eqnarray*}
and
$b^{-1}_{\sigma(i)}b^{\prime-1}_{\sigma'\sigma(i)}b_{\sigma'\sigma(i)}b'_{\sigma^{-1}\sigma'\sigma(i)}\in
\bar{B}. $  By Proposition  \ref {14.7.2}(iii),
$\sigma^{-1}\sigma^{\prime-1}\sigma\sigma'\in S'_m=A_m. $ Thus
$$((b_i)_i,  \sigma)^{-1}((b'_i)_i,  \sigma')^{-1}((b_i)_i,  \sigma)((b'_i)_i,  \sigma')\in \bar{B}\rtimes
A_m,  $$ which implies $W'\subset \bar{B}\rtimes A_m.$

Conversely,  for any  $ b=((b_i)_i,  1)\in \bar{B}$,  then
$\prod\limits_{i=1}^mb_i=1$. Let $\sigma=(1,  2,  \cdots,  m),
b'_1=1,  b'_i=\prod\limits_{j=1}^{i-1}b_j\;(2\leq i\leq m)$,   See
\begin{eqnarray*}
((b'_i)_i,  1)^{-1}(1,  \sigma)^{-1}((b'_i)_i,  1)(1,  \sigma)
&=&((b^{\prime-1}_i)_i,  1)(1,  \sigma^{-1})((b'_i)_i,  \sigma)\\
&=&((b^{\prime-1}_i,  b'_{\sigma(i)})_i,  1)\\
&=&((\prod\limits_{j=1}^{i-1}b^{-1}_j\prod\limits_{j=1}^{i}b_j)_i,  1)\\
&=&((b_i)_i,  1)=b.
\end{eqnarray*}
Therefore,   $b$ is a commutator of $((b'_i)_i,  1)$ and $(1,
\sigma)$,  i.e. $b\in W'.$ This show $\bar{B}\subset W'.$

Furthermore,  for any  $ \sigma,  \sigma'\in S_m,  $,  we have
$$(1,  \sigma)^{-1}(1,  \sigma')^{-1}(1,  \sigma)(1,  \sigma')=(1,  \sigma^{-1}\sigma^{\prime-1}\sigma\sigma'), $$
which implies $(1,  \sigma^{-1}\sigma^{\prime-1}\sigma\sigma')\in
W'$ and $S'_m=A_m\subset W'.$

Consequently,   $\bar{B}\rtimes A_m \subset W'. $$\Box$

\begin{Proposition}\label{14.8.17}
 (i) Let $B$ and  $\bar{B}$ be the same as in Lemma  \ref {14.8.15}. Then $B/\bar{B}\cong C_l;$

(ii) Let $W$ and  $W'$ be the same as in Lemma \ref {14.8.15}. Then
$$ W/W'\cong\left\{
\begin{array}{lll}
C_l\times C_2& m\geq2\\
C_l &m =1\\
1 & m =0
\end{array}
\right. ; $$

(iii) Let $\sigma$ be the same as in Proposition \ref{14.8.12}. Then
$$Z_{\sigma}/Z'_{\sigma}\cong\left(\prod\limits_{\lambda_i=1}C_{i}\right)
\left(\prod\limits_{\lambda_i\geq2}(C_{i}\times C_2)\right).$$
\end{Proposition}

{\bf Proof. } (i) It is clear that  $\nu : B\rightarrow C_l$ by
sending $(b_i)_i$ to $ \prod\limits_ib_i$  is a group homomorphism
with Ker$\nu=\bar{B}$. Thus (i) holds.

(ii) By (i) and Proposition \ref{14.7.2}(iii),  we have to show
$W/W'\cong (B/\bar{B})\times (S_{m}/A_{m}). $ Define
\begin{align*}
\psi:W &\rightarrow (B/\bar{B})\times (S_{m}/A_{m})\\
(b,  \sigma) &\mapsto (\bar{b},  \bar{\sigma}),
\end{align*}
where $\bar{b}$ and $\bar{\sigma}$ are the images of $b$ and
$\sigma$ under canonical epimorphisms $B\rightarrow B/\bar{B}$ and
$S_m\rightarrow S_{m}/A_{m}$,  respectively.  Therefore $$ \psi((b',
\sigma')(b,  \sigma))=(\overline{b'\sigma'(b)},
\overline{\sigma'\sigma})=(\bar{b'}\overline{\sigma'(b)},
\bar{\sigma'}\bar{\sigma}). $$  Let $b=(b_i)_i$. Since
$\nu(\sigma'(b))=\nu((b_{\sigma^{\prime-1}(i)})_i)=\prod\limits_ib_{\sigma^{\prime-1}(i)}=\prod\limits_ib_i=\nu(b)$,
we have $\overline{\sigma'(b)}=\bar{b}$. Consequently,
$$\psi((\sigma',  b')(\sigma,  b))=(\bar{b'}\overline{\sigma(b)},  \bar{\sigma'}\bar{\sigma})=(\bar{b'}\bar{b},  \bar{\sigma'}
\bar{\sigma},  )=\psi((b',  \sigma'))\psi((b,  \sigma)). $$ This
show that $\psi$ is a group homomorphism. Obviously,   it is
surjective and Ker$\psi=\bar{B}\rtimes A_{m}=W'$. we complete the
proof of (ii).

(iii) By Proposition \ref{14.8.12}(i) and Lemma \ref{14.8.14}(ii),
$Z_{\sigma}/Z'_{\sigma}=(\prod\limits_i
Z_{\sigma_i})/(\prod\limits_i Z'_{\sigma_i})\cong\prod\limits_i
(Z_{\sigma_i}/Z'_{\sigma_i}).$ Applying (ii),  we complete the proof
of (iii). $\Box$

\begin{Corollary}\label{14.8.18}
If $\sigma \in C$ is a permutation  of type
$1^{\lambda_1}2^{\lambda_2}\cdots n^{\lambda_n}$,  then $ \mid \!
\widehat {Z_{\sigma }}\! \mid
=\left(\prod\limits_{\lambda_i=1}i\right)\left(\prod\limits_{\lambda_i\geq2}2i\right)$,
written as  $\gamma _C$.
\end{Corollary}

{\bf Proof. } It follows from Corollary \ref{14.7.10} and
Proposition \ref{14.8.17}(iii). $\Box$

{\bf Remark.} Obviously,  $\gamma_{_C}$ only depends on the
conjugate class $C$．

\subsection { Fixed Point Set $F_{(\phi_{g_{_C}},  \phi_C)_C}$  }
  Let $(\phi_{g_{_C}}, \phi_{_C})_C\in M$. If  $\left(\left\{\chi^{(i)}_{C}\right\}_i,  u(C)\right)_C \in
F_{(\phi_{g_{_C}},  \phi_C)_C},  $ then,  according to (\ref
{e14.8.1}),  we have
\begin{eqnarray*}
(\phi_{g_{_C}}, \phi_C)_C\circ\left(\left\{\chi^{(i)}_{C}\right\}_i,
u(C)\right)_C&=&\left(\left\{\chi_C^{(\phi_C^{-1}(i))}\phi^{-1}_{g_{_C}}\right\}_i,
\phi_{g_{_C}}u(C)
\right)_C\\
&=&\left(\left\{\chi^{(i)}_{C}\right\}_i,  u(C)\right)_C.
\end{eqnarray*}
Consequently,
\begin{eqnarray} \label {e14.8.2}
\phi_{g_{_C}}(u(C))=u(C)
\end{eqnarray}
and
\begin{eqnarray} \label {e14.8.3}
\chi_C^{(\phi_C^{-1}(i))}\phi^{-1}_{g_{_C}}=\chi^{(i)}_{C}
\end{eqnarray}
for any $ C \in\mathcal {K}_r(G)$,  $i\in I_C(r)$.

 (\ref
{e14.8.2}) implies $u(C)\in Z_{g_{_C}}$ and $g_{_C}\in Z_{u(C)}$.
Considering $u(C)\in C$,  we have $u(C)\in Z_{g_{_C}}\cap C$.
Conversely,  if  $\sigma\in Z_{g_{_C}}\cap C$,  then $u(C)=\sigma$
satisfies (\ref {e14.8.2}). Consequently,  the number of $u(C)$'s
which satisfy (\ref {e14.8.2}) is  $\mid Z_{g_{_C}}\cap C \mid, $
written
$$\beta_{g_{_C}}:= \mid Z_{g_{_C}}\cap C \mid .$$

(\ref {e14.8.3}) is equivalent to
\begin{eqnarray} \label {e14.8.4}
\chi_C^{(i)}\phi^{-1}_{g_{_C}}=\chi^{(\phi_C(i))}_{C},
\end{eqnarray}
where $\phi_C\in S_{I_C(r)}=S_{r_{_C}}$. It is clear that (\ref
{e14.8.4}) holds if and only if the following (\ref {e14.8.4'})
holds for every  cycle,  such as  $\tau=(i_1,  \cdots,  i_r)$,   of
$\sigma$:
\begin{equation}\tag{\ref{e14.8.4}$'$}\label{e14.8.4'}
\chi_C^{(i)}\phi^{-1}_{g_{_C}}=\chi^{(\tau(i))}_{C}, \ i=i_1,
\ldots, i_r.
\end{equation}

By (\ref{e14.8.4'}),  we have
\begin{eqnarray}\label {e14.8.5}
\chi^{(i_2)}_{C}&=&\chi_C^{(i_1)}\phi^{-1}_{g_{_C}},  \nonumber\\
\chi^{(i_3)}_{C}&=&\chi_C^{(i_2)}\phi^{-1}_{g_{_C}}=\chi_C^{(i_1)}(\phi^{-1}_{g_{_C}})^2,  \nonumber\\
&\cdots&\nonumber\\
\chi^{(i_r)}_{C}&=&\chi_C^{(i_{r-1})}\phi^{-1}_{g_{_C}}=\chi_C^{(i_1)}(\phi^{-1}_{g_{_C}})^{r-1},  \nonumber\\
\chi^{(i_1)}_{C}&=&\chi_C^{(i_{r})}\phi^{-1}_{g_{_C}}=\chi_C^{(i_1)}(\phi^{-1}_{g_{_C}})^{r}.
\end{eqnarray}

We can obtain $\chi^{(i_2)}_{C},  \cdots,  \chi^{(i_r)}_{C}$ when we
obtain $\chi_C^{(i_1)}$. Notice
$(\phi^{-1}_{g_{_C}})^{r}=(\phi^{r}_{g_{_C}})^{-1}=(\phi_{g^r_{_C}})^{-1}$.
Therefore  (\ref {e14.8.5}) is equivalent to $$
\chi^{(i_1)}_{C}\phi_{g^r_{_C}}=\chi^{(i_1)}_{C},  $$ i.e.
\begin{eqnarray} \label{e14.8.6}
\chi^{(i_1)}_{C}\left(g^r_{_C}h(g^r_{_C})^{-1}\right)=\chi^{(i_1)}_{C}(h)
\end{eqnarray} for any $
h\in Z_{u(C)}$.  Since $g_{_C}\in Z_{u(C)}$ implies $g_{_C}^r\in
Z_{u(C)}$,  we have $$
\chi^{(i_1)}_{C}\left(g^r_{_C}h(g^r_{_C})^{-1}\right)=
\chi^{(i_1)}_{C}\left(g^r_{_C}\right)\chi^{(i_1)}_{C}\left(h\right)\chi^{(i_1)}_{C}\left((g^r_{_C})^{-1}\right)
=\chi^{(i_1)}_{C}\left(h\right).
$$
This shows (\ref {e14.8.6}) holds for any characters in $Z_{u(C)}$,
By Corollary \ref{14.8.18},  there exist $\gamma_{_C}$ characters
$\chi_C^{(i_1)}$'s such that  (\ref {e14.8.5}) holds. Assume that
$\phi_C$  is  written as multiplication of $k_{\phi_C}$ independent
cycles. Thus given a $u(C)$,  there exist $\gamma_C^{k_{\phi_C}}$
elements $\{\chi^{(i)}_C\}_{i}$'s such that  (\ref {e14.8.4}) holds.
Consequently,  for any $C\in \mathcal {K}_r(G)$,  there exist
$\gamma_C^{k_{\phi_C}}\beta_{g_{_C}}$ distinct elements
$\left(\left\{\chi^{(i)}_{C}\right\}_i, u(C)\right)$'s  such that
both (\ref {e14.8.2})and (\ref {e14.8.3}) hold.

 Consequently,  we have
\begin{equation}\label{e14.8.77}
  \mid  F_{(\phi_{g_{_C}},
\phi_C)_C} \mid  =\prod\limits_{ C\in \mathcal
{K}_r(G)}\left(\gamma_C^{k_{\phi_C}}\beta_{g_{_C}}\right)
\end{equation}
for any  $\phi_{g_{_C}},  \phi_C)_C\in M.$
\subsection{The Formula Computing the Number of Isomorphic Classes}

Applying Burnside's Lemma and (\ref{e14.8.77}),  we have
\begin{eqnarray}\label {e14.8.7} {\cal N} (G, r)
&=&\frac{1}{ \mid  M \mid  }\sum\limits_{(\phi_{g_{_C}},
\phi_C)_C\in
M} \mid  F_{(\phi_{g_{_C}},  \phi_C)_C} \mid  \nonumber\\
&=&\frac{1}{ \mid  M \mid  }\sum\limits_{(\phi_{g_{_C}},
\phi_C)_C\in M}\prod\limits_{ C\in \mathcal
{K}_r(G)}\left(\gamma_C^{k_{\phi_C}}\beta_{g_{_C}}\right)\nonumber\\
&=&\frac{1}{ \mid  M \mid  }\prod\limits_{ C\in \mathcal
{K}_r(G)}\sum\limits_{\phi_{g_{_C}}\in \mathrm{Inn}G\atop \phi_C\in
S_{r_C}}\left(\gamma_C^{k_{\phi_C}}\beta_{g_{_C}}\right)\nonumber\\
&=&\frac{1}{ \mid  M \mid  }\prod\limits_{ C\in \mathcal
{K}_r(G)}\sum\limits_{\phi_{g_{_C}}\in
\mathrm{Inn}G}\sum\limits_{\phi_C\in
S_{r_C}}\left(\gamma_C^{k_{\phi_C}}\beta_{g_{_C}}\right)\nonumber\\
&=&\frac{1}{ \mid  M \mid  }\prod\limits_{ C\in \mathcal
{K}_r(G)}\left(\sum\limits_{\phi_{g_{_C}}\in
\mathrm{Inn}G}\beta_{g_{_C}}\sum\limits_{\phi_C\in
S_{r_C}}\gamma_C^{k_{\phi_C}}\right)
\end{eqnarray}
where  $\beta_{g_{_C}}$ and $k_{\phi_C}$ are the same as above.  $
\mid  M \mid =\prod\limits_{C\in \mathcal {K}_r(G)}(n!r_{_C}!)$ when
$n \not=2, 6$;
 $ \mid  M \mid  =\prod\limits_{C\in \mathcal {K}_r(G)}(r_{_C}!)$
 when $n=2$.

Next we simplify formula (\ref {e14.8.7}). In  (\ref {e14.8.7}),
since $\gamma_C$ depends only on conjugate class,  it can be
computed by means of Corollary \ref{14.8.18}. Notice that
$k_{\phi_C}$ denotes the number of independent cycles in independent
cycle decomposition
 of $\phi _C$ in $S_{r_{_C}}$.

\begin{Definition}\label{14.8.19} (See \cite[pages 292-295]{Ri} )
Let
$$
[x]_n=x(x-1)(x-2)\cdots (x-n+1)\qquad n=1,  2,  \cdots.
$$ Obviously,  $[x]_n$ is a polynomial  with degree $n$. Let  $(-1)^{n-k}s(n,  k)$ denote the coefficient
of $x^k$ in $[x]_n$ and  $s(n,  k)$ is called the 1st Stirling
number. That is,
$$
[x]_n=\sum\limits_{k=1}^n(-1)^{n-k}s(n,  k)x^k.
$$
\end{Definition}

\cite[Theorem 8.2.9]{Ri} obtained the meaning about the 1st Stirling
number:
\begin{Lemma}\label{14.8.20} (\cite[Theorem 8.2.9]{Ri})
There exactly exist $s(n,  k)$ permutations in $S_n$,   of  which
the numbers of independent cycles in independent cycle decomposition
 are $k$. That is,  $s(n,  k)$ = $\mid \! \{ \tau \in S_n \mid$ the number of  independent cycles in independent cycle
 decomposition of $\tau$ is $k \}\! \mid$. \end{Lemma}

About $\beta_{g}$,  we have

\begin{Lemma}\label{14.8.21}
Let $G=S_n$.  Then $\sum\limits_{g\in G}\beta_{g}=\sum\limits_{g\in
G} \mid  Z_g\cap C \mid  = \mid G \mid =n!$  for any $C\in\mathcal
{K}(G)$.
\end{Lemma}

{\bf Proof. } For any  $C\in \mathcal {K}(G),  $ we have
\begin{eqnarray*}
\sum\limits_{g\in G} \mid  Z_g\cap C \mid  &=& \sum\limits_{g\in
G}\left |  \bigcup\limits_{h\in C}\left(Z_g\cap
\{h\}\right)\right | \\
&=&\sum\limits_{g\in G}\sum\limits_{h\in
C}\left | Z_g\cap\{h\}\right |  \\
&=&\sum\limits_{h\in C}\sum\limits_{g\in
G}\left | Z_g\cap\{h\}\right | \\
&=&\sum\limits_{h\in C} \mid  Z_h \mid  \\
&=&n!  \ \ ( \mbox { \hbox { by Corollary }\ref {14.7.15}(ii)} ). \
\ \ \Box
\end{eqnarray*}

\begin{Theorem}\label{14.8.22}
Let $G=S_n$ with $n\neq6$. Then
\begin{eqnarray} \label {e14.8.8}
{\cal N} (G,  r)=\prod\limits_{C\in \mathcal
{K}_r(G)}{\gamma_C+r_C-1 \choose r_C}. \end  {eqnarray}
\end{Theorem}

{\bf Proof. } We shall show in following  two cases.

(i)  $n \not= 2,  6$. By Proposition \ref{14.7.2}(ii),  Inn$G=G$.
Applying (\ref{e14.8.7}),  we have
\begin{eqnarray*}
{\cal N} (S_n,   r) &=&\frac{1}{\prod\limits_{ C\in \mathcal
{K}_r(G)}(n!r_C!)}\prod\limits_{ C\in \mathcal
{K}_r(G)}\left(\sum\limits_{g_{_C}\in
G}\beta_{g_{_C}}\sum\limits_{\phi_C\in
S_{r_C}}\gamma_C^{k_{\phi_C}}\right)\\
&=&\prod\limits_{ C\in \mathcal
{K}_r(G)}\left(\frac{1}{n!}\sum\limits_{g_{_C}\in
G}\beta_{g_{_C}}\right)\prod\limits_{ C\in \mathcal
{K}_r(G)}\left(\frac{1}{r_C!}\sum\limits_{\phi_C\in
S_{r_C}}\gamma_C^{k_{\phi_C}}\right)\\
&=&\prod\limits_{ C\in \mathcal
{K}_r(G)}\left(\frac{1}{n!}n!\right)\prod\limits_{ C\in \mathcal
{K}_r(G)}\left(\frac{1}{r_C!}\sum\limits_{k=1}^{r_{_C}}s(r_{_C},
k)\gamma_C^k\right)
\ \ ( \hbox {by Lemma  \ref {14.8.21},   \ref {14.8.20}} ) \\
&=&\prod\limits_{ C\in \mathcal
{K}_r(G)}\left(\frac{1}{r_C!}(-1)^{r_{_C}}
\sum\limits_{k=1}^{r_{_C}}(-1)^{r_{_C}-k}s(r_{_C},  k)(-\gamma_C)^k\right)\ \ \\
&=&\prod\limits_{ C\in \mathcal
{K}_r(G)}\left(\frac{1}{r_C!}(-1)^{r_{_C}}(-\gamma_C)(-\gamma_C-1)\cdots(-\gamma_C-r_{_C}+1)\right)(\hbox{by Definition \ref{14.8.19}})\\
&=&\prod\limits_{ C\in \mathcal
{K}_r(G)}\frac{(\gamma_C+r_{_C}-1)\cdots(\gamma_C+1)\gamma_C}{r_{_C}!}\\
&=&\prod\limits_{C\in \mathcal {K}_r(G)}{\gamma_C+r_C-1 \choose
r_C}. \ \Box
\end{eqnarray*}

(ii) $n =2$. In this case we have  Inn$G=\{1\}$, $\phi_{g_{_C}}=1$
and  $\beta_{g_{_C}}=1$ in (\ref{e14.8.7}). Thus
\begin{eqnarray*}
{\cal N} (S_n,   r) &=&\frac{1}{\prod\limits_{ C\in \mathcal
{K}_r(G)}(r_C!)}\prod\limits_{ C\in \mathcal
{K}_r(G)}\sum\limits_{\phi_C\in
S_{r_C}}\left(\gamma_C^{k_{\phi_C}}\right)\\
&=&\prod\limits_{ C\in \mathcal
{K}_r(G)}\left(\frac{1}{r_C!}\sum\limits_{\phi_C\in
S_{r_C}}\gamma_C^{k_{\phi_C}}\right)\ \ (\hbox{ similar to the computation of (i)})\\
&=&\prod\limits_{C\in \mathcal {K}_r(G)}{\gamma_C+r_C-1 \choose
r_C}. \ \Box
\end{eqnarray*}

\begin{Corollary}\label{14.8.23}
If $r=\sum\limits_{C\in \mathcal {K}(G)}C$,  i.e. $r_{_C}=1 $ for
any $ C\in \mathcal {K}(G)$,  then
\begin{eqnarray} \label {e14.8.9}
{\cal N} (G,  r)=\prod\limits_{C\in \mathcal {K}(G)}\gamma_C.
\end{eqnarray}
\end{Corollary}

{\bf Proof. } Considering $\mathcal {K}_r(G)=\mathcal {K}(G)$ and
applying Theorem \ref{14.8.22},  we can complete the proof. $\Box$

\begin{Corollary}\label{14.8.24}
If $C_0=\{(1)\}$ and $\mathcal {K}_r(G)=\{C_0\}, $ then
$$ {\cal N} (G,  r)=\left\{
\begin{array}{lll}
1& n=1\\
r_{_{C_0}}+1 & n\geq2
\end{array}
\right..$$
\end{Corollary}

{\bf Proof. } By Corollary \ref{14.8.18},  we have
$\gamma_{_{C_0}}=1 $ when $n=1$;  $\gamma_{_{C_0}}=2 $ when
$n\geq2$. Applying Theorem \ref{14.8.22},  we can complete the
proof. $\Box$

\subsection{Represetative}
 We have obtained  the number of isomorphic classes,  now  we give the
 representatives

\begin{Definition}\label{14.7.23}
Given an $RSC (G,   r,   \overrightarrow{\chi },   u)$,  written
$\widehat {Z_{u(C)}}=\{ \xi _{u(C)}^{(i)} \mid i = 1,  2,  \cdots,
\gamma _C\}$  and \ $n_i : =$ $ \mid \{ j \mid \chi _C ^{(j)} = \xi
_{u(C)}^{(i)} \} \! \mid, $ \ for any $C\in \mathcal {K}_r(G)$ \ \
and \ \ \ $1\le i \le \gamma _{u(C)}\  $ ,  then    $\{(n_1, n_2,
\cdots,  n_{\gamma _C})\}_{C\in\mathcal {K}_r(G)}$ is called the
type of $RSC (G,  r,  \overrightarrow{\chi },  $ $ u)$.
\end{Definition}

\begin {Lemma} \label {14.8.30}
 $RSC (G,   r,   \overrightarrow {\chi},   u)$ and  $RSC (G,   r,
\overrightarrow {\chi'},   u)$ are isomorphic if and only if they
have the same type.

\end {Lemma}
{\bf Proof. } If  $RSC (G,   r,   \overrightarrow {\chi},   u)$ and
$RSC (G,   r,   \overrightarrow {\chi'},   u)$ have the same type,
then there exists a bijective map  $\phi _C \in S_{r_C}$ such that
 $\chi'{} ^{(\phi _C(i))}_C = \chi
^{(i)}_C$ for any  $i \in I_C(r)$,   $C \in {\cal K}_r (G)$. Let
$\phi = id _G$ and   $g_{_C} = 1$. It follows from Corollary
\ref{14.8.6} that they are isomorphic.

Conversely,  if  $RSC (G,   r,   \overrightarrow {\chi},   u)$ and
$RSC (G,  r,  \overrightarrow {\chi'},   u)$ are isomorphic, then,
by Corollary \ref{14.8.6},   there exist
 $\phi \in \mathrm{Aut }G$ and    $\phi _C \in S_{r_C}$,  as well as
 $g_{_C}\in G$ such that  $u(C)=\phi_{g_{_C}}(u(C))$ and
$\chi_{C}^{\prime(\phi_C(i))}\phi_{g_{_C}}(h)=\chi_C^{(i)}(h)$ for
any
 $C \in {\cal K}_r(G)$,   $h \in Z_{u(C)}$,    $i\in I_C(r)$.
Thus
  $g_C \in  Z_{u(C)}$ and  $\chi_{C}^{\prime(\phi_C(i))}\phi_{g_{_C}}(h)$ $=\chi '{}_{C}^{(\phi_C(i))}(g_C h g_C
^{-1})$ $=\chi_{C}^{\prime(\phi_C(i))}(h)=\chi_C^{(i)}(h) $ for any
$h \in Z_{u(C)}$,  i.e. $\chi_{C}^{\prime(\phi_C(i))}=\chi_C^{(i)}.$
This implies that $RSC (G,   r,  \overrightarrow {\chi},  u)$  and
$RSC (G,  r, \overrightarrow {\chi'},  u)$ have the same  type.
$\Box$

\begin{Theorem}\label{14.8.31}
Let $G=S_n(n\neq6)$ and  $u_0$ be a map from ${\cal K}(G)
\rightarrow G$ with  $u_0(C) \in C$  for any $ C\in {\cal K}(G)$.
Let $\bar {\Omega }(G,   r,   u_0)$ denote the set consisting of all
elements with  distinct type in  $\{ (G,  r, \overrightarrow{\chi },
u_0) \mid (G,  r,  \overrightarrow{\chi },  u_0)$ is an  RSC  $\}$.
Then $\bar {\Omega }(G,
 r,   u_0)$ becomes the representative set of  $\Omega (G,   r)$.
\end{Theorem}

{\bf Proof. } According to Theorem \ref {14.8.22} and \cite  {cao},
the number of elements in  $\bar {\Omega }(G,   r,   u_0)$  is the
same as the number of isomorphic classes in  $\Omega (G,  r)$.
Applying Lemma \ref {14.8.30} we complete the proof. $\Box$

For example,  applying Corollary \ref{14.8.23} we compute ${\cal
N}(G, r)$ for $ G=S_3,  S_4,  S_5$ and $r=\sum\limits_{C\in \mathcal
{K}(G)}C$.

\begin{Example}\label{14.7.16}
(i) Let  $G = S_3$ and $r=\sum\limits_{C\in \mathcal {K}(G)}C$. Then
\begin{center}
\begin{tabular}{c | c | c |  c}
type  &$1^3$ & $1^12^1$ & $3^1$\\
\hline
 $\gamma_C$ & 2 & 2 & 3\\
\end{tabular}
\end{center}

By (\ref {e14.8.9}),  ${\cal N}(S_3,  r)=12.$ We explicitly   write
the type of representative elements as follows:\\$\{(1,  0)_{1^3},
(1, 0)_{1^12^1},  (1,  0,  0)_{3^1}\},  \{(0,  1)_{1^3},  (1,
0)_{1^12^1},  (1,  0,  0)_{3^1}\},  \{(1,  0)_{1^3},  (0,
1)_{1^12^1},  (1,  0,  0)_{3^1}\},
\\\{(0,  1)_{1^3},  (0,  1)_{1^12^1},  (1,  0,  0)_{3^1}\},
\{(1,  0)_{1^3},  (1,  0)_{1^12^1},  (0,  1,  0)_{3^1}\},  \{(1,  0)_{1^3},  (1,  0)_{1^12^1},  (0,  0,  1)_{3^1}\},  \\
\{(0,  1)_{1^3},  (1,  0)_{1^12^1},  (0,  1,  0)_{3^1}\},  \{(0,
1)_{1^3},  (1,  0)_{1^12^1},  (0,  0,  1)_{3^1}\},  \{(1, 0)_{1^3},
(0,  1)_{1^12^1},  (0,  1,  0)_{3^1}\},  \\\{(1, 0)_{1^3},  (0,
1)_{1^12^1},  (0,  0,  1)_{3^1}\},  \{(0, 1)_{1^3},  (0,
1)_{1^12^1},  (0,  1,  0)_{3^1}\},  \{(0, 1)_{1^3},  (0,
1)_{1^12^1},  (0,  0,  1)_{3^1}\}. $

(ii)  Let  $G = S_4$ and $r=\sum\limits_{C\in \mathcal {K}(G)}C$.
Then

\begin{center}
\begin{tabular}{c |  c |  c |  c |  c |  c}
type  & $1^4$ & $1^13^1$ & $1^22^1$&$2^2$&$4^1$\\
\hline
 $\gamma_C$ & 2 & 3 & 4& 4&4\\
\end{tabular}
\end{center}
By (\ref {e14.8.9}),  ${\cal N}(S_4,   r)=384. $

(iii)  Let  $G = S_4$ and $r=\sum\limits_{C\in \mathcal {K}(G)}C$.
Then
\begin{center}
\begin{tabular}{c | c |  c |  c | c | c |  c | c}
type  & $1^5$ & $1^14^1$ & $1^23^1$&$1^32^1$&$1^12^2$&$2^13^1$&$5^1$\\
\hline
 $\gamma_C$ & 2 & 4& 6& 4&4&6&5\\
\end{tabular}
\end{center}
By (\ref {e14.8.9}),   ${\cal N}(S_5,   r)=23040. $
\end{Example}

\chapter {Appendix}\label {c13}

The proof of Theorem \ref {14.2}: \  We define two functors $W$ and
$V$ as follows.
$$\begin{array}{rrcl}
W:& ^{kG}_{kG}{\mathcal M}^{kG}_{kG} & \rightarrow&
\prod_{C \in {\mathcal K}(G)}{\mathcal M}_{kZ_{u(C)}}\\
&B&\mapsto &\{^{u(C)}\! B^1\}_{C\in{\mathcal K}(G)}\\
&f& \mapsto& W(f),\\
\end{array}$$
where for any morphism $f :B\rightarrow B'$ in $^{kG}_{kG}{\mathcal
M}^{kG}_{kG}$, $W(f)=\{^{u(C)}\! f^1\}_{C\in{\mathcal K}(G)}$.
$$\begin{array}{rrcl}
V:&\prod_{C\in{\mathcal K}(G)}{\mathcal M}_{kZ_{u(C)}}
&\rightarrow&^{kG}_{kG}{\mathcal M}^{kG}_{kG}\\
&\{M(C)\}_{C\in{\mathcal K}(G)}&\mapsto& \bigoplus_{y
= xg_\theta ^{-1}u(C) g_\theta, \ x, y\in G}x\otimes M(C)\otimes_{kZ_{u(C)}}g_{\theta},\\
&f&\mapsto&V(f),\\
\end{array}$$
where for any morphism $f=\{f_C\}_{C\in{\mathcal K}(G)}:
\{M(C)\}_{C\in{\mathcal K}(G)}\rightarrow \{N(C)\}_{C\in{\mathcal
K}(G)}$, $V(f)(x\otimes  m\otimes_{kZ_{u(C)}}g_{\theta})=x\otimes
f_C(m)\otimes_{kZ_{(u(C))}}g_{\theta}$ for any $m\in M_{C}$, $x, y
\in G$ with $x^{-1}y=g^{-1}_{\theta}u(C)g_{\theta}$. That is,
$^yV(f)^x={\rm id}\otimes f_C\otimes{\rm id}$.

We will show that $W$ and $V$ are mutually inverse functors by three
steps.

(i) Let $B$ be a $kG$-Hopf bimodule. Then it is easy to see that
${}^{u(C)}\!  B^1$ is a right $kZ_{u(C)}$-module with the module
action given by
\begin {eqnarray}\label {th1e1'''}
 b \lhd  h = h^{-1} \cdot b \cdot h,\ \ b\in\ ^{u(C)}\! B^1,\ h\in Z_{u(C)}.
 \end {eqnarray}
Hence $W(B)$ is an object in $\prod_{C \in {\mathcal K}(G)}
{\mathcal M}_{kZ_{u(C)}}$. Moreover, one can see that $W(f)$ is a
morphism from $W(B)$ to $W(B')$ in $\prod_{C \in {\mathcal K}(G)}
{\mathcal M}_{kZ_{u(C)}}$ if $f: B\rightarrow B'$ is a morphism in
$^{kG}_{kG}{\mathcal M}^{kG}_{kG}$. Now it is straightforward to
check that $W$ is a functor from $^{kG}_{kG}{\mathcal M}^{kG}_{kG}$
to $\prod_{C \in {\mathcal K}(G)} {\mathcal M}_{kZ_{u(C)}}$.

(ii) Let $M =\{ M(C)\}_{C \in {\mathcal K}(G)}\in \prod_{C \in
{\mathcal K}(G)} {\mathcal M}_{kZ_{u(C)}}$, we define
\begin {eqnarray}\label {e1.2'''} \begin {array}{llll}
^yV(M)^x &=& x \otimes M(C)\otimes _{kZ_{u(C)}} g_{\theta},\\
h \cdot (x \otimes m \otimes _{kZ_{u(C)}} g_{\theta})  &=&
hx \otimes m \otimes _{kZ_{u(C)}} g_{\theta},\\
 (x \otimes m \otimes _{kZ_{u(C)}} g_{\theta}) \cdot h  &=&
xh \otimes m \otimes _{kZ_{u(C)}} g_{\theta}h=xh \otimes (m\lhd
\zeta_{\theta}(h))\otimes_{kZ_{u(C)}}g_{\theta'},
\end {array}
\end {eqnarray}
where $h, x, y \in G$ with $x^{-1}y\in C$ and the relation
(\ref{e0.2}) and (\ref {e0.3}), $m\in M(C)$. It is easy to check
that $V(M)$ is a $kG$-Hopf bimodule. If $f:\{M(C)\}_{C \in {\mathcal
K}(G)} \rightarrow\{N(C)\}_{C \in {\mathcal K}(G)}$ is a morphism in
$\prod_{C \in {\mathcal K}(G)} {\mathcal M}_{kZ_{u(C)}}$, then one
can see that $V(f): V(M)\rightarrow V(N)$ is a morphism in
$^{kG}_{kG}{\mathcal M}^{kG}_{kG}$. Now let $f:\{M(C)\}_{C \in
{\mathcal K}(G)} \rightarrow \{N(C)\}_{C \in {\mathcal K}(G)}$ and
$f':\{N(C)\} \rightarrow \{ L(C)\}_{C \in {\mathcal K}(G)}$ be two
morphisms in
 $\prod_{C \in {\mathcal K}(G)} {\mathcal M}_{kZ_{u(C)}}$.
Then for any $x, y \in G$ with
$x^{-1}y=g^{-1}_{\theta}u(C)g_{\theta}$ and $m\in M(C)$, we have
 \begin {eqnarray*}
V(f'f)(x\otimes m\otimes_{kZ_{u(C)} }g_{\theta})&=&x\otimes
(f'f)_C(m)\otimes_{kZ_{u(C)}}
g_{\theta}\\
&=& x\otimes f'_C(f_C(m))\otimes _{kZ_{u(C)}}g_{\theta}\\
&=& V(f')(V(f)(x\otimes m\otimes_{kZ_{u(C)}} g_{\theta})) .
\end {eqnarray*}
Therefore, $V(f'f)=V(f')V(f)$. It follows that $V$ is a functor from
$\prod_{C \in {\mathcal K}(G)} {\mathcal M}_{kZ_{u(C)}}$ to
$^{kG}_{kG}{\mathcal M}^{kG}_{kG}$.

(iii) We now construct two natural isomorphisms $\varphi $ and
$\psi$ as follows. Let $B \in\ ^{kG}_{kG}{\mathcal M}^{kG}_{kG}$.
Define a $k$-linear map $\varphi_B$ from $B$ to $VW(B)=\bigoplus_{y
= xg_\theta ^{-1}u(C) g_\theta, \ x, y\in
G}x\otimes\hspace{0.1cm}^{u(C)}\! B^1\otimes_{kZ_{u(C)}}g_{\theta}$
 by
 \begin {eqnarray} \label {e1.6'''}\varphi_B(b) &=& x \otimes(g_{\theta} x^{-1})\cdot
b\cdot g_{\theta}^{-1}
  \otimes _{kZ_{u(C)}}g_{\theta},
  \end {eqnarray}
  for any $x, y \in G$ satisfying Eq.(\ref{e0.2}), $b\in {}^yB^x$. Now we have
  \begin {eqnarray*} \delta ^-(g_{\theta} x^{-1}\cdot b\cdot g_{\theta}^{-1})&=&
  (g_{\theta}x^{-1})\cdot(\delta ^-(b))\cdot g_{\theta}^{-1} \ \ \ \
  ( B\hbox { is a }kG\hbox {-Hopf } \hbox { bimodule.})\\
&=&(g_{\theta}x^{-1})\cdot(y\otimes b)\cdot g_{\theta}^{-1} \\
&=& g_{\theta}x^{-1}yg_{\theta}^{-1}\otimes
g_{\theta}x^{-1}\cdot b\cdot g_{\theta}^{-1}\\
&=& u(C)\otimes g_{\theta}x^{-1}\cdot b\cdot
g_{\theta}^{-1}\hspace{0.5cm}(\mbox{by Eq.}(\ref{e0.2}))
\end {eqnarray*}
and $$\delta^+(g_{\theta}x^{-1}\cdot b\cdot
g_{\theta}^{-1})=g_{\theta}x^{-1}\cdot b\cdot g_{\theta}^{-1}\otimes
1.$$ It follows that $g_{\theta}x^{-1}\cdot b\cdot
g_{\theta}^{-1}\in {\ }^{u(C)}\! B^1$, and hence $x\otimes
g_{\theta}x^{-1}\cdot b\cdot g_{\theta}^{-1} \otimes
_{kZ_{u(C)}}g_{\theta}\in {\ }^y VW(B)^x$. This shows that
$\varphi_B$ is a $kG$-bicomodule map.

Next, let $x, y\in G$ with $x^{-1}y=g_{\theta}^{-1}u(C)g_{\theta}$,
$C\in{\mathcal K}(G)$, and let $b\in{ }^yB^x$ and $h\in G$. Then
$h\cdot b\in{ }^{hy}B^{hx}$ and
$(hx)^{-1}(hy)=x^{-1}y=g_{\theta}^{-1}u(C)g_{\theta}$. Hence
$$\begin{array}{rcl}
\varphi_B(h\cdot b)&=&hx\otimes g_{\theta}(hx)^{-1}\cdot(h\cdot
b)\cdot g_{\theta}^{-1}
\otimes_{kZ_{u(C)}}g_{\theta}\\
&=&hx \otimes g_{\theta}x^{-1}\cdot b\cdot g_{\theta}^{-1}
  \otimes _{kZ_{u(C)}}g_{\theta}\\
&=&h\cdot(x \otimes g_{\theta}x^{-1}\cdot b\cdot g_{\theta}^{-1}
\otimes _{kZ_{u(C)}}g_{\theta})\\
&=&h\cdot\varphi_B(b).
\end{array}$$
This shows that $\varphi_B$ is a left $kG$-module map. On the other
hand, $b\cdot h\in{ }^{yh}B^{xh}$ and
$(xh)^{-1}(yh)=h^{-1}x^{-1}yh=h^{-1}g_{\theta}^{-1}u(C)g_{\theta}h$.
By Eq.(\ref{e0.3}), one gets
$g_{\theta}h=\zeta_{\theta}(h)g_{\theta'}$ with
$\zeta_{\theta}(h)\in Z_{u(C)}$ and $\theta'\in\Theta_C$. Hence
$(xh)^{-1}(yh)=g_{\theta'}^{-1}\zeta_{\theta}(h)^{-1}u(C)\zeta_{\theta}(h)g_{\theta'}
=g_{\theta'}^{-1}u(C)g_{\theta'}$. Thus we have
$$\begin{array}{rcl}
\varphi_B(b\cdot h)&=&xh\otimes g_{\theta'}(xh)^{-1}\cdot(b\cdot
h)\cdot g_{\theta'}^{-1}
\otimes_{kZ_{u(C)}}g_{\theta'}\\
&=&xh\otimes g_{\theta'}h^{-1}x^{-1}\cdot b\cdot hg_{\theta'}^{-1}
\otimes_{kZ_{u(C)}}g_{\theta'}\\
&=&xh\otimes\zeta _{\theta }(h)^{-1}g_{\theta}x^{-1}\cdot b\cdot
g_{\theta}^{-1}\zeta_{\theta }(h)
\otimes_{kZ_{u(C)}}g_{\theta'}\\
&=&xh\otimes((g_{\theta}x^{-1}\cdot b\cdot
g_{\theta}^{-1})\lhd\zeta_{\theta }(h))
\otimes_{kZ_{u(C)}}g_{\theta'} \ \ \ ( \hbox {by } (\ref {th1e1}))\\
&=&(x\otimes g_{\theta}x^{-1}\cdot b\cdot g_{\theta}^{-1}
\otimes_{kZ_{u(C)}}g_{\theta})\cdot h\\
&=&\varphi_B(b)\cdot h.
\end{array}$$
It follows that $\varphi_B$ is a right $kG$-module map, and hence
$\varphi_B$ is a $kG$-Hopf bimodule map.

Suppose that $\varphi_B(b)=0$ with $b\in{ }^yB^x$. Then $x\otimes
g_{\theta}x^{-1}\cdot b\cdot g_{\theta}^{-1}\otimes
_{kZ_{u(C)}}g_{\theta}=0$, and hence $g_{\theta}x^{-1}\cdot b\cdot
g_{\theta}^{-1}\otimes_{kZ_{u(C)}}g_{\theta}=0$. Since
$M(C)\otimes_{kZ_{u(C)}}kZ_{u(C)}g_{\theta}\cong
M(C)\otimes_{kZ_{u(C)}}kZ_{u(C)}\cong M(C)$ as $k$-vector spaces for
any right $kZ_{u(C)}$-module $M(C)$, we have $g_{\theta}x^{-1}\cdot
b\cdot g_{\theta}^{-1}=0$. This implies $b=0$. It follows that
$\varphi_B$ is injective.

Let $x\otimes b\otimes _{kZ_{u(C)}}g_{\theta}\in\ ^yVW(B)^x$ with
$b\in\ ^{u(C)}\! B^1$ and $x^{-1}y=g_{\theta}^{-1}u(C)g_{\theta}$.
Put $b':= xg_{\theta}^{-1}\cdot b\cdot g_{\theta}$. Then
\begin {eqnarray*}
\delta^-(b')&=&xg_{\theta}^{-1}\cdot\delta^-(b)\cdot g_{\theta}\\
&=&xg_{\theta}^{-1}\cdot(u(C)\otimes b)\cdot g_{\theta}\\
&=& y\otimes b',\\
\delta ^+(b')&=&b'\otimes x.
\end {eqnarray*}
Hence $b'\in {\ }^yB^x $ and $\varphi_B(b')=x\otimes b\otimes
_{kZ_{u(C)}}g_{\theta}$. This shows that $\varphi_B$ is surjective.

If $f: B\rightarrow B'$ is a $kG$-Hopf bimodule map, then it is easy
to see that $\varphi_{B'}\circ f=(VW)(f)\circ\varphi_B$. Thus we
have proven that $\varphi$ is a natural equivalence from the
identity functor $id_{{\mathcal B}(kG)}$ to the composition functor
$VW$.

Now for any $M=\{ M(C)\}_{C \in {\mathcal K}(G)}\in\prod_{C \in
{\mathcal K}(G)}
 {\mathcal M }_{kZ_{u(C)}}$, it follows from Eq.(\ref {e0.2}) that ${}^{u(C)}\! V(M)^1 = 1\otimes
M(C)\otimes _{kZ_{u(C)}} 1$. Hence
 $WV(M)= \{ 1 \otimes M(C)\otimes _{kZ_{u(C)}} 1\}_{C \in {\mathcal K}(G)}$. Define a map $\psi _M$
from $M$ to  $WV(M)$ by
$$(\psi _M )_C (m ) = 1 \otimes m \otimes _{kZ_{u(C)}}1,$$
where $C \in {\mathcal K}(G),m \in M(C)$. Clearly, $(\psi_M)_C$ is
$k$-linear. Let $m\in M(C), h\in Z_{u(C)}$. Then
\begin {eqnarray*}
((\psi _M)_C(m ))  \lhd  h &=& (1 \otimes m \otimes _{kZ_{u(C)}}1 )\lhd h\\
 &=&h^{-1} \cdot (1 \otimes m \otimes _{kZ_{u(C)}}1)\cdot h\\
 &=&1 \otimes m \lhd  h \otimes _{kZ_{u(C)}}1) \\
 &=&(\psi _M)_C(m \lhd h).
\end {eqnarray*}
That is, $(\psi_M)_C$ is a $kZ_{u(C)}$-module map. Obviously,
$(\psi_M)_C$ is bijective and $\psi$ is a natural transformation.
Thus $\psi$ is a natural equivalence from the identity functor
$id_{\prod_{C\in{\mathcal K}(G)}{\mathcal M}_{kZ_{u(C)}}}$ to the
composition functor $WV$.\  \ $\Box$

The proof of Lemma  \ref {14.1.2}: (i) We first show that every
irreducible representation of $G$ is one dimensional. In fact, this
is a well-known conclusion. For  the completeness, we give a proof.
Assume that $\rho $ is an irreducible representation of $G$
corresponding to $V$. For any $g \in G,$ since $\rho (g^m) = id _V$,
all eigenvalues of $\rho (g)$ is $m$-th roots of 1, which belong to
$k$. Say that $\lambda _g$ is an eigenvalue of $\rho (g)$. It is
clear that $Ker (\rho (g) - \lambda _g id _V): = W$ is a
$kG$-submodule of $V$ and $W\not=0$ since $\lambda _g$ is an
eigenvalue of $\rho (g)$. Therefore, $W = V$, which implies $\rho
(g) = \lambda _g id _V.$ Considering $V$ is irreducible, we have
that $dim \ V =1.$

Next we show that  the order $n$ of $G$ is not divisible by
char$(k)$ = $p$. Assume that $\beta\in k$ is a primitive $m$-root of
1. If $p \mid n$, then $p \mid m$ since there exists an element with
order $p$. Say $m = pt$. See  $( \beta ^m -1) = (\beta ^t - 1)^p =
0$, which implies $\beta ^t =1.$ We get a contradiction.

Using  Maschke's Theorem, we have that every $kG$-module is a
pointed module.

(ii) Assume $M$ is a $kG$- Hopf bimodule.  By Theorem \ref {14.2}, $
 M \cong \bigoplus _{y = g_\theta ^{-1} u(C) g_\theta, \  x, y \in G}\   x \otimes M(C)\otimes _{kZ_{u(C)}}
g_{\theta}$ as $kG$-Hopf bimodules with $M(C) = ^{u(C)}\! M ^1$ for
any $C \in {\cal K}(G).$ It follows from (i) that  $M(C)$ is a right
pointed $kG$-module. Therefore,  $M$ is a PM $kG$-Hopf bimodule.
$\Box$

The proof of Lemma  \ref {14.1.3}:  It is clear that If $\xi$ is a
graded linear map from graded vector space $V=\bigoplus_{i\geq
0}V_{(i)}$ to graded vector space $V'=\bigoplus_{i\geq 0}V_{(i)}'$,
then \begin {eqnarray}\label {le13e1} \xi \iota _i \pi_i = \iota _i
\pi_i \xi, \ \ \xi \iota _i = \iota _i \pi_i \xi \iota _i, \ \ \pi_i
\xi = \pi_i \xi \iota _i \pi_i. \end {eqnarray}for $i =1, 2,
\cdots.$

Observe that $(M,\alpha^-,\alpha^+)$ is a $B$-bimodule. Since
$\Delta _H$ and $\mu _H$ are graded, we have the following equation
by (\ref {le13e1}):
 \begin {eqnarray}\label {e2.1.1}
\Delta\iota_1\pi_1&=&(\iota_1\pi_1\otimes\iota_0\pi_0)\Delta+
(\iota_0\pi_0\otimes\iota_1\pi_1)\Delta.
 \end {eqnarray}
 Since $\pi_1\iota_1={\rm id}$, it
follows from Eq.(\ref{e2.1.1}) that
 \begin {eqnarray}\label {e2.1.1'}
\Delta\iota_1&=&(\iota_1\pi_1\otimes\iota_0\pi_0)\Delta\iota_1+
(\iota_0\pi_0\otimes\iota_1\pi_1)\Delta\iota_1.
\end {eqnarray}
By a straightforward computation, one can show the following two
equations:
\begin {eqnarray}\label {e2.1.9}
\delta^-\pi _1=(\pi_0\otimes\pi_1)\Delta, & \ \ &
 \delta^+\pi_1=(\pi_1\otimes\pi_0)\Delta.
 \end {eqnarray}
 Now we have
 \begin {eqnarray*} ({\rm id}\otimes\delta^-)\delta^-&=&
({\rm id}\otimes\delta^-)(\pi_0\otimes\pi_1)\Delta\iota_1 \\
&=&(\pi_0\otimes\pi_0\otimes\pi_1)({\rm
id}\otimes\Delta)\Delta\iota_1 \ \ (\hbox { by Eq.}(\ref {e2.1.9})),
\end {eqnarray*}
and
\begin {eqnarray*} (\Delta_B\otimes{\rm id})\delta^-&=&
(\Delta_B\otimes{\rm id})(\pi_0\otimes\pi_1)\Delta\iota_1\\
&=&(\pi_0\otimes\pi_0\otimes\pi_1)(\Delta\otimes{\rm
id})\Delta\iota_1.
\end {eqnarray*}
This shows that $({\rm
id}\otimes\delta^-)\delta^-=(\Delta_B\otimes{\rm id})\delta^-$.
Clearly, $(\varepsilon\otimes{\rm id})\delta^-={\rm id}$. Hence $(M,
\delta^-)$ is a left $B$-comodule. Similarly, one can show that $(M,
\delta^+)$ is a right $B$-comodule and $(M, \delta^-, \delta^+)$ is
a $B$-bicomodule. Next, we have
\begin {eqnarray*}
\delta^-\alpha^+&=&(\pi_0\otimes\pi_1)\Delta\iota_1\pi_1\mu(\iota_1\otimes\iota_0)\\
&=&(\pi_0\otimes\pi_1)\Delta\mu(\iota_1\otimes\iota_0)       \ \ ( \hbox {   by  Eq.}(\ref{e2.1.1}))\\
&=&(\pi_0\mu\otimes\pi_1\mu)({\rm id}\otimes\tau\otimes{\rm
id})(\Delta\iota_1\otimes\Delta\iota_0)\\
&=&(\mu_B\otimes\pi_1)(\pi_0\otimes\pi_0\otimes\mu)({\rm
id}\otimes\tau\otimes{\rm
id})(\Delta\otimes\iota_0\otimes\iota_0)(\iota_1\otimes\Delta_B)\\
&=&(\mu_B\otimes\pi_1)({\rm id}\otimes{\rm id}\otimes\mu)({\rm
id}\otimes{\rm id}\otimes
\iota_1\otimes\iota_0)({\rm id}\otimes\tau\otimes{\rm id})(\pi_0\otimes\pi_1\otimes{\rm id}\otimes{\rm id})\\
&{\ }&(\Delta\iota_1\otimes\Delta_B)\ \ ( \hbox {    by Eq.}(\ref {e2.1.1'}))\\
&=&(\mu\otimes\alpha^+)({\rm id}\otimes\tau\otimes{\rm
id})(\delta^-\otimes\Delta_B),
   \end {eqnarray*}
where $\tau $ is the usual twist map. This shows that $\delta^-$ is
a right $B$-module map. Similarly, one can show that $\delta^-$ is a
left $B$-module map, and hence a $B$-bimodule map. A similar
argument shows that $\delta^+$ is also a $B$-bimodule map. Hence
$(M, \alpha^-, \alpha^+, \delta^-, \delta^+)$ is a $B$-Hopf
bimodule. $\Box$\

Proof of Lemma \ref {14.1.4}: (i) Suppose that the tensor algebra
$T_B(M)$ of $M$ over $B$ admits a graded Hopf algebra structure.
Obviously, $B$ is a Hopf algebra. It follows from Lemma \ref
{14.1.3} that $M$ is a $B$-Hopf bimodule.

Conversely, suppose that $B$ admits a Hopf algebra structure and $M$
admits a $B$-Hopf bimodule structure. By \cite[Section 1.4]{Ni78},
$T_B(M)$ is a bialgebra with the counit
$\varepsilon=\varepsilon_B\pi_0$ and the comultiplication
$\Delta=(\iota_0\otimes\iota_0)\Delta_B+\sum_{n>0}\mu^{n-1}T_n
(\Delta _M)$, where $\varepsilon_B$ is the counit of $B$,
$\Delta_M=(\iota_0\otimes\iota_1)\delta _M^-+
(\iota_1\otimes\iota_0)\delta_M^+$. It is easy to see that $T_B(M)$
is a graded bialgebra. Then it follows from \cite [Proposition
1.5.1]{Ni78} that $T_B(M)$ is a graded Hopf algebra.

 (ii)  Assume that the cotensor coalgebra $T_B^c(M)$ of $M$ over $B$
admits a graded Hopf algebra structure. Clearly, $B$ is a Hopf
algebra. By Lemma \ref {14.1.3} one gets that $M$ is a $B$-Hopf
bimodule.

Conversely, assume that $B$ admits a Hopf algebra structure and
 $M$ admits a $B$-Hopf bimodule structure. By \cite [Section 1.4] {Ni78},
 $T_B^c(M)$ is a bialgebra containing $B$ as a subbialgebra, and the multiplication of
 $T_B^c(M)$ is $\mu = \mu_B (\pi  _0 \otimes \pi  _0) + \sum _{n>0}
 T_n ^c(\mu _M) \Delta _{n-1}$, where $\mu _M = \alpha   _M^-
 (\pi _0 \otimes \pi _1) + \alpha   _M^+(\pi _1 \otimes \pi _0)$. Since $\Delta $ preserves the grading,
$\mu$ is a graded map. Thus $T_B^c(M)$ is a graded bialgebra. By
\cite [Proposition 1.5.1]{Ni78}, $T_B^c(M)$ is a graded Hopf
algebra. $\Box$\

The proof of Lemma \ref  {14.1.5}:  (i) Since $\phi$ is a Hopf
algebra map, $\iota_0\phi: B\rightarrow T_{B'}(M')$ is an algebra
map. Since $\psi$ is a $B$-bimodule map from $M$ to
$_{\phi}M'_{\phi}$, $\iota_1\psi$ is a $B$-bimodule map from $M$ to
$_{\iota_0\phi}T_{B'}(M')_{\iota_0\phi}$. It follows that $T_B
(\iota_0\phi , \iota_1\psi)$ is an algebra homomorphism from
$T_B(M)$ to $T_{B'}(M')$. Clearly, $T_B (\iota_0\phi, \iota_1\psi)$
preserves the gradation  of $T_B(M)$ and $T_{B'}(M')$.

Let $\Phi:=T_B(\iota_0\phi,\iota_1\psi)$. Then both
$\Delta_{T_{B'}(M')}\Phi$ and $(\Phi\otimes\Phi)\Delta_{T_B(M)}$ are
graded algebra maps from $T_B(M)$ to $T_{B'}(M')\otimes T_{B'}(M')$.
Now we show that
$\Delta_{T_{B'}(M')}\Phi=(\Phi\otimes\Phi)\Delta_{T_B(M)}$. By the
universal property of $T_B(M)$, we only have to prove that
$\Delta_{T_{B'}(M')}\Phi\iota_n=(\Phi\otimes\Phi)\Delta_{T_B(M)}\iota_n$
for $n=1, 2$. Since $\phi$ is a coalgebra map, one can see that
$\Delta_{T_{B'}(M')}\Phi\iota_0=(\Phi\otimes\Phi)\Delta_{T_B(M)}\iota_0$.
Note that $\psi$ is a $B'$-bicomodule map from $^{\phi}M^{\phi}$ to
$M'$. Hence by the proof of Lemma \ref{14.1.4} we have
\begin{eqnarray*}
\Delta_{T_{B'}(M')}\Phi\iota_1&=&\Delta_{T_{B'}(M')}\iota_1\psi\\
&=&(\iota_0\otimes\iota_1)\delta^-_{M'}\psi+(\iota_1\otimes\iota_0)\delta^+_{M'}\psi\\
&=&(\iota_0\otimes\iota_1)(\phi\otimes \psi)\delta^-_M+(\iota_1\otimes\iota_0)(\psi\otimes\phi)\delta^+_M\\
&=&(\iota_0\phi\otimes\iota_1\psi)\delta^-_M+(\iota_1\psi\otimes\iota_0\phi)\delta^+_M\\
&=&(\Phi\otimes\Phi)(\iota_0\otimes\iota_1)\delta^-_M+(\Phi\otimes\Phi)(\iota_1\otimes\iota_0)\delta^+_M\\
&=&(\Phi\otimes\Phi)[(\iota_0\otimes\iota_1)\delta^-_M+(\iota_1\otimes\iota_0)\delta^+_M]\\
&=&(\Phi\otimes\Phi)\Delta_{T_B(M)}\iota_1.
\end{eqnarray*}
Furthermore, it is easy to see that
$\varepsilon_{T_{B'}(M')}\Phi=\varepsilon_{T_B(M)}$. It follows that
$\Phi$ is a coalgebra map, and hence a Hopf algebra map since any
bialgebra map between two Hopf algebras is a Hopf algebra map.

(ii)  Since $\phi$ is a Hopf algebra map, $\phi\pi_0:
T^c_B(M)\rightarrow B'$ is e a coalgebra map. Since $\psi$ is a
$B'$-bicomodule map from $^{\phi}M^{\phi}$ to $M'$, $\psi\pi_1$ is a
$B'$-bicomodule map from $^{\phi\pi_0}T^c_B(M)^{\phi\pi_0}$ to $M'$.
We also have ${\rm corad}(T^c_B(M))\subseteq\iota_0(B)$, and hence
$(\psi\pi_1)({\rm corad}(T^c_B(M)))=0$. It follows that
$T^c_{B'}(\phi\pi_0, \psi\pi_1)$ is a coalgebra map from $T^c_B(M)$
to $T^c_{B'}(M')$. Clearly, $T^c_{B'}(\phi\pi_0, \psi\pi_1)$ is a
graded map.

Put $\Psi:=T_{B'}^c(\phi\pi _0, \psi\pi_1)$. Then both
$\Psi\mu_{T^c_B(M)}$ and $\mu_{T^c_{B'}(M')}(\Psi\otimes\Psi)$ are
graded coalgebra maps from $T^c_B(M)\otimes T^c_B(M)$ to
$T^c_{B'}(M')$. Since $T^c_B(M)\otimes T^c_B(M)$ is a graded
coalgebra, ${\rm corad}(T^c_B(M)\otimes
T^c_B(M))\subseteq(T^c_B(M)\otimes
T^c_B(M))_0=\iota_0(B)\otimes\iota_0(B)$. It follows that
$(\pi_1\Psi\mu_{T^c_B(M)})({\rm corad}(T^c_B(M)\otimes T^c_B(M)))=0$
and $(\pi_1\mu_{T^c_{B'}(M')}(\Psi\otimes\Psi))({\rm
corad}(T^c_B(M)\otimes T^c_B(M))) =0$. Thus by the universal
property of $T^c_{B'}(M')$, in order to prove
$\Psi\mu_{T^c_B(M)}=\mu_{T^c_{B'}(M')}(\Psi\otimes\Psi)$, we only
need to show
$\pi_n\Psi\mu_{T^c_B(M)}=\pi_n\mu_{T^c_{B'}(M')}(\Psi\otimes\Psi)$
for $n=1, 2$. However, this follows from a straightforward
computation dual to Part (i). Furthermore, one can see $\Psi(1)=1$.
Hence $\Psi$ is an algebra map, and so a Hopf algebra map.  $\Box$\

The Proof of Lemma  {1.6}:   (i) $\Rightarrow$ (ii) Assume that
$\phi: B\rightarrow B'$ is a Hopf algebra isomorphism and $\psi:
M\rightarrow\ ^{\phi^{-1}}_{\phi}{M'}_{\phi}^{\phi^{-1}}$ is a
$B$-Hopf bimodule isomorphism. Then $\phi^{-1}: B'\rightarrow B$ is
also a Hopf algebra isomorphism and $\psi^{-1}: M'\rightarrow\
_{\phi^{-1}}^{\phi}M^{\phi}_{\phi^{-1}}$ is a $B'$-Hopf bimodule
isomorphism. By Lemma \ref {14.1.5}(i), we have two graded Hopf
algebra homomorphisms $T_B(\iota_0\phi, \iota_1\psi):
T_B(M)\rightarrow T_{B'}(M')$ and $T_{B'}(\iota_0\phi^{-1},
\iota_1\psi^{-1}): T_{B'}(M')\rightarrow T_B(M)$. Then we have
$$\begin{array}{rcl}
T_{B'}(\iota_0\phi^{-1}, \iota_1\psi^{-1})T_B(\iota_0\phi,
\iota_1\psi)\iota_0&=&T_{B'}(\iota_0\phi^{-1},
\iota_1\psi^{-1})\iota_0\phi\\
&=&\iota_0\phi^{-1}\phi=\iota_0\\
\end{array}$$
and
$$\begin{array}{rcl}
T_{B'}(\iota_0\phi^{-1}, \iota_1\psi^{-1})T_B(\iota_0\phi,
\iota_1\psi)\iota_1&=&T_{B'}(\iota_0\phi^{-1},
\iota_1\psi^{-1})\iota_1\psi\\
&=&\iota_1\psi^{-1}\psi=\iota_1.
\end{array}$$
It follows from the universal property of $T_B(M)$ that
$T_{B'}(\iota_0\phi^{-1}, \iota_1\psi^{-1})T_B(\iota_0\phi,
\iota_1\psi)={\rm id}_{T_B(M)}$. Similarly, one can show that
$T_B(\iota_0\phi, \iota_1\psi)T_{B'}(\iota_0\phi^{-1},
\iota_1\psi^{-1})={\rm id}_{T_{B'}(M')}$. Hence $T_B(M)\cong
T_{B'}(M')$ as graded Hopf algebras.

(i) $\Rightarrow $ (iii) It is similar to the proof above.

(ii) $\Rightarrow $ (i) Assume that $f: T_B(M)\rightarrow
T_{B'}(M')$ is a graded Hopf algebra isomorphism. Put
$\phi:=\pi_0f\iota_0$ and $\psi:=\pi_1f\iota_1$. Then it is easy to
check that $\phi$ is a Hopf algebra isomorphism from $B$ to $B'$ and
$\psi$ is a $B$-Hopf bimodule isomorphism from $M$ to
$^{\phi^{-1}}_{\phi}{M'}^{\phi^{-1}}_{\phi}$.

(iv) $\Rightarrow$ (i) Assume that $\varphi $ is a graded Hopf
algebra isomorphism from $B[M]$ to $B'[M']$. Let $\phi := \pi
_0\varphi \iota _0$ and $\psi := \pi_1 \varphi \iota _1$. Obviously,
$\phi$ is a Hopf algebra isomorphism from $B$ to $B'$ and $\psi$ is
a bijection from $M$ to $M'$. Now we show that $\psi$ is a $B$-Hopf
bimodule isomorphism from $M$ to $^{\phi ^{-1}} _\phi M'{}^{\phi
^{-1}}_\phi$. Assume that $\alpha ^-, \alpha ^+, \delta ^- $ and
$\delta ^+$  are the module operations and comodule operations of
$M$; $\alpha'{} ^-, \alpha '{}^+, \delta '{}^- $ and $\delta'{} ^+$
are the module operations and comodule operations of $M'$.

See that
\begin {eqnarray*}
\alpha '{} ^{-} (\phi \otimes \psi ) &=& \pi_1 \mu (\iota _0 \pi_0
\varphi \iota _0 \otimes \iota _1\pi _1\varphi \iota _1) \\
&=& \pi_1 \mu (\varphi \iota _0 \otimes \varphi \iota _1) \ \ \
\mbox
{and } \\
\psi \alpha ^- &=& \pi_1 \varphi \iota _1 \pi_1 \mu (\iota _0\otimes
\iota _1) \\
&=& \pi_1\varphi \mu (\iota _0 \otimes \iota _1) \\
&= &\pi_1 \mu (\varphi \iota _0 \otimes \varphi \iota _1).
\end {eqnarray*}
Thus $\psi$ is a left $B$-module homomorphism from $M$ to $_
{\phi}M'$. Similarly, we have that $\psi$ is a right $B$-module
homomorphism from $M$ to $M' {}_\phi.$
 See that
\begin {eqnarray*}
(\phi ^{-1} \otimes id _{M'} )\delta'{} ^-\psi &=&  (\phi ^{-1}
\otimes id _{M'}) (\pi _0
\otimes \pi _1)\Delta \iota _1 \pi _1 \varphi \iota _1 \\
&=&  (\phi ^{-1} \pi_0\otimes  \pi _1)\Delta  \varphi \iota _1
\ \ \  \mbox {and }\\
(id _B \otimes \psi )\delta ^- &=& (\phi ^{-1} \pi _0 \varphi \iota
_0 \pi_0 \otimes \pi _1 \varphi \iota _0\pi_0)
 \Delta \iota _1 \\
&=&  (\phi ^{-1} \pi_0\otimes  \pi _1)\Delta  \varphi \iota _1 \ .
\end {eqnarray*}
Thus $\psi$ is a left $B$-comodule homomorphism from $M$ to $^
{\phi^{-1}}M'$. Similarly, we have that $\psi$ is a right
$B$-comodule homomorphism from $M$ to $M' {}^ {\phi^{-1}}$.

 (iii) $\Rightarrow$ (iv) Obvious.  $\Box$

The proof of Theorem \ref {14.3}:  (i) $\Rightarrow$ (ii). Assume
that $(Q,G,r)$ is a Hopf quiver, where $r$ is a ramification of $G$.
For any $C \in {\mathcal K} (G)$, let $M(C)$ be a trivial right
$kZ_{u(C)}$-module with ${ dim} M(C)=r_C$. That is, $m\cdot h=m$ for
any $m\in M(C)$ and $h\in Z_{u(C)}$. From the proof of Theorem \ref
{14.2}, $M:=\bigoplus_{y = xg_\theta ^{-1}u(C) g_\theta, \ x, y\in
G} (x\otimes M(C)\otimes_{kZ_{u(C)}}g_{\theta})$ is a $kG$-Hopf
bimodule. Let $\{\xi^{(i)}_C\mid i\in I_C(r)\}$ be a basis of $M(C)$
over $k$. Let $\phi: kQ_1 \rightarrow M$ be a $k$-linear map defined
by
$$\phi(a^{(i)}_{x,y})= x \otimes \xi^{(i)}_C \otimes _{kZ_{u(C)}} g_{\theta},$$
where $x, y \in G$ with $x^{-1}y = g_{\theta}u(C)g_{\theta}^{-1}$,
$C\in{\mathcal K}_r(G)$ and $i\in I_C(r)$. It is easy to see that
$\phi$ is a $k$-linear isomorphism and a $kG$-bicomodule
homomorphism from $kQ_1^c$ to $M$. It follows that the arrow
comodule $kQ_1^c$ admits a $kG$-Hopf bimodule structure.

(ii) $\Rightarrow$ (i). Assume that the arrow comodule $kQ_1^c$
admits a $kG$-Hopf bimodule
 structure. Then by Theorem \ref{14.2} and its proof, there is a
 vector space isomorphism
  $^y(kQ_1^c)^x \cong  x \otimes {}  ^{u(C)}\! (kQ_1^c)^1
 \otimes _{kZ_{u(C)}} g_{\theta}$  for any $x, y\in G$ with
 $y=xg_{\theta}u(C)g_{\theta}^{-1}$. Hence for any $x, y \in G$ with
 $y=xg_{\theta}u(C)g_{\theta}^{-1}$, we have ${ dim}^y(kQ_1^c)^x={ dim} ^{u(C)}\! (kQ_1^c)^1$.
 This shows that $Q$ is a Hopf quiver.

(ii) $\Leftrightarrow $ (iii) It follows from Lemma \ref {14.1.7}.
$\Box$

The proof of Lemma \ref {14.2.7}:  We first  show that
$^{\phi^{-1}}_{\phi}{\mathcal B}(V)$ is a Nichols algebra of
$^{\phi^{-1}}_{\phi}V$ by following steps. Let $R:=
^{\phi^{-1}}_{\phi}{\mathcal B}(V) $.

(i) Obviously,  $R_0 =k,$  $R_1 = P(R)$ and  $R$ is generated by
$R_1$ as algebras. Since $C(x\otimes y ) = C'(x'\otimes y')$ for any
$x' , y' \in R,$  $ \ x, y \in {\cal B}(V)$ with $x'=x, y'=y$, where
$C$ and $C'$ denote the braidings in $^B_B {\cal YD}$ and
$^{B'}_{B'} {\cal YD}$, respectively, we have that $R$ is a graded
braided Hopf algebra in $^{B'}_{B'} {\cal YD}$.

(ii) $R_1 =  ^{\phi ^{-1}} _\phi  {\cal B} (V)_1 $ $ \cong ^{\phi
^{-1}} _\phi V $ as YD $B'$-modules, since ${\cal B} (V)_1 \cong V $
as YD $B$-modules.

Consequently, $^{\phi^{-1}}_{\phi}{\mathcal B}(V)$ is a Nichols
algebra of $^{\phi^{-1}}_{\phi}V$.  By \cite [Proposition 2.2 (iv)
]{AS02}, $R \cong {\cal B} ( ^{\phi^{-1}}_{\phi}V)$ as graded
braided Hopf algebra in $^{B'}_{B'}{\cal YD}$. $\Box$

\begin{Proposition}\label {14.1.12} If $(kQ_1^c, G, r , \overrightarrow{\chi},
u )$ is a $kG$-Hopf bimodule, then
  $(kG)^*$-coactions on the $(kG)^*$-Hopf bimodule $(kQ_1^a,
G, r , \overrightarrow{\chi}, u)$ are given by
$$\delta^-(a^{(i)}_{y,x})=\sum _{h \in G} p _h \otimes  a
^{(i)} _{h^{-1}y, h^{-1}x},\ \ \delta^+(a^{(i)}_{y,x})=\sum_{h\in G}
\chi _C^{(i)}(\zeta_{\theta}(h^{-1})^{-1})a^{(i)}_{yh^{-1},
xh^{-1}}\otimes
 p_h$$
where $x, y\in G$ with $x^{-1}y=g^{-1}_{\theta}u(C)g_{\theta}$,
$\zeta_{\theta}$ is given by {\rm(\ref{e0.3})}, $C\in{\mathcal
K}_r(G)$, $i\in I_C(r)$ and $p_h\in(kG)^*$ is defined by
$p_h(g)=\delta_{h,g}$ for all $g, h\in G$.
\end{Proposition}

{\bf Proof.}  Let $(kQ_1^c, G, r , \overrightarrow{\chi}, u ) =
(kQ_1^c, \alpha^-, \alpha ^+)$ and $(kQ_1^a, G, r ,
\overrightarrow{\chi}, u )$ $  = (kQ_1^a, \delta ^-, \delta ^+)$.
Since $(kQ_1^c, kQ_1^a)$ is  an arrow dual pairing, $\xi _{kQ_1^a}$
is a $(kG)^*$- Hopf bimodule isomorphism from $(kQ_1^a, \delta ^-,
\delta ^+)$ to $((kQ_1^c{})^*, \alpha^{-*}, \alpha ^{+*})$ (see
Lemma \ref {14.1.7}). Now let $h, x, y, v, w\in G$ with
$x^{-1}y=g^{-1}_{\theta}u(C)g_{\theta}$, where $C\in{\mathcal
K}_r(G)$, $\theta\in\Theta_C$. By Proposition \ref {14.1.10}, we
have
$$\begin{array}{rcl}
\langle\delta^-((a^{(i)}_{y,x})^*), h\otimes a^{(j)}_{w,v}\rangle
&=&\langle(a^{(i)}_{y,x})^*, h\cdot a^{(j)}_{w,v}\rangle\\
&=&\langle(a^{(i)}_{y,x})^*, a^{(j)}_{hw,hv}\rangle\\
&=&\delta_{x,hv}\delta_{y,hw}\delta_{i,j}\\
&=&\delta_{h^{-1}x,v}\delta_{h^{-1}y,w}\delta_{i,j}.
\end{array}$$
This shows that $\delta^-((a^{(i)}_{y,x})^*)=\sum_{h\in
G}p_h\otimes(a^{(i)}_{h^{-1}y,h^{-1}x})^*$. Hence the left
$(kG)^*$-coaction of the $(kG)^*$-Hopf bimodule $(kQ_1^a; G, r ,
\overrightarrow{\chi}, u)$ is given by
$$\delta^-(a^{(i)}_{y,x})=\sum_{h\in
G}p_h\otimes a^{(i)}_{h^{-1}y,h^{-1}x}.$$

In order to show the second equation, we only need to show that
\begin {eqnarray}\label {2.9e1} <\delta^+((a^{(i)}_{y,x})^*), a_{wv}^{(j)}\otimes h >=\sum_{g\in G} <\chi
_C^{(i)}(\zeta_{\theta}(g^{-1})^{-1})(a^{(i)}_{yg^{-1},
xg^{-1}})^*\otimes
 p_g, \ a_{wv}^{(j)}\otimes h >
 \end {eqnarray} for any $h, v, w \in G, D \in {\cal K} (G), j \in
 I_D(r)$ and $v^{-1}w \in D, $ since $\xi _{kQ_1^a}$
is a $(kG)^*$-Hopf bimodule isomorphism. If $v \not= x h^{-1}$ or $w
\not= y h^{-1}$, then there exists $\beta \in k$ such that
$a_{w,v}^{(j)} \cdot h = \beta a_{wh,vh}^{(j)}$ by Proposition \ref
{14.1.10}. By simple computation, both the right hand side and the
left hand side of (\ref {2.9e1}) are zero.

Now assume $v = x h^{-1}$ and $w = y h^{-1}$.  Assume $g_\theta
h^{-1} $ $= \zeta _\theta (h^{-1}) g_{\theta '}. $ Thus $v^{-1}w=$ $
 hg^{-1}_{\theta}u(C)g_{\theta}h^{-1}  $ $ = $ $
g^{-1}_{\theta'}\zeta_{\theta}(h^{-1})^{-1}u(C)\zeta_{\theta}(h^{-1})g_{\theta'}
=$ $g^{-1}_{\theta'}u(C)g_{\theta'}$.  Since $g_{\theta '} h = \zeta
_{\theta '} (h^{-1})^{-1} g_\theta $,  $ \zeta _{\theta' } (h) =
\zeta _{\theta } (h^{-1})^{-1}.$
 Now we have
$$\begin{array}{rcl} \hbox { the left hand side of } (\ref {2.9e1}) &=&
 \langle \delta^+((a^{(i)}_{y,x})^*), a^{(j)}_{w,v}\otimes h \rangle\\
&=&\langle (a^{(i)}_{y,x})^*, a^{(j)}_{w,v}\cdot h\rangle\\
&=&\chi_C^{(j)}(\zeta_{\theta '}(h))\langle(a^{(i)}_{y,x})^*, a^{(j)}_{wh,vh}\rangle\\
&=&\chi_C^{(j)}(\zeta_{\theta'}(h))\delta_{i,j}\\
&=&\chi_C^{(i)}(\zeta_{\theta}(h^{-1})^{-1})\delta_{i,j}\\
&= &\hbox { the right hand side of } (\ref {2.9e1}).\ \Box
\end{array}$$

Now, we give  an interesting quantum combinatoric formula by means
of multiple Taft algebras. Let $n$ be a positive integer and $S_n$
be the symmetric group on the set $\{1,2,\cdots,n\}$. For any
permutation $\sigma\in S_n$, let $\tau(\sigma)$ denote the number of
reverse order of $\sigma $, i.e., $\tau(\sigma )=|\{(i,j)\mid 1\leq
i<j\leq n, \sigma(i)>\sigma(j)\}|$. For any $0\not= q\in k$, let
$S_n(q):=\sum_{\sigma\in S_n}q^{\tau(\sigma)}$.

\begin{Lemma}\label{14.3.10} In $kQ^c(G, r , \overrightarrow{\chi},
u )$, assume ${g} \in {\mathcal K}_r (G)$ and  $j\in I_{\{g\}}(r)$.
Let $q:= \chi ^{(j)}_{\{g\}} (g)$.
 If $i_1, i_2, \cdots, i_m$ be non-negative integers, then
$$\begin{array}{rcl}
a^{(j)}_{g^{i_m +1},g^{i_m }}\cdot a^{(j)}_{g^{i_{m-1}
+1},g^{i_{m-1}}}\cdot\cdots\cdot a^{(j)}_{g^{i_1 +1},g^{i_1 }}
&=&q^{\beta_m+\frac{m(m-1)}{2}}S_m(q^{-1})P^{(j)}_{g^{\alpha_m}}(g,m)
\end{array}$$
where $\alpha _m = i_1 + i_2 + \cdots + i_m $, $P^{(j)}_h(g,m) =$ \
$ a^{(j)}_{g^mh,g^{m-1}h}a^{(j)}_{g^{m-1}h, g^{m-2}h}\cdots
a^{(j)}_{gh,h}$, $\beta_1=0$ and $\beta_m=\sum_{j
=1}^{m-1}(i_1+i_2+\cdots+i_j )$ if $m>1$.

\end{Lemma}
{\bf Proof.} Now let $a_l=a^{(j)}_{g^{i_l+1},g^{i_l}}$ for
$l=1,2,\cdots,m$ and consider the sequence of $m$ arrows $A=(a_m,
\cdots, a_1)$. We shall use the notations of \cite[p.247]{CR02}. For
any $\sigma\in S_m$ and $1\leq l\leq m$, let $m(\sigma,l)=|\{i\mid
i<l,\sigma(i)<\sigma(l)\}|$. Then
$$m(\sigma,l)=(l-1)-|\{i\mid
i<l,\sigma(i)>\sigma(l)\}|,$$ and hence
$$\begin{array}{rcl}
\sum_{l=1}^mm(\sigma,l)&=&\frac{m(m-1)}{2}-\sum_{l=1}^m|\{i\mid
i<l,\sigma(i)>\sigma(l)\}|\\
&=&\frac{m(m-1)}{2}-\tau(\sigma).\\
\end{array}$$
Now it follows from \cite[Proposition 3.13]{CR02} that
$$\begin{array}{rcl}
a_m\cdot a_{m-1}\cdot\cdots\cdot a_1&=& \sum_{\sigma\in S_m}
A_m^\sigma\cdots  A_2^\sigma A_1^\sigma\\
&=&\sum_{\sigma \in S_m}(\prod_{l=1}^m q^{i_1+i_2+\cdots
+i_{\sigma(l)-1}+m(\sigma,l)})P^{(j)}_{g^{\alpha_m}}(g,m)\\
&=&q^{\beta_m+\frac{m(m-1)}{2}}S_m(q^{-1})
P^{(j)}_{g^{\alpha_m}}(g,m).\ \ \Box
\end{array}$$\\

 Let $G\cong{\bf  Z}$ be the infinite cyclic group with
generator $g$ and let $r$ be a ramification of $G$ given by
$r_{\{g\}}=1$ and $r_{\{g^n\}}=0$ if $n\not=1$. Let $0\not=q\in k$.
Define $\chi_{g} ^{(1)}\in\widehat G$ by $\chi_{g}^{(1)}(g)=q$. Then
$(G, r, \overrightarrow{\chi}, u)$ is an $RSC$. Using  Lemma
\ref{14.3.4} and Lemma \ref {14.3.10}, one gets the following
result.

\begin{Example}\label{14.3.11}
{\rm(i)}\ For any  \ $ 0\not=q\in k$,
$(m)_q!=q^{\frac{m(m-1)}{2}}S_m(q^{-1}) $, \ and \ $ S_m(q) =$ $
(m)_{\frac {1}{q}}!q^{\frac{m(m-1)}{2}}$, where $m$ is a positive
integer.

{\rm(ii)} \  Assume $q$ is a primitive $n$-th root of unity with
$n>1$. Then
$S_m(q)=0$ if $m\geq n$, and $S_m(q)\not=0$ if $0<m <n$.\\
\end{Example}

%%\chapter {Appendix }\label {c13}

Let $R$ and $B$  be algebras over field $k$ and $H$ be a Hopf
algebra over $k$.

If $B$  is a right $H$-comodule algebra with $R \subseteq B$ and
$B^{coH} = R$, then $B$  is called an $H$-extension of $A$.
Furthermore, if there exists an convolution-invertible
 right $H$-comodule homomorphism
$\gamma $   from $H$ to $B$, then $B$ is called an $H$-cleft
extension.

\begin {Theorem} \label {12.2.1} (see \cite [theorem 7.2.2]{Mo93})
An $H$-extension $R\subseteq  B$ is an $H$-cleft  iff
$$B \cong R\# _\sigma H \ \ \ \ \ \ \hbox { as algebras. }$$
\end {Theorem}
{\bf Sketch of proof. }  If $R \subseteq B$  is an $H$-cleft, then
there is an convolution-invertible right $H$-comodule homomorphism
$\gamma $  from $H$  to $B$. Set
$$\alpha (h, a) = \sum \gamma (h_1) a \gamma ^{-1} (h_2)$$
and
$$\sigma (h, k)= \sum \gamma (h_1)\gamma (k_1)\gamma ^{-1} (h_2k_2)$$
for any $h, k \in H, a \in R.$ We can show that the conditions in
Theorem \ref {3.1.1} are satisfied. Thus $R \#_\sigma H$  is an
algebra. We can check that $\Phi $  is an algebra isomorphism from
$R\#_\sigma H$  to $B$  by sending $a\# h$  to $a \gamma (h)$ for
any $a \in R, h\in H$.

Conversely, if $R\#_\sigma H$  is  an algebra, we define that $B=
R\#_\sigma H$  and $\gamma (h) = 1\#h$ for any $h\in H.$  It is
clear that the convolution-inverse of $\gamma $  is
$$\gamma ^{-1} (h) = \sum \sigma ^{-1}(S(h_2), h_3) \#S(h_1)$$
for any $h \in H.$  We can show that $R\subseteq B $  is an
$H$-cleft extension.
\begin{picture}(8,8)\put(0,0){\line(0,1){8}}\put(8,8){\line(0,-1){8}}\put(0,0){\line(1,0){8}}\put(8,8){\line(-1,0){8}}\end{picture}

\begin {Theorem} \label {12.2.2} (see \cite [theorem 7.2.2]{MS95})
Let $H$ be a finite-dimensional, semisimple and either  commutative
or cocommutative.  Let $A\# _\sigma H$  be a crossed product. Assume
that $R$  is $H$-prime. Then
 $R\# _\sigma H$  has only finite many minimal prime ideals and
their intersection  is zero.

\end {Theorem}

\chapter* {Index}

$\begin {array}{ll} A^{op}, \  A^{cop} \ \ \ref {s13} \ \ \ \ \  \ \
\ \ \ \ \ \ \ \ \ \ \ \ \ \ \ \ \ \  \ \ \ \ \ \ \ \ \ \ \ \ \ \
 \ &A^* , A^{\hat *} \ \ \ \ref {s13}, \ref {s15}, Ch. \ref {c1} \\
ad   \ \ \  Ch.  \ref {c4}\ \ \ \ \  \ \ \ \ \ \ \ \ \ \ \ \ \ \ \
&A \# H \ \ \ \ref {s5} \\
%\end {array} $
A \# _\sigma H, A _{\alpha , \sigma }  \# ^ {\psi , Q} H, Ch.\ref
{c3} &
 A ^{\phi , P} \# _{\beta , \mu } H,
A ^\phi _{\alpha } \# H, \  Ch. \ref {c3} \\
 A \bowtie H , \ A
\stackrel {c} {\bowtie } H, \ A ^{\phi }\bowtie ^{\psi } H \ \ \
Ch.\ref {c3} \ \ \ \ \ \ \ \ \ \ \ \ \ \ \ \ \ \ \ \        \ &A
^{\phi }_\alpha \bowtie ^{\psi }_\beta  H, \ A \stackrel {b}{\bowtie
} H \ \ \  \ref
{s3} \\
A \bowtie ^R H \ \ \ref {s14} &
  A  {\bowtie }_\tau H \ \  \ref {s13} \\
  A \bowtie _r H  \ \ref {s14} &C,  C^r, C^R , \ Ch. \ref {c1}\\
CYBE \ \ \ Ch. \ref {c12} &C(M) \ \ \  \ref {s17} \\
 {C_0} \ \ \ \ref {s17}
&\Delta  , \epsilon, \ \ \ Ch. \ref {c1}\\
\eta, \ \ \  Ch. \ref {c1}  &f*g ,\ \ \ Ch. \ref {c1}\\
    d_V, \ \ \ Ch. \ref {c1} & b_V, \ \ \ Ch. \ref {c1}\\
    H_4 , \ \ \ \ref {s9} &i_A , \ i_H , \ \ \ Ch. \ref {c3}\\
 \int_H^l, \ \int _H^r, \ \ \ Ch. \ref {c1}
      &rgD, wD, lpd, \ \ \ Ch. \ref {c9} \\
    m \ \ \ Ch. \ref {c1}
    &M_0, \ \ \  \ref {s17}\\
    M^{({\infty})} \ \ \  \ref {s17}     &M^H ,  M^{coH}, \ \ \ Ch. \ref {c1}\\
   O(H, \bar \Delta ), \ \ \ Ch. \ref {c1}
    &\pi  _A , \pi _H, \ \ \ Ch. \ref {c3}\\
     QYBE \ \ \ Ch. \ref {c12}   & r_{H}, r_{Hb}, r_{bH}, r_{Hj} ,
      r_{jH} \ \ \  \ref {s24},\ref {s25},\ref {s26}
 \\
 S \ \ \ Ch.\ref {c1}  &\Delta (h) = \sum  (h_1 \otimes h_2)
  \ \ \ Ch.\ref {c1} \\
 \rightharpoonup , \ \ \ \ref {s15}    &I ,  \ \ \ Ch. \ref {c1} \\
\hbox {action ,
adjoint action,}  \ \ Ch.\ref {c4}  & \hbox {adjoint coaction,} \ \ Ch.\ref {c4} \\
 \hbox {algebra }, \ \  Ch. \ref {c1}        &\hbox {almost cocommutative, }
 \ \ \ Ch.\ref {c6} \\
 \hbox { almost commutative, } \ Ch.\ref {c6}
& \hbox { antipode} , Ch. \ref {c1}\\
\hbox {anyonic Hopf algebra,} \ref {s15}
 &\hbox { anyonic algebra, \ \ \ \ref {s1}}\\
{\ } \hbox {bialgebra}, Ch. \ref {c1} &  {\ } \\
\end {array}$
\newpage

$\begin {array}{ll} \hbox { biproduct, bicrossproduct }, \ref {s5}
 &{ \ } \hbox {Braiding}, Ch. \ref {c1}\\
\hbox {braided group,} Ch. \ref {c1}   &\hbox { braided Hopf algebra} , Ch. \ref {c1} \\
{ \ }\hbox {Braided tensor category }, Ch. \ref {c1}
&\hbox {braided group analogue}, \ref {s8} \\
{\ }\hbox {classical Yang-Baxter equation}, 260
&{\ } \hbox {coalgebra }, Ch. \ref {c1}\\
\hbox {cobound Lie bialgebra}, Ch. \ref {c12}
&{\ } \hbox {2-cocycle, 2-cycle },Ch. \ref {c4}\\
\hbox {coevaluation } ,Ch. \ref {c1}
&{\ } \hbox {cocommutative, commutative ,} Ch. \ref {c4} \\
\hbox { comodule, comultiplication }, Ch. \ref {c1}
&{\ } \hbox {convolution }, Ch. \ref {c1}\\
\hbox {component}, \ref {s17} &{\ } \hbox {
coquasitriangular }, Ch. \ref {c1}\\
\hbox {coradical }, \ref {s17}
&{\ } \hbox {counit }, Ch. \ref {c1} \\
 \hbox { double bicrossproduct,} Ch.\ref {c3} &\hbox { cross (co)product}, Ch.\ref {c3}\\
{\ }\hbox {Drinfeld double }, \ref {s13}
&{\ }\hbox {evaluation}, Ch. \ref {c1}\\
\epsilon  \hbox {  Lie algebra }, \ref {s1}
&{\ } \hbox {full} \  H \hbox {-comodule}, \ref {s17}\\
\hbox {graded algebra, graded coalgebra}, \ref {s2}
 &{\ } \hbox {   Hopf algebra }, Ch. \ref {c1}\\
 \hbox {   Hopf module }, \ref {s1}
 &{\ } \hbox {    Hopf module algebra }, \ref {s21}\\
H\hbox {-Hopf module}, \ref {s1} &{\ }    H\hbox {-module
(co)algebra},
     \ref {s21}\\
H\hbox {-}m \hbox {-}\hbox {sequence}, Ch.\ref {c10}
&{\ }    H\hbox {-} m \hbox {-nilpotent },Ch.\ref {c10}\\
    H\hbox {-radical}, Ch.\ref {c9},
   &{\ }
  \hbox {inner action} , Ch.\ref {c4}\\
   &{\ } \hbox {integral}, Ch.\ref {c1} \\
  \hbox {left dual}, Ch.\ref {c1}
 &{\ } \hbox {Lie colour algebra, Lie H-algebra}, \ref {s1} \\
  \hbox {Maschke's theorem} , \ref {s1}
&{\ } \hbox { module algebra,
  module coalgebra}, Ch. \ref {c1} \\
  \hbox {  quantum Yang-Baxter equation ,} \ref {s11}
  &{\ }  \hbox { quantun group,}    Ch. \ref {c1}\\
\hbox {  quasitriangular,} Ch. \ref {c1} &\hbox {quantum commutative ,}  \ref {s7}  \\
    \hbox {relative irreducible}, \ref {s17}
    &{\ } R\hbox {-matrix}, Ch. \ref {c6}  \\
\hbox {  restricted Lie algebra } , \ref {s20}
 &{\ } \hbox {smash coproduct,  smash product}, Ch. \ref {c4} \\
\hbox {  strongly symmetric element }, \ref {s31}
&{\ } \hbox {  symmetric  tensor category }, Ch. \ref {c1}\\
\hbox { superspace,  super-Hopf algebra }, \ref {s1}
&{\ } \alpha , \beta\hbox {-symmetric} ,  \ref {s33} \\
 \alpha , \beta\hbox {- skew symmetric} , \ref {s31}
&{\ } \hbox {skew pairing,
skew copairing} , Ch. \ref {c6} \\
\hbox {skew group ring}, Ch.\ref {c9}
&{\ } \hbox {triangular Hopf algebra},\ref {s1}  \\
\hbox {tensor category}, Ch. \ref {c1}
&{\ } \hbox {twist map}, Ch. \ref {c1}\\
\hbox {twisted} \  H\hbox {-module}, Ch.\ref {c4}
&{\ } \hbox {wedge product}, \ref {s17}  \\
\hbox {weak action} , 79 &{\ } \hbox {weak} \  R\hbox {-matrix},
    \hbox {weak}\  r\hbox {-comatrix}, Ch.\ref {c6} \\
    \hbox {weak-closed}, \ref {s17}
    &{\ }  \hbox {w-relational hereditary}, \ref {s17}\\
    Yetter-Drinfeld module, Ch. \ref {c1}
       \end {array}$


\providecommand{\bysame}{\leavevmode\hbox
to3em{\hrulefill}\thinspace}



\begin{thebibliography}{150}



\bibitem{Ab80} E. Abe. Hopf Algebra. Cambridge University Press, 1980.



 \bibitem {AF74} F. W. Anderson and K. P. Fuller. Rings
and Categories of Modules. Springer-Verlag. New York, 1974.

\bibitem {ARS95} M. Auslander, I. Reiten and S.O. Smal$\phi$, Representation
theory of Artin algebras, Cambridge University Press, 1995.



\bibitem {AS98a} N. Andruskewisch and H.J.Schneider,
Hopf algebras of order $p^2$ and braided Hopf algebras, J. Alg. {\bf
199} (1998), 430--454.
\bibitem {AS98b} N. Andruskiewitsch and H. J. Schneider,
Lifting of quantum linear spaces and pointed Hopf algebras of order
$p^3$,  J. Alg. {\bf 209} (1998), 645--691

\bibitem {AS00} N. Andruskiewitsch and H. J. Schneider,
 Finite quantum groups and Cartan matrices, Adv. Math.,  {\bf 154} (2000), 1--45.


 \bibitem {AS02} N. Andruskewisch and H.J.Schneider, Pointed Hopf algebras,
new directions in Hopf algebras, edited by S. Montgomery and H.J.
Schneider, Cambradge University Press, 2002.


\bibitem {AZ93}  M. Aslam and A. M. Zaidi, The matrix equation in radicals,
Studia Sci.Math.Hungar., 28 (1993), 447--452.



\bibitem {BCM86} R. J. Blattner, M. Cohen and S. Montgomery, Crossed products and inner
  actions of  Hopf algebras, Transactions of the AMS., {\bf 298} (1986)2, 671--711.

\bibitem {BD82} A. A. Belavin and V. G. Drinfel'd.
Solutions of the classical Yang--Baxter equations for simple Lie
algebras. Functional Anal. Appl, {\bf 16} (1982)3, 159--180.


\bibitem {BD98} Y.Bespalov, B.Drabant, Hopf  (bi-)modules and crossed modules
in braided monoidal categories, J. Pure and Applied Algebra,
123(1998), 105-129.

\bibitem {BD99} Yuri Bespalov, Bernhard Drabant,  Cross Product Bialgebras - Part I,
   J ALGEBRA, 219 (1999), 466-505.

\bibitem {Be97}  Yu. Bespalov. Crossed Modules and quantum groups in braided
categories. Applied categorical structures, {\bf 5}(1997), 155--204.


\bibitem {BFM96}  Y. Bahturin, D. Fischman and S. Montgomery.
On the generalized Lie structure of associative algebras. Israel J.
of Math., {\bf 96}(1996) , 27--48.


\bibitem {BFM01} Y. Bahturin, D. Fischman and  S. Montgomery,
Bicharacter, twistings and Scheunert's theorem for Hopf algebra, J.
Alg. {\bf 236} (2001), 246-276.
\bibitem {BMZP92} Y. Bahturin, D. Mikhalev, M. Zaicev and V. Petrogradsky,
Infinite dimensional Lie superalgebras, Walter de Gruyter Publ.
Berlin, New York, 1992.

\bibitem {BM90}  E. Beggs and S. Majid.
Matched pairs of topological Lie algebras corresponding to Lie
bialgebra structures on $diff(S^1)$  and $diff(R).$ Ann. Inst. H.
Poincare, Phy., Theor. {\bf 53}(1990)1 ,15--34.




\bibitem {BM92} J. Bergen  and S. Montgomery. Ideal and
quotients in crossed products of Hopf algebras. J. Algebra {\bf 125}
(1992), 374--396.
\bibitem {BM89} R. J. Blattner and S. Montgomery.  Crossed products and Galois
  extensions of Hopf algebras. Pacific Journal of Math. {\bf 127} (1989), 27--55.



\bibitem {BM85}  R. J. Blattner and   S. Montgomery. A duality theorem for Hopf
module algebras.  J. algebra, {\bf 95} (1985), 153--172.

\bibitem {CHYZ04} X. W. Chen, H. L. Huang, Y. Ye and P. Zhang, Monomial Hopf
algebras,  J. Alg. {\bf 275} (2004), 212--232.


\bibitem{char} Charles W. Curtis,  Irving Reiner. Representation theory of finite group and associative
algebras, New York:John Wiley\&Sons,  1962.


\bibitem {Ch92}  William Chin. Crossed products
of semisimple cocommutive of  Hopf algebras.
  Proceeding of AMS, {\bf 116} (1992)2, 321--327.

\bibitem {Ch91}  William Chin. Crossed products and generalized inner actions of  Hopf
  algebras. Pacific Journal of Math., {\bf 150} (1991)2, 241--259.

\bibitem {Ch00}  H.X.Chen, Quantum double in monoidal categories,
Comm.Algebra,
 {\bf 28} (2000)5, 2303--2328.



\bibitem {Ch98}  Huixiang Chen.
Quasitriangular structures of bicrossed coproducts.
  J. algebra, {\bf 204} (1998), 504--531.

\bibitem {CM97}  W. Chin and S. Montgomery, Basic coalgebras, modular
interfaces,  AMS/IP Stud. Adv. Math., {\bf 4}, Amer. Math. Soc.,
Providence, RI, 1997, pp.41--47.

\bibitem {CM84a}  M. Cohen and S. Montgomery.
Smash products and group acting. Trans.  Amer. Math. Soc., {\bf 281}
(1984)1, 237--258.




\bibitem {CM84b}  M. Cohen and S. Montgomery.
 Group--graded rings, smash products,
  and group actions. Trans. Amer. Math. Soc., {\bf 282} (1984)1, 237--258.



\bibitem {CNS75} L. Corwin, Y. Ne'eman, and S. Strernberg.
 Graded Lie algebras
in mathematics and physics (Bose-fermi symmetry). Rev. Mod. Phys.,
{\bf 47} (1975), 573--604.




\bibitem {CW94}  M. Cohen and S. Westreich.
From supersymmetry to quantum commutativity. J. algebra,  {\bf 168}
(1994), 1--27.


\bibitem {CS96} Chuanren Cai and Jianhuan Sun.
 Hopf -Jacobson radical of Hopf
module algebras. Chinese Science Bulletin, {\bf 41}(1996)41,
621--625.




\bibitem {Ch00}  H.X.Chen, Quantum double in monoidal categories,
Comm.Algebra,
 {\bf 28} (2000)5, 2303--2328.



\bibitem {CZ93} Weixin Chen and Shouchuan Zhang.
 The module theoretic characterization
  of special and supernilpotent radicals for $\Gamma$-rings.
   Math.Japonic,
   {\bf 38}(1993)3, 541--547.


\bibitem{cao} Ru-Cheng Cao,   Combination Mathematics,  South China University of technology Press,
 2002.


\bibitem {CR02} C. Cibils and M. Rosso,  Hopf quivers,  J. Alg.,  {\bf  254}
(2002), 241-251.

\bibitem {CR97} C. Cibils and M. Rosso, Algebres des chemins quantiques,
Adv. Math.,   {\bf 125} (1997), 171--199.


  \bibitem {Di65} N. Divinsky.
   Rings and Radicals. Allen, London, 1965.
 \bibitem {Do93} Y. Doi.
  Braided bialgebras and quadratic bialgebras.
 Communications in algebra,
{\bf 5}  (1993)21, 1731--1749.



\bibitem {DNR01} S.Dascalescu, C.Nastasecu and S. Raianu,
Hopf algebras: an introduction,  Marcel Dekker Inc. , 2001.

  \bibitem {Dr92}  V. G. Drinfel'd.
  On some unsolved problems in Quantum group theory. In ``Proceedings,
  Leningrad 1990, p. 394''.  Lecture Notes in Mathematics,
   Vol. 1510, Springer-Verlag. Berlin, 1992.

\bibitem {Dr86} V. G. Drinfel'd.  Quantum groups. In ``Proceedings International
Congress of Mathematicians, August 3-11, 1986, Berkeley, CA" pp.
798--820, Amer. Math. Soc., Providence, RI, 1987.

 \bibitem {DT94}  Y. Doi and M. Takeuchi. Multiplication algebra by two-cocycles
 --the quantum version--. Communications in algebra,
{\bf 14}(1994)22, 5715--5731.


\bibitem {EG00} Pavel Etingof, Shlomo Gelaki,
    The classification of finite-dimensional triangular Hopf algebras
over an algebraically closed field of characteristic 0, Internat.
Math. Res. Notices, 5(2000), 223-234.


\bibitem {Fa73} C.Faith, Algebra : Rings, modules and categries, Springer-Verlag, 1973.


\bibitem {FZ92} I.  Frenkel and Y.  Zhu,   Vertex operator algebras associated to representations of affine and Virasoro
algebras,   Duke Math.  J. ,   {\bf 66} (1992),   123-168.



\bibitem {Ke99}  T.Kerler, Bridged links and tangle presentations of
cobordism categories, Adv. Math. {\bf 141} (1999), 207-- 281.

 \bibitem {Fi75}  J.R. Fisher. The Jacobson radicals
 for Hopf module algebras.
 J. algebra, {\bf 25}(1975), 217--221.

\bibitem {GM03} X. Gomez and S. Majid,  Braided Lie algebras and bicovariant
differential calculi over coquasitriangular Hopf algebras, J. Alg.
{\bf 261}(2003), 334--388.

\bibitem {GRR95}  D. Gurevich, A. Radul and V. Rubtsov,
Noncommutative differential geometry related to the Yang-Baxter
equation, Zap. Nauchn. Sem. S.-Peterburg Otdel. Mat. Inst. Steklov.
(POMI) {\bf 199 } (1992); translation in J. Math. Sci. {\bf 77 }
(1995), 3051--3062.

\bibitem {Gu86}  D. I. Gurevich, The Yang-Baxter equation and the
generalization of formal Lie theory, Dokl. Akad. Nauk SSSR, {\bf
288} (1986), 797--801.


 \bibitem {Hi92} {  J.Hietarinta.
  all solutions to the contant quantum Yang-Baxter equation in
  two dimensions.
  Phy.Letters A, } {\bf 165} (1992), 245--251.


\bibitem {He91} M.  A.  Hennings,   Hopf algebras and regular isotopy
invariants for link diagrams,   Pro.  Cambridge Phil.  Soc {\bf 109
} (1991),   59-77.


\bibitem{HR74} R. G. Heyneman and D. E. Radford. Reflexivity and coalgebras of
 finite type.
J. Algebra {\bf 28 }(1974), 215--246.

 \bibitem {HX97}  G.Q. Hu and Y. H. Xu.
Quantum commutative algebras and dualities.
   Science in China, ser A,  {\bf 26} (1996) 10, 892--900.


\bibitem {Hu74} T.  W.  Hungerford,   Algebra,
GTM {\bf 73},    New York:Springer-Verlag,  1974.

\bibitem{Hu} Guan-Zhang Hu,   Applied Modern Algebras,   Qinghua  University Press,  Beijing,   2005.
\bibitem {Ka97} L.  Kauffman,   invariants of links and 3-manifolds via Hopf algebras,   in Geometry and Physics,
Marcel Dekker Lecture Notes in Pure and Appl.  Math.  vol {\bf 184},
1997,  471-479.

\bibitem {Ja62}  {  N. Jacobson. Lie Algebras. Interscience publishers a
division of John Wiley and Sons , New York,} 1962.


\bibitem {Ji89} {  M. Jimbo. Introduction to the Yang-Baxter eqation in:
Braid group, knot theory and statistical mechanics. } Eds. C. N.
Yang and M. L. Ge (World scientific , Singapore, 1989).


\bibitem {JS86}  A. Joyal and R. Street. Braided Monoidal Categories.
Math. reports 86008, Macquaries university, 1986.
 \bibitem {JS91}  A. Joyal and R. Street. The geometry of tensor calculus
 I. Adv. Math.
 {\bf 88}(1991), 55--112.

 \bibitem {JS93}  A. Joyal and R. Street. Braided tensor categories.
  Adv. Math.
 {\bf 102}(1993), 20--78.


\bibitem {Ka95} C. Kassel.  Quantum  Groups. Graduate Texts in
Mathematics 155, Springer-Verlag, 1995.





\bibitem {Ka77} V. G. Kac. Lie superalgebras. Adv. in Math.,
{\bf 26}(1977) , 8--96.

 \bibitem {Kh99} V. K. Kharchenko,
An existence condition for multilinear quantum operations, J. Alg.
{\bf 217} (1999), 188--228.

 \bibitem {La71}  R. G. Larson. characters of Hopf algebras. J. algebra
 {\bf 17} (1971), 352--368.



 \bibitem {Li87}  Shaoxue Liu. Baer radical and
 Levitzki radical for additive category.
  J. Bejing Normal University  {\bf 4} (1987), 13--27(in Chinese).
 \bibitem{Li83}   Shaoxue Liu. Rings and Algebras. Science Press,
  1983 ( in Chinese).



 \bibitem{LR88}  R. G. Larson and D. E. Radford.
 Finite dimensional cosemisimple  Hopf algebras in characteristic 0 are
 semisimple. J. algebra
 {\bf 117} (1988), 267--289.

\bibitem{LR87}  R. G. Larson and D. E. Radford. Semisimple cosemisimple
 Hopf algebras.  Amer. J. Math.
 {\bf 109} (1987), 187--195.


\bibitem {Lu93} G. Lusztig, Introduction to  Quantum groups, Progress
Math. {\bf 110} , Berlin, 1993.

\bibitem {Ly95} V. Lyubashenko.
Tangles and Hopf algebras in braided tensor categories. J.pure and
applied algebra, {\bf 98} (1995), 245--278.









\bibitem {Ma90a} {  S. Majid.
Physics for algebraists: Non-commutative and non-cocommutative Hopf
algebras by a bicrossproduct construction. J.Algebra }  {\bf 130}
(1990), 17--64.

  \bibitem {Ma90b} S. Majid.  Quasitriangular Hopf algebras and Yang-Baxter
  equations. Int. J. Mod. Phys. A, {\bf 5} (1990), 1--91.


\bibitem {Ma93a} S. Majid.  Beyond supersymmetry and quantum symmetry
(an introduction to braided groups and braided matries). (Proceeding
of the 5th Nakai workshop, Tianjin, China), Quantum Groups,
Integrable Statistical Models and Knot Theory, 1992, edited by M-L.
Ge and H. J. de Vega, World Scientific, Singapore, 231--282, (1993).

  \bibitem {Ma93b} S. Majid. Free braided differential calculus, braided binomial
theorem, and the braided exponential map. J. Math. Phys., {\bf 34},
1993, 4843--4856 .


  \bibitem {Ma93c} S. Majid.  Braided groups. J. Pure and Applied Algebra,
 {\bf 86}  (1993), 187--221.


  \bibitem {Ma94a} S. Majid.  Bicrossproduct structure of the quantum Weyl
  group.
 J. algebra, {\bf 163}(1994), 68--87.

  \bibitem {Ma94b} S. Majid.  Crossproducts by braided groups and
  bosonization. J. algebra, {\bf 165}(1994), 165--190.



\bibitem {Ma94c} S. Majid, Quantum and braided Lie algebras, J. Geom. Phys.
{\bf 13} (1994), 307--356.



  \bibitem {Ma95a} S. Majid. Algebras and Hopf algebras
  in braided categories.
Lecture notes in pure and applied mathematics advances in Hopf
algebras, Vol. 158, edited by J. Bergen and S. Montgomery,   1995.

  \bibitem {Ma95b} S. Majid, Foundations of  Quantum Group Theory,  Cambradge University Press, 1995.



\bibitem {Mi94} {  W. Michaelis. A class of infinite-dimensional Lie bialgebras containing the Virasoro
 algebra. Advances in mathematics, } {\bf 107}  (1994), 365--392.



\bibitem {Mi80}  {  W. Michaelis.
 Lie coalgebras. Advances in mathematics, } { \bf 38}  (1980), 1--54.

\bibitem {Mi85}  {  W. Michaelis. The dual Poincare-Birkhoff-Witt Theorem,
 Advances in mathematics. } { \bf 57 }  (1985), 93--162.



\bibitem {Mo93}  S. Montgomery. Hopf algebras and their actions on rings. CBMS
  Number 82, Published by AMS, 1993.
\bibitem {MR94} S. Majid and M.J. Rodriguez Piaza. Random walk and the heat
equation on superspace and anyspace. J. Math. Phys. {\bf 35}(1994)
7, 3753--3760.


 \bibitem {MR87} J. C. McCommell and J. C. Robson. Noncommutative Noetherian
 Rings. John Wiley $\&$ Sons, New York, 1987.




 \bibitem {MS95}  S. Montgomery and H. J. Schneider. Hopf crossed products
  rings
of quotients and prime ideals.
  Advances in Mathematics,
 {\bf 112} 1995, 1--55.
\bibitem {Ni78}  W. Nichols,  Bialgebras of type one,
 Comm. Alg. {\bf 6} (1978), 1521--1552.


\bibitem {NO82}  C. Nastasescu and F. van Oystaeyen.
Graded Ring Theory. North-Holland Publishing Company, 1982.

\bibitem {NZ89}  W. D. Nichels and M. B. Zoeller. A Hopf algebras
freeness theorem. Amer. J. Math., {bf 111 } (1989), 381--385.

\bibitem {OZ04} F. Van Oystaeyen and P. Zhang, Quiver Hopf algebras, J. Alg.,
{\bf 280} (2004), 577--589.


\bibitem {Pa77}   D. S. Passman. The Algebraic Structure of Group Rings.
 John Wiley and Sons, New York,
1977.

\bibitem {Pe79}  G. K. Pedersen. $C^*$-Algebras and Their Automorphism
Groups. Academic Press INC, (London), LTD,  1979.


 \bibitem {Po97} H. C. Pop. Quantum group construction in
 a symmetric monoidal
 category.
 Communications in algebra,
{\bf 25},  1997, 117--158.

\bibitem {Pa98} B. Pareigis, On Lie algebras in the category of
Yetter-Drinfeld modules. Appl. Categ. Structures,  {\bf 6} (1998),
151--175.

 \bibitem{Ra85}   D.E. Radford. The structure of Hopf algebras
 with a projection.
 J. Alg., {\bf 92} , 1985, 322--347.


 \bibitem{Ra93}   D. E. Radford.
 Minimal quasitriangular Hopf algebras.
 J. algebra
 {\bf 157} (1993), 281--315.



     \bibitem{Ra76}   D. E. Radford.
 The order of antipode of a finite dimensional   Hopf algebras
 is finite. Amer. J. Math. {\bf 98} (1976), 333--355.




 \bibitem{Ra91}   D. E. Radford.
 On the quasitriangular structures of a semisimple Hopf algebra.
J. algebra,
 {\bf 141} (1991), 354--358.
\bibitem{Ra94a} D.  E.  Radford,  On Kauffman's knot invariants arising
 from finite-dimensional Hopf algebras,  in Advances in Hopf
 Algebras.  Marcel Dekker Lecture Notes in Pure and Appl.  Math,
 {\bf
 158} (1994),  205-266.


 \bibitem {Ra94}   D.E. Radford  The trace function and   Hopf algebras.
 J. algebra,
{\bf 163} (1994), 583--622.


\bibitem{RR04} D. Robles-Llana and M. Rocek, Quivers, quotients, and
duality, e-print hep-th/0405230.

 \bibitem {Ro79} J. J. Rotman. An introduction to homological algebras. Academic
 press, New York, 1979.

\bibitem{RT93} D.E.Radford, J.Towber,
 Yetter-Drinfeld categories associated to an arbitrary bialgebra,
 J. Pure and Applied Algebra,  87(1993) 259-279.

\bibitem{RT90} N. Yu. Reshetikhin and V.G. Turaev, Ribbon graphs and
their invariants derived from quantum groups, Commun. Math. Phys.,
{\bf 127} (1990), 1--26.

 \bibitem {Sa73}  A. D. Sands. Radicals and Morita contexts.
 J. algebra, {\bf 24}(1973), 335--345.


\bibitem {Sc79}  M. Scheunert. Generalized Lie algebras.
 J. Math. Phys., {\bf 20} (1979), 712--720.


\bibitem{suzuki} Michio Suzuki,   Group Theory I,  New
York:Springer-Verlag,  1978.
 \bibitem {SM78} T. Shudo And H. Miyamito. On the decomposition of
 coalgebras. Hiroshima
Math. J., 8(1978), 499--504.

\bibitem {St81} S. Stratila. Modular Theory in Operator Algebras.
Editura Academiei and Abacus Press, 1981.


\bibitem {Sw69a}   M. E. Sweedler. Hopf Algebras. Benjamin, New York, 1969.
 \bibitem {Sw69b}   M. E. Sweedler. Integrals for Hopf algebras. Ann. of Math. ,
{\bf 89} (1969), 323--335.

\bibitem {Sz82}   F. A. Szasz. Radicals of rings. John Wiley and Sons, New York,
1982.

 \bibitem{Ta72}   E. J. Taft. Reflexivity of algebras and coalgebras. J. Amer. Math.,
 {\bf 94} (1972)4, 1111--1130.

  \bibitem {Ta93} {  E. J. Taft.
Witt and Virasoro  algebras as Lie bialgebras. J.Pure  Appl.
algebra,
 } {\bf 87}  (1993), 301--312.



    \bibitem {To98} Wenting Tong. An Introduction to  Homological Algebras.
 Chinese Education Press,
 1998.
\bibitem {Wo89} S. L. Woronowicz, Differential calculus on
 compactmatrix pseudogroups (quantum groups), Commun. Math. Phys.,
 122 (1989), 125-170.

\bibitem {XF92} Y. H. Xu and Y. Fong. On the decomposition of coalgebras.
 In
``Proceedings International Colloquium on Words, Languages and
Combinatorics- Kyoto,August 28--31, 1990". pp. 504--522, World
Scientific , Singapore, 1992.


\bibitem{Xu95} Quan-Yan Xiang,  Modern algebras,  WuHan  University Press,  WuHan,  1995.



\bibitem{Ri} R.  A.  Brualdi,   Introductory Combinatorics,   Prentice Hall,
2004.


\bibitem {XSF94} Y. H. Xu,  K. P. Shum and Y. Fong. A decomposition Theory
of Comodules. J. Algebra, 170(1994), 880--896.

\bibitem {YG89}  C. N. Yang and M. L. Ge.
Braid group, Knot theory and Statistical Mechanics.
 World scientific , Singapore, 1989.



  \bibitem {Yi84} Zhong Yi. Homological dimension of skew group rings and
 crossed products. Journal of algebra, {\bf 164} 1984, 101--123.


\bibitem {ZC91} Shouchuan Zhang and Weixin Chen. The general theory of radicals and
  Baer radical for $\Gamma$-rings. J. Zhejiang University, {\bf 25} (1991),
719--724(in Chinese).


 \bibitem{Zh97a}
 Shouchuan Zhang. Relation between decomposition
of comodules and coalgebras.  Science in China,
 {\bf 40 }(1997), 25--29.


\bibitem {Zh97b} Shouchuan Zhang. The radicals of Hopf module algebras.
Chinese Ann. Mathematics, Ser B, 18(1997)4, 495--502.



\bibitem {Zh98a} Shouchuan Zhang. The Baer and Jacobson radicals of
crossed products. Acta math. Hungar., 78(1998), 11--24.


\bibitem {Zh98b} Shouchuan Zhang.
Classical Yang-Baxter equation and low dimensional triangular Lie
bialgebras. Physics Letters A, {\bf 246 } ( 1998 ) , 71--81.






\bibitem {ZC99}Zhang, Shouchuan,  Chen,  Hui-Xiang:    The double
bicrossproducts in braided tensor categories, { Communications in
Algebra} {\bf 29}(2001)(1), 31--66.



\bibitem {Zh03} Shouchuan Zhang,
Duality theorem and Drinfeld double in braided tensor categories,
Algebra colloq. {\bf 10 } (2003)2,  127-134.


\bibitem {Zh99d} Shouchuan Zhang. The radicals of Hopf module algebras.
Chinese Ann. Mathematics, Ser B, {\bf 18}(1997)4, 495--502.


\bibitem {Zh99e} Shouchuan Zhang.
The $H$-von Neumann regular radical and  the Jacobson radical of
twisted graded algebra. Acta Math. Hungar.,   Acta Math. Hungar.,
86(2000)4, 319-333.
\bibitem {Zh99f} Shouchuan Zhang.
The strongly symmetric elements and Yang-Baxter equations. Physics
Letters A,  261(1999)5-6, 275-283.
\bibitem {Zh99g} Shouchuan Zhang. The Homological dimension of crossed products.
 Chinese Math. Ann.  22A(2001), 767--772.




\bibitem {Zh93} S.C. Zhang, The Baer radical of generalized matrix rings,
in Proc. of the Sixth SIAM Conf. on Parallel Processing for
Scientific Computing, pp.546--551, Norfolk, Virginia, 1993. Eds:
R.F. Sincovec, D.E. Keyes, M.R. Leuze, L.R. Petzold, D.A.  Reed.
Also in  math.RA/0403280.


\bibitem{ZZ04} Zhang S C,  Zhang Y Z. Braided $m$-Lie algebras.
Letters in Mathematical Physics,  2004,  70: 155-167. Also in
math.RA/0308095.

\bibitem{ZZC07} Shouchuan Zhang, Y-Z Zhang, H.X. Chen, Classification of PM Quiver
Hopf Algebras,  Journal of Algebra and Its Applications, 6(2007)4,
1-32。

\bibitem{Zh82} Yuan-Da Zhang,  The Structure of Finite Groups (I),  Science Press,  Beijing,  1982.


\bibitem {ZWJ98} W. Zhao, S. H.  Wang and Z. Jiao. On the Hopf algebra structure
over double crossproduct.    Communications in algebra, {\bf
26}(1998)2,  467--476.



\bibitem {ZX94} X. G. Zou and Y. H. Xu. The decomposition  of comodules.
Science in China (Series A), 37(1994)8, 946--953.

\bibitem{Zh05} X. Zhu, Finite representations of a quiver arising from
string theory, e-print math.AG/0507316.
\end {thebibliography}
\end {document}